**Aleks Kleyn**

# Introduction into Noncommutative Algebra

## *Volume 1*
## *Division Algebra*

Aleks_Kleyn@MailAPS.org
http://AleksKleyn.dyndns-home.com:4080/
http://arxiv.org/a/kleyn_a_1
http://AleksKleyn.blogspot.com/

2020 *Mathematics Subject Classification.* Primary 16D10; Secondary 16K20,16K40

ABSTRACT. I dedicated the volume 1 of monograph 'Introduction into Noncommutative Algebra' to studying of algebra over commutative ring. The main topics that I covered in this volume:
- definition of module and algebra over commutative ring
- linear map of algebra over commutative ring
- vector space over associative division algebra
- system of linear equations over associative division algebra
- homomorphism of vector space over associative division algebra and covariance theory
- eigenvalue and eigenvector over associative division algebra

'

# Contents

# Preface

## 1.1. About this Book

I planned to write several different papers dedicated to linear algebra. I considered these papers in the following sections in this chapter. However at some point I realized that I am ready to assemble these papers in one book where I will consider all problems that I studied in the field of linear algebra. However, answer to one question raises new questions that must be answered also. For instance, since I want to consider Eigenvalue problem, I must to consider the theory of algebraic equations.

The list of subjects was gradually increasing. i realized that I want give a broad overview of the theory of non-commutative algebra. I understand that I will not give complete description of the subject. This is wide field. However list of questions that I want to discuss requires to write few volumes.

## 1.2. Linear Map

I will start by the following statement.

Theorem 1.2.1. *Since division $D$-algebra $A$ satisfies the following conditions*

1.2.1.1: *$D$-algebra $A$ is commutative and associative*

1.2.1.2: *the field $D$ is the center of $D$-algebra $A$*

*then*

$$(1.2.1) \qquad\qquad \mathcal{L}(D; A \to A) = A$$

We will consider the proof of the theorem 1.2.1 in volume 3. If at least one condition of the theorem 1.2.1 is not true, then the the equality (1.2.1) is imposible.

- The set of linear maps $\mathcal{L}(R; C \to C)$ of complex field is $C$-vector space of dimension 2. This is because the product in real field is commutative and complex field $C$ is vector space of dimension 2 over real field $R$.
- The set of linear maps $\mathcal{L}(R; H \to H)$ of quaternion algebra is $H$-vector space of dimension 4. This is because the product in quaternion algebra is non-commutative

During my research in quaternion algebra, I considered set of maps with certain symmetry in papers [10, 11]. It was start point of new research. The desire to find analogues of statements from theory of functions of complex variable was one of the reasons of this research.

I have made several attempts to find a set of maps that would allow me to find a structure of the set $\mathcal{L}(R; H \to H)$ similar to the structure of the set $\mathcal{L}(R; C \to C)$. However, I do not consider any of these attempts successful. However, during this research, I realized that the set $\mathcal{L}(R; H \to H)$ is left or right $H$-vector space of dimension 4.

I concidered the set of linear maps $\mathcal{L}(R; O \to O)$ of octonion algebra the similar way. However I realized that similar statement is true for any division $D$-algebra.

It does not mean that consideration of the algebra $\mathcal{L}(D; A \to A)$ as $A$-vector space can substitute representation of map as tensor.

## 1.3. System of Linear Equations

Let $V$ be left vector space over $D$-algebra $A$. Let $\overline{\overline{e}}$ be basis of left $A$-vector space $V$. In this case, we can consider two different problems, whose condition can be written as system of linear equations.

- Let $V$ be left $A$-vector space of columns. Then we write vectors of basis $\overline{\overline{e}}$ as row of matrix

$$e = \begin{pmatrix} e_1 & ... & e_n \end{pmatrix}$$

and coordinates of vector $\overline{w}$ with respect to basis $\overline{\overline{e}}$ as column of matrix

$$w = \begin{pmatrix} w^1 \\ ... \\ w^n \end{pmatrix}$$

We represent the set of vectors $v_1, ..., v_m$ as row of matrix

$$v = \begin{pmatrix} v_1 & ... & v_m \end{pmatrix}$$

We write coordinates of vector $v_i$ with respect to basis $\overline{\overline{e}}$ as column of matrix

$$v_i = \begin{pmatrix} v_i^1 \\ ... \\ v_i^n \end{pmatrix}$$

To answer the question whether it is possible to represent the vector $\overline{w}$ as linear combination of the set of vectors $v_1, ..., v_m$, we must solve the system of linear equations

$$c^1 v_1^1 + ... + c^m v_m^1 = w^1$$
$$... = ...$$
$$c^1 v_1^n + ... + c^m v_m^n = w^n$$

with respect to the set of unknown scalars $c^1, ..., c^m$.

- Let $A$ be associative division $D$-algebra. Let

$$f : V \to W$$

linear map of $A$-vector space $V$ into $A$-vector space $W$. The map of coordinate matrix of vector $v \in V$ with respect to the basis $\overline{\overline{e}}_V$ into coordinate matrix of vector $f \circ v \in W$ with respect to the basis $\overline{\overline{e}}_W$ has the following form

$$f_{11}^1 v^1 f_{12}^1 + \ldots + f_{m1}^1 v^m f_{m2}^1 = (f \circ v)^1$$
$$\ldots = \ldots$$
$$f_{11}^n v^1 f_{12}^n + \ldots + f_{m1}^n v^m f_{m2}^n = (f \circ v)^n$$

and does not depend on whether we consider left $A$-vector space or right $A$-vector space. To find coordinate matrix of vector $v \in V$, $f \circ v = w$, we must solve the system of linear equations

(1.3.1)
$$f_{11}^1 v^1 f_{12}^1 + \ldots + f_{m1}^1 v^m f_{m2}^1 = w^1$$
$$\ldots = \ldots$$
$$f_{11}^n v^1 f_{12}^n + \ldots + f_{m1}^n v^m f_{m2}^n = w^n$$

## 1.4. Eigenvalue of Matrix

### 1.4.1. Eigenvalue in commutative algebra.
I will start from the definition of eigenvalue of matrix in commutative algebra. [1.1]

Let $a$ be $n \times n$ matrix whose entries belong to the field $F$. The polynomial

$$\det(a - bE_n)$$

where $E_n$ is $n \times n$ unit matrix is called characteristic polynomial of the matrix $a$. Roots of characteristic polynomial are called eigenvalues of matrix $a$. Therefore, eigenvalues of matrix $a$ satisfy the equation

$$\det(a - bE_n) = 0$$

and the system of linear equations

(1.4.1)
$$(a - bE_n)x = 0$$

has non-trivial solution if $b$ is eigenvalue of the matrix $a$.

If matrices $a$ and $a_1$ are similar

(1.4.2)
$$a_1 = f^{-1}af$$

then the equality

(1.4.3)
$$a_1 - bE_n = f^{-1}af - bE_n = f^{-1}af - f^{-1}bE_n f = f^{-1}(a - bE_n)f$$

implies that matricess $a$ and $a_1$ have the same characteristic polynomial, and, therefore, the same characteristic roots.

Let $V$ be vector space over field $F$. Let $\overline{\overline{e}}$, $\overline{\overline{e}}_1$ be bases of vector space $V$ and $f$ be coordinate matrix of the basis $\overline{\overline{e}}_1$ with respect to the basis $\overline{\overline{e}}$

(1.4.4)
$$\overline{\overline{e}}_1 = f\overline{\overline{e}}$$

Let $a$ be the matrix of the endomorphism $\overline{a}$ with respect to the basis $\overline{\overline{e}}$. Then the matrix $a_1$ defined by the equality (1.4.2) is the matrix of the endomorphism $\overline{a}$ with respect to the basis $\overline{\overline{e}}_1$. If we identify endomorphism and its matrix with respect to given basis, then we will say that eigenvalues of endomorphism are eigenvalues of its matrix with respect to given basis. Therefore, eigenvalues of endomorphism do not depend on choice of a basis.

---

[1.1] See also the definition on the page [3]-200.

Let

$$(1.4.5) \qquad \overline{x} = x^i e_i$$

Then the equality (1.4.1) gets the following form

$$(1.4.6) \qquad \overline{a}\overline{x} = b\overline{x}$$

$$(1.4.7) \qquad ax = bx$$

with respect to the basis $\overline{\overline{e}}$. Let $f$ be coordinate matrix of the basis $\overline{\overline{e}}_1$ with respect to the basis $\overline{\overline{e}}$ $\boxed{(1.4.4) \quad \overline{\overline{e}}_1 = f\overline{\overline{e}}}$ . Then coordinate matrix $x_1$ of the vector $\overline{x}$ with respect to the basis $\overline{\overline{e}}_1$ has the following form

$$(1.4.8) \qquad x_1 = f^{-1}x$$

The equality

$$(1.4.9) \qquad a_1 x_1 = f^{-1}aff^{-1}x = f^{-1}ax = f^{-1}bx = bf^{-1}x = bx_1$$

follows from the equalities (1.4.2), (1.4.7), (1.4.8). Therefore, the equality (1.4.6) does not depend on choice of a basis, and we will say that eigenvector of endomorphism is eigenvector of its matrix with respect to given basis.

The set of eigenvalues of endomorphism $\overline{a}$ each eigenvalue being taken with their multiplicity in the characteristic polynomial is called the spectrum of endomorphism $\overline{a}$.

**1.4.2. Eigenvalue in non-commutative algebra.** Let $A$ be non-commutative division $D$-algebra $A$. If we consider vector space over $D$-algebra $A$, then we want to keep the definition of eigenvalue and eigenvector considered in case of vector space over commutative algebra.

Let $V$ be left $A$-vector space of columns. Since the product is non-commutative, the question arises. Should we use the eqaulity

$$(1.4.10) \qquad \overline{a} \circ \overline{x} = b\overline{x}$$

or the eqaulity

$$(1.4.11) \qquad \overline{a} \circ \overline{x} = \overline{x}b$$

The equality (1.4.10) looks natural because, on the right side of the equality, we wrote left-side product of vector over scalar. However, if $cb \neq bc$, $b$, $c \in A$, then

$$(1.4.12) \qquad \overline{a} \circ (c\overline{x}) \neq b(c\overline{x})$$

Statement 1.4.1. *The statement* (1.4.12) *implies that the set of vectors* $\overline{x}$ *which satisfy the equality* (1.4.10), *do not generate subspace of vector space.*      ⊙

If we write down the equality (1.4.10) using matrices, then we get the equality

$$(1.4.13) \qquad x^*{}_*a = bx$$

The equality (1.4.13) is system of linear equations with respect to the set of unknown $x$. However, we cannot say that this system of linear equations is singular and has infinitely many solutions. However this statement is an obvious consequence of the statement 1.4.1.

The impression is that the equality (1.4.10) not suitable for definition of eigenvalue. However, we will return to this definition in the section 1.5.

Consider the equality (1.4.11). Let $a$ be matrix of endomorphism $\overline{\overline{a}}$ with respect to the basis $\overline{\overline{e}}$ of left $A$-vector space of columns. Let $a$ be matrix of vector $\overline{x}$ with respect to the basis $\overline{\overline{e}}$. The system of linear equations

$$(1.4.14) \qquad x^*{}_*a = xb$$

follows from the equation (1.4.11). The system of linear equations (1.4.14) has non-trivial solution when the matrix $a - bE_n$ is $^*{}_*$-singular. In such case, the system of linear equations (1.4.14) has infinitely many solutions.

At first sight, considered model corresponds to the model considered in the subsection 1.4.1. However, if we choose another basis $\overline{\overline{e}}_1$

$$(1.4.15) \qquad \overline{\overline{e}}_1 = g^*{}_*\overline{\overline{e}}$$

then the matrix $a$ will change

$$(1.4.16) \qquad a_1 = g^*{}_*a^*{}_*g^{-1^*{}_*}$$

$$(1.4.17) \qquad a = g^{-1^*{}_*{}^*}{}_*a_1^*{}_*g$$

and the matrix $x$ will change

$$(1.4.18) \qquad x_1 = x^*{}_*g^{-1^*{}_*}$$

$$(1.4.19) \qquad x = x_1^*{}_*g$$

The equality

$$(1.4.20) \qquad x_1^*{}_*g^*{}_*g^{-1^*{}_*{}^*}{}_*a_1^*{}_*g = x_1^*{}_*gb$$

follows from equalities (1.4.14), (1.4.17), (1.4.19). From the equality (1.4.20), it follows that column vector $x_1$ is not solution of the system of linear equations (1.4.14). This contradicts the statement that the equality (1.4.11) does not depend on choice of basis.

The equality (1.4.11) requires more careful analysis. On the right side of the equality (1.4.11) I wrote down the operation which is not defined in left $A$-vector space, namely, right-sided product of vector over scalar. However we can write the system of linear equations (1.4.14) differently

$$(1.4.21) \qquad x^*{}_*a = x^*{}_*(Eb)$$

where $E$ is unit matrix. On the right side of the equality (1.4.21) as well as on the left side, I wrote down $^*{}_*$-product of column vector over matrix of endomorphism. Therefore, left and right sides have the same transformation when changing the basis.

Similar to vector space over field, the endomorphism $\overline{\overline{e}}b$ which has the matrix $Eb$ with respect to the basis $\overline{\overline{e}}$ is called similarity transformation with respect to the basis $\overline{\overline{e}}$. Therefore, we should write down the equality (1.4.11) as

$$(1.4.22) \qquad a \circ x = \overline{\overline{e}}b \circ x$$

The equality (1.4.22) is the definition of eigenvalue $b$ and eigenvector $x$ of endomorphism $a$. If we consider matrix of endomorphism $a$ with respect to the basis $\overline{\overline{e}}$ of left $A$-vector space of columns, then the equality (1.4.14) is the definition of $^*{}_*$-eigenvalue $b$ and eigenvector $x$ corresponding to $^*{}_*$-eigenvalue $b$.

## 1.5. Left and Right Eigenvalues

To better understand essence of eigenvalues, consider the following theorem.

THEOREM 1.5.1. *Let $D$ be field. Let $V$ be $D$-vector space. If endomorphism has $k$ different eigenvalues* [1.2] *$b(1)$, ..., $b(k)$, then eigenvectors $v(1)$, ..., $v(k)$ which correspond to different eigenvalues are linearly independent.*

PROOF. We will prove the theorem by induction with respect to $k$.
Since eigenvector is non-zero, the theorem is true for $k = 1$.

STATEMENT 1.5.2. *Let the theorem be true for $k = m - 1$.*                    ⊙

Let $b(1)$, ..., $b(m-1)$, $b(m)$ be different eigenvalues and $v(1)$, ..., $v(m-1)$, $v(m)$ corresponding eigenvectors

$$(1.5.1) \qquad f \circ v(i) = b(i)v(i)$$

$$i \neq j \Rightarrow b(i) \neq b(j)$$

Let the statement 1.5.3 be true.

STATEMENT 1.5.3. *There exists linear dependence*

$$(1.5.2) \qquad a(1)v(1) + ... + a(m-1)v(m-1) + a(m)v(m) = 0$$

*where $a(1) \neq 0$.*                                                        ⊙

The equality

$$(1.5.3) \quad a(1)(f \circ v(1)) + ... + a(m-1)(f \circ v(m-1)) + a(m)(f \circ v(m)) = 0$$

follows from the equality (1.5.2). The equality

$$(1.5.4) \qquad \begin{aligned} &a(1)b(1)v(1) + ... \\ &+ a(m-1)b(m-1)v(m-1) + a(m)b(m)v(m) = 0 \end{aligned}$$

follows from the equality (1.5.1), (1.5.3). The equality

$$(1.5.5) \qquad \begin{aligned} &a(1)(b(1) - b(m))v(1) + ... \\ &+ a(m-1)(b(m-1) - b(m))v(m-1) = 0 \end{aligned}$$

follows from the equality (1.5.2), (1.5.4). Since

$$a(1)(b(1) - b(m)) \neq 0$$

then the equality (1.5.5) is linear dependence between vectors $v(1)$, ..., $v(m-1)$. This statement contradicts the statement 1.5.2. Therefore, the statement 1.5.3 is not true and vectors $v(1)$, ..., $v(m-1)$, $v(m)$ are linearly independent.    □

Let dimension of $D$-vector space $V$ be $n$. From the theorem 1.5.1, it follows that if $k = n$, then the set of corresponding eigenvectors is the basis of $D$-vector space $V$. [1.3]

---

[1.2] See also the theorem on the page [3]-203.

[1.3] If multiplicity of eigenvalue is greater than 1, then dimension of $D$-vector space generated by corresponding eigenvectors maybe less than multiplicity of eigenvalue.

Let $\overline{\overline{e}}$ be a basis of $D$-vector space $V$. Let $a$ be matrix of endomorphism $\overline{a}$ with respect to the basis $\overline{\overline{e}}$ and

$$x(k) = \begin{pmatrix} x(k)^1 \\ ... \\ x(k)^n \end{pmatrix}$$

be coordinates of eigenvector $x(k)$ with respect to the basis $\overline{\overline{e}}$. Gather a set of coordinates of eigenvectors $v(1), ..., v(n-1), v(n)$ into the matrix

$$x = \left( \begin{pmatrix} x(1)^1 \\ ... \\ x(1)^n \end{pmatrix} \quad ... \quad \begin{pmatrix} x(n)^1 \\ ... \\ x(n)^n \end{pmatrix} \right) = \begin{pmatrix} x(1)^1 & ... & x(n)^1 \\ ... & ... & ... \\ x(1)^n & ... & x(n)^n \end{pmatrix}$$

Then we can write down the equality

$$\boxed{(1.4.6) \quad \overline{a}\overline{x} = b\overline{x}}$$

using product of matrices

$$(1.5.6) \quad \begin{aligned} & \begin{pmatrix} a_1^1 & ... & a_n^1 \\ ... & ... & ... \\ a_1^n & ... & a_n^n \end{pmatrix} \begin{pmatrix} x(1)^1 & ... & x(n)^1 \\ ... & ... & ... \\ x(1)^n & ... & x(n)^n \end{pmatrix} \\ & = \begin{pmatrix} x(1)^1 & ... & x(n)^1 \\ ... & ... & ... \\ x(1)^n & ... & x(n)^n \end{pmatrix} \begin{pmatrix} b(1) & ... & 0 \\ ... & ... & ... \\ 0 & ... & b(n) \end{pmatrix} \end{aligned}$$

Since the set of vectors $v(1), ..., v(n-1), v(n)$ is basis of $D$-vector space $V$, then the matrix $x$ is non-singular. Therefore, the equality

$$(1.5.7) \quad \begin{aligned} & \begin{pmatrix} x(1)^1 & ... & x(n)^1 \\ ... & ... & ... \\ x(1)^n & ... & x(n)^n \end{pmatrix}^{-1} \begin{pmatrix} a_1^1 & ... & a_n^1 \\ ... & ... & ... \\ a_1^n & ... & a_n^n \end{pmatrix} \begin{pmatrix} x(1)^1 & ... & x(n)^1 \\ ... & ... & ... \\ x(1)^n & ... & x(n)^n \end{pmatrix} \\ & = \begin{pmatrix} b(1) & ... & 0 \\ ... & ... & ... \\ 0 & ... & b(n) \end{pmatrix} \end{aligned}$$

follows from the equality (1.5.6). The equality (1.5.7) implies that the matrix $a$ is similar to diagonal matrix

$$b = \operatorname{diag}(b(1), ..., b(n))$$

Let $V$ be left $A$-vector space of columns. Then the equality (1.5.6) gets form

(1.5.8)
$$\begin{pmatrix} x(1)^1 & ... & x(n)^1 \\ ... & ... & ... \\ x(1)^n & ... & x(n)^n \end{pmatrix} {}^*_* \begin{pmatrix} a_1^1 & ... & a_n^1 \\ ... & ... & ... \\ a_1^n & ... & a_n^n \end{pmatrix}$$
$$= \begin{pmatrix} b(1) & ... & 0 \\ ... & ... & ... \\ 0 & ... & b(n) \end{pmatrix} {}^*_* \begin{pmatrix} x(1)^1 & ... & x(n)^1 \\ ... & ... & ... \\ x(1)^n & ... & x(n)^n \end{pmatrix}$$

If we select a column $x(k)$ in the equality (1.5.8), then we get the equality

(1.5.9)
$$\begin{pmatrix} x(k)^1 \\ ... \\ x(k)^n \end{pmatrix} {}^*_* \begin{pmatrix} a_1^1 & ... & a_n^1 \\ ... & ... & ... \\ a_1^n & ... & a_n^n \end{pmatrix} = b(k) \begin{pmatrix} x(k)^1 \\ ... \\ x(k)^n \end{pmatrix}$$

The equality (1.5.9) coincides with the equality

$$(1.4.10) \quad \overline{a} \circ \overline{x} = b\overline{x}$$

up to notation.

Thus, we get the definition of left ${}^*_*$-eigenvalue which is different from ${}^*_*$-eigenvalue considered in the subsection 1.4.2. If we similarly consider right $A$-vector space of rows, then get the definition of right ${}^*_*$-eigenvalue. [1.4]

## 1.6. System of Differential Equations

In commutative algebra there is one more problem to solve which we need to find eigenvalue of matrix. We are looking for a solution of system of differential equations

(1.6.1)
$$\frac{dx^1}{dt} = a_1^1 x^1 + ... + a_n^1 x^n$$
$$......$$
$$\frac{dx^n}{dt} = a_1^n x^1 + ... + a_n^n x^n$$

as

$$\begin{pmatrix} x^1 \\ ... \\ x^n \end{pmatrix} = \begin{pmatrix} c^1 \\ ... \\ c^n \end{pmatrix} e^{bt}$$

---

[1.4] In the section [20]-8.2, professor Cohn considers definition of left and right eigenvalues of $n \times n$ matrix. Cohn traditionally considers product of row over column. Since there exist two products of matrices, I slightly changed Cohn's definition.

In this case, the number $b$ is eigenvalue of the matrix

$$a = \begin{pmatrix} a_1^1 & ... & a_n^1 \\ ... & ... & ... \\ a_1^n & ... & a_n^n \end{pmatrix}$$

and the vector

$$c = \begin{pmatrix} c^1 \\ ... \\ c^n \end{pmatrix}$$

is eigenvector of the matrix $a$ for eigenvalue $b$.

In case of non-commutative $D$-algebra A, we can write the system of differential equations (1.6.1) as follows

$$(1.6.2) \qquad \frac{dx}{dt} = a_* {}^* x$$

In such case, the solution of the system of differential equations (1.6.2) may have one of following representations.

- We can write solution as

$$(1.6.3) \qquad x = e^{bt} c = \begin{pmatrix} e^{bt} c^1 \\ ... \\ e^{bt} c^n \end{pmatrix} \qquad c = \begin{pmatrix} c^1 \\ ... \\ c^n \end{pmatrix}$$

  where $A$-number $b$ is $_*{}^*$-eigenvalue of the matrix $a$ and column vector $c$ is eigencolumn for $_*{}^*$-eigenvalue $b$.

- We can write solution as

$$(1.6.4) \qquad x = c e^{bt} = \begin{pmatrix} c^1 e^{bt} \\ ... \\ c^n e^{bt} \end{pmatrix} \qquad c = \begin{pmatrix} c^1 \\ ... \\ c^n \end{pmatrix}$$

  where $A$-number $b$ is right $_*{}^*$-eigenvalue of the matrix $a$ and column vector $c$ is eigencolumn for right $_*{}^*$-eigenvalue $b$.

In non-commutative $D$-algebra A, we can consider a more general form of the system of linear differential equations

$$(1.6.5) \qquad \begin{aligned} \frac{dx^1}{dt} &= a_{1\,s0}^1 x^1 a_{1\,s1}^1 + ... + a_{n\,s0}^1 x^n a_{n\,s1}^1 \\ &\quad ...... \\ \frac{dx^n}{dt} &= a_{1\,s0}^n x^1 a_{1\,s1}^n + ... + a_{n\,s0}^n x^n a_{n\,s1}^n \end{aligned}$$

We can represent the system of differential equations (1.6.5) using product of matrices

$$(1.6.6) \quad \begin{pmatrix} \dfrac{dx^1}{dt} \\ ... \\ \dfrac{dx^n}{dt} \end{pmatrix} = \begin{pmatrix} a^1_{1\,s0} \otimes a^1_{1\,s1} & ... & a^1_{n\,s0} \otimes a^1_{n\,s1} \\ ... & ... & ... \\ a^n_{1\,s0} \otimes a^n_{1\,s1} & ... & a^n_{n\,s0} \otimes a^n_{n\,s1} \end{pmatrix} \circ^\circ \begin{pmatrix} x^1 \\ ... \\ x^n \end{pmatrix}$$

The solution of the system of differential equations (1.6.6) may have one of following representations.

- We can write solution as

$$(1.6.7) \quad x = Ce^{bt}c = Ce^{bt} \begin{pmatrix} c^1 \\ ... \\ c^n \end{pmatrix}$$

  where
  - $A$-number $b$ is left $_\circ{}^\circ$-eigenvalue of the matrix

$$a = \begin{pmatrix} a^1_{1\,s0} \otimes a^1_{1\,s1} & ... & a^1_{n\,s0} \otimes a^1_{n\,s1} \\ ... & ... & ... \\ a^n_{1\,s0} \otimes a^n_{1\,s1} & ... & a^n_{n\,s0} \otimes a^n_{n\,s1} \end{pmatrix}$$

    and the column $c$ is eigencolumn of the matrix $a$ corresponding to the left $_\circ{}^\circ$-eigenvalue $b$.
  - $A$-number $b$ and $A$-number $C$ satisfy the condition

$$(1.6.8) \quad b \in \bigcap_{i=1}^{n} \bigcap_{j=1}^{n} Z(A, a^i_{j\,s1})$$

$$(1.6.9) \quad C \in \bigcap_{i=1}^{n} \bigcap_{j=1}^{n} Z(A, a^i_{j\,s1})$$

- We can write solution as

$$(1.6.10) \quad x = ce^{bt}C = \begin{pmatrix} c^1 \\ ... \\ c^n \end{pmatrix} e^{bt}C$$

  where
  - $A$-number $b$ is right $_\circ{}^\circ$-eigenvalue of the matrix

$$a = \begin{pmatrix} a^1_{1\,s0} \otimes a^1_{1\,s1} & ... & a^1_{n\,s0} \otimes a^1_{n\,s1} \\ ... & ... & ... \\ a^n_{1\,s0} \otimes a^n_{1\,s1} & ... & a^n_{n\,s0} \otimes a^n_{n\,s1} \end{pmatrix}$$

    and the column $c$ is eigencolumn of the matrix $a$ corresponding to the right $_\circ{}^\circ$-eigenvalue $b$.

— $A$-number $b$ and $A$-number $C$ satisfy the condition

$$(1.6.11) \qquad b \in \bigcap_{i=1}^{n} \bigcap_{j=1}^{n} Z(A, a_{j\,s2}^{i})$$

$$(1.6.12) \qquad C \in \bigcap_{i=1}^{n} \bigcap_{j=1}^{n} Z(A, a_{j\,s2}^{i})$$

## 1.7. Polynomial

In commutative algebra polynomials have a lot of common with number theory: divisibility, greatest common factor, least common multiple, factorization into irreducible factors. When we study polynomials in non-commutative algebra, we want to preserve statements from commutative algebra as much as possible.

In the book [20] on the page 48, professor Cohn wrote that the study of polynomials over skew field is difficult problem. For this reason professor Ore considered skew right polynomial ring [1.5] where Ore introduced the transformation of monomial

$$(1.7.1) \qquad at = ta^{\alpha} + a^{\delta}$$

The map

$$\alpha : a \in A \to a^{\alpha} \in A$$

is endomorphism of $D$-algebra $A$ and the map

$$\delta : a \in A \to a^{\delta} \in A$$

is $\alpha$-derivation of $D$-algebra $A$

$$(1.7.2) \qquad (a+b)^{\delta} = a^{\delta} + b^{\delta} \qquad (ab)^{\delta} = a^{\delta}b^{\alpha} + ab^{\delta}$$

In such case all coefficients of polynomial can be written on the right. Such polynomials are called right-sided polynomials.

The same way we may consider the transformation of monomial

$$(1.7.3) \qquad ta = a^{\alpha}t + a^{\delta}$$

In such case all coefficients of polynomial can be written on the left. Such polynomials are called left-sided polynomials. [1.6]

The point of view proposed by Ore allows us to get qualitative picture of the theory of polynomials in non-commutative algebra. However, there exist algebras where this point of view is not complete.

THEOREM 1.7.1. *Let $A$ be division $D$-algebra. Let $a \in A \otimes A$. Let the equation*

$$(1.7.4) \qquad a \circ x + b = 0$$

*have infinitely many roots, but not any $A$-number is the root of the equation* (1.7.4). *Then ring of polynomials $A[x]$ is not Ore ring.*

PROOF. Let the transformation $\boxed{(1.7.1) \quad at = ta^{\alpha} + a^{\delta}}$ be defined in $D$ algebra $A$. Let

$$(1.7.5) \qquad a = a_{s0} \otimes a_{s1}$$

---

[1.5] See also the definition on pages [20]-(48 - 49).

[1.6] See also the definition on the page [16]-3.

The equality

$$(1.7.6) \qquad (a_{s0} \otimes a_{s1}) \circ x + b = a_{s0} x a_{s1} + b = (x a_{s0}^\alpha + a_{s0}^\delta) a_{s1} + b$$
$$= x a_{s0}^\alpha a_{s1} + a_{s0}^\delta a_{s1} + b = 0$$

follows from equalities (1.7.1), (1.7.4), (1.7.5). Let

$$A = a_{s0}^\alpha a_{s1} \qquad B = a_{s0}^\delta a_{s1} + b$$

The equation

$$(1.7.7) \qquad xA + B = 0$$

has the same roots as the equation (1.7.4).

- If $A \neq 0$, then the equation (1.7.7) has single root.
- If $A = 0$, $B \neq 0$, then the equation (1.7.7) has no root.
- If $A = 0$, $B = 0$, then any $A$-number is root of the equation (1.7.7).

All three options contradict the condition of the theorem. Therefore, the transformation $\boxed{(1.7.1) \quad at = t a^\alpha + a^\delta}$ is not defined in $D$ algebra $A$. $\qquad\square$

---

THEOREM 1.7.2. *General solution of the linear equation*

$$(1.7.8) \qquad ix - xi = k$$

*in quaternion algebra is*

$$(1.7.9) \qquad x = C_1 + C_2 i + \frac{1}{2} j \qquad C_1, C_2 \in R$$

PROOF. The tensor

$$i \otimes 1 - 1 \otimes i \in H^{2\otimes}$$

corresponds to the matrix

$$(1.7.10) \qquad \begin{pmatrix} 0 & -1 & 0 & 0 \\ 1 & 0 & 0 & 0 \\ 0 & 0 & 0 & -1 \\ 0 & 0 & 1 & 0 \end{pmatrix} - \begin{pmatrix} 0 & -1 & 0 & 0 \\ 1 & 0 & 0 & 0 \\ 0 & 0 & 0 & 1 \\ 0 & 0 & -1 & 0 \end{pmatrix} = \begin{pmatrix} 0 & 0 & 0 & 0 \\ 0 & 0 & 0 & 0 \\ 0 & 0 & 0 & -2 \\ 0 & 0 & 2 & 0 \end{pmatrix}$$

$\qquad\square$

---

THEOREM 1.7.3. *The ring of polynomials $H[x]$ in quaternion algebra is not Ore ring.*

PROOF. The theorem follows from theorems 1.7.1, 1.7.2. $\qquad\square$

In the section [19]-6.4, professor Cohn considered the definition of left-sided polynomial regardless of the definition of Ore ring.

Let $A$ be associative algebra.

On the page [19]-148, professor Cohn considered the set $A[x]^*$ of left-sided polynomials where left-sided polynomial is defined as formal expression

$$(1.7.11) \qquad f = a_0 + a_1 x + \ldots + a_n x^n$$

$$(1.7.12) \qquad g = b_0 + b_1 x + \ldots + b_m x^m$$

where $a_0, ..., a_n, b_0, ..., b_m$, are $A$-numbers and $x$ is variable.

Left-sided polynomials $f$ and $g$ are equal, if $n = m, a_0 = b_0, ..., a_n = b_n$.

We also define on the set $A[x]^*$ sum [1.7] of left-sided polynomials $f$ and $g$

$$(1.7.13) \qquad f + g = (a_0 + b_0) + (a_1 + b_1)x + ... + (a_n + b_n)x^n$$

and product

$$(1.7.14) \quad f * g = a_0 b_0 + (a_0 b_1 + a_1 b_0)x + (a_0 b_2 + a_1 b_1 + a_2 b_0)x^2 + ... + a_n b_m x^{n+m}$$

It is easy to verify that the set $A[x]^*$ is $D$-algebra. Evidently left-sided polynomial $f$ is polynomial $f$. However, in general, product of left-sided polynomials $f$ and $g$ does not equal to product of polynomials $f$ and $g$.

EXAMPLE 1.7.4. *Let*

$$(1.7.15) \qquad f = 1 + ix$$

$$(1.7.16) \qquad g = k + jx$$

*Equalities*

$$(1.7.17) \qquad f * g = k + (j + ik)x + ijx^2 = k + kx^2$$

$$(1.7.18) \quad fg = (1 + ix)(k + jx) = (k + jx) + ix(k + jx) = k + jx + ixk + ixjx$$

*follow from the equalities* (1.7.14), (1.7.15), (1.7.16)*. Let $x = k$. The equality*

$$(1.7.19) \qquad f(k) = 1 + ik = 1 + j$$

*follows from the equality* (1.7.15)*. The equality*

$$(1.7.20) \qquad g(k) = k + jk = k - i = k - i$$

*follows from the equality* (1.7.16)*. The equality*

$$(1.7.21) \quad f(k)g(k) = (1+j)(k-i) = k-i+j(k-i) = k-i+jk-ji = k-i+i+k = 2k$$

*follows from the equalities* (1.7.19), (1.7.20)*. The equality*

$$(1.7.22) \qquad f * g(k) = k + kk^2 = 0$$

*follows from the equality* (1.7.17)*. The equality*

$$(1.7.23) \qquad fg(k) = k + jk + ik^2 + ikjk = k + i - i + k = 2k$$

*follows from the equality* (1.7.18)*. Therefore,*

$$f * g(k) \neq fg(k) = f(k)g(k)$$

$\square$

## 1.8. Division in Non-commutative Algebra

To answer the question about dividing polynomials in non-commutative algebra, we have to go beyond algebra of polynomials.

Everything is simple in division algebra. The equation $ax = b$ has unique solution as soon as $a \neq 0$. However, we considered this question in the proof of the theorem 1.7.1. Division algebra is small county in vast kingdom of Algebra.

Consider any $D$-algebra $A$. Suppose we claim that $A$-number $b$ divides $A$-number $a$. What does this statement mean?

---

[1.7] If $m < n$, for instance, then we set $b_{m+1} = ... = b_n = 0$.

It would seem that there is simple answer.

**Definition 1.8.1.** *There exists A-number c such that*

$$a = bc$$

*In such case, A-number b is called left divisor of A-number a.* □

**Definition 1.8.2.** *There exists A-number c such that*

$$a = cb$$

*In such case, A-number b is called right divisor of A-number a.* □

But should we distinguish between left and right divisors? According to the definition, $D$-algebra $A$ is $D$-module $A$, equipped by **bilinear map**

$$f : A \times A \to A$$

called **product**. But if we assume that first argument of the map $f$ is given $A$-number, then bilinear map becomes linear map of second argument.

Therefore, we will say that $A$-number $b$ divides $A$-number $a$, if there exists linear map

$$c : A \to A$$

such that

$$a = c \circ b$$

For instance, we can rewrite definitions 1.8.1, 1.8.2 in the following way:

**Definition 1.8.3.** *There exists linear map $1 \otimes c$ such that $(1 \otimes c) \circ b = a$* □

**Definition 1.8.4.** *There exists linear map $c \otimes 1$ such that $(c \otimes 1) \circ b = a$* □

We will also say that $A$-number $b$ divides $A$-number $a$ with remainder, if there exists linear map

$$c : A \to A$$

and $A$-number $d$ such that

$$a = c \circ b + d$$

Division algorithm, algorithms to search greatest common factor or least common multiple are standard procedure.

The solution of one question very often raises new questions. Among the issues that will need to be resolved is the expansion of $A$-number using its factors.

For instance, what will be the result of division of $A$-number $a$ by $A$-numbers $b_1$ and $b_2$. To answer this question, we divide the task into several steps. [1.8]

The quotient of $A$-number $a$ divided by $A$-number $b_1$ has the following form

(1.8.1)                             $a = c_1 \circ b_1$

where linear map $c_1$ has the following form

(1.8.2)                             $c_1 = c_{10s} \otimes c_{11s}$

---

[1.8] To simplify the reasoning, we assume that at each stage we have division without reminder. We also consider associative $D$-algebra.

We do not have an algorith of dividing a tensor $c_1$ by $A$-number $b_2$. However the equality

$$(1.8.3) \qquad a = c_{10s}b_1c_{11s}$$

follows from equalities (1.8.1), (1.8.2). The equality (1.8.3) implies that if we want to divide $A$-number $a$ by $A$-numbers $b_1$ and $b_2$, then we have to divide every $A$-number $c_{10s}$, $c_{11s}$ by $A$-number $b_2$ and substitute the result into equality (1.8.3)

$$(1.8.4) \qquad a = (c_{10s2} \circ b_2)b_1c_{11s} + c_{10s}b_1(c_{11s2} \circ b_2)$$

However, the expression on the right side of the equality (1.8.4) is bilinear map of $A$-numbers $b_1$ and $b_2$

$$(1.8.5) \qquad \begin{aligned} &a = c \circ (b_1, b_2) \\ &c = c_{10s20t} \otimes_2 c_{10s20t} \otimes_1 c_{11s} + c_{10s} \otimes_1 c_{11s20t} \otimes_2 c_{11s21t} \end{aligned}$$

Gradually increasing the number of factors, we come to the conclusion that if $A$-number $a$ is product of $A$-numbers $b_1$, ..., $b_n$, then $A$-number $a$ is polylinear map of $A$-numbers $b_1$, ..., $b_n$.

Now we can return to polynomials. In order to find factorization of a polynomial into irreducible factors, we have to answer some questions. How to find the set of roots? How big is the set of factorizations, if the set of roots is infinite? How we can find factor if this factor does not have roots?

### 1.9. Inverse Matrix

This topic came up in the process of writing the book.

In this book, I considered definition of left and right eigenvalues of $n \times n$ matrix $a_2$ in case when the matrix $a_2$ has $n$ different eigenvalues.

In such case, the set of eigenvectors has two equivalent definitions

$$(18.1.1) \quad u_{2*}{}^*a_{2*}{}^*u_2^{-1}{}^* = a_1$$

and

$$(18.1.2) \quad u_{2*}{}^*a_2 = a_{1*}{}^*u_2$$

In such case, the set of eigenvectors has two equivalent definitions

$$(18.1.4) \quad u_2^{-1}{}^*{}_*{}^*a_{2*}{}^*u_2 = a_1$$

and

$$(18.1.5) \quad a_{2*}{}^*u_2 = u_{2*}{}^*a_1$$

Matrices $u_2$, $a_1$ are also $n \times n$ matrices. In such case, sets of left and right $_*{}^*$-eigenvalues are the same.

However, professor Cohn pointed out that there are cases when the set of keft $_*{}^*$-eigenvalues is different from the set of right $_*{}^*$-eigenvalues. To find such matrix, I decided to consider the case when the matrix $a_2$ has $k$, $k < n$, different eigenvalues. In this case I have to confine myself to the definition which is similar to the definition

$$(18.1.2) \quad u_{2*}{}^*a_2 = a_{1*}{}^*u_2$$

In such case, the matrix $u_2$ is $n \times k$ matrix and the matrix $a_1$ is $k \times k$ matrix. Trying to answer the question can I write a definition similar to the definition (18.1.1), I realized that I need to consider the concept of right $_*{}^*$-inverse matrix.

### 1.10. Non-commutative Sum

Since our childhood, we know that sum does not depend on order of summands. So, when I learned about non-commutative addition ([24]), I was little bit confused. However I was lucky to see geometry where sum of vectors is non-commutative.

At school I liked many subjects: physics, chemistry, drafting. I did not think about mathematics seriously. I did not like arithmetic, but I liked algebra and geometry. I did not yet know that all roads lead to mathematics.

Drafting was my favorite subject. My teacher, Gabis Sergey Aleksandrovich, gave me problem book in drafting; I solved every problem from this problem book. After this I introduced my problem and solved it. To my surprise, the first theorem which we proved on geometry class next year was exactly the solution to this problem.

I knew that it will be hard to enter college after school. So I decided to enter technical school after eight class. Bela Markovna offered me few lessons during summer; she wanted to see how far I was prepared for the exam on math. I remember one lesson. When we considered line and parabola graphs, Bela Markovna said: it is easy to draw straight line or parabola; however, what can we do with cardiogram? I decided to find an analytical expression for cardiogram.

My friend offered me to study calculus together during summer. Boys who studied at a technical school lived in our building; they gave me calculus textbook for summer time. At this moment my friend lost interest to his idea; and I began to study calculus independently.

The first chapter of the textbook was dedicated to analytic geometry. So I realized that this is the tool that I need to analyze curves of cardiogram. analytic geometry determined my future. I realized that I will prepare myself for math department of university.

I have read a lot of books dedicated to calculus. However, at some point of time, I got the impression that one author copied off from another. I started to look for different books. I do not remember where did I get the first volume of Fihtengolts. But this book impressed me very much. Closer to spring, I felt that the book becomes difficult; however I continued to read this book. In March, I had the flu with an unusually high temperature. But I doubt that it was a flu, because it became easier to read the first volume of Fihtengolts.

The study of calculus reflected in my attempts to restore functional dependence on graph. I did not yet understand that not every map could be represented analytically. I drew graphs, made tables of dependence of function increment on increment of argument. However soon this idea bothered me; and I changed direction of my learning.

When I started to read the book dedicated to tensor calculus, I discovered that simple theorems of linear algebra substituted complex calculations of analytic geometry related to classification of second-order curves. I lost my interest to analytic geometry after this.

Almost half a century has passed. And suddenly, analytic geometry brought me two pleasant surprises. I discovered that analytic geometry has interesting structure from the point of view of the theory of representation of universal algebra. I identify analytic and affine geometries; however, I think that this is not big mistake.

The second surprise is related to the fact that I realized that I have model of affine geometry in Riemann space. Since I studied metric-affine manifold for a very long time, then there was one step left to see model of affine geometry where sum is non-commutative. Henceforth non-commutative addition ceased to be abstract concept for me.

# Representation of Universal Algebra

This chapter contains definitions and theorems which are necessary for an understanding of the text of this book. So the reader may read the statements from this chapter in process of reading the main text of the book.

## 2.1. Universal Algebra

DEFINITION 2.1.1. *For the map*

$$f : A \to B$$

*the set*

$$\operatorname{Im} f = \{f(a) : a \in A\}$$

*is called* **image of map** *$f$.* $\square$

DEFINITION 2.1.2. *For any sets* [2.1]*$A$, $B$,* **Cartesian power** *$B^A$ is the set of maps*

$$f : A \to B$$

$\square$

DEFINITION 2.1.3. *For any $n \ge 0$, a map* [2.2]

$$\omega : A^n \to A$$

*is called $n$-**ary operation on set** $A$ or just **operation on set** $A$. For any $a_1$, ..., $a_n \in A$, we use either notation $\omega(a_1, ..., a_n)$, $a_1...a_n\omega$ to denote image of map $\omega$.* $\square$

REMARK 2.1.4. *According to definitions 2.1.2, 2.1.3, $n$-ari operation $\omega \in A^{A^n}$.* $\square$

DEFINITION 2.1.5. **Operator domain** *is the set of operators* [2.3] *$\Omega$ with a map*

$$a : \Omega \to N$$

*If $\omega \in \Omega$, then $a(\omega)$ is called the **arity** of operator $\omega$. If $a(\omega) = n$, then operator $\omega$ is called $n$-ary. We use notation*

$$\Omega(n) = \{\omega \in \Omega : a(\omega) = n\}$$

*for the set of $n$-ary operators.* $\square$

DEFINITION 2.1.6. *Let $A$ be a set. Let $\Omega$ be an operator domain.* [2.4]*The family of maps*

$$\Omega(n) \to A^{A^n} \quad n \in N$$

---

*is called $\Omega$-algebra structure on $A$. The set $A$ with $\Omega$-algebra structure is called $\Omega$-algebra* $A_\Omega$ *or* **universal algebra**. *The set $A$ is called* **carrier of $\Omega$-algebra**.

<div style="text-align: right">□</div>

The operator domain $\Omega$ describes a set of $\Omega$-algebras. An element of the set $\Omega$ is called operator, because an operation assumes certain set. According to the remark 2.1.4 and the definition 2.1.6, for each operator $\omega \in \Omega(n)$, we match $n$-ary operation $\omega$ on $A$.

DEFINITION 2.1.7. *Let $A$, $B$ be $\Omega$-algebras and $\omega \in \Omega(n)$. The map* [2.5]

$$f : A \to B$$

**is compatible with operation** $\omega$, *if, for all $a_1$, ..., $a_n \in A$,*

$$f(a_1...a_n\omega) = f(a_1)...f(a_n)\omega$$

*The map $f$ is called* **homomorphism** *from $\Omega$-algebra $A$ to $\Omega$-algebra $B$, if $f$ is compatible with each $\omega \in \Omega$.*

<div style="text-align: right">□</div>

DEFINITION 2.1.8. *Homomorphism*

$$f : A \to B$$

*is called* [2.6] **isomorphism** *between $A$ and $B$, if correspondence $f^{-1}$ is homomorphism. If there is an isomorphism between $A$ and $B$, then we say that $A$ and $B$ are isomorphic and write* $A \cong B$. *An injective homomorphis is called* **monomorphism**. *A surjective homomorphis is called* **epimorphism**.

<div style="text-align: right">□</div>

DEFINITION 2.1.9. *A homomorphism*

$$f : A \to A$$

*in which source and target are the same algebra is called* **endomorphism**. *We use notation* $\mathrm{End}(\Omega; A)$ *for the set of endomorphisms of $\Omega$-algebra $A$. Isomorphism*

$$f : A \to A$$

*is called* **automorphism**.

<div style="text-align: right">□</div>

THEOREM 2.1.10. *Let the map*

$$f : A \to B$$

*be homomorphism of $\Omega$-algebra $A$ into $\Omega$-algebra $B$. Let the map*

$$g : B \to C$$

*be homomorphism of $\Omega$-algebra $B$ into $\Omega$-algebra $C$. Then the map* [2.7]

$$h = g \circ f : A \to C$$

*is homomorphism of $\Omega$-algebra $A$ into $\Omega$-algebra $C$.*

---

[2.4] I follow the definition (2), page [18]-48.
[2.5] I follow the definition on page [18]-49.
[2.6] I follow the definition on page [18]-49.

PROOF. Let $\omega \in \Omega_n$ be n-ary operation. Then

$$(2.1.1) \qquad f(a_1...a_n\omega) = f(a_1)...f(a_n)\omega$$

for any $a_1, ..., a_n \in A$, and

$$(2.1.2) \qquad g(b_1...b_n\omega) = g(b_1)...g(b_n)\omega$$

for any $b_1, ..., b_n \in B$. The equality

$$(2.1.3) \quad \begin{aligned} (g \circ f)(a_1...a_n\omega) &= g(f(a_1...a_n\omega)) \\ &= g(f(a_1)...f(a_n)\omega) \\ &= g(f(a_1))...g(f(a_n))\omega \\ &= (g \circ f)(a_1)...(g \circ f)(a_n)\omega \end{aligned}$$

follows from equalities (2.1.1), (2.1.2). The theorem follows from the equality (2.1.3) for any operation $\omega$. $\square$

DEFINITION 2.1.11. *If there is a monomorphism from $\Omega$-algebra $A$ to $\Omega$-algebra $B$, then we say that $A$ **can be embedded in** $B$.* $\square$

DEFINITION 2.1.12. *If there is an epimorphism from $A$ to $B$, then $B$ is called* **homomorphic image** *of algebra $A$.* $\square$

## 2.2. Cartesian Product of Universal Algebras

DEFINITION 2.2.1. *Let $\mathcal{C}$ be a category. An object $P$ of category $\mathcal{C}$ is called* **universally attracting** [2.8]*if, for any object $R$ of category $\mathcal{C}$, there exists unique morphism*

$$f : R \to P$$

$\square$

DEFINITION 2.2.2. *Let $\mathcal{C}$ be a category. An object $P$ of category $\mathcal{C}$ is called* **universally repelling** [2.9]*if, for any object $R$ of category $\mathcal{C}$, there exists unique morphism*

$$f : P \to R$$

$\square$

DEFINITION 2.2.3. *Let $\mathcal{A}$ be a category. Let $\{B_i, i \in I\}$ be the set of objects of $\mathcal{A}$. Let $\mathcal{A}(B_i, i \in I)$ be a category whose objects are tuples $(P, f)$ where $P$ is object of category $\mathcal{A}$ and $f$ is set of morphisms*

$$\{f_i : P \to B_i, i \in I\}$$

---

[2.7] I follow the proposition [18]-3.2 on page [18]-57.

[2.8] See also the definition in [1], pages 57.

[2.9] See also the definition in [1], pages 57.

*Universally attracting object of category* $\mathcal{A}(B_i, i \in I)$

$$P = \prod_{i \in I} B_i$$

*is called a* **product of set of objects** $\{B_i, i \in I\}$ **in category** $\mathcal{A}$.[2.10]

*If* $|I| = n$, *then we also will use notation*

$$P = \prod_{i=1}^{n} B_i = B_1 \times ... \times B_n$$

*for product of set of objects* $\{B_i, i \in I\}$ *in* $\mathcal{A}$. $\qquad\qquad\square$

THEOREM 2.2.4. *Let* $\mathcal{A}$ *be a category. Let* $\{B_i, i \in I\}$ *be the set of objects of* $\mathcal{A}$. *The product of set of objects* $\{B_i, i \in I\}$ *in category* $\mathcal{A}$ *is an object*

$$P = \prod_{i \in I} B_i$$

*and set of morphisms*

$$\{f_i : P \to B_i, i \in I\}$$

*such that for any object* $R$ *and set of morphisms*

$$\{g_i : R \to B_i, i \in I\}$$

*there exists a unique morphism*

$$h : R \to P$$

*such that diagram*

$$f_i \circ h = g_i$$

*is commutative for all* $i \in I$.

PROOF. The theorem follows from definitions 2.2.1, 2.2.3. $\qquad\qquad\square$

EXAMPLE 2.2.5. *Let* $\mathcal{S}$ *be the category of sets.*[2.11] *According to the definition 2.2.3, Cartesian product*

$$A = \prod_{i \in I} A_i$$

*of family of sets* $(A_i, i \in I)$ *and family of projections on the i-th factor*

$$p_i : A \to A_i$$

*are product in the category* $\mathcal{S}$. $\qquad\qquad\square$

THEOREM 2.2.6. *The product exists in the category* $\mathcal{A}$ *of* $\Omega$-*algebras. Let* $\Omega$-*algebra* $A$ *and family of morphisms*

$$p_i : A \to A_i \quad i \in I$$

---

[2.10] I made definition according to the definition on page [1]-58.

[2.11] See also the example in [1], page 59.

be product in the category $\mathcal{A}$. Then

2.2.6.1: *The set $A$ is Cartesian product of family of sets  $(A_i, i \in I)$*

2.2.6.2: *The homomorphism of $\Omega$-algebra*

$$p_i : A \to A_i$$

*is projection on $i$-th factor.*

2.2.6.3: *We can represent any $A$-number $a$ as tuple  $(p_i(a), i \in I)$  of $A_i$-numbers.*

2.2.6.4: *Let  $\omega \in \Omega$  be $n$-ary operation. Then operation $\omega$ is defined component-wise*

(2.2.1) $$a_1...a_n\omega = (a_{1i}...a_{ni}\omega, i \in I)$$

*where  $a_1 = (a_{1i}, i \in I),\ ...,\ a_n = (a_{ni}, i \in I)$ .*

PROOF. Let

$$A = \prod_{i \in I} A_i$$

be Cartesian product of family of sets  $(A_i, i \in I)$  and, for each $i \in I$, the map

$$p_i : A \to A_i$$

be projection on the $i$-th factor. Consider the diagram of morphisms in category of sets $\mathcal{S}$

(2.2.2)
$$\begin{array}{ccc} A & \xrightarrow{p_i} & A_i \\ \uparrow\scriptstyle\omega & \nearrow\scriptstyle g_i & \\ A^n & & \end{array} \qquad p_i \circ \omega = g_i$$

where the map $g_i$ is defined by the equality

$$g_i(a_1, ..., a_n) = p_i(a_1)...p_i(a_n)\omega$$

According to the definition   2.2.3,  the map $\omega$ is defined uniquely from the set of diagrams (2.2.2)

(2.2.3) $$a_1...a_n\omega = (p_i(a_1)...p_i(a_n)\omega, i \in I)$$

The equality (2.2.1) follows from the equality (2.2.3). □

DEFINITION 2.2.7. *If $\Omega$-algebra $A$ and family of morphisms*

$$p_i : A \to A_i \quad i \in I$$

*is product in the category $\mathcal{A}$, then $\Omega$-algebra $A$ is called* **direct** *or* **Cartesian product of $\Omega$-algebras**  $(A_i, i \in I)$. □

THEOREM 2.2.8. *Let set $A$ be Cartesian product of sets  $(A_i, i \in I)$  and set $B$ be Cartesian product of sets  $(B_i, i \in I)$.  For each $i \in I$, let*

$$f_i : A_i \to B_i$$

*be the map from the set $A_i$ into the set $B_i$. For each $i \in I$, consider commutative diagram*

(2.2.4)

$$
\begin{array}{ccc}
B & \xrightarrow{\;p_i'\;} & B_i \\[2pt]
\big\uparrow{\scriptstyle f} & & \big\uparrow{\scriptstyle f_i} \\[2pt]
A & \xrightarrow{\;p_i\;} & A_i
\end{array}
$$

*where maps $p_i$, $p_i'$ are projection on the $i$-th factor. The set of commutative diagrams (2.2.4) uniquely defines map*

$$f : A \to B$$

$$f(a_i, i \in I) = (f_i(a_i), i \in I)$$

PROOF. For each $i \in I$, consider commutative diagram

(2.2.5)

Let $a \in A$. According to the statement 2.2.6.3, we can represent $A$-number $a$ as tuple of $A_i$-numbers

(2.2.6)          $a = (a_i, i \in I) \quad a_i = p_i(a) \in A_i$

Let

(2.2.7)                          $b = f(a) \in B$

According to the statement 2.2.6.3, we can represent $B$-number $b$ as tuple of $B_i$-numbers

(2.2.8)          $b = (b_i, i \in I) \quad b_i = p_i'(b) \in B_i$

From commutativity of diagram (1) and from equalities (2.2.7), (2.2.8), it follows that

(2.2.9)                          $b_i = g_i(b)$

From commutativity of diagram (2) and from the equality (2.2.6), it follows that

$$b_i = f_i(a_i)$$

<div align="right">□</div>

THEOREM 2.2.9. *Let $\Omega$-algebra $A$ be Cartesian product of $\Omega$-algebras $(A_i, i \in I)$ and $\Omega$-algebra $B$ be Cartesian product of $\Omega$-algebras $(B_i, i \in I)$. For each $i \in I$, let the map*

$$f_i : A_i \to B_i$$

*be homomorphism of $\Omega$-algebra. Then the map*

$$f : A \to B$$

*defined by the equality*

(2.2.10) $$f(a_i, i \in I) = (f_i(a_i), i \in I)$$

*is homomorphism of $\Omega$-algebra.*

PROOF. Let $\omega \in \Omega$ be n-ary operation. Let $a_1 = (a_{1i}, i \in I)$, ..., $a_n = (a_{ni}, i \in I)$ , $b_1 = (b_{1i}, i \in I)$, ..., $b_n = (b_{ni}, i \in I)$ . From equalities (2.2.1), (2.2.10), it follows that

$$
\begin{aligned}
f(a_1...a_n\omega) &= f(a_{1i}...a_{ni}\omega, i \in I) \\
&= (f_i(a_{1i}...a_{ni}\omega), i \in I) \\
&= ((f_i(a_{1i}))...(f_i(a_{ni})), i \in I) \\
&= (b_{1i}...b_{ni}\omega, i \in I)
\end{aligned}
$$

$$f(a_1)...f(a_n)\omega = b_1...b_n\omega = (b_{1i}...b_{ni}\omega, i \in I)$$

$\square$

DEFINITION 2.2.10. *Equivalence on $\Omega$-algebra $A$, which is subalgebra of $\Omega$-algebra $A^2$, is called* **congruence** [2.12] *on $A$.* $\square$

THEOREM 2.2.11 (first isomorphism theorem). *Let*

$$f : A \to B$$

*be homomorphism of $\Omega$-algebras with kernel s. Then there is decomposition*

$$A/\ker f \xrightarrow{q} f(A) \qquad f = p \circ q \circ r$$

2.2.11.1: *The* **kernel of homomorphism** $\ker f = f \circ f^{-1}$ *is a congruence on $\Omega$-algebra $A$.*

2.2.11.2: *The set $A/\ker f$ is $\Omega$-algebra.*

2.2.11.3: *The map*

$$p : a \in A \to a^{\ker f} \in A/\ker f$$

*is epimorphism and is called* **natural homomorphism**.

2.2.11.4: *The map*

$$q : p(a) \in A/\ker f \to f(a) \in f(A)$$

*is the isomorphism*

2.2.11.5: *The map*

$$r : f(a) \in f(A) \to f(a) \in B$$

*is the monomorphism*

---

[2.12] I follow the definition on page [18]-57.

PROOF. The statement 2.2.11.1 follows from the proposition II.3.4 ([18], page 58). Statements 2.2.11.2, 2.2.11.3 follow from the theorem II.3.5 ([18], page 58) and from the following definition. Statements 2.2.11.4, 2.2.11.5 follow from the theorem II.3.7 ([18], page 60). □

DEFINITION 2.2.12. *Let $\mathcal{A}$ be a category. Let $\{B_i, i \in I\}$ be the set of objects of $\mathcal{A}$. Let $\mathcal{A}[B_i, i \in I]$ be a category whose objects are tuples $(P, f)$ where $P$ is object of category $\mathcal{A}$ and $f$ is set of morphisms*

$$\{f_i : B_i \to P, i \in I\}$$

*Universally repelling object of category $\mathcal{A}[B_i, i \in I]$*

$$P = \coprod_{i \in I} B_i$$

*is called a* **coproduct of set of objects** $\{B_i, i \in I\}$ **in category** $\mathcal{A}$.[2.13]

*If $|I| = n$, then we also will use notation*

$$P = \coprod_{i=1}^{n} B_i = B_1 \coprod \cdots \coprod B_n$$

*for coproduct of set of objects $\{B_i, i \in I\}$ in $\mathcal{A}$.* □

THEOREM 2.2.13. *Let $\mathcal{A}$ be a category. Let $\{B_i, i \in I\}$ be the set of objects of $\mathcal{A}$. Object*

$$P = \coprod_{i \in I} B_i$$

*and set of morphisms*

$$\{f_i : B_i \to P, i \in I\}$$

*is called a* **coproduct of set of objects** $\{B_i, i \in I\}$ **in category** $\mathcal{A}$[2.14]*if for any object $R$ and set of morphisms*

$$\{g_i : B_i \to R, i \in I\}$$

*there exists a unique morphism*

$$h : P \to R$$

*such that diagram*

$$\begin{array}{ccc} P & \xleftarrow{f_i} & B_i \\ {\scriptstyle h}\downarrow & \swarrow{\scriptstyle g_i} & \\ R & & \end{array} \qquad h \circ f_i = g_i$$

*is commutative for all $i \in I$.*

PROOF. The theorem follows from definitions 2.2.2, 2.2.12. □

---

[2.13]  I made definition according to the definition on page [1]-59.

[2.14]  I made definition according to [1], page 59.

## 2.3. Representation of Universal Algebra

DEFINITION 2.3.1. *Let the set $A_2$ be $\Omega_2$-algebra. Let the set of transformations* End$(\Omega_2, A_2)$ *be $\Omega_1$-algebra. The homomorphism*

$$f : A_1 \to \text{End}(\Omega_2; A_2)$$

*of $\Omega_1$-algebra $A_1$ into $\Omega_1$-algebra* End$(\Omega_2, A_2)$ *is called* **representation of $\Omega_1$-algebra** $A_1$ *or* **$A_1$-representation** *in $\Omega_2$-algebra $A_2$.* $\square$

We also use notation

$$f : A_1 \longrightarrow_* A_2$$

to denote the representation of $\Omega_1$-algebra $A_1$ in $\Omega_2$-algebra $A_2$.

DEFINITION 2.3.2. *Let the map*

$$f : A_1 \to \text{End}(\Omega_2; A_2)$$

*be an isomorphism of the $\Omega_1$-algebra $A_1$ into* End$(\Omega_2, A_2)$. *Then the representation*

$$f : A_1 \longrightarrow_* A_2$$

*of the $\Omega_1$-algebra $A_1$ is called* **effective**.[2.15] $\square$

DEFINITION 2.3.3. *The representation*

$$g : A_1 \longrightarrow_* A_2$$

*of the $\Omega_1$-algebra $A_1$ is called* **free**,[2.16] *if the statement*

$$f(a_1)(a_2) = f(b_1)(a_2)$$

*for any $a_2 \in A_2$ implies that $a = b$.* $\square$

THEOREM 2.3.4. *Free representation is effective.*

PROOF. The theorem follows from the theorem [15]-3.1.6. $\square$

REMARK 2.3.5. *Representation of rotation group in affine space is effective. However this representation is not free, since origin is fixed point of every transformation.* $\square$

DEFINITION 2.3.6. *The representation*

$$g : A_1 \longrightarrow_* A_2$$

*of $\Omega_1$-algebra is called* **transitive**[2.17] *if, for any $a$, $b \in 2$, there exists such $g$ that*

$$a = f(g)(b)$$

---

[2.15] See similar definition of effective representation of group in [21], page 16, [22], page 111, [19], page 51 (Cohn calls such representation faithful).
[2.16] See similar definition of free representation of group in [21], page 16.

*The representation of $\Omega_1$-algebra is called **single transitive** if it is transitive and free.*                                                                  □

THEOREM 2.3.7. *Representation is single transitive iff for any $a, b \in A_2$ exists one and only one $g \in A_1$ such that $a = f(g)(b)$.*

PROOF. The theorem follows from definitions 2.3.3 and 2.3.6.              □

THEOREM 2.3.8. *Let*

$$f : A_1 \longrightarrow A_2$$

*be a single transitive representation of $\Omega_1$-algebra $A_1$ in $\Omega_2$-algebra $A_2$. There is the structure of $\Omega_1$-algebra on the set $A_2$.*

PROOF. The theorem follows from the theorem [15]-3.1.10.                □

DEFINITION 2.3.9. *Let*

$$f : A_1 \longrightarrow A_2$$

*be representation of $\Omega_1$-algebra $A_1$ in $\Omega_2$-algebra $A_2$ and*

$$g : B_1 \longrightarrow B_2$$

*be representation of $\Omega_1$-algebra $B_1$ in $\Omega_2$-algebra $B_2$. For $i = 1, 2$, let the map*

$$r_i : A_i \to B_i$$

*be homomorphism of $\Omega_i$-algebra. The tuple of maps $r = (r_1, r_2)$ such, that*

$$(2.3.1) \qquad r_2 \circ f(a) = g(r_1(a)) \circ r_2$$

*is called **morphism of representations from $f$ into $g$**. We also say that **morphism of representations of $\Omega_1$-algebra in $\Omega_2$-algebra** is defined.*          □

REMARK 2.3.10. *We may consider a pair of maps $r_1, r_2$ as map*

$$F : A_1 \cup A_2 \to B_1 \cup B_2$$

*such that*

$$F(A_1) = B_1 \qquad F(A_2) = B_2$$

*Therefore, hereinafter the tuple of maps $r = (r_1, r_2)$ also is called map and we will use map*

$$r : f \to g$$

*Let $a = (a_1, a_2)$ be tuple of A-numbers. We will use notation*

$$r(a) = (r_1(a_1), r_2(a_2))$$

*for image of tuple of A-numbers with respect to morphism of representations $r$.*          □

REMARK 2.3.11. *Consider morphism of representations*

$$(r_1 : A_1 \to B_1, \ r_2 : A_2 \to B_2)$$

---

2.17 See similar definition of transitive representation of group in [22], page 110, [19], page 51.

*We denote elements of the set $B_1$ by letter using pattern $b \in B_1$. However if we want to show that $b$ is image of element $a \in A_1$, we use notation $r_1(a)$. Thus equation*

$$r_1(a) = r_1(a)$$

*means that $r_1(a)$ (in left part of equation) is image $a \in A_1$ (in right part of equation). Using such considerations, we denote element of set $B_2$ as $r_2(m)$. We will follow this convention when we consider correspondences between homomorphisms of $\Omega_1$-algebra and maps between sets where we defined corresponding representations.*          □

THEOREM 2.3.12. *Let*

$$f : A_1 \longrightarrow\!* A_2$$

*be representation of $\Omega_1$-algebra $A_1$ in $\Omega_2$-algebra $A_2$ and*

$$g : B_1 \longrightarrow\!* B_2$$

*be representation of $\Omega_1$-algebra $B_1$ in $\Omega_2$-algebra $B_2$. The map*

$$(r_1 : A_1 \to B_1, \ r_2 : A_2 \to B_2)$$

*is morphism of representations iff*

$$(2.3.2) \qquad r_2(f(a)(m)) = g(r_1(a))(r_2(m))$$

PROOF. The theorem follows from the theorem [15]-3.2.5.          □

REMARK 2.3.13. *There are two ways to interpret* (2.3.2)

- *Let transformation $f(a)$ map $m \in A_2$ into $f(a)(m)$. Then transformation $g(r_1(a))$ maps $r_2(m) \in B_2$ into $r_2(f(a)(m))$.*
- *We represent morphism of representations from $f$ into $g$ using diagram*

(2.3.3)

*From* (2.3.1), *it follows that diagram* (1) *is commutative.*

*We also use diagram*

(2.3.4)

*instead of diagram* (2.3.3).          □

DEFINITION 2.3.14. *If representation $f$ and $g$ coincide, then morphism of representations $r = (r_1, r_2)$ is called* **morphism of representation** $f$.    □

DEFINITION 2.3.15. *Let*

$$f : A_1 \longrightarrow_* A_2$$

*be representation of $\Omega_1$-algebra $A_1$ in $\Omega_2$-algebra $A_2$ and*

$$g : A_1 \longrightarrow_* B_2$$

*be representation of $\Omega_1$-algebra $A_1$ in $\Omega_2$-algebra $B_2$. Let*

$$(\mathrm{id} : A_1 \to A_1, \ r_2 : A_2 \to B_2)$$

*be morphism of representations. In this case we identify morphism* $(\mathrm{id}, r_2)$ *of representations of $\Omega_1$-algebra and corresponding homomorphism $r_2$ of $\Omega_2$-algebra and the homomorphism $r_2$ is called* **reduced morphism of representations**. *We will use diagram*

(2.3.5)

*to represent reduced morphism $r_2$ of representations of $\Omega_1$-algebra. From diagram it follows*

(2.3.6)
$$r_2 \circ f(a) = g(a) \circ r_2$$

*We also use diagram*

*instead of diagram* (2.3.5).    □

THEOREM 2.3.16. *Let*

$$f : A_1 \longrightarrow_* A_2$$

*be representation of $\Omega_1$-algebra $A_1$ in $\Omega_2$-algebra $A_2$ and*

$$g : A_1 \longrightarrow_* B_2$$

*be representation of $\Omega_1$-algebra $A_1$ in $\Omega_2$-algebra $B_2$. The map*

$$r_2 : A_2 \to B_2$$

*is reduced morphism of representations iff*

$$(2.3.7) \qquad\qquad r_2(f(a)(m)) = g(a)(r_2(m))$$

PROOF. The theorem follows from the theorem [15]-3.4.3. □

DEFINITION 2.3.17. *Let $A_1$ be $\Omega_1$-algebra. Let $\mathcal{A}_2$ be category of $\Omega_2$-algebras. We define* **category $A_1(\mathcal{A}_2)$ of representations** *of $\Omega_1$-algebra $A_1$ in $\Omega_2$-algebra. Representations of $\Omega_1$-algebra $A_1$ in $\Omega_2$-algebra are objects of this category. Reduced morphisms of corresponding representations are morphisms of this category.* □

THEOREM 2.3.18. *In category $A_1(\mathcal{A}_2)$ there exists product of effective representations of $\Omega_1$-algebra $A_1$ in $\Omega_2$-algebra and the product is effective representation of $\Omega_1$-algebra $A_1$.*

PROOF. The theorem follows from the theorem [15]-4.4.2. □

## 2.4. Ω-Group

DEFINITION 2.4.1. *Let sum be defined in $\Omega$-algebra $A$. A map*

$$f : A \to A$$

*of $\Omega_1$-algebra $A$ is called* **additive map** *if*

$$f(a + b) = f(a) + f(b)$$

□

DEFINITION 2.4.2. *A map*

$$g : A^n \to A$$

*is called* **polyadditive map** *if for any $i$, $i = 1, ..., n$,*

$$f(a_1, ..., a_i + b_i, ..., a_n) = f(a_1, ..., a_i, ..., a_n) + f(a_1, ..., b_i, ..., a_n)$$

□

DEFINITION 2.4.3. *Let sum which is not necessarily commutative be defined in $\Omega_1$-algebra $A$. We use the symbol $+$ to denote sum. Let*

$$\Omega = \Omega_1 \setminus \{+\}$$

*If $\Omega_1$-algebra $A$ is group relative to sum and any operation $\omega \in \Omega$ is polyadditive map, then $\Omega_1$-algebra $A$ is called* **$\Omega$-group**. *If $\Omega$-group $A$ is associative group relative to sum, then $\Omega_1$-algebra $A$ is called* **associative $\Omega$-group**. *If $\Omega$-group $A$ is Abelian group relative to sum, then $\Omega_1$-algebra $A$ is called* **Abelian $\Omega$-group**. □

THEOREM 2.4.4. *Let $\omega \in \Omega(n)$ be polyadditive map. The operation $\omega$ is* **distributive** *over addition*

$$a_1...(a_i + b_i)...a_n\omega = a_1...a_i...a_n\omega + a_1...b_i...a_n\omega \quad i = 1, ..., n$$

Proof. The theorem follows from definitions 2.4.2, 2.4.3.                    □

Definition 2.4.5. *Let product*

$$c_1 = a_1 * b_1$$

*be operation of* $\Omega_1$-*algebra* $A$. *Let* $\Omega = \Omega_1 \setminus \{*\}$. *If* $\Omega_1$-*algebra* $A$ *is group with respect to product and, for any operation* $\omega \in \Omega(n)$, *the product is distributive over the operation* $\omega$

$$a * (b_1...b_n\omega) = (a * b_1)...(a * b_n)\omega$$

$$(b_1...b_n\omega) * a = (b_1 * a)...(b_n * a)\omega$$

*then* $\Omega_1$-*algebra* $A$ *is called* **multiplicative** $\Omega$-**group**.                    □

Definition 2.4.6. *If*

(2.4.1)                    $$a * b = b * a$$

*then multiplicative* $\Omega$-*group is called* **Abelian**.                    □

Definition 2.4.7. *Let sum*

$$c_1 = a_1 + b_1$$

*and product*

$$c_1 = a_1 * b_1$$

*be operations of* $\Omega_1$-*algebra* $A$. *Let* $\Omega = \Omega_1 \setminus \{+, *\}$. *If* $\Omega_1$-*algebra* $A$ *is Abelian* $\Omega \cup \{*\}$-*group and multiplicative* $\Omega \cup \{+\}$-*group, then* $\Omega_1$-*algebra* $A$ *is called* $\Omega$-**ring**.                    □

Theorem 2.4.8. *The product in* $\Omega$-*ring* **is distributive** *over addition*

$$a * (b_1 + b_2) = a * b_1 + a * b_2$$

$$(b_1 + b_2) * a = b_1 * a + b_2 * a$$

Proof. The theorem follows from the definitions 2.4.3, 2.4.5, 2.4.7.                    □

## 2.5. Representation of Multiplicative $\Omega$-Group

We assume that transformations of representation

$$g : A_1 \longrightarrow A_2$$

of multiplicative $\Omega$-group $A_1$ in $\Omega_2$-algebra $A_2$ may act on $A_2$-numbers either on the left or on the right. We must consider an endomorphism of $\Omega_2$-algebra $A_2$ as operator. In such case, the following conditions must be satisfied.

Definition 2.5.1. *Let* $\mathrm{End}(\Omega_2, A_2)$ *be a multiplicative* $\Omega$-*group with product*[2.18]

$$(f, g) \to f \bullet g$$

*Let an endomorphism* $f$ *act on* $A_2$-*number* $a$ *on the left. We will use notation*

(2.5.1)                    $$f(a_2) = f \bullet a_2$$

*Let $A_1$ be multiplicative $\Omega$-group with product*

$$(a, b) \to a * b$$

*We call a homomorphism of multiplicative $\Omega$-group*

(2.5.2) $$f : A_1 \to \text{End}(\Omega_2, A_2)$$

**left-side representation** *of multiplicative $\Omega$-group $A_1$ or* **left-side $A_1$-representation** *in $\Omega_2$-algebra $A_2$ if the map $f$ holds*

(2.5.3) $$f(a_1 * b_1) \bullet a_2 = (f(a_1) \bullet f(b_1)) \bullet a_2$$

*We identify an $A_1$-number $a_1$ and its image $f(a_1)$ and write left-side transformation caused by $A_1$-number $a_1$ as*

$$a_2' = f(a_1) \bullet a_2 = a_1 * a_2$$

*In this case, the equality* (2.5.3) *gets following form*

(2.5.4) $$f(a_1 * b_1) \bullet a_2 = (a_1 * b_1) * a_2$$

*The map*

$$(a_1, a_2) \in A_1 \times A_2 \to a_1 * a_2 \in A_2$$

*generated by left-side representation $f$ is called* **left-side product** *of $A_2$-number $a_2$ over $A_1$-number $a_1$.* $\square$

From the equality (2.5.4), it follows that

(2.5.5) $$(a_1 * b_1) * a_2 = a_1 * (b_1 * a_2)$$

We consider the equality (2.5.5) as **associative law**.

EXAMPLE 2.5.2. *Let*

$$f : A_2 \to A_2$$
$$g : A_2 \to A_2$$

*be endomorphisms of $\Omega_2$-algebra $A_2$. Let product in multiplicative $\Omega$-group $\text{End}(\Omega_2, A_2)$ is composition of endomorphisms. Since the product of maps $f$ and $g$ is defined in the same order as these maps act on $A_2$-number, then we consider the equality*

(2.5.6) $$(f \circ g) \circ a = f \circ (g \circ a)$$

*as associative law. This allows writing of equality* (2.5.6) *without using of brackets*

$$f \circ g \circ a = f \circ (g \circ a) = (f \circ g) \circ a$$

*as well it allows writing of equality* (2.5.3) *in the following form*

(2.5.7) $$f(a_1 * b_1) \circ a_2 = f(a_1) \circ f(b_1) \circ a_2$$

$\square$

---

[2.18] Very often a product in multiplicative $\Omega$-group $\text{End}(\Omega_2, A_2)$ is superposition of endomorphisms

$$f \bullet g = f \circ g$$

However, a product in multiplicative $\Omega$-group $\text{End}(\Omega_2, A_2)$ may be different from superposition of endomorphisms.

REMARK 2.5.3. *Let the map*

$$f : A_1 \longrightarrow\!\!\!*\!\!\!\longrightarrow A_2$$

*be the left-side representation of multiplicative $\Omega$-group $A_1$ in $\Omega_2$-algebra $A_2$. Let the map*

$$g : B_1 \longrightarrow\!\!\!*\!\!\!\longrightarrow B_2$$

*be the left-side representation of multiplicative $\Omega$-group $B_1$ in $\Omega_2$-algebra $B_2$. Let the map*

$$(r_1 : A_1 \to B_1, \ r_2 : A_2 \to B_2)$$

*be morphism of representations. We use notation*

$$r_2(a_2) = r_2 \circ a_2$$

*for image of $A_2$-number $a_2$ with respect to the map $r_2$. Then we can write the equality (2.3.2) in the following form*

$$r_2 \circ (a_1 * a_2) = r_1(a_1) * (r_2 \circ a_2)$$

$\square$

DEFINITION 2.5.4. *Let* $\operatorname{End}(\Omega_2, A_2)$ *be a multiplicative $\Omega$-group with product*[2.19]

$$(f, g) \to f \bullet g$$

*Let an endomorphism $f$ act on $A_2$-number $a$ on the right. We will use notation*

$$(2.5.8) \qquad\qquad f(a_2) = a_2 \bullet f$$

*Let $A_1$ be multiplicative $\Omega$-group with product*

$$(a, b) \to a * b$$

*We call a homomorphism of multiplicative $\Omega$-group*

$$(2.5.9) \qquad\qquad f : A_1 \to \operatorname{End}(\Omega_2, A_2)$$

**right-side representation** *of multiplicative $\Omega$-group $A_1$ or* **right-side $A_1$-representation** *in $\Omega_2$-algebra $A_2$ if the map $f$ holds*

$$(2.5.10) \qquad\qquad a_2 \bullet f(a_1 * b_1) = a_2 \bullet (f(a_1) \bullet f(b_1))$$

*We identify an $A_1$-number $a_1$ and its image $f(a_1)$ and write right-side transformation caused by $A_1$-number $a_1$ as*

$$a_2' = a_2 \bullet f(a_1) = a_2 * a_1$$

*In this case, the equality (2.5.10) gets following form*

$$(2.5.11) \qquad\qquad a_2 \bullet f(a_1 * b_1) = a_2 * (a_1 * b_1)$$

*The map*

$$(a_1, a_2) \in A_1 \times A_2 \to a_2 * a_1 \in A_2$$

*generated by right-side representation $f$ is called* **right-side product** *of $A_2$-number $a_2$ over $A_1$-number $a_1$.* $\square$

From the equality (2.5.11), it follows that

(2.5.12) $$a_2 * (a_1 * b_1) = (a_2 * a_1) * b_1$$

We consider the equality (2.5.12) as **associative law**.

EXAMPLE 2.5.5. *Let*

$$f : A_2 \to A_2$$
$$g : A_2 \to A_2$$

*be endomorphisms of $\Omega_2$-algebra $A_2$. Let product in multiplicative $\Omega$-group* $\mathrm{End}(\Omega_2, A_2)$ *is composition of endomorphisms. Since the product of maps $f$ and $g$ is defined in the same order as these maps act on $A_2$-number, then we consider the equality*

(2.5.13) $$a \circ (g \circ f) = (a \circ g) \circ f$$

*as associative law. This allows writing of equality (2.5.13) without using of brackets*

$$a \circ g \circ f = (a \circ g) \circ f = a \circ (g \circ f)$$

*as well it allows writing of equality (2.5.10) in the following form*

(2.5.14) $$a_2 \circ f(a_1 * b_1) = a_2 \circ f(a_1) \circ f(b_1)$$

$\square$

REMARK 2.5.6. *Let the map*

$$f : A_1 \overset{*}{\longrightarrow} A_2$$

*be the right-side representation of multiplicative $\Omega$-group $A_1$ in $\Omega_2$-algebra $A_2$. Let the map*

$$g : B_1 \overset{*}{\longrightarrow} B_2$$

*be the right-side representation of multiplicative $\Omega$-group $B_1$ in $\Omega_2$-algebra $B_2$. Let the map*

$$(r_1 : A_1 \to B_1, \ r_2 : A_2 \to B_2)$$

*be morphism of representations. We use notation*

$$r_2(a_2) = r_2 \circ a_2$$

*for image of $A_2$-number $a_2$ with respect to the map $r_2$. Then we can write the equality (2.3.2) in the following form*

$$r_2 \circ (a_2 * a_1) = (r_2 \circ a_2) * r_1(a_1)$$

$\square$

If multiplicative $\Omega$-group $A_1$ is Abelian, then there is no difference between left-side and right-side representations.

---

[2.19]  Very often a product in multiplicative $\Omega$-group $\mathrm{End}(\Omega_2, A_2)$ is superposition of endomorphisms

$$f \bullet g = f \circ g$$

However, a product in multiplicative $\Omega$-group $\mathrm{End}(\Omega_2, A_2)$ may be different from superposition of endomorphisms.

DEFINITION 2.5.7. *Let $A_1$ be Abelian multiplicative $\Omega$-group. We call a homomorphism of multiplicative $\Omega$-group*

$$(2.5.15) \qquad\qquad f : A_1 \to \operatorname{End}(\Omega_2, A_2)$$

**representation** *of multiplicative $\Omega$-group $A_1$ or $A_1$-**representation** in $\Omega_2$-algebra $A_2$ if the map $f$ holds*

$$(2.5.16) \qquad\qquad f(a_1 * b_1) \bullet a_2 = (f(a_1) \bullet f(b_1)) \bullet a_2$$

$\square$

Usually we identify a representation of the Abelian multiplicative $\Omega$-group $A_1$ and a left-side representation of the multiplicative $\Omega$-group $A_1$. However, if it is necessary for us, we identify a representation of the Abelian multiplicative $\Omega$-group $A_1$ and a right-side representation of the multiplicative $\Omega$-group $A_1$.

DEFINITION 2.5.8. *Left-side representation*

$$f : A_1 \longrightarrow_* A_2$$

*is called* **covariant** *if the equality*

$$a_1 * (b_1 * a_2) = (a_1 * b_1) * a_2$$

*is true.*

$\square$

DEFINITION 2.5.9. *Left-side representation*

$$f : A_1 \longrightarrow_* A_2$$

*is called* **contravariant** *if the equality*

$$(2.5.17) \qquad\qquad a_1^{-1} * (b_1^{-1} * a_2) = (b_1 * a_1)^{-1} * a_2$$

*is true.*

$\square$

THEOREM 2.5.10. *The product*

$$(a, b) \to a * b$$

*in multiplicative $\Omega$-group $A$ determines two different representations.*

- *the* **left shift**

$$a' = L(b) \circ a = b * a$$

  *is left-side representation of multiplicative $\Omega$-group $A$ in $\Omega$-algebra $A$*

$$(2.5.18) \qquad\qquad L(c * b) = L(c) \circ L(b)$$

- *the* **right shift**

$$a' = a \circ R(b) = a * b$$

  *is right-side representation of multiplicative $\Omega$-group $A$ in $\Omega$-algebra $A$*

$$(2.5.19) \qquad\qquad R(b * c) = R(b) \circ R(c)$$

PROOF. The theorem follows from the theorem [15]-5.2.1.                    $\square$

DEFINITION 2.5.11. *Let $A_1$ be $\Omega$-groupoid with product*

$$(a, b) \to a * b$$

*Let the map*

$$f : A_1 \longrightarrow_* \longrightarrow A_2$$

*be the left-side representation of $\Omega$-groupoid $A_1$ in $\Omega_2$-algebra $A_2$. For any $a_2 \in A_2$, we define* **orbit of representation** *of the $\Omega$-groupoid $A_1$ as set*

$$A_1 * a_2 = \{b_2 = a_1 * a_2 : a_1 \in A_1\}$$

$\square$

DEFINITION 2.5.12. *Let $A_1$ be $\Omega$-groupoid with product*

$$(a, b) \to a * b$$

*Let the map*

$$f : A_1 \longrightarrow_* \longrightarrow A_2$$

*be the right-side representation of $\Omega$-groupoid $A_1$ in $\Omega_2$-algebra $A_2$. For any $a_2 \in A_2$, we define* **orbit of representation** *of the $\Omega$-groupoid $A_1$ as set*

$$a_2 * A_1 = \{b_2 = a_2 * a_1 : a_1 \in A_1\}$$

$\square$

THEOREM 2.5.13. *Let the map*

$$f : A_1 \longrightarrow_* \longrightarrow A_2$$

*be the left-side representation of multiplicative $\Omega$-group $A_1$. Let*

$$(2.5.20) \qquad b_2 \in A_1 * a_2$$

*Then*

$$(2.5.21) \qquad A_1 * a_2 = A_1 * b_2$$

PROOF. The theorem follows from the theorem [15]-5.3.7.    $\square$

Since left-side representation and right-side representation are based on homomorphism of $\Omega$-group, then the following statement is true.

THEOREM 2.5.14 (duality principle for representation of multiplicative $\Omega$-group). *Any statement which holds for left-side representation of multiplicative $\Omega$-group $A_1$ holds also for right-side representation of multiplicative $\Omega$-group $A_1$, if we will use right-side product over $A_1$-number $a_1$ instead of left-side product over $A_1$-number $a_1$.*    $\square$

THEOREM 2.5.15. *If there exists single transitive representation*

$$f : A_1 \longrightarrow_* \longrightarrow A_2$$

*of multiplicative $\Omega$-group $A_1$ in $\Omega_2$-algebra $A_2$, then we can uniquely define coordinates on $A_2$ using $A_1$-numbers.*

> If $f$ is left-side single transitive representation then $f(a)$ is equivalent to the left shift $L(a)$ on the group $A_1$. If $f$ is right-side single transitive representation then $f(a)$ is equivalent to the right shift $R(a)$ on the group $A_1$.

Proof. The theorem follows from the theorem [15]-5.5.4.     □

Theorem 2.5.16. *Left and right shifts on multiplicative $\Omega$-group $A_1$ are commuting.*

Proof. The theorem follows from the associativity of product on multiplicative $\Omega$-group $A_1$

$$(L(a) \circ c) \circ R(b) = (a * c) * b = a * (c * b) = L(a) \circ (c \circ R(b))$$

□

Theorem 2.5.17. *Let $A_1$ be multiplicative $\Omega$-group. Let left-side $A_1$-representation $f$ on $\Omega_2$-algebra $A_2$ be single transitive. Then we can uniquely define a single transitive right-side $A_1$-representation $h$ on $\Omega_2$-algebra $A_2$ such that diagram*

$$
\begin{array}{ccc}
A_2 & \xrightarrow{\ h(a_1)\ } & A_2 \\
\downarrow{\scriptstyle f(b_1)} & & \downarrow{\scriptstyle f(b_1)} \\
A_2 & \xrightarrow{\ h(a_1)\ } & A_2
\end{array}
$$

*is commutative for any $a_1$, $b_1 \in A_1$.*

Proof. We use group coordinates for $A_2$-numbers $a_2$. Then according to theorem 2.5.15 we can write the left shift $L(a_1)$ instead of the transformation $f(a_1)$.

Let $a_2$, $b_2 \in A_2$. Then we can find one and only one $a_1 \in A_1$ such that

$$b_2 = a_2 * a_1 = a_2 \circ R(a_1)$$

We assume

$$h(a) = R(a)$$

For some $b_1 \in A_1$, we have

$$c_2 = f(b_1) \bullet a_2 = L(b_1) \circ a_2 \quad d_2 = f(b_1) \bullet b_2 = L(b_1) \circ b_2$$

According to the theorem 2.5.16, the diagram

$$(2.5.22)
\begin{array}{ccc}
a_2 & \xrightarrow{\ h(a_1)=R(a_1)\ } & b_2 \\
\downarrow{\scriptstyle f(b_1)=L(b_1)} & & \downarrow{\scriptstyle f(b_1)=L(b_1)} \\
c_2 & \xrightarrow{\ h(a_1)=R(a_1)\ } & d_2
\end{array}
$$

is commutative.

Changing $b_1$ we get that $c_2$ is an arbitrary $A_2$-number.

We see from the diagram that if $a_2 = b_2$ then $c_2 = d_2$ and therefore $h(e) = \delta$. On other hand if $a_2 \neq b_2$ then $c_2 \neq d_2$ because the left-side $A_1$-representation $f$ is single transitive. Therefore the right-side $A_1$-representation $h$ is effective.

In the same way we can show that for given $c_2$ we can find $a_1$ such that $d_2 = c_2 \bullet h(a_1)$. Therefore the right-side $A_1$-representation $h$ is single transitive.

In general the product of transformations of the left-side $A_1$-representation $f$ is not commutative and therefore the right-side $A_1$-representation $h$ is different from the left-side $A_1$-representation $f$. In the same way we can create a left-side $A_1$-representation $f$ using the right-side $A_1$-representation $h$.          □

Representations $f$ and $h$ are called **twin representations** of the multiplicative $\Omega$-group $A_1$.

REMARK 2.5.18. *It is clear that transformations $L(a)$ and $R(a)$ are different until the multiplicative $\Omega$-group $A_1$ is nonabelian. However they both are maps onto. Theorem 2.5.17 states that if both right and left shift presentations exist on the set $A_2$, then we can define two commuting representations on the set $A_2$. The right shift or the left shift only cannot represent both types of representation. To understand why it is so let us change diagram (2.5.22) and assume*

$$h(a_1) \bullet a_2 = L(a_1) \circ a_2 = b_2$$

*instead of*

$$a_2 \bullet h(a_1) = a_2 \circ R(a_1) = b_2$$

*and let us see what expression $h(a_1)$ has at the point $c_2$. The diagram*

$$
\begin{array}{ccc}
a_2 & \xrightarrow{\ h(a_1)=L(a_1)\ } & b_2 \\
{\scriptstyle f(b_1)=L(b_1)}\big\downarrow & & \big\downarrow{\scriptstyle f(b_1)=L(b_1)} \\
c_2 & \xrightarrow[\ h(a_1)\ ]{} & d_2
\end{array}
$$

*is equivalent to the diagram*

$$
\begin{array}{ccc}
a_2 & \xrightarrow{\ h(a_1)=L(a_1)\ } & b_2 \\
{\scriptstyle (f(b_1))^{-1}=L(b_1^{-1})}\big\uparrow & & \big\downarrow{\scriptstyle f(b_1)=L(b_1)} \\
c_2 & \xrightarrow[\ h(a_1)\ ]{} & d_2
\end{array}
$$

*and we have $d_2 = b_1 b_2 = b_1 a_1 a_2 = b_1 a_1 b_1^{-1} c_2$. Therefore*

$$h(a_1) \bullet c_2 = (b_1 a_1 b_1^{-1}) c_2$$

*We see that the representation of $h$ depends on its argument.*          □

## 2.6. Tensor Product of Representations

DEFINITION 2.6.1. *Let $A$, $B_1$, ..., $B_n$, $B$ be universal algebras. Let, for any $k$, $k = 1, ..., n$,*

$$f_k : A \longrightarrow\!\!\!* B_k$$

*be effective representation of $\Omega_1$-algebra $A$ in $\Omega_2$-algebra $B_k$. Let*

$$f : A \longrightarrow\!\!\!* B$$

*be effective representation of $\Omega_1$-algebra $A$ in $\Omega_2$-algebra $B$. The map*

$$r_2 : B_1 \times ... \times B_n \to B$$

*is called* **reduced polymorphism of representations** $f_1$, ..., $f_n$ *into representation* $f$, *if, for any* $k$, $k = 1$, ..., $n$, *provided that all variables except the variable* $x_k \in B_k$ *have given value, the map* $r_2$ *is a reduced morphism of representation* $f_k$ *into representation* $f$.

*If* $f_1 = ... = f_n$, *then we say that the map* $r_2$ *is reduced polymorphism of representation* $f_1$ *into representation* $f$.

*If* $f_1 = ... = f_n = f$, *then we say that the map* $r_2$ *is reduced polymorphism of representation* $f$.                    □

DEFINITION 2.6.2. *Let* $A$ *be Abelian multiplicative* $\Omega_1$-*group. Let* $A_1$, ..., $A_n$ *be* $\Omega_2$-*algebras.* [2.20] *Let, for any* $k$, $k = 1$, ..., $n$,

$$f_k : A \dashrightarrow A_k$$

*be effective representation of multiplicative* $\Omega_1$-*group* $A$ *in* $\Omega_2$-*algebra* $A_k$. *Consider category* $\mathcal{A}$ *whose objects are reduced polymorphisms of representations* $f_1$, ..., $f_n$

$$r_1 : A_1 \times ... \times A_n \longrightarrow S_1 \qquad r_2 : A_1 \times ... \times A_n \longrightarrow S_2$$

*where* $S_1$, $S_2$ *are* $\Omega_2$-*algebras and*

$$g_1 : A \dashrightarrow S_1 \qquad g_2 : A \dashrightarrow S_2$$

*are effective representations of multiplicative* $\Omega_1$-*group* $A$. *We define morphism* $r_1 \to r_2$ *to be reduced morphism of representations* $h : S_1 \to S_2$ *making following diagram commutative*

*Universal object* $B_1 \otimes ... \otimes B_n$ *of category* $\mathcal{A}$ *is called* **tensor product** *of representations* $A_1$, ..., $A_n$.                    □

THEOREM 2.6.3. *Since there exists polymorphism of representations, then there exists tensor product of representations.*

PROOF. The theorem follows from the theorem [15]-4.7.5.                    □

THEOREM 2.6.4. *Let*

$$b_k \in B_k \quad k = 1, ..., n \quad b_{i\cdot 1}, ..., b_{i\cdot p} \in B_i \quad \omega \in \Omega_2(p) \quad a \in A$$

*Tensor product is distributive over operation* $\omega$

(2.6.1)
$$b_1 \otimes ... \otimes (b_{i\cdot 1}...b_{i\cdot p}\omega) \otimes ... \otimes b_n$$
$$= (b_1 \otimes ... \otimes b_{i\cdot 1} \otimes ... \otimes b_n)...(b_1 \otimes ... \otimes b_{i\cdot p} \otimes ... \otimes b_n)\omega$$

---

*The representation of multiplicative $\Omega_1$-group $A$ in tensor product is defined by equality*

$$(2.6.2) \qquad b_1 \otimes ... \otimes (f_i(a) \circ b_i) \otimes ... \otimes b_n = f(a) \circ (b_1 \otimes ... \otimes b_i \otimes ... \otimes b_n)$$

PROOF. The theorem follows from the theorem [15]-4.7.11. □

THEOREM 2.6.5. *Let $B_1$, ..., $B_n$ be $\Omega_2$-algebras. Let*

$$f : A_1 \times ... \times A_n \to B_1 \otimes ... \otimes B_n$$

*be reduced polymorphism defined by equality*

$$(2.6.3) \qquad f \circ (b_1, ..., b_n) = b_1 \otimes ... \otimes b_n$$

*Let*

$$g : A_1 \times ... \times A_n \to V$$

*be reduced polymorphism into $\Omega$-algebra $V$. There exists morphism of representations*

$$h : B_1 \otimes ... \otimes B_n \to V$$

*such that the diagram*

*is commutative.*

PROOF. The theorem follows from the theorem [15]-4.7.10. □

THEOREM 2.6.6. *The map*

$$(x_1, ..., x_n) \in B_1 \times ... \times B_n \to x_1 \otimes ... \otimes x_n \in B_1 \otimes ... \otimes B_n$$

*is polymorphism.*

PROOF. The theorem follows from the theorem [15]-4.7.9. □

## 2.7. Basis of Representation of Universal Algebra

DEFINITION 2.7.1. *Let*

$$f : A_1 \longrightarrow\!\!\!* \longrightarrow A_2$$

*be representation of $\Omega_1$-algebra $A_1$ in $\Omega_2$-algebra $A_2$. The set $B_2 \subset A_2$ is called* **stable set of representation** *$f$, if $f(a)(m) \in B_2$ for each $a \in A_1$, $m \in B_2$.* □

THEOREM 2.7.2. *Let*

$$f : A_1 \longrightarrow\!\!\!* \longrightarrow A_2$$

*be representation of $\Omega_1$-algebra $A_1$ in $\Omega_2$-algebra $A_2$. Let set $B_2 \subset A_2$ be subalgebra of $\Omega_2$-algebra $A_2$ and stable set of representation $f$. Then there exists representation*

$$f_{B_2} : A_1 \relbar* \longrightarrow B_2$$

*such that $f_{B_2}(a) = f(a)|_{B_2}$. Representation $f_{B_2}$ is called* **subrepresentation** *of representation $f$.*

PROOF. The theorem follows from the theorem [15]-6.1.2.  $\square$

THEOREM 2.7.3. *The set[2.21] $\mathcal{B}_f$ of all subrepresentations of representation $f$ generates a closure system on $\Omega_2$-algebra $A_2$ and therefore is a complete lattice.*

PROOF. The theorem follows from the theorem [15]-6.1.3.  $\square$

We denote the corresponding closure operator by $J[f]$. Thus $J[f, X]$ is the intersection of all subalgebras of $\Omega_2$-algebra $A_2$ containing $X$ and stable with respect to representation $f$.

DEFINITION 2.7.4. $J[f, X]$ *is called* **subrepresentation** *generated by set $X$, and $X$ is a generating set of subrepresentation $J[f, X]$. In particular, a* **generating set** *of representation $f$ is a subset $X \subset A_2$ such that $J[f, X] = A_2$.*  $\square$

THEOREM 2.7.5. *Let[2.22]*

$$g : A_1 \longrightarrow* \longrightarrow A_2$$

*be representation of $\Omega_1$-algebra $A_1$ in $\Omega_2$-algebra $A_2$. Let $X \subset A_2$. Define a subset $X_k \subset A_2$ by induction on $k$.*

2.7.5.1: $X_0 = X$

2.7.5.2: $x \in X_k => x \in X_{k+1}$

2.7.5.3: $x_1 \in X_k, \ ..., \ x_n \in X_k, \ \omega \in \Omega_2(n) => x_1...x_n\omega \in X_{k+1}$

2.7.5.4: $x \in X_k, \ a \in A => f(a)(x) \in X_{k+1}$

*Then*

$$(2.7.1) \qquad \bigcup_{k=0}^{\infty} X_k = J[f, X]$$

PROOF. The theorem follows from the theorem [15]-6.1.4.  $\square$

DEFINITION 2.7.6. *Let $X \subset A_2$. For each $m \in J[f, X]$ there exists $\Omega_2$-**word** defined according to following rules.* $w[f, X, m]$

2.7.6.1: *If $m \in X$, then $m$ is $\Omega_2$-word.*

2.7.6.2: *If $m_1, \ ..., \ m_n$ are $\Omega_2$-words and $\omega \in \Omega_2(n)$, then $m_1...m_n\omega$ is $\Omega_2$-word.*

2.7.6.3: *If $m$ is $\Omega_2$-word and $a \in A_1$, then $f(a)(m)$ is $\Omega_2$-word.*

---

[2.21]  This definition is similar to definition of the lattice of subalgebras ([18], p. 79, 80).

[2.22]  The statement of theorem is similar to the statement of theorem 5.1, [18], page 79.

We will identify an element $m \in J[f, X]$ and corresponding it $\Omega_2$-word using equation

$$m = w[f, X, m]$$

Similarly, for an arbitrary set $B \subset J[f, X]$ we consider the set of $\Omega_2$-words[2.23]

$$w[f, X, B] = \{w[f, X, m] : m \in B\}$$

We also use notation

$$w[f, X, B] = (w[f, X, m], m \in B)$$

Denote $w[f, X]$ **the set of $\Omega_2$-words of representation** $J[f, X]$.  □

---

THEOREM 2.7.7. *Let*

$$f : A_1 \longrightarrow_* A_2$$

*be representation of $\Omega_1$-algebra $A_1$ in $\Omega_2$-algebra $A_2$. Let*

$$g : A_1 \longrightarrow_* B_2$$

*be representation of $\Omega_1$-algebra $A_1$ in $\Omega_2$-algebra $B_2$. Let $X$ be the generating set of representation $f$. Let*

$$R : A_2 \to B_2$$

*be reduced morphism of representation[2.24] and $X' = R(X)$. Reduced morphism $R$ of representation generates the map of $\Omega_2$-words*

$$w[f \to g, X, R] : w[f, X] \to w[g, X']$$

*such that*

2.7.7.1: *If* $m \in X$, $m' = R(m)$, *then*

$$w[f \to g, X, R](m) = m'$$

2.7.7.2: *If*

$$m_1, ..., m_n \in w[f, X]$$

$$m'_1 = w[f \to g, X, R](m_1) \quad ... \quad m'_n = w[f \to g, X, R](m_n)$$

*then for operation $\omega \in \Omega_2(n)$ holds*

$$w[f \to g, X, R](m_1...m_n\omega) = m'_1...m'_n\omega$$

2.7.7.3: *If*

$$m \in w[f, X] \quad m' = w[f \to g, X, R](m) \quad a \in A_1$$

*then*

$$w[f \to g, X, R](f(a)(m)) = g(a)(m')$$

PROOF. The theorem follows from the theorem [15]-6.1.10.  □

---

DEFINITION 2.7.8. *Let*

$$f : A_1 \longrightarrow_* A_2$$

---

[2.23] The expression $w[f, X, m]$ is a special case of the expression $w[f, X, B]$, namely

$$w[f, X, \{m\}] = \{w[f, X, m]\}$$

[2.24] I considered morphism of representation in the theorem [15]-8.1.7.

*be representation of $\Omega_1$-algebra $A_1$ in $\Omega_2$-algebra $A_2$ and*

$$Gen[f] = \{X \subseteq A_2 : J[f, X] = A_2\}$$

*If, for the set $X \subset A_2$, it is true that $X \in Gen[f]$, then for any set $Y$, $X \subset Y \subset A_2$, also it is true that $Y \in Gen[f]$. If there exists minimal set $X \in Gen[f]$, then the set $X$ is called **quasibasis** of representation $f$.* □

THEOREM 2.7.9. *If the set $X$ is the quasibasis of the representation $f$, then, for any $m \in X$, the set $X \setminus \{m\}$ is not generating set of the representation $f$.*

PROOF. The theorem follows from the theorem [15]-6.2.2. □

REMARK 2.7.10. *The proof of the theorem 2.7.9 gives us effective method for constructing the quasibasis of the representation $f$. Choosing an arbitrary generating set, step by step, we remove from set those elements which have coordinates relative to other elements of the set. If the generating set of the representation is infinite, then this construction may not have the last step. If the representation has finite generating set, then we need a finite number of steps to construct a quasibasis of this representation.* □

We introduced $\Omega_2$-word of $x \in A_2$ relative generating set $X$ in the definition 2.7.6. From the theorem 2.7.9, it follows that if the generating set $X$ is not an quasibasis, then a choice of $\Omega_2$-word relative generating set $X$ is ambiguous. However, even if the generating set $X$ is an quasibasis, then a representation of $m \in A_2$ in form of $\Omega_2$-word is ambiguous.

REMARK 2.7.11. *There are three reasons of ambiguity in notation of $\Omega_2$-word.*

2.7.11.1: *In $\Omega_i$-algebra $A_i$, $i = 1, 2$, equalities may be defined. For instance, if $e$ is unit of multiplicative group $A_i$, then the equality*

$$ae = a$$

*is true for any $a \in A_i$.*

2.7.11.2: *Ambiguity of choice of $\Omega_2$-word may be associated with properties of representation. For instance, if $m_1$, ..., $m_n$ are $\Omega_2$-words, $\omega \in \Omega_2(n)$ and $a \in A_1$, then* [2.25]

$$(2.7.2) \qquad f(a)(m_1...m_n\omega) = (f(a)(m_1))...(f(a)(m_n))\omega$$

*At the same time, if $\omega$ is operation of $\Omega_1$-algebra $A_1$ and operation of $\Omega_2$-algebra $A_2$, then we require that $\Omega_2$-words $f(a_1...a_n\omega)(x)$ and $(f(a_1)(x))...(f(a_n)(x))\omega$ describe the same element of $\Omega_2$-algebra $A_2$.* [2.26]

$$(2.7.3) \qquad f(a_1...a_n\omega)(x) = (f(a_1)(x))...(f(a_n)(x))\omega$$

---

[2.25]For instance, let $\{e_1, e_2\}$ be the basis of vector space over field $k$. The equation (2.7.2) has the form of distributive law

$$a(b^1 e_1 + b^2 e_2) = (ab^1)e_1 + (ab^2)e_2$$

[2.26]For vector space, this requirement has the form of distributive law

$$(a + b)e_1 = ae_1 + be_1$$

2.7.11.3: *Equalities like* (2.7.2), (2.7.3) *persist under morphism of representation. Therefore we can ignore this form of ambiguity of $\Omega_2$-word. However, a fundamentally different form of ambiguity is possible. We can see an example of such ambiguity in theorems [15]-9.9.23, [15]-9.9.24.*

*So we see that we can define different equivalence relations on the set of $\Omega_2$-words.* [2.27] *Our goal is to find a maximum equivalence on the set of $\Omega_2$-words which persist under morphism of representation.* □

DEFINITION 2.7.12. *Generating set $X$ of representation $f$ generates equivalence*

$$\rho[f, X] = \{(w[f, X, m], w_1[f, X, m]) : m \in A_2\}$$

*on the set of $\Omega_2$-words.* □

THEOREM 2.7.13. *Let $X$ be quasibasis of the representation*

$$f : A_1 \overset{*}{\longrightarrow} A_2$$

*For an generating set $X$ of representation $f$, there exists equivalence*

$$\lambda[f, X] \subseteq w[f, X] \times w[f, X]$$

*which is generated exclusively by the following statements.*

2.7.13.1: *If in $\Omega_2$-algebra $A_2$ there is an equality*

$$w_1[f, X, m] = w_2[f, X, m]$$

*defining structure of $\Omega_2$-algebra, then*

$$(w_1[f, X, m], w_2[f, X, m]) \in \lambda[f, X]$$

2.7.13.2: *If in $\Omega_1$-algebra $A_1$ there is an equality*

$$w_1[f, X, m] = w_2[f, X, m]$$

*defining structure of $\Omega_1$-algebra, then*

$$(f(w_1)(w[f, X, m]), f(w_2)(w[f, X, m])) \in \lambda[f, X]$$

2.7.13.3: *For any operation $\omega \in \Omega_1(n)$,*

$$(f(a_{11}...a_{1n}\omega)(a_2), (f(a_{11})...f(a_{1n})\omega)(a_2)) \in \lambda[f, X]$$

2.7.13.4: *For any operation $\omega \in \Omega_2(n)$,*

$$(f(a_1)(a_{21}...a_{2n}\omega), f(a_1)(a_{21})...f(a_1)(a_{2n})\omega) \in \lambda[f, X]$$

2.7.13.5: *Let $\omega \in \Omega_1(n) \cap \Omega_2(n)$. If the representation $f$ satisfies equality* [2.28]

$$f(a_{11}...a_{1n}\omega)(a_2) = (f(a_{11})(a_2))...(f(a_{1n})(a_2))\omega$$

*then we can assume that the following equality is true*

$$(f(a_{11}...a_{1n}\omega)(a_2), (f(a_{11})(a_2))...(f(a_{1n})(a_2))\omega) \in \lambda[f, X]$$

---

[2.27] Evidently each of the equalities (2.7.2), (2.7.3) generates some equivalence relation.

[2.28] Consider a representation of commutative ring $D$ in $D$-algebra $A$. We will use notation

$$f(a)(v) = av$$

The operations of addition and multiplication are defined in both algebras. However the equality

$$f(a + b)(v) = f(a)(v) + f(b)(v)$$

PROOF. The theorem follows from the theorem [15]-6.2.5.          □

DEFINITION 2.7.14. *Quasibasis* $\overline{\overline{e}}$ *of the representation* $f$ *such that*

$$\rho[f, \overline{\overline{e}}] = \lambda[f, \overline{\overline{e}}]$$

*is called* **basis of representation** $f$.          □

THEOREM 2.7.15. *Let* $\overline{\overline{e}}$ *be the basis of the representation* $f$. *Let*

$$R_1 : \overline{\overline{e}} \to \overline{\overline{e}}'$$

*be arbitrary map of the set* $\overline{\overline{e}}$. *Consider the map of* $\Omega_2$-*words*

$$w[f \to g, \overline{\overline{e}}, \overline{\overline{e}}', R_1] : w[f, \overline{\overline{e}}] \to w[g, \overline{\overline{e}}']$$

*that satisfies conditions* 2.7.7.1, 2.7.7.2, 2.7.7.3 *and such that*

$$e \in \overline{\overline{e}} => w[f \to g, \overline{\overline{e}}, \overline{\overline{e}}', R_1](e) = R_1(e)$$

*There exists unique endomorphism of representation* $f$ [2.29]

$$r_2 : A_2 \to A_2$$

*defined by rule*

$$R(m) = w[f \to g, \overline{\overline{e}}, \overline{\overline{e}}', R_1](w[f, \overline{\overline{e}}, m])$$

PROOF. The theorem follows from the theorem [15]-6.2.10.          □

## 2.8. Diagram of Representations of Universal Algebras

DEFINITION 2.8.1. **Diagram** $(f, A)$ **of representations of universal algebras** *is oriented graph such that*

2.8.1.1: *the vertex of* $A_k$, $k = 1, ..., n$, *is* $\Omega_k$-*algebra;*

2.8.1.2: *the edge* $f_{kl}$ *is representation of* $\Omega_k$-*algebra* $A_k$ *in* $\Omega_l$-*algebra* $A_l$;

*We require that this graph is connected graph and does not have loops. Let* $A_{[0]}$ *be set of initial vertices of the graph. Let* $A_{[k]}$ *be set of vertices of the graph for which the maximum path from the initial vertices is* $k$.          □

REMARK 2.8.2. *Since different vertices of the graph can be the same algebra, then we denote* $A = (A_{(1)} \ ... \ A_{(n)})$ *the set of universal algebras which are distinct. From the equality*

$$A = (A_{(1)} \ ... \ A_{(n)}) = (A_1 \ ... \ A_n)$$

*it follows that, for any index* $(i)$, *there exists at least one index* $i$ *such that* $A_{(i)} = A_i$. *If there are two sets of sets* $A = (A_{(1)} \ ... \ A_{(n)})$, $B = (B_{(1)} \ ... \ B_{(n)})$ *a nd there is a map*

$$h_{(i)} : A_{(i)} \to B_{(i)}$$

*for an index* $(i)$, *then also there is a map*

$$h_i : A_i \to B_i$$

_______________

is true, and the equality

$$f(ab)(v) = f(a)(v)f(b)(v)$$

is wrong.

[2.29] This statement is similar to the theorem [1]-4.1, p. 135.

*for any index $i$ such that  $A_{(i)} = A_i$  and in this case $h_i = h_{(i)}$.*                    □

DEFINITION 2.8.3. *Diagram $(f, A)$ of representations of universal algebras is called* **commutative** *when diagram meets the following requirement. for each pair of representations*

$$f_{ik} : A_i \longrightarrow_* A_k$$
$$f_{jk} : A_j \longrightarrow_* A_k$$

*the following equality is true*[2.30]

$$(2.8.1) \qquad f_{ik}(a_i)(f_{jk}(a_j)(a_k)) = f_{jk}(a_j)(f_{ik}(a_i)(a_k))$$

□

Consider the theorem 2.8.4 for the purposes of illustration of the definition 2.8.1.

THEOREM 2.8.4. *Let*

$$f_{ij} : A_i \longrightarrow_* A_j$$

*be representation of $\Omega_i$-algebra $A_i$ in $\Omega_j$-algebra $A_j$. Let*

$$f_{jk} : A_j \longrightarrow_* A_k$$

*be representation of $\Omega_j$-algebra $A_j$ in $\Omega_k$-algebra $A_k$. We represent the fragment*[2.31]

$$A_i \xrightarrow{f_{ij}} A_j \xrightarrow{f_{jk}} A_k$$

*of the diagram of representations using the diagram*

$$(2.8.2)$$

*The map*

$$f_{ijk} : A_i \to \mathrm{End}(\Omega_j, \mathrm{End}(\Omega_k, A_k))$$

*is defined by the equality*

$$(2.8.3) \qquad f_{ijk}(a_i)(f_{jk}(a_j)) = f_{jk}(f_{ij}(a_i)(a_j))$$

----

[2.30]  Metaphorically speaking, representations $f_{ik}$ and $f_{jk}$ are transparent to each other.

*where $a_i \in A_i$, $a_j \in A_j$. If the representation $f_{jk}$ is effective and the representation $f_{ij}$ is free, then the map $f_{ijk}$ is free representation*

$$f_{ijk} : A_i \longrightarrow \mathrm{End}(\Omega_k, A_k)$$

*of $\Omega_i$-algebra $A_i$ in $\Omega_j$-algebra $\mathrm{End}(\Omega_k, A_k)$.*

PROOF. The theorem follows from the theorem [15]-7.1.6. $\qquad\square$

DEFINITION 2.8.5. *Let $(f, A)$ be the diagram of representations where $A = (A_{(1)} \ \dots \ A_{(n)})$ is the set of universal algebras. Let $(B, g)$ be the diagram of representations where $B = (B_{(1)} \ \dots \ B_{(n)})$ is the set of universal algebras. The set of maps $h = (h_{(1)} \ \dots \ h_{(n)})$*

$$h_{(i)} : A_{(i)} \to B_{(i)}$$

*is called* **morphism from diagram of representations** $(f, A)$ **into diagram of representations** $(B, g)$, *if for any indexes $(i)$, $(j)$, $i$, $j$ such that $A_{(i)} = A_i$, $A_{(j)} = A_j$ and for any representation*

$$f_{ji} : A_j \longrightarrow A_i$$

*the tuple of maps $(h_j \ h_i)$ is morphism of representations from $f_{ji}$ into $g_{ji}$.* $\qquad\square$

For any representation $f_{ij}$, $i = 1, \dots, n$, $j = 1, \dots, n$, we have diagram

(2.8.4)

Equalities

(2.8.5) $$h_j \circ f_{ij}(a_i) = g_{ij}(h_i(a_i)) \circ h_j$$

(2.8.6) $$h_j(f_{ij}(a_i)(a_j)) = g_{ij}(h_i(a_i))(h_j(a_j))$$

express commutativity of diagram (1).

---

[2.31] The theorem 2.8.4 states that transformations in diagram of representations are coordinated.

# Abelian Group

## 3.1. Group

DEFINITION 3.1.1. *Let binary operation* [3.1]

$$(a, b) \in A \times A \to ab \in A$$

*which we call product be defined on the set A. The set A is called multiplicative monoid if*

3.1.1.1: *the product is* **associative**

$$(ab)c = a(bc)$$

3.1.1.2: *the product has unit element e*

$$ea = ae = a$$

*If product is* **commutative**

$$(3.1.1) \qquad ab = ba$$

*then multiplicative monoid A is called commutative or Abelian.* □

DEFINITION 3.1.2. *Let binary operation* [3.1]

$$(a, b) \in A \times A \to a + b \in A$$

*which we call sum be defined on the set A. The set A is called additive monoid if*

3.1.2.1: *the sum is* **associative**

$$(a + b) + c = a + (b + c)$$

3.1.2.2: *the sum has zero element 0*

$$0 + a = a + 0 = a$$

*If sum is* **commutative**

$$(3.1.2) \qquad a + b = b + a$$

*then additive monoid A is called commutative or Abelian.* □

THEOREM 3.1.3. *Unit element of multiplicative monoid A is unique.*

PROOF. Let $e_1$, $e_2$ be unit elements of multiplicative monoid $A$. The equality

$$(3.1.3) \qquad e_1 = e_1 e_2 = e_2$$

follows from the statement 3.1.1.2. The theorem follows from the equality (3.1.3). □

THEOREM 3.1.4. *Zero element of additive monoid A is unique.*

PROOF. Let $0_1$, $0_2$ be zero elements of additive monoid $A$. The equality

$$(3.1.4) \qquad 0_1 = 0_1 + 0_2 = 0_2$$

follows from the statement 3.1.2.2. The theorem follows from the equality (3.1.4). □

REMARK 3.1.5. *According to definitions 2.1.6, 3.1.1 and the theorem 3.1.3, multiplicative monoid is universal algebra with one binary operation (multipli-*

REMARK 3.1.6. *According to definitions 2.1.6, 3.1.2 and the theorem 3.1.4, additive monoid is universal algebra with one binary operation (addition) and one 0-ary operation (zero element).* □

---

[3.1] See also the definition of monoid on the page [1]-3.

*cation) and one* 0-*ary operation (unit element).*                                               □

---

DEFINITION 3.1.7. *The multiplicative monoid* [3.2] *$A$ is called multiplicative group if for any $A$-number $x$ there exists $A$-number $y$ such that*

$$xy = yx = e$$

*$A$-number $y$ is called inverse for $A$-number $x$. If product is commutative, then multiplicative group $A$ is called commutative or Abelian.*      □

DEFINITION 3.1.8. *The additive monoid* [3.2] *$A$ is called additive group if for any $A$-number $x$ there exists $A$-number $y$ such that*

$$x + y = y + x = 0$$

*$A$-number $y$ is called inverse for $A$-number $x$. If sum is commutative, then additive group $A$ is called commutative or Abelian.*      □

---

Definitions of multiplicative and additive monoids differ only in form of notation of operation. Therefore, if choice of operation is clear from the context, the corresponding algebra will be called monoid. Siimilar remark is true for groups. Definitions of multiplicative and additive groups differ only in form of notation of operation. Therefore, if choice of operation is clear from the context, the corresponding algebra will be called group. Unless otherwise stated, we usually assume that additive group is Abelian.

---

THEOREM 3.1.9. *Let $A$ be multiplicative group. $A$-number $y$, which is inverse for $A$-number $a$, is unique. We use notation $y = a^{-1}$*

$$(3.1.5) \qquad aa^{-1} = e$$

PROOF. Let $y_1$, $y_2$ be inverse numbers for $A$-number $a$. The equality

$$(3.1.6) \qquad \begin{aligned} y_1 &= y_1 e = y_1(ay_2) \\ &= (y_1a)y_2 = ey_2 = y_2 \end{aligned}$$

follows from the statements 3.1.1.1, 3.1.1.2 and from the definition 3.1.7. The theorem follows from the equality (3.1.6).      □

THEOREM 3.1.10. *Let $A$ be additive group. $A$-number $y$, which is inverse for $A$-number $a$, is unique. We use notation $y = -a$*

$$(3.1.7) \qquad a + (-a) = 0$$

PROOF. Let $y_1$, $y_2$ be inverse numbers for $A$-number $a$. The equality

$$(3.1.8) \qquad \begin{aligned} y_1 &= y_1 + 0 = y_1 + (a + y_2) \\ &= (y_1 + a) + y_2 = 0 + y_2 \\ &= y_2 \end{aligned}$$

follows from the statements 3.1.2.1, 3.1.2.2 and from the definition 3.1.8. The theorem follows from the equality (3.1.8).      □

---

THEOREM 3.1.11. *If $A$ is multiplicative group, then*

$$(3.1.9) \qquad (a^{-1})^{-1} = a$$

PROOF. According to the definition (3.1.5), $A$-number $y$, which is inverse for

THEOREM 3.1.12. *If $A$ is additive group, then*

$$(3.1.11) \qquad -(-a) = a$$

PROOF. According to the definition (3.1.7), $A$-number $y$, which is inverse for

---

$A$-number $a^{-1}$, is $(a^{-1})^{-1}$

(3.1.10)          $(a^{-1})^{-1}a^{-1} = e$

According to the theorem 3.1.9, the equality (3.1.9) follows from equalities (3.1.5), (3.1.10).          □

$A$-number $-a$, is $-(-a)$

(3.1.12)          $(-(-a)) + (-a) = 0$

According to the theorem 3.1.10, the equality (3.1.11) follows from equalities (3.1.7), (3.1.12).

THEOREM 3.1.13. *Let $A$ be multiplicative group with unit element $e$. $A$-number, which is inverse for unit element is unit element*

(3.1.13)          $e^{-1} = e$

THEOREM 3.1.14. *Let $A$ be additive group with zero element $0$. $A$-number, which is inverse for zero element is zero element*

(3.1.15)          $-0 = 0$

PROOF. The equality

(3.1.14)          $e^{-1} = e^{-1}e = e$

follows from theorems 3.1.3, 3.1.9. The equality (3.1.13) follows from the equality (3.1.14).          □

PROOF. The equality

(3.1.16)          $-0 = (-0) + 0 = 0$

follows from theorems 3.1.4, 3.1.10. The equality (3.1.15) follows from the equality (3.1.16).          □

DEFINITION 3.1.15. *Let $H$ be subgroup of multiplicative group $G$. The set*

(3.1.17)          $aH = \{ab : b \in H\}$

*is called* **left coset**. *The set*

$$Ha = \{ba : b \in H\}$$

*is called* **right coset**.          □

DEFINITION 3.1.16. *Let $H$ be subgroup of additive group $G$. The set*

(3.1.18)          $a + H = \{a + b : b \in H\}$

*is called* **left coset**. *The set*

$$H + a = \{b + a : b \in H\}$$

*is called* **right coset**.          □

THEOREM 3.1.17. *Let $H$ be subgroup of multiplicative group $G$. If*

(3.1.19)          $aH \cap bH \neq \emptyset$

*then*

$$aH = bH$$

THEOREM 3.1.18. *Let $H$ be subgroup of additive group $G$. If*

(3.1.21)          $a + H \cap b + H \neq \emptyset$

*then*

$$a + H = b + H$$

PROOF. From the statement (3.1.19) and from the definition (3.1.17), it follows that there exist $H$-numbers $c$, $d$ such that

(3.1.20)          $ac = bd$

3.1.17.1: From the theorem 3.1.9, it follows that there exists $H$-number $d^{-1}$. Therefore,

$$cd^{-1} \in H$$

3.1.17.2: If $c \in H$, then from the definition (3.1.17), it follows that $cH = H$.

PROOF. From the statement (3.1.21) and from the definition (3.1.18), it follows that there exist $H$-numbers $c$, $d$ such that

(3.1.22)          $a + c = b + d$

3.1.18.1: From the theorem 3.1.10, it follows that there exists $H$-number $-d$. Therefore,

$$c + (-d) \in H$$

3.1.18.2: If $c \in H$, then from the definition (3.1.18), it follows that $c + H = H$.

3.1.17.3: From statements 3.1.17.1, 3.1.17.2, it follows that

$$(cd^{-1})H = H$$

3.1.17.4: From the statement 3.1.1.1 and from the definition 3.1.7, it follows that

$$(ab)H = a(bH)$$

3.1.17.5: The equality

$$a(cd^{-1}) = (ac)d^{-1} = (bd)d^{-1}$$
$$= b(dd^{-1}) = be = b$$

follows from statements 3.1.1.1, 3.1.1.2 and from equalities (3.1.5), (3.1.20).

The equality

$$bH = a(cd^{-1})H = aH$$

follows from statements 3.1.17.3, 3.1.17.5.                              □

---

3.1.18.3: From statements 3.1.18.1, 3.1.18.2, it follows that

$$(c + (-d)) + H = H$$

3.1.18.4: From the statement 3.1.2.1 and from the definition 3.1.8, it follows that

$$(a + b) + H = a + (b + H)$$

3.1.18.5: The equality

$$a + (c + (-d)) = (a + c) + (-d)$$
$$= (b + d) + (-d)$$
$$= b + (d + (-d))$$
$$= b + 0 = b$$

follows from statements 3.1.2.1, 3.1.2.2 and from equalities (3.1.7), (3.1.22).

The equality

$$b + H = (a + (c + (-d))) + H = a + H$$

follows from statements 3.1.18.3, 3.1.18.5.                              □

---

**Definition 3.1.19.** *Subgroup $H$ of multiplicative group $A$ such that*

$$aH = Ha$$

*for any $A$-number is called* **normal**. [3.3]

□

**Definition 3.1.20.** *Subgroup $H$ of additive group $A$ such that*

$$a + H = H + a$$

*for any $A$-number is called* **normal**. [3.3]

□

---

According to the theorem 3.1.17, normal subgroup $H$ of multiplicative group $A$ generates equivalence relation mod $H$ on the set $A$

$$(3.1.23) \qquad a \equiv b(\mathrm{mod}\, H) \Leftrightarrow aH = bH$$

We use notation $A/H$ for the set of equivalence classes of equivalence relation mod $H$.

According to the theorem 3.1.18, normal subgroup $H$ of multiplicative group $A$ generates equivalence relation mod $H$ on the set $A$

$$(3.1.24)$$
$$a \equiv b(\mathrm{mod}\, H) \Leftrightarrow a + H = b + H$$

We use notation $A/H$ for the set of equivalence classes of equivalence relation mod $H$.

---

**Theorem 3.1.21.** *Let $H$ be normal subgroup of multiplicative group $A$. The product*

$$(3.1.25) \qquad (aH)(bH) = abH$$

*generates group on the set $A/H$.*

**Theorem 3.1.22.** *Let $H$ be normal subgroup of additive group $A$. The sum*

$$(3.1.27) \quad (a + H) + (b + H) = a + b + H$$

*generates group on the set $A/H$.*

---

[3.3]  See also the definition on the page [1]-14.

PROOF. The equality

$$(aH)(bH)$$
$$= a(Hb)H$$
(3.1.26) $$= a(bH)H$$
$$= (ab)(HH)$$
$$= abH$$

follows from the definition 3.1.19 and from statements 3.1.1.1, 3.1.17.4.

The equality (3.1.25) follows from the equality (3.1.26).

Let $c \in aH$, $d \in bH$. Then the equality

$$(ab)H = (aH)(bH)$$
$$= (cH)(dH)$$
$$= (cd)H$$

follows from the definition (3.1.25) and from the theorem 3.1.17. Therefore, the product (3.1.25) does not depend on the choice of $A$-numbers $c \in aH$, $d \in bH$.

Since $eH = H$, then, from the equality (3.1.25), it follows that the product (3.1.25) satisfies the definition of multiplicative group. □

PROOF. The equality

$$(a + H) + (b + H)$$
$$= a + (H + b) + H$$
(3.1.28) $$= a + (b + H) + H$$
$$= (a + b) + (H + H)$$
$$= a + b + H$$

follows from the definition 3.1.20 and from statements 3.1.2.1, 3.1.18.4.

The equality (3.1.27) follows from the equality (3.1.28).

Let $c \in a + H$, $d \in b + H$. Then the equality

$$(a + b) + H = (a + H) + (b + H)$$
$$= (c + H) + (d + H)$$
$$= (c + d) + H$$

follows from the definition (3.1.27) and from the theorem 3.1.18. Therefore, the product (3.1.27) does not depend on the choice of $A$-numbers $c \in a + H$, $d \in b + H$.

Since $e + H = H$, then, from the equality (3.1.27), it follows that the product (3.1.27) satisfies the definition of additive group. □

.

DEFINITION 3.1.23. *Let $H$ be normal subgroup of multiplicative group $A$. Group $A/H$ is called* **factor group** *of group $A$ by $H$ and the map*

$$f : a \in A \to aH \in A/H$$

*is called* **canonical map**. □

DEFINITION 3.1.24. *Let $H$ be normal subgroup of additive group $A$. Group $A/H$ is called* **factor group** *of group $A$ by $H$ and the map*

$$f : a \in A \to a + H \in A/H$$

*is called* **canonical map**. □

## 3.2. Ring

DEFINITION 3.2.1. *Let two binary operations* [3.4] *be defined on the set $A$*

- *sum*

$$(a, b) \in A \times A \to a + b \in A$$

- *product*

$$(a, b) \in A \times A \to ab \in A$$

3.2.1.1: *The set $A$ is Abelian additive group with respect to sum and $0$ is zero element.*

3.2.1.2: *The set $A$ is multiplicative monoid with respect to product.*

3.2.1.3: *Product is* **distributive** *over sum*

$$(3.2.1) \qquad a(b+c) = ab + ac$$

$$(3.2.2) \qquad (a+b)c = ac + bc$$

*If product is commutative, then ring is called commutative.* □

DEFINITION 3.2.2. *If multiplicative monoid of the ring $D$ has unit element* 1, *then the ring $D$ is called* **unital ring** [3.5] □

THEOREM 3.2.3. *Let $A$ be a ring. Then*

$$(3.2.3) \qquad 0a = 0$$

*for any $A$-number $a$.*

PROOF. The equality

$$(3.2.4) \qquad 0a + ba = (0+b)a = ba$$

follows from statements 3.2.1.2, 3.2.1.3. The theorem follows from the equality (3.2.4), from the theorem 3.1.4 and from the statement 3.1.2.2. □

THEOREM 3.2.4. *Let $A$ be a ring. Then*

$$(3.2.5) \qquad (-a)b = -(ab)$$

$$(3.2.6) \qquad a(-b) = -(ab)$$

$$(3.2.7) \qquad (-a)(-b) = ab$$

*for any $A$-numbers $a$, $b$. If ring has unit element, then*

$$(3.2.8) \qquad (-1)a = -a$$

*for any $A$-number $a$.*

PROOF. Equalities

$$(3.2.9) \qquad ab + (-a)b = (a+(-a))b = 0b = 0$$

$$(3.2.10) \qquad ab + a(-b) = a(b+(-b)) = a0 = 0$$

follow from statements 3.2.1.1, 3.2.1.3 and from equalities (3.1.7), (3.2.3). The equality (3.2.5) follows from the theorem 3.1.10 and from the equality (3.2.9). The equality (3.2.6) follows from the theorem 3.1.10 and from the equality (3.2.10).

The equality

$$(3.2.11) \qquad (-a)(-b) = -(a(-b)) = -(-(ab))$$

follows from equalities (3.2.5), (3.2.6). The equality (3.2.7) follows from equalities (3.1.11), (3.2.11).

---

[3.4] See also the definition of ring on the page [1]-83.

[3.5] See similar statement on the page [7]-52.

If ring has unit element, then the equality

(3.2.12) $$b + (-1)b = 1b + (-1)b = (1 + (-1))b = 0b = 0$$

follows from statements 3.2.1.1, 3.2.1.2, 3.2.1.3 and from equalities (3.1.7), (3.2.3).

□

DEFINITION 3.2.5. *A* **left ideal** [3.6] *$D_1$ in a ring $D$ is a subgroup of the additive group of $D$ such that*

$$DD_1 \subseteq D_1$$

□

DEFINITION 3.2.7. *A* **right ideal** [3.6] *$D_1$ in a ring $D$ is a subgroup of the additive group of $D$ such that*

$$D_1D \subseteq D_1$$

□

THEOREM 3.2.6. *If $D_1$ is a left ideal in unital ring $D$ has unit, then*

$$DD_1 = D_1$$

PROOF. The theorem follows from the definition 3.2.5 and equality

$$eD_1 = D_1$$

□

THEOREM 3.2.8. *If $D_1$ is a right ideal in unital ring $D$ has unit, then*

$$D_1D = D_1$$

PROOF. The theorem follows from the definition 3.2.7 and equality

$$D_1e = D_1$$

□

DEFINITION 3.2.9. **Ideal** [3.6] *is a subgroup of the additive group of $D$ which is both left and right ideal.*

□

If $D$ is ring, then the set $\{0\}$ and the set $A$ are ideals which we call trivial ideals.

DEFINITION 3.2.10. *Let $A$ be commutative ring. If the set $A' = A \setminus \{0\}$ is multiplicative group with respect to product, then the set $A$ is called field.*

□

## 3.3. Group-Homomorphism

THEOREM 3.3.1. *The map*

$$f : A \to B$$

*is homomorphism of multiplicative monoid $A$ with unit element $e_1$ into multiplicative monoid $B$ with unit element $e_2$ if*

(3.3.1) $$f(ab) = f(a)f(b)$$

(3.3.2) $$f(e_1) = e_2$$

*for any $A$-numbers $a$, $b$.*

THEOREM 3.3.2. *The map*

$$f : A \to B$$

*is homomorphism of additive monoid $A$ with zero element $0_1$ into additive monoid $B$ with zero element $0_2$ if*

(3.3.3) $$f(a + b) = f(a) + f(b)$$

(3.3.4) $$f(0_1) = 0_2$$

*for any $A$-numbers $a$, $b$.*

---

[3.6] See also the definition on the page [1]-86.

PROOF. The theorem follows from definitions 2.1.6, 2.1.7, 3.1.1 and from the remark 3.1.5.                □

THEOREM 3.3.3. *Let maps*

$$f : A \to B$$
$$g : A \to B$$

*be homomorphisms of multiplicative Abelian group $A$ with unit element $e_1$ into multiplicative Abelian group $B$ with unit element $e_2$. Then the map*

$$h : A \to B$$

*defined by the equality*

$$h(x) = (fg)(x) = f(x)g(x)$$

*is homomorphism of multiplicative Abelian group $A$ into multiplicative Abelian group $B$.*

PROOF. The theorem follows from the equalities

$$h(xy) = f(xy)g(xy)$$
$$= f(x)f(y)g(x)g(y)$$
$$= f(x)g(x)f(y)g(y)$$
$$= h(x)h(y)$$

$$h(e_1) = (fg)(e_1)$$
$$= f(e_1)g(e_1)$$
$$= e_2 e_2 = e_2$$

and the theorem 3.3.1.                □

THEOREM 3.3.5. *Let*

$$f : A \to B$$

*be homomorphism* [3.7] *of multiplicative group $A$ with unit element $e_1$ into multiplicative group $B$ with unit element $e_2$. Then*

(3.3.5)                $f(a^{-1}) = f(a)^{-1}$

PROOF. The equality

(3.3.6)    $e_2 = f(e_1) = f(aa^{-1})$
$$= f(a)f(a^{-1})$$

PROOF. The theorem follows from definitions 2.1.6, 2.1.7, 3.1.2 and from the remark 3.1.6.                □

THEOREM 3.3.4. *Let maps*

$$f : A \to B$$
$$g : A \to B$$

*be homomorphisms of additive Abelian group $A$ with zero element $0_1$ into additive Abelian group $B$ with zero element $0_2$. Then the map*

$$h : A \to B$$

*defined by the equality*

$$h(x) = (f + g)(x) = f(x) + g(x)$$

*is homomorphism of additive Abelian group $A$ into additive Abelian group $B$.*

PROOF. The theorem follows from the equalities

$$h(x + y) = f(x + y) + g(x + y)$$
$$= f(x) + f(y) + g(x) + g(y)$$
$$= f(x) + g(x) + f(y) + g(y)$$
$$= h(x) + h(y)$$

$$h(0_1) = (f + g)(0_1)$$
$$= f(0_1) + g(0_1)$$
$$= 0_2 + 0_2 = 0_2$$

and the theorem 3.3.2.                □

THEOREM 3.3.6. *Let*

$$f : A \to B$$

*be homomorphism* [3.7] *of additive group $A$ with zero element $0_1$ into additive group $B$ with zero element $0_2$. Then*

(3.3.7)                $f(-a) = -f(a)$

PROOF. The equality

(3.3.8)    $0_2 = f(0_1) = f(a + (-a))$
$$= f(a) + f(-a)$$

---

[3.7] The map $f$ also is called group-homomorphism.

follows from equalities (3.1.5), (3.3.1), (3.3.2). The equality (3.3.5) follows from the equality (3.3.6) and from the theorem 3.1.9.    □

follows from equalities (3.1.7), (3.3.3), (3.3.4). The equality (3.3.7) follows from the equality (3.3.8) and from the theorem 3.1.10.    □

---

THEOREM 3.3.7. *Image of group-homomorphism*

$$f : A \to B$$

*is subgroup of group $B$.*

PROOF.

Let $A$ and $B$ be multiplicative groups. According to definitions 2.1.6, 3.1.1 and the theorem 3.3.1, the set $\operatorname{Im} f$ is monoid. According to the theorem 3.3.5, the set $\operatorname{Im} f$ is group.

Let $A$ and $B$ be additive groups. According to definitions 2.1.6, 3.1.2 and the theorem 3.3.2, the set $\operatorname{Im} f$ is monoid. According to the theorem 3.3.6, the set $\operatorname{Im} f$ is group.    □

There exist homomorphism of additive group into multiplicative. For instance, the map $y = e^x$ is homomorphism of additive group of real numbers into multiplicative group of positive real numbers.

---

THEOREM 3.3.8. *Let*

$$f : A \to B$$

*be homomorphism of multiplicative group $A$ with unit element $e_1$ into multiplicative group $B$ with unit element $e_2$. The set*

(3.3.9)    $\ker f = \{a \in A : f(a) = e_2\}$

*is called* **kernel of group-homomorphism** $f$. *Kernel of group-homomorphism is normal subgroup of the group $A$.*

THEOREM 3.3.9. *Let*

$$f : A \to B$$

*be homomorphism of additive group $A$ with zero element $0_1$ into additive group $B$ with zero element $0_2$. The set*

(3.3.12)    $\ker f = \{a \in A : f(a) = 0_2\}$

*is called* **kernel of group-homomorphism** $f$. *Kernel of group-homomorphism is normal subgroup of the group $A$.*

PROOF. Let $A$ be multiplicative group.

3.3.8.1: Let $a$, $b \in \ker f$. The equality

$$f(ab) = f(a)f(b) = e_2 e_2 = e_2$$

follows from equalities (3.3.1), (3.3.9) and from the statement 3.1.1.2. According to the definition (3.3.9), $ab \in \ker f$.

3.3.8.2: According to equalities (3.3.2), (3.3.9), $e_1 \in \ker f$.

3.3.8.3: According to statements 3.3.8.1, 3.3.8.2, the set $\ker f$ is multiplicative monoid.

3.3.8.4: Let $a \in \ker f$. According to the theorem 3.1.9, there exists $A$-

PROOF. Let $A$ be additive group.

3.3.9.1: Let $a$, $b \in \ker f$. The equality

$$f(a + b) = f(a) + f(b) = 0_2 + 0_2 = 0_2$$

follows from equalities (3.3.3), (3.3.12) and from the statement 3.1.2.2. According to the definition (3.3.12), $ab \in \ker f$.

3.3.9.2: According to equalities (3.3.4), (3.3.12), $0_1 \in \ker f$.

3.3.9.3: According to statements 3.3.9.1, 3.3.9.2, the set $\ker f$ is additive monoid.

3.3.9.4: Let $a \in \ker f$. According to the theorem 3.1.10, there exists $A$-

number $a^{-1}$. The equality

$$f(a^{-1}) = f(a)^{-1} = e_2^{-1} = e_2$$

follows from equalities (3.1.13), (3.3.5), (3.3.9). According to the definition (3.3.9), $a^{-1} \in \ker f$.

According to statements 3.3.8.3, 3.3.8.4, the set $\ker f$ is multiplicative group.

The equality

$$(3.3.10) \quad \begin{aligned} & f(a(\ker f)a^{-1}) \\ &= f(a)f(\ker f)f(a^{-1}) \\ &= f(a)e_2(f(a))^{-1} \\ &= f(a)(f(a))^{-1} = e_2 \end{aligned}$$

follows from equalities (3.3.1), (3.3.5), (3.3.9). The equality

$$a(\ker f)a^{-1} = \ker f$$

follows from the equality (3.3.10). Therefore,

$$(3.3.11) \qquad a(\ker f) = (\ker f)a$$

The theorem follows from the equality (3.3.11) and from the definition 3.1.19.

□

number $-a$. The equality

$$f(-a) = -f(a) = -0_2 = 0_2$$

follows from equalities (3.1.15), (3.3.7), (3.3.12). According to the definition (3.3.12), $-a \in \ker f$.

According to statements 3.3.9.3, 3.3.9.4, the set $\ker f$ is additive group.

The equality

$$(3.3.13) \quad \begin{aligned} & f(a + (\ker f) - a) \\ &= f(a) + f(\ker f) - f(a) \\ &= f(a) + 0_2 - f(a) \\ &= f(a) - f(a) = 0_2 \end{aligned}$$

follows from equalities (3.3.3), (3.3.7), (3.3.12). The equality

$$a + (\ker f) + (-a) = \ker f$$

follows from the equality (3.3.13). Therefore,

$$(3.3.14) \qquad a + (\ker f) = (\ker f) + a$$

The theorem follows from the equality (3.3.14) and from the definition 3.1.20.

□

<div style="background-color:#eef">

THEOREM 3.3.10. *Let*

$$f : A \to B$$

*be homomorphism of multiplicative group $A$ with unit element $e_1$ into multiplicative group $B$ with unit element $e_2$. If*

$$(3.3.15) \qquad \ker f = \{e_1\}$$

*then the kernel of the homomorphism $f$ is called* **trivial***, and homomorphism $f$ is monomorphism.*

THEOREM 3.3.11. *Let*

$$f : A \to B$$

*be homomorphism of additive group $A$ with zero element $0_1$ into additive group $B$ with zero element $0_2$. If*

$$(3.3.16) \qquad \ker f = \{0_1\}$$

*then the kernel of the homomorphism $f$ is called* **trivial***, and homomorphism $f$ is monomorphism.*

</div>

PROOF. Let there exist $A$-numbers $a$, $b$ such that

$$(3.3.17) \qquad f(a) = f(b)$$

PROOF. Let there exist $A$-numbers $a$, $b$ such that

$$(3.3.20) \qquad f(a) = f(b)$$

The equality

$$(3.3.18) \quad \begin{aligned} f(ab^{-1}) &= f(a)f(b^{-1}) \\ &= f(a)f(b)^{-1} \\ &= f(b)f(b)^{-1} \\ &= e_2 \end{aligned}$$

follows from equalities (3.3.1), (3.3.5), (3.3.17) and from the statement 3.1.1.2. The equality

$$(3.3.19) \quad ab^{-1} = e_1$$

follows from equalities (3.3.9), (3.3.15), (3.3.18). From equalities (3.1.9), (3.3.19) and from the theorem 3.1.9, it follows that $a = b$. □

The equality

$$(3.3.21) \quad \begin{aligned} f(a+(-b)) &= f(a) + f(-b) \\ &= f(a) + (-f(b)) \\ &= f(b) + (-f(b)) \\ &= 0_2 \end{aligned}$$

follows from equalities (3.3.3), (3.3.7), (3.3.20) and from the statement 3.1.2.2. The equality

$$(3.3.22) \quad a+(-b) = 0_1$$

follows from equalities (3.3.12), (3.3.16), (3.3.21). From equalities (3.1.11), (3.3.22) and from the theorem 3.1.10, it follows that $a = b$. □

**THEOREM 3.3.12.** *Let $B$ be normal subgroup of multiplicative group $A$ with unit element $e$. There exist multiplicative group $C$ and homomorphism*

$$f : A \to C$$

*such that*

$$(3.3.23) \quad \ker f = B$$

PROOF. Let

$$f(x) = xB$$

Since the set $eB = B$ is unit element of multiplicative group $A/B$, then

$$(3.3.24) \quad B \subseteq \ker f$$

If $f(x) = B$, then $xB = B$. According to the theorem 3.1.17, $x \in B$. Therefore

$$(3.3.25) \quad \ker f \subseteq B$$

The theorem follows from statements (3.3.24), (3.3.25). □

**THEOREM 3.3.13.** *Let $B$ be normal subgroup of additive group $A$ with zero element $0$. There exist additive group $C$ and homomorphism*

$$f : A \to C$$

*such that*

$$(3.3.26) \quad \ker f = B$$

PROOF. Let

$$f(x) = x + B$$

Since the set $0 + B = B$ is zero element of additive group $A/B$, then

$$(3.3.27) \quad B \subseteq \ker f$$

If $f(x) = B$, then $x + B = B$. According to the theorem 3.1.18, $x \in B$. Therefore

$$(3.3.28) \quad \ker f \subseteq B$$

The theorem follows from statements (3.3.27), (3.3.28). □

**DEFINITION 3.3.14.** *Let $A_1$, ..., $A_n$ be groups. A sequence of homomorphisms*

$$A_1 \xrightarrow{f_1} A_2 \xrightarrow{f_2} ... \xrightarrow{f_{n-1}} A_n$$

*is called* **exact** *if* $\forall i,\ i = 1, ...., n-2,\ \operatorname{Im} f_i = \ker f_{i+1}$. □

**EXAMPLE 3.3.15.** *Let*

$$(3.3.29) \quad e \longrightarrow A \xrightarrow{f} B \xrightarrow{g} C \xrightarrow{h} e$$

*be exact sequence of homomorphisms of multiplicative groups.*

According to the equality (3.3.2) and the definition *3.3.14*, $\ker f = \{e\}$. According to the theorem *3.3.11*, the map $f$ is monomorphism.

From the diagram (3.3.29), it follows that $\ker h = C$. According to the definition *3.3.14*, $\operatorname{Im} g = C$. According to the definition *2.1.8*, homomorphism $g$ is epimorphism.                                                                 □

EXAMPLE 3.3.16. *According to the example 3.3.15, if the sequence of homomorphisms of additive groups*

$$e \longrightarrow A \xrightarrow{\ f\ } B \longrightarrow e$$

*is exact, then homomorphism $f$ is isomorphism.*                                  □

EXAMPLE 3.3.17. *Joining examples 3.3.15, 3.3.16, we get diagram*

*where rows and column are exact sequence of homomorphisms and homomorphism $f$ is isomorphism.*                                                            □

EXAMPLE 3.3.18. *Let*

$$(3.3.30) \qquad 0 \longrightarrow A \xrightarrow{\ f\ } B \xrightarrow{\ g\ } C \xrightarrow{\ h\ } 0$$

*be exact sequence of homomorphisms of additive groups.*

According to the equality (3.3.4) and the definition *3.3.14*, $\ker f = \{0\}$. According to the theorem *3.3.11*, the map $f$ is monomorphism.

From the diagram (3.3.30), it follows that $\ker h = C$. According to the definition *3.3.14*, $\operatorname{Im} g = C$. According to the definition *2.1.8*, homomorphism $g$ is epimorphism.                                                    □

EXAMPLE 3.3.19. *According to the example 3.3.18, if the sequence of homomorphisms of additive groups*

$$0 \longrightarrow A \xrightarrow{\ f\ } B \longrightarrow 0$$

*is exact, then homomorphism $f$ is isomorphism.*                                  □

EXAMPLE 3.3.20. *Joining examples 3.3.18, 3.3.19, we get diagram*

*where rows and column are exact sequence of homomorphisms and homomorphism f is isomorphism.*  □

## 3.4. Basis of Additive Abelian Group

In this section, $G$ is additive Abelian group.

DEFINITION 3.4.1. *The action of ring of rational integers $Z$ in additive Abelian group $G$ is defined using following rules*

$$(3.4.1) \qquad 0g = 0$$
$$(3.4.2) \qquad (n+1)g = ng + g$$
$$(3.4.3) \qquad (n-1)g = ng - g$$

□

THEOREM 3.4.2. *The action of ring of rational integers $Z$ in additive Abelian group $G$ defined in the definition 3.4.1 is representation. The following equalities are true*

$$(3.4.4) \qquad 1a = a$$
$$(3.4.5) \qquad (nm)a = n(ma)$$
$$(3.4.6) \qquad (m+n)a = ma + na$$
$$(3.4.7) \qquad (m-n)a = ma - na$$
$$(3.4.8) \qquad n(a+b) = na + nb$$

PROOF. The equality (3.4.4) follows from the equality (3.4.1) and from the equality (3.4.2) when $n = 0$.

From the equality (3.4.1), it follows that the equality (3.4.6) is true when $n = 0$.

- Let the equality (3.4.6) is true when $n = k \geq 0$. Then

$$(m+k)a = ma + ka$$

The equality

$$(m + (k+1))a = ((m+k) + 1)a = (m+k)a + a = ma + ka + a$$
$$= ma + (k+1)a$$

follows from the equality (3.4.2). Therefore, the equality (3.4.6) is true when $n = k+1$. According to mathematical induction, the equality (3.4.6) is true for any $n \geq 0$.

- Let the equality (3.4.6) is true when $n = k \leq 0$. Then

$$(m+k)a = ma + ka$$

The equality

$$(m + (k-1))a = ((m+k) - 1)a = (m+k)a - a = ma + ka - a$$
$$= ma + (k-1)a$$

follows from the equality (3.4.3). Therefore, the equality (3.4.6) is true when $n = k-1$. According to mathematical induction, the equality (3.4.6) is true for any $n \leq 0$.

- Therefore, the equality (3.4.6) is true for any $n \in Z$.

The equality

$$(3.4.9) \qquad\qquad (k+n)a - na = ka$$

follows from the equality (3.4.6). The equality (3.4.7) follows from the equality (3.4.9), if we assume $m = k + n$, $k = m - n$.

From the equality (3.4.1), it follows that the equality (3.4.5) is true when $n = 0$.

- Let the equality (3.4.5) is true when $n = k \geq 0$. Then

$$(km)a = k(ma)$$

The equality

$$((k+1)m)a = (km + m)a = (km)a + ma = k(ma) + ma$$
$$= (k+1)(ma)$$

follows from equalities (3.4.2), (3.4.6). Therefore, the equality (3.4.5) is true when $n = k + 1$. According to mathematical induction, the equality (3.4.5) is true for any $n \geq 0$.

- Let the equality (3.4.6) is true when $n = k \leq 0$. Then

$$(km)a = k(ma)$$

The equality

$$((k-1)m)a = (km - m)a = (km)a - ma = k(ma) - ma$$
$$= (k-1)(ma)$$

follows from equalities (3.4.3), (3.4.7). Therefore, the equality (3.4.5) is true when $n = k - 1$. According to mathematical induction, the equality (3.4.5) is true for any $n \leq 0$.

- Therefore, the equality (3.4.5) is true for any $n \in Z$.

From the equality (3.4.1), it follows that the equality (3.4.8) is true when $n = 0$.

- Let the equality (3.4.8) is true when $n = k \geq 0$. Then

$$k(a + b) = ka + kb$$

The equality

$$(k + 1)(a + b) = k(a + b) + a + b = ka + kb + a + b$$
$$= ka + a + kb + b$$
$$= (k + 1)a + (k + 1)b$$

follows from the equality (3.4.2). Therefore, the equality (3.4.8) is true when $n = k+1$. According to mathematical induction, the equality (3.4.8) is true for any $n \geq 0$.

- Let the equality (3.4.6) is true when $n = k \leq 0$. Then

$$k(a + b) = ka + kb$$

The equality

$$(k - 1)(a + b) = k(a + b) - (a + b) = ka + kb - a - b$$
$$= ka - a + kb - b$$
$$= (k - 1)a + (k - 1)b$$

follows from the equality (3.4.3). Therefore, the equality (3.4.8) is true when $n = k-1$. According to mathematical induction, the equality (3.4.8) is true for any $n \leq 0$.

- Therefore, the equality (3.4.8) is true for any $n \in Z$.

From the equality (3.4.8), it follows that the map

$$\varphi(n) : a \in G \to na \in G$$

is an endomorphism of Abelian group $G$. From equalities (3.4.6), (3.4.5), it follows that the map

$$\varphi : Z \to \operatorname{End}(Ab, G)$$

is a homomorphism of the ring $Z$. According to the definition 2.3.1, the map $\varphi$ is representation of ring of rational integers $Z$ in Abelian group $G$. $\square$

THEOREM 3.4.3. *The set of $G$-numbers generated by the set $S = \{s_i : i \in I\}$ has form*

$$(3.4.10) \qquad J(S) = \left\{ g : g = \sum_{i \in I} g^i s_i, g^i \in Z \right\}$$

*where the set $\{i \in I : g^i \neq 0\}$ is finite.*

PROOF. We prove the theorem by induction based on the theorems [18]-5.1, page 79 and 2.7.5.

For any $s_k \in S$, let $g^i = \delta_k^i$. Then

$$(3.4.11) \qquad s_k = \sum_{i \in I} g^i s_i$$

$s_k \in J(S)$ follows from (3.4.10), (3.4.11).

Let  $g_1$, $g_2 \in X_k \subseteq J(S)$.  Since $G$ is Abelian group, then, according to the statement 2.7.5.3, $g_1 + g_2 \in J(S)$. According to the equality (3.4.10), there exist $Z$-numbers $g_1^i$, $g_2^i$, $i \in I$, such that

$$(3.4.12) \qquad\qquad g_1 = \sum_{i \in I} g_1^i s_i \qquad g_2 = \sum_{i \in I} g_2^i s_i$$

where sets

$$(3.4.13) \qquad\qquad H_1 = \{i \in I : g_1^i \neq 0\} \qquad H_2 = \{i \in I : g_2^i \neq 0\}$$

are finite. From the equality (3.4.12), it follows that

$$(3.4.14) \qquad g_1 + g_2 = \sum_{i \in I} g_1^i s_i + \sum_{i \in I} g_2^i s_i = \sum_{i \in I} (g_1^i s_i + g_2^i s_i)$$

The equality

$$(3.4.15) \qquad\qquad g_1 + g_2 = \sum_{i \in I} (g_1^i + g_2^i) s_i$$

follows from equalities (3.4.6), (3.4.14). From the equality (3.4.13), it follows that the set

$$\{i \in I : g_1^i + g_2^i \neq 0\} \subseteq H_1 \cup H_2$$

is finite. From the equality (3.4.15), it follows that  $g_1 + g_2 \in J(S)$.          □

CONVENTION 3.4.4. *We will use Einstein summation convention in which repeated index (one above and one below) implies summation with respect to repeated index. In this case we assume that we know the set of summation index and do not use summation symbol*

$$g^i s_i = \sum_{i \in I} g^i s_i$$

□

THEOREM 3.4.5. *The homomorphism*

$$f : G \to G'$$

*holds the equality*

$$(3.4.16) \qquad\qquad f(na) = nf(a)$$

PROOF. From the equality (3.4.1), it follows that the equality (3.4.16) is true when $n = 0$.

• Let the equality (3.4.16) is true when $n = k \geq 0$. Then

$$f(ka) = kf(a)$$

The equality

$$f((k+1)a) = f(ka + a) = f(ka) + f(a) = kf(a) + f(a)$$
$$= (k+1)f(a)$$

follows from the equality (3.4.2). Therefore, the equality (3.4.16) is true when $n = k + 1$. According to mathematical induction, the equality (3.4.16) is true for any $n \geq 0$.

- Let the equality (3.4.6) is true when $n = k \le 0$. Then

$$f(ka) = kf(a)$$

The equality

$$f((k-1)a) = f(ka - a) = f(ka) - f(a) = kf(a) - f(a)$$
$$= (k-1)f(a)$$

follows from the equality (3.4.3). Therefore, the equality (3.4.16) is true when $n = k - 1$. According to mathematical induction, the equality (3.4.16) is true for any $n \le 0$.

- Therefore, the equality (3.4.16) is true for any $n \in Z$.

$\square$

DEFINITION 3.4.6. *Let* $S = \{s_i : i \in I\}$ *be set of G-numbers. The expression* $g^i s_i$ *is called* **linear combination** *of G-numbers* $s_i$. *A G-number* $g = g^i s_i$ *is called* **linearly dependent** *on G-numbers* $s_i$. $\square$

The following definition follows from the theorems 2.7.5, 3.4.3 and from the definition 2.7.4.

DEFINITION 3.4.7. $J(S)$ *is called* **subgroup generated by set** $S$, *and* $S$ *is a* **generating set** *of subgroup* $J(S)$. *In particular, a* **generating set** *of Abelian group is a subset* $X \subset G$ *such that* $J(X) = G$. $\square$

The following definition follows from the theorems 2.7.5, 3.4.3 and from the definition 2.7.14.

DEFINITION 3.4.8. *If the set* $X \subset G$ *is generating set of Abelian group* $G$, *then any set* $Y$, $X \subset Y \subset G$ *also is generating set of Abelian group* $G$. *If there exists minimal set* $X$ *generating the Abelian group* $G$, *then the set* $X$ *is called* **quasibasis** *of Abelian group* $G$. $\square$

DEFINITION 3.4.9. *The set of G-numbers*[3.8]$g_i$, $i \in I$, *is* **linearly independent** *if* $w = 0$ *follows from the equation*

$$w^i g_i = 0$$

*Otherwise the set of G-numbers* $g_i$, $i \in I$, *is* **linearly dependent**. $\square$

DEFINITION 3.4.10. *Let the set of G-numbers* $S = \{s_i : i \in I\}$ *be quasibasis. If the set of G-numbers* $S$ *is linearly independent, then quasibasis* $S$ *is called* **basis** *of Abelian group* $G$. $\square$

THEOREM 3.4.11. *Product in the ring* $D$ *is bilinear over the ring* $Z$

$$(3.4.17) \qquad\qquad (nd_1)d_2 = n(d_1 d_2)$$

---

[3.8] I follow to the definition in [1], page 130.

(3.4.18)                                    $$d_1(nd_2) = n(d_1 d_2)$$

Proof. From the equality (3.4.1), it follows that the equality (3.4.17) is true when $n = 0$.

- Let the equality (3.4.17) is true when $n = k \geq 0$. Then

$$(kd_1)d_2 = k(d_1 d_2)$$

  The equality

  $$((k+1)d_1)d_2 = (kd_1 + d_1)d_2 = (kd_1)d_2 + d_1 d_2 = k(d_1 d_2) + d_1 d_2$$
  $$= (k+1)(d_1 d_2)$$

  follows from equalities (3.4.2), (3.2.2). Therefore, the equality (3.4.17) is true when $n = k + 1$. According to mathematical induction, the equality (3.4.17) is true for any $n \geq 0$.

- Let the equality (3.4.6) is true when $n = k \leq 0$. Then

$$(kd_1)d_2 = k(d_1 d_2)$$

  The equality

  $$((k-1)d_1)d_2 = (kd_1 - d_1)d_2 = (kd_1)d_2 - d_1 d_2 = k(d_1 d_2) - d_1 d_2$$
  $$= (k-1)(d_1 d_2)$$

  follows from equalities (3.4.2), (3.2.2). Therefore, the equality (3.4.17) is true when $n = k - 1$. According to mathematical induction, the equality (3.4.17) is true for any $n \leq 0$.

- Therefore, the equality (3.4.17) is true for any $n \in Z$.

From the equality (3.4.1), it follows that the equality (3.4.18) is true when $n = 0$.

- Let the equality (3.4.18) is true when $n = k \geq 0$. Then

$$d_1(kd_2) = k(d_1 d_2)$$

  The equality

  $$d_1((k+1)d_2) = d_1(kd_2 + d_2) = d_1(kd_2) + d_1 d_2 = k(d_1 d_2) + d_1 d_2$$
  $$= (k+1)(d_1 d_2)$$

  follows from equalities (3.4.2), (3.2.1). Therefore, the equality (3.4.18) is true when $n = k + 1$. According to mathematical induction, the equality (3.4.18) is true for any $n \geq 0$.

- Let the equality (3.4.6) is true when $n = k \leq 0$. Then

$$d_1(kd_2) = k(d_1 d_2)$$

  The equality

  $$d_1((k-1)d_2) = d_1(kd_2 - d_2) = d_1(kd_2) - d_1 d_2 = k(d_1 d_2) - d_1 d_2$$
  $$= (k-1)(d_1 d_2)$$

  follows from equalities (3.4.2), (3.2.1). Therefore, the equality (3.4.18) is true when $n = k - 1$. According to mathematical induction, the equality (3.4.18) is true for any $n \leq 0$.

- Therefore, the equality (3.4.18) is true for any $n \in Z$.

□

Theorem 3.4.12. *Let $D$ be ring with unit $1_D$.*

3.4.12.1: *The map*

$$f : n \in Z \to n 1_D \in D$$

*is homomorphism of the ring $Z$ into the ring $D$ and allows to identify $Z$-number $n$ and $D$-number $n 1_D$.*

3.4.12.2: *The ring $Z / \ker f$ is subring of ring $D$.*

Proof. Let $n_1, n_2 \in Z$. The equality

$$(3.4.19) \qquad n_1 1_D + n_2 1_D = (n_1 + n_2) 1_D$$

follows from the equality (3.4.6). The equality

$$(3.4.20) \qquad \begin{aligned} (n_1 1_D)(n_2 1_D) &= n_1 (1_D (n_2 1_D)) = n_1 (n_2 (1_D 1_D)) \\ &= n_1 (n_2 1_D) = (n_1 n_2) 1_D \end{aligned}$$

follows from equalities (3.4.17), (3.4.18). The statement 3.4.12.1 follows from equalities (3.4.19), (3.4.20). The statement 3.4.12.2 follows from the statement 3.4.12.1.

□

## 3.5. Basis of Multiplicative Abelian Group

In this section, $G$ is multiplicative Abelian group.

Definition 3.5.1. *The action of ring of rational integers $Z$ in multiplicative Abelian group $G$ is defined using following rules*

$$(3.5.1) \qquad g^0 = e$$

$$(3.5.2) \qquad g^{n+1} = g^n g$$

$$(3.5.3) \qquad g^{n-1} = g^n g^{-1}$$

□

Theorem 3.5.2. *The action of ring of rational integers $Z$ in multiplicative Abelian group $G$ defined in the definition 3.5.1 is representation. The following equalities are true*

$$(3.5.4) \qquad a^1 = a$$

$$(3.5.5) \qquad a^{nm} = (a^n)^m$$

$$(3.5.6) \qquad a^{m+n} = a^m a^n$$

$$(3.5.7) \qquad a^{m-n} = a^m a^{-n}$$

$$(3.5.8) \qquad (ab)^n = a^n b^n$$

Proof. The equality (3.5.4) follows from the equality (3.5.1) and from the equality (3.5.2) when $n = 0$.

From the equality (3.5.1), it follows that the equality (3.5.6) is true when $n = 0$.

- Let the equality $(3.5.6)$ is true when $n = k \geq 0$. Then

$$a^{m+k} = a^m a^k$$

The equality

$$a^{m+(k+1)} = a^{(m+k)+1} = a^{m+k}a = a^m a^k a$$
$$= a^m a^{k+1}$$

follows from the equality $(3.5.2)$. Therefore, the equality $(3.5.6)$ is true when $n = k+1$. According to mathematical induction, the equality $(3.5.6)$ is true for any $n \geq 0$.

- Let the equality $(3.5.6)$ is true when $n = k \leq 0$. Then

$$a^{m+k} = a^m a^k$$

The equality

$$a^{m+(k-1)} = a^{(m+k)-1} = a^{m+k}a^{-1} = a^m a^k a^{-1}$$
$$= a^m a^{k-1}$$

follows from the equality $(3.5.3)$. Therefore, the equality $(3.5.6)$ is true when $n = k-1$. According to mathematical induction, the equality $(3.5.6)$ is true for any $n \leq 0$.

- Therefore, the equality $(3.5.6)$ is true for any $n \in Z$.

The equality

$$(3.5.9) \qquad\qquad a^{k+n}a^{-n} = a^k$$

follows from the equality $(3.5.6)$. The equality $(3.5.7)$ follows from the equality $(3.5.9)$, if we assume $m = k + n$, $k = m - n$.

From the equality $(3.5.1)$, it follows that the equality $(3.5.5)$ is true when $n = 0$.

- Let the equality $(3.5.5)$ is true when $n = k \geq 0$. Then

$$a^{km} = (a^m)^k$$

The equality

$$a^{(k+1)m} = a^{km+m} = a^{km}a^m = (a^m)^k a^m$$
$$= (a^m)^{k+1}$$

follows from equalities $(3.5.2)$, $(3.5.6)$. Therefore, the equality $(3.5.5)$ is true when $n = k + 1$. According to mathematical induction, the equality $(3.5.5)$ is true for any $n \geq 0$.

- Let the equality $(3.5.6)$ is true when $n = k \leq 0$. Then

$$a^{km} = (a^m)^k$$

The equality

$$a^{(k-1)m} = a^{km-m} = a^{km}a^{-m} = (a^m)^k a^{-m}$$
$$= (a^m)^{k-1}$$

follows from equalities $(3.5.3)$, $(3.5.7)$. Therefore, the equality $(3.5.5)$ is true when $n = k - 1$. According to mathematical induction, the equality $(3.5.5)$ is true for any $n \leq 0$.

- Therefore, the equality $(3.5.5)$ is true for any $n \in Z$.

From the equality (3.5.1), it follows that the equality (3.5.8) is true when $n = 0$.

- Let the equality (3.5.8) is true when $n = k \geq 0$. Then

$$(ab)^k = a^k b^k$$

The equality

$$\begin{aligned}
(ab)^{k+1} = (ab)^k ab &= a^k b^k ab \\
&= a^k a b^k b \\
&= a^{k+1} b^{k+1}
\end{aligned}$$

follows from the equality (3.5.2). Therefore, the equality (3.5.8) is true when $n = k + 1$. According to mathematical induction, the equality (3.5.8) is true for any $n \geq 0$.

- Let the equality (3.5.6) is true when $n = k \leq 0$. Then

$$(ab)^k = a^k b^k$$

The equality

$$\begin{aligned}
(ab)^{k-1} = (ab)^k (ab)^{-1} &= a^k b^k a^{-1} b^{-1} \\
&= a^k a^{-1} b^k b^{-1} \\
&= a^{k-1} b^{k-1}
\end{aligned}$$

follows from the equality (3.5.3). Therefore, the equality (3.5.8) is true when $n = k - 1$. According to mathematical induction, the equality (3.5.8) is true for any $n \leq 0$.

- Therefore, the equality (3.5.8) is true for any $n \in Z$.

From the equality (3.5.8), it follows that the map

$$\varphi(n) : a \in G \to a^n \in G$$

is an endomorphism of Abelian group $G$. From equalities (3.5.6), (3.5.5), it follows that the map

$$\varphi : Z \to \operatorname{End}(Ab, G)$$

is a homomorphism of the ring $Z$. According to the definition 2.3.1, the map $\varphi$ is representation of ring of rational integers $Z$ in Abelian group $G$. □

THEOREM 3.5.3. *The set of $G$-numbers generated by the set $S = \{s_i : i \in I\}$ has form*

$$(3.5.10) \qquad J(S) = \left\{ g : g = \prod_{i \in I} s_i^{g^i}, g^i \in Z \right\}$$

*where the set $\{i \in I : g^i \neq 0\}$ is finite.*

PROOF. We prove the theorem by induction based on the theorems [18]-5.1, page 79 and 2.7.5.

For any $s_k \in S$, let $g^i = \delta_k^i$. Then

$$(3.5.11) \qquad s_k = \prod_{i \in I} s_i^{g^i}$$

$s_k \in J(S)$ follows from (3.5.10), (3.5.11).

Let $g_1$, $g_2 \in X_k \subseteq J(S)$. Since $G$ is Abelian group, then, according to the statement 2.7.5.3, $g_1 g_2 \in J(S)$. According to the equality (3.5.10), there exist $Z$-numbers $g_1^i$, $g_2^i$, $i \in I$, such that

$$(3.5.12) \qquad g_1 = \prod_{i \in I} s_i^{g_1^i} \qquad g_2 = \prod_{i \in I} s_i^{g_2^i}$$

where sets

$$(3.5.13) \qquad H_1 = \{i \in I : g_1^i \neq 0\} \qquad H_2 = \{i \in I : g_2^i \neq 0\}$$

are finite. From the equality (3.5.12), it follows that

$$(3.5.14) \qquad g_1 g_2 = \prod_{i \in I} s_i^{g_1^i} \prod_{i \in I} s_i^{g_2^i} = \prod_{i \in I} (s_i^{g_1^i} s_i^{g_2^i})$$

The equality

$$(3.5.15) \qquad g_1 g_2 = \prod_{i \in I} s_i^{g_1^i + g_2^i}$$

follows from equalities (3.5.6), (3.5.14). From the equality (3.5.13), it follows that the set

$$\{i \in I : g_1^i g_2^i \neq 0\} \subseteq H_1 \cup H_2$$

is finite. From the equality (3.5.15), it follows that $g_1 g_2 \in J(S)$. $\qquad\square$

CONVENTION 3.5.4. *We will use product convention in which repeated index (one at the base and one at the exponent) implies product with respect to repeated index. In this case we assume that we know the set of product index and do not use product symbol*

$$s_i^{g^i} = \prod_{i \in I} s_i^{g^i}$$

$\qquad\square$

THEOREM 3.5.5. *The homomorphism*

$$f : G \to G'$$

*holds the equality*

$$(3.5.16) \qquad f(a^n) = f(a)^n$$

PROOF. From the equality (3.5.1), it follows that the equality (3.5.16) is true when $n = 0$.

• Let the equality (3.5.16) is true when $n = k \geq 0$. Then

$$f(a^k) = f(a)^k$$

The equality

$$f(a^{k+1}) = f(a^k a) = f(a^k) f(a) = f(a)^k f(a)$$
$$= f(a)^{k+1}$$

follows from the equality (3.5.2). Therefore, the equality (3.5.16) is true when $n = k + 1$. According to mathematical induction, the equality (3.5.16) is true for any $n \geq 0$.

- Let the equality (3.5.6) is true when $n = k \le 0$. Then

$$f(a^k) = f(a)^k$$

The equality

$$f(a^{k-1}) = f(a^k a^{-1}) = f(a^k) f(a)^{-1} = f(a)^k f(a)^{-1}$$
$$= f(a)^{k-1}$$

follows from the equality (3.5.3). Therefore, the equality (3.5.16) is true when $n = k - 1$. According to mathematical induction, the equality (3.5.16) is true for any $n \le 0$.

- Therefore, the equality (3.5.16) is true for any $n \in Z$.

$\square$

DEFINITION 3.5.6. *Let $S = \{s_i : i \in I\}$ be set of $G$-numbers. The expression $s_i^{g^i}$ is called* **linear combination** *of $G$-numbers $s_i$. A $G$-number $g = s_i^{g^i}$ is called* **linearly dependent** *on $G$-numbers $s_i$.* $\square$

The following definition follows from the theorems 2.7.5, 3.5.3 and from the definition 2.7.4.

DEFINITION 3.5.7. *$J(S)$ is called* **subgroup generated by set** *$S$, and $S$ is a* **generating set** *of subgroup $J(S)$. In particular, a* **generating set** *of Abelian group is a subset $X \subset G$ such that $J(X) = G$.* $\square$

The following definition follows from the theorems 2.7.5, 3.5.3 and from the definition 2.7.14.

DEFINITION 3.5.8. *If the set $X \subset G$ is generating set of Abelian group $G$, then any set $Y$, $X \subset Y \subset G$ also is generating set of Abelian group $G$. If there exists minimal set $X$ generating the Abelian group $G$, then the set $X$ is called* **quasibasis** *of Abelian group $G$.* $\square$

DEFINITION 3.5.9. *The set of $G$-numbers[3.9] $g_i$, $i \in I$, is* **linearly independent** *if $w = e$ follows from the equation*

$$g_i^{w^i} = e$$

*Otherwise the set of $G$-numbers $g_i$, $i \in I$, is* **linearly dependent**. $\square$

DEFINITION 3.5.10. *Let the set of $G$-numbers $S = \{s_i : i \in I\}$ be quasibasis. If the set of $G$-numbers $S$ is linearly independent, then quasibasis $S$ is called* **basis** *of Abelian group $G$.* $\square$

## 3.6. Free Abelian Group

---

[3.9] I follow to the definition in [1], page 130.

[3.10] See also the definition in [1], pages 37.

DEFINITION 3.6.1. *Let $S$ be a set and [3.10] $Ab(S, *)$ be category objects of which are maps*

$$f : S \to G$$

*of the set $S$ into multiplicative Abelian groups. If*

$$f : S \to G$$
$$f' : S \to G'$$

*are objects of the category $Ab(S, *)$, then morphism from $f$ to $f'$ is the homomorphism of groups*

$$g : G \to G'$$

*such that the diagram*

*is commutative.*  □

DEFINITION 3.6.2. *Let $S$ be a set and [3.10] $Ab(S, +)$ be category objects of which are maps*

$$f : S \to G$$

*of the set $S$ into additive Abelian groups. If*

$$f : S \to G$$
$$f' : S \to G'$$

*are objects of the category $Ab(S, +)$, then morphism from $f$ to $f'$ is the homomorphism of groups*

$$g : G \to G'$$

*such that the diagram*

*is commutative.*  □

THEOREM 3.6.3. *There exists [3.11] universally repelling object $G$ of the category $Ab(S, *)$*

3.6.3.1: *$S \subseteq G$.*

3.6.3.2: *The set $S$ generates multiplicative Abelian group $G$.*

*Abelian group $G$ is called* **free Abelian group** *generated by the set $S$.*

THEOREM 3.6.6. *There exists [3.11] universally repelling object $G$ of the category $Ab(S, +)$*

3.6.6.1: *$S \subseteq G$.*

3.6.6.2: *The set $S$ generates additive Abelian group $G$.*

*Abelian group $G$ is called* **free Abelian group** *generated by the set $S$.*

PROOF. Let $G$ be set of maps

$$h : S \to Z$$

such that the set

$$H = \{x \in S : h(x) \neq 0\}$$

is finite.

PROOF. Let $G$ be set of maps

$$h : S \to Z$$

such that the set

$$H = \{x \in S : h(x) \neq 0\}$$

is finite.

LEMMA 3.6.4. *We introduce the following operation on the set $G$*

$$(h_1 h_2)(x) = h_1(x) + h_2(x)$$

*The set $G$ is multiplicative Abelian group.*

LEMMA 3.6.7. *We introduce the following operation on the set $G$*

$$(h_1 + h_2)(x) = h_1(x) + h_2(x)$$

*The set $G$ is additive Abelian group.*

---

[3.11]  See also similar statement in [1], pages 38.

PROOF. Let $h_1$, $h_2 \in G$. According to the definition, sets

$$H_1 = \{x \in S : h_1(x) \neq 0\}$$
$$H_2 = \{x \in S : h_2(x) \neq 0\}$$

are finite. Then the set

$$\{x \in S : (h_1 h_2)(x) \neq 0\} \subseteq H_1 \cup H_2$$

is finite. Therefore, $h_1 h_2 \in G$. Properties of Abelian group are evident. ⊙

Let $k \in Z$, $x \in S$. Consider the map $h = k * x$ defined by the equality

$$k * x(y) = \begin{cases} k & y = x \\ 0 & y \neq x \end{cases}$$

The map

$$f : x \in S \to 1 * x \in G$$

is injective and allows us to identify $S$ and image $f(S) \subseteq G$ (the statement 3.6.3.1) and to use notation

$$(1 * x)^k = k * x$$

According to the definition of the group $G$, we can write any map $h \in G$ as

$$(3.6.1) \qquad h = \prod_{x \in S} k_x * x$$

where $k_x \in Z$, $x \in S$, and the set $\{x : k_x \neq 0\}$ is finite. The statement 3.6.3.2 follows from the theorem 3.5.3.

LEMMA 3.6.5. *The representation* (3.6.1) *of the map* $h$ *is unique.*

PROOF. If the map $h$ admits representations

$$(3.6.2) \qquad \prod_{x \in S} k_{1x} * x = \prod_{x \in S} k_{2x} * x$$

then the equality

$$(3.6.3) \qquad \prod_{x \in S} (k_{1x} - k_{2x}) * x = e$$

follows from the equality (3.6.2). $k_{1x} = k_{2x}$, $x \in S$, follows from the equality (3.6.3). ⊙

Let

$$f' : S \to G'$$

PROOF. Let $h_1$, $h_2 \in G$. According to the definition, sets

$$H_1 = \{x \in S : h_1(x) \neq 0\}$$
$$H_2 = \{x \in S : h_2(x) \neq 0\}$$

are finite. Then the set

$$\{x \in S : (h_1 + h_2)(x) \neq 0\} \subseteq H_1 \cup H_2$$

is finite. Therefore, $h_1 + h_2 \in G$. Properties of Abelian group are evident. ⊙

Let $k \in Z$, $x \in S$. Consider the map $h = k * x$ defined by the equality

$$k * x(y) = \begin{cases} k & y = x \\ 0 & y \neq x \end{cases}$$

The map

$$f : x \in S \to 1 * x \in G$$

is injective and allows us to identify $S$ and image $f(S) \subseteq G$ (the statement 3.6.6.1) and to use notation

$$k(1 * x) = k * x$$

According to the definition of the group $G$, we can write any map $h \in G$ as

$$(3.6.6) \qquad h = \sum_{x \in S} k_x * x$$

where $k_x \in Z$, $x \in S$, and the set $\{x : k_x \neq 0\}$ is finite. The statement 3.6.6.2 follows from the theorem 3.4.3.

LEMMA 3.6.8. *The representation* (3.6.6) *of the map* $h$ *is unique.*

PROOF. If the map $h$ admits representations

$$(3.6.7) \qquad \sum_{x \in S} k_{1x} * x = \sum_{x \in S} k_{2x} x$$

then the equality

$$(3.6.8) \qquad \sum_{x \in S} (k_{1x} - k_{2x}) * x = 0$$

follows from the equality (3.6.7). $k_{1x} = k_{2x}$, $x \in S$, follows from the equality (3.6.8). ⊙

Let

$$f' : S \to G'$$

be the map of the set $S$ into Abelian group $G'$. We define the homomorphism

$$g : G \to G'$$

with request that following diagram is commutative

(3.6.4)

From the diagram (3.6.4), it follows that

(3.6.5) $$g(1 * x) = f'(x)$$

According to the theorem 3.5.5, since the map $g$ is homomorphism of Abelian group, then the equality

$$g(h) = g\left(\prod_{x \in S} k_x * x\right) = \prod_{x \in S} g(k_x * x)$$
$$= \prod_{x \in S} g(1 * x)^{k_x} = \prod_{x \in S} f'(x)^{k_x}$$

for any map $h \in G$ follows from equalities (3.6.5), (3.6.1). Therefore, homomorphism $g$ is defined uniquely. According to definitions 2.2.2, 3.6.1, Abelian group $G$ is free group. □

be the map of the set $S$ into Abelian group $G'$. We define the homomorphism

$$g : G \to G'$$

with request that following diagram is commutative

(3.6.9)

From the diagram (3.6.9), it follows that

(3.6.10) $$g(1 * x) = f'(x)$$

According to the theorem 3.4.5, since the map $g$ is homomorphism of Abelian group, then the equality

$$g(h) = g\left(\sum_{x \in S} k_x * x\right) = \sum_{x \in S} g(k_x * x)$$
$$= \sum_{x \in S} k_x g(1 * x) = \sum_{x \in S} k_x f'(x)$$

for any map $h \in G$ follows from equalities (3.6.10), (3.6.6). Therefore, homomorphism $g$ is defined uniquely. According to definitions 2.2.2, 3.6.2, Abelian group $G$ is free group. □

THEOREM 3.6.9. *Multiplicative Abelian group $G$ is free iff Abelian group $G$ has a basis $S = \{s_i : i \in I\}$ and the set of $G$-numbers $S$ is linear independent.*

PROOF. According to the definitions 3.6.1 and proof of the theorem 3.6.3, there exist the set $S \subseteq G$ such that following diagram is commutative

□

THEOREM 3.6.10. *Additive Abelian group $G$ is free iff Abelian group $G$ has a basis $S = \{s_i : i \in I\}$ and the set of $G$-numbers $S$ is linear independent.*

PROOF. According to the definitions 3.6.2 and proof of the theorem 3.6.6, there exist the set $S \subseteq G$ such that following diagram is commutative

□

THEOREM 3.6.11. *Let $A$ be Abelian additive monoid. Consider category $\mathcal{K}(A)$ of homomorphisms of monoid $A$ into Abelian groups. Let $B$ and $C$ be Abelian groups. Then, for homomorphisms*

$$f : A \to B$$

$$g : A \to C$$

*morphism of category $\mathcal{K}(A)$ is homomorphism*

$$h : B \to C$$

*such that the diagram*

$$A \xrightarrow{\ f\ } B$$
$$\searrow{\scriptstyle g} \quad \swarrow{\scriptstyle h}$$
$$C$$

*is commutative. There exists Abelian group $K(A)$ and universally repelling homomorphism*

$$\gamma : A \to K(A)$$

*The group $K(A)$ is called* **Grothendieck group** *and the map $\gamma$ is called* **canonical map***.*

REMARK 3.6.12. *According to the definition 2.2.2, if $B$ is Abelian group and*

$$f : A \to B$$

*is homomorphism, then there exists unique homomorphism*

$$g : K(A) \to B$$

*such that the following diagram*

$$A \xrightarrow{\ \gamma\ } K(A)$$
$$\searrow{\scriptstyle f} \quad \swarrow{\scriptstyle g}$$
$$C$$

*is commutative.* ⊙

PROOF. Let $F(A)$ be free Abelian group generated by the set $A$. We denote the generator of the group $F(A)$, corresponding to $A$-number $a$ by $[a]$. Let $B$ be the subgroup generated by all $F(A)$-numbers of type

$$[a + b] - [a] - [b] \qquad a, b \in A$$

Let $K(A) = F(A)/B$. We define the homomorphism $\gamma$ such way that the following diagram

$$A \xrightarrow{\ \gamma\ } K(A)$$
$$\searrow{\scriptstyle f} \quad \nearrow{\scriptstyle g}$$
$$F(A)$$

$$f : a \in A \to [a] \in F(A)$$
$$g : a \in F(A) \to aB \in K(A)$$

is commutative. □

DEFINITION 3.6.13. *The* **cancellation law** *holds in additive monoid $A$, if, for every $A$-numbers $a$, $b$, $c$, the equality $a + c = b + c$ implies $a = b$.* □

THEOREM 3.6.14. *If the cancellation law holds in additive monoid $A$, then the canonical map*

$$\gamma : A \to K(A)$$

*is injective.*

PROOF. Consider tuples $(a, b)$, $a$, $b \in A$. The tuple $(a, b)$ is equivalent to the tuple $(a', b')$ when

$$a + b' = a' + b$$

We define addition by the equality

$$(a, b) + (c, d) = (a + c, b + d)$$

Equivalence classes of tuples form group, whose zero element is the class of $(0, 0)$. The tuple $(b, a)$ is negative for the tuple $(a, b)$. According to cancelation law, the homomorphism

$$a \to \text{class of the tuple } (0, a)$$

is injective.                                                                     □

EXAMPLE 3.6.15. *In the set of numbers*

$$A = \{0, 1, ...\}$$

*we introduce addition by the equality*

$$(3.6.11) \qquad\qquad a + b = \max(a, b)$$

LEMMA 3.6.16. *Addition* (3.6.11) *is commutative.*

PROOF. *The lemma follows from the definition* (3.6.11).                          ⊙

LEMMA 3.6.17. *Addition* (3.6.11) *is associative.*

PROOF. *To prove the equality*

$$(3.6.12) \qquad\qquad a + (b + c) = (a + b) + c$$

*it is enough to consider following relations.*

- *If* $a \le \max(b, c)$, *then or* $a \le b$, *either* $a \le c$.
  - *Let* $b \le c$. *Then*

$$(3.6.13) \qquad\qquad a + (b + c) = a + c = c$$

$$(3.6.14) \qquad\qquad (a + b) + c = \max(a, b) + c = c$$

  *The equality* (3.6.12) *follows from equalities* (3.6.13), (3.6.14).
  - *Let* $c \le b$. *Then*

$$(3.6.15) \qquad\qquad a + (b + c) = a + b = b$$

$$(3.6.16) \qquad\qquad (a + b) + c = b + c = b$$

  *The equality* (3.6.12) *follows from equalities* (3.6.15), (3.6.16). *Therefore*

$$(3.6.17) \qquad\qquad a + (b + c) = a$$

$$(3.6.18) \qquad\qquad (a + b) + c = a + c = a$$

- *If $a > \max(b, c)$, then $a > b$, $a > c$. The equality (3.6.12) follows from equalities (3.6.17), (3.6.18).*

*The lemma follows from the equality (3.6.12).* ⊙

LEMMA 3.6.18. *$A$-number $0$ is zero element of addition (3.6.11).*

PROOF. *Since $\forall a \in A$, $0 \le a$, the lemma follows from the definition (3.6.11).* ⊙

THEOREM 3.6.19. *The set $A$ equiped by addition (3.6.11), is Abelian additive monoid.*

PROOF. *The theorem follows from lemmas 3.6.16, 3.6.17, 3.6.18 and from the definition 3.1.1.* ⊙

THEOREM 3.6.20. *The cancellation law does not hold in additive monoid $A$.*

PROOF. *To prove the theorem, it is enough to see that*

$$2 + 9 = 9$$
$$4 + 9 = 9$$

*However, $2 \ne 4$.* ⊙ □

## 3.7. Direct Sum of Abelian Groups

DEFINITION 3.7.1. *Coproduct in category of multiplicative Abelian groups $Ab(*)$ is called* **direct sum**. [3.12] *We will use notation $A \oplus B$ for direct sum of Abelian groups $A$ and $B$.* □

DEFINITION 3.7.2. *Coproduct in category of additive Abelian groups $Ab(+)$ is called* **direct sum**. [3.12] *We will use notation $A \oplus B$ for direct sum of Abelian groups $A$ and $B$.* □

THEOREM 3.7.3. *Let $\{A_i, i \in I\}$ be set of multiplicative Abelian groups. Let*

$$A \subseteq \prod_{i \in I} A_i$$

*be such set that*

$$A = \{(x_i, i \in I) : |\{i : x_i \ne e\}| < \infty\}$$

*Then* [3.13]

$$A = \bigoplus_{i \in I} A_i$$

PROOF. *According to the theorem 2.2.6, there exists Abelian group $B =$*

THEOREM 3.7.4. *Let $\{A_i, i \in I\}$ be set of additive Abelian groups. Let*

$$A \subseteq \prod_{i \in I} A_i$$

*be such set that*

$$A = \{(x_i, i \in I) : |\{i : x_i \ne 0\}| < \infty\}$$

*Then* [3.13]

$$A = \bigoplus_{i \in I} A_i$$

PROOF. *According to the theorem 2.2.6, there exists Abelian group $B =$*

---

[3.12]  See also definition in [1], pages 36, 37.
[3.13]  See also proposition [1]-7.1, page 37.

$\prod A_i$. According to the statement 2.2.6.3, $B$-number $a$ can be represented as tuple $\{a_i, i \in I\}$ $A_i$-numbers. According to the statement 2.2.6.4, the product in group $B$ is defined by the equality

$$(3.7.1) \qquad (a_i, i \in I)(b_i, i \in I)$$
$$= (a_i b_i, i \in I)$$

According to construction, $A$ is subgroup of Abelian group $B$. The map

$$\lambda_j : A_j \to A$$

defined by the equality

$$(3.7.2) \qquad \lambda_j(x) = (\delta_j^i x, i \in I)$$

is an injective homomorphism.

Let

$$\{f_i : A_i \to C, i \in I\}$$

be set of homomorphisms into multiplicative Abelian group $C$. We define the map

$$f : A \to C$$

by the equality

$$(3.7.3) \qquad f\left(\bigoplus_{i \in I} x_i\right) = \prod_{i \in I} f_i(x_i)$$

The product in the right side of the equality (3.7.3) is finite, since all factors, except for a finite number, equal $e$. From the equality

$$f_i(x_i y_i) = f_i(x_i) f_i(y_i)$$

and the equality (3.7.3), it follows that

$$f(\bigoplus_{i \in I}(x_i y_i)) = \prod_{i \in I} f_i(x_i y_i)$$
$$= \prod_{i \in I} (f_i(x_i) f_i(y_i))$$
$$= \prod_{i \in I} f_i(x_i) \prod_{i \in I} f_i(y_i)$$
$$= f(\bigoplus_{i \in I} x_i) f(\bigoplus_{i \in I} y_i)$$

Therefore, the map $f$ is homomorphism of Abelian group. The equality

$$f \circ \lambda_j(x) = \prod_{i \in I} f_i(\delta_j^i x) = f_j(x)$$

$\prod A_i$. According to the statement 2.2.6.3, $B$-number $a$ can be represented as tuple $\{a_i, i \in I\}$ $A_i$-numbers. According to the statement 2.2.6.4, the sum in group $B$ is defined by the equality

$$(3.7.4) \qquad (a_i, i \in I) + (b_i, i \in I)$$
$$= (a_i + b_i, i \in I)$$

According to construction, $A$ is subgroup of Abelian group $B$. The map

$$\lambda_j : A_j \to A$$

defined by the equality

$$(3.7.5) \qquad \lambda_j(x) = (\delta_j^i x, i \in I)$$

is an injective homomorphism.

Let

$$\{f_i : A_i \to C, i \in I\}$$

be set of homomorphisms into additive Abelian group $C$. We define the map

$$f : A \to C$$

by the equality

$$(3.7.6) \qquad f\left(\bigoplus_{i \in I} x_i\right) = \sum_{i \in I} f_i(x_i)$$

The sum in the right side of the equality (3.7.6) is finite, since all terms, except for a finite number, equal 0. From the equality

$$f_i(x_i + y_i) = f_i(x_i) + f_i(y_i)$$

and the equality (3.7.6), it follows that

$$f(\bigoplus_{i \in I}(x_i + y_i)) = \sum_{i \in I} f_i(x_i + y_i)$$
$$= \sum_{i \in I} (f_i(x_i) + f_i(y_i))$$
$$= \sum_{i \in I} f_i(x_i) + \sum_{i \in I} f_i(y_i)$$
$$= f(\bigoplus_{i \in I} x_i) + f(\bigoplus_{i \in I} y_i)$$

Therefore, the map $f$ is homomorphism of Abelian group. The equality

$$f \circ \lambda_j(x) = \sum_{i \in I} f_i(\delta_j^i x) = f_j(x)$$

follows from equalities (3.7.5), (3.7.6). Since the map $\lambda_i$ is injective, then the follows from equalities (3.7.2), (3.7.3). Since the map $\lambda_i$ is injective, then the map $f$ is unique. Therefore, the theorem folows from definitions 2.2.12, 3.7.1. □

map $f$ is unique. Therefore, the theorem folows from definitions 2.2.12, 3.7.2. □

**Theorem** 3.7.5. *Let $B$, $C$ be subgroups of multiplicative Abelian group $A$. Let $A = BC$, $B \cap C = e$. Then* [3.14] *$A = B \oplus C$.*

**Theorem** 3.7.6. *Let $B$, $C$ be subgroups of additive Abelian group $A$. Let $A = B + C$, $B \cap C = 0$. Then* [3.14] *$A = B \oplus C$.*

**Proof.** The map

$$(3.7.7) \quad f : (x, y) \in B \times C \to xy \in A$$

is homomorphism of groups. This statement follows from the equality

$$f(x_1 x_2, y_1 y_2) = x_1 x_2 y_1 y_2$$
$$= x_1 y_1 x_2 y_2$$
$$= f(x_1, y_1) f(x_2, y_2)$$

Let $f(x, y) = e$. From the equality (3.7.7), it follows that

$$(3.7.8) \quad xy = e \quad x \in B \quad y \in C$$

From the equality (3.7.8), it follows that $x = y^{-1} \in C$. Therefore, $x = y = e$. Therefore, the map $f$ is isomorphism of groups.

According to the theorem 3.7.3,

$$(3.7.9) \quad B \times C = B \oplus C$$

The theorem follows from the equality (3.7.9). □

**Proof.** The map

$$(3.7.10) \quad f : (x, y) \in B \times C \to x + y \in A$$

is homomorphism of groups. This statement follows from the equality

$$f(x_1 + x_2, y_1 + y_2) = x_1 + x_2 + y_1 + y_2$$
$$= x_1 + y_1 + x_2 + y_2$$
$$= f(x_1, y_1) + f(x_2, y_2)$$

Let $f(x, y) = 0$. From the equality (3.7.10), it follows that

$$(3.7.11) \quad x + y = 0 \quad x \in B \quad y \in C$$

From the equality (3.7.11), it follows that $x = -y \in C$. Therefore, $x = y = 0$. Therefore, the map $f$ is isomorphism of groups.

According to the theorem 3.7.4,

$$(3.7.12) \quad B \times C = B \oplus C$$

The theorem follows from the equality (3.7.12). □

**Example** 3.7.7. *Let Abelian group $H$ be subgroup of multiplicative Abelian group $G$. We want to find $H \oplus G$. We will use the theorem 3.7.3. Let $A = H \times G$. We define the multiplication by the equality*

$$(a, b)(c, d) = (ac, bd)$$

*Therefore, if $a \in H$, then tuples $(a, 0)$ and $(0, a)$ represent different $(H \oplus G)$-numbers.* □

**Example** 3.7.8. *Let Abelian group $H$ be subgroup of additive Abelian group $G$. We want to find $H \oplus G$. We will use the theorem 3.7.4. Let $A = H \times G$. We define the addition by the equality*

$$(a, b) + (c, d) = (a + c, b + d)$$

*Therefore, if $a \in H$, then tuples $(a, 0)$ and $(0, a)$ represent different $(H \oplus G)$-numbers.* □

---

[3.14] See also remark in [1], page 37.

THEOREM 3.7.9. *Let*

$$f : A \to B$$

*be a surjective homomorphism of Abelian groups and* $C = \ker f$. *If $B$ is free group, then there exists a subgroup $D$ of the group $A$ such that the restriction $f$ to $D$ induces isomorphism $D$ with $B$ and $A = C \oplus D$.*

PROOF.

Let $A$, $B$ be multiplicative Abelian groups. Let $(b_i, i \in I)$ be basis in group $B$. For each $i \in I$, let $a_i$ be $A$-number such that

$$(3.7.13) \qquad f(a_i) = b_i$$

Let $D$ be the subgroup of $A$ generated by $A$-numbers $a_i$, $i \in I$. Let

$$(3.7.14) \qquad \prod_{i \in I} a_i^{n^i} = e$$

with integers $n^i$ where the set $\{i \in I : n^i \neq 0\}$ is finite. The equality

$$(3.7.15) \qquad e = \prod_{i \in I} f(a_i)^{n^i} = \prod_{i \in I} b_i^{n^i}$$

follows from equalities (3.7.13), (3.7.14) and theorems 3.3.1, 3.5.5. Since $(b_i, i \in I)$ is the basis, then according to definitions 3.5.9, 3.5.10, the statement $n_i = 0$, $i \in I$, follows from the equality (3.7.15). Therefore, the set $(a_i, i \in I)$ is the basis of the subgroup $D$.

Let $a \in D$ and $f(a) = e$. Since the set $(a_i, i \in I)$ is the basis of the subgroup $D$, then there exist integers $n^i$ such that

$$(3.7.16) \qquad a = a_i^{n^i}$$

and the set $\{i \in I : n^i \neq 0\}$ is finite. Therefore, the equality

$$(3.7.17) \qquad e = f(a) = f(a_i)^{n^i} = b_i^{n^i}$$

follows from the equality (3.7.16) and theorems 3.3.1, 3.5.5. Since the set $(b_i, i \in I)$ is the basis of the group $B$, $n^i = 0$, $i \in I$, follows from the equality (3.7.17). Therefore, $a = e$ and $D \cap C = \{e\}$.

Let $a \in A$. Since $f(x) \in B$ and the set $(b_i, i \in I)$ is the basis of the group $B$, then there exist integers $n^i$ such that

$$(3.7.18) \qquad f(a) = b_i^{n^i}$$

Let $A$, $B$ be additive Abelian groups. Let $(b_i, i \in I)$ be basis in group $B$. For each $i \in I$, let $a_i$ be $A$-number such that

$$(3.7.21) \qquad f(a_i) = b_i$$

Let $D$ be the subgroup of $A$ generated by $A$-numbers $a_i$, $i \in I$. Let

$$(3.7.22) \qquad \sum_{i \in I} n^i a_i = 0$$

with integers $n^i$ where the set $\{i \in I : n^i \neq 0\}$ is finite. The equality

$$(3.7.23) \qquad 0 = \sum_{i \in I} n^i f(a_i) = \sum_{i \in I} n^i b_i$$

follows from equalities (3.7.21), (3.7.22) and theorems 3.3.2, 3.4.5. Since $(b_i, i \in I)$ is the basis, then according to definitions 3.4.9, 3.4.10, the statement $n_i = 0$, $i \in I$, follows from the equality (3.7.23). Therefore, the set $(a_i, i \in I)$ is the basis of the subgroup $D$.

Let $a \in D$ and $f(a) = 0$. Since the set $(a_i, i \in I)$ is the basis of the subgroup $D$, then there exist integers $n^i$ such that

$$(3.7.24) \qquad a = n^i a_i$$

and the set $\{i \in I : n^i \neq 0\}$ is finite. Therefore, the equality

$$(3.7.25) \qquad 0 = f(a) = n^i f(a_i) = n^i b_i$$

follows from the equality (3.7.24) and theorems 3.3.2, 3.4.5. Since the set $(b_i, i \in I)$ is the basis of the group $B$, $n^i = 0$, $i \in I$, follows from the equality (3.7.25). Therefore, $a = 0$ and $D \cap C = \{0\}$.

Let $a \in A$. Since $f(x) \in B$ and the set $(b_i, i \in I)$ is the basis of the group $B$, then there exist integers $n^i$ such that

$$(3.7.26) \qquad f(a) = n^i b_i$$

and the set $\{i \in I : n^i \neq 0\}$ is finite. The equality

$$(3.7.19) \qquad \begin{aligned} f(xa_i^{-n^i}) &= f(x)f(a_i)^{-n^i} \\ &= f(x)b_i^{-n^i} = e \end{aligned}$$

follows from equalities (3.7.13), (3.7.18). From the equality (3.7.19), it follows that

$$(3.7.20) \qquad b = aa_i^{-n^i} \in C = \ker f$$

From the statement (3.7.20), it follows that $a \in C + D$. According to the theorem 3.7.5, $A = C \oplus D$.

and the set $\{i \in I : n^i \neq 0\}$ is finite. The equality

$$(3.7.27) \qquad \begin{aligned} f(x - n^i a_i) &= f(x) - n^i f(a_i) \\ &= f(x) - n^i b_i = 0 \end{aligned}$$

follows from equalities (3.7.21), (3.7.26). From the equality (3.7.27), it follows that

$$(3.7.28) \qquad b = a - n^i a_i \in C = \ker f$$

From the statement (3.7.28), it follows that $a \in C + D$. According to the theorem 3.7.6, $A = C \oplus D$.

$\square$

**Theorem 3.7.10.** *Let $A$ be subgroup of Abelian group $B$. If group $B/A$ is free, then*

$$(3.7.29) \qquad\qquad B = A \oplus B/A$$

**Proof.**

Let $B$ be multiplicative Abelian group. According to examples 3.3.15, 3.3.17, the sequence of homomorphisms

$$e \longrightarrow A \xrightarrow{\ f\ } B \xrightarrow{\ g\ } B/A \longrightarrow e$$

is exact. According to the example 3.3.15, the homomorphism $f$ is monomorphism. According to the definition 3.3.14, $A = \ker g$. According to the example 3.3.15, the homomorphism $g$ is epimorphism. According to the theorem 3.7.9, there exists subgroup $D$ of the group $B$ which is isomorphic to group $B/A$ and such that

$$(3.7.30) \qquad\qquad B = A \oplus D$$

If we identify $D$-number $d$ and the set $dA$, then the equality (3.7.29) follows from the equality (3.7.30).

Let $B$ be additive Abelian group. According to examples 3.3.18, 3.3.20, the sequence of homomorphisms

$$0 \longrightarrow A \xrightarrow{\ f\ } B \xrightarrow{\ g\ } B/A \longrightarrow 0$$

is exact. According to the example 3.3.18, the homomorphism $f$ is monomorphism. According to the definition 3.3.14, $A = \ker g$. According to the example 3.3.18, the homomorphism $g$ is epimorphism. According to the theorem 3.7.9, there exists subgroup $D$ of the group $B$ which is isomorphic to group $B/A$ and such that

$$(3.7.31) \qquad\qquad B = A \oplus D$$

If we identify $D$-number $d$ and the set $dA$, then the equality (3.7.29) follows from the equality (3.7.31). $\square$

**Theorem 3.7.11.** *Let*

$$f : A \to B$$

*be homomorphism of Abelian groups. If* $\operatorname{Im} f$ *is free group, then*

(3.7.32) $$A = \ker f \oplus \operatorname{Im} f$$

PROOF. Since $A/\ker f = \operatorname{Im} f$, the equality (3.7.32) follows from the equality (3.7.29). $\square$

## 3.8. Unital Extension of Ring

Let $D$ be ring.

LEMMA 3.8.1. *Let* $n \in Z$. *The map*

(3.8.1) $$n : d \in D \to nd \in D$$

*is endomorphism of Abelian group* $D$.

PROOF. The statement

$$nd \in D$$

follows from   the definition   3.4.1.      According to   the theorem 3.4.2 and the definition 2.3.1,   the map (3.8.1) is endomorphism of Abelian group $D$. $\square$

LEMMA 3.8.2. *Let* $a \in D$. *The map*

(3.8.2) $$a : d \in D \to ad \in D$$

*is endomorphism of Abelian group* $D$.

PROOF. The statement

$$ad \in D$$

follows from   the definition   3.2.1.      According to   the statement 3.2.1.3 and theorems 3.2.3, 3.3.2,   the map (3.8.2) is endomorphism of Abelian group $D$. $\square$

LEMMA 3.8.3. *Let* $n \in Z$, $a \in D$. *The map*

(3.8.3) $$a \oplus n : d \in D \to (a \oplus n)d \in D$$

*defined by the equality*

(3.8.4) $$(a \oplus n)d = ad + nd$$

*is endomorphism of Abelian group* $D$.

PROOF. Since $D$ is Abelian group, then the statement

$$nd + ad \in D \quad n \in Z \quad a \in D$$

follows from lemmas 3.8.1, 3.8.2. Therefore, for any $Z$-number $n$ and for any $D$-number $a$, we defined the map (3.8.3). According to lemmas 3.8.1, 3.8.2 and the theorem 3.3.4, the map (3.8.3) is endomorphism of Abelian group $V$. $\square$

THEOREM 3.8.4. *The set of maps* $a \oplus n$ *generates* [3.15] *ring* $D_{(1)}$ *where the sum is defined by the equality*

(3.8.5) $$(a_1 \oplus n_1) + (a_2 \oplus n_2) = (a_1 + a_2) \oplus (n_1 + n_2)$$

*and the product is defined by the equality*

(3.8.6)          $(a_1 \oplus n_1) \circ (a_2 \oplus n_2) = (a_1 a_2 + n_2 a_1 + n_1 a_2) \oplus (n_1 n_2)$

3.8.4.1: *If ring $D$ has unit, then $Z \subseteq D$, $D_{(1)} = D$.*
3.8.4.2: *If ring $D$ is ideal of $Z$, then $D \subseteq Z$, $D_{(1)} = Z$.*
3.8.4.3: *Otherwise $D_{(1)} = D \oplus Z$.*
3.8.4.4: *The ring $D$ is ideal of ring $D_{(1)}$.*
*The ring $D_{(1)}$ is called* **unital extension** *of the ring $D$.*

PROOF.     According to the theorem 3.3.4, the equality

$$((a_1 \oplus n_1) + (a_2 \oplus n_2))d$$
$$= (a_1 \oplus n_1)d + (a_2 \oplus n_2)d$$
$$= a_1 d + n_1 d + a_2 d + n_2 d$$
(3.8.7)
$$= a_1 d + a_2 d + n_1 d + n_2 d$$
$$= (a_1 + a_2)d + (n_1 + n_2)d$$
$$= (a_1 + a_2) \oplus (n_1 + n_2)d$$

follows from the equality (3.8.4). The equality (3.8.5) follows from the equality (3.8.7).

According to the theorem 2.1.10, the equality

$$((a_1 \oplus n_1) \circ (a_2 \oplus n_2))d$$
$$= (a_1 \oplus n_1)((a_2 \oplus n_2)d)$$
$$= (a_1 \oplus n_1)(a_2 d + n_2 d)$$
(3.8.8)
$$= (a_1 \oplus n_1)(a_2 d)$$
$$+ (a_1 \oplus n_1)(n_2 d)$$
$$= a_1(a_2 d) + n_1(a_2 d)$$
$$+ a_1(n_2 d) + n_1(n_2 d)$$

follows from the equality (3.8.4). The equality

$$((a_1 \oplus n_1) \circ (a_2 \oplus n_2))d$$
(3.8.9)
$$= (a_1 a_2)d + (n_1 a_2)d$$
$$+ (n_2 a_1)d + (n_1 n_2)d$$

follows from    equalities (3.4.5),    (3.8.8),    from statements 3.1.1.1, 3.2.1.2, from the theorem 3.4.5,   from the lemma 3.8.2. The equality

$$((a_1 \oplus n_1) \circ (a_2 \oplus n_2))d$$
$$= (a_1 a_2 + n_1 a_2 + n_2 a_1 + n_1 n_2)d$$
$$= ((a_1 a_2 + n_1 a_2 + n_2 a_1) \oplus n_1 n_2)d$$

and therefore equality (3.8.6), follows from equalities (3.8.4), (3.8.9) and from lemmas 3.8.1, 3.8.2.

Let ring $D$ have unit $e$. Then we can identify $n \in Z$ and $ne \in D$. The equality

(3.8.10)               $(a \oplus n)d = ad + nd = (a + ne)d$

---

3.15   See the definition of unital extension also on the pages [7]-52, [8]-64.

follows from the equality

$$(3.8.4) \quad (a \oplus n)d = ad + nd$$

The statement 3.8.4.1 follows from the equality (3.8.10).

Let ring $D$ be ideal of $Z$. Then $D \subseteq Z$. The equality

$$(3.8.11) \quad\quad (a \oplus n)d = ad + nd = (a + n)d$$

follows from the equality (3.8.4). The statement 3.8.4.2 follows from the equality (3.8.11).

The statement 3.8.4.4 follows from the lemma 3.8.3.                    $\square$

# $D$-Module

## 4.1. Module over Commutative Ring

According to the theorem 3.3.4, considering the representation

$$f : D \longrightarrow_* V$$

of the ring $D$ in the Abelian group $V$, we assume

$$(4.1.1) \qquad (f(a) + f(b))(x) = f(a)(x) + f(b)(x)$$

According to the definition 2.3.1, the map $f$ is homomorphism of ring $D$. Therefore

$$(4.1.2) \qquad f(a + b) = f(a) + f(b)$$

The equalty

$$(4.1.3) \qquad f(a + b)(x) = f(a)(x) + f(b)(x)$$

follows from equalities (4.1.1), (4.1.2).

**Theorem 4.1.1.** *Representation*

$$f : D \longrightarrow_* V$$

*of the ring $D$ in an Abelian group $V$ satisfies the equality*

$$(4.1.4) \qquad f(0) = \overline{0}$$

*where*

$$\overline{0} : V \to V$$

*is map such that*

$$(4.1.5) \qquad \overline{0} \circ v = 0$$

Proof. The equality

$$(4.1.6) \qquad f(a)(x) = f(a + 0)(x) = f(a)(x) + f(0)(x)$$

follows from the equality (4.1.3). The equality

$$(4.1.7) \qquad f(0)(x) = 0$$

follows from the equality (4.1.6). The equality (4.1.5) follows from the equalities (4.1.4), (4.1.7). $\qquad \square$

**Theorem 4.1.2.** *Representation*

$$f : D \longrightarrow_* V$$

*of the ring $D$ in an Abelian group $V$ is* **effective** *iff $a = 0$ follows from equation $f(a) = 0$.*

Proof. Suppose $a$, $b \in D$ cause the same transformation. Then

$$(4.1.8) \qquad f(a)(m) = f(b)(m)$$

for any $m \in V$. From equalities (4.1.3), (4.1.8), it follows that

$$(4.1.9) \qquad f(a - b)(m) = 0$$

The equality

$$f(a - b) = \overline{0}$$

follows from equalities (4.1.5), (4.1.9). Therefore, the representation $f$ is effective iff $a = b$.  □

DEFINITION 4.1.3. *Let $D$ be commutative ring. Effective representation*

$$(4.1.10) \qquad f : D \overset{*}{\longrightarrow} V \qquad f(d) : v \to d\, v$$

*of ring $D$ in Abelian group $V$ is called* **module** *over ring $D$. We will also say that Abelian group $V$ is $D$-**module**. $V$-number is called* **vector***. Bilinear map*

$$(4.1.11) \qquad (a, v) \in D \times V \to av \in V$$

*generated by representation*

$$(4.1.12) \qquad (a, v) \to av$$

*is called product of vector over scalar.*  □

DEFINITION 4.1.4. *Let $D$ be field. Effective representation*

$$(4.1.13) \qquad f : D \overset{*}{\longrightarrow} V \qquad f(d) : v \to d\, v$$

*of field $D$ in Abelian group $V$ is called* **vector space** *over field $D$. We will also say that Abelian group $V$ is $D$-**vector space**. $V$-number is called* **vector***. Bilinear map*

$$(4.1.14) \qquad (a, v) \in D \times V \to av \in V$$

*generated by representation*

$$(4.1.15) \qquad (a, v) \to av$$

*is called product of vector over scalar.*  □

THEOREM 4.1.5. *The following diagram of representations describes $D$-module $V$*

$$(4.1.16)$$

*The diagram of representations (4.1.16) holds* **commutativity of representations** *of ring of rational integers $Z$ and commutative ring $D$ in Abelian group $V$*

$$(4.1.17) \qquad a(nv) = n(av)$$

PROOF. The diagram of representations (4.1.16) follows from the definition 4.1.3 and from the theorem 3.4.2. Since transformation $g_2(a)$ is endomorphism of $Z$-module $V$, we obtain the equality (4.1.17).  □

DEFINITION 4.1.6. *Subrepresentation of $D$-module $V$ is called* **submodule** *of $D$-module $V$.*  □

**Theorem 4.1.7.** *Let $V$ be $D$-module. For any $D_{(1)}$-number the map*

$$a \oplus n : v \in V \to (a \oplus n)v \in V \tag{4.1.18}$$

$$(a \oplus n)v = av + nv \quad a \oplus n \in D_{(1)} \tag{4.1.19}$$

*is endomorphism of Abelian group $V$. The set of transormations (3.8.3) is representation of ring $D_{(1)}$ in Abelian group $V$. We use the notation $D_{(1)}v$ for the set of vectors generated by vector $v$.*

Proof. According to the theorem 4.1.5, expressions $av$, $nv$ are $V$-numbers. Since $V$ is Abelian group, expression (4.1.19) is also $V$-number. $\square$

**Theorem 4.1.8.** *Let $V$ be $D$-module. Following conditions hold for $V$-numbers:*

4.1.8.1: **commutative law**

$$v + w = w + v \tag{4.1.20}$$

4.1.8.2: **associative law**

$$(pq)v = p(qv) \tag{4.1.21}$$

4.1.8.3: **distributive law**

$$p(v + w) = pv + pw \tag{4.1.22}$$

$$(p + q)v = pv + qv \tag{4.1.23}$$

4.1.8.4: **unitarity law**

$$(0 \oplus 1)v = v \tag{4.1.24}$$

*If ring has unit, then the equality (4.1.24) has form*

$$1v = v \tag{4.1.25}$$

*for any $p$, $q \in D_{(1)}$, $v$, $w \in V$.*

Proof. The equality (4.1.20) follows from the definition 4.1.3.
The equality (4.1.22) follows from theorems 3.3.2, 4.1.7.
The equality (4.1.23) follows from the theorem 3.3.4.
The equality (4.1.21) follows from the theorem 2.1.10.
According to theorems 4.1.5, 4.1.7, the equality (4.1.24) follows from the equality

$$(0 \oplus 1)v = 0v + 1v = v$$

and from the equality (3.4.4). If ring has unit, then the equality (4.1.25) follows from the theorem 3.4.12. $\square$

**Theorem 4.1.9.** *Let $V$ be $D$-module. For any set of $V$-numbers*

$$v = (v(i) \in V, i \in I) \tag{4.1.26}$$

*vector generated by the diagram of representations (4.1.16) has the following form*[4.1]

$$J(v) = \{w : w = c(i)v(i), c(i) \in D_{(1)}, |\{i : c(i) \neq 0\}| < \infty\} \tag{4.1.27}$$

*where ring $D_{(1)}$ is unital extension of the ring $D$.*

PROOF. We prove the theorem by induction based on the theorem 2.7.5.
For any $v(k) \in J(v),$ let

$$c(i) = \begin{cases} 0 \oplus 1, & i = k \\ 0 \end{cases}$$

Then the equality

$$(4.1.28) \qquad v(k) = \sum_{i \in I} c(i)v(i) \quad c(i) \in D_{(1)}$$

follows from the equality (4.1.19). From equalities (4.1.27), (4.1.28), it follows that
the theorem is true on the set $X_0 = v \subseteq J(v)$.

Let $X_{k-1} \subseteq J(v)$. According to the definition 4.1.3 and theorems 2.7.5,
4.1.7, if $w \in X_k$, then or $w = w_1 + w_2$, $w_1$, $w_2 \in X_{k-1}$, either $w = aw_1$, $a \in D_{(1)}$,
$w_1 \in X_{k-1}$.

LEMMA 4.1.10. *Let $w = w_1 + w_2$, $w_1$, $w_2 \in X_{k-1}$. Then $w \in J(v)$.*

PROOF. According to the equality (4.1.27), there exist $D_{(1)}$-numbers
$w_1(i)$, $w_2(i)$, $i \in I$, such that

$$(4.1.29) \qquad w_1 = \sum_{i \in I} w_1(i)v(i) \quad w_2 = \sum_{i \in I} w_2(i)v(i)$$

where sets

$$(4.1.30) \qquad H_1 = \{i \in I : w_1(i) \neq 0\} \quad H_2 = \{i \in I : w_2(i) \neq 0\}$$

are finite. Since $V$ is $D$-module, then from equalities (4.1.23), (4.1.29), it
follows that

$$(4.1.31) \qquad \begin{aligned} w_1 + w_2 &= \sum_{i \in I} w_1(i)v(i) + \sum_{i \in I} w_2(i)v(i) = \sum_{i \in I}(w_1(i)v(i) + w_2(i)v(i)) \\ &= \sum_{i \in I}(w_1(i) + w_2(i))v(i) \end{aligned}$$

From the equality (4.1.30), it follows that the set

$$\{i \in I : w_1(i) + w_2(i) \neq 0\} \subseteq H_1 \cup H_2$$

is finite.                                                                      ⊙

LEMMA 4.1.11. *Let $w = aw_1$, $a \in D_{(1)}$, $w_1 \in X_{k-1}$. Then $w \in J(v)$.*

PROOF. According to the statement 2.7.5.4, for any $D_{(1)}$-number $a$,

$$(4.1.32) \qquad aw_1 \in X_k$$

According to the equality (4.1.27), there exist $D_{(1)}$-numbers $w(i)$, $i \in I$,
such that

$$(4.1.33) \qquad w_1 = \sum_{i \in I} w_1(i)v(i)$$

---

4.1 For a set $A$, we denote by $|A|$ the cardinal number of the set $A$. The notation $|A| < \infty$ means
that the set $A$ is finite.

where

(4.1.34)                    $|\{i \in I : w_1(i) \neq 0\}| < \infty$

From the equality (4.1.33) it follows that

(4.1.35)        $aw_1 = a\sum_{i \in I} w_1(i)v(i) = \sum_{i \in I} a(w_1(i)v(i)) = \sum_{i \in I}(aw_1(i))v(i)$

From the statement (4.1.34), it follows that the set $\{i \in I : aw(i) \neq 0\}$ is finite.                                                                          $\odot$

From lemmas 4.1.10, 4.1.11, it follows that $X_k \subseteq J(v)$.        $\square$

Convention 4.1.12. *We will use summation convention in which repeated index in linear combination implies summation with respect to repeated index. In this case we assume that we know the set of summation index and do not use summation symbol*

$$c(i)v(i) = \sum_{i \in I} c(i)v(i)$$

*If needed to clearly show a set of indices, I will do it.*        $\square$

Definition 4.1.13. *Let $V$ be module. Let*

$$v = (v(i) \in V, i \in I)$$

*be set of vectors. The expression $c(i)v(i)$, $c(i) \in D_{(1)}$, is called* **linear combination** *of vectors $v(i)$. A vector $w = c(i)v(i)$ is called* **linearly dependent** *on vectors $v(i)$.*        $\square$

Theorem 4.1.14. *Let $v = (v_i \in V, i \in I)$ be set of vectors of $D$-module $V$. If vectors $v_i$, $i \in I$, belongs submodule $V'$ of $D$-module $V$, then linear combination of vectors $v_i$, $i \in I$, belongs submodule $V'$.*

Proof. The theorem follows from the theorem 4.1.9 and definitions 4.1.13, 4.1.6.        $\square$

Theorem 4.1.15. *Let $D$ be field. Since the equation*

(4.1.36)                        $w(i)v(i) = 0$

*implies existence of index $i = j$ such that $w(j) \neq 0$, then the vector $v(j)$ linearly depends on rest of vectors $v$.*

Proof. The theorem follows from the equality

$$v(j) = \sum_{i \in I \setminus \{j\}} \frac{w(i)}{w(j)} v(i)$$

and from the definition 4.1.13.                                                  $\square$

It is evident that for any set of vectors $v(i)$

$$w(i) = 0 \Rightarrow w(i)v(i) = 0$$

DEFINITION 4.1.16. *The set of vectors* $^{4.2}v(i)$, $i \in I$, *of D-module V is* **linearly independent** *if* $c(i) = 0$, $i \in I$, *follows from the equation*

$$(4.1.37) \qquad\qquad c(i)v(i) = 0$$

*Otherwise the set of vectors* $v(i)$, $i \in I$, *is* **linearly dependent**.  $\square$

The following definition follows from the theorems 2.7.5, 4.1.9 and from the definition 2.7.4.

DEFINITION 4.1.17. $J(v)$ *is called* **submodule** *generated by set v, and v is a* **generating set** *of submodule* $J(v)$. *In particular, a* **generating set** *of D-module V is a subset* $X \subset V$ *such that* $J(X) = V$.  $\square$

The following definition follows from the theorems 2.7.5, 4.1.9 and from the definition 2.7.14.

DEFINITION 4.1.18. *If the set* $X \subset V$ *is generating set of D-module V, then any set Y, $X \subset Y \subset V$ also is generating set of D-module V. If there exists minimal set X generating the D-module V, then the set X is called* **quasi-basis** *of D-module V.*  $\square$

DEFINITION 4.1.19. *Let* $\overline{\overline{e}}$ *be the quasibasis of D-module V and vector* $\overline{v} \in V$ *has expansion*

$$(4.1.38) \qquad\qquad \overline{v} = v(i)e(i)$$

*with respect to the quasibasis* $\overline{\overline{e}}$. $D_{(1)}$*-numbers* $v(i)$ *are called* **coordinates** *of vector* $\overline{v}$ *with respect to the quasibasis* $\overline{\overline{e}}$. *Matrix of* $D_{(1)}$*-numbers* $v = (v(i), i \in I)$ *is called* **coordinate matrix of vector** $\overline{v}$ *in quasibasis* $\overline{\overline{e}}$.  $\square$

THEOREM 4.1.20. *The set of vectors* $\overline{\overline{e}} = (e(i), i \in I)$ *is quasi-basis of D-module V, if following statements are true.*

4.1.20.1: *Arbitrary vector* $v \in V$ *is linear combination of vectors of the set* $\overline{\overline{e}}$.

4.1.20.2: *Vector* $e(i)$ *cannot be represented as a linear combination of the remaining vectors of the set* $\overline{\overline{e}}$.

PROOF. According to the statement 4.1.20.1, the theorem 4.1.9 and the definition 4.1.13, the set $\overline{\overline{e}}$ generates *D*-module *V* (the definition 4.1.17). According to the statement 4.1.20.2, the set $\overline{\overline{e}}$ is minimal set generating *D*-module *V*. According to the definitions 4.1.18, the set $\overline{\overline{e}}$ is a quasi-basis of *D*-module *V*.  $\square$

THEOREM 4.1.21. *Let D be commutative ring. Let* $\overline{\overline{e}}$ *be quasi-basis of D-module V. Let*

$$(4.1.39) \qquad\qquad c(i)e(i) = 0$$

*be linear dependence of vectors of the quasi-basis* $\overline{\overline{e}}$. *Then*

---

[4.2] I follow to the definition on page [1]-130.

4.1.21.1: $D_{(1)}$-number $c(i)$, $i \in I$, does not have inverse element in ring $D_{(1)}$.
4.1.21.2: The set $V'$ of matrices $c = (c(i), i \in I)$ generates $D$-module $V'$.

PROOF. Let there exist matrix $c = (c(i), i \in I)$ such that the equality (4.1.39) is true and there exist index $i = j$ such that $c(j) \neq 0$. If we assume that $D$-number $c(j)$ has inverse one, then the equality

$$e(j) = \sum_{i \in I \setminus \{j\}} \frac{c(i)}{c(j)} e(i)$$

follows from the equality (4.1.39). Therefore, the vector $e(j)$ is linear combination of other vectors of the set $\overline{\overline{e}}$ and the set $\overline{\overline{e}}$ is not quasi-basis. Therefore, our assumption is false, and $D_{(1)}$-number $c(j)$ does not have inverse.

Let matrices $b = (b(i), i \in I) \in D'$, $c = (c(i), i \in I) \in D'$. From equalities

$$b(i)e(i) = 0$$

$$c(i)e(i) = 0$$

it follows that

$$(b(i) + c(i))e(i) = 0$$

Therefore, the set $D'$ is Abelian group.

Let matrix $c = (c(i), i \in I) \in D'$ and $a \in D$. From the equality

$$c(i)e(i) = 0$$

it follows that

$$(ac(i))e(i) = 0$$

Therefore, Abelian group $D'$ is $D$-module. □

THEOREM 4.1.22. Let $D$-module $V$ have the quasi-basis $\overline{\overline{e}}$ such that in the equality

(4.1.40) $$c(i)e(i) = 0$$

there exists index $i = j$ such that $c(j) \neq 0$. Then
4.1.22.1: The matrix $c = (c(i), i \in I)$ determines coordinates of vector $0 \in V$ with respect to quasi-basis $\overline{\overline{e}}$.
4.1.22.2: Coordinates of vector $\overline{v}$ with respect to quasi-basis $\overline{\overline{e}}$ are uniquely determined up to a choice of coordinates of vector $0 \in V$.

PROOF. The statement 4.1.22.1 follows from the equality (4.1.40) and from the definition 4.1.19.

Let vector $\overline{v}$ have expansion

(4.1.41) $$\overline{v} = v(i)e(i)$$

with respect to quasi-basis $\overline{\overline{e}}$. The equality

(4.1.42) $$\overline{v} = \overline{v} + 0 = v(i)e(i) + c(i)e(i) = (v(i) + c(i))e(i)$$

follows from equalities (4.1.40), (4.1.41). The statement 4.1.22.2 follows from equalities (4.1.41), (4.1.42) and from the definition 4.1.19. □

THEOREM 4.1.23. *Let the set of vectors of quasi-basis $\overline{\overline{e}}$ of D-module $V$ be linear independent. Then the quasi-basis $\overline{\overline{e}}$ is basis of D-module $V$.*

PROOF. The theorem follows from the definition 2.7.14 and the theorem 4.1.22. □

DEFINITION 4.1.24. *The D-module $V$ is* **free *D*-module**,[4.3]*if D-module $V$ has basis.* □

THEOREM 4.1.25. *The D-vector space is free D-module.*

PROOF. Let the set of vectors $e(i)$, $i \in I$, be linear dependent. Then the equation

$$w(i)e(i) = 0$$

implies existence of index $i = j$ such that $w(j) \neq 0$. According to the theorem 4.1.15, the vector $e(j)$ linearly depends on rest of vectors of the set $\overline{\overline{e}}$. According to the definition 4.1.18, the set of vectors $e(i)$, $i \in I$, is not a basis for *D*-vector space $V$.

Therefore, if the set of vectors $e(i)$, $i \in I$, is a basis, then these vectors are linearly independent. Since an arbitrary vector $v \in V$ is linear combination of vectors $e(i)$, $i \in I$, then the set of vectors $v$, $e(i)$, $i \in I$, is not linearly independent. □

THEOREM 4.1.26. *Coordinates of vector $v \in V$ relative to basis $\overline{\overline{e}}$ of free D-module $V$ are uniquely defined. The equality*

(4.1.43)                    $ve = we \qquad v(i)e(i) = w(i)e(i)$

*implies the equality*

$$v = w \qquad v(i) = w(i)$$

PROOF. The theorem follows from the theorem 4.1.22 and from definitions 4.1.16, 4.1.24. □

EXAMPLE 4.1.27. *From theorems 3.4.2, 3.6.10 and the definition 4.1.3, it follows that free Abelian group $G$ is module over ring of integers $Z$.* □

## 4.2. *D*-Module Type

**4.2.1. Concept of Generalized Index.** The only advantage of expression $c(i)v(i)$, $c(i) \in D_{(1)}$, of linear combination of vectors

$$v = (v(i) \in V, i \in I)$$

is its versatility. However to see vector expression more clearly, it is convenient to gather the coordinates of a vector in table called matrix. Therefore, instead of argument written in brackets, we will use indices which say how coordinates of vector are organized in matrix.

---

[4.3]  I follow to the definition in [1], page 135.

Organization of coordinates of vector in matrix is called *D*-module type.

In this section, we consider column vector and row vector. it is evident that there exist other forms of representation of vector. For instance, we can represent coordinates of vector as $n \times m$ matrix or as triangular matrix. Format of representation depends on considered problem.

If a statement depends on the format of representation of vector, we will specify either type of *D*-module, or format of representation of basis. Both ways of specifying the type of *D*-module are equivalent.

If we consider expression $i$ in the notation $v(i)$ of coordinates of vector as generalized index, then we see that behavior of vector does not depend on *D*-module type. We can extend the definition of generalized index and consider index $i$ in expression $a^i$ as generalized index. We will use the symbol $\cdot$ in front of a generalized index when we need to describe its structure. we put the sign $'-'$ in place of the index whose position was changed. For instance, if an original term was $a_{ij}$ I will use notation $a_{i\cdot-}^{\phantom{i\cdot}j}$ instead of notation $a_i^j$. This allows us to write the equality

$$(4.2.1) \qquad b'^i = b'^{\cdot-}_{\phantom{.}k} = f^{\cdot-\phantom{.}l}_{\phantom{.}k\cdot-}b'^{\cdot-}_{\phantom{.}l} = f^i_j b^j.$$

The equality (4.2.1) assumes the equality of indices $^i = {^{\cdot-}_{\phantom{.}k}}$.

Even though the structure of a generalized index is arbitrary we assume that there exists a one-to-one map of the interval of positive integers $1, ..., n$ to the range of index. Let $I$ be the range of the index $i$. We denote the power of this set by symbol $|I|$ and assume that $|I| = n$. If we want to enumerate elements $a_i$ we use notation $a_1, ..., a_n$.

In geometry, there is a tradition to write coordinates of vector with index at the top

$$v(i) = v^i$$

and to enumerate the vectors of basis using index at the bottom

$$e(i) = e_i$$

We also use the convention 4.2.1.

Convention 4.2.1. *We will use Einstein summation convention in which repeated index (one above and one below) implies summation with respect to repeated index. In this case we assume that we know the set of summation index and do not use summation symbol*

$$c^i v_i = \sum_{i \in I} c^i v_i$$

*If needed to clearly show a set of indices, I will do it.* □

If we write the system of linear equations

$$a_1^1 x^1 + ... + a_n^1 x^n = 0$$
$$(4.2.2) \qquad\qquad ......$$
$$a_1^n x^1 + ... + a_n^n x^n = 0$$

using matrices

$$(4.2.3) \qquad \begin{pmatrix} a_1^1 & ... & a_n^1 \\ ... & ... & ... \\ a_1^n & ... & a_n^n \end{pmatrix} \begin{pmatrix} a^i \\ ... \\ a^n \end{pmatrix} = \begin{pmatrix} 0 \\ ... \\ 0 \end{pmatrix}$$

then we can identify the solution of the system of linear equations (4.2.2) and column of matrix. At the same time, the set of solutions of the system of linear equations (4.2.2) is *D*-module.

**4.2.2. Notation for Matrix Entry.** Representation of coordinates of a vector as a matrix allows making a notation more compact. The question of the presentation of vector as a row or a column of the matrix is just a question of convention. We extend the concept of generalized index to entries of the matrix. A matrix is a two dimensional table, the rows and columns of which are enumerated by generalized indeices. Since we use generalized index, we cannot tell whether index *a* of matrix enumerates rows or columns until we know the structure of index. However as we can see below the form of presentation of matrix is not important for us. To make sure that notation offered below is consistent with the traditional we will assume that the matrix is presented in the form

$$a = \begin{pmatrix} a_1^1 & ... & a_n^1 \\ ... & ... & ... \\ a_1^m & ... & a_n^m \end{pmatrix}$$

The upper index enumerates rows and the lower index enumerates columns.

DEFINITION 4.2.2. *I use the following names and notation for different* **submatrices** *of the matrix* $a$

- **column of matrix** *with the index* $i$

$$a_i = \begin{pmatrix} a_i^1 \\ ... \\ a_i^m \end{pmatrix}$$

  *The upper index enumerates entries of the column and the lower index enumerates columns.*

  $a_T$ : *the submatrix obtained from the matrix* $a$ *by selecting columns with an index from the set* $T$

  $a_{[i]}$ : *the submatrix obtained from the matrix* $a$ *by deleting column* $a_i$

  $a_{[T]}$ : *the submatrix obtained from the matrix* $a$ *by deleting columns with an index from the set* $T$

- **row of matrix** *with the index* $j$

$$a^j = \begin{pmatrix} a_1^j & ... & a_n^j \end{pmatrix}$$

  *The lower index enumerates entries of the row and the upper index enumerates rows.*

$a^S$    : *the submatrix obtained from the matrix a by selecting rows with an index from the set S*

$a^{[j]}$   : *the submatrix obtained from the matrix a by deleting row $a^j$*

$a^{[S]}$   : *the submatrix obtained from the matrix a by deleting rows with an index from the set S*

$\square$

REMARK 4.2.3. *We will combine the notation of indices. Thus $a^j_i$ is $1 \times 1$ submatrix. The same time this is the notation for a matrix entry. This allows an identifying of $1 \times 1$ matrix and its entry. The index a is number of column of matrix and the index b is number of row of matrix.* $\square$

Let $I$, $|I| = n$, be a set of indices. We introduce the **Kronecker symbol**

$$\delta^i_j = \begin{cases} 1 & i = j \\ 0 & i \neq j \end{cases} \quad i, j \in I$$

Let $k$, $l$, $k \leq l$, be integers. We will consider the set

$$[k, l] = \{ i : k \leq i \leq l \}$$

### 4.2.3. *D*-Module of Columns.

DEFINITION 4.2.4. *We represented the set of vectors $v(1) = v_1, ..., v(m) = v_m$ as row of matrix*

$$(4.2.4) \qquad v = \begin{pmatrix} v_1 & ... & v_m \end{pmatrix}$$

*and the set of $D_{(1)}$-nummbers $c(1) = c^1, ..., c(m) = c^m$ as column of matrix*

$$(4.2.5) \qquad c = \begin{pmatrix} c^1 \\ ... \\ c^m \end{pmatrix}$$

*Corresponding representation of D-module V is called D-**module of columns**, and V-number is called **column vector**.* $\square$

THEOREM 4.2.5. *We can represent linear combination*

$$\sum_{i=1}^m c(i)v(i) = c^i v_i$$

*of vectors $v_1, ..., v_m$ as product of matrices*

$$(4.2.6) \qquad cv = \begin{pmatrix} c^1 \\ ... \\ c^m \end{pmatrix} \begin{pmatrix} v_1 & ... & v_m \end{pmatrix} = c^i v_i$$

PROOF. The theorem follows from definitions 4.1.13, 4.2.4. $\square$

THEOREM 4.2.6. *If we write vectors of basis $\overline{\overline{e}}$ as row of matrix*

$$(4.2.7) \qquad\qquad e = \begin{pmatrix} e_1 & \dots & e_n \end{pmatrix}$$

*and coordinates of vector $\overline{w} = w^i e_i$ with respect to basis $\overline{\overline{e}}$ as column of matrix*

$$(4.2.8) \qquad\qquad w = \begin{pmatrix} w^1 \\ \dots \\ w^n \end{pmatrix}$$

*then we can represent the vector $\overline{w}$ as product of matrices*

$$(4.2.9) \qquad \overline{w} = we = \begin{pmatrix} w^1 \\ \dots \\ w^n \end{pmatrix} \begin{pmatrix} e_1 & \dots & e_n \end{pmatrix} = w^i e_i$$

PROOF. The theorem follows from the theorem 4.2.5.  □

THEOREM 4.2.7. *Coordinates of vector $v \in V$ relative to basis $\overline{\overline{e}}$ of free D-module $V$ are uniquely defined. The equality*

$$(4.2.10) \qquad\qquad ve = we \qquad v^i e_i = w^i e_i$$

*implies the equality*

$$v = w \qquad v^i = w^i$$

PROOF. The theorem follows from the theorem 4.1.26.  □

THEOREM 4.2.8. *Let the set of vectors $\overline{\overline{e}}_A = (e_{Ak}, k \in K)$ be a basis of D-algebra of columns $A$. Let $\overline{\overline{e}}_V = (e_{Vi}, i \in I)$ be a basis of $A$-module of columns $V$. Then the set of $V$-numbers*

$$(4.2.11) \qquad\qquad \overline{\overline{e}}_{VA} = (e_{VAki} = e_{Ak}e_{Vi})$$

*is the basis of D-module $V$. The basis $\overline{\overline{e}}_{VA}$ is called extension of the basis $(\overline{\overline{e}}_A, \overline{\overline{e}}_V)$.*

LEMMA 4.2.9. *The set of $V$-numbers $e_{2ik}$ generates D-module $V$.*

PROOF. Let

$$(4.2.12) \qquad\qquad a = a^i e_{1i}$$

be any $V$-number. For any $A$-number $a^i$, according to the definition 5.1.1 and the convention 5.2.1, there exist $D$-numbers $a^{ik}$ such that

$$(4.2.13) \qquad\qquad a^i = a^{ik} e_k$$

The equality

$$(4.2.14) \qquad\qquad a = a^{ik} e_k e_{1i} = a^{ik} e_{2ik}$$

follows from equalities (4.2.11), (4.2.12), (4.2.13). According to the theorem 4.1.9, the lemma follows from the equality (4.2.14).  ⊙

Lemma 4.2.10. *The set of V-numbers $e_{2ik}$ is linearly independent over the ring D.*

Proof. Let

$$(4.2.15) \qquad a^{ik}e_{2ik} = 0$$

The equality

$$(4.2.16) \qquad a^{ik}e_k e_{1i} = 0$$

follows from equalities (4.2.11), (4.2.15). According to the theorem 4.2.7, the equality

$$(4.2.17) \qquad a^{ik}e_k = 0 \quad k \in K$$

follows from the equality (4.2.16). According to the theorem 4.2.7, the equality

$$(4.2.18) \qquad a^{ik} = 0 \quad k \in K \quad i \in I$$

follows from the equality (4.2.17). Therefore, the set of *V*-numbers $e_{2ik}$ is linearly independent over the ring *D*. ⊙

Proof of theorem 4.2.8. The theorem follows from the theorem 4.2.7 and from lemmas 4.2.9, 4.2.10. □

### 4.2.4. *D*-Module of Rows.

Definition 4.2.11. *We represented the set of vectors $v(1) = v^1$, ..., $v(m) = v^m$ as column of matrix*

$$(4.2.19) \qquad v = \begin{pmatrix} v^1 \\ ... \\ v^m \end{pmatrix}$$

*and the set of $D_{(1)}$-nummbers $c(1) = c_1$, ..., $c(m) = c_m$ as row of matrix*

$$(4.2.20) \qquad c = \begin{pmatrix} c_1 & ... & c_m \end{pmatrix}$$

*Corresponding representation of D-module V is called D-**module of rows**, and V-number is called **row vector**.* □

Theorem 4.2.12. *We can represent linear combination*

$$\sum_{i=1}^{m} c(i)v(i) = c_i v^i$$

*of vectors $v^1$, ..., $v^m$ as product of matrices*

$$(4.2.21) \qquad vc = \begin{pmatrix} v^1 \\ ... \\ v^m \end{pmatrix} \begin{pmatrix} c_1 & ... & c_m \end{pmatrix} = c_i v^i$$

Proof. The theorem follows from definitions 4.1.13, 4.2.11. □

THEOREM 4.2.13. *If we write vectors of basis $\overline{\overline{e}}$ as column of matrix*

$$(4.2.22) \qquad e = \begin{pmatrix} e^1 \\ ... \\ e^n \end{pmatrix}$$

*and coordinates of vector $\overline{w} = w_i e^i$ with respect to basis $\overline{\overline{e}}$ as row of matrix*

$$(4.2.23) \qquad w = \begin{pmatrix} w_1 & ... & w_n \end{pmatrix}$$

*then we can represent the vector $\overline{w}$ as product of matrices*

$$(4.2.24) \qquad \overline{w} = ew = \begin{pmatrix} e^1 \\ ... \\ e^n \end{pmatrix} \begin{pmatrix} w_1 & ... & w_n \end{pmatrix} = w_i e^i$$

PROOF. The theorem follows from the theorem 4.2.12. □

THEOREM 4.2.14. *Coordinates of vector $v \in V$ relative to basis $\overline{\overline{e}}$ of free D-module $V$ are uniquely defined. The equality*

$$(4.2.25) \qquad ve = we \qquad v_i e^i = w_i e^i$$

*implies the equality*

$$v = w \qquad v_i = w_i$$

PROOF. The theorem follows from the theorem 4.1.26. □

THEOREM 4.2.15. *Let the set of vectors $\overline{\overline{e}}_A = (e^k_A, k \in K)$ be a basis of D-algebra of rows $A$. Let $\overline{\overline{e}}_V = (e^i_V, i \in I)$ be a basis of $A$-module of rows $V$. Then the set of $V$-numbers*

$$(4.2.26) \qquad \overline{\overline{e}}_{VA} = (e^{ki}_{VA} = e^k_A e^i_V)$$

*is the basis of D-module $V$. The basis $\overline{\overline{e}}_{VA}$ is called extension of the basis $(\overline{\overline{e}}_A, \overline{\overline{e}}_V)$.*

LEMMA 4.2.16. *The set of $V$-numbers $e^{ik}_2$ generates D-module $V$.*

PROOF. Let

$$(4.2.27) \qquad a = a_i e^i_1$$

be any $V$-number. For any $A$-number $a_i$, according to the definition 5.1.1 and the convention 5.2.1, there exist $D$-numbers $a_{ik}$ such that

$$(4.2.28) \qquad a_i = a_{ik} e^k$$

The equality

$$(4.2.29) \qquad a = a_{ik} e^k e^i_1 = a_{ik} e^{ik}_2$$

follows from equalities (4.2.26), (4.2.27), (4.2.28). According to the theorem 4.1.9, the lemma follows from the equality (4.2.29).                                    ⊙

Lemma 4.2.17. *The set of $V$-numbers $e_2^{ik}$ is linearly independent over the ring $D$.*

Proof. Let

$$(4.2.30) \qquad a_{ik}e_2^{ik} = 0$$

The equality

$$(4.2.31) \qquad a_{ik}e^k e_1^i = 0$$

follows from equalities (4.2.26), (4.2.30). According to the theorem 4.2.14, the equality

$$(4.2.32) \qquad a_{ik}e^k = 0 \quad k \in K$$

follows from the equality (4.2.31). According to the theorem 4.2.14, the equality

$$(4.2.33) \qquad a_{ik} = 0 \quad k \in K \quad i \in I$$

follows from the equality (4.2.32). Therefore, the set of $V$-numbers $e_2^{ik}$ is linearly independent over the ring $D$. ⊙

Proof of theorem 4.2.15. The theorem follows from the theorem 4.2.14 and from lemmas 4.2.16, 4.2.17. □

## 4.3. Linear Map of $D$-Module

### 4.3.1. General Definition.

Definition 4.3.1. *Morphism of representations*

$$(4.3.1) \qquad (h : D_1 \to D_2 \quad f : V_1 \to V_2)$$

*of $D_1$-module $V_1$ into $D_2$-module $V_2$ is called homomorphism or* **linear map** *of $D_1$-module $V_1$ into $D_2$-module $V_2$. Let us denote $\mathcal{L}(D_1 \to D_2; V_1 \to V_2)$ set of linear maps of $D_1$-module $V_1$ into $D_2$-module $V_2$.* □

If the map (4.3.1) is linear map of $D_1$-module $V_1$ into $D_2$-module $V_2$, then, according to the definition 2.5.1, I use notation

$$f \circ a = f(a)$$

for image of the map $f$.

Theorem 4.3.2. *Linear map*

$$(4.3.1) \qquad (h : D_1 \to D_2 \quad f : V_1 \to V_2)$$

*of $D_1$-module $V_1$ into $D_2$-module $V_2$ satisfies to equalities* [4.4]

$$(4.3.2) \qquad h(p + q) = h(p) + h(q)$$

$$(4.3.3) \qquad h(pq) = h(p)h(q)$$

$$(4.3.4) \qquad f \circ (a + b) = f \circ a + f \circ b$$

$$(4.3.5) \qquad f \circ (pa) = h(p)(f \circ a)$$

$$p, q \in D_1 \quad v, w \in V_1$$

PROOF.    According to definitions 2.3.9, the map $h$ is homomorphism of ring $D_1$ into ring $D_2$. Equalities (4.3.2), (4.3.3). follow from this statement.    From the definition 2.3.9, it follows that the map $f$ is a homomorphism of the Abelian group $V_1$ into the Abelian group $V_2$ (the equality (4.3.4)).    The equality (4.3.5) follows from the equality (2.3.6).                    □

---

THEOREM 4.3.3. *Linear map*[4.5]

(4.3.6)    $(h : D_1 \to D_2 \quad \overline{f} : V_1 \to V_2)$

*of $D_1$-module of columns $A_1$ into $D_2$-module of columns $A_2$ has presentation*

(4.3.7)    $\overline{f} \circ (e_1 a) = e_2 f h(a)$

(4.3.8)    $\overline{f} \circ (a^i e_{1i}) = h(a^i) f_i^k e_{2k}$

(4.3.9)    $b = f h(a)$

*relative to selected bases. Here*

- $a$ *is coordinate matrix of $V_1$- number $\overline{a}$ relative the basis $\overline{\overline{e}}_1$*

(4.3.10)    $\overline{a} = e_1 a$

- $h(a) = (h(a_i), i \in I)$ *is a matrix of $D_2$-numbers.*
- $b$ *is coordinate matrix of $V_2$- number*

    $\overline{b} = \overline{f} \circ \overline{a}$

    *relative the basis $\overline{\overline{e}}_2$*

(4.3.11)    $\overline{b} = e_2 b$

- $f$ *is coordinate matrix of set of $V_2$-numbers $(\overline{f} \circ e_{1i}, i \in I)$ relative the basis $\overline{\overline{e}}_2$.*

THEOREM 4.3.4. *Linear map*[4.6]

(4.3.12)    $(h : D_1 \to D_2 \quad \overline{f} : V_1 \to V_2)$

*of $D_1$-module of rows $A_1$ into $D_2$-module of rows $A_2$ has presentation*

(4.3.13)    $\overline{f} \circ (a e_1) = h(a) f e_2$

(4.3.14)    $\overline{f} \circ (a_i e_1^i) = h(a_i) f_k^i e_2^k$

(4.3.15)    $b = h(a) f$

*relative to selected bases. Here*

- $a$ *is coordinate matrix of $V_1$- number $\overline{a}$ relative the basis $\overline{\overline{e}}_1$*

(4.3.16)    $\overline{a} = a e_1$

- $h(a) = (h(a^i), i \in I)$ *is a matrix of $D_2$-numbers.*
- $b$ *is coordinate matrix of $V_2$- number*

    $\overline{b} = \overline{f} \circ \overline{a}$

    *relative the basis $\overline{\overline{e}}_2$*

(4.3.17)    $\overline{b} = b e_2$

- $f$ *is coordinate matrix of set of $V_2$-numbers $(\overline{f} \circ e_1^i, i \in I)$ relative the basis $\overline{\overline{e}}_2$.*

---

PROOF OF THEOREM 4.3.3.    Since

(4.3.6)    $(h : D_1 \to D_2 \quad \overline{f} : V_1 \to V_2)$

is a linear map, then the equality

(4.3.18)    $\overline{b} = \overline{f} \circ \overline{a} = \overline{f} \circ (e_1 a) = (\overline{f} \circ e_1) h(a)$

follows from equalities

(4.3.5)   $f \circ (pa) = h(p)(f \circ a)$         (4.3.10)   $\overline{a} = e_1 a$

---

[4.4]    In some books (for instance, on page [1]-119) the theorem 4.3.2 is considered as a definition.

[4.5]    In theorems 4.3.3, 4.3.5, we use the following convention. Let $\overline{\overline{e}}_1 = (e_{1i}, i \in I)$ be a basis of $D_1$-module of columns $V_1$. Let $\overline{\overline{e}}_2 = (e_{2j}, j \in J)$ be a basis of $D_2$-module of columns $V_2$.

[4.6]    In theorems 4.3.4, 4.3.6, we use the following convention. Let $\overline{\overline{e}}_1 = (e_1^i, i \in I)$ be a basis of $D_1$-module of rows $V_1$.    Let $\overline{\overline{e}}_2 = (e_2^j, j \in J)$ be a basis of $D_2$-module of rows $V_2$.

$V_2$-number $\overline{f} \circ e_{1i}$ has expansion

$$(4.3.19) \qquad \overline{f} \circ e_{1i} = e_2 f_i = e_{2j} f_i^j$$

relative to basis $\overline{\overline{e}}_2$. Combining (4.3.18) and (4.3.19), we get

$$(4.3.20) \qquad \overline{b} = e_2 f h(a)$$

(4.3.9) follows from comparison of

$$\boxed{(4.3.11) \quad \overline{b} = e_2 b}$$

and (4.3.20) and the theorem 4.2.7.                     □

PROOF OF THEOREM 4.3.4.    Since

$$\boxed{(4.3.12) \quad (h : D_1 \to D_2 \quad \overline{f} : V_1 \to V_2)}$$

is a linear map, then the equality

$$(4.3.21) \qquad \overline{b} = \overline{f} \circ \overline{a} = \overline{f} \circ (ae_1) = h(a)(\overline{f} \circ e_1)$$

follows from equalities

$$\boxed{(4.3.5) \quad f \circ (pa) = h(p)(f \circ a)} \qquad \boxed{(4.3.16) \quad \overline{a} = ae_1}$$

$V_2$-number $\overline{f} \circ e_1^i$ has expansion

$$(4.3.22) \qquad \overline{f} \circ e_1^i = f_i e_2 = f_j^i e_2^j$$

relative to basis $\overline{\overline{e}}_2$. Combining (4.3.21) and (4.3.22), we get

$$(4.3.23) \qquad \overline{b} = h(a) f e_2$$

(4.3.15) follows from comparison of

$$\boxed{(4.3.17) \quad \overline{b} = be_2}$$

and (4.3.23) and the theorem 4.2.14.                     □

THEOREM 4.3.5. *Let the map*

$$h : D_1 \to D_2$$

*be homomorphism of the ring $D_1$ into the ring $D_2$.     Let*

$$f = (f_j^i, i \in I, j \in J)$$

*be matrix of $D_2$-numbers.    The map*[4.5]

$$\boxed{(4.3.6) \quad (h : D_1 \to D_2 \quad \overline{f} : V_1 \to V_2)}$$

*defined by the equality*

$$\boxed{(4.3.7) \quad \overline{f} \circ (e_1 a) = e_2 f h(a)}$$

*is homomorphism of $D_1$-module of columns $V_1$ into $D_2$-module of columns $V_2$.   The homomorphism (4.3.6) which has the given  matrix $f$  is unique.*

THEOREM 4.3.6. *Let the map*

$$h : D_1 \to D_2$$

*be homomorphism of the ring $D_1$ into the ring $D_2$.     Let*

$$f = (f_i^j, i \in I, j \in J)$$

*be matrix of $D_2$-numbers.    The map*[4.6]

$$\boxed{(4.3.12) \quad (h : D_1 \to D_2 \quad \overline{f} : V_1 \to V_2)}$$

*defined by the equality*

$$\boxed{(4.3.13) \, \overline{f} \circ (ae_1) = h(a) f e_2}$$

*is homomorphism of $D_1$-module of rows $V_1$  into $D_2$-module of rows $V_2$.   The homomorphism (4.3.12) which has the given matrix $f$  is unique.*

Proof. The equality

$$
\begin{aligned}
&\overline{f} \circ (e_1(v+w)) \\
={}& e_2 f h(v+w) \\
(4.3.24)\qquad ={}& e_2 f(h(v)+h(w)) \\
={}& e_2 f h(v) + e_2 f h(w) \\
={}& \overline{f} \circ (e_1 v) + \overline{f} \circ (e_1 w)
\end{aligned}
$$

follows from equalities

$$\boxed{(4.3.4)\quad f \circ (a+b) = f \circ a + f \circ b}$$

$$\boxed{(4.3.5)\quad f \circ (pa) = h(p)(f \circ a)}$$

From the equality (4.3.24), it follows that the map $\overline{f}$ is homomorphism of Abelian group. The equality

$$
\begin{aligned}
&\overline{f} \circ (e_1(av)) \\
(4.3.25)\qquad ={}& e_2 f h(av) \\
={}& h(a)(e_2 f h(v)) \\
={}& h(a)(\overline{f} \circ (e_1 v))
\end{aligned}
$$

follows from the equality (4.3.5). From the equality (4.3.25), and definitions 2.8.5, 4.3.1, it follows that the map

$$\boxed{(4.3.6)\quad (h: D_1 \to D_2 \quad \overline{f}: V_1 \to V_2)}$$

is homomorphism of $D_1$-module of columns $V_1$ into $D_2$-module of columns $V_2$.

Let $f$ be matrix of homomorphisms $\overline{f}$, $\overline{g}$ relative to bases $\overline{\overline{e}}_1$, $\overline{\overline{e}}_2$. The equality

$$\overline{f} \circ (e_1 v) = e_2 f h(v) = \overline{g} \circ (e_1 v)$$

follows from the theorem 4.3.2. Therefore, $\overline{f} = \overline{g}$.  □

Proof. The equality

$$
\begin{aligned}
&\overline{f} \circ ((v+w)e_1) \\
={}& h(v+w) f e_2 \\
(4.3.26)\qquad ={}& (h(v)+h(w)) f e_2 \\
={}& h(v) f e_2 + h(w) f e_2 \\
={}& \overline{f} \circ (v e_1) + \overline{f} \circ (w e_1)
\end{aligned}
$$

follows from equalities

$$\boxed{(4.3.4)\quad f \circ (a+b) = f \circ a + f \circ b}$$

$$\boxed{(4.3.5)\quad f \circ (pa) = h(p)(f \circ a)}$$

From the equality (4.3.26), it follows that the map $\overline{f}$ is homomorphism of Abelian group. The equality

$$
\begin{aligned}
&\overline{f} \circ ((av)e_1) \\
(4.3.27)\qquad ={}& h(av) f e_2 \\
={}& h(a)(h(v) f e_2) \\
={}& h(a)(\overline{f} \circ (v e_1))
\end{aligned}
$$

follows from the equality (4.3.5). From the equality (4.3.27), and definitions 2.8.5, 4.3.1, it follows that the map

$$\boxed{(4.3.12)\quad (h: D_1 \to D_2 \quad \overline{f}: V_1 \to V_2)}$$

is homomorphism of $D_1$-module of rows $V_1$ into $D_2$-module of rows $V_2$.

Let $f$ be matrix of homomorphisms $\overline{f}$, $\overline{g}$ relative to bases $\overline{\overline{e}}_1$, $\overline{\overline{e}}_2$. The equality

$$\overline{f} \circ (v e_1) = h(v) f e_2 = \overline{g} \circ (v e_1)$$

follows from the theorem 4.3.2. Therefore, $\overline{f} = \overline{g}$.  □

On the basis of theorems 4.3.3, 4.3.5, as well on the basis of theorems 4.3.4, 4.3.6, we identify the linear map

$$(h: D_1 \to D_2 \quad \overline{f}: V_1 \to V_2)$$

and the matrix $f$ of its presentation.

Theorem 4.3.7. *The map*

$$(I_{(1)}: d \in D \to d+0 \in D_{(1)} \quad \text{id}: v \in V \to v \in V)$$

*is homomorphism of $D$-module $V$ into $D_{(1)}$-module $V$.*

Proof. The equality

$$(4.3.28) \qquad (d_1 + 0) + (d_2 + 0) = (d_1 + d_2) + (0 + 0) = (d_1 + d_2) + 0$$

follows from the equality (3.8.5). The equality

$$(4.3.29) \qquad (d_1 + 0)(d_2 + 0) = (d_1 d_2 + 0 d_1 + 0 d_2) + (0\,0) = (d_1 d_2) + 0$$

follows from the equality (3.8.6). Equalities (4.3.28), (4.3.29) imply that the map $I_{(1)}$ is homomorphism of ring $D$ into ring $D_{(1)}$.

The equality

$$(4.3.30) \qquad (d + 0)v = dv + 0v = dv + 0 = dv$$

follows from the equality

$$\boxed{(4.1.23) \quad (p + q)v = pv + qv}$$

The theorem follows from the equality (4.3.30) and the definition 2.3.9. □

Convention 4.3.8. *According to theorems 4.1.9, 4.3.7, the set of words generating $D$-module $V$ is the same as the set of words generating $D_{(1)}$-module $V$. Therefore, without loss of generality, we will assume that ring $D$ has unit.* □

### 4.3.2. Linear Map When Rings $D_1 = D_2 = D$.

Definition 4.3.9. *Morphism of representations*

$$(4.3.31) \qquad f : V_1 \to V_2$$

*of $D$-module $V_1$ into $D$-module $V_2$ is called homomorphism or* **linear map** *of $D$-module $V_1$ into $D$-module $V_2$. Let us denote* $\mathcal{L}(D; V_1 \to V_2)$ *set of linear maps of $D$-module $V_1$ into $D$-module $V_2$.* □

If the map (4.3.31) is linear map of $D$-module $V_1$ into $D$-module $V_2$, then, according to the definition 2.5.1, I use notation

$$f \circ a = f(a)$$

for image of the map $f$.

Theorem 4.3.10. *Linear map*

$$\boxed{(4.3.31) \quad f : V_1 \to V_2}$$

*of $D$-module $V_1$ into $D$-module $V_2$ satisfies to equalities*[4.7]

$$(4.3.32) \qquad f \circ (a + b) = f \circ a + f \circ b$$

$$(4.3.33) \qquad f \circ (pa) = p(f \circ a)$$

$$p \in D \quad v, w \in V_1$$

---

[4.7] In some books (for instance, on page [1]-119) the theorem 4.3.10 is considered as a definition.

Proof.            From the definition 2.3.9, it follows that the map $f$ is a homomorphism of the Abelian group $V_1$ into the Abelian group $V_2$ (the equality (4.3.32)). The equality (4.3.33) follows from the equality (2.3.6).                □

| | |
|---|---|
| Theorem 4.3.11. *Linear map*[4.8] | Theorem 4.3.12. *Linear map*[4.9] |
| (4.3.34)       $\overline{f} : V_1 \to V_2$ | (4.3.40)       $\overline{f} : V_1 \to V_2$ |
| *of D-module of columns $A_1$ into D-module of columns $A_2$ has presentation* | *of D-module of rows $A_1$ into D-module of rows $A_2$ has presentation* |
| (4.3.35)       $\overline{f} \circ (e_1 a) = e_2 f a$ | (4.3.41)       $\overline{f} \circ (a e_1) = a f e_2$ |
| (4.3.36)       $\overline{f} \circ (a^i e_{1i}) = a^i f_i^k e_{2k}$ | (4.3.42)       $\overline{f} \circ (a_i e_1^i) = a_i f_k^i e_2^k$ |
| (4.3.37)       $b = f a$ | (4.3.43)       $b = a f$ |
| *relative to selected bases. Here* | *relative to selected bases. Here* |
| • *a is coordinate matrix of $V_1$-number $\overline{a}$ relative the basis $\overline{\overline{e}}_1$* | • *a is coordinate matrix of $V_1$-number $\overline{a}$ relative the basis $\overline{\overline{e}}_1$* |
| (4.3.38)       $\overline{a} = e_1 a$ | (4.3.44)       $\overline{a} = a e_1$ |
| • *b is coordinate matrix of $V_2$-number* $$\overline{b} = \overline{f} \circ \overline{a}$$ *relative the basis* $\overline{\overline{e}}_2$ | • *b is coordinate matrix of $V_2$-number* $$\overline{b} = \overline{f} \circ \overline{a}$$ *relative the basis* $\overline{\overline{e}}_2$ |
| (4.3.39)       $\overline{b} = e_2 b$ | (4.3.45)       $\overline{b} = b e_2$ |
| • *f is coordinate matrix of set of $V_2$-numbers $(\overline{f} \circ e_{1i}, i \in I)$ relative the basis $\overline{\overline{e}}_2$.* | • *f is coordinate matrix of set of $V_2$-numbers $(\overline{f} \circ e_1^i, i \in I)$ relative the basis $\overline{\overline{e}}_2$.* |

Proof of theorem 4.3.11.      Since

$\boxed{(4.3.34)\quad \overline{f} : V_1 \to V_2}$

is a linear map, then the equality

(4.3.46)              $\overline{b} = \overline{f} \circ \overline{a} = \overline{f} \circ (e_1 a) = (\overline{f} \circ e_1) a$

follows from equalities

$\boxed{(4.3.33)\quad f \circ (pa) = p(f \circ a)}$     $\boxed{(4.3.38)\quad \overline{a} = e_1 a}$

$V_2$-number $\overline{f} \circ e_{1i}$ has expansion

(4.3.47)                    $\overline{f} \circ e_{1i} = e_2 f_i = e_{2j} f_i^j$

relative to basis $\overline{\overline{e}}_2$. Combining (4.3.46) and (4.3.47), we get

(4.3.48)                    $\overline{b} = e_2 f a$

---

[4.8]  In theorems 4.3.11, 4.3.13, we use the following convention. Let $\overline{\overline{e}}_1 = (e_{1i}, i \in I)$ be a basis of $D$-module of columns $V_1$. Let $\overline{\overline{e}}_2 = (e_{2j}, j \in J)$ be a basis of $D$-module of columns $V_2$.

[4.9]  In theorems 4.3.12, 4.3.14, we use the following convention. Let $\overline{\overline{e}}_1 = (e_1^i, i \in I)$ be a basis of $D$-module of rows $V_1$. Let $\overline{\overline{e}}_2 = (e_2^j, j \in J)$ be a basis of $D$-module of rows $V_2$.

(4.3.37) follows from comparison of

$$(4.3.39) \quad \overline{b} = e_2 b$$

and (4.3.48) and the theorem 4.2.7.          □

PROOF OF THEOREM 4.3.12.     Since

$$(4.3.40) \quad \overline{f} : V_1 \to V_2$$

is a linear map, then the equality

$$(4.3.49) \qquad \overline{b} = \overline{f} \circ \overline{a} = \overline{f} \circ (ae_1) = a(\overline{f} \circ e_1)$$

follows from equalities

$$(4.3.33) \quad f \circ (pa) = p(f \circ a) \qquad (4.3.44) \quad \overline{a} = ae_1$$

$V_2$-number $\overline{f} \circ e_1^i$ has expansion

$$(4.3.50) \qquad \overline{f} \circ e_1^i = f_i e_2 = f_j^i e_2^j$$

relative to basis $\overline{\overline{e}}_2$. Combining (4.3.49) and (4.3.50), we get

$$(4.3.51) \qquad \overline{b} = afe_2$$

(4.3.43) follows from comparison of

$$(4.3.45) \quad \overline{b} = be_2$$

and (4.3.51) and the theorem 4.2.14.          □

THEOREM 4.3.13. *Let*

$$f = (f_j^i, i \in I, j \in J)$$

*be matrix of $D$-numbers.     The map*[4.8]

$$(4.3.34) \quad \overline{f} : V_1 \to V_2 \quad \text{defined by the}$$

*equality*

$$(4.3.35) \quad \overline{f} \circ (e_1 a) = e_2 f a$$

*is homomorphism of $D$-module of columns $V_1$ into $D$-module of columns $V_2$. The homomorphism (4.3.34) which has the given matrix $f$ is unique.*

THEOREM 4.3.14. *Let*

$$f = (f_i^j, i \in I, j \in J)$$

*be matrix of $D$-numbers.     The map*[4.9]

$$(4.3.40) \quad \overline{f} : V_1 \to V_2 \quad \text{defined by the}$$

*equality*

$$(4.3.41) \quad \overline{f} \circ (ae_1) = afe_2$$

*is homomorphism of $D$-module of rows $V_1$ into $D$-module of rows $V_2$. The homomorphism (4.3.40) which has the given matrix $f$ is unique.*

PROOF. The equality

$$(4.3.52) \quad \begin{aligned} &\overline{f} \circ (e_1(v+w)) \\ &= e_2 f(v+w) \\ &= e_2 fv + e_2 fw \\ &= \overline{f} \circ (e_1 v) + \overline{f} \circ (e_1 w) \end{aligned}$$

follows from equalities

$$(4.3.32) \quad f \circ (a+b) = f \circ a + f \circ b$$

$$(4.3.33) \quad f \circ (pa) = p(f \circ a)$$

From the equality (4.3.52), it follows that the map $\overline{f}$ is homomorphism of Abelian group. The equality

PROOF. The equality

$$(4.3.54) \quad \begin{aligned} &\overline{f} \circ ((v+w)e_1) \\ &= (v+w)fe_2 \\ &= vfe_2 + wfe_2 \\ &= \overline{f} \circ (ve_1) + \overline{f} \circ (we_1) \end{aligned}$$

follows from equalities

$$(4.3.32) \quad f \circ (a+b) = f \circ a + f \circ b$$

$$(4.3.33) \quad f \circ (pa) = p(f \circ a)$$

From the equality (4.3.54), it follows that the map $\overline{f}$ is homomorphism of Abelian group. The equality

$$(4.3.53) \quad \begin{aligned} &\overline{f} \circ (e_1(av)) \\ &= e_2 f(av) \\ &= a(e_2 f v) \\ &= a(\overline{f} \circ (e_1 v)) \end{aligned}$$

$$(4.3.55) \quad \begin{aligned} &\overline{f} \circ ((av)e_1) \\ &= (av)f e_2 \\ &= a(v f e_2) \\ &= a(\overline{f} \circ (v e_1)) \end{aligned}$$

follows from the equality (4.3.33). From the equality (4.3.53), and definitions 2.8.5, 4.3.9, it follows that the map

$$\boxed{(4.3.34) \quad \overline{f} : V_1 \to V_2}$$

is homomorphism of $D$-module of columns $V_1$ into $D$-module of columns $V_2$.

Let $f$ be matrix of homomorphisms $\overline{f}$, $\overline{g}$ relative to bases $\overline{\overline{e}}_1$, $\overline{\overline{e}}_2$. The equality

$$\overline{f} \circ (e_1 v) = e_2 f v = \overline{g} \circ (e_1 v)$$

follows from the theorem 4.3.10. Therefore, $\overline{f} = \overline{g}$. $\qquad \square$

follows from the equality (4.3.33). From the equality (4.3.55), and definitions 2.8.5, 4.3.9, it follows that the map

$$\boxed{(4.3.40) \quad \overline{f} : V_1 \to V_2}$$

is homomorphism of $D$-module of rows $V_1$ into $D$-module of rows $V_2$.

Let $f$ be matrix of homomorphisms $\overline{f}$, $\overline{g}$ relative to bases $\overline{\overline{e}}_1$, $\overline{\overline{e}}_2$. The equality

$$\overline{f} \circ (v e_1) = v f e_2 = \overline{g} \circ (v e_1)$$

follows from the theorem 4.3.10. Therefore, $\overline{f} = \overline{g}$. $\qquad \square$

On the basis of theorems 4.3.11, 4.3.13, as well on the basis of theorems 4.3.12, 4.3.14, we identify the linear map

$$\overline{f} : V_1 \to V_2$$

and the matrix $f$ of its presentation.

THEOREM 4.3.15. *Let map*

$$f : A_1 \to A_2$$

*be linear map of $D$-module $A_1$ into $D$-module $A_2$. Then*

$$f \circ 0 = 0$$

PROOF. The theorem follows from the equality

$$f \circ (a + 0) = f \circ a + f \circ 0$$

$$\square$$

DEFINITION 4.3.16. *Linear map*

$$f : V \to V$$

*of $D$-module $V$ is called* **endomorphism** *of $D$-module $V$. We use notation* $\mathrm{End}(D, V)$ *set of endomorphisms of $D$-module $V$.* $\qquad \square$

## 4.4. Polylinear Map of $D$-Module

DEFINITION 4.4.1. *Let $D$ be the commutative ring. Reduced polymorphism of $D$-modules $A_1$, ..., $A_n$ into $D$-module $S$*

$$f : A_1 \times ... \times A_n \to S$$

*is called* **polylinear map** *of $D$-modules $A_1$, ..., $A_n$ into $D$-module $S$.*

*We denote $\mathcal{L}(D; A_1 \times ... \times A_n \to S)$ the set of polylinear maps of $D$-modules $A_1$, ..., $A_n$ into $D$-module $S$. Let us denote $\mathcal{L}(D; A^n \to S)$ set of $n$-linear maps of $D$-module $A$ $(A_1 = ... = A_n = A)$ into $D$-module $S$.*          $\square$

THEOREM 4.4.2. *Let $D$ be the commutative ring. The polylinear map of $D$-modules $A_1$, ..., $A_n$ into $D$-module $S$*

$$f : A_1 \times ... \times A_n \to S$$

*satisfies to equalities*

$$f \circ (a_1, ..., a_i + b_i, ..., a_n) = f \circ (a_1, ..., a_i, ..., a_n) + f \circ (a_1, ..., b_i, ..., a_n)$$

$$f \circ (a_1, ..., pa_i, ..., a_n) = pf \circ (a_1, ..., a_i, ..., a_n)$$

$$1 \le i \le n \quad a_i, b_i \in A_i \quad p \in D$$

PROOF. The theorem follows from definitions 2.6.1, 4.3.9, 4.4.1 and from the theorem 4.3.10.          $\square$

THEOREM 4.4.3. *Let $D$ be the commutative ring. Let $A_1$, ..., $A_n$, $S$ be $D$-modules. The map*

$$(4.4.1) \qquad f + g : A_1 \times ... \times A_n \to S \quad f, g \in \mathcal{L}(D; A_1 \times ... \times A_n \to S)$$

*defined by the equality*

$$(4.4.2) \qquad (f + g) \circ (a_1, ..., a_n) = f \circ (a_1, ..., a_n) + g \circ (a_1, ..., a_n)$$

*is called* **sum of polylinear maps** *$f$ and $g$ and is polylinear map. The set $\mathcal{L}(D; A_1 \times ... \times A_n \to S)$ is an Abelian group relative sum of maps.*

PROOF. According to the theorem 4.4.2

$$(4.4.3) \qquad f \circ (a_1, ..., a_i + b_i, ..., a_n) = f \circ (a_1, ..., a_i, ..., a_n) + f \circ (a_1, ..., b_i, ..., a_n)$$

$$(4.4.4) \qquad f \circ (a_1, ..., pa_i, ..., a_n) = pf \circ (a_1, ..., a_i, ..., a_n)$$

$$(4.4.5) \qquad g \circ (a_1, ..., a_i + b_i, ..., a_n) = g \circ (a_1, ..., a_i, ..., a_n) + g \circ (a_1, ..., b_i, ..., a_n)$$

$$(4.4.6) \qquad g \circ (a_1, ..., pa_i, ..., a_n) = pg \circ (a_1, ..., a_i, ..., a_n)$$

The equality

$$(4.4.7) \qquad \begin{aligned} & (f + g) \circ (x_1, ..., x_i + y_i, ..., x_n) \\ = & f \circ (x_1, ..., x_i + y_i, ..., x_n) + g \circ (x_1, ..., x_i + y_i, ..., x_n) \\ = & f \circ (x_1, ..., x_i, ..., x_n) + f \circ (x_1, ..., y_i, ..., x_n) \\ & + g \circ (x_1, ..., x_i, ..., x_n) + g \circ (x_1, ..., y_i, ..., x_n) \\ = & (f + g) \circ (x_1, ..., x_i, ..., x_n) + (f + g) \circ (x_1, ..., y_i, ..., x_n) \end{aligned}$$

follows from the equalities (4.4.2), (4.4.3), (4.4.5). The equality

$$
\begin{aligned}
&(f+g)\circ(x_1,...,px_i,...,x_n)\\
=&f\circ(x_1,...,px_i,...,x_n)+g\circ(x_1,...,px_i,...,x_n)\\
=&pf\circ(x_1,...,x_i,...,x_n)+pg\circ(x_1,...,x_i,...,x_n)\\
=&p(f\circ(x_1,...,x_i,...,x_n)+g\circ(x_1,...,x_i,...,x_n))\\
=&p(f+g)\circ(x_1,...,x_i,...,x_n)
\end{aligned}
$$

(4.4.8)

follows from the equalities (4.4.2), (4.4.4), (4.4.6). From equalities (4.4.7), (4.4.8) and from the theorem 4.4.2, it follows that the map (4.4.1) is linear map of $D$-modules.

Let $f$, $g$, $h \in \mathcal{L}(D; A_1 \times ... \times A_2 \to S)$. For any $a=(a_1,...,a_n)$, $a_1 \in A_1$, ..., $a_n \in A_n$,

$$
\begin{aligned}
(f+g)\circ a &= f\circ a+g\circ a=g\circ a+f\circ a\\
&=(g+f)\circ a\\
((f+g)+h)\circ a &=(f+g)\circ a+h\circ a=(f\circ a+g\circ a)+h\circ a\\
&=f\circ a+(g\circ a+h\circ a)=f\circ a+(g+h)\circ a\\
&=(f+(g+h))\circ a
\end{aligned}
$$

Therefore, sum of polylinear maps is commutative and associative.

From the equality (4.4.2), it follows that the map

$$0: v \in A_1 \times ... \times A_n \to 0 \in S$$

is zero of addition

$$(0+f)\circ(a_1,...,a_n)=0\circ(a_1,...,a_n)+f\circ(a_1,...,a_n)=f\circ(a_1,...,a_n)$$

From the equality (4.4.2), it follows that the map

$$-f:(a_1,...,a_n)\in A_1 \times ... \times A_n \to -(f\circ(a_1,...,a_n))\in S$$

is map inversed to map $f$

$$f+(-f)=0$$

because

$$
\begin{aligned}
(f+(-f))\circ(a_1,...,a_n) &= f\circ(a_1,...,a_n)+(-f)\circ(a_1,...,a_n)\\
&=f\circ(a_1,...,a_n)-f\circ(a_1,...,a_n)\\
&=0=0\circ(a_1,...,a_n)
\end{aligned}
$$

From the equality

$$
\begin{aligned}
(f+g)\circ(a_1,...,a_n) &= f\circ(a_1,...,a_n)+g\circ(a_1,...,a_n)\\
&=g\circ(a_1,...,a_n)+f\circ(a_1,...,a_n)\\
&=(g+f)\circ(a_1,...,a_n)
\end{aligned}
$$

it follows that sum of maps is commutative. Therefore, the set $\mathcal{L}(D; A_1 \times ... \times A_n \to S)$ is an Abelian group. $\square$

Corollary 4.4.4. *Let $A_1$, $A_2$ be D-modules. The map*

$$(4.4.9) \qquad f + g : A_1 \to A_2 \quad f, g \in \mathcal{L}(D; A_1 \to A_2)$$

*defined by equation*

$$(4.4.10) \qquad (f + g) \circ x = f \circ x + g \circ x$$

*is called* **sum of maps** *f and g and is linear map. The set $\mathcal{L}(D; A_1 \to A_2)$ is an Abelian group relative sum of maps.* □

Theorem 4.4.5. *Let D be the commutative ring. Let $A_1$, ..., $A_n$, S be D-modules. The map*

$$(4.4.11) \qquad d\,f : A_1 \times ... \times A_n \to S \quad d \in D \quad f \in \mathcal{L}(D; A_1 \times ... \times A_n \to S)$$

*defined by equality*

$$(4.4.12) \qquad (d\,f) \circ (a_1, ..., a_n) = d(f \circ (a_1, ..., a_n))$$

*is polylinear map and is called* **product of map** *f* **over scalar** *d. The representation*

$$(4.4.13) \qquad a : f \in \mathcal{L}(D; A_1 \times ... \times A_n \to S) \to af \in \mathcal{L}(D; A_1 \times ... \times A_n \to S)$$

*of ring D in Abelian group $\mathcal{L}(D; A_1 \times ... \times A_n \to S)$ generates structure of D-module.*

Proof. According to the theorem 4.4.2

$$(4.4.14) \quad f \circ (a_1, ..., a_i + b_i, ..., a_n) = f \circ (a_1, ..., a_i, ..., a_n) + f \circ (a_1, ..., b_i, ..., a_n)$$

$$(4.4.15) \qquad f \circ (a_1, ..., pa_i, ..., a_n) = pf \circ (a_1, ..., a_i, ..., a_n)$$

The equality

$$(4.4.16) \qquad \begin{aligned} &(pf) \circ (x_1, ..., x_i + y_i, ..., x_n) \\ &= p\ f \circ (x_1, ..., x_i + y_i, ..., x_n) \\ &= p\ (f \circ (x_1, ..., x_i, ..., x_n) + f \circ (x_1, ..., y_i, ..., x_n)) \\ &= p(f \circ (x_1, ..., x_i, ..., x_n)) + p(f \circ (x_1, ..., y_i, ..., x_n)) \\ &= (pf) \circ (x_1, ..., x_i, ..., x_n) + (pf) \circ (x_1, ..., y_i, ..., x_n) \end{aligned}$$

follows from equalities (4.4.12), (4.4.14). The equality

$$(4.4.17) \qquad \begin{aligned} &(pf) \circ (x_1, ..., qx_i, ..., x_n) \\ &= p(f \circ (x_1, ..., qx_i, ..., x_n)) = pq(f \circ (x_1, ..., x_i, ..., x_n)) \\ &= qp(f \circ (x_1, ..., x_n)) = q(pf) \circ (x_1, ..., x_n) \end{aligned}$$

follows from equalities (4.4.12), (4.4.15). From equalities (4.4.16), (4.4.17) and from the theorem 4.4.2, it follows that the map (4.4.11) is polylinear map of *D*-modules.

The equality

$$(4.4.18) \qquad (p + q)f = pf + qf$$

follows from the equality

$$((p+q)f) \circ (x_1, ..., x_n) = (p+q)(f \circ (x_1, ..., x_n))$$
$$= p(f \circ (x_1, ..., x_n)) + q(f \circ (x_1, ..., x_n))$$
$$= (pf) \circ (x_1, ..., x_n) + (qf) \circ (x_1, ..., x_n)$$

The equality

$$(4.4.19) \qquad\qquad p(qf) = (pq)f$$

follows from the equality

$$(p(qf)) \circ (x_1, ..., x_n) = p (qf) \circ (x_1, ..., x_n) = p (q f \circ (x_1, ..., x_n))$$
$$= (pq) f \circ (x_1, ..., x_n) = ((pq)f) \circ (x_1, ..., x_n)$$

From equalities (4.4.18) (4.4.19) it follows that the map (4.4.13) is representation of ring $D$ in Abelian group $\mathcal{L}(D; A_1 \times ... \times A_n \to S)$. Since specified representation is effective, then, according to the definition 4.1.3 and the theorem 4.4.3, Abelian group $\mathcal{L}(D; A_1 \to A_2)$ is $D$-module. □

Corollary 4.4.6. *Let $A_1$, $A_2$ be $D$-modules. The map*

$$(4.4.20) \qquad\qquad d\,f : A_1 \to A_2 \quad d \in D \quad f \in \mathcal{L}(D; A_1 \to A_2)$$

*defined by the equality*

$$(4.4.21) \qquad\qquad (d\,f) \circ x = d(f \circ x)$$

*is linear map and is called* **product of map** $f$ **over scalar** $d$. *The representation*

$$(4.4.22) \qquad a : f \in \mathcal{L}(D; A_1 \to A_2) \to af \in \mathcal{L}(D; A_1 \to A_2)$$

*of ring $D$ in Abelian group $\mathcal{L}(D; A_1 \to A_2)$ generates structure of $D$-module.* □

Theorem 4.4.7. *The set* $\mathrm{End}(D, V)$ *of endomorphisms of $D$-module $V$ is $D$-module.*

Proof. The theorem follows from the definition 4.3.16 and the corollary 4.4.6.
□

## 4.5. $D$-module $\mathcal{L}(D; A \to B)$

Theorem 4.5.1.

$$(4.5.1) \qquad \mathcal{L}(D; A^p \to \mathcal{L}(D; A^q \to B)) = \mathcal{L}(D; A^{p+q} \to B)$$

Proof. □

Theorem 4.5.2. *Let*

$$\overline{\overline{e}} = \{e_i : i \in I\}$$

*be basis of $D$-module $A$. The set*

$$(4.5.2) \qquad \overline{\overline{h}} = \{h^i \in \mathcal{L}(D; A \to D) : i \in I, h^i \circ e_j = \delta^i_j\}$$

*is the basis of $D$-module $\mathcal{L}(D; A \to D)$.*

Proof.

Lemma 4.5.3. *Maps $h^i$ are linear independent.*

Proof. Let there exist $D$-numbers $c_i$ such that

$$c_i h^i = 0$$

Then for any $A$-number $e_j$

$$0 = c_i h^i \circ e_j = c_i \delta^i_j = c_j$$

The lemma follows from the definition 4.1.16.    ⊙

Lemma 4.5.4. *The map $f \in \mathcal{L}(D; A \to D)$, is linear composition of maps $h^i$.*

Proof. For any $A$-number $a$,

(4.5.3) $$a = a^i e_i$$

the equality

(4.5.4) $$h^i \circ a = h^i \circ (a^j e_j) = a^j (h^i \circ e_j) = a^j \delta^i_j = a^i$$

follows from equalities (4.5.2), (4.5.3) and from the theorem 4.3.10. The equality

(4.5.5) $$f \circ a = f \circ (a^i e_i) = a^i (f \circ e_i) = (f \circ e_i)(h^i \circ a)$$

follows from the equality (4.5.4). The equality

$$f = (f \circ e_i) h^i$$

follows from equalities (4.4.10), (4.4.21), (4.5.5).    ⊙

The theorem follows from lemmas 4.5.3, 4.5.4 and from the definition 4.1.24.  □

Theorem 4.5.5. *Let $D$ be commutative ring. Let*

$$\overline{\overline{e}}_i = \{ e_{i \cdot i} : i \in I_i \}$$

*be basis of $D$-module $A_i$, $i = 1, \ldots, n$. Let*

$$\overline{\overline{e}}_B = \{ e_{B \cdot i} : i \in I \}$$

*be basis of $D$-module $B$. The set*

(4.5.6) $$\overline{\overline{h}} = \{ h_i^{i_1 \ldots i_n} \in \mathcal{L}(D; A_1 \times \ldots \times A_n \to B) : i \in I, i_i \in I_i, i = 1, \ldots, n,$$
$$h_i^{i_1 \ldots i_n} \circ (e_{1 \cdot j_1}, \ldots, e_{n \cdot j_n}) = \delta^{i_1}_{j_1} \ldots \delta^{i_n}_{j_n} e_{B \cdot i} \}$$

*is the basis of $D$-module $\mathcal{L}(D; A_1 \times \ldots \times A_n \to B)$.*

Proof.

Lemma 4.5.6. *Maps $h_i^{i_1 \ldots i_n}$ are linear independent.*

Proof. Let there exist $D$-numbers $c^i_{i_1 \ldots i_n}$ such that

$$c^i_{i_1 \ldots i_n} h_i^{i_1 \ldots i_n} = 0$$

Then for any set of indices $j_1, \ldots, j_n$

$$0 = c^i_{i_1 \ldots i_n} h_i^{i_1 \ldots i_n} \circ (e_{1 \cdot j_1}, \ldots, e_{n \cdot j_n}) = c^i_{i_1 \ldots i_n} \delta^{i_1}_{j_1} \ldots \delta^{i_n}_{j_n} e_{B \cdot i} = c^i_{j_1 \ldots j_n} e_{B \cdot i}$$

Therefore, $c^i_{j_1 \ldots j_n} = 0$. The lemma follows from the definition 4.1.16.    ⊙

LEMMA 4.5.7. *The map* $f \in \mathcal{L}(D; A_1 \times ... \times A_n \to B)$ *is linear composition of maps* $h_i^{i_1 ... i_n}$.

PROOF. For any $A_1$-number $a_1$

$$(4.5.7) \qquad\qquad a_1 = a_1^{i_1} e_{1 \cdot i_1}$$

..., for any $A_n$-number $a_n$

$$(4.5.8) \qquad\qquad a_n = a_n^{i_n} e_{n \cdot i_n}$$

the equality

$$
\begin{aligned}
(4.5.9) \qquad h_i^{i_1 ... i_n} \circ (a_1, ..., a_n) &= h_i^{i_1 ... i_n} \circ (a_1^{j_1} e_{1 \cdot j_1}, ..., a_n^{j_n} e_{n \cdot j_n}) \\
&= a_1^{j_1} ... a_n^{j_n} (h_i^{i_1 ... i_n} \circ (e_{1 \cdot j_1}, ..., e_{n \cdot j_n})) \\
&= a_1^{j_1} ... a_n^{j_n} \delta_{j_1}^{i_1} ... \delta_{j_n}^{i_n} e_{B \cdot i} \\
&= a_1^{i_1} ... a_n^{i_n} e_{B \cdot i}
\end{aligned}
$$

follows from equalities $(4.5.6)$, $(4.5.7)$, $(4.5.8)$ and from the theorem $4.3.10$. The equality

$$
\begin{aligned}
(4.5.10) \qquad f \circ (a_1, ..., a_n) &= f \circ (a_1^{j_1} e_{1 \cdot j_1}, ..., a_n^{j_n} e_{n \cdot j_n}) \\
&= a_1^{j_1} ... a_n^{j_n} f \circ (e_{1 \cdot j_1} ... e_{n \cdot j_n})
\end{aligned}
$$

follows from equalities $(4.5.7)$, $(4.5.8)$. Since

$$f \circ (e_{1 \cdot j_1} ... e_{n \cdot j_n}) \in B$$

then

$$(4.5.11) \qquad\qquad f \circ (e_{1 \cdot j_1} ... e_{n \cdot j_n}) = f_{j_1 ... j_n}^i e_i$$

The equality

$$
\begin{aligned}
(4.5.12) \qquad f \circ (a_1, ..., a_n) &= a_1^{i_1} ... a_n^{i_n} f_{i_1 ... i_n}^i e_i \\
&= f_{i_1 ... i_n}^i h_i^{i_1 ... i_n} \circ (a_1, ..., a_n)
\end{aligned}
$$

follows from equalities $(4.5.9)$, $(4.5.10)$, $(4.5.11)$. The equality

$$f = f_{i_1 ... i_n}^i h_i^{i_1 ... i_n}$$

follows from equalities $(4.4.2)$, $(4.4.12)$, $(4.5.12)$. $\qquad\qquad\qquad\odot$

The theorem follows from lemmas $4.5.6$, $4.5.7$ and from the definition $4.1.24$. $\quad\square$

THEOREM 4.5.8. *Let* $A_1$, ..., $A_n$, $B$ *be free modules over commutative ring* $D$. *D-module* $\mathcal{L}(D; A_1 \times ... \times A_n \to B)$ *is free D-module.*

PROOF. The theorem follows from the theorem $4.5.5$. $\qquad\qquad\qquad\qquad\square$

## 4.6. Tensor Product of *D*-Modules

THEOREM 4.6.1. *The commutative ring* $D$ *is Abelian multiplicative* $\Omega$-*group.*

PROOF. Let the product $\circ$ in the ring $D$ be defined according to rule

$$a \circ b = ab$$

Since the product in the ring is distributive over addition, the theorem follows from definitions 2.4.5, 2.4.6. □

THEOREM 4.6.2. *Let $A_1$, ..., $A_n$ be modules over commutative ring $D$. There exists and unique* **tensor product**

$$f : A_1 \times ... \times A_n \rightarrow A_1 \otimes ... \otimes A_n$$

*of $D$-modules $A_1$, ..., $A_n$. We use notation*

$$f \circ (a_1, ..., a_n) = a_1 \otimes ... \otimes a_n$$

*for the image of the map $f$.*

PROOF. The theorem follows from the definition 4.1.3 and from theorems 2.6.3, 4.6.1. □

THEOREM 4.6.3. *Let $D$ be the commutative ring. Let $A_1$, ..., $A_n$ be $D$-modules. Tensor product is distributive over sum*

$$(4.6.1) \quad \begin{aligned} & a_1 \otimes ... \otimes (a_i + b_i) \otimes ... \otimes a_n \\ = & a_1 \otimes ... \otimes a_i \otimes ... \otimes a_n + a_1 \otimes ... \otimes b_i \otimes ... \otimes a_n \\ & a_i, b_i \in A_i \end{aligned}$$

*The representation of the ring $D$ in tensor product is defined by equality*

$$(4.6.2) \quad \begin{aligned} & a_1 \otimes ... \otimes (ca_i) \otimes ... \otimes a_n = c(a_1 \otimes ... \otimes a_i \otimes ... \otimes a_n) \\ & a_i \in A_i \quad c \in D \end{aligned}$$

PROOF. The equality (4.6.1) follows from the equality (2.6.1). The equality (4.6.2) follows from the equality (2.6.2). □

THEOREM 4.6.4. *Let $A_1$, ..., $A_n$ be modules over commutative ring $D$. Let*

$$f : A_1 \times ... \times A_n \rightarrow A_1 \otimes ... \otimes A_n$$

*be polylinear map defined by the equality*

$$(4.6.3) \quad f \circ (d_1, ..., d_n) = d_1 \otimes ... \otimes d_n$$

*Let*

$$g : A_1 \times ... \times A_n \rightarrow V$$

*be polylinear map into $D$-module $V$. There exists a linear map*

$$h : A_1 \otimes ... \otimes A_n \rightarrow V$$

*such that the diagram*

$$(4.6.4)$$

*is commutative. The map h is defined by the equality*

$$(4.6.5) \qquad h \circ (a_1 \otimes ... \otimes a_n) = g \circ (a_1, ..., a_n)$$

PROOF. The theorem follows from the theorem 2.6.5 and from definitions 4.3.9, 4.4.1. □

THEOREM 4.6.5. *The map*

$$(v_1, ..., v_n) \in V_1 \times ... \times V_n \to v_1 \otimes ... \otimes v_n \in V_1 \otimes ... \otimes V_n$$

*is polylinear map.*

PROOF. The theorem follows from the theorem 2.6.6 and from the definition 4.4.1. □

THEOREM 4.6.6. *Tensor product $A_1 \otimes ... \otimes A_n$ of free finite dimensional modules $A_1, ..., A_n$ over the commutative ring $D$ is free finite dimensional module.*

*Let $\overline{\overline{e}}_i$ be the basis of module $A_i$ over ring $D$. We can represent any tensor $a \in A_1 \otimes ... \otimes A_n$ in the following form*

$$(4.6.6) \qquad a = a^{i_1 \cdots i_n} e_{1 \cdot i_1} \otimes ... \otimes e_{n \cdot i_n}$$

*The expression $a^{i_1 \cdots i_n}$ is defined uniquely and is called* **standard component of tensor**.

PROOF. Vector $a_i \in A_i$ has expansion

$$a_i = a_i^k \overline{e}_{i \cdot k}$$

relative to basis $\overline{\overline{e}}_i$. From equalities (4.6.1), (4.6.2), it follows

$$a_1 \otimes ... \otimes a_n = a_1^{i_1} ... a_n^{i_n} e_{1 \cdot i_1} \otimes ... \otimes e_{n \cdot i_n}$$

Since set of tensors $a_1 \otimes ... \otimes a_n$ is the generating set of module $A_1 \otimes ... \otimes A_n$, then we can write tensor $a \in A_1 \otimes ... \otimes A_n$ in form

$$(4.6.7) \qquad a = a^s a_{s \cdot 1}^{i_1} ... a_{s \cdot n}^{i_n} e_{1 \cdot i_1} \otimes ... \otimes e_{n \cdot i_n}$$

where $a^s, a_{s \cdot 1}^{i_1}, ..., a_{s \cdot n}^{i_n} \in F$. Let

$$a^s a_{s \cdot 1}^{i_1} ... a_{s \cdot n}^{i_n} = a^{i_1 \cdots i_n}$$

Then equality (4.6.7) has form (4.6.6).

Therefore, set of tensors $e_{1 \cdot i_1} \otimes ... \otimes e_{n \cdot i_n}$ is the generating set of module $A_1 \otimes ... \otimes A_n$. Since the dimension of module $A_i$, $i = 1, ..., n$, is finite, then the set of tensors $e_{1 \cdot i_1} \otimes ... \otimes e_{n \cdot i_n}$ is finite. Therefore, the set of tensors $e_{1 \cdot i_1} \otimes ... \otimes e_{n \cdot i_n}$ contains a basis of module $A_1 \otimes ... \otimes A_n$, and the module $A_1 \otimes ... \otimes A_n$ is free module over the ring $D$. □

## 4.7. Direct Sum of *D*-Modules

Let *D*-MODULE be category of *D*-modules morphisms of which are homomorphisms of *D*-module.

DEFINITION 4.7.1. *Coproduct in category D-*MODULE *is called* **direct sum**. [4.10] *We will use notation $V \oplus W$ for direct sum of D-modules $V$ and $W$.* □

Theorem 4.7.2. *In category $D$-Module there exists direct sum of $D$-modules.*[4.11] *Let $\{V_i, i \in I\}$ be set of $D$-modules. Then the representation*

$$D \xrightarrow{\ *\ } \bigoplus_{i \in I} V_i \qquad d\left(\bigoplus_{i \in I} a_i\right) = \bigoplus_{i \in I} da_i$$

*of the ring $D$ in direct sum of Abelian groups*

$$V = \bigoplus_{i \in I} V_i$$

*is direct sum of $D$-modules*

$$V = \bigoplus_{i \in I} V_i$$

Proof. Let $V_i$, $i \in I$, be set of $D$-modules. Let

$$V = \bigoplus_{i \in I} V_i$$

be direct sum of Abelian groups. Let $v = (v_i, i \in I) \in V$. We define action of $D$-number $a$ over $V$-number $v$ according to the equality

$$(4.7.1) \qquad a(v_i, i \in I) = (av_i, i \in I)$$

According to the theorem 4.1.8,

$$(4.7.2) \qquad (ab)v_i = a(bv_i)$$

$$(4.7.3) \qquad a(v_i + w_i) = av_i + aw_i$$

$$(4.7.4) \qquad (a+b)v_i = av_i + bv_i$$

$$(4.7.5) \qquad 1v_i = v_i$$

for any $a$, $b \in D$, $v_i$, $w_i \in V_i$, $i \in I$. Equalities

$$
\begin{aligned}
(4.7.6) \qquad (ab)v &= (ab)(v_i, i \in I) = ((ab)v_i, i \in I) \\
&= (a(bv_i), i \in I) = a(bv_i, i \in I) = a(b(v_i, i \in I)) \\
&= a(bv)
\end{aligned}
$$

$$
\begin{aligned}
(4.7.7) \qquad a(v+w) &= a((v_i, i \in I) + (w_i, i \in I)) = a(v_i + w_i, i \in I) \\
&= (a(v_i + w_i), i \in I) = (av_i + aw_i, i \in I) \\
&= (av_i, i \in I) + (aw_i, i \in I) = a(v_i, i \in I) + a(w_i, i \in I) \\
&= av + aw
\end{aligned}
$$

$$
\begin{aligned}
(4.7.8) \qquad (a+b)v &= (a+b)(v_i, i \in I) = ((a+b)v_i, i \in I) \\
&= (av_i + bv_i, i \in I) = (av_i, i \in I) + (bv_i, i \in I) \\
&= a(v_i, i \in I) + b(v_i, i \in I) = av + bv
\end{aligned}
$$

$$1v = 1(v_i, i \in I) = (1v_i, i \in I) = (v_i, i \in I) = v$$

---

[4.10] See also the definition of direct sum of modules in [1], page 128. On the same page, Lang proves the existence of direct sum of modules.

[4.11] See also proposition on the page [1]-1.1

follow from equalities   (3.7.4), (4.7.1), (4.7.2), (4.7.3), (4.7.4), (4.7.5)  for any   $a$,
$b \in D$, $v = (v_i, i \in I)$, $w = (w_i, i \in I) \in V$.

From the equality  (4.7.7), it follows that the map

$$a : v \to av$$

defined by the equality  (4.7.1) is endomorphism of Abelian group $V$.  From equalities
(4.7.6), (4.7.8), it follows that the map  (4.7.1) is homomorphism of ring $D$ into set
of endomorphisms of Abelian group $V$.  Therefore, the map  (4.7.1) is representation
of ring $D$ in Abelian group $V$.

Let $a_1$, $a_2$ be *D*-numbers such that

$$(4.7.9) \qquad\qquad a_1(v_i, i \in I) = a_2(v_i, i \in I)$$

for any $V$-number   $v = (v_i, i \in I)$.  The equality

$$(4.7.10) \qquad\qquad (a_1 v_i, i \in I) = (a_2 v_i, i \in I)$$

follows from the equality  (4.7.9).  From the equality  (4.7.10), it follows that

$$(4.7.11) \qquad\qquad \forall i \in I : a_1 v_i = a_2 v_i$$

Since for any $i \in I$, $V_i$ is $D$-module, then corresponding representation is effective
according to the definition   4.1.3.    Therefore, the equality $a_1 = a_2$ follows from
the statement  (4.7.11) and the map  (4.7.1) is effective representation of the ring $D$
in Abelian group $V$.  Therefore, Abelian group $V$ is $D$-module.

Let

$$\{f_i : V_i \to W, i \in I\}$$

be set of linear maps into $D$-module $W$.  We define the map

$$f : V \to W$$

by the equality

$$(4.7.12) \qquad\qquad f\left(\bigoplus_{i \in I} x_i\right) = \sum_{i \in I} f_i(x_i)$$

The sum in the right side of the equality  (4.7.12) is finite, since all summands,
except for a finite number, equal 0.  From the equality

$$f_i(x_i + y_i) = f_i(x_i) + f_i(y_i)$$

and the equality  (4.7.12), it follows that

$$\begin{aligned}
f(\bigoplus_{i \in I}(x_i + y_i)) &= \sum_{i \in I} f_i(x_i + y_i) \\
&= \sum_{i \in I} (f_i(x_i) + f_i(y_i)) \\
&= \sum_{i \in I} f_i(x_i) + \sum_{i \in I} f_i(y_i) \\
&= f(\bigoplus_{i \in I} x_i) + f(\bigoplus_{i \in I} y_i)
\end{aligned}$$

From the equality

$$f_i(dx_i) = df_i(x_i)$$

and the equality (4.7.12), it follows that

$$f((dx_i, i \in I)) = \sum_{i \in I} f_i(dx_i) = \sum_{i \in I} df_i(x_i) = d \sum_{i \in I} f_i(x_i)$$
$$= df((x_i, i \in I))$$

Therefore, the map $f$ is linear map. The equality

$$f \circ \lambda_j(x) = \sum_{i \in I} f_i(\delta^i_j x) = f_j(x)$$

follows from equalities (3.7.5), (4.7.12). Since the map $\lambda_i$ is injective, then the map $f$ is unique. Therefore, the theorem follows from definitions 2.2.12, 3.7.2. $\square$

THEOREM 4.7.3. *Direct sum of $D$-modules $V_1$, ..., $V_n$ coincides with their Cartesian product*

$$V_1 \oplus ... \oplus V_n = V_1 \times ... \times V_n$$

PROOF. The theorem follows from the theorem 4.7.2. $\square$

THEOREM 4.7.4. *Let $\overline{\overline{e}}_1 = (e_{1i}, i \in I)$ be a basis of $D$-module of columns $V$. Then we can represent $D$-module $V$ as direct sum*

$$(4.7.13) \qquad V = \bigoplus_{i \in I} De_i$$

PROOF. The map

$$f : v^i e_i \in V \to (v^i, i \in I) \in \bigoplus_{i \in I} De_i$$

is isomorphism. $\square$

THEOREM 4.7.5. *Let $\overline{\overline{e}}_1 = (e_1^i, i \in I)$ be a basis of $D$-module of rows $V$. Then we can represent $D$-module $V$ as direct sum*

$$(4.7.14) \qquad V = \bigoplus_{i \in I} De^i$$

PROOF. The map

$$f : v_i e^i \in V \to (v_i, i \in I) \in \bigoplus_{i \in I} De^i$$

is isomorphism. $\square$

THEOREM 4.7.6. *Let $V^1$, ..., $V^n$ be $D$-modules and*

$$V = V^1 \oplus ... \oplus V^n$$

*Let us represent a-number*

$$a = a^1 \oplus ... \oplus a^n$$

*as column vector*

$$a = \begin{pmatrix} a^1 \\ ... \\ a^n \end{pmatrix}$$

*Let us represent a linear map*

$$f : V \to W$$

*as row vector*

$$f = \begin{pmatrix} f_1 & ... & f_n \end{pmatrix}$$
$$f_i : V^i \to W$$

*Then we can represent value of the map $f$ in $V$-number $a$ as product of matrices*

$$(4.7.15) \qquad f \circ a = \begin{pmatrix} f_1 & \dots & f_n \end{pmatrix} \circ \begin{pmatrix} a^1 \\ \dots \\ a^n \end{pmatrix} = f_i \circ a^i$$

PROOF. The theorem follows from the definition (4.7.12). □

THEOREM 4.7.7. *Let $W^1$, ..., $W^m$ be D-modules and*

$$W = W^1 \oplus \dots \oplus W^m$$

*Let us represent $W$-number*

$$b = b^1 \oplus \dots \oplus b^m$$

*as column vector*

$$b = \begin{pmatrix} b^1 \\ \dots \\ b^m \end{pmatrix}$$

*Then the linear map*

$$f : V \to W$$

*has representation as column vector of maps*

$$f = \begin{pmatrix} f^1 \\ \dots \\ f^m \end{pmatrix}$$

*such way that, if $b = f \circ a$, then*

$$\begin{pmatrix} b^1 \\ \dots \\ b^m \end{pmatrix} = \begin{pmatrix} f^1 \\ \dots \\ f^m \end{pmatrix} \circ a = \begin{pmatrix} f^1 \circ a \\ \dots \\ f^m \circ a \end{pmatrix}$$

PROOF. The theorem follows from the theorem 2.2.8. □

THEOREM 4.7.8. *Let $V^1$, ..., $V^n$, $W^1$, ..., $W^m$ be D-modules and*

$$V = V^1 \oplus \dots \oplus V^n$$
$$W = W^1 \oplus \dots \oplus W^m$$

*Let us represent $V$-number*

$$a = a^1 \oplus \dots \oplus a^n$$

*as column vector*

$$(4.7.16) \qquad a = \begin{pmatrix} a^1 \\ \dots \\ a^n \end{pmatrix}$$

*Let us represent W-number*

$$b = b^1 \oplus ... \oplus b^m$$

*as column vector*

$$(4.7.17) \qquad b = \begin{pmatrix} b^1 \\ ... \\ b^m \end{pmatrix}$$

*Then the linear map*

$$f : V \to W$$

*has representation as a matrix of maps*

$$(4.7.18) \qquad f = \begin{pmatrix} f_1^1 & ... & f_n^1 \\ ... & ... & ... \\ f_1^m & ... & f_n^m \end{pmatrix}$$

*such way that, if $b = f \circ a$, then*

$$(4.7.19) \qquad \begin{pmatrix} b^1 \\ ... \\ b^m \end{pmatrix} = \begin{pmatrix} f_1^1 & ... & f_n^1 \\ ... & ... & ... \\ f_1^m & ... & f_n^m \end{pmatrix} \circ \begin{pmatrix} a^1 \\ ... \\ a^n \end{pmatrix} = \begin{pmatrix} f_1^1 \circ a^i \\ ... \\ f_i^m \circ a^i \end{pmatrix}$$

*The map*

$$f_j^i : V^j \to W^i$$

*is a linear map and is called* **partial linear map**.

PROOF. According to the theorem 4.7.7, there exists the set of linear maps

$$f^i : V \to W^i$$

such that

$$(4.7.20) \qquad \begin{pmatrix} b^1 \\ ... \\ b^m \end{pmatrix} = \begin{pmatrix} f^1 \\ ... \\ f^m \end{pmatrix} \circ a = \begin{pmatrix} f^1 \circ a \\ ... \\ f^m \circ a \end{pmatrix}$$

According to the theorem 4.7.6, for every $i$, there exists the set of linear maps

$$f_j^i : V^j \to W^i$$

such that

$$(4.7.21) \qquad f^i \circ a = \begin{pmatrix} f_1^i & ... & f_n^i \end{pmatrix} \circ \begin{pmatrix} a^1 \\ ... \\ a^n \end{pmatrix} = f_j^i \circ a^j$$

If we identify matrices

$$\begin{pmatrix} \begin{pmatrix} f_1^1 & ... & f_n^1 \end{pmatrix} \\ ... \\ \begin{pmatrix} f_1^m & ... & f_n^m \end{pmatrix} \end{pmatrix} = \begin{pmatrix} f_1^1 & ... & f_n^1 \\ ... & ... & ... \\ f_1^m & ... & f_n^m \end{pmatrix}$$

then the equality (4.7.19) follows from equalities (4.7.20), (4.7.21). $\square$

Let $B^1$, ..., $B^m$ be $D$-algebras. Then we can represent linear map $f_j^i$ using $B^i \otimes B^i$-number.

THEOREM 4.7.9. *Let $D$-module $V$ be free additive group. Let $U$ be submodule of $D$-module $V$. There exists submodule $W$ such that*

$$V = U \oplus W$$

PROOF. According to definitions 4.1.3, 4.1.6, $U$ is subgroup of additive group $V$. According to theorems 3.7.9, 3.7.10, there exists subgroup $W$ of additive group $V$ such that group $W$ is isomorphic to group $V/U$ and

$$V = U \oplus W$$

According to the theorem 4.7.2, action of $D$-number $a$ on any $V$-number

$$v = u \oplus w$$

has the form

$$av = (au) \oplus (aw)$$

Therefore, $W$ is submodule of $D$-module $V$. $\square$

## 4.8. Kernel of Linear Map

THEOREM 4.8.1. *Let the map*

$$(4.3.12) \quad (h : D_1 \to D_2 \quad f : V_1 \to V_2)$$

*be linear map of $D_1$-module $V_1$ into $D_2$-module $V_2$. The set*

$$(4.8.1) \qquad \ker(h, f) = \{v \in V : f \circ v = 0\}$$

*is called* **kernel of linear map** $(h, f)$. *Kernel of linear map $(h, f)$ is submodule of the module $V_1$.*

PROOF. Let $v, w \in V_1$. Then

$$(4.8.2) \qquad f \circ v = 0 \quad f \circ w = 0$$

The equality

$$f \circ (v + w) = f \circ v + f \circ w = 0$$

follows from equalities (4.3.4), (4.8.2). Therefore, the set $\ker(h, f)$ is subgroup of additive group $V_1$. The equality

$$f \circ (pv) = h(p)(f \circ v) = 0$$

follows from equalities (4.3.5), (4.8.2). Therefore, we have representation of ring $D_1$ on the set $\ker(h, f)$. According to the definition 4.1.6, the set $\ker(h, f)$ is submodule of $D_1$-module $V_1$. $\square$

DEFINITION 4.8.2. *Let* $V_i$, $i = 1, ..., n$, *be* $D_i$*-modules. A sequence of homomorphisms*

$$V_1 \xrightarrow{f_1} V_2 \xrightarrow{f_2} ... \xrightarrow{f_{n-1}} V_n$$

$$D_1 \xrightarrow{h_1} D_2 \xrightarrow{h_2} ... \xrightarrow{h_{n-1}} D_n$$

*is called* **exact** *if* $\forall i$, $i = 1, ..., n-2$, $\operatorname{Im} f_i = \ker(h_{i+1}, f_{i+1})$.  □

THEOREM 4.8.3. *Let the map*

$$(4.3.40) \quad f : V_1 \to V_2$$

*be linear map of* $D$*-module* $V_1$ *into* $D$*-module* $V_2$. *The set*

$$(4.8.3) \qquad \ker f = \{v \in V : f \circ v = 0\}$$

*is called* **kernel of linear map** $f$. *Kernel of linear map* $f$ *is submodule of the module* $V_1$.

PROOF. Let $v, w \in V_1$. Then

$$(4.8.4) \qquad f \circ v = 0 \quad f \circ w = 0$$

The equality

$$f \circ (v + w) = f \circ v + f \circ w = 0$$

follows from equalities (4.3.32), (4.8.4). Therefore, the set $\ker f$ is subgroup of additive group $V_1$. The equality

$$f \circ (pv) = p(f \circ v) = 0$$

follows from equalities (4.3.5), (4.8.4). Therefore, we have representation of ring $D$ on the set $\ker f$. According to the definition 4.1.6, the set $\ker f$ is submodule of $D$-module $V_1$.  □

DEFINITION 4.8.4. *Let* $V_i$, $i = 1, ..., n$, *be* $D$*-modules. A sequence of homomorphisms*

$$V_1 \xrightarrow{f_1} V_2 \xrightarrow{f_2} ... \xrightarrow{f_{n-1}} V_n$$

*is called* **exact** *if* $\forall i$, $i = 1, ..., n-2$, $\operatorname{Im} f_i = \ker f_{i+1}$.  □

THEOREM 4.8.5. *Let* $D$*-module* $V_1$ *be free additive group. Let* $U$ *be submodule of* $D$*-module* $V_1$. *There exist homomorphism*

$$(4.3.40) \quad f : V_1 \to V_2$$

*such that*

$$\ker f = U$$

PROOF. According to the theorem 4.7.9, there exists submodule $W_1$ such that

$$V_1 = U \oplus W$$

The map

$$f : u \oplus w \in U \oplus W \to w \in W$$

is homomotphism and

$$\ker f = U$$

$\square$

# $D$-Algebra

## 5.1. Algebra over Commutative Ring

DEFINITION 5.1.1. *Let $D$ be commutative ring. $D$-module $A$ is called **algebra over ring** $D$ or $D$-**algebra**, if we defined product*[5.1]*in $A$*

$$(5.1.1) \qquad v\,w = C \circ (v, w)$$

*where $C$ is bilinear map*

$$C : A \times A \to A$$

*If $A$ is free $D$-module, then $A$ is called **free algebra over ring** $D$.* □

THEOREM 5.1.2. *The multiplication in the algebra $A$ is distributive over addition*

$$(5.1.2) \qquad (a + b)c = ac + bc$$

$$(5.1.3) \qquad a(b + c) = ab + ac$$

PROOF. The statement of the theorem follows from the chain of equations

$$(a + b)c = f \circ (a + b, c) = f \circ (a, c) + f \circ (b, c) = ac + bc$$
$$a(b + c) = f \circ (a, b + c) = f \circ (a, b) + f \circ (a, c) = ab + ac$$

□

DEFINITION 5.1.3. *If the product in $D$-algebra $A$ has unit element, then $D$-algebra $A$ is called **unital algebra**[5.2]* □

The multiplication in algebra can be neither commutative nor associative. Following definitions are based on definitions given in [23], p. 13.

DEFINITION 5.1.4. *The **commutator***

$$[a, b] = ab - ba$$

*measures commutativity in $D$-algebra $A$. $D$-algebra $A$ is called **commutative**, if*

$$[a, b] = 0$$

□

DEFINITION 5.1.5. *The **associator***

$$(5.1.4) \qquad (a, b, c) = (ab)c - a(bc)$$

*measures associativity in $D$-algebra $A$. $D$-algebra $A$ is called **associative**, if*

$$(a, b, c) = 0$$

---

[5.1] I follow the definition given in [23], page 1, [17], page 4. The statement which is true for any $D$-module, is true also for $D$-algebra.

[5.2] See the definition of unital $D$-algebra also on the pages [7]-137.

$\square$

Theorem 5.1.6. *Let $A$ be algebra over commutative ring $D$.* [5.3]

(5.1.5)        $a(b,c,d) + (a,b,c)d = (ab,c,d) - (a,bc,d) + (a,b,cd)$

*for any $a$, $b$, $c$, $d \in A$.*

Proof. The equation (5.1.5) follows from the chain of equations

$$a(b,c,d) + (a,b,c)d = a((bc)d - b(cd)) + ((ab)c - a(bc))d$$
$$= a((bc)d) - a(b(cd)) + ((ab)c)d - (a(bc))d$$
$$= ((ab)c)d - (ab)(cd) + (ab)(cd)$$
$$+ a((bc)d) - a(b(cd)) - (a(bc))d$$
$$= (ab,c,d) - (a(bc))d + a((bc)d) + (ab)(cd) - a(b(cd))$$
$$= (ab,c,d) - (a,(bc),d) + (a,b,cd)$$

$\square$

Theorem 5.1.7. *Let $\overline{\overline{e}}$ be basis of $D$-algebra $A$. Then*

(5.1.6)        $(a,b,c) = A^k_{ijl} a^i b^j c^l e_k$

*where* **coordinates of associator** *are defined by the equality*

(5.1.7)        $A^k_{ijl} = C^k_{pl} C^p_{ij} - C^k_{ip} C^p_{jl}$

Proof. The equality

(5.1.8)        $(a,b,c) = (C^k_{pl}(ab)^p c^l - C^k_{ip} a^i (bc)^p) e_k$
$$= (C^k_{pl} C^p_{ij} - C^k_{ip} C^p_{jl}) a^i b^j c^l e_k$$

follows from the equality (5.1.4). The equality (5.1.7) follows from the equality (5.1.8). $\square$

Definition 5.1.8. *The set* [5.4]

$$N(A) = \{a \in A : \forall b, c \in A, (a,b,c) = (b,a,c) = (b,c,a) = 0\}$$

*is called the* **nucleus of an** $D$**-algebra** $A$. $\square$

Definition 5.1.9. *The set* [5.5]

$$Z(A) = \{a \in A : a \in N(A), \forall b \in A, ab = ba\}$$

*is called the* **center of an** $D$**-algebra** $A_1$. $\square$

---

[5.3] The statement of the theorem is based on the equation [23]-(2.4).
[5.4] The definition is based on the similar definition in [23], p. 13.
[5.5] The definition is based on the similar definition in [23], page 14.
[5.6] See also the definition on the page [1]-86.

DEFINITION 5.1.10. *A* **left ideal** [5.6] $A_1$ *in a D-algebra A is a subgroup of the additive group of A such that*

$$AA_1 \subseteq A_1$$

□

DEFINITION 5.1.12. *A* **right ideal** [5.6] $A_1$ *in a D-algebra A is a subgroup of the additive group of A such that*

$$A_1A \subseteq A_1$$

□

THEOREM 5.1.11. *If $A_1$ is a left ideal in unital D-algebra A has unit, then*

$$AA_1 = A_1$$

PROOF. The theorem follows from the definition 5.1.10 and equality

$$eA_1 = A_1$$

□

THEOREM 5.1.13. *If $A_1$ is a right ideal in unital D-algebra A has unit, then*

$$A_1A = A_1$$

PROOF. The theorem follows from the definition 5.1.12 and equality

$$A_1e = A_1$$

□

DEFINITION 5.1.14. **Ideal** [5.6] *is a subgroup of the additive group of A which is both left and right ideal.* □

If $A$ is $D$-algebra, then the set $\{0\}$ and the set $A$ are ideals which we call trivial ideals.

THEOREM 5.1.15. *Let D be commutative ring. Let A be associative D-algebra with unit $1_A$.*

5.1.15.1: *The map*

$$f : d \in D \to d1_A \in A$$

*is homomorphism of the ring D into the center of the algebra A and allows to identify D-number d and A-number $d\,1_A$.*

5.1.15.2: *The ring D is subalgebra of algebra A.*

5.1.15.3: $D \subset Z(A)$.

PROOF. Let $d_1, d_2 \in D$. According to the statement 4.1.8.3,

$$(5.1.9) \qquad d_1 1_A + d_2 1_A = (d_1 + d_2) 1_A$$

According to the definition 5.1.1,

$$(5.1.10) \qquad (d_1 1_A)(d_2 1_A) = d_1(1_A(d_2 1_A)) = d_1(d_2(1_A 1_A))$$
$$= d_1(d_2 1_A) = (d_1 d_2) 1_A$$

The statement 5.1.15.1 follows from equalities (5.1.9), (5.1.10). The statement 5.1.15.2 follows from the statement 5.1.15.1.

According to the definition 5.1.1,

$$(5.1.11) \qquad a(d1_A) = d(a1_A) = d(1_A a) = (d1_A)a$$

The statement 5.1.15.3 follows from the equation (5.1.11). □

LEMMA 5.1.16. *Let* $d \in D$. *The map*

(5.1.12)        $d : b \in A \rightarrow db \in A$

*is endomorphism of Abelian group* $A$.

PROOF. The statement

$$db \in A$$

follows from definitions 4.1.3, 5.1.1. According to the equality (4.1.22) the map (5.1.12) is endomorphism of Abelian group $A$.                    □

LEMMA 5.1.17. *Let* $a \in A$. *The map*

(5.1.13)        $a : b \in A \rightarrow ab \in A$

*is endomorphism of Abelian group* $A$.

PROOF. The statement

$$ab \in A$$

follows from the definition 5.1.1. According to the theorem 5.1.2 the map (5.1.13) is endomorphism of Abelian group $A$.                    □

LEMMA 5.1.20. *Let* $d \in D$, $a \in A$. *The map*

(5.1.16)    $a \oplus d : b \in A \rightarrow (a \oplus d)b \in A$

*defined by the equality*

(5.1.17)        $(a \oplus d)b = ab + db$

*is endomorphism of Abelian group* $A$.

PROOF. Since $A$ is Abelian group, then the statement

$$db + ab \in A \quad d \in D \quad a \in A$$

follows from lemmas 5.1.16, 5.1.17. Therefore, for any $D$-number $d$ and for any $D$-number $a$, we defined the map (5.1.16). According to lemmas 5.1.16, 5.1.17 and the theorem 3.3.4, the map (5.1.16) is endomorphism of Abelian group $V$.    □

LEMMA 5.1.18. *Let* $d \in D$. *The map*

(5.1.14)        $d : b \in A \rightarrow bd \in A$

*is endomorphism of Abelian group* $A$.

PROOF. The statement

$$bd \in A$$

follows from definitions 4.1.3, 5.1.1. According to the equality (4.1.22) the map (5.1.14) is endomorphism of Abelian group $A$.                    □

LEMMA 5.1.19. *Let* $a \in A$. *The map*

(5.1.15)        $a : b \in A \rightarrow ba \in A$

*is endomorphism of Abelian group* $A$.

PROOF. The statement

$$ba \in A$$

follows from the definition 5.1.1. According to the theorem 5.1.2 the map (5.1.15) is endomorphism of Abelian group $A$.                    □

LEMMA 5.1.21. *Let* $d \in D$, $a \in A$. *The map*

(5.1.18)    $a \oplus d : b \in A \rightarrow b(a \oplus d) \in A$

*defined by the equality*

(5.1.19)        $b(a \oplus d) = ba + bd$

*is endomorphism of Abelian group* $A$.

PROOF. Since $A$ is Abelian group, then the statement

$$bd + ba \in A \quad d \in D \quad a \in A$$

follows from lemmas 5.1.18, 5.1.19. Therefore, for any $D$-number $d$ and for any $D$-number $a$, we defined the map (5.1.18). According to lemmas 5.1.18, 5.1.19 and the theorem 3.3.4, the map (5.1.18) is endomorphism of Abelian group $V$.    □

THEOREM 5.1.22. *Let* $A$ *be associative* $D$-*algebra.    The set of maps* $a \oplus d$

*generates* [5.7] *$D$-algebra $A_{(1)}$ where the sum is defined by the equality*

$$(5.1.20) \qquad (a_1 \oplus d_1) + (a_2 \oplus d_2) = (a_1 + a_2) \oplus (d_1 + d_2)$$

*and the product is defined by the equality*

$$(5.1.21) \qquad (a_1 \oplus d_1) \circ (a_2 \oplus d_2) = (a_1 a_2 + d_2 a_1 + d_1 a_2) \oplus (d_1 d_2)$$

5.1.22.1: *If $D$-algebra $A$ has unit, then $D \subseteq A$, $A_{(1)} = A$.*

5.1.22.2: *If $D$-algebra $D$ is ideal of $D$, then $A \subseteq D$, $A_{(1)} = D$.*

5.1.22.3: *Otherwise $A_{(1)} = A \oplus D$.*

5.1.22.4: *The $D$-algebra $A$ is ideal of $D$-algebra $A_{(1)}$.*

*The $D$-algebra $A_{(1)}$ is called* **unital extension** *of the $D$-algebra $D$.*

PROOF. Expression for operations on the set $A_{(1)}$ does not depend on whether the endomorphism $a \oplus d$ acts on $A$-numbers on the left or on the right. However, the text of the proof is slightly different. Therefore, we will consider both versions of the proof.

According to the theorem 3.3.4, the equality

$$(5.1.22) \quad \begin{aligned} &((a_1 \oplus d_1) + (a_2 \oplus d_2))b \\ =&(a_1 \oplus d_1)b + (a_2 \oplus d_2)b \\ =&a_1 b + d_1 b + a_2 b + d_2 b \\ =&a_1 b + a_2 b + d_1 b + d_2 b \\ =&(a_1 + a_2)b + (d_1 + d_2)b \\ =&(a_1 + a_2) \oplus (d_1 + d_2)b \end{aligned}$$

follows from the equality (5.1.17). The equality (5.1.20) follows from the equality (5.1.22).

According to the theorem 2.1.10, the equality

$$(5.1.23) \quad \begin{aligned} &((a_1 \oplus d_1) \circ (a_2 \oplus d_2))b \\ =&(a_1 \oplus d_1)((a_2 \oplus d_2)b) \\ =&(a_1 \oplus d_1)(a_2 b + d_2 b) \\ =&(a_1 \oplus d_1)(a_2 b) \\ +&(a_1 \oplus d_1)(d_2 b) \\ =&a_1(a_2 b) + d_1(a_2 b) \\ +&a_1(d_2 b) + d_1(d_2 b) \end{aligned}$$

follows from the equality (5.1.17). The equality

$$(5.1.24) \quad \begin{aligned} &((a_1 \oplus d_1) \circ (a_2 \oplus d_2))b \\ =&(a_1 a_2)b + (d_1 a_2)b \\ +&(d_2 a_1)b + (d_1 d_2)b \end{aligned}$$

According to the theorem 3.3.4, the equality

$$(5.1.27) \quad \begin{aligned} &b((a_1 \oplus d_1) + (a_2 \oplus d_2)) \\ =&b(a_1 \oplus d_1) + b(a_2 \oplus d_2) \\ =&ba_1 + bd_1 + ba_2 + bd_2 \\ =&ba_1 + ba_2 + bd_1 + bd_2 \\ =&b(a_1 + a_2) + b(d_1 + d_2) \\ =&b(a_1 + a_2) \oplus (d_1 + d_2) \end{aligned}$$

follows from the equality (5.1.19). The equality (5.1.20) follows from the equality (5.1.27).

According to the theorem 2.1.10, the equality

$$(5.1.28) \quad \begin{aligned} &b((a_2 \oplus d_2) \circ (a_1 \oplus d_1)) \\ =&(b(a_2 \oplus d_2))(a_1 \oplus d_1) \\ =&(ba_2 + bd_2)(a_1 \oplus d_1) \\ =&(ba_2)(a_1 \oplus d_1) \\ +&(bd_2)(a_1 \oplus d_1) \\ =&(ba_2)a_1 + (ba_2)d_1 \\ +&(bd_2)a_1 + (bd_2)d_1 \end{aligned}$$

follows from the equality (5.1.19). The equality

$$(5.1.29) \quad \begin{aligned} &b((a_2 \oplus d_2) \circ (a_1 \oplus d_1)) \\ =&b(a_2 a_1) + b(a_2 d_1) \\ +&b(a_1 d_2) + b(d_2 d_1) \end{aligned}$$

---

[5.7] See the definition of unital extension also on the pages [7]-52, [8]-64.

follows from   from definitions 5.1.1, 5.1.5, from the equality  (5.1.23),   from the lemma 5.1.17. The equality

$$((a_1 \oplus d_1) \circ (a_2 \oplus d_2))b$$
$$=(a_1 a_2 + d_1 a_2 + d_2 a_1 + d_1 d_2)b$$
$$=((a_1 a_2 + d_1 a_2 + d_2 a_1) \oplus d_1 d_2)b$$

and therefore equality (5.1.21), follows from equalities (5.1.17), (5.1.24) and from lemmas 5.1.16, 5.1.17.

Let $D$-algebra $A$ have unit $e$. Then we can identify $d \in D$ and $de \in A$. The equality

$$(5.1.25) \quad (a \oplus d)b = ab + db = (a + de)b$$

follows from the equality

$$\boxed{(5.1.17) \quad (a \oplus d)b = ab + db}$$

The statement 5.1.22.1 follows from the equality (5.1.25).

Let $D$-algebra $A$ be ideal of $D$. Then $A \subseteq D$. The equality

$$(5.1.26) \quad (a \oplus d)b = ab + db = (a + d)b$$

follows from the equality (5.1.17). The statement 5.1.22.2 follows from the equality (5.1.26).

The statement 5.1.22.4 follows from the definition 5.1.14 and the lemmas 5.1.20, 5.1.21.                                    □

follows from   from definitions 5.1.1, 5.1.5, from the equality  (5.1.28),   from the lemma 5.1.19. The equality

$$b((a_2 \oplus d_2) \circ (a_1 \oplus d_1))$$
$$=b(a_2 a_1 + a_2 d_1 + a_1 d_2 + d_2 d_1)$$
$$=b((a_2 a_1 + a_2 d_1 + a_1 d_2) \oplus d_2 d_1)$$

and therefore equality (5.1.21), follows from equalities (5.1.19), (5.1.29) and from lemmas 5.1.18, 5.1.19.

Let $D$-algebra $A$ have unit $e$. Then we can identify $d \in D$ and $ed \in A$. The equality

$$(5.1.30) \quad b(a \oplus d) = ba + bd = b(a + ed)$$

follows from the equality

$$\boxed{(5.1.19) \quad b(a \oplus d) = ba + bd}$$

The statement 5.1.22.1 follows from the equality (5.1.30).

Let $D$-algebra $A$ be ideal of $D$. Then $A \subseteq D$. The equality

$$(5.1.31) \quad b(a \oplus d) = ba + bd = b(a + d)$$

follows from the equality (5.1.19). The statement 5.1.22.2 follows from the equality (5.1.31).

THEOREM 5.1.23. *Let $A$ be non-commutative $D$-algebra. For any $b \in A$, there exists subalgebra $Z(A, b)$ of $D$-algebra $A$ such that*

$$(5.1.32) \qquad\qquad c \in Z(A, b) \Leftrightarrow cb = bc$$

*$D$-algebra $Z(A, b)$ is called* **center** *of $A$-**number** $b$.*

PROOF. The subset $Z(A, b)$ is not empty because $0 \in Z(A, b)$, $b \in Z(A, b)$. Let $d_1$, $d_2 \in D$, $c_1$, $c_2 \in Z(A, b)$. Then

$$(d_1 c_1 + d_2 c_2)b = d_1 c_1 b + d_2 c_2 b = b d_1 c_1 + b d_2 c_2$$
$$= b(d_1 c_1 + d_2 c_2)$$

And this implies $d_1 c_1 + d_2 c_2 \in Z(A, b)$. Therefore, the set $Z(A, b)$ is $D$-module.

Let $c_1$, $c_2 \in Z(A, b)$. Then

$$(5.1.33) \qquad b(c_1 c_2) = (b c_1)c_2 = (c_1 b)c_2 = c_1(b c_2) = c_1(c_2 b) = (c_1 c_2)b$$

From the equality (5.1.33), it follows that $c_1 c_2 \in Z(A, b)$. Therefore, $D$-module $Z(A, b)$ is $D$-algebra.                                    □

LEMMA 5.1.24. *Let $A$ be non-commutative $D$-algebra. For any $a \in A$, if $p$, $c \in Z(A, a)$, then*

$$(5.1.34) \qquad pc^n \in Z(A, a)$$

PROOF. According to the theorem 5.1.23,

$$p = pc^0 \in Z(A, a)$$

Let the statement (5.1.34) is true for $n = k - 1$, $k > 0$,

$$(5.1.35) \qquad pc^{k-1} \in Z(A, a)$$

Then the statement

$$pc^k = pc^{k-1}c \in Z(A, a)$$

follows from the statement (5.1.35) and the theorem 5.1.23. According to mathematical induction, the statement (5.1.34) is true for any $n \geq 0$. $\qquad\square$

THEOREM 5.1.25. *Let $A$ be non-commutative $D$-algebra. For any $a \in A$, if $c \in Z(A, a)$, then*

$$(5.1.36) \qquad p(c) \in Z(A, a)$$

*for any polynomial*

$$(5.1.37) \qquad p(x) = p_0 + p_1 x + ... + p_n x^n$$
$$p_0, \ ..., \ p_n \in D$$

PROOF. The theorem follows from the lemma 5.1.24 and the theorem 5.1.23. $\qquad\square$

THEOREM 5.1.26. *Since $a \in Z(A, b)$, then $b \in Z(A, a)$.*

PROOF. Let $a \in Z(A, b)$. The equality

$$(5.1.38) \qquad ab = ba$$

follows from the statement (5.1.32). The theorem follows from the equality (5.1.38). $\qquad\square$

THEOREM 5.1.27. *Since $a \in Z(A, b)$ and $c \in Z(A)$, then $a \in Z(A, bc)$.*

PROOF. The equality

$$(5.1.39) \qquad a(bc) = (ab)c = (ba)c = b(ac) = b(ca) = (bc)a$$

follows from definitions 5.1.8, 5.1.9 and from the theorem 5.1.23. The theorem follows from the equality (5.1.39) and from the theorem 5.1.23. $\qquad\square$

DEFINITION 5.1.28. *$D$-algebra $A$ is called **division algebra**, if for any $A$-number $a \neq 0$ there exists $A$-number $a^{-1}$.* $\qquad\square$

THEOREM 5.1.29. *Let $D$-algebra $A$ be division algebra. $D$-algebra $A$ has unit.*

PROOF. The theorem follows from the statement that the equation

$$ax = a$$

has solution for any $a \in A$.                                                    □

THEOREM 5.1.30. *Let $D$-algebra $A$ be division algebra.   The ring $D$ is the field and subset of the center of $D$-algebra $A$.*

PROOF. Let $e$ be unit of $D$-algebra $A$. Since the map

$$d \in D \rightarrow d\, e \in A$$

is embedding of the ring $D$ into $D$-algebra $A$, then the ring $D$ is subset of the center of $D$-algebra $A$.                                                    □

THEOREM 5.1.31. *Let $A$ be associative division $D$-algebra.  The statement*

(5.1.40)                              $ab = ac \quad a \neq 0$

*implies $b = c$.*

PROOF. The equality

(5.1.41)          $b = (a^{-1}a)b = a^{-1}(ab) = a^{-1}(ac) = (a^{-1}a)c = c$

follows from the statement (5.1.41).                                               □

DEFINITION 5.1.32. *$A$-numbers $a$ and $b$ are called similar, if there exists $A$-number $c$ such that*

(5.1.42)                                 $a = c^{-1}bc$

                                                                                  □

THEOREM 5.1.33. *Let $A$ be algebra over commutative ring $D$. For any set $B$, the set of maps $A^B$ is Abelian group with respect to operation*

$$(f + g)(x) = f(x) + g(x)$$

PROOF. The theorem follows from definitions   4.1.3,  5.1.1.            □

## 5.2. $D$-Algebra Type

CONVENTION 5.2.1. *According to definitions 4.1.24, 5.1.1, free $D$-algebra $A$ is $D$-module $A$, which has a basis $\overline{\overline{e}}$.  However, generally speaking, $\overline{\overline{e}}$ is not a basis of $D$-algebra $A$, because there is product in $D$-algebra $A$. For instance, in quaternion algebra, any quaternion*

$$a = a^0 + a^1 i + a^2 j + a^3 k$$

*can be represented as*

$$a = a^0 + a^1 i + a^2 j + a^3 ij$$

*However, for most problems using a basis of $D$-module $A$ is easier than using a basis of $D$-algebra $A$.  For instance, it is easier to determine coordinates of $H$-number with respect to the basis $(1, i, j, k)$ of $R$-vector space $H$, than to determine coordinates of $H$-number with respect to the basis $(1, i, j)$ of $R$-algebra $H$.  Therefore*

*the phrase "we consider basis of D-algebra A" means that we consider D-algebra A and basis of D-module A.* □

CONVENTION 5.2.2. *Let A be free algebra with finite or countable basis. Considering expansion of element of algebra A relative basis $\overline{\overline{e}}$ we use the same root letter to denote this element and its coordinates. In expression $a^2$, it is not clear whether this is component of expansion of element a relative a basis, or this is operation $a^2 = aa$. To make text clearer we use separate color for index of element of algebra. For instance,*

$$a = a^i e_i$$

□

*D-algebra A inherits type of D-module A.*

DEFINITION 5.2.3. *D-algebra A is called D-algebra of columns, if D-module A is D-module of columns.* □

DEFINITION 5.2.5. *D-algebra A is called D-algebra of rows, if D-module A is D-module of rows.* □

CONVENTION 5.2.4. *Let $\overline{\overline{e}}$ be the basis of D-algebra of columns A. If D-algebra A has unit, then we assume that $e_0$ is the unit of D-algebra A.* □

CONVENTION 5.2.6. *Let $\overline{\overline{e}}$ be the basis of D-algebra of rows A. If D-algebra A has unit, then we assume that $e^0$ is the unit of D-algebra A.* □

THEOREM 5.2.7. *Let the set of vectors $\overline{\overline{e}} = (e_i, i \in I)$ be a basis of D-algebra of columns A. Let*

$$a = a^i e_i \quad b = b^i e_i \quad a, b \in A$$

*We can get the product of a, b according to rule*

$$(5.2.1) \qquad (ab)^k = C_{ij}^k a^i b^j$$

*where $C_{ij}^k$ are **structural constants** of algebra A over ring D. The product of basis vectors in the algebra A is defined according to rule*

$$(5.2.2) \qquad e_i e_j = C_{ij}^k e_k$$

PROOF. The equality (5.2.2) is corollary of the statement that $\overline{\overline{e}}$ is the basis of D-algebra A. Since the product in the algebra is a bilinear map, then we can write the product of a and b as

$$(5.2.3) \qquad ab = a^i b^j e_i e_j$$

THEOREM 5.2.8. *Let the set of vectors $\overline{\overline{e}} = (e^i, i \in I)$ be a basis of D-algebra of rows A. Let*

$$a = a_i e^i \quad b = b_i e^i \quad a, b \in A$$

*We can get the product of a, b according to rule*

$$(5.2.5) \qquad (ab)_k = C_k^{ij} a_i b_j$$

*where $C_k^{ij}$ are **structural constants** of algebra A over ring D. The product of basis vectors in the algebra A is defined according to rule*

$$(5.2.6) \qquad e^i e^j = C_k^{ij} e^k$$

PROOF. The equality (5.2.6) is corollary of the statement that $\overline{\overline{e}}$ is the basis of D-algebra A. Since the product in the algebra is a bilinear map, then we can write the product of a and b as

$$(5.2.7) \qquad ab = a_i b_j e^i e^j$$

From equalities (5.2.2), (5.2.3), it follows that

$$(5.2.4) \qquad ab = a^i b^j C_{ij}^k e_k$$

Since $\overline{\overline{e}}$ is a basis of the algebra $A$, then the equality (5.2.1) follows from the equality (5.2.4). $\qquad\square$

From equalities (5.2.6), (5.2.7), it follows that

$$(5.2.8) \qquad ab = a_i b_j C_k^{ij} e^k$$

Since $\overline{\overline{e}}$ is a basis of the algebra $A$, then the equality (5.2.5) follows from the equality (5.2.8). $\qquad\square$

**Theorem 5.2.9.** *Since the algebra $A$ is commutative, then*

$$(5.2.9) \qquad C_{ij}^p = C_{ji}^p$$

*Since the algebra $A$ is associative, then*

$$(5.2.10) \qquad C_{ij}^p C_{pk}^q = C_{ip}^q C_{jk}^p$$

**Theorem 5.2.10.** *Since the algebra $A$ is commutative, then*

$$(5.2.11) \qquad C_p^{ij} = C_p^{ji}$$

*Since the algebra $A$ is associative, then*

$$(5.2.12) \qquad C_p^{ij} C_q^{pk} = C_q^{ip} C_p^{jk}$$

**Proof.** For commutative algebra, the equation (5.2.9) follows from equation

$$e_i e_j = e_j e_i$$

For associative algebra, the equation (5.2.10) follows from equation

$$(e_i e_j) e_k = e_i (e_j e_k)$$

$\qquad\square$

**Proof.** For commutative algebra, the equation (5.2.11) follows from equation

$$e^i e^j = e^j e^i$$

For associative algebra, the equation (5.2.12) follows from equation

$$(e^i e^j) e^k = e^i (e^j e^k)$$

$\qquad\square$

## 5.3. Diagram of Representations

**Theorem 5.3.1.** *The representation*

$$(5.3.1) \qquad g_{23} : A \dashrightarrow A$$

*of $D$-module $A$ in $D$-module $A$ is equivalent to structure of $D$-algebra $A$.*

**Theorem 5.3.2.** *Let $D$ be commutative ring and $A$ be Abelian group. The diagram of representations*

$$g_{12}(d) : v \to d\, v$$
$$g_{23}(v) : w \to C \circ (v, w)$$
$$C \in \mathcal{L}(D; A^2 \to A)$$

*generates the structure of $D$-algebra $A$.*

**Proof.** The structure of $D$-module $A$ is generated by effective representation

$$g_{12} : D \dashrightarrow A$$

of ring $D$ in Abelian group $A$.

LEMMA 5.3.3. *Let the structure of $D$-algebra $A$ defined in $D$-module $A$, be generated by product*

$$v\,w = C \circ (v, w)$$

**Left shift of $D$-module** $A$ *defined by equation*

(5.3.2) $$l \circ v : w \in A \to v\,w \in A$$

*generates the representation*

$$A \xrightarrow{g_{23}} A \qquad \begin{array}{l} g_{23} : v \to l \circ v \\ g_{23} \circ v : w \to (l \circ v) \circ w \end{array}$$

*of $D$-module $A$ in $D$-module $A$*

PROOF. According to definitions 5.1.1 and 4.4.1, left shift of $D$-module $A$ is linear map. According to the definition 4.3.9, the map $l(v)$ is endomorphism of $D$-module $A$. The equation

(5.3.3) $$(l \circ (v_1 + v_2)) \circ w = (v_1 + v_2)w = v_1 w + v_2 w = (l \circ v_1) \circ w + (l \circ v_2) \circ w$$

follows from the definition 4.4.1 and from the equation (5.3.2).     According to the corollary 4.4.4, the equation

(5.3.4) $$l \circ (v_1 + v_2) = l \circ v_1 + l \circ v_2$$

follows from equation (5.3.3). The equation

(5.3.5) $$(l \circ (dv)) \circ w = (dv)w = d(vw) = d((l \circ v) \circ w)$$

follows from the definition 4.4.1 and from the equation (5.3.2).   According to the corollary 4.4.4, the equation

(5.3.6) $$l \circ (dv) = d(l \circ v)$$

follows from equation (5.3.5). The lemma follows from equalities (5.3.4), (5.3.6). ⊙

LEMMA 5.3.4. *The representation*

$$A \xrightarrow{g_{23}} A \qquad \begin{array}{l} g_{23} : v \to l \circ v \\ g_{23} \circ v : w \to (l \circ v) \circ w \end{array}$$

*of $D$-module $A$ in $D$-module $A$ determines the product in $D$-module $A$ according to rule*

$$ab = (g_{23} \circ a) \circ b$$

PROOF. Since map $g_{23} \circ v$ is endomorphism of $D$-module $A$, then

(5.3.7) $$\begin{array}{l} (g_{23} \circ v)(w_1 + w_2) = (g_{23} \circ v) \circ w_1 + (g_{23} \circ v) \circ w_2 \\ (g_{23} \circ v) \circ (dw) = d((g_{23} \circ v) \circ w) \end{array}$$

Since the map $g_{23}$ is linear map

$$g_{23} : A \to \mathcal{L}(D; A \to A)$$

then,   according to corollarys 4.4.4, 4.4.6,

(5.3.8) $$(g_{23} \circ (v_1 + v_2)) \circ w = (g_{23} \circ v_1 + g_{23} \circ v_2)(w) = (g_{23} \circ v_1) \circ w + (g_{23} \circ v_2) \circ w$$

(5.3.9)         $(g_{23} \circ (d\, v)) \circ w = (d\, (g_{23} \circ v)) \circ w = d\, ((g_{23} \circ v) \circ w)$

From equations (5.3.7), (5.3.8), (5.3.9) and the definition 4.4.1, it follows that the map $g_{23}$ is bilinear map. Therefore, the map $g_{23}$ determines the product in $D$-module $A$ according to rule

$$ab = (g_{23} \circ a) \circ b$$

$\odot$

The theorem follows from lemmas 5.3.3, 5.3.4.                                     $\square$

## 5.4. Homomorphism of $D$-Algebra

### 5.4.1. General Definition.

Definition 5.4.1. *Let* $A_i$, $i = 1, 2$, *be algebra over commutative ring* $D_i$. *Let diagram of representations*

(5.4.1)
$$\begin{array}{c} A_1 \xrightarrow{\;h_{1.23}\;} A_1 \\ {}_{h_{1.12}}\!\!\!\searrow \quad \nearrow_{h_{1.12}} \\ D_1 \end{array}$$
$\qquad h_{1.12}(d) : v \to d\, v$
$\qquad h_{1.23}(v) : w \to C_1 \circ (v, w)$
$\qquad C_1 \in \mathcal{L}(D_1; A_1^2 \to A_1)$

*describe* $D_1$-*algebra* $A_1$. *Let diagram of representations*

(5.4.2)
$$\begin{array}{c} A_2 \xrightarrow{\;h_{2.23}\;} A_2 \\ {}_{h_{2.12}}\!\!\!\searrow \quad \nearrow_{h_{2.12}} \\ D_2 \end{array}$$
$\qquad h_{2.12}(d) : v \to d\, v$
$\qquad h_{2.23}(v) : w \to C_2 \circ (v, w)$
$\qquad C_2 \in \mathcal{L}(D_2; A_2^2 \to A_2)$

*describe* $D_2$-*algebra* $A_2$. *Morphism*

(5.4.3)                     $(h : D_1 \to D_2 \quad f : A_1 \to A_2)$

*of diagram of representations* (5.4.1) *into diagram of representations* (5.4.2) *is called* **homomorphism** *of* $D_1$-*algebra* $A_1$ *into* $D_2$-*algebra* $A_2$. *Let us denote* $\mathrm{Hom}(D_1 \to D_2; A_1 \to A_2)$ *set of homomorphisms of* $D_1$-*algebra* $A_1$ *into* $D_2$-*algebra* $A_2$.                     $\square$

Theorem 5.4.2. *Homomorphism*

$$\boxed{(5.4.3) \quad (h : D_1 \to D_2 \quad f : A_1 \to A_2)}$$

*of* $D_1$-*algebra* $A_1$ *into* $D_2$-*algebra* $A_2$ *is linear map* (5.4.3) *of* $D_1$-*module* $A_1$ *into* $D_2$-*module* $A_2$ *such that*

(5.4.4)                         $f \circ (ab) = (f \circ a)(f \circ b)$

*and satisfies following equalities*

(5.4.5)                         $h(p + q) = h(p) + h(q)$

(5.4.6)                         $h(pq) = h(p)h(q)$

(5.4.7)                         $f \circ (a + b) = f \circ a + f \circ b$

$$\boxed{(5.4.4) \quad f \circ (ab) = (f \circ a)(f \circ b)}$$

(5.4.8)                         $f \circ (pa) = h(p)(f \circ a)$

$$p, q \in D_1 \quad a, b \in A_1$$

Proof. According to definitions 2.3.9, 2.8.5, the map (5.4.3) is morphism of representation $h_{1.12}$ describing $D_1$-module $A_1$ into representation $h_{2.12}$ describing $D_2$-module $A_2$. Therefore, the map $f$ is linear map of $D_1$-module $A_1$ into $D_2$-module $A_2$ and equalities (5.4.5), (5.4.6), (5.4.7), (5.4.8) follow from the theorem 4.3.2.

According to definitions 2.3.9, 2.8.5, the map $f$ is morphism of representation $h_{1.23}$ into representation $h_{2.23}$. The equality

$$(5.4.9) \qquad f \circ (h_{1.23}(a)(b)) = h_{2.23}(f \circ a)(f \circ b)$$

follows from the equality (2.3.2). From equalities (5.4.1), (5.4.2), (5.4.9), it follows that

$$(5.4.10) \qquad f \circ (C_1 \circ (a, b)) = C_2 \circ (f \circ a, f \circ b)$$

The equality (5.4.4) follows from equalities (5.1.1), (5.4.10). $\qquad\square$

| | |
|---|---|
| Theorem 5.4.3. *The homomorphism*[5.8] | Theorem 5.4.4. *The homomorphism*[5.9] |
| $$(h : D_1 \to D_2 \quad \overline{f} : A_1 \to A_2)$$ | $$(h : D_1 \to D_2 \quad \overline{f} : A_1 \to A_2)$$ |
| *of $D_1$-algebra of columns $A_1$ into $D_2$-algebra of columns $A_2$ has presentation* | *of $D_1$-algebra of rows $A_1$ into $D_2$-algebra of rows $A_2$ has presentation* |
| $$(5.4.11) \qquad b = f h(a)$$ | $$(5.4.22) \qquad b = h(a) f$$ |
| $$(5.4.12) \qquad \overline{f} \circ (a^i e_{1i}) = h(a^i) f_i^k e_{2k}$$ | $$(5.4.23) \qquad \overline{f} \circ (a_i e_1^i) = h(a_i) f_k^i e_2^k$$ |
| $$(5.4.13) \qquad \overline{f} \circ (e_1 a) = e_2 f h(a)$$ | $$(5.4.24) \qquad \overline{f} \circ (a e_1) = h(a) f e_2$$ |
| *relative to selected bases. Here* | *relative to selected bases. Here* |
| • *$a$ is coordinate matrix of $A_1$-number $\overline{a}$ relative the basis $\overline{\overline{e}}_1$* | • *$a$ is coordinate matrix of $A_1$-number $\overline{a}$ relative the basis $\overline{\overline{e}}_1$* |
| $$(5.4.14) \qquad \overline{a} = e_1 a$$ | $$(5.4.25) \qquad \overline{a} = a e_1$$ |
| • *$h(a) = (h(a_i), i \in I)$ is a matrix of $D_2$-numbers.* | • *$h(a) = (h(a^i), i \in I)$ is a matrix of $D_2$-numbers.* |
| • *$b$ is coordinate matrix of $A_2$-number* | • *$b$ is coordinate matrix of $A_2$-number* |
| $$\overline{b} = \overline{f} \circ \overline{a}$$ | $$\overline{b} = \overline{f} \circ \overline{a}$$ |
| *relative the basis $\overline{\overline{e}}_2$* | *relative the basis $\overline{\overline{e}}_2$* |
| $$(5.4.15) \qquad \overline{b} = e_2 b$$ | $$(5.4.26) \qquad \overline{b} = b e_2$$ |

---

[5.8] In theorems 5.4.3, 5.4.5, we use the following convention. Let the set of vectors $\overline{\overline{e}}_1 = (e_{1i}, i \in I)$ be a basis and $C_{1\,ij}^{\phantom{1}\phantom{ij}k}$, $k$, $i, j \in I$, be structural constants of $D_1$-algebra of columns $A_1$. Let the set of vectors $\overline{\overline{e}}_2 = (e_{2j}, j \in J)$ be a basis and $C_{2\,ij}^{\phantom{2}\phantom{ij}k}$, $k$, $i, j \in J$, be structural constants of $D_2$-algebra of columns $A_2$.

[5.9] In theorems 5.4.4, 5.4.6, we use the following convention. Let the set of vectors $\overline{\overline{e}}_1 = (e_1^i, i \in I)$ be a basis and $C_{1\,k}^{\phantom{1}\phantom{k}ij}$, $k$, $i, j \in I$, be structural constants of $D_1$-algebra of rows $A_1$. Let the set of vectors $\overline{\overline{e}}_2 = (e_2^j, j \in J)$ be a basis and $C_{2\,k}^{\phantom{2}\phantom{k}ij}$, $k$, $i, j \in J$, be structural constants of $D_2$-algebra of rows $A_2$.

- $f$ is coordinate matrix of set of $A_2$-numbers ($\overline{f} \circ e_{1i}, i \in I$) relative the basis $\overline{\overline{e}}_2$.
- There is relation between the matrix of homomorphism and structural constants

$$(5.4.16) \qquad h(C_{1\,ij}^{\ k})f_k^l = f_i^p f_j^q C_{2\,pq}^{\ l}$$

- $f$ is coordinate matrix of set of $A_2$-numbers ($\overline{f} \circ e_1^i, i \in I$) relative the basis $\overline{\overline{e}}_2$.
- There is relation between the matrix of homomorphism and structural constants

$$(5.4.27) \qquad h(C_{1\,k}^{\ ij})f_l^k = f_p^i f_q^j C_{2\,l}^{\ pq}$$

Proof. Equalities

$$\boxed{(5.4.11) \quad b = fh(a)}$$

$$\boxed{(5.4.12) \quad \overline{f} \circ (a^i e_{1i}) = h(a^i)f_i^k e_{2k}}$$

$$\boxed{(5.4.13) \quad \overline{f} \circ (e_1 a) = e_2 fh(a)}$$

follow from theorems 4.3.3, 5.4.2.

The equality

$$(5.4.17) \qquad \overline{f} \circ (e_{1i}e_{1j}) = h(C_{1\,ij}^{\ k})\overline{f} \circ e_{1k}$$

follows from equalities

$$\boxed{(5.2.2) \quad e_i e_j = C_{ij}^{\ k} e_k}$$

$$\boxed{(5.4.8) \quad f \circ (pa) = h(p)(f \circ a)}$$

The equality

$$(5.4.18) \qquad \overline{f} \circ (e_{1i}e_{1j}) = h(C_{1\,ij}^{\ k})f_k^l e_{2l}$$

follows from the equality

$$\boxed{(4.3.19) \quad \overline{f} \circ e_{1i} = e_2 f_i = e_{2j}f_i^j}$$

and the equality (5.4.17). The equality

$$(5.4.19) \qquad \begin{aligned} &\overline{f} \circ (e_{1i}e_{1j}) \\ &= (\overline{f} \circ e_{1i})(\overline{f} \circ e_{1j}) \\ &= f_i^p e_{2p} f_j^q e_{2q} \end{aligned}$$

follows from the equality

$$\boxed{(5.4.4) \quad f \circ (ab) = (f \circ a)(f \circ b)}$$

The equality

$$(5.4.20) \qquad \overline{f} \circ (e_{1i}e_{2j}) = f_i^p f_j^q C_{2\,pq}^{\ l} e_{2l}$$

follows from the equality

$$\boxed{(5.2.2) \quad e_i e_j = C_{ij}^{\ k} e_k}$$

and the equality (5.4.19). The equality

$$(5.4.21) \qquad h(C_{1\,ij}^{\ k})f_k^l e_{2l} = f_i^p f_j^q C_{2\,pq}^{\ l} e_{2l}$$

follows from equalities (5.4.18), (5.4.20). The equality

$$\boxed{(5.4.16) \quad h(C_{1\,ij}^{\ k})f_k^l = f_i^p f_j^q C_{2\,pq}^{\ l}}$$

follows from the equality (5.4.21), and and the theorem 4.2.7.  □

Proof. Equalities

$$\boxed{(5.4.22) \quad b = h(a)f}$$

$$\boxed{(5.4.23) \quad \overline{f} \circ (a_i e_1^i) = h(a_i)f_k^i e_2^k}$$

$$\boxed{(5.4.24) \quad \overline{f} \circ (ae_1) = h(a)fe_2}$$

follow from theorems 4.3.4, 5.4.2.

The equality

$$(5.4.28) \qquad \overline{f} \circ (e_1^i e_1^j) = h(C_{1\,k}^{\ ij})\overline{f} \circ e_1^k$$

follows from equalities

$$\boxed{(5.2.6) \quad e^i e^j = C_k^{\ ij} e^k}$$

$$\boxed{(5.4.8) \quad f \circ (pa) = h(p)(f \circ a)}$$

The equality

$$(5.4.29) \qquad \overline{f} \circ (e_1^i e_1^j) = h(C_{1\,k}^{\ ij})f_k^l e_2^l$$

follows from the equality

$$\boxed{(4.3.22) \quad \overline{f} \circ e_1^i = f_i e_2 = f_2^i e^j}$$

and the equality (5.4.28). The equality

$$(5.4.30) \qquad \begin{aligned} &\overline{f} \circ (e_1^i e_1^j) \\ &= (\overline{f} \circ e_1^i)(\overline{f} \circ e_1^j) \\ &= f_p^i e_2^p f_q^j e_2^q \end{aligned}$$

follows from the equality

$$\boxed{(5.4.4) \quad f \circ (ab) = (f \circ a)(f \circ b)}$$

The equality

$$(5.4.31) \qquad \overline{f} \circ (e_1^i e_2^j) = f_p^i f_q^j C_{2\,l}^{\ pq} e_2^l$$

follows from the equality

$$\boxed{(5.2.6) \quad e^i e^j = C_k^{\ ij} e^k}$$

and the equality (5.4.30). The equality

$$(5.4.32) \qquad h(C_{1\,k}^{\ ij})f_l^k e_2^l = f_p^i f_q^j C_{2\,l}^{\ pq} e_2^l$$

follows from equalities (5.4.29), (5.4.31). The equality

$$\boxed{(5.4.27) \quad h(C_{1\,k}^{\ ij})f_l^k = f_p^i f_q^j C_{2\,l}^{\ pq}}$$

follows from the equality (5.4.32), and and the theorem 4.2.14.     □

**THEOREM 5.4.5.** *Let the map*

$$h : D_1 \to D_2$$

*be homomorphism of the ring $D_1$ into the ring $D_2$.     Let*

$$f = (f_j^i, i \in \mathbf{I}, j \in \mathbf{J})$$

*be matrix of $D_2$-numbers which satisfy the equality*

(5.4.16)     $h(C_{1j}^{\ k})f_k^l = f_i^p f_j^q C_{2pq}^{\ \ l}$

*The map* [5.8]

(4.3.6)     $(h : D_1 \to D_2 \quad \overline{f} : V_1 \to V_2)$

*defined by the equality*

(4.3.7)     $\overline{f} \circ (e_1 a) = e_2 f h(a)$

*is homomorphism of $D_1$-algebra of columns $A_1$ into $D_2$-algebra of columns $A_2$.   The homomorphism (4.3.6) which has the given matrix $f$ is unique.*

PROOF. According to the theorem 4.3.5, the map (4.3.6) is homomorphism of $D_1$-module $A_1$ into $D_2$-module $A_2$ and the homomorphism (4.3.6) is unique.

Let  $a = a^i e_{1i}$,   $b = b^i e_{1i}$.  The equality

(5.4.33)  $\begin{aligned}&h(C_{1ij}^{\ \ k})f_k^l h(a^i) h(b^j) e_{2l} \\ &= f_i^p f_j^q C_{2pq}^{\ \ l} h(a^i) h(b^j) e_{2l}\end{aligned}$

follows from the equality (5.4.16).   Since the map $h$ is homomorphism of the ring $D_1$ into the ring $D_2$, then

(5.4.34)  $\begin{aligned}&h(C_{1ij}^{\ \ k})f_k^l h(a^i) h(b^j) e_{2l} \\ &= h(C_{1ij}^{\ \ k} a^i b^j) f_k^l e_{2l}\end{aligned}$

The equality

(5.4.35)  $\begin{aligned}&h(C_{1ij}^{\ \ k})f_k^l h(a^i) h(b^j) e_{2l} \\ &= h((ab)^k) f_k^l e_{2l}\end{aligned}$

follows from the equality (5.4.34) and the equality

(5.2.1)    $(ab)^k = C_{ij}^k a^i b^j$

**THEOREM 5.4.6.** *Let the map*

$$h : D_1 \to D_2$$

*be homomorphism of the ring $D_1$ into the ring $D_2$.     Let*

$$f = (f_i^j, i \in \mathbf{I}, j \in \mathbf{J})$$

*be matrix of $D_2$-numbers which satisfy the equality*

(5.4.27)    $h(C_{1k}^{\ ij})f_l^k = f_p^i f_q^j C_{2l}^{\ pq}$

*The map* [5.9]

(4.3.12)    $(h : D_1 \to D_2 \quad \overline{f} : V_1 \to V_2)$

*defined by the equality*

(4.3.13)    $\overline{f} \circ (ae_1) = h(a) f e_2$

*is homomorphism of $D_1$-algebra of rows $A_1$ into $D_2$-algebra of rows $A_2$.   The homomorphism (4.3.12) which has the given matrix $f$ is unique.*

PROOF. According to the theorem 4.3.6, the map (4.3.12) is homomorphism of $D_1$-module $A_1$ into $D_2$-module $A_2$ and the homomorphism (4.3.12) is unique.

Let  $a = a_i e_1^i$,   $b = b_i e_1^i$.  The equality

(5.4.39)  $\begin{aligned}&h(C_{1k}^{\ ij})f_l^k h(a_i) h(b_j) e_2^l \\ &= f_p^i f_q^j C_{2l}^{\ pq} h(a_i) h(b_j) e_2^l\end{aligned}$

follows from the equality (5.4.27).   Since the map $h$ is homomorphism of the ring $D_1$ into the ring $D_2$, then

(5.4.40)  $\begin{aligned}&h(C_{1k}^{\ ij})f_l^k h(a_i) h(b_j) e_2^l \\ &= h(C_{1k}^{\ ij} a_i b_j) f_l^k e_2^l\end{aligned}$

The equality

(5.4.41)  $\begin{aligned}&h(C_{1k}^{\ ij})f_l^k h(a_i) h(b_j) e_2^l \\ &= h((ab)_k) f_l^k e_2^l\end{aligned}$

follows from the equality (5.4.40) and the equality

(5.2.5)    $(ab)_k = C_k^{ij} a_i b_j$

The equality

$$(5.4.36) \qquad h(C_1{}^k_{ij}) f^l_k h(a^i) h(b^j) e_{2l}$$
$$= f \circ (ab)$$

follows from the equality (5.4.35) and the equality

$$\boxed{(5.4.12) \quad \overline{f} \circ (a^i e_{1i}) = h(a^i) f^k_i e_{2k}}$$

The equality

$$(5.4.37) \qquad f^p_i f^q_j C_{2pq}^l h(a^i) h(b^j) e_{2l}$$
$$= (h(a^i) f^p_i e_{2p})(h(b^j) f^q_j e_{2q})$$

follows from the equality

$$\boxed{(5.2.2) \quad e_i e_j = C^k_{ij} e_k}$$

The equality

$$(5.4.38) \qquad f^p_i f^q_j C_{2pq}^l h(a^i) h(b^j) e_{2l}$$
$$= (f \circ a)(f \circ b)$$

follows from the equality (5.4.37) and the equality

$$\boxed{(5.4.12) \quad \overline{f} \circ (a^i e_{1i}) = h(a^i) f^k_i e_{2k}}$$

The equality

$$\boxed{(5.4.4) \quad f \circ (ab) = (f \circ a)(f \circ b)}$$

follows from equalities (5.4.33), (5.4.36), (5.4.38). According to the theorem 5.4.2, the map (4.3.6) is homomorphism of $D_1$-algebra $A_1$ into $D_2$-algebra $A_2$. □

The equality

$$(5.4.42) \qquad h(C_1{}^{ij}_k) f^k_l h(a_i) h(b_j) e^l_2$$
$$= f \circ (ab)$$

follows from the equality (5.4.41) and the equality

$$\boxed{(5.4.23) \quad \overline{f} \circ (a_i e^i_1) = h(a_i) f^i_k e^k_2}$$

The equality

$$(5.4.43) \qquad f^i_p f^j_q C_{2l}^{pq} h(a_i) h(b_j) e^l_2$$
$$= (h(a_i) f^p_i e^p_2)(h(b_j) f^j_q e^q_2)$$

follows from the equality

$$\boxed{(5.2.6) \quad e^i e^j = C^{ij}_k e^k}$$

The equality

$$(5.4.44) \qquad f^i_p f^j_q C_{2l}^{pq} h(a_i) h(b_j) e^l_2$$
$$= (f \circ a)(f \circ b)$$

follows from the equality (5.4.43) and the equality

$$\boxed{(5.4.23) \quad \overline{f} \circ (a_i e^i_1) = h(a_i) f^i_k e^k_2}$$

The equality

$$\boxed{(5.4.4) \quad f \circ (ab) = (f \circ a)(f \circ b)}$$

follows from equalities (5.4.39), (5.4.42), (5.4.44). According to the theorem 5.4.2, the map (4.3.12) is homomorphism of $D_1$-algebra $A_1$ into $D_2$-algebra $A_2$. □

## 5.4.2. Homomorphism When Rings $D_1 = D_2 = D$.

**Definition 5.4.7.** *Let $A_i$, $i = 1, 2$, be algebra over commutative ring $D$. Let diagram of representations*

$$(5.4.45)$$

$$h_{1.12}(d) : v \to d\,v$$
$$h_{1.23}(v) : w \to C_1 \circ (v, w)$$
$$C_1 \in \mathcal{L}(D; A_1^2 \to A_1)$$

*describe $D_1$-algebra $A_1$. Let diagram of representations*

$$(5.4.46)$$

$$h_{2.12}(d) : v \to d\,v$$
$$h_{2.23}(v) : w \to C_2 \circ (v, w)$$
$$C_2 \in \mathcal{L}(D; A_2^2 \to A_2)$$

*describe $D_2$-algebra $A_2$. Morphism*

$$(5.4.47) \qquad\qquad f : A_1 \to A_2$$

*of diagram of representations* (5.4.45) *into diagram of representations* (5.4.46) *is called* **homomorphism** *of $D$-algebra $A_1$ into $D$-algebra $A_2$. Let us denote* $\mathrm{Hom}(D; A_1 \to A_2)$ *set of homomorphisms of $D$-algebra $A_1$ into $D$-algebra $A_2$.*

$\square$

THEOREM 5.4.8. *Homomorphism*

$$(5.4.47) \quad f : A_1 \to A_2$$

*of $D$-algebra $A_1$ into $D$-algebra $A_2$ is linear map* (5.4.47) *of $D$-module $A_1$ into $D$-module $A_2$ such that*

$$(5.4.48) \qquad\qquad f \circ (ab) = (f \circ a)(f \circ b)$$

*and satisfies following equalities*

$$(5.4.49) \qquad\qquad f \circ (a + b) = f \circ a + f \circ b$$

$$(5.4.48) \quad f \circ (ab) = (f \circ a)(f \circ b)$$

$$(5.4.50) \qquad\qquad f \circ (pa) = p(f \circ a)$$
$$p \in D \quad a, b \in A_1$$

PROOF. According to definitions 2.3.9, 2.8.5, the map (5.4.47) is morphism of representation $h_{1.12}$ describing $D$-module $A_1$ into representation $h_{2.12}$ describing $D$-module $A_2$. Therefore, the map $f$ is linear map of $D$-module $A_1$ into $D$-module $A_2$ and equalities (5.4.49), (5.4.50) follow from the theorem 4.3.10.

According to definitions 2.3.9, 2.8.5, the map $f$ is morphism of representation $h_{1.23}$ into representation $h_{2.23}$. The equality

$$(5.4.51) \qquad f \circ (h_{1.23}(a)(b)) = h_{2.23}(f \circ a)(f \circ b)$$

follows from the equality (2.3.2). From equalities (5.4.45), (5.4.46), (5.4.51), it follows that

$$(5.4.52) \qquad f \circ (C_1 \circ (a, b)) = C_2 \circ (f \circ a, f \circ b)$$

The equality (5.4.48) follows from equalities (5.1.1), (5.4.52).           $\square$

THEOREM 5.4.9. *The homomor­phism* [5.10]

$$\overline{f} : A_1 \to A_2$$

*of $D$-algebra of columns $A_1$ into $D$-alge­bra of columns $A_2$ has presentation*

$$(5.4.53) \qquad b = fa$$

THEOREM 5.4.10. *The homomor­phism* [5.11]

$$\overline{f} : A_1 \to A_2$$

*of $D$-algebra of rows $A_1$ into $D$-algebra of rows $A_2$ has presentation*

$$(5.4.64) \qquad b = af$$

---

[5.10] In theorems 5.4.9, 5.4.11, we use the following convention. Let the set of vectors $\overline{\overline{e}}_1 = (e_{1i}, i \in I)$ be a basis and $C_{1\ ij}^{\ \ k}$, $k, i, j \in I$, be structural constants of $D$-algebra of columns $A_1$. Let the set of vectors $\overline{\overline{e}}_2 = (e_{2j}, j \in J)$ be a basis and $C_{2\ ij}^{\ \ k}$, $k, i, j \in J$, be structural constants of $D$-algebra of columns $A_2$.

[5.11] In theorems 5.4.10, 5.4.12, we use the following convention. Let the set of vectors $\overline{\overline{e}}_1 = (e_1^i, i \in I)$ be a basis and $C_{1\ k}^{\ ij}$, $k, i, j \in I$, be structural constants of $D$-algebra of rows $A_1$. Let the set of vectors $\overline{\overline{e}}_2 = (e_2^j, j \in J)$ be a basis and $C_{2\ k}^{\ ij}$, $k, i, j \in J$, be structural constants of $D$-algebra of rows $A_2$.

(5.4.54)    $\overline{f} \circ (a^i e_{1i}) = a^i f_i^k e_{2k}$

(5.4.55)    $\overline{f} \circ (e_1 a) = e_2 f a$

*relative to selected bases. Here*

- *$a$ is coordinate matrix of $A_1$-number $\overline{a}$ relative the basis $\overline{\overline{e}}_1$*

(5.4.56)    $\overline{a} = e_1 a$

- *$b$ is coordinate matrix of $A_2$-number*

    $\overline{b} = \overline{f} \circ \overline{a}$

    *relative the basis $\overline{\overline{e}}_2$*

(5.4.57)    $\overline{b} = e_2 b$

- *$f$ is coordinate matrix of set of $A_2$-numbers $(\overline{f} \circ e_{1i}, i \in {\color{red}I})$ relative the basis $\overline{\overline{e}}_2$.*
- *There is relation between the matrix of homomorphism and structural constants*

(5.4.58)    $C_{1\,ij}^{\ \ k} f_k^l = f_i^p f_j^q C_{2\,pq}^{\ \ l}$

PROOF. Equalities

$\boxed{(5.4.53)\quad b = fa}$

$\boxed{(5.4.54)\quad \overline{f} \circ (a^i e_{1i}) = a^i f_i^k e_{2k}}$

$\boxed{(5.4.55)\quad \overline{f} \circ (e_1 a) = e_2 f a}$

follow from theorems {\color{red}4.3.11}, {\color{red}5.4.8}.
    The equality

(5.4.59)    $\overline{f} \circ (e_{1i} e_{1j}) = C_{1\,ij}^{\ \ k} \overline{f} \circ e_{1k}$

follows from equalities

$\boxed{(5.2.2)\quad e_i e_j = C_{ij}^k e_k}$

$\boxed{(5.4.50)\quad f \circ (pa) = p(f \circ a)}$

The equality

(5.4.60)    $\overline{f} \circ (e_{1i} e_{1j}) = C_{1\,ij}^{\ \ k} f_k^l e_{2l}$

follows from the equality

$\boxed{(4.3.47)\quad \overline{f} \circ e_{1i} = e_2 f_i = e_{2j} f_i^j}$

and the equality {\color{red}(5.4.59)}. The equality

$$\begin{aligned}
&\overline{f} \circ (e_{1i} e_{1j}) \\
\text{(5.4.61)}\quad &= (\overline{f} \circ e_{1i})(\overline{f} \circ e_{1j}) \\
&= f_i^p e_{2p} f_j^q e_{2q}
\end{aligned}$$

(5.4.65)    $\overline{f} \circ (a_i e_1^i) = a_i f_k^i e_2^k$

(5.4.66)    $\overline{f} \circ (a e_1) = a f e_2$

*relative to selected bases. Here*

- *$a$ is coordinate matrix of $A_1$-number $\overline{a}$ relative the basis $\overline{\overline{e}}^1$*

(5.4.67)    $\overline{a} = a e_1$

- *$b$ is coordinate matrix of $A_2$-number*

    $\overline{b} = \overline{f} \circ \overline{a}$

    *relative the basis $\overline{\overline{e}}^2$*

(5.4.68)    $\overline{b} = b e_2$

- *$f$ is coordinate matrix of set of $A_2$-numbers $(\overline{f} \circ e_1^i, i \in {\color{red}I})$ relative the basis $\overline{\overline{e}}^2$.*
- *There is relation between the matrix of homomorphism and structural constants*

(5.4.69)    $C_{1\,k}^{\ \ ij} f_l^k = f_p^i f_j^j C_{2l}^{\ pq}$

PROOF. Equalities

$\boxed{(5.4.64)\quad b = af}$

$\boxed{(5.4.65)\quad \overline{f} \circ (a_i e_1^i) = a_i f_k^i e_2^k}$

$\boxed{(5.4.66)\quad \overline{f} \circ (a e_1) = a f e_2}$

follow from theorems {\color{red}4.3.12}, {\color{red}5.4.8}.
    The equality

(5.4.70)    $\overline{f} \circ (e_1^i e_1^j) = C_{1\,k}^{\ \ ij} \overline{f} \circ e_1^k$

follows from equalities

$\boxed{(5.2.6)\quad e^i e^j = C_k^{ij} e^k}$

$\boxed{(5.4.50)\quad f \circ (pa) = p(f \circ a)}$

The equality

(5.4.71)    $\overline{f} \circ (e_1^i e_1^j) = C_{1\,k}^{\ \ ij} f_l^k e_2^l$

follows from the equality

$\boxed{(4.3.50)\quad \overline{f} \circ e_1^i = f_i e_2 = f_i^j e_2^j}$

and the equality {\color{red}(5.4.70)}. The equality

$$\begin{aligned}
&\overline{f} \circ (e_1^i e_1^j) \\
\text{(5.4.72)}\quad &= (\overline{f} \circ e_1^i)(\overline{f} \circ e_1^j) \\
&= f_p^i e_2^p f_q^j e_2^q
\end{aligned}$$

follows from the equality

$$(5.4.48) \quad f \circ (ab) = (f \circ a)(f \circ b)$$

The equality

$$(5.4.62) \quad \overline{f} \circ (e_{1i} e_{2j}) = f_i^p f_j^q C_{2\,pq}^{\phantom{pq}l} e_{2l}$$

follows from the equality

$$(5.2.2) \quad e_i e_j = C_{ij}^k e_k$$

and the equality (5.4.61). The equality

$$(5.4.63) \quad C_{1\,ij}^{\phantom{ij}k} f_k^l e_{2l} = f_i^p f_j^q C_{2\,pq}^{\phantom{pq}l} e_{2l}$$

follows from equalities (5.4.60), (5.4.62). The equality

$$(5.4.58) \quad C_{1\,ij}^{\phantom{ij}k} f_k^l = f_i^p f_j^q C_{2\,pq}^{\phantom{pq}l}$$

follows from the equality (5.4.63), and and the theorem 4.2.7.     □

**Theorem 5.4.11.** *Let*

$$f = (f_j^i, i \in \textcolor{red}{I}, j \in \textcolor{red}{J})$$

*be matrix of $D$-numbers which satisfy the equality*

$$(5.4.58) \quad C_{1\,ij}^{\phantom{ij}k} f_k^l = f_i^p f_j^q C_{2\,pq}^{\phantom{pq}l}$$

*The map* [5.10]

$$(4.3.34) \quad \overline{f} : V_1 \to V_2 \quad \text{defined by the}$$

*equality*

$$(4.3.35) \quad \overline{f} \circ (e_{1}a) = e_{2}fa$$

*is homomorphism of $D$-algebra of columns $A_1$ into $D$-algebra of columns $A_2$. The homomorphism (4.3.34) which has the given matrix $f$ is unique.*

**Proof.** According to the theorem 4.3.13, the map (4.3.34) is homomorphism of $D$-module $A_1$ into $D$-module $A_2$ and the homomorphism (4.3.34) is unique.

Let $a = a^i e_{1i}$, $b = b^i e_{1i}$. The equality

$$(5.4.75) \quad \begin{aligned} & C_{1\,ij}^{\phantom{ij}k} f_k^l a^i b^j e_{2l} \\ & = f_i^p f_j^q C_{2\,pq}^{\phantom{pq}l} a^i b^j e_{2l} \end{aligned}$$

follows from the equality (5.4.58). The equality

$$(5.4.76) \quad \begin{aligned} & C_{1\,ij}^{\phantom{ij}k} f_k^l a^i b^j e_{2l} \\ & = (ab)^k f_k^l e_{2l} \end{aligned}$$

follows from the equality

follows from the equality

$$(5.4.48) \quad f \circ (ab) = (f \circ a)(f \circ b)$$

The equality

$$(5.4.73) \quad \overline{f} \circ (e_1^i e_2^j) = f_p^i f_q^j C_{2l}^{\phantom{2l}pq} e_2^l$$

follows from the equality

$$(5.2.6) \quad e^i e^j = C_k^{ij} e^k$$

and the equality (5.4.72). The equality

$$(5.4.74) \quad C_{1k}^{\phantom{1k}ij} f_l^k e_2^l = f_p^i f_q^j C_{2l}^{\phantom{2l}pq} e_2^l$$

follows from equalities (5.4.71), (5.4.73). The equality

$$(5.4.69) \quad C_{1k}^{\phantom{1k}ij} f_l^k = f_p^i f_q^j C_{2l}^{\phantom{2l}pq}$$

follows from the equality (5.4.74), and and the theorem 4.2.14.     □

**Theorem 5.4.12.** *Let*

$$f = (f_i^j, i \in \textcolor{red}{I}, j \in \textcolor{red}{J})$$

*be matrix of $D$-numbers which satisfy the equality*

$$(5.4.69) \quad C_{1k}^{\phantom{1k}ij} f_l^k = f_p^i f_q^j C_{2l}^{\phantom{2l}pq}$$

*The map* [5.11]

$$(4.3.40) \quad \overline{f} : V_1 \to V_2 \quad \text{defined by the}$$

*equality*

$$(4.3.41) \quad \overline{f} \circ (ae_1) = af e_2$$

*is homomorphism of $D$-algebra of rows $A_1$ into $D$-algebra of rows $A_2$. The homomorphism (4.3.40) which has the given matrix $f$ is unique.*

**Proof.** According to the theorem 4.3.14, the map (4.3.40) is homomorphism of $D$-module $A_1$ into $D$-module $A_2$ and the homomorphism (4.3.40) is unique.

Let $a = a_i e_1^i$, $b = b_i e_1^i$. The equality

$$(5.4.80) \quad \begin{aligned} & C_{1k}^{\phantom{1k}ij} f_l^k a_i b_j e_2^l \\ & = f_p^i f_q^j C_{2l}^{\phantom{2l}pq} a_i b_j e_2^l \end{aligned}$$

follows from the equality (5.4.69). The equality

$$(5.4.81) \quad \begin{aligned} & C_{1k}^{\phantom{1k}ij} f_l^k a_i b_j e_2^l \\ & = (ab)_k f_l^k e_2^l \end{aligned}$$

$$(5.2.1) \quad (ab)^k = C_{ij}^k a^i b^j$$

The equality

$$(5.4.77) \qquad C_{1\,ij}^{\ \ k} f_k^l a^i b^j e_{2l}$$
$$= f \circ (ab)$$

follows from the equality $(5.4.76)$ and the equality

$$(5.4.54) \quad \overline{f} \circ (a^i e_{1i}) = a^i f_i^k e_{2k}$$

The equality

$$(5.4.78) \qquad f_i^p f_j^q C_{2\,pq}^{\ \ l} a^i b^j e_{2l}$$
$$= (a^i f_i^p e_{2p})(b^j f_j^q e_{2q})$$

follows from the equality

$$(5.2.2) \quad e_i e_j = C_{ij}^k e_k$$

The equality

$$(5.4.79) \qquad f_i^p f_j^q C_{2\,pq}^{\ \ l} a^i b^j e_{2l}$$
$$= (f \circ a)(f \circ b)$$

follows from the equality $(5.4.78)$ and the equality

$$(5.4.54) \quad \overline{f} \circ (a^i e_{1i}) = a^i f_i^k e_{2k}$$

The equality

$$(5.4.48) \quad f \circ (ab) = (f \circ a)(f \circ b)$$

follows from equalities $(5.4.75)$, $(5.4.77)$, $(5.4.79)$. According to the theorem $5.4.8$, the map $(4.3.34)$ is homomorphism of *D*-algebra $A_1$ into *D*-algebra $A_2$. $\quad\square$

follows from the equality

$$(5.2.5) \quad (ab)_k = C_k^{ij} a_i b_j$$

The equality

$$(5.4.82) \qquad C_{1\,k}^{\ \ ij} f_l^k a_i b_j e_2^l$$
$$= f \circ (ab)$$

follows from the equality $(5.4.81)$ and the equality

$$(5.4.65) \quad \overline{f} \circ (a_i e_1^i) = a_i f_k^i e_2^k$$

The equality

$$(5.4.83) \qquad f_p^i f_q^j C_{2\,l}^{\ \ pq} a_i b_j e_2^l$$
$$= (a_i f_p^i e_2^p)(b_j f_q^j e_2^q)$$

follows from the equality

$$(5.2.6) \quad e^i e^j = C_k^{ij} e^k$$

The equality

$$(5.4.84) \qquad f_p^i f_q^j C_{2\,l}^{\ \ pq} a_i b_j e_2^l$$
$$= (f \circ a)(f \circ b)$$

follows from the equality $(5.4.83)$ and the equality

$$(5.4.65) \quad \overline{f} \circ (a_i e_1^i) = a_i f_k^i e_2^k$$

The equality

$$(5.4.48) \quad f \circ (ab) = (f \circ a)(f \circ b)$$

follows from equalities $(5.4.80)$, $(5.4.82)$, $(5.4.84)$. According to the theorem $5.4.8$, the map $(4.3.40)$ is homomorphism of *D*-algebra $A_1$ into *D*-algebra $A_2$. $\quad\square$

# Linear Map of Algebra

## 6.1. Linear Map of $D$-Algebra

DEFINITION 6.1.1. *Let $A_1$ and $A_2$ be algebras over commutative ring $D$. The linear map of the $D$-module $A_1$ into the $D$-module $A_2$ is called* **linear map** *of $D$-algebra $A_1$ into $D$-algebra $A_2$.*

*Let us denote* $\mathcal{L}(D; A_1 \to A_2)$ *set of linear maps of $D$-algebra $A_1$ into $D$-algebra $A_2$.* □

DEFINITION 6.1.2. *Let $A_1$, ..., $A_n$, $S$ be $D$-algebras. Polylinear map*

$$f : A_1 \times ... \times A_n \to S$$

*of $D$-modules $A_1$, ..., $A_n$ into $D$-module $S$ is called* **polylinear map** *of $D$-algebras $A_1$, ..., $A_n$ into $D$-algebra $S$. Let us denote* $\mathcal{L}(D; A_1 \times ... \times A_n \to S)$ *set of polylinear maps of $D$-modules $A_1$, ..., $A_n$ into $D$-module $S$. Let us denote* $\mathcal{L}(D; A^n \to S)$ *set of $n$-linear maps of $D$-module $A$ ($A_1 = ... = A_n = A$) into $D$-module $S$.* □

THEOREM 6.1.3. *Let $A_1$, ..., $A_n$ be $D$-algebras. Tensor product $A_1 \otimes ... \otimes A_n$ of $D$-modules $A_1$, ..., $A_n$ is $D$-algebra, if we define product by the equality*

$$(a_1, ..., a_n) * (b_1, ..., b_n) = (a_1 b_1) \otimes ... \otimes (a_n b_n)$$

PROOF. According to the definition 5.1.1 and to the theorem 4.6.2, tensor product $A_1 \otimes ... \otimes A_n$ of $D$-algebras $A_1$, ..., $A_n$ is $D$-module.

Consider the map

$$(6.1.1) \qquad * : (A_1 \times ... \times A_n) \times (A_1 \times ... \times A_n) \to A_1 \otimes ... \otimes A_n$$

defined by the equality

$$(a_1, ..., a_n) * (b_1, ..., b_n) = (a_1 b_1) \otimes ... \otimes (a_n b_n)$$

For given values of variables $b_1$, ..., $b_n$, the map (6.1.1) is polylinear map with respect to variables $a_1$, ..., $a_n$. According to the theorem 4.6.4, there exists a linear map

$$(6.1.2) \qquad *(b_1, ..., b_n) : A_1 \otimes ... \otimes A_n \to A_1 \otimes ... \otimes A_n$$

defined by the equality

$$(6.1.3) \qquad (a_1 \otimes ... \otimes a_n) * (b_1, ..., b_n) = (a_1 b_1) \otimes ... \otimes (a_n b_n)$$

Since we can present any tensor $a \in A_1 \otimes ... \otimes A_n$ as sum of tensors $a_1 \otimes ... \otimes a_n$, then, for given tensor $a \in A_1 \otimes ... \otimes A_n$, the map (6.1.2) is polylinear map of variables $b_1$, ..., $b_n$. According to the theorem 4.6.4, there exists a linear map

$$(6.1.4) \qquad *(a) : A_1 \otimes ... \otimes A_n \to A_1 \otimes ... \otimes A_n$$

defined by the equality

$$(6.1.5) \qquad (a_1 \otimes ... \otimes a_n) * (b_1 \otimes ... \otimes b_n) = (a_1 b_1) \otimes ... \otimes (a_n b_n)$$

Therefore, the equality (6.1.5) defines bilinear map

$$(6.1.6) \qquad *: (A_1 \otimes ... \otimes A_n) \times (A_1 \otimes ... \otimes A_n) \to A_1 \otimes ... \otimes A_n$$

Bilinear map (6.1.6) generates the product in $D$-module $A_1 \otimes ... \otimes A_n$.    □

In case of tensor product of $D$-algebras $A_1$, $A_2$ we consider product defined by the equality

$$(6.1.7) \qquad (a_1 \otimes a_2) \circ (b_1 \otimes b_2) = (a_1 b_1) \otimes (b_2 a_2)$$

THEOREM 6.1.4. *Let $\overline{e}_i$ be the basis of the algebra $A_i$ over the ring $D$. Let $B_{i\cdot kl}^{\;\;\;j}$ be structural constants of the algebra $A_i$ relative the basis $\overline{e}_i$. Structural constants of the tensor product $A_1 \otimes ... \otimes A_n$ relative to the basis $e_{1 \cdot i_l} \otimes ... \otimes e_{n \cdot i_n}$ have form*

$$(6.1.8) \qquad C_{\cdot k_l \dots k_n, l_l \dots l_n}^{\cdot j_l \dots j_n} = C_{1 \cdot k_l l_l}^{\;\; j_l} \dots C_{n \cdot k_n l_n}^{\;\; j_n}$$

PROOF. Direct multiplication of tensors $e_{1 \cdot i_l} \otimes ... \otimes e_{n \cdot i_n}$ has form

$$
\begin{aligned}
(6.1.9) \qquad & (e_{1 \cdot k_l} \otimes ... \otimes e_{n \cdot k_n})(e_{1 \cdot l_l} \otimes ... \otimes e_{n \cdot l_n}) \\
=& (\overline{e}_{1 \cdot k_l} \overline{e}_{1 \cdot l_l}) \otimes ... \otimes (\overline{e}_{n \cdot k_n} \overline{e}_{n \cdot l_n}) \\
=& (\overline{e}_{1 \cdot k_l} \overline{e}_{1 \cdot l_l}) \otimes ... \otimes (\overline{e}_{n \cdot k_n} \overline{e}_{n \cdot l_n}) \\
=& (C_{1 \cdot k_l l_l}^{\;\; j_l} \overline{e}_{1 \cdot j_l}) \otimes ... \otimes (C_{n \cdot k_n l_n}^{\;\; j_n} \overline{e}_{n \cdot j_n}) \\
=& C_{1 \cdot k_l l_l}^{\;\; j_l} \dots C_{n \cdot k_n l_n}^{\;\; j_n} \overline{e}_{1 \cdot j_l} \otimes ... \otimes \overline{e}_{n \cdot j_n}
\end{aligned}
$$

According to the definition of structural constants

$$(6.1.10) \quad (e_{1 \cdot k_l} \otimes ... \otimes e_{n \cdot k_n})(e_{1 \cdot l_l} \otimes ... \otimes e_{n \cdot l_n}) = C_{\cdot k_l \dots k_n, l_l \dots l_n}^{\cdot j_l \dots j_n}(e_{1 \cdot j_l} \otimes ... \otimes e_{n \cdot j_n})$$

The equality (6.1.8) follows from comparison (6.1.9), (6.1.10).

From the chain of equalities

$$
\begin{aligned}
& (a_1 \otimes ... \otimes a_n)(b_1 \otimes ... \otimes b_n) \\
=& (a_1^{k_l} \overline{e}_{1 \cdot k_l} \otimes ... \otimes a_n^{k_n} \overline{e}_{n \cdot k_n})(b_1^{l_l} \overline{e}_{1 \cdot l_l} \otimes ... \otimes b_n^{l_n} \overline{e}_{n \cdot l_n}) \\
=& a_1^{k_l} \dots a_n^{k_n} b_1^{l_l} \dots b_n^{l_n} (\overline{e}_{1 \cdot k_l} \otimes ... \otimes \overline{e}_{n \cdot k_n})(\overline{e}_{1 \cdot l_l} \otimes ... \otimes \overline{e}_{n \cdot l_n}) \\
=& a_1^{k_l} \dots a_n^{k_n} b_1^{l_l} \dots b_n^{l_n} C_{\cdot k_l \dots k_n, l_l \dots l_n}^{\cdot j_l \dots j_n} (e_{1 \cdot j_l} \otimes ... \otimes e_{n \cdot j_n}) \\
=& a_1^{k_l} \dots a_n^{k_n} b_1^{l_l} \dots b_n^{l_n} C_{1 \cdot k_l l_l}^{\;\; j_l} \dots C_{n \cdot k_n l_n}^{\;\; j_n} (e_{1 \cdot j_l} \otimes ... \otimes e_{n \cdot j_n}) \\
=& (a_1^{k_l} b_1^{l_l} \; C_{1 \cdot k_l l_l}^{\;\; j_l} \overline{e}_{1 \cdot j_l}) \otimes ... \otimes (a_n^{k_n} b_n^{l_n} C_{n \cdot k_n l_n}^{\;\; j_n} \overline{e}_{n \cdot j_n}) \\
=& (a_1 b_1) \otimes ... \otimes (a_n b_n)
\end{aligned}
$$

it follows that definition of product (6.1.10) with structural constants (6.1.8) agreed with the definition of product (6.1.5).    □

THEOREM 6.1.5. *For tensors $a$, $b \in A_1 \otimes ... \otimes A_n$, standard components of product satisfy to equality*

$$(6.1.11) \qquad (ab)^{j_l \dots j_n} = C_{\cdot k_l \dots k_n, l_l \dots l_n}^{\cdot j_l \dots j_n} a^{k_l \dots k_n} b^{l_l \dots l_n}$$

PROOF. According to the definition

$$(6.1.12) \qquad ab = (ab)^{j_l \dots j_n} e_{1 \cdot j_l} \otimes ... \otimes e_{n \cdot j_n}$$

At the same time

$$
\begin{aligned}
(6.1.13) \qquad ab =& a^{k_l \dots k_n} e_{1 \cdot k_l} \otimes ... \otimes e_{n \cdot k_n} b^{k_l \dots k_n} e_{1 \cdot l_l} \otimes ... \otimes e_{n \cdot l_n} \\
=& a^{k_l \dots k_n} b^{k_l \dots k_n} C_{\cdot k_l \dots k_n, l_l \dots l_n}^{\cdot j_l \dots j_n} e_{1 \cdot j_l} \otimes ... \otimes e_{n \cdot j_n}
\end{aligned}
$$

The equality (6.1.11) follows from equalities (6.1.12), (6.1.13).    □

THEOREM 6.1.6. *If the algebra $A_i$, $i = 1, ..., n$, is associative, then the tensor product $A_1 \otimes ... \otimes A_n$ is associative algebra.*

PROOF. Since

$$((e_{1 \cdot i_l} \otimes ... \otimes e_{n \cdot i_n})(e_{1 \cdot j_l} \otimes ... \otimes e_{n \cdot j_n}))(e_{1 \cdot k_l} \otimes ... \otimes e_{n \cdot k_n})$$

$$=((\overline{e}_{1 \cdot i_l} \overline{e}_{1 \cdot j_l}) \otimes ... \otimes (\overline{e}_{n \cdot i_n} \overline{e}_{1 \cdot j_n}))(e_{1 \cdot k_l} \otimes ... \otimes e_{n \cdot k_n})$$

$$=((\overline{e}_{1 \cdot i_l} \overline{e}_{1 \cdot j_l})\overline{e}_{1 \cdot k_l}) \otimes ... \otimes ((\overline{e}_{n \cdot i_n} \overline{e}_{1 \cdot j_n})\overline{e}_{1 \cdot k_n})$$

$$=(\overline{e}_{1 \cdot i_l}(\overline{e}_{1 \cdot j_l} \overline{e}_{1 \cdot k_l})) \otimes ... \otimes (\overline{e}_{n \cdot i_n}(\overline{e}_{1 \cdot j_n} \overline{e}_{1 \cdot k_n}))$$

$$=(e_{1 \cdot i_l} \otimes ... \otimes e_{n \cdot i_n})((\overline{e}_{1 \cdot j_l} \overline{e}_{1 \cdot k_l}) \otimes ... \otimes (\overline{e}_{1 \cdot j_n} \overline{e}_{1 \cdot k_n}))$$

$$=(e_{1 \cdot i_l} \otimes ... \otimes e_{n \cdot i_n})((e_{1 \cdot j_l} \otimes ... \otimes e_{n \cdot j_n})(e_{1 \cdot k_l} \otimes ... \otimes e_{n \cdot k_n}))$$

then

$$(ab)c = a^{i_l...i_n} b^{j_l...j_n} c^{k_l...k_n}$$

$$((e_{1 \cdot i_l} \otimes ... \otimes e_{n \cdot i_n})(e_{1 \cdot j_l} \otimes ... \otimes e_{n \cdot j_n}))(e_{1 \cdot k_l} \otimes ... \otimes e_{n \cdot k_n})$$

$$= a^{i_l...i_n} b^{j_l...j_n} c^{k_l...k_n}$$

$$(e_{1 \cdot i_l} \otimes ... \otimes e_{n \cdot i_n})((e_{1 \cdot j_l} \otimes ... \otimes e_{n \cdot j_n})(e_{1 \cdot k_l} \otimes ... \otimes e_{n \cdot k_n}))$$

$$= a(bc)$$

□

THEOREM 6.1.7. *Let $A$ be algebra over commutative ring $D$. There exists a linear map*

$$h : a \otimes b \in A \otimes A \to ab \in A$$

PROOF. The theorem is corollary of the definition 5.1.1 and the theorem 4.6.4.    □

THEOREM 6.1.8. *Let map*

$$f : A_1 \to A_2$$

*be linear map of $D$-algebra $A_1$ into $D$-algebra $A_2$. Then maps $a \circ f$, $b \star f$, $a$, $b \in A_2$, defined by equalities*

$$(a \circ f) \circ x = a(f \circ x)$$

$$(b \star f) \circ x = (f \circ x)b$$

*are linear.*

PROOF. Statement of theorem follows from chains of equalities

$$(af) \circ (x + y) = a(f \circ (x + y)) = a(f \circ x + f \circ y) = a(f \circ x) + a(f \circ y)$$
$$= (af) \circ x + (af) \circ y$$
$$(af) \circ (px) = a(f \circ (px)) = ap(f \circ x) = pa(f \circ x)$$
$$= p((af) \circ x)$$
$$(fb) \circ (x + y) = (f \circ (x + y))b = (f \circ x + f \circ y)\ b = (f \circ x)b + (f \circ y)b$$
$$= (fb) \circ x + (fb) \circ y$$
$$(fb) \circ (px) = (f \circ (px))b = p(f \circ x)b$$
$$= p((fb) \circ x)$$

$\square$

## 6.2. Algebra $\mathcal{L}(D; A \to A)$

THEOREM 6.2.1. *Let $A$, $B$, $C$ be algebras over commutative ring $D$. Let $f$ be linear map from $D$-algebra $A$ into $D$-algebra $B$. Let $g$ be linear map from $D$-algebra $B$ into $D$-algebra $C$. The map $g \circ f$ defined by diagram*

(6.2.1)

*is linear map from $D$-algebra $A$ into $D$-algebra $C$.*

PROOF. The proof of the theorem follows from chains of equalities

$$(g \circ f) \circ (a + b) = g \circ (f \circ (a + b)) = g \circ (f \circ a + f \circ b)$$
$$= g \circ (f \circ a) + g \circ (f \circ b) = (g \circ f) \circ a + (g \circ f) \circ b$$
$$(g \circ f) \circ (pa) = g \circ (f \circ (pa)) = g \circ (p\ f \circ a) = p\ g \circ (f \circ a)$$
$$= p\ (g \circ f) \circ a$$

$\square$

THEOREM 6.2.2. *Let $A$, $B$, $C$ be algebras over the commutative ring $D$. Let $f$ be a linear map from $D$-algebra $A$ into $D$-algebra $B$. The map $f$ generates a linear map*

(6.2.2)        $$f^* : g \in \mathcal{L}(D; B \to C) \to g \circ f \in \mathcal{L}(D; A \to C)$$

(6.2.3)

Proof. The proof of the theorem follows from chains of equalities [6.1]

$$((g_1 + g_2) \circ f) \circ a = (g_1 + g_2) \circ (f \circ a) = g_1 \circ (f \circ a) + g_2 \circ (f \circ a)$$
$$= (g_1 \circ f) \circ a + (g_2 \circ f) \circ a$$
$$= (g_1 \circ f + g_2 \circ f) \circ a$$
$$((pg) \circ f) \circ a = (pg) \circ (f \circ a) = p \, g \circ (f \circ a) = p \, (g \circ f) \circ a$$
$$= (p(g \circ f)) \circ a$$

$\square$

Theorem 6.2.3. *Let $A$, $B$, $C$ be algebras over the commutative ring $D$. Let $g$ be a linear map from $D$-algebra $B$ into $D$-algebra $C$. The map $g$ generates a linear map*

(6.2.6) $$g_* : f \in \mathcal{L}(D; A \to B) \to g \circ f \in \mathcal{L}(D; A \to C)$$

(6.2.7)

Proof. The proof of the theorem follows from chains of equalities [6.2]

$$(g \circ (f_1 + f_2)) \circ a = g \circ ((f_1 + f_2) \circ a) = g \circ (f_1 \circ a + f_2 \circ a)$$
$$= g \circ (f_1 \circ a) + g \circ (f_2 \circ a) = (g \circ f_1) \circ a + (g \circ f_2) \circ a$$
$$= (g \circ f_1 + g \circ f_2) \circ a$$
$$(g \circ (pf)) \circ a = g \circ ((pf) \circ a) = g \circ (p \, (f \circ a)) = p \, g \circ (f \circ a)$$
$$= p \, (g \circ f) \circ a = (p(g \circ f)) \circ a$$

$\square$

Theorem 6.2.4. *Let $A$, $B$, $C$ be algebras over the commutative ring $D$. The map*

(6.2.10) $$\circ : (g, f) \in \mathcal{L}(D; B \to C) \times \mathcal{L}(D; A \to B) \to g \circ f \in \mathcal{L}(D; A \to C)$$

*is bilinear map.*

Proof. The theorem follows from theorems 6.2.2, 6.2.3. $\square$

Theorem 6.2.5. *Let $A$ be algebra over commutative ring $D$. $D$-module $\mathcal{L}(D; A \to A)$ equiped by product*

(6.2.11) $$\circ : (g, f) \in \mathcal{L}(D; A \to A) \times \mathcal{L}(D; A \to A) \to g \circ f \in \mathcal{L}(D; A \to A)$$

---

[6.1]We use following definitions of operations over maps

(6.2.4) $$(f + g) \circ a = f \circ a + g \circ a$$
(6.2.5) $$(pf) \circ a = p \, f \circ a$$

[6.2]We use following definitions of operations over maps

(6.2.8) $$(f + g) \circ a = f \circ a + g \circ a$$
(6.2.9) $$(pf) \circ a = p \, f \circ a$$

*is algebra over $D$.*

PROOF. The theorem follows from definition 5.1.1 and theorem 6.2.4.  □

### 6.3. Linear Map into Associative Algebra

THEOREM 6.3.1. *Consider $D$-algebras $A_1$ and $A_2$. For given map $f \in \mathcal{L}(D; A_1 \to A_2)$, the map*

$$g : A_2 \times A_2 \to \mathcal{L}(D; A_1 \to A_2)$$
$$g(a, b) \circ f = afb$$

*is bilinear map.*

PROOF. The statement of theorem follows from chains of equalities

$$((a_1 + a_2)fb) \circ x = (a_1 + a_2) \ f \circ x \ b = a_1 \ f \circ x \ b + a_2 \ f \circ x \ b$$
$$= (a_1 fb) \circ x + (a_2 fb) \circ x = (a_1 fb + a_2 fb) \circ x$$
$$((pa)fb) \circ x = (pa) \ f \circ x \ b = p(a \ f \circ x \ b) = p((afb) \circ x) = (p(afb)) \circ x$$
$$(af(b_1 + b_2)) \circ x = a \ f \circ x \ (b_1 + b_2) = a \ f \circ x \ b_1 + a \ f \circ x \ b_2$$
$$= (afb_1) \circ x + (afb_2) \circ x = (afb_1 + afb_2) \circ x$$
$$(af(pb)) \circ x = a \ f \circ x \ (pb) = p(a \ f \circ x \ b) = p((afb) \circ x) = (p(afb)) \circ x$$

□

THEOREM 6.3.2. *For given map $f \in \mathcal{L}(D; A_1 \to A_2)$ of $D$-algebra $A_1$ into $D$-algebra $A_2$, there exists linear map*

$$h : A_2 \otimes A_2 \to \mathcal{L}(D; A_1 \to A_2)$$

*defined by the equality*

(6.3.1)                               $$(a \otimes b) \circ f = b \star (a \circ f)$$
$$((a \otimes b) \circ f) \circ x = a(f \circ x)b$$

PROOF. The statement of the theorem is corollary of theorems 4.6.4, 6.3.1.  □

THEOREM 6.3.3. *Let $A_1$ and $A_2$ be $D$-algebras. Let product in algebra $A \otimes A$ be defined according to rule*

(6.3.2)                         $$(p_0 \otimes p_1) \circ (q_0 \otimes q_1) = (p_0 q_0) \otimes (q_1 p_1)$$

*A linear map*

(6.3.3)                         $$h : A_2 \otimes A_2 \dashrightarrow \mathcal{L}(D; A_1 \to A_2)$$

*defined by the equality*

(6.3.4)           $$(a \otimes b) \circ f = afb \quad a, b \in A_2 \quad f \in \mathcal{L}(D; A_1 \to A_2)$$

*is representation [6.3] of algebra $A_2 \otimes A_2$ in module $\mathcal{L}(D; A_1 \to A_2)$.*

PROOF. According to theorem 6.1.8, map (6.3.4) is transformation of module $\mathcal{L}(D; A_1 \rightarrow A_2)$. For a given tensor $c \in A_2 \otimes A_2$, a transformation $h(c)$ is a linear transformation of module $\mathcal{L}(D; A_1 \rightarrow A_2)$, because

$$((a \otimes b) \circ (f_1 + f_2)) \circ x = (a(f_1 + f_2)b) \circ x = a((f_1 + f_2) \circ x)b$$
$$= a(f_1 \circ x + f_2 \circ x)b = a(f_1 \circ x)b + a(f_2 \circ x)b$$
$$= (af_1 b) \circ x + (af_2 b) \circ x$$
$$= (a \otimes b) \circ f_1 \circ x + (a \otimes b) \circ f_2 \circ x$$
$$= ((a \otimes b) \circ f_1 + (a \otimes b) \circ f_2) \circ x$$
$$((a \otimes b) \circ (pf)) \circ x = (a(pf)b) \circ x = a((pf) \circ x)b$$
$$= a(p \ f \circ x)b = pa(f \circ x)b$$
$$= p \ (afb) \circ x = p \ ((a \otimes b) \circ f) \circ x$$
$$= (p((a \otimes b) \circ f)) \circ x$$

According to theorem 6.3.2, map (6.3.4) is linear map.

Let $f \in \mathcal{L}(D; A_1 \rightarrow A_2)$, $a \otimes b$, $c \otimes d \in A_2 \otimes A_2$. According to the theorem 6.3.2

$$(a \otimes b) \circ f = afb \in \mathcal{L}(D; A_1 \rightarrow A_2)$$

Therefore, according to theorem 6.3.2

$$(c \otimes d) \circ ((a \otimes b) \circ f) = c(afb)d$$

Since the product in algebra $A_2$ is associative, then

$$(c \otimes d) \circ ((a \otimes b) \circ f) = c(afb)d = (ca)f(bd) = (ca \otimes bd) \circ f$$

Therefore, since we define the product in algebra $A_2 \otimes A_2$ according to equality (6.3.2), then the map (6.3.3) is morphism of algebras. According to the definition 2.5.1, map (6.3.4) is a representation of the algebra $A_2 \otimes A_2$ in the module $\mathcal{L}(D; A_1 \rightarrow A_2)$. $\qquad \square$

THEOREM 6.3.4. *Let $A$ be $D$-algebra. Let product in $D$-module $A \otimes A$ be defined according to rule*

$$(p_0 \otimes p_1) \circ (q_0 \otimes q_1) = (p_0 q_0) \otimes (q_1 p_1)$$

*A representation*

$$h : A \otimes A \longrightarrow\!\!* \ \mathcal{L}(D; A \rightarrow A) \qquad h(p) : g \rightarrow p \circ g$$

*of $D$-algebra $A \otimes A$ in module $\mathcal{L}(D; A \rightarrow A)$ defined by the equality*

$$(a \otimes b) \circ g = agb \qquad a, b \in A \qquad g \in \mathcal{L}(D; A \rightarrow A)$$

*allows us to identify tensor $d \in A \times A$ and linear map $d \circ \delta \in \mathcal{L}(D; A \rightarrow A)$ where $\delta \in \mathcal{L}(D; A \rightarrow A)$ is identity map. Linear map $(a \otimes b) \circ \delta$ has form*

$$(6.3.5) \qquad\qquad (a \otimes b) \circ c = acb$$

PROOF. According to the theorem 6.3.2, the map $f \in \mathcal{L}(D; A \rightarrow A)$ and the tensor $d \in A \otimes A$ generate the map

$$(6.3.6) \qquad\qquad x \rightarrow (d \circ f) \circ x$$

---

[6.3] See the definition of representation of $\Omega$-algebra in the definition 2.5.1.

If we assume $f = \delta$, $d = a \otimes b$, then the equality (6.3.6) gets form

$$(6.3.7) \qquad ((a \otimes b) \circ \delta) \circ x = (a\delta b) \circ x = a\,(\delta \circ x)\,b = axb$$

If we assume

$$(6.3.8) \qquad ((a \otimes b) \circ \delta) \circ x = (a \otimes b) \circ (\delta \circ x) = (a \otimes b) \circ x$$

then comparison of equalities (6.3.7) and (6.3.8) gives a basis to identify the action of the tensor $a \otimes b$ and transformation $(a \otimes b) \circ \delta$. □

From the theorem 6.3.4, it follows that we can consider the map (6.3.4) as the product of maps $a \otimes b$ and $f$. The tensor $a \in A_2 \otimes A_2$ is **nonsingular**, if there exists the tensor $b \in A_2 \otimes A_2$ such that $a \circ b = 1 \otimes 1$.

DEFINITION 6.3.5. Consider[6.4] the representation of algebra $A_2 \otimes A_2$ in the module $\mathcal{L}(D; A_1 \to A_2)$. The set

$$(A_2 \otimes A_2) \circ f = \{g = d \circ f : d \in A_2 \otimes A_2\}$$

is called **orbit of linear map** $f \in \mathcal{L}(D; A_1 \to A_2)$. □

THEOREM 6.3.6. Consider $D$-algebra $A_1$ and associative $D$-algebra $A_2$. Consider the representation of algebra $A_2 \otimes A_2$ in the module $\mathcal{L}(D; A_1 \to A_2)$. The map

$$h : A_1 \to A_2$$

generated by the map

$$f : A_1 \to A_2$$

has form

$$(6.3.9) \qquad h = (a_{s\cdot 0} \otimes a_{s\cdot 1}) \circ f = a_{s\cdot 0} f a_{s\cdot 1}$$

PROOF. We can represent any tensor $a \in A_2 \otimes A_2$ in the form

$$a = a_{s\cdot 0} \otimes a_{s\cdot 1}$$

According to the theorem 6.3.3, the map (6.3.4) is linear. This proofs the statement of the theorem. □

DEFINITION 6.3.7. Expression $a_{s\cdot p}$, $p = 0,\ 1,\ $ in equality (6.3.9) is called **component of linear map** $h$. □

THEOREM 6.3.8. Let $A_2$ be algebra with unit $e$. Let $a \in A_2 \otimes A_2$ be a nonsingular tensor. Orbits of linear maps $f \in \mathcal{L}(D; A_1 \to A_2)$ and $g = a \circ f$ coincide

$$(6.3.10) \qquad (A_2 \otimes A_2) \circ f = (A_2 \otimes A_2) \circ g$$

PROOF. If $h \in (A_2 \otimes A_2) \circ g$, then there exists $b \in A_2 \otimes A_2$ such that $h = b \circ g$. In that case

$$(6.3.11) \qquad h = b \circ (a \circ f) = (b \circ a) \circ f$$

Therefore, $h \in (A_2 \otimes A_2) \circ f$,

$$(6.3.12) \qquad (A_2 \otimes A_2) \circ g \subset (A_2 \otimes A_2) \circ f$$

---

[6.4]The definition is made by analogy with the definition 2.5.11.

Since $a$ is nonsingular tensor, then

(6.3.13)
$$f = a^{-1} \circ g$$

If $h \in (A_2 \otimes A_2) \circ f$, then there exists $b \in A_2 \otimes A_2$ such that

(6.3.14)
$$h = b \circ f$$

From equalities (6.3.13), (6.3.14), it follows that

$$h = b \circ (a^{-1} \circ g) = (b \circ a^{-1}) \circ g$$

Therefore, $h \in (A_2 \otimes A_2) \circ g$,

(6.3.15)
$$(A_2 \otimes A_2) \circ f \subset (A_2 \otimes A_2) \circ g$$

(6.3.10) follows from equalities (6.3.12), (6.3.15).                        □

From the theorem 6.3.8, it also follows that if $g = a \circ f$ and $a \in A_2 \otimes A_2$ is a singular tensor, then relationship (6.3.12) is true. However, the main result of the theorem 6.3.8 is that the representations of the algebra $A_2 \otimes A_2$ in module $\mathcal{L}(D; A_1 \to A_2)$ generates an equivalence in the module $\mathcal{L}(D; A_1 \to A_2)$. If we successfully choose the representatives of each equivalence class, then the resulting set will be generating set of considered representation.[6.5]

## 6.4. Polylinear Map into Associative Algebra

THEOREM 6.4.1. *Let $A_1$, ..., $A_n$, $A$ be associative $D$-algebras. Let*

$$f_i \in \mathcal{L}(D; A_i \to A) \quad i = 1, ..., n$$

$$a_j \in A \quad j = 0, ..., n$$

*For given transposition $\sigma$ of $n$ variables, the map*

(6.4.1)
$$((a_0, ..., a_n, \sigma) \circ (f_1, ..., f_n)) \circ (x_1, ..., x_n)$$
$$= (a_0 \sigma(f_1) a_1 ... a_{n-1} \sigma(f_n) a_n) \circ (x_1, ..., x_n)$$
$$= a_0 \sigma(f_1 \circ x_1) a_1 ... a_{n-1} \sigma(f_n \circ x_n) a_n$$

*is $n$-linear map into algebra $A$.*

---

[6.5]Generating set of representation is defined in definition 2.7.4.

PROOF. The statement of theorem follows from chains of equalities

$$((a_0, ..., a_n, \sigma) \circ (f_1, ..., f_n)) \circ (x_1, ..., x_i + y_i, ..., x_n)$$

$$= a_0 \sigma(f_1 \circ x_1) a_1 ... \sigma(f_i \circ (x_i + y_i)) ... a_{n-1} \sigma(f_n \circ x_n) a_n$$

$$= a_0 \sigma(f_1 \circ x_1) a_1 ... \sigma(f_i \circ x_i + f_i \circ y_i) ... a_{n-1} \sigma(f_n \circ x_n) a_n$$

$$= a_0 \sigma(f_1 \circ x_1) a_1 ... \sigma(f_i \circ x_i) ... a_{n-1} \sigma(f_n \circ x_n) a_n$$

$$+ a_0 \sigma(f_1 \circ x_1) a_1 ... \sigma(f_i \circ y_i) ... a_{n-1} \sigma(f_n \circ x_n) a_n$$

$$= ((a_0, ..., a_n, \sigma) \circ (f_1, ..., f_n)) \circ (x_1, ..., x_i, ..., x_n)$$

$$+ ((a_0, ..., a_n, \sigma) \circ (f_1, ..., f_n)) \circ (x_1, ..., y_i, ..., x_n)$$

$$((a_0, ..., a_n, \sigma) \circ (f_1, ..., f_n)) \circ (x_1, ..., p x_i, ..., x_n)$$

$$= a_0 \sigma(f_1 \circ x_1) a_1 ... \sigma(f_i \circ (p x_i)) ... a_{n-1} \sigma(f_n \circ x_n) a_n$$

$$= a_0 \sigma(f_1 \circ x_1) a_1 ... \sigma(p(f_i \circ x_i)) ... a_{n-1} \sigma(f_n \circ x_n) a_n$$

$$= p(a_0 \sigma(f_1 \circ x_1) a_1 ... \sigma(f_i \circ x_i) ... a_{n-1} \sigma(f_n \circ x_n) a_n)$$

$$= p(((a_0, ..., a_n, \sigma) \circ (f_1, ..., f_n)) \circ (x_1, ..., x_i, ..., x_n))$$

□

In the equality (6.4.1), as well as in other expressions of polylinear map, we have convention that map $f_i$ has variable $x_i$ as argument.

THEOREM 6.4.2. *Let $A_1$, ..., $A_n$, $A$ be associative $D$-algebras. For given set of maps*

$$f_i \in \mathcal{L}(D; A_i \to A) \quad i = 1, ..., n$$

*the map*

$$h : A^{n+1} \to \mathcal{L}(D; A_1 \times ... \times A_n \to A)$$

*defined by equality*

$$(a_0, ..., a_n, \sigma) \circ (f_1, ..., f_n) = a_0 \sigma(f_1) a_1 ... a_{n-1} \sigma(f_n) a_n$$

*is $n + 1$-linear map into $D$-module $\mathcal{L}(D; A_1 \times ... \times A_n \to A)$.*

PROOF. The statement of theorem follows from chains of equalities

$$((a_0, ..., a_i + b_i, ... a_n, \sigma) \circ (f_1, ..., f_n)) \circ (x_1, ..., x_n)$$

$$= a_0 \sigma(f_1 \circ x_1) a_1 ... (a_i + b_i) ... a_{n-1} \sigma(f_n \circ x_n) a_n$$

$$= a_0 \sigma(f_1 \circ x_1) a_1 ... a_i ... a_{n-1} \sigma(f_n \circ x_n) a_n + a_0 \sigma(f_1 \circ x_1) a_1 ... b_i ... a_{n-1} \sigma(f_n \circ x_n) a_n$$

$$= ((a_0, ..., a_i, ..., a_n, \sigma) \circ (f_1, ..., f_n)) \circ (x_1, ..., x_n)$$

$$+ ((a_0, ..., b_i, ..., a_n, \sigma) \circ (f_1, ..., f_n)) \circ (x_1, ..., x_n)$$

$$= ((a_0, ..., a_i, ..., a_n, \sigma) \circ (f_1, ..., f_n) + (a_0, ..., b_i, ..., a_n, \sigma) \circ (f_1, ..., f_n)) \circ (x_1, ..., x_n)$$

$$((a_0, ..., p a_i, ... a_n, \sigma) \circ (f_1, ..., f_n)) \circ (x_1, ..., x_n)$$

$$= a_0 \sigma(f_1 \circ x_1) a_1 ... p a_i ... a_{n-1} \sigma(f_n \circ x_n) a_n$$

$$= p(a_0 \sigma(f_1 \circ x_1) a_1 ... a_i ... a_{n-1} \sigma(f_n \circ x_n) a_n)$$

$$= p(((a_0, ..., a_i, ..., a_n, \sigma) \circ (f_1, ..., f_n)) \circ (x_1, ..., x_n))$$

$$= (p((a_0, ..., a_i, ..., a_n, \sigma) \circ (f_1, ..., f_n))) \circ (x_1, ..., x_n)$$

□

Theorem 6.4.3. *Let $A_1$, ..., $A_n$, $A$ be associative $D$-algebras. For given set of maps*

$$f_i \in \mathcal{L}(D; A_i \to A) \quad i = 1, ..., n$$

*there exists linear map*

$$h : A^{n+1 \otimes} \times S_n \to \mathcal{L}(D; A_1 \times ... \times A_n \to A)$$

*defined by the equality*

(6.4.2)
$$(a_0 \otimes ... \otimes a_n, \sigma) \circ (f_1, ..., f_n) = (a_0, ..., a_n, \sigma) \circ (f_1, ..., f_n)$$
$$= a_0 \sigma(f_1) a_1 ... a_{n-1} \sigma(f_n) a_n$$

Proof. The statement of the theorem is corollary of theorems 4.6.4, 6.4.2. □

Theorem 6.4.4. *Let $A_1$, ..., $A_n$, $A$ be associative $D$-algebras. For given tensor $a \in A^{n+1 \otimes}$ and given transposition $\sigma \in S_n$ the map*

$$h : \prod_{i=1}^{n} \mathcal{L}(D; A_i \to A) \to \mathcal{L}(D; A_1 \times ... \times A_n \to A)$$

*defined by equality*

$$(a_0 \otimes ... \otimes a_n, \sigma) \circ (f_1, ..., f_n) = a_0 \sigma(f_1) a_1 ... a_{n-1} \sigma(f_n) a_n$$

*is $n$-linear map into $D$-module $\mathcal{L}(D; A_1 \times ... \times A_n \to A)$.*

Proof. The statement of theorem follows from chains of equalities

$$((a_0 \otimes ... \otimes a_n, \sigma) \circ (f_1, ..., f_i + g_i, ..., f_n)) \circ (x_1, ..., x_n)$$
$$= (a_0 \sigma(f_1) a_1 ... \sigma(f_i + g_i) ... a_{n-1} \sigma(f_n) a_n) \circ (x_1, ..., x_n)$$
$$= a_0 \sigma(f_1 \circ x_1) a_1 ... \sigma((f_i + g_i) \circ x_i) ... a_{n-1} \sigma(f_n \circ x_n) a_n$$
$$= a_0 \sigma(f_1 \circ x_1) a_1 ... \sigma(f_i \circ x_i + g_i \circ x_i) ... a_{n-1} \sigma(f_n \circ x_n) a_n$$
$$= a_0 \sigma(f_1 \circ x_1) a_1 ... \sigma(f_i \circ x_i) ... a_{n-1} \sigma(f_n \circ x_n) a_n$$
$$+ a_0 \sigma(f_1 \circ x_1) a_1 ... \sigma(g_i \circ x_i) ... a_{n-1} \sigma(f_n \circ x_n) a_n$$
$$= (a_0 \sigma(f_1) a_1 ... \sigma(f_i) ... a_{n-1} \sigma(f_n) a_n) \circ (x_1, ..., x_n)$$
$$+ (a_0 \sigma(f_1) a_1 ... \sigma(g_i) ... a_{n-1} \sigma(f_n) a_n) \circ (x_1, ..., x_n)$$
$$= ((a_0 \otimes ... \otimes a_n, \sigma) \circ (f_1, ..., f_i, ..., f_n)) \circ (x_1, ..., x_n)$$
$$+ ((a_0 \otimes ... \otimes a_n, \sigma) \circ (f_1, ..., g_i, ..., f_n)) \circ (x_1, ..., x_n)$$
$$= ((a_0 \otimes ... \otimes a_n, \sigma) \circ (f_1, ..., f_i, ..., f_n)$$
$$+ (a_0 \otimes ... \otimes a_n, \sigma) \circ (f_1, ..., g_i, ..., f_n)) \circ (x_1, ..., x_n)$$

$$((a_0 \otimes ... \otimes a_n, \sigma) \circ (f_1, ..., pf_i, ..., f_n)) \circ (x_1, ..., x_n)$$
$$= (a_0 \sigma(f_1) a_1, ... \sigma(pf_i) ... a_{n-1} \sigma(f_n) a_n) \circ (x_1, ..., x_n)$$
$$= a_0 \sigma(f_1 \circ x_1) a_1, ... \sigma((pf_i) \circ x_i) ... a_{n-1} \sigma(f_n \circ x_n) a_n$$
$$= a_0 \sigma(f_1 \circ x_1) a_1, ... \sigma(p(f_i \circ x_i)) ... a_{n-1} \sigma(f_n \circ x_n) a_n$$
$$= p(a_0 \sigma(f_1 \circ x_1) a_1, ... \sigma(f_i \circ x_i) ... a_{n-1} \sigma(f_n \circ x_n) a_n)$$
$$= p(((a_0 \otimes ... \otimes a_n, \sigma) \circ (f_1, ..., f_i, ..., f_n)) \circ (x_1, ..., x_n))$$
$$= (p((a_0 \otimes ... \otimes a_n, \sigma) \circ (f_1, ..., f_i, ..., f_n))) \circ (x_1, ..., x_n)$$

□

THEOREM 6.4.5. *Let $A_1$, ..., $A_n$, $A$ be associative D-algebras. For given tensor* $a \in A^{n+1\otimes}$ *and given transposition $\sigma \in S_n$ there exists linear map*

$$h : \mathcal{L}(D; A_1 \to A) \otimes ... \otimes \mathcal{L}(D; A_n \to A) \to \mathcal{L}(D; A_1 \times ... \times A_n \to A)$$

*defined by the equality*

(6.4.3)        $(a_0 \otimes ... \otimes a_n, \sigma) \circ (f_1 \otimes ... \otimes f_n) = (a_0 \otimes ... \otimes a_n, \sigma) \circ (f_1, ..., f_n)$

PROOF. The statement of the theorem is corollary of theorems 4.6.4, 6.4.4.   $\square$

# Vector Space over Division Algebra

## 7.1. Effective Representation of $D$-Algebra

According to the theorem 3.3.4, considering the representation

$$f : A \longrightarrow_* V$$

of the $D$-algebra $A$ in the Abelian group $V$, we assume

$$(7.1.1) \qquad (f(a) + f(b))(x) = f(a)(x) + f(b)(x)$$

According to the definition 2.3.1, the map $f$ is homomorphism of $D$-algebra $A$. Therefore

$$(7.1.2) \qquad f(a + b) = f(a) + f(b)$$

The equality

$$(7.1.3) \qquad f(a + b)(x) = f(a)(x) + f(b)(x)$$

follows from equalities (7.1.1), (7.1.2).

THEOREM 7.1.1. *Representation*

$$f : A \longrightarrow_* V$$

*of the D-algebra A in an Abelian group V satisfies the equality*

$$(7.1.4) \qquad f(0) = \overline{0}$$

*where*

$$\overline{0} : V \to V$$

*is map such that*

$$(7.1.5) \qquad \overline{0} \circ v = 0$$

PROOF. The equality

$$(7.1.6) \qquad f(a)(x) = f(a + 0)(x) = f(a)(x) + f(0)(x)$$

follows from the equality (7.1.3). The equality

$$(7.1.7) \qquad f(0)(x) = 0$$

follows from the equality (7.1.6). The equality (7.1.5) follows from the equalities (7.1.4), (7.1.7). $\qquad \square$

THEOREM 7.1.2. *Representation*

$$f : A \longrightarrow_* V$$

*of the $D$-algebra $A$ in an Abelian group $V$ is **effective** iff $a = 0$ follows from
equation $f(a) = 0$.*

PROOF. Suppose $a$, $b \in A$ cause the same transformation. Then

$$(7.1.8) \qquad f(a)(m) = f(b)(m)$$

for any $m \in V$. From equalities (7.1.3), (7.1.8), it follows that

$$(7.1.9) \qquad f(a - b)(m) = 0$$

The equality

$$f(a - b) = \overline{0}$$

follows from equalities (7.1.5), (7.1.9). Therefore, the representation $f$ is effective
iff $a = b$.                                                                           □

## 7.2. Left Vector Space over Division Algebra

DEFINITION 7.2.1. *Let $A$ be associative division $D$-algebra. Let $V$ be $D$-vector
space. Let in $D$-vector space $\operatorname{End}(D, V)$ the product of endomorphisms is defined
as composition of maps. Let there exist homomorphism*

$$(7.2.1) \qquad g_{34} : A \to \operatorname{End}(D, V)$$

*of $D$-algebra $A$ into $D$-algebra $\operatorname{End}(D, V)$.*
*Effective left-side representation*

$$(7.2.2) \qquad g_{34} : A \twoheadrightarrow V \quad g_{34}(a) : v \in V \to av \in V \quad a \in A$$

*of $D$-algebra $A$ in $D$-vector space $V$ is called **left vector space** over $D$-algebra
$A$. We will also say that $D$-vector space $V$ is **left $A$-vector space**. $V$-number is
called **vector**. Bilinear map*

$$(7.2.3) \qquad (a, v) \in A \times V \to av \in V$$

*generated by left-side representation*

$$(7.2.4) \qquad (a, v) \to av$$

*is called left-side product of vector over scalar.*                                   □

THEOREM 7.2.2. *Let $D$ be commutative ring with unit. Let $A$ be associative
division $D$-algebra. The representation generated by left shift*

$$f : A \twoheadrightarrow A \quad f(a) : b \in A \to ab \in A$$

*is left vector space $\mathcal{L}A$ over $D$-algebra $A$.*

PROOF. According to definitions 7.2.1, 5.1.1, $D$-algebra $A$ is Abelian group.
From the equality

$$a(b + c) = ab + ac$$

and definitions 2.1.7, 2.1.9, it follows that the map

$$(7.2.5) \qquad f(a) : b \in A \to ab \in A$$

is endomorphism of Abelian group $A$. From equalities

$$(a_1 + a_2)b = a_1 b + a_2 b$$

$$(a_2 a_1)b = a_2(a_1 b)$$

and definitions 2.1.7, 2.3.1, it follows that the map $f$ is representation of $D$-algebra $A$ in Abelian group $A$. The theorem follows from the definition 7.2.1. □

**Theorem 7.2.3.** *Let $D$ be commutative ring with unit. Let $A$ be associative division $D$-algebra. For any set $B$, the representation*

$$f : A \longrightarrow A^B \quad f(a) : g \in A^B \to ag \in A^B \quad (ag)(b) = ag(b)$$

*is left vector space $\mathcal{L}A^B$ over $D$-algebra $A$.*

**Proof.** According to the theorem 5.1.33, the set of maps $A^B$ is Abelian group. From the equality

$$a(h(b) + g(b)) = ah(b) + ag(b)$$

and definitions 2.1.7, 2.1.9, it follows that the map $f(a)$ is endomorphism of Abelian group $A^B$. From equalities

$$(a_1 + a_2)h(b) = a_1 h(b) + a_2 h(b)$$

$$(a_2 a_1)h(b) = a_2(a_1 h(b))$$

and definitions 2.1.7, 2.3.1, it follows that the map $f$ is representation of $D$-algebra $A$ in Abelian group $A^B$. The theorem follows from the definition 7.2.1. □

**Theorem 7.2.4.** *The following diagram of representations describes left $A$-vector space $V$*

$$(7.2.6)$$

$$g_{12}(d) : a \to d\,a$$
$$g_{23}(v) : w \to C(w,v)$$
$$C \in \mathcal{L}(A^2 \to A)$$
$$g_{34}(a) : v \to a\,v$$
$$g_{14}(d) : v \to d\,v$$

*The diagram of representations* (7.2.6) *holds* **commutativity of representations** *of commutative ring $D$ and $D$-algebra $A$ in Abelian group $V$*

$$(7.2.7) \qquad a(dv) = d(av)$$

**Proof.** The diagram of representations (7.2.6) follows from the definition 7.2.1 and from the theorem 5.3.2. Since left-side transformation $g_{3,4}(a)$ is endomorphism of $D$-vector space $V$, we obtain the equality (7.2.7). □

**Definition 7.2.5.** *Subrepresentation of left $A$-vector space $V$ is called* **subspace** *of left $A$-vector space $V$.* □

**Theorem 7.2.6.** *Let $V$ be left $A$-vector space. Following conditions hold for $V$-numbers:*

7.2.6.1: **commutative law**

(7.2.8) $$v + w = w + v$$

7.2.6.2: **associative law**

(7.2.9) $$(pq)v = p(qv)$$

7.2.6.3: **distributive law**

(7.2.10) $$p(v + w) = pv + pw$$

(7.2.11) $$(p + q)v = pv + qv$$

7.2.6.4: **unitarity law**

(7.2.12) $$1v = v$$

*for any* $p, q \in A, v, w \in V$.

PROOF. The equality (7.2.8) follows from the definition 7.2.1.
The equality (7.2.10) follows from the theorem 3.3.2.
The equality (7.2.11) follows from the theorem 3.3.4.
The equality (7.2.9) follows from the theorem 2.1.10.
The equality (7.2.12) follows fro theorems 3.4.12, 7.2.4 because representation $g_{3,4}$ is left-side representation of the multiplicative group of $D$-algebra $A$,           □

## 7.3. Basis of Left $A$-Vector Space

THEOREM 7.3.1. *Let $V$ be left $A$-vector space. For any set of $V$-numbers*

(7.3.1) $$v = (v(i) \in V, i \in I)$$

*vector generated by the diagram of representations* (7.2.6) *has the following form* [7.1]

(7.3.2) $$J(v) = \{w : w = c(i)v(i), c(i) \in A, |\{i : c(i) \neq 0\}| < \infty\}$$

PROOF. We prove the theorem by induction based on the theorem 2.7.5.
For any $v(k) \in J(v)$, let

$$c(i) = \begin{cases} 1, & i = k \\ 0 \end{cases}$$

Then

(7.3.3) $$v(k) = \sum_{i \in I} c(i)v(i) \quad c(i) \in A$$

From equalities (7.3.2), (7.3.3), it follows that the theorem is true on the set $X_0 = v \subseteq J(v)$.
Let $X_{k-1} \subseteq J(v)$. According to the definition 7.2.1 and the theorem 2.7.5, if $w \in X_k$, then or $w = w_1 + w_2$, $w_1, w_2 \in X_{k-1}$, either $w = aw_1$, $a \in A$, $w_1 \in X_{k-1}$.

---

[7.1] For a set $A$, we denote by $|A|$ the cardinal number of the set $A$. The notation $|A| < \infty$ means that the set $A$ is finite.

Lemma 7.3.2. *Let* $w = w_1 + w_2$, $w_1$, $w_2 \in X_{k-1}$. *Then* $w \in J(v)$.

Proof. According to the equality (7.3.2), there exist $A$-numbers $w_1(i)$, $w_2(i)$, $i \in I$, such that

$$(7.3.4) \qquad w_1 = \sum_{i \in I} w_1(i)v(i) \qquad w_2 = \sum_{i \in I} w_2(i)v(i)$$

where sets

$$(7.3.5) \qquad H_1 = \{i \in I : w_1(i) \neq 0\} \qquad H_2 = \{i \in I : w_2(i) \neq 0\}$$

are finite. Since $V$ is left $A$-vector space, then from equalities (7.2.11), (7.3.4), it follows that

$$(7.3.6) \qquad \begin{aligned} w_1 + w_2 &= \sum_{i \in I} w_1(i)v(i) + \sum_{i \in I} w_2(i)v(i) = \sum_{i \in I}(w_1(i)v(i) + w_2(i)v(i)) \\ &= \sum_{i \in I}(w_1(i) + w_2(i))v(i) \end{aligned}$$

From the equality (7.3.5), it follows that the set

$$\{i \in I : w_1(i) + w_2(i) \neq 0\} \subseteq H_1 \cup H_2$$

is finite. $\odot$

Lemma 7.3.3. *Let* $w = aw_1$, $a \in A$, $w_1 \in X_{k-1}$. *Then* $w \in J(v)$.

Proof. According to the statement 2.7.5.4, for any $A$-number $a$,

$$(7.3.7) \qquad aw_1 \in X_k$$

According to the equality (7.3.2), there exist $A$-numbers $w(i)$, $i \in I$, such that

$$(7.3.8) \qquad w_1 = \sum_{i \in I} w_1(i)v(i)$$

where

$$(7.3.9) \qquad |\{i \in I : w_1(i) \neq 0\}| < \infty$$

From the equality (7.3.8) it follows that

$$(7.3.10) \qquad aw_1 = a\sum_{i \in I} w_1(i)v(i) = \sum_{i \in I} a(w_1(i)v(i)) = \sum_{i \in I}(aw_1(i))v(i)$$

From the statement (7.3.9), it follows that the set $\{i \in I : aw(i) \neq 0\}$ is finite. $\odot$

From lemmas 7.3.2, 7.3.3, it follows that $X_k \subseteq J(v)$. $\square$

Convention 7.3.4. *We will use summation convention in which repeated index in linear combination implies summation with respect to repeated index. In this case we assume that we know the set of summation index and do not use summation symbol*

$$c(i)v(i) = \sum_{i \in I} c(i)v(i)$$

*If needed to clearly show a set of indices, I will do it.*                    □

**Definition 7.3.5.** *Let $V$ be left-side vector space. Let*

$$v = (v(i) \in V, i \in I)$$

*be set of vectors. The expression $c(i)v(i)$, $c(i) \in A$, is called* **linear combination** *of vectors $v(i)$. A vector $w = c(i)v(i)$ is called* **linearly dependent** *on vectors $v(i)$.*                    □

**Convention 7.3.6.** *If it is necessary to explicitly show that we multiply vector $v(i)$ over $A$-number $c(i)$ on the left, then we will call the expression $c(i)v(i)$* **left linear combination***. We will use this convention to similar terms. For instance, we will say that a vector $w = c(i)v(i)$ linearly depends on vectors $v(i)$, $i \in I$, on the left.*                    □

**Theorem 7.3.7.** *Let $v = (v_i \in V, i \in I)$ be set of vectors of left $A$-vector space $V$. If vectors $v_i$, $i \in I$, belongs subspace $V'$ of left $A$-vector space $V$, then linear combination of vectors $v_i$, $i \in I$, belongs subspace $V'$.*

Proof. The theorem follows from the theorem 7.3.1 and definitions 7.3.5, 7.2.5.                    □

**Definition 7.3.8.** *$J(v)$ is called vector subspace generated by set $v$, and $v$ is a* **generating set** *of vector subspace $J(v)$. In particular, a* **generating set** *of left $A$-vector space $V$ is a subset $X \subset V$ such that $J(X) = V$.*                    □

**Definition 7.3.9.** *If the set $X \subset V$ is generating set of left $A$-vector space $V$, then any set $Y$, $X \subset Y \subset V$ also is generating set of left $A$-vector space $V$. If there exists minimal set $X$ generating the left $A$-vector space $V$, then the set $X$ is called* **quasi-basis** *of left $A$-vector space $V$.*                    □

**Definition 7.3.10.** *Let $\overline{\overline{e}}$ be the quasibasis of left $A$-vector space $V$ and vector $\overline{v} \in V$ has expansion*

$$(7.3.11) \qquad\qquad \overline{v} = v(i)e(i)$$

*with respect to the quasibasis $\overline{\overline{e}}$. $A$-numbers $v(i)$ are called* **coordinates** *of vector $\overline{v}$ with respect to the quasibasis $\overline{\overline{e}}$. Matrix of $A$-numbers $v = (v(i), i \in I)$ is called* **coordinate matrix of vector** $\overline{v}$ *in quasibasis $\overline{\overline{e}}$.*                    □

**Theorem 7.3.11.** *Let $A$ be associative division $D$-algebra. Since the equation*

$$(7.3.12) \qquad\qquad w(i)v(i) = 0$$

*implies existence of index $i = j$ such that $w(j) \neq 0$, then the vector $v(j)$ linearly depends on rest of vectors $v$.*

Proof. The theorem follows from the equality

$$v(j) = \sum_{i \in I \setminus \{j\}} (w(j))^{-1} w(i) v(i)$$

and from the definition 7.3.5.                                                                □

It is evident that for any set of vectors $v(i)$

$$w(i) = 0 \Rightarrow w(i)v(i) = 0$$

Definition 7.3.12. *The set of vectors* $^{7.2}v(i)$, $i \in I$, *of left $A$-vector space* $V$ *is* **linearly independent** *if* $c(i) = 0$, $i \in I$, *follows from the equation*

(7.3.13)                                      $$c(i)v(i) = 0$$

*Otherwise the set of vectors* $v(i)$, $i \in I$, *is* **linearly dependent**.             □

Theorem 7.3.13. *Let $A$ be associative division $D$-algebra. The set of vectors* $\overline{e} = (e(i), i \in I)$ *is a* **basis of left $A$-vector space** $V$ *if vectors $e(i)$ are linearly independent and any vector $v \in V$ linearly depends on vectors $e(i)$.*

Proof. Let the set of vectors $e(i)$, $i \in I$, be linear dependent. Then the equation

$$w(i)e(i) = 0$$

implies existence of index $i = j$ such that $w(j) \neq 0$. According to the theorem 7.3.11, the vector $e(j)$ linearly depends on rest of vectors of the set $\overline{e}$. According to the definition 7.3.9, the set of vectors $e(i)$, $i \in I$, is not a basis for left $A$-vector space $V$.

Therefore, if the set of vectors $e(i)$, $i \in I$, is a basis, then these vectors are linearly independent. Since an arbitrary vector $v \in V$ is linear combination of vectors $e(i)$, $i \in I$, then the set of vectors $v$, $e(i)$, $i \in I$, is not linearly independent.                                                                □

Theorem 7.3.14. *Coordinates of vector $v \in V$ relative to basis $\overline{\overline{e}}$ of left $A$-vector space $V$ are uniquely defined. The equality*

(7.3.14)                          $$ve = we \qquad v(i)e(i) = w(i)e(i)$$

*implies the equality*

$$v = w \qquad v(i) = w(i)$$

Proof. According to the theorem 7.3.13, the system of vectors $\overline{v}$, $e(i)$, $i \in I$, is linearly dependent and in equality

(7.3.15)                                  $$b\overline{v} + c(i)e(i) = 0$$

at least $b$ is different from 0. Then equality

(7.3.16)                                  $$\overline{v} = (-b^{-1}c(i))e(i)$$

follows from (7.3.15). The equality

(7.3.17)                          $$\overline{v} = v(i)e(i) \qquad v(i) = -b^{-1}c(i)$$

---

$^{7.2}$ I follow to the definition on page [1]-130.

follows from the equality (7.3.16).

Assume we get another expansion

(7.3.18)                         $\overline{v} = v'(i)e(i)$

We subtract (7.3.17) from (7.3.18) and get

(7.3.19)                         $0 = (v'(i) - v(i))e(i)$

Since vectors $e(i)$ are linearly independent, then the equality

(7.3.20)                    $v'(i) - v(i) = 0 \quad i \in I$

follows from the equality (7.3.19). Therefore, the theorem follows from the equality
(7.3.20).                                                                          $\square$

THEOREM 7.3.15. *We can consider division $D$-algebra $A$ as left $A$-vector space
with basis $\overline{\overline{e}} = \{1\}$.*

PROOF. The theorem follows from the equality (7.2.12).                          $\square$

## 7.4. Right Vector Space over Division Algebra

DEFINITION 7.4.1. *Let $A$ be associative division $D$-algebra. Let $V$ be $D$-vector
space. Let in $D$-vector space $\operatorname{End}(D, V)$ the product of endomorphisms is defined
as composition of maps. Let there exist homomorphism*

(7.4.1)                         $g_{34} : A \to \operatorname{End}(D, V)$

*of $D$-algebra $A$ into $D$-algebra $\operatorname{End}(D, V)$.*

*Effective right-side representation*

(7.4.2)        $g_{34} : A \twoheadrightarrow V \quad g_{34}(a) : v \in V \to va \in V \quad a \in A$

*of $D$-algebra $A$ in $D$-vector space $V$ is called* **right vector space** *over $D$-algebra
$A$. We will also say that $D$-vector space $V$ is* **right $A$-vector space**. *$V$-number
is called* **vector**. *Bilinear map*

(7.4.3)                         $(v, a) \in V \times A \to va \in V$

*generated by right-side representation*

(7.4.4)                         $(v, a) \to va$

*is called right-side product of vector over scalar.*                           $\square$

THEOREM 7.4.2. *Let $D$ be commutative ring with unit. Let $A$ be associative
division $D$-algebra. The representation generated by right shift*

$$f : A \twoheadrightarrow A \quad f(a) : b \in A \to ba \in A$$

*is right vector space $\mathcal{R}A$ over $D$-algebra $A$.*

PROOF. According to definitions 7.4.1, 5.1.1, $D$-algebra $A$ is Abelian group.
From the equality

$$(b + c)a = ba + ca$$

and definitions 2.1.7, 2.1.9, it follows that the map

(7.4.5) $$f(a) : b \in A \to ba \in A$$

is endomorphism of Abelian group $A$. From equalities

$$b(a_1 + a_2) = ba_1 + ba_2$$
$$b(a_1 a_2) = (ba_1)a_2$$

and definitions 2.1.7, 2.3.1, it follows that the map $f$ is representation of $D$-algebra $A$ in Abelian group $A$. The theorem follows from the definition 7.4.1. □

THEOREM 7.4.3. *Let $D$ be commutative ring with unit. Let $A$ be associative division $D$-algebra. For any set $B$, the representation*

$$f : A \longrightarrow A^B \quad f(a) : g \in A^B \to ga \in A^B \quad (ga)(b) = g(b)a$$

*is right vector space $\mathcal{R}A^B$ over $D$-algebra $A$.*

PROOF. According to the theorem 5.1.33, the set of maps $A^B$ is Abelian group. From the equality

$$(h(b) + g(b))a = h(b)a + g(b)a$$

and definitions 2.1.7, 2.1.9, it follows that the map $f(a)$ is endomorphism of Abelian group $A^B$. From equalities

$$h(b)(a_1 + a_2) = h(b)a_1 + h(b)a_2$$
$$h(b)(a_1 a_2) = (h(b)a_1)a_2$$

and definitions 2.1.7, 2.3.1, it follows that the map $f$ is representation of $D$-algebra $A$ in Abelian group $A^B$. The theorem follows from the definition 7.4.1. □

THEOREM 7.4.4. *The following diagram of representations describes right $A$-vector space $V$*

(7.4.6)

$$g_{12}(d) : a \to d\,a$$
$$g_{23}(v) : w \to C(w,v)$$
$$C \in \mathcal{L}(A^2 \to A)$$
$$g_{34}(a) : v \to v\,a$$
$$g_{14}(d) : v \to v\,d$$

*The diagram of representations* (7.4.6) *holds* **commutativity of representations of commutative ring $D$ and $D$-algebra $A$ in Abelian group $V$**

(7.4.7) $$(vd)a = (va)d$$

PROOF. The diagram of representations (7.4.6) follows from the definition 7.4.1 and from the theorem 5.3.2. Since right-side transformation $g_{3,4}(a)$ is endomorphism of $D$-vector space $V$, we obtain the equality (7.4.7). □

DEFINITION 7.4.5. *Subrepresentation of right $A$-vector space $V$ is called* **subspace** *of right $A$-vector space $V$.* □

THEOREM 7.4.6. *Let $V$ be right $A$-vector space. Following conditions hold for $V$-numbers:*

7.4.6.1: **commutative law**

$$(7.4.8) \qquad\qquad v + w = w + v$$

7.4.6.2: **associative law**

$$(7.4.9) \qquad\qquad v(pq) = (vp)q$$

7.4.6.3: **distributive law**

$$(7.4.10) \qquad\qquad (v + w)p = vp + wp$$

$$(7.4.11) \qquad\qquad v(p + q) = vp + vq$$

7.4.6.4: **unitarity law**

$$(7.4.12) \qquad\qquad v1 = v$$

*for any $p, q \in A$, $v, w \in V$.*

PROOF. The equality (7.4.8) follows from the definition 7.4.1.
The equality (7.4.10) follows from the theorem 3.3.2.
The equality (7.4.11) follows from the theorem 3.3.4.
The equality (7.4.9) follows from the theorem 2.1.10.
The equality (7.4.12) follows fro theorems 3.4.12, 7.4.4 because representation $g_{3,4}$ is right-side representation of the multiplicative group of $D$-algebra $A$, $\qquad \square$

## 7.5. Basis of Right $A$-Vector Space

THEOREM 7.5.1. *Let $V$ be right $A$-vector space. For any set of $V$-numbers*

$$(7.5.1) \qquad\qquad v = (v(i) \in V, i \in I)$$

*vector generated by the diagram of representations (7.4.6) has the following form* [7.3]

$$(7.5.2) \qquad J(v) = \{w : w = v(i)c(i), c(i) \in A, |\{i : c(i) \neq 0\}| < \infty\}$$

PROOF. We prove the theorem by induction based on the theorem 2.7.5.
For any $v(k) \in J(v)$, let

$$c(i) = \begin{cases} 1, & i = k \\ 0 \end{cases}$$

Then

$$(7.5.3) \qquad\qquad v(k) = \sum_{i \in I} v(i)c(i) \quad c(i) \in A$$

From equalities (7.5.2), (7.5.3), it follows that the theorem is true on the set $X_0 = v \subseteq J(v)$.

---

[7.3] For a set $A$, we denote by $|A|$ the cardinal number of the set $A$. The notation $|A| < \infty$ means that the set $A$ is finite.

Let $X_{k-1} \subseteq J(v)$. According to the definition 7.4.1 and the theorem 2.7.5, if $w \in X_k$, then or $w = w_1 + w_2$, $w_1$, $w_2 \in X_{k-1}$, either $w = w_1 a$, $a \in A$, $w_1 \in X_{k-1}$.

LEMMA 7.5.2. *Let $w = w_1 + w_2$, $w_1$, $w_2 \in X_{k-1}$. Then $w \in J(v)$.*

PROOF. According to the equality (7.5.2), there exist $A$-numbers $w_1(i)$, $w_2(i)$, $i \in I$, such that

$$(7.5.4) \qquad w_1 = \sum_{i \in I} v(i) w_1(i) \quad w_2 = \sum_{i \in I} v(i) w_2(i)$$

where sets

$$(7.5.5) \qquad H_1 = \{ i \in I : w_1(i) \neq 0 \} \quad H_2 = \{ i \in I : w_2(i) \neq 0 \}$$

are finite. Since $V$ is right $A$-vector space, then from equalities (7.4.11), (7.5.4), it follows that

$$
\begin{aligned}
(7.5.6) \quad w_1 + w_2 &= \sum_{i \in I} v(i) w_1(i) + \sum_{i \in I} v(i) w_2(i) = \sum_{i \in I} (v(i) w_1(i) + v(i) w_2(i)) \\
&= \sum_{i \in I} v(i) (w_1(i) + w_2(i))
\end{aligned}
$$

From the equality (7.5.5), it follows that the set

$$\{ i \in I : w_1(i) + w_2(i) \neq 0 \} \subseteq H_1 \cup H_2$$

is finite.                                                                                      $\odot$

LEMMA 7.5.3. *Let $w = w_1 a$, $a \in A$, $w_1 \in X_{k-1}$. Then $w \in J(v)$.*

PROOF. According to the statement 2.7.5.4, for any $A$-number $a$,

$$(7.5.7) \qquad w_1 a \in X_k$$

According to the equality (7.5.2), there exist $A$-numbers $w(i)$, $i \in I$, such that

$$(7.5.8) \qquad w_1 = \sum_{i \in I} v(i) w_1(i)$$

where

$$(7.5.9) \qquad |\{ i \in I : w_1(i) \neq 0 \}| < \infty$$

From the equality (7.5.8) it follows that

$$(7.5.10) \qquad w a_1 = \left( \sum_{i \in I} v(i) w_1(i) \right) a = \sum_{i \in I} (v(i) w_1(i)) a = \sum_{i \in I} v(i) (w_1(i) a)$$

From the statement (7.5.9), it follows that the set $\{ i \in I : w(i) a \neq 0 \}$ is finite.                                                                                      $\odot$

From lemmas 7.5.2, 7.5.3, it follows that $X_k \subseteq J(v)$.                    $\square$

CONVENTION 7.5.4. *We will use summation convention in which repeated index in linear combination implies summation with respect to repeated index. In this case we assume that we know the set of summation index and do not use summation symbol*

$$v(i)c(i) = \sum_{i \in I} v(i)c(i)$$

*If needed to clearly show a set of indices, I will do it.*                    □

DEFINITION 7.5.5. *Let $V$ be right-side vector space. Let*

$$v = (v(i) \in V, i \in I)$$

*be set of vectors. The expression $v(i)c(i)$, $c(i) \in A$, is called* **linear combination** *of vectors $v(i)$. A vector $w = v(i)c(i)$ is called* **linearly dependent** *on vectors $v(i)$.*                    □

CONVENTION 7.5.6. *If it is necessary to explicitly show that we multiply vector $v(i)$ over $A$-number $c(i)$ on the right, then we will call the expression $v(i)c(i)$* **right linear combination**. *We will use this convention to similar terms. For instance, we will say that a vector $w = v(i)c(i)$ linearly depends on vectors $v(i)$, $i \in I$, on the right.*                    □

THEOREM 7.5.7. *Let $v = (v_i \in V, i \in I)$ be set of vectors of right $A$-vector space $V$. If vectors $v_i$, $i \in I$, belongs subspace $V'$ of right $A$-vector space $V$, then linear combination of vectors $v_i$, $i \in I$, belongs subspace $V'$.*

PROOF. The theorem follows from the theorem 7.5.1 and definitions 7.5.5, 7.4.5.                    □

DEFINITION 7.5.8. *$J(v)$ is called vector subspace generated by set $v$, and $v$ is a* **generating set** *of vector subspace $J(v)$. In particular, a* **generating set** *of right $A$-vector space $V$ is a subset $X \subset V$ such that $J(X) = V$.*                    □

DEFINITION 7.5.9. *If the set $X \subset V$ is generating set of right $A$-vector space $V$, then any set $Y$, $X \subset Y \subset V$ also is generating set of right $A$-vector space $V$. If there exists minimal set $X$ generating the right $A$-vector space $V$, then the set $X$ is called* **quasi-basis** *of right $A$-vector space $V$.*                    □

DEFINITION 7.5.10. *Let $\overline{\overline{e}}$ be the quasibasis of right $A$-vector space $V$ and vector $\overline{v} \in V$ has expansion*

(7.5.11)                    $\overline{v} = e(i)v(i)$

*with respect to the quasibasis $\overline{\overline{e}}$. $A$-numbers $v(i)$ are called* **coordinates** *of vector $\overline{v}$ with respect to the quasibasis $\overline{\overline{e}}$. Matrix of $A$-numbers $v = (v(i), i \in I)$ is called* **coordinate matrix of vector** $\overline{v}$ *in quasibasis $\overline{\overline{e}}$.*                    □

THEOREM 7.5.11. *Let $A$ be associative division $D$-algebra. Since the equation*

(7.5.12) $$v(i)w(i) = 0$$

*implies existence of index $i = j$ such that $w(j) \neq 0$, then the vector $v(j)$ linearly depends on rest of vectors $v$.*

PROOF. The theorem follows from the equality

$$v(j) = \sum_{i \in I \setminus \{j\}} v(i)w(i)(w(j))^{-1}$$

and from the definition 7.5.5.                                                □

It is evident that for any set of vectors $v(i)$

$$w(i) = 0 \Rightarrow v(i)w(i) = 0$$

DEFINITION 7.5.12. *The set of vectors [7.4] $v(i)$, $i \in I$, of right $A$-vector space $V$ is* **linearly independent** *if $c(i) = 0$, $i \in I$, follows from the equation*

(7.5.13) $$v(i)c(i) = 0$$

*Otherwise the set of vectors $v(i)$, $i \in I$, is* **linearly dependent**.                □

THEOREM 7.5.13. *Let $A$ be associative division $D$-algebra. The set of vectors $\overline{\overline{e}} = (e(i), i \in I)$ is a* **basis of right $A$-vector space** *$V$ if vectors $e(i)$ are linearly independent and any vector $v \in V$ linearly depends on vectors $e(i)$.*

PROOF. Let the set of vectors $e(i)$, $i \in I$, be linear dependent. Then the equation

$$e(i)w(i) = 0$$

implies existence of index $i = j$ such that $w(j) \neq 0$. According to the theorem 7.5.11, the vector $e(j)$ linearly depends on rest of vectors of the set $\overline{\overline{e}}$. According to the definition 7.5.9, the set of vectors $e(i)$, $i \in I$, is not a basis for right $A$-vector space $V$.

Therefore, if the set of vectors $e(i)$, $i \in I$, is a basis, then these vectors are linearly independent. Since an arbitrary vector $v \in V$ is linear combination of vectors $e(i)$, $i \in I$, then the set of vectors $v$, $e(i)$, $i \in I$, is not linearly independent.                                                □

THEOREM 7.5.14. *Coordinates of vector $v \in V$ relative to basis $\overline{\overline{e}}$ of right $A$-vector space $V$ are uniquely defined. The equality*

(7.5.14) $$ev = ew \qquad e(i)v(i) = e(i)w(i)$$

*implies the equality*

$$v = w \qquad v(i) = w(i)$$

----

[7.4] I follow to the definition on page [1]-130.

PROOF. According to the theorem 7.5.13, the system of vectors $\overline{v}, e(i), i \in I$, is linearly dependent and in equality

$$(7.5.15) \qquad \overline{v}b + e(i)c(i) = 0$$

at least $b$ is different from 0. Then equality

$$(7.5.16) \qquad \overline{v} = e(i)(-c(i)b^{-1})$$

follows from (7.5.15). The equality

$$(7.5.17) \qquad \overline{v} = e(i)v(i) \quad v(i) = -c(i)b^{-1}$$

follows from the equality (7.5.16).

Assume we get another expansion

$$(7.5.18) \qquad \overline{v} = e(i)v'(i)$$

We subtract (7.5.17) from (7.5.18) and get

$$(7.5.19) \qquad 0 = e(i)(v'(i) - v(i))$$

Since vectors $e(i)$ are linearly independent, then the equality

$$(7.5.20) \qquad v'(i) - v(i) = 0 \quad i \in I$$

follows from the equality (7.5.19). Therefore, the theorem follows from the equality (7.5.20). $\square$

# Vector Space Type

## 8.1. Biring

We consider matrices whose entries belong to associative division $D$-algebra $A$.

According to the custom the product of matrices $a$ and $b$ is defined as product of rows of the matrix $a$ and columns of the matrix $b$. In non-commutative algebra, this product is not enough to solve some problems.

EXAMPLE 8.1.1. *We represent the set of vectors* $e(1) = e_1$, ..., $e(n) = e_n$ *of basis* $\overline{e}$ *of left vector space* $V$ *over* $D$-*algebra* $A$ *(see the definition 7.2.1 and the theorem 7.3.13) as row of matrix*

$$(8.1.1) \qquad e = \begin{pmatrix} e_1 & ... & e_n \end{pmatrix}$$

*We represent coordinates* $v(1) = v^1$, ..., $v(n) = v^n$ *of vector* $\overline{v} = v^i e_i$ *as column of matrix*

$$(8.1.2) \qquad v = \begin{pmatrix} v^1 \\ ... \\ v^n \end{pmatrix}$$

*Therefore, we can represent the vector* $\overline{v}$ *as product of matrices*

$(8.1.3)$

$$v = \begin{pmatrix} e_1 & ... & e_n \end{pmatrix} \begin{pmatrix} v^1 \\ ... \\ v^n \end{pmatrix} = e_i v^i$$

□

EXAMPLE 8.1.2. *We represent the set of vectors* $e(1) = e_1$, ..., $e(n) = e_n$ *of basis* $\overline{e}$ *of right vector space* $V$ *over* $D$-*algebra* $A$ *(see the definition 7.4.1 and the theorem 7.3.13) as row of matrix*

$$(8.1.4) \qquad e = \begin{pmatrix} e_1 & ... & e_n \end{pmatrix}$$

*We represent coordinates* $v(1) = v^1$, ..., $v(n) = v^n$ *of vector* $\overline{v} = v^i e_i$ *as column of matrix*

$$(8.1.5) \qquad v = \begin{pmatrix} v^1 \\ ... \\ v^n \end{pmatrix}$$

*However, we cannot represent the vector* $\overline{v}$ *as product of matrices*

$(8.1.6)$

$$v = \begin{pmatrix} v^1 \\ ... \\ v^n \end{pmatrix} \qquad e = \begin{pmatrix} e_1 & ... & e_n \end{pmatrix}$$

*because this product is not defined.*

□

From examples 8.1.1, 8.1.2 it follows that we cannot confine ourselves to traditional product of matrices and we need to define two products of matrices. To distinguish between these products we introduced a new notation.

DEFINITION 8.1.3. *Let the nubmer of columns of the matrix* $a$ *equal the number*

*of rows of the matrix* $b$. $_*{}^*$**-product** *of matrices* $a$ *and* $b$ *has form*

$$(8.1.7) \qquad a_*{}^*b = \left(a_k^i b_j^k\right)$$

$$(8.1.8) \qquad (a_*{}^*b)_j^i = a_k^i b_j^k$$

$$(8.1.9) \qquad \begin{pmatrix} a_1^1 & ... & a_p^1 \\ ... & ... & ... \\ a_1^n & ... & a_p^n \end{pmatrix} {}_*{}^* \begin{pmatrix} b_1^1 & ... & b_m^1 \\ ... & ... & ... \\ b_1^p & ... & b_m^p \end{pmatrix} = \begin{pmatrix} a_k^1 b_1^k & ... & a_k^1 b_m^k \\ ... & ... & ... \\ a_k^n b_1^k & ... & a_k^n b_m^k \end{pmatrix}$$

$$= \begin{pmatrix} (a_*{}^*b)_1^1 & ... & (a_*{}^*b)_m^1 \\ ... & ... & ... \\ (a_*{}^*b)_1^n & ... & (a_*{}^*b)_m^n \end{pmatrix}$$

$_*{}^*$*-product can be expressed as product of a row of the matrix* $a$ *over a column of the matrix* $b$. $\qquad\square$

DEFINITION 8.1.4. *Let the nubmer of rows of the matrix* $a$ *equal the number of columns of the matrix* $b$. $^*{}_*$**-product** *of matrices* $a$ *and* $b$ *has form*

$$(8.1.10) \qquad a^*{}_*b = \left(a_i^k b_k^j\right)$$

$$(8.1.11) \qquad (a^*{}_*b)_j^i = a_k^i b_k^j$$

$$(8.1.12) \qquad \begin{pmatrix} a_1^1 & ... & a_m^1 \\ ... & ... & ... \\ a_1^p & ... & a_m^p \end{pmatrix} {}^*{}_* \begin{pmatrix} b_1^1 & ... & b_p^1 \\ ... & ... & ... \\ b_1^n & ... & b_p^n \end{pmatrix} = \begin{pmatrix} a_1^k b_k^1 & ... & a_m^k b_k^1 \\ ... & ... & ... \\ a_1^k b_k^n & ... & a_m^k b_k^n \end{pmatrix}$$

$$= \begin{pmatrix} (a^*{}_*b)_1^1 & ... & (a^*{}_*b)_m^1 \\ ... & ... & ... \\ (a^*{}_*b)_1^n & ... & (a^*{}_*b)_m^n \end{pmatrix}$$

$^*{}_*$*-product can be expressed as product of a column of the matrix* $a$ *over a row of the matrix* $b$. $\qquad\square$

We also consider following operations on the set of matrices.

DEFINITION 8.1.5. *The transpose* $a^T$ *of the matrix* $a$ *exchanges rows and columns*

$$(8.1.13) \qquad (a^T)_j^i = a_i^j$$

$\qquad\square$

THEOREM 8.1.6. *Double transpose is original matrix*

$$(8.1.14) \qquad (a^T)^T = a$$

Proof. The theorem follows from the definition 8.1.5. □

Definition 8.1.7. *The sum of matrices $a$ and $b$ is defined by the equality*

$$(8.1.15) \qquad (a+b)^i_j = a^i_j + b^i_j$$

□

Remark 8.1.8. *We will use symbol $_*{}^*$- or $^*{}_*$- in name of properties of each product and in the notation. We can read the symbol $_*{}^*$ as rc-product (product of row over column) and the symbol $^*{}_*$ as cr-product (product of column over row). In order to keep this notation consistent with the existing one we assume that we have in mind $_*{}^*$-product when no clear notation is present. This rule we extend to following terminology. To draw symbol of product, we put two symbols of product in the place of index which participate in sum. For instance, if product of A-numbers has form $a \circ b$, then $_*{}^*$-product of matrices $a$ and $b$ has form $a_\circ{}^\circ b$ and $^*{}_*$-product of matrices $a$ and $b$ has form $a^\circ{}_\circ b$.* □

Theorem 8.1.9. *The $_*{}^*$-product in $\Omega$-ring is **distributive** over addition*

$$(8.1.16) \quad a_*{}^*(b_1 + b_2) = a_*{}^*b_1 + a_*{}^*b_2$$

$$(8.1.17) \quad (b_1 + b_2)_*{}^*a = b_1{}_*{}^*a + b_2{}_*{}^*a$$

Proof. The equality

$$
\begin{aligned}
&(a_*{}^*(b_1+b_2))^i_j\\
&= a^i_k(b_1+b_2)^k_j\\
(8.1.20) \quad &= a^i_k(b_1{}^k_j + b_2{}^k_j)\\
&= a^i_k b_1{}^k_j + a^i_k b_2{}^k_j\\
&= (a_*{}^*b_1)^i_j + (a_*{}^*b_2)^i_j
\end{aligned}
$$

follows from equalities (5.1.3), (8.1.7), (8.1.15). The equality (8.1.16) follows from the equality (8.1.20). The equality

$$
\begin{aligned}
&((b_1+b_2)_*{}^*a)^i_j\\
&= (b_1+b_2)^i_k a^k_j\\
(8.1.21) \quad &= (b_1{}^i_k + b_2{}^i_k)a^k_j\\
&= b_1{}^i_k a^k_j + b_2{}^i_k a^k_j\\
&= (b_1{}_*{}^*a)^i_j + (b_2{}_*{}^*a)^i_j
\end{aligned}
$$

follows from equalities (5.1.2), (8.1.7), (8.1.15). The equality (8.1.17) follows from the equality (8.1.21). □

Theorem 8.1.10. *The $^*{}_*$-product in $\Omega$-ring is **distributive** over addition*

$$(8.1.18) \quad a^*{}_*(b_1 + b_2) = a^*{}_*b_1 + a^*{}_*b_2$$

$$(8.1.19) \quad (b_1 + b_2)^*{}_*a = b_1{}^*{}_*a + b_2{}^*{}_*a$$

Proof. The equality

$$
\begin{aligned}
&(a^*{}_*(b_1+b_2))^i_j\\
&= a^k_j(b_1+b_2)^i_k\\
(8.1.22) \quad &= a^k_j(b_1{}^i_k + b_2{}^i_k)\\
&= a^k_j b_1{}^i_k + a^k_j b_2{}^i_k\\
&= (a^*{}_*b_1)^i_j + (a^*{}_*b_2)^i_j
\end{aligned}
$$

follows from equalities (5.1.3), (8.1.10), (8.1.15). The equality (8.1.18) follows from the equality (8.1.22). The equality

$$
\begin{aligned}
&((b_1+b_2)^*{}_*a)^i_j\\
&= (b_1+b_2)^k_j a^i_k\\
(8.1.23) \quad &= (b_1{}^k_j + b_2{}^k_j)a^i_k\\
&= b_1{}^k_j a^i_k + b_2{}^k_j a^i_k\\
&= (b_1{}^*{}_*a)^i_j + (b_2{}^*{}_*a)^i_j
\end{aligned}
$$

follows from equalities (5.1.2), (8.1.10), (8.1.15). The equality (8.1.19) follows from the equality (8.1.23). □

$$(8.1.24) \qquad (a_*{}^*b)^T = a^{T*}{}_*b^T$$

$$(8.1.25) \qquad (a^*{}_*b)^T = (a^T)_*{}^*(b^T)$$

Proof of theorem 8.1.11.    The chain of equalities

$$(8.1.26) \qquad ((a_*{}^*b)^T)^j_i = (a_*{}^*b)^i_j = a^i_k b^k_j = (a^T)^k_i (b^T)^j_k = ((a^T)^*{}_*(b^T))^j_i$$

follows from (8.1.13), (8.1.7) and (8.1.10).    The equality (8.1.24) follows from (8.1.26). □

Proof of theorem 8.1.12.    We can prove (8.1.25) in case of matrices the same way as we proved (8.1.24). However it is more important for us to show that (8.1.25) follows directly from (8.1.24).

Applying (8.1.14) to each term in left side of (8.1.25) we get

$$(8.1.27) \qquad (a^*{}_*b)^T = ((a^T)^{T*}{}_*(b^T)^T)^T$$

From (8.1.27) and (8.1.24) it follows that

$$(8.1.28) \qquad (a^*{}_*b)^T = ((a^T{}_*{}^*b^T)^T)^T$$

(8.1.25) follows from (8.1.28) and (8.1.14).                        □

Definition 8.1.13.    *Let*

$$a = \begin{pmatrix} a^1_1 & ... & a^1_n \\ ... & ... & ... \\ a^n_1 & ... & a^n_n \end{pmatrix}$$

*be a matrix of A-numbers. We call matrix*[8.1]

$$(8.1.29) \qquad \mathcal{H}a = \begin{pmatrix} (a^1_1)^{-1} & ... & (a^n_1)^{-1} \\ ... & ... & ... \\ (a^1_n)^{-1} & ... & (a^n_n)^{-1} \end{pmatrix}$$

$$(8.1.30) \qquad (\mathcal{H}a)^i_j = ((a^T)^i_j)^{-1} = (a^j_i)^{-1}$$

**Hadamard inverse of matrix** *a ([6]-page 4).*                □

Set of $n \times n$ matrices is closed relative to $_*{}^*$-product and $^*{}_*$-product as well relative to sum.

Definition 8.1.14.    *The set of matrices $\mathcal{A}$ is a* **biring** *if we defined on $\mathcal{A}$ an unary operation, say transpose, and three binary operations, say $_*{}^*$-product, $^*{}_*$-product and sum, such that*

- $_*{}^*$-product and sum define structure of ring on $\mathcal{A}$
- $^*{}_*$-product and sum define structure of ring on $\mathcal{A}$

---

[8.1]  The notation in the equality (8.1.30) means that we exchange rows and columns in Hadamard inverse. We can formally write right site of the equality (8.1.30) in the following form

$$(a^j_i)^{-1} = \frac{1}{a^i_j}$$

- *both products have common identity*

$$E_n = (\delta_j^i) \quad 1 \le i, j \le n$$

$\square$

THEOREM 8.1.15. *Transpose of identity is identity*

(8.1.31) $$E_n^T = E_n$$

PROOF. The theorem follows from definitions 8.1.5, 8.1.14. $\square$

DEFINITION 8.1.16. *We introduce* $*^*$-**power** *of the* $n \times n$ *matrix* $a$ *using recursive definition*

(8.1.32) $$a^{0*^*} = E_n$$

(8.1.33) $$a^{k*^*} = a^{k-1*^*} {}_*^* a$$

$\square$

DEFINITION 8.1.17. *We introduce* ${}_*^*$-**power** *of the* $n \times n$ *matrix* $a$ *using recursive definition*

(8.1.34) $$a^{0^*{}_*} = E_n$$

(8.1.35) $$a^{k^*{}_*} = a^{k-1^*{}_*} {}^*_* a$$

$\square$

THEOREM 8.1.18.

(8.1.36) $$a^{1*^*} = a$$

THEOREM 8.1.19.

(8.1.38) $$a^{1^*{}_*} = a$$

PROOF. The equality

(8.1.37) $$a^{1*^*} = a^{0*^*} {}_*^* a = E_n {}_*^* a = a$$

follows from equalities (8.1.32), (8.1.33). The equality (8.1.36) follows from the equality (8.1.37). $\square$

PROOF. The equality

(8.1.39) $$a^{1^*{}_*} = a^{0^*{}_*} {}^*_* a = E_n {}^*_* a = a$$

follows from equalities (8.1.34), (8.1.35). The equality (8.1.38) follows from the equality (8.1.39). $\square$

THEOREM 8.1.20.

(8.1.40) $$(a^T)^{n*^*} = (a^{n*^*})^T$$

THEOREM 8.1.21.

(8.1.45) $$(a^T)^{n^*{}_*} = (a^{n^*{}_*})^T$$

PROOF. We proceed by induction on $n$.

For $n = 0$ the statement immediately follows from equalities

| (8.1.32) | $a^{0*^*} = E_n$ |
|---|---|
| (8.1.34) | $a^{0^*{}_*} = E_n$ |
| (8.1.31) | $E_n^T = E_n$ |

Suppose the statement of theorem holds when $n = k - 1$

(8.1.41) $$(a^T)^{k-1*^*} = (a^{k-1*^*})^T$$

It follows from

| (8.1.33) | $a^{k*^*} = a^{k-1*^*} {}_*^* a$ |
|---|---|

PROOF. We proceed by induction on $n$.

For $n = 0$ the statement immediately follows from equalities

| (8.1.32) | $a^{0*^*} = E_n$ |
|---|---|
| (8.1.34) | $a^{0^*{}_*} = E_n$ |
| (8.1.31) | $E_n^T = E_n$ |

Suppose the statement of theorem holds when $n = k - 1$

(8.1.46) $$(a^T)^{k-1^*{}_*} = (a^{k-1^*{}_*})^T$$

It follows from

| (8.1.35) | $a^{k^*{}_*} = a^{k-1^*{}_*} {}^*_* a$ |
|---|---|

that

$$(8.1.42) \qquad (a^T)^{k_*{}^*} = (a^T)^{k-1_*{}^*}{}_* a^T$$

It follows from (8.1.42) and (8.1.41) that

$$(8.1.43) \qquad (a^T)^{k_*{}^*} = (a^{k-1^*{}^*})^T{}_* a^T$$

It follows from (8.1.43) and (8.1.25) that

$$(8.1.44) \qquad (a^T)^{k_*{}^*} = (a^{k-1^*{}^*}{}_* a)^T$$

(8.1.40) follows from (8.1.42) and (8.1.35). $\square$

that

$$(8.1.47) \qquad (a^T)^{k^*{}^*} = (a^T)^{k-1^*{}^*}{}^* a^T$$

It follows from (8.1.47) and (8.1.46) that

$$(8.1.48) \qquad (a^T)^{k^*{}^*} = (a^{k-1_*{}^*})^T{}^* a^T$$

It follows from (8.1.48) and (8.1.25) that

$$(8.1.49) \qquad (a^T)^{k^*{}^*} = (a^{k-1^*{}_*}{}^* a)^T$$

(8.1.45) follows from (8.1.47) and (8.1.33). $\square$

---

**Definition 8.1.22.** *Let a be $n \times n$ matrix. If there exists $n \times n$ matrix b such that*

$$(8.1.50) \qquad a_*{}^* b = E_n$$

*then the matrix $b = a^{-1_*{}^*}$ is called ${}_*{}^*$-* **inverse element** *of the matrix a*

$$(8.1.51) \qquad a_*{}^* a^{-1_*{}^*} = E_n$$

$\square$

**Definition 8.1.24.** *Let a be $n \times n$ matrix. If there exists $n \times n$ matrix b such that*

$$(8.1.52) \qquad a^*{}_* b = E_n$$

*then the matrix $b = a^{-1^*{}_*}$ is called ${}^*{}_*$-* **inverse element** *of the matrix a*

$$(8.1.53) \qquad a^*{}_* a^{-1^*{}_*} = E_n$$

$\square$

---

**Definition 8.1.23.** *If $n \times n$ matrix a has ${}_*{}^*$-inverse matrix we call such matrix ${}_*{}^*$-* **nonsingular matrix**. *Otherwise, we call such matrix ${}_*{}^*$-* **singular matrix**. $\square$

**Definition 8.1.25.** *If $n \times n$ matrix a has ${}^*{}_*$-inverse matrix we call such matrix ${}^*{}_*$-* **nonsingular matrix**. *Otherwise, we call such matrix ${}^*{}_*$-* **singular matrix**. $\square$

---

**Theorem 8.1.26.** *Let $n \times n$ matrix a have ${}_*{}^*$-inverse matrix. Then transpose matrix $a^T$ has ${}^*{}_*$-inverse matrix and these matrices satisfy the equality*

$$(8.1.54) \qquad (a^T)^{-1^*{}_*} = (a^{-1_*{}^*})^T$$

**Theorem 8.1.27.** *Let $n \times n$ matrix a have ${}^*{}_*$-inverse matrix. Then transpose matrix $a^T$ has ${}_*{}^*$-inverse matrix and these matrices satisfy the equality*

$$(8.1.56) \qquad (a^T)^{-1_*{}^*} = (a^{-1^*{}_*})^T$$

**Proof.** If we get transpose of both side (8.1.51) and apply

| (8.1.31)   $E_n^T = E_n$ |
|---|

we get

$$(a_*{}^* a^{-1_*{}^*})^T = E_n^T = E_n$$

Applying

| (8.1.24)   $(a_*{}^* b)^T = a^T{}^*{}_* b^T$ |
|---|

we get

$$(8.1.55) \qquad E_n = a^T{}^*{}_* (a^{-1_*{}^*})^T$$

(8.1.54) follows from comparison (8.1.53) and (8.1.55). $\square$

**Proof.** If we get transpose of both side (8.1.53) and apply

| (8.1.31)   $E_n^T = E_n$ |
|---|

we get

$$(a^*{}_* a^{-1^*{}_*})^T = E_n^T = E_n$$

Applying

| (8.1.25)   $(a^*{}_* b)^T = (a^T){}_*{}^* (b^T)$ |
|---|

we get

$$(8.1.57) \qquad E_n = a^T{}_*{}^* (a^{-1^*{}_*})^T$$

(8.1.56) follows from comparison (8.1.53) and (8.1.57). $\square$

Theorems 8.1.11, 8.1.12, 8.1.20, 8.1.21, 8.1.26, 8.1.27 show that some kind

of duality exists between $_*^*$-product and $^*_*$-product. We can combine these statements.

---

THEOREM 8.1.28. **Duality principle for biring**. *Let $\mathfrak{A}$ be true statement about biring $\mathcal{A}$. If we exchange the same time*

- $a \in \mathcal{A}$ *and* $a^T$
- $_*^*$-*product and* $^*_*$-*product*

*then we soon get true statement.*

---

THEOREM 8.1.29. **Duality principle for biring of matrices**. *Let $\mathcal{A}$ be biring of matrices. Let $\mathfrak{A}$ be true statement about matrices. If we exchange the same time*

- *columns and rows of all matrices*
- $_*^*$-*product and* $^*_*$-*product*

*then we soon get true statement.*

PROOF. This is the immediate consequence of the theorem 8.1.28. □

REMARK 8.1.30. *We execute operations in expression*

$$(8.1.58) \qquad a_*^*b_*^*c$$

*from left to right. However we can execute product from right to left. In such case the expression* (8.1.58) *gets the form*

$$c^*_*b^*_*a$$

*We follow the rule that to write power and indices from right of expression. For instance, let original expression be like*

$$a^{-1^*}{}_*^*b^i$$

*Then expression which we read from right to left is like*

$$b^{i^*}{}_*a^{-1^*}{}_*$$

*Suppose we established the order in which we write indices. Then we state that we read an expression from top to bottom reading first upper indices, then lower ones. We can read this expression from bottom to top. We extend this rule stating that we read symbols of operation in the same order as indices. For instance, if we read expression*

$$a^i{}_*^*b^{-1^*} = c^i$$

*from bottom to top, then we can write this expression as*

$$a_i{}^*_*b^{-1^*} = c_i$$

*According to the duality principle if we can prove one statement then we can prove other as well.* □

---

| THEOREM 8.1.31. *Let matrix $a$ have $_*^*$-inverse matrix. Then for any matrices $b$ and $c$ the equality* | THEOREM 8.1.32. *Let matrix $a$ have $^*_*$-inverse matrix. Then for any matrices $b$ and $c$ the equality* |
|---|---|
| $(8.1.59) \qquad b = c$ | $(8.1.61) \qquad b = c$ |

*follows from the equality*

$$(8.1.60) \qquad b_{*}{}^{*}a = c_{*}{}^{*}a$$

Proof. The equality (8.1.59) follows from the equality (8.1.60) if we multiply both parts of the equality (8.1.60) over $a^{-1_{*}{}^{*}}$. $\qquad \square$

Theorem 8.1.33. *Let $n \times n$ matrices $a$, $b$ have $_{*}^{*}$-inverse matrix. Then the matrix $a_{*}{}^{*}b$ also has $_{*}^{*}$-inverse matrix and*

$$(8.1.63) \quad (a_{*}{}^{*}b)^{-1_{*}{}^{*}} = b^{-1_{*}{}^{*}}{}_{*}{}^{*}a^{-1_{*}{}^{*}}$$

Proof. The equality

$$(8.1.64) \quad \begin{aligned} &(a_{*}{}^{*}b)_{*}{}^{*}(b^{-1_{*}{}^{*}}{}_{*}{}^{*}a^{-1_{*}{}^{*}}) \\ &= a_{*}{}^{*}(b_{*}{}^{*}b^{-1_{*}{}^{*}})_{*}{}^{*}a^{-1_{*}{}^{*}} \\ &= a_{*}{}^{*}E_{n*}{}^{*}a^{-1_{*}{}^{*}} = a_{*}{}^{*}a^{-1_{*}{}^{*}} \\ &= E_{n} \end{aligned}$$

follows from the equality

$$\boxed{(8.1.51) \quad a_{*}{}^{*}a^{-1_{*}{}^{*}} = E_{n}}$$

and from the definition 8.1.14. According to the definition 8.1.22 and the equality (8.1.64), the matrix $a_{*}{}^{*}b$ has $_{*}^{*}$-inverse matrix

$$(8.1.65) \qquad (a_{*}{}^{*}b)_{*}{}^{*}(a_{*}{}^{*}b)^{-1_{*}{}^{*}} = E_{n}$$

The equality (8.1.63) follows from equalities (8.1.64), (8.1.65). $\qquad \square$

*follows from the equality*

$$(8.1.62) \qquad b^{*}{}_{*}a = c^{*}{}_{*}a$$

Proof. The equality (8.1.61) follows from the equality (8.1.62) if we multiply both parts of the equality (8.1.62) over $a^{-1^{*}{}_{*}}$. $\qquad \square$

Theorem 8.1.34. *Let $n \times n$ matrices $a$, $b$ have $^{*}_{*}$-inverse matrix. Then the matrix $a^{*}{}_{*}b$ also has $^{*}_{*}$-inverse matrix and*

$$(8.1.66) \quad (a^{*}{}_{*}b)^{-1^{*}{}_{*}} = b^{-1^{*}{}_{*}}{}^{*}{}_{*}a^{-1^{*}{}_{*}}$$

Proof. The equality

$$(8.1.67) \quad \begin{aligned} &(a^{*}{}_{*}b)^{*}{}_{*}(b^{-1^{*}{}_{*}}{}^{*}{}_{*}a^{-1^{*}{}_{*}}) \\ &= a^{*}{}_{*}(b^{*}{}_{*}b^{-1^{*}{}_{*}})^{*}{}_{*}a^{-1^{*}{}_{*}} \\ &= a^{*}{}_{*}E_{n}{}^{*}{}_{*}a^{-1^{*}{}_{*}} = a^{*}{}_{*}a^{-1^{*}{}_{*}} \\ &= E_{n} \end{aligned}$$

follows from the equality

$$\boxed{(8.1.53) \quad a^{*}{}_{*}a^{-1^{*}{}_{*}} = E_{n}}$$

and from the definition 8.1.14. According to the definition 8.1.24 and the equality (8.1.67), the matrix $a^{*}{}_{*}b$ has $^{*}_{*}$-inverse matrix

$$(8.1.68) \qquad (a^{*}{}_{*}b)^{*}{}_{*}(a^{*}{}_{*}b)^{-1^{*}{}_{*}} = E_{n}$$

The equality (8.1.66) follows from equalities (8.1.67), (8.1.68). $\qquad \square$

## 8.2. $_{*}^{*}$--Quasideterminant

Theorem 8.2.1. *Let $n \times n$ matrix $a$ have $_{*}^{*}$-inverse matrix.[8.2] Then $k \times k$ submatrix of $_{*}^{*}$-inverse matrix satisfies to the equality*

$$(8.2.1) \qquad \left(\left(a^{-1_{*}{}^{*}}\right)_{J}^{I}\right)^{-1_{*}{}^{*}} = a_{I}^{J} - a_{[I]}^{J}{}_{*}{}^{*}\left(a_{[I]}^{[J]}\right)^{-1_{*}{}^{*}}{}_{*}{}^{*}a_{I}^{[J]}$$

$$|I| = |J| = k < n.$$

Proof. We can represent the definition

$$\boxed{(8.1.51) \quad a_{*}{}^{*}a^{-1_{*}{}^{*}} = E_{n}}$$

---

[8.2] This statement and its proof is based on statement 1.2.1 on page [5]-8 for matrix over free division ring.

of ₊*-inverse matrix the following way

$$(8.2.2) \qquad \begin{pmatrix} a_I^J & a_{[I]}^J \\ a_I^{[J]} & a_{[I]}^{[J]} \end{pmatrix} {}_{*}{}^{*} \begin{pmatrix} (a^{-1*})_J^I & (a^{-1*})_{[J]}^I \\ (a^{-1*})_J^{[I]} & (a^{-1*})_{[J]}^{[I]} \end{pmatrix} = \begin{pmatrix} E_k & 0 \\ 0 & E_{n-k} \end{pmatrix}$$

The equality $(a^{[J]}{}_{*}{}^{*}(a^{-1*}))_J = 0$

$$(8.2.3) \qquad a_I^{[J]}{}_{*}{}^{*}\left(a^{-1*}\right)_J^I + a_{[I]}^{[J]}{}_{*}{}^{*}(a^{-1*})_J^{[I]} = 0$$

and the equality $(a^J{}_{*}{}^{*}(a^{-1*}))_J = E_k$

$$(8.2.4) \qquad a_I^J{}_{*}{}^{*}(a^{-1*})_J^I + a_{[I]}^J{}_{*}{}^{*}(a^{-1*})_J^{[I]} = E_k$$

follow from the equality (8.2.2) We multiply the equality (8.2.3) by $\left(a_{[I]}^{[J]}\right)^{-1*}$

$$(8.2.5) \qquad \left(a_{[I]}^{[J]}\right)^{-1*}{}_{*}{}^{*}a_I^{[J]}{}_{*}{}^{*}\left(a^{-1*}\right)_J^I + (a^{-1*})_J^{[I]} = 0$$

The equality

$$(8.2.6) \qquad a_I^J{}_{*}{}^{*}(a^{-1*})_J^I - a_{[I]}^J{}_{*}{}^{*}\left(a_{[I]}^{[J]}\right)^{-1*}{}_{*}{}^{*}a_I^{[J]}{}_{*}{}^{*}(a^{-1*})_J^I = E_k$$

follows from equalities (8.2.5), (8.2.4). The equality (8.2.1) follows from (8.2.6) when we both sides of the equality (8.2.6) multiply by $\left((a^{-1*})_J^I\right)^{-1*}$.    □

**Theorem 8.2.2.** *Let $n \times n$ matrix $a$ have a $_{*}{}^{*}$-inverse matrix. Then entries of the $_{*}{}^{*}$-inverse matrix satisfy to the equality*

$$(8.2.7) \qquad (a^{-1*})_j^i = \left(a_i^j - a_{[i]}^j{}_{*}{}^{*}\left(a_{[i]}^{[j]}\right)^{-1*}{}_{*}{}^{*}a_i^{[j]}\right)^{-1}$$

$$(8.2.8) \qquad (\mathcal{H}a^{-1*})_i^j = a_i^j - a_{[i]}^j{}_{*}{}^{*}\left(a_{[i]}^{[j]}\right)^{-1*}{}_{*}{}^{*}a_i^{[j]}$$

Proof. The equality (8.2.7) follows from the equality (8.2.1) when we assume $I = i$, $J = j$. The equality (8.2.8) follows from equalities (8.1.30), (8.2.7).    □

**Definition 8.2.3.** $\binom{j}{i}$-$_{*}{}^{*}$-**quasideterminant** *of $n \times n$ matrix $a$ is formal expression*[8.3]

$$(8.2.9) \qquad \det({}_{*}{}^{*})_i^j a = a_i^j - a_{[i]}^j{}_{*}{}^{*}\left(a_{[i]}^{[j]}\right)^{-1*}{}_{*}{}^{*}a_i^{[j]} = ((a^{-1*})_j^i)^{-1}$$

*We consider $\binom{j}{i}$-$_{*}{}^{*}$-quasideterminant as an entry of the matrix*

$$(8.2.10) \qquad \begin{aligned} \det({}_{*}{}^{*})\, a &= \begin{pmatrix} \det({}_{*}{}^{*})_1^1 a & ... & \det({}_{*}{}^{*})_n^1 a \\ ... & ... & ... \\ \det({}_{*}{}^{*})_1^n a & ... & \det({}_{*}{}^{*})_n^n a \end{pmatrix} \\ &= \begin{pmatrix} ((a^{-1*})_1^1)^{-1} & ... & ((a^{-1*})_1^n)^{-1} \\ ... & ... & ... \\ ((a^{-1*})_n^1)^{-1} & ... & ((a^{-1*})_n^n)^{-1} \end{pmatrix} \end{aligned}$$

*which is called $_*{}^*$-**quasideterminant**.* □

**Theorem 8.2.4.** *Expression for $_*{}^*$-inverse matrix has form*

$$(8.2.11) \qquad a^{-1_*{}^*} = \mathcal{H}\det(_*{}^*)a$$

$$\det(_*{}^*)\,a = \mathcal{H}a^{-1_*{}^*}$$

**Proof.** The equality (8.2.11) follows from the equality

$$(8.2.12) \qquad \det(_*{}^*)_i^j\,a = (\mathcal{H}a^{-1_*{}^*})_i^j = ((a^{-1_*{}^*})_j^i)^{-1}$$

□

**Theorem 8.2.5.** *Expression for $(_i^j)$-$_*{}^*$-quasideterminant has the following form*

$$(8.2.13) \qquad \det(_*{}^*)_i^j\,a = a_i^j - a_{[i]}^j{}^* \mathcal{H}\det(_*{}^*)_{[i]}^{[j]}\,a_*{}^* a_i^{[j]}$$

**Proof.** The equality (8.2.13) follows from the equalities (8.2.8), (8.2.9). □

**Theorem 8.2.6.**

$$(8.2.14) \qquad \begin{pmatrix} \det(_*{}^*)_1^1 a & ... & \det(_*{}^*)_n^1 a \\ ... & ... & ... \\ \det(_*{}^*)_1^n a & ... & \det(_*{}^*)_n^n a \end{pmatrix} = \begin{pmatrix} ((a^{-1_*{}^*})_1^1)^{-1} & ... & ((a^{-1_*{}^*})_1^n)^{-1} \\ ... & ... & ... \\ ((a^{-1_*{}^*})_n^1)^{-1} & ... & ((a^{-1_*{}^*})_n^n)^{-1} \end{pmatrix}$$

$$(8.2.15) \qquad \begin{pmatrix} (\det(_*{}^*)_1^1 a)^{-1} & ... & (\det(_*{}^*)_1^n a)^{-1} \\ ... & ... & ... \\ (\det(_*{}^*)_n^1 a)^{-1} & ... & (\det(_*{}^*)_n^n a)^{-1} \end{pmatrix} = \begin{pmatrix} (a^{-1_*{}^*})_1^1 & ... & (a^{-1_*{}^*})_n^1 \\ ... & ... & ... \\ (a^{-1_*{}^*})_1^n & ... & (a^{-1_*{}^*})_n^n \end{pmatrix}$$

**Proof.** The equality (8.2.14) follows from the equality (8.2.9). The equality (8.2.15) follows from the equality (8.2.14). □

**Theorem 8.2.7.** *Consider matrix*

$$a = \begin{pmatrix} a_1^1 & a_2^1 \\ a_1^2 & a_2^2 \end{pmatrix}$$

*Then*

$$(8.2.16) \qquad \det(_*{}^*)a = \begin{pmatrix} a_1^1 - a_2^1(a_2^2)^{-1}a_1^2 & a_2^1 - a_1^1(a_1^2)^{-1}a_2^2 \\ a_1^2 - a_2^2(a_2^1)^{-1}a_1^1 & a_2^2 - a_1^2(a_1^1)^{-1}a_2^1 \end{pmatrix}$$

$$(8.2.17) \qquad a^{-1_*{}^*} = \begin{pmatrix} (a_1^1 - a_2^1(a_2^2)^{-1}a_1^2)^{-1} & (a_2^1 - a_2^2(a_2^1)^{-1}a_1^1)^{-1} \\ (a_1^2 - a_1^1(a_1^2)^{-1}a_2^2)^{-1} & (a_2^2 - a_1^2(a_1^1)^{-1}a_2^1)^{-1} \end{pmatrix}$$

---

<sup>8.3</sup> In the definition 8.2.3, I follow the definition [5]-1.2.2 on page=9.

Proof. According to the equality (8.2.9)

$$(8.2.18) \qquad \det({}_*{}^*)^1_1 a = a^1_1 - a^1_2 (a^2_2)^{-1} a^2_1$$

$$(8.2.19) \qquad \det({}_*{}^*)^2_1 a = a^1_2 - a^1_1 (a^2_1)^{-1} a^2_2$$

$$(8.2.20) \qquad \det({}_*{}^*)^1_2 a = a^2_1 - a^2_2 (a^1_2)^{-1} a^1_1$$

$$(8.2.21) \qquad \det({}_*{}^*)^2_2 a = a^2_2 - a^2_1 (a^1_1)^{-1} a^1_2$$

The equality (8.2.16) follows from equalities (8.2.18), (8.2.19), (8.2.20), (8.2.21). □

Theorem 8.2.8.

$$(8.2.22) \qquad (ma)^{-1}{}_*{}^* = a^{-1}{}_*{}^* m^{-1}$$

$$(8.2.23) \qquad (am)^{-1}{}_*{}^* = m^{-1} a^{-1}{}_*{}^*$$

Proof. To prove the equality (8.2.22) we proceed by induction on size of the matrix.

Since

$$(ma)^{-1}{}_*{}^* = \left((ma)^{-1}\right) = \left(a^{-1} m^{-1}\right) = \left(a^{-1}\right) m^{-1} = a^{-1}{}_*{}^* m^{-1}$$

the statement is evident for $1 \times 1$ matrix.

Let the statement holds for $(n-1) \times (n-1)$ matrix. Then from the equality (8.2.1), it follows that

$$(((ma)^{-1}{}_*{}^*)^I_J)^{-1}{}_*{}^* = (ma)^J_I - (ma)^J_{[I]}{}_*{}^* \left((ma)^{[J]}_{[I]}\right)^{-1}{}_*{}^* {}_*{}^* (ma)^{[J]}_I$$

$$= ma^J_I - m\ a^J_{[I]}{}_*{}^* \left(a^{[J]}_{[I]}\right)^{-1}{}_*{}^* m^{-1}{}_*{}^* m\ a^{[J]}_I$$

$$= m\ a^J_I - m\ a^J_{[I]}{}_*{}^* \left(a^{[J]}_{[I]}\right)^{-1}{}_*{}^* {}_*{}^* a^{[J]}_I$$

$$(8.2.24) \qquad (((ma)^{-1}{}_*{}^*)^I_J)^{-1}{}_*{}^* = m\ (a^{-1}{}_*{}^*)^I_J$$

The equality (8.2.22) follows from the equality (8.2.24). In the same manner we prove the equality (8.2.23). □

Theorem 8.2.9. *Let*

$$(8.2.25) \qquad a = \begin{pmatrix} 1 & 0 \\ 0 & 1 \end{pmatrix}$$

*Then*

$$(8.2.26) \qquad a^{-1}{}_*{}^* = \begin{pmatrix} 1 & 0 \\ 0 & 1 \end{pmatrix}$$

Proof. It is clear from (8.2.18), and (8.2.21) that $\left(a^{-1}{}_*{}^*\right)^1_1 = 1$ and $\left(a^{-1}{}_*{}^*\right)^2_2 = 1$. However expression for $\left(a^{-1}{}_*{}^*\right)^2_1$ and $\left(a^{-1}{}_*{}^*\right)^1_2$ cannot be defined from (8.2.20) and (8.2.19) since $a_1^2 = a_2^1 = 0$. We can transform these expressions. For instance

$$\begin{aligned}
(a^{-1*^*})_1^2 &= (a_2^I - a_I^I (a_I^2)^{-1} a_2^2)^{-1} \\
&= (a_I^I ((a_I^I)^{-1} a_2^I - (a_I^I)^{-1} a_2^2))^{-1} \\
&= ((a_I^2)^{-1} a_I^I (a_I^2 (a_I^I)^{-1} a_2^I - a_2^2))^{-1} \\
&= (a_I^I (a_I^2 (a_I^I)^{-1} a_2^I - a_2^2))^{-1} a_I^2
\end{aligned}$$

It follows immediately that $(a^{-1*^*})_1^2 = 0$. In the same manner we can find that $(a^{-1*^*})_2^I = 0$. This completes the proof of (8.2.26).           $\square$

## 8.3.   $*_*$--Quasideterminant

THEOREM 8.3.1. Let $n \times n$ matrix $a$ have $*_*$-inverse matrix.[8.4] Then $k \times k$ submatrix of $*_*$-inverse matrix satisfies to the equality

(8.3.1)          $$\left( (a^{-1*^*})_I^J \right)^{-1*^*} = a_J^I - a_J^{[I]} *_* \left( a_{[J]}^{[I]} \right)^{-1*^*} *_* a_{[J]}^I$$

$|I| = |J| = k < n$.

PROOF. We can represent the definition

$$\boxed{(8.1.53) \quad a^* *_* a^{-1*^*} = E_n}$$

of $*_*$-inverse matrix the following way

(8.3.2)          $$\begin{pmatrix} a_J^I & a_{[J]}^I \\ a_J^{[I]} & a_{[J]}^{[I]} \end{pmatrix} *_* \begin{pmatrix} (a^{-1*^*})_I^J & (a^{-1*^*})_{[I]}^J \\ (a^{-1*^*})_I^{[J]} & (a^{-1*^*})_{[I]}^{[J]} \end{pmatrix} = \begin{pmatrix} E_k & 0 \\ 0 & E_{n-k} \end{pmatrix}$$

The equality $(a_{[J]} *_* (a^{-1*^*})^J = 0)$

(8.3.3)          $$a_{[J]}^I *_* (a^{-1*^*})_I^J + a_{[J]}^{[I]} *_* (a^{-1*^*})_{[I]}^J = 0$$

and the equality $(a_{.J} *_* (a^{-1*^*})^J = E_k)$

(8.3.4)          $$a_J^I *_* (a^{-1*^*})_I^J + a_J^{[I]} *_* (a^{-1*^*})_{[I]}^J = E_k$$

follow from the equality (8.3.2) We multiply the equality (8.3.3) by $\left( a_{[J]}^{[I]} \right)^{-1*^*}$

(8.3.5)          $$\left( a_{[J]}^{[I]} \right)^{-1*^*} *_* a_{[J]}^I *_* (a^{-1*^*})_I^J + (a^{-1*^*})_{[I]}^J = 0$$

The equality

(8.3.6)          $$a_J^I *_* (a^{-1*^*})_I^J - a_J^{[I]} *_* \left( a_{[J]}^{[I]} \right)^{-1*^*} *_* a_{[J]}^I *_* (a^{-1*^*})_I^J = E_k$$

follows from equalities (8.3.5), (8.3.4). The equality (8.3.1) follows from (8.3.6) when we both sides of the equality (8.3.6) multiply by $\left( (a^{-1*^*})_I^J \right)^{-1*^*}$.           $\square$

---

[8.4]   This statement and its proof is based on statement 1.2.1 on page [5]-8 for matrix over free division ring.

THEOREM 8.3.2. *Let $n \times n$ matrix $a$ have a $_*^*$-inverse matrix. Then entries of the $_*^*$-inverse matrix satisfy to the equality*

$$(8.3.7) \qquad (a^{-1^*_*})^j_i = \left( a^i_j - a^{[i]}_j {}_*{}^* \left( a^{[i]}_{[j]} \right)^{-1^*_*} {}^*_* a^i_{[j]} \right)^{-1}$$

$$(8.3.8) \qquad (\mathcal{H}a^{-1^*_*})^i_j = a^i_j - a^{[i]}_j {}_*{}^* \left( a^{[i]}_{[j]} \right)^{-1^*_*} {}^*_* a^i_{[j]}$$

PROOF. The equality (8.3.7) follows from the equality (8.3.1) when we assume $I = i$, $J = j$. The equality (8.3.8) follows from equalities (8.1.30), (8.3.7). $\square$

DEFINITION 8.3.3. $\binom{i}{j}$-$^*_*$-**quasideterminant** *of $n \times n$ matrix $a$ is formal expression* [8.5]

$$(8.3.9) \qquad \det({}^*_*)^i_j\, a = a^i_j - a^{[i]}_j {}_*{}^* \left( a^{[i]}_{[j]} \right)^{-1^*_*} {}^*_* a^i_{[j]} = ((a^{-1^*_*})^i_j)^{-1}$$

*We consider $\binom{i}{j}$-$^*_*$-quasideterminant as an entry of the matrix*

$$(8.3.10) \qquad
\begin{aligned}
\det({}^*_*)\, a &= \begin{pmatrix} \det({}^*_*)^1_1 a & ... & \det({}^*_*)^1_n a \\ ... & ... & ... \\ \det({}^*_*)^n_1 a & ... & \det({}^*_*)^n_n a \end{pmatrix} \\[2mm]
&= \begin{pmatrix} ((a^{-1^*_*})^1_1)^{-1} & ... & ((a^{-1^*_*})^n_1)^{-1} \\ ... & ... & ... \\ ((a^{-1^*_*})^1_n)^{-1} & ... & ((a^{-1^*_*})^n_n)^{-1} \end{pmatrix}
\end{aligned}$$

*which is called $^*_*$-**quasideterminant**.* $\square$

THEOREM 8.3.4. *Expression for $^*_*$-inverse matrix has form*

$$(8.3.11) \qquad a^{-1^*_*} = \mathcal{H} \det({}^*_*)a$$
$$\det({}^*_*)\, a = \mathcal{H}a^{-1^*_*}$$

PROOF. The equality (8.3.11) follows from the equality

$$(8.3.12) \qquad \det({}^*_*)^i_j\, a = (\mathcal{H}a^{-1^*_*})^i_j = ((a^{-1^*_*})^j_i)^{-1}$$

$\square$

THEOREM 8.3.5. *Expression for $\binom{i}{j}$-$^*_*$-quasideterminant has the following form*

$$(8.3.13) \qquad \det({}^*_*)^i_j\, a = a^i_j - a^{[i]}_j {}_*{}^* \mathcal{H}\det({}^*_*)^{[i]}_{[j]}\, a^* {}_* a^i_{[j]}$$

PROOF. The equality (8.3.13) follows from the equalities (8.3.8), (8.3.9). $\square$

---

[8.5] In the definition 8.3.3, I follow the definition [5]-1.2.2 on page=9.

Theorem 8.3.6.

(8.3.14)
$$\begin{pmatrix} \det({}^*{}_*)^1_1 a & ... & \det({}^*{}_*)^1_n a \\ ... & ... & ... \\ \det({}^*{}_*)^n_1 a & ... & \det({}^*{}_*)^n_n a \end{pmatrix} = \begin{pmatrix} ((a^{-1^*{}_*})^1_1)^{-1} & ... & ((a^{-1^*{}_*})^n_1)^{-1} \\ ... & ... & ... \\ ((a^{-1^*{}_*})^1_n)^{-1} & ... & ((a^{-1^*{}_*})^n_n)^{-1} \end{pmatrix}$$

(8.3.15)
$$\begin{pmatrix} (\det({}^*{}_*)^1_1 a)^{-1} & ... & (\det({}^*{}_*)^n_1 a)^{-1} \\ ... & ... & ... \\ (\det({}^*{}_*)^1_n a)^{-1} & ... & (\det({}^*{}_*)^n_n a)^{-1} \end{pmatrix} = \begin{pmatrix} (a^{-1^*{}_*})^1_1 & ... & (a^{-1^*{}_*})^n_1 \\ ... & ... & ... \\ (a^{-1^*{}_*})^1_n & ... & (a^{-1^*{}_*})^n_n \end{pmatrix}$$

Proof. The equality (8.3.14) follows from the equality (8.3.9). The equality (8.3.15) follows from the equality (8.3.14). □

Theorem 8.3.7. *Consider matrix*

$$a = \begin{pmatrix} a^1_1 & a^1_2 \\ a^2_1 & a^2_2 \end{pmatrix}$$

*Then*

(8.3.16)
$$\det({}^*{}_*)a = \begin{pmatrix} a^1_1 - a^2_1(a^2_2)^{-1}a^1_2 & a^1_2 - a^2_2(a^2_1)^{-1}a^1_1 \\ a^2_1 - a^1_1(a^1_2)^{-1}a^2_2 & a^2_2 - a^1_2(a^1_1)^{-1}a^2_1 \end{pmatrix}$$

(8.3.17)
$$a^{-1^*{}_*} = \begin{pmatrix} (a^1_1 - a^2_1(a^2_2)^{-1}a^1_2)^{-1} & (a^2_1 - a^1_1(a^1_2)^{-1}a^2_2)^{-1} \\ (a^1_2 - a^2_2(a^2_1)^{-1}a^1_1)^{-1} & (a^2_2 - a^1_2(a^1_1)^{-1}a^2_1)^{-1} \end{pmatrix}$$

Proof. According to the equality (8.3.9)

(8.3.18) $$\det({}^*{}_*)^1_1 a = a^1_1 - a^2_1(a^2_2)^{-1}a^1_2$$

(8.3.19) $$\det({}^*{}_*)^2_1 a = a^2_1 - a^2_2(a^2_1)^{-1}a^1_1$$

(8.3.20) $$\det({}^*{}_*)^1_2 a = a^1_2 - a^1_1(a^1_2)^{-1}a^2_2$$

(8.3.21) $$\det({}^*{}_*)^2_2 a = a^2_2 - a^1_2(a^1_1)^{-1}a^2_1$$

The equality (8.3.16) follows from equalities (8.3.18), (8.3.19), (8.3.20), (8.3.21). □

Theorem 8.3.8.

(8.3.22) $$(ma)^{-1^*{}_*} = a^{-1^*{}_*}m^{-1}$$

(8.3.23) $$(am)^{-1^*{}_*} = m^{-1}a^{-1^*{}_*}$$

Proof. To prove the equality (8.3.22) we proceed by induction on size of the matrix.

Since
$$(ma)^{-1^{*_*}} = \left((ma)^{-1}\right) = \left(a^{-1}m^{-1}\right) = \left(a^{-1}\right)m^{-1} = a^{-1^{*_*}}m^{-1}$$
the statement is evident for $1 \times 1$ matrix.

Let the statement holds for $(n-1) \times (n-1)$ matrix. Then from the equality (8.3.1), it follows that

$$(((ma)^{-1^{*_*}})_I^J)^{-1^{*_*}} = (ma)_J^I - (ma)_J^{[I]}{}_* \left((ma)_{[J]}^{[I]}\right)^{-1^{*_*}}{}_*(ma)_{[J]}^I$$

$$= m a_J^I - m\, a_J^{[I]}{}_* \left(a_{[J]}^{[I]}\right)^{-1^{*_*}} m^{-1}{}_* m\, a_{[J]}^I$$

$$= m\, a_J^I - m\, a_J^{[I]}{}_* \left(a_{[J]}^{[I]}\right)^{-1^{*_*}}{}_*a_{[J]}^I$$

$$(8.3.24) \qquad (((ma)^{-1^{*_*}})_I^J)^{-1^{*_*}} = m\, (a^{-1^{*_*}})_I^J$$

The equality (8.3.22) follows from the equality (8.3.24). In the same manner we prove the equality (8.3.23). $\square$

**THEOREM 8.3.9.** *Let*

$$(8.3.25) \qquad a = \begin{pmatrix} 1 & 0 \\ 0 & 1 \end{pmatrix}$$

*Then*

$$(8.3.26) \qquad a^{-1^{*_*}} = \begin{pmatrix} 1 & 0 \\ 0 & 1 \end{pmatrix}$$

**PROOF.** It is clear from (8.3.18), and (8.3.21) that $(a^{-1^{*_*}})_1^1 = 1$ and $(a^{-1^{*_*}})_2^2 = 1$. However expression for $(a^{-1^{*_*}})_1^2$ and $(a^{-1^{*_*}})_2^1$ cannot be defined from (8.3.20) and (8.3.19) since $a_1^2 = a_2^1 = 0$. We can transform these expressions. For instance

$$(a^{-1^{*_*}})_1^2 = (a_2^1 - a_2^2(a_1^2)^{-1}a_1^1)^{-1}$$

$$= (a_1^1((a_1^1)^{-1}a_2^1 - (a_1^2)^{-1}a_2^2))^{-1}$$

$$= ((a_1^2)^{-1}a_1^1(a_2^2(a_1^1)^{-1}a_2^1 - a_2^2))^{-1}$$

$$= (a_1^1(a_2^2(a_1^1)^{-1}a_2^1 - a_2^2))^{-1}a_1^2$$

It follows immediately that $(a^{-1^{*_*}})_1^2 = 0$. In the same manner we can find that $(a^{-1^{*_*}})_2^1 = 0$. This completes the proof of (8.3.26). $\square$

**THEOREM 8.3.10.**

$$(8.3.27) \qquad \det(_*^*)_j^i\, a^T = \det(^*_*)_i^j\, a$$

**PROOF.** According to (8.2.9) and (8.1.29)

$$\det(_*^*)_j^i\, a^T = (((a^T)^{-1^{*_*}})_{\,i}^{\,-\,\,j}{}_{\,\cdot\,-})^{-1}$$

Using theorem 8.1.27 we get

$$\det({}_*{}^*)^i_j\, a^T = (((a^{-1}{}^*{}^*)^T)^{\cdot\,-\,j}_{i\,\cdot\,-})^{-1}$$

Using (8.1.13) we get

(8.3.28)
$$\det({}_*{}^*)^i_j\, a^T = ((a^{-1}{}^*{}^*)^{\,\cdot\,i\,\cdot\,-}_{-\,\cdot\,j})^{-1}$$

Using (8.3.28), (8.3.9) we get (8.3.27).                                   □

The theorem 8.3.10 extends the duality principle stated in the theorem 8.1.29 to statements on quasideterminants and tells us that the same expression is ${}_*{}^*$-quasideterminant of matrix $a$ and ${}^*{}_*$-quasideterminant of matrix $a^T$.

> THEOREM 8.3.11. **duality principle**. *Let $\mathfrak{A}$ be true statement about matrix biring. If we exchange the same time*
> - *row and column*
> - *${}_*{}^*$-quasideterminant and ${}^*{}_*$-quasideterminant*
>
> *then we soon get true statement.*

## 8.4. Vector Space Type

### 8.4.1. Left $A$-Vector Space of Columns.

> DEFINITION 8.4.1. *We represented the set of vectors $v(1) = v_1, ..., v(m) = v_m$ as row of matrix*
>
> (8.4.1)
> $$v = \begin{pmatrix} v_1 & ... & v_m \end{pmatrix}$$
>
> *and the set of $A$-nummbers $c(1) = c^1, ..., c(m) = c^m$ as column of matrix*
>
> (8.4.2)
> $$c = \begin{pmatrix} c^1 \\ ... \\ c^m \end{pmatrix}$$
>
> *Corresponding representation of left $A$-vector space $V$ is called* **left $A$-vector space of columns**, *and $V$-number is called* **column vector**.                □

> THEOREM 8.4.2. *We can represent linear combination*
> $$\sum_{i=1}^m c(i)v(i) = c^i v_i$$
> *of vectors $v_1, ..., v_m$ as ${}^*{}_*$-product of matrices*
>
> (8.4.3)
> $$c\,{}^*{}_*\,v = \begin{pmatrix} c^1 \\ ... \\ c^m \end{pmatrix} {}^*{}_* \begin{pmatrix} v_1 & ... & v_m \end{pmatrix} = c^i v_i$$

PROOF. The theorem follows from definitions 7.3.5, 8.4.1.                □

THEOREM 8.4.3. *If we write vectors of basis $\overline{\overline{e}}$ as row of matrix*

$$(8.4.4) \qquad e = \begin{pmatrix} e_1 & ... & e_n \end{pmatrix}$$

*and coordinates of vector $\overline{w} = w^i e_i$ with respect to basis $\overline{\overline{e}}$ as column of matrix*

$$(8.4.5) \qquad w = \begin{pmatrix} w^1 \\ ... \\ w^n \end{pmatrix}$$

*then we can represent the vector $\overline{w}$ as $^*{}_*$-product of matrices*

$$(8.4.6) \qquad \overline{w} = w^*{}_* e = \begin{pmatrix} w^1 \\ ... \\ w^n \end{pmatrix} {}^*{}_* \begin{pmatrix} e_1 & ... & e_n \end{pmatrix} = w^i e_i$$

PROOF. The theorem follows from the theorem 8.4.2.                        □

THEOREM 8.4.4. *Coordinates of vector $v \in V$ relative to basis $\overline{\overline{e}}$ of left $A$-vector space $V$ are uniquely defined. The equality*

$$(8.4.7) \qquad v^*{}_* e = w^*{}_* e \qquad v^i e_i = w^i e_i$$

*implies the equality*

$$v = w \qquad v^i = w^i$$

PROOF. The theorem follows from the theorem 7.3.14.                       □

THEOREM 8.4.5. *Let the set of vectors $\overline{\overline{e}}_A = (e_{Ak}, k \in K)$ be a basis of $D$-algebra of columns $A$. Let $\overline{\overline{e}}_V = (e_{Vi}, i \in I)$ be a basis of left $A$-vector space of columns $V$. Then the set of $V$-numbers*

$$(8.4.8) \qquad \overline{\overline{e}}_{VA} = (e_{VAki} = e_{Ak} e_{Vi})$$

*is the basis of $D$-vector space $V$. The basis $\overline{\overline{e}}_{VA}$ is called extension of the basis $(\overline{\overline{e}}_A, \overline{\overline{e}}_V)$.*

LEMMA 8.4.6. *The set of $V$-numbers $e_{2ik}$ generates $D$-vector space $V$.*

PROOF. Let

$$(8.4.9) \qquad a = a^i e_{1i}$$

be any $V$-number. For any $A$-number $a^i$, according to the definition 5.1.1 and the convention 5.2.1, there exist $D$-numbers $a^{ik}$ such that

$$(8.4.10) \qquad a^i = a^{ik} e_k$$

The equality

$$(8.4.11) \qquad a = a^{ik} e_k e_{1i} = a^{ik} e_{2ik}$$

follows from equalities $(8.4.8)$, $(8.4.9)$, $(8.4.10)$. According to the theorem $4.1.9$, the lemma follows from the equality $(8.4.11)$.                                    $\odot$

Lemma 8.4.7. *The set of $V$-numbers $e_{2ik}$ is linearly independent over the ring $D$.*

Proof. Let

$$(8.4.12) \qquad\qquad a^{ik} e_{2ik} = 0$$

The equality

$$(8.4.13) \qquad\qquad a^{ik} e_k e_{1i} = 0$$

follows from equalities $(8.4.8)$, $(8.4.12)$. According to the theorem $4.4.4$, the equality

$$(8.4.14) \qquad\qquad a^{ik} e_k = 0 \quad k \in K$$

follows from the equality $(8.4.13)$. According to the theorem $4.2.7$, the equality

$$(8.4.15) \qquad\qquad a^{ik} = 0 \quad k \in K \quad i \in I$$

follows from the equality $(8.4.14)$. Therefore, the set of $V$-numbers $e_{2ik}$ is linearly independent over the ring $D$.                                    $\odot$

Proof of theorem 8.4.5. The theorem follows from the theorem $4.2.7$ and from lemmas $8.4.6$, $8.4.7$.                                    $\square$

### 8.4.2. Left $A$-Vector Space of Rows.

Definition 8.4.8. *We represented the set of vectors $v(1) = v^1$, ..., $v(m) = v^m$ as column of matrix*

$$(8.4.16) \qquad\qquad v = \begin{pmatrix} v^1 \\ ... \\ v^m \end{pmatrix}$$

*and the set of $A$-nummbers $c(1) = c_1$, ..., $c(m) = c_m$ as row of matrix*

$$(8.4.17) \qquad\qquad c = \begin{pmatrix} c_1 & ... & c_m \end{pmatrix}$$

*Corresponding representation of left $A$-vector space $V$ is called* **left $A$-vector space of rows***, and $V$-number is called* **row vector***.*                                    $\square$

Theorem 8.4.9. *We can represent linear combination*

$$\sum_{i=1}^{m} c(i) v(i) = c_i v^i$$

*of vectors $v^1$, ..., $v^m$ as $_*{}^*$-product of matrices*

$$(8.4.18) \qquad\qquad c_*{}^* v = \begin{pmatrix} c_1 & ... & c_m \end{pmatrix} {}_*{}^* \begin{pmatrix} v^1 \\ ... \\ v^m \end{pmatrix} = c_i v^i$$

Proof. The theorem follows from definitions 7.3.5, 8.4.8.                    □

Theorem 8.4.10. *If we write vectors of basis $\overline{\overline{e}}$ as column of matrix*

$$(8.4.19) \qquad e = \begin{pmatrix} e^1 \\ ... \\ e^n \end{pmatrix}$$

*and coordinates of vector $\overline{w} = w_i e^i$ with respect to basis $\overline{\overline{e}}$ as row of matrix*

$$(8.4.20) \qquad w = \begin{pmatrix} w_1 & ... & w_n \end{pmatrix}$$

*then we can represent the vector $\overline{w}$ as $_*{}^*$-product of matrices*

$$(8.4.21) \qquad \overline{w} = w_*{}^* e = \begin{pmatrix} w_1 & ... & w_n \end{pmatrix} {}_*{}^* \begin{pmatrix} e^1 \\ ... \\ e^n \end{pmatrix} = w_i e^i$$

Proof. The theorem follows from the theorem 8.4.9.                    □

Theorem 8.4.11. *Coordinates of vector $v \in V$ relative to basis $\overline{\overline{e}}$ of left $A$-vector space $V$ are uniquely defined. The equality*

$$(8.4.22) \qquad v_*{}^* e = w_*{}^* e \qquad v_i e^i = w_i e^i$$

*implies the equality*

$$v = w \qquad v_i = w_i$$

Proof. The theorem follows from the theorem 7.3.14.                    □

Theorem 8.4.12. *Let the set of vectors $\overline{\overline{e}}_A = (e_A^k, k \in K)$ be a basis of $D$-algebra of rows $A$. Let $\overline{\overline{e}}_V = (e_V^i, i \in I)$ be a basis of left $A$-vector space of rows $V$. Then the set of $V$-numbers*

$$(8.4.23) \qquad \overline{\overline{e}}_{VA} = (e_{VA}^{ki} = e_A^k e_V^i)$$

*is the basis of $D$-vector space $V$. The basis $\overline{\overline{e}}_{VA}$ is called extension of the basis $(\overline{\overline{e}}_A, \overline{\overline{e}}_V)$.*

Lemma 8.4.13. *The set of $V$-numbers $e_2^{ik}$ generates $D$-vector space $V$.*

Proof. Let

$$(8.4.24) \qquad a = a_i e_1^i$$

be any $V$-number. For any $A$-number $a_i$, according to the definition 5.1.1 and the convention 5.2.1, there exist $D$-numbers $a_{ik}$ such that

$$(8.4.25) \qquad a_i = a_{ik} e^k$$

The equality

$$(8.4.26) \qquad a = a_{ik} e^k e_1^i = a_{ik} e_2^{ik}$$

follows from equalities ($8.4.23$), ($8.4.24$), ($8.4.25$). According to the theorem $4.1.9$, the lemma follows from the equality ($8.4.26$).                    ⊙

**Lemma 8.4.14.** *The set of $V$-numbers $e_2^{ik}$ is linearly independent over the ring $D$.*

**Proof.**  Let

$$(8.4.27) \qquad\qquad a_{ik}e_2^{ik} = 0$$

The equality

$$(8.4.28) \qquad\qquad a_{ik}e^k e_1^i = 0$$

follows from equalities ($8.4.23$), ($8.4.27$). According to the theorem $8.4.11$, the equality

$$(8.4.29) \qquad\qquad a_{ik}e^k = 0 \quad k \in K$$

follows from the equality ($8.4.28$). According to the theorem $4.2.14$, the equality

$$(8.4.30) \qquad\qquad a_{ik} = 0 \quad k \in K \quad i \in I$$

follows from the equality ($8.4.29$). Therefore, the set of $V$-numbers $e_2^{ik}$ is linearly independent over the ring $D$.                    ⊙

**Proof of theorem 8.4.12.**    The theorem follows from the theorem $4.2.14$ and from lemmas $8.4.13$, $8.4.14$.                    □

### 8.4.3. Right $A$-Vector Space of Columns.

**Definition 8.4.15.** *We represented the set of vectors  $v(1) = v_1$, ..., $v(m) = v_m$  as row of matrix*

$$(8.4.31) \qquad\qquad v = \begin{pmatrix} v_1 & ... & v_m \end{pmatrix}$$

*and the set of $A$-nummbers  $c(1) = c^1$, ..., $c(m) = c^m$   as column of matrix*

$$(8.4.32) \qquad\qquad c = \begin{pmatrix} c^1 \\ ... \\ c^m \end{pmatrix}$$

*Corresponding representation of right $A$-vector space $V$  is called* **right $A$-vector space of columns***, and $V$-number is called* **column vector***.*                    □

**Theorem 8.4.16.** *We can represent linear combination*

$$\sum_{i=1}^{m} v(i)c(i) = v_i c^i$$

*of vectors  $v_1$, ..., $v_m$   as $_*{}^*$-product of matrices*

$$(8.4.33) \qquad v_*{}^* c = \begin{pmatrix} v_1 & ... & v_m \end{pmatrix} {}_*{}^* \begin{pmatrix} c^1 \\ ... \\ c^m \end{pmatrix} = v_i c^i$$

Proof. The theorem follows from definitions 7.5.5, 8.4.15. □

Theorem 8.4.17. *If we write vectors of basis $\overline{\overline{e}}$ as row of matrix*

$$(8.4.34) \qquad e = \begin{pmatrix} e_1 & \dots & e_n \end{pmatrix}$$

*and coordinates of vector $\overline{w} = e_i w^i$ with respect to basis $\overline{\overline{e}}$ as column of matrix*

$$(8.4.35) \qquad w = \begin{pmatrix} w^1 \\ \dots \\ w^n \end{pmatrix}$$

*then we can represent the vector $\overline{w}$ as $_*{}^*$-product of matrices*

$$(8.4.36) \qquad \overline{w} = e_*{}^* w = \begin{pmatrix} e_1 & \dots & e_n \end{pmatrix} _*{}^* \begin{pmatrix} w^1 \\ \dots \\ w^n \end{pmatrix} = e_i w^i$$

Proof. The theorem follows from the theorem 8.4.16. □

Theorem 8.4.18. *Coordinates of vector $v \in V$ relative to basis $\overline{\overline{e}}$ of right $A$-vector space $V$ are uniquely defined. The equality*

$$(8.4.37) \qquad e_*{}^* v = e_*{}^* w \qquad e_i v^i = e_i w^i$$

*implies the equality*

$$v = w \qquad v^i = w^i$$

Proof. The theorem follows from the theorem 7.5.14. □

Theorem 8.4.19. *Let the set of vectors $\overline{\overline{e}}_A = (e_{Ak}, k \in K)$ be a basis of $D$-algebra of columns $A$. Let $\overline{\overline{e}}_V = (e_{Vi}, i \in I)$ be a basis of right $A$-vector space of columns $V$. Then the set of $V$-numbers*

$$(8.4.38) \qquad \overline{\overline{e}}_{VA} = (e_{VAki} = e_{Vi} e_{Ak})$$

*is the basis of $D$-vector space $V$. The basis $\overline{\overline{e}}_{VA}$ is called extension of the basis $(\overline{\overline{e}}_A, \overline{\overline{e}}_V)$.*

Lemma 8.4.20. *The set of $V$-numbers $e_{2ik}$ generates $D$-vector space $V$.*

Proof. Let

$$(8.4.39) \qquad a = e_{1i} a^i$$

be any $V$-number. For any $A$-number $a^i$, according to the definition 5.1.1 and the convention 5.2.1, there exist $D$-numbers $a^{ik}$ such that

$$(8.4.40) \qquad a^i = e_k a^{ik}$$

The equality

$$(8.4.41) \qquad a = e_{1i} e_k a^{ik} = e_{2ik} a^{ik}$$

follows from equalities ($8.4.38$), ($8.4.39$), ($8.4.40$). According to the theorem $4.1.9$, the lemma follows from the equality ($8.4.41$). ⊙

Lemma 8.4.21. *The set of $V$-numbers $e_{2ik}$ is linearly independent over the ring $D$.*

Proof.  Let

$$(8.4.42) \qquad\qquad e_{2ik}a^{ik} = 0$$

The equality

$$(8.4.43) \qquad\qquad e_{1i}e_k a^{ik} = 0$$

follows from equalities ($8.4.38$), ($8.4.42$). According to the theorem $8.4.18$, the equality

$$(8.4.44) \qquad\qquad e_k a^{ik} = 0 \quad k \in K$$

follows from the equality ($8.4.43$). According to the theorem $4.2.7$, the equality

$$(8.4.45) \qquad\qquad a^{ik} = 0 \quad k \in K \quad i \in I$$

follows from the equality ($8.4.44$). Therefore, the set of $V$-numbers $e_{2ik}$ is linearly independent over the ring $D$. ⊙

Proof of theorem 8.4.19.    The theorem follows from the theorem $4.2.7$ and from lemmas $8.4.20$, $8.4.21$. □

### 8.4.4. Right $A$-Vector Space of Rows.

Definition 8.4.22. *We represented the set of vectors $v(1) = v^1$, ..., $v(m) = v^m$ as column of matrix*

$$(8.4.46) \qquad\qquad v = \begin{pmatrix} v^1 \\ ... \\ v^m \end{pmatrix}$$

*and the set of $A$-nummbers $c(1) = c_1$, ..., $c(m) = c_m$ as row of matrix*

$$(8.4.47) \qquad\qquad c = \begin{pmatrix} c_1 & ... & c_m \end{pmatrix}$$

*Corresponding representation of right $A$-vector space $V$ is called* **right $A$-vector space of rows**, *and $V$-number is called* **row vector**. □

Theorem 8.4.23. *We can represent linear combination*

$$\sum_{i=1}^m v(i)c(i) = v^i c_i$$

*of vectors $v^1$, ..., $v^m$ as $*_*$-product of matrices*

$$(8.4.48) \qquad v *_* c = \begin{pmatrix} v^1 \\ ... \\ v^m \end{pmatrix} *_* \begin{pmatrix} c_1 & ... & c_m \end{pmatrix} = v^i c_i$$

PROOF. The theorem follows from definitions 7.5.5, 8.4.22.  □

THEOREM 8.4.24. *If we write vectors of basis $\overline{\overline{e}}$ as column of matrix*

$$(8.4.49) \qquad e = \begin{pmatrix} e^1 \\ ... \\ e^n \end{pmatrix}$$

*and coordinates of vector $\overline{w} = e^i w_i$ with respect to basis $\overline{\overline{e}}$ as row of matrix*

$$(8.4.50) \qquad w = \begin{pmatrix} w_1 & ... & w_n \end{pmatrix}$$

*then we can represent the vector $\overline{w}$ as $*_*$-product of matrices*

$$(8.4.51) \qquad \overline{w} = e *_* w = \begin{pmatrix} e^1 \\ ... \\ e^n \end{pmatrix} *_* \begin{pmatrix} w_1 & ... & w_n \end{pmatrix} = e^i w_i$$

PROOF. The theorem follows from the theorem 8.4.23.  □

THEOREM 8.4.25. *Coordinates of vector $v \in V$ relative to basis $\overline{\overline{e}}$ of right $A$-vector space $V$ are uniquely defined. The equality*

$$(8.4.52) \qquad e *_* v = e *_* w \qquad e^i v_i = e^i w_i$$

*implies the equality*

$$v = w \qquad v_i = w_i$$

PROOF. The theorem follows from the theorem 7.5.14.  □

THEOREM 8.4.26. *Let the set of vectors $\overline{\overline{e}}_A = (e_A^k, k \in K)$ be a basis of $D$-algebra of rows $A$. Let $\overline{\overline{e}}_V = (e_V^i, i \in I)$ be a basis of right $A$-vector space of rows $V$. Then the set of $V$-numbers*

$$(8.4.53) \qquad \overline{\overline{e}}_{VA} = (e_{VA}^{ki} = e_V^i e_A^k)$$

*is the basis of $D$-vector space $V$. The basis $\overline{\overline{e}}_{VA}$ is called extension of the basis $(\overline{\overline{e}}_A, \overline{\overline{e}}_V)$.*

LEMMA 8.4.27. *The set of $V$-numbers $e_2^{ik}$ generates $D$-vector space $V$.*

PROOF. Let

$$(8.4.54) \qquad a = e_1^i a_i$$

be any $V$-number. For any $A$-number $a_i$, according to the definition 5.1.1 and the convention 5.2.1, there exist $D$-numbers $a_{ik}$ such that

$$(8.4.55) \qquad a_i = e^k a_{ik}$$

The equality

$$(8.4.56) \qquad a = e_1^i e^k a_{ik} = e_2^{ik} a_{ik}$$

follows from equalities (8.4.53), (8.4.54), (8.4.55). According to the theorem 4.1.9, the lemma follows from the equality (8.4.56).                                                   ⊙

LEMMA 8.4.28. *The set of V-numbers $e_2^{ik}$ is linearly independent over the ring D.*

PROOF. Let

$$(8.4.57) \qquad\qquad e_2^{ik} a_{ik} = 0$$

The equality

$$(8.4.58) \qquad\qquad e_1^i e^k a_{ik} = 0$$

follows from equalities (8.4.53), (8.4.57). According to the theorem 8.4.25, the equality

$$(8.4.59) \qquad\qquad e^k a_{ik} = 0 \quad k \in K$$

follows from the equality (8.4.58). According to the theorem 4.2.14, the equality

$$(8.4.60) \qquad\qquad a_{ik} = 0 \quad k \in K \quad i \in I$$

follows from the equality (8.4.59). Therefore, the set of $V$-numbers $e_2^{ik}$ is linearly independent over the ring $D$.                                                         ⊙

PROOF OF THEOREM 8.4.26.    The theorem follows from the theorem 4.2.14 and from lemmas 8.4.27, 8.4.28.                                                            □

CHAPTER 9

# System of Linear Equations in $A$-Vector Space

### 9.1. Linear Span in Left $A$-Vector Space

#### 9.1.1. Left $A$-vector space of columns.

DEFINITION 9.1.1. *Let $V$ be a left $A$-vector space of columns and* $\{\overline{a}_i \in V, i \in I\}$ *be set of vectors.* **Linear span** *in left $A$-vector space of columns is set* $span(\overline{a}_i, i \in I)$ *of vectors linearly dependent on vectors $\overline{a}_i$.* □

THEOREM 9.1.2. *Let $V$ be a left $A$-vector space of columns and row*

$$\overline{a} = \begin{pmatrix} \overline{a}_1 & ... & \overline{a}_m \end{pmatrix} = (\overline{a}_i, i \in I)$$

*be set of vectors. The vector $\overline{b}$ linearly dependends on vectors $\overline{a}_i$*

$$\overline{b} \in span(\overline{a}_i, i \in I)$$

*if there exist column of A-numbers*

$$x = \begin{pmatrix} x^1 \\ ... \\ x^m \end{pmatrix}$$

*such that the following equality is true*

(9.1.1) $$\overline{b} = x^* {}_* \overline{a} = x^i \overline{a}_i$$

PROOF. The equality (9.1.1) follows from the definition 9.1.1 and the theorem 8.4.2. □

THEOREM 9.1.3. $span(\overline{a}_i, i \in I)$ *is subspace of left $A$-vector space of columns.*

PROOF. Suppose

$$\overline{b} \in \mathrm{span}(\overline{a}_i, i \in I)$$
$$\overline{c} \in \mathrm{span}(\overline{a}_i, i \in I)$$

According to the theorem 9.1.2

(9.1.2) $$\overline{b} = b^* {}_* \overline{a} = b^i \overline{a}_i$$

(9.1.3) $$\overline{c} = c^* {}_* \overline{a} = c^i \overline{a}_i$$

Then

$$\overline{b} + \overline{c} = a^* {}_* b + a^* {}_* c = a^* {}_* (b+c) \in \mathrm{span}(\overline{a}_i, i \in I)$$
$$\overline{b}k = (a^* {}_* b)k = a^* {}_* (bk) \in \mathrm{span}(\overline{a}_i, i \in I)$$

This proves the statement. □

THEOREM 9.1.4. *Let* $\overline{\overline{e}} = (e_j, j \in J)$ *be a basis. Let vectors* $\overline{b}$, $\overline{a}_i$ *have expansion*

(9.1.4)
$$\overline{b} = b^*{}_* e = b^j e_j \qquad b = \begin{pmatrix} b^1 \\ ... \\ b^n \end{pmatrix}$$

(9.1.5)
$$\overline{a_i} = a_i{}^*{}_* e = a_i^j e_j \qquad a_i = \begin{pmatrix} a_i^1 \\ ... \\ a_i^n \end{pmatrix}$$

*The vector* $\overline{b}$ *linearly dependends on vectors* $\overline{a}_i$

$$\overline{b} \in span(\overline{a}_i, i \in I)$$

*if there exist column of A-numbers*

(9.1.6)
$$x = \begin{pmatrix} x^1 \\ ... \\ x^m \end{pmatrix}$$

*which satisfy the* **system of linear equations**

(9.1.7)
$$x^1 a_1^1 + ... + x^m a_m^1 = b^1$$
$$......$$
$$x^1 a_1^n + ... + x^m a_m^n = b^n$$

PROOF. The equality

(9.1.8)
$$b^j e_j = x^i (a_i^j e_j) = (x^i a_i^j) e_j$$

follows from equalities (9.1.1), (9.1.4), (9.1.5). The equality

(9.1.9)
$$b^j = x^i a_i^j$$

follows from the equality (9.1.8) and the theorem 8.4.4. The equality (9.1.7) follows from the equality (9.1.9). $\square$

THEOREM 9.1.5. *Consider the matrix*

(9.1.10)
$$a = \left( \begin{pmatrix} a_1^1 \\ ... \\ a_1^n \end{pmatrix} \quad ... \quad \begin{pmatrix} a_m^1 \\ ... \\ a_m^n \end{pmatrix} \right) = \begin{pmatrix} a_1^1 & ... & a_m^1 \\ ... & ... & ... \\ a_1^n & ... & a_m^n \end{pmatrix}$$

*Then the system of linear equations (9.1.7) can be represented as*

(9.1.11)
$$\begin{pmatrix} x^1 \\ ... \\ x^m \end{pmatrix} {}^*{}_* \begin{pmatrix} a_1^1 & ... & a_m^1 \\ ... & ... & ... \\ a_1^n & ... & a_m^n \end{pmatrix} = \begin{pmatrix} b^1 \\ ... \\ b^n \end{pmatrix}$$

(9.1.12) $$x^*{}_*a = b$$

Proof. The equality (9.1.11) follows from equalities (9.1.4), (9.1.6), (9.1.7), (9.1.10) and from the definition 8.1.3. ☐

Definition 9.1.6. *Let $a$ be $^*{}_*$-nonsingular matrix. We call appropriate system of linear equations*

(9.1.12)   $x^*{}_*a = b$

**nonsingular system of linear equations.** ☐

Theorem 9.1.7. *Solution of nonsingular system of linear equations*

(9.1.12)   $x^*{}_*a = b$

*is determined uniquely and can be presented in either form* [9.1]

(9.1.13) $$x = b^*{}_*a^{-1^*{}_*}$$

(9.1.14) $$x = b^*{}_*\mathcal{H}\det({}^*{}_*)a$$

Proof. Multiplying both sides of the equation (9.1.12) from left by $a^{-1^*{}_*}$, we get (9.1.13). Using definition

(8.3.11)   $a^{-1^*{}_*} = \mathcal{H}\det({}^*{}_*)a$

we get (9.1.14). Based on the theorem 8.1.31, the solution is unique. ☐

**9.1.2. Left $A$-vector space of rows.**

Definition 9.1.8. *Let $V$ be a left $A$-vector space of rows and $\{\overline{a}^i \in V, i \in I\}$ be set of vectors.* **Linear span** *in left $A$-vector space of rows is set* $span(\overline{a}^i, i \in I)$ *of vectors linearly dependent on vectors $\overline{a}^i$.* ☐

Theorem 9.1.9. *Let $V$ be a left $A$-vector space of rows and column*

$$\overline{a} = \begin{pmatrix} \overline{a}^i \\ ... \\ \overline{a}^m \end{pmatrix} = (\overline{a}^i, i \in I)$$

*be set of vectors. The vector $\overline{b}$ linearly dependends on vectors $\overline{a}^i$*

$$\overline{b} \in span(\overline{a}^i, i \in I)$$

*if there exist row of $A$-numbers*

$$x = \begin{pmatrix} x_1 & ... & x_m \end{pmatrix}$$

*such that the following equality is true*

(9.1.15) $$\overline{b} = x_*{}^*\overline{a} = x_i\overline{a}^i$$

---

PROOF. The equality (9.1.15) follows from the definition 9.1.8 and the theorem 8.4.9. □

THEOREM 9.1.10. $span(\overline{a}^i, i \in I)$ is subspace of left $A$-vector space of rows.

PROOF. Suppose

$$\overline{b} \in \mathrm{span}(\overline{a}^i, i \in I)$$

$$\overline{c} \in \mathrm{span}(\overline{a}^i, i \in I)$$

According to the theorem 9.1.9

$$(9.1.16) \qquad \overline{b} = b_*{}^* \overline{a} = b_i \overline{a}^i$$

$$(9.1.17) \qquad \overline{c} = c_*{}^* \overline{a} = c_i \overline{a}^i$$

Then

$$\overline{b} + \overline{c} = a_*{}^* b + a_*{}^* c = a_*{}^*(b+c) \in \mathrm{span}(\overline{a}^i, i \in I)$$

$$\overline{b}k = (a_*{}^* b)k = a_*{}^*(bk) \in \mathrm{span}(\overline{a}^i, i \in I)$$

This proves the statement. □

THEOREM 9.1.11. Let $\overline{\overline{e}} = (e^j, j \in J)$ be a basis. Let vectors $\overline{b}$, $\overline{a}^i$ have expansion

$$(9.1.18) \qquad \overline{b} = b_*{}^* e = b_j e^j \qquad b = \begin{pmatrix} b_1 & \dots & b_n \end{pmatrix}$$

$$(9.1.19) \qquad \overline{a^i} = a^i{}_*{}^* e = a^i_j e^j \qquad a^i = \begin{pmatrix} a^i_1 & \dots & a^i_n \end{pmatrix}$$

The vector $\overline{b}$ linearly dependends on vectors $\overline{a}^i$

$$\overline{b} \in span(\overline{a}^i, i \in I)$$

if there exist row of $A$-numbers

$$(9.1.20) \qquad x = \begin{pmatrix} x_1 & \dots & x_m \end{pmatrix}$$

which satisfy the **system of linear equations**

$$x_1 a^1_1 + \dots + x_m a^m_1 = b_1$$
$$(9.1.21) \qquad \dots\dots$$
$$x_1 a^1_n + \dots + x_m a^m_n = b_n$$

PROOF. The equality

$$(9.1.22) \qquad b_j e^j = x_i(a^i_j e^j) = (x_i a^i_j) e^j$$

follows from equalities (9.1.15), (9.1.18), (9.1.19). The equality

$$(9.1.23) \qquad b_j = x_i a^i_j$$

follows from the equality (9.1.22) and the theorem 8.4.11. The equality (9.1.21) follows from the equality (9.1.23). □

Theorem 9.1.12. *Consider the matrix*

$$(9.1.24) \qquad a = \begin{pmatrix} \begin{pmatrix} a_1^1 & ... & a_n^1 \end{pmatrix} \\ ... \\ \begin{pmatrix} a_1^m & ... & a_n^m \end{pmatrix} \end{pmatrix} = \begin{pmatrix} a_1^1 & ... & a_n^1 \\ ... & ... & ... \\ a_1^m & ... & a_n^m \end{pmatrix}$$

*Then the system of linear equations* (9.1.21) *can be represented as*

$$(9.1.25) \qquad \begin{pmatrix} x_1 & ... & x_m \end{pmatrix} {}_*{}^* \begin{pmatrix} a_1^1 & ... & a_n^1 \\ ... & ... & ... \\ a_1^m & ... & a_n^m \end{pmatrix} = \begin{pmatrix} b^1 \\ ... \\ b^n \end{pmatrix}$$

$$(9.1.26) \qquad x_*{}^* a = b$$

Proof. The equality (9.1.25) follows from equalities (9.1.18), (9.1.20), (9.1.21), (9.1.24) and from the definition 8.1.4. □

Definition 9.1.13. *Let $a$ be $_*{}^*$-nonsingular matrix. We call appropriate system of linear equations*

(9.1.26)  $x_*{}^* a = b$

**nonsingular system of linear equations**. □

Theorem 9.1.14. *Solution of nonsingular system of linear equations*

(9.1.26)  $x_*{}^* a = b$

*is determined uniquely and can be presented in either form* [9.2]

$$(9.1.27) \qquad x = b_*{}^* a^{-1_*{}^*}$$

$$(9.1.28) \qquad x = b_*{}^* \mathcal{H} \det({}_*{}^*) a$$

Proof. Multiplying both sides of the equation (9.1.26) from left by $a^{-1_*{}^*}$, we get (9.1.27). Using definition

(8.2.11)  $a^{-1_*{}^*} = \mathcal{H} \det({}_*{}^*) a$

we get (9.1.28). Based on the theorem 8.1.31, the solution is unique. □

## 9.2. Linear Span in Right $A$-Vector Space

### 9.2.1. Right $A$-vector space of columns.

Definition 9.2.1. *Let $V$ be a right $A$-vector space of columns and $\{\overline{a}_i \in V, i \in I\}$ be set of vectors.* **Linear span** *in right $A$-vector space of columns is set* $span(\overline{a}_i, i \in I)$ *of vectors linearly dependent on vectors* $\overline{a}_i$. □

<hr>

[9.2]  We can see a solution of system (9.1.26) in theorem [5]-1.6.1. I repeat this statement because I slightly changed the notation.

Theorem 9.2.2. *Let $V$ be a right $A$-vector space of columns and row*

$$\overline{a} = \begin{pmatrix} \overline{a}_1 & ... & \overline{a}_m \end{pmatrix} = (\overline{a}_i, i \in I)$$

*be set of vectors. The vector $\overline{b}$ linearly dependends on vectors $\overline{a}_i$*

$$\overline{b} \in span(\overline{a}_i, i \in I)$$

*if there exist column of $A$-numbers*

$$x = \begin{pmatrix} x^1 \\ ... \\ x^m \end{pmatrix}$$

*such that the following equality is true*

(9.2.1) $$\overline{b} = \overline{a}_*{}^* x = \overline{a}_i x^i$$

Proof. The equality (9.2.1) follows from the definition 9.2.1 and the theorem 8.4.16. $\square$

Theorem 9.2.3. *$span(\overline{a}_i, i \in I)$ is subspace of right $A$-vector space of columns.*

Proof. Suppose

$$\overline{b} \in span(\overline{a}_i, i \in I)$$
$$\overline{c} \in span(\overline{a}_i, i \in I)$$

According to the theorem 9.2.2

(9.2.2) $$\overline{b} = \overline{a}_*{}^* b = \overline{a}_i b^i$$

(9.2.3) $$\overline{c} = \overline{a}_*{}^* c = \overline{a}_i c^i$$

Then

$$\overline{b} + \overline{c} = \overline{a}_*{}^* b + \overline{a}_*{}^* c = \overline{a}_*{}^*(b + c) \in span(\overline{a}_i, i \in I)$$
$$\overline{b}k = (\overline{a}_*{}^* b)k = \overline{a}_*{}^*(bk) \in span(\overline{a}_i, i \in I)$$

This proves the statement. $\square$

Theorem 9.2.4. *Let $\overline{\overline{e}} = (e_j, j \in J)$ be a basis. Let vectors $\overline{b}$, $\overline{a}_i$ have expansion*

(9.2.4) $$\overline{b} = e_*{}^* b = e_j b^j \qquad b = \begin{pmatrix} b^1 \\ ... \\ b^n \end{pmatrix}$$

(9.2.5) $$\overline{a_i} = e_*{}^* a_i = e_j a_i^j \qquad a_i = \begin{pmatrix} a_i^1 \\ ... \\ a_i^n \end{pmatrix}$$

*The vector $\overline{b}$ linearly dependends on vectors $\overline{a}_i$*

$$\overline{b} \in span(\overline{a}_i, i \in I)$$

*if there exist column of A-numbers*

$$(9.2.6) \qquad x = \begin{pmatrix} x^1 \\ ... \\ x^m \end{pmatrix}$$

*which satisfy the* **system of linear equations**

$$(9.2.7) \qquad \begin{array}{c} a_1^1 x^1 + ... + a_m^1 x^m = b^1 \\ ...... \\ a_1^n x^1 + ... + a_m^n x^m = b^n \end{array}$$

Proof. The equality

$$(9.2.8) \qquad e_j b^j = (e_j a_i^j) x^i = e_j (a_i^j x^i)$$

follows from equalities (9.2.1), (9.2.4), (9.2.5). The equality

$$(9.2.9) \qquad b^j = a_i^j x^i$$

follows from the equality (9.2.8) and the theorem 8.4.18. The equality (9.2.7) follows from the equality (9.2.9). □

Theorem 9.2.5. *Consider the matrix*

$$(9.2.10) \qquad a = \begin{pmatrix} \begin{pmatrix} a_1^1 \\ ... \\ a_1^n \end{pmatrix} & ... & \begin{pmatrix} a_m^1 \\ ... \\ a_m^n \end{pmatrix} \end{pmatrix} = \begin{pmatrix} a_1^1 & ... & a_m^1 \\ ... & ... & ... \\ a_1^n & ... & a_m^n \end{pmatrix}$$

*Then the system of linear equations* (9.2.7) *can be represented as*

$$(9.2.11) \qquad \begin{pmatrix} a_1^1 & ... & a_m^1 \\ ... & ... & ... \\ a_1^n & ... & a_m^n \end{pmatrix} {}_*^* \begin{pmatrix} x^1 \\ ... \\ x^m \end{pmatrix} = \begin{pmatrix} b^1 \\ ... \\ b^n \end{pmatrix}$$

$$(9.2.12) \qquad a_*^* x = b$$

Proof. The equality (9.2.11) follows from equalities (9.2.4), (9.2.6), (9.2.7), (9.2.10) and from the definition 8.1.3. □

Definition 9.2.6. *Let a be* $_*^*$-*nonsingular matrix. We call appropriate system of linear equations*

$$\boxed{(9.2.12) \quad a_*^* x = b}$$

**nonsingular system of linear equations**. □

THEOREM 9.2.7. *Solution of nonsingular system of linear equations*

$$(9.2.12) \quad a_*{}^* x = b$$

*is determined uniquely and can be presented in either form* [9.3]

$$(9.2.13) \qquad\qquad x = a^{-1*}{}_*{}^* b$$

$$(9.2.14) \qquad\qquad x = \mathcal{H} \det({}_*{}^*) a_*{}^* b$$

PROOF. Multiplying both sides of the equation (9.2.12) from left by $a^{-1*}$, we get (9.2.13). Using definition

$$(8.2.11) \quad a^{-1*} = \mathcal{H} \det({}_*{}^*) a$$

we get (9.2.14). Based on the theorem 8.1.31, the solution is unique.    $\square$

**9.2.2. Right $A$-vector space of rows.**

DEFINITION 9.2.8. *Let $V$ be a right $A$-vector space of rows and $\{\overline{a}^i \in V, i \in I\}$ be set of vectors.* **Linear span** *in right $A$-vector space of rows is set* $span(\overline{a}^i, i \in I)$ *of vectors linearly dependent on vectors* $\overline{a}^i$.    $\square$

THEOREM 9.2.9. *Let $V$ be a right $A$-vector space of rows and column*

$$\overline{a} = \begin{pmatrix} \overline{a}^1 \\ ... \\ \overline{a}^m \end{pmatrix} = (\overline{a}^i, i \in I)$$

*be set of vectors. The vector $\overline{b}$ linearly dependends on vectors $\overline{a}^i$*

$$\overline{b} \in span(\overline{a}^i, i \in I)$$

*if there exist row of $A$-numbers*

$$x = \begin{pmatrix} x_1 & ... & x_m \end{pmatrix}$$

*such that the following equality is true*

$$(9.2.15) \qquad\qquad \overline{b} = \overline{a}^*{}_* x = \overline{a}^i x_i$$

PROOF. The equality (9.2.15) follows from the definition 9.2.8 and the theorem 8.4.23.    $\square$

THEOREM 9.2.10. *$span(\overline{a}^i, i \in I)$ is subspace of right $A$-vector space of rows.*

PROOF. Suppose

$$\overline{b} \in \mathrm{span}(\overline{a}^i, i \in I)$$

$$\overline{c} \in \mathrm{span}(\overline{a}^i, i \in I)$$

---

[9.3] We can see a solution of system (9.2.12) in theorem [5]-1.6.1. I repeat this statement because I slightly changed the notation.

According to the theorem 9.2.9

$$(9.2.16) \qquad \overline{b} = \overline{a}^*{}_*b = \overline{a}^i b_i$$

$$(9.2.17) \qquad \overline{c} = \overline{a}^*{}_*c = \overline{a}^i c_i$$

Then

$$\overline{b} + \overline{c} = \overline{a}^*{}_*b + \overline{a}^*{}_*c = \overline{a}^*{}_*(b+c) \in \text{span}(\overline{a}^i, i \in I)$$

$$\overline{b}k = (\overline{a}^*{}_*b)k = \overline{a}^*{}_*(bk) \in \text{span}(\overline{a}^i, i \in I)$$

This proves the statement. $\qquad\square$

**THEOREM 9.2.11.** *Let* $\overline{\overline{e}} = (e^j, j \in J)$ *be a basis. Let vectors* $\overline{b}, \overline{a}^i$ *have expansion*

$$(9.2.18) \qquad \overline{b} = e^*{}_*b = e^j b_j \qquad b = \begin{pmatrix} b_1 & \dots & b_n \end{pmatrix}$$

$$(9.2.19) \qquad \overline{a^i} = e^*{}_*a^i = e^j a^i_j \qquad a^i = \begin{pmatrix} a^i_1 & \dots & a^i_n \end{pmatrix}$$

*The vector* $\overline{b}$ *linearly dependends on vectors* $\overline{a}^i$

$$\overline{b} \in span(\overline{a}^i, i \in I)$$

*if there exist row of A-numbers*

$$(9.2.20) \qquad x = \begin{pmatrix} x_1 & \dots & x_m \end{pmatrix}$$

*which satisfy the* **system of linear equations**

$$(9.2.21) \qquad \begin{matrix} a^1_1 x_1 + \dots + a^m_1 x_m = b_1 \\ \dots\dots \\ a^1_n x_1 + \dots + a^m_n x_m = b_n \end{matrix}$$

PROOF. The equality

$$(9.2.22) \qquad e^j b_j = (e^j a^i_j) x_i = e^j (a^i_j x_i)$$

follows from equalities (9.2.15), (9.2.18), (9.2.19). The equality

$$(9.2.23) \qquad b_j = a^i_j x_i$$

follows from the equality (9.2.22) and the theorem 8.4.25. The equality (9.2.21) follows from the equality (9.2.23). $\qquad\square$

**THEOREM 9.2.12.** *Consider the matrix*

$$(9.2.24) \qquad a = \begin{pmatrix} \begin{pmatrix} a^1_1 & \dots & a^1_n \end{pmatrix} \\ \dots \\ \begin{pmatrix} a^m_1 & \dots & a^m_n \end{pmatrix} \end{pmatrix} = \begin{pmatrix} a^1_1 & \dots & a^1_n \\ \dots & \dots & \dots \\ a^m_1 & \dots & a^m_n \end{pmatrix}$$

*Then the system of linear equations* (9.2.21) *can be represented as*

$$(9.2.25) \qquad \begin{pmatrix} a^1_1 & \dots & a^1_n \\ \dots & \dots & \dots \\ a^m_1 & \dots & a^m_n \end{pmatrix}^*{}_* \begin{pmatrix} x_1 & \dots & x_m \end{pmatrix} = \begin{pmatrix} b_1 & \dots & b_n \end{pmatrix}$$

$$(9.2.26) \qquad\qquad a^*{}_*x = b$$

Proof. The equality (9.2.25) follows from equalities (9.2.18), (9.2.20), (9.2.21), (9.2.24) and from the definition 8.1.4. ∎

Definition 9.2.13. *Let $a$ be $^*{}_*$-nonsingular matrix. We call appropriate system of linear equations*

$$(9.2.26) \quad a^*{}_*x = b$$

**nonsingular system of linear equations**. ∎

Theorem 9.2.14. *Solution of nonsingular system of linear equations*

$$(9.2.26) \quad a^*{}_*x = b$$

*is determined uniquely and can be presented in either form*[9.4]

$$(9.2.27) \qquad\qquad x = a^{-1^*}{}_*{}^*{}_*b$$

$$(9.2.28) \qquad\qquad x = \mathcal{H}\det({}^*{}_*)a^*{}_*b$$

Proof. Multiplying both sides of the equation (9.2.26) from left by $a^{-1^*}{}_*$, we get (9.2.27). Using definition

$$(8.3.11) \quad a^{-1^*}{}_* = \mathcal{H}\det({}^*{}_*)a$$

we get (9.2.28). Based on the theorem 8.1.31, the solution is unique. ∎

### 9.3. Rank of Matrix

In this section, we will make the following assumption.

- $i \in M$, $|M| = m$, $j \in N$, $|N| = n$.
- $a = (a_j^i)$ is an arbitrary matrix.
- $k$, $s \in S \supseteq M$, $l$, $t \in T \supseteq N$, $k = |S| = |T|$.
- $p \in M \setminus S$, $r \in N \setminus T$.

Matrix $a_T^S$ is $k \times k$ submatrix.

Definition 9.3.1. *If submatrix $a_T^S$ is $_*{}^*$-nonsingular matrix then we say that $_*{}^*$-rank of matrix $a$ is not less then $k$. $_*{}^*$-**rank of matrix** $a$, $\operatorname{rank}({}_*{}^*)a$, is the maximal value of $k$. We call an appropriate submatrix the $_*{}^*$-**major submatrix**.* ∎

Definition 9.3.2. *If submatrix $a_T^S$ is $^*{}_*$-nonsingular matrix then we say that $^*{}_*$-rank of matrix $a$ is not less then $k$. $^*{}_*$-**rank of matrix** $a$, $\operatorname{rank}({}^*{}_*)a$, is the maximal value of $k$. We call an appropriate submatrix the $^*{}_*$-**major submatrix**.* ∎

Theorem 9.3.3. *Let matrix $a$ be $_*{}^*$-singular matrix and submatrix $a_T^S$ be $_*{}^*$-major submatrix. Then*

$$\det({}_*{}^*)_r^p\, a_{T \cup \{r\}}^{S \cup \{p\}} = 0$$

Theorem 9.3.4. *Let matrix $a$ be $^*{}_*$-singular matrix and submatrix $a_T^S$ be $^*{}_*$-major submatrix. Then*

$$\det({}^*{}_*)_r^p\, a_{T \cup \{r\}}^{S \cup \{p\}} = 0$$

---

[9.4] We can see a solution of system (9.2.26) in theorem [5]-1.6.1. I repeat this statement because I slightly changed the notation.

Proof of theorem 9.3.3.          To understand why submatrix

$$(9.3.1) \qquad b = a^{S \cup \{p\}}_{T \cup \{r\}}$$

does not have $_*{}^*$-inverse matrix, [9.5] we assume that there exists $_*{}^*$-inverse matrix $b^{-1_*{}^*}$. We write down the system of linear equations

$$(8.2.3) \quad a^{[J]}_{I}{}_*{}^*(a^{-1_*{}^*})^{I}_{J} + a^{[J]}_{[I]}{}_*{}^*(a^{-1_*{}^*})^{[I]}_{J} = 0$$

$$(8.2.4) \quad a^{J}_{I}{}_*{}^*(a^{-1_*{}^*})^{I}_{J} + a^{J}_{[I]}{}_*{}^*(a^{-1_*{}^*})^{[I]}_{J} = E_k$$

in case $I = \{r\}$, $J = \{p\}$  (then $[I] = T$, $[J] = S$ )

$$(9.3.2) \qquad b^{S}_{r} \ b^{-1_*{}^{*r}}_{\quad p} + b^{S}_{T}{}_*{}^*b^{-1_*{}^*}{}^{T}_{p} = 0$$

$$(9.3.3) \qquad b^{p}_{r} \ b^{-1_*{}^{*r}}_{\quad p} + b^{p}_{T}{}_*{}^*b^{-1_*{}^*}{}^{T}_{p} = 1$$

and will try to solve this system. We multiply (9.3.2) by $(b^{S}_{T})^{-1_*{}^*}$

$$(9.3.4) \qquad (b^{S}_{T})^{-1_*{}^*}{}_*{}^*b^{S}_{r} \ b^{-1_*{}^{*r}}_{\quad p} + b^{-1_*{}^*}{}^{T}_{p} = 0$$

Now we can substitute (9.3.4) into (9.3.3)

$$(9.3.5) \qquad b^{p}_{r} \ b^{-1_*{}^{*r}}_{\quad p} - b^{p}_{T}{}_*{}^*(b^{S}_{T})^{-1_*{}^*}{}_*{}^*b^{S}_{r} \ b^{-1_*{}^{*r}}_{\quad p} = 1$$

From (9.3.5) it follows that

$$(9.3.6) \qquad (b^{p}_{r} - b^{p}_{T}{}_*{}^*(b^{S}_{T})^{-1_*{}^*}{}_*{}^*b^{S}_{r}) \ b^{-1_*{}^{*r}}_{\quad p} = 1$$

The equality

$$(9.3.7) \qquad (\det(_*{}^*)^{p}_{r} b)b^{-1_*{}^{*r}}_{\quad p} = 1$$

follows from the equality

$$(8.2.9) \quad \det(_*{}^*)^{j}_{i} \, a = a^{j}_{i} - a^{j}_{[i]}{}_*{}^*\left(a^{[j]}_{[i]}\right)^{-1_*{}^*}{}_*{}^*a^{[j]}_{i} = ((a^{-1_*{}^*})^{i}_{j})^{-1}$$

and from the equality (9.3.6). Thus we proved that quasideterminant $\det(_*{}^*)^{p}_{r} b$ is defined and its equation to 0 is necessary and sufficient condition that the matrix $b$ is singular.          □

Proof of theorem 9.3.4.          To understand why submatrix

$$(9.3.8) \qquad b = a^{S \cup \{p\}}_{T \cup \{r\}}$$

does not have $^*{}_*$-inverse matrix, [9.6] we assume that there exists $^*{}_*$-inverse matrix $b^{-1^*{}_*}$. We write down the system of linear equations

$$(8.3.3) \quad a^{I}_{[J]}{}^*{}_*(a^{-1^*{}_*})^{J}_{I} + a^{[I]}_{[J]}{}^*{}_*(a^{-1^*{}_*})^{J}_{[I]} = 0$$

---

[9.5]It is natural to expect relationship between $_*{}^*$-singularity of the matrix and its $_*{}^*$-quasideterminant similar to relationship which is known in commutative case. However $_*{}^*$-quasideterminant is defined not always. For instance, it is not defined when $_*{}^*$-inverse matrix has too much elements equal 0. As it follows from this theorem, the $_*{}^*$-quasideterminant is undefined also in case when $_*{}^*$-rank of the matrix is less then $n - 1$.

[9.6]It is natural to expect relationship between $^*{}_*$-singularity of the matrix and its $^*{}_*$-quasideterminant similar to relationship which is known in commutative case. However $^*{}_*$-quasideterminant is defined not always. For instance, it is not defined when $^*{}_*$-inverse matrix has too much elements equal 0. As it follows from this theorem, the $^*{}_*$-quasideterminant is undefined also in case when $^*{}_*$-rank of the matrix is less then $n - 1$.

(8.3.4)    $a_J^I {}_*{}^*(a^{-1}{}^*)_I^J + a_J^{[I]}{}_*{}^*(a^{-1}{}^*)_{[I]}^J = E_k$

in case $I = \{p\}$, $J = \{r\}$  (then $[I] = S$, $[J] = T$ )

(9.3.9)    $b_T^p \ b^{-1}{}^*{}_*{}_p^r + b_T^{S}{}_*{}^* b^{-1}{}^*{}_*{}_S^r = 0$

(9.3.10)    $b_r^p \ b^{-1}{}^*{}_*{}_p^r + b_r^{S}{}_*{}^* b^{-1}{}^*{}_*{}_S^r = 1$

and will try to solve this system. We multiply (9.3.9) by $(b_T^S)^{-1}{}^*{}^*$

(9.3.11)    $(b_T^S)^{-1}{}^*{}_*{}^* b_T^p \ b^{-1}{}^*{}_*{}_p^r + b^{-1}{}^*{}_*{}_S^r = 0$

Now we can substitute (9.3.11) into (9.3.10)

(9.3.12)    $b_r^p \ b^{-1}{}^*{}_*{}_p^r - b_r^{S}{}_*{}^* (b_T^S)^{-1}{}^*{}_*{}^* b_T^p \ b^{-1}{}^*{}_*{}_p^r = 1$

From (9.3.12) it follows that

(9.3.13)    $(b_r^p - b_r^{S}{}_*{}^* (b_T^S)^{-1}{}^*{}_*{}^* b_T^p) \ b^{-1}{}^*{}_*{}_p^r = 1$

The equality

(9.3.14)    $(\det({}^*_*)_r^p \ b) b^{-1}{}^*{}_*{}_p^r = 1$

follows from the equality

(8.3.9)    $\det({}^*_*)_j^i \ a = a_j^i - a_j^{[i]}{}_*{}^* \left(a_{[j]}^{[i]}\right)^{-1}{}^*{}_*{}^* a_{[j]}^i = ((a^{-1}{}^*)_j^i)^{-1}$

from the equality (9.3.13). Thus we proved that quasideterminant $\det({}^*_*)_r^p \ b$ is defined and its equation to 0 is necessary and sufficient condition that the matrix $b$ is singular.    □

THEOREM 9.3.5. *Let $a$ be a matrix,*

$$\mathrm{rank}({}^*_*)a = k < n$$

*and $a_T^S$ is ${}^*_*$-major submatrix. Then column $a_r$ is a left linear composition of columns $a_T$*

(9.3.15)    $a_{N \setminus T} = R^*{}_* a_T$

(9.3.16)    $a_r = R_r^*{}_* a_T$

(9.3.17)    $a_r^t = R_r^t a_t^t$

PROOF. If the matrix $a$ has $k$ rows, then assuming that column $a_r$ is a left linear combination (9.3.16) of columns $a_T$ with coefficients $R_r^t$, we get system of left linear equations (9.3.17). We assume that $x^t = R_r^t$ are unknown variables in the system of left linear equations (9.3.17). According to the theorem 9.1.7, the system of left linear equations (9.3.17) has a unique solution and this solution is nontrivial because all ${}^*_*$-

THEOREM 9.3.6. *Let $a$ be a matrix,*

$$\mathrm{rank}({}_*{}^*)a = k < n$$

*and $a_T^S$ is ${}_*{}^*$-major submatrix. Then column $a_r$ is a right linear composition of columns $a_T$*

(9.3.23)    $a_{N \setminus T} = a_T{}_*{}^* R$

(9.3.24)    $a_r = a_T{}_*{}^* R_r$

(9.3.25)    $a_r^a = a_t^a \ R_r^t$

PROOF. If the matrix $a$ has $k$ rows, then assuming that column $a_r$ is a right linear combination (9.3.24) of columns $a_T$ with coefficients $R_r^t$, we get system of right linear equations (9.3.25). We assume that $x^t = R_r^t$ are unknown variables in the system of right linear equations (9.3.25). According to the theorem 9.2.7, the system of right linear equations (9.3.25) has a unique solution and this solution is nontrivial because all quasideterminants are different from $0$.

It remains to prove this statement in case when a number of rows of the matrix $a$ is more than $k$. I get column $a_r$ and row $a^p$. According to assumption, submatrix $a_{T \cup \{r\}}^{S \cup \{p\}}$ is a $_*^*$-singular matrix and its $_*^*$-quasideterminant

$$(9.3.18) \qquad \det({}^*_*)_r^p \, a_{T \cup \{r\}}^{S \cup \{p\}} = 0$$

According to the equality (8.3.9) the equality (9.3.18) has form

$$a_r^p - a_r^{S*}*((a_T^S)^{-1^*})^**a_T^p = 0$$

Matrix

$$(9.3.19) \qquad R_r = a_r^{S*}*((a_T^S)^{-1^*})$$

does not depend on $p$, Therefore, for any $p \in M \setminus S$

$$(9.3.20) \qquad a_r^p = R_r * a_T^p$$

To prove the equality (9.3.17), it is necessary to prove the equality

$$(9.3.21) \qquad a_r^k = R_r * a_T^k$$

From the equality

$$((a_T^S)^{-1^*})^**a_T^k = \delta_S^k$$

it follows that

$$(9.3.22) \qquad \begin{aligned} a_r^k &= a_r^{S*}*\delta_S^k \\ &= a_r^{S*}*((a_T^S)^{-1^*})^**a_T^k \end{aligned}$$

The equality (9.3.21) follows from equalities (9.3.19), (9.3.22). The equality (9.3.17) follows from equalities (9.3.20), (9.3.21).
□

---

It remains to prove this statement in case when a number of rows of the matrix $a$ is more than $k$. I get column $a_r$ and row $a^p$. According to assumption, submatrix $a_{T \cup \{r\}}^{S \cup \{p\}}$ is a $_*^*$-singular matrix and its $_*^*$-quasideterminant

$$(9.3.26) \qquad \det(_*^*)_r^p \, a_{T \cup \{r\}}^{S \cup \{p\}} = 0$$

According to the equality (8.2.9) the equality (9.3.26) has form

$$a_r^p - a_T^p*((a_T^S)^{-1^*})_**a_r^S = 0$$

Matrix

$$(9.3.27) \qquad R_r = ((a_T^S)^{-1^*})_**a_r^S$$

does not depend on $p$, Therefore, for any $p \in M \setminus S$

$$(9.3.28) \qquad a_r^p = a_T^p * R_r$$

To prove the equality (9.3.25), it is necessary to prove the equality

$$(9.3.29) \qquad a_r^k = a_T^k * R_r$$

From the equality

$$a_{T*}^k*((a_T^S)^{-1^*}) = \delta_S^k$$

it follows that

$$(9.3.30) \qquad \begin{aligned} a_r^k &= \delta_{S*}^k*a_r^S \\ &= a_{T*}^k*((a_T^S)^{-1^*})_**a_r^S \end{aligned}$$

The equality (9.3.29) follows from equalities (9.3.27), (9.3.30). The equality (9.3.25) follows from equalities (9.3.28), (9.3.29).
□

---

THEOREM 9.3.7. *Let $a$ be a matrix,*

$$\operatorname{rank}(^*_*)a = k < m$$

*and $a_T^S$ is $_*^*$-major submatrix. Then row $a^p$ is a right linear composition of rows $a^S$*

$$(9.3.31) \qquad a^{M \setminus S} = a^{S}*R$$

$$(9.3.32) \qquad a^p = a^{S}*R^p$$

$$(9.3.33) \qquad a_p^b = a_s^b R_p^s$$

PROOF. If the matrix $a$ has $k$ columns,

---

THEOREM 9.3.8. *Let $a$ be a matrix,*

$$\operatorname{rank}(_*^*)a = k < m$$

*and $a_T^S$ is $_*^*$-major submatrix. Then row $a^p$ is a left linear composition of rows $a^S$*

$$(9.3.39) \qquad a^{M \setminus S} = R_* * a^S$$

$$(9.3.40) \qquad a^p = R^p * a^S$$

$$(9.3.41) \qquad a_b^p = R_s^p a_b^s$$

PROOF. If the matrix $a$ has $k$ columns, then assuming that row $a^p$ is a left lin-then assuming that row $a^p$ is a right linear combination (9.3.32) of rows $a^S$ with coefficients $R^p_s$, we get system of right linear equations (9.3.33). We assume that $x_s = R^p_s$ are unknown variables in the system of right linear equations (9.3.33). According to the theorem 9.2.14, the system of right linear equations (9.3.33) has a unique solution and this solution is nontrivial because all $_*{}^*$-quasideterminants are different from 0.

It remains to prove this statement in case when a number of columns of the matrix $a$ is more than $k$. I get row $a^p$ and column $a_r$. According to assumption, submatrix $a^{S \cup \{p\}}_{T \cup \{r\}}$ is a $_*{}^*$-singular matrix and its $_*{}^*$-quasideterminant

$$(9.3.34) \qquad \det({}_*{}^*)^p_r\, a^{S \cup \{p\}}_{T \cup \{r\}} = 0$$

According to the equality (8.3.9) the equality (9.3.34) has form

$$a^p_r - a^{S}_r{}_*{}^*((a^S_T)^{-1^{*^*}})^*{}_* a^p_T = 0$$

Matrix

$$(9.3.35) \qquad R^p = ((a^S_T)^{-1^{*^*}})^*{}_* a^p_T$$

does not depend on $r$, Therefore, for any $r \in N \setminus T$

$$(9.3.36) \qquad a^p_r = a^{S}_r{}_*{}^* R^p$$

To prove the equality (9.3.33), it is necessary to prove the equality

$$(9.3.37) \qquad a^p_l = a^{S}_l{}_*{}^* R^p$$

From the equality

$$a^{S}_l{}_*((a^S_T)^{-1^{*^*}}) = \delta^T_l$$

it follows that

$$(9.3.38) \qquad a^p_l = \delta^{T}_l{}_*{}^* a^p_T$$
$$= a^{S}_l{}_*{}^*((a^S_T)^{-1^{*^*}})^*{}_* a^p_T$$

The equality (9.3.37) follows from equalities (9.3.35), (9.3.38). The equality (9.3.33) follows from equalities (9.3.36), (9.3.37). □

ear combination (9.3.40) of rows $a^S$ with coefficients $R^p_s$, we get system of left linear equations (9.3.41). We assume that $x_s = R^p_s$ are unknown variables in the system of left linear equations (9.3.41). According to the theorem 9.1.14, the system of left linear equations (9.3.41) has a unique solution and this solution is nontrivial because all $_*{}^*$-quasideterminants are different from 0.

It remains to prove this statement in case when a number of columns of the matrix $a$ is more than $k$. I get row $a^p$ and column $a_r$. According to assumption, submatrix $a^{S \cup \{p\}}_{T \cup \{r\}}$ is a $_*{}^*$-singular matrix and its $_*{}^*$-quasideterminant

$$(9.3.42) \qquad \det({}_*{}^*)^p_r\, a^{S \cup \{p\}}_{T \cup \{r\}} = 0$$

According to the equality (8.2.9) the equality (9.3.42) has form

$$a^p_r - a^p_T{}_*{}^*((a^S_T)^{-1^{*^*}})_*{}^* a^S_r = 0$$

Matrix

$$(9.3.43) \qquad R^p = a^p_T{}_*{}^*((a^S_T)^{-1^{*^*}})$$

does not depend on $r$, Therefore, for any $r \in N \setminus T$

$$(9.3.44) \qquad a^p_r = R^p{}_*{}^* a^S_r$$

To prove the equality (9.3.41), it is necessary to prove the equality

$$(9.3.45) \qquad a^p_l = R^p{}_*{}^* a^S_l$$

From the equality

$$((a^S_T)^{-1^{*^*}})_*{}^* a^S_l = \delta^T_l$$

it follows that

$$(9.3.46) \qquad a^p_l = a^p_T{}_*{}^* \delta^T_l$$
$$= a^p_T{}_*{}^*((a^S_T)^{-1^{*^*}})_*{}^* a^S_l$$

The equality (9.3.45) follows from equalities (9.3.43), (9.3.46). The equality (9.3.41) follows from equalities (9.3.44), (9.3.45). □

THEOREM 9.3.9. *Let matrix $a$ have $n$ columns. If*

$$\mathrm{rank}({}_*{}^*)a = k < n$$

THEOREM 9.3.10. *Let matrix $a$ have $n$ columns. If*

$$\mathrm{rank}({}^*{}_*)a = k < n$$

*then columns of the matrix are right linearly dependent*

$$a_*{}^*\lambda = 0$$

PROOF. Let column $a_r$ be a right linear composition of rows (9.3.24). We assume $\lambda^r = -1$, $\lambda^t = R_r^t$ and the rest $\lambda^c = 0$. □

THEOREM 9.3.11. *Let matrix $a$ have $m$ rows. If*

$$\mathrm{rank}(_*{}^*)a = k < m$$

*then rows of the matrix are left linearly dependent*

$$\lambda_*{}^*a = 0$$

PROOF. Let row $a^p$ be a left linear composition of cols (9.3.40). We assume $\lambda_p = -1$, $\lambda_s = R_s^p$ and the rest $\lambda_c = 0$. □

THEOREM 9.3.13. *Let $V$ be right $A$-vector space of columns. Let $(\overline{a}^i, i \in M, |M| = m)$ be set of vectors linearly independent from left. Then $_*{}^*$-rank of their coordinate matrix equal $m$.*

PROOF. Let $\overline{e}$ be the basis of right $A$-vector space $V$. According to the theorem 8.4.17, the coordinate matrix of set of vectors $(\overline{a}^i)$ relative basis $\overline{e}$ consists from columns which are coordinate matrices of vectors $\overline{a}^i$ relative the basis $\overline{e}$. Therefore $_*{}^*$-rank of this matrix cannot be more then $m$.

Let $_*{}^*$-rank of the coordinate matrix be less then $m$. According to the theorem 9.3.9, columns of matrix are right linear dependent

$$(9.3.47) \qquad a_*{}^*\lambda = 0$$

Suppose $c = a_*{}^*\lambda$. From the equation (9.3.47) it follows that right linear composition of rows

$$\overline{e}_*{}^*c = 0$$

of vectors of basis equal 0. This contradicts to statement that vectors $\overline{e}$ form

*then columns of the matrix are left linearly dependent*

$$\lambda^*{}_*a = 0$$

PROOF. Let column $a_r$ be a left linear composition of rows (9.3.16). We assume $\lambda^r = -1$, $\lambda^t = R_r^t$ and the rest $\lambda^c = 0$. □

THEOREM 9.3.12. *Let matrix $a$ have $m$ rows. If*

$$\mathrm{rank}(^*{}_*)a = k < m$$

*then rows of the matrix are right linearly dependent*

$$a^*{}_*\lambda = 0$$

PROOF. Let row $a^p$ be a right linear composition of cols (9.3.32). We assume $\lambda_p = -1$, $\lambda_s = R_s^p$ and the rest $\lambda_c = 0$. □

THEOREM 9.3.14. *Let $V$ be left $A$-vector space of columns. Let $(\overline{a}^i, i \in M, |M| = m)$ be set of vectors linearly independent from left. Then $^*{}_*$-rank of their coordinate matrix equal $m$.*

PROOF. Let $\overline{\overline{e}}$ be the basis of left $A$-vector space $V$. According to the theorem 8.4.3, the coordinate matrix of set of vectors $(\overline{a}^i)$ relative basis $\overline{\overline{e}}$ consists from columns which are coordinate matrices of vectors $\overline{a}^i$ relative the basis $\overline{\overline{e}}$. Therefore $^*{}_*$-rank of this matrix cannot be more then $m$.

Let $^*{}_*$-rank of the coordinate matrix be less then $m$. According to the theorem 9.3.10, columns of matrix are left linear dependent

$$(9.3.48) \qquad \lambda^*{}_*a = 0$$

Suppose $c = \lambda^*{}_*a$. From the equation (9.3.48) it follows that left linear composition of rows

$$c^*{}_*\overline{\overline{e}} = 0$$

of vectors of basis equal 0. This contradicts to statement that vectors $\overline{e}$ form basis. We proved statement of theorem. $\square$

basis. We proved statement of theorem. $\square$

THEOREM 9.3.15. *Let $V$ be right $A$-vector space of rows. Let $(\overline{a}^i, i \in M, |M| = m)$ be set of vectors linearly independent from left. Then $^*_*$-rank of their coordinate matrix equal $m$.*

PROOF. Let $\overline{e}$ be the basis of right $A$-vector space $V$. According to the theorem 8.4.24, the coordinate matrix of set of vectors $(\overline{a}^i)$ relative basis $\overline{e}$ consists from rows which are coordinate matrices of vectors $\overline{a}^i$ relative the basis $\overline{e}$. Therefore $^*_*$-rank of this matrix cannot be more then $m$.

Let $^*_*$-rank of the coordinate matrix be less then $m$. According to the theorem 9.3.12, rows of matrix are right linear dependent

(9.3.49)        $a^*_* \lambda = 0$

Suppose $c = a^*_* \lambda$. From the equation (9.3.49) it follows that right linear composition of rows

$$\overline{e}^*_* c = 0$$

of vectors of basis equal 0. This contradicts to statement that vectors $\overline{e}$ form basis. We proved statement of theorem. $\square$

THEOREM 9.3.16. *Let $V$ be left $A$-vector space of rows. Let $(\overline{a}^i, i \in M, |M| = m)$ be set of vectors linearly independent from left. Then $_*^*$-rank of their coordinate matrix equal $m$.*

PROOF. Let $\overline{\overline{e}}$ be the basis of left $A$-vector space $V$. According to the theorem 8.4.10, the coordinate matrix of set of vectors $(\overline{a}^i)$ relative basis $\overline{\overline{e}}$ consists from rows which are coordinate matrices of vectors $\overline{a}^i$ relative the basis $\overline{\overline{e}}$. Therefore $_*^*$-rank of this matrix cannot be more then $m$.

Let $_*^*$-rank of the coordinate matrix be less then $m$. According to the theorem 9.3.11, rows of matrix are left linear dependent

(9.3.50)        $\lambda_*^* a = 0$

Suppose $c = \lambda_*^* a$. From the equation (9.3.50) it follows that left linear composition of rows

$$c_*^* \overline{e} = 0$$

of vectors of basis equal 0. This contradicts to statement that vectors $\overline{e}$ form basis. We proved statement of theorem. $\square$

## 9.4. System of Linear Equations $a_*^* x = b$

DEFINITION 9.4.1. *Let $a$ be a matrix of system of left linear equations*

$$x_1 a_1^1 + ... + x_m a_1^m = b_1$$

(9.4.1)                ......

$$x_1 a_n^1 + ... + x_m a_n^m = b_n$$

*We call matrix*

(9.4.2)        $$\begin{pmatrix} a_i^j \\ b_i \end{pmatrix} = \begin{pmatrix} a_1^1 & ... & a_n^1 \\ ... & ... & ... \\ a_1^m & ... & a_n^m \\ b_1 & ... & b_n \end{pmatrix}$$

*an* **augmented matrix** *of the system of left linear equations* (9.4.1).        $\square$

**Definition 9.4.2.** *Let a be a matrix of system of right linear equations*

$$(9.4.3) \quad \begin{array}{llll} a_1^1 x^1 & +... & +a_n^1 x^n & = b^1 \\ ... & ... & ... & ... \\ a_1^m x^1 & +... & +a_n^m x^n & = b^m \end{array}$$

*We call matrix*

$$(9.4.4) \qquad \begin{pmatrix} a_i^j & b^j \end{pmatrix} = \begin{pmatrix} a_1^1 & ... & a_n^1 & b^1 \\ ... & ... & ... & ... \\ a_1^m & ... & a_n^m & b^m \end{pmatrix}$$

*an* **augmented matrix** *of the system of right linear equations* (9.4.3). □

**Theorem 9.4.3.** *System of right linear equations* (9.4.3) *has a solution iff*

$$(9.4.5) \qquad \mathrm{rank}(_{\ast}{}^{\ast})(a_i^j) = \mathrm{rank}(_{\ast}{}^{\ast}) \begin{pmatrix} a_i^j & b^j \end{pmatrix}$$

**Proof.** Let $a_T^S$ be $_{\ast}{}^{\ast}$-major submatrix of matrix $a$.

Let a system of right linear equations (9.4.3) have solution $x^i = d^i$. Then

$$(9.4.6) \qquad a_{\ast}{}^{\ast}d = b$$

Equation (9.4.6) can be rewritten in form

$$(9.4.7) \qquad a_{T\ast}{}^{\ast}d^T + a_{N\setminus T\ast}{}^{\ast}d^{N\setminus T} = b$$

Substituting (9.3.23)] into (9.4.7) we get

$$(9.4.8) \qquad a_{T\ast}{}^{\ast}d^T + a_{T\ast}{}^{\ast}r_{\ast}{}^{\ast}d^{N\setminus T} = b$$

From (9.4.8) it follows that $A$-column $b$ is a $_{\ast}{}^{\ast}D$-linear combination of $A$-columns $a_T$

$$a_{T\ast}{}^{\ast}(d^T + r_{\ast}{}^{\ast}d^{N\setminus T}) = b$$

This holds equation (9.4.5).

It remains to prove that an existence of solution of system of right linear equations (9.4.3) follows from (9.4.5). Holding (9.4.5) means that $a_T^S$ is $_{\ast}{}^{\ast}$-major submatrix of augmented matrix as well. From theorem 9.3.6 it follows that column $b$ is a right linear composition of columns $a_T$

$$b = a_{T\ast}{}^{\ast}R^T$$

Assigning $R^r = 0$ we get

$$b = a_{\ast}{}^{\ast}R$$

Therefore, we found at least one solution of system of right linear equations (9.4.3).

□

**Theorem 9.4.4.** *Suppose* (9.4.3) *is a system of right linear equations satisfying* (9.4.5). *If* $\mathrm{rank}(_{\ast}{}^{\ast})a = k \le m,$ *then solution of the system depends on arbitrary values of* $m - k$ *variables not included into* $_{\ast}{}^{\ast}$-*major submatrix.*

Proof. Let $a_T^S$ be $_*{}^*$-major submatrix of matrix $a$. Suppose

(9.4.9)                                    $a^p{}_*{}^* x = b^p$

is an equation with number $p$. Applying the theorem 9.3.8 to augmented matrix (9.4.2) we get

(9.4.10)                                   $a^p = R^p{}_*{}^* a^S$

(9.4.11)                                   $b^p = R^p{}_*{}^* b^S$

Substituting (9.4.10) and (9.4.11) into (9.4.9) we get

(9.4.12)                                   $R^p{}_*{}^* a^S{}_*{}^* x = R^p{}_*{}^* b^S$

(9.4.12) means that we can exclude equation (9.4.9) from system (9.4.3) and the new system is equivalent to the old one. Therefore, a number of equations can be reduced to $k$.

At this point, we have two choices. If the number of variables is also $k$ then according to theorem 9.2.14 the system has unique solution

$$\boxed{(9.2.28) \quad x = \mathcal{H}\det({}_*{}_*)a^*{}_* b}$$

If the number of variables $m > k$ then we can move $m - k$ variables that are not included into $_*{}^*$-major submatrix in right side. Giving arbitrary values to these variables, we determine value of the right side and for this value we get a unique solution according to the theorem 9.2.14.                            $\square$

Corollary 9.4.5. *System of right linear equations (9.4.3) has a unique solution iff its matrix is nonsingular.*                            $\square$

Theorem 9.4.6. *Solutions of a homogenous system of right linear equations*

(9.4.13)                                   $a_*{}^* x = 0$

*form a right $A$-vector space of columns.*

Proof. Let $X$ be set of solutions of system of right $A$-linear equations (9.4.13). Suppose $x = (x^a) \in X$ and $y = (y^a) \in X$ Then

$$a_j^i x^j = 0$$
$$a_j^i y^j = 0$$

Therefore

$$[a_j^i(x^j + y^j) = [a_j^i x^j + [a_j^i y^j = 0$$
$$x + y = (x^j + y^j) \in X$$

The same way we see

$$a_j^i(x^j b) = (a_j^i x^j)b = 0$$
$$xb = (x^j b) \in X$$

According to the definition           7.4.1, $X$ is a right $A$-vector space of columns.   $\square$

## 9.5. Nonsingular Matrix

Suppose $A$ is $n \times n$ matrix. Theorems 9.3.11 and 9.3.9 tell us that if $\text{rank}(_*{}^*)a < n$, then rows are linearly dependent from left and columns are linearly dependent from right.[9.7]

THEOREM 9.5.1. *Let a be $n \times n$ matrix and column $a_r$ be a left linear combination of other columns. Then $\text{rank}(_*{}^*)a < n$.*

PROOF. The statement that column $a_r$ is a right linear combination of other columns means that system of inear equations

$$a_r = a_{[r]*}{}^* \lambda$$

has at least one solution. According theorem 9.4.3

$$\text{rank}(_*{}^*)a = \text{rank}(_*{}^*)a_{[r]}$$

Since a number of columns is less then $n$ we get $\text{rank}(_*{}^*)a_{[r]} < n$.  □

THEOREM 9.5.2. *Let a be $n \times n$ matrix and row $a^p$ be a right linear combination of other rows. Then $\text{rank}(_*{}^*)a < n$.*

PROOF. Proof of statement is similar to proof of the theorem 9.5.1.  □

THEOREM 9.5.3. *Suppose a and b are $n \times n$ matrices and*

$$(9.5.1) \qquad c = a_*{}^* b$$

*$c$ is $_*{}^*$-singular matrix iff either matrix a or matrix b is $_*{}^*$-singular matrix.*

PROOF. Suppose matrix $b$ is $_*{}^*$-singular. According to the theorem 9.3.6, columns of matrix $b$ are right linearly dependent. Therefore

$$(9.5.2) \qquad 0 = b_*{}^* \lambda$$

where $\lambda \neq 0$. From (9.5.1) and (9.5.2) it follows that

$$c_*{}^* \lambda = a_*{}^* b_*{}^* \lambda = 0$$

According theorem 9.5.1 matrix $c$ is $_*{}^*$-singular.

Suppose matrix $b$ is not $_*{}^*$-singular, but matrix $a$ is $_*{}^*$-singular According to the theorem 9.3.8, columns of matrix $a$ are right linearly dependent. Therefore

$$(9.5.3) \qquad 0 = a_*{}^* \mu$$

where $\mu \neq 0$. According to the theorem 9.2.14, the system of linear equations

$$b_*{}^* \lambda = \mu$$

has only solution where $\lambda \neq = 0$. Therefore

$$c_*{}^* \lambda = a_*{}^* b_*{}^* \lambda = a_*{}^* \mu = 0$$

According to theorem 9.5.1 matrix $c$ is $_*{}^*$-singular.

Suppose matrix $c$ is $_*{}^*$-singular matrix. According to the theorem 9.3.8 columnss of matrix $c$ are right linearly dependent. Therefore

$$(9.5.4) \qquad 0 = c_*{}^* \lambda$$

---

[9.7]This statement is similar to proposition [6]-1.2.5.

where $\lambda \neq 0$. From (9.5.1) and (9.5.4) it follows that

$$0 = a_* {}^* b_* {}^* \lambda$$

If

$$0 = b_* {}^* \lambda$$

satisfied then matrix $b$ is $_*{}^*$-singular. Suppose that matrix $b$ is not $_*{}^*$-singular. Let us introduce

$$\mu = b_* {}^* \lambda$$

where $\mu \neq 0$. Then

(9.5.5)                         $$0 = a_* {}^* \mu$$

From (9.5.5) it follows that matrix $a$ is $_*{}^*$-singular.                    □

According to the theorem 8.3.11, we can write similar statements for left linear combination of rows or right linear combination of columns and $^*{}_*$-quasideterminant.

THEOREM 9.5.4. *Let $a$ be $n \times n$ matrix and column $a_r$ be a right linear combination of other columns. Then* $\mathrm{rank}(^*{}_*)a < n$.                    □

THEOREM 9.5.5. *Let $a$ be $n \times n$ matrix and row $a^p$ be a left linear combination of other rows. Then* $\mathrm{rank}(^*{}_*)a < n$.                    □

THEOREM 9.5.6. *Suppose $a$ and $b$ are $n \times n$ matrices and*

(9.5.6)                         $$c = a^* {}_* b$$

*$c$ is $^*{}_*$-singular matrix iff either matrix $a$ or matrix $b$ is $^*{}_*$-singular matrix.*                    □

DEFINITION 9.5.7. $_*{}^*$-**matrix group** $GL(n_*{}^*A)$ *is a group of $_*{}^*$-nonsingular matrices where we define $_*{}^*$-product of matrices (8.1.7) and $_*{}^*$-inverse matrix $a^{-1_*{}^*}$.*                    □

DEFINITION 9.5.8. $^*{}_*$-**matrix group** $GL(n^*{}_*A)$ *is a group of $^*{}_*$-nonsingular matrices where we define $^*{}_*$-product of matrices (8.1.10) and $^*{}_*$-inverse matrix $a^{-1^*{}_*}$.*                    □

THEOREM 9.5.9.

$$GL(n_*{}^*A) \neq GL(n^*{}_*A)$$

REMARK 9.5.10. *From theorems 8.2.9, 8.3.9, it follows that there are matrices which are $_*{}^*$-nonsingular and $^*{}_*$-nonsingular. Theorem 9.5.9 implies that sets of $_*{}^*$-nonsingular matrices and $^*{}_*$-nonsingular matrices are not identical. For instance, there exists such $_*{}^*$-nonsingular matrix which is $^*{}_*$-singular matrix.*                    □

PROOF. It is enough to prove this statement for $n = 2$. Assume every $^*{}_*$-singular matrix

(9.5.7)                         $$a = \begin{pmatrix} a_1^1 & a_2^1 \\ a_1^2 & a_2^2 \end{pmatrix}$$

is $_*{}^*$-singular matrix. It follows from theorems 9.3.8, 9.3.6, that $_*{}^*$-singular matrix satisfies to condition

$$(9.5.8) \qquad\qquad a_2^1 = ba_1^1$$

$$(9.5.9) \qquad\qquad a_2^2 = ba_1^2$$

$$(9.5.10) \qquad\qquad a_1^2 = a_1^1 c$$

$$(9.5.11) \qquad\qquad a_2^2 = a_2^1 c$$

If we substitute (9.5.10) into (9.5.9) we get

$$a_2^2 = b\ a_1^1\ c$$

$b$ and $c$ are arbitrary $A$-numbers and $_*{}^*$-singular matrix matrix (9.5.7) has form $(d = a_1^1)$

$$(9.5.12) \qquad\qquad a = \begin{pmatrix} d & dc \\ bd & bdc \end{pmatrix}$$

The similar way we can show that $^*{}_*$-singular matrix has form

$$(9.5.13) \qquad\qquad a = \begin{pmatrix} d & c'd \\ db' & c'db' \end{pmatrix}$$

From assumption it follows that (9.5.13) and (9.5.12) represent the same matrix. Comparing (9.5.13) and (9.5.12) we get that for every $d, c \in D$ exists such $c' \in D$ which does not depend on $d$ and satisfies equation

$$dc = c'd$$

This contradicts the fact that $A$ is division algebra. $\qquad\qquad\square$

EXAMPLE 9.5.11. *Since we get quaternion algebra, we assume $b = 1 + k$, $c = j$, $d = k$. Then we get*

$$a = \begin{pmatrix} k & kj \\ (1+k)k & (1+k)kj \end{pmatrix} = \begin{pmatrix} k & -i \\ k-1 & -i-j \end{pmatrix}$$

$$\det(_*{}^*)_2^2 a = a_2^2 - a_2^1(a_1^1)^{-1}a_1^2$$

$$= -i - j - (k-1)(k)^{-1}(-i) = -i - j - (k-1)(-k)(-i)$$

$$= -i - j - kki + ki = -i - j + i + j = 0$$

$$\det(_*{}^*)_1^1 a = a_1^1 - a_1^2(a_2^2)^{-1}a_2^1$$

$$= k - (k-1)(-i-j)^{-1}(-i) = k - (k-1)\frac{1}{2}(i+j)(-i)$$

$$= k + \frac{1}{2}((k-1)i + (k-1)j)i = k + \frac{1}{2}(ki - i + kj - j)i$$

$$= k + \frac{1}{2}(j - i - i - j)i = k - ii = k + 1$$

$$\det({}_*^*)^2_1 a = a^2_1 - a^1_1(a^1_2)^{-1}a^2_2$$

$$= k - 1 - k(-i)^{-1}(-i-j) = k - 1 + ki(i+j)$$

$$= k - 1 + j(i+j) = k - 1 + ji + jj$$

$$= k - 1 - k - 1 = -2$$

$$\det({}_*^*)^1_2 a = a^1_2 - a^2_2(a^2_1)^{-1}a^1_1$$

$$= (-i) - (-i-j)(k-1)^{-1}k = -i + (i+j)\frac{1}{2}(-k-1)k$$

$$= -i - \frac{1}{2}(i+j)(k+1)k = -i - \frac{1}{2}(ik+i+jk+j)k$$

$$= -i - \frac{1}{2}(-j+i+i+j)k = -i - ik = -i + j$$

$$\det({}^*_*)^2_2 a = a^2_2 - a^1_2(a^1_1)^{-1}a^2_1$$

$$= -i - j - (-i)(k)^{-1}(k-1) = -i - j + i(-k)(k-1)$$

$$= -i - j + j(k-1) = -i - j + jk - j = -i - j + i - j = -2j$$

*The system of linear equations*

(9.5.14)
$$\begin{pmatrix} k & -i \\ k-1 & -i-j \end{pmatrix} {}^*_* \begin{pmatrix} x^1 \\ x^2 \end{pmatrix} = \begin{pmatrix} b^1 \\ b^2 \end{pmatrix}$$

*has $_*^*$-singular matrix. We can write the system of linear equations* (9.5.14) *in the form*

$$\begin{cases} k\ x^1 - & i\ x^2 = b^1 \\ (k-1)\ x^1 - (i+j)\ x^2 = b^2 \end{cases}$$

*The system of linear equations*

(9.5.15)
$$\begin{pmatrix} k & -i \\ k-1 & -i-j \end{pmatrix} {}^*_* \begin{pmatrix} x_1 & x_2 \end{pmatrix} = \begin{pmatrix} b_1 & b_2 \end{pmatrix}$$

*has ${}^*_*$-nonsingular matrix. We can write the system of linear equations* (9.5.15) *in the form*

$$\begin{cases} k\ x_1 - & i\ x_1 \\ + (k-1)\ x_2 & - (i+j)\ x_2 \\ = \quad b_1 & = \quad b_2 \end{cases} \qquad \begin{cases} k\ x_1 + (k-1)\ x_2 = b_1 \\ -i\ x_1 - (i+j)\ x_2 = b_2 \end{cases}$$

*The system of linear equations*

(9.5.16)
$$\begin{pmatrix} x_1 & x_2 \end{pmatrix} {}^*_* \begin{pmatrix} k & -i \\ k-1 & -i-j \end{pmatrix} = \begin{pmatrix} b_1 & b_2 \end{pmatrix}$$

has $_*^*$-singular matrix. We can write the system of linear equations (9.5.16) in the form

$$\begin{cases} x_1 k & -x_1 i \\ + x_2(k-1) & -x_2(i+j) \\ = b_1 & = b_2 \end{cases} \qquad \begin{cases} x_1 k & +x_2(k-1) = b_1 \\ -x_1 i & -x_2(i+j) = b_2 \end{cases}$$

*The system of linear equations*

$$(9.5.17) \qquad \begin{pmatrix} x^1 \\ x^2 \end{pmatrix} {}_*^* \begin{pmatrix} k & -i \\ k-1 & -i-j \end{pmatrix} = \begin{pmatrix} b^1 \\ b^2 \end{pmatrix}$$

has $^*_*$-nonsingular matrix. We can write the system of linear equations (9.5.17) in the form

$$\begin{cases} x^1 \; k & - \; x^2 \; i & = b^1 \\ x^1 \; (k-1) - x^2 \; (i+j) & = b^2 \end{cases}$$

<div style="text-align: right;">□</div>

## 9.6. Inverse Matrix

**Definition 9.6.1.** *Let a be $n \times m$ matrix*

$$a = \begin{pmatrix} a_1^1 & ... & a_n^1 \\ ... & ... & ... \\ a_1^m & ... & a_n^m \end{pmatrix}$$

*and*

$$k = \min(n, m)$$

*The matrix a is called left $^*_*$-regular if there exists $m \times n$ matrix b*

$$b = \begin{pmatrix} b_1^1 & ... & b_m^1 \\ ... & ... & ... \\ b_1^n & ... & b_m^n \end{pmatrix}$$

*such that $^*_*$-product*

$$(9.6.1) \qquad c = b^*_* a$$

*has submatrix*

$$c_J^I = E_k$$

*and other elements of the matrix c are 0. The matrix b is called left $^*_*$-inverse for the matrix a.* □

**Definition 9.6.2.** *Let a be $n \times m$ matrix*

$$a = \begin{pmatrix} a_1^1 & ... & a_n^1 \\ ... & ... & ... \\ a_1^m & ... & a_n^m \end{pmatrix}$$

*and*

$$k = \min(n, m)$$

*The matrix a is called right $^*_*$-regular if there exists $m \times n$ matrix b*

$$b = \begin{pmatrix} b_1^1 & ... & b_m^1 \\ ... & ... & ... \\ b_1^n & ... & b_m^n \end{pmatrix}$$

*such that $_*^*$-product*

$$(9.6.2) \qquad c = a_*^* b$$

*has submatrix*

$$c_J^I = E_k$$

*and other elements of the matrix c are 0. The matrix b is called right $_*^*$-inverse for the matrix a.* □

**Definition 9.6.3.** *Let $a$ be $n \times m$ matrix*

$$a = \begin{pmatrix} a_1^1 & \dots & a_n^1 \\ \dots & \dots & \dots \\ a_1^m & \dots & a_n^m \end{pmatrix}$$

*and*

$$k = \min(n, m)$$

*The matrix $a$ is called left $_*^*$-regular if there exists $m \times n$ matrix $b$*

$$b = \begin{pmatrix} b_1^1 & \dots & b_m^1 \\ \dots & \dots & \dots \\ b_1^n & \dots & b_m^n \end{pmatrix}$$

*such that $_*^*$-product*

(9.6.3) $\qquad c = b_* {}^* a$

*has submatrix*

$$c_J^I = E_k$$

*and other elements of the matrix $c$ are 0. The matrix $b$ is called left $_*^*$-inverse for the matrix $a$.* □

**Definition 9.6.4.** *Let $a$ be $n \times m$ matrix*

$$a = \begin{pmatrix} a_1^1 & \dots & a_n^1 \\ \dots & \dots & \dots \\ a_1^m & \dots & a_n^m \end{pmatrix}$$

*and*

$$k = \min(n, m)$$

*The matrix $a$ is called right $_*^*$-regular if there exists $m \times n$ matrix $b$*

$$b = \begin{pmatrix} b_1^1 & \dots & b_m^1 \\ \dots & \dots & \dots \\ b_1^n & \dots & b_m^n \end{pmatrix}$$

*such that $^*_*$-product*

(9.6.4) $\qquad c = a^* {}_* b$

*has submatrix*

$$c_J^I = E_k$$

*and other elements of the matrix $c$ are 0. The matrix $b$ is called right $_*^*$-inverse for the matrix $a$.* □

**Theorem 9.6.5.** *Let*

$$k = \min(n, m)$$

*The matrix $a$ has left $^*_*$-inverse matrix $b$ iff*

(9.6.5) $\qquad \operatorname{rank}(^*_*)a = k$

**Proof.** Let

$$J = [1, k]$$

Let $I$ be set of indices of left linear independent rows of matrix $a$. If $|I| < k$, then we will complement the set $I$ with indices of the remaining rows of matrix $a$ in such a way as to get $|I| = k$. According to the definition 8.1.4, the equality

(9.6.6) $\qquad c_J^I = b_J {}^*_* a^I$

follows from the equality

$\boxed{(9.6.1) \quad c = b^*_* a}$

According to definitions 8.1.4, 9.6.1,

**Theorem 9.6.6.** *Let*

$$k = \min(n, m)$$

*The matrix $a$ has right $^*_*$-inverse matrix $b$ iff*

(9.6.8) $\qquad \operatorname{rank}(_*^*)a = k$

**Proof.** Let

$$J = [1, k]$$

Let $I$ be set of indices of right linear independent rows of matrix $a$. If $|I| < k$, then we will complement the set $I$ with indices of the remaining rows of matrix $a$ in such a way as to get $|I| = k$. According to the definition 8.1.3, the equality

(9.6.9) $\qquad c_J^I = a^I {}^*_* b_J$

follows from the equality

$\boxed{(9.6.2) \quad c = a_* {}^* b}$

According to definitions 8.1.3, 9.6.2,

For each index $j \in J$, the system of linear equations

$$b_j^1 a_1^{i_1} + ... + b_j^n a_n^{i_1} = \delta_j^1$$

(9.6.7)                ......

$$b_j^1 a_1^{i_k} + ... + b_j^n a_n^{i_k} = \delta_j^k$$

follows from the equality (9.6.6). We assume that unknown variables are $x^i = b_j^i$.

If we consider all systems of linear equations (9.6.7), then their augmented matrix will have the form

$$\begin{pmatrix} a_1^{i_1} & ... & a_n^{i_1} & \delta_1^1 & ... & \delta_k^1 \\ ... & ... & ... & ... & ... & ... \\ a_1^{i_k} & ... & a_n^{i_k} & \delta_1^k & ... & \delta_k^k \end{pmatrix}$$

Rank of augmented matrix is $k$. According to the theorem 9.4.3, system of linear equations (9.6.7) have solution iff condition (9.6.5). is met. In particular, the set of rows $a^i$, $i \in I$, is left linear independent.          □

---

**THEOREM 9.6.7.** *Let*

$$k = \min(n, m)$$

*The matrix $a$ has left $_*{}^*$-inverse matrix $b$ iff*

(9.6.11)          $\operatorname{rank}(_*{}^*)a = k$

PROOF. Let

$$I = [1, k]$$

Let $J$ be set of indices of left linear independent columns of matrix $a$. If $|J| < k$, then we will complement the set $J$ with indices of the remaining columns of matrix $a$ in such a way as to get $|J| = k$. According to the definition 8.1.3, the equality

(9.6.12)          $c_J^I = b^I {}_*{}^* a_J$

follows from the equality

(9.6.3)   $c = b_* a$

According to definitions 8.1.3, 9.6.3,

---

For each index $j \in J$, the system of linear equations

$$a_1^{i_1} b_j + ... + a_n^{i_1} b_j^n = \delta_j^1$$

(9.6.10)                ......

$$a_1^{i_k} b_j + ... + a_n^{i_k} b_j^n = \delta_j^k$$

follows from the equality (9.6.9). We assume that unknown variables are $x^i = b_j^i$.

If we consider all systems of linear equations (9.6.10), then their augmented matrix will have the form

$$\begin{pmatrix} a_1^{i_1} & ... & a_n^{i_1} & \delta_1^1 & ... & \delta_k^1 \\ ... & ... & ... & ... & ... & ... \\ a_1^{i_k} & ... & a_n^{i_k} & \delta_1^k & ... & \delta_k^k \end{pmatrix}$$

Rank of augmented matrix is $k$. According to the theorem 9.4.3, system of linear equations (9.6.10) have solution iff condition (9.6.8). is met. In particular, the set of rows $a^i$, $i \in I$, is right linear independent.          □

---

**THEOREM 9.6.8.** *Let*

$$k = \min(n, m)$$

*The matrix $a$ has right $^*{}_*$-inverse matrix $b$ iff*

(9.6.14)          $\operatorname{rank}(^*{}_*)a = k$

PROOF. Let

$$I = [1, k]$$

Let $J$ be set of indices of right linear independent columns of matrix $a$. If $|J| < k$, then we will complement the set $J$ with indices of the remaining columns of matrix $a$ in such a way as to get $|J| = k$. According to the definition 8.1.4, the equality

(9.6.15)          $c_J^I = a_J {}^*{}_* b^I$

follows from the equality

(9.6.4)   $c = a^*{}_* b$

According to definitions 8.1.4, 9.6.4,

For each index $i \in I$, the system of linear equations

$$(9.6.13)\quad \begin{aligned} b_1^i a_{j_1}^1 + ... + b_m^i a_{j_1}^m &= \delta_1^i \\ &...... \\ b_1^i a_{j_k}^1 + ... + b_m^i a_{j_k}^m &= \delta_k^i \end{aligned}$$

follows from the equality (9.6.12). We assume that unknown variables are $x_j = b_j^i$.

If we consider all systems of linear equations (9.6.13), then their augmented matrix will have the form

$$\begin{pmatrix} a_{j_1}^1 & ... & a_{j_k}^1 \\ ... & ... & ... \\ a_{j_1}^m & ... & a_{j_k}^m \\ \delta_1^1 & ... & \delta_k^1 \\ ... & ... & ... \\ \delta_1^k & ... & \delta_k^k \end{pmatrix}$$

Rank of augmented matrix is $k$. According to the theorem 9.4.3, system of linear equations (9.6.13) have solution iff condition (9.6.11). is met. In particular, the set of columns $a_j$, $j \in J$, is left linear independent. □

For each index $i \in I$, the system of linear equations

$$(9.6.16)\quad \begin{aligned} a_{j_1}^1 b_1^i + ... + a_{j_1}^m b_m^i &= \delta_1^i \\ &...... \\ a_{j_k}^1 b_1^i + ... + a_{j_k}^m b_m^i &= \delta_k^i \end{aligned}$$

follows from the equality (9.6.15). We assume that unknown variables are $x_j = b_j^i$.

If we consider all systems of linear equations (9.6.16), then their augmented matrix will have the form

$$\begin{pmatrix} a_{j_1}^1 & ... & a_{j_k}^1 \\ ... & ... & ... \\ a_{j_1}^m & ... & a_{j_k}^m \\ \delta_1^1 & ... & \delta_k^1 \\ ... & ... & ... \\ \delta_1^k & ... & \delta_k^k \end{pmatrix}$$

Rank of augmented matrix is $k$. According to the theorem 9.4.3, system of linear equations (9.6.16) have solution iff condition (9.6.14). is met. In particular, the set of columns $a_j$, $j \in J$, is right linear independent. □

If $m = n$, then left $^*_*$-inverse matrix is unique. In general, this is not true when $m \neq n$. Consider the system of linear equations (9.6.7).

- If $m < n$, then the system of linear equations (9.6.7) contains more unknowns than the number of equations. A solution of the system of linear equations (9.6.7) is unique up to choice of values of redundant variables.
- If $m > n$, then the system of linear equations (9.6.7) contains $k$ unknowns and $k$ equations. Therefore, the system of linear equations (9.6.7) has a unique solution. Howeve, a choice of the set $I$ is not unique.

DEFINITION 9.6.9. *The set* $left(a^{-1^*_*})$ *is the set of left* $^*_*$-*inverse matrices for the matrix* $a$. □

DEFINITION 9.6.10. *The set* $right(a^{-1^*_*})$ *is the set of right* $^*_*$-*inverse matrices for the matrix* $a$. □

DEFINITION 9.6.11. *The set* $left(a^{-1^*_*})$ *is the set of left* $^*_*$-*inverse matrices for the matrix* $a$. □

DEFINITION 9.6.12. *The set* $right(a^{-1^*_*})$ *is the set of right* $^*_*$-*inverse matrices for the matrix* $a$. □

THEOREM 9.6.13. $a \in left(b^{-1^*_*})$

THEOREM 9.6.14. $a \in right(b^{-1^*_*})$

*iff* $b \in right(a^{-1^{*}*})$.

PROOF. The theorem follows from definitions 9.6.1, 9.6.4, 9.6.9, 9.6.12. □

*iff* $b \in left(a^{-1}{}_{*}{}^{*})$.

PROOF. The theorem follows from definitions 9.6.2, 9.6.3, 9.6.10, 9.6.11. □

# Homomorphism of Left $A$-Vector Space

## 10.1. General Definition

Let $A_i$, $i = 1$, $2$, be division algebra over commutative ring $D_i$. Let $V_i$, $i = 1$, $2$, be left $A_i$-vector space.

DEFINITION 10.1.1. *Let diagram of representations*

$$g_{1.12}(d) : a \to d\,a$$
$$g_{1.23}(v) : w \to C_1(w, v)$$

(10.1.1)

$$A_1 \xrightarrow{g_{1.23}} A_1 \xrightarrow{g_{1.34}} V_1$$

$$C_1 \in \mathcal{L}(A_1^2 \to A_1)$$
$$g_{1.34}(a) : v \to a\,v$$
$$g_{1.14}(d) : v \to d\,v$$

*describe left $A_1$-vector space $V_1$. Let diagram of representations*

$$g_{2.12}(d) : a \to d\,a$$
$$g_{2.23}(v) : w \to C_2(w, v)$$

(10.1.2)

$$A_2 \xrightarrow{g_{2.23}} A_2 \xrightarrow{g_{2.34}} V_2$$

$$C_2 \in \mathcal{L}(A_2^2 \to A_2)$$
$$g_{2.34}(a) : v \to a\,v$$
$$g_{2.14}(d) : v \to d\,v$$

*describe left $A_2$-vector space $V_2$. Morphism*

(10.1.3) $$(h : D_1 \to D_2 \quad g : A_1 \to A_2 \quad f : V_1 \to V_2)$$

*of diagram of representations* (10.1.1) *into diagram of representations* (10.1.2) *is called* **homomorphism** *of left $A_1$-vector space $V_1$ into left $A_2$-vector space $V_2$. Let us denote* $\mathrm{Hom}(D_1 \to D_2; A_1 \to A_2; *V_1 \to *V_2)$ *set of homomorphisms of left $A_1$-vector space $V_1$ into left $A_2$-vector space $V_2$.* □

We will use notation

$$f \circ a = f(a)$$

for image of homomorphism $f$.

THEOREM 10.1.2. *The homomorphism*

$$(h : D_1 \to D_2 \quad g : A_1 \to A_2 \quad f : V_1 \to V_2)$$
(10.1.3)

*of left $A_1$-vector space $V_1$ into left $A_2$-vector space $V_2$ satisfies following equalities*

(10.1.4) $$h(p + q) = h(p) + h(q)$$
(10.1.5) $$h(pq) = h(p)h(q)$$
(10.1.6) $$g \circ (a + b) = g \circ a + g \circ b$$
(10.1.7) $$g \circ (ab) = (g \circ a)(g \circ b)$$

(10.1.8)                                   $g \circ (pa) = h(p)(g \circ a)$

(10.1.9)                                   $f \circ (u + v) = f \circ u + f \circ v$

(10.1.10)                                  $f \circ (av) = (g \circ a)(f \circ v)$

$$p, q \in D_1 \quad a, b \in A_1 \quad u, v \in V_1$$

PROOF.    According to definitions    2.3.9, 2.8.5, 7.2.1, 11 the map

$$(h : D_1 \to D_2 \quad g : A_1 \to A_2)$$

is homomorphism of $D_1$-algebra $A_1$ into $D_2$-algebra $A_2$. Therefore, equalities
(10.1.4), (10.1.5), (10.1.6), (10.1.7), (10.1.8)    follow from the theorem 4.3.2.
The equality (10.1.9) follows from the definition 10.1.1, since, accorfing to the def-
inition 2.3.9, the map $f$ is homomorpism of Abelian group. The equality (10.1.10)
follows from the equality (2.3.6) because the map

$$(g : A_1 \to A_2 \quad f : V_1 \to V_2)$$

is morphism of representation $g_{1.34}$ into representation $g_{2.34}$.                □

DEFINITION 10.1.3. *Homomorphism* [10.1]

$$(h : D_1 \to D_2 \quad g : A_1 \to A_2 \quad f : V_1 \to V_2)$$

*is called* **isomorphism** *between left $A_1$-vector space $V_1$ and left $A_2$-vector space $V_2$,
if there exists the map*

$$(h^{-1} : D_2 \to D_1 \quad g^{-1} : A_2 \to A_1 \quad f^{-1} : V_2 \to V_1)$$

*which is homomorphism.*                                                          □

### 10.1.1. Homomorphism of Left $A$-Vector Space of columns.

THEOREM 10.1.4. *The homomorphism* [10.2]

| (10.1.3)    $(h : D_1 \to D_2 \quad \overline{g} : A_1 \to A_2 \quad \overline{f} : V_1 \to V_2)$ |
|---|

*of left $A_1$-vector space of columns $V_1$ into left $A_2$-vector space of columns $V_2$ has
presentation*

(10.1.11)                                  $b = gh(a)$

(10.1.12)                                  $\overline{g} \circ (a^i e_{A_1 i}) = h(a^i) g_i^k e_{A_2 k}$

(10.1.13)                                  $\overline{g} \circ (e_{A_1} a) = e_{A_2} gh(a)$

(10.1.14)                                  $w = (\overline{g} \circ v)^* {}_* f$

(10.1.15)                                  $\overline{f} \circ (v^i e_{V_1 i}) = (\overline{g} \circ v^i) f_i^k e_{V_2 k}$

---

[10.1] I follow the definition on page [18]-49.

[10.2] In theorems 10.1.4, 10.1.5, we use the following convention.    Let $i = 1, 2$.    Let the
set of vectors $\overline{e}_{A_i} = (e_{A_i k}, k \in K_i)$ be a basis and $C_{i j}^{k}$, $k, i, j \in K_i$, be structural constants
of $D_i$-algebra of columns $A_i$.    Let the set of vectors $\overline{e}_{V_i} = (e_{V_i i}, i \in I_i)$ be a basis of left $A_i$-
vector space $V_i$.

(10.1.16) $$\overline{f} \circ (v^* {}_* e_{V_1}) = (\overline{g} \circ v)^* {}_* (f^* {}_* e_{V_2})$$

*relative to selected bases. Here*

- *$a$ is coordinate matrix of $A_1$-number $\overline{a}$ relative the basis $\overline{\overline{e}}_{A_1}$*

$$\overline{a} = e_{A_1} a$$

- *$h(a) = (h(a_k), k \in K)$ is a matrix of $D_2$-numbers.*
- *$b$ is coordinate matrix of $A_2$-number*

$$\overline{b} = \overline{g} \circ \overline{a}$$

  *relative the basis $\overline{\overline{e}}_{A_2}$*

$$\overline{b} = e_{A_2} b$$

- *$g$ is coordinate matrix of set of $A_2$-numbers $(\overline{g} \circ e_{A_1 k}, k \in K)$ relative the basis $\overline{\overline{e}}_{A_2}$.*
- *There is relation between the matrix of homomorphism and structural constants*

(10.1.17) $$h(C_1{}_{ij}^{\ k}) g_k^l = g_i^p g_j^q C_2{}_{pq}^{\ l}$$

- *$v$ is coordinate matrix of $V_1$-number $\overline{v}$ relative the basis $\overline{\overline{e}}_{V_1}$*

$$\overline{v} = v^* {}_* e_{V_1}$$

- *$g(v) = (g(v_i), i \in I)$ is a matrix of $A_2$-numbers.*
- *$w$ is coordinate matrix of $V_2$-number*

$$\overline{w} = \overline{f} \circ \overline{v}$$

  *relative the basis $\overline{\overline{e}}_{V_2}$*

$$\overline{w} = w^* {}_* e_{V_2}$$

- *$f$ is coordinate matrix of set of $V_2$-numbers $(\overline{f} \circ e_{V_1 i}, i \in I)$ relative the basis $\overline{\overline{e}}_{V_2}$.*

*For given homomorphism*

> (10.1.3)   $(h : D_1 \to D_2$   $\overline{g} : A_1 \to A_2$   $\overline{f} : V_1 \to V_2)$

*the set of matrices $(g, f)$ is unique and is called **coordinates of homomorphism**
(10.1.3)   relative bases $(\overline{\overline{e}}_{A_1}, \overline{\overline{e}}_{V_1})$,   $(\overline{\overline{e}}_{A_2}, \overline{\overline{e}}_{V_2})$.*

PROOF.    According to definitions   2.3.9, 2.8.5,  5.4.1, 10.1.1,  the map

$$(h : D_1 \to D_2 \quad \overline{g} : A_1 \to A_2)$$

is homomorphism of $D_1$-algebra $A_1$ into $D_2$-algebra $A_2$. Therefore, equalities
(10.1.11), (10.1.12), (10.1.13), (10.1.17)  follow from equalities

> (5.4.11)   $b = f h(a)$

> (5.4.12)   $\overline{f} \circ (a^i e_{1i}) = h(a^i) f_i^k e_{2k}$

> (5.4.13)   $\overline{f} \circ (e_1 a) = e_2 f h(a)$

> (5.4.16)   $h(C_1{}_{ij}^{\ k}) f_k^l = f_i^p f_j^q C_2{}_{pq}^{\ l}$

From the theorem 4.3.3, it follows that the matrix $g$ is unique.

Vector  $\overline{v} \in V_1$  has expansion

(10.1.18) $$\overline{v} = v^* {}_* e_{V_1}$$

relative to the basis $\overline{\overline{e}}_{V_1}$.   Vector $\overline{w} \in V_2$ has expansion

(10.1.19) $$\overline{w} = w^*{}_* e_{V_2}$$

relative to the basis $\overline{\overline{e}}_{V_2}$.   Since $\overline{f}$ is a homomorphism, then the equality

(10.1.20) $$\overline{w} = \overline{f} \circ \overline{v} = \overline{f} \circ (v^*{}_* e_{V_1})$$
$$= (\overline{g} \circ v)^*{}_* (\overline{f} \circ e_{V_1})$$

follows from equalities (10.1.9), (10.1.10). $V_2$-number $\overline{f} \circ e_{V_1 i}$ has expansion

(10.1.21) $$\overline{f} \circ e_{V_1 i} = f_i{}^*{}_* e_{V_2} = f_i^j e_{V_2 j}$$

relative to basis $\overline{\overline{e}}_{V_2}$.   The equality

(10.1.22) $$\overline{w} = (\overline{g} \circ v)^*{}_* f^*{}_* e_{V_2}$$

follows from equalities (10.1.20), (10.1.21).   The equality (10.1.14) follows from comparison of (10.1.19) and (10.1.22) and the theorem 8.4.4.   From the equality (10.1.21) and from the theorem 8.4.4, it follows that the matrix $f$ is unique.   $\square$

**Theorem 10.1.5.** *Let the map*

$$h : D_1 \to D_2$$

*be homomorphism of the ring $D_1$ into the ring $D_2$.*        *Let*

$$g = (g_l^k, k \in K, l \in L)$$

*be matrix of $D_2$-numbers which satisfy the equality*

(10.1.17)   $h(C_{1\,ij}^{\ \ k}) g_k^l = g_i^p g_j^q C_{2\,pq}^{\ \ l}$

*Let*

$$f = (f_j^i, i \in I, j \in J)$$

*be matrix of $A_2$-numbers.*        *The map*[10.2]

(10.1.3)   $(h : D_1 \to D_2 \quad \overline{g} : A_1 \to A_2 \quad \overline{f} : V_1 \to V_2)$

*defined by equalities*

(10.1.13)   $\overline{g} \circ (e_{A_1}, a) = e_{A_2} g h(a)$

(10.1.16)   $\overline{f} \circ (v^*{}_* e_{V_1}) = (\overline{g} \circ v)^*{}_* (f^*{}_* e_{V_2})$

*is homomorphism of left $A_1$-vector space of columns $V_1$ into left $A_2$-vector space of columns $V_2$.*        *The homomorphism*

(10.1.3)   $(h : D_1 \to D_2 \quad \overline{g} : A_1 \to A_2 \quad \overline{f} : V_1 \to V_2)$

*which has the given  set of matrices $(g, f)$  is unique.*

**Proof.**   According to the theorem 5.4.5, the map

(10.1.23) $$(h : D_1 \to D_2 \quad \overline{g} : A_1 \to A_2)$$

is homomorphism of $D_1$-algebra $A_1$ into $D_2$-algebra $A_2$ and the homomorphism (10.1.23) is unique.

The equality

$$
\begin{aligned}
&\overline{f} \circ ((v+w)^*{}_* e_{V_1}) \\
&= (g \circ (v+w))^*{}_* f^*{}_* e_{V_2} \\
&= (g \circ v)^*{}_* f^*{}_* e_{V_2} \\
&\quad + (g \circ w)^*{}_* f^*{}_* e_{V_2} \\
&= \overline{f} \circ (v^*{}_* e_{V_1}) + \overline{f} \circ (w^*{}_* e_{V_1})
\end{aligned}
$$

(10.1.24)

follows from equalities

| (8.1.19)   $(b_1 + b_2)^*{}_* a = b_1{}^*{}_* a + b_2{}^*{}_* a$ | (10.1.6)   $g \circ (a+b) = g \circ a + g \circ b$ |
|---|---|
| (10.1.10)   $f \circ (av) = (g \circ a)(f \circ v)$ | (10.1.16)   $\overline{f} \circ (v^*{}_* e_{V_1}) = (\overline{g} \circ v)^*{}_* (f^*{}_* e_{V_2})$ |

From the equality (10.1.24), it follows that the map $\overline{f}$ is homomorphism of Abelian group. The equality

$$
\begin{aligned}
&\overline{f} \circ ((av)^*{}_* e_{V_1}) \\
&= (\overline{g} \circ (av))^*{}_* f^*{}_* e_{V_2} \\
&= (\overline{g} \circ a)((\overline{g} \circ v)^*{}_* f^*{}_* e_{V_2}) \\
&= (\overline{g} \circ a)(\overline{f} \circ (v^*{}_* e_{V_1}))
\end{aligned}
$$

(10.1.25)

follows from equalities

| (10.1.10)   $f \circ (av) = (g \circ a)(f \circ v)$ |
|---|
| (5.4.4)   $g \circ (ab) = (g \circ a)(g \circ b)$ |
| (10.1.16)   $\overline{f} \circ (v^*{}_* e_{V_1}) = (\overline{g} \circ v)^*{}_* (f^*{}_* e_{V_2})$ |

From the equality (10.1.25), and definitions 2.8.5, 10.1.1, it follows that the map

| (10.1.3)   $(h : D_1 \to D_2 \quad \overline{g} : A_1 \to A_2 \quad \overline{f} : V_1 \to V_2)$ |
|---|

is homomorphism of $A_1$-vector space of columns $V_1$ into $A_2$-vector space of columns $V_2$.

Let $f$ be matrix of homomorphisms $\overline{f}$, $\overline{g}$ relative to bases $\overline{\overline{e}}_V$, $\overline{\overline{e}}_W$. The equality

$$
\overline{f} \circ (v^*{}_* e_V) = v^*{}_* f^*{}_* e_W = \overline{g} \circ (v^*{}_* e_V)
$$

follows from the theorem 10.3.6. Therefore, $\overline{f} = \overline{g}$.          □

On the basis of theorems 10.1.4, 10.1.5 we identify the homomorphism

| (10.1.3)   $(h : D_1 \to D_2 \quad \overline{g} : A_1 \to A_2 \quad \overline{f} : V_1 \to V_2)$ |
|---|

of left $A$-vector space $V$ of columns and coordinates of its presentation

| (10.1.13)   $\overline{g} \circ (e_{A_1} a) = e_{A_2} gh(a)$ |
|---|

| (10.1.16)   $\overline{f} \circ (v^*{}_* e_{V_1}) = (\overline{g} \circ v)^*{}_* (f^*{}_* e_{V_2})$ |
|---|

## 10.1.2. Homomorphism of Left $A$-Vector Space of rows.

[10.3] In theorems 10.1.6, 10.1.7, we use the following convention.  Let  $i = 1, 2$.     Let the set of vectors $\overline{\overline{e}}_{A_i} = (e^k_{A_i}, k \in K_i)$  be a basis and  $C^{ij}_{ik}$,   $k, i, j \in K_i$,  be structural constants of $D_i$-algebra of rows $A_i$.    Let the set of vectors $\overline{\overline{e}}_{V_i} = (e^i_{V_i}, i \in I_i)$  be a basis of left $A_i$-vector

THEOREM 10.1.6. *The homomorphism* [10.3]

$$(10.1.3) \quad (h : D_1 \to D_2 \quad \overline{g} : A_1 \to A_2 \quad \overline{f} : V_1 \to V_2)$$

*of left $A_1$-vector space of rows $V_1$ into left $A_2$-vector space of rows $V_2$ has presentation*

$$(10.1.26) \qquad\qquad b = h(a)g$$

$$(10.1.27) \qquad\qquad \overline{g} \circ (a_i e^i_{A_1}) = h(a_i)g^i_k e^k_{A_2}$$

$$(10.1.28) \qquad\qquad \overline{g} \circ (a e_{A_1}) = h(a) g e_{A_2}$$

$$(10.1.29) \qquad\qquad w = (\overline{g} \circ v)_* {}^* f$$

$$(10.1.30) \qquad\qquad \overline{f} \circ (v_i e^i_{V_1}) = (\overline{g} \circ v_i) f^i_k e^k_{V_2}$$

$$(10.1.31) \qquad\qquad \overline{f} \circ (v_* {}^* e_{V_1}) = (\overline{g} \circ v)_* f_* {}^* e_{V_2}$$

*relative to selected bases. Here*

- *a is coordinate matrix of $A_1$-number $\overline{a}$ relative the basis $\overline{\overline{e}}_{A_1}$*

$$\overline{a} = a e_{A_1}$$

- $h(a) = (h(a^k), k \in K)$ *is a matrix of $D_2$-numbers.*
- *b is coordinate matrix of $A_2$-number*

$$\overline{b} = \overline{g} \circ \overline{a}$$

  *relative the basis $\overline{\overline{e}}_{A_2}$*

$$\overline{b} = b e_{A_2}$$

- *g is coordinate matrix of set of $A_2$-numbers $(\overline{g} \circ e^k_{A_1}, k \in K)$ relative the basis $\overline{\overline{e}}_{A_2}$.*
- *There is relation between the matrix of homomorphism and structural constants*

$$(10.1.32) \qquad\qquad h(C_{1\,k}^{\;ij})g^k_l = g^i_p g^j_q C_{2\,l}^{\;pq}$$

- *v is coordinate matrix of $V_1$-number $\overline{v}$ relative the basis $\overline{\overline{e}}_{V_1}$*

$$\overline{v} = v_* {}^* e_{V_1}$$

- $g(v) = (g(v^i), i \in I)$ *is a matrix of $A_2$-numbers.*
- *w is coordinate matrix of $V_2$-number*

$$\overline{w} = \overline{f} \circ \overline{v}$$

  *relative the basis $\overline{\overline{e}}_{V_2}$*

$$\overline{w} = w_* {}^* e_{V_2}$$

- *f is coordinate matrix of set of $V_2$-numbers $(\overline{f} \circ e^i_{V_1}, i \in I)$ relative the basis $\overline{\overline{e}}_{V_2}$.*

*For given homomorphism*

$$(10.1.3) \quad (h : D_1 \to D_2 \quad \overline{g} : A_1 \to A_2 \quad \overline{f} : V_1 \to V_2)$$

*the set of matrices $(g, f)$ is unique and is called* **coordinates of homomorphism** (10.1.3) *relative bases* $(\overline{\overline{e}}_{A_1}, \overline{\overline{e}}_{V_1}), \quad (\overline{\overline{e}}_{A_2}, \overline{\overline{e}}_{V_2}).$

---

space $V_i$.

PROOF.    According to definitions   2.3.9, 2.8.5, 5.4.1, 10.1.1,  the map

$$(h : D_1 \to D_2 \quad \overline{g} : A_1 \to A_2)$$

is homomorphism of $D_1$-algebra $A_1$ into $D_2$-algebra $A_2$.  Therefore, equalities
(10.1.26), (10.1.27), (10.1.28), (10.1.32)  follow from equalities

| (5.4.22)  $b = h(a)f$ |

| (5.4.23)  $\overline{f} \circ (a_i e_1^i) = h(a_i) f_k^i e_2^k$ |

| (5.4.24)  $\overline{f} \circ (a e_1) = h(a) f e_2$ |

| (5.4.27)  $h(C_{1\ k}^{\ ij}) f_l^k = f_p^i f_q^j C_{2\ l}^{\ pq}$ |

From the theorem 4.3.4, it follows that the matrix $g$ is unique.

Vector  $\overline{v} \in V_1$  has expansion

$$(10.1.33) \qquad\qquad \overline{v} = v_*{}^* e_{V_1}$$

relative to the basis  $\overline{\overline{e}}_{V_1}$.   Vector  $\overline{w} \in V_2$  has expansion

$$(10.1.34) \qquad\qquad \overline{w} = w_*{}^* e_{V_2}$$

relative to the basis  $\overline{\overline{e}}_{V_2}$.   Since $\overline{f}$ is a homomorphism, then the equality

$$(10.1.35) \qquad \begin{aligned} \overline{w} = \overline{f} \circ \overline{v} &= \overline{f} \circ (v_*{}^* e_{V_1}) \\ &= (\overline{g} \circ v)_*{}^* (\overline{f} \circ e_{V_1}) \end{aligned}$$

follows from equalities (10.1.9), (10.1.10).  $V_2$-number  $\overline{f} \circ e_{V_1}^i$  has expansion

$$(10.1.36) \qquad \overline{f} \circ e_{V_1}^i = f^i{}_*{}^* e_{V_2} = f_j^i e_{V_2}^j$$

relative to basis  $\overline{\overline{e}}_{V_2}$.   The equality

$$(10.1.37) \qquad\qquad \overline{w} = (\overline{g} \circ v)_*{}^* f_*{}^* e_{V_2}$$

follows from equalities (10.1.35), (10.1.36).  The equality (10.1.29) follows from
comparison of (10.1.34) and (10.1.37) and the theorem 8.4.11.  From the equality
(10.1.36) and from the theorem 8.4.11, it follows that the matrix $f$ is unique.   □

THEOREM 10.1.7.  *Let the map*

$$h : D_1 \to D_2$$

*be homomorphism of the ring $D_1$ into the ring $D_2$.*         *Let*

$$g = (g_k^l, k \in K, l \in L)$$

*be matrix of $D_2$-numbers which satisfy the equality*

| (10.1.32)  $h(C_{1\ k}^{\ ij}) g_l^k = g_p^i g_q^j C_{2\ l}^{\ pq}$ |

*Let*

$$f = (f_i^j, i \in I, j \in J)$$

*be matrix of $A_2$-numbers.*         *The map* [10.3]

| (10.1.3)  $(h : D_1 \to D_2 \quad \overline{g} : A_1 \to A_2 \quad \overline{f} : V_1 \to V_2)$ |

*defined by equalities*

| (10.1.28)  $\overline{g} \circ (a e_{A_1}) = h(a) g e_{A_2}$ |

| (10.1.31)  $\overline{f} \circ (v_*{}^* e_{V_1}) = (\overline{g} \circ v)_*{}^* f_*{}^* e_{V_2}$ |

*is homomorphism of left $A_1$-vector space of rows $V_1$ into left $A_2$-vector space of rows $V_2$. The homomorphism*

$$(10.1.3) \quad (h: D_1 \to D_2 \quad \overline{g}: A_1 \to A_2 \quad \overline{f}: V_1 \to V_2)$$

*which has the given set of matrices $(g, f)$ is unique.*

Proof.    According to the theorem 5.4.6, the map

$$(10.1.38) \qquad (h: D_1 \to D_2 \quad \overline{g}: A_1 \to A_2)$$

is homomorphism of $D_1$-algebra $A_1$ into $D_2$-algebra $A_2$ and the homomorphism (10.1.38) is unique.

The equality

$$
\begin{aligned}
&\overline{f} \circ ((v+w)_*{}^* e_{V_1}) \\
&= (g \circ (v+w))_*{}^* f_*{}^* e_{V_2} \\
&= (g \circ v)_*{}^* f_*{}^* e_{V_2} \\
&+ (g \circ w)_*{}^* f_*{}^* e_{V_2} \\
&= \overline{f} \circ (v_*{}^* e_{V_1}) + \overline{f} \circ (w_*{}^* e_{V_1})
\end{aligned}
$$

(10.1.39)

follows from equalities

$$(8.1.17) \quad (b_1 + b_2)_*{}^* a = b_{1*}{}^* a + b_{2*}{}^* a \qquad (10.1.6) \quad g \circ (a+b) = g \circ a + g \circ b$$

$$(10.1.10) \quad f \circ (av) = (g \circ a)(f \circ v) \qquad (10.1.31) \quad \overline{f} \circ (v_*{}^* e_{V_1}) = (\overline{g} \circ v)_*{}^* f_*{}^* e_{V_2}$$

From the equality (10.1.39), it follows that the map $\overline{f}$ is homomorphism of Abelian group. The equality

$$
\begin{aligned}
&\overline{f} \circ ((av)_*{}^* e_{V_1}) \\
&= (\overline{g} \circ (av))_*{}^* f_*{}^* e_{V_2} \\
&= (\overline{g} \circ a)((\overline{g} \circ v)_*{}^* f_*{}^* e_{V_2}) \\
&= (\overline{g} \circ a)(\overline{f} \circ (e_{V_1}{}^* {}_* v))
\end{aligned}
$$

(10.1.40)

follows from equalities

$$(10.1.10) \quad f \circ (av) = (g \circ a)(f \circ v)$$

$$(5.4.4) \quad g \circ (ab) = (g \circ a)(g \circ b)$$

$$(10.1.31) \quad \overline{f} \circ (v_*{}^* e_{V_1}) = (\overline{g} \circ v)_*{}^* f_*{}^* e_{V_2}$$

From the equality (10.1.40), and definitions 2.8.5, 10.1.1, it follows that the map

$$(10.1.3) \quad (h: D_1 \to D_2 \quad \overline{g}: A_1 \to A_2 \quad \overline{f}: V_1 \to V_2)$$

is homomorphism of $A_1$-vector space of rows $V_1$ into $A_2$-vector space of rows $V_2$.

Let $f$ be matrix of homomorphisms $\overline{f}$, $\overline{g}$ relative to bases $\overline{\overline{e}}_V$, $\overline{\overline{e}}_W$. The equality

$$\overline{f} \circ (v_*{}^* e_V) = v_*{}^* f_*{}^* e_W = \overline{g} \circ (v_*{}^* e_V)$$

follows from the theorem 10.3.7. Therefore, $\overline{f} = \overline{g}$.    □

On the basis of theorems 10.1.6, 10.1.7 we identify the homomorphism

$$(10.1.3) \quad (h: D_1 \to D_2 \quad \overline{g}: A_1 \to A_2 \quad \overline{f}: V_1 \to V_2)$$

of left $A$-vector space $V$ of rows and coordinates of its presentation

$$(10.1.28) \quad \overline{g} \circ (a e_{A_1}) = h(a) g e_{A_2}$$

$$(10.1.31) \quad \overline{f} \circ (v_* {}^* e_{V_1}) = (\overline{g} \circ v)_* {}^* f_* {}^* e_{V_2}$$

## 10.2. Homomorphism When Rings $D_1 = D_2 = D$

Let $A_i$, $i = 1,\ 2$, be division algebra over commutative ring $D$. Let $V_i$, $i = 1,\ 2$, be left $A_i$-vector space.

DEFINITION 10.2.1. *Let diagram of representations*

$$(10.2.1)$$

$$g_{1.12}(d) : a \to d\,a$$
$$g_{1.23}(v) : w \to C_1(w, v)$$
$$C_1 \in \mathcal{L}(A_1^2 \to A_1)$$
$$g_{1.34}(a) : v \to a\,v$$
$$g_{1.14}(d) : v \to d\,v$$

*describe left $A_1$-vector space $V_1$. Let diagram of representations*

$$(10.2.2)$$

$$g_{2.12}(d) : a \to d\,a$$
$$g_{2.23}(v) : w \to C_2(w, v)$$
$$C_2 \in \mathcal{L}(A_2^2 \to A_2)$$
$$g_{2.34}(a) : v \to a\,v$$
$$g_{2.14}(d) : v \to d\,v$$

*describe left $A_2$-vector space $V_2$. Morphism*

$$(10.2.3) \quad (g : A_1 \to A_2 \quad f : V_1 \to V_2)$$

*of diagram of representations* (10.2.1) *into diagram of representations* (10.2.2) *is called* **homomorphism** *of left $A_1$-vector space $V_1$ into left $A_2$-vector space $V_2$. Let us denote* $\mathrm{Hom}(D; A_1 \to A_2; *V_1 \to *V_2)$ *set of homomorphisms of left $A_1$-vector space $V_1$ into left $A_2$-vector space $V_2$.* $\qquad\square$

We will use notation

$$f \circ a = f(a)$$

for image of homomorphism $f$.

THEOREM 10.2.2. *The homomorphism*

$$(10.2.3) \quad (g : A_1 \to A_2 \quad f : V_1 \to V_2)$$

*of left $A_1$-vector space $V_1$ into left $A_2$-vector space $V_2$ satisfies following equalities*

$$(10.2.4) \quad g \circ (a + b) = g \circ a + g \circ b$$

$$(10.2.5) \quad g \circ (ab) = (g \circ a)(g \circ b)$$

$$(10.2.6) \quad g \circ (pa) = p(g \circ a)$$

$$(10.2.7) \qquad f \circ (u + v) = f \circ u + f \circ v$$

$$(10.2.8) \qquad f \circ (av) = (g \circ a)(f \circ v)$$

$$p \in D \quad a, b \in A_1 \quad u, v \in V_1$$

PROOF. According to definitions 2.3.9, 2.8.5, 7.2.1, 1 the map

$$g : A_1 \to A_2$$

is homomorphism of $D$-algebra $A_1$ into $D$-algebra $A_2$. Therefore, equalities (10.2.4), (10.2.5), (10.2.6) follow from the theorem 4.3.10. The equality (10.2.7) follows from the definition 10.2.1, since, accorfing to the definition 2.3.9, the map $f$ is homomorpism of Abelian group. The equality (10.2.8) follows from the equality (2.3.6) because the map

$$(g : A_1 \to A_2 \quad f : V_1 \to V_2)$$

is morphism of representation $g_{1.34}$ into representation $g_{2.34}$. □

DEFINITION 10.2.3. *Homomorphism* [10.4]

$$(g : A_1 \to A_2 \quad f : V_1 \to V_2)$$

*is called* **isomorphism** *between left $A_1$-vector space $V_1$ and left $A_2$-vector space $V_2$, if there exists the map*

$$(g^{-1} : A_2 \to A_1 \quad f^{-1} : V_2 \to V_1)$$

*which is homomorphism.* □

THEOREM 10.2.4. *The homomorphism* [10.5]

$$(10.2.3) \qquad (\overline{g} : A_1 \to A_2 \quad \overline{f} : V_1 \to V_2)$$

*of left $A_1$-vector space of columns $V_1$ into left $A_2$-vector space of columns $V_2$ has presentation*

$$(10.2.9) \qquad b = ga$$

$$(10.2.10) \qquad \overline{g} \circ (a^i e_{A_1 i}) = a^i g_i^k e_{A_2 k}$$

$$(10.2.11) \qquad \overline{g} \circ (e_{A_1} a) = e_{A_2} ga$$

THEOREM 10.2.5. *The homomorphism* [10.6]

$$(10.2.3) \qquad (\overline{g} : A_1 \to A_2 \quad \overline{f} : V_1 \to V_2)$$

*of left $A_1$-vector space of rows $V_1$ into left $A_2$-vector space of rows $V_2$ has presentation*

$$(10.2.16) \qquad b = ag$$

$$(10.2.17) \qquad \overline{g} \circ (a_i e_{A_1}^i) = a_i g_k^i e_{A_2}^k$$

$$(10.2.18) \qquad \overline{g} \circ (ae_{A_1}) = age_{A_2}$$

---

[10.4] I follow the definition on page [18]-49.

[10.5] In theorems 10.2.4, 10.2.6, we use the following convention. Let $i = 1, 2$. Let the set of vectors $\overline{\overline{e}}_{A_i} = (e_{A_i k}, k \in K_i)$ be a basis and $C_{ij}^k$, $k, i, j \in K_i$, be structural constants of $D$-algebra of columns $A_i$. Let the set of vectors $\overline{\overline{e}}_{V_i} = (e_{V_i i}, i \in I_i)$ be a basis of left $A_i$-vector space $V_i$.

[10.6] In theorems 10.2.5, 10.2.7, we use the following convention. Let $i = 1, 2$. Let the set of vectors $\overline{\overline{e}}_{A_i} = (e_{A_i}^k, k \in K_i)$ be a basis and $C_{ik}^{ij}$, $k, i, j \in K_i$, be structural constants of $D$-algebra of rows $A_i$. Let the set of vectors $\overline{\overline{e}}_{V_i} = (e_{V_i}^i, i \in I_i)$ be a basis of left $A_i$-vector space $V_i$.

$$(10.2.12) \qquad w = (\overline{g} \circ v)^* {}_* f$$

$$(10.2.13) \quad \overline{f} \circ (v^i e_{V_1 i}) = (\overline{g} \circ v^i) f_i^k e_{V_2 k}$$

$$(10.2.14)$$
$$\overline{f} \circ (v^* {}_* e_{V_1}) = (\overline{g} \circ v)^* {}_* (f^* {}_* e_{V_2})$$

*relative to selected bases. Here*

- *$a$ is coordinate matrix of $A_1$-number $\overline{a}$ relative the basis $\overline{\overline{e}}_{A_1}$*

$$\overline{a} = e_{A_1} a$$

- *$b$ is coordinate matrix of $A_2$-number*

$$\overline{b} = \overline{g} \circ \overline{a}$$

  *relative the basis $\overline{\overline{e}}_{A_2}$*

$$\overline{b} = e_{A_2} b$$

- *$g$ is coordinate matrix of set of $A_2$-numbers $(\overline{g} \circ e_{A_1 k}, k \in K)$ relative the basis $\overline{\overline{e}}_{A_2}$.*
- *There is relation between the matrix of homomorphism and structural constants*

$$(10.2.15) \qquad C_{1\,ij}^{\ \ k} g_k^l = g_i^p g_j^q C_{2\,pq}^{\ \ l}$$

- *$v$ is coordinate matrix of $V_1$-number $\overline{v}$ relative the basis $\overline{\overline{e}}_{V_1}$*

$$\overline{v} = v^* {}_* e_{V_1}$$

- *$g(v) = (g(v_i), i \in I)$ is a matrix of $A_2$-numbers.*
- *$w$ is coordinate matrix of $V_2$-number*

$$\overline{w} = \overline{f} \circ \overline{v}$$

  *relative the basis $\overline{\overline{e}}_{V_2}$*

$$\overline{w} = w^* {}_* e_{V_2}$$

- *$f$ is coordinate matrix of set of $V_2$-numbers $(\overline{f} \circ e_{V_1 i}, i \in I)$ relative the basis $\overline{\overline{e}}_{V_2}$.*

*For given homomorphism*

$$(10.2.3) \quad (\overline{g} : A_1 \to A_2 \quad \overline{f} : V_1 \to V_2)$$

*the set of matrices $(g, f)$ is unique and is called* **coordinates of homomorphism** (10.2.3) *relative bases* $(\overline{e}_{A_1}, \overline{\overline{e}}_{V_1}), \quad (\overline{e}_{A_2}, \overline{\overline{e}}_{V_2}).$

$$(10.2.19) \qquad w = (\overline{g} \circ v)_* {}^* f$$

$$(10.2.20) \quad \overline{f} \circ (v_i e_{V_1}^i) = (\overline{g} \circ v_i) f_k^i e_{V_2}^k$$

$$(10.2.21) \quad \overline{f} \circ (v_* {}^* e_{V_1}) = (\overline{g} \circ v)_* {}^* f_* {}^* e_{V_2}$$

*relative to selected bases. Here*

- *$a$ is coordinate matrix of $A_1$-number $\overline{a}$ relative the basis $\overline{\overline{e}}_{A_1}$*

$$\overline{a} = a e_{A_1}$$

- *$b$ is coordinate matrix of $A_2$-number*

$$\overline{b} = \overline{g} \circ \overline{a}$$

  *relative the basis $\overline{\overline{e}}_{A_2}$*

$$\overline{b} = b e_{A_2}$$

- *$g$ is coordinate matrix of set of $A_2$-numbers $(\overline{g} \circ e_{A_1}^k, k \in K)$ relative the basis $\overline{\overline{e}}_{A_2}$.*
- *There is relation between the matrix of homomorphism and structural constants*

$$(10.2.22) \qquad C_{1k}^{\ ij} g_l^k = g_p^i g_q^j C_{2l}^{\ pq}$$

- *$v$ is coordinate matrix of $V_1$-number $\overline{v}$ relative the basis $\overline{\overline{e}}_{V_1}$*

$$\overline{v} = v_* {}^* e_{V_1}$$

- *$g(v) = (g(v^i), i \in I)$ is a matrix of $A_2$-numbers.*
- *$w$ is coordinate matrix of $V_2$-number*

$$\overline{w} = \overline{f} \circ \overline{v}$$

  *relative the basis $\overline{\overline{e}}_{V_2}$*

$$\overline{w} = w_* {}^* e_{V_2}$$

- *$f$ is coordinate matrix of set of $V_2$-numbers $(\overline{f} \circ e_{V_1}^i, i \in I)$ relative the basis $\overline{\overline{e}}_{V_2}$.*

*For given homomorphism*

$$(10.2.3) \quad (\overline{g} : A_1 \to A_2 \quad \overline{f} : V_1 \to V_2)$$

*the set of matrices $(g, f)$ is unique and is called* **coordinates of homomorphism** (10.2.3) *relative bases* $(\overline{e}_{A_1}, \overline{\overline{e}}_{V_1}), \quad (\overline{e}_{A_2}, \overline{\overline{e}}_{V_2}).$

Proof.   According to definitions 2.3.9, 2.8.5, 5.4.7, 10.2.1, the map

$$\overline{g} : A_1 \to A_2$$

is homomorphism of $D$-algebra $A_1$ into $D$-algebra $A_2$. Therefore, equalities (10.2.9), (10.2.10), (10.2.11), (10.2.15) follow from equalities

| (5.4.53)   $b = fa$ |

| (5.4.54)   $\overline{f} \circ (a^i e_{1i}) = a^i f_i^k e_k$ |

| (5.4.55)   $\overline{f} \circ (e_1 a) = e_2 fa$ |

| (5.4.58)   $C_{1ij}^{\ k} f_k^l = f_i^p f_j^q C_{2pq}^{\ l}$ |

From the theorem 4.3.11, it follows that the matrix $g$ is unique.

Vector $\overline{v} \in V_1$ has expansion

$$(10.2.23) \qquad \overline{v} = v^* {}_* e_{V_1}$$

relative to the basis $\overline{\overline{e}}_{V_1}$. Vector $\overline{w} \in V_2$ has expansion

$$(10.2.24) \qquad \overline{w} = w^* {}_* e_{V_2}$$

relative to the basis $\overline{\overline{e}}_{V_2}$. Since $\overline{f}$ is a homomorphism, then the equality

$$(10.2.25) \qquad \begin{aligned} \overline{w} &= \overline{f} \circ \overline{v} = \overline{f} \circ (v^* {}_* e_{V_1}) \\ &= (\overline{g} \circ v)^* {}_* (\overline{f} \circ e_{V_1}) \end{aligned}$$

follows from equalities (10.2.7), (10.2.8). $V_2$-number $\overline{f} \circ e_{V_1 i}$ has expansion

$$(10.2.26) \quad \overline{f} \circ e_{V_1 i} = f_i^* {}_* e_{V_2} = f_i^j e_{V_2 j}$$

relative to basis $\overline{\overline{e}}_{V_2}$. The equality

$$(10.2.27) \qquad \overline{w} = (\overline{g} \circ v)^* {}_* f^* {}_* e_{V_2}$$

follows from equalities (10.2.25), (10.2.26). The equality (10.2.12) follows from comparison of (10.2.24) and (10.2.27) and the theorem 8.4.4. From the equality (10.2.26) and from the theorem 8.4.4, it follows that the matrix $f$ is unique.   □

Proof.   According to definitions 2.3.9, 2.8.5, 5.4.7, 10.2.1, the map

$$\overline{g} : A_1 \to A_2$$

is homomorphism of $D$-algebra $A_1$ into $D$-algebra $A_2$. Therefore, equalities (10.2.16), (10.2.17), (10.2.18), (10.2.22) follow from equalities

| (5.4.64)   $b = af$ |

| (5.4.65)   $\overline{f} \circ (a_i e_1^i) = a_i f_k^i e^k$ |

| (5.4.66)   $\overline{f} \circ (ae_1) = af e_2$ |

| (5.4.69)   $C_{1k}^{\ ij} f_l^k = f_p^i f_q^j C_{2l}^{\ pq}$ |

From the theorem 4.3.12, it follows that the matrix $g$ is unique.

Vector $\overline{v} \in V_1$ has expansion

$$(10.2.28) \qquad \overline{v} = v_* {}^* e_{V_1}$$

relative to the basis $\overline{\overline{e}}_{V_1}$. Vector $\overline{w} \in V_2$ has expansion

$$(10.2.29) \qquad \overline{w} = w_* {}^* e_{V_2}$$

relative to the basis $\overline{\overline{e}}_{V_2}$. Since $\overline{f}$ is a homomorphism, then the equality

$$(10.2.30) \qquad \begin{aligned} \overline{w} &= \overline{f} \circ \overline{v} = \overline{f} \circ (v_* {}^* e_{V_1}) \\ &= (\overline{g} \circ v)_* {}^* (\overline{f} \circ e_{V_1}) \end{aligned}$$

follows from equalities (10.2.7), (10.2.8). $V_2$-number $\overline{f} \circ e_{V_1}^i$ has expansion

$$(10.2.31) \quad \overline{f} \circ e_{V_1}^i = f^i {}_* {}^* e_{V_2} = f_j^i e_{V_2}^j$$

relative to basis $\overline{\overline{e}}_{V_2}$. The equality

$$(10.2.32) \qquad \overline{w} = (\overline{g} \circ v)_* {}^* f_* {}^* e_{V_2}$$

follows from equalities (10.2.30), (10.2.31). The equality (10.2.19) follows from comparison of (10.2.29) and (10.2.32) and the theorem 8.4.11. From the equality (10.2.31) and from the theorem 8.4.11, it follows that the matrix $f$ is unique.   □

Theorem 10.2.6. *Let*

$$g = (g_l^k, k \in K, l \in L)$$

*be matrix of $D$-numbers which satisfy the equality*

| (10.2.15)   $C_{1ij}^{\ k} g_k^l = g_i^p g_j^q C_{2pq}^{\ l}$ |

Theorem 10.2.7. *Let*

$$g = (g_k^l, k \in K, l \in L)$$

*be matrix of $D$-numbers which satisfy the equality*

| (10.2.22)   $C_{1k}^{\ ij} g_l^k = g_p^i g_q^j C_{2l}^{\ pq}$ |

*Let*
$$f = (f_j^i, i \in I, j \in J)$$
*be matrix of $A_2$-numbers.        The map* 10.5

| (10.2.3)    $(\overline{g} : A_1 \to A_2 \quad \overline{f} : V_1 \to V_2)$ |

*defined by equalities*

| (10.2.11)    $\overline{f} \circ (e_{A_1} a) = e_{A_2} f a$ |

| (10.2.14)    $\overline{f} \circ (v^*{}_* e_{V_1}) = (\overline{g} \circ v)^*{}_* (f^*{}_* e_{V_2})$ |

*is homomorphism of left $A_1$-vector space of columns $V_1$ into left $A_2$-vector space of columns $V_2$.   The homomorphism*

| (10.2.3)    $(\overline{g} : A_1 \to A_2 \quad \overline{f} : V_1 \to V_2)$ |

*which has the given   set of matrices $(g, f)$  is unique.*

PROOF.   According to the theorem 5.4.11, the map
$$(10.2.33) \qquad \overline{g} : A_1 \to A_2$$
is homomorphism of $D$-algebra $A_1$ into $D$-algebra $A_2$ and the homomorphism (10.2.33) is unique.

The equality
$$\begin{aligned}
& \overline{f} \circ ((v+w)^*{}_* e_{V_1}) \\
&= (g \circ (v+w))^*{}_* f^*{}_* e_{V_2} \\
(10.2.34) &= (g \circ v)^*{}_* f^*{}_* e_{V_2} \\
&\quad + (g \circ w)^*{}_* f^*{}_* e_{V_2} \\
&= \overline{f} \circ (v^*{}_* e_{V_1}) + \overline{f} \circ (w^*{}_* e_{V_1})
\end{aligned}$$
follows from equalities

| (8.1.19)    $(b_1 + b_2)^*{}_* a = b_1{}^*{}_* a + b_2{}^*{}_* a$ |
| (10.1.6)    $g \circ (a+b) = g \circ a + g \circ b$ |
| (10.1.10)    $f \circ (av) = (g \circ a)(f \circ v)$ |
| (10.1.16)    $\overline{f} \circ (v^*{}_* e_{V_1}) = (\overline{g} \circ v)^*{}_* (f^*{}_* e_{V_2})$ |

From the equality (10.2.34), it follows that the map $\overline{f}$ is homomorphism of Abelian group. The equality

*Let*
$$f = (f_i^j, i \in I, j \in J)$$
*be matrix of $A_2$-numbers.        The map* 10.6

| (10.2.3)    $(\overline{g} : A_1 \to A_2 \quad \overline{f} : V_1 \to V_2)$ |

*defined by equalities*

| (10.2.18)    $\overline{f} \circ (ae_{A_1}) = af e_{A_2}$ |

| (10.2.21)    $\overline{f} \circ (v_* {}^* e_{V_1}) = (\overline{g} \circ v)_* {}^* f_* {}^* e_{V_2}$ |

*is homomorphism of left $A_1$-vector space of rows $V_1$ into left $A_2$-vector space of rows $V_2$.    The homomorphism*

| (10.2.3)    $(\overline{g} : A_1 \to A_2 \quad \overline{f} : V_1 \to V_2)$ |

*which has the given   set of matrices $(g, f)$  is unique.*

PROOF.   According to the theorem 5.4.12, the map
$$(10.2.36) \qquad \overline{g} : A_1 \to A_2$$
is homomorphism of $D$-algebra $A_1$ into $D$-algebra $A_2$ and the homomorphism (10.2.36) is unique.

The equality
$$\begin{aligned}
& \overline{f} \circ ((v+w)_* {}^* e_{V_1}) \\
&= (g \circ (v+w))_* {}^* f_* {}^* e_{V_2} \\
(10.2.37) &= (g \circ v)_* {}^* f_* {}^* e_{V_2} \\
&\quad + (g \circ w)_* {}^* f_* {}^* e_{V_2} \\
&= \overline{f} \circ (v_* {}^* e_{V_1}) + \overline{f} \circ (w_* {}^* e_{V_1})
\end{aligned}$$
follows from equalities

| (8.1.17)    $(b_1 + b_2)_* {}^* a = b_1{}_* {}^* a + b_2{}_* {}^* a$ |
| (10.1.6)    $g \circ (a+b) = g \circ a + g \circ b$ |
| (10.1.10)    $f \circ (av) = (g \circ a)(f \circ v)$ |
| (10.1.31)    $\overline{f} \circ (v_* {}^* e_{V_1}) = (\overline{g} \circ v)_* {}^* f_* {}^* e_{V_2}$ |

From the equality (10.2.37), it follows that the map $\overline{f}$ is homomorphism of Abelian group. The equality

$$(10.2.35) \quad \begin{aligned} &\overline{f} \circ ((av)^* {}_* e_{V_1}) \\ &= (\overline{g} \circ (av))^* {}_* f^* {}_* e_{V_2} \\ &= (\overline{g} \circ a)((\overline{g} \circ v)^* {}_* f^* {}_* e_{V_2}) \\ &= (\overline{g} \circ a)(\overline{f} \circ (v^* {}_* e_{V_1})) \end{aligned}$$

follows from equalities

| (10.1.10) | $f \circ (av) = (g \circ a)(f \circ v)$ |
|---|---|

| (5.4.4) | $g \circ (ab) = (g \circ a)(g \circ b)$ |
|---|---|

| (10.1.16) | $\overline{f} \circ (v^* {}_* e_{V_1}) = (\overline{g} \circ v)^* {}_* (f^* {}_* e_{V_2})$ |
|---|---|

From the equality (10.2.35), and definitions 2.8.5, 10.2.1, it follows that the map

| (10.2.3) | $(\overline{g} : A_1 \to A_2 \quad \overline{f} : V_1 \to V_2)$ |
|---|---|

is homomorphism of $A_1$-vector space of columns $V_1$ into $A_2$-vector space of columns $V_2$.

Let $f$ be matrix of homomorphisms $\overline{f}$, $\overline{g}$ relative to bases $\overline{\overline{e}}_V$, $\overline{\overline{e}}_W$. The equality

$$\overline{f} \circ (v^* {}_* e_V) = v^* {}_* f^* {}_* e_W = \overline{g} \circ (v^* {}_* e_V)$$

follows from the theorem 10.3.6. Therefore, $\overline{f} = \overline{g}$.    □

On the basis of theorems 10.2.4, 10.2.6 we identify the homomorphism

| (10.2.3) | $(\overline{g} : A_1 \to A_2 \quad \overline{f} : V_1 \to V_2)$ |
|---|---|

of left $A$-vector space $V$ of columns and coordinates of its presentation

| (10.2.11) | $\overline{f} \circ (e_{A_1} a) = e_{A_2} f a$ |
|---|---|

| (10.2.14) | $\overline{f} \circ (v^* {}_* e_{V_1}) = (\overline{g} \circ v)^* {}_* (f^* {}_* e_{V_2})$ |
|---|---|

$$(10.2.38) \quad \begin{aligned} &\overline{f} \circ ((av)_* {}^* e_{V_1}) \\ &= (\overline{g} \circ (av))_* {}^* f_* {}^* e_{V_2} \\ &= (\overline{g} \circ a)((\overline{g} \circ v)_* {}^* f_* {}^* e_{V_2}) \\ &= (\overline{g} \circ a)(\overline{f} \circ (e_{V_1} {}_* {}^* v)) \end{aligned}$$

follows from equalities

| (10.1.10) | $f \circ (av) = (g \circ a)(f \circ v)$ |
|---|---|

| (5.4.4) | $g \circ (ab) = (g \circ a)(g \circ b)$ |
|---|---|

| (10.1.31) | $\overline{f} \circ (v_* {}^* e_{V_1}) = (\overline{g} \circ v)_* {}^* f_* {}^* e_{V_2}$ |
|---|---|

From the equality (10.2.38), and definitions 2.8.5, 10.2.1, it follows that the map

| (10.2.3) | $(\overline{g} : A_1 \to A_2 \quad \overline{f} : V_1 \to V_2)$ |
|---|---|

is homomorphism of $A_1$-vector space of rows $V_1$ into $A_2$-vector space of rows $V_2$.

Let $f$ be matrix of homomorphisms $\overline{f}$, $\overline{g}$ relative to bases $\overline{\overline{e}}_V$, $\overline{\overline{e}}_W$. The equality

$$\overline{f} \circ (v_* {}^* e_V) = v_* {}^* f_* {}^* e_W = \overline{g} \circ (v_* {}^* e_V)$$

follows from the theorem 10.3.7. Therefore, $\overline{f} = \overline{g}$.    □

On the basis of theorems 10.2.5, 10.2.7 we identify the homomorphism

| (10.2.3) | $(\overline{g} : A_1 \to A_2 \quad \overline{f} : V_1 \to V_2)$ |
|---|---|

of left $A$-vector space $V$ of rows and coordinates of its presentation

| (10.2.18) | $\overline{f} \circ (a e_{A_1}) = a f e_{A_2}$ |
|---|---|

| (10.2.21) | $\overline{f} \circ (v_* {}^* e_{V_1}) = (\overline{g} \circ v)_* {}^* f_* {}^* e_{V_2}$ |
|---|---|

## 10.3. Homomorphism When $D$-algebras $A_1 = A_2 = A$

Let $A$ be division algebra over commutative ring $D$. Let $V_i$, $i = 1, 2$, be left $A$-vector space.

DEFINITION 10.3.1. *Let diagram of representations*

(10.3.1)

$$A \xrightarrow{g_{1.23}} A \xrightarrow{g_{1.34}} V_1$$

with $g_{1.12}$, $g_{1.12}$, $g_{1.14}$, $D$

$g_{1.12}(d) : a \to d\,a$

$g_{1.23}(v) : w \to C(w, v)$

$C \in \mathcal{L}(A^2 \to A)$

$g_{1.34}(a) : v \to a\,v$

$g_{1.14}(d) : v \to d\,v$

*describe left $A$-vector space $V_1$. Let diagram of representations*

(10.3.2)

$$A \xrightarrow{g_{2.23}} A \xrightarrow{g_{2.34}} V_2$$

with $g_{2.12}$, $g_{2.12}$, $g_{2.14}$, $D$

$g_{2.12}(d) : a \to d\,a$

$g_{2.23}(v) : w \to C(w, v)$

$C \in \mathcal{L}(A^2 \to A)$

$g_{2.34}(a) : v \to a\,v$

$g_{2.14}(d) : v \to d\,v$

*describe left $A$-vector space $V_2$. Morphism*

(10.3.3)
$$f : V_1 \to V_2$$

*of diagram of representations* (10.3.1) *into diagram of representations* (10.3.2) *is called* **homomorphism** *of left $A$-vector space $V_1$ into left $A$-vector space $V_2$. Let us denote* $\mathrm{Hom}(D; A; *V_1 \to *V_2)$ *set of homomorphisms of left $A$-vector space $V_1$ into left $A$-vector space $V_2$.* □

We will use notation
$$f \circ a = f(a)$$
for image of homomorphism $f$.

THEOREM 10.3.2. *The homomorphism*

(10.3.3)   $f : V_1 \to V_2$

*of left $A$-vector space $V_1$ into left $A$-vector space $V_2$ satisfies following equalities*

(10.3.4)
$$f \circ (u + v) = f \circ u + f \circ v$$

(10.3.5)
$$f \circ (av) = a(f \circ v)$$
$$a \in A \quad u, v \in V_1$$

PROOF.     The equality (10.3.4) follows from the definition 10.3.1, since, accorfing to the definition 2.3.9, the map $f$ is homomorpism of Abelian group. The equality (10.3.5) follows from the equality (2.3.6) because the map

$$f : V_1 \to V_2$$

is morphism of representation $g_{1.34}$ into representation $g_{2.34}$.     □

---

10.7   I follow the definition on page [18]-49.

DEFINITION 10.3.3. *Homomorphism* [10.7]

$$f : V_1 \to V_2$$

*is called* **isomorphism** *between left $A$-vector space $V_1$ and left $A$-vector space $V_2$, if there exists the map*

$$f^{-1} : V_2 \to V_1$$

*which is homomorphism.* □

DEFINITION 10.3.4. *A homomorphism* [10.7]

$$f : V \to V$$

*in which source and target are the same left $A$-vector space is called* **endomorphism**. *Endomorphism*

$$f : V \to V$$

*of left $A$-vector space $V$ is called* **automorphism**, *if there exists the map $f^{-1}$ which is endomorphism.* □

THEOREM 10.3.5. *The set $GL(V)$ of automorphisms of left $A$-vector space $V$ is group.*

THEOREM 10.3.6. *The homomorphism* [10.8]

(10.3.3)  $\boxed{\overline{f} : V_1 \to V_2}$

*of left $A$-vector space of columns $V_1$ into left $A$-vector space of columns $V_2$ has presentation*

(10.3.6)  $\qquad w = v^* {}_* f$

(10.3.7)  $\qquad \overline{f} \circ (v^i e_{V_1 i}) = v^i f_i^k e_{V_2 k}$

(10.3.8)  $\qquad \overline{f} \circ (v^* {}_* e_{V_1}) = v^* {}_* f^* {}_* e_{V_2}$

*relative to selected bases. Here*

- $v$ *is coordinate matrix of $V_1$-number $\overline{v}$ relative to the basis $\overline{\overline{e}}_{V_1}$*

  $$\overline{v} = v^* {}_* e_{V_1}$$

- $w$ *is coordinate matrix of $V_2$-number*

  $$\overline{w} = \overline{f} \circ \overline{v}$$

THEOREM 10.3.7. *The homomorphism* [10.9]

(10.3.3)  $\boxed{\overline{f} : V_1 \to V_2}$

*of left $A$-vector space of rows $V_1$ into left $A$-vector space of rows $V_2$ has presentation*

(10.3.9)  $\qquad w = v_* {}^* f$

(10.3.10)  $\qquad \overline{f} \circ (v_i e_{V_1}^i) = v_i f_k^i e_{V_2}^k$

(10.3.11)  $\qquad \overline{f} \circ (v_* {}^* e_{V_1}) = v_* {}^* f_* {}^* e_{V_2}$

*relative to selected bases. Here*

- $v$ *is coordinate matrix of $V_1$-number $\overline{v}$ relative to the basis $\overline{\overline{e}}_{V_1}$*

  $$\overline{v} = v_* {}^* e_{V_1}$$

- $w$ *is coordinate matrix of $V_2$-number*

  $$\overline{w} = \overline{f} \circ \overline{v}$$

---

[10.8] In theorems 10.3.6, 10.3.8, we use the following convention. Let $i = 1, 2$. Let the set of vectors $\overline{\overline{e}}_{V_i} = (e_{V_i i}, i \in \mathcal{I}_i)$ be a basis of left $A$-vector space $V_i$.

[10.9] In theorems 10.3.7, 10.3.9, we use the following convention. Let $i = 1, 2$. Let the set of vectors $\overline{\overline{e}}_{V_i} = (e_{V_i}^i, i \in \mathcal{I}_i)$ be a basis of left $A$-vector space $V_i$.

*relative the basis* $\overline{\overline{e}}_{V_2}$

$$\overline{w} = w^* {}_* e_{V_2}$$

- $f$ *is coordinate matrix of set of $V_2$-numbers* $(\overline{f} \circ e_{V_1 i}, i \in I)$ *relative the basis* $\overline{\overline{e}}_{V_2}$.

*For given homomorphism* $\overline{f}$, *the matrix* $f$ *is unique and is called* **matrix of homomorphism** $\overline{f}$ *relative bases* $\overline{\overline{e}}_1$, $\overline{\overline{e}}_2$.

PROOF. Vector $\overline{v} \in V_1$ has expansion

(10.3.12) $\qquad \overline{v} = v^* {}_* e_{V_1}$

relative to the basis $\overline{\overline{e}}_{V_1}$. Vector $\overline{w} \in V_2$ has expansion

(10.3.13) $\qquad \overline{w} = w^* {}_* e_{V_2}$

relative to the basis $\overline{\overline{e}}_{V_2}$. Since $\overline{f}$ is a homomorphism, then the equality

(10.3.14) $\quad \begin{aligned} \overline{w} &= \overline{f} \circ \overline{v} = \overline{f} \circ (v^* {}_* e_{V_1}) \\ &= v^* {}_* (\overline{f} \circ e_{V_1}) \end{aligned}$

follows from equalities (10.3.4), (10.3.5). $V_2$-number $\overline{f} \circ e_{V_1 i}$ has expansion

(10.3.15) $\quad \overline{f} \circ e_{V_1 i} = f_i{}^* {}_* e_{V_2} = f_i^j e_{V_2 j}$

relative to basis $\overline{\overline{e}}_{V_2}$. The equality

(10.3.16) $\qquad \overline{w} = v^* {}_* f^* {}_* e_{V_2}$

follows from equalities (10.3.14), (10.3.15). The equality (10.3.6) follows from comparison of (10.3.13) and (10.3.16) and the theorem 8.4.4. From the equality (10.3.15) and from the theorem 8.4.4, it follows that the matrix $f$ is unique. $\quad \square$

---

*relative the basis* $\overline{\overline{e}}_{V_2}$

$$\overline{w} = w_* {}^* e_{V_2}$$

- $f$ *is coordinate matrix of set of $V_2$-numbers* $(\overline{f} \circ e^i_{V_1}, i \in I)$ *relative the basis* $\overline{\overline{e}}_{V_2}$.

*For given homomorphism* $\overline{f}$, *the matrix* $f$ *is unique and is called* **matrix of homomorphism** $\overline{f}$ *relative bases* $\overline{\overline{e}}_1$, $\overline{\overline{e}}_2$.

PROOF. Vector $\overline{v} \in V_1$ has expansion

(10.3.17) $\qquad \overline{v} = v_* {}^* e_{V_1}$

relative to the basis $\overline{\overline{e}}_{V_1}$. Vector $\overline{w} \in V_2$ has expansion

(10.3.18) $\qquad \overline{w} = w_* {}^* e_{V_2}$

relative to the basis $\overline{\overline{e}}_{V_2}$. Since $\overline{f}$ is a homomorphism, then the equality

(10.3.19) $\quad \begin{aligned} \overline{w} &= \overline{f} \circ \overline{v} = \overline{f} \circ (v_* {}^* e_{V_1}) \\ &= v_* {}^* (\overline{f} \circ e_{V_1}) \end{aligned}$

follows from equalities (10.3.4), (10.3.5). $V_2$-number $\overline{f} \circ e^i_{V_1}$ has expansion

(10.3.20) $\quad \overline{f} \circ e^i_{V_1} = f^i {}_* {}^* e_{V_2} = f^i_j e^j_{V_2}$

relative to basis $\overline{\overline{e}}_{V_2}$. The equality

(10.3.21) $\qquad \overline{w} = v_* {}^* f_* {}^* e_{V_2}$

follows from equalities (10.3.19), (10.3.20). The equality (10.3.9) follows from comparison of (10.3.18) and (10.3.21) and the theorem 8.4.11. From the equality (10.3.20) and from the theorem 8.4.11, it follows that the matrix $f$ is unique. $\quad \square$

---

THEOREM 10.3.8. *Let*

$$f = (f^i_j, i \in I, j \in J)$$

*be matrix of $A$-numbers. The map*[10.8]

| (10.3.3) | $\overline{f} : V_1 \to V_2$ |
|---|---|

*defined by the equality*

| (10.3.8) | $\overline{f} \circ (v^* {}_* e_{V_1}) = v^* {}_* f^* {}_* e_{V_2}$ |
|---|---|

*is homomorphism of left $A$-vector space*

---

THEOREM 10.3.9. *Let*

$$f = (f^j_i, i \in I, j \in J)$$

*be matrix of $A$-numbers. The map*[10.9]

| (10.3.3) | $\overline{f} : V_1 \to V_2$ |
|---|---|

*defined by the equality*

| (10.3.11) | $\overline{f} \circ (v_* {}^* e_{V_1}) = v_* {}^* f_* {}^* e_{V_2}$ |
|---|---|

*is homomorphism of left $A$-vector space*

of columns $V_1$ into left $A$-vector space of columns $V_2$. The homomorphism

(10.3.3) $\quad \overline{f} : V_1 \to V_2$

which has the given matrix $f$ is unique.

Proof. The equality

$$\overline{f} \circ ((v + w)^*{}_*e_{V_1})$$
$$= (v + w)^*{}_*f^*{}_*e_{V_2}$$
(10.3.22)
$$= v^*{}_*f^*{}_*e_{V_2} + w^*{}_*f^*{}_*e_{V_2}$$
$$= \overline{f} \circ (v^*{}_*e_{V_1}) + \overline{f} \circ (w^*{}_*e_{V_1})$$

follows from equalities

(8.1.19) $\quad (b_1 + b_2)^*{}_*a = b_1{}^*{}_*a + b_2{}^*{}_*a$

(10.2.8) $\quad f \circ (av) = a(f \circ v)$

(10.3.8) $\quad \overline{f} \circ (v^*{}_*e_{V_1}) = v^*{}_*f^*{}_*e_{V_2}$

From the equality (10.3.22), it follows that the map $\overline{f}$ is homomorphism of Abelian group. The equality

$$\overline{f} \circ ((av)^*{}_*e_{V_1})$$
$$= (av)^*{}_*f^*{}_*e_{V_2}$$
(10.3.23)
$$= a(v^*{}_*f^*{}_*e_{V_2})$$
$$= a(\overline{f} \circ (v^*{}_*e_{V_1}))$$

follows from equalities

(10.2.8) $\quad f \circ (av) = a(f \circ v)$

(10.3.8) $\quad \overline{f} \circ (v^*{}_*e_{V_1}) = v^*{}_*f^*{}_*e_{V_2}$

From the equality (10.3.23), and definitions 2.8.5, 10.3.1, it follows that the map

(10.3.3) $\quad \overline{f} : V_1 \to V_2$

is homomorphism of $A$-vector space of columns $V_1$ into $A$-vector space of columns $V_2$.

Let $f$ be matrix of homomorphisms $\overline{f}$, $\overline{g}$ relative to bases $\overline{e}_V$, $\overline{e}_W$. The equality

$\overline{f} \circ (v^*{}_*e_V) = v^*{}_*f^*{}_*e_W = \overline{g} \circ (v^*{}_*e_V)$

follows from the theorem 10.3.6. Therefore, $\overline{f} = \overline{g}$. □

On the basis of theorems 10.3.6, 10.3.8 we identify the homomorphism

of rows $V_1$ into left $A$-vector space of rows $V_2$. The homomorphism

(10.3.3) $\quad \overline{f} : V_1 \to V_2$

which has the given matrix $f$ is unique.

Proof. The equality

$$\overline{f} \circ ((v + w)_*{}^*e_{V_2})$$
$$= (v + w)_*f_*{}^*e_{V_2}$$
(10.3.24)
$$= v_*f_*{}^*e_{V_2} + w_*{}^*f_*{}^*e_{V_2}$$
$$= \overline{f} \circ (v_*{}^*e_{V_1}) + \overline{f} \circ (w_*{}^*e_{V_1})$$

follows from equalities

(8.1.17) $\quad (b_1 + b_2)_*{}^*a = b_1{}_*{}^*a + b_2{}_*{}^*a$

(10.2.8) $\quad f \circ (av) = a(f \circ v)$

(10.3.11) $\quad \overline{f} \circ (v_*{}^*e_{V_1}) = v_*f_*{}^*e_{V_2}$

From the equality (10.3.24), it follows that the map $\overline{f}$ is homomorphism of Abelian group. The equality

$$\overline{f} \circ ((av)_*{}^*e_{V_1})$$
$$= (av)_*{}^*f_*{}^*e_{V_2}$$
(10.3.25)
$$= a(v_*f_*{}^*e_{V_2})$$
$$= a(\overline{f} \circ (v_*{}^*e_{V_1}))$$

follows from equalities

(10.2.8) $\quad f \circ (av) = a(f \circ v)$

(10.3.11) $\quad \overline{f} \circ (v_*{}^*e_{V_1}) = v_*f_*{}^*e_{V_2}$

From the equality (10.3.25), and definitions 2.8.5, 10.3.1, it follows that the map

(10.3.3) $\quad \overline{f} : V_1 \to V_2$

is homomorphism of $A$-vector space of rows $V_1$ into $A$-vector space of rows $V_2$.

Let $f$ be matrix of homomorphisms $\overline{f}$, $\overline{g}$ relative to bases $\overline{e}_V$, $\overline{e}_W$. The equality

$\overline{f} \circ (v_*{}^*e_V) = v_*f_*{}^*e_W = \overline{g} \circ (v_*{}^*e_V)$

follows from the theorem 10.3.7. Therefore, $\overline{f} = \overline{g}$. □

On the basis of theorems 10.3.7, 10.3.9 we identify the homomorphism

| $(10.3.3)$   $\overline{f} : V_1 \to V_2$ |
|---|

of left $A$-vector space $V$ of columns and coordinates of its presentation

| $(10.3.8)$   $\overline{f} \circ (v^* {}_* e_{V_1}) = v^* {}_* f^* {}_* e_{V_2}$ |
|---|

| $(10.3.3)$   $\overline{f} : V_1 \to V_2$ |
|---|

of left $A$-vector space $V$ of rows and coordinates of its presentation

| $(10.3.11)$   $\overline{f} \circ (v_* {}^* e_{V_1}) = v_* {}^* f_* {}^* e_{V_2}$ |
|---|

## 10.4. Sum of Homomorphisms of $A$-Vector Space

THEOREM 10.4.1. *Let $V$, $W$ be left $A$-vector spaces and maps $g : V \to W$, $h : V \to W$ be homomorphisms of left $A$-vector space. Let the map $f : V \to W$ be defined by the equality*

$$(10.4.1) \qquad \forall v \in V : f \circ v = g \circ v + h \circ v$$

*The map $f$ is homomorphism of left $A$-vector space and is called* **sum**

$$f = g + h$$

**of homomorphisms** $g$ *and* $h$.

PROOF. According to the theorem 7.2.6 (the equality (7.2.8)), the equality

$$(10.4.2) \qquad \begin{aligned} & f \circ (u + v) \\ = {} & g \circ (u + v) + h \circ (u + v) \\ = {} & g \circ u + g \circ v \\ & + h \circ u + h \circ v \\ = {} & f \circ u + f \circ v \end{aligned}$$

follows from the equality

| $(10.3.4)$   $f \circ (u + v) = f \circ u + f \circ v$ |
|---|

and the equality (10.4.1). According to the theorem 7.2.6 (the equality (7.2.10)), the equality

$$(10.4.3) \qquad \begin{aligned} & f \circ (av) \\ = {} & g \circ (av) + h \circ (av) \\ = {} & a(g \circ v) + a(h \circ v) \\ = {} & a(g \circ v + h \circ v) \\ = {} & a(f \circ v) \end{aligned}$$

follows from the equality

| $(10.3.5)$   $f \circ (av) = a(f \circ v)$ |
|---|

and the equality (10.4.1).

The theorem follows from the theorem 10.3.2 and equalities (10.4.2), (10.4.3).
□

THEOREM 10.4.2. *Let $V$, $W$ be left $A$-vector spaces of columns. Let $\overline{\overline{e}}_V$ be a basis of left $A$-vector space $V$. Let $\overline{\overline{e}}_W$ be a basis of left $A$-vector space $W$. The homomorphism $\overline{f} : V \to W$ is sum of ho-momorphism*

THEOREM 10.4.3. *Let $V$, $W$ be left $A$-vector spaces of rows. Let $\overline{\overline{e}}_V$ be a basis of left $A$-vector space $V$. Let $\overline{\overline{e}}_W$ be a basis of left $A$-vector space $W$. The ho-momorphism $\overline{f} : V \to W$ is sum of ho-* *momorphisms $\overline{g} : V \to W$, $\overline{h} : V \to W$ iff matrix $f$ of homomorphism $\overline{f}$ relative to bases $\overline{\overline{e}}_V$, $\overline{\overline{e}}_W$ is equal to the sum of the matrix $g$ of the homomorphism $\overline{g}$ relative to bases $\overline{\overline{e}}_V$, $\overline{\overline{e}}_W$ and the matrix $h$ of the homomorphism $\overline{h}$ relative to bases $\overline{\overline{e}}_V$, $\overline{\overline{e}}_W$.*

PROOF. According to theorem 10.3.6

$$(10.4.4) \qquad \overline{f} \circ (a^*{}_*e_U) = a^*{}_* f^*{}_* e_V$$

$$(10.4.5) \qquad \overline{g} \circ (a^*{}_*e_U) = a^*{}_* g^*{}_* e_V$$

$$(10.4.6) \qquad \overline{h} \circ (a^*{}_*e_U) = a^*{}_* h^*{}_* e_V$$

According to the theorem 10.4.1, the equality

$$
\begin{aligned}
(10.4.7) \quad & \overline{f} \circ (v^*{}_*e_U) \\
&= \overline{g} \circ (v^*{}_*e_U) + \overline{h} \circ (v^*{}_*e_U) \\
&= v^*{}_* g^*{}_* e_V + v^*{}_* h^*{}_* e_V \\
&= v^*{}_* (g+h)^*{}_* e_V
\end{aligned}
$$

follows from equalities (8.1.18), (8.1.19), (10.4.5), (10.4.6). The equality

$$f = g + h$$

follows from equalities (10.4.4), (10.4.7) and the theorem 10.3.6.

Let

$$f = g + h$$

According to the theorem 10.3.8, there exist homomorphisms $\overline{f} : U \to V$, $\overline{g} : U \to V$, $\overline{h} : U \to V$ such that

$$
\begin{aligned}
(10.4.8) \quad & \overline{f} \circ (v^*{}_*e_U) \\
&= v^*{}_* f^*{}_* e_V \\
&= v^*{}_* (g+h)^*{}_* e_V \\
&= v^*{}_* g^*{}_* e_V + v^*{}_* h^*{}_* e_V \\
&= \overline{g} \circ (v^*{}_*e_U) + \overline{h} \circ (v^*{}_*e_U)
\end{aligned}
$$

From the equality (10.4.8) and the theorem 10.4.1, it follows that the homomorphism $\overline{f}$ is sum of homomorphisms $\overline{g}$, $\overline{h}$. ☐

*momorphisms $\overline{g} : V \to W$, $\overline{h} : V \to W$ iff matrix $f$ of homomorphism $\overline{f}$ relative to bases $\overline{\overline{e}}_V$, $\overline{\overline{e}}_W$ is equal to the sum of the matrix $g$ of the homomorphism $\overline{g}$ relative to bases $\overline{\overline{e}}_V$, $\overline{\overline{e}}_W$ and the matrix $h$ of the homomorphism $\overline{h}$ relative to bases $\overline{\overline{e}}_V$, $\overline{\overline{e}}_W$.*

PROOF. According to theorem 10.3.7

$$(10.4.9) \qquad \overline{f} \circ (a_*{}^*e_U) = a_*{}^* f_*{}^* e_V$$

$$(10.4.10) \qquad \overline{g} \circ (a_*{}^*e_U) = a_*{}^* g_*{}^* e_V$$

$$(10.4.11) \qquad \overline{h} \circ (a_*{}^*e_U) = a_*{}^* h_*{}^* e_V$$

According to the theorem 10.4.1, the equality

$$
\begin{aligned}
(10.4.12) \quad & \overline{f} \circ (v_*{}^*e_U) \\
&= \overline{g} \circ (v_*{}^*e_U) + \overline{h} \circ (v_*{}^*e_U) \\
&= v_*{}^* g_*{}^* e_V + v_*{}^* h_*{}^* e_V \\
&= v_*{}^* (g+h)_*{}^* e_V
\end{aligned}
$$

follows from equalities (8.1.16), (8.1.17), (10.4.10), (10.4.11). The equality

$$f = g + h$$

follows from equalities (10.4.9), (10.4.12) and the theorem 10.3.7.

Let

$$f = g + h$$

According to the theorem 10.3.9, there exist homomorphisms $\overline{f} : U \to V$, $\overline{g} : U \to V$, $\overline{h} : U \to V$ such that

$$
\begin{aligned}
(10.4.13) \quad & \overline{f} \circ (v_*{}^*e_U) \\
&= v_*{}^* f_*{}^* e_V \\
&= v_*{}^* (g+h)_*{}^* e_V \\
&= v_*{}^* g_*{}^* e_V + v_*{}^* h_*{}^* e_V \\
&= \overline{g} \circ (v_*{}^*e_U) + \overline{h} \circ (v^*{}_*e_U)
\end{aligned}
$$

From the equality (10.4.13) and the theorem 10.4.1, it follows that the homomorphism $\overline{f}$ is sum of homomorphisms $\overline{g}$, $\overline{h}$. ☐

## 10.5. Product of Homomorphisms of $A$-Vector Space

THEOREM 10.5.1. *Let $U$, $V$, $W$ be left $A$-vector spaces. Let diagram of maps*

(10.5.1)

$$U \xrightarrow{\ f\ } W$$
$$\searrow g \qquad h \nearrow$$
$$V$$

(10.5.2)
$$\forall u \in U : f \circ u = (h \circ g) \circ u$$
$$= h \circ (g \circ u)$$

*be commutative diagram where maps $g$, $h$ are homomorphisms of left $A$-vector space. The map $f = h \circ g$ is homomorphism of left $A$-vector space and is called **product of homomorphisms** $h$, $g$.*

PROOF. The equality (10.5.2) follows from the commutativity of the diagram (10.5.1). The equality

(10.5.3)
$$f \circ (u + v)$$
$$= h \circ (g \circ (u + v))$$
$$= h \circ (g \circ u + g \circ v))$$
$$= h \circ (g \circ u) + h \circ (g \circ v)$$
$$= f \circ u + f \circ v$$

follows from the equality

$$\boxed{(10.3.4) \quad f \circ (u + v) = f \circ u + f \circ v}$$

and the equality (10.5.2). The equality

(10.5.4)
$$f \circ (av) = h \circ (g \circ (av))$$
$$= h \circ (a(g \circ v))$$
$$= a(h \circ (g \circ v))$$
$$= a(f \circ v)$$

follows from the equality

$$\boxed{(10.3.5) \quad f \circ (av) = a(f \circ v)}$$

and the equality (10.4.1). The theorem follows from the theorem 10.3.2 and equalities (10.5.3), (10.5.4). $\qquad\square$

THEOREM 10.5.2. *Let $U$, $V$, $W$ be left $A$-vector spaces of columns. Let $\overline{\overline{e}}_U$ be a basis of left $A$-vector space $U$. Let $\overline{\overline{e}}_V$ be a basis of left $A$-vector space $V$. Let $\overline{\overline{e}}_W$ be a basis of left $A$-vector space $W$. The homomorphism $\overline{f} : U \to W$ is product of homomorphisms $\overline{h} : V \to W$, $\overline{g} : U \to V$ iff matrix $f$ of homomorphism $\overline{f}$ relative to bases $\overline{\overline{e}}_U$, $\overline{\overline{e}}_W$ is equal to the $^*{}_*$-product of matrix $g$ of the*

THEOREM 10.5.3. *Let $U$, $V$, $W$ be left $A$-vector spaces of rows. Let $\overline{\overline{e}}_U$ be a basis of left $A$-vector space $U$. Let $\overline{\overline{e}}_V$ be a basis of left $A$-vector space $V$. Let $\overline{\overline{e}}_W$ be a basis of left $A$-vector space $W$. The homomorphism $\overline{f} : U \to W$ is product of homomorphisms $\overline{h} : V \to W$, $\overline{g} : U \to V$ iff matrix $f$ of homomorphism $\overline{f}$ relative to bases $\overline{\overline{e}}_U$, $\overline{\overline{e}}_W$ is equal to the $_*{}^*$-product of matrix $h$ of the* *homomorphism $\overline{g}$ relative to bases $\overline{\overline{e}}_U$, $\overline{\overline{e}}_V$ over matrix $h$ of the homomorphism $\overline{h}$ relative to bases $\overline{\overline{e}}_V$, $\overline{\overline{e}}_W$*

(10.5.5)    $f = g^*{}_* h$

Proof. According to the theorem 10.3.6

(10.5.7)   $\overline{f} \circ (u^*{}_* e_U) = u^*{}_* f^*{}_* e_W$

(10.5.8)   $\overline{g} \circ (u^*{}_* e_U) = u^*{}_* g^*{}_* e_V$

(10.5.9)   $\overline{h} \circ (v^*{}_* e_V) = v^*{}_* h^*{}_* e_W$

According to the theorem 10.5.1, the equality

$$
\begin{aligned}
(10.5.10) \quad & \overline{f} \circ (u^*{}_* e_U) \\
&= \overline{h} \circ (\overline{g} \circ (u^*{}_* e_U)) \\
&= \overline{h} \circ (u^*{}_* g^*{}_* e_V) \\
&= u^*{}_* g^*{}_* h^*{}_* e_W
\end{aligned}
$$

follows from equalities (10.5.8), (10.5.9). The equality (10.5.5) follows from equalities (10.5.7), (10.5.10) and the theorem 10.3.6.

Let $f = g^*{}_* h$. According to the theorem 10.3.8, there exist homomorphisms $\overline{f} : U \to V$, $\overline{g} : U \to V$, $\overline{h} : U \to V$ such that

$$
\begin{aligned}
& \overline{f} \circ (u^*{}_* e_U) = u^*{}_* f^*{}_* e_W \\
&= u^*{}_* g^*{}_* h^*{}_* e_W \\
(10.5.11) \quad &= \overline{h} \circ (u^*{}_* g^*{}_* e_V) \\
&= \overline{h} \circ (\overline{g} \circ (u^*{}_* e_U)) \\
&= (\overline{h} \circ \overline{g}) \circ (u^*{}_* e_U)
\end{aligned}
$$

From the equality (10.5.11) and the theorem 10.4.1, it follows that the homomorphism $\overline{f}$ is product of homomorphisms $\overline{h}$, $\overline{g}$. $\qquad\square$

*homomorphism $\overline{h}$ relative to bases $\overline{\overline{e}}_V$, $\overline{\overline{e}}_W$ over matrix $g$ of the homomorphism $\overline{g}$ relative to bases $\overline{\overline{e}}_U$, $\overline{\overline{e}}_V$*

(10.5.6)    $f = h_*{}^* g$

Proof. According to the theorem 10.3.7

(10.5.12)   $\overline{f} \circ (u_*{}^* e_U) = u_*{}^* f_*{}^* e_W$

(10.5.13)   $\overline{g} \circ (u_*{}^* e_U) = u_*{}^* g_*{}^* e_V$

(10.5.14)   $\overline{h} \circ (v_*{}^* e_V) = v_*{}^* h_*{}^* e_W$

According to the theorem 10.5.1, the equality

$$
\begin{aligned}
(10.5.15) \quad & \overline{f} \circ (u_*{}^* e_U) \\
&= \overline{h} \circ (\overline{g} \circ (u_*{}^* e_U)) \\
&= \overline{h} \circ (u_*{}^* g_*{}^* e_V) \\
&= u_*{}^* g_*{}^* h_*{}^* e_W
\end{aligned}
$$

follows from equalities (10.5.13), (10.5.14). The equality (10.5.6) follows from equalities (10.5.12), (10.5.15) and the theorem 10.3.7.

Let $f = h_*{}^* g$. According to the theorem 10.3.9, there exist homomorphisms $\overline{f} : U \to V$, $\overline{g} : U \to V$, $\overline{h} : U \to V$ such that

$$
\begin{aligned}
& \overline{f} \circ (u_*{}^* e_U) = u_*{}^* f_*{}^* e_W \\
&= u_*{}^* g_*{}^* h_*{}^* e_W \\
(10.5.16) \quad &= \overline{h} \circ (u_*{}^* g_*{}^* e_V) \\
&= \overline{h} \circ (\overline{g} \circ (u_*{}^* e_U)) \\
&= (\overline{h} \circ \overline{g}) \circ (u_*{}^* e_U)
\end{aligned}
$$

From the equality (10.5.16) and the theorem 10.4.1, it follows that the homomorphism $\overline{f}$ is product of homomorphisms $\overline{h}$, $\overline{g}$. $\qquad\square$

## 10.6. Direct Sum of Left $A$-Vector Space

Let Left $A$-Vector be category of left $A$-vector spaces morphisms of which are homomorphisms of $A$-vector space.

Definition 10.6.1. *Coproduct in category* Left $A$-Vector *is called* **direct sum**.[10.10] *We will use notation* $V \oplus W$ *for direct sum of left $A$-vector spaces $V$ and $W$.* $\qquad\square$

THEOREM 10.6.2. *In category* LEFT $A$-VECTOR *there exists direct sum of left $A$-vector spaces.* [10.11] *Let $\{V_i, i \in I\}$ be set of left $A$-vector spaces. Then the representation*

$$A \xrightarrow{\quad\quad} \bigoplus_{i \in I} V_i \qquad a\left(\bigoplus_{i \in I} v_i\right) = \bigoplus_{i \in I} a v_i$$

*of the $D$-algebra $A$ in direct sum of Abelian groups*

$$V = \bigoplus_{i \in I} V_i$$

*is direct sum of left $A$-vector spaces*

$$V = \bigoplus_{i \in I} V_i$$

PROOF. Let $V_i$, $i \in I$, be set of left $A$-vector spaces. Let

$$V = \bigoplus_{i \in I} V_i$$

be direct sum of Abelian groups. Let $v = (v_i, i \in I) \in V$. We define action of $A$-number $a$ over $V$-number $v$ according to the equality

$$(10.6.1) \qquad a(v_i, i \in I) = (a v_i, i \in I)$$

According to the theorem 7.2.6,

$$(10.6.2) \qquad (ab)v_i = a(b v_i)$$

$$(10.6.3) \qquad a(v_i + w_i) = a v_i + a w_i$$

$$(10.6.4) \qquad (a + b)v_i = a v_i + b v_i$$

$$(10.6.5) \qquad 1 v_i = v_i$$

for any $a, b \in A$, $v_i, w_i \in V_i$, $i \in I$. Equalities

$$(10.6.6) \quad \begin{aligned} (ab)v &= (ab)(v_i, i \in I) = ((ab)v_i, i \in I) \\ &= (a(b v_i), i \in I) = a(b v_i, i \in I) = a(b(v_i, i \in I)) \\ &= a(bv) \end{aligned}$$

$$(10.6.7) \quad \begin{aligned} a(v + w) &= a((v_i, i \in I) + (w_i, i \in I)) = a(v_i + w_i, i \in I) \\ &= (a(v_i + w_i), i \in I) = (a v_i + a w_i, i \in I) \\ &= (a v_i, i \in I) + (a w_i, i \in I) = a(v_i, i \in I) + a(w_i, i \in I) \\ &= av + aw \end{aligned}$$

$$(10.6.8) \quad \begin{aligned} (a + b)v &= (a + b)(v_i, i \in I) = ((a + b)v_i, i \in I) \\ &= (a v_i + b v_i, i \in I) = (a v_i, i \in I) + (b v_i, i \in I) \\ &= a(v_i, i \in I) + b(v_i, i \in I) = av + bv \end{aligned}$$

$$1v = 1(v_i, i \in I) = (1 v_i, i \in I) = (v_i, i \in I) = v$$

---

[10.10]   See also the definition of direct sum of modules in [1], page 128. On the same page, Lang proves the existence of direct sum of modules.
[10.11]   See also proposition on the page [1]-1.1

follow from equalities (3.7.4), (10.6.1), (10.6.2), (10.6.3), (10.6.4), (10.6.5) for any $a$, $b \in A$, $v = (v_i, i \in I)$, $w = (w_i, i \in I) \in V$.

From the equality (10.6.7), it follows that the map

$$a : v \to av$$

defined by the equality (10.6.1) is endomorphism of Abelian group $V$. From equalities (10.6.6), (10.6.8), it follows that the map (10.6.1) is homomorphism of $D$-algebra $A$ into set of endomorphisms of Abelian group $V$. Therefore, the map (10.6.1) is representation of $D$-algebra $A$ in Abelian group $V$.

Let $a_1$, $a_2$ be $A$-numbers such that

$$(10.6.9) \qquad a_1(v_i, i \in I) = a_2(v_i, i \in I)$$

for any $V$-number  $v = (v_i, i \in I)$.  The equality

$$(10.6.10) \qquad (a_1 v_i, i \in I) = (a_2 v_i, i \in I)$$

follows from the equality (10.6.9). From the equality (10.6.10), it follows that

$$(10.6.11) \qquad \forall i \in I : a_1 v_i = a_2 v_i$$

Since for any $i \in I$, $V_i$ is $A$-vector space, then corresponding representation is effective according to the definition 7.2.1 Therefore, the equality $a_1 = a_2$ follows from the statement (10.6.11) and the map (10.6.1) is effective representation of the $D$-algebra $A$ in Abelian group $V$. Therefore, Abelian group $V$ is left $A$-vector space.

Let

$$\{f_i : V_i \to W, i \in I\}$$

be set of linear maps into left $A$-vector space $W$. We define the map

$$f : V \to W$$

by the equality

$$(10.6.12) \qquad f\left(\bigoplus_{i \in I} x_i\right) = \sum_{i \in I} f_i(x_i)$$

The sum in the right side of the equality (10.6.12) is finite, since all summands, except for a finite number, equal 0. From the equality

$$f_i(x_i + y_i) = f_i(x_i) + f_i(y_i)$$

and the equality (10.6.12), it follows that

$$\begin{aligned} f(\bigoplus_{i \in I}(x_i + y_i)) &= \sum_{i \in I} f_i(x_i + y_i) \\ &= \sum_{i \in I}(f_i(x_i) + f_i(y_i)) \\ &= \sum_{i \in I} f_i(x_i) + \sum_{i \in I} f_i(y_i) \\ &= f(\bigoplus_{i \in I} x_i) + f(\bigoplus_{i \in I} y_i) \end{aligned}$$

From the equality

$$f_i(dx_i) = df_i(x_i)$$

and the equality (10.6.12), it follows that

$$f((dx_i, i \in I)) = \sum_{i \in I} f_i(dx_i) = \sum_{i \in I} df_i(x_i) = d \sum_{i \in I} f_i(x_i)$$
$$= df((x_i, i \in I))$$

Therefore, the map $f$ is linear map. The equality

$$f \circ \lambda_j(x) = \sum_{i \in I} f_i(\delta^i_j x) = f_j(x)$$

follows from equalities (3.7.5), (10.6.12). Since the map $\lambda_i$ is injective, then the map $f$ is unique. Therefore, the theorem folows from definitions 2.2.12, 3.7.2. □

DEFINITION 10.6.3. *Let $D$ be commutative ring with unit. Let $A$ be associative division $D$-algebra. For any integer $n > 0$, we introduce left vector space $\mathcal{L}^n A$ using definition by induction*

$$\mathcal{L}^1 A = \mathcal{L} A$$
$$\mathcal{L}^k A = \mathcal{L}^{k-1} A \oplus \mathcal{L} A$$

□

DEFINITION 10.6.4. *Let $D$ be commutative ring with unit. Let $A$ be associative division $D$-algebra. For any set $B$ and any integer $n > 0$, we introduce left vector space $\mathcal{L}^n A^B$ using definition by induction*

$$\mathcal{L}^1 A^B = \mathcal{L} A^B$$
$$\mathcal{L}^k A^B = \mathcal{L}^{k-1} A^B \oplus \mathcal{L} A^B$$

□

THEOREM 10.6.5. *Direct sum of left $A$-vector spaces $V_1, ..., V_n$ coincides with their Cartesian product*

$$V_1 \oplus ... \oplus V_n = V_1 \times ... \times V_n$$

PROOF. The theorem follows from the theorem 10.6.2. □

THEOREM 10.6.6. *Let $\overline{e}_1 = (e_{1i}, i \in I)$ be a basis of left $A$-vector space of columns $V$. Then we can represent $A$-module $V$ as direct sum*

$$(10.6.13) \qquad V = \bigoplus_{i \in I} A e_i$$

PROOF. The map

$$f : v^i e_i \in V \to (v^i, i \in I) \in \bigoplus_{i \in I} A e_i$$

is isomorphism. □

THEOREM 10.6.7. *Let $\overline{e}_1 = (e^i_1, i \in I)$ be a basis of left $A$-vector space of rows $V$. Then we can represent $A$-module $V$ as direct sum*

$$(10.6.14) \qquad V = \bigoplus_{i \in I} A e^i$$

PROOF. The map

$$f : v_i e^i \in V \to (v_i, i \in I) \in \bigoplus_{i \in I} A e^i$$

is isomorphism. □

# Homomorphism of Right $A$-Vector Space

## 11.1. General Definition

Let $A_i$, $i = 1$, $2$, be division algebra over commutative ring $D_i$. Let $V_i$, $i = 1$, $2$, be right $A_i$-vector space.

DEFINITION 11.1.1. *Let diagram of representations*

$$(11.1.1)$$

$$g_{1.12}(d) : a \to d\,a$$
$$g_{1.23}(v) : w \to C_1(w, v)$$
$$C_1 \in \mathcal{L}(A_1^2 \to A_1)$$
$$g_{1.34}(a) : v \to v\,a$$
$$g_{1.14}(d) : v \to v\,d$$

*describe right $A_1$-vector space $V_1$. Let diagram of representations*

$$(11.1.2)$$

$$g_{2.12}(d) : a \to d\,a$$
$$g_{2.23}(v) : w \to C_2(w, v)$$
$$C_2 \in \mathcal{L}(A_2^2 \to A_2)$$
$$g_{2.34}(a) : v \to v\,a$$
$$g_{2.14}(d) : v \to v\,d$$

*describe right $A_2$-vector space $V_2$. Morphism*

$$(11.1.3) \qquad (h : D_1 \to D_2 \quad g : A_1 \to A_2 \quad f : V_1 \to V_2)$$

*of diagram of representations* (11.1.1) *into diagram of representations* (11.1.2) *is called* **homomorphism** *of right $A_1$-vector space $V_1$ into right $A_2$-vector space $V_2$. Let us denote* $\mathrm{Hom}(D_1 \to D_2; A_1 \to A_2; V_1* \to V_2*)$ *set of homomorphisms of right $A_1$-vector space $V_1$ into right $A_2$-vector space $V_2$.* □

We will use notation

$$f \circ a = f(a)$$

for image of homomorphism $f$.

THEOREM 11.1.2. *The homomorphism*

$$(11.1.3) \qquad (h : D_1 \to D_2 \quad g : A_1 \to A_2 \quad f : V_1 \to V_2)$$

*of right $A_1$-vector space $V_1$ into right $A_2$-vector space $V_2$ satisfies following equalities*

$$(11.1.4) \qquad h(p + q) = h(p) + h(q)$$

$$(11.1.5) \qquad h(pq) = h(p)h(q)$$

$$(11.1.6) \qquad g \circ (a + b) = g \circ a + g \circ b$$

$$(11.1.7) \qquad g \circ (ab) = (g \circ a)(g \circ b)$$

(11.1.8)                         $g \circ (ap) = (g \circ a)h(p)$

(11.1.9)                         $f \circ (u + v) = f \circ u + f \circ v$

(11.1.10)                        $f \circ (va) = (f \circ v)(g \circ a)$

$$p, q \in D_1 \quad a, b \in A_1 \quad u, v \in V_1$$

PROOF.    According to definitions    2.3.9, 2.8.5, 7.4.1, 11 the map

$$(h : D_1 \to D_2 \quad g : A_1 \to A_2)$$

is homomorphism of $D_1$-algebra $A_1$ into $D_2$-algebra $A_2$.    Therefore, equalities
(11.1.4), (11.1.5),    (11.1.6), (11.1.7), (11.1.8)      follow from the theorem 4.3.2.
The equality (11.1.9) follows from the definition 11.1.1, since, accorfing to the def-
inition 2.3.9, the map $f$ is homomorpism of Abelian group. The equality (11.1.10)
follows from the equality (2.3.6) because the map

$$(g : A_1 \to A_2 \quad f : V_1 \to V_2)$$

is morphism of representation $g_{1.34}$ into representation $g_{2.34}$.                         □

DEFINITION 11.1.3.  *Homomorphism* [11.1]

$$(h : D_1 \to D_2 \quad g : A_1 \to A_2 \quad f : V_1 \to V_2)$$

*is called* **isomorphism** *between right $A_1$-vector space $V_1$ and right $A_2$-vector space*
*$V_2$, if there exists the map*

$$(h^{-1} : D_2 \to D_1 \quad g^{-1} : A_2 \to A_1 \quad f^{-1} : V_2 \to V_1)$$

*which is homomorphism.*                                                □

## 11.1.1. Homomorphism of Right $A$-Vector Space of columns.

THEOREM 11.1.4.  *The homomorphism* [11.2]

$$\boxed{(11.1.3) \quad (h : D_1 \to D_2 \quad \overline{g} : A_1 \to A_2 \quad \overline{f} : V_1 \to V_2)}$$

*of right $A_1$-vector space of columns $V_1$ into right $A_2$-vector space of columns $V_2$ has*
*presentation*

(11.1.11)                        $b = gh(a)$

(11.1.12)                        $\overline{g} \circ (a^i e_{A_1 i}) = h(a^i)g_i^k e_{A_2 k}$

(11.1.13)                        $\overline{g} \circ (e_{A_1} a) = e_{A_2} gh(a)$

(11.1.14)                        $w = f_*{}^*(\overline{g} \circ v)$

(11.1.15)                        $\overline{f} \circ (e_{V_1 i} v^i) = e_{V_2 k} f_i^k (\overline{g} \circ v^i)$

---

[11.1]  I follow the definition on page [18]-49.

[11.2]  In theorems 11.1.4, 11.1.5, we use the following convention.    Let $i = 1, 2$.    Let the
set of vectors $\overline{\overline{e}}_{A_i} = (e_{A_i k}, k \in K_i)$ be a basis and $C_{i \cdot ij}^{\cdot \cdot k}$, $k$, $i$, $j \in K_i$, be structural constants
of $D_i$-algebra of columns $A_i$.    Let the set of vectors $\overline{\overline{e}}_{V_i} = (e_{V_i i}, i \in I_i)$ be a basis of right $A_i$-
vector space $V_i$.

(11.1.16) $$\overline{f} \circ (e_{V_1 *}{}^* v) = e_{V_2 *}{}^* f_*{}^* (\overline{g} \circ v)$$

*relative to selected bases. Here*

- *a is coordinate matrix of $A_1$-number $\overline{a}$ relative the basis $\overline{\overline{e}}_{A_1}$*

$$\overline{a} = e_{A_1} a$$

- *$h(a) = (h(a_k), k \in K)$ is a matrix of $D_2$-numbers.*
- *b is coordinate matrix of $A_2$-number*

$$\overline{b} = \overline{g} \circ \overline{a}$$

  *relative the basis $\overline{\overline{e}}_{A_2}$*

$$\overline{b} = e_{A_2} b$$

- *g is coordinate matrix of set of $A_2$-numbers $(\overline{g} \circ e_{A_1 k}, k \in K)$ relative the basis $\overline{\overline{e}}_{A_2}$.*
- *There is relation between the matrix of homomorphism and structural constants*

(11.1.17) $$h(C_{1ij}{}^k) g_k^l = g_i^p g_j^q C_{2pq}{}^l$$

- *v is coordinate matrix of $V_1$-number $\overline{v}$ relative the basis $\overline{\overline{e}}_{V_1}$*

$$\overline{v} = e_{V_1 *}{}^* v$$

- *$g(v) = (g(v_i), i \in I)$ is a matrix of $A_2$-numbers.*
- *w is coordinate matrix of $V_2$-number*

$$\overline{w} = \overline{f} \circ \overline{v}$$

  *relative the basis $\overline{\overline{e}}_{V_2}$*

$$\overline{w} = e_{V_2 *}{}^* w$$

- *f is coordinate matrix of set of $V_2$-numbers $(\overline{f} \circ e_{V_1 i}, i \in I)$ relative the basis $\overline{\overline{e}}_{V_2}$.*

*For given homomorphism*

> (11.1.3)   $(h : D_1 \to D_2 \quad \overline{g} : A_1 \to A_2 \quad \overline{f} : V_1 \to V_2)$

*the set of matrices $(g, f)$ is unique and is called* **coordinates of homomorphism** *(11.1.3) relative bases $(\overline{\overline{e}}_{A_1}, \overline{\overline{e}}_{V_1})$,   $(\overline{\overline{e}}_{A_2}, \overline{\overline{e}}_{V_2})$.*

PROOF.    According to definitions   2.3.9, 2.8.5,   5.4.1, 11.1.1,   the map

$$(h : D_1 \to D_2 \quad \overline{g} : A_1 \to A_2)$$

is homomorphism of $D_1$-algebra $A_1$ into $D_2$-algebra $A_2$. Therefore, equalities (11.1.11), (11.1.12), (11.1.13), (11.1.17) follow from equalities

> (5.4.11)   $b = f h(a)$

> (5.4.12)   $\overline{f} \circ (a^i e_{1i}) = h(a^i) f_i^k e_{2k}$

> (5.4.13)   $\overline{f} \circ (e_1 a) = e_2 f h(a)$

> (5.4.16)   $h(C_{1ij}{}^k) f_k^l = f_i^p f_j^q C_{2pq}{}^l$

From the theorem 4.3.3, it follows that the matrix $g$ is unique.

Vector $\overline{v} \in V_1$ has expansion

(11.1.18) $$\overline{v} = e_{V_1 *}{}^* v$$

relative to the basis $\overline{\overline{e}}_{V_1}$. Vector $\overline{w} \in V_2$ has expansion

$$(11.1.19) \qquad \overline{w} = e_{V_2*}{}^* w$$

relative to the basis $\overline{\overline{e}}_{V_2}$. Since $\overline{f}$ is a homomorphism, then the equality

$$(11.1.20) \qquad \overline{w} = \overline{f} \circ \overline{v} = \overline{f} \circ (e_{V_1*}{}^* v)$$
$$= (\overline{f} \circ e_{V_1})_*{}^* (\overline{g} \circ v)$$

follows from equalities (11.1.9), (11.1.10). $V_2$-number $\overline{f} \circ e_{V_1 i}$ has expansion

$$(11.1.21) \qquad \overline{f} \circ e_{V_1 i} = e_{V_2*}{}^* f_i = e_{V_2}^j f_i^j$$

relative to basis $\overline{\overline{e}}_{V_2}$. The equality

$$(11.1.22) \qquad \overline{w} = e_{V_2*}{}^* f_*{}^* (\overline{g} \circ v)$$

follows from equalities (11.1.20), (11.1.21). The equality (11.1.14) follows from comparison of (11.1.19) and (11.1.22) and the theorem 8.4.18. From the equality (11.1.21) and from the theorem 8.4.18, it follows that the matrix $f$ is unique. $\square$

THEOREM 11.1.5. *Let the map*

$$h : D_1 \to D_2$$

*be homomorphism of the ring $D_1$ into the ring $D_2$.* *Let*

$$g = (g_l^k, k \in K, l \in L)$$

*be matrix of $D_2$-numbers which satisfy the equality*

$$(11.1.17) \quad h(C_{1\,ij}^{\ k}) g_k^l = g_i^p g_j^q C_{2\,pq}^{\ l}$$

*Let*

$$f = (f_j^i, i \in I, j \in J)$$

*be matrix of $A_2$-numbers.* *The map* [11.2]

$$(11.1.3) \quad (h : D_1 \to D_2 \quad \overline{g} : A_1 \to A_2 \quad \overline{f} : V_1 \to V_2)$$

*defined by equalities*

$$(11.1.13) \quad \overline{g} \circ (e_{A_1} a) = e_{A_2} g h(a)$$

$$(11.1.16) \quad \overline{f} \circ (e_{V_1*}{}^* v) = e_{V_2*}{}^* f_*{}^* (\overline{g} \circ v)$$

*is homomorphism of right $A_1$-vector space of columns $V_1$ into right $A_2$-vector space of columns $V_2$.* *The homomorphism*

$$(11.1.3) \quad (h : D_1 \to D_2 \quad \overline{g} : A_1 \to A_2 \quad \overline{f} : V_1 \to V_2)$$

*which has the given set of matrices $(g, f)$ is unique.*

PROOF. According to the theorem 5.4.5, the map

$$(11.1.23) \qquad (h : D_1 \to D_2 \quad \overline{g} : A_1 \to A_2)$$

is homomorphism of $D_1$-algebra $A_1$ into $D_2$-algebra $A_2$ and the homomorphism (11.1.23) is unique.

The equality

$$(11.1.24) \quad \begin{aligned} &\overline{f} \circ (e_{V_1*}{}^*(v+w)) \\ &= e_{V_2*}{}^* f_*{}^* (g \circ (v+w)) \\ &= e_{V_2*}{}^* f_*{}^* (g \circ v) \\ &+ e_{V_2*}{}^* f_*{}^* (g \circ w) \\ &= \overline{f} \circ (e_{V_1*}{}^* v) + \overline{f} \circ (e_{V_1*}{}^* w) \end{aligned}$$

follows from equalities

| | |
|---|---|
| (8.1.16)  $a_*{}^*(b_1 + b_2) = a_*{}^* b_1 + a_*{}^* b_2$ | (11.1.6)  $g \circ (a+b) = g \circ a + g \circ b$ |
| (11.1.10)  $f \circ (va) = (f \circ v)(g \circ a)$ | (11.1.16)  $\overline{f} \circ (e_{V_1*}{}^* v) = e_{V_2*}{}^* f_*{}^* (\overline{g} \circ v)$ |

From the equality (11.1.24), it follows that the map $\overline{f}$ is homomorphism of Abelian group. The equality

$$(11.1.25) \quad \begin{aligned} &\overline{f} \circ (e_{V_1*}{}^*(va)) \\ &= e_{V_2*}{}^* f_*{}^* (\overline{g} \circ (va)) \\ &= (e_{V_2*}{}^* f_*{}^* (\overline{g} \circ v))(\overline{g} \circ a) \\ &= (\overline{f} \circ (e_{V_1*}{}^* v))(\overline{g} \circ a) \end{aligned}$$

follows from equalities

| |
|---|
| (11.1.10)  $f \circ (va) = (f \circ v)(g \circ a)$ |
| (5.4.4)  $g \circ (ab) = (g \circ a)(g \circ b)$ |
| (11.1.16)  $\overline{f} \circ (e_{V_1*}{}^* v) = e_{V_2*}{}^* f_*{}^* (\overline{g} \circ v)$ |

From the equality (11.1.25), and definitions 2.8.5, 11.1.1, it follows that the map

| |
|---|
| (11.1.3)  $(h : D_1 \to D_2 \quad \overline{g} : A_1 \to A_2 \quad \overline{f} : V_1 \to V_2)$ |

is homomorphism of $A_1$-vector space of columns $V_1$ into $A_2$-vector space of columns $V_2$.

Let $f$ be matrix of homomorphisms $\overline{f}$, $\overline{g}$ relative to bases $\overline{\overline{e}}_V$, $\overline{\overline{e}}_W$. The equality

$$\overline{f} \circ (e_{V*}{}^* v) = e_{W*}{}^* f_*{}^* v = \overline{g} \circ (e_{V*}{}^* v)$$

follows from the theorem 11.3.6. Therefore, $\overline{f} = \overline{g}$. $\square$

On the basis of theorems 11.1.4, 11.1.5 we identify the homomorphism

| |
|---|
| (11.1.3)  $(h : D_1 \to D_2 \quad \overline{g} : A_1 \to A_2 \quad \overline{f} : V_1 \to V_2)$ |

of right $A$-vector space $V$ of columns and coordinates of its presentation

| |
|---|
| (11.1.13)  $\overline{g} \circ (e_{A_1} a) = e_{A_2} gh(a)$ |
| (11.1.16)  $\overline{f} \circ (e_{V_1*}{}^* v) = e_{V_2*}{}^* f_*{}^* (\overline{g} \circ v)$ |

## 11.1.2. Homomorphism of Right $A$-Vector Space of rows.

[11.3] In theorems 11.1.6, 11.1.7, we use the following convention. Let $i = 1, 2$. Let the set of vectors $\overline{\overline{e}}_{A_i} = (e_{A_i}^k, k \in K_i)$ be a basis and $C_{ik}^{\cdot j}$, $k, i, j \in K_i$, be structural constants of $D_i$-algebra of rows $A_i$. Let the set of vectors $\overline{\overline{e}}_{V_i} = (e_{V_i}^i, i \in I_i)$ be a basis of right $A_i$-

Theorem 11.1.6. *The homomorphism* [11.3]

$$(11.1.3) \quad (h: D_1 \to D_2 \quad \overline{g}: A_1 \to A_2 \quad \overline{f}: V_1 \to V_2)$$

*of right $A_1$-vector space of rows $V_1$ into right $A_2$-vector space of rows $V_2$ has presentation*

$$(11.1.26) \qquad b = h(a)g$$

$$(11.1.27) \qquad \overline{g} \circ (a_i e_{A_1}^i) = h(a_i)g_k^i e_{A_2}^k$$

$$(11.1.28) \qquad \overline{g} \circ (a e_{A_1}) = h(a)g e_{A_2}$$

$$(11.1.29) \qquad w = f^*_*(\overline{g} \circ v)$$

$$(11.1.30) \qquad \overline{f} \circ (e_{V_1}^i v_i) = e_{V_2}^k f_k^i(\overline{g} \circ v_i)$$

$$(11.1.31) \qquad \overline{f} \circ (e_{V_1}{}^*_* v) = e_{V_2}{}^*_* f^*_*(\overline{g} \circ v)$$

*relative to selected bases. Here*

- *$a$ is coordinate matrix of $A_1$-number $\overline{a}$ relative the basis $\overline{\overline{e}}_{A_1}$*

$$\overline{a} = a e_{A_1}$$

- *$h(a) = (h(a^k), k \in K)$ is a matrix of $D_2$-numbers.*
- *$b$ is coordinate matrix of $A_2$-number*

$$\overline{b} = \overline{g} \circ \overline{a}$$

  *relative the basis $\overline{\overline{e}}_{A_2}$*

$$\overline{b} = b e_{A_2}$$

- *$g$ is coordinate matrix of set of $A_2$-numbers $(\overline{g} \circ e_{A_1}^k, k \in K)$ relative the basis $\overline{\overline{e}}_{A_2}$.*
- *There is relation between the matrix of homomorphism and structural constants*

$$(11.1.32) \qquad h(C_{1k}^{ij})g_l^k = g_p^i g_q^j C_{2l}^{pq}$$

- *$v$ is coordinate matrix of $V_1$-number $\overline{v}$ relative the basis $\overline{\overline{e}}_{V_1}$*

$$\overline{v} = e_{V_1}{}^*_* v$$

- *$g(v) = (g(v^i), i \in I)$ is a matrix of $A_2$-numbers.*
- *$w$ is coordinate matrix of $V_2$-number*

$$\overline{w} = \overline{f} \circ \overline{v}$$

  *relative the basis $\overline{\overline{e}}_{V_2}$*

$$\overline{w} = e_{V_2}{}^*_* w$$

- *$f$ is coordinate matrix of set of $V_2$-numbers $(\overline{f} \circ e_{V_1}^i, i \in I)$ relative the basis $\overline{\overline{e}}_{V_2}$.*

*For given homomorphism*

$$(11.1.3) \quad (h: D_1 \to D_2 \quad \overline{g}: A_1 \to A_2 \quad \overline{f}: V_1 \to V_2)$$

*the set of matrices $(g, f)$ is unique and is called* **coordinates of homomorphism** (11.1.3) *relative bases $(\overline{\overline{e}}_{A_1}, \overline{\overline{e}}_{V_1})$, $(\overline{\overline{e}}_{A_2}, \overline{\overline{e}}_{V_2})$.*

---

vector space $V_i$.

PROOF.    According to definitions    2.3.9, 2.8.5, 5.4.1, 11.1.1,  the map

$$(h : D_1 \to D_2 \quad \overline{g} : A_1 \to A_2)$$

is homomorphism of $D_1$-algebra $A_1$ into $D_2$-algebra $A_2$. Therefore, equalities (11.1.26), (11.1.27), (11.1.28), (11.1.32)  follow from equalities

$$\boxed{(5.4.22) \quad b = h(a)f}$$

$$\boxed{(5.4.23) \quad \overline{f} \circ (a_i e_1^i) = h(a_i) f_k^i e_2^k}$$

$$\boxed{(5.4.24) \quad \overline{f} \circ (a e_1) = h(a) f e_2}$$

$$\boxed{(5.4.27) \quad h(C_1{}_k^{ij}) f_l^k = f_p^i f_q^j C_2{}_l^{pq}}$$

From the theorem 4.3.4, it follows that the matrix $g$ is unique.

Vector  $\overline{v} \in V_1$  has expansion

$$(11.1.33) \qquad\qquad \overline{v} = e_{V_1}{}^*{}_* v$$

relative to the basis  $\overline{\overline{e}}_{V_1}$.  Vector  $\overline{w} \in V_2$  has expansion

$$(11.1.34) \qquad\qquad \overline{w} = e_{V_2}{}^*{}_* w$$

relative to the basis  $\overline{\overline{e}}_{V_2}$.  Since  $\overline{f}$  is a homomorphism, then the equality

$$(11.1.35) \qquad \begin{aligned} \overline{w} &= \overline{f} \circ \overline{v} = \overline{f} \circ (e_{V_1}{}^*{}_* v) \\ &= (\overline{f} \circ e_{V_1})^*{}_* (\overline{g} \circ v) \end{aligned}$$

follows from equalities (11.1.9), (11.1.10). $V_2$-number  $\overline{f} \circ e_{V_1}^i$  has expansion

$$(11.1.36) \qquad \overline{f} \circ e_{V_1}^i = e_{V_2}{}^*{}_* f^i = e_{V_2 j} f_j^i$$

relative to basis  $\overline{\overline{e}}_{V_2}$.  The equality

$$(11.1.37) \qquad\qquad \overline{w} = e_{V_2}{}^*{}_* f^*{}_* (\overline{g} \circ v)$$

follows from equalities (11.1.35), (11.1.36).  The equality (11.1.29) follows from comparison of (11.1.34) and (11.1.37) and the theorem 8.4.25. From the equality (11.1.36) and from the theorem 8.4.25, it follows that the matrix $f$ is unique.    □

THEOREM 11.1.7.  *Let the map*

$$h : D_1 \to D_2$$

*be homomorphism of the ring $D_1$ into the ring $D_2$.*          *Let*

$$g = (g_k^l, k \in K, l \in L)$$

*be matrix of $D_2$-numbers which satisfy the equality*

$$\boxed{(11.1.32) \quad h(C_1{}_k^{ij}) g_l^k = g_p^i g_q^j C_2{}_l^{pq}}$$

*Let*

$$f = (f_i^j, i \in I, j \in J)$$

*be matrix of $A_2$-numbers.          The map* [11.3]

$$\boxed{(11.1.3) \quad (h : D_1 \to D_2 \quad \overline{g} : A_1 \to A_2 \quad \overline{f} : V_1 \to V_2)}$$

*defined by equalities*

$$\boxed{(11.1.28) \quad \overline{g} \circ (a e_{A_1}) = h(a) g e_{A_2}}$$

$$\boxed{(11.1.31) \quad \overline{f} \circ (e_{V_1}{}^*{}_* v) = e_{V_2}{}^*{}_* f^*{}_* (\overline{g} \circ v)}$$

*is homomorphism of right $A_1$-vector space of rows $V_1$ into right $A_2$-vector space of rows $V_2$. The homomorphism*

(11.1.3)     $(h : D_1 \rightarrow D_2 \quad \overline{g} : A_1 \rightarrow A_2 \quad \overline{f} : V_1 \rightarrow V_2)$

*which has the given set of matrices $(g, f)$ is unique.*

PROOF.     According to the theorem 5.4.6, the map

(11.1.38)                    $(h : D_1 \rightarrow D_2 \quad \overline{g} : A_1 \rightarrow A_2)$

is homomorphism of $D_1$-algebra $A_1$ into $D_2$-algebra $A_2$ and the homomorphism (11.1.38) is unique.

The equality

$$
\begin{aligned}
&\overline{f} \circ (e_{V_1}{}^*{}_* (v + w)) \\
&= e_{V_2}{}^*{}_* f^*{}_* (g \circ (v + w)) \\
&= e_{V_2}{}^*{}_* f^*{}_* (g \circ v) \\
&\quad + e_{V_2}{}^*{}_* f^*{}_* (g \circ w) \\
&= \overline{f} \circ (e_{V_1}{}^*{}_* v) + \overline{f} \circ (e_{V_1}{}^*{}_* w)
\end{aligned}
$$

(11.1.39)

follows from equalities

(8.1.18)   $a^*{}_*(b_1 + b_2) = a^*{}_* b_1 + a^*{}_* b_2$          (11.1.6)   $g \circ (a + b) = g \circ a + g \circ b$

(11.1.10)   $f \circ (va) = (f \circ v)(g \circ a)$          (11.1.31)   $\overline{f} \circ (e_{V_1}{}^*{}_* v) = e_{V_2}{}^*{}_* f^*{}_* (\overline{g} \circ v)$

From the equality (11.1.39), it follows that the map $\overline{f}$ is homomorphism of Abelian group. The equality

$$
\begin{aligned}
&\overline{f} \circ (e_{V_1}{}^*{}_* (va)) \\
&= e_{V_2}{}^*{}_* f^*{}_* (\overline{g} \circ (va)) \\
&= (e_{V_2}{}^*{}_* f^*{}_* (\overline{g} \circ v))(\overline{g} \circ a) \\
&= (\overline{f} \circ (e_{V_1}{}^*{}_* v))(\overline{g} \circ a)
\end{aligned}
$$

(11.1.40)

follows from equalities

(11.1.10)   $f \circ (va) = (f \circ v)(g \circ a)$

(5.4.4)   $g \circ (ab) = (g \circ a)(g \circ b)$

(11.1.31)   $\overline{f} \circ (e_{V_1}{}^*{}_* v) = e_{V_2}{}^*{}_* f^*{}_* (\overline{g} \circ v)$

From the equality (11.1.40), and definitions 2.8.5, 11.1.1, it follows that the map

(11.1.3)     $(h : D_1 \rightarrow D_2 \quad \overline{g} : A_1 \rightarrow A_2 \quad \overline{f} : V_1 \rightarrow V_2)$

is homomorphism of $A_1$-vector space of rows $V_1$ into $A_2$-vector space of rows $V_2$.

Let $f$ be matrix of homomorphisms $\overline{f}$, $\overline{g}$ relative to bases $\overline{\overline{e}}_V$, $\overline{\overline{e}}_W$. The equality

$$\overline{f} \circ (e_V{}^*{}_* v) = e_W{}^*{}_* f^*{}_* v = \overline{g} \circ (e_V{}^*{}_* v)$$

follows from the theorem 11.3.7. Therefore, $\overline{f} = \overline{g}$.                    □

On the basis of theorems 11.1.6, 11.1.7 we identify the homomorphism

(11.1.3)     $(h : D_1 \rightarrow D_2 \quad \overline{g} : A_1 \rightarrow A_2 \quad \overline{f} : V_1 \rightarrow V_2)$

of right $A$-vector space $V$ of rows and coordinates of its presentation

$$(11.1.28) \quad \overline{g} \circ (a e_{A_1}) = h(a) g e_{A_2}$$

$$(11.1.31) \quad \overline{f} \circ (e_{V_1}{}^*{}_* v) = e_{V_2}{}^*{}_* f^*{}_* (\overline{g} \circ v)$$

## 11.2. Homomorphism When Rings $D_1 = D_2 = D$

Let $A_i$, $i = 1, 2$, be division algebra over commutative ring $D$. Let $V_i$, $i = 1, 2$, be right $A_i$-vector space.

DEFINITION 11.2.1. *Let diagram of representations*

$$(11.2.1)$$

$$g_{1.12}(d) : a \to d\,a$$
$$g_{1.23}(v) : w \to C_1(w, v)$$
$$C_1 \in \mathcal{L}(A_1^2 \to A_1)$$
$$g_{1.34}(a) : v \to v\,a$$
$$g_{1.14}(d) : v \to v\,d$$

*describe right $A_1$-vector space $V_1$. Let diagram of representations*

$$(11.2.2)$$

$$g_{2.12}(d) : a \to d\,a$$
$$g_{2.23}(v) : w \to C_2(w, v)$$
$$C_2 \in \mathcal{L}(A_2^2 \to A_2)$$
$$g_{2.34}(a) : v \to v\,a$$
$$g_{2.14}(d) : v \to v\,d$$

*describe right $A_2$-vector space $V_2$. Morphism*

$$(11.2.3) \quad (g : A_1 \to A_2 \quad f : V_1 \to V_2)$$

*of diagram of representations* (11.2.1) *into diagram of representations* (11.2.2) *is called* **homomorphism** *of right $A_1$-vector space $V_1$ into right $A_2$-vector space $V_2$. Let us denote* $\mathrm{Hom}(D; A_1 \to A_2; V_1* \to V_2*)$ *set of homomorphisms of right $A_1$-vector space $V_1$ into right $A_2$-vector space $V_2$.* □

We will use notation

$$f \circ a = f(a)$$

for image of homomorphism $f$.

THEOREM 11.2.2. *The homomorphism*

$$(11.2.3) \quad (g : A_1 \to A_2 \quad f : V_1 \to V_2)$$

*of right $A_1$-vector space $V_1$ into right $A_2$-vector space $V_2$ satisfies following equalities*

$$(11.2.4) \quad g \circ (a + b) = g \circ a + g \circ b$$

$$(11.2.5) \quad g \circ (ab) = (g \circ a)(g \circ b)$$

$$(11.2.6) \qquad g \circ (ap) = (g \circ a)p$$

$$(11.2.7) \qquad f \circ (u + v) = f \circ u + f \circ v$$

$$(11.2.8) \qquad f \circ (va) = (f \circ v)(g \circ a)$$

$$p \in D \quad a, b \in A_1 \quad u, v \in V_1$$

Proof.    According to definitions    2.3.9, 2.8.5,   7.4.1,   1 the map

$$g : A_1 \rightarrow A_2$$

is homomorphism of $D$-algebra $A_1$ into $D$-algebra $A_2$. Therefore, equalities (11.2.4), (11.2.5), (11.2.6) follow from the theorem 4.3.10.    The equality (11.2.7) follows from the definition 11.2.1, since, accorfing to the definition 2.3.9, the map $f$ is homomorpism of Abelian group. The equality (11.2.8) follows from the equality (2.3.6) because the map

$$(g : A_1 \rightarrow A_2 \quad f : V_1 \rightarrow V_2)$$

is morphism of representation $g_{1.34}$ into representation $g_{2.34}$.    □

Definition 11.2.3. *Homomorphism* 11.4

$$(g : A_1 \rightarrow A_2 \quad f : V_1 \rightarrow V_2)$$

*is called* **isomorphism** *between right $A_1$-vector space $V_1$ and right $A_2$-vector space $V_2$, if there exists the map*

$$(g^{-1} : A_2 \rightarrow A_1 \quad f^{-1} : V_2 \rightarrow V_1)$$

*which is homomorphism.*    □

Theorem 11.2.4. *The homomorphism* 11.5

$$(11.2.3) \quad (\overline{g} : A_1 \rightarrow A_2 \quad \overline{f} : V_1 \rightarrow V_2)$$

*of right $A_1$-vector space of columns $V_1$ into right $A_2$-vector space of columns $V_2$ has presentation*

$$(11.2.9) \qquad b = ga$$

$$(11.2.10) \qquad \overline{g} \circ (a^i e_{A_1 i}) = a^i g_i^k e_{A_2 k}$$

Theorem 11.2.5. *The homomorphism* 11.6

$$(11.2.3) \quad (\overline{g} : A_1 \rightarrow A_2 \quad \overline{f} : V_1 \rightarrow V_2)$$

*of right $A_1$-vector space of rows $V_1$ into right $A_2$-vector space of rows $V_2$ has presentation*

$$(11.2.16) \qquad b = ag$$

$$(11.2.17) \qquad \overline{g} \circ (a_i e_{A_1}^i) = a_i g_k^i e_{A_2}^k$$

---

11.4  I follow the definition on page [18]-49.

11.5  In theorems 11.2.4, 11.2.6, we use the following convention.    Let $i = 1, 2$.    Let the set of vectors $\overline{\overline{e}}_{A_i} = (e_{A_i k}, k \in K_i)$ be a basis and $C_{ij}^k,$ $k, i, j \in K_i,$ be structural constants of $D$-algebra of columns $A_i$.    Let the set of vectors $\overline{\overline{e}}_{V_i} = (e_{V_i i}, i \in I_i)$ be a basis of right $A_i$-vector space $V_i$.

11.6  In theorems 11.2.5, 11.2.7, we use the following convention.    Let $i = 1, 2$.    Let the set of vectors $\overline{\overline{e}}_{A_i} = (e_{A_i}^k, k \in K_i)$ be a basis and $C_{ik}^{ij},$ $k, i, j \in K_i,$ be structural constants of $D$-algebra of rows $A_i$.    Let the set of vectors $\overline{\overline{e}}_{V_i} = (e_{V_i}^i, i \in I_i)$ be a basis of right $A_i$-vector space $V_i$.

(11.2.11) $\qquad \overline{g} \circ (e_{A_1} a) = e_{A_2} ga$

(11.2.12) $\qquad w = f_*{}^* (\overline{g} \circ v)$

(11.2.13) $\qquad \overline{f} \circ (e_{V_1 i} v^i) = e_{V_2 k} f_i^k (\overline{g} \circ v^i)$

(11.2.14) $\qquad \overline{f} \circ (e_{V_1 *} {}^* v) = e_{V_2 *} {}^* f_*{}^* (\overline{g} \circ v)$

*relative to selected bases. Here*

- *a is coordinate matrix of $A_1$-number $\overline{a}$ relative the basis $\overline{\overline{e}}_{A_1}$*

  $$\overline{a} = e_{A_1} a$$

- *b is coordinate matrix of $A_2$-number*

  $$\overline{b} = \overline{g} \circ \overline{a}$$

  *relative the basis $\overline{\overline{e}}_{A_2}$*

  $$\overline{b} = e_{A_2} b$$

- *g is coordinate matrix of set of $A_2$-numbers $(\overline{g} \circ e_{A_1 k}, k \in K)$ relative the basis $\overline{\overline{e}}_{A_2}$.*

- *There is relation between the matrix of homomorphism and structural constants*

(11.2.15) $\qquad C_{1 ij}^{k} g_k^l = g_i^p g_j^q C_{2 pq}^{l}$

- *v is coordinate matrix of $V_1$-number $\overline{v}$ relative the basis $\overline{\overline{e}}_{V_1}$*

  $$\overline{v} = e_{V_1 *} {}^* v$$

- *$g(v) = (g(v_i), i \in I)$ is a matrix of $A_2$-numbers.*

- *w is coordinate matrix of $V_2$-number*

  $$\overline{w} = \overline{f} \circ \overline{v}$$

  *relative the basis $\overline{\overline{e}}_{V_2}$*

  $$\overline{w} = e_{V_2 *} {}^* w$$

- *f is coordinate matrix of set of $V_2$-numbers $(\overline{f} \circ e_{V_1 i}, i \in I)$ relative the basis $\overline{\overline{e}}_{V_2}$.*

*For given homomorphism*

(11.2.3) $\qquad (\overline{g} : A_1 \rightarrow A_2 \quad \overline{f} : V_1 \rightarrow V_2)$

*the set of matrices $(g, f)$ is unique*

(11.2.18) $\qquad \overline{g} \circ (a e_{A_1}) = a g e_{A_2}$

(11.2.19) $\qquad w = f^* {}_* (\overline{g} \circ v)$

(11.2.20) $\qquad \overline{f} \circ (e_{V_1}^i v_i) = e_{V_2}^k f_k^i (\overline{g} \circ v_i)$

(11.2.21) $\qquad \overline{f} \circ (e_{V_1 *}^{\quad *} v) = e_{V_2 *}^{\quad *} f^*{}_* (\overline{g} \circ v)$

*relative to selected bases. Here*

- *a is coordinate matrix of $A_1$-number $\overline{a}$ relative the basis $\overline{\overline{e}}_{A_1}$*

  $$\overline{a} = a e_{A_1}$$

- *b is coordinate matrix of $A_2$-number*

  $$\overline{b} = \overline{g} \circ \overline{a}$$

  *relative the basis $\overline{\overline{e}}_{A_2}$*

  $$\overline{b} = b e_{A_2}$$

- *g is coordinate matrix of set of $A_2$-numbers $(\overline{g} \circ e_{A_1}^k, k \in K)$ relative the basis $\overline{\overline{e}}_{A_2}$.*

- *There is relation between the matrix of homomorphism and structural constants*

(11.2.22) $\qquad C_{1 k}^{ij} g_l^k = g_p^i g_q^j C_{2 l}^{pq}$

- *v is coordinate matrix of $V_1$-number $\overline{v}$ relative the basis $\overline{\overline{e}}_{V_1}$*

  $$\overline{v} = e_{V_1 *}^{\quad *} v$$

- *$g(v) = (g(v^i), i \in I)$ is a matrix of $A_2$-numbers.*

- *w is coordinate matrix of $V_2$-number*

  $$\overline{w} = \overline{f} \circ \overline{v}$$

  *relative the basis $\overline{\overline{e}}_{V_2}$*

  $$\overline{w} = e_{V_2 *}^{\quad *} w$$

- *f is coordinate matrix of set of $V_2$-numbers $(\overline{f} \circ e_{V_1}^i, i \in I)$ relative the basis $\overline{\overline{e}}_{V_2}$.*

*For given homomorphism*

(11.2.3) $\qquad (\overline{g} : A_1 \rightarrow A_2 \quad \overline{f} : V_1 \rightarrow V_2)$

*the set of matrices $(g, f)$ is unique*

*and is called* **coordinates of homomorphism** (11.2.3) *relative bases* $(\overline{\overline{e}}_{A_1}, \overline{\overline{e}}_{V_1})$, $(\overline{\overline{e}}_{A_2}, \overline{\overline{e}}_{V_2})$.

PROOF. According to definitions 2.3.9, 2.8.5, 5.4.7, 11.2.1, the map

$$\overline{g} : A_1 \to A_2$$

is homomorphism of $D$-algebra $A_1$ into $D$-algebra $A_2$. Therefore, equalities (11.2.9), (11.2.10), (11.2.11), (11.2.15) follow from equalities

$$\boxed{(5.4.53) \quad b = fa}$$

$$\boxed{(5.4.54) \quad \overline{f} \circ (a^i e_{1 i}) = a^i f_i^k e_k}$$

$$\boxed{(5.4.55) \quad \overline{f} \circ (e_1 a) = e_2 fa}$$

$$\boxed{(5.4.58) \quad C_{1 ij}^{\ k} f_k^l = f_i^p f_j^q C_{2 pq}^{\ \ l}}$$

From the theorem 4.3.11, it follows that the matrix $g$ is unique.

Vector $\overline{v} \in V_1$ has expansion

$$(11.2.23) \qquad \overline{v} = e_{V_1 *}{}^* v$$

relative to the basis $\overline{\overline{e}}_{V_1}$. Vector $\overline{w} \in V_2$ has expansion

$$(11.2.24) \qquad \overline{w} = e_{V_2 *}{}^* w$$

relative to the basis $\overline{\overline{e}}_{V_2}$. Since $\overline{f}$ is a homomorphism, then the equality

$$(11.2.25) \quad \begin{aligned} \overline{w} = \overline{f} \circ \overline{v} &= \overline{f} \circ (e_{V_1 *}{}^* v) \\ &= (\overline{f} \circ e_{V_1})_*{}^* (\overline{g} \circ v) \end{aligned}$$

follows from equalities (11.2.7), (11.2.8). $V_2$-number $\overline{f} \circ e_{V_1 i}$ has expansion

$$(11.2.26) \quad \overline{f} \circ e_{V_1 i} = e_{V_2 *}{}^* f_i = e_{V_2}^j f_i^j$$

relative to basis $\overline{\overline{e}}_{V_2}$. The equality

$$(11.2.27) \qquad \overline{w} = e_{V_2 *}{}^* f_*{}^* (\overline{g} \circ v)$$

follows from equalities (11.2.25), (11.2.26). The equality (11.2.12) follows from comparison of (11.2.24) and (11.2.27) and the theorem 8.4.18. From the equality (11.2.26) and from the theorem 8.4.18, it follows that the matrix $f$ is unique. $\square$

---

*and is called* **coordinates of homomorphism** (11.2.3) *relative bases* $(\overline{\overline{e}}_{A_1}, \overline{\overline{e}}_{V_1})$, $(\overline{\overline{e}}_{A_2}, \overline{\overline{e}}_{V_2})$.

PROOF. According to definitions 2.3.9, 2.8.5, 5.4.7, 11.2.1, the map

$$\overline{g} : A_1 \to A_2$$

is homomorphism of $D$-algebra $A_1$ into $D$-algebra $A_2$. Therefore, equalities (11.2.16), (11.2.17), (11.2.18), (11.2.22) follow from equalities

$$\boxed{(5.4.64) \quad b = af}$$

$$\boxed{(5.4.65) \quad \overline{f} \circ (a_i e_1^i) = a_i f_k^i e^k}$$

$$\boxed{(5.4.66) \quad \overline{f} \circ (ae_1) = afe_2}$$

$$\boxed{(5.4.69) \quad C_{1 k}^{\ ij} f_l^k = f_p^i f_q^j C_{2 l}^{\ pq}}$$

From the theorem 4.3.12, it follows that the matrix $g$ is unique.

Vector $\overline{v} \in V_1$ has expansion

$$(11.2.28) \qquad \overline{v} = e_{V_1}{}^*{}_* v$$

relative to the basis $\overline{\overline{e}}_{V_1}$. Vector $\overline{w} \in V_2$ has expansion

$$(11.2.29) \qquad \overline{w} = e_{V_2}{}^*{}_* w$$

relative to the basis $\overline{\overline{e}}_{V_2}$. Since $\overline{f}$ is a homomorphism, then the equality

$$(11.2.30) \quad \begin{aligned} \overline{w} = \overline{f} \circ \overline{v} &= \overline{f} \circ (e_{V_1}{}^*{}_* v) \\ &= (\overline{f} \circ e_{V_1})^*{}_* (\overline{g} \circ v) \end{aligned}$$

follows from equalities (11.2.7), (11.2.8). $V_2$-number $\overline{f} \circ e_{V_1}^i$ has expansion

$$(11.2.31) \quad \overline{f} \circ e_{V_1}^i = e_{V_2}{}^*{}_* f^i = e_{V_2}^j f_j^i$$

relative to basis $\overline{\overline{e}}_{V_2}$. The equality

$$(11.2.32) \qquad \overline{w} = e_{V_2}{}^*{}_* f^*{}_* (\overline{g} \circ v)$$

follows from equalities (11.2.30), (11.2.31). The equality (11.2.19) follows from comparison of (11.2.29) and (11.2.32) and the theorem 8.4.25. From the equality (11.2.31) and from the theorem 8.4.25, it follows that the matrix $f$ is unique. $\square$

---

THEOREM 11.2.6. *Let*
$$g = (g_l^k, k \in K, l \in L)$$

---

THEOREM 11.2.7. *Let*
$$g = (g_k^l, k \in K, l \in L)$$

be matrix of $D$-numbers which satisfy the equality

$$(11.2.15) \quad C_1{}^k_{ij}g^l_k = g^p_i g^q_j C_{2}{}^l_{pq}$$

Let
$$f = (f^i_j, i \in I, j \in J)$$
be matrix of $A_2$-numbers.     The map[11.5]

$$(11.2.3) \quad (\overline{g} : A_1 \to A_2 \quad \overline{f} : V_1 \to V_2)$$

defined by equalities

$$(11.2.11) \quad \overline{f} \circ (e_{A_1}a) = e_{A_2}fa$$

$$(11.2.14) \quad \overline{f} \circ (e_{V_1*}{}^*v) = e_{V_2*}{}^*f_*{}^*(\overline{g} \circ v)$$

is homomorphism of right $A_1$-vector space of columns $V_1$ into right $A_2$-vector space of columns $V_2$.     The homomorphism

$$(11.2.3) \quad (\overline{g} : A_1 \to A_2 \quad \overline{f} : V_1 \to V_2)$$

which has the given set of matrices $(g, f)$ is unique.

Proof.     According to the theorem 5.4.11, the map

$$(11.2.33) \quad \overline{g} : A_1 \to A_2$$

is homomorphism of $D$-algebra $A_1$ into $D$-algebra $A_2$ and the homomorphism (11.2.33) is unique.

The equality
$$\overline{f} \circ (e_{V_1*}{}^*(v+w))$$
$$= e_{V_2*}{}^*f_*{}^*(g \circ (v+w))$$
$$(11.2.34) = e_{V_2*}{}^*f_*{}^*(g \circ v)$$
$$+ e_{V_2*}{}^*f_*{}^*(g \circ w)$$
$$= \overline{f} \circ (e_{V_1*}{}^*v) + \overline{f} \circ (e_{V_1*}{}^*w)$$

follows from equalities

$$(8.1.16) \quad a_*{}^*(b_1 + b_2) = a_*{}^*b_1 + a_*{}^*b_2$$

$$(11.1.6) \quad g \circ (a + b) = g \circ a + g \circ b$$

$$(11.1.10) \quad f \circ (va) = (f \circ v)(g \circ a)$$

$$(11.1.16) \quad \overline{f} \circ (e_{V_1*}{}^*v) = e_{V_2*}{}^*f_*{}^*(\overline{g} \circ v)$$

From the equality (11.2.34), it fol-

---

be matrix of $D$-numbers which satisfy the equality

$$(11.2.22) \quad C_{1k}{}^{ij}g^k_l = g^i_p g^j_q C_{2l}{}^{pq}$$

Let
$$f = (f^j_i, i \in I, j \in J)$$
be matrix of $A_2$-numbers.     The map[11.6]

$$(11.2.3) \quad (\overline{g} : A_1 \to A_2 \quad \overline{f} : V_1 \to V_2)$$

defined by equalities

$$(11.2.18) \quad \overline{f} \circ (ae_{A_1}) = afe_{A_2}$$

$$(11.2.21) \quad \overline{f} \circ (e_{V_1}{}^*{}_*v) = e_{V_2}{}^*{}_*f^*{}_*(\overline{g} \circ v)$$

is homomorphism of right $A_1$-vector space of rows $V_1$ into right $A_2$-vector space of rows $V_2$.     The homomorphism

$$(11.2.3) \quad (\overline{g} : A_1 \to A_2 \quad \overline{f} : V_1 \to V_2)$$

which has the given     set of matrices $(g, f)$ is unique.

Proof.     According to the theorem 5.4.12, the map

$$(11.2.36) \quad \overline{g} : A_1 \to A_2$$

is homomorphism of $D$-algebra $A_1$ into $D$-algebra $A_2$ and the homomorphism (11.2.36) is unique.

The equality
$$\overline{f} \circ (e_{V_1}{}^*{}_*(v+w))$$
$$= e_{V_2}{}^*{}_*f^*{}_*(g \circ (v+w))$$
$$(11.2.37) = e_{V_2}{}^*{}_*f^*{}_*(g \circ v)$$
$$+ e_{V_2}{}^*{}_*f^*{}_*(g \circ w)$$
$$= \overline{f} \circ (e_{V_1}{}^*{}_*v) + \overline{f} \circ (e_{V_1}{}^*{}_*w)$$

follows from equalities

$$(8.1.18) \quad a^*{}_*(b_1 + b_2) = a^*{}_*b_1 + a^*{}_*b_2$$

$$(11.1.6) \quad g \circ (a + b) = g \circ a + g \circ b$$

$$(11.1.10) \quad f \circ (va) = (f \circ v)(g \circ a)$$

$$(11.1.31) \quad \overline{f} \circ (e_{V_1}{}^*{}_*v) = e_{V_2}{}^*{}_*f^*{}_*(\overline{g} \circ v)$$

From the equality (11.2.37), it fol-

lows that the map $\overline{f}$ is homomorphism of Abelian group. The equality

$$
\begin{aligned}
(11.2.35) \quad & \overline{f} \circ (e_{V_1*}{}^*(va)) \\
= & e_{V_2*}{}^* f_*{}^*(\overline{g} \circ (va)) \\
= & (e_{V_2*}{}^* f_*{}^*(\overline{g} \circ v))(\overline{g} \circ a) \\
= & (\overline{f} \circ (e_{V_1*}{}^* v))(\overline{g} \circ a)
\end{aligned}
$$

follows from equalities

| | |
|---|---|
| (11.1.10) | $f \circ (va) = (f \circ v)(g \circ a)$ |
| (5.4.4) | $g \circ (ab) = (g \circ a)(g \circ b)$ |
| (11.1.16) | $\overline{f} \circ (e_{V_1*}{}^* v) = e_{V_2*}{}^* f_*{}^* (\overline{g} \circ v)$ |

From the equality (11.2.35), and definitions 2.8.5, 11.2.1, it follows that the map

| | |
|---|---|
| (11.2.3) | $(\overline{g} : A_1 \to A_2 \quad \overline{f} : V_1 \to V_2)$ |

is homomorphism of $A_1$-vector space of columns $V_1$ into $A_2$-vector space of columns $V_2$.

Let $f$ be matrix of homomorphisms $\overline{f}$, $\overline{g}$ relative to bases $\overline{\overline{e}}_V$, $\overline{\overline{e}}_W$. The equality

$$\overline{f} \circ (e_{V*}{}^* v) = e_{W*}{}^* f_*{}^* v = \overline{g} \circ (e_{V*}{}^* v)$$

follows from the theorem 11.3.6. Therefore, $\overline{f} = \overline{g}$. $\qquad \square$

On the basis of theorems 11.2.4, 11.2.6 we identify the homomorphism

| | |
|---|---|
| (11.2.3) | $(\overline{g} : A_1 \to A_2 \quad \overline{f} : V_1 \to V_2)$ |

of right $A$-vector space $V$ of columns and coordinates of its presentation

| | |
|---|---|
| (11.2.11) | $\overline{f} \circ (e_{A_1} a) = e_{A_2} fa$ |
| (11.2.14) | $\overline{f} \circ (e_{V_1*}{}^* v) = e_{V_2*}{}^* f_*{}^* (\overline{g} \circ v)$ |

lows that the map $\overline{f}$ is homomorphism of Abelian group. The equality

$$
\begin{aligned}
(11.2.38) \quad & \overline{f} \circ (e_{V_1}{}^*(va)) \\
= & e_{V_2}{}^* f^*{}_*(\overline{g} \circ (va)) \\
= & (e_{V_2}{}^* f^*{}_*(\overline{g} \circ v))(\overline{g} \circ a) \\
= & (\overline{f} \circ (e_{V_1}{}^*_* v))(\overline{g} \circ a)
\end{aligned}
$$

follows from equalities

| | |
|---|---|
| (11.1.10) | $f \circ (va) = (f \circ v)(g \circ a)$ |
| (5.4.4) | $g \circ (ab) = (g \circ a)(g \circ b)$ |
| (11.1.31) | $\overline{f} \circ (e_{V_1}{}^*_* v) = e_{V_2}{}^* f^*{}_* (\overline{g} \circ v)$ |

From the equality (11.2.38), and definitions 2.8.5, 11.2.1, it follows that the map

| | |
|---|---|
| (11.2.3) | $(\overline{g} : A_1 \to A_2 \quad \overline{f} : V_1 \to V_2)$ |

is homomorphism of $A_1$-vector space of rows $V_1$ into $A_2$-vector space of rows $V_2$.

Let $f$ be matrix of homomorphisms $\overline{f}$, $\overline{g}$ relative to bases $\overline{\overline{e}}_V$, $\overline{\overline{e}}_W$. The equality

$$\overline{f} \circ (e_V{}^*_* v) = e_W{}^*_* f^*{}_* v = \overline{g} \circ (e_V{}^*_* v)$$

follows from the theorem 11.3.7. Therefore, $\overline{f} = \overline{g}$. $\qquad \square$

On the basis of theorems 11.2.5, 11.2.7 we identify the homomorphism

| | |
|---|---|
| (11.2.3) | $(\overline{g} : A_1 \to A_2 \quad \overline{f} : V_1 \to V_2)$ |

of right $A$-vector space $V$ of rows and coordinates of its presentation

| | |
|---|---|
| (11.2.18) | $\overline{f} \circ (ae_{A_1}) = af e_{A_2}$ |
| (11.2.21) | $\overline{f} \circ (e_{V_1}{}^*_* v) = e_{V_2}{}^* f^*{}_* (\overline{g} \circ v)$ |

## 11.3. Homomorphism When $D$-algebras $A_1 = A_2 = A$

Let $A$ be division algebra over commutative ring $D$. Let $V_i$, $i = 1, 2$, be right $A$-vector space.

Definition 11.3.1. *Let diagram of representations*

(11.3.1)

$$A \xrightarrow{g_{1.23}} A \xrightarrow{g_{1.34}} V_1$$

$$g_{1.12} \qquad g_{1.12} \qquad g_{1.14}$$

$$D$$

$$g_{1.12}(d) : a \to d \, a$$

$$g_{1.23}(v) : w \to C(w, v)$$

$$C \in \mathcal{L}(A^2 \to A)$$

$$g_{1.34}(a) : v \to v \, a$$

$$g_{1.14}(d) : v \to v \, d$$

*describe right $A$-vector space $V_1$. Let diagram of representations*

(11.3.2)

$$A \xrightarrow{g_{2.23}} A \xrightarrow{g_{2.34}} V_2$$

$$g_{2.12} \qquad g_{2.12} \qquad g_{2.14}$$

$$D$$

$$g_{2.12}(d) : a \to d \, a$$

$$g_{2.23}(v) : w \to C(w, v)$$

$$C \in \mathcal{L}(A^2 \to A)$$

$$g_{2.34}(a) : v \to v \, a$$

$$g_{2.14}(d) : v \to v \, d$$

*describe right $A$-vector space $V_2$. Morphism*

(11.3.3)
$$f : V_1 \to V_2$$

*of diagram of representations* (11.3.1) *into diagram of representations* (11.3.2) *is called* **homomorphism** *of right $A$-vector space $V_1$ into right $A$-vector space $V_2$. Let us denote* $\mathrm{Hom}(D; A; V_1* \to V_2*)$ *set of homomorphisms of right $A$-vector space $V_1$ into right $A$-vector space $V_2$.*  □

We will use notation

$$f \circ a = f(a)$$

for image of homomorphism $f$.

Theorem 11.3.2. *The homomorphism*

$$\boxed{(11.3.3) \quad f : V_1 \to V_2}$$

*of right $A$-vector space $V_1$ into right $A$-vector space $V_2$ satisfies following equalities*

(11.3.4)
$$f \circ (u + v) = f \circ u + f \circ v$$

(11.3.5)
$$f \circ (va) = (f \circ v)a$$

$$a \in A \quad u, v \in V_1$$

Proof.    The equality (11.3.4) follows from the definition 11.3.1, since, accorfing to the definition 2.3.9, the map $f$ is homomorpism of Abelian group. The equality (11.3.5) follows from the equality (2.3.6) because the map

$$f : V_1 \to V_2$$

is morphism of representation $g_{1.34}$ into representation $g_{2.34}$.    □

---

Definition 11.3.3. *Homomorphism* [11.7]

$$f : V_1 \to V_2$$

*is called* **isomorphism** *between right $A$-vector space $V_1$ and right $A$-vector space $V_2$, if there exists the map*

$$f^{-1} : V_2 \to V_1$$

*which is homomorphism.*                                                $\square$

Definition 11.3.4. *A homomorphism* [11.7]

$$f : V \to V$$

*in which source and target are the same right $A$-vector space is called* **endomorphism**. *Endomorphism*

$$f : V \to V$$

*of right $A$-vector space $V$ is called* **automorphism**, *if there exists the map $f^{-1}$ which is endomorphism.*                $\square$

Theorem 11.3.5. *The set $GL(V)$ of automorphisms of right $A$-vector space $V$ is group.*

Theorem 11.3.6. *The homomorphism* [11.8]

$$(11.3.3) \quad \boxed{\overline{f} : V_1 \to V_2}$$

*of right $A$-vector space of columns $V_1$ into right $A$-vector space of columns $V_2$ has presentation*

$$(11.3.6) \qquad w = f_*{}^* v$$

$$(11.3.7) \qquad \overline{f} \circ (v i^{e_{V_1 i}}) = e_{V_2 k} f_i^k v^i$$

$$(11.3.8) \qquad \overline{f} \circ (e_{V_1 *}{}^* v) = e_{V_2 *}{}^* f_*{}^* v$$

*relative to selected bases. Here*

- $v$ *is coordinate matrix of $V_1$-number $\overline{v}$ relative the basis $\overline{\overline{e}}_{V_1}$*

$$\overline{v} = e_{V_1 *}{}^* v$$

- $w$ *is coordinate matrix of $V_2$-number*

$$\overline{w} = \overline{f} \circ \overline{v}$$

Theorem 11.3.7. *The homomorphism* [11.9]

$$(11.3.3) \quad \boxed{\overline{f} : V_1 \to V_2}$$

*of right $A$-vector space of rows $V_1$ into right $A$-vector space of rows $V_2$ has presentation*

$$(11.3.9) \qquad w = f^*{}_* v$$

$$(11.3.10) \qquad \overline{f} \circ (v i_{e_{V_1}^i}) = e_{V_2}^k f_k^i v_i$$

$$(11.3.11) \qquad \overline{f} \circ (e_{V_1}{}^*{}_* v) = e_{V_2}{}^*{}^* f^*{}_* v$$

*relative to selected bases. Here*

- $v$ *is coordinate matrix of $V_1$-number $\overline{v}$ relative the basis $\overline{\overline{e}}_{V_1}$*

$$\overline{v} = e_{V_1}{}^*{}_* v$$

- $w$ *is coordinate matrix of $V_2$-number*

$$\overline{w} = \overline{f} \circ \overline{v}$$

---

[11.8]  In theorems 11.3.6, 11.3.8, we use the following convention.    Let  $i = 1, 2$.    Let the set of vectors $\overline{\overline{e}}_{V_i} = (e_{V_i i}, i \in I_i)$ be a basis of right $A$-vector space $V_i$.

[11.9]  In theorems 11.3.7, 11.3.9, we use the following convention.    Let  $i = 1, 2$.    Let the set of vectors $\overline{\overline{e}}_{V_i} = (e_{V_i}^i, i \in I_i)$ be a basis of right $A$-vector space $V_i$.

*relative the basis* $\overline{\overline{e}}_{V_2}$

$$\overline{w} = e_{V_2 *}{}^* w$$

- *f is coordinate matrix of set of $V_2$-numbers $(\overline{f} \circ e_{V_1 i}, i \in I)$ relative the basis* $\overline{\overline{e}}_{V_2}$.

*For given homomorphism $\overline{f}$, the matrix $f$ is unique and is called* **matrix of homomorphism** $\overline{f}$ *relative bases* $\overline{\overline{e}}_1$, $\overline{\overline{e}}_2$.

PROOF. Vector $\overline{v} \in V_1$ has expansion

$$(11.3.12) \qquad \overline{v} = e_{V_1 *}{}^* v$$

relative to the basis $\overline{\overline{e}}_{V_1}$. Vector $\overline{w} \in V_2$ has expansion

$$(11.3.13) \qquad \overline{w} = e_{V_2 *}{}^* w$$

relative to the basis $\overline{\overline{e}}_{V_2}$. Since $\overline{f}$ is a homomorphism, then the equality

$$(11.3.14) \qquad \begin{aligned} \overline{w} &= \overline{f} \circ \overline{v} = \overline{f} \circ (e_{V_1 *}{}^* v) \\ &= (\overline{f} \circ e_{V_1})_* {}^* v \end{aligned}$$

follows from equalities (11.3.4), (11.3.5). $V_2$-number $\overline{f} \circ e_{V_1 i}$ has expansion

$$(11.3.15) \qquad \overline{f} \circ e_{V_1 i} = e_{V_2 *}{}^* f_i = e_{V_2}^j f_i^j$$

relative to basis $\overline{\overline{e}}_{V_2}$. The equality

$$(11.3.16) \qquad \overline{w} = e_{V_2 *}{}^* f_* {}^* v$$

follows from equalities (11.3.14), (11.3.15). The equality (11.3.6) follows from comparison of (11.3.13) and (11.3.16) and the theorem 8.4.18. From the equality (11.3.15) and from the theorem 8.4.18, it follows that the matrix $f$ is unique. $\square$

THEOREM 11.3.8. *Let*

$$f = (f_j^i, i \in I, j \in J)$$

*be matrix of $A$-numbers. The map* [11.8]

$$\boxed{(11.3.3) \quad \overline{f} : V_1 \to V_2}$$

*defined by the equality*

$$\boxed{(11.3.8) \quad \overline{f} \circ (e_{V_1 *}{}^* v) = e_{V_2 *}{}^* f_* {}^* v}$$

*is homomorphism of right $A$-vector space of columns $V_1$ into right $A$-vector*

---

*relative the basis* $\overline{\overline{e}}_{V_2}$

$$\overline{w} = e_{V_2}{}^*{}_* w$$

- *f is coordinate matrix of set of $V_2$-numbers $(\overline{f} \circ e_{V_1}^i, i \in I)$ relative the basis* $\overline{\overline{e}}_{V_2}$.

*For given homomorphism $\overline{f}$, the matrix $f$ is unique and is called* **matrix of homomorphism** $\overline{f}$ *relative bases* $\overline{\overline{e}}_1$, $\overline{\overline{e}}_2$.

PROOF. Vector $\overline{v} \in V_1$ has expansion

$$(11.3.17) \qquad \overline{v} = e_{V_1}{}^*{}_* v$$

relative to the basis $\overline{\overline{e}}_{V_1}$. Vector $\overline{w} \in V_2$ has expansion

$$(11.3.18) \qquad \overline{w} = e_{V_2}{}^*{}_* w$$

relative to the basis $\overline{\overline{e}}_{V_2}$. Since $\overline{f}$ is a homomorphism, then the equality

$$(11.3.19) \qquad \begin{aligned} \overline{w} &= \overline{f} \circ \overline{v} = \overline{f} \circ (e_{V_1}{}^*{}_* v) \\ &= (\overline{f} \circ e_{V_1})^*{}_* v \end{aligned}$$

follows from equalities (11.3.4), (11.3.5). $V_2$-number $\overline{f} \circ e_{V_1}^i$ has expansion

$$(11.3.20) \qquad \overline{f} \circ e_{V_1}^i = e_{V_2}{}^*{}_* f^i = e_{V_2 j} f_j^i$$

relative to basis $\overline{\overline{e}}_{V_2}$. The equality

$$(11.3.21) \qquad \overline{w} = e_{V_2}{}^*{}_* f^*{}_* v$$

follows from equalities (11.3.19), (11.3.20). The equality (11.3.9) follows from comparison of (11.3.18) and (11.3.21) and the theorem 8.4.25. From the equality (11.3.20) and from the theorem 8.4.25, it follows that the matrix $f$ is unique. $\square$

THEOREM 11.3.9. *Let*

$$f = (f_i^j, i \in I, j \in J)$$

*be matrix of $A$-numbers. The map* [11.9]

$$\boxed{(11.3.3) \quad \overline{f} : V_1 \to V_2}$$

*defined by the equality*

$$\boxed{(11.3.11) \quad \overline{f} \circ (e_{V_1}{}^*{}_* v) = e_{V_2}{}^*{}_* f^*{}_* v}$$

*is homomorphism of right $A$-vector space of rows $V_1$ into right $A$-vector*

*space of columns $V_2$.    The homomorphism*

(11.3.3)    $\overline{f}: V_1 \to V_2$

*which has the given matrix $f$ is unique.*

Proof.    The equality

$$\overline{f} \circ (e_{V_1 *}{}^*(v+w))$$
$$= e_{V_2 *} f_*{}^*(v+w)$$
$$= e_{V_2 *} f_*{}^* v + e_{V_2 *} f_*{}^* w$$
$$= \overline{f} \circ (e_{V_1 *}{}^* v) + \overline{f} \circ (e_{V_1 *}{}^* w)$$

(11.3.22)

follows from equalities

(8.1.16)    $a_*{}^*(b_1 + b_2) = a_*{}^* b_1 + a_*{}^* b_2$

(11.2.8)    $f \circ (va) = (f \circ v)a$

(11.3.8)    $\overline{f} \circ (e_{V_1 *}{}^* v) = e_{V_2 *} f_*{}^* v$

From the equality (11.3.22), it follows that the map $\overline{f}$ is homomorphism of Abelian group. The equality

$$\overline{f} \circ (e_{V_1 *}{}^*(va))$$
$$= e_{V_2 *} f_*{}^*(va)$$
$$= (e_{V_2 *} f_*{}^* v)a$$
$$= (\overline{f} \circ (e_{V_1 *}{}^* v))a$$

(11.3.23)

follows from equalities

(11.2.8)    $f \circ (va) = (f \circ v)a$

(11.3.8)    $\overline{f} \circ (e_{V_1 *}{}^* v) = e_{V_2 *} f_*{}^* v$

From the equality (11.3.23), and definitions 2.8.5, 11.3.1, it follows that the map

(11.3.3)    $\overline{f}: V_1 \to V_2$

is homomorphism of $A$-vector space of columns $V_1$ into $A$-vector space of columns $V_2$.

Let $f$ be matrix of homomorphisms $\overline{f}$, $\overline{g}$ relative to bases $\overline{\overline{e}}_V$, $\overline{\overline{e}}_W$. The equality

$$\overline{f} \circ (e_V *{}^* v) = e_W *{}^* f_*{}^* v = \overline{g} \circ (e_V *{}^* v)$$

follows from the theorem 11.3.6. Therefore, $\overline{f} = \overline{g}$.    □

On the basis of theorems 11.3.6, 11.3.8 we identify the homomorphism

*space of rows $V_2$.    The homomorphism*

(11.3.3)    $\overline{f}: V_1 \to V_2$

*which has the given matrix $f$ is unique.*

Proof.    The equality

$$\overline{f} \circ (e_{V_1}{}^*{}_*(v+w))$$
$$= e_{V_2}{}^*{}_* f^*{}_*(v+w)$$
$$= e_{V_2}{}^*{}_* f^*{}_* v + e_{V_2}{}^*{}_* f^*{}_* w$$
$$= \overline{f} \circ (e_{V_1}{}^*{}_* v) + \overline{f} \circ (e_{V_1}{}^*{}_* w)$$

(11.3.24)

follows from equalities

(8.1.18)    $a^*{}_*(b_1 + b_2) = a^*{}_* b_1 + a^*{}_* b_2$

(11.2.8)    $f \circ (va) = (f \circ v)a$

(11.3.11)    $\overline{f} \circ (e_{V_1}{}^*{}_* v) = e_{V_2}{}^*{}_* f^*{}_* v$

From the equality (11.3.24), it follows that the map $\overline{f}$ is homomorphism of Abelian group. The equality

$$\overline{f} \circ (e_{V_1}{}^*{}_*(va))$$
$$= e_{V_2}{}^*{}_* f^*{}_*(va)$$
$$= (e_{V_2}{}^*{}_* f^*{}_* v)a$$
$$= (\overline{f} \circ (e_{V_1}{}^*{}_* v))a$$

(11.3.25)

follows from equalities

(11.2.8)    $f \circ (va) = (f \circ v)a$

(11.3.11)    $\overline{f} \circ (e_{V_1}{}^*{}_* v) = e_{V_2}{}^*{}_* f^*{}_* v$

From the equality (11.3.25), and definitions 2.8.5, 11.3.1, it follows that the map

(11.3.3)    $\overline{f}: V_1 \to V_2$

is homomorphism of $A$-vector space of rows $V_1$ into $A$-vector space of rows $V_2$.

Let $f$ be matrix of homomorphisms $\overline{f}$, $\overline{g}$ relative to bases $\overline{\overline{e}}_V$, $\overline{\overline{e}}_W$. The equality

$$\overline{f} \circ (e_V *{}^* v) = e_W *{}^* f^*{}_* v = \overline{g} \circ (e_V *{}^* v)$$

follows from the theorem 11.3.7. Therefore, $\overline{f} = \overline{g}$.    □

On the basis of theorems 11.3.7, 11.3.9 we identify the homomorphism

| $(11.3.3)$   $\overline{f} : V_1 \to V_2$ | $(11.3.3)$   $\overline{f} : V_1 \to V_2$ |
|---|---|
| of right $A$-vector space $V$ of columns and coordinates of its presentation | of right $A$-vector space $V$ of rows and coordinates of its presentation |
| $(11.3.8)$   $\overline{f} \circ (e_{V_1} {}_* {}^* v) = e_{V_2} {}^* {}_* f_* {}^* v$ | $(11.3.11)$   $\overline{f} \circ (e_{V_1} {}^* {}_* v) = e_{V_2} {}^* {}_* f^* {}_* v$ |

## 11.4. Sum of Homomorphisms of $A$-Vector Space

THEOREM 11.4.1. *Let $V$, $W$ be right $A$-vector spaces and maps $g : V \to W$, $h : V \to W$ be homomorphisms of right $A$-vector space. Let the map $f : V \to W$ be defined by the equality*

$$(11.4.1) \qquad \forall v \in V : f \circ v = g \circ v + h \circ v$$

*The map $f$ is homomorphism of right $A$-vector space and is called* **sum**

$$f = g + h$$

**of homomorphisms** *$g$ and $h$.*

PROOF. According to the theorem 7.4.6 (the equality (7.4.8)), the equality

$$(11.4.2) \qquad \begin{aligned} &f \circ (u + v) \\ &= g \circ (u + v) + h \circ (u + v) \\ &= g \circ u + g \circ v \\ &\quad + h \circ u + h \circ v \\ &= f \circ u + f \circ v \end{aligned}$$

follows from the equality

| $(11.3.4)$   $f \circ (u + v) = f \circ u + f \circ v$ |
|---|

and the equality (11.4.1). According to the theorem 7.4.6 (the equality (7.4.10)), the equality

$$(11.4.3) \qquad \begin{aligned} &f \circ (va) \\ &= g \circ (va) + h \circ (va) \\ &= (g \circ v)a + (h \circ v)a \\ &= (g \circ v + h \circ v)a \\ &= (f \circ v)a \end{aligned}$$

follows from the equality

| $(11.3.5)$   $f \circ (va) = (f \circ v)a$ |
|---|

and the equality (11.4.1).

The theorem follows from the theorem 11.3.2 and equalities (11.4.2), (11.4.3). $\qquad\square$

THEOREM 11.4.2. *Let $V$, $W$ be right $A$-vector spaces of columns. Let $\overline{\overline{e}}_V$ be a basis of right $A$-vector space $V$. Let $\overline{\overline{e}}_W$ be a basis of right $A$-vector space $W$. The homomorphism $\overline{f} : V \to W$ is*

THEOREM 11.4.3. *Let $V$, $W$ be right $A$-vector spaces of rows. Let $\overline{\overline{e}}_V$ be a basis of right $A$-vector space $V$. Let $\overline{\overline{e}}_W$ be a basis of right $A$-vector space $W$. The homomorphism $\overline{f} : V \to W$ is* *sum of homomorphisms $\overline{g} : V \to W$, $\overline{h} : V \to W$ iff matrix $f$ of homomorphism $\overline{f}$ relative to bases $\overline{\overline{e}}_V$, $\overline{\overline{e}}_W$ is equal to the sum of the matrix $g$ of the homomorphism $\overline{g}$ relative to bases $\overline{\overline{e}}_V$, $\overline{\overline{e}}_W$ and the matrix $h$ of the homomorphism $\overline{h}$ relative to bases $\overline{\overline{e}}_V$, $\overline{\overline{e}}_W$.*

PROOF. According to theorem 11.3.6

$$(11.4.4) \qquad \overline{f} \circ (e_{U*}{}^*a) = e_{V*}{}^*f_*{}^*a$$

$$(11.4.5) \qquad \overline{g} \circ (e_{U*}{}^*a) = e_{V*}{}^*g_*{}^*a$$

$$(11.4.6) \qquad \overline{h} \circ (e_{U*}{}^*a) = e_{V*}{}^*h_*{}^*a$$

According to the theorem 11.4.1, the equality

$$
\begin{aligned}
(11.4.7) \quad & \overline{f} \circ (e_{U*}{}^*v) \\
& = \overline{g} \circ (e_{U*}{}^*v) + \overline{h} \circ (e_{U*}{}^*v) \\
& = e_{V*}{}^*g_*{}^*v + e_{V*}{}^*h_*{}^*v \\
& = e_{V*}{}^*(g + h)_*{}^*v
\end{aligned}
$$

follows from equalities (8.1.16), (8.1.17), (11.4.5), (11.4.6). The equality

$$f = g + h$$

follows from equalities (11.4.4), (11.4.7) and the theorem 11.3.6.

Let

$$f = g + h$$

According to the theorem 11.3.8, there exist homomorphisms $\overline{f} : U \to V$, $\overline{g} : U \to V$, $\overline{h} : U \to V$ such that

$$
\begin{aligned}
(11.4.8) \quad & \overline{f} \circ (e_{U*}{}^*v) \\
& = e_{V*}{}^*f_*{}^*v \\
& = e_{V*}{}^*(g + h)_*{}^*v \\
& = e_{V*}{}^*g_*{}^*v + e_{V*}{}^*h_*{}^*v \\
& = \overline{g} \circ (e_{U*}{}^*v) + \overline{h} \circ (e_{U*}{}^*v)
\end{aligned}
$$

From the equality (11.4.8) and the theorem 11.4.1, it follows that the homomorphism $\overline{f}$ is sum of homomorphisms $\overline{g}$, $\overline{h}$.  □

*sum of homomorphisms $\overline{g} : V \to W$, $\overline{h} : V \to W$ iff matrix $f$ of homomorphism $\overline{f}$ relative to bases $\overline{\overline{e}}_V$, $\overline{\overline{e}}_W$ is equal to the sum of the matrix $g$ of the homomorphism $\overline{g}$ relative to bases $\overline{\overline{e}}_V$, $\overline{\overline{e}}_W$ and the matrix $h$ of the homomorphism $\overline{h}$ relative to bases $\overline{\overline{e}}_V$, $\overline{\overline{e}}_W$.*

PROOF. According to theorem 11.3.7

$$(11.4.9) \qquad \overline{f} \circ (e_{U}{}^*{}_*a) = e_{V}{}^*{}_*f^*{}_*a$$

$$(11.4.10) \qquad \overline{g} \circ (e_{U}{}^*{}_*a) = e_{V}{}^*{}_*g^*{}_*a$$

$$(11.4.11) \qquad \overline{h} \circ (e_{U}{}^*{}_*a) = e_{V}{}^*{}_*h^*{}_*a$$

According to the theorem 11.4.1, the equality

$$
\begin{aligned}
(11.4.12) \quad & \overline{f} \circ (e_{U}{}^*{}_*v) \\
& = \overline{g} \circ (e_{U}{}^*{}_*v) + \overline{h} \circ (e_{U}{}^*{}_*v) \\
& = e_{V}{}^*{}_*g^*{}_*v + e_{V}{}^*{}_*h^*{}_*v \\
& = e_{V}{}^*{}_*(g + h)^*{}_*v
\end{aligned}
$$

follows from equalities (8.1.18), (8.1.19), (11.4.10), (11.4.11). The equality

$$f = g + h$$

follows from equalities (11.4.9), (11.4.12) and the theorem 11.3.7.

Let

$$f = g + h$$

According to the theorem 11.3.9, there exist homomorphisms $\overline{f} : U \to V$, $\overline{g} : U \to V$, $\overline{h} : U \to V$ such that

$$
\begin{aligned}
(11.4.13) \quad & \overline{f} \circ (e_{U}{}^*{}_*v) \\
& = e_{V}{}^*{}_*f^*{}_*v \\
& = e_{V}{}^*{}_*(g + h)^*{}_*v \\
& = e_{V}{}^*{}_*g^*{}_*v + e_{V}{}^*{}_*h^*{}_*v \\
& = \overline{g} \circ (e_{U}{}^*{}_*v) + \overline{h} \circ (e_{U}{}^*{}_*v)
\end{aligned}
$$

From the equality (11.4.13) and the theorem 11.4.1, it follows that the homomorphism $\overline{f}$ is sum of homomorphisms $\overline{g}$, $\overline{h}$.  □

## 11.5. Product of Homomorphisms of $A$-Vector Space

**Theorem 11.5.1.** *Let $U$, $V$, $W$ be right $A$-vector spaces. Let diagram of maps*

(11.5.1)

$$U \xrightarrow{\phantom{xx} f \phantom{xx}} W$$
$$g \searrow \quad \nearrow h$$
$$V$$

(11.5.2)
$$\forall u \in U : f \circ u = (h \circ g) \circ u$$
$$= h \circ (g \circ u)$$

*be commutative diagram where maps $g$, $h$ are homomorphisms of right $A$-vector space. The map $f = h \circ g$ is homomorphism of right $A$-vector space and is called* **product of homomorphisms** $h$, $g$.

Proof. The equality (11.5.2) follows from the commutativity of the diagram (11.5.1). The equality

(11.5.3)
$$f \circ (u + v)$$
$$= h \circ (g \circ (u + v))$$
$$= h \circ (g \circ u + g \circ v))$$
$$= h \circ (g \circ u) + h \circ (g \circ v)$$
$$= f \circ u + f \circ v$$

follows from the equality

(11.3.4)   $f \circ (u + v) = f \circ u + f \circ v$

and the equality (11.5.2). The equality

(11.5.4)
$$f \circ (va) = h \circ (g \circ (va))$$
$$= h \circ ((g \circ v)a)$$
$$= (h \circ (g \circ v))a$$
$$= (f \circ v)a$$

follows from the equality

(11.3.5)   $f \circ (va) = (f \circ v)a$

and the equality (11.4.1). The theorem follows from the theorem 11.3.2 and equalities (11.5.3), (11.5.4).                                                     $\square$

**Theorem 11.5.2.** *Let $U$, $V$, $W$ be right $A$-vector spaces of columns. Let $\overline{\overline{e}}_U$ be a basis of right $A$-vector space $U$. Let $\overline{\overline{e}}_V$ be a basis of right $A$-vector space $V$. Let $\overline{\overline{e}}_W$ be a basis of right $A$-vector space $W$. The homomorphism $\overline{f} : U \to W$ is product of homomorphisms $\overline{h} : V \to W$, $\overline{g} : U \to V$ iff matrix $f$ of homomorphism $\overline{f}$ relative to bases $\overline{\overline{e}}_U$, $\overline{\overline{e}}_W$ is equal to the $_*{}^*$-product of matrix $g$ of the*

**Theorem 11.5.3.** *Let $U$, $V$, $W$ be right $A$-vector spaces of rows. Let $\overline{\overline{e}}_U$ be a basis of right $A$-vector space $U$. Let $\overline{\overline{e}}_V$ be a basis of right $A$-vector space $V$. Let $\overline{\overline{e}}_W$ be a basis of right $A$-vector space $W$. The homomorphism $\overline{f} : U \to W$ is product of homomorphisms $\overline{h} : V \to W$, $\overline{g} : U \to V$ iff matrix $f$ of homomorphism $\overline{f}$ relative to bases $\overline{\overline{e}}_U$, $\overline{\overline{e}}_W$ is equal to the $^*{}_*$-product of matrix $h$ of the* *homomorphism $\overline{\overline{g}}$ relative to bases $\overline{\overline{e}}_U$, $\overline{\overline{e}}_V$ over matrix $h$ of the homomorphism $\overline{\overline{h}}$ relative to bases $\overline{\overline{e}}_V$, $\overline{\overline{e}}_W$*

$$(11.5.5) \qquad f = g_{*} {}_{*} h$$

*homomorphism $\overline{\overline{h}}$ relative to bases $\overline{\overline{e}}_V$, $\overline{\overline{e}}_W$ over matrix $g$ of the homomorphism $\overline{\overline{g}}$ relative to bases $\overline{\overline{e}}_U$, $\overline{\overline{e}}_V$*

$$(11.5.6) \qquad f = h^{*}{}_{*} g$$

PROOF. According to the theorem 11.3.6

$$(11.5.7) \qquad \overline{f} \circ (e_{U\,*}{}^{*} u) = e_{W\,*}{}^{*} f_{*}{}^{*} u$$

$$(11.5.8) \qquad \overline{g} \circ (e_{U\,*}{}^{*} u) = e_{V\,*}{}^{*} g_{*}{}^{*} u$$

$$(11.5.9) \qquad \overline{h} \circ (e_{V\,*}{}^{*} v) = e_{W\,*}{}^{*} h_{*}{}^{*} v$$

According to the theorem 11.5.1, the equality

$$
(11.5.10)
\begin{aligned}
& \overline{f} \circ (e_{U\,*}{}^{*} u) \\
&= \overline{h} \circ (\overline{g} \circ (e_{U\,*}{}^{*} u)) \\
&= \overline{h} \circ (u_{*}{}^{*} g_{*}{}^{*} e_{V}) \\
&= e_{W\,*}{}^{*} h_{*}{}^{*} g_{*}{}^{*} u
\end{aligned}
$$

follows from equalities (11.5.8), (11.5.9). The equality (11.5.5) follows from equalities (11.5.7), (11.5.10) and the theorem 11.3.6.

Let $f = g_{*}{}_{*} h$. According to the theorem 11.3.8, there exist homomorphisms $\overline{f} : U \to V$, $\overline{g} : U \to V$, $\overline{h} : U \to V$ such that

$$
(11.5.11)
\begin{aligned}
& \overline{f} \circ (e_{U\,*}{}^{*} u) = e_{W\,*}{}^{*} f_{*}{}^{*} u \\
&= e_{W\,*}{}^{*} h_{*}{}^{*} g_{*}{}^{*} u \\
&= \overline{h} \circ (e_{V\,*}{}^{*} g_{*}{}^{*} u) \\
&= \overline{h} \circ (\overline{g} \circ (e_{U\,*}{}^{*} u)) \\
&= (\overline{h} \circ \overline{g}) \circ (e_{U\,*}{}^{*} u)
\end{aligned}
$$

From the equality (11.5.11) and the theorem 11.4.1, it follows that the homomorphism $\overline{f}$ is product of homomorphisms $\overline{h}$, $\overline{g}$. □

PROOF. According to the theorem 11.3.7

$$(11.5.12) \qquad \overline{f} \circ (e_{U}{}^{*}{}_{*} u) = e_{W}{}^{*}{}_{*} f^{*}{}_{*} u$$

$$(11.5.13) \qquad \overline{g} \circ (e_{U}{}^{*}{}_{*} u) = e_{V}{}^{*}{}_{*} g^{*}{}_{*} u$$

$$(11.5.14) \qquad \overline{h} \circ (e_{V}{}^{*}{}_{*} v) = e_{W}{}^{*}{}_{*} h^{*}{}_{*} v$$

According to the theorem 11.5.1, the equality

$$
(11.5.15)
\begin{aligned}
& \overline{f} \circ (e_{U}{}^{*}{}_{*} u) \\
&= \overline{h} \circ (\overline{g} \circ (e_{U}{}^{*}{}_{*} u)) \\
&= \overline{h} \circ (u^{*}{}_{*} g^{*}{}_{*} e_{V}) \\
&= u^{*}{}_{*} g^{*}{}_{*} h^{*}{}_{*} e_{W}
\end{aligned}
$$

follows from equalities (11.5.13), (11.5.14). The equality (11.5.6) follows from equalities (11.5.12), (11.5.15) and the theorem 11.3.7.

Let $f = h^{*}{}_{*} g$. According to the theorem 11.3.9, there exist homomorphisms $\overline{f} : U \to V$, $\overline{g} : U \to V$, $\overline{h} : U \to V$ such that

$$
(11.5.16)
\begin{aligned}
& \overline{f} \circ (e_{U}{}^{*}{}_{*} u) = e_{W}{}^{*}{}_{*} f^{*}{}_{*} u \\
&= e_{W}{}^{*}{}_{*} h^{*}{}_{*} g^{*}{}_{*} u \\
&= \overline{h} \circ (e_{V}{}^{*}{}_{*} g^{*}{}_{*} u) \\
&= \overline{h} \circ (\overline{g} \circ (e_{U}{}^{*}{}_{*} u)) \\
&= (\overline{h} \circ \overline{g}) \circ (e_{U}{}^{*}{}_{*} u)
\end{aligned}
$$

From the equality (11.5.16) and the theorem 11.4.1, it follows that the homomorphism $\overline{f}$ is product of homomorphisms $\overline{h}$, $\overline{g}$. □

## 11.6. Direct Sum of Right $A$-Vector Space

Let RIGHT $A$-VECTOR be category of right $A$-vector spaces morphisms of which are homomorphisms of $A$-vector space.

DEFINITION 11.6.1. *Coproduct in category* RIGHT $A$-VECTOR *is called* **direct sum**.[11.10] *We will use notation* $V \oplus W$ *for direct sum of right $A$-vector spaces $V$ and $W$.* □

THEOREM 11.6.2. *In category* RIGHT $A$-VECTOR *there exists direct sum of right $A$-vector spaces.*[11.11]*Let* $\{V_i, i \in I\}$ *be set of right $A$-vector spaces. Then the representation*

$$A \dashrightarrow \bigoplus_{i \in I} V_i \qquad \left(\bigoplus_{i \in I} v_i\right) a = \bigoplus_{i \in I} v_i a$$

*of the $D$-algebra $A$ in direct sum of Abelian groups*

$$V = \bigoplus_{i \in I} V_i$$

*is direct sum of right $A$-vector spaces*

$$V = \bigoplus_{i \in I} V_i$$

PROOF. Let $V_i$, $i \in I$, be set of right $A$-vector spaces. Let

$$V = \bigoplus_{i \in I} V_i$$

be direct sum of Abelian groups. Let $v = (v_i, i \in I) \in V$. We define action of $A$-number $a$ over $V$-number $v$ according to the equality

$$(11.6.1) \qquad a(v_i, i \in I) = (av_i, i \in I)$$

According to the theorem 7.4.6,

$$(11.6.2) \qquad v_i(ba) = (v_i b)a$$

$$(11.6.3) \qquad (v_i + w_i)a = v_i a + w_i a$$

$$(11.6.4) \qquad v_i(a + b) = v_i a + v_i b$$

$$(11.6.5) \qquad v_i 1 = v_i$$

for any $a$, $b \in A$, $v_i$, $w_i \in V_i$, $i \in I$. Equalities

$$(11.6.6) \qquad \begin{aligned} v(ba) &= (v_i, i \in I)(ba) = (v_i(ba), i \in I) \\ &= ((v_i b)a, i \in I) = (v_i b, i \in I)a = ((v_i, i \in I)b)a \\ &= (vb)a \end{aligned}$$

$$(11.6.7) \qquad \begin{aligned} (v + w)a &= ((v_i, i \in I) + (w_i, i \in I))a = (v_i + w_i, i \in I)a \\ &= ((v_i + w_i)a, i \in I) = (v_i a + w_i a, i \in I) \\ &= (v_i a, i \in I) + (w_i a, i \in I) = (v_i, i \in I)a + (w_i, i \in I)a \\ &= va + wa \end{aligned}$$

$$(11.6.8) \qquad \begin{aligned} v(a + b) &= (v_i, i \in I)(a + b) = (v_i(a + b), i \in I) \\ &= (v_i a + v_i b, i \in I) = (v_i a, i \in I) + (v_i b, i \in I) \\ &= (v_i, i \in I)a + (v_i, i \in I)b = va + vb \\ v1 &= (v_i, i \in I)1 = (v_i 1, i \in I) = (v_i, i \in I) = v \end{aligned}$$

---

[11.10] See also the definition of direct sum of modules in [1], page 128. On the same page, Lang proves the existence of direct sum of modules.
[11.11] See also proposition on the page [1]-1.1

follow from equalities  (3.7.4),  (11.6.1),  (11.6.2),  (11.6.3),  (11.6.4),  (11.6.5)  for any $a$, $b \in A$, $v = (v_i, i \in I)$, $w = (w_i, i \in I) \in V$.

From the equality (11.6.7), it follows that the map

$$a : v \to va$$

defined by the equality (11.6.1) is endomorphism of Abelian group $V$. From equalities (11.6.6), (11.6.8), it follows that the map (11.6.1) is homomorphism of $D$-algebra $A$ into set of endomorphisms of Abelian group $V$. Therefore, the map (11.6.1) is representation of $D$-algebra $A$ in Abelian group $V$.

Let $a_1$, $a_2$ be $A$-numbers such that

$$(11.6.9) \qquad\qquad a_1(v_i, i \in I) = a_2(v_i, i \in I)$$

for any $V$-number  $v = (v_i, i \in I)$.  The equality

$$(11.6.10) \qquad\qquad (a_1 v_i, i \in I) = (a_2 v_i, i \in I)$$

follows from the equality (11.6.9). From the equality (11.6.10), it follows that

$$(11.6.11) \qquad\qquad \forall i \in I : a_1 v_i = a_2 v_i$$

Since for any $i \in I$, $V_i$ is $A$-vector space, then corresponding representation is effective according to the definition  7.4.1  Therefore, the equality $a_1 = a_2$ follows from the statement (11.6.11) and the map (11.6.1) is effective representation of the $D$-algebra $A$ in Abelian group $V$. Therefore, Abelian group $V$ is right $A$-vector space.

Let

$$\{f_i : V_i \to W, i \in I\}$$

be set of linear maps into right $A$-vector space $W$. We define the map

$$f : V \to W$$

by the equality

$$(11.6.12) \qquad\qquad f\left(\bigoplus_{i \in I} x_i\right) = \sum_{i \in I} f_i(x_i)$$

The sum in the right side of the equality (11.6.12) is finite, since all summands, except for a finite number, equal 0. From the equality

$$f_i(x_i + y_i) = f_i(x_i) + f_i(y_i)$$

and the equality (11.6.12), it follows that

$$\begin{aligned}
f(\bigoplus_{i \in I}(x_i + y_i)) &= \sum_{i \in I} f_i(x_i + y_i) \\
&= \sum_{i \in I}(f_i(x_i) + f_i(y_i)) \\
&= \sum_{i \in I} f_i(x_i) + \sum_{i \in I} f_i(y_i) \\
&= f(\bigoplus_{i \in I} x_i) + f(\bigoplus_{i \in I} y_i)
\end{aligned}$$

From the equality

$$f_i(dx_i) = df_i(x_i)$$

and the equality (11.6.12), it follows that

$$f((dx_i, i \in I)) = \sum_{i \in I} f_i(dx_i) = \sum_{i \in I} df_i(x_i) = d \sum_{i \in I} f_i(x_i)$$
$$= df((x_i, i \in I))$$

Therefore, the map $f$ is linear map. The equality

$$f \circ \lambda_j(x) = \sum_{i \in I} f_i(\delta^i_j x) = f_j(x)$$

follows from equalities (3.7.5), (11.6.12). Since the map $\lambda_i$ is injective, then the map $f$ is unique. Therefore, the theorem folows from definitions 2.2.12, 3.7.2. $\square$

DEFINITION 11.6.3. *Let $D$ be commutative ring with unit. Let $A$ be associative division $D$-algebra. For any integer $n > 0$, we introduce right vector space $\mathcal{R}^n A$ using definition by induction*

$$\mathcal{R}^1 A = \mathcal{R} A$$
$$\mathcal{R}^k A = \mathcal{R}^{k-1} A \oplus \mathcal{R} A$$

$\square$

DEFINITION 11.6.4. *Let $D$ be commutative ring with unit. Let $A$ be associative division $D$-algebra. For any set $B$ and any integer $n > 0$, we introduce right vector space $\mathcal{R}^n A^B$ using definition by induction*

$$\mathcal{R}^1 A^B = \mathcal{R} A^B$$
$$\mathcal{R}^k A^B = \mathcal{R}^{k-1} A^B \oplus \mathcal{R} A^B$$

$\square$

THEOREM 11.6.5. *Direct sum of right $A$-vector spaces $V_1$, ..., $V_n$ coincides with their Cartesian product*

$$V_1 \oplus ... \oplus V_n = V_1 \times ... \times V_n$$

PROOF. The theorem follows from the theorem 11.6.2. $\square$

THEOREM 11.6.6. *Let $\overline{e}_1 = (e_{1i}, i \in I)$ be a basis of right $A$-vector space of columns $V$. Then we can represent $A$-module $V$ as direct sum*

$$(11.6.13) \qquad V = \bigoplus_{i \in I} e_i A$$

PROOF. The map

$$f : e_i v^i \in V \to (v^i, i \in I) \in \bigoplus_{i \in I} e_i A$$

is isomorphism. $\square$

THEOREM 11.6.7. *Let $\overline{e}_1 = (e_1^i, i \in I)$ be a basis of right $A$-vector space of rows $V$. Then we can represent $A$-module $V$ as direct sum*

$$(11.6.14) \qquad V = \bigoplus_{i \in I} e^i A$$

PROOF. The map

$$f : e^i v_i \in V \to (v_i, i \in I) \in \bigoplus_{i \in I} e^i A$$

is isomorphism. $\square$

CHAPTER 12

# Basis Manifold for $A$-Vector Space

## 12.1. Dimension of Left $A$-Vector Space

THEOREM 12.1.1. *Let $V$ be a left $A$-vector space of columns. Let the set of vectors $\overline{\overline{e}} = (e_i, i \in I)$ be a basis of left $A$-vector space $V$. Let the set of vectors $\overline{\overline{e}} = (e_j, j \in J)$ be a basis of left $A$-vector space $V$. If $|I|$ and $|J|$ are finite numbers then $|I| = |J|$.*

PROOF. Let $|I| = m$ and $|J| = n$. Let

(12.1.1) $\qquad m < n$

Because $\overline{\overline{e}}_1$ is a basis, any vector $e_{2j}$, $j \in J$, has expansion

$$e_{2j} = a_j{}^*{}_* e_1$$

Because $\overline{\overline{e}}_2$ is a basis,

(12.1.2) $\qquad \lambda = 0$

should follow from

(12.1.3) $\qquad \lambda^*{}_* e_2 = \lambda^*{}_* a^*{}_* e_1 = 0$

Because $\overline{\overline{e}}_1$ is a basis we get

(12.1.4) $\qquad \lambda^*{}_* a = 0$

According to (12.1.1), $\operatorname{rank}(^*{}_*)a \leq m$ and system of linear equations (12.1.4) has more variables then equations. According to the theorem 9.4.4, we get the statement $\lambda \neq 0$ which contradicts statement (12.1.2). Therefore, statement $m < n$ is not valid. In the same manner we can prove that the statement $n < m$ is not valid. This completes the proof of the theorem. $\square$

THEOREM 12.1.2. *Let $V$ be a left $A$-vector space of rows. Let the set of vectors $\overline{\overline{e}} = (e^i, i \in I)$ be a basis of left $A$-vector space $V$. Let the set of vectors $\overline{\overline{e}} = (e^j, j \in J)$ be a basis of left $A$-vector space $V$. If $|I|$ and $|J|$ are finite numbers then $|I| = |J|$.*

PROOF. Let $|I| = m$ and $|J| = n$. Let

(12.1.5) $\qquad m < n$

Because $\overline{\overline{e}}_1$ is a basis, any vector $e_2^j$, $j \in J$, has expansion

$$e_2^j = a^j{}_*{}^* e_1$$

Because $\overline{\overline{e}}_2$ is a basis,

(12.1.6) $\qquad \lambda = 0$

should follow from

(12.1.7) $\qquad \lambda_*{}^* e_2 = \lambda_*{}^* a_*{}^* e_1 = 0$

Because $\overline{\overline{e}}_1$ is a basis we get

(12.1.8) $\qquad \lambda_*{}^* a = 0$

According to (12.1.5), $\operatorname{rank}(_*{}^*)a \leq m$ and system of linear equations (12.1.8) has more variables then equations. According to the theorem 9.4.4, we get the statement $\lambda \neq 0$ which contradicts statement (12.1.6). Therefore, statement $m < n$ is not valid. In the same manner we can prove that the statement $n < m$ is not valid. This completes the proof of the theorem. $\square$

DEFINITION 12.1.3. *We call **dimension of left $A$-vector space** the number of vectors in a basis.* $\square$

THEOREM 12.1.4. *If left $A$-vector space $V$ has finite dimension, then the statement*

(12.1.9) $\quad \forall v \in V, av = bv \quad a, b \in A$

THEOREM 12.1.5. *If left $A$-vector space $V$ has finite dimension, then the statement*

(12.1.10) $\quad \forall v \in V, av = bv \quad a, b \in A$

*implies $a = b$.*

PROOF. Let the set of vectors $\overline{\overline{e}} = (e_j, j \in J)$ be the basis of left $A$-vector space $V$. According to theorems 7.2.6, 8.4.4, the equality

$$(12.1.11) \qquad av^i = bv^i$$

follows from the equality (12.1.9). The theorem follows from the equality (12.1.11), if we assume $v = e_1$. $\square$

*implies $a = b$.*

PROOF. Let the set of vectors $\overline{\overline{e}} = (e_j, j \in J)$ be the basis of left $A$-vector space $V$. According to theorems 7.2.6, 8.4.11, the equality

$$(12.1.12) \qquad av_i = bv_i$$

follows from the equality (12.1.10). The theorem follows from the equality (12.1.12), if we assume $v = e_1$. $\square$

THEOREM 12.1.6. *If left $A$-vector space $V$ has finite dimension, then corresponding representation is free.*

PROOF. The theorem follows from the definition 2.3.3 and from theorems 12.1.4, 12.1.5. $\square$

THEOREM 12.1.7. *Let $V$ be a left $A$-vector space of columns. The coordinate matrix of basis $\overline{\overline{g}}$ relative basis $\overline{\overline{e}}$ of left $A$-vector space $V$ is $_*$-nonsingular matrix.*

PROOF. According to the theorem 9.3.5, $_*$-rank of the coordinate matrix of basis $\overline{\overline{g}}$ relative basis $\overline{\overline{e}}$ equal to the dimension of left $A$-space of columns. This proves the statement of the theorem. $\square$

THEOREM 12.1.8. *Let $V$ be a left $A$-vector space of rows. The coordinate matrix of basis $\overline{\overline{g}}$ relative basis $\overline{\overline{e}}$ of left $A$-vector space $V$ is $_*$-nonsingular matrix.*

PROOF. According to the theorem 9.3.8, $_*$-rank of the coordinate matrix of basis $\overline{\overline{g}}$ relative basis $\overline{\overline{e}}$ equal to the dimension of left $A$-space of columns. This proves the statement of the theorem. $\square$

THEOREM 12.1.9. *Let $V$ be a left $A$-vector space of columns. Let $\overline{\overline{e}}$ be a basis of left $A$-vector space $V$. Then any automorphism $\overline{f}$ of left $A$-vector space $V$ has form*

$$(12.1.13) \qquad v' = v *_* f$$

*where $f$ is a $_*$-nonsingular matrix.*

PROOF. The equality (12.1.13) follows from theorem 10.3.6. Because $\overline{f}$ is an isomorphism, for each vector $\overline{v}'$ there exist one and only one vector $\overline{v}$ such that $\overline{v}' = \overline{f} \circ \overline{v}$. Therefore, system of linear equations (12.1.13) has a unique solution. According to corollary 9.4.5 matrix $f$ is a $_*$-nonsingular matrix. $\square$

THEOREM 12.1.10. *Let $V$ be a left $A$-vector space of rows. Let $\overline{\overline{e}}$ be a basis of left $A$-vector space $V$. Then any automorphism $\overline{f}$ of left $A$-vector space $V$ has form*

$$(12.1.14) \qquad v' = v_* {}^* f$$

*where $f$ is a $_*{}^*$-nonsingular matrix.*

PROOF. The equality (12.1.14) follows from theorem 10.3.7. Because $\overline{f}$ is an isomorphism, for each vector $\overline{v}'$ there exist one and only one vector $\overline{v}$ such that $\overline{v}' = \overline{f} \circ \overline{v}$. Therefore, system of linear equations (12.1.14) has a unique solution. According to corollary 9.4.5 matrix $f$ is a $_*{}^*$-nonsingular matrix. $\square$

THEOREM 12.1.11. *Matrices of au-*

THEOREM 12.1.12. *Matrices of au-*

*tomorphisms of left $A$-space $V$ of columns form a group $GL(n^*{}_*A)$.*

Proof. If we have two automorphisms $\overline{f}$ and $\overline{g}$ then we can write

$$v' = v^*{}_* f$$
$$v'' = v'^*{}_* g = v^*{}_* f^*{}_* g$$

Therefore, the resulting automorphism has matrix $f^*{}_* g$. □

Theorem 12.1.13. *Let $\overline{\overline{e}}_{V1}$, $\overline{\overline{e}}_{V2}$ be bases of left $A$-vector space $V$ of columns. Let $\overline{\overline{e}}_{W1}$, $\overline{\overline{e}}_{W2}$ be bases of left $A$-vector space $W$ of columns. Let $a_1$ be matrix of homomorphism*

(12.1.15)          $a : V \to W$

*relative to bases $\overline{\overline{e}}_{V1}$, $\overline{\overline{e}}_{W1}$ and $a_2$ be matrix of homomorphism (12.1.15) relative to bases $\overline{\overline{e}}_{V2}$, $\overline{\overline{e}}_{W2}$. Suppose the basis $\overline{\overline{e}}_{V1}$ has coordinate matrix $b$ relative the basis $\overline{\overline{e}}_{V2}$*

$$\overline{\overline{e}}_{V1} = b^*{}_* \overline{\overline{e}}_{V2}$$

*and $\overline{\overline{e}}_{W1}$ has coordinate matrix $c$ relative the basis $\overline{\overline{e}}_{W2}$*

(12.1.16)          $\overline{\overline{e}}_{W1} = c^*{}_* \overline{\overline{e}}_{W2}$

*Then there is relationship between matrices $a_1$ and $a_2$*

$$a_1 = b^*{}_* a_2{}^*{}_* c^{-1^*{}^*}$$

Proof. Vector $\overline{v} \in V$ has expansion

$$\overline{v} = v^*{}_* e_{V1} = v^*{}_* b^*{}_* e_{V2}$$

relative to bases $\overline{\overline{e}}_{V1}$, $\overline{\overline{e}}_{V2}$. Since $a$ is homomorphism, we can write it as

(12.1.19)          $\overline{w} = v^*{}_* a_1{}^*{}_* e_{W1}$

relative to bases $\overline{\overline{e}}_{V1}$, $\overline{\overline{e}}_{W1}$ and as

(12.1.20)          $\overline{w} = v^*{}_* b^*{}_* a_2{}^*{}_* e_{W2}$

relative to bases $\overline{\overline{e}}_{V2}$, $\overline{\overline{e}}_{W2}$. According to the theorem 12.1.7, matrix $c$ has $^*{}_*$-inverse and from equation (12.1.16) it follows that

(12.1.21)          $\overline{\overline{e}}_{W2} = c^{-1^*{}^*}{}_* \overline{\overline{e}}_{W1}$

Theorem 12.1.14. *Let $\overline{\overline{e}}_{V1}$, $\overline{\overline{e}}_{V2}$ be bases of left $A$-vector space $V$ of columns. Let $\overline{\overline{e}}_{W1}$, $\overline{\overline{e}}_{W2}$ be bases of left $A$-vector space $W$ of columns. Let $a_1$ be matrix of homomorphism*

(12.1.17)          $a : V \to W$

*relative to bases $\overline{\overline{e}}_{V1}$, $\overline{\overline{e}}_{W1}$ and $a_2$ be matrix of homomorphism (12.1.17) relative to bases $\overline{\overline{e}}_{V2}$, $\overline{\overline{e}}_{W2}$. Suppose the basis $\overline{\overline{e}}_{V1}$ has coordinate matrix $b$ relative the basis $\overline{\overline{e}}_{V2}$*

$$\overline{\overline{e}}_{V1} = b_*{}^* \overline{\overline{e}}_{V2}$$

*and $\overline{\overline{e}}_{W1}$ has coordinate matrix $c$ relative the basis $\overline{\overline{e}}_{W2}$*

(12.1.18)          $\overline{\overline{e}}_{W1} = c_*{}^* \overline{\overline{e}}_{W2}$

*Then there is relationship between matrices $a_1$ and $a_2$*

$$a_1 = b_*{}^* a_2{}_*{}^* c^{-1_*{}^*}$$

Proof. Vector $\overline{v} \in V$ has expansion

$$\overline{v} = v_*{}^* e_{V1} = v_*{}^* b_*{}^* e_{V2}$$

relative to bases $\overline{\overline{e}}_{V1}$, $\overline{\overline{e}}_{V2}$. Since $a$ is homomorphism, we can write it as

(12.1.24)          $\overline{w} = v_*{}^* a_1{}_*{}^* e_{W1}$

relative to bases $\overline{\overline{e}}_{V1}$, $\overline{\overline{e}}_{W1}$ and as

(12.1.25)          $\overline{w} = v_*{}^* b_*{}^* a_2{}_*{}^* e_{W2}$

relative to bases $\overline{\overline{e}}_{V2}$, $\overline{\overline{e}}_{W2}$. According to the theorem 12.1.8, matrix $c$ has $_*{}^*$-inverse and from equation (12.1.18) it follows that

(12.1.26)          $\overline{\overline{e}}_{W2} = c^{-1_*{}^*}{}^* \overline{\overline{e}}_{W1}$

The right column top:

*tomorphisms of left $A$-space $V$ of rows form a group $GL(n_*{}^*A)$.*

Proof. If we have two automorphisms $\overline{f}$ and $\overline{g}$ then we can write

$$v' = v_*{}^* f$$
$$v'' = v'_*{}^* g = v_*{}^* f_*{}^* g$$

Therefore, the resulting automorphism has matrix $f_*{}^* g$. □

The equality

(12.1.22)  $\overline{w} = v^*{}_* b^*{}_* a_2{}^*{}_* c^{-1^*}{}^*{}_* e_{W1}$

follows from equalities (12.1.20), (12.1.21). From the theorem 8.4.4 and comparison of equations (12.1.19) and (12.1.22) it follows that

(12.1.23)  $v^*{}_* a_1 = v^*{}_* b^*{}_* a_2{}^*{}_* c^{-1^*}{}$

Since vector $a$ is arbitrary vector, the statement of theorem follows from the theorem 12.2.1 and equation (12.1.23).

□

The equality

(12.1.27)  $\overline{w} = v_*{}^* b_*{}^* a_{2*}{}^* c^{-1_*}{}^*{}^*{}_* e_{W1}$

follows from equalities (12.1.25), (12.1.26). From the theorem 8.4.11 and comparison of equations (12.1.24) and (12.1.27) it follows that

(12.1.28)  $v^*{}_* a_1 = v_*{}^* b_*{}^* a_{2*}{}^* c^{-1_*}{}^*$

Since vector $a$ is arbitrary vector, the statement of theorem follows from the theorem 12.2.2 and equation (12.1.28).

□

THEOREM 12.1.15. *Let $\overline{a}$ be automorphism of left $A$-vector space $V$ of columns. Let $a_1$ be matrix of this automorphism defined relative to basis $\overline{\overline{e}}_1$ and $a_2$ be matrix of the same automorphism defined relative to basis $\overline{\overline{e}}_2$. Suppose the basis $\overline{\overline{e}}_1$ has coordinate matrix $b$ relative the basis $\overline{\overline{e}}_2$*

$$\overline{\overline{e}}_1 = b^*{}_* \overline{\overline{e}}_2$$

*Then there is relationship between matrices $a_1$ and $a_2$*

$$a_1 = b^*{}_* a_2{}^*{}_* c^{-1^*}{}$$

PROOF. Statement follows from theorem 12.1.13, because in this case $c = b$.

□

THEOREM 12.1.16. *Let $\overline{a}$ be automorphism of left $A$-vector space $V$ of columns. Let $a_1$ be matrix of this automorphism defined relative to basis $\overline{\overline{e}}_1$ and $a_2$ be matrix of the same automorphism defined relative to basis $\overline{\overline{e}}_2$. Suppose the basis $\overline{\overline{e}}_1$ has coordinate matrix $b$ relative the basis $\overline{\overline{e}}_2$*

$$\overline{\overline{e}}_1 = b_*{}^* \overline{\overline{e}}_2$$

*Then there is relationship between matrices $a_1$ and $a_2$*

$$a_1 = b_*{}^* a_{2*}{}^* c^{-1_*}{}^*$$

PROOF. Statement follows from theorem 12.1.14, because in this case $c = b$.

□

## 12.2. Basis Manifold for Left $A$-Vector Space

Let $V$ be a left $A$-vector space of columns. We proved in the theorem 12.1.9 that, if we select a basis of left $A$-vector space $V$, then we can identify any automorphism $\overline{f}$ of left $A$-vector space $V$ with $*_*$-nonsingular matrix $f$. Corresponding transformation of coordinates of vector

$$v' = v^*{}_* f$$

is called **linear transformation**.

Let $V$ be a left $A$-vector space of rows. We proved in the theorem 12.1.10 that, if we select a basis of left $A$-vector space $V$, then we can identify any automorphism $\overline{f}$ of left $A$-vector space $V$ with $*_*$-nonsingular matrix $f$. Corresponding transformation of coordinates of vector

$$v' = v_*{}^* f$$

is called **linear transformation**.

THEOREM 12.2.1. *Let $V$ be a left $A$-vector space of columns and $\overline{\overline{e}}$ be basis of left $A$-vector space $V$. Automorphisms of left $A$-vector space of columns form a right-side linear effective $GL(\mathbf{n}^*{}_* A)$-*

THEOREM 12.2.2. *Let $V$ be a left $A$-vector space of rows and $\overline{\overline{e}}$ be basis of left $A$-vector space $V$. Automorphisms of left $A$-vector space of rows form a right-side linear effective $GL(\mathbf{n}_*{}^* A)$-represen-*

*representation.*

PROOF. Let $a$, $b$ be matrices of automorphisms $\overline{a}$ and $\overline{b}$ with respect to basis $\overline{\overline{e}}$. According to the theorem 10.3.6, coordinate transformation has following form

$$(12.2.1) \qquad v' = v^*{}_*a$$

$$(12.2.2) \qquad v'' = v'^*{}_*b$$

The equality

$$v'' = v^*{}_*a^*{}_*b$$

follows from equalities (12.2.1), (12.2.2). According to the theorem 10.5.2, the product of automorphisms $\overline{a}$ and $\overline{b}$ has matrix $a^*{}_*b$. Therefore, automorphisms of left $A$-vector space of columns form a right-side linear $GL(n^*{}_*A)$-representation.

It remains to prove that the kernel of inefficiency consists only of identity. Identity transformation satisfies to equation

$$v^i = v^j a^i_j$$

Choosing values of coordinates as $v^i = \delta^i_k$ where we selected $k$ we get

$$(12.2.3) \qquad \delta^i_k = \delta^j_k a^i_j$$

From (12.2.3) it follows

$$\delta^i_k = a^i_k$$

Since $k$ is arbitrary, we get the conclusion $a = \delta$. $\qquad\square$

THEOREM 12.2.3. *Automorphism a acting on each vector of basis of left $A$-vector space of columns maps a basis into another basis.*

PROOF. Let $\overline{\overline{e}}$ be basis of left $A$-vector space $V$ of columns. According to theorem 12.1.9, vector $e_i$ maps into a vector $e'_i$

$$(12.2.7) \qquad e'_i = e_i{}^*{}_*a$$

Let vectors $e'_i$ be linearly dependent. Then $\lambda \neq 0$ in the equality

$$(12.2.8) \qquad \lambda^*{}_*e' = 0$$

*tation.*

PROOF. Let $a$, $b$ be matrices of automorphisms $\overline{a}$ and $\overline{b}$ with respect to basis $\overline{\overline{e}}$. According to the theorem 10.3.7, coordinate transformation has following form

$$(12.2.4) \qquad v' = v_*{}^*a$$

$$(12.2.5) \qquad v'' = v'_*{}^*b$$

The equality

$$v'' = v_*{}^*a_*{}^*b$$

follows from equalities (12.2.4), (12.2.5). According to the theorem 10.5.3, the product of automorphisms $\overline{a}$ and $\overline{b}$ has matrix $a_*{}^*b$. Therefore, automorphisms of left $A$-vector space of rows form a right-side linear $GL(n_*{}^*A)$-representation.

It remains to prove that the kernel of inefficiency consists only of identity. Identity transformation satisfies to equation

$$v_i = v_j a^j_i$$

Choosing values of coordinates as $v_i = \delta^k_i$ where we selected $k$ we get

$$(12.2.6) \qquad \delta^k_i = \delta^k_j a^j_i$$

From (12.2.6) it follows

$$\delta^k_i = a^k_i$$

Since $k$ is arbitrary, we get the conclusion $a = \delta$. $\qquad\square$

THEOREM 12.2.4. *Automorphism a acting on each vector of basis of left $A$-vector space of rows maps a basis into another basis.*

PROOF. Let $\overline{\overline{e}}$ be basis of left $A$-vector space $V$ of rows. According to theorem 12.1.10, vector $e^i$ maps into a vector $e'^i$

$$(12.2.9) \qquad e'^i = e^i{}_*{}^*a$$

Let vectors $e'^i$ be linearly dependent. Then $\lambda \neq 0$ in the equality

$$(12.2.10) \qquad \lambda_*{}^*e' = 0$$

From equations (12.2.7) and (12.2.8) it follows that

$$\lambda^*{}_\star e'^*{}_\star a^{-1^*}{}_\star = \lambda^*{}_\star e = 0$$

and $\lambda \neq 0$. This contradicts to the statement that vectors $e_i$ are linearly independent. Therefore vectors $e'_i$ are linearly independent and form basis.    □

From equations (12.2.9) and (12.2.10) it follows that

$$\lambda_\star{}^* e'_\star{}^* a^{-1^*}{}_\star = \lambda_\star{}^* e = 0$$

and $\lambda \neq 0$. This contradicts to the statement that vectors $e^i$ are linearly independent. Therefore vectors $e'^i$ are linearly independent and form basis.    □

Thus we can extend a right-side linear $GL(V_\star)$-representation in left $A$-vector space $V_\star$ to the set of bases of left $A$-vector space $V$.       Transformation of this right-side representation on the set of bases of left $A$-vector space $V$ is called **active transformation** because the homomorphism of the left $A$-vector space induced this transformation (See also definition in the section [4]-14.1-3 as well the definition on the page [2]-214).

According to definition we write the action of the transformation $a \in GL(V_\star)$ on the basis $\overline{\overline{e}}$ as $\overline{\overline{e}}^*{}_\star a$.   Consider the equality

(12.2.11)     $$v^*{}_\star \overline{\overline{e}}^*{}_\star a = v^*{}_\star \overline{\overline{e}}^*{}_\star a$$

The expression $\overline{\overline{e}}^*{}_\star a$ on the left side of the equality (12.2.11) is image of basis $\overline{\overline{e}}$ with respect to active transformation $a$. The expression $v^*{}_\star \overline{\overline{e}}$ on the right side of the equality (12.2.11) is expansion of vector $\overline{v}$ with respect to basis $\overline{\overline{e}}$. Therefore, the expression on the right side of the equality (12.2.11) is image of vector $\overline{v}$ with respect to endomorphism $a$ and the expression on the left side of the equality

$$\boxed{(12.2.11) \quad v^*{}_\star \overline{\overline{e}}^*{}_\star a = v^*{}_\star \overline{\overline{e}}^*{}_\star a}$$

is expansion of image of vector $\overline{v}$ with respect to image of basis $\overline{\overline{e}}$. Therefore, from the equality (12.2.11) it follows that endomorphism $a$ of left $A$-vector space and corresponding active transformation $a$ act synchronously and coordinates $a \circ \overline{v}$ of image of the vector $\overline{v}$ with respect to the image $\overline{\overline{e}}^*{}_\star a$ of the basis $\overline{\overline{e}}$ are the same as coordinates of the vector $\overline{v}$ with respect to the basis $\overline{\overline{e}}$.

According to definition we write the action of the transformation $a \in GL(V_\star)$ on the basis $\overline{\overline{e}}$ as $\overline{\overline{e}}_\star{}^* a$.   Consider the equality

(12.2.12)     $$v_\star{}^* \overline{\overline{e}}_\star{}^* a = v_\star{}^* \overline{\overline{e}}_\star{}^* a$$

The expression $\overline{\overline{e}}_\star{}^* a$ on the left side of the equality (12.2.12) is image of basis $\overline{\overline{e}}$ with respect to active transformation $a$. The expression $v_\star{}^* \overline{\overline{e}}$ on the right side of the equality (12.2.12) is expansion of vector $\overline{v}$ with respect to basis $\overline{\overline{e}}$. Therefore, the expression on the right side of the equality (12.2.12) is image of vector $\overline{v}$ with respect to endomorphism $a$ and the expression on the left side of the equality

$$\boxed{(12.2.12) \quad v_\star{}^* \overline{\overline{e}}_\star{}^* a = v_\star{}^* \overline{\overline{e}}_\star{}^* a}$$

is expansion of image of vector $\overline{v}$ with respect to image of basis $\overline{\overline{e}}$. Therefore, from the equality (12.2.12) it follows that endomorphism $a$ of left $A$-vector space and corresponding active transformation $a$ act synchronously and coordinates $a \circ \overline{v}$ of image of the vector $\overline{v}$ with respect to the image $\overline{\overline{e}}_\star{}^* a$ of the basis $\overline{\overline{e}}$ are the same as coordinates of the vector $\overline{v}$ with respect to the basis $\overline{\overline{e}}$.

THEOREM 12.2.5. *Active right-side $GL(n^*{}_\star A)$-representation on the set of bases is single transitive representation. The set of bases identified with tragectory $\overline{\overline{e}}^*{}_\star GL(n^*{}_\star A)(V)$ of active right-side $GL(n^*{}_\star A)$-representation is called the* **basis manifold** *of left $A$-vector*

THEOREM 12.2.6. *Active right-side $GL(n_\star{}^* A)$-representation on the set of bases is single transitive representation. The set of bases identified with tragectory $\overline{\overline{e}}_\star{}^* GL(n_\star{}^* A)(V)$ of active right-side $GL(n_\star{}^* A)$-representation is called the* **basis manifold** *of left $A$-vector*

| | |
|---|---|
| *space $V$.* | *space $V$.* |

To prove theorems 12.2.5, 12.2.6, it is sufficient to show that at least one transformation of right-side representation is defined for any two bases and this transformation is unique.

| | |
|---|---|
| Proof. Homomorphism $a$ operating on basis $\overline{\overline{e}}$ has form | Proof. Homomorphism $a$ operating on basis $\overline{\overline{e}}$ has form |
| $$e_i' = e_i{}^*{}_* a$$ | $$e'^i = e^i{}^*{}_* a$$ |
| where $e_i'$ is coordinate matrix of vector $\overline{e}_i'$ relative basis $\overline{\overline{h}}$ and $e_i$ is coordinate matrix of vector $\overline{e}_i$ relative basis $\overline{\overline{h}}$. Therefore, coordinate matrix of image of basis equal to ${}^*{}_*$-product of coordinate matrix of original basis over matrix of automorphism | where $e'^i$ is coordinate matrix of vector $\overline{e}'^i$ relative basis $\overline{h}$ and $e^i$ is coordinate matrix of vector $\overline{e}^i$ relative basis $\overline{h}$. Therefore, coordinate matrix of image of basis equal to ${}_*{}^*$-product of coordinate matrix of original basis over matrix of automorphism |
| $$e' = e^*{}_* g$$ | $$e' = e_*{}^* g$$ |
| According to the theorem 12.1.7, matrices $g$ and $e$ are nonsingular. Therefore, matrix | According to the theorem 12.1.8, matrices $g$ and $e$ are nonsingular. Therefore, matrix |
| $$a = e^{-1}{}^*{}_*{}^* e'$$ | $$a = e^{-1}{}_*{}^*{}_* e'$$ |
| is the matrix of automorphism mapping basis $\overline{\overline{e}}$ to basis $\overline{\overline{e}}'$. | is the matrix of automorphism mapping basis $\overline{\overline{e}}$ to basis $\overline{\overline{e}}'$. |
| Suppose elements $g_1$, $g_2$ of group $G$ and basis $\overline{\overline{e}}$ satisfy equation | Suppose elements $g_1$, $g_2$ of group $G$ and basis $\overline{\overline{e}}$ satisfy equation |
| $$e^*{}_* g_1 = e^*{}_* g_2$$ | $$e_*{}^* g_1 = e_*{}^* g_2$$ |
| According to the theorems 12.1.7 and 8.1.32, we get $g_1 = g_2$. This proves statement of theorem.                      □ | According to the theorems 12.1.8 and 8.1.32, we get $g_1 = g_2$. This proves statement of theorem.                      □ |

Let us define an additional structure on left $A$-vector space $V$. Then not every automorphism keeps properties of the selected structure. For imstance, if we introduce norm in left $A$-vector space $V$, then we are interested in automorphisms which preserve the norm of the vector.

| | |
|---|---|
| Definition 12.2.7. *Normal subgroup $G(V)$ of the group $GL(V)$ such that subgroup $G(V)$ generates automorphisms which hold properties of the selected structure is called* **symmetry group**. | Definition 12.2.8. *Normal subgroup $G(V)$ of the group $GL(V)$ such that subgroup $G(V)$ generates automorphisms which hold properties of the selected structure is called* **symmetry group**. |
| *Without loss of generality we identify element $g$ of group $G(V)$ with corresponding transformation of representation and write its action on vector $v \in V$ as $v^*{}_* g$.*                      □ | *Without loss of generality we identify element $g$ of group $G(V)$ with corresponding transformation of representation and write its action on vector $v \in V$ as $v_*{}^* g$.*                      □ |

The right-side representation of group $G$ in left $A$-vector space is called **linear $G$-representation**.

---

DEFINITION 12.2.9. *The right-side representation of group $G(V)$ in the set of bases of left $A$-vector space $V$ is called* **active left-side $G$-representation**.          □

---

If the basis $\overline{\overline{e}}$ is given, then we can identify the automorphism of left $A$-vector space $V$ and its coordinates with respect to the basis $\overline{\overline{e}}$. The set $G(V_*)$ of coordinates of automorphisms with respect to the basis $\overline{\overline{e}}$ is group isomorphic to the group $G(V)$.

Not every two bases can be mapped by a transformation from the symmetry group because not every nonsingular linear transformation belongs to the representation of group $G(V)$. Therefore, we can represent the set of bases as a union of orbits of group $G(V)$. In particular, if the basis $\overline{\overline{e}} \in G(V)$, then the group $G(V)$ is orbit of the basis $\overline{\overline{e}}$.

---

DEFINITION 12.2.10. *We call orbit $\overline{\overline{e}}{}^*_* G(V)$ of the selected basis $\overline{\overline{e}}$ the* **basis manifold** *of left $A$-vector space $V$ of columns.*          □

DEFINITION 12.2.11. *We call orbit $\overline{\overline{e}}{}_*^* G(V)$ of the selected basis $\overline{\overline{e}}$ the* **basis manifold** *of left $A$-vector space $V$ of rows.*          □

---

THEOREM 12.2.12. *Active right-side $G(V)$-representation on basis manifold is single transitive representation.*

---

PROOF. The theorem follows from the theorem 12.2.5 and the definition 12.2.10, as well the theorem follows from the theorem 12.2.6 and the definition 12.2.11.          □

---

THEOREM 12.2.13. *On the basis manifold $\overline{\overline{e}}{}^*_* G(V)$ of left $A$-vector space of columns, there exists single transitive left-side $G(V)$-representation, commuting with active.*

THEOREM 12.2.14. *On the basis manifold $\overline{\overline{e}}{}_*^* G(V)$ of left $A$-vector space of rows, there exists single transitive left-side $G(V)$-representation, commuting with active.*

---

PROOF. Theorem 12.2.12 means that the basis manifold $\overline{\overline{e}}{}^*_* G(V)$ is a homogenous space of group $G$. The theorem follows from the theorem 2.5.17.          □

PROOF. Theorem 12.2.12 means that the basis manifold $\overline{\overline{e}}{}_*^* G(V)$ is a homogenous space of group $G$. The theorem follows from the theorem 2.5.17.          □

---

As we see from remark 2.5.18 transformation of left-side $G(V)$-representation is different from an active transformation and cannot be reduced to transformation of space $V$.

---

DEFINITION 12.2.15. *A transformation of left-side $G(V)$-representation is called* **passive transformation** *of basis manifold $\overline{\overline{e}}{}^*_* G(V)$ of left $A$-vector*

DEFINITION 12.2.16. *A transformation of left-side $G(V)$-representation is called* **passive transformation** *of basis manifold $\overline{\overline{e}}{}_*^* G(V)$ of left $A$-vector space* *space of columns, and the left-side $G(V)$-representation is called* **passive left-side $G(V)$-representation**. *According to the definition we write the passive transformation of basis $\overline{\overline{e}}$ defined by element $a \in G(V)$ as* $a^*{}_*\overline{\overline{e}}$. □

*of rows, and the left-side $G(V)$-representation is called* **passive left-side $G(V)$-representation**. *According to the definition we write the passive transformation of basis $\overline{\overline{e}}$ defined by element $a \in G(V)$ as* $a_*{}^*\overline{\overline{e}}$. □

THEOREM 12.2.17. *The coordinate matrix of basis $\overline{\overline{e}}'$ relative basis $\overline{\overline{e}}$ of left $A$-vector space $V$ of columns is identical with the matrix of passive transformation mapping basis $\overline{\overline{e}}$ to basis $\overline{\overline{e}}'$.*

THEOREM 12.2.18. *The coordinate matrix of basis $\overline{\overline{e}}'$ relative basis $\overline{\overline{e}}$ of left $A$-vector space $V$ of rows is identical with the matrix of passive transformation mapping basis $\overline{\overline{e}}$ to basis $\overline{\overline{e}}'$.*

PROOF. According to the theorem 8.4.3, the coordinate matrix of basis $\overline{\overline{e}}'$ relative basis $\overline{\overline{e}}$ consist from rowswhich are coordinate matrices of vectors $\overline{e}'_i$ relative the basis $\overline{\overline{e}}$. Therefore,

$$\overline{e}'_i = e'_i{}^*{}_*\overline{\overline{e}}$$

At the same time the passive transformation $a$ mapping one basis to another has a form

$$\overline{e}'_i = a_i{}^*{}_*\overline{\overline{e}}$$

According to the theorem 8.4.4,

$$e'_i = a_i$$

for any $i$. This proves the theorem. □

PROOF. According to the theorem 8.4.10, the coordinate matrix of basis $\overline{\overline{e}}'$ relative basis $\overline{\overline{e}}$ consist from colswhich are coordinate matrices of vectors $\overline{e}'_i$ relative the basis $\overline{\overline{e}}$. Therefore,

$$\overline{e}'^i = e'^i{}_*{}^*\overline{\overline{e}}$$

At the same time the passive transformation $a$ mapping one basis to another has a form

$$\overline{e}'^i = a^i{}_*{}^*\overline{\overline{e}}$$

According to the theorem 8.4.11,

$$e'^i = a^i$$

for any $i$. This proves the theorem. □

REMARK 12.2.19. *According to theorems 2.5.15, 12.2.17, we can identify a basis $\overline{\overline{e}}'$ of the basis manifold $\overline{\overline{e}}^*{}_*G(V)$ of left $A$-vector space of columns and the matrix $g$ of coordinates of the basis $\overline{\overline{e}}'$ with respect to the basis $\overline{\overline{e}}$*

$$(12.2.13) \qquad \overline{\overline{e}}' = g^*{}_*\overline{\overline{e}}$$

*Since $\overline{\overline{e}}' = \overline{\overline{e}}$, then the equality (12.2.13) gets form*

$$(12.2.14) \qquad \overline{\overline{e}} = E_n{}^*{}_*\overline{\overline{e}}$$

*Based on the equality (12.2.14), we will use notation $\overline{\overline{E_n}}{}^*{}_*G(V)$ for basis manifold $\overline{\overline{e}}^*{}_*G(V)$.* □

REMARK 12.2.20. *According to theorems 2.5.15, 12.2.18, we can identify a basis $\overline{\overline{e}}'$ of the basis manifold $\overline{\overline{e}}_*{}^*G(V)$ of left $A$-vector space of columns and the matrix $g$ of coordinates of the basis $\overline{\overline{e}}'$ with respect to the basis $\overline{\overline{e}}$*

$$(12.2.15) \qquad \overline{\overline{e}}' = g_*{}^*\overline{\overline{e}}$$

*Since $\overline{\overline{e}}' = \overline{\overline{e}}$, then the equality (12.2.15) gets form*

$$(12.2.16) \qquad \overline{\overline{e}} = E_n{}_*{}^*\overline{\overline{e}}$$

*Based on the equality (12.2.16), we will use notation $\overline{\overline{E_n}}{}_*{}^*G(V)$ for basis manifold $\overline{\overline{e}}_*{}^*G(V)$.* □

### 12.3. Passive Transformation in Left $A$-Vector Space

An active transformation changes bases and vectors uniformly and coordinates of vector relative basis do not change. A passive transformation changes only the basis and it leads to transformation of coordinates of vector relative to basis.

**Theorem 12.3.1.** *Let $V$ be a left $A$-vector space of columns. Let passive transformation $g \in G(V_*)$ map the basis $\overline{\overline{e}}_1$ into the basis $\overline{\overline{e}}_2$*

$$(12.3.1) \qquad \overline{\overline{e}}_2 = g^* {}_* \overline{\overline{e}}_1$$

*Let*

$$(12.3.2) \qquad v_i = \begin{pmatrix} v_i^1 \\ ... \\ v_i^n \end{pmatrix}$$

*be matrix of coordinates of the vector $\overline{v}$ with respect to the basis $\overline{\overline{e}}_i$, $i = 1, 2$. Coordinate transformations*

$$(12.3.3) \qquad v_1 = v_2 {}^* {}_* g$$

$$(12.3.4) \qquad v_2 = v_1 {}^* {}_* g^{-1^*}$$

*do not depend on vector $\overline{v}$ or basis $\overline{\overline{e}}$, but is defined only by coordinates of vector $\overline{v}$ relative to basis $\overline{\overline{e}}$.*

**Proof.** According to the theorem 8.4.3,

$$(12.3.8) \qquad \overline{v} = v_1 {}^* {}_* \overline{\overline{e}}_1 = v_2 {}^* {}_* \overline{\overline{e}}_2$$

The equality

$$(12.3.9) \qquad v_1 {}^* {}_* \overline{\overline{e}}_1 = v_2 {}^* {}_* g^* {}_* \overline{\overline{e}}_1$$

follows from equalities (12.3.1), (12.3.8). According to the theorem 8.4.4, the coordinate transformation (12.3.3) follows from the equality (12.3.9). Because $g$ is $_*$-nonsingular matrix, coordinate transformation (12.3.4) follows from the equality (12.3.3). □

**Theorem 12.3.3.** *Coordinate transformations*

$$\boxed{(12.3.4) \qquad v_2 = v_1 {}^* {}_* g^{-1^*}}$$

*form effective linear right-side contravariant $G$-representation which is called* **coordinate representation in**

**Theorem 12.3.2.** *Let $V$ be a left $A$-vector space of rows. Let passive transformation $g \in G(V_*)$ map the basis $\overline{\overline{e}}_1$ into the basis $\overline{\overline{e}}_2$*

$$(12.3.5) \qquad \overline{\overline{e}}_2 = g_* {}^* \overline{\overline{e}}_1$$

*Let*

$$v_i = \begin{pmatrix} v_{i\,1} & ... & v_{i\,n} \end{pmatrix}$$

*be matrix of coordinates of the vector $\overline{v}$ with respect to the basis $\overline{\overline{e}}_i$, $i = 1, 2$. Coordinate transformations*

$$(12.3.6) \qquad v_1 = v_2 {}_* {}^* g$$

$$(12.3.7) \qquad v_2 = v_1 {}_* {}^* g^{-1^*}$$

*do not depend on vector $\overline{v}$ or basis $\overline{\overline{e}}$, but is defined only by coordinates of vector $\overline{v}$ relative to basis $\overline{\overline{e}}$.*

**Proof.** According to the theorem 8.4.10,

$$(12.3.10) \qquad \overline{v} = v_1 {}_* {}^* \overline{\overline{e}}_1 = v_2 {}_* {}^* \overline{\overline{e}}_2$$

The equality

$$(12.3.11) \qquad v_1 {}_* {}^* \overline{\overline{e}}_1 = v_2 {}_* {}^* g_* {}^* \overline{\overline{e}}_1$$

follows from equalities (12.3.5), (12.3.10). According to the theorem 8.4.11, the coordinate transformation (12.3.6) follows from the equality (12.3.11). Because $g$ is $_*$-nonsingular matrix, coordinate transformation (12.3.7) follows from the equality (12.3.6). □

**Theorem 12.3.4.** *Coordinate transformations*

$$\boxed{(12.3.7) \qquad v_2 = v_1 {}_* {}^* g^{-1^*}}$$

*form effective linear right-side contravariant $G$-representation which is called* **coordinate representation in**

**left $A$-vector space** *of columns.*

Proof.    Let

$$(12.3.12) \qquad v_i = \begin{pmatrix} v_i^1 \\ ... \\ v_i^n \end{pmatrix}$$

be matrix of coordinates of the vector $\overline{v}$ with respect to the basis $\overline{\overline{e}}_i$, $i = 1, 2, 3$. Then

$$(12.3.13) \qquad \overline{v} = v_1{}^*{}_* \overline{\overline{e}}_1 = v_2{}^*{}_* \overline{\overline{e}}_2$$
$$= v_3{}^*{}_* \overline{\overline{e}}_3$$

Let passive transformation $g_1$ map the basis $\overline{\overline{e}}_1 \in \overline{e}^*{}_* G(V)$ into the basis $\overline{\overline{e}}_2 \in \overline{\overline{e}}^*{}_* G(V)$

$$(12.3.14) \qquad \overline{\overline{e}}_2 = g_1{}^*{}_* \overline{\overline{e}}_1$$

Let passive transformation $g_2$ map the basis $\overline{\overline{e}}_2 \in \overline{e}^*{}_* G(V)$ into the basis $\overline{\overline{e}}_3 \in \overline{\overline{e}}^*{}_* G(V)$

$$(12.3.15) \qquad \overline{\overline{e}}_3 = g_2{}^*{}_* \overline{\overline{e}}_2$$

The transformation

$$(12.3.16) \qquad \overline{\overline{e}}_3 = (g_2{}^*{}_* g_1){}^*{}_* \overline{\overline{e}}_1$$

is the product of transformations (12.3.14), (12.3.15). The coordinate transformation

$$(12.3.17) \qquad v_2 = v_1{}^*{}_* g_1^{-1}{}^*{}_*$$

corresponds to passive transformation $g_1$ and follows from equalities (12.3.13), (12.3.14) and the theorem 12.3.1. Similarly, according to the theorem 12.3.1, the coordinate transformation

$$(12.3.18) \qquad v_3 = v_2{}^*{}_* g_2^{-1}{}^*{}_*$$

corresponds to passive transformation $g_2$ and follows from equalities (12.3.13), (12.3.15). The product of coordinate transformations (12.3.17), (12.3.18) has form

$$(12.3.19) \qquad v_3 = v_1{}^*{}_* g_1^{-1}{}^*{}_*{}^*{}_* g_2^{-1}{}^*{}_*$$

and is coordinate transformation corresponding to passive transformation (12.3.16). The equality

$$(12.3.20) \qquad v_3 = v_1{}^*{}_* (g_2{}^*{}_* g_1)^{-1}{}^*{}_*$$

**left $A$-vector space** *of rows.*

Proof.    Let

$$v_i = \begin{pmatrix} v_{i\,1} & ... & v_{i\,n} \end{pmatrix}$$

be matrix of coordinates of the vector $\overline{v}$ with respect to the basis $\overline{\overline{e}}_i$, $i = 1, 2, 3$. Then

$$(12.3.21) \qquad \overline{v} = v_{1*}{}^* \overline{\overline{e}}_1 = v_{2*}{}^* \overline{\overline{e}}_2$$
$$= v_{3*}{}^* \overline{\overline{e}}_3$$

Let passive transformation $g_1$ map the basis $\overline{\overline{e}}_1 \in \overline{e}_*{}^* G(V)$ into the basis $\overline{\overline{e}}_2 \in \overline{\overline{e}}_*{}^* G(V)$

$$(12.3.22) \qquad \overline{\overline{e}}_2 = g_{1*}{}^* \overline{\overline{e}}_1$$

Let passive transformation $g_2$ map the basis $\overline{\overline{e}}_2 \in \overline{e}_*{}^* G(V)$ into the basis $\overline{\overline{e}}_3 \in \overline{\overline{e}}_*{}^* G(V)$

$$(12.3.23) \qquad \overline{\overline{e}}_3 = g_{2*}{}^* \overline{\overline{e}}_2$$

The transformation

$$(12.3.24) \qquad \overline{\overline{e}}_3 = (g_{2*}{}^* g_1)_*{}^* \overline{\overline{e}}_1$$

is the product of transformations (12.3.22), (12.3.23). The coordinate transformation

$$(12.3.25) \qquad v_2 = v_{1*}{}^* g_1^{-1}{}_*{}^*$$

corresponds to passive transformation $g_1$ and follows from equalities (12.3.21), (12.3.22) and the theorem 12.3.2. Similarly, according to the theorem 12.3.2, the coordinate transformation

$$(12.3.26) \qquad v_3 = v_{2*}{}^* g_2^{-1}{}_*{}^*$$

corresponds to passive transformation $g_2$ and follows from equalities (12.3.21), (12.3.23). The product of coordinate transformations (12.3.17), (12.3.18) has form

$$(12.3.27) \qquad v_3 = v_{1*}{}^* g_1^{-1}{}_*{}^*{}_*{}^* g_2^{-1}{}_*{}^*$$

and is coordinate transformation corresponding to passive transformation (12.3.24). The equality

$$(12.3.28) \qquad v_3 = v_{1*}{}^* (g_{2*}{}^* g_1)^{-1}{}_*{}^*$$

follows from equalities (8.1.66), (12.3.19). According to the definition 2.5.9, coordinate transformations form linear right-side contravariant $G$-representation. Suppose coordinate transformation does not change vectors $\delta_k$. Then unit of group $G$ corresponds to it because representation is single transitive. Therefore, coordinate representation is effective.  □

follows from equalities (8.1.66), (12.3.27). According to the definition 2.5.9, coordinate transformations form linear right-side contravariant $G$-representation. Suppose coordinate transformation does not change vectors $\delta_k$. Then unit of group $G$ corresponds to it because representation is single transitive. Therefore, coordinate representation is effective.  □

**Theorem 12.3.5.** *Let $V$ be a left $A$-vector space of columns and $\overline{\overline{e}}_1$, $\overline{\overline{e}}_2$ be bases in left $A$-vector space $V$. Let $g$ be passive transformation of basis $\overline{\overline{e}}_1$ into basis $\overline{\overline{e}}_2$*

$$\overline{\overline{e}}_2 = g^*{}_*\overline{\overline{e}}_1$$

*Let $f$ be endomorphism of left $A$-vector space $V$. Let $f_i$, $i = 1, 2$, be the matrix of endomorphism $f$ with respect to the basis $\overline{\overline{e}}_i$. Then*

$$(12.3.29) \qquad f_2 = g^*{}_*f_1{}^*{}_*g^{-1}{}^{*}{}^{*}$$

**Theorem 12.3.6.** *Let $V$ be a left $A$-vector space of rows and $\overline{\overline{e}}_1$, $\overline{\overline{e}}_2$ be bases in left $A$-vector space $V$. Let $g$ be passive transformation of basis $\overline{\overline{e}}_1$ into basis $\overline{\overline{e}}_2$*

$$\overline{\overline{e}}_2 = g_*{}^*\overline{\overline{e}}_1$$

*Let $f$ be endomorphism of left $A$-vector space $V$. Let $f_i$, $i = 1, 2$, be the matrix of endomorphism $f$ with respect to the basis $\overline{\overline{e}}_i$. Then*

$$(12.3.30) \qquad f_2 = g_*{}^*f_1{}_*{}^*g^{-1}{}_*{}^{*}$$

**Proof of the Theorem 12.3.5.** Let endomorphism $f$ maps a vector $v$ into a vector $w$

$$(12.3.31) \qquad w = f \circ v$$

Let $v_i$, $i = 1, 2$, be the matrix of vector $v$ with respect to the basis $\overline{\overline{e}}_i$. Then

$$(12.3.32) \qquad v_1 = v_2{}^*{}_*g$$

follows from the theorem 12.3.1.

Let $w_i$, $i = 1, 2$, be the matrix of vector $w$ with respect to the basis $\overline{\overline{e}}_i$. Then

$$(12.3.33) \qquad w_1 = w_2{}^*{}_*g$$

follows from the theorem 12.3.1.

According to the theorem 10.3.6, equalities

$$(12.3.34) \qquad w_1 = v_1{}^*{}_*f_1$$

$$(12.3.35) \qquad w_2 = v_2{}^*{}_*f_2$$

follow from the equality (12.3.31). The equality

$$(12.3.36) \qquad w_2{}^*{}_*g = v_1{}^*{}_*f_1 = v_2{}^*{}_*g^*{}_*f_1$$

follows from equalities (12.3.32), (12.3.34), (12.3.33). Since $g$ is regular matrix, then the equality

$$(12.3.37) \qquad w_2 = v_1{}^*{}_*f_1 = v_2{}^*{}_*g^*{}_*f_1{}^*{}_*g^{-1}{}^{*}$$

follows from the equality (12.3.36). The equality

$$(12.3.38) \qquad v_2{}^*{}_*f_2 = v_2{}^*{}_*g^*{}_*f_1{}^*{}_*g^{-1}{}^{*}$$

follows from equalities (12.3.35), (12.3.37). The equality (12.3.29) follows from the equality (12.3.38).                                                                    □

PROOF OF THE THEOREM 12.3.6. Let endomorphism $f$ maps a vector $v$ into a vector $w$

$$(12.3.39) \qquad\qquad w = f \circ v$$

Let $v_i$, $i = 1, 2$, be the matrix of vector $v$ with respect to the basis $\overline{\overline{e}}_i$. Then

$$(12.3.40) \qquad v_1 = v_{2*}{}^* g$$

follows from the theorem 12.3.2.

Let $w_i$, $i = 1, 2$, be the matrix of vector $w$ with respect to the basis $\overline{\overline{e}}_i$. Then

$$(12.3.41) \qquad w_1 = w_{2*}{}^* g$$

follows from the theorem 12.3.2.

According to the theorem 10.3.6, equalities

$$(12.3.42) \qquad\qquad w_1 = v_{1*}{}^* f_1$$

$$(12.3.43) \qquad\qquad w_2 = v_{2*}{}^* f_2$$

follow from the equality (12.3.39). The equality

$$(12.3.44) \qquad\qquad w_{2*}{}^* g = v_{1*}{}^* f_1 = v_{2*}{}^* g_*{}^* f_1$$

follows from equalities (12.3.40), (12.3.42), (12.3.41). Since $g$ is regular matrix, then the equality

$$(12.3.45) \qquad\qquad w_2 = v_{1*}{}^* f_1 = v_{2*}{}^* g_*{}^* f_{1*}{}^* g^{-1*}$$

follows from the equality (12.3.44). The equality

$$(12.3.46) \qquad\qquad v_{2*}{}^* f_2 = v_{2*}{}^* g_*{}^* f_{1*}{}^* g^{-1*}$$

follows from equalities (12.3.43), (12.3.45). The equality (12.3.30) follows from the equality (12.3.46).                                                                    □

THEOREM 12.3.7. *Let $V$ be a left $A$-vector space of columns. Let endomorphism $\overline{f}$ of left $A$-vector space $V$ is sum of endomorphisms $\overline{f}_1$, $\overline{f}_2$. The matrix $f$ of endomorphism $\overline{f}$ equal to the sum of matrices $f_1$, $f_2$ of endomorphisms $\overline{f}_1$, $\overline{f}_2$ and this equality does not depend on the choice of a basis.*

THEOREM 12.3.8. *Let $V$ be a left $A$-vector space of rows. Let endomorphism $\overline{f}$ of left $A$-vector space $V$ is sum of endomorphisms $\overline{f}_1$, $\overline{f}_2$. The matrix $f$ of endomorphism $\overline{f}$ equal to the sum of matrices $f_1$, $f_2$ of endomorphisms $\overline{f}_1$, $\overline{f}_2$ and this equality does not depend on the choice of a basis.*

PROOF. Let $\overline{\overline{e}}$, $\overline{\overline{e}}'$ be bases of left $A$-vector space $V$ and $g$ be passive transformation of basis $\overline{\overline{e}}$ into basis $\overline{\overline{e}}'$

$$e' = e^*{}_* g$$

Let $f$ be the matrix of endomorphism $\overline{f}$ with respect to the basis $\overline{\overline{e}}$. Let $f'$ be the matrix of endomorphism $\overline{f}$ with respect to the basis $\overline{\overline{e}}'$. Let $f_i$, $i = 1, 2$, be the matrix of endomorphism $\overline{f}_i$ with

PROOF. Let $\overline{\overline{e}}$, $\overline{\overline{e}}'$ be bases of left $A$-vector space $V$ and $g$ be passive transformation of basis $\overline{\overline{e}}$ into basis $\overline{\overline{e}}'$

$$e' = e_*{}^* g$$

Let $f$ be the matrix of endomorphism $\overline{f}$ with respect to the basis $\overline{\overline{e}}$. Let $f'$ be the matrix of endomorphism $\overline{f}$ with respect to the basis $\overline{\overline{e}}'$. Let $f_i$, $i = 1, 2$, be the matrix of endomorphism $\overline{f}_i$ with respect to the basis $\overline{\overline{e}}$. Let $f'_i$, $i = 1, 2$, be the matrix of endomorphism $\overline{f}_i$ with respect to the basis $\overline{\overline{e}}'$.    According to the theorem 10.4.2,

(12.3.47) $$f = f_1 + f_2$$

(12.3.48) $$f' = f'_1 + f'_2$$

According to the theorem 12.3.5,

(12.3.49) $$f' = g^*{}_*f^*{}_*g^{-1^*{}^*}$$

(12.3.50) $$f'_1 = g^*{}_*f_1{}^*{}_*g^{-1^*{}^*}$$

(12.3.51) $$f'_2 = g^*{}_*f_2{}^*{}_*g^{-1^*{}^*}$$

The equality

$$
\begin{aligned}
f' &= g^*{}_*f^*{}_*g^{-1^*{}^*}\\
&= g^*{}_*(f_1+f_2)^*{}_*g^{-1^*{}^*}\\
(12.3.52)\quad &= g^*{}_*f_1{}^*{}_*g^{-1^*{}^*}\\
&+ g^*{}_*f_2{}^*{}_*g^{-1^*{}^*}\\
&= f'_1 + f'_2
\end{aligned}
$$

follows from equalities (12.3.47), (12.3.49), (12.3.50), (12.3.51). From equalities (12.3.48), (12.3.52) it follows that the equality (12.3.47) is covariant with respect to choice of basis of left $A$-vector space $V$.  □

**Theorem 12.3.9.** *Let $V$ be a left $A$-vector space of columns. Let endomorphism $\overline{f}$ of left $A$-vector space $V$ is product of endomorphisms $\overline{f}_1$, $\overline{f}_2$. The matrix $f$ of endomorphism $\overline{f}$ equal to the $_*{}^*$-product of matrices $f_1$, $f_2$ of endomorphisms $\overline{f}_1$, $\overline{f}_2$ and this equality does not depend on the choice of a basis.*

**Proof.** Let $\overline{\overline{e}}$, $\overline{\overline{e}}'$ be bases of left $A$-vector space $V$ and $g$ be passive transformation of basis $\overline{\overline{e}}$ into basis $\overline{\overline{e}}'$

$$e' = e^*{}_*g$$

Let $f$ be the matrix of endomorphism $\overline{f}$ with respect to the basis $\overline{\overline{e}}$. Let $f'$ be the matrix of endomorphism $\overline{f}$ with respect to the basis $\overline{\overline{e}}'$. Let $f_i$, $i = 1, 2$, be the matrix of endomorphism $\overline{f}_i$ with respect to the basis $\overline{\overline{e}}$. Let $f'_i$, $i = 1, 2$, respect to the basis $\overline{\overline{e}}$. Let $f'_i$, $i = 1, 2$, be the matrix of endomorphism $\overline{f}_i$ with respect to the basis $\overline{\overline{e}}'$.    According to the theorem 10.4.3,

(12.3.53) $$f = f_1 + f_2$$

(12.3.54) $$f' = f'_1 + f'_2$$

According to the theorem 12.3.6,

(12.3.55) $$f' = g_*{}^*f_*{}^*g^{-1_*{}^*}$$

(12.3.56) $$f'_1 = g_*{}^*f_1{}_*{}^*g^{-1_*{}^*}$$

(12.3.57) $$f'_2 = g_*{}^*f_2{}_*{}^*g^{-1_*{}^*}$$

The equality

$$
\begin{aligned}
f' &= g_*{}^*f_*{}^*g^{-1_*{}^*}\\
&= g_*{}^*(f_1+f_2)_*{}^*g^{-1_*{}^*}\\
(12.3.58)\quad &= g_*{}^*f_1{}_*{}^*g^{-1_*{}^*}\\
&+ g_*{}^*f_2{}_*{}^*g^{-1_*{}^*}\\
&= f'_1 + f'_2
\end{aligned}
$$

follows from equalities (12.3.53), (12.3.55), (12.3.56), (12.3.57). From equalities (12.3.54), (12.3.58) it follows that the equality (12.3.53) is covariant with respect to choice of basis of left $A$-vector space $V$.  □

**Theorem 12.3.10.** *Let $V$ be a left $A$-vector space of rows. Let endomorphism $\overline{f}$ of left $A$-vector space $V$ is product of endomorphisms $\overline{f}_1$, $\overline{f}_2$. The matrix $f$ of endomorphism $\overline{f}$ equal to the $_*{}^*$-product of matrices $f_1$, $f_2$ of endomorphisms $\overline{f}_1$, $\overline{f}_2$ and this equality does not depend on the choice of a basis.*

**Proof.** Let $\overline{\overline{e}}$, $\overline{\overline{e}}'$ be bases of left $A$-vector space $V$ and $g$ be passive transformation of basis $\overline{\overline{e}}$ into basis $\overline{\overline{e}}'$

$$e' = e_*{}^*g$$

Let $f$ be the matrix of endomorphism $\overline{f}$ with respect to the basis $\overline{\overline{e}}$. Let $f'$ be the matrix of endomorphism $\overline{f}$ with respect to the basis $\overline{\overline{e}}'$. Let $f_i$, $i = 1, 2$, be the matrix of endomorphism $\overline{f}_i$ with respect to the basis $\overline{\overline{e}}$. Let $f'_i$, $i = 1, 2$, be the matrix of endomorphism $\overline{f}_i$ with respect to the basis $\overline{\overline{e}}'$. According to the theorem 10.5.2,

$$(12.3.59) \qquad f = f_{1*}{}^* f_2$$

$$(12.3.60) \qquad f' = f'_{1*}{}^* f'_2$$

According to the theorem 12.3.5,

$$(12.3.61) \qquad f' = g^*{}_* f^*{}_* g^{-1^*{}^*}$$

$$(12.3.62) \qquad f'_1 = g^*{}_* f_1{}^*{}_* g^{-1^*{}^*}$$

$$(12.3.63) \qquad f'_2 = g^*{}_* f_2{}^*{}_* g^{-1^*{}^*}$$

The equality

$$
\begin{aligned}
f' &= g^*{}_* f^*{}_* g^{-1^*{}^*}\\
&= g^*{}_* f_1{}^*{}_* f_2{}^*{}_* g^{-1^*{}^*}\\
(12.3.64) \qquad &= g^*{}_* f_1{}^*{}_* g^{-1^*{}^*}\\
&\quad {}_*{}^* g^*{}_* f_2{}^*{}_* g^{-1^*{}^*}\\
&= f'_1{}^*{}_* f'_2
\end{aligned}
$$

follows from equalities (12.3.59), (12.3.61), (12.3.62), (12.3.63). From equalities (12.3.60), (12.3.64) it follows that the equality (12.3.59) is covariant with respect to choice of basis of left $A$-vector space $V$. $\qquad\square$

be the matrix of endomorphism $\overline{f}_i$ with respect to the basis $\overline{\overline{e}}'$. According to the theorem 10.5.3,

$$(12.3.65) \qquad f = f_1{}^*{}_* f_2$$

$$(12.3.66) \qquad f' = f'_1{}^*{}_* f'_2$$

According to the theorem 12.3.6,

$$(12.3.67) \qquad f' = g_*{}^* f_*{}^* g^{-1_*{}^*}$$

$$(12.3.68) \qquad f'_1 = g_*{}^* f_{1*}{}^* g^{-1_*{}^*}$$

$$(12.3.69) \qquad f'_2 = g_*{}^* f_{2*}{}^* g^{-1_*{}^*}$$

The equality

$$
\begin{aligned}
f' &= g_*{}^* f_*{}^* g^{-1_*{}^*}\\
&= g_*{}^* f_{1*}{}^* f_{2*}{}^* g^{-1_*{}^*}\\
(12.3.70) \qquad &= g_*{}^* f_{1*}{}^* g^{-1_*{}^*}\\
&\quad {}_*{}^* g_*{}^* f_{2*}{}^* g^{-1_*{}^*}\\
&= f'_{1*}{}^* f'_2
\end{aligned}
$$

follows from equalities (12.3.65), (12.3.67), (12.3.68), (12.3.69). From equalities (12.3.66), (12.3.70) it follows that the equality (12.3.65) is covariant with respect to choice of basis of left $A$-vector space $V$. $\qquad\square$

## 12.4. Dimension of Right $A$-Vector Space

**Theorem 12.4.1.** *Let $V$ be a right $A$-vector space of columns. Let the set of vectors $\overline{\overline{e}} = (e_i, i \in I)$ be a basis of right $A$-vector space $V$. Let the set of vectors $\overline{\overline{e}} = (e_j, j \in J)$ be a basis of right $A$-vector space $V$. If $|I|$ and $|J|$ are finite numbers then $|I| = |J|$.*

**Proof.** Let $|I| = m$ and $|J| = n$. Let

$$(12.4.1) \qquad m < n$$

Because $\overline{e}_1$ is a basis, any vector $e_{2j}$, $j \in J$, has expansion

$$e_{2j} = e_{1*}{}^* a_j$$

Because $\overline{\overline{e}}_2$ is a basis,

$$(12.4.2) \qquad \lambda = 0$$

**Theorem 12.4.2.** *Let $V$ be a right $A$-vector space of rows. Let the set of vectors $\overline{\overline{e}} = (e^i, i \in I)$ be a basis of right $A$-vector space $V$. Let the set of vectors $\overline{\overline{e}} = (e^j, j \in J)$ be a basis of right $A$-vector space $V$. If $|I|$ and $|J|$ are finite numbers then $|I| = |J|$.*

**Proof.** Let $|I| = m$ and $|J| = n$. Let

$$(12.4.5) \qquad m < n$$

Because $\overline{e}_1$ is a basis, any vector $e_2^j$, $j \in J$, has expansion

$$e_2^j = e_1{}^*{}_* a^j$$

Because $\overline{\overline{e}}_2$ is a basis,

$$(12.4.6) \qquad \lambda = 0$$

should follow from

$$(12.4.3) \qquad e_{2*}{}^*\lambda = e_{1*}{}^*a_*{}^*\lambda = 0$$

Because $\overline{\overline{e}}_1$ is a basis we get

$$(12.4.4) \qquad a_*{}^*\lambda = 0$$

According to (12.4.1), $\operatorname{rank}({}_*{}^*)a \le m$ and system of linear equations (12.4.4) has more variables then equations. According to the theorem 9.4.4, we get the statement $\lambda \ne 0$ which contradicts statement (12.4.2). Therefore, statement $m < n$ is not valid. In the same manner we can prove that the statement $n < m$ is not valid. This completes the proof of the theorem. □

should follow from

$$(12.4.7) \qquad e_2{}^*{}_*\lambda = e_1{}^*{}_*a^*{}_*\lambda = 0$$

Because $\overline{\overline{e}}_1$ is a basis we get

$$(12.4.8) \qquad a^*{}_*\lambda = 0$$

According to (12.4.5), $\operatorname{rank}({}^*{}_*)a \le m$ and system of linear equations (12.4.8) has more variables then equations. According to the theorem 9.4.4, we get the statement $\lambda \ne 0$ which contradicts statement (12.4.6). Therefore, statement $m < n$ is not valid. In the same manner we can prove that the statement $n < m$ is not valid. This completes the proof of the theorem. □

DEFINITION 12.4.3. *We call* **dimension of right $A$-vector space** *the number of vectors in a basis.* □

THEOREM 12.4.4. *If right $A$-vector space $V$ has finite dimension, then the statement*

$$(12.4.9) \qquad \forall v \in V, va = vb \quad a, b \in A$$

*implies $a = b$.*

PROOF. Let the set of vectors $\overline{\overline{e}} = (e_j, j \in J)$ be the basis of left $A$-vector space $V$. According to theorems 7.4.6, 8.4.18, the equality

$$(12.4.11) \qquad v^i a = v^i b$$

follows from the equality (12.4.9). The theorem follows from the equality (12.4.11), if we assume $v = e_1$. □

THEOREM 12.4.5. *If right $A$-vector space $V$ has finite dimension, then the statement*

$$(12.4.10) \quad \forall v \in V, va = vb \quad a, b \in A$$

*implies $a = b$.*

PROOF. Let the set of vectors $\overline{\overline{e}} = (e_j, j \in J)$ be the basis of left $A$-vector space $V$. According to theorems 7.4.6, 8.4.25, the equality

$$(12.4.12) \qquad v_i a = v_i b$$

follows from the equality (12.4.10). The theorem follows from the equality (12.4.12), if we assume $v = e_1$. □

THEOREM 12.4.6. *If right $A$-vector space $V$ has finite dimension, then corresponding representation is free.*

PROOF. The theorem follows from the definition 2.3.3 and from theorems 12.4.4, 12.4.5. □

THEOREM 12.4.7. *Let $V$ be a right $A$-vector space of columns. The coordinate matrix of basis $\overline{\overline{g}}$ relative basis $\overline{\overline{e}}$ of right $A$-vector space $V$ is ${}_*{}^*$-nonsingular matrix.*

THEOREM 12.4.8. *Let $V$ be a right $A$-vector space of rows. The coordinate matrix of basis $\overline{\overline{g}}$ relative basis $\overline{\overline{e}}$ of right $A$-vector space $V$ is ${}^*{}_*$-nonsingular matrix.*

PROOF. According to the theorem 9.3.6, $^*{}_*$-rank of the coordinate matrix of basis $\overline{g}$ relative basis $\overline{e}$ equal to the dimension of left $A$-space of columns. This proves the statement of the theorem.

$\square$

PROOF. According to the theorem 9.3.7, $^*{}_*$-rank of the coordinate matrix of basis $\overline{g}$ relative basis $\overline{e}$ equal to the dimension of left $A$-space of columns. This proves the statement of the theorem.

$\square$

THEOREM 12.4.9. *Let $V$ be a right $A$-vector space of columns. Let $\overline{e}$ be a basis of right $A$-vector space $V$. Then any automorphism $\overline{f}$ of right $A$-vector space $V$ has form*

$$(12.4.13) \qquad v' = f_*{}^* v$$

*where $f$ is a $_*{}^*$-nonsingular matrix.*

THEOREM 12.4.10. *Let $V$ be a right $A$-vector space of rows. Let $\overline{e}$ be a basis of right $A$-vector space $V$. Then any automorphism $\overline{f}$ of right $A$-vector space $V$ has form*

$$(12.4.14) \qquad v' = f^*{}_* v$$

*where $f$ is a $^*{}_*$-nonsingular matrix.*

PROOF. The equality (12.4.13) follows from theorem 11.3.6. Because $\overline{f}$ is an isomorphism, for each vector $\overline{v}'$ there exist one and only one vector $\overline{v}$ such that $\overline{v}' = \overline{f} \circ \overline{v}$. Therefore, system of linear equations (12.4.13) has a unique solution. According to corollary 9.4.5 matrix $f$ is a $_*{}^*$-nonsingular matrix. $\square$

PROOF. The equality (12.4.14) follows from theorem 11.3.7. Because $\overline{f}$ is an isomorphism, for each vector $\overline{v}'$ there exist one and only one vector $\overline{v}$ such that $\overline{v}' = \overline{f} \circ \overline{v}$. Therefore, system of linear equations (12.4.14) has a unique solution. According to corollary 9.4.5 matrix $f$ is a $^*{}_*$-nonsingular matrix. $\square$

THEOREM 12.4.11. *Matrices of automorphisms of right $A$-space $V$ of columns form a group $GL(n_*{}^* A)$.*

THEOREM 12.4.12. *Matrices of automorphisms of right $A$-space $V$ of rows form a group $GL(n^*{}_* A)$.*

PROOF. If we have two automorphisms $\overline{f}$ and $\overline{g}$ then we can write

$$v' = f_*{}^* v$$

$$v'' = g_*{}^* v' = g_*{}^* f_*{}^* v$$

Therefore, the resulting automorphism has matrix $f_*{}^* g$. $\square$

PROOF. If we have two automorphisms $\overline{f}$ and $\overline{g}$ then we can write

$$v' = f^*{}_* v$$

$$v'' = g^*{}_* v' = g^*{}_* f^*{}_* v$$

Therefore, the resulting automorphism has matrix $f^*{}_* g$. $\square$

THEOREM 12.4.13. *Let $\overline{\overline{e}}_{V1}$, $\overline{\overline{e}}_{V2}$ be bases of right $A$-vector space $V$ of columns. Let $\overline{\overline{e}}_{W1}$, $\overline{\overline{e}}_{W2}$ be bases of right $A$-vector space $W$ of columns. Let $a_1$ be matrix of homomorphism*

$$(12.4.15) \qquad a : V \to W$$

*relative to bases $\overline{\overline{e}}_{V1}$, $\overline{\overline{e}}_{W1}$ and $a_2$ be matrix of homomorphism (12.4.15) relative to bases $\overline{\overline{e}}_{V2}$, $\overline{\overline{e}}_{W2}$. Suppose the basis $\overline{\overline{e}}_{V1}$ has coordinate matrix $b$ relative*

THEOREM 12.4.14. *Let $\overline{\overline{e}}_{V1}$, $\overline{\overline{e}}_{V2}$ be bases of right $A$-vector space $V$ of columns. Let $\overline{\overline{e}}_{W1}$, $\overline{\overline{e}}_{W2}$ be bases of right $A$-vector space $W$ of columns. Let $a_1$ be matrix of homomorphism*

$$(12.4.17) \qquad a : V \to W$$

*relative to bases $\overline{\overline{e}}_{V1}$, $\overline{\overline{e}}_{W1}$ and $a_2$ be matrix of homomorphism (12.4.17) relative to bases $\overline{\overline{e}}_{V2}$, $\overline{\overline{e}}_{W2}$. Suppose the basis $\overline{\overline{e}}_{V1}$ has coordinate matrix $b$ relative* the basis $\overline{\overline{e}}_{V2}$

$$\overline{\overline{e}}_{V1} = \overline{\overline{e}}_{V2*}{}^*b$$

and $\overline{\overline{e}}_{W1}$ has coordinate matrix c relative the basis $\overline{\overline{e}}_{W2}$

(12.4.16) $\qquad \overline{\overline{e}}_{W1} = \overline{\overline{e}}_{W2*}{}^*c$

Then there is relationship between matrices $a_1$ and $a_2$

$$a_1 = c^{-1*}{}_*{}^*a_{2*}{}^*b$$

PROOF. Vector $\overline{v} \in V$ has expansion

$$\overline{v} = e_{V1}{}^*{}_*v = e_{V2*}{}^*b_*{}^*v$$

relative to bases $\overline{\overline{e}}_{V1}$ , $\overline{\overline{e}}_{V2}$ . Since $a$ is homomorphism, we can write it as

(12.4.19) $\qquad \overline{w} = e_{W1*}{}^*a_{1*}{}^*v$

relative to bases $\overline{\overline{e}}_{V1}$, $\overline{\overline{e}}_{W1}$ and as

(12.4.20) $\qquad \overline{w} = e_{W2*}{}^*a_{2*}{}^*b_*{}^*v$

relative to bases $\overline{\overline{e}}_{V2}$, $\overline{\overline{e}}_{W2}$. According to the theorem 12.4.7, matrix $c$ has $^*_*$-inverse and from equation (12.4.16) it follows that

(12.4.21) $\qquad \overline{\overline{e}}_{W2} = e_{W1*}{}^*\overline{\overline{c}}^{-1*}{}^*$

The equality

(12.4.22) $\qquad \overline{w} = e_{W1*}{}^*c^{-1*}{}_*{}^*a_{2*}{}^*b_*{}^*v$

follows from equalities (12.4.20), (12.4.21). From the theorem 8.4.18 and comparison of equations (12.4.19) and (12.4.22) it follows that

(12.4.23) $\qquad v^*{}_*a_1 = c^{-1*}{}_*{}^*a_{2*}{}^*b_*{}^*v$

Since vector $a$ is arbitrary vector, the statement of theorem follows from the theorem 12.5.1 and equation (12.4.23).

<div align="right">□</div>

THEOREM 12.4.15. *Let $\overline{a}$ be automorphism of left $A$-vector space $V$ of columns. Let $a_1$ be matrix of this automorphism defined relative to basis $\overline{\overline{e}}_1$ and $a_2$ be matrix of the same automorphism defined relative to basis $\overline{\overline{e}}_2$. Suppose the basis $\overline{\overline{e}}_1$ has coordinate matrix b relative*

---

the basis $\overline{\overline{e}}_{V2}$

$$\overline{\overline{e}}_{V1} = \overline{\overline{e}}_{V2}{}^*{}_*b$$

and $\overline{\overline{e}}_{W1}$ has coordinate matrix c relative the basis $\overline{\overline{e}}_{W2}$

(12.4.18) $\qquad \overline{\overline{e}}_{W1} = \overline{\overline{e}}_{W2}{}^*{}_*c$

Then there is relationship between matrices $a_1$ and $a_2$

$$a_1 = c^{-1*}{}^*{}_*a_2{}^*{}_*b$$

PROOF. Vector $\overline{v} \in V$ has expansion

$$\overline{v} = e_{V1}{}^*{}_*v = e_{V2}{}^*{}_*b_*{}^*v$$

relative to bases $\overline{\overline{e}}_{V1}$ , $\overline{\overline{e}}_{V2}$ . Since $a$ is homomorphism, we can write it as

(12.4.24) $\qquad \overline{w} = e_{W1}{}^*{}_*a_1{}^*{}_*v$

relative to bases $\overline{\overline{e}}_{V1}$, $\overline{\overline{e}}_{W1}$ and as

(12.4.25) $\qquad \overline{w} = e_{W2}{}^*{}_*a_2{}^*{}_*b_*{}^*v$

relative to bases $\overline{\overline{e}}_{V2}$, $\overline{\overline{e}}_{W2}$. According to the theorem 12.4.8, matrix $c$ has $^*_*$-inverse and from equation (12.4.18) it follows that

(12.4.26) $\qquad \overline{\overline{e}}_{W2} = e_{W1}{}^*{}_*\overline{\overline{c}}^{-1*}{}^*$

The equality

(12.4.27) $\qquad \overline{w} = e_{W1}{}^*{}_*c^{-1*}{}^*{}_*a_2{}^*{}_*b_*{}^*v$

follows from equalities (12.4.25), (12.4.26). From the theorem 8.4.25 and comparison of equations (12.4.24) and (12.4.27) it follows that

(12.4.28) $\qquad v^*{}_*a_1 = c^{-1*}{}^*{}_*a_2{}^*{}_*b_*{}^*v$

Since vector $a$ is arbitrary vector, the statement of theorem follows from the theorem 12.5.2 and equation (12.4.28).

<div align="right">□</div>

THEOREM 12.4.16. *Let $\overline{a}$ be automorphism of left $A$-vector space $V$ of columns. Let $a_1$ be matrix of this automorphism defined relative to basis $\overline{\overline{e}}_1$ and $a_2$ be matrix of the same automorphism defined relative to basis $\overline{\overline{e}}_2$. Suppose the basis $\overline{\overline{e}}_1$ has coordinate matrix b relative* *the basis $\overline{\overline{e}}_2$*

$$\overline{\overline{e}}_1 = \overline{\overline{e}}_{2*}{}^*b$$

*Then there is relationship between matrices $a_1$ and $a_2$*

$$a_1 = c^{-1}{}_*{}^*{}_*{}^*a_{2*}{}^*b$$

PROOF. Statement follows from theorem 12.4.13, because in this case $c = b$. □

---

*the basis $\overline{\overline{e}}_2$*

$$\overline{\overline{e}}_1 = \overline{\overline{e}}_2{}^*{}_*b$$

*Then there is relationship between matrices $a_1$ and $a_2$*

$$a_1 = c^{-1}{}^*{}_*{}^*{}_*a_2{}^*{}_*b$$

PROOF. Statement follows from theorem 12.4.14, because in this case $c = b$. □

## 12.5. Basis Manifold for Right $A$-Vector Space

Let $V$ be a right $A$-vector space of columns. We proved in the theorem 12.4.9 that, if we select a basis of right $A$-vector space $V$, then we can identify any automorphism $\overline{f}$ of right $A$-vector space $V$ with $_*{}^*$-nonsingular matrix $f$. Corresponding transformation of coordinates of vector

$$v' = f_*{}^*v$$

is called **linear transformation**.

Let $V$ be a right $A$-vector space of rows. We proved in the theorem 12.4.10 that, if we select a basis of right $A$-vector space $V$, then we can identify any automorphism $\overline{f}$ of right $A$-vector space $V$ with $^*_*$-nonsingular matrix $f$. Corresponding transformation of coordinates of vector

$$v' = f^*{}_*v$$

is called **linear transformation**.

THEOREM 12.5.1. *Let $V$ be a right $A$-vector space of columns and $\overline{\overline{e}}$ be basis of right $A$-vector space $V$. Automorphisms of right $A$-vector space of columns form a left-side linear effective $GL(n_*{}^*A)$-representation.*

THEOREM 12.5.2. *Let $V$ be a right $A$-vector space of rows and $\overline{\overline{e}}$ be basis of right $A$-vector space $V$. Automorphisms of right $A$-vector space of rows form a left-side linear effective $GL(n^*_*A)$-representation.*

PROOF. Let $a$, $b$ be matrices of automorphisms $\overline{a}$ and $\overline{b}$ with respect to basis $\overline{\overline{e}}$. According to the theorem 11.3.6, coordinate transformation has following form

(12.5.1)        $v' = a_*{}^*v$

(12.5.2)        $v'' = b_*{}^*v'$

The equality

$$v'' = b_*{}^*a_*{}^*v$$

follows from equalities (12.5.1), (12.5.2). According to the theorem 11.5.2, the product of automorphisms $\overline{a}$ and $\overline{b}$ has matrix $b_*{}^*a$. Therefore, automorphisms of right $A$-vector space of columns form a left-side linear $GL(n_*{}^*A)$-representation.

PROOF. Let $a$, $b$ be matrices of automorphisms $\overline{a}$ and $\overline{b}$ with respect to basis $\overline{\overline{e}}$. According to the theorem 11.3.7, coordinate transformation has following form

(12.5.4)        $v' = a^*{}_*v$

(12.5.5)        $v'' = b^*{}_*v'$

The equality

$$v'' = b^*{}_*a^*{}_*v$$

follows from equalities (12.5.4), (12.5.5). According to the theorem 11.5.3, the product of automorphisms $\overline{a}$ and $\overline{b}$ has matrix $b^*{}_*a$. Therefore, automorphisms of right $A$-vector space of rows form a left-side linear $GL(n^*_*A)$-representation.

It remains to prove that the kernel of inefficiency consists only of identity. Identity transformation satisfies to equation

$$v^i = a^i_j v^j$$

Choosing values of coordinates as $v^i = \delta^i_k$ where we selected $k$ we get

(12.5.3) $$\delta^i_k = a^i_j \delta^j_k$$

From (12.5.3) it follows

$$\delta^i_k = a^i_k$$

Since $k$ is arbitrary, we get the conclusion $a = \delta$. $\qquad\square$

It remains to prove that the kernel of inefficiency consists only of identity. Identity transformation satisfies to equation

$$v_i = a^j_i v_j$$

Choosing values of coordinates as $v_i = \delta^k_i$ where we selected $k$ we get

(12.5.6) $$\delta^k_i = a^j_i \delta^k_j$$

From (12.5.6) it follows

$$\delta^k_i = a^k_i$$

Since $k$ is arbitrary, we get the conclusion $a = \delta$. $\qquad\square$

**Theorem 12.5.3.** *Automorphism $a$ acting on each vector of basis of right $A$-vector space of columns maps a basis into another basis.*

**Theorem 12.5.4.** *Automorphism $a$ acting on each vector of basis of right $A$-vector space of rows maps a basis into another basis.*

**Proof.** Let $\overline{\overline{e}}$ be basis of right $A$-vector space $V$ of columns. According to theorem 12.4.9, vector $e_i$ maps into a vector $e'_i$

(12.5.7) $$e'_i = a_*{}^* e_i$$

Let vectors $e'_i$ be linearly dependent. Then $\lambda \neq 0$ in the equality

(12.5.8) $$e'_*{}^* \lambda = 0$$

From equations (12.5.7) and (12.5.8) it follows that

$$a^{-1}{}^*{}_* e'_*{}^* \lambda = e_*{}^* \lambda = 0$$

and $\lambda \neq 0$. This contradicts to the statement that vectors $e_i$ are linearly independent. Therefore vectors $e'_i$ are linearly independent and form basis. $\quad\square$

**Proof.** Let $\overline{\overline{e}}$ be basis of right $A$-vector space $V$ of rows. According to theorem 12.4.10, vector $e^i$ maps into a vector $e'^i$

(12.5.9) $$e'^i = a^*{}_* e^i$$

Let vectors $e'^i$ be linearly dependent. Then $\lambda \neq 0$ in the equality

(12.5.10) $$e'^*{}_* \lambda = 0$$

From equations (12.5.9) and (12.5.10) it follows that

$$a^{-1}{}^*{}_* e'^*{}_* \lambda = e^*{}_* \lambda = 0$$

and $\lambda \neq 0$. This contradicts to the statement that vectors $e^i$ are linearly independent. Therefore vectors $e'^i$ are linearly independent and form basis. $\quad\square$

Thus we can extend a left-side linear $GL(V_*)$-representation in right $A$-vector space $V_*$ to the set of bases of right $A$-vector space $V$. Transformation of this left-side representation on the set of bases of right $A$-vector space $V$ is called **active transformation** because the homomorphism of the right $A$-vector space induced this transformation (See also definition in the section [4]-14.1-3 as well the definition on the page [2]-214).

According to definition we write the action of the transformation $a \in GL(V_*)$ on the basis $\overline{\overline{e}}$ as $a_*{}^* \overline{\overline{e}}$. Consider the equality

(12.5.11) $$a_*{}^* \overline{\overline{e}}_*{}^* v = a_*{}^* \overline{\overline{e}}_*{}^* v$$

According to definition we write the action of the transformation $a \in GL(V_*)$ on the basis $\overline{\overline{e}}$ as $a^*{}_* \overline{\overline{e}}$. Consider the equality

(12.5.12) $$a^*{}_* \overline{\overline{e}}^*{}_* v = a^*{}_* \overline{\overline{e}}^*{}_* v$$

The expression $a_*{}^*\overline{\overline{e}}$ on the left side of the equality (12.5.11) is image of basis $\overline{\overline{e}}$ with respect to active transformation $a$. The expression $\overline{\overline{e}}_*{}^*v$ on the right side of the equality (12.5.11) is expansion of vector $\overline{v}$ with respect to basis $\overline{\overline{e}}$. Therefore, the expression on the right side of the equality (12.5.11) is image of vector $\overline{v}$ with respect to endomorphism $a$ and the expression on the left side of the equality

$$(12.5.11) \qquad a_*{}^*\overline{\overline{e}}_*{}^*v = a_*{}^*\overline{\overline{e}}_*{}^*v$$

is expansion of vector $\overline{v}$ with respect to image of basis $\overline{\overline{e}}$. Therefore, from the equality (12.5.11) it follows that endomorphism $a$ of right $A$-vector space and corresponding active transformation $a$ act synchronously and coordinates $a \circ \overline{v}$ of image of the vector $\overline{v}$ with respect to the image $a_*{}^*\overline{\overline{e}}$ of the basis $\overline{\overline{e}}$ are the same as coordinates of the vector $\overline{v}$ with respect to the basis $\overline{\overline{e}}$.

The expression $a^*{}_*\overline{\overline{e}}$ on the left side of the equality (12.5.12) is image of basis $\overline{\overline{e}}$ with respect to active transformation $a$. The expression $\overline{\overline{e}}^*{}_*v$ on the right side of the equality (12.5.12) is expansion of vector $\overline{v}$ with respect to basis $\overline{\overline{e}}$. Therefore, the expression on the right side of the equality (12.5.12) is image of vector $\overline{v}$ with respect to endomorphism $a$ and the expression on the left side of the equality

$$(12.5.12) \qquad a^*{}_*\overline{\overline{e}}^*{}_*v = a^*{}_*\overline{\overline{e}}^*{}_*v$$

is expansion of vector $\overline{v}$ with respect to image of basis $\overline{\overline{e}}$. Therefore, from the equality (12.5.12) it follows that endomorphism $a$ of right $A$-vector space and corresponding active transformation $a$ act synchronously and coordinates $a \circ \overline{v}$ of image of the vector $\overline{v}$ with respect to the image $a^*{}_*\overline{\overline{e}}$ of the basis $\overline{\overline{e}}$ are the same as coordinates of the vector $\overline{v}$ with respect to the basis $\overline{\overline{e}}$.

**Theorem 12.5.5.** *Active left-side $GL(n_*{}^*A)$-representation on the set of bases is single transitive representation. The set of bases identified with tragectory $GL(n_*{}^*A)(V)_*{}^*\overline{\overline{e}}$ of active left-side $GL(n_*{}^*A)$-representation is called the* **basis manifold** *of right $A$-vector space $V$.*

**Theorem 12.5.6.** *Active left-side $GL(n^*{}_*A)$-representation on the set of bases is single transitive representation. The set of bases identified with tragectory $GL(n^*{}_*A)(V)^*{}_*\overline{\overline{e}}$ of active left-side $GL(n^*{}_*A)$-representation is called the* **basis manifold** *of right $A$-vector space $V$.*

To prove theorems 12.5.5, 12.5.6, it is sufficient to show that at least one transformation of left-side representation is defined for any two bases and this transformation is unique.

**Proof.** Homomorphism $a$ operating on basis $\overline{\overline{e}}$ has form

$$e'_i = a_*{}^*e_i$$

where $e'_i$ is coordinate matrix of vector $\overline{e}'_i$ relative basis $\overline{\overline{h}}$ and $e_i$ is coordinate matrix of vector $\overline{e}_i$ relative basis $\overline{\overline{h}}$. Therefore, coordinate matrix of image of basis equal to $_*{}^*$-product of matrix of automorphism over coordinate matrix of original basis

$$e' = g_*{}^*e$$

According to the theorem 12.4.7, matri-

**Proof.** Homomorphism $a$ operating on basis $\overline{\overline{e}}$ has form

$$e'^i = a^*{}_*e^i$$

where $e'^i$ is coordinate matrix of vector $\overline{e}'^i$ relative basis $\overline{\overline{h}}$ and $e^i$ is coordinate matrix of vector $\overline{e}^i$ relative basis $\overline{\overline{h}}$. Therefore, coordinate matrix of image of basis equal to $^*{}_*$-product of matrix of automorphism over coordinate matrix of original basis

$$e' = g^*{}_*e$$

According to the theorem 12.4.8, matrices $g$ and $e$ are nonsingular. Therefore,

ces $g$ and $e$ are nonsingular. Therefore, matrix

$$a = {e'}_*{}^* e^{-1}{}_*{}^*$$

is the matrix of automorphism mapping basis $\overline{\overline{e}}$ to basis $\overline{\overline{e}}'$.

Suppose elements $g_1$, $g_2$ of group $G$ and basis $\overline{\overline{e}}$ satisfy equation

$$g_1{}_*{}^* e = g_2{}_*{}^* e$$

According to the theorems 12.4.7 and 8.1.32, we get $g_1 = g_2$. This proves statement of theorem. $\qquad\square$

matrix

$$a = {e'}^*{}_* e^{-1}{}^*{}_*$$

is the matrix of automorphism mapping basis $\overline{\overline{e}}$ to basis $\overline{\overline{e}}'$.

Suppose elements $g_1$, $g_2$ of group $G$ and basis $\overline{\overline{e}}$ satisfy equation

$$g_1{}^*{}_* e = g_2{}^*{}_* e$$

According to the theorems 12.4.8 and 8.1.32, we get $g_1 = g_2$. This proves statement of theorem. $\qquad\square$

Let us define an additional structure on right $A$-vector space $V$. Then not every automorphism keeps properties of the selected structure. For instance, if we introduce norm in right $A$-vector space $V$, then we are interested in automorphisms which preserve the norm of the vector.

DEFINITION 12.5.7. *Normal subgroup $G(V)$ of the group $GL(V)$ such that subgroup $G(V)$ generates automorphisms which hold properties of the selected structure is called* **symmetry group**.

*Without loss of generality we identify element $g$ of group $G(V)$ with corresponding transformation of representation and write its action on vector $v \in V$ as $v_*{}^* g$.* $\qquad\square$

DEFINITION 12.5.8. *Normal subgroup $G(V)$ of the group $GL(V)$ such that subgroup $G(V)$ generates automorphisms which hold properties of the selected structure is called* **symmetry group**.

*Without loss of generality we identify element $g$ of group $G(V)$ with corresponding transformation of representation and write its action on vector $v \in V$ as $v^*{}_* g$.* $\qquad\square$

The left-side representation of group $G$ in right $A$-vector space is called **linear $G$-representation**.

DEFINITION 12.5.9. *The left-side representation of group $G(V)$ in the set of bases of right $A$-vector space $V$ is called* **active right-side $G$-representation**. $\qquad\square$

If the basis $\overline{\overline{e}}$ is given, then we can identify the automorphism of right $A$-vector space $V$ and its coordinates with respect to the basis $\overline{\overline{e}}$. The set $G(V_*)$ of coordinates of automorphisms with respect to the basis $\overline{\overline{e}}$ is group isomorphic to the group $G(V)$.

Not every two bases can be mapped by a transformation from the symmetry group because not every nonsingular linear transformation belongs to the representation of group $G(V)$. Therefore, we can represent the set of bases as a union of orbits of group $G(V)$. In particular, if the basis $\overline{\overline{e}} \in G(V)$, then the group $G(V)$ is orbit of the basis $\overline{\overline{e}}$.

DEFINITION 12.5.10. *We call orbit $G(V)_*{}^* \overline{\overline{e}}$ of the selected basis $\overline{\overline{e}}$ the* **basis**

DEFINITION 12.5.11. *We call orbit $G(V)^*{}_* \overline{\overline{e}}$ of the selected basis $\overline{\overline{e}}$ the* **basis**

**manifold** *of right $A$-vector space $V$ of columns.* □

**manifold** *of right $A$-vector space $V$ of rows.* □

THEOREM 12.5.12. *Active left-side $G(V)$-representation on basis manifold is single transitive representation.*

PROOF. The theorem follows from the theorem 12.5.5 and the definition 12.5.10, as well the theorem follows from the theorem 12.5.6 and the definition 12.5.11. □

THEOREM 12.5.13. *On the basis manifold $G(V)_{*}{}^{*}\overline{\overline{e}}$ of right $A$-vector space of columns, there exists single transitive right-side $G(V)$-representation, commuting with active.*

THEOREM 12.5.14. *On the basis manifold $G(V)^{*}{}_{*}\overline{\overline{e}}$ of right $A$-vector space of rows, there exists single transitive right-side $G(V)$-representation, commuting with active.*

PROOF. Theorem 12.5.12 means that the basis manifold $G(V)_{*}{}^{*}\overline{\overline{e}}$ is a homogenous space of group $G$. The theorem follows from the theorem 2.5.17. □

PROOF. Theorem 12.5.12 means that the basis manifold $G(V)^{*}{}_{*}\overline{\overline{e}}$ is a homogenous space of group $G$. The theorem follows from the theorem 2.5.17. □

As we see from remark 2.5.18 transformation of right-side $G(V)$-representation is different from an active transformation and cannot be reduced to transformation of space $V$.

DEFINITION 12.5.15. *A transformation of right-side $G(V)$-representation is called* **passive transformation** *of basis manifold $G(V)_{*}{}^{*}\overline{\overline{e}}$ of right $A$-vector space of columns, and the right-side $G(V)$-representation is called* **passive right-side $G(V)$-representation**. *According to the definition we write the passive transformation of basis $\overline{\overline{e}}$ defined by element $a \in G(V)$ as $\overline{\overline{e}}_{*}{}^{*}a$.* □

DEFINITION 12.5.16. *A transformation of right-side $G(V)$-representation is called* **passive transformation** *of basis manifold $G(V)^{*}{}_{*}\overline{\overline{e}}$ of right $A$-vector space of rows, and the right-side $G(V)$-representation is called* **passive right-side $G(V)$-representation**. *According to the definition we write the passive transformation of basis $\overline{\overline{e}}$ defined by element $a \in G(V)$ as $\overline{\overline{e}}^{*}{}_{*}a$.* □

THEOREM 12.5.17. *The coordinate matrix of basis $\overline{\overline{e}}'$ relative basis $\overline{\overline{e}}$ of right $A$-vector space $V$ of columns is identical with the matrix of passive transformation mapping basis $\overline{\overline{e}}$ to basis $\overline{\overline{e}}'$.*

THEOREM 12.5.18. *The coordinate matrix of basis $\overline{\overline{e}}'$ relative basis $\overline{\overline{e}}$ of right $A$-vector space $V$ of rows is identical with the matrix of passive transformation mapping basis $\overline{\overline{e}}$ to basis $\overline{\overline{e}}'$.*

PROOF. According to the theorem 8.4.17, the coordinate matrix of basis $\overline{\overline{e}}'$ relative basis $\overline{\overline{e}}$ consist from rows which are coordinate matrices of vectors $\overline{e}'_i$ rel-

PROOF. According to the theorem 8.4.24, the coordinate matrix of basis $\overline{\overline{e}}'$ relative basis $\overline{\overline{e}}$ consist from cols which are coordinate matrices of vectors $\overline{e}'_i$ rel-ative the basis $\overline{\overline{e}}$. Therefore,

$$\overline{e}'_i = \overline{e}_{*}{}^* e'_i$$

At the same time the passive transformation $a$ mapping one basis to another has a form

$$\overline{e}'_i = \overline{e}_{*}{}^* a_i$$

According to the theorem 8.4.18,

$$e'_i = a_i$$

for any $i$. This proves the theorem.    □

REMARK 12.5.19. *According to theorems 2.5.15, 12.5.17, we can identify a basis $\overline{\overline{e}}'$ of the basis manifold $G(V)_{*}{}^*\overline{\overline{e}}$ of right $A$-vector space of columns and the matrix $g$ of coordinates of the basis $\overline{\overline{e}}'$ with respect to the basis $\overline{\overline{e}}$*

$$(12.5.13) \qquad \overline{\overline{e}}' = \overline{\overline{e}}_{*}{}^* g$$

*Since $\overline{\overline{e}}' = \overline{\overline{e}}$, then the equality (12.5.13) gets form*

$$(12.5.14) \qquad \overline{\overline{e}} = \overline{\overline{e}}_{*}{}^* E_n$$

*Based on the equality (12.5.14), we will use notation $G(V)_{*}{}^{*\overline{\overline{E_n}}}$ for basis manifold $G(V)_{*}{}^*\overline{\overline{e}}$.*    □

ative the basis $\overline{\overline{e}}$. Therefore,

$$\overline{e}'^i = \overline{e}^*{}_* e'^i$$

At the same time the passive transformation $a$ mapping one basis to another has a form

$$\overline{e}'^i = \overline{e}^*{}_* a^i$$

According to the theorem 8.4.25,

$$e'^i = a^i$$

for any $i$. This proves the theorem.    □

REMARK 12.5.20. *According to theorems 2.5.15, 12.5.18, we can identify a basis $\overline{\overline{e}}'$ of the basis manifold $G(V)^*{}_*\overline{\overline{e}}$ of right $A$-vector space of columns and the matrix $g$ of coordinates of the basis $\overline{\overline{e}}'$ with respect to the basis $\overline{\overline{e}}$*

$$(12.5.15) \qquad \overline{\overline{e}}' = \overline{\overline{e}}^*{}_* g$$

*Since $\overline{\overline{e}}' = \overline{\overline{e}}$, then the equality (12.5.15) gets form*

$$(12.5.16) \qquad \overline{\overline{e}} = \overline{\overline{e}}^*{}_* E_n$$

*Based on the equality (12.5.16), we will use notation $G(V)^*{}_*\overline{\overline{E_n}}$ for basis manifold $G(V)^*{}_*\overline{\overline{e}}$.*    □

## 12.6. Passive Transformation in Right $A$-Vector Space

An active transformation changes bases and vectors uniformly and coordinates of vector relative basis do not change. A passive transformation changes only the basis and it leads to transformation of coordinates of vector relative to basis.

THEOREM 12.6.1. *Let $V$ be a right $A$-vector space of columns. Let passive transformation $g \in G(V_{*})$ map the basis $\overline{\overline{e}}_1$ into the basis $\overline{\overline{e}}_2$*

$$(12.6.1) \qquad \overline{\overline{e}}_2 = \overline{\overline{e}}_{1*}{}^* g$$

*Let*

$$(12.6.2) \qquad v_i = \begin{pmatrix} v_i^1 \\ ... \\ v_i^n \end{pmatrix}$$

*be matrix of coordinates of the vector $\overline{v}$ with respect to the basis $\overline{\overline{e}}_i$, $i = 1, 2$. Coordinate transformations*

$$(12.6.3) \qquad v_1 = g_{*}{}^* v_2$$

THEOREM 12.6.2. *Let $V$ be a right $A$-vector space of rows. Let passive transformation $g \in G(V_{*})$ map the basis $\overline{\overline{e}}_1$ into the basis $\overline{\overline{e}}_2$*

$$(12.6.5) \qquad \overline{\overline{e}}_2 = \overline{\overline{e}}_1{}^*{}_* g$$

*Let*

$$v_i = \begin{pmatrix} v_{i\,1} & ... & v_{i\,n} \end{pmatrix}$$

*be matrix of coordinates of the vector $\overline{v}$ with respect to the basis $\overline{\overline{e}}_i$, $i = 1, 2$. Coordinate transformations*

$$(12.6.6) \qquad v_1 = g^*{}_* v_2$$

$$(12.6.4) \qquad v_2 = g^{-1*^*}{}_*v_1$$

*do not depend on vector $\overline{v}$ or basis $\overline{\overline{e}}$, but is defined only by coordinates of vector $\overline{v}$ relative to basis $\overline{\overline{e}}$.*

Proof. According to the theorem 8.4.17,

$$(12.6.8) \qquad \overline{v} = \overline{\overline{e}}_{1*}{}^*v_1 = \overline{\overline{e}}_{2*}{}^*v_2$$

The equality

$$(12.6.9) \qquad \overline{\overline{e}}_{1*}{}^*v_1 = \overline{\overline{e}}_{1*}{}^*g_*{}^*v_2$$

follows from equalities (12.6.1), (12.6.8). According to the theorem 8.4.18, the coordinate transformation (12.6.3) follows from the equality (12.6.9). Because $g$ is $_*$*-nonsingular matrix, coordinate transformation (12.6.4) follows from the equality (12.6.3). □

Theorem 12.6.3. *Coordinate transformations*

$$\boxed{(12.6.4) \quad v_2 = g^{-1*^*}{}_*v_1}$$

*form effective linear left-side contravariant $G$-representation which is called* **coordinate representation in right $A$-vector space** *of columns.*

Proof. Let

$$(12.6.12) \qquad v_i = \begin{pmatrix} v_i^1 \\ ... \\ v_i^n \end{pmatrix}$$

be matrix of coordinates of the vector $\overline{v}$ with respect to the basis $\overline{\overline{e}}_i$, $i = 1, 2, 3$. Then

$$(12.6.13) \qquad \overline{v} = \overline{\overline{e}}_{1*}{}^*v_1 = \overline{\overline{e}}_{2*}{}^*v_2$$
$$= \overline{\overline{e}}_{3*}{}^*v_3$$

Let passive transformation $g_1$ map the basis $\overline{\overline{e}}_1 \in G(V)_*{}^*\overline{e}$ into the basis $\overline{\overline{e}}_2 \in G(V)_*{}^*\overline{e}$

$$(12.6.14) \qquad \overline{\overline{e}}_2 = \overline{\overline{e}}_{1*}{}^*g_1$$

Let passive transformation $g_2$ map the basis $\overline{\overline{e}}_2 \in G(V)_*{}^*\overline{e}$ into the basis $\overline{\overline{e}}_3 \in$

$$(12.6.7) \qquad v_2 = g^{-1*^*}{}_*v_1$$

*do not depend on vector $\overline{v}$ or basis $\overline{\overline{e}}$, but is defined only by coordinates of vector $\overline{v}$ relative to basis $\overline{\overline{e}}$.*

Proof. According to the theorem 8.4.24,

$$(12.6.10) \qquad \overline{v} = \overline{\overline{e}}_1{}^*{}_*v_1 = \overline{\overline{e}}_2{}^*{}_*v_2$$

The equality

$$(12.6.11) \qquad \overline{\overline{e}}_1{}^*{}_*v_1 = \overline{\overline{e}}_1{}^*{}_*g^*{}_*v_2$$

follows from equalities (12.6.5), (12.6.10). According to the theorem 8.4.25, the coordinate transformation (12.6.6) follows from the equality (12.6.11). Because $g$ is $_*$*-nonsingular matrix, coordinate transformation (12.6.7) follows from the equality (12.6.6). □

Theorem 12.6.4. *Coordinate transformations*

$$\boxed{(12.6.7) \quad v_2 = g^{-1*^*}{}_*v_1}$$

*form effective linear left-side contravariant $G$-representation which is called* **coordinate representation in right $A$-vector space** *of rows.*

Proof. Let

$$v_i = \begin{pmatrix} v_{i\,1} & ... & v_{i\,n} \end{pmatrix}$$

be matrix of coordinates of the vector $\overline{v}$ with respect to the basis $\overline{\overline{e}}_i$, $i = 1, 2, 3$. Then

$$(12.6.21) \qquad \overline{v} = \overline{\overline{e}}_1{}^*{}_*v_1 = \overline{\overline{e}}_2{}^*{}_*v_2$$
$$= \overline{\overline{e}}_3{}^*{}_*v_3$$

Let passive transformation $g_1$ map the basis $\overline{\overline{e}}_1 \in G(V)^*{}_*\overline{e}$ into the basis $\overline{\overline{e}}_2 \in G(V)^*{}_*\overline{e}$

$$(12.6.22) \qquad \overline{\overline{e}}_2 = \overline{\overline{e}}_1{}^*{}_*g_1$$

Let passive transformation $g_2$ map the basis $\overline{\overline{e}}_2 \in G(V)^*{}_*\overline{e}$ into the basis $\overline{\overline{e}}_3 \in G(V)^*{}_*\overline{e}$

$$(12.6.23) \qquad \overline{\overline{e}}_3 = \overline{\overline{e}}_2{}^*{}_*g_2$$

$G(V)_*{}^*\overline{\overline{e}}$

(12.6.15) $$\overline{\overline{e}}_3 = \overline{\overline{e}}_{2*}{}^* g_2$$

The transformation

(12.6.16) $$\overline{\overline{e}}_3 = \overline{\overline{e}}_{1*}{}^*(g_{1*}{}^* g_2)$$

is the product of transformations (12.6.14), (12.6.15). The coordinate transformation

(12.6.17) $$v_2 = g_1^{-1}{}^*{}_* v_1$$

corresponds to passive transformation $g_1$ and follows from equalities (12.6.13), (12.6.14) and the theorem 12.6.1. Similarly, according to the theorem 12.6.1, the coordinate transformation

(12.6.18) $$v_3 = g_2^{-1}{}^*{}_* v_2$$

corresponds to passive transformation $g_2$ and follows from equalities (12.6.13), (12.6.15). The product of coordinate transformations (12.3.17), (12.3.18) has form

(12.6.19) $$v_3 = g_2^{-1}{}^*{}_* g_1^{-1}{}^*{}_* v_1$$

and is coordinate transformation corresponding to passive transformation (12.6.16). The equality

(12.6.20) $$v_3 = (g_{2*}{}^* g_1)^{-1}{}^*{}_* v_1$$

follows from equalities (8.1.66), (12.6.19). According to the definition 2.5.9, coordinate transformations form linear right-side contravariant $G$-representation. Suppose coordinate transformation does not change vectors $\delta_k$. Then unit of group $G$ corresponds to it because representation is single transitive. Therefore, coordinate representation is effective. □

The transformation

(12.6.24) $$\overline{\overline{e}}_3 = \overline{\overline{e}}_{1*}{}^*(g_{1*}{}^* g_2)$$

is the product of transformations (12.6.22), (12.6.23). The coordinate transformation

(12.6.25) $$v_2 = g_1^{-1}{}^*{}_* v_1$$

corresponds to passive transformation $g_1$ and follows from equalities (12.6.21), (12.6.22) and the theorem 12.6.2. Similarly, according to the theorem 12.6.2, the coordinate transformation

(12.6.26) $$v_3 = g_2^{-1}{}^*{}_* v_2$$

corresponds to passive transformation $g_2$ and follows from equalities (12.6.21), (12.6.23). The product of coordinate transformations (12.3.17), (12.3.18) has form

(12.6.27) $$v_3 = g_2^{-1}{}^*{}_* g_1^{-1}{}^*{}_* v_1$$

and is coordinate transformation corresponding to passive transformation (12.6.24). The equality

(12.6.28) $$v_3 = (g_{2*}{}^* g_1)^{-1}{}^*{}_* v_1$$

follows from equalities (8.1.66), (12.6.27). According to the definition 2.5.9, coordinate transformations form linear right-side contravariant $G$-representation. Suppose coordinate transformation does not change vectors $\delta_k$. Then unit of group $G$ corresponds to it because representation is single transitive. Therefore, coordinate representation is effective. □

THEOREM 12.6.5. *Let $V$ be a right $A$-vector space of columns and $\overline{\overline{e}}_1$, $\overline{\overline{e}}_2$ be bases in right $A$-vector space $V$. Let $g$ be passive transformation of basis $\overline{\overline{e}}_1$ into basis $\overline{\overline{e}}_2$*

$$\overline{\overline{e}}_2 = \overline{\overline{e}}_{1*}{}^* g$$

*Let $f$ be endomorphism of right $A$-vector space $V$. Let $f_i$, $i = 1, 2$, be the matrix*

THEOREM 12.6.6. *Let $V$ be a right $A$-vector space of rows and $\overline{\overline{e}}_1$, $\overline{\overline{e}}_2$ be bases in right $A$-vector space $V$. Let $g$ be passive transformation of basis $\overline{\overline{e}}_1$ into basis $\overline{\overline{e}}_2$*

$$\overline{\overline{e}}_2 = \overline{\overline{e}}_{1*}{}^* g$$

*Let $f$ be endomorphism of right $A$-vector space $V$. Let $f_i$, $i = 1, 2$, be the matrix of endomorphism $f$ with respect to the*

*of endomorphism $f$ with respect to the basis $\overline{\overline{e}}_i$. Then*

$$(12.6.29) \qquad f_2 = g^{-1^*}{}_*^* f_{1*}^* g$$

*basis $\overline{\overline{e}}_i$. Then*

$$(12.6.30) \qquad f_2 = g^{-1^*}{}_*^* f_{1*}^* g$$

PROOF OF THE THEOREM 12.6.5. Let endomorphism $f$ maps a vector $v$ into a vector $w$

$$(12.6.31) \qquad w = f \circ v$$

Let $v_i$, $i = 1, 2$, be the matrix of vector $v$ with respect to the basis $\overline{\overline{e}}_i$. Then

$$(12.6.32) \qquad v_1 = g_*^* v_2$$

follows from the theorem 12.6.1.

Let $w_i$, $i = 1, 2$, be the matrix of vector $w$ with respect to the basis $\overline{\overline{e}}_i$. Then

$$(12.6.33) \qquad w_1 = g_*^* w_2$$

follows from the theorem 12.6.1.

According to the theorem 11.3.6, equalities

$$(12.6.34) \qquad w_1 = f_{1*}^* v_1$$

$$(12.6.35) \qquad w_2 = f_{2*}^* v_2$$

follow from the equality (12.6.31). The equality

$$(12.6.36) \qquad g_*^* w_2 = f_{1*}^* v_1 = f_{1*}^* g_*^* v_2$$

follows from equalities (12.6.32), (12.6.34), (12.6.33). Since $g$ is regular matrix, then the equality

$$(12.6.37) \qquad w_2 = f_{1*}^* v_1 = g^{-1^*}{}_*^* f_{1*}^* g_*^* v_2$$

follows from the equality (12.6.36). The equality

$$(12.6.38) \qquad f_{2*}^* v_2 = g^{-1^*}{}_*^* f_{1*}^* g_*^* v_2$$

follows from equalities (12.6.35), (12.6.37). The equality (12.6.29) follows from the equality (12.6.38). $\qquad\square$

PROOF OF THE THEOREM 12.6.6. Let endomorphism $f$ maps a vector $v$ into a vector $w$

$$(12.6.39) \qquad w = f \circ v$$

Let $v_i$, $i = 1, 2$, be the matrix of vector $v$ with respect to the basis $\overline{\overline{e}}_i$. Then

$$(12.6.40) \qquad v_1 = g^*{}_* v_2$$

follows from the theorem 12.6.2.

Let $w_i$, $i = 1, 2$, be the matrix of vector $w$ with respect to the basis $\overline{\overline{e}}_i$. Then

$$(12.6.41) \qquad w_1 = g^*{}_* w_2$$

follows from the theorem 12.6.2.

According to the theorem 11.3.6, equalities

$$(12.6.42) \qquad w_1 = f_1^*{}_* v_1$$

$$(12.6.43) \qquad w_2 = f_2^*{}_* v_2$$

follow from the equality (12.6.39). The equality

$$(12.6.44) \qquad g^*{}_* w_2 = f_1{}^*{}_* v_1 = f_1{}^*{}_* g^*{}_* v_2$$

follows from equalities (12.6.40), (12.6.42), (12.6.41). Since $g$ is regular matrix, then the equality

$$(12.6.45) \qquad w_2 = f_1{}^*{}_* v_1 = g^{-1^*}{}^*{}_* f_1{}^*{}_* g^*{}_* v_2$$

follows from the equality (12.6.44). The equality

$$(12.6.46) \qquad f_2{}^*{}_* v_2 = g^{-1^*}{}^*{}_* f_1{}^*{}_* g^*{}_* v_2$$

follows from equalities (12.6.43), (12.6.45). The equality (12.6.30) follows from the equality (12.6.46). $\qquad\square$

---

THEOREM 12.6.7. *Let $V$ be a right $A$-vector space of columns. Let endomorphism $\overline{f}$ of right $A$-vector space $V$ is sum of endomorphisms $\overline{f}_1$, $\overline{f}_2$. The matrix $f$ of endomorphism $\overline{f}$ equal to the sum of matrices $f_1$, $f_2$ of endomorphisms $\overline{f}_1$, $\overline{f}_2$ and this equality does not depend on the choice of a basis.*

PROOF. Let $\overline{\overline{e}}$, $\overline{\overline{e}}'$ be bases of right $A$-vector space $V$ and $g$ be passive transformation of basis $\overline{\overline{e}}$ into basis $\overline{\overline{e}}'$

$$e' = g_*{}^* e$$

Let $f$ be the matrix of endomorphism $\overline{f}$ with respect to the basis $\overline{\overline{e}}$. Let $f'$ be the matrix of endomorphism $\overline{f}$ with respect to the basis $\overline{\overline{e}}'$. Let $f_i$, $i = 1, 2$, be the matrix of endomorphism $\overline{f}_i$ with respect to the basis $\overline{\overline{e}}$. Let $f'_i$, $i = 1, 2$, be the matrix of endomorphism $\overline{f}_i$ with respect to the basis $\overline{\overline{e}}'$. According to the theorem 11.4.2,

$$(12.6.47) \qquad f = f_1 + f_2$$

$$(12.6.48) \qquad f' = f'_1 + f'_2$$

According to the theorem 12.6.5,

$$(12.6.49) \qquad f' = g^{-1^*}{}^*{}_* f_*{}^* g$$

$$(12.6.50) \qquad f'_1 = g^{-1^*}{}^*{}_* f_1{}_*{}^* g$$

$$(12.6.51) \qquad f'_2 = g^{-1^*}{}^*{}_* f_2{}_*{}^* g$$

---

THEOREM 12.6.8. *Let $V$ be a right $A$-vector space of rows. Let endomorphism $\overline{f}$ of right $A$-vector space $V$ is sum of endomorphisms $\overline{f}_1$, $\overline{f}_2$. The matrix $f$ of endomorphism $\overline{f}$ equal to the sum of matrices $f_1$, $f_2$ of endomorphisms $\overline{f}_1$, $\overline{f}_2$ and this equality does not depend on the choice of a basis.*

PROOF. Let $\overline{e}$, $\overline{e}'$ be bases of right $A$-vector space $V$ and $g$ be passive transformation of basis $\overline{\overline{e}}$ into basis $\overline{\overline{e}}'$

$$e' = g^*{}_* e$$

Let $f$ be the matrix of endomorphism $\overline{f}$ with respect to the basis $\overline{\overline{e}}$. Let $f'$ be the matrix of endomorphism $\overline{f}$ with respect to the basis $\overline{\overline{e}}'$. Let $f_i$, $i = 1, 2$, be the matrix of endomorphism $\overline{f}_i$ with respect to the basis $\overline{\overline{e}}$. Let $f'_i$, $i = 1, 2$, be the matrix of endomorphism $\overline{f}_i$ with respect to the basis $\overline{\overline{e}}'$. According to the theorem 11.4.3,

$$(12.6.53) \qquad f = f_1 + f_2$$

$$(12.6.54) \qquad f' = f'_1 + f'_2$$

According to the theorem 12.6.6,

$$(12.6.55) \qquad f' = g^{-1^*}{}^*{}_* f^*{}_* g$$

$$(12.6.56) \qquad f'_1 = g^{-1^*}{}^*{}_* f_1{}^*{}_* g$$

$$(12.6.57) \qquad f'_2 = g^{-1^*}{}^*{}_* f_2{}^*{}_* g$$

The equality

$$
\begin{aligned}
f' &= g^{-1\cdot^*}{}_*{}^* f_*{}^* g \\
&= g^{-1\cdot^*}{}_*{}^* (f_1 + f_2)_*{}^* g \\
(12.6.52) \qquad &= g^{-1\cdot^*}{}_*{}^* f_{1*}{}^* g \\
&\quad + g^{-1\cdot^*}{}_*{}^* f_{2*}{}^* g \\
&= f_1' + f_2'
\end{aligned}
$$

follows from equalities (12.6.47), (12.6.49), (12.6.50), (12.6.51). From equalities (12.6.48), (12.6.52) it follows that the equality (12.6.47) is covariant with respect to choice of basis of right $A$-vector space $V$. □

---

**THEOREM 12.6.9.** *Let $V$ be a right $A$-vector space of columns. Let endomorphism $\overline{f}$ of right $A$-vector space $V$ is product of endomorphisms $\overline{f}_1$, $\overline{f}_2$. The matrix $f$ of endomorphism $\overline{f}$ equal to the $_*{}^*$-product of matrices $f_2$, $f_1$ of endomorphisms $\overline{f}_2$, $\overline{f}_1$ and this equality does not depend on the choice of a basis.*

PROOF. Let $\overline{\overline{e}}$, $\overline{\overline{e}}'$ be bases of right $A$-vector space $V$ and $g$ be passive transformation of basis $\overline{\overline{e}}$ into basis $\overline{\overline{e}}'$

$$
e' = g_*{}^* e
$$

Let $f$ be the matrix of endomorphism $\overline{f}$ with respect to the basis $\overline{\overline{e}}$. Let $f'$ be the matrix of endomorphism $\overline{f}$ with respect to the basis $\overline{\overline{e}}'$. Let $f_i$, $i = 1, 2$, be the matrix of endomorphism $\overline{f}_i$ with respect to the basis $\overline{\overline{e}}$. Let $f_i'$, $i = 1, 2$, be the matrix of endomorphism $\overline{f}_i$ with respect to the basis $\overline{\overline{e}}'$. According to the theorem 11.5.2,

$$
(12.6.59) \qquad f = f_{1*}{}^* f_2
$$

$$
(12.6.60) \qquad f' = f_{1*}'{}^* f_2'
$$

According to the theorem 12.6.5,

$$
(12.6.61) \qquad f' = g^{-1\cdot^*}{}_*{}^* f_*{}^* g
$$

$$
(12.6.62) \qquad f_1' = g^{-1\cdot^*}{}_*{}^* f_{1*}{}^* g
$$

$$
(12.6.63) \qquad f_2' = g^{-1\cdot^*}{}_*{}^* f_{2*}{}^* g
$$

---

The equality

$$
\begin{aligned}
f' &= g^{-1^*\cdot^*}{}_* f^*{}_* g \\
&= g^{-1^*\cdot^*}{}_* (f_1 + f_2)^*{}_* g \\
(12.6.58) \qquad &= g^{-1^*\cdot^*}{}_* f_1^*{}_* g \\
&\quad + g^{-1^*\cdot^*}{}_* f_2^*{}_* g \\
&= f_1' + f_2'
\end{aligned}
$$

follows from equalities (12.6.53), (12.6.55), (12.6.56), (12.6.57). From equalities (12.6.54), (12.6.58) it follows that the equality (12.6.53) is covariant with respect to choice of basis of right $A$-vector space $V$. □

---

**THEOREM 12.6.10.** *Let $V$ be a right $A$-vector space of rows. Let endomorphism $\overline{f}$ of right $A$-vector space $V$ is product of endomorphisms $\overline{f}_1$, $\overline{f}_2$. The matrix $f$ of endomorphism $\overline{f}$ equal to the $^*{}_*$-product of matrices $f_2$, $f_1$ of endomorphisms $\overline{f}_2$, $\overline{f}_1$ and this equality does not depend on the choice of a basis.*

PROOF. Let $\overline{\overline{e}}$, $\overline{\overline{e}}'$ be bases of right $A$-vector space $V$ and $g$ be passive transformation of basis $\overline{\overline{e}}$ into basis $\overline{\overline{e}}'$

$$
e' = g^*{}_* e
$$

Let $f$ be the matrix of endomorphism $\overline{f}$ with respect to the basis $\overline{\overline{e}}$. Let $f'$ be the matrix of endomorphism $\overline{f}$ with respect to the basis $\overline{\overline{e}}'$. Let $f_i$, $i = 1, 2$, be the matrix of endomorphism $\overline{f}_i$ with respect to the basis $\overline{\overline{e}}$. Let $f_i'$, $i = 1, 2$, be the matrix of endomorphism $\overline{f}_i$ with respect to the basis $\overline{\overline{e}}'$. According to the theorem 11.5.3,

$$
(12.6.65) \qquad f = f_{1*}{}^* f_2
$$

$$
(12.6.66) \qquad f' = f_{1*}'{}^* f_2'
$$

According to the theorem 12.6.6,

$$
(12.6.67) \qquad f' = g^{-1^*\cdot^*}{}_* f^*{}_* g
$$

$$
(12.6.68) \qquad f_1' = g^{-1^*\cdot^*}{}_* f_1^*{}_* g
$$

$$
(12.6.69) \qquad f_2' = g^{-1^*\cdot^*}{}_* f_2^*{}_* g
$$

The equality

$$
\begin{aligned}
f' &= g^{-1\ast\ast}{}_{\ast}f_{\ast}{}^{\ast}g \\
&= g^{-1\ast\ast}{}_{\ast}f_{1\ast}{}^{\ast}f_{2\ast}{}^{\ast}g \\
(12.6.64) \quad &= g^{-1\ast\ast}{}_{\ast}f_{1\ast}{}^{\ast}g \\
&{}_{\ast}{}^{\ast}g^{-1\ast\ast}{}_{\ast}f_{2\ast}{}^{\ast}g \\
&= f'_{1\ast}{}^{\ast}f'_{2}
\end{aligned}
$$

follows from equalities (12.6.59), (12.6.61), (12.6.62), (12.6.63). From equalities (12.6.60), (12.6.64) it follows that the equality (12.6.59) is covariant with respect to choice of basis of right $A$-vector space $V$.  □

The equality

$$
\begin{aligned}
f' &= g^{-1\ast\ast}{}_{\ast}f^{\ast}{}_{\ast}g \\
&= g^{-1\ast\ast}{}_{\ast}f_{1}{}^{\ast}{}_{\ast}f_{2}{}^{\ast}{}_{\ast}g \\
(12.6.70) \quad &= g^{-1\ast\ast}{}_{\ast}f_{1}{}^{\ast}{}_{\ast}g \\
&{}^{\ast}{}_{\ast}g^{-1\ast\ast}{}_{\ast}f_{2}{}^{\ast}{}_{\ast}g \\
&= f'^{\ast}_{1\ \ast}f'_{2}
\end{aligned}
$$

follows from equalities (12.6.65), (12.6.67), (12.6.68), (12.6.69). From equalities (12.6.66), (12.6.70) it follows that the equality (12.6.65) is covariant with respect to choice of basis of right $A$-vector space $V$.  □

# Covariance in $A$-Vector Space

## 13.1. Geometric Object of Left $A$-Vector Space

Let $V$, $W$ be left $A$-vector spaces of columns and $G(V_*)$ be symmetry group of left $A$-vector space $V$. Homomorphism

$$(13.1.1) \qquad F : G(V_*) \to GL(W_*)$$

maps passive transformation $g \in G(V_*)$

$$(13.1.2) \qquad \overline{\overline{e}}_{V2} = g^*{}_*\overline{\overline{e}}_{V1}$$

of left $A$-vector spaces $V$ into passive transformation $F(g) \in GL(W_*)$

$$(13.1.3) \qquad \overline{\overline{e}}_{W2} = F(g)^*{}_*\overline{\overline{e}}_{W1}$$

of left $A$-vector spaces $W$.

Then coordinate transformation in left $A$-vector space $W$ gets form

$$(13.1.4) \qquad w_2 = w_1{}^*{}_*F(g)^{-1^*{}_*}$$

Therefore, the map

$$(13.1.5) \qquad F_1(g) = F(g)^{-1^*{}_*}$$

is right-side representation

$$F_1 : G(V_*) \xrightarrow{\quad * \quad} W_*$$

of group $G(V_*)$ in the set $W_*$.

Let $V$, $W$ be left $A$-vector spaces of rows and $G(V_*)$ be symmetry group of left $A$-vector space $V$. Homomorphism

$$(13.1.6) \qquad F : G(V_*) \to GL(W_*)$$

maps passive transformation $g \in G(V_*)$

$$(13.1.7) \qquad \overline{\overline{e}}_{V2} = g_*{}^*\overline{\overline{e}}_{V1}$$

of left $A$-vector spaces $V$ into passive transformation $F(g) \in GL(W_*)$

$$(13.1.8) \qquad \overline{\overline{e}}_{W2} = F(g)_*{}^*\overline{\overline{e}}_{W1}$$

of left $A$-vector spaces $W$.

Then coordinate transformation in left $A$-vector space $W$ gets form

$$(13.1.9) \qquad w_2 = w_1{}_*{}^*F(g)^{-1^*}$$

Therefore, the map

$$(13.1.10) \qquad F_1(g) = F(g)^{-1^*}$$

is right-side representation

$$F_1 : G(V_*) \xrightarrow{\quad * \quad} W_*$$

of group $G(V_*)$ in the set $W_*$.

**DEFINITION 13.1.1.** *Orbit*

$$(13.1.11) \quad \begin{aligned} &\color{purple}{\mathcal{O}(V, W, \overline{\overline{e}}_V, w)} \\ &= (w^*{}_*F(G)^{-1^*{}_*}, G^*{}_*\overline{\overline{e}}_V) \end{aligned}$$

*of representation $F_1$ is called* **coordinate manifold of geometric object** *in left $A$-vector space $V$ of columns. For any basis*

$$\boxed{(13.1.2) \quad \overline{\overline{e}}_{V2} = g^*{}_*\overline{\overline{e}}_{V1}}$$

*corresponding point*

$$\boxed{(13.1.4) \quad w_2 = w_1{}^*{}_*F(g)^{-1^*{}_*}}$$

*of orbit defines* **coordinates of geometric object** *in coordinate left $A$-vector space relative basis $\overline{\overline{e}}_{V2}$.* □

**DEFINITION 13.1.2.** *Orbit*

$$(13.1.12) \quad \begin{aligned} &\color{purple}{\mathcal{O}(V, W, \overline{\overline{e}}_V, w)} \\ &= (w_*{}^*F(G)^{-1^*}, G_*{}^*\overline{\overline{e}}_V) \end{aligned}$$

*of representation $F_1$ is called* **coordinate manifold of geometric object** *in left $A$-vector space $V$ of rows. For any basis*

$$\boxed{(13.1.7) \quad \overline{\overline{e}}_{V2} = g_*{}^*\overline{\overline{e}}_{V1}}$$

*corresponding point*

$$\boxed{(13.1.9) \quad w_2 = w_1{}_*{}^*F(g)^{-1^*}}$$

*of orbit defines* **coordinates of geometric object** *in coordinate left $A$-vector space relative basis $\overline{\overline{e}}_{V2}$.* □

**DEFINITION 13.1.3.** *Let us say the coordinates $w_1$ of vector $\overline{w}$ with respect*

**DEFINITION 13.1.4.** *Let us say the coordinates $w_1$ of vector $\overline{w}$ with respect*

*to the basis  $\overline{\overline{e}}_{W1}$  are given.  The set of vectors*

(13.1.13)            $\overline{w}_2 = w_2{}^*{}_* e_{W2}$

*is called* **geometric object** *defined in left $A$-vector space $V$ of columns.  For any basis  $\overline{\overline{e}}_{W2}$,  corresponding point*

$$\boxed{(13.1.4) \quad w_2 = w_1{}^*{}_* F(g)^{-1}{}^*{}_*}$$

*of coordinate manifold defines the vector*

(13.1.14)            $\overline{w}_2 = w_2{}^*{}_* e_{W2}$

*which is called* **representative of geometric object** *in left $A$-vector space $V$ in basis  $\overline{\overline{e}}_{V2}$.* □

*to the basis  $\overline{\overline{e}}_{W1}$  are given.  The set of vectors*

(13.1.15)            $\overline{w}_2 = w_{2*}{}^* e_{W2}$

*is called* **geometric object** *defined in left $A$-vector space $V$ of rows.  For any basis  $\overline{\overline{e}}_{W2}$,  corresponding point*

$$\boxed{(13.1.9) \quad w_2 = w_{1*}{}^* F(g)^{-1}{}_*{}^*}$$

*of coordinate manifold defines the vector*

(13.1.16)            $\overline{w}_2 = w_{2*}{}^* e_{W2}$

*which is called* **representative of geometric object** *in left $A$-vector space $V$ in basis  $\overline{\overline{e}}_{V2}$.* □

Since a geometric object is an orbit of representation, we see that according to theorem 2.5.13 the definition of the geometric object is a proper definition.

We also say that $\overline{w}$ is a **geometric object of type** $F$.

Definitions 13.1.1, 13.1.2 introduce a geometric object in coordinate space. We assume in definitions 13.1.3, 13.1.4 that we selected a basis of vector space $W$. This allows using a representative of the geometric object instead of its coordinates.

THEOREM 13.1.5. **(Principle of covariance).**  *Representative of geometric object does not depend on selection of basis  $\overline{\overline{e}}_{V2}$.*

PROOF. To define representative of geometric object, we need to select basis $\overline{\overline{e}}_V$, basis $\overline{\overline{e}}_W$ and coordinates of geometric object $w^\alpha$.  Corresponding representative of geometric object has form

$$\overline{w} = w^*{}_* e_W$$

Suppose we map basis $\overline{\overline{e}}_V$ to basis $\overline{\overline{e}}'_V$ by passive transformation

$$\overline{\overline{e}}'_V = g^*{}_* \overline{\overline{e}}_V$$

According to construction this generates passive transformation

$$\boxed{(13.1.3) \quad \overline{\overline{e}}_{W2} = F(g)^*{}_* \overline{\overline{e}}_{W1}}$$

and coordinate transformation

$$\boxed{(13.1.4) \quad w_2 = w_1{}^*{}_* F(g)^{-1}{}^*{}_*}$$

Corresponding representative of geometric object has form

$$\overline{w}' = w'^*{}_* e'_W$$
$$= w^*{}_* F(a)^{-1}{}^*{}_*{}^*{}_* F(a)^*{}_* e_W$$
$$= w^*{}_* e_W = \overline{w}$$

PROOF. To define representative of geometric object, we need to select basis $\overline{\overline{e}}_V$, basis $\overline{\overline{e}}_W$ and coordinates of geometric object $w^\alpha$.  Corresponding representative of geometric object has form

$$\overline{w} = w_*{}^* e_W$$

Suppose we map basis $\overline{\overline{e}}_V$ to basis $\overline{\overline{e}}'_V$ by passive transformation

$$\overline{\overline{e}}'_V = g_*{}^* \overline{\overline{e}}_V$$

According to construction this generates passive transformation

$$\boxed{(13.1.8) \quad \overline{\overline{e}}_{W2} = F(g)_*{}^* \overline{\overline{e}}_{W1}}$$

and coordinate transformation

$$\boxed{(13.1.9) \quad w_2 = w_{1*}{}^* F(g)^{-1}{}_*{}^*}$$

Corresponding representative of geometric object has form

$$\overline{w}' = w'_*{}^* e'_W$$
$$= w_*{}^* F(a)^{-1}{}_*{}^*{}_*{}^* F(a)_*{}^* e_W$$
$$= w_*{}^* e_W = \overline{w}$$

Therefore representative of geometric object is invariant relative selection of basis.  □

Therefore representative of geometric object is invariant relative selection of basis.  □

DEFINITION 13.1.6. *Let $V$ be a left $A$-vector space of columns and*

$$\overline{w}_1 = w_1{}^*{}_*e_W$$

$$\overline{w}_2 = w_2{}^*{}_*e_W$$

*be geometric objects of the same type defined in left $A$-vector space $V$. Geometric object*

$$\overline{w} = (w_1 + w_2)^*{}_*e_W$$

*is called* **sum**

$$\overline{w} = \overline{w}_1 + \overline{w}_2$$

**of geometric objects** $\overline{w}_1$ *and* $\overline{w}_2$.  □

DEFINITION 13.1.7. *Let $V$ be a left $A$-vector space of rows and*

$$\overline{w}_1 = w_{1*}{}^*e_W$$

$$\overline{w}_2 = w_{2*}{}^*e_W$$

*be geometric objects of the same type defined in left $A$-vector space $V$. Geometric object*

$$\overline{w} = (w_1 + w_2)_*{}^*e_W$$

*is called* **sum**

$$\overline{w} = \overline{w}_1 + \overline{w}_2$$

**of geometric objects** $\overline{w}_1$ *and* $\overline{w}_2$.  □

DEFINITION 13.1.8. *Let $V$ be a left $A$-vector space of columns and*

$$\overline{w}_1 = w_1{}^*{}_*e_W$$

*be geometric object defined in left $A$-vector space $V$. Geometric object*

$$\overline{w}_2 = (kw_1)^*{}_*e_W$$

*is called* **product**

$$\overline{w}_2 = k\overline{w}_1$$

**of geometric object** $\overline{w}_1$ **and constant** $k \in A$.  □

DEFINITION 13.1.9. *Let $V$ be a left $A$-vector space of rows and*

$$\overline{w}_1 = w_{1*}{}^*e_W$$

*be geometric object defined in left $A$-vector space $V$. Geometric object*

$$\overline{w}_2 = (kw_1)_*{}^*e_W$$

*is called* **product**

$$\overline{w}_2 = k\overline{w}_1$$

**of geometric object** $\overline{w}_1$ **and constant** $k \in A$.  □

THEOREM 13.1.10. *Geometric objects of type $F$ defined in left $A$-vector space $V$ of columns form left $A$-vector space of columns.*

THEOREM 13.1.11. *Geometric objects of type $F$ defined in left $A$-vector space $V$ of rows form left $A$-vector space of rows.*

PROOF. The statement of theorems 13.1.10, 13.1.11 follows from immediate verification of the properties of vector space.  □

## 13.2. Geometric Object of Right $A$-Vector Space

Let $V$, $W$ be right $A$-vector spaces of columns and $G(V_*)$ be symmetry group of right $A$-vector space $V$. Homomorphism

$$(13.2.1) \qquad F : G(V_*) \to GL(W_*)$$

Let $V$, $W$ be right $A$-vector spaces of rows and $G(V_*)$ be symmetry group of right $A$-vector space $V$. Homomorphism

$$(13.2.6) \qquad F : G(V_*) \to GL(W_*)$$

maps passive transformation $g \in G(V_*)$

$$(13.2.2) \qquad \overline{\overline{e}}_{V2} = \overline{\overline{e}}_{V1*}{}^*g$$

of right $A$-vector spaces $V$ into passive transformation $F(g) \in GL(W_*)$

$$(13.2.3) \qquad \overline{\overline{e}}_{W2} = \overline{\overline{e}}_{W1*}{}^*F(g)$$

of right $A$-vector spaces $W$.

Then coordinate transformation in right $A$-vector space $W$ gets form

$$(13.2.4) \qquad w_2 = F(g)^{-1*}{}_*{}^*w_1$$

Therefore, the map

$$(13.2.5) \qquad F_1(g) = F(g)^{-1*}{}^*$$

is left-side representation

$$F_1 : G(V_*) \longrightarrow{}^* W_*$$

of group $G(V_*)$ in the set $W_*$.

maps passive transformation $g \in G(V_*)$

$$(13.2.7) \qquad \overline{\overline{e}}_{V2} = \overline{\overline{e}}_{V1}{}^*{}_*g$$

of right $A$-vector spaces $V$ into passive transformation $F(g) \in GL(W_*)$

$$(13.2.8) \qquad \overline{\overline{e}}_{W2} = \overline{\overline{e}}_{W1}{}^*{}_*F(g)$$

of right $A$-vector spaces $W$.

Then coordinate transformation in right $A$-vector space $W$ gets form

$$(13.2.9) \qquad w_2 = F(g)^{-1*}{}^*{}_*w_1$$

Therefore, the map

$$(13.2.10) \qquad F_1(g) = F(g)^{-1*}{}_*$$

is left-side representation

$$F_1 : G(V_*) \longrightarrow{}^* W_*$$

of group $G(V_*)$ in the set $W_*$.

DEFINITION 13.2.1. *Orbit*

$$(13.2.11) \qquad \begin{aligned} &\mathcal{O}(V, W, \overline{\overline{e}}_V, w) \\ &= (F(G)^{-1*}{}_*{}^*w, \overline{\overline{e}}_{V*}{}^*G) \end{aligned}$$

*of representation $F_1$ is called* **coordinate manifold of geometric object** *in right $A$-vector space $V$ of columns. For any basis*

$$\boxed{(13.2.2) \quad \overline{\overline{e}}_{V2} = \overline{\overline{e}}_{V1*}{}^*g}$$

*corresponding point*

$$\boxed{(13.2.4) \quad w_2 = F(g)^{-1*}{}_*{}^*w_1}$$

*of orbit defines* **coordinates of geometric object** *in coordinate right $A$-vector space relative basis* $\overline{\overline{e}}_{V2}$ . $\square$

DEFINITION 13.2.2. *Orbit*

$$(13.2.12) \qquad \begin{aligned} &\mathcal{O}(V, W, \overline{\overline{e}}_V, w) \\ &= (F(G)^{-1*}{}^*{}_*w, \overline{\overline{e}}_{V*}{}^*G) \end{aligned}$$

*of representation $F_1$ is called* **coordinate manifold of geometric object** *in right $A$-vector space $V$ of rows. For any basis*

$$\boxed{(13.2.7) \quad \overline{\overline{e}}_{V2} = \overline{\overline{e}}_{V1}{}^*{}_*g}$$

*corresponding point*

$$\boxed{(13.2.9) \quad w_2 = F(g)^{-1*}{}^*{}_*w_1}$$

*of orbit defines* **coordinates of geometric object** *in coordinate right $A$-vector space relative basis* $\overline{\overline{e}}_{V2}$ . $\square$

DEFINITION 13.2.3. *Let us say the coordinates $w_1$ of vector $\overline{w}$ with respect to the basis* $\overline{\overline{e}}_{W1}$ *are given. The set of vectors*

$$(13.2.13) \qquad \overline{w}_2 = e_{W2*}{}^*w_2$$

*is called* **geometric object** *defined in right $A$-vector space $V$ of columns. For any basis* $\overline{\overline{e}}_{W2}$, *corresponding point*

$$\boxed{(13.2.4) \quad w_2 = F(g)^{-1*}{}_*{}^*w_1}$$

*of coordinate manifold defines the vector*

$$(13.2.14) \qquad \overline{w}_2 = e_{W2*}{}^*w_2$$

DEFINITION 13.2.4. *Let us say the coordinates $w_1$ of vector $\overline{w}$ with respect to the basis* $\overline{\overline{e}}_{W1}$ *are given. The set of vectors*

$$(13.2.15) \qquad \overline{w}_2 = e_{W2}{}^*{}_*w_2$$

*is called* **geometric object** *defined in right $A$-vector space $V$ of rows. For any basis* $\overline{\overline{e}}_{W2}$, *corresponding point*

$$\boxed{(13.2.9) \quad w_2 = F(g)^{-1*}{}^*{}_*w_1}$$

*of coordinate manifold defines the vector*

$$(13.2.16) \qquad \overline{w}_2 = e_{W2}{}^*{}_*w_2$$

*which is called* **representative of geometric object** *in right $A$-vector space $V$ in basis* $\overline{\overline{e}}_{V2}$. $\square$

Since a geometric object is an orbit of representation, we see that according to theorem 2.5.13 the definition of the geometric object is a proper definition.

We also say that $\overline{w}$ is a **geometric object of type** $F$.

Definitions 13.2.1, 13.2.2 introduce a geometric object in coordinate space. We assume in definitions 13.2.3, 13.2.4 that we selected a basis of vector space $W$. This allows using a representative of the geometric object instead of its coordinates.

**THEOREM 13.2.5. (Principle of covariance).** *Representative of geometric object does not depend on selection of basis* $\overline{\overline{e}}_{V2}$.

PROOF. To define representative of geometric object, we need to select basis $\overline{\overline{e}}_V$, basis $\overline{\overline{e}}_W$ and coordinates of geometric object $w^\alpha$. Corresponding representative of geometric object has form

$$\overline{w} = e_{W*}{}^* w$$

Suppose we map basis $\overline{\overline{e}}_V$ to basis $\overline{\overline{e}}'_V$ by passive transformation

$$\overline{\overline{e}}'_V = \overline{\overline{e}}_{V*}{}^* g$$

According construction this generates passive transformation

$$\boxed{(13.2.3) \quad \overline{\overline{e}}_{W2} = \overline{\overline{e}}_{W1*}{}^* F(g)}$$

and coordinate transformation

$$\boxed{(13.2.4) \quad w_2 = F(g)^{-1}{}^*{}_*{}^* w_1}$$

Corresponding representative of geometric object has form

$$\overline{w}' = e'_{W*}{}^* w'$$
$$= e_{W*}{}^* F(a)_*{}^* F(a)^{-1}{}^*{}_*{}^* w$$
$$= e_{W*}{}^* w = \overline{w}$$

Therefore representative of geometric object is invariant relative selection of basis. $\square$

PROOF. To define representative of geometric object, we need to select basis $\overline{\overline{e}}_V$, basis $\overline{\overline{e}}_W$ and coordinates of geometric object $w^\alpha$. Corresponding representative of geometric object has form

$$\overline{w} = e_{W}{}^*{}_* w$$

Suppose we map basis $\overline{\overline{e}}_V$ to basis $\overline{\overline{e}}'_V$ by passive transformation

$$\overline{\overline{e}}'_V = \overline{\overline{e}}_V{}^*{}_* g$$

According construction this generates passive transformation

$$\boxed{(13.2.8) \quad \overline{\overline{e}}_{W2} = \overline{\overline{e}}_{W1}{}^*{}_* F(g)}$$

and coordinate transformation

$$\boxed{(13.2.9) \quad w_2 = F(g)^{-1}{}^*{}_*{}^* w_1}$$

Corresponding representative of geometric object has form

$$\overline{w}' = e'_{W}{}^*{}_* w'$$
$$= e_{W}{}^*{}_* F(a)^*{}_* F(a)^{-1}{}^*{}^*{}_* w$$
$$= e_{W}{}^*{}_* w = \overline{w}$$

Therefore representative of geometric object is invariant relative selection of basis. $\square$

**DEFINITION 13.2.6.** *Let $V$ be a right $A$-vector space of columns and*

$$\overline{w}_1 = e_{W*}{}^* w_1$$
$$\overline{w}_2 = e_{W*}{}^* w_2$$

*be geometric objects of the same type defined in right $A$-vector space $V$. Geomet-*

**DEFINITION 13.2.7.** *Let $V$ be a right $A$-vector space of rows and*

$$\overline{w}_1 = e_{W}{}^*{}_* w_1$$
$$\overline{w}_2 = e_{W}{}^*{}_* w_2$$

*be geometric objects of the same type defined in right $A$-vector space $V$. Geomet-*

| *ric object* | *ric object* |
| --- | --- |
| $$\overline{w} = e_{W*}{}^*(w_1 + w_2)$$ | $$\overline{w} = e_{W}{}^*{}_*(w_1 + w_2)$$ |
| *is called* **sum** | *is called* **sum** |
| $$\overline{w} = \overline{w}_1 + \overline{w}_2$$ | $$\overline{w} = \overline{w}_1 + \overline{w}_2$$ |
| **of geometric objects** $\overline{w}_1$ *and* $\overline{w}_2$.  □ | **of geometric objects** $\overline{w}_1$ *and* $\overline{w}_2$.  □ |

| DEFINITION 13.2.8. *Let $V$ be a right $A$-vector space of columns and* | DEFINITION 13.2.9. *Let $V$ be a right $A$-vector space of rows and* |
| --- | --- |
| $$\overline{w}_1 = e_{W*}{}^*w_1$$ | $$\overline{w}_1 = e_{W}{}^*{}_*w_1$$ |
| *be geometric object defined in right $A$-vector space $V$. Geometric object* | *be geometric object defined in right $A$-vector space $V$. Geometric object* |
| $$\overline{w}_2 = e_{W*}{}^*(w_1 k)$$ | $$\overline{w}_2 = e_{W}{}^*{}_*(w_1 k)$$ |
| *is called* **product** | *is called* **product** |
| $$\overline{w}_2 = \overline{w}_1 k$$ | $$\overline{w}_2 = \overline{w}_1 k$$ |
| **of geometric object** $\overline{w}_1$ **and constant** $k \in A$.  □ | **of geometric object** $\overline{w}_1$ **and constant** $k \in A$.  □ |

| THEOREM 13.2.10. *Geometric objects of type $F$ defined in right $A$-vector space $V$ of columns form right $A$-vector space of columns.* | THEOREM 13.2.11. *Geometric objects of type $F$ defined in right $A$-vector space $V$ of rows form right $A$-vector space of rows.* |
| --- | --- |

PROOF. The statement of theorems 13.2.10, 13.2.11 follows from immediate verification of the properties of vector space.  □

## 13.3. Summary

The question how large a diversity of geometric objects is well studied in case of vector spaces. However it is not such obvious in case of left or right $A$-vector spaces. As we can see from table 13.3.1 left $A$-vector space of columns and right $A$-vector space of rows have common symmetry group $GL(n^*{}_*A)$. This allows study a geometric object in left $A$-vector space of columns when we study passive representation in right $A$-vector space of rows. Can we study the same time a geometric object in left $A$-vector space of rows. At first glance, the answer is negative based theorem 9.5.9. However, equation (8.1.24) determines required linear map between $GL(n_*{}^*A)$ and $GL(n^*{}_*A)$.

THEOREM 13.3.1. *Basis manifold of left $A$-vector space $V$ and basis manifold of right $A$-vector space $V$ are different*

$$GL(n_*{}^*A)_*{}^*E_n \ne E_n{}^*{}_*GL(n^*{}_*A)$$

PROOF.

| LEMMA 13.3.2. *We can identify basis manifold* $GL(n_*{}^*A)(V)_*{}^*\overline{E_n}$ *and* | LEMMA 13.3.3. *We can identify basis manifold* $\overline{E_n}{}^*{}_*GL(n^*{}_*A)(V)$ *and* |
| --- | --- |

TABLE 13.3.1. Active and Passive Representations
in $A$-vector space

| Vector space | Group of Representation | Active Representation | Passive Representation |
|---|---|---|---|
| left $A$-vector space of rows | $GL(n_{*}{}^{*}A)$ | right-side | left-side |
| left $A$-vector space of columns | $GL(n^{*}{}_{*}A)$ | right-side | left-side |
| right $A$-vector space of rows | $GL(n^{*}{}_{*}A)$ | left-side | right-side |
| right $A$-vector space of columns | $GL(n_{*}{}^{*}A)$ | left-side | right-side |

*the set of $_{*}{}^{*}$-regular matrices $GL(n_{*}{}^{*}A)$.*

PROOF. According to the remark
12.5.20, we can identify a basis $\overline{\overline{e}}$ of right
$A$-vector space $V$ and the matrix $e$ of
coordinates of basis $\overline{\overline{e}}$ withb respect to
the basis $\overline{\overline{E_n}}$

$$\overline{\overline{e}} = e^{*}{}_{*}\overline{\overline{E_n}}$$

⊙

*the set of $^{*}{}_{*}$-regular matrices $GL(n^{*}{}_{*}A)$.*

PROOF. According to the remark
12.2.20, we can identify a basis $\overline{\overline{e}}$ of left
$A$-vector space $V$ and the matrix $e$ of
coordinates of basis $\overline{\overline{e}}$ withb respect to
the basis $\overline{\overline{E_n}}$

$$\overline{\overline{e}} = \overline{\overline{E_{n*}}}{}^{*}e$$

⊙

From lemmas 13.3.2, 13.3.3, it follows that if the set of vectors $(f_i, i \in I, |I| = n)$ is basis of left $A$-vector space of columns and basis of right $A$-vector space of columns, then their coordinate matrix

$$(13.3.1) \qquad f = \begin{pmatrix} f_1^1 & ... & f_1^n \\ ... & ... & ... \\ f_n^1 & ... & f_n^n \end{pmatrix}$$

is $_{*}{}^{*}$-nonsingular and $^{*}{}_{*}$-nonsingular matrix. The statement of the theorem follows
from the theorem 9.5.9.                                                    □

# Linear Map of $A$-Vector Space

## 14.1. Linear Map

DEFINITION 14.1.1. *Let* $A_i, i = 1, 2,$ *be division algebra over commutative ring* $D_i$. *Let* $V_i, i = 1, 2,$ *be* $A_i$-*vector space. Morphism of diagram of representations*

*into diagram of representations*

*is called* **linear map** *of* $A_1$-*vector spaces* $V_1$ *into* $A_2$-*vector spaces* $V_2$. *Let us denote* $\mathcal{L}(D_1 \to D_2; A_1 \to A_2; V_1 \to V_2)$ *set of linear maps of* $A_1$-*vector spaces module* $V_1$ *into* $A_2$-*vector spaces* $V_2$. $\qquad \square$

From the definition 14.1.1, it follows that linear map of $A_1$-vector spaces $V_1$ into $A_2$-vector spaces $V_2$ is linear map of $D_1$-vector spaces $V_1$ into $D_2$-vector spaces $V_2$.

THEOREM 14.1.2. *Let* $V^1, ..., V^n$ *be* $A$-*vector spaces and*

$$V = V^1 \oplus ... \oplus V^n$$

*Let us represent* $v$-*number*

$$v = v^1 \oplus ... \oplus v^n$$

*as column vector*

$$v = \begin{pmatrix} v^1 \\ ... \\ v^n \end{pmatrix}$$

*Let us represent a linear map*

$$f : V \to W$$

*as row vector*

$$f = \begin{pmatrix} f_1 & ... & f_n \end{pmatrix}$$

$$f_i : V^i \to W$$

*Then we can represent value of the map* $f$ *in* $V$-*number* $v$ *as product of matrices*

$$(14.1.1) \qquad f \circ v = \begin{pmatrix} f_1 & ... & f_n \end{pmatrix} \circ \begin{pmatrix} v^1 \\ ... \\ v^n \end{pmatrix} = f_i \circ v^i$$

PROOF. The theorem follows from definitions    (10.6.12),    (11.6.12).    □

THEOREM 14.1.3. *Let  $W^1, ..., W^m$  be A-vector spaces and*

$$W = W^1 \oplus ... \oplus W^m$$

*Let us represent W-number*

$$w = w^1 \oplus ... \oplus w^m$$

*as column vector*

$$w = \begin{pmatrix} w^1 \\ ... \\ w^m \end{pmatrix}$$

*Then the linear map*

$$f : V \to W$$

*has representation as column vector of maps*

$$f = \begin{pmatrix} f^1 \\ ... \\ f^m \end{pmatrix}$$

*such way that, if  $w = f \circ v$,  then*

$$\begin{pmatrix} w^1 \\ ... \\ w^m \end{pmatrix} = \begin{pmatrix} f^1 \\ ... \\ f^m \end{pmatrix} \circ v = \begin{pmatrix} f^1 \circ v \\ ... \\ f^m \circ v \end{pmatrix}$$

PROOF. The theorem follows from the theorem 2.2.8.    □

THEOREM 14.1.4. *Let  $V^1, ..., V^n$,   $W^1, ..., W^m$  be A-vector spaces and*

$$V = V^1 \oplus ... \oplus V^n$$

$$W = W^1 \oplus ... \oplus W^m$$

*Let us represent V-number*

$$v = v^1 \oplus ... \oplus v^n$$

*as column vector*

(14.1.2)
$$v = \begin{pmatrix} v^1 \\ ... \\ v^n \end{pmatrix}$$

*Let us represent W-number*

$$w = w^1 \oplus ... \oplus w^m$$

*as column vector*

$$(14.1.3) \qquad w = \begin{pmatrix} w^1 \\ ... \\ w^m \end{pmatrix}$$

*Then the linear map*

$$f : V \to W$$

*has representation as a matrix of maps*

$$(14.1.4) \qquad f = \begin{pmatrix} f_1^1 & ... & f_n^1 \\ ... & ... & ... \\ f_1^m & ... & f_n^m \end{pmatrix}$$

*such way that, if  $w = f \circ v$,  then*

$$(14.1.5) \qquad \begin{pmatrix} w^1 \\ ... \\ w^m \end{pmatrix} = \begin{pmatrix} f_1^1 & ... & f_n^1 \\ ... & ... & ... \\ f_1^m & ... & f_n^m \end{pmatrix} \circ \begin{pmatrix} v^1 \\ ... \\ v^n \end{pmatrix} = \begin{pmatrix} f_i^1 \circ v^i \\ ... \\ f_i^m \circ v^i \end{pmatrix}$$

*The map*

$$f_j^i : V^j \to W^i$$

*is a linear map and is called* **partial linear map**.

Proof. According to the theorem 14.1.3, there exists the set of linear maps

$$f^i : V \to W^i$$

such that

$$(14.1.6) \qquad \begin{pmatrix} w^1 \\ ... \\ w^m \end{pmatrix} = \begin{pmatrix} f^1 \\ ... \\ f^m \end{pmatrix} \circ v = \begin{pmatrix} f^1 \circ v \\ ... \\ f^m \circ v \end{pmatrix}$$

According to the theorem 14.1.2, for every $i$, there exists the set of linear maps

$$f_j^i : V^j \to W^i$$

such that

$$(14.1.7) \qquad f^i \circ v = \begin{pmatrix} f_1^i & ... & f_n^i \end{pmatrix} \circ \begin{pmatrix} v^1 \\ ... \\ v^n \end{pmatrix} = f_j^i \circ v^j$$

If we identify matrices

$$\begin{pmatrix} \begin{pmatrix} f_1^1 & ... & f_n^1 \end{pmatrix} \\ ... \\ \begin{pmatrix} f_1^m & ... & f_n^m \end{pmatrix} \end{pmatrix} = \begin{pmatrix} f_1^1 & ... & f_n^1 \\ ... & ... & ... \\ f_1^m & ... & f_n^m \end{pmatrix}$$

then the equality (14.1.5) follows from equalities (14.1.6), (14.1.7). $\qquad\square$

## 14.2. Matrix of Linear Map

DEFINITION 14.2.1. *Let* $A_i$, $i = 1, ..., n$, *be the set of $D$-algebras. The matrix*

$$f = \begin{pmatrix} f_1^1 & ... & f_n^1 \\ ... & ... & ... \\ f_1^m & ... & f_n^m \end{pmatrix}$$

*whose elements are maps*

$$f_j^i : A_i \to A_j$$

*is called* **matrix of maps**. $\qquad\square$

Consider following operations on the set of matrices of maps

- addition

$$\begin{pmatrix} f_1^1 & ... & f_n^1 \\ ... & ... & ... \\ f_1^m & ... & f_n^m \end{pmatrix} + \begin{pmatrix} g_1^1 & ... & g_n^1 \\ ... & ... & ... \\ g_1^m & ... & g_n^m \end{pmatrix} = \begin{pmatrix} f_1^1 + g_1^1 & ... & f_n^1 + g_n^1 \\ ... & ... & ... \\ f_1^m + g_1^m & ... & f_n^m + g_n^m \end{pmatrix}$$

- $_\circ{}^\circ$-product

$$\begin{pmatrix} f_1^1 & ... & f_k^1 \\ ... & ... & ... \\ f_1^m & ... & f_k^m \end{pmatrix} {}_\circ{}^\circ \begin{pmatrix} g_1^1 & ... & g_n^1 \\ ... & ... & ... \\ g_1^k & ... & g_n^k \end{pmatrix} = \begin{pmatrix} f_i^1 \circ g_1^i & ... & f_i^1 \circ g_n^i \\ ... & ... & ... \\ f_i^m \circ g_1^i & ... & f_i^m \circ g_n^i \end{pmatrix}$$

- ${}^\circ{}_\circ$-product

$$\begin{pmatrix} f_1^1 & ... & f_n^1 \\ ... & ... & ... \\ f_1^k & ... & f_n^k \end{pmatrix} {}^\circ{}_\circ \begin{pmatrix} g_1^1 & ... & g_k^1 \\ ... & ... & ... \\ g_1^m & ... & g_k^m \end{pmatrix} = \begin{pmatrix} f_1^i \circ g_i^1 & ... & f_n^i \circ g_i^1 \\ ... & ... & ... \\ f_1^i \circ g_i^m & ... & f_n^i \circ g_i^m \end{pmatrix}$$

# Linear Map of Left $A$-Vector Space

## 15.1. General Definition

Let $A_i$, $i = 1, 2$, be division algebra over commutative ring $D_i$. Let $V_i$, $i = 1, 2$, be left $A_i$-vector space.

> **Theorem 15.1.1.** *Linear map*
>
> (15.1.1) $$(h : D_1 \to D_2 \quad g : A_1 \to A_2 \quad f : V_1 \to V_2)$$
>
> *of left $A_1$-vector spaces $V_1$ into left $A_2$-vector spaces $V_2$ satisfies to equalities*
>
> (15.1.2) $$h(p + q) = h(p) + h(q)$$
>
> (15.1.3) $$h(pq) = h(p)h(q)$$
>
> (15.1.4) $$g \circ (a + b) = g \circ a + g \circ b$$
>
> (15.1.5) $$g \circ (ab) = (g \circ a)(g \circ b)$$
>
> (15.1.6) $$g \circ (pa) = h(p)(g \circ a)$$
>
> (15.1.7) $$f \circ (u + v) = f \circ u + f \circ v$$
>
> (15.1.8) $$f \circ (pv) = h(p)(f \circ v)$$
>
> $$p, q \in D_1 \quad a, b \in A_1 \quad u, v \in V_1$$

Proof. According to definitions 2.3.9, 2.8.5, 14.1.1, 11 the map

$$(h : D_1 \to D_2 \quad g : A_1 \to A_2)$$

is homomorphism of $D_1$-algebra $A_1$ into $D_2$-algebra $A_2$. Therefore, equalities (15.1.2), (15.1.3), (15.1.4), (15.1.5), (15.1.6) follow from the theorem 4.3.2. The equality (15.1.7) follows from the definition 10.1.1, since, accorfing to the definition 2.3.9, the map $f$ is homomorpism of Abelian group. The equality (15.1.8) follows from the equality (2.3.6) because the map

$$(h : D_1 \to D_2 \quad f : V_1 \to V_2)$$

is morphism of representation $g_{1.34}$ into representation $g_{2.34}$. $\square$

---

[15.1] In theorem 15.1.2, we use the following convention. Let $i = 1, 2$. Let the set of vectors $\overline{\overline{e}}_{A_i} = (e_{A_i k}, k \in K_i)$ be a basis and $C_{i \ ij}^{\ k}$, $k, i, j \in K_i$, be structural constants of $D_i$-algebra of columns $A_i$. Let the set of vectors $\overline{\overline{e}}_{V_i} = (e_{V_i i}, i \in I_i)$ be a basis of left $A_i$-vector space $V_i$.

Theorem 15.1.2. *Linear map* 15.1

(15.1.1)   $(h : D_1 \to D_2 \quad \overline{g} : A_1 \to A_2 \quad \overline{f} : V_1 \to V_2)$

*has presentation*

(15.1.9)   $$b = gh(a)$$

(15.1.10)   $$\overline{g} \circ (a^i e_{A_1 i}) = h(a^i) g_i^k e_{A_2 k}$$

(15.1.11)   $$\overline{g} \circ (e_{A_1} a) = e_{A_2} gh(a)$$

(15.1.12)   $$f \circ (v^*_{\;*} e_{V_1}) = (f_i^j \circ g \circ v^i) e_{V_2 \cdot j}$$

$$w = f_{\;\circ}{}^\circ (g \circ v)$$

*relative to selected bases. Here*

- *$a$ is coordinate matrix of $A_1$-number $\overline{a}$ relative the basis $\overline{\overline{e}}_{A_1}$*

$$\overline{a} = e_{A_1} a$$

- *$h(a) = (h(a_k), k \in K)$   is a matrix of $D_2$-numbers.*
- *$b$ is coordinate matrix of $A_2$-number*

$$\overline{b} = \overline{g} \circ \overline{a}$$

  *relative the basis $\overline{\overline{e}}_{A_2}$*

$$\overline{b} = e_{A_2} b$$

- *$g$ is coordinate matrix of set of $A_2$-numbers  $(\overline{g} \circ e_{A_1 k}, k \in K)$   relative the basis $\overline{\overline{e}}_{A_2}$.*
- *There is relation between the matrix of homomorphism and structural constants*

(15.1.13)   $$h(C1_{ij}^{\;k}) g_k^l = g_i^p g_j^q C2_{pq}^{\;\;l}$$

- *$v$ is coordinate matrix of $V_1$-number $\overline{v}$ relative the basis $\overline{\overline{e}}_{V_1}$*

$$\overline{v} = v^*_{\;*} e_{V_1}$$

- *$g(v) = (g(v_i), i \in I)$   is a matrix of $A_2$-numbers.*
- *$w$ is coordinate matrix of $V_2$-number*

$$\overline{w} = \overline{f} \circ \overline{v}$$

  *relative the basis $\overline{\overline{e}}_{V_2}$*

$$\overline{w} = w^*_{\;*} e_{V_2}$$

- *$f$ is coordinate matrix of set of $V_2$-numbers  $(\overline{f} \circ e_{V_1 i}, i \in I)$   relative the basis $\overline{\overline{e}}_{V_2}$.*

*The map*

$$f_j^i : A_2(g \circ e_{V_1 \cdot j}) \to A_2 e_i$$

*is a linear map and is called* **partial linear map**.

Proof.   According to definitions 2.3.9, 2.8.5, 14.1.1,  the map

(15.1.14)   $$(h : D_1 \to D_2 \quad \overline{g} : A_1 \to A_2)$$

is homomorphism of $D_1$-algebra $A_1$ into $D_2$-algebra $A_2$.   Therefore, equalities (15.1.9), (15.1.10), (15.1.11), (15.1.13)  follow from equalities

(5.4.11)   $b = fh(a)$

$$(5.4.12) \quad \overline{f} \circ (a^i e_{1i}) = h(a^i) f_i^k e_{2k}$$

$$(5.4.13) \quad \overline{f} \circ (e_1 a) = e_2 f h(a)$$

$$(5.4.16) \quad h(C_{1\,ij}^{\ \ k}) f_k^l = f_i^p f_j^q C_{2\,pq}^{\ \ l}$$

From the theorem 4.3.3, it follows that the matrix $g$ is unique.

$A_1$-module $V_1$ is direct sum

$$V_1 = \bigoplus_{k \in K} A_1 e_{V_1 \cdot k}$$

$A_2$-module $V_2$ is direct sum

$$V_2 = \bigoplus_{l \in L} A_2 e_{V_2 \cdot l}$$

The equality (15.1.12) follows from the theorem 14.1.4. $\qquad \square$

## 15.2. Linear Map When Rings $D_1 = D_2 = D$

Let $A_i$, $i = 1,\ 2$, be division algebra over commutative ring $D$. Let $V_i$, $i = 1,\ 2$, be left $A_i$-vector space.

THEOREM 15.2.1. *Linear map*

$$(15.2.1) \qquad (g : A_1 \to A_2 \quad f : V_1 \to V_2)$$

*of left $A_1$-vector spaces $V_1$ into left $A_2$-vector spaces $V_2$ satisfies to equalities*

$$(15.2.2) \qquad g \circ (a + b) = g \circ a + g \circ b$$

$$(15.2.3) \qquad g \circ (ab) = (g \circ a)(g \circ b)$$

$$(15.2.4) \qquad g \circ (pa) = p(g \circ a)$$

$$(15.2.5) \qquad f \circ (u + v) = f \circ u + f \circ v$$

$$(15.2.6) \qquad f \circ (pv) = p(f \circ v)$$

$$p \in D \quad a, b \in A_1 \quad u, v \in V_1$$

PROOF. According to definitions 2.3.9, 2.8.5, 14.1.1, 1 the map

$$g : A_1 \to A_2$$

is homomorphism of $D$-algebra $A_1$ into $D$-algebra $A_2$. Therefore, equalities (15.2.2), (15.2.3), (15.2.4) follow from the theorem 4.3.10. The equality (15.2.5) follows from the definition 10.2.1, since, accorfing to the definition 2.3.9, the map $f$ is homomorpism of Abelian group. The equality (15.2.6) follows from the equality (2.3.6) because the map

$$f : V_1 \to V_2$$

is morphism of representation $g_{1.34}$ into representation $g_{2.34}$. $\qquad \square$

---

<sup>15.2</sup> In theorem 15.2.2, we use the following convention. Let $i = 1,\ 2$. Let the set of vectors $\overline{\overline{e}}_{A_i} = (e_{A_i k}, k \in K_i)$ be a basis and $C_{1\,ij}^{\ \ k}$, $k,\ i,\ j \in K_i$, be structural constants of $D$-algebra of columns $A_i$. Let the set of vectors $\overline{\overline{e}}_{V_i} = (e_{V_i i}, i \in I_i)$ be a basis of left $A_i$-vector space $V_i$.

THEOREM 15.2.2. *Linear map* 15.2

$$(15.2.1) \quad (\overline{g} : A_1 \to A_2 \quad \overline{f} : V_1 \to V_2)$$

*has presentation*

$$(15.2.7) \qquad b = ga$$

$$(15.2.8) \qquad \overline{g} \circ (a^i e_{A_1 i}) = a^i g_i^k e_{A_2 k}$$

$$(15.2.9) \qquad \overline{g} \circ (e_{A_1} a) = e_{A_2} ga$$

$$(15.2.10) \qquad f \circ (v^* {}_* e_{V_1}) = (f_i^j \circ g \circ v^i) e_{V_2 \cdot j}$$

$$w = f_* \circ (g \circ v)$$

*relative to selected bases. Here*

- *$a$ is coordinate matrix of $A_1$-number $\overline{a}$ relative the basis $\overline{\overline{e}}_{A_1}$*

$$\overline{a} = e_{A_1} a$$

- *$b$ is coordinate matrix of $A_2$-number*

$$\overline{b} = \overline{g} \circ \overline{a}$$

  *relative the basis $\overline{\overline{e}}_{A_2}$*

$$\overline{b} = e_{A_2} b$$

- *$g$ is coordinate matrix of set of $A_2$-numbers $(\overline{g} \circ e_{A_1 k}, k \in K)$ relative the basis $\overline{\overline{e}}_{A_2}$.*
- *There is relation between the matrix of homomorphism and structural constants*

$$(15.2.11) \qquad C_{1 \, ij}^{\ \ k} g_k^l = g_i^p g_j^q C_{2 \, pq}^{\ \ l}$$

- *$v$ is coordinate matrix of $V_1$-number $\overline{v}$ relative the basis $\overline{\overline{e}}_{V_1}$*

$$\overline{v} = v^* {}_* e_{V_1}$$

- *$g(v) = (g(v_i), i \in I)$ is a matrix of $A_2$-numbers.*
- *$w$ is coordinate matrix of $V_2$-number*

$$\overline{w} = \overline{f} \circ \overline{v}$$

  *relative the basis $\overline{\overline{e}}_{V_2}$*

$$\overline{w} = w^* {}_* e_{V_2}$$

- *$f$ is coordinate matrix of set of $V_2$-numbers $(\overline{f} \circ e_{V_1 i}, i \in I)$ relative the basis $\overline{\overline{e}}_{V_2}$.*

*The map*

$$f_j^i : A_2(g \circ e_{V_1 \cdot j}) \to A_2 e_i$$

*is a linear map and is called **partial linear map**.*

PROOF.　According to definitions 2.3.9, 2.8.5, 14.1.1, the map

$$(15.2.12) \qquad \overline{g} : A_1 \to A_2$$

is homomorphism of $D$-algebra $A_1$ into $D$-algebra $A_2$. Therefore, equalities (15.2.7), (15.2.8), (15.2.9), (15.2.11) follow from equalities

$$(5.4.53) \quad b = fa$$

$$\boxed{(5.4.54) \quad \overline{f} \circ (a^i e_i) = a^i f_i^k e_k}$$

$$\boxed{(5.4.55) \quad \overline{f} \circ (e_1 a) = e_2 f a}$$

$$\boxed{(5.4.58) \quad C_1{}_{ij}^{\phantom{ij}k} f_k^l = f_i^p f_j^q C_2{}_{pq}^{\phantom{pq}l}}$$

From the theorem 4.3.11, it follows that the matrix $g$ is unique.

$A_1$-module $V_1$ is direct sum

$$V_1 = \bigoplus_{k \in K} A_1 e_{V_1 \cdot k}$$

$A_2$-module $V_2$ is direct sum

$$V_2 = \bigoplus_{l \in L} A_2 e_{V_2 \cdot l}$$

The equality (15.2.10) follows from the theorem 14.1.4.                    $\square$

## 15.3. Linear Map When $D$-algebras  $A_1 = A_2 = A$

Let $A$ be division algebra over commutative ring $D$.    Let $V_i$, $i = 1$, $2$,  be left $A$-vector space.

**Theorem 15.3.1.** *Linear map*

$$(15.3.1) \qquad\qquad f : V_1 \to V_2$$

*of left $A$-vector spaces $V_1$ into left $A$-vector spaces $V_2$ satisfies to equalities*

$$(15.3.2) \qquad\qquad f \circ (u + v) = f \circ u + f \circ v$$

$$(15.3.3) \qquad\qquad f \circ (av) = a(f \circ v)$$

$$a \in A \quad u, v \in V_1$$

**Proof.**      The equality (15.3.2) follows from the definition 10.3.1, since, accorfing to the definition 2.3.9, the map $f$ is homomorpism of Abelian group. The equality (15.3.3) follows from the equality (2.3.6) because the map

$$f : V_1 \to V_2$$

is morphism of representation $g_{1.34}$ into representation $g_{2.34}$.                    $\square$

**Theorem 15.3.2.** *Linear map* **15.3**

$$\boxed{(15.3.1) \quad \overline{f} : V_1 \to V_2}$$

*has presentation*

$$(15.3.4) \qquad\qquad f \circ (v^* {}_* e_{V_1}) = (f_i^j \circ v^i) e_{V_2 \cdot j}$$

$$w = f_\circ{}^\circ v$$

*relative to selected bases. Here*

- *$v$ is coordinate matrix of $V_1$-number $\overline{v}$ relative the basis  $\overline{\overline{e}}_{V_1}$*

$$\overline{v} = v^* {}_* e_{V_1}$$

---

**15.3**   In theorem 15.3.2, we use the following convention.    Let  $i = 1$, $2$.    Let the set of vectors $\overline{\overline{e}}_{V_i} = (e_{V_i i}, i \in I_i)$  be a basis of left $A$-vector space $V_i$.

- $w$ is coordinate matrix of $V_2$-number

$$\overline{w} = \overline{f} \circ \overline{v}$$

  relative the basis  $\overline{\overline{e}}_{V_2}$

$$\overline{w} = w^*{}_* e_{V_2}$$

- $f$ is coordinate matrix of set of $V_2$-numbers  $\left(\overline{f} \circ e_{V_1 i}, i \in I\right)$  relative the basis  $\overline{\overline{e}}_{V_2}$.

*The map*

$$f_j^i : A_2(g \circ e_{V_1 \cdot j}) \to A_2 e_i$$

*is a linear map and is called* **partial linear map**.

Proof.
$A_1$-module $V_1$ is direct sum

$$V_1 = \bigoplus_{k \in K} A_1 e_{V_1 \cdot k}$$

$A_2$-module $V_2$ is direct sum

$$V_2 = \bigoplus_{l \in L} A_2 e_{V_2 \cdot l}$$

The equality (15.3.4) follows from the theorem 14.1.4.                 $\square$

# Linear Map of Right $A$-Vector Space

## 16.1. General Definition

Let $A_i$, $i = 1$, $2$, be division algebra over commutative ring $D_i$. Let $V_i$, $i = 1$, $2$, be right $A_i$-vector space.

**Theorem 16.1.1.** *Linear map*

(16.1.1) $$(h : D_1 \rightarrow D_2 \quad g : A_1 \rightarrow A_2 \quad f : V_1 \rightarrow V_2)$$

*of right $A_1$-vector spaces $V_1$ into right $A_2$-vector spaces $V_2$ satisfies to equalities*

(16.1.2) $$h(p + q) = h(p) + h(q)$$

(16.1.3) $$h(pq) = h(p)h(q)$$

(16.1.4) $$g \circ (a + b) = g \circ a + g \circ b$$

(16.1.5) $$g \circ (ab) = (g \circ a)(g \circ b)$$

(16.1.6) $$g \circ (ap) = (g \circ a)h(p)$$

(16.1.7) $$f \circ (u + v) = f \circ u + f \circ v$$

(16.1.8) $$f \circ (vp) = (f \circ v)h(p)$$

$$p, q \in D_1 \quad a, b \in A_1 \quad u, v \in V_1$$

Proof. According to definitions 2.3.9, 2.8.5, 14.1.1, 11 the map

$$(h : D_1 \rightarrow D_2 \quad g : A_1 \rightarrow A_2)$$

is homomorphism of $D_1$-algebra $A_1$ into $D_2$-algebra $A_2$. Therefore, equalities (16.1.2), (16.1.3), (16.1.4), (16.1.5), (16.1.6) follow from the theorem 4.3.2. The equality (16.1.7) follows from the definition 11.1.1, since, accorfing to the definition 2.3.9, the map $f$ is homomorpism of Abelian group. The equality (16.1.8) follows from the equality (2.3.6) because the map

$$(h : D_1 \rightarrow D_2 \quad f : V_1 \rightarrow V_2)$$

is morphism of representation $g_{1.34}$ into representation $g_{2.34}$. $\qquad \square$

---

Theorem 16.1.2. *Linear map* 16.1

(16.1.1)    $(h : D_1 \to D_2 \quad \overline{g} : A_1 \to A_2 \quad \overline{f} : V_1 \to V_2)$

*has presentation*

(16.1.9)                $b = gh(a)$

(16.1.10)            $\overline{g} \circ (a^i e_{A_1 i}) = h(a^i) g_i^k e_{A_2 k}$

(16.1.11)            $\overline{g} \circ (e_{A_1} a) = e_{A_2} gh(a)$

(16.1.12)            $f \circ (e_{V_1 *}{}^* v) = e_{V_2 \cdot j} (f_i^j \circ g \circ v^i)$

                        $w = f_\circ{}^\circ (g \circ v)$

*relative to selected bases. Here*

- *a is coordinate matrix of $A_1$-number $\overline{a}$ relative the basis $\overline{\overline{e}}_{A_1}$*

                        $\overline{a} = e_{A_1} a$

- $h(a) = (h(a_k), k \in K)$  *is a matrix of $D_2$-numbers.*
- *b is coordinate matrix of $A_2$-number*

                        $\overline{b} = \overline{g} \circ \overline{a}$

  *relative the basis $\overline{\overline{e}}_{A_2}$*

                        $\overline{b} = e_{A_2} b$

- *g is coordinate matrix of set of $A_2$-numbers $(\overline{g} \circ e_{A_1 k}, k \in K)$ relative the basis $\overline{\overline{e}}_{A_2}$.*
- *There is relation between the matrix of homomorphism and structural constants*

(16.1.13)            $h(C_1{}_{ij}^k) g_k^l = g_i^p g_j^q C_2{}_{pq}^l$

- *v is coordinate matrix of $V_1$-number $\overline{v}$ relative the basis $\overline{\overline{e}}_{V_1}$*

                        $\overline{v} = e_{V_1 *}{}^* v$

- $g(v) = (g(v_i), i \in I)$  *is a matrix of $A_2$-numbers.*
- *w is coordinate matrix of $V_2$-number*

                        $\overline{w} = \overline{f} \circ \overline{v}$

  *relative the basis $\overline{\overline{e}}_{V_2}$*

                        $\overline{w} = e_{V_2 *}{}^* w$

- *f is coordinate matrix of set of $V_2$-numbers $(\overline{f} \circ e_{V_1 i}, i \in I)$ relative the basis $\overline{\overline{e}}_{V_2}$.*

*The map*

                        $f_j^i : (g \circ e_{V_1 \cdot j}) A_2 \to e_i A_2$

*is a linear map and is called* **partial linear map**.

Proof.    According to definitions  2.3.9, 2.8.5, 14.1.1,  the map

(16.1.14)                $(h : D_1 \to D_2 \quad \overline{g} : A_1 \to A_2)$

*is homomorphism of $D_1$-algebra $A_1$ into $D_2$-algebra $A_2$. Therefore, equalities (16.1.9), (16.1.10), (16.1.11), (16.1.13)  follow from equalities*

(5.4.11)    $b = fh(a)$

$(5.4.12)$   $\overline{f} \circ (a^i e_{1i}) = h(a^i) f_i^k e_{2k}$

$(5.4.13)$   $\overline{f} \circ (e_1 a) = e_2 f h(a)$

$(5.4.16)$   $h(C_{1ij}^{\ k}) f_k^l = f_i^p f_j^q C_{2pq}^{\ l}$

From the theorem 4.3.3, it follows that the matrix $g$ is unique.

$A_1$-module $V_1$ is direct sum

$$V_1 = \bigoplus_{k \in K} A_1 e_{V_1 \cdot k}$$

$A_2$-module $V_2$ is direct sum

$$V_2 = \bigoplus_{l \in L} A_2 e_{V_2 \cdot l}$$

The equality (16.1.12) follows from the theorem 14.1.4.                    $\square$

## 16.2. Linear Map When Rings  $D_1 = D_2 = D$

Let  $A_i$, $i = 1, 2$,  be division algebra over commutative ring $D$.   Let  $V_i$, $i = 1, 2$,  be right $A_i$-vector space.

THEOREM 16.2.1.  *Linear map*

$(16.2.1)$                $(g : A_1 \to A_2 \quad f : V_1 \to V_2)$

*of right $A_1$-vector spaces $V_1$ into right $A_2$-vector spaces $V_2$ satisfies to equalities*

$(16.2.2)$                $g \circ (a + b) = g \circ a + g \circ b$

$(16.2.3)$                $g \circ (ab) = (g \circ a)(g \circ b)$

$(16.2.4)$                $g \circ (ap) = (g \circ a)p$

$(16.2.5)$                $f \circ (u + v) = f \circ u + f \circ v$

$(16.2.6)$                $f \circ (vp) = (f \circ v)p$

$$p \in D \quad a, b \in A_1 \quad u, v \in V_1$$

PROOF.   According to definitions     2.3.9, 2.8.5,  14.1.1,  1 the map

$$g : A_1 \to A_2$$

is homomorphism of $D$-algebra $A_1$ into $D$-algebra $A_2$.  Therefore, equalities  (16.2.2), (16.2.3), (16.2.4)  follow from the theorem 4.3.10.   The equality (16.2.5) follows from the definition 11.2.1, since, accorfing to the definition 2.3.9, the map $f$ is homomorpism of Abelian group.  The equality (16.2.6) follows from the equality (2.3.6) because the map

$$f : V_1 \to V_2$$

is morphism of representation $g_{1.34}$ into representation $g_{2.34}$.                    $\square$

---

[16.2]   In theorem 16.2.2, we use the following convention.   Let  $i = 1, 2$.   Let the set of vectors $\overline{e}_{A_i} = (e_{A_i k}, k \in K_i)$  be a basis and  $C_{1ij}^{\ k}$,   $k, i, j \in K_i$,  be structural constants of $D$-algebra of columns $A_i$.   Let the set of vectors $\overline{e}_{V_i} = (e_{V_i i}, i \in I_i)$  be a basis of right $A_i$-vector space $V_i$.

THEOREM 16.2.2. *Linear map* **16.2**

$$\boxed{(16.2.1) \quad (\overline{g} : A_1 \to A_2 \quad \overline{f} : V_1 \to V_2)}$$

*has presentation*

$$(16.2.7) \qquad\qquad b = ga$$

$$(16.2.8) \qquad\qquad \overline{g} \circ (a^i e_{A_1 i}) = a^i g_i^k e_{A_2 k}$$

$$(16.2.9) \qquad\qquad \overline{g} \circ (e_{A_1} a) = e_{A_2} ga$$

$$(16.2.10) \qquad\qquad f \circ (e_{V_1*}{}^* v) = e_{V_2 \cdot j} (f_i^j \circ g \circ v^i)$$

$$w = f_\circ{}^\circ (g \circ v)$$

*relative to selected bases. Here*

- $a$ *is coordinate matrix of $A_1$-number $\overline{a}$ relative the basis* $\overline{\overline{e}}_{A_1}$

$$\overline{a} = e_{A_1} a$$

- $b$ *is coordinate matrix of $A_2$-number*

$$\overline{b} = \overline{g} \circ \overline{a}$$

  *relative the basis* $\overline{\overline{e}}_{A_2}$

$$\overline{b} = e_{A_2} b$$

- $g$ *is coordinate matrix of set of $A_2$-numbers* $(\overline{g} \circ e_{A_1 k}, k \in K)$ *relative the basis* $\overline{\overline{e}}_{A_2}$.
- *There is relation between the matrix of homomorphism and structural constants*

$$(16.2.11) \qquad\qquad C_1{}_{ij}^{\phantom{ij}k} g_k^l = g_i^p g_j^q C_2{}_{pq}^{\phantom{pq}l}$$

- $v$ *is coordinate matrix of $V_1$-number $\overline{v}$ relative the basis* $\overline{\overline{e}}_{V_1}$

$$\overline{v} = e_{V_1*}{}^* v$$

- $g(v) = (g(v_i), i \in I)$ *is a matrix of $A_2$-numbers.*
- $w$ *is coordinate matrix of $V_2$-number*

$$\overline{w} = \overline{f} \circ \overline{v}$$

  *relative the basis* $\overline{\overline{e}}_{V_2}$

$$\overline{w} = e_{V_2*}{}^* w$$

- $f$ *is coordinate matrix of set of $V_2$-numbers* $(\overline{f} \circ e_{V_1 i}, i \in I)$ *relative the basis* $\overline{\overline{e}}_{V_2}$.

*The map*

$$f_j^i : (g \circ e_{V_1 \cdot j}) A_2 \to e_i A_2$$

*is a linear map and is called* **partial linear map**.

PROOF. According to definitions 2.3.9, 2.8.5, 14.1.1, the map

$$(16.2.12) \qquad\qquad \overline{g} : A_1 \to A_2$$

*is homomorphism of $D$-algebra $A_1$ into $D$-algebra $A_2$.* Therefore, equalities (16.2.7), (16.2.8), (16.2.9), (16.2.11) follow from equalities

$$\boxed{(5.4.53) \quad b = fa}$$

$$\boxed{(5.4.54) \quad \overline{f} \circ (a^i e_{1i}) = a^i f_i^k e_k}$$

$$\boxed{(5.4.55) \quad \overline{f} \circ (e_1 a) = e_2 f a}$$

$$\boxed{(5.4.58) \quad C_1{}_{ij}^{\ k} f_k^l = f_i^p f_j^q C_2{}_{pq}^{\ l}}$$

From the theorem 4.3.11, it follows that the matrix $g$ is unique.

$A_1$-module $V_1$ is direct sum

$$V_1 = \bigoplus_{k \in K} A_1 e_{V_1 \cdot k}$$

$A_2$-module $V_2$ is direct sum

$$V_2 = \bigoplus_{l \in L} A_2 e_{V_2 \cdot l}$$

The equality (16.2.10) follows from the theorem 14.1.4. $\qquad\square$

## 16.3. Linear Map When $D$-algebras $A_1 = A_2 = A$

Let $A$ be division algebra over commutative ring $D$. Let $V_i$, $i = 1, 2$, be right $A$-vector space.

THEOREM 16.3.1. *Linear map*

$$(16.3.1) \qquad f : V_1 \to V_2$$

*of right $A$-vector spaces $V_1$ into right $A$-vector spaces $V_2$ satisfies to equalities*

$$(16.3.2) \qquad f \circ (u + v) = f \circ u + f \circ v$$

$$(16.3.3) \qquad f \circ (va) = (f \circ v)a$$

$$a \in A \quad u, v \in V_1$$

PROOF. The equality (16.3.2) follows from the definition 11.3.1, since, accorfing to the definition 2.3.9, the map $f$ is homomorpism of Abelian group. The equality (16.3.3) follows from the equality (2.3.6) because the map

$$f : V_1 \to V_2$$

is morphism of representation $g_{1.34}$ into representation $g_{2.34}$. $\qquad\square$

THEOREM 16.3.2. *Linear map* [16.3]

$$\boxed{(16.3.1) \quad \overline{f} : V_1 \to V_2}$$

*has presentation*

$$(16.3.4) \qquad f \circ (e_{V_1 *}{}^* v) = e_{V_2 \cdot j}(f_i^j \circ v^i)$$

$$w = f_\circ{}^\circ v$$

*relative to selected bases. Here*

- *$v$ is coordinate matrix of $V_1$-number $\overline{v}$ relative the basis $\overline{\overline{e}}_{V_1}$*

$$\overline{v} = e_{V_1 *}{}^* v$$

- $w$ is coordinate matrix of $V_2$-number

$$\overline{w} = \overline{f} \circ \overline{v}$$

  relative the basis $\overline{\overline{e}}_{V_2}$

$$\overline{w} = e_{V_2 *}{}^* w$$

- $f$ is coordinate matrix of set of $V_2$-numbers $\left( \overline{f} \circ e_{V_1 i}, i \in I \right)$ relative the basis $\overline{\overline{e}}_{V_2}$.

*The map*

$$f_j^i : (g \circ e_{V_1 \cdot j}) A_2 \to e_i A_2$$

*is a linear map and is called* **partial linear map**.

PROOF.
$A_1$-module $V_1$ is direct sum

$$V_1 = \bigoplus_{k \in K} A_1 e_{V_1 \cdot k}$$

$A_2$-module $V_2$ is direct sum

$$V_2 = \bigoplus_{l \in L} A_2 e_{V_2 \cdot l}$$

The equality (16.3.4) follows from the theorem 14.1.4.  $\square$

CHAPTER 17

# Similarity Transformation

Let $A$ be associative division $D$-algebra.

### 17.1. Similarity Transformation in Left $A$-Vector Space

**Theorem 17.1.1.** *Let $V$ be a left $A$-vector space of columns and $\overline{\overline{e}}$ be basis of left $A$-vector space $V$. Let $\dim V = n$ and $E_n$ be $n \times n$ identity matrix. For any $A$-number $a$, there exist endomorphism $\overline{\overline{e}}a$ of left $A$-vector space $V$ which has the matrix $E_n a$ with respecpect to the basis $\overline{\overline{e}}$.*

**Theorem 17.1.2.** *Let $V$ be a left $A$-vector space of rows and $\overline{\overline{e}}$ be basis of left $A$-vector space $V$. Let $\dim V = n$ and $E_n$ be $n \times n$ identity matrix. For any $A$-number $a$, there exist endomorphism $\overline{\overline{e}}a$ of left $A$-vector space $V$ which has the matrix $E_n a$ with respecpect to the basis $\overline{\overline{e}}$.*

**Proof.** Let $\overline{\overline{e}}$ be the basis of left $A$-vector space of columns. According to the theorem 10.3.8 and the definition 10.3.4, the map

$$\overline{\overline{e}}a : v^*{}_*e \in V \to v^*{}_*(E_n a)^*{}_* e \in V$$

is endomorphism of left $A$-vector space $V$. □

**Proof.** Let $\overline{\overline{e}}$ be the basis of left $A$-vector space of columns. According to the theorem 10.3.9 and the definition 10.3.4, the map

$$\overline{\overline{e}}a : v_*{}^*e \in V \to v_*{}^*(E_n a)_*{}^* e \in V$$

is endomorphism of left $A$-vector space $V$. □

**Definition 17.1.3.** *The endomorphism $\overline{\overline{e}}a$ of left $A$-vector space $V$ is called* **similarity transformation** *with respect to the basis $\overline{\overline{e}}$.* □

**Theorem 17.1.4.** *Let $V$ be a left $A$-vector space of columns and $\overline{\overline{e}}_1$, $\overline{\overline{e}}_2$ be bases of left $A$-vector space $V$. Let $g$ be passive transformation of basis $\overline{\overline{e}}_1$ into basis $\overline{\overline{e}}_2$*

$$\overline{\overline{e}}_2 = g^*{}_* \overline{\overline{e}}_1$$

*The similarity transformation $\overline{\overline{e}_1 a}$ has the matrix*

$$(17.1.1) \qquad g a^*{}_* g^{-1*}{}_*$$

*with respect to the basis $\overline{\overline{e}}_2$.*

**Theorem 17.1.5.** *Let $V$ be a left $A$-vector space of rows and $\overline{\overline{e}}_1$, $\overline{\overline{e}}_2$ be bases of left $A$-vector space $V$. Let $g$ be passive transformation of basis $\overline{\overline{e}}_1$ into basis $\overline{\overline{e}}_2$*

$$\overline{\overline{e}}_2 = g_*{}^* \overline{\overline{e}}_1$$

*The similarity transformation $\overline{\overline{e}_1 a}$ has the matrix*

$$(17.1.2) \qquad g a_*{}^* g^{-1}{}_*{}^*$$

*with respect to the basis $\overline{\overline{e}}_2$.*

Proof. According to the theorem 12.3.5, the similarity transformation $\overline{\overline{e}}_1 a$ has the matrix

$$(17.1.3) \qquad g^* {}_* E_n a^* {}_* g^{-1^* {}_*}$$

with respect to the basis $\overline{\overline{e}}_2$. The theorem follows from the expression (17.1.3). □

Proof. According to the theorem 12.3.6, the similarity transformation $\overline{\overline{e}}_1 a$ has the matrix

$$(17.1.4) \qquad g_* {}^* E_n a_* {}^* g^{-1_* {}^*}$$

with respect to the basis $\overline{\overline{e}}_2$. The theorem follows from the expression (17.1.4). □

Theorem 17.1.6. *Let $V$ be left $A$-vector space of columns and $\overline{\overline{e}}$ be basis of left $A$-vector space $V$. The set of similarity transformations $\overline{\overline{e}}a$ of left $A$-vector space $V$ is right-side effective representation of $D$-algebra $A$ in left $A$-vector space $V$.*

Proof. Consider the set of matrices $E_n a$ of similarity transformations $\overline{\overline{e}}a$ whith respecpect to the basis $\overline{\overline{e}}$. Let $V^*$ be the set of coordinates of $V$-numbers whith respecpect to the basis $\overline{\overline{e}}$. According to theorems 13.1.10, 13.1.11, the set $V^*$ is left $A$-vector space. According to the theorem 17.1.1, the map $\overline{\overline{e}}a$ is the endomorphism of left $A$-vector space of columns.

Lemma 17.1.7. *Let $V$ be a left $A$-vector space of columns. The map*

$$(17.1.5) \qquad a \in A \to E_n a \in \operatorname{End}(A, V^*)$$

*is homomorphism of $D$-algebra $A$.*

Lemma 17.1.8. *Let $V$ be a left $A$-vector space of rows. The map*

$$(17.1.6) \qquad a \in A \to E_n a \in \operatorname{End}(A, V^*)$$

*is homomorphism of $D$-algebra $A$.*

Proof. According to the theorem 5.4.8, the lemma follows from equalities

$$E_n(a + b) = E_n a + E_n b$$
$$E_n(ba) = b(E_n a)$$
$$E_n(ab) = (E_n a)^* {}_* (E_n b)$$

⊙

According to the definition 2.5.1 and the lemma 17.1.7, the map (17.1.5) is right-side representation of $D$-algebra $A$ in left $A$-vector space $V^*$. According to the theorem 5.1.31, the equality

$$E_n a = E_n b$$

implies $a = b$. According to the definition 2.3.2, the representation (17.1.5) is effective.

According to theorems 12.3.7, 12.3.9, the representation (17.1.5) is covariant with respect to choice of basis of left $A$-vector space $V$. According to the theorem 17.1.1, the set of maps $\overline{\overline{e}}a$ is right-side representation of $D$-algebra $A$ in left $A$-vector space $V$.

Proof. According to the theorem 5.4.8, the lemma follows from equalities

$$E_n(a + b) = E_n a + E_n b$$
$$E_n(ba) = b(E_n a)$$
$$E_n(ab) = (E_n a)_* {}^* (E_n b)$$

⊙

According to the definition 2.5.1 and the lemma 17.1.8, the map (17.1.6) is right-side representation of $D$-algebra $A$ in left $A$-vector space $V^*$. According to the theorem 5.1.31, the equality

$$E_n a = E_n b$$

implies $a = b$. According to the definition 2.3.2, the representation (17.1.6) is effective.

According to theorems 12.3.8, 12.3.10, the representation (17.1.6) is covariant with respect to choice of basis of left $A$-vector space $V$. According to the theorem 17.1.2, the set of maps $\overline{\overline{e}}a$ is right-side representation of $D$-algebra $A$ in left $A$-vector space $V$.

## 17.2. Similarity Transformation in Right $A$-Vector Space

**Theorem 17.2.1.** *Let $V$ be a right $A$-vector space of columns and $\overline{\overline{e}}$ be basis of right $A$-vector space $V$. Let $\dim V = n$ and $E_n$ be $n \times n$ identity matrix. For any $A$-number $a$, there exist endomorphism $\overline{a\overline{\overline{e}}}$ of right $A$-vector space $V$ which has the matrix $aE_n$ whith respecpect to the basis $\overline{\overline{e}}$.*

Proof. Let $\overline{\overline{e}}$ be the basis of right $A$-vector space of columns. According to the theorem 11.3.8 and the definition 11.3.4, the map

$$\overline{a\overline{\overline{e}}} : e_*{}^*v \in V \to e_*{}^*(aE_n)_*{}^*v \in V$$

is endomorphism of right $A$-vector space $V$. □

**Theorem 17.2.2.** *Let $V$ be a right $A$-vector space of rows and $\overline{\overline{e}}$ be basis of right $A$-vector space $V$. Let $\dim V = n$ and $E_n$ be $n \times n$ identity matrix. For any $A$-number $a$, there exist endomorphism $\overline{a\overline{\overline{e}}}$ of right $A$-vector space $V$ which has the matrix $aE_n$ with respecpect to the basis $\overline{\overline{e}}$.*

Proof. Let $\overline{\overline{e}}$ be the basis of right $A$-vector space of columns. According to the theorem 11.3.9 and the definition 11.3.4, the map

$$\overline{a\overline{\overline{e}}} : e^*{}_*v \in V \to e^*{}_*(aE_n)^*{}_*v \in V$$

is endomorphism of right $A$-vector space $V$. □

**Definition 17.2.3.** *The endomorphism $\overline{a\overline{\overline{e}}}$ of right $A$-vector space $V$ is called* **similarity transformation** *with respect to the basis $\overline{\overline{e}}$.* □

**Theorem 17.2.4.** *Let $V$ be a right $A$-vector space of columns and $\overline{\overline{e}}_1, \overline{\overline{e}}_2$ be bases of right $A$-vector space $V$. Let $g$ be passive transformation of basis $\overline{\overline{e}}_1$ into basis $\overline{\overline{e}}_2$*

$$\overline{\overline{e}}_2 = \overline{\overline{e}}_{1*}{}^*g$$

*The similarity transformation $\overline{a\overline{\overline{e}}}_1$ has the matrix*

$$(17.2.1) \qquad g^{-1}{}_*{}^*{}_*{}^*ag$$

*with respect to the basis $\overline{\overline{e}}_2$.*

Proof. According to the theorem 12.6.5, the similarity transformation $\overline{a\overline{\overline{e}}}_1$ has the matrix

$$(17.2.3) \qquad g^{-1}{}_*{}^*{}_*{}^*aE_{n*}{}^*g$$

with respect to the basis $\overline{\overline{e}}_2$. The theorem follows from the expression (17.2.3). □

**Theorem 17.2.5.** *Let $V$ be a right $A$-vector space of rows and $\overline{\overline{e}}_1, \overline{\overline{e}}_2$ be bases of right $A$-vector space $V$. Let $g$ be passive transformation of basis $\overline{\overline{e}}_1$ into basis $\overline{\overline{e}}_2$*

$$\overline{\overline{e}}_2 = \overline{\overline{e}}_1{}^*{}_*g$$

*The similarity transformation $\overline{a\overline{\overline{e}}}_1$ has the matrix*

$$(17.2.2) \qquad g^{-1^*}{}^*{}_*ag$$

*with respect to the basis $\overline{\overline{e}}_2$.*

Proof. According to the theorem 12.6.6, the similarity transformation $\overline{a\overline{\overline{e}}}_1$ has the matrix

$$(17.2.4) \qquad g^{-1^*}{}^*{}_*aE_n{}^*{}_*g$$

with respect to the basis $\overline{\overline{e}}_2$. The theorem follows from the expression (17.2.4). □

Theorem 17.2.6. *Let $V$ be right $A$-vector space of columns and $\overline{\overline{e}}$ be basis of right $A$-vector space $V$. The set of similarity transformations $\overline{a\overline{\overline{e}}}$ of right $A$-vector space $V$ is left-side effective representation of $D$-algebra $A$ in right $A$-vector space $V$.*

Proof. Consider the set of matrices $aE_n$ of similarity transformations $\overline{a\overline{\overline{e}}}$ whith respecpect to the basis $\overline{\overline{e}}$. Let $V^*$ be the set of coordinates of $V$-numbers whith respecpect to the basis $\overline{\overline{e}}$. According to theorems 13.2.10, 13.2.11, the set $V^*$ is right $A$-vector space. According to the theorem 17.2.1, the map $\overline{a\overline{\overline{e}}}$ is the endomorphism of right $A$-vector space of columns.

| | |
|---|---|
| Lemma 17.2.7. *Let $V$ be a right $A$-vector space of columns. The map* $$(17.2.5) \quad a \in A \to aE_n \in \operatorname{End}(A, V^*)$$ *is homomorphism of $D$-algebra $A$.* | Lemma 17.2.8. *Let $V$ be a right $A$-vector space of rows. The map* $$(17.2.6) \quad a \in A \to aE_n \in \operatorname{End}(A, V^*)$$ *is homomorphism of $D$-algebra $A$.* |

Proof. According to the theorem 5.4.8, the lemma follows from equalities

$$(a+b)E_n = aE_n + bE_n$$
$$(ba)E_n = b(aE_n)$$
$$((ab)E_n) = (aE_n)_*(bE_n)$$

According to the definition 2.5.4 and the lemma 17.2.7, the map (17.2.5) is left-side representation of $D$-algebra $A$ in right $A$-vector space $V^*$. According to the theorem 5.1.31, the equality

$$aE_n = bE_n$$

implies $a = b$. According to the definition 2.3.2, the representation (17.2.5) is effective.

According to theorems 12.6.7, 12.6.9, the representation (17.2.5) is covariant with respect to choice of basis of right $A$-vector space $V$. According to the theorem 17.2.1, the set of maps $\overline{a\overline{\overline{e}}}$ is left-side representation of $D$-algebra $A$ in right $A$-vector space $V$.

Proof. According to the theorem 5.4.8, the lemma follows from equalities

$$(a+b)E_n = aE_n + bE_n$$
$$(ba)E_n = b(aE_n)$$
$$((ab)E_n) = (aE_n)^*{}_*(bE_n)$$

According to the definition 2.5.4 and the lemma 17.2.8, the map (17.2.6) is left-side representation of $D$-algebra $A$ in right $A$-vector space $V^*$. According to the theorem 5.1.31, the equality

$$aE_n = bE_n$$

implies $a = b$. According to the definition 2.3.2, the representation (17.2.6) is effective.

According to theorems 12.6.8, 12.6.10, the representation (17.2.6) is covariant with respect to choice of basis of right $A$-vector space $V$. According to the theorem 17.2.2, the set of maps $\overline{a\overline{\overline{e}}}$ is left-side representation of $D$-algebra $A$ in right $A$-vector space $V$.

$\square$

## 17.3. Twin Representations of Division $D$-Algebra

| | |
|---|---|
| Definition 17.3.1. *Let $V$ be a right $A$-vector space and $\overline{\overline{e}}$ be basis of right $A$-* | Definition 17.3.2. *Let $V$ be a left $A$-vector space and $\overline{\overline{e}}$ be basis of left $A$-* |

vector space $V$. *Bilinear map*

(17.3.1)     $(a, v) \in A \times V \to av \in V$

*generated by left-side representation*

(17.3.2)     $(a, v) \to \overline{a\overline{\overline{e}}} \circ v$

*is called* **left-side product** *of vector over scalar.*                    □

---

THEOREM 17.3.3. $\forall a, b \in A, v \in V$

(17.3.5)          $(av)b = a(vb)$

PROOF. According to the theorem 17.2.1 and to the definition 17.3.1, the equality

(17.3.6)  $a(vb) = \overline{a\overline{\overline{e}}} \circ (vb)$

$= (\overline{a\overline{\overline{e}}} \circ v)b = (av)b$

follows from the equality

$\boxed{(11.3.5) \quad f \circ (va) = (f \circ v)a}$

The equality (17.3.5) follows from the equality (17.3.6).                    □

Vectors $v$ and $av$ are called left **colinear**. The statement that vectors $v$ and $w$ are left colinear vectors does not imply that vectors $v$ and $w$ are right linearly dependent.

---

THEOREM 17.3.5. *Let $V$ be a right $A$-vector space*

$$f : A \longrightarrow V$$

*An effective left-side representation*

$$g(\overline{\overline{e}}) : A \longrightarrow V$$

*of $D$-algebra $A$ in right $A$-vector space $V$ depends on the basis $\overline{\overline{e}}$ of right $A$-vector space $V$. The following diagram*

(17.3.9)

*is commutative for any $a, b \in A$.* (17.3.5).

---

vector space $V$. *Bilinear map*

(17.3.3)     $(v, a) \in V \times A \to va \in V$

*generated by right-side representation*

(17.3.4)     $(v, a) \to \overline{\overline{\overline{e}a}} \circ v$

*is called* **right-side product** *of vector over scalar.*                    □

---

THEOREM 17.3.4. $\forall a, b \in A, v \in V$

(17.3.7)          $(av)b = a(vb)$

PROOF. According to the theorem 17.1.1 and to the definition 17.3.2, the equality

(17.3.8)  $(av)b = \overline{\overline{\overline{e}}b} \circ (av)$

$= a(\overline{\overline{\overline{e}}b} \circ v) = a(vb)$

follows from the equality

$\boxed{(10.3.5) \quad f \circ (av) = a(f \circ v)}$

The equality (17.3.7) follows from the equality (17.3.8).                    □

Vectors $v$ and $va$ are called right **colinear**. The statement that vectors $v$ and $w$ are right colinear vectors does not imply that vectors $v$ and $w$ are left linearly dependent.

---

THEOREM 17.3.6. *Let $V$ be a left $A$-vector space*

$$f : A \longrightarrow V$$

*An effective right-side representation*

$$g(\overline{\overline{e}}) : A \longrightarrow V$$

*of $D$-algebra $A$ in left $A$-vector space $V$ depends on the basis $\overline{\overline{e}}$ of left $A$-vector space $V$. The following diagram*

(17.3.10)

*is commutative for any $a, b \in A$.* (17.3.7).

Proof. Commutativity of the diagram (17.3.9) follows from the equality $\square$

Proof. Commutativity of the diagram (17.3.10) follows from the equality $\square$

Definition 17.3.7. *We call representations $f$ and $g(\overline{\overline{e}})$ considered in the theorem 17.3.5* **twin representations** *of the $D$-algebra $A$. The equality (17.3.5) represents* **associative law** *for twin representations. This allows us writing of such expression without using of brackets.* $\square$

Definition 17.3.8. *We call representations $f$ and $g(\overline{\overline{e}})$ considered in the theorem 17.3.6* **twin representations** *of the $D$-algebra $A$. The equality (17.3.7) represents* **associative law** *for twin representations. This allows us writing of such expression without using of brackets.* $\square$

## 17.4. $A \otimes A$-Module

Theorem 17.4.1. *Let $f$, $g(\overline{\overline{e}})$ be twin representations of $D$-algebra $A$ in $D$-module $V$ and $\overline{\overline{e}}$ be the basis of right $A$-vector space $V$. Let product in $D$-module $A \otimes A$ be defined according to rule*

$$(p_0 \otimes p_1) \circ (q_0 \otimes q_1) = (p_0 q_0) \otimes (q_1 p_1)$$

*Then representation*

$$h : A \otimes A \longrightarrow\!\!\!\!\!* V$$

*of $D$-algebra $A \otimes A$ in $D$-module $V$ defined by the equality*

$$h(a \otimes b) \circ v = avb$$

*is left $A \otimes A$-module.*

Theorem 17.4.2. *Let $f$, $g(\overline{\overline{e}})$ be twin representations of $D$-algebra $A$ in $D$-module $V$ and $\overline{\overline{e}}$ be the basis of left $A$-vector space $V$. Let product in $D$-module $A \otimes A$ be defined according to rule*

$$(p_0 \otimes p_1) \circ (q_0 \otimes q_1) = (p_0 q_0) \otimes (q_1 p_1)$$

*Then representation*

$$h : A \otimes A \longrightarrow\!\!\!\!\!* V$$

*of $D$-algebra $A \otimes A$ in $D$-module $V$ defined by the equality*

$$h(a \otimes b) \circ v = avb$$

*is left $A \otimes A$-module.*

Proof. According to theorems 17.2.6, 17.1.6, for given $v \in V$, the map

$$p : (a, b) \in A \times A \to avb \in V$$

is bilinear map. According to the theorem 4.6.4, for given $v \in V$, there exist map

$$h(a \otimes b) = p(a, b)$$

$$h : a \otimes b \in A \otimes A \to avb \in V$$

Since $v \in V$ is arbitrary, then for any $a$, $b$  $h(a \otimes b)$  is a map

$$h(a \otimes b) : v \in V \to avb \in V$$

According to theorems 17.2.1, 17.1.1, (as well 17.2.2, 17.1.2)

$$h(a \otimes b) \in \mathcal{L}(D; V \to V)$$

Let we define product

$$(a_0 \otimes a_1) \circ (b_0 \otimes b_1) = (a_0 b_0) \otimes (b_1 a_1)$$

in $D$-module $A \otimes A$. Then $D$-module $A \otimes A$ is $D$-algebra. From equalities

$$(a_0 \otimes a_1 + b_0 \otimes b_1) \circ v = a_0 v a_1 + b_0 v b_1 = (a_0 \otimes a_1) \circ v + (b_0 \otimes b_1) \circ v$$

$$(p a_0 \otimes a_1) \circ v = p a_0 v a_1 = (p(a_0 \otimes a_1)) \circ v$$

$$(a_0 \otimes a_1) \circ (b_0 \otimes b_1) \circ v = (a_0 \otimes a_1) \circ (b_0 v b_1) = a_0 b_0 v b_1 a_1 = (a_0 b_0 \otimes (b))$$

it follows that the map

$$h : A \otimes A \to \mathcal{L}(D, V \to V)$$

is homomorphism. Therefore, the map $h$ is representation.    $\square$

According to theorems 17.4.2, 17.4.1, there exist three different representations of $D$-algebra $A$ in $D$-module $V$. The choice of the representation depends on the problem that we consider.

For instance, the set of solutions of system of linear equations

$$5i x^1 + 6j x^2 + 3k x^3 + 7x^4 = i + k$$

$$6k x^1 + 7j x^2 + 9k x^3 + j x^4 = j + 2k$$

in quaternion algebra is right $H$-vector space. In general, left linear combination of solutions is not solution. However there are problems where we forced to reject simple model and consider both structures of vector space at the same time.

Let $V$ be left $A \otimes A$-module. If $v$, $w \in V$, then for any $a$, $b$, $c$, $d \in A$

(17.4.1)                                $avb + cwd \in V$

The expression (17.4.1) is called **linear combination of vectors** of left $A \otimes A$-module $V$.

| | |
|---|---|
| **Definition 17.4.3.** *Let $V$ be left $A \otimes A$-module of columns. The set of vectors $a_i$, $i \in I$, of left $A \otimes A$-module $V$ is* **linearly independent** *if $b^i c^i = 0$, $i = 1, ..., n$, follows from the equality*<br><br>$$b^i a_i c^i = 0$$<br><br>*Otherwise the set of vectors $a_i$ is* **linearly dependent**.    $\square$ | **Definition 17.4.4.** *Let $V$ be left $A \otimes A$-module of rows. The set of vectors $a^i$, $i \in I$, of left $A \otimes A$-module $V$ is* **linearly independent** *if $b_i c_i = 0$, $i = 1, ..., n$, follows from the equality*<br><br>$$b_i a^i c_i = 0$$<br><br>*Otherwise the set of vectors $a_i$ is* **linearly dependent**.    $\square$ |
| **Definition 17.4.5.** *We call set of vectors*<br><br>$$\overline{\overline{e}} = (e_i, i \in I)$$<br><br>**basis** *for $A \otimes A$-module $V$ if the set of vectors $e_i$ is linearly independent and adding to this set any other vector we get a new set which is linearly dependent.*    $\square$ | **Definition 17.4.6.** *We call set of vectors*<br><br>$$\overline{\overline{e}} = (e^i, i \in I)$$<br><br>**basis** *for $A \otimes A$-module $V$ if the set of vectors $e^i$ is linearly independent and adding to this set any other vector we get a new set which is linearly dependent.*    $\square$ |

## 17.5. Eigenvalue of Endomorphism

| | |
|---|---|
| **Theorem 17.5.1.** *Let $V$ be a left $A$-vector space and $\overline{\overline{e}}$ be basis of left $A$-* | **Theorem 17.5.2.** *Let $V$ be a right $A$-vector space and $\overline{\overline{e}}$ be basis of right* |

*vector space $V$. If $b \notin Z(A)$, then there is no endomorphism $\overline{f}$ such that*

$$(17.5.1) \qquad \overline{f} \circ v = bv$$

Proof. Let there exist endomorphism $\overline{f}$ such that the equality (17.5.1) is true. The equality

$$
\begin{aligned}
(17.5.2) \qquad (ba)v &= b(av) = \overline{f} \circ (av) \\
&= a(\overline{f} \circ v) = a(bv) \\
&= (ab)v
\end{aligned}
$$

follows from equalities

$$\boxed{(7.2.9) \quad (pq)v = p(qv)}$$
$$\boxed{(10.3.5) \quad f \circ (av) = a(f \circ v)}$$

and from the equality (17.5.1). According to the theorem 12.1.4, the equality

$$(17.5.3) \qquad ab = ba$$

for any $a \in A$ follows from the equality (17.5.2). The equality (17.5.3) for any $a \in A$ contradicts to the statement $b \notin Z(A)$. Therefore, considered endomorphism $\overline{f}$ does not exist. $\qquad \square$

According to the theorem 17.5.1, as opposed to the right-side product of vector over scalar, left-side product of vector over scalar is not endomorphism of left $A$-vector space. Comparing the definition 17.1.3 and the theorem 17.5.1, we see that similarity transformation in non-commutative algebra and similarity transformation in commutative algebra play the same role.

---

DEFINITION 17.5.3. *Let $V$ be a left $A$-vector space and $\overline{\overline{e}}$ be basis of left $A$-vector space $V$. The vector $v \in V$ is called* **eigenvector** *of the endomorphism*

$$\overline{f} : V \to V$$

*with respect to the basis $\overline{\overline{e}}$, if there exists $b \in A$ such that*

$$(17.5.7) \qquad \overline{f} \circ v = \overline{\overline{e}} b \circ v$$

*$A$-number $b$ is called* **eigenvalue** *of the endomorphism $f$ with respect to the basis $\overline{\overline{e}}$.* $\qquad \square$

---

*$A$-vector space $V$. If $b \notin Z(A)$, then there is no endomorphism $\overline{f}$ such that*

$$(17.5.4) \qquad \overline{f} \circ v = vb$$

Proof. Let there exist endomorphism $\overline{f}$ such that the equality (17.5.4) is true. The equality

$$
\begin{aligned}
(17.5.5) \qquad v(ab) &= (va)b = \overline{f} \circ (va) \\
&= (\overline{f} \circ v)a = (vb)a \\
&= v(ba)
\end{aligned}
$$

follows from equalities

$$\boxed{(7.4.9) \quad v(pq) = (vp)q}$$
$$\boxed{(11.3.5) \quad f \circ (va) = (f \circ v)a}$$

and from the equality (17.5.4). According to the theorem 12.4.4, the equality

$$(17.5.6) \qquad ab = ba$$

for any $a \in A$ follows from the equality (17.5.5). The equality (17.5.6) for any $a \in A$ contradicts to the statement $b \notin Z(A)$. Therefore, considered endomorphism $\overline{f}$ does not exist. $\qquad \square$

According to the theorem 17.5.2, as opposed to the left-side product of vector over scalar, right-side product of vector over scalar is not endomorphism of right $A$-vector space. Comparing the definition 17.2.3 and the theorem 17.5.2, we see that similarity transformation in non-commutative algebra and similarity transformation in commutative algebra play the same role.

---

DEFINITION 17.5.4. *Let $V$ be a right $A$-vector space and $\overline{\overline{e}}$ be basis of right $A$-vector space $V$. The vector $v \in V$ is called* **eigenvector** *of the endomorphism*

$$\overline{f} : V \to V$$

*with respect to the basis $\overline{\overline{e}}$, if there exists $b \in A$ such that*

$$(17.5.8) \qquad \overline{f} \circ v = b\overline{\overline{e}} \circ v$$

*$A$-number $b$ is called* **eigenvalue** *of the endomorphism $f$ with respect to the basis $\overline{\overline{e}}$.* $\qquad \square$

THEOREM 17.5.5. *Let $V$ be a left $A$-vector space of columns and $\overline{\overline{e}}_1, \overline{\overline{e}}_2$ be bases of left $A$-vector space $V$. Let $g$ be passive transformation of basis $\overline{\overline{e}}_1$ into basis $\overline{\overline{e}}_2$*

$$\overline{\overline{e}}_2 = g^*{}_*\overline{\overline{e}}_1$$

*Let $f$ be the matrix of endomorphism $\overline{f}$ and $v$ be the matrix of vector $\overline{v}$ with respect to the basis $\overline{\overline{e}}_2$. The vector $\overline{v}$ is eigenvector of the endomorphism $\overline{f}$ with respect to the basis $\overline{\overline{e}}_1$ iff there exists $A$-number $b$ such that the system of linear equations*

$$(17.5.9) \qquad v^*{}_*(f - gb^*{}_*g^{-1^*{}^*}) = 0$$

*has non-trivial solution.*

PROOF. According to theorems 10.3.6, 17.1.4, the equality

$$(17.5.10) \qquad v^*{}_*f^*{}_*e = \overline{f} \circ \overline{v} = b\overline{\overline{e}} \circ \overline{v}$$
$$= v^*{}_*gb^*{}_*g^{-1^*{}^*}{}_*e$$

follows from the equality (17.5.7). The equality (17.5.9) follows from the equality (17.5.10) and from the theorem 8.4.4. Since the case $\overline{v} = 0$ is not interesting for us, the system of linear equations (17.5.9) should have non-trivial solution. □

THEOREM 17.5.6. *Let $V$ be a left $A$-vector space of columns and $\overline{\overline{e}}_1, \overline{\overline{e}}_2$ be bases of left $A$-vector space $V$. Let $g$ be passive transformation of basis $\overline{\overline{e}}_1$ into basis $\overline{\overline{e}}_2$*

$$\overline{\overline{e}}_2 = g^*{}_*\overline{\overline{e}}_1$$

*Let $f$ be the matrix of endomorphism $\overline{f}$ with respect to the basis $\overline{\overline{e}}_2$. $A$-number $b$ is eigenvalue of endomorphism $\overline{f}$ iff the matrix*

$$(17.5.11) \qquad f - gb^*{}_*g^{-1^*}$$

*is $^*{}_*$-singular.*

PROOF. The theorem follows from theorems 9.4.3, 17.5.5. □

THEOREM 17.5.7. *Let $V$ be a right $A$-vector space of columns and $\overline{\overline{e}}_1, \overline{\overline{e}}_2$ be bases of right $A$-vector space $V$. Let $g$ be passive transformation of basis $\overline{\overline{e}}_1$ into basis $\overline{\overline{e}}_2$*

$$\overline{\overline{e}}_2 = \overline{\overline{e}}_{1*}{}^*g$$

*Let $f$ be the matrix of endomorphism $\overline{f}$ and $v$ be the matrix of vector $\overline{v}$ with respect to the basis $\overline{\overline{e}}_2$. The vector $\overline{v}$ is eigenvector of the endomorphism $\overline{f}$ with respect to the basis $\overline{\overline{e}}_1$ iff there exists $A$-number $b$ such that the system of linear equations*

$$(17.5.12) \qquad (f - g^{-1^*{}^*}{}_*bg)_*{}^*v = 0$$

*has non-trivial solution.*

PROOF. According to theorems 11.3.6, 17.2.4, the equality

$$(17.5.13) \qquad e_*{}^*f_*{}^*v = \overline{f} \circ \overline{v} = b\overline{\overline{e}} \circ \overline{v}$$
$$= e_*{}^*g^{-1^*{}^*}{}_*bg_*{}^*v$$

follows from the equality (17.5.8). The equality (17.5.12) follows from the equality (17.5.13) and from the theorem 8.4.18. Since the case $\overline{v} = 0$ is not interesting for us, the system of linear equations (17.5.12) should have non-trivial solution. □

THEOREM 17.5.8. *Let $V$ be a right $A$-vector space of columns and $\overline{\overline{e}}_1, \overline{\overline{e}}_2$ be bases of right $A$-vector space $V$. Let $g$ be passive transformation of basis $\overline{\overline{e}}_1$ into basis $\overline{\overline{e}}_2$*

$$\overline{\overline{e}}_2 = \overline{\overline{e}}_{1*}{}^*g$$

*Let $f$ be the matrix of endomorphism $\overline{f}$ with respect to the basis $\overline{\overline{e}}_2$. $A$-number $b$ is eigenvalue of endomorphism $\overline{f}$ iff the matrix*

$$(17.5.14) \qquad f - g^{-1^*{}^*}{}_*{}^*bg$$

*is $_*{}^*$-singular.*

PROOF. The theorem follows from theorems 9.4.3, 17.5.7. □

THEOREM 17.5.9. *Let $V$ be a left $A$-vector space of rows and $\overline{\overline{e}}_1, \overline{\overline{e}}_2$ be bases of left $A$-vector space $V$. Let $g$ be passive transformation of basis $\overline{\overline{e}}_1$ into basis $\overline{\overline{e}}_2$*

$$\overline{\overline{e}}_2 = g_{*}{}^{*}\overline{\overline{e}}_1$$

*Let $f$ be the matrix of endomorphism $\overline{f}$ and $v$ be the matrix of vector $\overline{v}$ with respect to the basis $\overline{\overline{e}}_2$. The vector $\overline{v}$ is eigenvector of the endomorphism $\overline{f}$ with respect to the basis $\overline{\overline{e}}_1$ iff there exists $A$-number $b$ such that the system of linear equations*

$$(17.5.15) \qquad v_{*}{}^{*}(f - gb_{*}{}^{*}g^{-1*}{}^{*}) = 0$$

*has non-trivial solution.*

PROOF. According to theorems 10.3.7, 17.1.5, the equality

$$(17.5.16) \qquad \begin{aligned} v_{*}{}^{*}f_{*}{}^{*}e &= \overline{f} \circ \overline{v} = \overline{b\overline{e}} \circ \overline{v} \\ &= v_{*}{}^{*}gb_{*}{}^{*}g^{-1*}{}_{*}{}^{*}e \end{aligned}$$

follows from the equality (17.5.7). The equality (17.5.15) follows from the equality (17.5.16) and from the theorem 8.4.11. Since the case $\overline{v} = 0$ is not interesting for us, the system of linear equations (17.5.15) should have non-trivial solution. □

THEOREM 17.5.10. *Let $V$ be a left $A$-vector space of rows and $\overline{\overline{e}}_1, \overline{\overline{e}}_2$ be bases of left $A$-vector space $V$. Let $g$ be passive transformation of basis $\overline{\overline{e}}_1$ into basis $\overline{\overline{e}}_2$*

$$\overline{\overline{e}}_2 = g_{*}{}^{*}\overline{\overline{e}}_1$$

*Let $f$ be the matrix of endomorphism $\overline{f}$ with respect to the basis $\overline{\overline{e}}_2$. $A$-number $b$ is eigenvalue of endomorphism $\overline{f}$ iff the matrix*

$$(17.5.17) \qquad f - gb_{*}{}^{*}g^{-1*}{}^{*}$$

*is $_{*}{}^{*}$-singular.*

PROOF. The theorem follows from theorems 9.4.3, 17.5.9. □

THEOREM 17.5.11. *Let $V$ be a right $A$-vector space of rows and $\overline{\overline{e}}_1, \overline{\overline{e}}_2$ be bases of right $A$-vector space $V$. Let $g$ be passive transformation of basis $\overline{\overline{e}}_1$ into basis $\overline{\overline{e}}_2$*

$$\overline{\overline{e}}_2 = \overline{\overline{e}}_1{}^{*}{}_{*}g$$

*Let $f$ be the matrix of endomorphism $\overline{f}$ and $v$ be the matrix of vector $\overline{v}$ with respect to the basis $\overline{\overline{e}}_2$. The vector $\overline{v}$ is eigenvector of the endomorphism $\overline{f}$ with respect to the basis $\overline{\overline{e}}_1$ iff there exists $A$-number $b$ such that the system of linear equations*

$$(17.5.18) \qquad (f - g^{-1*}{}^{*}{}_{*}bg)^{*}{}_{*}v = 0$$

*has non-trivial solution.*

PROOF. According to theorems 11.3.7, 17.2.5, the equality

$$(17.5.19) \qquad \begin{aligned} e^{*}{}_{*}f^{*}{}_{*}v &= \overline{f} \circ \overline{v} = \overline{\overline{e}b} \circ \overline{v} \\ &= e^{*}{}_{*}g^{-1*}{}^{*}{}_{*}bg^{*}{}_{*}v \end{aligned}$$

follows from the equality (17.5.8). The equality (17.5.18) follows from the equality (17.5.19) and from the theorem 8.4.25. Since the case $\overline{v} = 0$ is not interesting for us, the system of linear equations (17.5.18) should have non-trivial solution. □

THEOREM 17.5.12. *Let $V$ be a right $A$-vector space of rows and $\overline{\overline{e}}_1, \overline{\overline{e}}_2$ be bases of right $A$-vector space $V$. Let $g$ be passive transformation of basis $\overline{\overline{e}}_1$ into basis $\overline{\overline{e}}_2$*

$$\overline{\overline{e}}_2 = \overline{\overline{e}}_1{}^{*}{}_{*}g$$

*Let $f$ be the matrix of endomorphism $\overline{f}$ with respect to the basis $\overline{\overline{e}}_2$. $A$-number $b$ is eigenvalue of endomorphism $\overline{f}$ iff the matrix*

$$(17.5.20) \qquad f - g^{-1*}{}^{*}{}_{*}bg$$

*is $^{*}{}_{*}$-singular.*

PROOF. The theorem follows from theorems 9.4.3, 17.5.11. □

THEOREM 17.5.13. *Let $V$ be a left $A$-vector space of columns and $\overline{\overline{e}}_1$, $\overline{\overline{e}}_2$ be bases of left $A$-vector space $V$.*

STATEMENT 17.5.14. *Eigenvalue b of the endomorphism $\overline{f}$ with respect to the basis $\overline{\overline{e}}_1$ does not depend on the choice of the basis $\overline{\overline{e}}_2$.* ⊙

STATEMENT 17.5.15. *Eigenvector $\overline{v}$ of the endomorphism $\overline{f}$ corresponding to eigenvalue b does not depend on the choice of the basis $\overline{\overline{e}}_2$.* ⊙

PROOF. Let $f_i$, $i = 1, 2$, be the matrix of endomorphism $f$ with respect to the basis $\overline{\overline{e}}_i$. Let $b$ be eigenvalue of the endomorphism $\overline{f}$ with respect to the basis $\overline{\overline{e}}_1$ and $\overline{v}$ be eigenvector of the endomorphism $\overline{f}$ corresponding to eigenvalue $b$. Let $v_i$, $i = 1, 2$, be the matrix of vector $\overline{v}$ with respect to the basis $\overline{\overline{e}}_i$.

Suppose first that $\overline{\overline{e}}_1 = \overline{\overline{e}}_2$. Then $E_n$ is passive transformation of basis $\overline{\overline{e}}_1$ into basis $\overline{\overline{e}}_2$

$$\overline{\overline{e}}_2 = E_n{}^*{}_*\overline{\overline{e}}_1$$

According to the theorem 17.1.4, the similarity transformation $\overline{\overline{e}}_1 b$ has the matrix

$$E_n b^*{}_* E_n^{-1*}{}^* = E_n b$$

with respect to the basis $\overline{\overline{e}}_2$. Therefore, according to the theorem 17.5.6, we proved following statements.

LEMMA 17.5.16. *The matrix $f_1 - E_n b$ is $_*{}^*$-singular matrix.* ⊙

LEMMA 17.5.17. *The column vector $v_1$ satisfies to the system of linear equations*

$$(17.5.21) \qquad v_1{}^*{}_*(f_1 - E_n b) = 0$$

⊙

Suppose that $\overline{\overline{e}}_1 \neq \overline{\overline{e}}_2$. According to the theorem 12.2.13, there exists

THEOREM 17.5.20. *Let $V$ be a right $A$-vector space of columns and $\overline{\overline{e}}_1$, $\overline{\overline{e}}_2$ be bases of right $A$-vector space $V$.*

STATEMENT 17.5.21. *Eigenvalue b of the endomorphism $\overline{f}$ with respect to the basis $\overline{\overline{e}}_1$ does not depend on the choice of the basis $\overline{\overline{e}}_2$.* ⊙

STATEMENT 17.5.22. *Eigenvector $\overline{v}$ of the endomorphism $\overline{f}$ corresponding to eigenvalue b does not depend on the choice of the basis $\overline{\overline{e}}_2$.* ⊙

PROOF. Let $f_i$, $i = 1, 2$, be the matrix of endomorphism $f$ with respect to the basis $\overline{\overline{e}}_i$. Let $b$ be eigenvalue of the endomorphism $\overline{f}$ with respect to the basis $\overline{\overline{e}}_1$ and $\overline{v}$ be eigenvector of the endomorphism $\overline{f}$ corresponding to eigenvalue $b$. Let $v_i$, $i = 1, 2$, be the matrix of vector $\overline{v}$ with respect to the basis $\overline{\overline{e}}_i$.

Suppose first that $\overline{\overline{e}}_1 = \overline{\overline{e}}_2$. Then $E_n$ is passive transformation of basis $\overline{\overline{e}}_1$ into basis $\overline{\overline{e}}_2$

$$\overline{\overline{e}}_2 = \overline{\overline{e}}_1{}_*{}^* E_n$$

According to the theorem 17.2.4, the similarity transformation $b\overline{\overline{e}}_1$ has the matrix

$$E_n^{-1*}{}_*{}^* b E_n = b E_n$$

with respect to the basis $\overline{\overline{e}}_2$. Therefore, according to the theorem 17.5.8, we proved following statements.

LEMMA 17.5.23. *The matrix $f_1 - b E_n$ is $_*{}^*$-singular matrix.* ⊙

LEMMA 17.5.24. *The column vector $v_1$ satisfies to the system of linear equations*

$$(17.5.25) \qquad (f_1 - b E_n)_*{}^* v_1 = 0$$

⊙

Suppose that $\overline{\overline{e}}_1 \neq \overline{\overline{e}}_2$. According to the theorem 12.5.13, there exists unique passive transformation

$$\overline{\overline{e}}_2 = g^*{}_* \overline{\overline{e}}_1$$

and $g$ is $^*{}_*$-non-singular matrix. According to the theorem 12.3.5,

(17.5.22)        $f_2 = g^*{}_* f_1{}^*{}_* g^{-1^*}{}_*$

According to the theorem 17.5.6, if $b$ is eigenvalue of endomorphism $\overline{f}$, then we have to prove the lemma 17.5.18.

LEMMA 17.5.18. *The matrix*

$$f_2 - gb^*{}_* g^{-1^*}{}_*$$

*is* $^*{}_*$-*singular matrix.*        ⊙

Since the equality

$$f_2 - gb^*{}_* g^{-1^*}{}_*$$
$$= g^*{}_* f_1{}^*{}_* g^{-1^*}{}_* - gb^*{}_* g^{-1^*}{}_*$$
$$= g^*{}_*(f_1 - E_n b)^*{}_* g^{-1^*}{}_*$$

follows from the equality (17.5.22), then the lemma 17.5.18 follows from the lemma 17.5.16. Therefore, we proved the statement 17.5.14.

According to the theorem 17.5.5, if $\overline{v}$ is eigenvector of endomorphism $\overline{f}$, then we have to prove the lemma 17.5.19.

LEMMA 17.5.19. *The column vector* $v_2$ *satisfies to the system of linear equations*

(17.5.23)    $v_2{}^*{}_*(f_2 - gb^*{}_* g^{-1^*}{}_*) = 0$

PROOF. According to the theorem 12.3.1,

(17.5.24)        $v_2 = v_1{}^*{}_* g^{-1^*}{}_*$

Since the equality

$$v_2{}^*{}_*(f_2 - gb^*{}_* g^{-1^*}{}_*)$$
$$= v_2{}^*{}_*(g^*{}_* f_1{}^*{}_* g^{-1^*}{}_* - gb^*{}_* g^{-1^*}{}_*)$$
$$= v_1{}^*{}_* g^{-1^*}{}_*{}^* g^*{}_*(f_1 - E_n b)^*{}_* g^{-1^*}{}_*$$

follows from equalities (17.5.22), (17.5.24), then the lemma 17.5.19 follows from the lemma 17.5.17.        ⊙

Therefore, we proved the statement 17.5.15.        □

---

unique passive transformation

$$\overline{\overline{e}}_2 = \overline{\overline{e}}_{1*}{}^* g$$

and $g$ is $_*{}^*$-non-singular matrix. According to the theorem 12.6.5,

(17.5.26)        $f_2 = g^{-1*}{}^*{}_*{}^* f_{1*}{}^* g$

According to the theorem 17.5.8, if $b$ is eigenvalue of endomorphism $\overline{f}$, then we have to prove the lemma 17.5.25.

LEMMA 17.5.25. *The matrix*

$$f_2 - g^{-1*}{}^*{}_*{}^* bg$$

*is* $_*{}^*$-*singular matrix.*        ⊙

Since the equality

$$f_2 - g^{-1*}{}^*{}_*{}^* bg$$
$$= g^{-1*}{}^*{}_*{}^* f_{1*}{}^* g - g^{-1*}{}^*{}_*{}^* bg$$
$$= g^{-1*}{}^*{}_*{}^*(f_1 - bE_n)_*{}^* g$$

follows from the equality (17.5.26), then the lemma 17.5.25 follows from the lemma 17.5.23. Therefore, we proved the statement 17.5.21.

According to the theorem 17.5.7, if $\overline{v}$ is eigenvector of endomorphism $\overline{f}$, then we have to prove the lemma 17.5.26.

LEMMA 17.5.26. *The column vector* $v_2$ *satisfies to the system of linear equations*

(17.5.27)    $(f_2 - g^{-1*}{}^*{}_*{}^* bg)_*{}^* v_2 = 0$

PROOF. According to the theorem 12.6.1,

(17.5.28)        $v_2 = g^{-1*}{}^*{}_*{}^* v_1$

Since the equality

$$(f_2 - g^{-1*}{}^*{}_*{}^* bg)_*{}^* v_2$$
$$= (g^{-1*}{}^*{}_*{}^* f_{1*}{}^* g - g^{-1*}{}^*{}_*{}^* bg)_*{}^* v_2$$
$$= g^{-1*}{}^*{}_*{}^*(f_1 - bE_n)_*{}^* g_*{}^* g^{-1*}{}^*{}_*{}^* v_1$$

follows from equalities (17.5.26), (17.5.28), then the lemma 17.5.26 follows from the lemma 17.5.24.        ⊙

Therefore, we proved the statement 17.5.22.        □

**Theorem 17.5.27.** *Let $V$ be a left $A$-vector space of rows and $\overline{\overline{e}}_1, \overline{\overline{e}}_2$ be bases of left $A$-vector space $V$.*

**Statement 17.5.28.** *Eigenvalue $b$ of the endomorphism $\overline{f}$ with respect to the basis $\overline{\overline{e}}_1$ does not depend on the choice of the basis $\overline{\overline{e}}_2$.* ⊙

**Statement 17.5.29.** *Eigenvector $\overline{v}$ of the endomorphism $\overline{f}$ corresponding to eigenvalue $b$ does not depend on the choice of the basis $\overline{\overline{e}}_2$.* ⊙

**Proof.** Let $f_i$, $i = 1, 2$, be the matrix of endomorphism $f$ with respect to the basis $\overline{\overline{e}}_i$. Let $b$ be eigenvalue of the endomorphism $\overline{f}$ with respect to the basis $\overline{\overline{e}}_1$ and $\overline{v}$ be eigenvector of the endomorphism $\overline{f}$ corresponding to eigenvalue $b$. Let $v_i$, $i = 1, 2$, be the matrix of vector $\overline{v}$ with respect to the basis $\overline{\overline{e}}_i$.

Suppose first that $\overline{\overline{e}}_1 = \overline{\overline{e}}_2$. Then $E_n$ is passive transformation of basis $\overline{\overline{e}}_1$ into basis $\overline{\overline{e}}_2$

$$\overline{\overline{e}}_2 = E_{n\,*}{}^*\overline{\overline{e}}_1$$

According to the theorem 17.1.5, the similarity transformation $\overline{\overline{\overline{e}}_1 b}$ has the matrix

$$E_n b_*{}^* E_n^{-1\,*}{}^* = E_n b$$

with respect to the basis $\overline{\overline{e}}_2$. Therefore, according to the theorem 17.5.10, we proved following statements.

**Lemma 17.5.33.** *The matrix $f_1 - E_n b$ is $_*{}^*$-singular matrix.* ⊙

**Lemma 17.5.34.** *The row vector $v_1$ satisfies to the system of linear equations*

$$(17.5.29) \qquad v_{1\,*}{}^*(f_1 - E_n b) = 0$$
⊙

Suppose that $\overline{\overline{e}}_1 \neq \overline{\overline{e}}_2$. According to the theorem 12.2.14, there exists unique passive transformation

$$\overline{\overline{e}}_2 = g_*{}^*\overline{\overline{e}}_1$$

**Theorem 17.5.30.** *Let $V$ be a right $A$-vector space of rows and $\overline{\overline{e}}_1, \overline{\overline{e}}_2$ be bases of right $A$-vector space $V$.*

**Statement 17.5.31.** *Eigenvalue $b$ of the endomorphism $\overline{f}$ with respect to the basis $\overline{\overline{e}}_1$ does not depend on the choice of the basis $\overline{\overline{e}}_2$.* ⊙

**Statement 17.5.32.** *Eigenvector $\overline{v}$ of the endomorphism $\overline{f}$ corresponding to eigenvalue $b$ does not depend on the choice of the basis $\overline{\overline{e}}_2$.* ⊙

**Proof.** Let $f_i$, $i = 1, 2$, be the matrix of endomorphism $f$ with respect to the basis $\overline{\overline{e}}_i$. Let $b$ be eigenvalue of the endomorphism $\overline{f}$ with respect to the basis $\overline{\overline{e}}_1$ and $\overline{v}$ be eigenvector of the endomorphism $\overline{f}$ corresponding to eigenvalue $b$. Let $v_i$, $i = 1, 2$, be the matrix of vector $\overline{v}$ with respect to the basis $\overline{\overline{e}}_i$.

Suppose first that $\overline{\overline{e}}_1 = \overline{\overline{e}}_2$. Then $E_n$ is passive transformation of basis $\overline{\overline{e}}_1$ into basis $\overline{\overline{e}}_2$

$$\overline{\overline{e}}_2 = \overline{\overline{e}}_1{}^*{}_* E_n$$

According to the theorem 17.2.5, the similarity transformation $\overline{b\overline{\overline{e}}_1}$ has the matrix

$$E_n^{-1\,*}{}^*{}_* b E_n = b E_n$$

with respect to the basis $\overline{\overline{e}}_2$. Therefore, according to the theorem 17.5.12, we proved following statements.

**Lemma 17.5.37.** *The matrix $f_1 - b E_n$ is $^*{}_*$-singular matrix.* ⊙

**Lemma 17.5.38.** *The row vector $v_1$ satisfies to the system of linear equations*

$$(17.5.33) \qquad (f_1 - b E_n)^*{}_* v_1 = 0$$
⊙

Suppose that $\overline{\overline{e}}_1 \neq \overline{\overline{e}}_2$. According to the theorem 12.5.14, there exists unique passive transformation

$$\overline{\overline{e}}_2 = \overline{\overline{e}}_1{}^*{}_* g$$

and $g$ is $_*{}^*$-non-singular matrix. According to the theorem 12.3.6,

(17.5.30)                $f_2 = g_*{}^* f_{1*}{}^* g^{-1}{}_*{}^*$

According to the theorem 17.5.10, if $b$ is eigenvalue of endomorphism $\overline{f}$, then we have to prove the lemma 17.5.35.

**Lemma 17.5.35.** *The matrix*
$$f_2 - g b_*{}^* g^{-1}{}_*{}^*$$
*is $_*{}^*$-singular matrix.*                                  ⊙

Since the equality
$$f_2 - g b_*{}^* g^{-1}{}_*{}^*$$
$$= g_*{}^* f_{1*}{}^* g^{-1}{}_*{}^* - g b_*{}^* g^{-1}{}_*{}^*$$
$$= g_*{}^* (f_1 - E_n b)_*{}^* g^{-1}{}_*{}^*$$
follows from the equality (17.5.30), then the lemma 17.5.35 follows from the lemma 17.5.33. Therefore, we proved the statement 17.5.28.

According to the theorem 17.5.9, if $\overline{v}$ is eigenvector of endomorphism $\overline{f}$, then we have to prove the lemma 17.5.36.

**Lemma 17.5.36.** *The row vector $v_2$ satisfies to the system of linear equations*

(17.5.31)      $v_{2*}{}^* (f_2 - g b_*{}^* g^{-1}{}_*{}^*) = 0$

**Proof.** According to the theorem 12.3.2,

(17.5.32)                $v_2 = v_{1*}{}^* g^{-1}{}_*{}^*$

Since the equality
$$v_{2*}{}^* (f_2 - g b_*{}^* g^{-1}{}_*{}^*)$$
$$= v_{2*}{}^* (g_*{}^* f_{1*}{}^* g^{-1}{}_*{}^* - g b_*{}^* g^{-1}{}_*{}^*)$$
$$= v_{1*}{}^* g^{-1}{}_*{}^* {}_*{}^* g_*{}^* (f_1 - E_n b)_*{}^* g^{-1}{}_*{}^*$$
follows from equalities (17.5.30), (17.5.32), then the lemma 17.5.36 follows from the lemma 17.5.34.          ⊙

Therefore, we proved the statement 17.5.29.                                  □

and $g$ is $^*{}_*$-non-singular matrix. According to the theorem 12.6.6,

(17.5.34)                $f_2 = g^{-1}{}^*{}_* {}_* f_{1}{}^*{}_* g$

According to the theorem 17.5.12, if $b$ is eigenvalue of endomorphism $\overline{f}$, then we have to prove the lemma 17.5.39.

**Lemma 17.5.39.** *The matrix*
$$f_2 - g^{-1}{}^*{}_* {}_* bg$$
*is $^*{}_*$-singular matrix.*                                  ⊙

Since the equality
$$f_2 - g^{-1}{}^*{}_* {}_* bg$$
$$= g^{-1}{}^*{}_* {}_* f_1{}^*{}_* g - g^{-1}{}^*{}_* {}_* bg$$
$$= g^{-1}{}^*{}_* {}_* (f_1 - b E_n)^*{}_* g$$
follows from the equality (17.5.34), then the lemma 17.5.39 follows from the lemma 17.5.37. Therefore, we proved the statement 17.5.31.

According to the theorem 17.5.11, if $\overline{v}$ is eigenvector of endomorphism $\overline{f}$, then we have to prove the lemma 17.5.40.

**Lemma 17.5.40.** *The row vector $v_2$ satisfies to the system of linear equations*

(17.5.35)      $(f_2 - g^{-1}{}^*{}_* {}_* bg)^*{}_* v_2 = 0$

**Proof.** According to the theorem 12.6.2,

(17.5.36)                $v_2 = g^{-1}{}^*{}_* {}_* v_1$

Since the equality
$$(f_2 - g^{-1}{}^*{}_* {}_* bg)^*{}_* v_2$$
$$= (g^{-1}{}^*{}_* {}_* f_1{}^*{}_* g - g^{-1}{}^*{}_* {}_* bg)^*{}_* v_2$$
$$= g^{-1}{}^*{}_* {}_* (f_1 - b E_n)^*{}_* g^*{}_* g^{-1}{}^*{}_* {}_* v_1$$
follows from equalities (17.5.34), (17.5.36), then the lemma 17.5.40 follows from the lemma 17.5.38.          ⊙

Therefore, we proved the statement 17.5.32.                                  □

## 17.6. Eigenvalue of Matrix

The definition 17.5.3 does not depend on structure of coordinates of left

The definition 17.5.4 does not depend on structure of coordinates of right

$A$-vector space $V$. As soon as we choose the basis $\overline{\overline{e}}$ of left $A$-vector space $V$, we see additional structure which is reflected in theorems 17.5.5, 17.5.9.

$A$-vector space $V$. As soon as we choose the basis $\overline{\overline{e}}$ of right $A$-vector space $V$, we see additional structure which is reflected in theorems 17.5.7, 17.5.11.

The same matrix may correspond to different endomorphisms depending on whether we consider left or right $A$-vector space. Therefore, we must consider all definitions of eigenvectors which may correspond to given matrix. From theorems 17.5.5, 17.5.7, 17.5.9, 17.5.11, it follows that for given matrix there exist four types of eigenvectors and four types of eigenvalues. Couple eigenvector and eigenvalue of a certain type is used in solving the problem where corresponding type of $A$-vector space arises. This statement is reflected in the names of eigenvector and eigenvalue.

Eigenvector of endomorphism depends on choice of two bases. However when we consider eigenvector of matrix, it is not the choice of bases that is important for us, but the matrix of passive transformation between these bases.

**DEFINITION 17.6.1.** *A-number $b$ is called* $^*_*$**-eigenvalue** *of the pair of matrices $(f, g)$ where $g$ is the* $^*_*$*-non-singular matrix if the matrix*

$$f - gb^*{}_*g^{-1^*}{}_*$$

*is* $^*_*$*-singular matrix.* □

**DEFINITION 17.6.2.** *A-number $b$ is called* $_*{}^*$**-eigenvalue** *of the pair of matrices $(f, g)$ where $g$ is the* $_*{}^*$*-non-singular matrix if the matrix*

$$f - g^{-1}{}_*{}^*{}_*bg$$

*is* $_*{}^*$*-singular matrix.* □

**THEOREM 17.6.3.** *Let $f$ be $n \times n$ matrix of A-numbers. Let $g$ be non-$^*_*$-singular $n \times n$ matrix of A-numbers. Then non-zero $^*_*$-eigenvalues of the pair of matrices $(f, g)$ are roots of any equation*

$$\det({}^*_*)^i_j \ (f - gb^*{}_*g^{-1^*}{}_*) = 0$$

**THEOREM 17.6.4.** *Let $f$ be $n \times n$ matrix of A-numbers. Let $g$ be non-$_*{}^*$-singular $n \times n$ matrix of A-numbers. Then non-zero $_*{}^*$-eigenvalues of the pair of matrices $(f, g)$ are roots of any equation*

$$\det({}_*{}^*)^i_j \ (f - g^{-1}{}_*{}^*{}_*bg) = 0$$

**PROOF.** The theorem follows from the theorem 9.3.4 and from the definition 17.6.1. □

**PROOF.** The theorem follows from the theorem 9.3.3 and from the definition 17.6.2. □

**DEFINITION 17.6.5.** *Let A-number $b$ be* $^*_*$*-eigenvalue of the pair of matrices $(f, g)$ where $g$ is the* $^*_*$*-non-singular matrix. A column $v$ is called* **eigencolumn** *of the pair of matrices $(f, g)$ corresponding to* $^*_*$*-eigenvalue $b$, if the following equality is true*

$$v^*{}_*f = v^*{}_*gb^*{}_*g^{-1^*}{}_*$$

□

**DEFINITION 17.6.7.** *Let A-number $b$ be* $_*{}^*$*-eigenvalue of the pair of matrices $(f, g)$ where $g$ is the* $_*{}^*$*-non-singular matrix. A column $v$ is called* **eigencolumn** *of the pair of matrices $(f, g)$ corresponding to* $_*{}^*$*-eigenvalue $b$, if the following equality is true*

$$f_*{}^*v = g^{-1}{}_*{}^*{}_*bg_*{}^*v$$

□

The definition 17.6.5 follows from the definition 17.6.1 and from the theorem

The definition 17.6.7 follows from the definition 17.6.2 and from the theorem

17.5.5.

THEOREM 17.6.6. *Let $A$-number $b$ be $^*_*$-eigenvalue of the pair of matrices $(f, g)$ where $g$ is the $^*_*$-non-singular matrix. The set of eigencolumns of the pair of matrices $(f, g)$ corresponding to $^*_*$-eigenvalue $b$ is left $A$-vector space of columns.*

PROOF. According to the definitions 17.6.1, 17.6.5, coordinates of eigencolumn are solution of the system of linear equations

$$(17.6.1) \qquad v^*{}_*(f - gb^*{}_*g^{-1^*}{}_*) = 0$$

with $^*_*$-singular matrix

$$f - gb^*{}_*g^{-1^*}{}_*$$

According to the theorem 9.4.6, the set of solutions of the system of linear equations (17.6.1) is left space of columns.  □

DEFINITION 17.6.9. *Let $A$-number $b$ be $^*_*$-eigenvalue of the pair of matrices $(f, g^{-1^*}{}_*)$ where $g$ is the $^*_*$-non-singular matrix. A row $v$ is called* **eigenrow** *of the pair of matrices $(f, g)$ corresponding to $^*_*$-eigenvalue $b$, if the following equality is true*

$$f^*{}_*v = g^{-1^*}{}_*{}^*{}_*bg^*{}_*v$$

□

Since

$$g^{-1^*}{}_*{}^*{}_*bg = (g^{-1^*}{}_*b)^*{}_*g$$
$$= (g^{-1^*}{}_*b)^*{}_*(g^{-1^*}{}_*)^{-1^*}{}_*$$

then the definition 17.6.9 follows from the definition 17.6.1 and from the theorem 17.5.11.

THEOREM 17.6.10. *Let $A$-number $b$ be $^*_*$-eigenvalue of the pair of matrices $(f, g^{-1^*}{}_*)$ where $g$ is the $^*_*$-non-singular matrix. The set of eigenrows of the pair of matrices $(f, g)$ corresponding to $^*_*$-eigenvalue $b$ is right $A$-vector space of rows.*

THEOREM 17.6.8. *Let $A$-number $b$ be $_*^*$-eigenvalue of the pair of matrices $(f, g)$ where $g$ is the $_*^*$-non-singular matrix. The set of eigencolumns of the pair of matrices $(f, g)$ corresponding to $_*^*$-eigenvalue $b$ is right $A$-vector space of columns.*

PROOF. According to the definitions 17.6.2, 17.6.7, coordinates of eigencolumn are solution of the system of linear equations

$$(17.6.2) \qquad (f - g^{-1_*^*}{}_*^*bg)_*^*v = 0$$

with $_*^*$-singular matrix

$$f - g^{-1_*^*}{}_*^*bg$$

According to the theorem 9.4.6, the set of solutions of the system of linear equations (17.6.2) is right space of columns.  □

DEFINITION 17.6.11. *Let $A$-number $b$ be $_*^*$-eigenvalue of the pair of matrices $(f, g^{-1_*^*})$ where $g$ is the $_*^*$-non-singular matrix. A row $v$ is called* **eigenrow** *of the pair of matrices $(f, g)$ corresponding to $_*^*$-eigenvalue $b$, if the following equality is true*

$$v_*^*f = v_*^*gb_*^*g^{-1_*^*}$$

□

Since

$$gb_*^*g^{-1_*^*} = g_*^*(bg^{-1_*^*})$$
$$= (g^{-1_*^*})^{-1_*^*}{}_*^*(bg^{-1_*^*})$$

then the definition 17.6.11 follows from the definition 17.6.2 and from the theorem 17.5.9.

THEOREM 17.6.12. *Let $A$-number $b$ be $_*^*$-eigenvalue of the pair of matrices $(f, g^{-1_*^*})$ where $g$ is the $_*^*$-non-singular matrix. The set of eigenrows of the pair of matrices $(f, g)$ corresponding to $_*^*$-eigenvalue $b$ is left $A$-vector space of rows.*

Proof. According to the definitions 17.6.1, 17.6.9, coordinates of eigenrow are solution of the system of linear equations

$$(17.6.3) \qquad (f - g^{-1^*}{}_{*}{}^{*}bg)^*{}_{*}v = 0$$

with $^*_{*}$-singular matrix

$$f - g^{-1^*}{}_{*}{}^{*}bg$$

According to the theorem 9.4.6, the set of solutions of the system of linear equations (17.6.3) is right space of rows. □

Proof. According to the definitions 17.6.2, 17.6.11, coordinates of eigenrow are solution of the system of linear equations

$$(17.6.4) \qquad v_{*}{}^{*}(f - gb_{*}{}^{*}g^{-1^*}) = 0$$

with $_{*}{}^{*}$-singular matrix

$$f - gb_{*}{}^{*}g^{-1^*}$$

According to the theorem 9.4.6, the set of solutions of the system of linear equations (17.6.4) is left space of rows. □

Remark 17.6.13. *Let*

$$(17.6.5) \qquad g^i_j \in Z(A) \qquad i = 1, ..., n \qquad j = 1, ..., n$$

*Then*

$$(17.6.6) \qquad g^{-1^*}{}_{*}{}^{*}bg = (g^{-1^*}{}_{*}{}^{*}g)b = E_n b$$

$$(17.6.7) \qquad gb^*{}_{*}g^{-1^*} = b(g^*{}_{*}g^{-1^*}) = bE_n$$

*From equalities* (17.6.6), (17.6.7), *it follows that, if the condition* (17.6.5), *is satisfied, then we can simplify definitions of eigenvalue and eigenvector of matrix and get definitions which correspond to definitions in commutative algebra.* □

Definition 17.6.14. *$A$-number $b$ is called $^*_{*}$-eigenvalue of the matrix $f$ if the matrix $f - bE_n$ is $^*_{*}$-singular matrix.* □

Definition 17.6.17. *$A$-number $b$ is called $_{*}{}^{*}$-eigenvalue of the matrix $f$ if the matrix $f - bE_n$ is $_{*}{}^{*}$-singular matrix.*

Theorem 17.6.15. *Let $f$ be $n \times n$ matrix of $A$-numbers and*

$$_{*}{}^{*}\operatorname{rank} a = k < n$$

*Then $0$ is $^*_{*}$-eigenvalue of multiplicity $n - k$.*

Theorem 17.6.18. *Let $f$ be $n \times n$ matrix of $A$-numbers and*

$$_{*}{}^{*}\operatorname{rank} a = k < n$$

*Then $0$ is $_{*}{}^{*}$-eigenvalue of multiplicity $n - k$.*

Proof. The theorem follows from the equality

$$a = a - 0E$$

from the definition 17.6.14 and the theorem 9.4.4. □

Proof. The theorem follows from the equality

$$a = a - 0E$$

from the definition 17.6.17 and the theorem 9.4.4. □

Theorem 17.6.16. *Let $f$ be $n \times n$ matrix of $A$-numbers. Then non-zero $^*_{*}$-eigenvalues of the matrix $f$ are roots of any equation*

$$\det({}^*_{*})^i_j (a - bE) = 0$$

Theorem 17.6.19. *Let $f$ be $n \times n$ matrix of $A$-numbers. Then non-zero $_{*}{}^{*}$-eigenvalues of the matrix $f$ are roots of any equation*

$$\det({}_{*}{}^{*})^i_j (a - bE) = 0$$

Proof. The theorem follows from the theorem 9.3.3 and from the definition 17.5.3. □

Proof. The theorem follows from the theorem 9.3.4 and from the definition 17.5.3. □

Definition 17.6.20. *Let $A$-number $b$ be $_*{}^*$-eigenvalue of the matrix $f$. A column $v$ is called* **eigencolumn** *of matrix $f$ corresponding to $_*{}^*$-eigenvalue $b$, if the following equality is true*

$$(17.6.8) \qquad f_*{}^*v = bv$$

□

Definition 17.6.21. *Let $A$-number $b$ be $^*{}_*$-eigenvalue of the matrix $f$. A column $v$ is called* **eigencolumn** *of matrix $f$ corresponding to $^*{}_*$-eigenvalue $b$, if the following equality is true*

$$(17.6.9) \qquad v^*{}_*f = vb$$

□

Definition 17.6.22. *Let $A$-number $b$ be $_*{}^*$-eigenvalue of the matrix $f$. A row $v$ is called* **eigenrow** *of matrix $f$ corresponding to $_*{}^*$-eigenvalue $b$, if the following equality is true*

$$(17.6.10) \qquad v_*{}^*f = vb$$

□

Definition 17.6.23. *Let $A$-number $b$ be $^*{}_*$-eigenvalue of the matrix $f$. A row $v$ is called* **eigenrow** *of matrix $f$ corresponding to $^*{}_*$-eigenvalue $b$, if the following equality is true*

$$(17.6.11) \qquad f^*{}_*v = bv$$

□

Theorem 17.6.24. *Let $A$-number $b$ be $_*{}^*$-eigenvalue of the matrix $a$. The set of eigencolumns of matrix $a$ corresponding to $_*{}^*$-eigenvalue $b$ is right $A$-vector space of columns.*

Proof. According to the definitions 17.6.14, 17.6.20, coordinates of eigencolumn are solution of the system of linear equations

$$(17.6.12) \qquad (f - bE_n)_*{}^*v = 0$$

with $_*{}^*$-singular matrix $f - bE_n$. According to the theorem 9.4.6, the set of solutions of the system of linear equations (17.6.12) is right $A$-vector space of columns. □

Theorem 17.6.25. *Let $A$-number $b$ be $^*{}_*$-eigenvalue of the matrix $a$. The set of eigencolumns of matrix $a$ corresponding to $^*{}_*$-eigenvalue $b$ is left $A$-vector space of columns.*

Proof. According to the definitions 17.6.17, 17.6.21, coordinates of eigencolumn are solution of the system of linear equations

$$(17.6.13) \qquad v^*{}_*(f - E_nb) = 0$$

with $^*{}_*$-singular matrix $f - E_nb$. According to the theorem 9.4.6, the set of solutions of the system of linear equations (17.6.13) is left $A$-vector space of columns. □

Theorem 17.6.26. *Let $A$-number $b$ be $_*{}^*$-eigenvalue of the matrix $a$. The set of eigenrows of matrix $a$ corresponding to $_*{}^*$-eigenvalue $b$ is left $A$-vector space of rows.*

Proof. According to the defi-

Theorem 17.6.27. *Let $A$-number $b$ be $^*{}_*$-eigenvalue of the matrix $a$. The set of eigenrows of matrix $a$ corresponding to $^*{}_*$-eigenvalue $b$ is right $A$-vector space of rows.*

Proof. According to the defi-

nitions 17.6.14, 17.6.22, coordinates of eigenrow are solution of the system of linear equations

$$(17.6.14) \qquad v_*{}^*(f - E_n b) = 0$$

with $_*^*$-singular matrix $f - E_n b$. According to the theorem 9.4.6, the set of solutions of the system of linear equations (17.6.14) is left $A$-vector space of rows. □

nitions 17.6.17, 17.6.23, coordinates of eigenrow are solution of the system of linear equations

$$(17.6.15) \qquad (f - b E_n)^*{}_* v = 0$$

with $^*_*$-singular matrix $f - b E_n$. According to the theorem 9.4.6, the set of solutions of the system of linear equations (17.6.15) is right $A$-vector space of rows. □

# Left and Right Eigenvalues

Let $A$ be associative division $D$-algebra.

## 18.1. Left and Right $_*{}^*$-Eigenvalues

**DEFINITION 18.1.1.** *Let* $a_2$ *be* $n \times n$ *matrix which is* $_*{}^*$-*similar to diagonal matrix* $a_1$

$$a_1 = \mathrm{diag}(b(1), ..., b(n))$$

*Thus, there exist non-*$_*{}^*$*-singular matrix* $u_2$ *such that*

(18.1.1) $\qquad u_{2*}{}^* a_{2*}{}^* u_2^{-1*}{}^* = a_1$

(18.1.2) $\qquad u_{2*}{}^* a_2 = a_{1*}{}^* u_2$

*The row* $u_2^i$ *of the matrix* $u_2$ *satisfies to the equality*

(18.1.3) $\qquad u_{2}^i{}_*{}^* a_2 = b(i) u_2^i$

*The A-number* $b(i)$ *is called left* $_*{}^*$-**eigenvalue** *and row vector* $u_2^i$ *is called* **eigenrow** *for left* $_*{}^*$-*eigenvalue* $b(i)$. □

**DEFINITION 18.1.2.** *Let* $a_2$ *be* $n \times n$ *matrix which is* $_*{}^*$-*similar to diagonal matrix* $a_1$

$$a_1 = \mathrm{diag}(b(1), ..., b(n))$$

*Thus, there exist non-*$_*{}^*$*-singular matrix* $u_2$ *such that*

(18.1.4) $\qquad u_2^{-1*}{}^* {}_*{}^* a_{2*}{}^* u_2 = a_1$

(18.1.5) $\qquad a_{2*}{}^* u_2 = u_{2*}{}^* a_1$

*The column* $u_{2i}$ *of the matrix* $u_2$ *satisfies to the equality*

(18.1.6) $\qquad a_{2*}{}^* u_{2i} = u_{2i} b(i)$

*The A-number* $b(i)$ *is called right* $_*{}^*$-**eigenvalue** *and column vector* $u_{2i}$ *is called* **eigencolumn** *for right* $_*{}^*$-*eigenvalue* $b(i)$. □

**DEFINITION 18.1.3.** *The set* $_*{}^*$-$\mathrm{spec}(a)$ *of all left and right* $_*{}^*$-*eigenvalues is called* $_*{}^*$-*spectrum of the matrix* $a$. □

**THEOREM 18.1.4.** *Let* $A$-*number* $b$ *be left* $_*{}^*$-*eigenvalue of the matrix* $a_2$.[18.1] *Then any* $A$-*number which* $A$-*conjugated with* $A$-*number* $b$, *is left* $_*{}^*$-*eigenvalue of the matrix* $a_2$.

PROOF. Let $v$ be eigenrow for left $_*{}^*$-eigenvalue $b$. Let $c \neq 0$ be $A$-number. The theorem follows from the equality

$$cv_*{}^* a_2 = cbv = cbc^{-1} \, cv$$

□

**THEOREM 18.1.5.** *Let* $A$-*number* $b$ *be right* $_*{}^*$-*eigenvalue of the matrix* $a_2$.[18.1] *Then any* $A$-*number which* $A$-*conjugated with* $A$-*number* $b$, *is right* $_*{}^*$-*eigenvalue of the matrix* $a_2$.

PROOF. Let $v$ be eigencolumn for right $_*{}^*$-eigenvalue $b$. Let $c \neq 0$ be $A$-number. The theorem follows from the equality

$$a_{2*}{}^* vc = vbc = vc \, c^{-1} bc$$

□

**THEOREM 18.1.6.** *Let* $v$ *be eigenrow for left* $_*{}^*$-*eigenvalue* $b$. *The vector* $vc$

**THEOREM 18.1.7.** *Let* $v$ *be eigencolumn for right* $_*{}^*$-*eigenvalue* $b$. *The vec-*

---

[18.1] See the similar statement and proof on the page [20]-375.

*is eigenrow for left $_*{}^*$-eigenvalue $b$ iff $A$-number $c$ commutes with all entries of the matrix $a_2$.*

PROOF. The vector $vc$ is eigenrow for left $_*{}^*$-eigenvalue $b$ iff the following equality is true

$$(18.1.7) \qquad bv_ic = v_ja_2{}^j_i c = v_j ca_{2i}{}^j$$

The equality (18.1.7) holds iff $A$-number $c$ commutes with all entries of the matrix $a_2$. □

According to theorems 18.1.4, 18.1.6, the set of eigenrows for given left $_*{}^*$-eigenvalue is not $A$-vector space.

---

THEOREM 18.1.8. *Eigenvector $u_2^i$ of the matrix $a_2$ for left $_*{}^*$-eigenvalue $b(i)$ does not depend on the choice of the basis $\overline{\overline{e}}_2$*

$$(18.1.9) \qquad e_1^i = u_{2*}{}^i {}^* \overline{\overline{e}}_2$$

PROOF. According to the theorem 10.3.9, we can consider the matrix $a_2$ as matrix of automorphism $a$ of left $A$-vector space with respect to the basis $\overline{\overline{e}}_2$. According to the theorem 12.3.6, the matrix $u_2$ is the matrix is the matrix of passive transformation

$$(18.1.10) \qquad \overline{\overline{e}}_1 = u_{2*}{}^* \overline{\overline{e}}_2$$

and the matrix $a_1$ is the matrix of automorphism $a$ of left $A$-vector space with respect to the basis $\overline{\overline{e}}_1$.

The coordinate matrix of the basis $\overline{\overline{e}}_1$ with respect to the basis $\overline{\overline{e}}_1$ is identity matrix $E_n$. Therefore, the vector $e_1^i$ is eigenrow for left $_*{}^*$-eigenvalue $b(i)$

$$(18.1.11) \qquad e_{1*}^i {}^* a_1 = b(i)e_1^i$$

The equality

$$(18.1.12) \qquad e_1^i {}_* {}^* u_{2*} {}^* a_{2*} {}^* u_2^{-1*}{}^* = b(i)e_1^i$$

follows from equalities (18.1.1), (18.1.11). The equality

$$(18.1.13) \qquad e_1^i {}_* {}^* u_{2*} {}^* a_2 = b(i)e_1^i {}_* {}^* u_2$$

---

*tor $cv$ is eigencolumn for right $_*{}^*$-eigenvalue $b$ iff $A$-number $c$ commutes with all entries of the matrix $a_2$.*

PROOF. The vector $cv$ is eigencolumn for right $_*{}^*$-eigenvalue $b$ iff the following equality is true

$$(18.1.8) \qquad cv^i b = ca_2{}^j_i v^i = a_j^i cv^j$$

The equality (18.1.8) holds iff $A$-number $c$ commutes with all entries of the matrix $a_2$. □

According to theorems 18.1.5, 18.1.7, the set of eigencolumns for given right $_*{}^*$-eigenvalue is not $A$-vector space.

---

THEOREM 18.1.9. *Eigenvector $u_{2i}$ of the matrix $a_2$ for right $_*{}^*$-eigenvalue $b(i)$ does not depend on the choice of the basis $\overline{\overline{e}}_2$*

$$(18.1.15) \qquad e_{1i} = \overline{\overline{e}}_{2*} {}^* u_{2i}$$

PROOF. According to the theorem 11.3.8, we can consider the matrix $a_2$ as matrix of automorphism $a$ of right $A$-vector space with respect to the basis $\overline{\overline{e}}_2$. According to the theorem 12.6.5, the matrix $u_2$ is the matrix is the matrix of passive transformation

$$(18.1.16) \qquad \overline{\overline{e}}_1 = \overline{\overline{e}}_{2*} {}^* u_2$$

and the matrix $a_1$ is the matrix of automorphism $a$ of right $A$-vector space with respect to the basis $\overline{\overline{e}}_1$.

The coordinate matrix of the basis $\overline{\overline{e}}_1$ with respect to the basis $\overline{\overline{e}}_1$ is identity matrix $E_n$. Therefore, the vector $e_{1i}$ is eigencolumn for right $_*{}^*$-eigenvalue $b(i)$

$$(18.1.17) \qquad a_{1*} {}^* e_{1i} = e_{1i} b(i)$$

The equality

$$(18.1.18) \qquad \begin{aligned} u_2^{-1*}{}^* {}_* {}^* a_{2*} {}^* u_{2*} {}^* e_{1i} \\ = e_{1i} b(i) \end{aligned}$$

follows from equalities (18.1.4), (18.1.17). The equality

$$(18.1.19) \qquad a_{2*} {}^* u_{2*} {}^* e_{1i} = u_{2*} {}^* e_{1i} b(i)$$

follows from the equality (18.1.12). The equality

$$(18.1.14) \qquad e_{1*}^i{}^* u_2 = u_2^i$$

follows from equalities (18.1.3), (18.1.13). The theorem follows from equalities (18.1.10), (18.1.14). □

The equality (18.1.9) follows from the equality (18.1.10). However, these equalities have different meaning. The equality (18.1.10) is representation of passive transformation. In the equality (18.1.9), we see the expansion of the eigenvector for given left $_*{}^*$-eigenvalue with respect to selected basis. The goal of the proof of the theorem 18.1.8 is to explain the meaning of the equality (18.1.9).

follows from the equality (18.1.18). The equality

$$(18.1.20) \qquad u_{2*}{}^* e_{1i} = u_{2i}$$

follows from equalities (18.1.6), (18.1.19). The theorem follows from equalities (18.1.16), (18.1.20). □

The equality (18.1.15) follows from the equality (18.1.16). However, these equalities have different meaning. The equality (18.1.16) is representation of passive transformation. In the equality (18.1.15), we see the expansion of the eigenvector for given right $_*{}^*$-eigenvalue with respect to selected basis. The goal of the proof of the theorem 18.1.9 is to explain the meaning of the equality (18.1.15).

THEOREM 18.1.10. *Any left $_*{}^*$-eigenvalue $b(i)$ is $_*{}^*$-eigenvalue. Eigenrow $u^i$ for left $_*{}^*$-eigenvalue $b(i)$ is eigenrow for $_*{}^*$-eigenvalue $b(i)$.*

THEOREM 18.1.11. *Any right $_*{}^*$-eigenvalue $b(i)$ is $_*{}^*$-eigenvalue. Eigencolumn $u_i$ for right $_*{}^*$-eigenvalue $b(i)$ is eigencolumn for $_*{}^*$-eigenvalue $b(i)$.*

PROOF. All entries of the row of the matrix

$$(18.1.21) \qquad a_1 - b(i)E_n$$

are equal to zero. Therefore, the matrix (18.1.21) is $_*{}^*$-singular and left $_*{}^*$-eigenvalue $b(i)$ is $_*{}^*$-eigenvalue. Since

$$(18.1.22) \qquad e_{1i}^j = \delta_j^i$$

then eigenvector $e_1^i$ for left $_*{}^*$-eigenvalue $b(i)$ is eigenvector for $_*{}^*$-eigenvalue $b(i)$. The equality

$$(18.1.23) \qquad e_{1*}^i{}^* a_1 = e_1^i b(i)$$

follows from the definition 17.6.22.

From the equality (18.1.1), it follows that we can present the matrix (18.1.21) as

$$u_{2*}{}^* a_{2*}{}^* u_2^{-1*}{}^* - b(i)E_n$$
$$= u_{2*}{}^* (a_2 - u_2^{-1*}{}^* {}_*{}^* b(i)u_2)_*{}^* u_2^{-1*}{}^*$$

Since the matrix (18.1.21) is $_*{}^*$-singular, then the matrix

$$(18.1.24) \qquad a_2 - u_2^{-1*}{}^* {}_*{}^* b(i)u_2$$

PROOF. All entries of the row of the matrix

$$(18.1.27) \qquad a_1 - b(i)E_n$$

are equal to zero. Therefore, the matrix (18.1.27) is $_*{}^*$-singular and right $_*{}^*$-eigenvalue $b(i)$ is $_*{}^*$-eigenvalue. Since

$$(18.1.28) \qquad e_{1i}^j = \delta_i^j$$

then eigenvector $e_{1i}$ for right $_*{}^*$-eigenvalue $b(i)$ is eigenvector for $_*{}^*$-eigenvalue $b(i)$. The equality

$$(18.1.29) \qquad a_{1*}{}^* e_{1i} = b(i)e_{1i}$$

follows from the definition 17.6.20.

From the equality (18.1.4), it follows that we can present the matrix (18.1.27) as

$$u_2^{-1*}{}^* {}_*{}^* a_{2*}{}^* u_2 - b(i)E_n$$
$$= u_2^{-1*}{}^* {}_*{}^* (a_2 - u_2 b(i)_*{}^* u_2^{-1*}{}^*)_*{}^* u_2$$

Since the matrix (18.1.27) is $_*{}^*$-singular, then the matrix

$$(18.1.30) \qquad a_2 - u_2 b(i)_*{}^* u_2^{-1*}{}^*$$

is $_*$-singular also. Therefore, and left $_*$-eigenvalue $b(i)$ is $_*$-eigenvalue of the pair of matrices $(a_2, u_2)$.

Equalities

$$(18.1.25) \quad e_{1*}^{i*} u_{2*}^{*} a_{2*}^{*} u_2^{-1*^*} = e_1^i b(i)$$

$$e_{1*}^{i*} u_{2*}^{*} a_2$$

$$(18.1.26) \quad = e_1^i b(i)_*^* u_2$$

$$= e_{1*}^{i*} u_{2*}^{*} u_2^{-1*^*} {}_*^* b(i) u_2$$

follow from equalities (18.1.1), (18.1.23). From equalities (18.1.9), (18.1.26) and from the theorem 18.1.8, it follows that eigenrow $u_2^i$ of the matrix $a_2$ for left $_*$-eigenvalue $b(i)$ is eigenrow $u_2^i$ of the pair of matrices $(a_2, u_2)$ corresponding to $_*$-eigenvalue $b(i)$. □

REMARK 18.1.12. *The equality*

$$(18.1.33) \quad u_{1*}^{i*} u_2^{-1*^*} {}_*^* b(i) u_2 = b(i) u_1^i$$

*follows from the theorem 18.1.10. The equality* (18.1.33) *seems unusual. However, if we slightly change it*

$$(18.1.34) \quad \begin{aligned} u_{1*}^{i*} u_2^{-1*^*} b(i) \\ = b(i) u_{1*}^{i*} u_2^{-1*^*} \end{aligned}$$

*then the equality* (18.1.34) *follows from the equality* (18.1.9). □

is $_*$-singular also. Therefore, and right $_*$-eigenvalue $b(i)$ is $_*$-eigenvalue of the pair of matrices $(a_2, u_2)$.

Equalities

$$(18.1.31)$$
$$u_2^{-1*^*} a_{2*}^{*} u_{2*}^{*} e_{1i} = b(i) e_{1i}$$

$$a_{2*}^{*} u_{2*}^{*} e_{1i}$$

$$(18.1.32) \quad = u_{2*}^{*} b(i) e_{1i}$$

$$= u_2 b(i)_*^* u_2^{-1*^*} {}_*^* u_{2*}^{*} e_{1i}$$

follow from equalities (18.1.4), (18.1.29). From equalities (18.1.15), (18.1.32) and from the theorem 18.1.9, it follows that eigencolumn $u_{2i}$ of the matrix $a_2$ for right $_*$-eigenvalue $b(i)$ is eigencolumn $u_{2i}$ of the pair of matrices $(a_2, u_2)$ corresponding to $_*$-eigenvalue $b(i)$. □

REMARK 18.1.13. *The equality*

$$(18.1.35)$$
$$u_2 b(i)_*^* u_2^{-1*^*} {}_*^* u_{1i} = u_{1i} b(i)$$

*follows from the theorem 18.1.11. The equality* (18.1.35) *seems unusual. However, if we slightly change it*

$$(18.1.36) \quad \begin{aligned} b(i) u_2^{-1*^*} {}_*^* u_{1i} \\ = u_2^{-1*^*} {}_*^* u_{1i} b(i) \end{aligned}$$

*then the equality* (18.1.36) *follows from the equality* (18.1.15). □

## 18.2. Left and Right $_*$-Eigenvalues

DEFINITION 18.2.1. *Let $a_2$ be $n \times n$ matrix which is $_*$-similar to diagonal matrix $a_1$*

$$a_1 = \mathrm{diag}(b(1), ..., b(n))$$

*Thus, there exist non-$_*$-singular matrix $u_2$ such that*

$$(18.2.1) \quad u_{2*}^{*} a_{2*}^{*} u_2^{-1*^*} = a_1$$

$$(18.2.2) \quad u_{2*}^{*} a_2 = a_{1*}^{*} u_2$$

*The column $u_{2i}$ of the matrix $u_2$ satisfies to the equality*

$$(18.2.3) \quad u_{2i*}^{*} a_2 = b(i) u_{2i}$$

DEFINITION 18.2.2. *Let $a_2$ be $n \times n$ matrix which is $_*$-similar to diagonal matrix $a_1$*

$$a_1 = \mathrm{diag}(b(1), ..., b(n))$$

*Thus, there exist non-$_*$-singular matrix $u_2$ such that*

$$(18.2.4) \quad u_2^{-1*^*} {}_*^* a_{2*}^{*} u_2 = a_1$$

$$(18.2.5) \quad a_{2*}^{*} u_2 = u_{2*}^{*} a_1$$

*The row $u_2^i$ of the matrix $u_2$ satisfies to the equality*

$$(18.2.6) \quad a_{2*}^{*} u_2^i = u_2^i b(i)$$

*The $A$-number $b(i)$ is called left $^*_*$-**eigenvalue** and column vector $u_{2i}$ is called* **eigencolumn** *for left $^*_*$-eigenvalue $b(i)$.*                    □

*The $A$-number $b(i)$ is called right $^*_*$-**eigenvalue** and row vector $u_2^i$ is called* **eigenrow** *for right $^*_*$-eigenvalue $b(i)$.*                    □

**DEFINITION 18.2.3.** *The set $^*_*$-$\mathrm{spec}(a)$ of all left and right $^*_*$-eigenvalues is called $^*_*$-spectrum of the matrix $a$.*                    □

**THEOREM 18.2.4.** *Let $A$-number $b$ be left $^*_*$-eigenvalue of the matrix $a_2$.[18.2] Then any $A$-number which $A$-conjugated with $A$-number $b$, is left $^*_*$-eigenvalue of the matrix $a_2$.*

**THEOREM 18.2.5.** *Let $A$-number $b$ be right $^*_*$-eigenvalue of the matrix $a_2$.[18.2] Then any $A$-number which $A$-conjugated with $A$-number $b$, is right $^*_*$-eigenvalue of the matrix $a_2$.*

PROOF. Let $v$ be eigencolumn for left $^*_*$-eigenvalue $b$. Let $c \neq 0$ be $A$-number. The theorem follows from the equality

$$cv^*_* a_2 = cbv = cbc^{-1}cv$$
□

PROOF. Let $v$ be eigenrow for right $^*_*$-eigenvalue $b$. Let $c \neq 0$ be $A$-number. The theorem follows from the equality

$$a_2{}^*_* vc = vbc = vc\,c^{-1}bc$$
□

**THEOREM 18.2.6.** *Let $v$ be eigencolumn for left $^*_*$-eigenvalue $b$. The vector $vc$ is eigencolumn for left $^*_*$-eigenvalue $b$ iff $A$-number $c$ commutes with all entries of the matrix $a_2$.*

**THEOREM 18.2.7.** *Let $v$ be eigenrow for right $^*_*$-eigenvalue $b$. The vector $cv$ is eigenrow for right $^*_*$-eigenvalue $b$ iff $A$-number $c$ commutes with all entries of the matrix $a_2$.*

PROOF. The vector $vc$ is eigencolumn for left $^*_*$-eigenvalue $b$ iff the following equality is true

$$(18.2.7) \qquad bv^i c = v^j a_2{}_j^i c = v^j c a_2{}_j^i$$

The equality (18.2.7) holds iff $A$-number $c$ commutes with all entries of the matrix $a_2$.                    □

PROOF. The vector $cv$ is eigenrow for right $^*_*$-eigenvalue $b$ iff the following equality is true

$$(18.2.8) \qquad c v_i b = c a_2{}_i^j v_j = a_2{}_i^j c v_j$$

The equality (18.2.8) holds iff $A$-number $c$ commutes with all entries of the matrix $a_2$.                    □

According to theorems 18.2.4, 18.2.6, the set of eigencolumns for given left $^*_*$-eigenvalue is not $A$-vector space.

According to theorems 18.2.5, 18.2.7, the set of eigenrows for given right $^*_*$-eigenvalue is not $A$-vector space.

**THEOREM 18.2.8.** *Eigenvector $u_{2i}$ of the matrix $a_2$ for left $^*_*$-eigenvalue $b(i)$ does not depend on the choice of the basis $\overline{\overline{e}}_2$*

$$(18.2.9) \qquad e_{1i} = u_{2i}{}^*_* \overline{\overline{e}}_2$$

**THEOREM 18.2.9.** *Eigenvector $u_2^i$ of the matrix $a_2$ for right $^*_*$-eigenvalue $b(i)$ does not depend on the choice of the basis $\overline{\overline{e}}_2$*

$$(18.2.15) \qquad e_1^i = \overline{\overline{e}}_2{}^*_* u_2^i$$

Proof. According to the theorem 10.3.8, we can consider the matrix $a_2$ as matrix of automorphism $a$ of left $A$-vector space with respect to the basis $\overline{\overline{e}}_2$. According to the theorem 12.3.5, the matrix $u_2$ is the matrix is the matrix of passive transformation

$$(18.2.10) \qquad \overline{\overline{e}}_1 = u_2{}^*{}_*\overline{\overline{e}}_2$$

and the matrix $a_1$ is the matrix of automorphism $a$ of left $A$-vector space with respect to the basis $\overline{\overline{e}}_1$.

The coordinate matrix of the basis $\overline{\overline{e}}_1$ with respect to the basis $\overline{\overline{e}}_1$ is identity matrix $E_n$. Therefore, the vector $e_{1i}$ is eigencolumn for left $*_*$-eigenvalue $b(i)$

$$(18.2.11) \qquad e_{1i}{}^*{}_*a_1 = b(i)e_{1i}$$

The equality
$$(18.2.12)$$
$$e_{1i}{}^*{}_*u_2{}^*{}_*a_2{}^*{}_*u_2^{-1^*} = b(i)e_{1i}$$

follows from equalities (18.2.1), (18.2.11). The equality

$$(18.2.13) \qquad e_{1i}{}^*{}_*u_2{}^*{}_*a_2 = b(i)e_{1i}{}^*{}_*u_2$$

follows from the equality (18.2.12). The equality

$$(18.2.14) \qquad e_{1i}{}^*{}_*u_2 = u_{2i}$$

follows from equalities (18.2.3), (18.2.13). The theorem follows from equalities (18.2.10), (18.2.14). □

The equality (18.2.9) follows from the equality (18.2.10). However, these equalities have different meaning. The equality (18.2.10) is representation of passive transformation. In the equality (18.2.9), we see the expansion of the eigenvector for given left $*_*$-eigenvalue with respect to selected basis. The goal of the proof of the theorem 18.2.8 is to explain the meaning of the equality (18.2.9).

Theorem 18.2.10. *Any left $*_*$-eigenvalue $b(i)$ is $*_*$-eigenvalue. Eigencolumn $u_i$ for left $*_*$-eigenvalue $b(i)$ is eigencolumn for $*_*$-eigenvalue $b(i)$.*

Proof. According to the theorem 11.3.9, we can consider the matrix $a_2$ as matrix of automorphism $a$ of right $A$-vector space with respect to the basis $\overline{\overline{e}}_2$. According to the theorem 12.6.6, the matrix $u_2$ is the matrix is the matrix of passive transformation

$$(18.2.16) \qquad \overline{\overline{e}}_1 = \overline{\overline{e}}_2{}^*{}_*u_2$$

and the matrix $a_1$ is the matrix of automorphism $a$ of right $A$-vector space with respect to the basis $\overline{\overline{e}}_1$.

The coordinate matrix of the basis $\overline{\overline{e}}_1$ with respect to the basis $\overline{\overline{e}}_1$ is identity matrix $E_n$. Therefore, the vector $e_1^i$ is eigenrow for right $*_*$-eigenvalue $b(i)$

$$(18.2.17) \qquad a_1{}^*{}_*e_1^i = e_1^i b(i)$$

The equality

$$(18.2.18) \quad u_2^{-1^*}{}^*{}_*a_2{}^*{}_*u_2{}^*{}_*e_1^i = e_1^i b(i)$$

follows from equalities (18.2.4), (18.2.17). The equality

$$(18.2.19) \qquad a_2{}^*{}_*u_2{}^*{}_*e_1^i = u_2{}^*{}_*e_1^i b(i)$$

follows from the equality (18.2.18). The equality

$$(18.2.20) \qquad u_2{}^*{}_*e_1^i = u_2^i$$

follows from equalities (18.2.6), (18.2.19). The theorem follows from equalities (18.2.16), (18.2.20). □

The equality (18.2.15) follows from the equality (18.2.16). However, these equalities have different meaning. The equality (18.2.16) is representation of passive transformation. In the equality (18.2.15), we see the expansion of the eigenvector for given right $*_*$-eigenvalue with respect to selected basis. The goal of the proof of the theorem 18.2.9 is to explain the meaning of the equality (18.2.15).

Theorem 18.2.11. *Any right $*_*$-eigenvalue $b(i)$ is $*_*$-eigenvalue. Eigenrow $u^i$ for right $*_*$-eigenvalue $b(i)$ is eigenrow for $*_*$-eigenvalue $b(i)$.*

Proof. All entries of the row of the matrix

$$(18.2.21) \qquad a_1 - b(i)E_n$$

are equal to zero. Therefore, the matrix (18.2.21) is ${}^*{}_*$-singular and left ${}^*{}_*$-eigenvalue $b(i)$ is ${}^*{}_*$-eigenvalue. Since

$$(18.2.22) \qquad e_{1i}^j = \delta_i^j$$

then eigenvector $e_{1i}$ for left ${}^*{}_*$-eigenvalue $b(i)$ is eigenvector for ${}^*{}_*$-eigenvalue $b(i)$. The equality

$$(18.2.23) \qquad e_{1i}{}^*{}_*a_1 = e_{1i}b(i)$$

follows from the definition 17.6.21.

From the equality (18.2.1), it follows that we can present the matrix (18.2.21) as

$$u_2{}^*{}_*a_2{}^*{}_*u_2^{-1*}{}^* - b(i)E_n$$
$$= u_2{}^*{}_*(a_2 - u_2^{-1*}{}^*{}_*b(i)u_2){}^*{}_*u_2^{-1*}{}^*$$

Since the matrix (18.2.21) is ${}^*{}_*$-singular, then the matrix

$$(18.2.24) \qquad a_2 - u_2^{-1*}{}^*{}_*b(i)u_2$$

is ${}^*{}_*$-singular also. Therefore, and left ${}^*{}_*$-eigenvalue $b(i)$ is ${}^*{}_*$-eigenvalue of the pair of matrices $(a_2, u_2)$.

Equalities
$$(18.2.25)$$
$$e_{1i}{}^*{}_*u_2{}^*{}_*a_2{}^*{}_*u_2^{-1*}{}^* = e_{1i}b(i)$$

$$e_{1i}{}^*{}_*u_2{}^*{}_*a_2$$
$$(18.2.26) \quad = e_{1i}b(i){}^*{}_*u_2$$
$$= e_{1i}{}^*{}_*u_2{}^*{}_*u_2^{-1*}{}^*{}_*b(i)u_2$$

follow from equalities (18.2.1), (18.2.23). From equalities (18.2.9), (18.2.26) and from the theorem 18.2.8, it follows that eigencolumn $u_{2i}$ of the matrix $a_2$ for left ${}^*{}_*$-eigenvalue $b(i)$ is eigencolumn $u_{2i}$ of the pair of matrices $(a_2, u_2)$ corresponding to ${}^*{}_*$-eigenvalue $b(i)$. □

Remark 18.2.12. *The equality*
$$(18.2.33)$$
$$u_{1i}{}^*{}_*u_2^{-1*}{}^*{}_*b(i)u_2 = b(i)u_{1i}$$

*follows from the theorem 18.2.10. The equality (18.2.33) seems unusual. How-*

Proof. All entries of the row of the matrix

$$(18.2.27) \qquad a_1 - b(i)E_n$$

are equal to zero. Therefore, the matrix (18.2.27) is ${}^*{}_*$-singular and right ${}^*{}_*$-eigenvalue $b(i)$ is ${}^*{}_*$-eigenvalue. Since

$$(18.2.28) \qquad e_1{}_i^j = \delta_j^i$$

then eigenvector $e_1^i$ for right ${}^*{}_*$-eigenvalue $b(i)$ is eigenvector for ${}^*{}_*$-eigenvalue $b(i)$. The equality

$$(18.2.29) \qquad a_1{}^*{}_*e_1^i = b(i)e_1^i$$

follows from the definition 17.6.23.

From the equality (18.2.4), it follows that we can present the matrix (18.2.27) as

$$u_2^{-1*}{}^*{}_*a_2{}^*{}_*u_2 - b(i)E_n$$
$$= u_2^{-1*}{}^*{}_*(a_2 - u_2b(i){}^*{}_*u_2^{-1*}{}^*){}^*{}_*u_2$$

Since the matrix (18.2.27) is ${}^*{}_*$-singular, then the matrix

$$(18.2.30) \qquad a_2 - u_2b(i){}^*{}_*u_2^{-1*}{}^*$$

is ${}^*{}_*$-singular also. Therefore, and right ${}^*{}_*$-eigenvalue $b(i)$ is ${}^*{}_*$-eigenvalue of the pair of matrices $(a_2, u_2)$.

Equalities
$$(18.2.31) \quad u_2^{-1*}{}^*{}_*a_2{}^*{}_*u_2{}^*{}_*e_1^i = b(i)e_1^i$$

$$a_2{}^*{}_*u_2{}^*{}_*e_1^i$$
$$(18.2.32) \quad = u_2{}^*{}_*b(i)e_1^i$$
$$= u_2b(i){}^*{}_*u_2^{-1*}{}^*{}_*u_2{}^*{}_*e_1^i$$

follow from equalities (18.2.4), (18.2.29). From equalities (18.2.15), (18.2.32) and from the theorem 18.2.9, it follows that eigenrow $u_2^i$ of the matrix $a_2$ for right ${}^*{}_*$-eigenvalue $b(i)$ is eigenrow $u_2^i$ of the pair of matrices $(a_2, u_2)$ corresponding to ${}^*{}_*$-eigenvalue $b(i)$. □

Remark 18.2.13. *The equality*
$$(18.2.35) \quad u_2b(i){}^*{}_*u_2^{-1*}{}^*{}_*u_1^i = u_1^ib(i)$$

*follows from the theorem 18.2.11. The equality (18.2.35) seems unusual. How-*

*ever, if we slightly change it*

$$(18.2.34) \qquad u_{1i}{}^{*}{}_{*}u_{2}^{-1^{*}}{}_{*}b(i)$$
$$= b(i)u_{1i}{}^{*}{}_{*}u_{2}^{-1^{*}}{}_{*}$$

*then the equality* (18.2.34) *follows from the equality* (18.2.9). □

*ever, if we slightly change it*

$$(18.2.36) \qquad b(i)u_{2}^{-1^{*}}{}_{*}{}_{*}u_{1}^{i}$$
$$= u_{2}^{-1^{*}}{}_{*}{}_{*}u_{1}^{i}b(i)$$

*then the equality* (18.2.36) *follows from the equality* (18.2.15). □

## 18.3. Converse Theorem Problem

QUESTION 18.3.1. *From the theorem* [18.3] *18.1.11, it follows that any right* $_{*}^{*}$*-eigenvalue $d(i)$ is $_{*}^{*}$-eigenvalue. Is the converse theorem true? Namely, is $_{*}^{*}$-eigenvalue a right $_{*}^{*}$-eigenvalue?* □

In general, it is not easy to get an answer to the question 18.3.1. What happens if the multiplicity of $_{*}^{*}$-eigenvalue is greater than 1? Can $n \times n$ matrix have less than $n$ eigenvalues, taking into account the multiplicity?

It may seem that, if $n \times n$ matrix has $n$ eigenvalues of the multiplicity 1, then this is most simple case. If $n \times n$ matrix with entries from commutative algebra has $n$ eigenvalues of the multiplicity 1, then the set of eigenvectors is the basis of vector space of dimension $n$. However this statement is not true when entries of the matrix belong to non-commutative algebra.

Concider the theorem 18.3.2.

THEOREM 18.3.2. *Let $A$ be associative division $D$-algebra. Let $V$ be $A$-vector space. If endomorphism has $k$ different eigenvalues* [18.4] *$b(1)$, ..., $b(k)$, then eigenvectors $v(1)$, ..., $v(k)$ which correspond to different eigenvalues are linearly independent.*

However, the theorem 18.3.2 is not correct. We will consider the proof separately for left and right $A$-vector spaces.

PROOF. We will prove the theorem by induction with respect to $k$.

Since eigenvector is non-zero, the theorem is true for $k = 1$.

STATEMENT 18.3.3. *Let the theorem be true for $k = m - 1$.* ⊙

According to theorems 17.5.13, 17.5.27, we can assume that $\overline{\overline{e}}_2 = \overline{\overline{e}}_1$. Let $b(1)$, ..., $b(m-1)$, $b(m)$ be different eigenvalues and $v(1)$, ..., $v(m-1)$, $v(m)$ corresponding eigenvectors

$$(18.3.1) \qquad f \circ v(i) = v(i)b(i)$$
$$i \neq j \Rightarrow b(i) \neq b(j)$$

Let the statement 18.3.4 be true.

PROOF. We will prove the theorem by induction with respect to $k$.

Since eigenvector is non-zero, the theorem is true for $k = 1$.

STATEMENT 18.3.5. *Let the theorem be true for $k = m - 1$.* ⊙

According to theorems 17.5.20, 17.5.30, we can assume that $\overline{\overline{e}}_2 = \overline{\overline{e}}_1$. Let $b(1)$, ..., $b(m-1)$, $b(m)$ be different eigenvalues and $v(1)$, ..., $v(m-1)$, $v(m)$ corresponding eigenvectors

$$(18.3.5) \qquad f \circ v(i) = b(i)v(i)$$
$$i \neq j \Rightarrow b(i) \neq b(j)$$

Let the statement 18.3.6 be true.

---

[18.3] Similar reasoning is true for left $_{*}^{*}$-eigenvalue, as well for left and right $^{*}_{*}$-eigenvalue.

[18.4] See also the theorem on the page [3]-203.

STATEMENT 18.3.4. *There exists linear dependence*

$$
\begin{aligned}
& a(1)v(1) + ... \\
(18.3.2) \quad & + a(m-1)v(m-1) \\
& + a(m)v(m) = 0
\end{aligned}
$$

*where* $a(1) \neq 0$.                    ⊙

The equality

$$
\begin{aligned}
& a(1)(f \circ v(1)) + ... \\
(18.3.3) \quad & + a(m-1)(f \circ v(m-1)) \\
& + a(m)(f \circ v(m)) = 0
\end{aligned}
$$

follows from the equality (18.3.2) and the theorem 10.3.2. The equality

$$
\begin{aligned}
& a(1)v(1)b(1) + ... \\
(18.3.4) \quad & + a(m-1) \\
& * v(m-1)b(m-1) \\
& + a(m)v(m)b(m) = 0
\end{aligned}
$$

follows from equalities (18.3.1), (18.3.3). Since the product is non-commutative, we cannot reduce expression (18.3.4) to linear dependence.                    □

STATEMENT 18.3.6. *There exists linear dependence*

$$
\begin{aligned}
& v(1)a(1) + ... \\
(18.3.6) \quad & + v(m-1)a(m-1) \\
& + v(m)a(m) = 0
\end{aligned}
$$

*where* $a(1) \neq 0$.                    ⊙

The equality

$$
\begin{aligned}
& (f \circ v(1))a(1) + ... \\
(18.3.7) \quad & + (f \circ v(m-1))a(m-1) \\
& + (f \circ v(m))a(m) = 0
\end{aligned}
$$

follows from the equality (18.3.6) and the theorem 11.3.2. The equality

$$
\begin{aligned}
& b(1)v(1)a(1) + ... \\
(18.3.8) \quad & + b(m-1)v(m-1) \\
& * a(m-1) \\
& + b(m)v(m)a(m) = 0
\end{aligned}
$$

follows from equalities (18.3.5), (18.3.7). Since the product is non-commutative, we cannot reduce expression (18.3.8) to linear dependence.                    □

REMARK 18.3.7. *According to the definition given above, sets of left and right $_*$\*-eigenvalues coincide. However, we considered the most simple case. Cohn considers definitions*

$$\boxed{(18.1.5) \quad a_{2*}{}^*u_2 = u_{2*}{}^*a_1} \qquad \boxed{(18.1.2) \quad u_{2*}{}^*a_2 = a_{1*}{}^*u_2}$$

*where $a_2$ is non-$_*$\*-singular matrix as starting point for research.*                    □

Let multiplicity of all left $_*$\*-eigenvalues of $n \times n$ matrix $a_2$ equal 1. Let $n \times n$ matrix $a_2$ do not have $n$ left $_*$\*-eigenvalues, but have $k$, $k < n$, left $_*$\*-eigenvalues. In such case, the matrix $u_2$ in the equality

$$\boxed{(18.1.2) \quad u_{2*}{}^*a_2 = a_{1*}{}^*u_2}$$

is $n \times k$ matrix and the matrix $a_1$ is $k \times k$ matrix. However, the equality (18.1.2) does not imply that

$$\operatorname{rank}(_*{}^*)u_2 = k$$

Let

$$\operatorname{rank}(_*{}^*)u_2 = k$$

Then there exists $k \times n$ matrix $w \in right(u_2^{-1}{}_*{}^*)$ such that

$$(18.3.9) \qquad u_{2*}{}^*w = E_k$$

Let multiplicity of all right $_*$\*-eigenvalues of $n \times n$ matrix $a_2$ equal 1. Let $n \times n$ matrix $a_2$ do not have $n$ right $_*$\*-eigenvalues, but have $k$, $k < n$, right $_*$\*-eigenvalues. In such case, the matrix $u_2$ in the equality

$$\boxed{(18.1.5) \quad a_{2*}{}^*u_2 = u_{2*}{}^*a_1}$$

is $k \times n$ matrix and the matrix $a_1$ is $k \times k$ matrix. However, the equality (18.1.5) does not imply that

$$\operatorname{rank}(_*{}^*)u_2 = k$$

Let

$$\operatorname{rank}(_*{}^*)u_2 = k$$

Then there exists $n \times k$ matrix $w \in left(u_2^{-1}{}_*{}^*)$ such that

$$(18.3.11) \qquad w_*{}^*u_2 = E_k$$

The equality

$$(18.3.10) \qquad u_{2*}{}^{*}a_{2*}{}^{*}w = a_1$$

follows from equalities (18.1.2), (18.3.9).

    Let

$$\mathrm{rank}({}_{*}{}^{*})u_2 < k$$

Then there exists linear dependence between left eigenvectors. In such case, sets of left and right ${}_{*}{}^{*}$-eigenvalues do not coincide.

The equality

$$(18.3.12) \qquad w_{*}{}^{*}a_{2*}{}^{*}u_2 = a_1$$

follows from equalities (18.1.5), (18.3.11).

    Let

$$\mathrm{rank}({}_{*}{}^{*})u_2 < k$$

Then there exists linear dependence between right eigenvectors. In such case, sets of left and right ${}_{*}{}^{*}$-eigenvalues do not coincide.

## 18.4. Eigenvalue of Matrix of Linear Map

    Let

$$a = \begin{pmatrix} a^{1}_{1\,s0} \otimes a^{1}_{1\,s1} & ... & a^{1}_{n\,s0} \otimes a^{1}_{n\,s1} \\ ... & ... & ... \\ a^{n}_{1\,s0} \otimes a^{n}_{1\,s1} & ... & a^{n}_{n\,s0} \otimes a^{n}_{n\,s1} \end{pmatrix}$$

be matrix of linear map

$$\overline{a} : V \to V$$

of $A$-vector space $V$.

### 18.4.1. ${}_{\circ}{}^{\circ}$-Eigenvalue of Matrix of Linear Map.

DEFINITION 18.4.1. *$A$-number $b$ is called left ${}_{\circ}{}^{\circ}$-eigenvalue of the matrix $a$ if there exists column vector $v$ such that*

$$(18.4.1) \qquad a_{\circ}{}^{\circ}v = bv$$

*The column vector $v$ is called* **eigencolumn** *for left ${}_{\circ}{}^{\circ}$-eigenvalue $b$.* □

DEFINITION 18.4.2. *$A$-number $b$ is called right ${}_{\circ}{}^{\circ}$-eigenvalue of the matrix $a$ if there exists column vector $v$ such that*

$$(18.4.2) \qquad a_{\circ}{}^{\circ}v = vb$$

*The column vector $v$ is called* **eigencolumn** *for right ${}_{\circ}{}^{\circ}$-eigenvalue $b$.* □

THEOREM 18.4.3. *Let entries of the matrix $a$ satisfy the equality*

$$(18.4.3) \qquad a^{i}_{j\,s0} \otimes a^{i}_{j\,s1} = 1 \otimes a^{i}_{1\,j}$$

*Then left ${}_{\circ}{}^{\circ}$-eigenvalue $b$ is left ${}^{*}_{*}$-eigenvalue of the matrix $a_1$.*

THEOREM 18.4.4. *Let entries of the matrix $a$ satisfy the equality*

$$(18.4.4) \qquad a^{i}_{j\,s0} \otimes a^{i}_{j\,s1} = a^{i}_{0\,j} \otimes 1$$

*Then right ${}_{\circ}{}^{\circ}$-eigenvalue $b$ is right ${}^{*}_{*}$-eigenvalue of the matrix $a_0$.*

PROOF OF THEOREM 18.4.3. The equality

$$(18.4.5) \qquad (a^{i}_{j\,s0} \otimes a^{i}_{j\,s1}) \circ v^{j} = (1 \otimes a^{i}_{1\,j}) \circ v^{j} = v^{j} a^{i}_{1\,j}$$

follows from the equality $\boxed{(18.4.3)\quad a_{j\,s0}^{i}\otimes a_{j\,s1}^{i}=1\otimes a_{1\,j}^{i}}$ .

The equality

(18.4.6)
$$\begin{pmatrix}v^1\\ ...\\ v^n\end{pmatrix}*_{*}\begin{pmatrix}a_{1\,1}^{1}&...&a_{1\,n}^{1}\\ ...&...&...\\ a_{1\,1}^{n}&...&a_{1\,n}^{n}\end{pmatrix}=\begin{pmatrix}a_{1\,s0}^{1}\otimes a_{1\,s1}^{1}&...&a_{n\,s0}^{1}\otimes a_{n\,s1}^{1}\\ ...&...&...\\ a_{1\,s0}^{n}\otimes a_{1\,s1}^{n}&...&a_{n\,s0}^{n}\otimes a_{n\,s1}^{n}\end{pmatrix}\circ^{\circ}\begin{pmatrix}v^1\\ ...\\ v^n\end{pmatrix}$$

follows from the equality (18.4.5). The equality

(18.4.7)
$$\begin{pmatrix}a_{1\,s0}^{1}\otimes a_{1\,s1}^{1}&...&a_{n\,s0}^{1}\otimes a_{n\,s1}^{1}\\ ...&...&...\\ a_{1\,s0}^{n}\otimes a_{1\,s1}^{n}&...&a_{n\,s0}^{n}\otimes a_{n\,s1}^{n}\end{pmatrix}\circ^{\circ}\begin{pmatrix}v^1\\ ...\\ v^n\end{pmatrix}=b\begin{pmatrix}v^1\\ ...\\ v^n\end{pmatrix}$$

follows from the definition 18.4.1 of left $_{\circ}^{\circ}$-eigenvalue of the matrix $a$. The equality

(18.4.8)
$$\begin{pmatrix}v^1\\ ...\\ v^n\end{pmatrix}*_{*}\begin{pmatrix}a_{1\,1}^{1}&...&a_{1\,n}^{1}\\ ...&...&...\\ a_{1\,1}^{n}&...&a_{1\,n}^{n}\end{pmatrix}=b\begin{pmatrix}v^1\\ ...\\ v^n\end{pmatrix}$$

follows from equalities (18.4.6), (18.4.7). The theorem follows from the equality (18.4.8) and from the definition 18.2.1 of left $^{*}_{*}$-eigenvalue of the matrix $a_1$.    $\square$

PROOF OF THEOREM 18.4.4.    The equality

(18.4.9)
$$(a_{j\,s0}^{i}\otimes a_{j\,s1}^{i})\circ v^j=(a_{0\,j}^{i}\otimes 1)\circ v^j=a_{0\,j}^{i}v^j$$

follows from the equality $\boxed{(18.4.4)\quad a_{j\,s0}^{i}\otimes a_{j\,s1}^{i}=a_{0\,j}^{i}\otimes 1}$ .

The equality

(18.4.10)
$$\begin{pmatrix}a_{0\,1}^{1}&...&a_{0\,n}^{1}\\ ...&...&...\\ a_{0\,1}^{n}&...&a_{0\,n}^{n}\end{pmatrix}*^{*}\begin{pmatrix}v^1\\ ...\\ v^n\end{pmatrix}=\begin{pmatrix}a_{1\,s0}^{1}\otimes a_{1\,s1}^{1}&...&a_{n\,s0}^{1}\otimes a_{n\,s1}^{1}\\ ...&...&...\\ a_{1\,s0}^{n}\otimes a_{1\,s1}^{n}&...&a_{n\,s0}^{n}\otimes a_{n\,s1}^{n}\end{pmatrix}\circ^{\circ}\begin{pmatrix}v^1\\ ...\\ v^n\end{pmatrix}$$

follows from the equality (18.4.9). The equality

(18.4.11)
$$\begin{pmatrix}a_{1\,s0}^{1}\otimes a_{1\,s1}^{1}&...&a_{n\,s0}^{1}\otimes a_{n\,s1}^{1}\\ ...&...&...\\ a_{1\,s0}^{n}\otimes a_{1\,s1}^{n}&...&a_{n\,s0}^{n}\otimes a_{n\,s1}^{n}\end{pmatrix}\circ^{\circ}\begin{pmatrix}v^1\\ ...\\ v^n\end{pmatrix}=\begin{pmatrix}v^1\\ ...\\ v^n\end{pmatrix}b$$

follows from the definition 18.4.2 of right $_{\circ}^{\circ}$-eigenvalue of the matrix $a$. The equality

(18.4.12)
$$\begin{pmatrix}a_{0\,1}^{1}&...&a_{0\,n}^{1}\\ ...&...&...\\ a_{0\,1}^{n}&...&a_{0\,n}^{n}\end{pmatrix}*^{*}\begin{pmatrix}v^1\\ ...\\ v^n\end{pmatrix}=\begin{pmatrix}v^1\\ ...\\ v^n\end{pmatrix}b$$

follows from equalities (18.4.10), (18.4.11). The theorem follows from the equality (18.4.12) and from the definition 18.1.2 of right $_{*}^{*}$-eigenvalue of the matrix $a_0$.    $\square$

THEOREM 18.4.5. *Let entries of the matrix a satisfy the equality*

$$(18.4.13) \qquad a^i_{j\,s0} \otimes a^i_{j\,s1} = a_{0\,j}^{\,i} \otimes 1$$

*Then left $_\circ{}^\circ$-eigenvalue $b$ is $_*{}^*$-eigenvalue of the matrix $a_0$.*

THEOREM 18.4.6. *Let entries of the matrix a satisfy the equality*

$$(18.4.14) \qquad a^i_{j\,s0} \otimes a^i_{j\,s1} = 1 \otimes a_{1\,j}^{\,i}$$

*Then right $_\circ{}^\circ$-eigenvalue $b$ is $^*{}_*$-eigenvalue of the matrix $a_1$.*

PROOF OF THEOREM 18.4.5. The equality

$$(18.4.15) \qquad (a^i_{j\,s0} \otimes a^i_{j\,s1}) \circ v^j = (a_{0\,j}^{\,i} \otimes 1) \circ v^j = a_{0\,j}^{\,i} v^j$$

follows from the equality $\boxed{(18.4.13) \quad a^i_{j\,s0} \otimes a^i_{j\,s1} = a_{0\,j}^{\,i} \otimes 1}$ . The equality

$$(18.4.16)$$

$$\begin{pmatrix} a_{0\,1}^{\,1} & ... & a_{0\,n}^{\,1} \\ ... & ... & ... \\ a_{0\,1}^{\,n} & ... & a_{0\,n}^{\,n} \end{pmatrix} *^* \begin{pmatrix} v^1 \\ ... \\ v^n \end{pmatrix} = \begin{pmatrix} a^1_{1\,s0} \otimes a^1_{1\,s1} & ... & a^1_{n\,s0} \otimes a^1_{n\,s1} \\ ... & ... & ... \\ a^n_{1\,s0} \otimes a^n_{1\,s1} & ... & a^n_{n\,s0} \otimes a^n_{n\,s1} \end{pmatrix} \circ^\circ \begin{pmatrix} v^1 \\ ... \\ v^n \end{pmatrix}$$

follows from the equality (18.4.15). The equality

$$(18.4.17) \qquad \begin{pmatrix} a^1_{1\,s0} \otimes a^1_{1\,s1} & ... & a^1_{n\,s0} \otimes a^1_{n\,s1} \\ ... & ... & ... \\ a^n_{1\,s0} \otimes a^n_{1\,s1} & ... & a^n_{n\,s0} \otimes a^n_{n\,s1} \end{pmatrix} \circ^\circ \begin{pmatrix} v^1 \\ ... \\ v^n \end{pmatrix} = b \begin{pmatrix} v^1 \\ ... \\ v^n \end{pmatrix}$$

follows from the definition 18.4.1 of left $_\circ{}^\circ$-eigenvalue of the matrix $a$. The equality

$$(18.4.18) \qquad \begin{pmatrix} a_{0\,1}^{\,1} & ... & a_{0\,n}^{\,1} \\ ... & ... & ... \\ a_{0\,1}^{\,n} & ... & a_{0\,n}^{\,n} \end{pmatrix} *^* \begin{pmatrix} v^1 \\ ... \\ v^n \end{pmatrix} = b \begin{pmatrix} v^1 \\ ... \\ v^n \end{pmatrix}$$

follows from equalities (18.4.16), (18.4.17). The theorem follows from the equality (18.4.18), from the definition 17.6.14 of $^*{}_*$-eigenvalue of the matrix $a_0$ and from the definition 17.6.20 of corresponding eigenvector $v$. $\square$

PROOF OF THEOREM 18.4.6. The equality

$$(18.4.19) \qquad (a^i_{j\,s0} \otimes a^i_{j\,s1}) \circ v^j = (1 \otimes a_{1\,j}^{\,i}) \circ v^j = v^j a_{1\,j}^{\,i}$$

follows from the equality $\boxed{(18.4.14) \quad a^i_{j\,s0} \otimes a^i_{j\,s1} = 1 \otimes a_{1\,j}^{\,i}}$ . The equality

$$(18.4.20)$$

$$\begin{pmatrix} v^1 \\ ... \\ v^n \end{pmatrix} *_* \begin{pmatrix} a_{1\,1}^{\,1} & ... & a_{1\,n}^{\,1} \\ ... & ... & ... \\ a_{1\,1}^{\,n} & ... & a_{1\,n}^{\,n} \end{pmatrix} = \begin{pmatrix} a^1_{1\,s0} \otimes a^1_{1\,s1} & ... & a^1_{n\,s0} \otimes a^1_{n\,s1} \\ ... & ... & ... \\ a^n_{1\,s0} \otimes a^n_{1\,s1} & ... & a^n_{n\,s0} \otimes a^n_{n\,s1} \end{pmatrix} \circ^\circ \begin{pmatrix} v^1 \\ ... \\ v^n \end{pmatrix}$$

follows from the equality (18.4.19). The equality

$$(18.4.21) \qquad \begin{pmatrix} a_{1s0}^1 \otimes a_{1s1}^1 & ... & a_{ns0}^1 \otimes a_{ns1}^1 \\ ... & ... & ... \\ a_{1s0}^n \otimes a_{1s1}^n & ... & a_{ns0}^n \otimes a_{ns1}^n \end{pmatrix} \circ^\circ \begin{pmatrix} v^1 \\ ... \\ v^n \end{pmatrix} = \begin{pmatrix} v^1 \\ ... \\ v^n \end{pmatrix} b$$

follows from the definition 18.4.2 of right $_\circ{}^\circ$-eigenvalue of the matrix $a$. The equality

$$(18.4.22) \qquad \begin{pmatrix} v^1 \\ ... \\ v^n \end{pmatrix} {}_*{}^* \begin{pmatrix} a_{11}^1 & ... & a_{1n}^1 \\ ... & ... & ... \\ a_{11}^n & ... & a_{1n}^n \end{pmatrix} = \begin{pmatrix} v^1 \\ ... \\ v^n \end{pmatrix} b$$

follows from equalities (18.4.20), (18.4.21). The theorem follows from the equality (18.4.22), from the definition 17.6.17 of $_*{}^*$-eigenvalue of the matrix $a_1$ and from the definition 17.6.21 of corresponding eigenvector $v$.  □

**Theorem** 18.4.7. *Let the column vector $v$ be eigencolumn for left $_\circ{}^\circ$-eigenvalue $b$ of the matrix $a$. Let A-number $b$ satisfy the condition*

$$(18.4.23) \qquad b \in \bigcap_{i=1}^{n} \bigcap_{j=1}^{n} Z(A, a_{js0}^i)$$

*Then for any polynomial*

$$(18.4.24) \qquad p(x) = p_0 + p_1 x + ... + p_n x^n$$
$$p_0, \ ..., \ p_n \in D$$

*the column vector*

$$cv = \begin{pmatrix} cv^1 \\ ... \\ cv^n \end{pmatrix}$$
$$c = p(b)$$

*is also eigencolumn for left $_\circ{}^\circ$-eigenvalue $b$.*

**Theorem** 18.4.8. *Let the column vector $v$ be eigencolumn for right $_\circ{}^\circ$-eigenvalue $b$ of the matrix $a$. Let A-number $b$ satisfy the condition*

$$(18.4.26) \qquad b \in \bigcap_{i=1}^{n} \bigcap_{j=1}^{n} Z(A, a_{js1}^i)$$

*Then for any polynomial*

$$(18.4.27) \qquad p(x) = p_0 + p_1 x + ... + p_n x^n$$
$$p_0, \ ..., \ p_n \in D$$

*the column vector*

$$vc = \begin{pmatrix} v^1 c \\ ... \\ v^n c \end{pmatrix}$$
$$c = p(b)$$

*is also eigencolumn for right $_\circ{}^\circ$-eigenvalue $b$.*

**Proof.** According to the theorem 5.1.25, the equality

$$(18.4.25) \qquad \begin{aligned} & (a_{js0}^i \otimes a_{js1}^i) \circ (cv^j) \\ = \ & a_{js0}^i cv^j a_{js1}^i \\ = \ & c(a_{js0}^i v^j a_{js1}^i) \\ = \ & cbv^i = b(cv^i) \end{aligned}$$

follows from the condition (18.4.23). According to the definition 18.4.1, the equality (18.4.25) implies that the column vec-

**Proof.** According to the theorem 5.1.25, the equality

$$(18.4.28) \qquad \begin{aligned} & (a_{js0}^i \otimes a_{js1}^i) \circ (v^j c) \\ = \ & a_{js0}^i v^j c a_{js1}^i \\ = \ & (a_{js0}^i v^j a_{js1}^i) c \\ = \ & v^i bc = (v^i c) b \end{aligned}$$

follows from the condition (18.4.26). According to the definition 18.4.2, the equality (18.4.28) implies that the column vec-tor $cv$ is eigencolumn for left $_\circ{}^\circ$-eigen­value $b$ of the matrix $a$.          □

tor $vc$ is eigencolumn for right $_\circ{}^\circ$-eigen­value $b$ of the matrix $a$.          □

### 18.4.2. $_\circ{}^\circ$-Eigenvalue of Matrix of Linear Map.

**DEFINITION 18.4.9.** *A-number $b$ is called left $_\circ{}^\circ$-**eigenvalue** of the matrix $a$ if there exists row vector $v$ such that*

$$(18.4.29) \qquad a^\circ{}_\circ v = bv$$

*The row vector $v$ is called **eigenrow** for left $_\circ{}^\circ$-eigenvalue $b$.*          □

**DEFINITION 18.4.10.** *A-number $b$ is called right $_\circ{}^\circ$-**eigenvalue** of the matrix $a$ if there exists row vector $v$ such that*

$$(18.4.30) \qquad a^\circ{}_\circ v = vb$$

*The row vector $v$ is called **eigenrow** for right $_\circ{}^\circ$-eigenvalue $b$.*          □

**THEOREM 18.4.11.** *Let entries of the matrix $a$ satisfy the equality*

$$(18.4.31) \qquad a^j_{i\,s0} \otimes a^j_{i\,s1} = 1 \otimes a^j_{1\,i}$$

*Then left $_\circ{}^\circ$-eigenvalue $b$ is left $_*{}^*$-eigen­value of the matrix $a_1$.*

**THEOREM 18.4.12.** *Let entries of the matrix $a$ satisfy the equality*

$$(18.4.32) \qquad a^j_{i\,s0} \otimes a^j_{i\,s1} = a^j_{0\,i} \otimes 1$$

*Then right $_\circ{}^\circ$-eigenvalue $b$ is right $_*{}^*$-eigenvalue of the matrix $a_0$.*

PROOF OF THEOREM 18.4.11.    The equality

$$(18.4.33) \qquad (a^j_{i\,s0} \otimes a^j_{i\,s1}) \circ v_i = (1 \otimes a^j_{1\,i}) \circ v_i = v_i a^j_{1\,i}$$

follows from the equality $\boxed{(18.4.31) \quad a^j_{i\,s0} \otimes a^j_{i\,s1} = 1 \otimes a^j_{1\,i}}$ .
The equality
$$(18.4.34)$$

$$\begin{pmatrix} v_1 & \dots & v_n \end{pmatrix} *{}* \begin{pmatrix} a_1{}^1_1 & \dots & a_1{}^1_n \\ \dots & \dots & \dots \\ a_1{}^n_1 & \dots & a_1{}^n_n \end{pmatrix} = \begin{pmatrix} a^1_{1\,s0} \otimes a^1_{1\,s1} & \dots & a^n_{1\,s0} \otimes a^n_{1\,s1} \\ \dots & \dots & \dots \\ a^1_{n\,s0} \otimes a^1_{n\,s1} & \dots & a^n_{n\,s0} \otimes a^n_{n\,s1} \end{pmatrix} \circ{}_\circ \begin{pmatrix} v_1 & \dots & v_n \end{pmatrix}$$

follows from the equality (18.4.33). The equality

$$(18.4.35) \qquad \begin{pmatrix} a^1_{1\,s0} \otimes a^1_{1\,s1} & \dots & a^n_{1\,s0} \otimes a^n_{1\,s1} \\ \dots & \dots & \dots \\ a^1_{n\,s0} \otimes a^1_{n\,s1} & \dots & a^n_{n\,s0} \otimes a^n_{n\,s1} \end{pmatrix} \circ{}_\circ \begin{pmatrix} v_1 & \dots & v_n \end{pmatrix} = b \begin{pmatrix} v_1 & \dots & v_n \end{pmatrix}$$

follows from the definition 18.4.9 of left $_\circ{}^\circ$-eigenvalue of the matrix $a$. The equality

$$(18.4.36) \qquad \begin{pmatrix} v_1 & \dots & v_n \end{pmatrix} *{}* \begin{pmatrix} a_1{}^1_1 & \dots & a_1{}^1_n \\ \dots & \dots & \dots \\ a_1{}^n_1 & \dots & a_1{}^n_n \end{pmatrix} = b \begin{pmatrix} v_1 & \dots & v_n \end{pmatrix}$$

follows from equalities (18.4.34), (18.4.35). The theorem follows from the equality (18.4.36) and from the definition 18.1.1 of left $_*{}^*$-eigenvalue of the matrix $a_1$.    □

PROOF OF THEOREM 18.4.12.    The equality

$$(18.4.37) \qquad (a^j_{i\,s0} \otimes a^j_{i\,s1}) \circ v_i = (a^j_{0\,i} \otimes 1) \circ v_i = a^j_{0\,i} v_i$$

follows from the equality $\boxed{(18.4.32) \quad a^j_{i\,s0} \otimes a^j_{i\,s1} = a^j_{0\,i} \otimes 1}$ .
The equality
(18.4.38)

$$\begin{pmatrix} a_0{}^1_1 & ... & a_0{}^1_n \\ ... & ... & ... \\ a_0{}^n_1 & ... & a_0{}^n_n \end{pmatrix} *_* \begin{pmatrix} v_1 & ... & v_n \end{pmatrix} = \begin{pmatrix} a^1_{1\,s0} \otimes a^1_{1\,s1} & ... & a^n_{1\,s0} \otimes a^n_{1\,s1} \\ ... & ... & ... \\ a^1_{n\,s0} \otimes a^1_{n\,s1} & ... & a^n_{n\,s0} \otimes a^n_{n\,s1} \end{pmatrix} \circ_\circ \begin{pmatrix} v_1 & ... & v_n \end{pmatrix}$$

follows from the equality (18.4.37). The equality

$$(18.4.39) \quad \begin{pmatrix} a^1_{1\,s0} \otimes a^1_{1\,s1} & ... & a^n_{1\,s0} \otimes a^n_{1\,s1} \\ ... & ... & ... \\ a^1_{n\,s0} \otimes a^1_{n\,s1} & ... & a^n_{n\,s0} \otimes a^n_{n\,s1} \end{pmatrix} \circ_\circ \begin{pmatrix} v_1 & ... & v_n \end{pmatrix} = \begin{pmatrix} v_1 & ... & v_n \end{pmatrix} b$$

follows from the definition 18.4.10 of right $\circ_\circ$-eigenvalue of the matrix $a$. The equality

$$(18.4.40) \quad \begin{pmatrix} a_0{}^1_1 & ... & a_0{}^1_n \\ ... & ... & ... \\ a_0{}^n_1 & ... & a_0{}^n_n \end{pmatrix} *_* \begin{pmatrix} v_1 & ... & v_n \end{pmatrix} = \begin{pmatrix} v_1 & ... & v_n \end{pmatrix} b$$

follows from equalities (18.4.38), (18.4.39). The theorem follows from the equality
(18.4.40) and from the definition 18.2.2 of right $*_*$-eigenvalue of the matrix $a_0$.   □

| THEOREM 18.4.13. *Let entries of the matrix $a$ satisfy the equality* | THEOREM 18.4.14. *Let entries of the matrix $a$ satisfy the equality* |
|---|---|
| $(18.4.41) \quad a^j_{i\,s0} \otimes a^j_{i\,s1} = a^j_{0\,i} \otimes 1$ | $(18.4.42) \quad a^j_{i\,s0} \otimes a^j_{i\,s1} = 1 \otimes a^j_{1\,i}$ |
| *Then left $\circ_\circ$-eigenvalue $b$ is $*_*$-eigenvalue of the matrix $a_0$.* | *Then right $\circ_\circ$-eigenvalue $b$ is $_**$-eigenvalue of the matrix $a_1$.* |

PROOF OF THEOREM 18.4.13.   The equality

$$(18.4.43) \quad (a^j_{i\,s0} \otimes a^j_{i\,s1}) \circ v_i = (a^j_{0\,i} \otimes 1) \circ v_i = a^j_{0\,i} v_i$$

follows from the equality $\boxed{(18.4.41) \quad a^j_{i\,s0} \otimes a^j_{i\,s1} = a^j_{0\,i} \otimes 1}$ .
The equality
(18.4.44)

$$\begin{pmatrix} a_0{}^1_1 & ... & a_0{}^1_n \\ ... & ... & ... \\ a_0{}^n_1 & ... & a_0{}^n_n \end{pmatrix} *_* \begin{pmatrix} v_1 & ... & v_n \end{pmatrix} = \begin{pmatrix} a^1_{1\,s0} \otimes a^1_{1\,s1} & ... & a^n_{1\,s0} \otimes a^n_{1\,s1} \\ ... & ... & ... \\ a^1_{n\,s0} \otimes a^1_{n\,s1} & ... & a^n_{n\,s0} \otimes a^n_{n\,s1} \end{pmatrix} \circ_\circ \begin{pmatrix} v_1 & ... & v_n \end{pmatrix}$$

follows from the equality (18.4.43). The equality

$$(18.4.45) \quad \begin{pmatrix} a^1_{1\,s0} \otimes a^1_{1\,s1} & ... & a^n_{1\,s0} \otimes a^n_{1\,s1} \\ ... & ... & ... \\ a^1_{n\,s0} \otimes a^1_{n\,s1} & ... & a^n_{n\,s0} \otimes a^n_{n\,s1} \end{pmatrix} \circ_\circ \begin{pmatrix} v_1 & ... & v_n \end{pmatrix} = b \begin{pmatrix} v_1 & ... & v_n \end{pmatrix}$$

follows from the definition 18.4.9 of left $^\circ_\circ$-eigenvalue of the matrix $a$. The equality

$$(18.4.46) \qquad \begin{pmatrix} a_0{}^1_1 & ... & a_0{}^1_n \\ ... & ... & ... \\ a_0{}^n_1 & ... & a_0{}^n_n \end{pmatrix} {}^*_* \begin{pmatrix} v_1 & ... & v_n \end{pmatrix} = b \begin{pmatrix} v_1 & ... & v_n \end{pmatrix}$$

follows from equalities (18.4.44), (18.4.45). The theorem follows from the equality (18.4.46), from the definition 17.6.17 of $_*^*$-eigenvalue of the matrix $a_0$ and from the definition 17.6.23 of corresponding eigenvector $v$. □

Proof of theorem 18.4.14. The equality

$$(18.4.47) \qquad (a^j_{i\,s0} \otimes a^j_{i\,s1}) \circ v_i = (1 \otimes a_1{}^j_i) \circ v_i = v_i a_1{}^j_i$$

follows from the equality $\boxed{(18.4.42) \quad a^j_{i\,s0} \otimes a^j_{i\,s1} = 1 \otimes a_1{}^j_i}$.
The equality
$(18.4.48)$

$$\begin{pmatrix} v_1 & ... & v_n \end{pmatrix} {}^*_* \begin{pmatrix} a_1{}^1_1 & ... & a_1{}^1_n \\ ... & ... & ... \\ a_1{}^n_1 & ... & a_1{}^n_n \end{pmatrix} = \begin{pmatrix} a_1{}^1_{s0} \otimes a_1{}^1_{s1} & ... & a_1{}^n_{s0} \otimes a_1{}^n_{s1} \\ ... & ... & ... \\ a_n{}^1_{s0} \otimes a_n{}^1_{s1} & ... & a_n{}^n_{s0} \otimes a_n{}^n_{s1} \end{pmatrix} {}^\circ_\circ \begin{pmatrix} v_1 & ... & v_n \end{pmatrix}$$

follows from the equality (18.4.47). The equality

$$(18.4.49) \qquad \begin{pmatrix} a_1{}^1_{s0} \otimes a_1{}^1_{s1} & ... & a_1{}^n_{s0} \otimes a_1{}^n_{s1} \\ ... & ... & ... \\ a_n{}^1_{s1} \otimes a_n{}^1_{s1} & ... & a_n{}^n_{s0} \otimes a_n{}^n_{s1} \end{pmatrix} {}^\circ_\circ \begin{pmatrix} v_1 & ... & v_n \end{pmatrix} = \begin{pmatrix} v_1 & ... & v_n \end{pmatrix} b$$

follows from the definition 18.4.10 of right $^\circ_\circ$-eigenvalue of the matrix $a$. The equality

$$(18.4.50) \qquad \begin{pmatrix} v_1 & ... & v_n \end{pmatrix} {}^*_* \begin{pmatrix} a_1{}^1_1 & ... & a_1{}^1_n \\ ... & ... & ... \\ a_1{}^n_1 & ... & a_1{}^n_n \end{pmatrix} = \begin{pmatrix} v_1 & ... & v_n \end{pmatrix} b$$

follows from equalities (18.4.48), (18.4.49). The theorem follows from the equality (18.4.50), from the definition 17.6.14 of $^*_*$-eigenvalue of the matrix $a_1$ and from the definition 17.6.22 of corresponding eigenvector $v$. □

---

Theorem 18.4.15. *Let the row vector $v$ be eigenrow for left $^\circ_\circ$-eigenvalue $b$ of the matrix $a$. Let A-number $b$ satisfy the condition*

$$(18.4.51) \qquad b \in \bigcap_{i=1}^{n} \bigcap_{j=1}^{n} Z(A, a^j_{i\,s0})$$

*Then for any polynomial*

$$(18.4.52) \qquad p(x) = p_0 + p_1 x + ... + p_n x^n$$
$$p_0, ..., p_n \in D$$

Theorem 18.4.16. *Let the row vector $v$ be eigenrow for right $^\circ_\circ$-eigenvalue $b$ of the matrix $a$. Let A-number $b$ satisfy the condition*

$$(18.4.54) \qquad b \in \bigcap_{i=1}^{n} \bigcap_{j=1}^{n} Z(A, a^j_{i\,s1})$$

*Then for any polynomial*

$$(18.4.55) \qquad p(x) = p_0 + p_1 x + ... + p_n x^n$$
$$p_0, ..., p_n \in D$$

*the row vector*

$$cv = \begin{pmatrix} cv_1 & ... & cv_n \end{pmatrix}$$

$$c = p(b)$$

*is also eigenrow for left $^\circ$-eigenvalue $b$.*

PROOF. According to the theorem 5.1.25, the equality

$$(18.4.53) \quad \begin{aligned} &(a^j_{i\,s0} \otimes a^j_{i\,s1}) \circ (cv_i) \\ &= a^j_{i\,s0} cv_i a^j_{i\,s1} \\ &= c(a^j_{i\,s0} v_i a^j_{i\,s1}) \\ &= cbv_j = b(cv_j) \end{aligned}$$

follows from the condition (18.4.51). According to the definition 18.4.9, the equality (18.4.53) implies that the row vector $cv$ is eigenrow for left $^\circ$-eigenvalue $b$ of the matrix $a$. □

*the row vector*

$$vc = \begin{pmatrix} v_1 c & ... & v_n c \end{pmatrix}$$

$$c = p(b)$$

*is also eigenrow for right $^\circ$-eigenvalue $b$.*

PROOF. According to the theorem 5.1.25, the equality

$$(18.4.56) \quad \begin{aligned} &(a^j_{i\,s0} \otimes a^j_{i\,s1}) \circ (v_i c) \\ &= a^j_{i\,s0} v_i c a^j_{i\,s1} \\ &= (a^j_{i\,s0} v_i a^j_{i\,s1}) c \\ &= v_j bc = (v_f c) b \end{aligned}$$

follows from the condition (18.4.54). According to the definition 18.4.10, the equality (18.4.56) implies that the row vector $vc$ is eigenrow for right $^\circ$-eigenvalue $b$ of the matrix $a$. □

# CHAPTER 20

# Index

# Special Symbols and Notations

**Александр Клейн**

# Введение в некоммутативную алгебру

## *Том 1*
## *Алгебра с делением*

Aleks_Kleyn@MailAPS.org
http://AleksKleyn.dyndns-home.com:4080/
http://arxiv.org/a/kleyn_a_1
http://AleksKleyn.blogspot.com/

2020 *Mathematics Subject Classification.* Primary 16D10; Secondary 16K20,16K40

Аннотация. Том 1 монографии "Введение в некоммутативную алгебру" посвящён изучению алгебры над коммутативным кольцом. Основные темы, которые я рассмотрел в этом томе:

- определение модуля и алгебры над коммутативным кольцом
- линейное отображение алгебры над коммутативным кольцом
- векторное пространство над ассоциативной алгброй с делением
- система линейных уравнений над ассоциативной алгброй с делением
- гомоморфизм векторного пространства над ассоциативной алгброй с делением и теорию ковариантности
- собственное значение и собственный вектор над ассоциативной алгброй с делением

'

# Оглавление

# Предисловие

## 1.1. Об этой книге

Я планировал написать несколько разных статей по линейной алгебре. Я рассмотрел эти статьи в последующих разделах этой главы. Но в какой-то момент я понял, что я готов собрать эти статьи в одну книгу и рассмотреть все задачи, которыми я занимался в области линейной алгебры. Однако ответ на один вопрос порождает новые вопросы, которые тоже должны быть отвечены. Например, если я хочу рассмотреть задачу о собственных значениях, я должен рассмотреть теорию алгебраических уравнений.

Список тем постепенно увеличивался. Мне стало ясно, что я хочу дать широкий обзор теории некоммутативной алгебры. Я понимаю что вряд ли дам полное описание, тема обширная. Но и те вопросы, что я хочу затронуть, займут несколько томов.

## 1.2. Линейное отображение

Я начну со следуещего утверждения.

**Теорема 1.2.1.** *Если $D$-алгебра $A$ с делением удовлетворяет следующим условиям*

1.2.1.1: *$D$-алгебра $A$ коммутативна и ассоциативна*

1.2.1.2: *поле $D$ является центром $D$-алгебры $A$*

*то*

$$(1.2.1) \qquad \mathcal{L}(D; A \to A) = A$$

Мы рассмотрим доказательство теоремы 1.2.1 в томе 3. Если хотя бы одно условие теоремы 1.2.1 нарушено, то равенство (1.2.1) невозможно.

- Множество линейных отображений $\mathcal{L}(R; C \to C)$ поля комплексных чисел является $C$-векторным пространством размерности 2. Это связано с тем, что произведение в поле комплексных чисел коммутативно и поле комплексных чисел $C$ является векторным пространством размерности 2 над полем действительных чисел $R$.
- Множество линейных отображений $\mathcal{L}(R; H \to H)$ алгебры кватернионов является $H$-векторным пространством размерности 4. Это связано с тем, что произведение в алгебре кватернионов не коммутативно

Изучая алгебру кватернионов, в статьях [10, 11] я рассмотрел множество отображений, обладающих некоторой симметрией. Это положило начало нового исследования. Одним из мотивов этого исследования было желание найти аналоги утверждениям из теории функции комплексного переменного.

Я сделал несколько попыток найти множество отображений, которые позволили бы мне найти структуру множества $\mathcal{L}(R; H \to H)$, подобную структуре множества $\mathcal{L}(R; C \to C)$. Но не считаю ни одну из этих попыток успешной. Тем не менее, в процессе этого исследования я понял, что множество $\mathcal{L}(R; H \to H)$ является левым или правым $H$-векторным пространством размерности 4.

Я рассмотрел множество линейных отображений $\mathcal{L}(R; O \to O)$ алгебры октонионов аналогичным образом. Однако я понял, что подобное утверждение верно для любой $D$-алгебры с делением.

Это не значит, что рассмотрение алгебры $\mathcal{L}(D; A \to A)$ как $A$-векторного пространства может заменить представление отображения как тензора.

### 1.3. Система линейных уравнений

Пусть $V$ - левое векторное пространство над $D$-алгеброй $A$. Пусть $\overline{\overline{e}}$ - базис левого $A$-векторного пространства $V$. В этом случае мы можем рассмотреть две различные задачи, условие которых можно записать в виде системы линейных уравнений.

- Пусть $V$ - левое $A$-векторное пространство столбцов. Тогда векторы базиса $\overline{\overline{e}}$ мы запишем в виде строки матрицы

$$e = \begin{pmatrix} e_1 & ... & e_n \end{pmatrix}$$

и координаты вектора $\overline{w}$ относительно базиса $\overline{\overline{e}}$ мы запишем в виде столбца матрицы

$$w = \begin{pmatrix} w^1 \\ ... \\ w^n \end{pmatrix}$$

Представим множество векторов $v_1, ..., v_m$ в виде строки матрицы

$$v = \begin{pmatrix} v_1 & ... & v_m \end{pmatrix}$$

Мы запишем координаты вектора $v_i$ относительно базиса $\overline{\overline{e}}$ мы запишем в виде столбца матрицы

$$v_i = \begin{pmatrix} v_i^1 \\ ... \\ v_i^n \end{pmatrix}$$

Чтобы ответить на вопрос, можно ли представить вектор $\overline{w}$ как линейную комбинацию множества векторов $v_1, ..., v_m$, мы должны решить систему линейных уравнений

$$c^1 v_1^1 + ... + c^m v_m^1 = w^1$$
$$... = ...$$
$$c^1 v_1^n + ... + c^m v_m^n = w^n$$

относительно множества неизвестных скаляров $c^1, ..., c^m$.

- Пусть $A$ - ассоциативная $D$-алгебра с делением.   Пусть

$$f : V \to W$$

линейное отображение $A$-векторного пространства $V$ в $A$-векторное пространство $W$. Отображение координатной матрицы вектора $v \in V$ относительно базиса $\overline{\overline{e}}_V$ в координатную матрицу вектора $f \circ v \in W$ относительно базиса $\overline{\overline{e}}_W$ имеет вид

$$f_{11}^1 v^1 f_{12}^1 + ... + f_{m1}^1 v^m f_{m2}^1 = (f \circ v)^1$$
$$... = ...$$
$$f_{11}^n v^1 f_{12}^n + ... + f_{m1}^n v^m f_{m2}^n = (f \circ v)^n$$

и не зависит от того, рассматриваем мы левое $A$-векторное пространство или правое $A$-векторное пространство. Чтобы найти координатную матрицу вектора $v \in V$, $f \circ v = w$,  мы должны решить систему линейных уравнений

(1.3.1)
$$f_{11}^1 v^1 f_{12}^1 + ... + f_{m1}^1 v^m f_{m2}^1 = w^1$$
$$... = ...$$
$$f_{11}^n v^1 f_{12}^n + ... + f_{m1}^n v^m f_{m2}^n = w^n$$

### 1.4. Собственное значение матрицы

**1.4.1. Собственное значение в коммутативной алгебре.** Я начну с определения собственного значения матрицы в коммутативной алгебре. [1.1]
Пусть $a$ - $n \times n$ матрица, элементы которой принадлежат полю $F$. Многочлен

$$\det(a - bE_n)$$

где $E_n$ - $n \times n$ единичная матрица, называется характерисическим многочленом матрицы $a$. Корни характеристического многочлена называются собственными значениями матрицы $a$. Следовательно, собственные значения матрицы $a$ удовлетворяют уравнению

$$\det(a - bE_n) = 0$$

и система линейных уравнений

(1.4.1)
$$(a - bE_n)x = 0$$

имеет нетривиальное решение, если $b$ - собственное значение матрицы $a$.

Если матрицы $a$ и $a_1$ подобны

(1.4.2)
$$a_1 = f^{-1}af$$

то из равенства

(1.4.3)
$$a_1 - bE_n = f^{-1}af - bE_n = f^{-1}af - f^{-1}bE_nf = f^{-1}(a - bE_n)f$$

следует, что матрицы $a$ and $a_1$ имеют одинаковые характеристические многочлены, и, следовательно, одинаковые характеристические корни.

Пусть $V$ - векторное пространство над полем $F$. Пусть $\overline{\overline{e}}$, $\overline{\overline{e}}_1$ - базисы векторного пространства $V$ и $f$ - координатная матрица базиса $\overline{\overline{e}}_1$ относительно базиса $\overline{\overline{e}}$

(1.4.4)
$$\overline{\overline{e}}_1 = f\overline{\overline{e}}$$

---

[1.1] Смотри также определение на странице [3]-207.

Пусть $a$ - матрица эндоморфизма $\overline{a}$ относительно базиса $\overline{e}$. Тогда матрица $a_1$, определённая равенством (1.4.2), является матрицей эндоморфизма $\overline{a}$ относительно базиса $\overline{e}_1$. Если мы отождествляем эндоморфизм и его матрицу относительно заданного базиса, то мы будем говорить, что собственные значения эндоморфизма - это собственные значения его матрицы относительно заданного базиса. Следовательно, собственные значения эндоморфизма не зависят от выбора базиса.

Если мы положим

$$(1.4.5) \qquad\qquad \overline{x} = x^i e_i$$

то равенство (1.4.1) принимает вид

$$(1.4.6) \qquad\qquad \overline{a}\overline{x} = b\overline{x}$$

$$(1.4.7) \qquad\qquad ax = bx$$

относительно базиса $\overline{e}$. Пусть $f$ - координатная матрица базиса $\overline{e}_1$ относительно базиса $\overline{e}$ $\boxed{(1.4.4) \quad \overline{e}_1 = f\overline{e}}$ . Тогда координатная матрица $x_1$ вектора $\overline{x}$ относительно базиса $\overline{e}_1$ имеет вид

$$(1.4.8) \qquad\qquad x_1 = f^{-1}x$$

Равенство

$$(1.4.9) \qquad a_1x_1 = f^{-1}aff^{-1}x = f^{-1}ax = f^{-1}bx = bf^{-1}x = bx_1$$

является следствием равенств (1.4.2), (1.4.7), (1.4.8). Следовательно, равенство (1.4.6) не зависит от выбора базиса, и мы будем говорить, что собственный вектор эндоморфизма - это собственный вектор его матрицы относительно заданного базиса.

Множество собственных значений эндоморфизма $\overline{a}$ с учётом их кратности в характеристическом многочлене называется спектром эндоморфизма $\overline{a}$.

**1.4.2. Собственное значение в некоммутативной алгебре.** Пусть $A$ - некоммутативная $D$-алгебра $A$ с делением. Если мы рассматриваем векторное пространство над $D$-алгеброй $A$, то мы хотим сохранить определение собственного числа и собственного вектора, рассмотренное в случае векторного пространства над коммутативной алгеброй.

Пусть $V$ - левое $A$-векторное пространство столбцов. Поскольку произведение некоммутативно, возникает вопрос. Должны мы пользоваться равенством

$$(1.4.10) \qquad\qquad \overline{a} \circ \overline{x} = b\overline{x}$$

или равенством

$$(1.4.11) \qquad\qquad \overline{a} \circ \overline{x} = \overline{x}b$$

Равенство (1.4.10) выглядит естественно, так как в правой части записанно левостороннее произведение вектора на скаляр. Однако, если $cb \neq bc$, $b$, $c \in A$, то

$$(1.4.12) \qquad\qquad \overline{a} \circ (c\overline{x}) \neq b(c\overline{x})$$

Утверждение 1.4.1. *Из утверждения (1.4.12) следует, что множество векторов $\overline{x}$, удовлетворяющих равенству (1.4.10), не порождает подпространство векторного пространства.* ⊙

Если мы запишем равенство (1.4.10) с помощью матриц, то мы получим равенство

$$(1.4.13) \qquad x^*{}_* a = bx$$

Равенство (1.4.13) является системой линейных уравнений относительно множества неизвестных $x$. Однако мы не можем сказать, что эта система линейных уравнений вырождена и имеет бесконечно много решений. Но это утверждение является очевидным следствием утверждения 1.4.1.

Создаётся впечатление, что равенство (1.4.10) не годится для определения собственного значения. Однако мы вернёмся к этому определению в разделе 1.5.

Рассмотрим равенство (1.4.11). Пусть $a$ - матрица эндоморфизма $\overline{a}$ относительно базиса $\overline{\overline{e}}$ левого $A$-векторного пространства столбцов. Пусть $a$ - матрица вектора $\overline{x}$ относительно базиса $\overline{\overline{e}}$. Система линейных уравнений

$$(1.4.14) \qquad x^*{}_* a = xb$$

является следствием уравнения (1.4.11). Система линейных уравнений (1.4.14) имеет нетривиальное решение, если матрица $a - bE_n$ $^*{}_*$-вырождена. В этом случае, система линейных уравнений (1.4.14) имеет бесконечно много решений.

На первый взгляд, рассмотренная модель соответствует модели, рассмотренной в разделе 1.4.1. Однако, если мы выберем другой базис $\overline{\overline{e}}_1$

$$(1.4.15) \qquad \overline{\overline{e}}_1 = g^*{}_* \overline{\overline{e}}$$

то изменятся матрица $a$

$$(1.4.16) \qquad a_1 = g^*{}_* a^*{}_* g^{-1^*}{}_*$$

$$(1.4.17) \qquad a = g^{-1^*}{}_*{}^*{}_* a_1{}^*{}_* g$$

и матрица $x$

$$(1.4.18) \qquad x_1 = x^*{}_* g^{-1^*}{}_*$$

$$(1.4.19) \qquad x = x_1{}^*{}_* g$$

Равенство

$$(1.4.20) \qquad x_1{}^*{}_* g^*{}_* g^{-1^*}{}_*{}^*{}_* a_1{}^*{}_* g = x_1{}^*{}_* gb$$

является следствием равенств (1.4.14), (1.4.17), (1.4.19). Из равенства (1.4.20) следует, что вектор-столбец $x_1$ не является решением системы линейных уравнений (1.4.14). Это противоречит утверждению, что равенство (1.4.11) не зависит от выбора базиса.

Равенство (1.4.11) требует более внимательного анализа. В правой части равенства (1.4.11) записана операция, которая не определена в левом $A$-векторном пространстве, а именно, правостороннее умножение вектора на скаляр. Однако мы можем записать систему линейных уравнений (1.4.14) иначе

$$(1.4.21) \qquad x^*{}_* a = x^*{}_* (Eb)$$

где $E$ - единичная матрица. В правой части равенства (1.4.21) также как и в левой части записано $^*{}_*$-произведение вектор-столбца на матрицу эндоморфизма. Следовательно, левая и правая части имеют один и тот же закон преобразования при замене базиса.

По аналогии с векторным пространством над полем, эндоморфизм $\overline{\overline{e}}b$ который имеет матрицу $Eb$ относительно базиса $\overline{\overline{e}}$, называется преобразованием подобия относительно базиса $\overline{e}$. Следовательно, мы должны записать равенство (1.4.11) в виде

$$(1.4.22) \qquad\qquad a \circ x = \overline{\overline{e}}b \circ x$$

Равенство (1.4.22) является определением собственного значения $b$ и собственного вектора $x$ эндоморфизма $a$. Если мы рассматриваем матрицу эндоморфизма $a$ относительно базиса $\overline{e}$ левого $A$-векторного пространства столбцов, то равенство (1.4.14) является определением $^*_*$-собственного значения $b$ и собственного вектора $x$, соответствующего $^*_*$-собственному значению $b$.

## 1.5. Левые и правые собственные значения

Чтобы лучше разобраться с природой собственных значений, рассмотрим следующую теорему.

> ТЕОРЕМА 1.5.1. *Пусть $D$ - поле $V$ - $D$-векторное пространство. Если эндоморфизм имеет $k$ различных собственных значений*[1.2] $b(1)$, ..., $b(k)$, *то собственные вектора $v(1)$, ..., $v(k)$, относящиеся к различным собственным значениям, линейно независимы.*

ДОКАЗАТЕЛЬСТВО. Мы будем доказывать теорему индукцией по $k$.
Поскольку собственный вектор отличен от нуля, теорема верна для $k = 1$.

УТВЕРЖДЕНИЕ 1.5.2. *Пусть теорема верна для $k = m - 1$.*            ⊙

Пусть $b(1)$, ..., $b(m-1)$, $b(m)$ - различные собственные значения и $v(1)$, ..., $v(m-1)$, $v(m)$ - соответствующие собственные векторы

$$(1.5.1) \qquad\qquad f \circ v(i) = b(i)v(i)$$

$$i \neq j \Rightarrow b(i) \neq b(j)$$

Пусть утверждение 1.5.3 верно.

УТВЕРЖДЕНИЕ 1.5.3. *Существует линейная зависимость*

$$(1.5.2) \qquad a(1)v(1) + ... + a(m-1)v(m-1) + a(m)v(m) = 0$$

*где $a(1) \neq 0$.*                                                  ⊙

Равенство

$$(1.5.3) \quad a(1)(f \circ v(1)) + ... + a(m-1)(f \circ v(m-1)) + a(m)(f \circ v(m)) = 0$$

является следствием равенства (1.5.2). Равенство

$$(1.5.4) \quad \begin{aligned} & a(1)b(1)v(1) + ... \\ & + a(m-1)b(m-1)v(m-1) + a(m)b(m)v(m) = 0 \end{aligned}$$

является следствием равенства (1.5.1), (1.5.3). Равенство

$$(1.5.5) \quad \begin{aligned} & a(1)(b(1) - b(m))v(1) + ... \\ & + a(m-1)(b(m-1) - b(m))v(m-1) = 0 \end{aligned}$$

---

[1.2] Смотри также теорему на странице [3]-210.

является следствием равенства (1.5.2), (1.5.4). Так как

$$a(1)(b(1) - b(m)) \neq 0$$

то равенство (1.5.5) является линейной зависимостью между векторами $v(1)$, ..., $v(m-1)$. Это утверждение противоречит утверждению 1.5.2. Следовательно, утверждение 1.5.3 не верно и векторы $v(1)$, ..., $v(m-1)$, $v(m)$ линейно независимы. $\square$

Пусть $D$-векторное пространство имеет размерность $n$. Из теоремы 1.5.1 следует, что если $k = n$, то множество соответствующих собственных векторов является базисом $D$-векторного пространства $V$.[1.3]

Пусть $\overline{e}$ - базис $D$-векторного пространства $V$. Пусть $a$ - матрица эндоморфизма $\overline{a}$ относительно базиса $\overline{\overline{e}}$ и

$$x(k) = \begin{pmatrix} x(k)^1 \\ ... \\ x(k)^n \end{pmatrix}$$

координаты собственного вектора $x(k)$ относительно базиса $\overline{\overline{e}}$. Соберём множество координат собственных векторов $v(1)$, ..., $v(n-1)$, $v(n)$ в матрицу

$$x = \begin{pmatrix} \begin{pmatrix} x(1)^1 \\ ... \\ x(1)^n \end{pmatrix} & ... & \begin{pmatrix} x(n)^1 \\ ... \\ x(n)^n \end{pmatrix} \end{pmatrix} = \begin{pmatrix} x(1)^1 & ... & x(n)^1 \\ ... & ... & ... \\ x(1)^n & ... & x(n)^n \end{pmatrix}$$

Тогда мы можем записать равенство

$$\boxed{(1.4.6) \quad \overline{ax} = b\overline{x}}$$

используя произведение матриц

$$(1.5.6) \quad \begin{aligned} & \begin{pmatrix} a_1^1 & ... & a_n^1 \\ ... & ... & ... \\ a_1^n & ... & a_n^n \end{pmatrix} \begin{pmatrix} x(1)^1 & ... & x(n)^1 \\ ... & ... & ... \\ x(1)^n & ... & x(n)^n \end{pmatrix} \\ & = \begin{pmatrix} x(1)^1 & ... & x(n)^1 \\ ... & ... & ... \\ x(1)^n & ... & x(n)^n \end{pmatrix} \begin{pmatrix} b(1) & ... & 0 \\ ... & ... & ... \\ 0 & ... & b(n) \end{pmatrix} \end{aligned}$$

Так как множество векторов $v(1)$, ..., $v(n-1)$, $v(n)$ является базисом $D$-векторного пространства $V$, то матрица $x$ - не вырождена. Следовательно,

---

[1.3] Если кратность собственного значения больше чем 1, то размерность $D$-векторного пространства, порождённого соответствующими собственными векторами, может быть меньше кратности собственного значения.

равенство

$$(1.5.7) \quad \begin{pmatrix} x(1)^1 & ... & x(n)^1 \\ ... & ... & ... \\ x(1)^n & ... & x(n)^n \end{pmatrix}^{-1} \begin{pmatrix} a_1^1 & ... & a_n^1 \\ ... & ... & ... \\ a_1^n & ... & a_n^n \end{pmatrix} \begin{pmatrix} x(1)^1 & ... & x(n)^1 \\ ... & ... & ... \\ x(1)^n & ... & x(n)^n \end{pmatrix}$$

$$= \begin{pmatrix} b(1) & ... & 0 \\ ... & ... & ... \\ 0 & ... & b(n) \end{pmatrix}$$

является следствием равенства (1.5.6). Из равенства (1.5.7) следует, что матрица $a$ подобна диагональной матрице

$$b = \mathrm{diag}(b(1), ..., b(n))$$

Пусть $V$ - левое $A$-векторное пространство столбцов. Тогда равенство (1.5.6) приобретает вид

$$(1.5.8) \quad \begin{pmatrix} x(1)^1 & ... & x(n)^1 \\ ... & ... & ... \\ x(1)^n & ... & x(n)^n \end{pmatrix} {}^*_* \begin{pmatrix} a_1^1 & ... & a_n^1 \\ ... & ... & ... \\ a_1^n & ... & a_n^n \end{pmatrix}$$

$$= \begin{pmatrix} b(1) & ... & 0 \\ ... & ... & ... \\ 0 & ... & b(n) \end{pmatrix} {}^*_* \begin{pmatrix} x(1)^1 & ... & x(n)^1 \\ ... & ... & ... \\ x(1)^n & ... & x(n)^n \end{pmatrix}$$

Если мы в равенстве (1.5.8) выделим столбец $x(k)$, то мы получим равенство

$$(1.5.9) \quad \begin{pmatrix} x(k)^1 \\ ... \\ x(k)^n \end{pmatrix} {}^*_* \begin{pmatrix} a_1^1 & ... & a_n^1 \\ ... & ... & ... \\ a_1^n & ... & a_n^n \end{pmatrix} = b(k) \begin{pmatrix} x(k)^1 \\ ... \\ x(k)^n \end{pmatrix}$$

Равенство (1.5.9) с точностью до обозначений совпадает с равенством

$$\boxed{(1.4.10) \quad \overline{a} \circ \overline{x} = b\overline{x}}$$

Таким образом, мы получаем определение левого ${}^*_*$-собственного значения, которое отличается от ${}^*_*$-собственного значения, рассмотренного в разделе 1.4.2. Если мы аналогичным образом рассмотрим правое $A$-векторное пространство строк, то мы получаем определение правого ${}^*_*$-собственного значения. [1.4]

---

[1.4] В разделе [20]-8.2, профессор Кон рассматривает определение левого и правого собственных значений $n \times n$ матрицы. Кон по традиции рассматривает произведение строки на столбец. Поскольку существует две операции произведения матриц, я слегка изменил определение Кона.

## 1.6. Система дифференциальных уравнений

В коммутативной алгебре есть ещё одна задача, для решения которой требуется найти собственное значение матрицы. Мы ищем решение системы дифференциальных уравнений

$$\frac{dx^1}{dt} = a^1_1 x^1 + ... + a^1_n x^n$$

(1.6.1)                    ......

$$\frac{dx^n}{dt} = a^n_1 x^1 + ... + a^n_n x^n$$

в виде

$$\begin{pmatrix} x^1 \\ ... \\ x^n \end{pmatrix} = \begin{pmatrix} c^1 \\ ... \\ c^n \end{pmatrix} e^{bt}$$

В этом случае число $b$ является собственным значением матрицы

$$a = \begin{pmatrix} a^1_1 & ... & a^1_n \\ ... & ... & ... \\ a^n_1 & ... & a^n_n \end{pmatrix}$$

и вектор

$$c = \begin{pmatrix} c^1 \\ ... \\ c^n \end{pmatrix}$$

является собственным вектором матрицы $a$, соответствующим собственному значению $b$.

В случае некоммутативной $D$-алгебры $A$, мы можем записать систему дифференциальных уравнений (1.6.1) в виде

(1.6.2)                    $$\frac{dx}{dt} = a_* {}^* x$$

В этом случае решение системы дифференциальных уравнений (1.6.2) может иметь одно из следующих представлений.

- Мы можем записать решение в виде

(1.6.3)          $$x = e^{bt} c = \begin{pmatrix} e^{bt} c^1 \\ ... \\ e^{bt} c^n \end{pmatrix} \qquad c = \begin{pmatrix} c^1 \\ ... \\ c^n \end{pmatrix}$$

где $A$-число $b$ является $_*{}^*$-собственным значением матрицы $a$ и вектор-столбец $c$ является собственным столбцом, соответствующим $_*{}^*$-собственному числу $b$.

- Мы можем записать решение в виде

$$(1.6.4) \qquad x = ce^{bt} = \begin{pmatrix} c^1 e^{bt} \\ ... \\ c^n e^{bt} \end{pmatrix} \qquad c = \begin{pmatrix} c^1 \\ ... \\ c^n \end{pmatrix}$$

где $A$-число $b$ является правым $_*{}^*$-собственным значением матрицы $a$ и вектор-столбец $c$ является собственным столбцом, соответствующим правому $_*{}^*$-собственному числу $b$.

В некоммутативной $D$-алгебры $A$, мы можем рассмотреть более общий вид системы линейных дифференциальных уравнений

$$(1.6.5) \qquad \begin{aligned} \frac{dx^1}{dt} &= a^1_{1\,s0} x^1 a^1_{1\,s1} + ... + a^1_{n\,s0} x^n a^1_{n\,s1} \\ &\quad ...... \\ \frac{dx^n}{dt} &= a^n_{1\,s0} x^1 a^n_{1\,s1} + ... + a^n_{n\,s0} x^n a^n_{n\,s1} \end{aligned}$$

Мы можем записать систему дифференциальных уравнений (1.6.5), пользуясь произведением матриц

$$(1.6.6) \qquad \begin{pmatrix} \dfrac{dx^1}{dt} \\ ... \\ \dfrac{dx^n}{dt} \end{pmatrix} = \begin{pmatrix} a^1_{1\,s0} \otimes a^1_{1\,s1} & ... & a^1_{n\,s0} \otimes a^1_{n\,s1} \\ ... & ... & ... \\ a^n_{1\,s0} \otimes a^n_{1\,s1} & ... & a^n_{n\,s0} \otimes a^n_{n\,s1} \end{pmatrix} \circ^\circ \begin{pmatrix} x^1 \\ ... \\ x^n \end{pmatrix}$$

Решение системы дифференциальных уравнений (1.6.6) может иметь одно из следующих представлений.

- Мы можем записать решение в виде

$$(1.6.7) \qquad x = Ce^{bt} c = Ce^{bt} \begin{pmatrix} c^1 \\ ... \\ c^n \end{pmatrix}$$

где
- $A$-число $b$ является левым $_\circ{}^\circ$-собственным значением матрицы

$$a = \begin{pmatrix} a^1_{1\,s0} \otimes a^1_{1\,s1} & ... & a^1_{n\,s0} \otimes a^1_{n\,s1} \\ ... & ... & ... \\ a^n_{1\,s0} \otimes a^n_{1\,s1} & ... & a^n_{n\,s0} \otimes a^n_{n\,s1} \end{pmatrix}$$

и столбец $c$ является собственным столбцом матрицы $a$, соответствующим левому $_\circ{}^\circ$-собственному значению $b$.
- $A$-число $b$ и $A$-число $C$ удовлетворяют условию

$$(1.6.8) \qquad b \in \bigcap_{i=1}^{n} \bigcap_{j=1}^{n} Z(A, a^i_{j\,s1})$$

$$(1.6.9) \qquad C \in \bigcap_{i=1}^{n} \bigcap_{j=1}^{n} Z(A, a_{j\,s1}^{i})$$

- Мы можем записать решение в виде

$$(1.6.10) \qquad x = c e^{bt} C = \begin{pmatrix} c^1 \\ ... \\ c^n \end{pmatrix} e^{bt} C$$

где
- $A$-число $b$ является правым $\circ^\circ$-собственным значением матрицы

$$a = \begin{pmatrix} a_{1\,s0}^{1} \otimes a_{1\,s1}^{1} & ... & a_{n\,s0}^{1} \otimes a_{n\,s1}^{1} \\ ... & ... & ... \\ a_{1\,s0}^{n} \otimes a_{1\,s1}^{n} & ... & a_{n\,s0}^{n} \otimes a_{n\,s1}^{n} \end{pmatrix}$$

и столбец $c$ является собственным столбцом матрицы $a$, соответствующим правому $\circ^\circ$-собственному значению $b$.
- $A$-число $b$ и $A$-число $C$ удовлетворяют условию

$$(1.6.11) \qquad b \in \bigcap_{i=1}^{n} \bigcap_{j=1}^{n} Z(A, a_{j\,s2}^{i})$$

$$(1.6.12) \qquad C \in \bigcap_{i=1}^{n} \bigcap_{j=1}^{n} Z(A, a_{j\,s2}^{i})$$

### 1.7. Многочлен

В коммутативной алгебре теория многочленов имеет много общего с теорией чисел: делимость, наибольший общий делитель, наименьшее общее кратное, разложение на неприводимые множители. Когда мы изучаем многочлены в некоммутативной алгебре, мы хотим максимально сохранить утверждения из коммутативной алгебры.

В книге [20] на странице 48, профессор Кон писал, что изучать многочлены над не коммутативным полем - сложная задача. Поэтому профессор Оре рассмотрел кольцо правых косых многочленов [1.5] где Оре определил преобразование одночлена

$$(1.7.1) \qquad at = ta^{\alpha} + a^{\delta}$$

Отображение

$$\alpha : a \in A \to a^{\alpha} \in A$$

является эндоморфизмом $D$-алгебры $A$ и отображение

$$\delta : a \in A \to a^{\delta} \in A$$

является $\alpha$-дифференцированием $D$-алгебры $A$

$$(1.7.2) \qquad (a + b)^{\delta} = a^{\delta} + b^{\delta} \qquad (ab)^{\delta} = a^{\delta} b^{\alpha} + a b^{\delta}$$

---

[1.5] Смотри также определение на страницах [20]-(48 - 49).

В этом случае все коэффициенты многочлена могут быть записаны справа. Подобные многочлены называются правосторонними многочленами.

Аналогично мы можем рассматривать преобразование одночлена

$$(1.7.3) \qquad ta = a^\alpha t + a^\delta$$

В этом случае все коэффициенты многочлена могут быть записаны слева. Подобные многочлены называются левосторонними многочленами. [1.6]

Точка зрения, предложенная Оре, позволила получить качественную картину теории полиномов в некоммутативной алгебре. Однако существуют алгебры где рассматриваемая картина не полна.

> **Теорема 1.7.1.** *Пусть $A$ - $D$-алгебра с делением. Пусть $a \in A \otimes A$. Пусть уравнение*
>
> $$(1.7.4) \qquad a \circ x + b = 0$$
>
> *имеет бесконечно много корней, но не любое $A$-число является корнем уравнения* (1.7.4). *Тогда кольцо многочленов $A[x]$ не является кольцом Оре.*

Доказательство. Допустим, что в $D$-алгебре $A$ определено преобразование $\boxed{(1.7.1) \quad at = ta^\alpha + a^\delta}$. Пусть

$$(1.7.5) \qquad a = a_{s0} \otimes a_{s1}$$

Равенство

$$(1.7.6) \qquad \begin{aligned} (a_{s0} \otimes a_{s1}) \circ x + b &= a_{s0} x a_{s1} + b = (x a_{s0}^\alpha + a_{s0}^\delta) a_{s1} + b \\ &= x a_{s0}^\alpha a_{s1} + a_{s0}^\delta a_{s1} + b = 0 \end{aligned}$$

является следствием равенств (1.7.1), (1.7.4), (1.7.5). Положим

$$A = a_{s0}^\alpha a_{s1} \qquad B = a_{s0}^\delta a_{s1} + b$$

Уравнение

$$(1.7.7) \qquad xA + B = 0$$

имеет теже корни, что и уравнение (1.7.4).

- Если $A \neq 0$, то уравнение (1.7.7) имеет единственный корень.
- Если $A = 0$, $B \neq 0$, то уравнение (1.7.7) не имеет корней.
- Если $A = 0$, $B = 0$, то любое $A$-число является корнем уравнения (1.7.7).

Все три варианта противоречат условию теоремы. Следовательно, преобразование $\boxed{(1.7.1) \quad at = ta^\alpha + a^\delta}$ не определено в $D$-алгебре $A$. $\qquad \square$

> **Теорема 1.7.2.** *Общее решение линейного уравнения*
>
> $$(1.7.8) \qquad ix - xi = k$$
>
> *в алгебре кватернионов имеет вид*
>
> $$(1.7.9) \qquad x = C_1 + C_2 i + \frac{1}{2} j \qquad C_1, C_2 \in R$$

---

[1.6] Смотри также определение на странице [16]-3.

Доказательство. Тензор

$$i \otimes 1 - 1 \otimes i \in H^{2\otimes}$$

соответствует матрице

$$(1.7.10) \quad \begin{pmatrix} 0 & -1 & 0 & 0 \\ 1 & 0 & 0 & 0 \\ 0 & 0 & 0 & -1 \\ 0 & 0 & 1 & 0 \end{pmatrix} - \begin{pmatrix} 0 & -1 & 0 & 0 \\ 1 & 0 & 0 & 0 \\ 0 & 0 & 0 & 1 \\ 0 & 0 & -1 & 0 \end{pmatrix} = \begin{pmatrix} 0 & 0 & 0 & 0 \\ 0 & 0 & 0 & 0 \\ 0 & 0 & 0 & -2 \\ 0 & 0 & 2 & 0 \end{pmatrix}$$

$\square$

Теорема 1.7.3. *Кольцо многочленов $H[x]$ в алгебре кватернионов не является кольцом Оре.*

Доказательство. Теорема является следствием теорем 1.7.1, 1.7.2.    $\square$

В разделе [19]-6.4, профессор Кон рассмотрел определение левостороннего многочлена, не опираясь на определение кольца Оре.

Пусть $A$ - ассоциативная алгебра.

На странице [19]-148, профессор Кон рассмотрел множество $A[x]^*$ левосторонних многочленов, где левосторонний многочлен определён как формальное выражение

$$(1.7.11) \quad f = a_0 + a_1 x + ... + a_n x^n$$

$$(1.7.12) \quad g = b_0 + b_1 x + ... + b_m x^m$$

где $a_0, ..., a_n, b_0, ..., b_m$, - $A$-числа и $x$ - переменная.

Левосторонние многочлены $f$ и $g$ равны, если $n = m$, $a_0 = b_0$, ..., $a_n = b_n$.

Мы также определим на множестве $A[x]^*$ сумму [1.7] левосторонних многочленов $f$ и $g$

$$(1.7.13) \quad f + g = (a_0 + b_0) + (a_1 + b_1) x + ... + (a_n + b_n) x^n$$

и произведение

$$(1.7.14) \quad f * g = a_0 b_0 + (a_0 b_1 + a_1 b_0) x + (a_0 b_2 + a_1 b_1 + a_2 b_0) x^2 + ... + a_n b_m x^{n+m}$$

Нетрудно убедиться, что множество $A[x]^*$ является $D$-алгеброй. Очевидно, что левосторонний многочлен $f$ является многочленом $f$. Однако произведение левосторонних многочленов $f$ и $g$, вообще говоря, не равно произведению многочленов $f$ и $g$.

Пример 1.7.4. *Положим*

$$(1.7.15) \quad f = 1 + ix$$

$$(1.7.16) \quad g = k + jx$$

*Равенства*

$$(1.7.17) \quad f * g = k + (j + ik) x + ij x^2 = k + kx^2$$

---

[1.7] Если $m < n$, например, то мы положим $b_{m+1} = ... = b_n = 0$.

$$(1.7.18) \quad fg = (1+ix)(k+jx) = (k+jx) + ix(k+jx) = k + jx + ixk + ixjx$$

*являются следствием равенств* (1.7.14), (1.7.15), (1.7.16). *Пусть* $x = k$. *Равенство*

$$(1.7.19) \qquad\qquad f(k) = 1 + ik = 1 + j$$

*является следствием равенства* (1.7.15). *Равенство*

$$(1.7.20) \qquad\qquad g(k) = k + jk = k - i = k - i$$

*является следствием равенства* (1.7.16). *Равенство*

$$(1.7.21)\ \ f(k)g(k) = (1+j)(k-i) = k-i+j(k-i) = k-i+jk-ji = k-i+i+k = 2k$$

*является следствием равенств* (1.7.19), (1.7.20). *Равенство*

$$(1.7.22) \qquad\qquad f * g(k) = k + kk^2 = 0$$

*является следствием равенства* (1.7.17). *Равенство*

$$(1.7.23) \qquad fg(k) = k + jk + ik^2 + ikjk = k + i - i + k = 2k$$

*является следствием равенства* (1.7.18). *Следовательно,*

$$f * g(k) \neq fg(k) = f(k)g(k)$$

$\square$

## 1.8. Деление в некоммутативной алгебре

Чтобы ответить на вопрос о делении многочленов в некоммутативной алгебре, мы должны выйти за рамки алгебры многочленов.

В алгебре с делением всё просто. Уравнение $ax = b$ имеет единственное решение, как только $a \neq 0$. Но не будем повторяться. Мы рассматривали этот вопрос в доказательстве теоремы 1.7.1. Алгебра с делением - это небольшое графство в обширном королевстве Алгебра.

Рассмотрим произвольную $D$-алгебру $A$. Допустим, мы утверждаем, что $A$-число $b$ делит $A$-число $a$. Что значит это утверждение?

Казалось бы, есть простой ответ.

**ОПРЕДЕЛЕНИЕ 1.8.1.** *Существует $A$-число $c$ такое, что*

$$a = bc$$

*В этом случае $A$-число $b$ называется левым делителем $A$-числа $a$.* $\square$

**ОПРЕДЕЛЕНИЕ 1.8.2.** *Существует $A$-число $c$ такое, что*

$$a = cb$$

*В этом случае $A$-число $b$ называется правым делителем $A$-числа $a$.* $\square$

Но должны ли мы различать левый и правый делители? Согласно определению, $D$-алгебра $A$ - это $D$-модуль $A$, в котором определено **билинейное отображение**

$$f : A \times A \to A$$

называемое **произведением**. Но если мы положим, что первый аргумент отображения $f$ является заданным $A$-числом, то билинейное отображение превращается в линейное отображение второго аргумента.

Следовательно, мы будем говорить, что $A$-число $b$ делит $A$-число $a$, если существует линейное отображение

$$c : A \to A$$

такое, что

$$a = c \circ b$$

Например, мы можем переписать определения 1.8.1, 1.8.2 следующим образом:

ОПРЕДЕЛЕНИЕ 1.8.3. *Существует линейное отображение* $1 \otimes c$ *такое, что* $(1 \otimes c) \circ b = a$ □

ОПРЕДЕЛЕНИЕ 1.8.4. *Существует линейное отображение* $c \otimes 1$ *такое, что* $(c \otimes 1) \circ b = a$ □

Мы также будем говорить, что $A$-число $b$ делит $A$-число $a$, если существует линейное отображение

$$c : A \to A$$

и $A$-число $d$ такие, что

$$a = c \circ b + d$$

Алгоритмы деления, поиска наибольшего общего делителя или наименьшего общего кратного являются стандартной процедурой.

Решение одного вопроса очень часто порождает новые вопросы. Среди вопросов, которые надо будет решить - это разложение $A$-числа на множители.

Например, какой будет результат деления $A$-числа $a$ на $A$-числа $b_1$ и $b_2$? Чтобы ответить на этот вопрос, разделим задачу на несколько шагов. [1.8]

Частное от деления $A$-числа $a$ на $A$-число $b_1$ имеет вид

$$(1.8.1) \qquad a = c_1 \circ b_1$$

где линейное отображение $c_1$ имеет вид

$$(1.8.2) \qquad c_1 = c_{10s} \otimes c_{11s}$$

У нас нет алгоритма деления тензора $c_1$ на $A$-число $b_2$. Однако равенство

$$(1.8.3) \qquad a = c_{10s} b_1 c_{11s}$$

является следствием равенств (1.8.1), (1.8.2). Из равенства (1.8.3) следует, что если мы хотим поделить $A$-число $a$ на $A$-числа $b_1$ и $b_2$, то мы должны поделить каждое $A$-число $c_{10s}$, $c_{11s}$ на $A$-число $b_2$ и результат подставить в равенство (1.8.3)

$$(1.8.4) \qquad a = (c_{10s2} \circ b_2) b_1 c_{11s} + c_{10s} b_1 (c_{11s2} \circ b_2)$$

---

[1.8] Чтобы упростить рассуждения, мы допускаем, что на каждом этапе мы имеем деление без остатка. Мы также рассматриваем ассоциативную $D$-алгебру.

Но выражение в правой части равенства ($1.8.4$) является билинейным отображением $A$-чисел $b_1$ и $b_2$

$$(1.8.5) \quad \begin{aligned} a &= c \circ (b_1, b_2) \\ c &= c_{10s20t} \otimes_2 c_{10s21t} \otimes_1 c_{11s} + c_{10s} \otimes_1 c_{11s20t} \otimes_2 c_{11s21t} \end{aligned}$$

Постепенно увеличивая число множителей, мы приходим к выводу, что если $A$-число $a$ является произведением $A$-чисел $b_1$, ..., $b_n$, то $A$-число $a$ является полилинейным отображением $A$-чисел $b_1$, ..., $b_n$.

Теперь мы можем вернуться к многочленам. Для того чтобы найти разложение многочлена на неприводимые множители, мы должны ответить на ряд вопросов. Как найти множество корней? Как велико множество разложений, если множество корней бесконечно? Если многочлен имеет неприводимый делитель, не имеющий корней, как найти этот делитель?

## 1.9. Обратная матрицы

Эта тема появилась в процессе написания книги.

В этой книге я рассмотрел определение левых и правых $_*$*-собственных чисел $n \times n$ матрицы $a_2$ в случае, когда матрица $a_2$ имеет $n$ различных собственных значений.

<table>
<tr><td>

В этом случае, множество собственных векторов имеет два эквивалентных определения

$$\boxed{(18.1.1) \quad u_{2*}{}^*a_{2*}{}^*u_2^{-1*}{}^* = a_1}$$
и
$$\boxed{(18.1.2) \quad u_{2*}{}^*a_2 = a_{1*}{}^*u_2}$$

</td><td>

В этом случае, множество собственных векторов имеет два эквивалентных определения

$$\boxed{(18.1.4) \quad u_2^{-1*}{}^*{}_*{}^*a_{2*}{}^*u_2 = a_1}$$
и
$$\boxed{(18.1.5) \quad a_{2*}{}^*u_2 = u_{2*}{}^*a_1}$$

</td></tr>
</table>

Матрицы $u_2$, $a_1$ также являются $n \times n$ матрицами. В этом случае, множества левых и правых $_*$*-собственных чисел совпадают.

Однако профессор Кон обратил внимание, что возможны случаи, когда множество левых $_*$*-собственных чисел отличается от множества правых $_*$*-собственных чисел. Чтобы найти такую матрицу, я решил рассмотреть случай, когда матрица $a_2$ имеет $k$, $k < n$, различных собственных значений. В этом случае я вынужден ограничиться определением, подобным определению

$$\boxed{(18.1.2) \quad u_{2*}{}^*a_2 = a_{1*}{}^*u_2}$$ В этом случае, матрица $u_2$ является $n \times k$ матрицей и матрица $a_1$ является $k \times k$ матрицей. Пытаясь ответить на вопрос могу ли я записать определение аналогичное определению ($18.1.1$), я понял, что мне надо рассмотреть понятие правой $_*$*-обратной матрицы.

## 1.10. Некоммутативная сумма

С детства мы знаем, что от перестановки мест слагаемых сумма не меняется. Поэтому, когда я услышал о некоммутативном сложении ([24]), я слегка растерялся. Однако мне недавно повезло встретить геометрию, где сумма векторов некоммутативна.

В школе я увлекался многими вещами: физикой, химией, черчением. Математика была немного в стороне. Я не любил арифметику, но мне нравились алгебра и геометрия. Я тогда ещё не знал, что все пути ведут в математику.

Черчение было моим любимым предметом. Мой учитель, Габис Сергей Александрович, подарил мне сборник задач по черчению, и я добросовестно все их перерешал. А потом придумал свою собственную задачу и решил её. К моему удивлению, первая теорема, которую мы доказывали на уроке геометрии на следующий год, была именно решением этой задачи.

Я понимал, что после школы мне будет сложно поступить в институт. Поэтому после восьмого класса я решил что я буду поступить в техникум. Бела Марковна предложила позаниматься со мной летом, посмотреть насколько я готов к экзамену по математике. Один из уроков остался в моей памяти. Когда мы рассматривали графики прямой и параболы, Бела Марковна сказала: нарисовать прямую или параболу легко, а как быть с кардиограммой? И я загорелся идеей найти аналитическое выражение для кардиограммы.

Мой товарищ предложил мне на протяжении лета вместе изучать математический анализ. В нашем доме жили мальчики, которые учились в техникуме, и они одолжили мне свой учебник по математическому анализу. К этому времени мой товарищ потерял интерес к своей идее, и я стал изучать математический анализ самостоятельно.

Первая глава учебника была посвящена аналитической геометрии. И я понял, что это то что мне нужно для анализа кривых кардиограммы. Аналитическая геометрия решила моё будущее. Я понял, что буду готовиться в университет на механико-математический факультет.

Я перечитал много книг по математическому анализу. Но в какой-то момент у меня сложилось впечатление что один автор списывает у другого. Я стал искать другие книги. Не помню откуда у меня появился первый том Фихтенгольца. Но эта книга произвела на меня неизгладимое впечатление. Ближе к весне я почувствовал, что книга становится трудной, но продолжал читать. В марте у меня был грипп с необычайно высокой температурой. Но я сомневаюсь, что это был грипп, так как после этого первый том Фихтенгольца стал для меня понятен.

Изучение математического анализа отразилось на моих попытках восстановить функциональную зависимость с графика. Я тогда ещё не понимал, что не всякое отображение может быть представлено аналитически. Я рисовал графики, составлял таблицы приращения функции в зависимости от приращения аргумента. Но вскоре эта затея мне надоела, и я сменил занятие.

Когда мне в руки попала книга по тензорному исчислению, я обнаружил, что сложные построения аналитической геометрии по классификации кривых второго порядка можно заменить простыми теоремами линейной алгебры. После этого, мой интерес к аналитической геометрии заметно стал меньше.

Прошло почти полвека. И неожиданно аналитическая геометрия преподнесла мне два приятных сюрприза. Сперва я обнаружил, что аналитическая геометрия является интересной структурой с точки зрения теории представлений универсальных алгебр. Точнее, я отождествляю аналитическую и аффинную геометрии, но я думаю, что это небольшая ошибка.

Второй сюрприз связан с тем, что я понял, что я могу построить модель аффинной геометрии в римановом пространстве. Так как я очень долго занимался метрико-аффинными многообразиями, то оставался один шаг, чтобы увидеть модель аффинной геометрии, где сумма некоммутативна. С этого момента некоммутативное сложение перестало быть для меня абстрактной концепцией.

Глава 2

# Представление универсальной алгебры

В этой главе собраны определения и теоремы, которые необходимы для понимания текста предлагаемой книги. Поэтому читатель может обращаться к утверждениям из этой главы по мере чтения основного текста книги.

## 2.1. Универсальная алгебра

Определение 2.1.1. *Для отображения*

$$f : A \to B$$

*множество*

$$\mathrm{Im}\, f = \{f(a) : a \in A\}$$

*называется* **образом отображения** $f$. □

Определение 2.1.2. *Для любых множеств*[2.1] $A$, $B$, **декартова степень** $B^A$ *- это множество отображений*

$$f : A \to B$$
□

Определение 2.1.3. *Пусть дано множество* $A$ *и целое число* $n \geq 0$. *Отображение*[2.2]

$$\omega : A^n \to A$$

*называется* $n$-**арной операцией на множестве** $A$ *или просто* **операцией на множестве** $A$. *Для любых* $a_1, ..., a_n \in A$, *мы пользуемся любой из форм записи* $\omega(a_1, ..., a_n)$, $a_1...a_n\omega$ *для обозначения образа отображения* $\omega$. □

Замечание 2.1.4. *Согласно определениям* 2.1.2, 2.1.3, $n$-*арная операция* $\omega \in A^{A^n}$. □

Определение 2.1.5. **Область операторов** *- это множество операторов*[2.3] $\Omega$ *вместе с отображением*

$$a : \Omega \to N$$

*Если* $\omega \in \Omega$, *то* $a(\omega)$ *называется* **арностью** *оператора* $\omega$. *Если* $a(\omega) = n$, *то оператор* $\omega$ *называется* $n$-*арным. Мы пользуемся обозначением*

$$\Omega(n) = \{\omega \in \Omega : a(\omega) = n\}$$

*для множества* $n$-*арных операторов.* □

---

[2.1] Я следую определению из примера (iV), [18], страницы 17, 18.
[2.2] Определение 2.1.3 опирается на определение в примере (vi), страница [18]-26.
[2.3] Я следую определению 1, страница [18]-62.

Определение 2.1.6. *Пусть $A$ - множество, а $\Omega$ - область операто-ров.*[2.4] *Семейство отображений*

$$\Omega(n) \to A^{A^n} \quad n \in N$$

*называется структурой $\Omega$-алгебры на $A$. Множество $A$ со структурой $\Omega$-алгебры называется $\Omega$-**алгеброй** $A_\Omega$ или **универсальной алгеброй**. Мно-жество $A$ называется **носителем** $\Omega$-алгебры.* $\square$

Область операторов $\Omega$ описывает множество $\Omega$-алгебр. Элемент множе-ства $\Omega$ называется оператором, так как операция предполагает некоторое мно-жество. Согласно замечанию 2.1.4 и определению 2.1.6, каждому оператору $\omega \in \Omega(n)$ сопоставляется $n$-арная операция $\omega$ на $A$.

Определение 2.1.7. *Пусть $A$, $B$ - $\Omega$-алгебры и $\omega \in \Omega(n)$. Отображе-ние*[2.5]

$$f : A \to B$$

**согласовано с операцией** $\omega$, *если, для любых $a_1, ..., a_n \in A$,*

$$f(a_1...a_n\omega) = f(a_1)...f(a_n)\omega$$

*Отображение $f$ называется **гомоморфизмом** $\Omega$-алгебры $A$ в $\Omega$-алгебру $B$, если $f$ согласовано с каждым $\omega \in \Omega$.* $\square$

Определение 2.1.8. *Гомоморфизм*

$$f : A \to B$$

*называется*[2.6] **изоморфизмом** *между $A$ и $B$, если соответствие $f^{-1}$ яв-ляется гомоморфизмом. Если существует изоморфизм между $A$ и $B$, то говорят, что $A$ и $B$ изоморфны, и пишут $A \cong B$. Инъективный гомомор-физм называется **мономорфизмом**. Сюръективный гомоморфизм называет-ся **эпиморфизмом**.* $\square$

Определение 2.1.9. *Гомоморфизм*

$$f : A \to A$$

*источником и целью которого является одна и таже алгебра, называется **эндоморфизмом**. Мы обозначим $End(\Omega; A)$ множество эндоморфизмов $\Omega$-алгебры $A$. Изоморфизм*

$$f : A \to A$$

*называется **автоморфизмом**.* $\square$

Теорема 2.1.10. *Пусть отображение*

$$f : A \to B$$

---

[2.4] Я следую определению 2, страница [18]-62.
[2.5] Я следую определению на странице [18]-63.
[2.6] Я следую определению на странице [18]-63.

*является гомоморфизмом $\Omega$-алгебры $A$ в $\Omega$-алгебру $B$. Пусть отображение*

$$g : B \to C$$

*является гомоморфизмом $\Omega$-алгебры $B$ в $\Omega$-алгебру $C$. Тогда отображение*[2.7]

$$h = g \circ f : A \to C$$

*является гомоморфизмом $\Omega$-алгебры $A$ в $\Omega$-алгебру $C$.*

ДОКАЗАТЕЛЬСТВО. Пусть $\omega \in \Omega_n$ - n-арная операция. Тогда

(2.1.1) $$f(a_1...a_n\omega) = f(a_1)...f(a_n)\omega$$

для любых $a_1, ..., a_n \in A$, и

(2.1.2) $$g(b_1...b_n\omega) = g(b_1)...g(b_n)\omega$$

для любых $b_1, ..., b_n \in B$. Равенство

$$
\begin{aligned}
(g \circ f)(a_1...a_n\omega) &= g(f(a_1...a_n\omega)) \\
&= g(f(a_1)...f(a_n)\omega) \\
&= g(f(a_1))...g(f(a_n))\omega \\
&= (g \circ f)(a_1)...(g \circ f)(a_n)\omega
\end{aligned}
$$

(2.1.3)

является следствием равенств (2.1.1), (2.1.2). Теорема является следствием равенства (2.1.3) для любой операции $\omega$.  □

ОПРЕДЕЛЕНИЕ 2.1.11. *Если существует мономорфизм $\Omega$-алгебры $A$ в $\Omega$-алгебру $B$, то говорят, что $A$ **может быть вложена в** $B$.*  □

ОПРЕДЕЛЕНИЕ 2.1.12. *Если существует эпиморфизм из $A$ в $B$, то $B$ называется **гомоморфным образом** алгебры $A$.*  □

## 2.2. Декартово произведение универсальных алгебр

ОПРЕДЕЛЕНИЕ 2.2.1. *Пусть $\mathcal{C}$ - некоторая категория. Объект $P$ категории $\mathcal{C}$ называется **универсально притягивающим**,[2.8]если для любого объекта $R$ категории $\mathcal{C}$ существует единственный морфизм*

$$f : R \to P$$

□

ОПРЕДЕЛЕНИЕ 2.2.2. *Пусть $\mathcal{C}$ - некоторая категория. Объект $P$ категории $\mathcal{C}$ называется **универсально отталкивающим**,[2.9]если для любого объекта $R$ категории $\mathcal{C}$ существует единственный морфизм*

$$f : P \to R$$

□

---

[2.7] Я следую предложению [18]-3.2 на странице [18]-71.
[2.8] Смотри также определение в [1], страница 47.

ОПРЕДЕЛЕНИЕ 2.2.3. *Пусть $\mathcal{A}$ - категория. Пусть $\{B_i, i \in I\}$ - множество объектов из $\mathcal{A}$. Пусть $\mathcal{A}(B_i, i \in I)$ - категория, объектами которой являются кортежи $(P, f)$, где $P$ - объект категории $\mathcal{A}$ и $f$ - множество морфизмов*

$$\{f_i : P \to B_i, i \in I\}$$

*Универсально притягивающий объект категории $\mathcal{A}(B_i, i \in I)$*

$$P = \prod_{i \in I} B_i$$

*называется **произведением множества объектов** $\{B_i, i \in I\}$ **в категории** $\mathcal{A}$.* [2.10]

*Если $|I| = n$, то для произведения множества объектов $\{B_i, i \in I\}$ в $\mathcal{A}$ мы так же будем пользоваться записью*

$$P = \prod_{i=1}^{n} B_i = B_1 \times ... \times B_n$$

□

ТЕОРЕМА 2.2.4. *Пусть $\mathcal{A}$ - категория. Пусть $\{B_i, i \in I\}$ - множество объектов из $\mathcal{A}$. Произведение множества объектов $\{B_i, i \in I\}$ в категории $\mathcal{A}$ - это объект*

$$P = \prod_{i \in I} B_i$$

*и множество морфизмов*

$$\{f_i : P \to B_i, i \in I\}$$

*такие, что для любого объекта $R$ и множества морфизмов*

$$\{g_i : R \to B_i, i \in I\}$$

*существует единственный морфизм*

$$h : R \to P$$

*такой, что диаграмма*

$$P \xrightarrow{f_i} B_i \qquad f_i \circ h = g_i$$
$$h \uparrow \quad \nearrow g_i$$
$$R$$

*коммутативна для всех $i \in I$.*

ДОКАЗАТЕЛЬСТВО. Теорема является следствием определений 2.2.1, 2.2.3.
□

ПРИМЕР 2.2.5. *Пусть $\mathcal{S}$ - категория множеств.* [2.11] *Согласно определению*

*2.2.3,  декартово произведение*

$$A = \prod_{i \in I} A_i$$

*семейства множеств  $(A_i, i \in I)$  и семейство проекций на i-й множитель*

$$p_i : A \to A_i$$

*являются произведением в категории $\mathcal{S}$.*          □

ТЕОРЕМА 2.2.6.  *Произведение существует в категории $\mathcal{A}$ $\Omega$-алгебр. Пусть $\Omega$-алгебра $A$ и семейство морфизмов*

$$p_i : A \to A_i    i \in I$$

*является произведением в категории $\mathcal{A}$. Тогда*

2.2.6.1: *Множество $A$ является декартовым произведением семейства множеств  $(A_i, i \in I)$*

2.2.6.2: *Гомоморфизм $\Omega$-алгебры*

$$p_i : A \to A_i$$

*является проекцией на i-й множитель.*

2.2.6.3: *Любое $A$-число $a$ может быть однозначно представлено в виде кортежа  $(p_i(a), i \in I)$  $A_i$-чисел.*

2.2.6.4: *Пусть  $\omega \in \Omega$  - n-арная операция. Тогда операция $\omega$ определена покомпонентно*

$$(2.2.1)          a_1...a_n\omega = (a_{1i}...a_{ni}\omega, i \in I)$$

*где  $a_1 = (a_{1i}, i \in I)$, ..., $a_n = (a_{ni}, i \in I)$ .*

ДОКАЗАТЕЛЬСТВО. Пусть

$$A = \prod_{i \in I} A_i$$

декартово произведение семейства множеств  $(A_i, i \in I)$  и, для каждого $i \in I$, отображение

$$p_i : A \to A_i$$

является проекцией на $i$-й множитель. Рассмотрим диаграмму морфизмов в категории множеств $\mathcal{S}$

$$(2.2.2)$$

$$p_i \circ \omega = g_i$$

где отображение $g_i$ определено равенством

$$g_i(a_1, ..., a_n) = p_i(a_1)...p_i(a_n)\omega$$

Согласно определению  2.2.3,  отображение $\omega$ определено однозначно из множества диаграмм (2.2.2)

$$(2.2.3)          a_1...a_n\omega = (p_i(a_1)...p_i(a_n)\omega, i \in I)$$

Равенство (2.2.1) является следствием равенства (2.2.3).          □

Определение 2.2.7. *Если $\Omega$-алгебра $A$ и семейство морфизмов*

$$p_i : A \to A_i \quad i \in I$$

*является произведением в категории $\mathcal{A}$, то $\Omega$-алгебра $A$ называется* **прямым** *или* **декартовым произведением** *$\Omega$-алгебр* $(A_i, i \in I)$.          $\square$

Теорема 2.2.8. *Пусть множество $A$ является декартовым произведением множеств $(A_i, i \in I)$ и множество $B$ является декартовым произведением множеств $(B_i, i \in I)$. Для каждого $i \in I$, пусть*

$$f_i : A_i \to B_i$$

*является отображением множества $A_i$ в множество $B_i$. Для каждого $i \in I$, рассмотрим коммутативную диаграмму*

(2.2.4)

$$
\begin{array}{ccc}
B & \xrightarrow{p'_i} & B_i \\
\uparrow{\scriptstyle f} & & \uparrow{\scriptstyle f_i} \\
A & \xrightarrow{p_i} & A_i
\end{array}
$$

*где отображения $p_i$, $p'_i$ являются проекцией на $i$-й множитель. Множество коммутативных диаграмм* (2.2.4) *однозначно определяет отображение*

$$f : A \to B$$

$$f(a_i, i \in I) = (f_i(a_i), i \in I)$$

Доказательство. Для каждого $i \in I$, рассмотрим коммутативную диаграмму

(2.2.5)

$$
\begin{array}{ccc}
B & \xrightarrow{p'_i} & B_i \\
\uparrow{\scriptstyle f} & \searrow{\scriptstyle g_i} & \uparrow{\scriptstyle f_i} \\
 & (1) & \\
 & (2) & \\
A & \xrightarrow{p_i} & A_i
\end{array}
$$

Пусть $a \in A$. Согласно утверждению 2.2.6.3, $A$-число $a$ может быть представлено в виде кортежа $A_i$-чисел

(2.2.6)          $$a = (a_i, i \in I) \quad a_i = p_i(a) \in A_i$$

Пусть

(2.2.7)          $$b = f(a) \in B$$

Согласно утверждению 2.2.6.3, $B$-число $b$ может быть представлено в виде кортежа $B_i$-чисел

(2.2.8)          $$b = (b_i, i \in I) \quad b_i = p'_i(b) \in B_i$$

Из коммутативности диаграммы (1) и из равенств (2.2.7), (2.2.8) следует, что
(2.2.9)                          $b_i = g_i(b)$
Из коммутативности диаграммы (2) и из равенства (2.2.6) следует, что
$$b_i = f_i(a_i)$$
□

Теорема 2.2.9. *Пусть $\Omega$-алгебра $A$ является декартовым произведением $\Omega$-алгебр $(A_i, i \in I)$ и $\Omega$-алгебра $B$ является декартовым произведением $\Omega$-алгебр $(B_i, i \in I)$. Для каждого $i \in I$, пусть отображение*
$$f_i : A_i \to B_i$$
*является гомоморфизмом $\Omega$-алгебры. Тогда отображение*
$$f : A \to B$$
*определённое равенством*
(2.2.10)                          $f(a_i, i \in I) = (f_i(a_i), i \in I)$
*является гомоморфизмом $\Omega$-алгебры.*

Доказательство. Пусть $\omega \in \Omega$ - n-арная операция. Пусть $a_1 = (a_{1i}, i \in I), ..., a_n = (a_{ni}, i \in I)$ и $b_1 = (b_{1i}, i \in I), ..., b_n = (b_{ni}, i \in I)$. Из равенств (2.2.1), (2.2.10) следует, что
$$f(a_1...a_n\omega) = f(a_{1i}...a_{ni}\omega, i \in I)$$
$$= (f_i(a_{1i}...a_{ni}\omega), i \in I)$$
$$= ((f_i(a_{1i}))...(f_i(a_{ni})), i \in I)$$
$$= (b_{1i}...b_{ni}\omega, i \in I)$$
$$f(a_1)...f(a_n)\omega = b_1...b_n\omega = (b_{1i}...b_{ni}\omega, i \in I)$$
□

Определение 2.2.10. *Эквивалентность на $\Omega$-алгебре $A$, которая является подалгеброй $\Omega$-алгебры $A^2$, называется* **конгруенцией** [2.12] *на $A$.* □

Теорема 2.2.11 (первая теорема об изоморфизмах). *Пусть*
$$f : A \to B$$
*гомоморфизм $\Omega$-алгебр с ядром $s$. Тогда отображение $f$ имеет разложение*

$$
\begin{array}{ccc}
A/\ker f & \xrightarrow{\ q\ } & f(A) \\
{\scriptstyle p}\big\uparrow & & \big\downarrow{\scriptstyle r} \\
A & \xrightarrow{\ f\ } & B
\end{array}
\qquad f = p \circ q \circ r
$$

2.2.11.1: **Ядро гомоморфизма** $\ker f = f \circ f^{-1}$ *является конгруэнцией на $\Omega$-алгебре $A$.*

2.2.11.2: *Множество $A/\ker f$ является $\Omega$-алгеброй.*

─────────

[2.12] Я следую определению на странице [18]-71.

2.2.11.3: *Отображение*

$$p : a \in A \to a^{\ker f} \in A/\ker f$$

*является эпиморфизмом и называется* **естественным гомоморфизмом**.

2.2.11.4: *Отображение*

$$q : p(a) \in A/\ker f \to f(a) \in f(A)$$

*является изоморфизмом.*

2.2.11.5: *Отображение*

$$r : f(a) \in f(A) \to f(a) \in B$$

*является мономорфизмом.*

Доказательство. Утверждение 2.2.11.1 является следствием предложения II.3.4 ([18], страница 72). Утверждения 2.2.11.2, 2.2.11.3 являются следствием теоремы II.3.5 ([18], страница 72) и последующего определения. Утверждения 2.2.11.4, 2.2.11.5 являются следствием теоремы II.3.7 ([18], страница 74). □

Определение 2.2.12. *Пусть* $\mathcal{A}$ *- категория. Пусть* $\{B_i, i \in I\}$ *- множество объектов из* $\mathcal{A}$. *Пусть* $\mathcal{A}[B_i, i \in I]$ *- категория, объектами которой являются кортежи* $(P, f)$, *где* $P$ *- объект категории* $\mathcal{A}$ *и* $f$ *- множество морфизмов*

$$\{f_i : B_i \to P, i \in I\}$$

*Универсально отталкивающий объект категории* $\mathcal{A}[B_i, i \in I]$

$$P = \coprod_{i \in I} B_i$$

*называется* **копроизведением множества объектов** $\{B_i, i \in I\}$ **в категории** $\mathcal{A}$.[2.13]

*Если* $|I| = n$, *то для копроизведения множества объектов* $\{B_i, i \in I\}$ *в* $\mathcal{A}$ *мы так же будем пользоваться записью*

$$P = \coprod_{i=1}^{n} B_i = B_1 \coprod \ldots \coprod B_n$$

□

Теорема 2.2.13. *Пусть* $\mathcal{A}$ *- категория. Пусть* $\{B_i, i \in I\}$ *- множество объектов из* $\mathcal{A}$. *Объект*

$$P = \coprod_{i \in I} B_i$$

*и множество морфизмов*

$$\{f_i : B_i \to P, i \in I\}$$

---

[2.13] Определение дано согласно определению на странице [1]-46.

*называется* **копроизведением множества объектов** $\{B_i, i \in I\}$ **в кате-**
**гории** $\mathcal{A}$, [2.14] *если для любого объекта $R$ и множество морфизмов*

$$\{g_i : B_i \to R, i \in I\}$$

*существует единственный морфизм*

$$h : P \to R$$

*такой, что диаграмма*

$$h \circ f_i = g_i$$

*коммутативна для всех $i \in I$.*

Доказательство. Теорема является следствием определений 2.2.2, 2.2.12.
□

## 2.3. Представление универсальной алгебры

Определение 2.3.1. *Пусть множество $A_2$ является $\Omega_2$-алгеброй. Пусть
на множестве преобразований* $\operatorname{End}(\Omega_2, A_2)$ *определена структура $\Omega_1$-алгеб-
ры. Гомоморфизм*

$$f : A_1 \to \operatorname{End}(\Omega_2; A_2)$$

$\Omega_1$-*алгебры $A_1$ в $\Omega_1$-алгебру* $\operatorname{End}(\Omega_2, A_2)$ *называется* **представлением** $\Omega_1$-
**алгебры** $A_1$ *или* $A_1$-**представлением** *в $\Omega_2$-алгебре $A_2$.* □

Мы будем также пользоваться записью

$$f : A_1 \longrightarrow^* A_2$$

для обозначения представления $\Omega_1$-алгебры $A_1$ в $\Omega_2$-алгебре $A_2$.

Определение 2.3.2. *Мы будем называть представление*

$$f : A_1 \longrightarrow^* A_2$$

$\Omega_1$-*алгебры $A_1$* **эффективным**, [2.15] *если отображение*

$$f : A_1 \to \operatorname{End}(\Omega_2; A_2)$$

*является изоморфизмом $\Omega_1$-алгебры $A_1$ в* $\operatorname{End}(\Omega_2, A_2)$. □

Определение 2.3.3. *Мы будем называть представление*

$$g : A_1 \longrightarrow^* A_2$$

---

[2.14] Определение дано согласно [1], страница 46.
[2.15] Аналогичное определение эффективного представления группы смотри в [21], страница
16, [22], страница 111, [19], страница 51 (Кон называет такое представление точным).

$\Omega_1$-алгебры $A_1$ **свободным**,[2.16]*если из утверждения*

$$f(a_1)(a_2) = f(b_1)(a_2)$$

*для любого* $a_2 \in A_2$ *следует, что* $a = b$.          □

---

ТЕОРЕМА 2.3.4. *Свободное представление эффективно.*

ДОКАЗАТЕЛЬСТВО. Теорема является следствием теоремы [15]-3.1.6.          □

ЗАМЕЧАНИЕ 2.3.5. *Представление группы вращений в аффинном пространстве эффективно, но не свободно, так как начало координат является неподвижной точкой любого преобразования.*          □

---

ОПРЕДЕЛЕНИЕ 2.3.6. *Мы будем называть представление*

$$g : A_1 \longrightarrow A_2$$

$\Omega_1$-алгебры **транзитивным**,[2.17]*если для любых* $a, b \in 2$ *существует такое* $g$, *что*

$$a = f(g)(b)$$

*Мы будем называть представление* $\Omega_1$-*алгебры* **однотранзитивным**, *если оно транзитивно и свободно.*          □

---

ТЕОРЕМА 2.3.7. *Представление однотранзитивно тогда и только тогда, когда для любых* $a, b \in A_2$ *существует одно и только одно* $g \in A_1$ *такое, что* $a = f(g)(b)$.

ДОКАЗАТЕЛЬСТВО. Следствие определений 2.3.3 и 2.3.6.          □

---

ТЕОРЕМА 2.3.8. *Пусть*

$$f : A_1 \longrightarrow A_2$$

*однотранзитивное представление* $\Omega_1$-*алгебры* $A_1$ *в* $\Omega_2$-*алгебре* $A_2$. *Существует структура* $\Omega_1$-*алгебры на множестве* $A_2$.

ДОКАЗАТЕЛЬСТВО. Теорема является следствием теоремы [15]-3.1.10.          □

---

ОПРЕДЕЛЕНИЕ 2.3.9. *Пусть*

$$f : A_1 \longrightarrow A_2$$

*представление* $\Omega_1$-*алгебры* $A_1$ *в* $\Omega_2$-*алгебре* $A_2$ *и*

$$g : B_1 \longrightarrow B_2$$

---

[2.16] Аналогичное определение свободного представления группы смотри в [21], страница 16.
[2.17] Аналогичное определение транзитивного представления группы смотри в [22], страница 110, [19], страница 51.

представление $\Omega_1$-алгебры $B_1$ в $\Omega_2$-алгебре $B_2$. Для $i = 1, 2,$ пусть отображение

$$r_i : A_i \to B_i$$

является гомоморфизмом $\Omega_i$-алгебры. Кортеж отображений $r = (r_1, r_2)$ таких, что

(2.3.1) $$r_2 \circ f(a) = g(r_1(a)) \circ r_2$$

называется **морфизмом представлений** из $f$ в $g$. Мы также будем говорить, что определён **морфизм представлений** $\Omega_1$-**алгебры в** $\Omega_2$-**алгебре**. □

Замечание 2.3.10. *Мы можем рассматривать пару отображений $r_1$, $r_2$ как отображение*

$$F : A_1 \cup A_2 \to B_1 \cup B_2$$

*такое, что*

$$F(A_1) = B_1 \qquad F(A_2) = B_2$$

*Поэтому в дальнейшем кортеж отображений $r = (r_1, r_2)$ мы будем также называть отображением и пользоваться записью*

$$r : f \to g$$

*Пусть $a = (a_1, a_2)$ - кортеж $A$-чисел. Мы будем пользоваться записью*

$$r(a) = (r_1(a_1), r_2(a_2))$$

*для образа кортежа $A$-чисел при морфизме представлений $r$.* □

Замечание 2.3.11. *Рассмотрим морфизм представлений*

$$(r_1 : A_1 \to B_1, \ r_2 : A_2 \to B_2)$$

*Мы можем обозначать элементы множества $B_1$, пользуясь буквой по образцу $b \in B_1$. Но если мы хотим показать, что $b$ является образом элемента $a \in A_1$, мы будем пользоваться обозначением $r_1(a)$. Таким образом, равенство*

$$r_1(a) = r_1(a)$$

*означает, что $r_1(a)$ (в левой части равенства) является образом $a \in A_1$ (в правой части равенства). Пользуясь подобными соображениями, мы будем обозначать элемент множества $B_2$ в виде $r_2(m)$. Мы будем следовать этому соглашению, изучая соотношения между гомоморфизмами $\Omega_1$-алгебр и отображениями между множествами, где определены соответствующие представления.* □

Теорема 2.3.12. *Пусть*

$$f : A_1 \longrightarrow_* A_2$$

*представление $\Omega_1$-алгебры $A_1$ в $\Omega_2$-алгебре $A_2$ и*

$$g : B_1 \longrightarrow_* B_2$$

*представление $\Omega_1$-алгебры $B_1$ в $\Omega_2$-алгебре $B_2$. Отображение*

$$(r_1 : A_1 \to B_1, \ r_2 : A_2 \to B_2)$$

*является морфизмом представлений тогда и только тогда, когда*

$$(2.3.2) \qquad r_2(f(a)(m)) = g(r_1(a))(r_2(m))$$

Доказательство. Теорема является следствием теоремы [15]-3.2.5.   $\square$

Замечание 2.3.13. *Мы можем интерпретировать* (2.3.2) *двумя способами*

- *Пусть преобразование* $f(a)$ *отображает* $m \in A_2$ *в* $f(a)(m)$. *Тогда преобразование* $g(r_1(a))$ *отображает* $r_2(m) \in B_2$ *в* $r_2(f(a)(m))$.
- *Мы можем представить морфизм представлений из* $f$ *в* $g$, *пользуясь диаграммой*

(2.3.3)

*Из* (2.3.1) *следует, что диаграмма* (1) *коммутативна.*

*Мы будем также пользоваться диаграммой*

(2.3.4)

*вместо диаграммы* (2.3.3).                                              $\square$

Определение 2.3.14. *Если представления* $f$ *и* $g$ *совпадают, то морфизм представлений* $r = (r_1, r_2)$ *называется* **морфизмом представления** $f$.   $\square$

Определение 2.3.15. *Пусть*

$$f : A_1 \longrightarrow\!\!\!* A_2$$

*представление* $\Omega_1$-*алгебры* $A_1$ *в* $\Omega_2$-*алгебре* $A_2$ *и*

$$g : A_1 \longrightarrow\!\!\!* B_2$$

*представление* $\Omega_1$-*алгебры* $A_1$ *в* $\Omega_2$-*алгебре* $B_2$. *Пусть*

$$(\mathrm{id} : A_1 \to A_1, \ r_2 : A_2 \to B_2)$$

*морфизм представлений. В этом случае мы можем отождествить морфизм* $(\mathrm{id}, r_2)$ *представлений* $\Omega_1$-*алгебры и соответствующий гомоморфизм* $r_2$ $\Omega_2$-алгебры и будем называть гомоморфизм $r_2$ **приведенным морфизмом представлений**. Мы будем пользоваться диаграммой

(2.3.5)

для представления приведенного морфизма $r_2$ представлений $\Omega_1$-алгебры. Из диаграммы следует

(2.3.6) $$r_2 \circ f(a) = g(a) \circ r_2$$

Мы будем также пользоваться диаграммой

вместо диаграммы (2.3.5).                                                   □

Теорема 2.3.16. *Пусть*

$$f : A_1 \longrightarrow A_2$$

*представление $\Omega_1$-алгебры $A_1$ в $\Omega_2$-алгебре $A_2$ и*

$$g : A_1 \longrightarrow B_2$$

*представление $\Omega_1$-алгебры $A_1$ в $\Omega_2$-алгебре $B_2$. Отображение*

$$r_2 : A_2 \to B_2$$

*является приведенным морфизмом представлений тогда и только тогда, когда*

(2.3.7) $$r_2(f(a)(m)) = g(a)(r_2(m))$$

Доказательство. Теорема является следствием теоремы [15]-3.4.3.   □

Определение 2.3.17. *Пусть $A_1$ - $\Omega_1$-алгебра. Пусть $\mathcal{A}_2$ - категория $\Omega_2$-алгебр. Мы определим* **категорию** $A_1(\mathcal{A}_2)$ **представлений** *$\Omega_1$-алгебры $A_1$ в $\Omega_2$-алгебре. Объектами этой категории являются представления $\Omega_1$-алгебры $A_1$ в $\Omega_2$-алгебре. Морфизмами этой категории являются приведенные морфизмы соответствующих представлений.*   □

Теорема 2.3.18. *В категории $A_1(\mathcal{A}_2)$ существует произведение эффективных представлений $\Omega_1$-алгебры $A_1$ в $\Omega_2$-алгебре и это произведение является эффективным представлением $\Omega_1$-алгебры $A_1$.*

Доказательство. Теорема является следствием теоремы [15]-4.4.2.   □

## 2.4. Ω-группа

Определение 2.4.1. *Пусть в $\Omega$-алгебре $A$ определена операция сложения. Отображение*

$$f : A \to A$$

*$\Omega_1$-алгебры $A$ называется* **аддитивным отображением**, *если*

$$f(a + b) = f(a) + f(b)$$

□

Определение 2.4.2. *Отображение*

$$g : A^n \to A$$

*называется* **полиаддитивным отображением**, *если для любого $i$, $i = 1, ..., n$,*

$$f(a_1, ..., a_i + b_i, ..., a_n) = f(a_1, ..., a_i, ..., a_n) + f(a_1, ..., b_i, ..., a_n)$$

□

Определение 2.4.3. *Пусть в $\Omega_1$-алгебре $A$ определена операция сложения, которая не обязательно коммутативна. Мы пользуемся символом $+$ для обозначения операции суммы. Положим*

$$\Omega = \Omega_1 \setminus \{+\}$$

*Если $\Omega_1$-алгебра $A$ является группой относительно операции сложения и любая операция $\omega \in \Omega$ является полиаддитивным отображением, то $\Omega_1$-алгебра $A$ называется* **$\Omega$-группой**. *Если $\Omega$-группа $A$ является ассоциативной группой относительно операции сложения, то $\Omega_1$-алгебра $A$ называется* **ассоциативной $\Omega$-группой**. *Если $\Omega$-группа $A$ является абелевой группой относительно операции сложения, то $\Omega_1$-алгебра $A$ называется* **абелевой $\Omega$-группой**.   □

Теорема 2.4.4. *Пусть $\omega \in \Omega(n)$ - полиаддитивное отображение. Операция $\omega$* **дистрибутивна** *относительно сложения*

$$a_1...(a_i + b_i)...a_n\omega = a_1...a_i...a_n\omega + a_1...b_i...a_n\omega \quad i = 1, ..., n$$

Доказательство. Теорема является следствием определений 2.4.2, 2.4.3.   □

Определение 2.4.5. *Пусть произведение*

$$c_1 = a_1 * b_1$$

*является операцией $\Omega_1$-алгебры $A$. Положим $\Omega = \Omega_1 \setminus \{*\}$. Если $\Omega_1$-алгебра $A$ является группой относительно произведения и для любой операции $\omega \in \Omega(n)$ умножение дистрибутивно относительно операция $\omega$*

$$a * (b_1...b_n\omega) = (a * b_1)...(a * b_n)\omega$$
$$(b_1...b_n\omega) * a = (b_1 * a)...(b_n * a)\omega$$

*то $\Omega_1$-алгебра $A$ называется **мультипликативной $\Omega$-группой**.* □

О‌пределение 2.4.6. *Если*

(2.4.1) $$a * b = b * a$$

*то мультипликативная $\Omega$-группа называется **абелевой**.* □

О‌пределение 2.4.7. *Пусть сложение*

$$c_1 = a_1 + b_1$$

*и произведение*

$$c_1 = a_1 * b_1$$

*являются операциями $\Omega_1$-алгебры $A$. Положим $\Omega = \Omega_1 \setminus \{+, *\}$. Если $\Omega_1$-алгебра $A$ является абелевой $\Omega \cup \{*\}$-группой и мультипликативной $\Omega \cup \{+\}$-группой, то $\Omega_1$-алгебра $A$ называется $\Omega$-**кольцом**.* □

Т‌еорема 2.4.8. *Произведение в $\Omega$-кольце **дистрибутивно** относительно сложения*

$$a * (b_1 + b_2) = a * b_1 + a * b_2$$
$$(b_1 + b_2) * a = b_1 * a + b_2 * a$$

Доказательство. Теорема является следствием определений 2.4.3, 2.4.5, 2.4.7. □

## 2.5. Представление мультипликативной $\Omega$-группы

Мы предполагаем, что преобразования представления

$$g : A_1 \longrightarrow A_2$$

мультипликативной $\Omega$-группы $A_1$ в $\Omega_2$-алгебре $A_2$ могут действовать на $A_2$-числа либо слева, либо справа. Мы должны рассматривать эндоморфизм $\Omega_2$-алгебры $A_2$ как оператор. При этом должны быть выполнены следующие условия.

О‌пределение 2.5.1. *Пусть* $\mathrm{End}(\Omega_2, A_2)$ *- мультипликативная $\Omega$-группа с произведением*[2.18]

$$(f, g) \to f \bullet g$$

*Пусть эндоморфизм $f$ действует на $A_2$-число а слева. Мы будем пользоваться записью*

(2.5.1) $$f(a_2) = f \bullet a_2$$

*Пусть $A_1$ - мультипликативная Ω-группа с произведением*

$$(a, b) \to a * b$$

*Мы будем называть гомоморфизм мультипликативной Ω-группы*

$$(2.5.2) \qquad f : A_1 \to \operatorname{End}(\Omega_2, A_2)$$

**левосторонним представлением** *мультипликативной Ω-группы $A_1$ или* **левосторонним** $A_1$-**представлением** *в $\Omega_2$-алгебре $A_2$, если отображение $f$ удовлетворяет условиям*

$$(2.5.3) \qquad f(a_1 * b_1) \bullet a_2 = (f(a_1) \bullet f(b_1)) \bullet a_2$$

*Мы будем отождествлять $A_1$-число $a_1$ с его образом $f(a_1)$ и записывать левостороннее преобразование, порождённое $A_1$-числом $a_1$, в форме*

$$a_2' = f(a_1) \bullet a_2 = a_1 * a_2$$

*В этом случае равенство (2.5.3) принимает вид*

$$(2.5.4) \qquad f(a_1 * b_1) \bullet a_2 = (a_1 * b_1) * a_2$$

*Отображение*

$$(a_1, a_2) \in A_1 \times A_2 \to a_1 * a_2 \in A_2$$

*порождённое левосторонним представлением $f$, называется* **левосторонним произведением** *$A_2$-числа $a_2$ на $A_1$-число $a_1$.*  $\square$

Из равенства (2.5.4) следует, что

$$(2.5.5) \qquad (a_1 * b_1) * a_2 = a_1 * (b_1 * a_2)$$

Мы можем рассматривать равенство (2.5.5) как **закон ассоциативности**,

Пример 2.5.2. *Пусть*

$$f : A_2 \to A_2$$

$$g : A_2 \to A_2$$

*эндоморфизмы $\Omega_2$-алгебры $A_2$. Пусть произведение в мультипликативной Ω-группе $\operatorname{End}(\Omega_2, A_2)$ является композицией эндоморфизмов. Так как произведение отображений $f$ и $g$ определено в том же порядке, как эти отображения действуют на $A_2$-число, то мы можем рассматривать равенство*

$$(2.5.6) \qquad (f \circ g) \circ a = f \circ (g \circ a)$$

*как закон ассоциативности, который позволяет записывать равенство (2.5.6) без использования скобок*

$$f \circ g \circ a = f \circ (g \circ a) = (f \circ g) \circ a$$

*а также записать равенство (2.5.3) в виде*

$$(2.5.7) \qquad f(a_1 * b_1) \circ a_2 = f(a_1) \circ f(b_1) \circ a_2$$

$\square$

---

2.18   Очень часто произведение в мультипликативной Ω-группе $\operatorname{End}(\Omega_2, A_2)$ совпадает с суперпозицией эндоморфизмов

$$f \bullet g = f \circ g$$

Однако произведение в мультипликативной Ω-группе $\operatorname{End}(\Omega_2, A_2)$ может отличаться от суперпозиции эндоморфизмов.

ЗАМЕЧАНИЕ 2.5.3. *Пусть отображение*

$$f : A_1 \xrightarrow{\;*\;} A_2$$

*является левосторонним представлением мультипликативной $\Omega$-группы $A_1$ в $\Omega_2$-алгебре $A_2$. Пусть отображение*

$$g : B_1 \xrightarrow{\;*\;} B_2$$

*является левосторонним представлением мультипликативной $\Omega$-группы $B_1$ в $\Omega_2$-алгебре $B_2$. Пусть отображение*

$$(r_1 : A_1 \to B_1, \; r_2 : A_2 \to B_2)$$

*является морфизмом представлений. Мы будем пользоваться записью*

$$r_2(a_2) = r_2 \circ a_2$$

*для образа $A_2$-числа $a_2$ при отображении $r_2$. Тогда мы можем записать равенство* (2.3.2) *следующим образом*

$$r_2 \circ (a_1 * a_2) = r_1(a_1) * (r_2 \circ a_2)$$

$\square$

ОПРЕДЕЛЕНИЕ 2.5.4. *Пусть* $\mathrm{End}(\Omega_2, A_2)$ *- мультипликативная $\Omega$-группа с произведением*[2.19]

$$(f, g) \to f \bullet g$$

*Пусть эндоморфизм $f$ действует на $A_2$-число $a$ справа. Мы будем пользоваться записью*

$$(2.5.8) \qquad\qquad f(a_2) = a_2 \bullet f$$

*Пусть $A_1$ - мультипликативная $\Omega$-группа с произведением*

$$(a, b) \to a * b$$

*Мы будем называть гомоморфизм мультипликативной $\Omega$-группы*

$$(2.5.9) \qquad\qquad f : A_1 \to \mathrm{End}(\Omega_2, A_2)$$

**правосторонним представлением** *мультипликативной $\Omega$-группы $A_1$ или* **правосторонним** *$A_1$-**представлением** в $\Omega_2$-алгебре $A_2$, если отображение $f$ удовлетворяет условиям*

$$(2.5.10) \qquad a_2 \bullet f(a_1 * b_1) = a_2 \bullet (f(a_1) \bullet f(b_1))$$

*Мы будем отождествлять $A_1$-число $a_1$ с его образом $f(a_1)$ и записывать правостороннее преобразование, порождённое $A_1$-числом $a_1$, в форме*

$$a_2' = a_2 \bullet f(a_1) = a_2 * a_1$$

*В этом случае равенство* (2.5.10) *принимает вид*

$$(2.5.11) \qquad a_2 \bullet f(a_1 * b_1) = a_2 * (a_1 * b_1)$$

*Отображение*

$$(a_1, a_2) \in A_1 \times A_2 \to a_2 * a_1 \in A_2$$

*порождённое правосторонним представлением $f$, называется* **правосторонним произведением** *$A_2$-числа $a_2$ на $A_1$-число $a_1$.*

$\square$

Из равенства (2.5.11) следует, что

(2.5.12) $$a_2 * (a_1 * b_1) = (a_2 * a_1) * b_1$$

Мы можем рассматривать равенство (2.5.12) как **закон ассоциативности**,

Пример 2.5.5. *Пусть*

$$f : A_2 \to A_2$$
$$g : A_2 \to A_2$$

*эндоморфизмы $\Omega_2$-алгебры $A_2$. Пусть произведение в мультипликативной $\Omega$-группе $\mathrm{End}(\Omega_2, A_2)$ является композицией эндоморфизмов. Так как произведение отображений $f$ и $g$ определено в том же порядке, как эти отображения действуют на $A_2$-число, то мы можем рассматривать равенство*

(2.5.13) $$a \circ (g \circ f) = (a \circ g) \circ f$$

*как закон ассоциативности, который позволяет записывать равенство (2.5.13) без использования скобок*

$$a \circ g \circ f = (a \circ g) \circ f = a \circ (g \circ f)$$

*а также записать равенство (2.5.10) в виде*

(2.5.14) $$a_2 \circ f(a_1 * b_1) = a_2 \circ f(a_1) \circ f(b_1)$$

<div align="right">□</div>

Замечание 2.5.6. *Пусть отображение*

$$f : A_1 \longrightarrow\!\!\!* A_2$$

*является правосторонним представлением мультипликативной $\Omega$-группы $A_1$ в $\Omega_2$-алгебре $A_2$. Пусть отображение*

$$g : B_1 \longrightarrow\!\!\!* B_2$$

*является правосторонним представлением мультипликативной $\Omega$-группы $B_1$ в $\Omega_2$-алгебре $B_2$. Пусть отображение*

$$(r_1 : A_1 \to B_1, \ r_2 : A_2 \to B_2)$$

*является морфизмом представлений. Мы будем пользоваться записью*

$$r_2(a_2) = r_2 \circ a_2$$

*для образа $A_2$-числа $a_2$ при отображении $r_2$. Тогда мы можем записать равенство (2.3.2) следующим образом*

$$r_2 \circ (a_2 * a_1) = (r_2 \circ a_2) * r_1(a_1)$$

<div align="right">□</div>

Если мультипликативная $\Omega$-группа $A_1$ - абелевая, то нет разницы между левосторонним и правосторонним представлениями.

---

2.19   Очень часто произведение в мультипликативной $\Omega$-группе $\mathrm{End}(\Omega_2, A_2)$ совпадает с суперпозицией эндоморфизмов

$$f \bullet g = f \circ g$$

Однако произведение в мультипликативной $\Omega$-группе $\mathrm{End}(\Omega_2, A_2)$ может отличаться от суперпозиции эндоморфизмов.

Определение 2.5.7. *Пусть $A_1$ - абелева мультипликативная $\Omega$-группа. Мы будем называть гомоморфизм мультипликативной $\Omega$-группы*

$$(2.5.15) \qquad f : A_1 \to \operatorname{End}(\Omega_2, A_2)$$

**представлением** *мультипликативной $\Omega$-группы $A_1$ или $A_1$-**представлением** в $\Omega_2$-алгебре $A_2$, если отображение $f$ удовлетворяет условиям*

$$(2.5.16) \qquad f(a_1 * b_1) \bullet a_2 = (f(a_1) \bullet f(b_1)) \bullet a_2$$

$\square$

Обычно мы отождествляем представление абелевой мультипликативной $\Omega$-группы $A_1$ и левостороннее представление мультипликативной $\Omega$-группы $A_1$. Однако, если это необходимо нам, мы можем отождествить представление абелевой мультипликативной $\Omega$-группы $A_1$ и правостороннее представление мультипликативной $\Omega$-группы $A_1$.

Определение 2.5.8. *Левостороннее представление*

$$f : A_1 \longrightarrow\!\!\!* A_2$$

*называется* **ковариантным**, *если равенство*

$$a_1 * (b_1 * a_2) = (a_1 * b_1) * a_2$$

*верно.*

$\square$

Определение 2.5.9. *Левостороннее представление*

$$f : A_1 \longrightarrow\!\!\!* A_2$$

*называется* **контравариантным**, *если равенство*

$$(2.5.17) \qquad a_1^{-1} * (b_1^{-1} * a_2) = (b_1 * a_1)^{-1} * a_2$$

*верно.*

$\square$

Теорема 2.5.10. *Произведение*

$$(a, b) \to a * b$$

*в мультипликативной $\Omega$-группе $A$ определяет два различных представления.*

- **Левый сдвиг**

$$a' = L(b) \circ a = b * a$$

  *является левосторонним представлением мультипликативной $\Omega$-группы $A$ в $\Omega$-алгебре $A$*

$$(2.5.18) \qquad L(c * b) = L(c) \circ L(b)$$

- **Правый сдвиг**

$$a' = a \circ R(b) = a * b$$

  *является правосторонним представлением мультипликативной $\Omega$-группы $A$ в $\Omega$-алгебре $A$*

$$(2.5.19) \qquad R(b * c) = R(b) \circ R(c)$$

Доказательство. Теорема является следствием теоремы [15]-5.2.1. □

---

ОПРЕДЕЛЕНИЕ 2.5.11. *Пусть $A_1$ является $\Omega$-группоидом с произведением*

$$(a, b) \to a * b$$

*Пусть отображение*

$$f : A_1 \xrightarrow{\;*\;} A_2$$

*является левосторонним представлением $\Omega$-группоида $A_1$ в $\Omega_2$-алгебре $A_2$. Для любого $a_2 \in A_2$, мы определим **орбиту представления** $\Omega$-группоида $A_1$ как множество*

$$A_1 * a_2 = \{b_2 = a_1 * a_2 : a_1 \in A_1\}$$

□

---

ОПРЕДЕЛЕНИЕ 2.5.12. *Пусть $A_1$ является $\Omega$-группоидом с произведением*

$$(a, b) \to a * b$$

*Пусть отображение*

$$f : A_1 \xrightarrow{\;*\;} A_2$$

*является правосторонним представлением $\Omega$-группоида $A_1$ в $\Omega_2$-алгебре $A_2$. Для любого $a_2 \in A_2$, мы определим **орбиту представления** $\Omega$-группоида $A_1$ как множество*

$$a_2 * A_1 = \{b_2 = a_2 * a_1 : a_1 \in A_1\}$$

□

---

ТЕОРЕМА 2.5.13. *Пусть отображение*

$$f : A_1 \xrightarrow{\;*\;} A_2$$

*является левосторонним представлением мультипликативной $\Omega$-группы $A_1$. Если*

$$(2.5.20) \qquad b_2 \in A_1 * a_2$$

*то*

$$(2.5.21) \qquad A_1 * a_2 = A_1 * b_2$$

---

Доказательство. Теорема является следствием теоремы [15]-5.3.7. □

Так как левостороннее представление и правостороннее представление опирается на гомоморфизм $\Omega$-группы, то верно следующее утверждение.

---

ТЕОРЕМА 2.5.14 (принцип двойственности для представления мультипликативной $\Omega$-группы). *Любое утверждение, справедливое для левостороннего представления мультипликативной $\Omega$-группы $A_1$, будет справедливо для правостороннего представления мультипликативной $\Omega$-группы $A_1$, если мы будем пользоваться правосторонним произведением на $A_1$-число $a_1$ вместо левостороннего произведения на $A_1$-число $a_1$.* □

Теорема 2.5.15. *Если существует однотранзитивное представление*

$$f : A_1 \longrightarrow A_2$$

*мультипликативной $\Omega$-группы $A_1$ в $\Omega_2$-алгебре $A_2$, то мы можем однозначно определить координаты на $A_2$, пользуясь $A_1$-числами.*

*Если $f$ - левостороннее однотранзитивное однотранзитивное представление, то $f(a)$ эквивалентно левому сдвигу $L(a)$ на группе $A_1$. Если $f$ - правостороннее однотранзитивное представление, то $f(a)$ эквивалентно правому сдвигу $R(a)$ на группе $A_1$.*

Доказательство. Теорема является следствием теоремы [15]-5.5.4.   $\square$

Теорема 2.5.16. *Правый и левый сдвиги на мультипликативной $\Omega$-группе $A_1$ перестановочны.*

Доказательство. Теорема является следствием ассоциативности произведения в мультипликативной $\Omega$-группе $A_1$

$$(L(a) \circ c) \circ R(b) = (a * c) * b = a * (c * b) = L(a) \circ (c \circ R(b))$$

$\square$

Теорема 2.5.17. *Пусть $A_1$ - мультипликативная $\Omega$-группа. Пусть левостороннее $A_1$-представление $f$ на $\Omega_2$-алгебре $A_2$ однотранзитивно. Тогда мы можем однозначно определить однотранзитивное правостороннее $A_1$-представление $h$ на $\Omega_2$-алгебре $A_2$ такое, что диаграмма*

$$
\begin{array}{ccc}
A_2 & \xrightarrow{\ h(a_1)\ } & A_2 \\
\downarrow{\scriptstyle f(b_1)} & & \downarrow{\scriptstyle f(b_1)} \\
A_2 & \xrightarrow{\ h(a_1)\ } & A_2
\end{array}
$$

*коммутативна для любых $a_1, b_1 \in A_1$.*

Доказательство. Мы будем пользоваться групповыми координатами для $A_2$-чисел $a_2$. Тогда согласно теореме 2.5.15 мы можем записать левый сдвиг $L(a_1)$ вместо преобразования $f(a_1)$.

Пусть $a_2, b_2 \in A_2$. Тогда мы можем найти одно и только одно $a_1 \in A_1$ такое, что

$$b_2 = a_2 * a_1 = a_2 \circ R(a_1)$$

Мы предположим

$$h(a) = R(a)$$

Существует $b_1 \in A_1$ такое, что

$$c_2 = f(b_1) \bullet a_2 = L(b_1) \circ a_2 \quad d_2 = f(b_1) \bullet b_2 = L(b_1) \circ b_2$$

Согласно теореме 2.5.16 диаграмма

(2.5.22)
$$
\begin{array}{ccc}
a_2 & \xrightarrow{\ h(a_1)=R(a_1)\ } & b_2 \\
{\scriptstyle f(b_1)=L(b_1)}\big\downarrow & & \big\downarrow{\scriptstyle f(b_1)=L(b_1)} \\
c_2 & \xrightarrow[\ h(a_1)=R(a_1)\ ]{} & d_2
\end{array}
$$

коммутативна.

Изменяя $b_1$ мы получим, что $c_2$ - это произвольное $A_2$-число.

Мы видим из диаграммы, что, если $a_2 = b_2$, то $c_2 = d_2$ и следовательно $h(e) = \delta$. С другой стороны, если $a_2 \neq b_2$, то $c_2 \neq d_2$ потому, что левостороннее $A_1$-представление $f$ однотранзитивно. Следовательно правостороннее $A_1$-представление $h$ эффективно.

Таким же образам мы можем показать, что для данного $c_2$ мы можем найти $a_1$ такое, что $d_2 = c_2 \bullet h(a_1)$. Следовательно правостороннее $A_1$-представление $h$ однотранзитивно.

В общем случае, произведение преобразований левостороннего $A_1$-представления $f$ не коммутативно и следовательно правосторонним $A_1$-представление $h$ отлично от левостороннего $A_1$-представления $f$. Таким же образом мы можем создать левостороннее $A_1$-представление $f$, пользуясь правосторонним $A_1$-представлением $h$. $\qquad\square$

Мы будем называть представления $f$ и $h$ **парными представлениями** мультипликативной Ω-группы $A_1$.

ЗАМЕЧАНИЕ 2.5.18. *Очевидно, что преобразования $L(a)$ и $R(a)$ отличаются, если мультипликативная Ω-группа $A_1$ неабелева. Тем не менее, они являются отображениями на. Теорема 2.5.17 утверждает, что, если оба представления правого и левого сдвига существуют на множестве $A_2$, то мы можем определить два перестановочных представления на множестве $A_2$. Только правый или левый сдвиг не может представлять оба типа представления. Чтобы понять почему это так, мы можем изменить диаграмму (2.5.22) и предположить*

$$ h(a_1) \bullet a_2 = L(a_1) \circ a_2 = b_2 $$

*вместо*

$$ a_2 \bullet h(a_1) = a_2 \circ R(a_1) = b_2 $$

*и проанализировать, какое выражение $h(a_1)$ имеет в точке $c_2$. Диаграмма*

$$
\begin{array}{ccc}
a_2 & \xrightarrow{\ h(a_1)=L(a_1)\ } & b_2 \\
{\scriptstyle f(b_1)=L(b_1)}\big\downarrow & & \big\downarrow{\scriptstyle f(b_1)=L(b_1)} \\
c_2 & \xrightarrow[\ h(a_1)\ ]{} & d_2
\end{array}
$$

*эквивалентна диаграмме*

$$
\begin{array}{ccc}
a_2 & \xrightarrow{\ h(a_1)=L(a_1)\ } & b_2 \\
{\scriptstyle (f(b_1))^{-1}=L(b_1^{-1})}\big\uparrow & & \big\downarrow{\scriptstyle f(b_1)=L(b_1)} \\
c_2 & \xrightarrow[\ h(a_1)\ ]{} & d_2
\end{array}
$$

*и мы имеем*  $d_2 = b_1 b_2 = b_1 a_1 a_2 = b_1 a_1 b_1^{-1} c_2$.  *Следовательно*

$$h(a_1) \bullet c_2 = (b_1 a_1 b_1^{-1}) c_2$$

*Мы видим, что представление* $h$ *зависит от его аргумента.*  □

## 2.6. Тензорное произведение представлений

ОПРЕДЕЛЕНИЕ 2.6.1. *Пусть* $A$, $B_1$, ..., $B_n$, $B$ *- универсальные алгебры. Пусть, для любого* $k$, $k = 1$, ..., $n$,

$$f_k : A \longrightarrow_* B_k$$

*эффективное представление* $\Omega_1$*-алгебры* $A$ *в* $\Omega_2$*-алгебре* $B_k$. *Пусть*

$$f : A \longrightarrow_* B$$

*эффективное представление* $\Omega_1$*-алгебры* $A$ *в* $\Omega_2$*-алгебре* $B$. *Отображение*

$$r_2 : B_1 \times ... \times B_n \to B$$

*называется* **приведенным полиморфизмом представлений** $f_1$, ..., $f_n$ *в представление* $f$, *если для любого* $k$, $k = 1$, ..., $n$, *при условии, что все переменные кроме переменной* $x_k \in B_k$ *имеют заданное значение, отображение* $r_2$ *является приведенным морфизмом представления* $f_k$ *в представление* $f$.

*Если* $f_1 = ... = f_n$, *то мы будем говорить, что отображение* $r_2$ *является приведенным полиморфизмом представления* $f_1$ *в представление* $f$.

*Если* $f_1 = ... = f_n = f$, *то мы будем говорить, что отображение* $r_2$ *является приведенным полиморфизмом представления* $f$.  □

ОПРЕДЕЛЕНИЕ 2.6.2. *Пусть* $A$ *является абелевой мультипликативной* $\Omega_1$*-группой. Пусть* $A_1$, ..., $A_n$ *-* $\Omega_2$*-алгебры.*[2.20] *Пусть для любого* $k$, $k = 1$, ..., $n$,

$$f_k : A \longrightarrow_* A_k$$

*эффективное представление мультипликативной* $\Omega_1$*-группы* $A$ *в* $\Omega_2$*-алгебре* $A_k$. *Рассмотрим категорию* $\mathcal{A}$ *объектами которой являются приведенные полиморфизмы представлений* $f_1$, ..., $f_n$

$$r_1 : A_1 \times ... \times A_n \longrightarrow S_1 \qquad r_2 : A_1 \times ... \times A_n \longrightarrow S_2$$

*где* $S_1$, $S_2$ *-* $\Omega_2$*-алгебры и*

$$g_1 : A \longrightarrow_* S_1 \qquad g_2 : A \longrightarrow_* S_2$$

*эффективные представления мультипликативной* $\Omega_1$*-группы* $A$. *Мы определим морфизм* $r_1 \to r_2$ *как приведенный морфизм представлений* $h : S_1 \to S_2$, *для которого коммутативна диаграмма*

*Универсальный объект* $B_1 \otimes ... \otimes B_n$ *категории* $\mathcal{A}$ *называется* **тензорным произведением** *представлений* $A_1, ..., A_n$. □

**Теорема 2.6.3.** *Если существует полиморфизм представлений, то тензорное произведение представлений существует.*

**Доказательство.** Теорема является следствием теоремы [15]-4.7.5. □

**Теорема 2.6.4.** *Пусть*

$$b_k \in B_k \quad k = 1, ..., n \quad b_{i\cdot1}, ..., b_{i\cdot p} \in B_i \quad \omega \in \Omega_2(p) \quad a \in A$$

*Тензорное произведение дистрибутивно относительно операции* $\omega$

$$(2.6.1) \qquad b_1 \otimes ... \otimes (b_{i\cdot1}...b_{i\cdot p}\omega) \otimes ... \otimes b_n$$
$$= (b_1 \otimes ... \otimes b_{i\cdot1} \otimes ... \otimes b_n)...(b_1 \otimes ... \otimes b_{i\cdot p} \otimes ... \otimes b_n)\omega$$

*Представление мультипликативной* $\Omega_1$-*группы* $A$ *в тензорном произведении определено равенством*

$$(2.6.2) \qquad b_1 \otimes ... \otimes (f_i(a) \circ b_i) \otimes ... \otimes b_n = f(a) \circ (b_1 \otimes ... \otimes b_i \otimes ... \otimes b_n)$$

**Доказательство.** Теорема является следствием теоремы [15]-4.7.11. □

**Теорема 2.6.5.** *Пусть* $B_1, ..., B_n$ - $\Omega_2$-*алгебры. Пусть*

$$f : A_1 \times ... \times A_n \to B_1 \otimes ... \otimes B_n$$

*приведенный полиморфизм, определённый равенством*

$$(2.6.3) \qquad f \circ (b_1, ..., b_n) = b_1 \otimes ... \otimes b_n$$

*Пусть*

$$g : A_1 \times ... \times A_n \to V$$

*приведенный полиморфизм в* $\Omega_2$-*алгебру* $V$. *Существует морфизм представлений*

$$h : B_1 \otimes ... \otimes B_n \to V$$

*такой, что диаграмма*

*коммутативна.*

**Доказательство.** Теорема является следствием теоремы [15]-4.7.10. □

---

[2.20] Я определяю тензорное произведение представлений универсальной алгебры по аналогии с определением в [1], с. 456 - 458.

Теорема 2.6.6. *Отображение*

$$(x_1, ..., x_n) \in B_1 \times ... \times B_n \to x_1 \otimes ... \otimes x_n \in B_1 \otimes ... \otimes B_n$$

*является полиморфизмом.*

Доказательство. Теорема является следствием теоремы [15]-4.7.9.  □

## 2.7. Базис представления универсальной алгебры

Определение 2.7.1. *Пусть*

$$f : A_1 \xrightarrow{\;\ast\;} A_2$$

*представление $\Omega_1$-алгебры $A_1$ в $\Omega_2$-алгебре $A_2$. Множество $B_2 \subset A_2$ называется* **стабильным множеством представления** $f$, *если $f(a)(m) \in B_2$ для любых $a \in A_1$, $m \in B_2$.*  □

Теорема 2.7.2. *Пусть*

$$f : A_1 \xrightarrow{\;\ast\;} A_2$$

*представление $\Omega_1$-алгебры $A_1$ в $\Omega_2$-алгебре $A_2$. Пусть множество $B_2 \subset A_2$ является подалгеброй $\Omega_2$-алгебры $A_2$ и стабильным множеством представления $f$. Тогда существует представление*

$$f_{B_2} : A_1 \xrightarrow{\;\ast\;} B_2$$

*такое, что $f_{B_2}(a) = f(a)|_{B_2}$. Представление $f_{B_2}$ называется* **подпредставлением** *представления $f$.*

Доказательство. Теорема является следствием теоремы [15]-6.1.2.  □

Теорема 2.7.3. *Множество[2.21] $\mathcal{B}_f$ всех подпредставлений представления $f$ порождает систему замыканий на $\Omega_2$-алгебре $A_2$ и, следовательно, является полной структурой.*

Доказательство. Теорема является следствием теоремы [15]-6.1.3.  □

Обозначим соответствующий оператор замыкания через $J[f]$. Таким образом, $J[f, X]$ является пересечением всех подалгебр $\Omega_2$-алгебры $A_2$, содержащих $X$ и стабильных относительно представления $f$.

Определение 2.7.4. $J[f, X]$ *называется* **подпредставлением**, *порождённым множеством $X$, а $X$ - множеством образующих подпредставления $J[f, X]$. В частности,* **множеством образующих** *представления $f$ будет такое подмножество $X \subset A_2$, что $J[f, X] = A_2$.*  □

Теорема 2.7.5. *Пусть[2.22]*

$$g : A_1 \xrightarrow{\;\ast\;} A_2$$

---

*представление $\Omega_1$-алгебры $A_1$ в $\Omega_2$-алгебре $A_2$. Пусть $X \subset A_2$. Определим подмножество $X_k \subset A_2$ индукцией по $k$.*

2.7.5.1: $X_0 = X$

2.7.5.2: $x \in X_k => x \in X_{k+1}$

2.7.5.3: $x_1 \in X_k, ..., x_n \in X_k, \omega \in \Omega_2(n) => x_1...x_n\omega \in X_{k+1}$

2.7.5.4: $x \in X_k, a \in A => f(a)(x) \in X_{k+1}$

*Тогда*

$$(2.7.1) \qquad\qquad \bigcup_{k=0}^{\infty} X_k = J[f, X]$$

Доказательство. Теорема является следствием теоремы [15]-6.1.4.    □

Определение 2.7.6. *Пусть $X \subset A_2$. Для любого $m \in J[f, X]$ существует $\Omega_2$-слово $w[f, X, m]$ , определённое согласно следующему правилам.*

2.7.6.1: *Если $m \in X$, то $m$ - $\Omega_2$-слово.*

2.7.6.2: *Если $m_1, ..., m_n$ - $\Omega_2$-слова и $\omega \in \Omega_2(n)$, то $m_1...m_n\omega$ - $\Omega_2$-слово.*

2.7.6.3: *Если $m$ - $\Omega_2$-слово и $a \in A_1$, то $f(a)(m)$ - $\Omega_2$-слово.*

*Мы будем отождествлять элемент $m \in J[f, X]$ и соответствующее ему $\Omega_2$-слово, выражая это равенством*

$$m = w[f, X, m]$$

*Аналогично, для произвольного множества $B \subset J[f, X]$ рассмотрим множество $\Omega_2$-слов[2.23]*

$$w[f, X, B] = \{w[f, X, m] : m \in B\}$$

*Мы будем также пользоваться записью*

$$w[f, X, B] = (w[f, X, m], m \in B)$$

*Обозначим $w[f, X]$ **множество $\Omega_2$-слов представления** $J[f, X]$.*    □

Теорема 2.7.7. *Пусть*

$$f : A_1 \longrightarrow\!\!\!* A_2$$

*представление $\Omega_1$-алгебры $A_1$ в $\Omega_2$-алгебре $A_2$. Пусть*

$$g : A_1 \longrightarrow\!\!\!* B_2$$

*представление $\Omega_1$-алгебры $A_1$ в $\Omega_2$-алгебре $B_2$. Пусть $X$ - множество образующих представления $f$. Пусть*

$$R : A_2 \to B_2$$

---

2.22  Утверждение теоремы аналогично утверждению теоремы 5.1, [18], страница 94.

2.23  Выражение $w[f, X, m]$ является частным случаем выражения $w[f, X, B]$, а именно

$$w[f, X, \{m\}] = \{w[f, X, m]\}$$

*приведенный морфизм представления*[2.24] и $X' = R(X)$. *Приведенный морфизм
$R$ представления порождает отображение $\Omega_2$-слов*

$$w[f \to g, X, R] : w[f, X] \to w[g, X']$$

*такое, что*

2.7.7.1: *Если* $m \in X$, $m' = R(m)$, *то*

$$w[f \to g, X, R](m) = m'$$

2.7.7.2: *Если*

$$m_1, ..., m_n \in w[f, X]$$

$$m'_1 = w[f \to g, X, R](m_1) \quad ... \quad m'_n = w[f \to g, X, R](m_n)$$

*то для операции* $\omega \in \Omega_2(n)$ *справедливо*

$$w[f \to g, X, R](m_1...m_n\omega) = m'_1...m'_n\omega$$

2.7.7.3: *Если*

$$m \in w[f, X] \quad m' = w[f \to g, X, R](m) \quad a \in A_1$$

*то*

$$w[f \to g, X, R](f(a)(m)) = g(a)(m')$$

Доказательство. Теорема является следствием теоремы [15]-6.1.10. □

Определение 2.7.8. *Пусть*

$$f : A_1 \longrightarrow\!\!\!\ast\!\!\!\longrightarrow A_2$$

*представление* $\Omega_1$-*алгебры* $A_1$ *в* $\Omega_2$-*алгебре* $A_2$ *и*

$$Gen[f] = \{X \subseteq A_2 : J[f, X] = A_2\}$$

*Если для множества* $X \subset A_2$ *верно* $X \in Gen[f]$, *то для любого множества
$Y$, $X \subset Y \subset A_2$, также верно* $Y \in Gen[f]$. *Если существует минималь-
ное множество* $X \in Gen[f]$, *то множество $X$ называется* **квазибазисом**
*представления* $f$. □

Теорема 2.7.9. *Если множество $X$ является квазибазисом представле-
ния $f$, то для любого $m \in X$ множество $X \setminus \{m\}$ не является множеством
образующих представления $f$.*

Доказательство. Теорема является следствием теоремы [15]-6.2.2. □

Замечание 2.7.10. *Доказательство теоремы 2.7.9 даёт нам эффектив-
ный метод построения квазибазиса представления $f$. Выбрав произвольное
множество образующих, мы шаг за шагом исключаем те элементы мно-
жества, которые имеют координаты относительно остальных элементов
множества. Если множество образующих представления бесконечно, то рас-
смотренная операция может не иметь последнего шага. Если представление
имеет конечное множество образующих, то за конечное число шагов мы мо-
жем построить квазибазис этого представления.* □

---

[2.24] Я рассмотрел морфизм представления в теореме [15]-8.1.7.

Мы ввели $\Omega_2$-слово элемента $x \in A_2$ относительно множества образующих $X$ в определении 2.7.6. Из теоремы 2.7.9 следует, что если множество образующих $X$ не является квазибазисом, то выбор $\Omega_2$-слова относительно множества образующих $X$ неоднозначен. Но даже если множество образующих $X$ является квазибазисом, то представление $m \in A_2$ в виде $\Omega_2$-слова неоднозначно.

Замечание 2.7.11. *Существует три источника неоднозначности в записи* $\Omega_2$*-слова.*

2.7.11.1: *В* $\Omega_i$*-алгебре* $A_i$, $i = 1, 2$, *могут быть определены равенства. Например, если* $e$ *- единица мультипликативной группы* $A_i$, *то верно равенство*

$$ae = a$$

*для любого* $a \in A_i$.

2.7.11.2: *Неоднозначность выбора* $\Omega_2$*-слова может быть связана со свойствами представления. Например, если* $m_1, ..., m_n$ *-* $\Omega_2$*-слова,* $\omega \in \Omega_2(n)$ *и* $a \in A_1$, *то*[2.25]

$$(2.7.2) \qquad f(a)(m_1...m_n\omega) = (f(a)(m_1))...(f(a)(m_n))\omega$$

*В тоже время, если* $\omega$ *является операцией* $\Omega_1$*-алгебры* $A_1$ *и операцией* $\Omega_2$*-алгебры* $A_2$, *то мы можем потребовать, что* $\Omega_2$*-слова* $f(a_1...a_n\omega)(x)$ *и* $(f(a_1)(x))...(f(a_n)(x))\omega$ *описывают один и тот же элемент* $\Omega_2$*-алгебры* $A_2$.[2.26]

$$(2.7.3) \qquad f(a_1...a_n\omega)(x) = (f(a_1)(x))...(f(a_n)(x))\omega$$

2.7.11.3: *Равенства вида* (2.7.2), (2.7.3) *сохраняются при морфизме представления. Поэтому мы можем игнорировать эту форму неоднозначности записи* $\Omega_2$*-слова. Однако возможна принципиально другая форма неоднозначности, пример которой можно найти в теоремах [15]-9.7.23, [15]-9.7.24.*

*Таким образом, мы видим, что на множестве* $\Omega_2$*-слов можно определить различные отношения эквивалентности.*[2.27] *Наша задача - найти максимальное отношение эквивалентности на множестве* $\Omega_2$*-слов, которое сохраняется при морфизме представления.* □

Определение 2.7.12. *Множество образующих* $X$ *представления* $f$ *порождает отношение эквивалентности*

$$\rho[f, X] = \{(w[f, X, m], w_1[f, X, m]) : m \in A_2\}$$

*на множестве* $\Omega_2$*-слов.* □

---

[2.25]Например, пусть $\{e_1, e_2\}$ - базис векторного пространства над полем $k$. Равенство (2.7.2) принимает форму закона дистрибутивности

$$a(b^1e_1 + b^2e_2) = (ab^1)e_1 + (ab^2)e_2$$

[2.26]Для векторного пространства это требование принимает форму закона дистрибутивности

$$(a + b)e_1 = ae_1 + be_1$$

[2.27]Очевидно, что каждое из равенств (2.7.2), (2.7.3) порождает некоторое отношение эквивалентности.

Теорема 2.7.13. *Пусть $X$ - квазибазис представления*

$$f : A_1 \longrightarrow A_2$$

*Для любого множества образующих $X$ представления $f$ существует отношение эквивалентности*

$$\lambda[f, X] \subseteq w[f, X] \times w[f, X]$$

*которое порождено исключительно следующими утверждениями.*

2.7.13.1: *Если в $\Omega_2$-алгебре $A_2$ существует равенство*

$$w_1[f, X, m] = w_2[f, X, m]$$

*определяющее структуру $\Omega_2$-алгебры, то*

$$(w_1[f, X, m], w_2[f, X, m]) \in \lambda[f, X]$$

2.7.13.2: *Если в $\Omega_1$-алгебре $A_1$ существует равенство*

$$w_1[f, X, m] = w_2[f, X, m]$$

*определяющее структуру $\Omega_1$-алгебры, то*

$$(f(w_1)(w[f, X, m]), f(w_2)(w[f, X, m])) \in \lambda[f, X]$$

2.7.13.3: *Для любой операции $\omega \in \Omega_1(n)$,*

$$(f(a_{11}...a_{1n}\omega)(a_2), (f(a_{11})...f(a_{1n})\omega)(a_2)) \in \lambda[f, X]$$

2.7.13.4: *Для любой операции $\omega \in \Omega_2(n)$,*

$$(f(a_1)(a_{21}...a_{2n}\omega), f(a_1)(a_{21})...f(a_1)(a_{2n})\omega) \in \lambda[f, X]$$

2.7.13.5: *Пусть $\omega \in \Omega_1(n) \cap \Omega_2(n)$. Если представление $f$ удовлетворяет равенству[2.28]*

$$f(a_{11}...a_{1n}\omega)(a_2) = (f(a_{11})(a_2))...(f(a_{1n})(a_2))\omega$$

*то мы можем предположить, что верно равенство*

$$(f(a_{11}...a_{1n}\omega)(a_2), (f(a_{11})(a_2))...(f(a_{1n})(a_2))\omega) \in \lambda[f, X]$$

Доказательство. Теорема является следствием теоремы [15]-6.2.5.    □

Определение 2.7.14. *Квазибазис $\overline{\overline{e}}$ представления $f$ такой, что*

$$\rho[f, \overline{\overline{e}}] = \lambda[f, \overline{\overline{e}}]$$

*называется* **базисом представления** *$f$.*    □

---

2.28   Рассмотрим представление коммутативного кольца $D$ в $D$-алгебре $A$. Мы будем пользоваться записью

$$f(a)(v) = av$$

В обеих алгебрах определены операции сложения и умножения. Однако равенство

$$f(a + b)(v) = f(a)(v) + f(b)(v)$$

верно, а равенство

$$f(ab)(v) = f(a)(v)f(b)(v)$$

является ошибочным.

Теорема 2.7.15. *Пусть $\overline{\overline{e}}$ - базис представления $f$. Пусть*

$$R_1 : \overline{\overline{e}} \to \overline{\overline{e}}'$$

*произвольное отображение множества $\overline{\overline{e}}$. Рассмотрим отображение $\Omega_2$-слов*

$$w[f \to g, \overline{\overline{e}}, \overline{\overline{e}}', R_1] : w[f, \overline{\overline{e}}] \to w[g, \overline{\overline{e}}']$$

*удовлетворяющее условиям* 2.7.7.1, 2.7.7.2, 2.7.7.3, *и такое, что*

$$e \in \overline{\overline{e}} => w[f \to g, \overline{\overline{e}}, \overline{\overline{e}}', R_1](e) = R_1(e)$$

*Существует единственный эндоморфизм представления $f$* [2.29]

$$r_2 : A_2 \to A_2$$

*определённый правилом*

$$R(m) = w[f \to g, \overline{\overline{e}}, \overline{\overline{e}}', R_1](w[f, \overline{\overline{e}}, m])$$

Доказательство. Теорема является следствием теоремы [15]-6.2.10.  □

## 2.8. Диаграмма представлений универсальных алгебр

Определение 2.8.1. **Диаграмма $(f, A)$ представлений универсальных алгебр** *- это такой ориентированный граф, что*

2.8.1.1: *вершина $A_k$, $k = 1$, ..., $n$, является $\Omega_k$-алгеброй;*
2.8.1.2: *ребро $f_{kl}$ является представлением $\Omega_k$-алгебры $A_k$ в $\Omega_l$-алгебре $A_l$;*
*Мы будем требовать, чтобы этот граф был связным и не содержал циклов. Мы будем полагать, что $A_{[0]}$ - это множество начальных вершин графа. Мы будем полагать, что $A_{[k]}$ - это множество вершин графа, для которых максимальный путь от начальных вершин равен $k$.*  □

Замечание 2.8.2. *Так как в разных вершинах графа может быть одна и та же алгебра, то мы обозначим $A = (A_{(1)} \ ... \ A_{(n)})$ множество универсальных алгебр, которые попарно различны. Из равенства*

$$A = (A_{(1)} \ ... \ A_{(n)}) = (A_1 \ ... \ A_n)$$

*следует, что для любого индекса $(i)$ существует по крайней мере один индекс $i$ такой, что $A_{(i)} = A_i$. Если даны два набора множеств $A = (A_{(1)} \ ... \ A_{(n)})$, $B = (B_{(1)} \ ... \ B_{(n)})$ и определено отображение*

$$h_{(i)} : A_{(i)} \to B_{(i)}$$

*для некоторого индекса $(i)$, то также определено отображение*

$$h_i : A_i \to B_i$$

*для любого индекса $i$ такого, что $A_{(i)} = A_i$ и в этом случае $h_i = h_{(i)}$.*  □

Определение 2.8.3. *Диаграмма $(f, A)$ представлений универсальных алгебр называется* **коммутативной**, *если выполнено следующее условие. для каждой пары представлений*

$$f_{ik} : A_i \longrightarrow A_k$$

---

[2.29] Это утверждение похоже на теорему [1]-1, с. 104.

$$f_{jk} : A_j \longrightarrow A_k$$

*следующее равенство верно*[2.30]

(2.8.1) $$f_{ik}(a_i)(f_{jk}(a_j)(a_k)) = f_{jk}(a_j)(f_{ik}(a_i)(a_k))$$

$\square$

Мы можем проиллюстрировать определение 2.8.1 с помощью теоремы 2.8.4.

Теорема 2.8.4. *Пусть*

$$f_{ij} : A_i \longrightarrow A_j$$

*представление $\Omega_i$-алгебры $A_i$ в $\Omega_j$-алгебре $A_j$. Пусть*

$$f_{jk} : A_j \longrightarrow A_k$$

*представление $\Omega_j$-алгебры $A_j$ в $\Omega_k$-алгебре $A_k$. Мы можем описать фрагмент*[2.31]

$$A_i \xrightarrow{f_{ij}} A_j \xrightarrow{f_{jk}} A_k$$

*диаграммы представлений с помощью диаграммы*

(2.8.2)

*Отображение*

$$f_{ijk} : A_i \to \mathrm{End}(\Omega_j, \mathrm{End}(\Omega_k, A_k))$$

*определено равенством*

(2.8.3) $$f_{ijk}(a_i)(f_{jk}(a_j)) = f_{jk}(f_{ij}(a_i)(a_j))$$

*где $a_i \in A_i$, $a_j \in A_j$. Если представление $f_{jk}$ эффективно и представление $f_{ij}$ свободно, то отображение $f_{ijk}$ является свободным представлением*

$$f_{ijk} : A_i \longrightarrow \mathrm{End}(\Omega_k, A_k)$$

*$\Omega_i$-алгебры $A_i$ в $\Omega_j$-алгебре $\mathrm{End}(\Omega_k, A_k)$.*

---

[2.30]  Образно говоря, представления $f_{ik}$ и $f_{jk}$ прозрачны друг для друга.

Доказательство. Теорема является следствием теоремы [15]-7.1.6. □

Определение 2.8.5. *Пусть $(f, A)$ - диаграмма представлений, где $A = (A_{(1)} \ ... \ A_{(n)})$ - множество универсальных алгебр. Пусть $(B, g)$ - диаграмма представлений, где $B = (B_{(1)} \ ... \ B_{(n)})$ - множество универсальных алгебр. Множество отображений $h = (h_{(1)} \ ... \ h_{(n)})$*

$$h_{(i)} : A_{(i)} \rightarrow B_{(i)}$$

*называется* **морфизмом из диаграммы представлений $(f, A)$ в диаграмму представлений $(B, g)$,** *если для любых индексов $(i)$, $(j)$, $i$, $j$ таких, что $A_{(i)} = A_i$, $A_{(j)} = A_j$, и для каждого представления*

$$f_{ji} : A_j \longrightarrow A_i$$

*пара отображений $(h_j \ \ h_i)$ является морфизмом представлений из $f_{ji}$ в $g_{ji}$.* □

Для любого представления $f_{ij}$, $i = 1, ..., n$, $j = 1, ..., n$, мы имеем диаграмму

(2.8.4)

Равенства

(2.8.5)
$$h_j \circ f_{ij}(a_i) = g_{ij}(h_i(a_i)) \circ h_j$$

(2.8.6)
$$h_{ij}(f_{ij}(a_i)(a_j)) = g_{ij}(h_i(a_i))(h_j(a_j))$$

выражают коммутативность диаграммы (1).

---

2.31  Теорема 2.8.4 утверждает, что преобразования в башне представлений согласованы.

# Абелева группа

## 3.1. Группа

ОПРЕДЕЛЕНИЕ 3.1.1. *Пусть на множестве $A$ определена бинарная операция* [3.1]

$$(a, b) \in A \times A \to ab \in A$$

*которую мы называем произведением. Множество $A$ называется мультипликативным моноидом, если*

3.1.1.1: *произведение* **ассоциативно**

$$(ab)c = a(bc)$$

3.1.1.2: *произведение имеет единичный элемент $e$*

$$ea = ae = a$$

*Если произведение* **коммутативно**

$$(3.1.1) \qquad ab = ba$$

*то мультипликативный моноид $A$ называется коммутативным или абелевым.* □

ОПРЕДЕЛЕНИЕ 3.1.2. *Пусть на множестве $A$ определена бинарная операция* [3.1]

$$(a, b) \in A \times A \to a + b \in A$$

*которую мы называем суммой. Множество $A$ называется аддитивным моноидом, если*

3.1.2.1: *сумма* **ассоциативна**

$$(a + b) + c = a + (b + c)$$

3.1.2.2: *сумма имеет нулевой элемент $0$*

$$0 + a = a + 0 = a$$

*Если сумма* **коммутативна**

$$(3.1.2) \qquad a + b = b + a$$

*то аддитивный моноид $A$ называется коммутативным или абелевым.* □

ТЕОРЕМА 3.1.3. *Единичный элемент мультипликативного моноида $A$ определён однозначно.*

ДОКАЗАТЕЛЬСТВО. Пусть $e_1$, $e_2$ - единичные элементы мультипликативного моноида $A$. Равенство

$$(3.1.3) \qquad e_1 = e_1 e_2 = e_2$$

следует из утверждения 3.1.1.2. Теорема следует из равенства (3.1.3). □

ЗАМЕЧАНИЕ 3.1.5. *Согласно определениям 2.1.6, 3.1.1 и теореме 3.1.3*

ТЕОРЕМА 3.1.4. *Нулевой элемент аддитивного моноида $A$ определён однозначно.*

ДОКАЗАТЕЛЬСТВО. Пусть $0_1$, $0_2$ - нулевые элементы аддитивного моноида $A$. Равенство

$$(3.1.4) \qquad 0_1 = 0_1 + 0_2 = 0_2$$

следует из утверждения 3.1.2.2. Теорема следует из равенства (3.1.4). □

ЗАМЕЧАНИЕ 3.1.6. *Согласно определениям 2.1.6, 3.1.2 и теореме 3.1.4*

---

[3.1] Смотри также определения моноида на странице [1]-17.

мультипликативный моноид - это универсальная алгебра с одной бинарной операцией (умножение) и одной 0-арной операцией (единичный элемент). □

аддитивный моноид - это универсальная алгебра с одной бинарной операцией (сложение) и одной 0-арной операцией (нулевой элемент). □

Определение 3.1.7. *Мультипликативный моноид* [3.2] *$A$ называется мультипликативной группой, если для каждого $A$-числа $x$ существует $A$-число $y$ такое, что*

$$xy = yx = e$$

*$A$-число $y$ называется обратным к $A$-числу $x$. Если произведение коммутативно, то мультипликативная группа $A$ называется коммутативной или абелевой.* □

Определение 3.1.8. *Аддитивный моноид* [3.2] *$A$ называется аддитивной группой, если для каждого $A$-числа $x$ существует $A$-число $y$ такое, что*

$$x + y = y + x = 0$$

*$A$-число $y$ называется обратным к $A$-числу $x$. Если сумма коммутативно, то аддитивная группа $A$ называется коммутативной или абелевой.* □

Определения мультипликативного и аддитивного моноида отличаются только формой записи операции. Поэтому, если выбор операции ясен из контекста, мы будем соответствующую алгебру называть моноидом. Аналогичное замечание верно для группы. Определения мультипликативной и аддитивной групп отличаются только формой записи операции. Поэтому, если выбор операции ясен из контекста, мы будем соответствующую алгебру называть группой. Если не оговоренно противное, обычно аддитивная группа предполагается абелевой.

Теорема 3.1.9. *Пусть $A$ - мультипликативная группа. $A$-число $y$, обратное к $A$-числу $a$, определено однозначно. Мы пользуемся записью $y = a^{-1}$*

$$(3.1.5) \qquad aa^{-1} = e$$

Доказательство. Пусть $y_1$, $y_2$ - обратные числа к $A$-числу $a$. Равенство

$$(3.1.6) \qquad \begin{aligned} y_1 &= y_1 e = y_1(ay_2) \\ &= (y_1 a)y_2 = ey_2 = y_2 \end{aligned}$$

является следствием утверждений 3.1.1.1, 3.1.1.2 и определения 3.1.7. Теорема является следствием равенства (3.1.6). □

Теорема 3.1.10. *Пусть $A$ - аддитивная группа. $A$-число $y$, обратное к $A$-числу $a$, определено однозначно. Мы пользуемся записью $y = -a$*

$$(3.1.7) \qquad a + (-a) = 0$$

Доказательство. Пусть $y_1$, $y_2$ - обратные числа к $A$-числу $a$. Равенство

$$(3.1.8) \qquad \begin{aligned} y_1 &= y_1 + 0 = y_1 + (a + y_2) \\ &= (y_1 + a) + y_2 = 0 + y_2 \\ &= y_2 \end{aligned}$$

является следствием утверждений 3.1.2.1, 3.1.2.2 и утверждения 3.1.8. Теорема является следствием равенства (3.1.8). □

---

[3.2] Смотри также определение группы на странице [1]-21.

ТЕОРЕМА 3.1.11. *Если $A$ - муль-типликативная группа, то*

$$(3.1.9) \qquad (a^{-1})^{-1} = a$$

ДОКАЗАТЕЛЬСТВО. Согласно определению (3.1.5), $A$-число $y$, обратное к $A$-числу $a^{-1}$, иммет вид $(a^{-1})^{-1}$

$$(3.1.10) \qquad (a^{-1})^{-1}a^{-1} = e$$

равенство (3.1.9) является следствием равенств (3.1.5), (3.1.10) и теоремы 3.1.9.                                    □

ТЕОРЕМА 3.1.12. *Если $A$ - аддитивная группа, то*

$$(3.1.11) \qquad -(-a) = a$$

ДОКАЗАТЕЛЬСТВО. Согласно определению (3.1.7), $A$-число $y$, обратное к $A$-числу $-a$, иммет вид $-(-a)$

$$(3.1.12) \qquad (-(-a)) + (-a) = 0$$

равенство (3.1.11) является следствием равенств (3.1.7), (3.1.12) и теоремы 3.1.10.                                    □

ТЕОРЕМА 3.1.13. *Пусть $A$ - муль-типликативная группа с единичным элементом $e$. $A$-число, обратное к единичному элементу, является единичным элементом*

$$(3.1.13) \qquad e^{-1} = e$$

ДОКАЗАТЕЛЬСТВО. Равенство

$$(3.1.14) \qquad e^{-1} = e^{-1}e = e$$

является следствием теорем 3.1.3, 3.1.9. Равенство (3.1.13) является следствием равенства (3.1.14).                 □

ТЕОРЕМА 3.1.14. *Пусть $A$ - аддитивная группа с нулевым элементом 0. $A$-число, обратное к нулевому элементу, является нулевым элементом*

$$(3.1.15) \qquad -0 = 0$$

ДОКАЗАТЕЛЬСТВО. Равенство

$$(3.1.16) \qquad -0 = (-0) + 0 = 0$$

является следствием теорем 3.1.4, 3.1.10. Равенство (3.1.15) является следствием равенства (3.1.16).                 □

ОПРЕДЕЛЕНИЕ 3.1.15. *Пусть $H$ - подгруппа мультипликативной группы $G$. Множество*

$$(3.1.17) \qquad aH = \{ab : b \in H\}$$

*называется* **левым смежным классом**. *Множество*

$$Ha = \{ba : b \in H\}$$

*называется* **правым смежным классом**.                               □

ОПРЕДЕЛЕНИЕ 3.1.16. *Пусть $H$ - подгруппа аддитивной группы $G$. Множество*

$$(3.1.18) \qquad a + H = \{a + b : b \in H\}$$

*называется* **левым смежным классом**. *Множество*

$$H + a = \{b + a : b \in H\}$$

*называется* **правым смежным классом**.                               □

ТЕОРЕМА 3.1.17. *Пусть $H$ - подгруппа мультипликативной группы $G$. Если*

$$(3.1.19) \qquad aH \cap bH \neq \emptyset$$

*то*

$$aH = bH$$

ДОКАЗАТЕЛЬСТВО. Из утверждения (3.1.19) и определения (3.1.17) сле-

ТЕОРЕМА 3.1.18. *Пусть $H$ - подгруппа аддитивной группы $G$. Если*

$$(3.1.21) \qquad a + H \cap b + H \neq \emptyset$$

*то*

$$a + H = b + H$$

ДОКАЗАТЕЛЬСТВО. Из утверждения (3.1.21) и определения (3.1.18) следует, что существуют $H$-числа $c$, $d$ та-дует, что существуют $H$-числа $c$, $d$ такие, что

(3.1.20) $\qquad ac = bd$

3.1.17.1: Согласно теореме 3.1.9, существует $H$-число $d^{-1}$. Следовательно,

$$cd^{-1} \in H$$

3.1.17.2: Если $c \in H$, то из определения (3.1.17) следует, что $cH = H$.

3.1.17.3: Из утверждений 3.1.17.1, 3.1.17.2 следует, что

$$(cd^{-1})H = H$$

3.1.17.4: Из утверждения 3.1.1.1 и определения 3.1.7 следует, что

$$(ab)H = a(bH)$$

3.1.17.5: Равенство

$$a(cd^{-1}) = (ac)d^{-1} = (bd)d^{-1}$$
$$= b(dd^{-1}) = be = b$$

является следствием утверждений 3.1.1.1, 3.1.1.2 и равенств (3.1.5), (3.1.20).

Равенство

$$bH = a(cd^{-1})H = aH$$

является следствием утверждений 3.1.17.3, 3.1.17.5. $\qquad\square$

кие, что

(3.1.22) $\qquad a + c = b + d$

3.1.18.1: Согласно теореме 3.1.10, существует $H$-число $-d$. Следовательно,

$$c + (-d) \in H$$

3.1.18.2: Если $c \in H$, то из определения (3.1.18) следует, что $c + H = H$.

3.1.18.3: Из утверждений 3.1.18.1, 3.1.18.2 следует, что

$$(c + (-d)) + H = H$$

3.1.18.4: Из утверждения 3.1.2.1 и определения 3.1.8 следует, что

$$(a + b) + H = a + (b + H)$$

3.1.18.5: Равенство

$$a + (c + (-d)) = (a + c) + (-d)$$
$$= (b + d) + (-d)$$
$$= b + (d + (-d))$$
$$= b + 0 = b$$

является следствием утверждений 3.1.2.1, 3.1.2.2 и равенств (3.1.7), (3.1.22).

Равенство

$$b + H = (a + (c + (-d))) + H = a + H$$

является следствием утверждений 3.1.18.3, 3.1.18.5. $\qquad\square$

---

ОПРЕДЕЛЕНИЕ 3.1.19. *Подгруппа $H$ мультипликативной группы $A$ такая, что*

$$aH = Ha$$

*для любого $A$-числа, называется* **нормальной**. 3.3 $\qquad\square$

ОПРЕДЕЛЕНИЕ 3.1.20. *Подгруппа $H$ аддитивной группы $A$ такая, что*

$$a + H = H + a$$

*для любого $A$-числа, называется* **нормальной**. 3.3 $\qquad\square$

---

Согласно теореме 3.1.17, нормальная подгруппа $H$ мультипликативной группы $A$ порождает отношение эквивалентности mod $H$ на множестве $A$

(3.1.23) $\quad a \equiv b (\mathrm{mod}\, H) \Leftrightarrow aH = bH$

Мы будем обозначать $A/H$ множество

Согласно теореме 3.1.18, нормальная подгруппа $H$ аддитивной группы $A$ порождает отношение эквивалентности mod $H$ на множестве $A$

(3.1.24)

$\quad a \equiv b (\mathrm{mod}\, H) \Leftrightarrow a + H = b + H$

Мы будем обозначать $A/H$ множество

---

классов эквивалентности отношения эквивалентности mod $H$.

классов эквивалентности отношения эквивалентности mod $H$.

**Теорема 3.1.21.** *Пусть $H$ - нормальная подгруппа мультипликативной группы $A$. Произведение*

$$(3.1.25) \qquad (aH)(bH) = abH$$

*порождает группу на множестве $A/H$.*

**Теорема 3.1.22.** *Пусть $H$ - нормальная подгруппа аддитивной группы $A$. Сумма*

$$(3.1.27) \quad (a+H)+(b+H) = a+b+H$$

*порождает группу на множестве $A/H$.*

Доказательство. Равенство

$$\begin{aligned}
(3.1.26) \qquad & (aH)(bH) \\
&= a(Hb)H \\
&= a(bH)H \\
&= (ab)(HH) \\
&= abH
\end{aligned}$$

следует из определения 3.1.19 и утверждений 3.1.1.1, 3.1.17.4.

Равенство (3.1.25) является следствием равенства (3.1.26).

Пусть $c \in aH$, $d \in bH$. Тогда равенство

$$\begin{aligned}
(ab)H &= (aH)(bH) \\
&= (cH)(dH) \\
&= (cd)H
\end{aligned}$$

следует из определения (3.1.25) и теоремы 3.1.17. Следовательно, произведение (3.1.25) не зависит от выбора $A$-чисел $c \in aH$, $d \in bH$.

Так как $eH = H$, то из равенства (3.1.25) следует, что произведение (3.1.25) удовлетворяет определению мультипликативной группы.                □
.

Доказательство. Равенство

$$\begin{aligned}
(3.1.28) \qquad & (a+H)+(b+H) \\
&= a+(H+b)+H \\
&= a+(b+H)+H \\
&= (a+b)+(H+H) \\
&= a+b+H
\end{aligned}$$

следует из определения 3.1.20 и утверждений 3.1.2.1, 3.1.18.4.

Равенство (3.1.27) является следствием равенства (3.1.28).

Пусть $c \in a+H$, $d \in b+H$. Тогда равенство

$$\begin{aligned}
(a+b)+H &= (a+H)+(b+H) \\
&= (c+H)+(d+H) \\
&= (c+d)+H
\end{aligned}$$

следует из определения (3.1.27) и теоремы 3.1.18. Следовательно, произведение (3.1.27) не зависит от выбора $A$-чисел $c \in a+H$, $d \in b+H$.

Так как $e + H = H$, то из равенства (3.1.27) следует, что произведение (3.1.27) удовлетворяет определению аддитивной группы.                □

**Определение 3.1.23.** *Пусть $H$ - нормальная подгруппа мультипликативной группы $A$. Группа $A/H$ называется* **факторгруппой** *группы $A$ по $H$ и отображение*

$$f : a \in A \to aH \in A/H$$

*называется* **каноническим отображением.**                □

**Определение 3.1.24.** *Пусть $H$ - нормальная подгруппа аддитивной группы $A$. Группа $A/H$ называется* **факторгруппой** *группы $A$ по $H$ и отображение*

$$f : a \in A \to a+H \in A/H$$

*называется* **каноническим отображением.**                □

## 3.2. Кольцо

Определение 3.2.1. *Пусть на множестве $A$ определены две бинарные операции* [3.4]
- *сумма*
$$(a, b) \in A \times A \to a + b \in A$$
- *произведение*
$$(a, b) \in A \times A \to ab \in A$$

3.2.1.1: *Множество $A$ является абелевой аддитивной группой относительно суммы и $0$ - нулевой элемент.*

3.2.1.2: *Множество $A$ является мультипликативным моноидом относительно произведения.*

3.2.1.3: *Произведение* **дистрибутивно** *относительно сложения*
$$(3.2.1) \qquad a(b + c) = ab + ac$$
$$(3.2.2) \qquad (a + b)c = ac + bc$$

*Если умножение коммутативно, то кольцо называется коммутативным.* □

Определение 3.2.2. *Если мультипликативный моноид кольца $D$ содержит единичный элемент $1$, то кольцо $D$ называется* **унитальным кольцом**. [3.5] □

Теорема 3.2.3. *Пусть $A$ - кольцо. Тогда*
$$(3.2.3) \qquad 0a = 0$$
*для любого $A$-числа $a$.*

Доказательство. Равенство
$$(3.2.4) \qquad 0a + ba = (0 + b)a = ba$$

является следствием утверждений 3.2.1.2, 3.2.1.3. Теорема является следствием равенства (3.2.4), теоремы 3.1.4 и утверждения 3.1.2.2. □

Теорема 3.2.4. *Пусть $A$ - кольцо. Тогда*
$$(3.2.5) \qquad (-a)b = -(ab)$$
$$(3.2.6) \qquad a(-b) = -(ab)$$
$$(3.2.7) \qquad (-a)(-b) = ab$$
*для любого $A$-чисел $a$, $b$. Если кольцо имеет единичный элемент, то*
$$(3.2.8) \qquad (-1)a = -a$$

---

[3.4] Смотри также определения кольца на странице [1]-73.

[3.5] Смотри аналогичное определение на странице [7]-52.

*для любого $A$-числа $a$.*

Доказательство. Равенства

$$(3.2.9) \qquad ab + (-a)b = (a + (-a))b = 0b = 0$$

$$(3.2.10) \qquad ab + a(-b) = a(b + (-b)) = a0 = 0$$

являются следствием утверждений 3.2.1.1, 3.2.1.3 и равенств (3.1.7), (3.2.3). Равенство (3.2.5) является следствием теоремы 3.1.10 и равенства (3.2.9). Равенство (3.2.6) является следствием теоремы 3.1.10 и равенства (3.2.10).

Равенство

$$(3.2.11) \qquad (-a)(-b) = -(a(-b)) = -(-(ab))$$

является следствием равенств (3.2.5), (3.2.6). Равенство (3.2.7) является следствием равенств (3.1.11), (3.2.11).

Если кольцо имеет единичный элемент, то равенство

$$(3.2.12) \qquad b + (-1)b = 1b + (-1)b = (1 + (-1))b = 0b = 0$$

является следствием утверждений 3.2.1.1, 3.2.1.2, 3.2.1.3 и равенств (3.1.7), (3.2.3). $\qquad\square$

ОПРЕДЕЛЕНИЕ 3.2.5. **Левый идеал** [3.6] $D_1$ *кольца $D$ - это подгруппа аддитивной группы $D$ такая, что*

$$DD_1 \subseteq D_1$$
$\qquad\square$

ОПРЕДЕЛЕНИЕ 3.2.7. **Правый идеал** [3.6] $D_1$ *кольца $D$ - это подгруппа аддитивной группы $D$ такая, что*

$$D_1 D \subseteq D_1$$
$\qquad\square$

ТЕОРЕМА 3.2.6. *Если $D_1$ является левым идеалом унитального кольца $D$, то*

$$DD_1 = D_1$$

ТЕОРЕМА 3.2.8. *Если $D_1$ является правым идеалом унитального кольца $D$, то*

$$D_1 D = D_1$$

Доказательство. Теорема является следствием определения 3.2.5 и равенства

$$eD_1 = D_1$$
$\qquad\square$

Доказательство. Теорема является следствием определения 3.2.7 и равенства

$$D_1 e = D_1$$
$\qquad\square$

ОПРЕДЕЛЕНИЕ 3.2.9. **Идеал** [3.6] *- это подгруппа аддитивной группы $D$, которая одновременно является левым и правым идеалами.* $\qquad\square$

Если $D$ - кольцо, то множество $\{0\}$ и само $D$ являются идеалами, которые мы называем тривиальными идеалами.

ОПРЕДЕЛЕНИЕ 3.2.10. *Пусть $A$ - коммутативное кольцо. Если множество $A' = A \setminus \{0\}$ является мультипликативной группой относительно*

---

[3.6] Смотри также определение на странице [1]-75.

*произведения, то множество A называется полем.* □

## 3.3. Гомоморфизм групп

**Теорема 3.3.1.** *Отображение*

$$f : A \to B$$

*является гомоморфизмом мультипликативного моноида A с единичным элементом $e_1$ в мультипликативный моноид B с единичным элементом $e_2$, если*

(3.3.1)          $f(ab) = f(a)f(b)$

(3.3.2)          $f(e_1) = e_2$

*для любых A-чисел a, b.*

Доказательство. Теорема является следствием определений 2.1.6, 2.1.7, 3.1.1 и замечания 3.1.5. □

**Теорема 3.3.2.** *Отображение*

$$f : A \to B$$

*является гомоморфизмом аддитивного моноида A с нулевым элементом $0_1$ в аддитивный моноид B с нулевым элементом $0_2$, если*

(3.3.3)          $f(a + b) = f(a) + f(b)$

(3.3.4)          $f(0_1) = 0_2$

*для любых A-чисел a, b.*

Доказательство. Теорема является следствием определений 2.1.6, 2.1.7, 3.1.2 и замечания 3.1.6. □

**Теорема 3.3.3.** *Пусть отображения*

$$f : A \to B$$
$$g : A \to B$$

*являются гомоморфизмами мультипликативной абелевой группы A с единичным элементом $e_1$ в мультипликативную абелеву группу B с единичным элементом $e_2$. Тогда отображение*

$$h : A \to B$$

*определённое равенством*

$$h(x) = (fg)(x) = f(x)g(x)$$

*является гомоморфизмом мультипликативной абелевой группы A в мультипликативный абелеву группу B.*

Доказательство. Теорема явля-

**Теорема 3.3.4.** *Пусть отображения*

$$f : A \to B$$
$$g : A \to B$$

*являются гомоморфизмами аддитивной абелевой группы A с нулевым элементом $0_1$ в аддитивную абелеву группу B с нулевым элементом $0_2$. Тогда отображение*

$$h : A \to B$$

*определённое равенством*

$$h(x) = (f + g)(x) = f(x) + g(x)$$

*является гомоморфизмом аддитивной абелевой группы A в аддитивный абелеву группу B.*

Доказательство. Теорема является следствием равенств

$$
\begin{aligned}
h(x + y) &= f(x + y) + g(x + y) \\
&= f(x) + f(y) + g(x) + g(y) \\
&= f(x) + g(x) + f(y) + g(y) \\
&= h(x) + h(y)
\end{aligned}
$$

ется следствием равенств

$$h(xy) = f(xy)g(xy)$$
$$= f(x)f(y)g(x)g(y)$$
$$= f(x)g(x)f(y)g(y)$$
$$= h(x)h(y)$$

$$h(e_1) = (fg)(e_1)$$
$$= f(e_1)g(e_1)$$
$$= e_2 + e_2 = e_2$$

и теоремы 3.3.1. □

$$h(0_1) = (f + g)(0_1)$$
$$= f(0_1) + g(0_1)$$
$$= 0_2 + 0_2 = 0_2$$

и теоремы 3.3.2. □

ТЕОРЕМА 3.3.5. *Пусть*

$$f : A \to B$$

*гомоморфизм* 3.7 *мультипликативной группы* $A$ *с единичным элементом* $e_1$ *в мультипликативную группу* $B$ *с единичным элементом* $e_2$. *Тогда*

(3.3.5) $$f(a^{-1}) = f(a)^{-1}$$

ДОКАЗАТЕЛЬСТВО. Равенство

(3.3.6) $$e_2 = f(e_1) = f(aa^{-1})$$
$$= f(a)f(a^{-1})$$

является следствием равенств (3.1.5), (3.3.1), (3.3.2). Равенство (3.3.5) является следствием равенства (3.3.6) и теоремы 3.1.9. □

ТЕОРЕМА 3.3.6. *Пусть*

$$f : A \to B$$

*гомоморфизм* 3.7 *аддитивной группы* $A$ *с нулевым элементом* $0_1$ *в аддитивную группу* $B$ *с нулевым элементом* $0_2$. *Тогда*

(3.3.7) $$f(-a) = -f(a)$$

ДОКАЗАТЕЛЬСТВО. Равенство

(3.3.8) $$0_2 = f(0_1) = f(a + (-a))$$
$$= f(a) + f(-a)$$

является следствием равенств (3.1.7), (3.3.3), (3.3.4). Равенство (3.3.7) является следствием равенства (3.3.8) и теоремы 3.1.10. □

ТЕОРЕМА 3.3.7. *Образ гомоморфизма групп*

$$f : A \to B$$

*является подгруппой группы* $B$.

ДОКАЗАТЕЛЬСТВО.

Пусть $A$ и $B$ - мультипликативные группы. Согласно определениям 2.1.6, 3.1.1 и теореме 3.3.1, множество Im $f$ является моноидом. Согласно теореме 3.3.5, множество Im $f$ является группой.

Пусть $A$ и $B$ - аддитивные группы. Согласно определениям 2.1.6, 3.1.2 и теореме 3.3.2, множество Im $f$ является моноидом. Согласно теореме 3.3.6, множество Im $f$ является группой. □

Существует гомоморфизм аддитивной группы в мультипликативную. Например, отображение $y = e^x$ является гомоморфизмом аддитивной группы действительных чисел в мультипликативную группу положительных действительных чисел.

---

3.7 Мы будем также пользоваться термином гомоморфизм групп для отображения $f$.

Теорема 3.3.8. *Пусть*

$$f : A \to B$$

*гомоморфизм мультипликативной группы $A$ с единичным элементом $e_1$ в мультипликативную группу $B$ с единичным элементом $e_2$. Множество*

(3.3.9)      $\ker f = \{a \in A : f(a) = e_2\}$

*называется* **ядром гомоморфизма** $f$ **групп**. *Ядро гомоморфизма групп является нормальной подгруппой группы $A$.*

Теорема 3.3.9. *Пусть*

$$f : A \to B$$

*гомоморфизм аддитивной группы $A$ с нулевым элементом $0_1$ в аддитивную группу $B$ с нулевым элементом $0_2$. Множество*

(3.3.12)      $\ker f = \{a \in A : f(a) = 0_2\}$

*называется* **ядром гомоморфизма** $f$ **групп**. *Ядро гомоморфизма групп является нормальной подгруппой группы $A$.*

Доказательство. Пусть $A$ - мультипликативная группа.

3.3.8.1: Пусть $a, b \in \ker f$. Равенство

$$f(ab) = f(a)f(b) = e_2 e_2 = e_2$$

является следствием равенств (3.3.1), (3.3.9) и утверждения 3.1.1.2. Согласно определению (3.3.9), $ab \in \ker f$.

3.3.8.2: Согласно равенствам (3.3.2), (3.3.9), $e_1 \in \ker f$.

3.3.8.3: Из утверждений 3.3.8.1, 3.3.8.2 следует, что множество $\ker f$ является мультипликативным моноидом.

3.3.8.4: Пусть $a \in \ker f$. Согласно теореме 3.1.9, существует $A$-число $a^{-1}$. Равенство

$$f(a^{-1}) = f(a)^{-1} = e_2^{-1} = e_2$$

является следствием равенств (3.1.13), (3.3.5), (3.3.9). Согласно определению (3.3.9), $a^{-1} \in \ker f$.

Согласно утверждениям 3.3.8.3, 3.3.8.4, множество $\ker f$ является мультипликативной группой.

Равенство

(3.3.10)
$$\begin{aligned} & f(a(\ker f)a^{-1}) \\ &= f(a)f(\ker f)f(a^{-1}) \\ &= f(a)e_2(f(a))^{-1} \\ &= f(a)(f(a))^{-1} = e_2 \end{aligned}$$

Доказательство. Пусть $A$ - аддитивная группа.

3.3.9.1: Пусть $a, b \in \ker f$. Равенство

$$f(a + b) = f(a) + f(b) = 0_2 + 0_2 = 0_2$$

является следствием равенств (3.3.3), (3.3.12) и утверждения 3.1.2.2. Согласно определению (3.3.12), $a + b \in \ker f$.

3.3.9.2: Согласно равенствам (3.3.4), (3.3.12), $0_1 \in \ker f$.

3.3.9.3: Из утверждений 3.3.9.1, 3.3.9.2 следует, что множество $\ker f$ является аддитивным моноидом.

3.3.9.4: Пусть $a \in \ker f$. Согласно теореме 3.1.10, существует $A$-число $-a$. Равенство

$$f(-a) = -f(a) = -0_2 = 0_2$$

является следствием равенств (3.1.15), (3.3.7), (3.3.12). Согласно определению (3.3.12), $-a \in \ker f$.

Согласно утверждениям 3.3.9.3, 3.3.9.4, множество $\ker f$ является аддитивной группой.

Равенство

(3.3.13)
$$\begin{aligned} & f(a + (\ker f) - a) \\ &= f(a) + f(\ker f) - f(a) \\ &= f(a) + 0_2 - f(a) \\ &= f(a) - f(a) = 0_2 \end{aligned}$$

является следствием равенств (3.3.1), (3.3.5), (3.3.9). Равенство

$$a(\ker f)a^{-1} = \ker f$$

является следствием равенства (3.3.10). Следовательно,

(3.3.11)        $a(\ker f) = (\ker f)a$

Теорема является следствием равенства (3.3.11) и определения 3.1.19.     □

является следствием равенств (3.3.3), (3.3.7), (3.3.12). Равенство

$$a + (\ker f) + (-a) = \ker f$$

является следствием равенства (3.3.13). Следовательно,

(3.3.14)        $a + (\ker f) = (\ker f) + a$

Теорема является следствием равенства (3.3.14) и определения 3.1.20.     □

---

ТЕОРЕМА 3.3.10. *Пусть*

$$f : A \to B$$

*гомоморфизм мультипликативной группы $A$ с единичным элементом $e_1$ в мультипликативную группу $B$ с единичным элементом $e_2$. Если*

(3.3.15)        $\ker f = \{e_1\}$

*то ядро гомоморфизма $f$ называется* **тривиальным**, *а гомоморфизм $f$ является мономорфизмом.*

ТЕОРЕМА 3.3.11. *Пусть*

$$f : A \to B$$

*гомоморфизм аддитивной группы $A$ с нулевым элементом $0_1$ в аддитивную группу $B$ с нулевым элементом $0_2$. Если*

(3.3.16)        $\ker f = \{0_1\}$

*то ядро гомоморфизма $f$ называется* **тривиальным**, *а гомоморфизм $f$ является мономорфизмом.*

---

Доказательство. Пусть существуют $A$-числа $a$, $b$ такие, что

(3.3.17)        $f(a) = f(b)$

Равенство

$$
\begin{aligned}
f(ab^{-1}) &= f(a)f(b^{-1}) \\
&= f(a)f(b)^{-1} \\
&= f(b)f(b)^{-1} \\
&= e_2
\end{aligned}
$$
(3.3.18)

является следствием равенств (3.3.1), (3.3.5), (3.3.17) и утверждения 3.1.1.2. Равенство

(3.3.19)        $ab^{-1} = e_1$

является следствием равенств (3.3.9), (3.3.15), (3.3.18). Из равенств (3.1.9), (3.3.19) и теоремы 3.1.9 следует, что $a = b$.     □

Доказательство. Пусть существуют $A$-числа $a$, $b$ такие, что

(3.3.20)        $f(a) = f(b)$

Равенство

$$
\begin{aligned}
f(a + (-b)) &= f(a) + f(-b) \\
&= f(a) + (-f(b)) \\
&= f(b) + (-f(b)) \\
&= 0_2
\end{aligned}
$$
(3.3.21)

является следствием равенств (3.3.3), (3.3.7), (3.3.20) и утверждения 3.1.2.2. Равенство

(3.3.22)        $a + (-b) = 0_1$

является следствием равенств (3.3.12), (3.3.16), (3.3.21). Из равенств (3.1.11), (3.3.22) и теоремы 3.1.10 следует, что $a = b$.     □

---

ТЕОРЕМА 3.3.12. *Пусть $B$ - нормальная подгруппа мультипликативной группы $A$ с единичным элементом $e$. Существует мультипликативная группа $C$ и гомомор-*

ТЕОРЕМА 3.3.13. *Пусть $B$ - нормальная подгруппа аддитивной группы $A$ с нулевым элементом $0$. Существует аддитивная группа $C$ и гомоморфизм*

$$f : A \to C$$

*физм*

$$f : A \to C$$

*такие, что*

(3.3.23) $\qquad \ker f = B$

ДОКАЗАТЕЛЬСТВО. Положим

$$f(x) = xB$$

Поскольку множество $eB = B$ является единичным элементом мультипликативной группы $A/B$, то

(3.3.24) $\qquad B \subseteq \ker f$

Если $f(x) = B$, то $xB = B$. Согласно теореме 3.1.17, $x \in B$. Следовательно,

(3.3.25) $\qquad \ker f \subseteq B$

Теорема является следствием утверждений (3.3.24), (3.3.25). □

*такие, что*

(3.3.26) $\qquad \ker f = B$

ДОКАЗАТЕЛЬСТВО. Положим

$$f(x) = x + B$$

Поскольку множество $0 + B = B$ является нулевым элементом аддитивной группы $A/B$, то

(3.3.27) $\qquad B \subseteq \ker f$

Если $f(x) = B$, то $x + B = B$. Согласно теореме 3.1.18, $x \in B$. Следовательно,

(3.3.28) $\qquad \ker f \subseteq B$

Теорема является следствием утверждений (3.3.27), (3.3.28). □

ОПРЕДЕЛЕНИЕ 3.3.14. *Пусть $A_1$, ..., $A_n$ - группы. Последовательность гомоморфизмов*

$$A_1 \xrightarrow{f_1} A_2 \xrightarrow{f_2} ... \xrightarrow{f_{n-1}} A_n$$

*называется* **точной**, *если* $\forall i, \ i = 1, ..., n-2, \ \operatorname{Im} f_i = \ker f_{i+1}$. □

ПРИМЕР 3.3.15. *Пусть*

(3.3.29) $\qquad e \longrightarrow A \xrightarrow{f} B \xrightarrow{g} C \xrightarrow{h} e$

*точная последовательность гомоморфизмов мультипликативных групп.*

*Согласно равенству* (3.3.2) *и определению 3.3.14*, $\ker f = \{e\}$. *Согласно теореме 3.3.11, отображение $f$ является мономорфизмом.*

*Из диаграммы* (3.3.29) *следует, что* $\ker h = C$. *Согласно определению 3.3.14*, $\operatorname{Im} g = C$. *Согласно определению 2.1.8, гомоморфизм $g$ является эпиморфизмом.* □

ПРИМЕР 3.3.16. *Согласно примеру 3.3.15, если последовательность гомоморфизмов аддитивных групп*

$$e \longrightarrow A \xrightarrow{f} B \longrightarrow e$$

*точна, то гомоморфизм $f$ является изоморфизмом.* □

ПРИМЕР 3.3.17. *Объединяя примеры 3.3.15, 3.3.16, мы получим диаграмму*

$$
\begin{array}{ccccccccc}
 & & & & & & e & & \\
 & & & & & & \downarrow & & \\
e & \longrightarrow & A & \longrightarrow & B & \longrightarrow & C & \longrightarrow & e \\
 & & & & & & \downarrow f & & \\
e & \longrightarrow & A & \longrightarrow & B & \longrightarrow & B/A & \longrightarrow & e \\
 & & & & & & \downarrow & & \\
 & & & & & & e & & 
\end{array}
$$

*где строки и столбец являются точными последовательностями гомоморфизмов и гомоморфизм f является изоморфизмом.* □

**Пример 3.3.18.** *Пусть*

(3.3.30)          $0 \longrightarrow A \xrightarrow{f} B \xrightarrow{g} C \xrightarrow{h} 0$

*точная последовательность гомоморфизмов аддитивных групп.*

*Согласно равенству* (3.3.4) *и определению* 3.3.14, $\ker f = \{0\}$. *Согласно теореме* 3.3.11, *отображение f является мономорфизмом.*

*Из диаграммы* (3.3.30) *следует, что* $\ker h = C$. *Согласно определению* 3.3.14, $\operatorname{Im} g = C$. *Согласно определению* 2.1.8, *гомоморфизм g является эпиморфизмом.*

**Пример 3.3.19.** *Согласно примеру* 3.3.18, *если последовательность гомоморфизмов аддитивных групп*

$$0 \longrightarrow A \xrightarrow{f} B \longrightarrow 0$$

*точна, то гомоморфизм f является изоморфизмом.* □

**Пример 3.3.20.** *Объединяя примеры* 3.3.18, 3.3.19, *мы получим диаграмму*

*где строки и столбец являются точными последовательностями гомоморфизмов и гомоморфизм f является изоморфизмом.* □

## 3.4. Базис аддитивной абелевой группы

В этом разделе $G$ - аддитивная абелева группа.

**Определение 3.4.1.** *Мы определим действие кольца целых чисел Z в аддитивной абелевой группе G согласно правилу*

(3.4.1)                    $0g = 0$

(3.4.2)                    $(n+1)g = ng + g$

(3.4.3)                    $(n-1)g = ng - g$

□

Теорема 3.4.2. *Действие кольца целых чисел $Z$ в аддитивной абелевой группе $G$, рассмотренное в определении* 3.4.1, *является представлением. Верны следующие равенства*

$$(3.4.4) \qquad 1a = a$$

$$(3.4.5) \qquad (nm)a = n(ma)$$

$$(3.4.6) \qquad (m+n)a = ma + na$$

$$(3.4.7) \qquad (m-n)a = ma - na$$

$$(3.4.8) \qquad n(a+b) = na + nb$$

Доказательство. Равенство (3.4.4) является следствием равенства (3.4.1) и равенства (3.4.2), когда $n = 0$.

Из равенства (3.4.1) следует, что равенство (3.4.6) верно, когда $n = 0$.

- Пусть равенство (3.4.6) верно, когда $n = k \geq 0$. Тогда

$$(m+k)a = ma + ka$$

Равенство

$$(m+(k+1))a = ((m+k)+1)a = (m+k)a + a = ma + ka + a$$
$$= ma + (k+1)a$$

является следствием равенства (3.4.2). Следовательно, равенство (3.4.6) верно, когда $n = k+1$. Согласно принципу математической индукции, равенство (3.4.6) верно для любого $n \geq 0$.

- Пусть равенство (3.4.6) верно, когда $n = k \leq 0$. Тогда

$$(m+k)a = ma + ka$$

Равенство

$$(m+(k-1))a = ((m+k)-1)a = (m+k)a - a = ma + ka - a$$
$$= ma + (k-1)a$$

является следствием равенства (3.4.3). Следовательно, равенство (3.4.6) верно, когда $n = k-1$. Согласно принципу математической индукции, равенство (3.4.6) верно для любого $n \leq 0$.

- Следовательно, равенство (3.4.6) верно для любого $n \in Z$.

Равенство

$$(3.4.9) \qquad (k+n)a - na = ka$$

является следствием равенства (3.4.6). Равенство (3.4.7) является следствием равенства (3.4.9), если мы положим $m = k+n$, $k = m-n$.

Из равенства (3.4.1) следует, что равенство (3.4.5) верно, когда $n = 0$.

- Пусть равенство (3.4.5) верно, когда $n = k \geq 0$. Тогда

$$(km)a = k(ma)$$

Равенство

$$((k+1)m)a = (km+m)a = (km)a + ma = k(ma) + ma$$
$$= (k+1)(ma)$$

является следствием равенств (3.4.2), (3.4.6). Следовательно, равенство (3.4.5) верно, когда $n = k+1$. Согласно принципу математической индукции, равенство (3.4.5) верно для любого $n \geq 0$.

- Пусть равенство (3.4.6) верно, когда $n = k \leq 0$. Тогда

$$(km)a = k(ma)$$

Равенство

$$((k-1)m)a = (km - m)a = (km)a - ma = k(ma) - ma$$
$$= (k-1)(ma)$$

является следствием равенств (3.4.3), (3.4.7). Следовательно, равенство (3.4.5) верно, когда $n = k-1$. Согласно принципу математической индукции, равенство (3.4.5) верно для любого $n \leq 0$.

- Следовательно, равенство (3.4.5) верно для любого $n \in Z$.

Из равенства (3.4.1) следует, что равенство (3.4.8) верно, когда $n = 0$.

- Пусть равенство (3.4.8) верно, когда $n = k \geq 0$. Тогда

$$k(a + b) = ka + kb$$

Равенство

$$(k+1)(a + b) = k(a + b) + a + b = ka + kb + a + b$$
$$= ka + a + kb + b$$
$$= (k+1)a + (k+1)b$$

является следствием равенства (3.4.2). Следовательно, равенство (3.4.8) верно, когда $n = k+1$. Согласно принципу математической индукции, равенство (3.4.8) верно для любого $n \geq 0$.

- Пусть равенство (3.4.6) верно, когда $n = k \leq 0$. Тогда

$$k(a + b) = ka + kb$$

Равенство

$$(k-1)(a + b) = k(a + b) - (a + b) = ka + kb - a - b$$
$$= ka - a + kb - b$$
$$= (k-1)a + (k-1)b$$

является следствием равенства (3.4.3). Следовательно, равенство (3.4.8) верно, когда $n = k-1$. Согласно принципу математической индукции, равенство (3.4.8) верно для любого $n \leq 0$.

- Следовательно, равенство (3.4.8) верно для любого $n \in Z$.

Из равенства (3.4.8) следует, что отображение

$$\varphi(n) : a \in G \to na \in G$$

является эндоморфизмом абелевой группы $G$. Из равенств (3.4.6), (3.4.5) следует, что отображение

$$\varphi : Z \to \text{End}(Ab, G)$$

является гомоморфизмом кольца $Z$. Согласно определению 2.3.1, отображение $\varphi$ является представлением кольца целых чисел $Z$ в абелевой группе $G$. $\quad\square$

Теорема 3.4.3. *Множество $G$-чисел, порождённое множеством* $S = \{s_i : i \in I\}$ *, имеет вид*

$$(3.4.10) \qquad J(S) = \left\{ g : g = \sum_{i \in I} g^i s_i, g^i \in Z \right\}$$

*где множество* $\{i \in I : g^i \neq 0\}$ *конечно.*

Доказательство. Мы докажем теорему по индукции, опираясь на теоремы [18]-5.1, страница 94, и 2.7.5.

Для произвольного $s_k \in S$, положим $g^i = \delta^i_k$. Тогда

$$(3.4.11) \qquad s_k = \sum_{i \in I} g^i s_i$$

$s_k \in J(S)$ следует из (3.4.10), (3.4.11).

Пусть $g_1$, $g_2 \in X_k \subseteq J(S)$. Так как $G$ является абелевой группой, то, согласно утверждению 2.7.5.3, $g_1 + g_2 \in J(S)$. Согласно равенству (3.4.10), существуют $Z$-числа $g^i_1, g^i_2, i \in I$, такие, что

$$(3.4.12) \qquad g_1 = \sum_{i \in I} g^i_1 s_i \quad g_2 = \sum_{i \in I} g^i_2 s_i$$

где множества

$$(3.4.13) \qquad H_1 = \{i \in I : g^i_1 \neq 0\} \quad H_2 = \{i \in I : g^i_2 \neq 0\}$$

конечны. Из равенства (3.4.12) следует, что

$$(3.4.14) \qquad g_1 + g_2 = \sum_{i \in I} g^i_1 s_i + \sum_{i \in I} g^i_2 s_i = \sum_{i \in I} (g^i_1 s_i + g^i_2 s_i)$$

Равенство

$$(3.4.15) \qquad g_1 + g_2 = \sum_{i \in I} (g^i_1 + g^i_2) s_i$$

является следствием равенств (3.4.6), (3.4.14). Из равенства (3.4.13) следует, что множество

$$\{i \in I : g^i_1 + g^i_2 \neq 0\} \subseteq H_1 \cup H_2$$

конечно. Из равенства (3.4.15) следует, что $g_1 + g_2 \in J(S)$. $\qquad \square$

Соглашение 3.4.4. *Мы будем пользоваться соглашением Эйнштейна о сумме, в котором повторяющийся индекс (один вверху и один внизу) подразумевает сумму по повторяющемуся индексу. В этом случае предполагается известным множество индекса суммирования и знак суммы опускается*

$$g^i s_i = \sum_{i \in I} g^i s_i$$

$\square$

Теорема 3.4.5. *Гомоморфизм*

$$f : G \to G'$$

*удовлетворяет равенству*

$$(3.4.16) \qquad f(na) = nf(a)$$

Доказательство. Из равенства (3.4.1) следует, что равенство (3.4.16) верно, когда $n = 0$.

- Пусть равенство (3.4.16) верно, когда $n = k \geq 0$. Тогда

$$f(ka) = kf(a)$$

Равенство

$$f((k+1)a) = f(ka + a) = f(ka) + f(a) = kf(a) + f(a)$$
$$= (k+1)f(a)$$

является следствием равенства (3.4.2). Следовательно, равенство (3.4.16) верно, когда $n = k + 1$. Согласно принципу математической индукции, равенство (3.4.16) верно для любого $n \geq 0$.

- Пусть равенство (3.4.6) верно, когда $n = k \leq 0$. Тогда

$$f(ka) = kf(a)$$

Равенство

$$f((k-1)a) = f(ka - a) = f(ka) - f(a) = kf(a) - f(a)$$
$$= (k-1)f(a)$$

является следствием равенства (3.4.3). Следовательно, равенство (3.4.16) верно, когда $n = k - 1$. Согласно принципу математической индукции, равенство (3.4.16) верно для любого $n \leq 0$.

- Следовательно, равенство (3.4.16) верно для любого $n \in Z$.

□

Определение 3.4.6. *Пусть* $S = \{s_i : i \in I\}$ *- множество $G$-чисел. Выражение* $g^i s_i$ *называется* **линейной комбинацией** *$G$-чисел $s_i$. $G$-число $g = g^i s_i$ называется* **линейно зависимым** *от $G$-чисел $s_i$.* □

Следующее определение является следствием теорем 2.7.5, 3.4.3 и определения 2.7.4.

Определение 3.4.7. *$J(S)$ называется* **подгруппой, порождённой множеством** *$S$, а $S$ -* **множеством образующих** *подгруппы $J(S)$. В частности,* **множеством образующих** *абелевой группы будет такое подмножество $X \subset G$, что* $J(X) = G$. □

Следующее определение является следствием теорем 2.7.5, 3.4.3 и определения 2.7.14.

Определение 3.4.8. *Если множество $X \subset G$ является множеством образующих абелевой группы $G$, то любое множество $Y$, $X \subset Y \subset G$ также является множеством образующих абелевой группы $G$. Если существует минимальное множество $X$, порождающее абелеву группу $G$, то такое множество $X$ называется* **квазибазисом** *абелевой группы $G$.* □

ОПРЕДЕЛЕНИЕ 3.4.9. *Множество $G$-чисел*[3.8]$g_i$, $i \in I$, **линейно независимо**, *если $w = 0$ следует из уравнения*

$$w^i g_i = 0$$

*В противном случае, множество $G$-чисел $g_i$, $i \in I$,* **линейно зависимо**. □

ОПРЕДЕЛЕНИЕ 3.4.10. *Пусть множество $G$-чисел $S = \{s_i : i \in I\}$ является квазибазисом. Если множество $G$-чисел $S$ линейно независимо, то квазибазис $S$ называется* **базисом** *абелевой группы $G$*. □

ТЕОРЕМА 3.4.11. *Произведение в кольце $D$ билинейно над кольцом $Z$*

(3.4.17)                          $$(nd_1)d_2 = n(d_1 d_2)$$

(3.4.18)                          $$d_1(nd_2) = n(d_1 d_2)$$

ДОКАЗАТЕЛЬСТВО. Из равенства (3.4.1) следует, что равенство (3.4.17) верно, когда $n = 0$.

- Пусть равенство (3.4.17) верно, когда $n = k \geq 0$. Тогда

$$(kd_1)d_2 = k(d_1 d_2)$$

Равенство

$$((k+1)d_1)d_2 = (kd_1 + d_1)d_2 = (kd_1)d_2 + d_1 d_2 = k(d_1 d_2) + d_1 d_2$$
$$= (k+1)(d_1 d_2)$$

является следствием равенств (3.4.2), (3.2.2). Следовательно, равенство (3.4.17) верно, когда $n = k + 1$. Согласно принципу математической индукции, равенство (3.4.17) верно для любого $n \geq 0$.

- Пусть равенство (3.4.6) верно, когда $n = k \leq 0$. Тогда

$$(kd_1)d_2 = k(d_1 d_2)$$

Равенство

$$((k-1)d_1)d_2 = (kd_1 - d_1)d_2 = (kd_1)d_2 - d_1 d_2 = k(d_1 d_2) - d_1 d_2$$
$$= (k-1)(d_1 d_2)$$

является следствием равенств (3.4.2), (3.2.2). Следовательно, равенство (3.4.17) верно, когда $n = k - 1$. Согласно принципу математической индукции, равенство (3.4.17) верно для любого $n \leq 0$.

- Следовательно, равенство (3.4.17) верно для любого $n \in Z$.

Из равенства (3.4.1) следует, что равенство (3.4.18) верно, когда $n = 0$.

- Пусть равенство (3.4.18) верно, когда $n = k \geq 0$. Тогда

$$d_1(kd_2) = k(d_1 d_2)$$

---

[3.8] Я следую определению в [1], страница 100.

Равенство

$$d_1((k+1)d_2) = d_1(kd_2 + d_2) = d_1(kd_2) + d_1d_2 = k(d_1d_2) + d_1d_2$$
$$= (k+1)(d_1d_2)$$

является следствием равенств ($3.4.2$), ($3.2.1$). Следовательно, равенство ($3.4.18$) верно, когда $n = k + 1$. Согласно принципу математической индукции, равенство ($3.4.18$) верно для любого $n \geq 0$.

- Пусть равенство ($3.4.6$) верно, когда $n = k \leq 0$. Тогда

$$d_1(kd_2) = k(d_1d_2)$$

Равенство

$$d_1((k-1)d_2) = d_1(kd_2 - d_2) = d_1(kd_2) - d_1d_2 = k(d_1d_2) - d_1d_2$$
$$= (k-1)(d_1d_2)$$

является следствием равенств ($3.4.2$), ($3.2.1$). Следовательно, равенство ($3.4.18$) верно, когда $n = k - 1$. Согласно принципу математической индукции, равенство ($3.4.18$) верно для любого $n \leq 0$.

- Следовательно, равенство ($3.4.18$) верно для любого $n \in Z$. $\square$

ТЕОРЕМА 3.4.12. *Пусть $D$ - кольцо с единицей $1_D$.*

3.4.12.1: *Отображение*

$$f : n \in Z \to n1_D \in D$$

*является гомоморфизмом кольца $Z$ в кольцо $D$ и позволяет отождествить $Z$-число $n$ и $D$-число $n1_D$.*

3.4.12.2: *Кольцо $Z/\ker f$ является подкольцом кольца $D$.*

ДОКАЗАТЕЛЬСТВО. Пусть $n_1, n_2 \in Z$. Равенство

$$(3.4.19) \qquad n_1 1_D + n_2 1_D = (n_1 + n_2) 1_D$$

является следствием равенства ($3.4.6$). Равенство

$$(3.4.20) \qquad (n_1 1_D)(n_2 1_D) = n_1(1_D(n_2 1_D)) = n_1(n_2(1_D 1_D))$$
$$= n_1(n_2 1_D) = (n_1 n_2) 1_D$$

является следствием равенств ($3.4.17$), ($3.4.18$). Утверждение $3.4.12.1$ является следствием равенств ($3.4.19$), ($3.4.20$). Утверждение $3.4.12.2$ является следствием утверждения $3.4.12.1$. $\square$

## 3.5. Базис мультипликативной абелевой группы

В этом разделе $G$ - мультипликативная абелева группа.

ОПРЕДЕЛЕНИЕ 3.5.1. *Мы определим действие кольца целых чисел $Z$ в мультипликативной абелевой группе $G$ согласно правилу*

$$(3.5.1) \qquad\qquad g^0 = e$$
$$(3.5.2) \qquad\qquad g^{n+1} = g^n g$$
$$(3.5.3) \qquad\qquad g^{n-1} = g^n g^{-1}$$

$\square$

Теорема 3.5.2. *Действие кольца целых чисел $Z$ в мультипликативной абелевой группе $G$, рассмотренное в определении* 3.5.1, *является представлением. Верны следующие равенства*

(3.5.4) $$a^1 = a$$
(3.5.5) $$a^{nm} = (a^n)^m$$
(3.5.6) $$a^{m+n} = a^m a^n$$
(3.5.7) $$a^{m-n} = a^m a^{-n}$$
(3.5.8) $$(ab)^n = a^n b^n$$

Доказательство. Равенство (3.5.4) является следствием равенства (3.5.1) и равенства (3.5.2), когда $n = 0$.

Из равенства (3.5.1) следует, что равенство (3.5.6) верно, когда $n = 0$.

- Пусть равенство (3.5.6) верно, когда $n = k \geq 0$. Тогда

$$a^{m+k} = a^m a^k$$

Равенство

$$a^{m+(k+1)} = a^{(m+k)+1} = a^{m+k} a = a^m a^k a$$
$$= a^m a^{k+1}$$

является следствием равенства (3.5.2). Следовательно, равенство (3.5.6) верно, когда $n = k+1$. Согласно принципу математической индукции, равенство (3.5.6) верно для любого $n \geq 0$.

- Пусть равенство (3.5.6) верно, когда $n = k \leq 0$. Тогда

$$a^{m+k} = a^m a^k$$

Равенство

$$a^{m+(k-1)} = a^{(m+k)-1} = a^{m+k} a^{-1} = a^m a^k a^{-1}$$
$$= a^m a^{k-1}$$

является следствием равенства (3.5.3). Следовательно, равенство (3.5.6) верно, когда $n = k-1$. Согласно принципу математической индукции, равенство (3.5.6) верно для любого $n \leq 0$.

- Следовательно, равенство (3.5.6) верно для любого $n \in Z$.

Равенство

(3.5.9) $$a^{k+n} a^{-n} = a^k$$

является следствием равенства (3.5.6). Равенство (3.5.7) является следствием равенства (3.5.9), если мы положим $m = k+n$, $k = m-n$.

Из равенства (3.5.1) следует, что равенство (3.5.5) верно, когда $n = 0$.

- Пусть равенство (3.5.5) верно, когда $n = k \geq 0$. Тогда

$$a^{km} = (a^m)^k$$

Равенство

$$a^{(k+1)m} = a^{km+m} = a^{km} a^m = (a^m)^k a^m$$
$$= (a^m)^{k+1}$$

является следствием равенств (3.5.2), (3.5.6). Следовательно, равенство (3.5.5) верно, когда $n = k+1$. Согласно принципу математической индукции, равенство (3.5.5) верно для любого $n \geq 0$.

• Пусть равенство (3.5.6) верно, когда $n = k \leq 0$. Тогда

$$a^{km} = (a^m)^k$$

Равенство

$$a^{(k-1)m} = a^{km-m} = a^{km}a^{-m} = (a^m)^k a^{-m}$$
$$= (a^m)^{k-1}$$

является следствием равенств (3.5.3), (3.5.7). Следовательно, равенство (3.5.5) верно, когда $n = k-1$. Согласно принципу математической индукции, равенство (3.5.5) верно для любого $n \leq 0$.

• Следовательно, равенство (3.5.5) верно для любого $n \in Z$.

Из равенства (3.5.1) следует, что равенство (3.5.8) верно, когда $n = 0$.

• Пусть равенство (3.5.8) верно, когда $n = k \geq 0$. Тогда

$$(ab)^k = a^k b^k$$

Равенство

$$(ab)^{k+1} = (ab)^k ab = a^k b^k ab$$
$$= a^k ab^k b$$
$$= a^{k+1} b^{k+1}$$

является следствием равенства (3.5.2). Следовательно, равенство (3.5.8) верно, когда $n = k+1$. Согласно принципу математической индукции, равенство (3.5.8) верно для любого $n \geq 0$.

• Пусть равенство (3.5.6) верно, когда $n = k \leq 0$. Тогда

$$(ab)^k = a^k b^k$$

Равенство

$$(ab)^{k-1} = (ab)^k (ab)^{-1} = a^k b^k a^{-1} b^{-1}$$
$$= a^k a^{-1} b^k b^{-1}$$
$$= a^{k-1} b^{k-1}$$

является следствием равенства (3.5.3). Следовательно, равенство (3.5.8) верно, когда $n = k-1$. Согласно принципу математической индукции, равенство (3.5.8) верно для любого $n \leq 0$.

• Следовательно, равенство (3.5.8) верно для любого $n \in Z$.

Из равенства (3.5.8) следует, что отображение

$$\varphi(n) : a \in G \to a^n \in G$$

является эндоморфизмом абелевой группы $G$. Из равенств (3.5.6), (3.5.5) следует, что отображение

$$\varphi : Z \to \mathrm{End}(Ab, G)$$

является гомоморфизмом кольца $Z$. Согласно определению 2.3.1, отображение $\varphi$ является представлением кольца целых чисел $Z$ в абелевой группе $G$.    □

Теорема 3.5.3. *Множество $G$-чисел, порождённое множеством* $S = \{s_i : i \in I\}$ , *имеет вид*

$$(3.5.10) \qquad J(S) = \left\{ g : g = \prod_{i \in I} s_i^{g^i}, g^i \in Z \right\}$$

*где множество* $\{i \in I : g^i \neq 0\}$ *конечно.*

Доказательство. Мы докажем теорему по индукции, опираясь на теоремы [18]-5.1, страница 94, и 2.7.5.

Для произвольного $s_k \in S$, положим $g^i = \delta_k^i$. Тогда

$$(3.5.11) \qquad s_k = \prod_{i \in I} s_i^{g^i}$$

$s_k \in J(S)$ следует из (3.5.10), (3.5.11).

Пусть $g_1, g_2 \in X_k \subseteq J(S)$. Так как $G$ является абелевой группой, то, согласно утверждению 2.7.5.3, $g_1 g_2 \in J(S)$. Согласно равенству (3.5.10), существуют $Z$-числа $g_1^i, g_2^i, i \in I$, такие, что

$$(3.5.12) \qquad g_1 = \prod_{i \in I} s_i^{g_1^i} \qquad g_2 = \prod_{i \in I} s_i^{g_2^i}$$

где множества

$$(3.5.13) \qquad H_1 = \{i \in I : g_1^i \neq 0\} \qquad H_2 = \{i \in I : g_2^i \neq 0\}$$

конечны. Из равенства (3.5.12) следует, что

$$(3.5.14) \qquad g_1 g_2 = \prod_{i \in I} s_i^{g_1^i} \prod_{i \in I} s_i^{g_2^i} = \prod_{i \in I} (s_i^{g_1^i} s_i^{g_2^i})$$

Равенство

$$(3.5.15) \qquad g_1 g_2 = \prod_{i \in I} s_i^{g_1^i + g_2^i}$$

является следствием равенств (3.5.6), (3.5.14). Из равенства (3.5.13) следует, что множество

$$\{i \in I : g_1^i g_2^i \neq 0\} \subseteq H_1 \cup H_2$$

конечно. Из равенства (3.5.15) следует, что $g_1 g_2 \in J(S)$.          $\square$

Соглашение 3.5.4. *Мы будем пользоваться соглашением о произведении, в котором повторяющийся индекс (один в основании и один в экспоненте) подразумевает произведение по повторяющемуся индексу. В этом случае предполагается известным множество индекса произведения и знак произведения опускается*

$$s_i^{g^i} = \prod_{i \in I} s_i^{g^i}$$

$\square$

Теорема 3.5.5. *Гомоморфизм*

$$f : G \to G'$$

*удовлетворяет равенству*

$$(3.5.16) \qquad\qquad f(a^n) = f(a)^n$$

Доказательство. Из равенства (3.5.1) следует, что равенство (3.5.16) верно, когда $n = 0$.

- Пусть равенство (3.5.16) верно, когда $n = k \geq 0$. Тогда

$$f(a^k) = f(a)^k$$

Равенство

$$f(a^{k+1}) = f(a^k a) = f(a^k)f(a) = f(a)^k f(a)$$
$$= f(a)^{k+1}$$

является следствием равенства (3.5.2). Следовательно, равенство (3.5.16) верно, когда $n = k + 1$. Согласно принципу математической индукции, равенство (3.5.16) верно для любого $n \geq 0$.

- Пусть равенство (3.5.6) верно, когда $n = k \leq 0$. Тогда

$$f(a^k) = f(a)^k$$

Равенство

$$f(a^{k-1}) = f(a^k a^{-1}) = f(a^k)f(a)^{-1} = f(a)^k f(a)^{-1}$$
$$= f(a)^{k-1}$$

является следствием равенства (3.5.3). Следовательно, равенство (3.5.16) верно, когда $n = k - 1$. Согласно принципу математической индукции, равенство (3.5.16) верно для любого $n \leq 0$.

- Следовательно, равенство (3.5.16) верно для любого $n \in Z$.

$\square$

Определение 3.5.6. *Пусть* $S = \{s_i : i \in I\}$ - *множество* $G$-*чисел. Выражение* $s_i^{g^i}$ *называется* **линейной комбинацией** $G$-*чисел* $s_i$. $G$-*число* $g = s_i^{g^i}$ *называется* **линейно зависимым** *от* $G$-*чисел* $s_i$. $\square$

Следующее определение является следствием теорем 2.7.5, 3.5.3 и определения 2.7.4.

Определение 3.5.7. $J(S)$ *называется* **подгруппой, порождённой множеством** $S$, *а* $S$ - **множеством образующих** *подгруппы* $J(S)$. *В частности,* **множеством образующих** *абелевой группы будет такое подмножество* $X \subset G$, *что* $J(X) = G$. $\square$

Следующее определение является следствием теорем 2.7.5, 3.5.3 и определения 2.7.14.

ОПРЕДЕЛЕНИЕ 3.5.8. *Если множество $X \subset G$ является множеством образующих абелевой группы $G$, то любое множество $Y$, $X \subset Y \subset G$ также является множеством образующих абелевой группы $G$. Если существует минимальное множество $X$, порождающее абелевую группу $G$, то такое множество $X$ называется* **квазибазисом** *абелевой группы $G$.* $\square$

ОПРЕДЕЛЕНИЕ 3.5.9. *Множество $G$-чисел[3.9]$g_i$, $i \in I$,* **линейно независимо**, *если $w = e$ следует из уравнения*

$$g_i^{w^i} = e$$

*В противном случае, множество $G$-чисел $g_i$, $i \in I$,* **линейно зависимо**. $\square$

ОПРЕДЕЛЕНИЕ 3.5.10. *Пусть множество $G$-чисел $S = \{s_i : i \in I\}$ является квазибазисом. Если множество $G$-чисел $S$ линейно независимо, то квазибазис $S$ называется* **базисом** *абелевой группы $G$.* $\square$

## 3.6. Свободная абелева группа

ОПРЕДЕЛЕНИЕ 3.6.1. *Пусть $S$ - множество и[3.10] $Ab(S, *)$ - категория, объектами которой являются отображения*

$$f : S \to G$$

*множества $S$ в мультипликативные абелевые группы. Если*

$$f : S \to G$$
$$f' : S \to G'$$

*объекты категории $Ab(S, *)$, то морфизм из $f$ в $f'$ - это гомоморфизм групп*

$$g : G \to G'$$

*такой, что диаграмма*

$$G \xrightarrow{g} G'$$

*коммутативна.* $\square$

ОПРЕДЕЛЕНИЕ 3.6.2. *Пусть $S$ - множество и[3.10] $Ab(S, +)$ - категория, объектами которой являются отображения*

$$f : S \to G$$

*множества $S$ в аддитивные абелевые группы. Если*

$$f : S \to G$$
$$f' : S \to G'$$

*объекты категории $Ab(S, +)$, то морфизм из $f$ в $f'$ - это гомоморфизм групп*

$$g : G \to G'$$

*такой, что диаграмма*

$$G \xrightarrow{g} G'$$

*коммутативна.* $\square$

---

[3.9]  Я следую определению в [1], страница 100.
[3.10]  Смотри также определение в [1], страницы 56, 57.
[3.11]  Смотри аналогичное утверждение в [1], страница 57.

Теорема 3.6.3. *Универсально отталкивающий объект $G$ категории $Ab(S, *)$ существует.* [3.11]

3.6.3.1: $S \subseteq G$.

3.6.3.2: *Множество $S$ порождает мультипликавную абелевую группу $G$.*

*Абелева группа $G$ называется* **free Abelian group**, *порождённой множеством $S$.*

Доказательство. Пусть $G$ - множество отображений

$$h : S \to Z$$

таких, что множество

$$H = \{x \in S : h(x) \neq 0\}$$

конечно.

Лемма 3.6.4. *Определим на множестве $G$ операцию*

$$(h_1 h_2)(x) = h_1(x) + h_2(x)$$

*Множество $G$ является мультипликативной абелевой группой.*

Доказательство. Пусть $h_1, h_2 \in G$. Согласно определению, множества

$$H_1 = \{x \in S : h_1(x) \neq 0\}$$
$$H_2 = \{x \in S : h_2(x) \neq 0\}$$

конечны. Тогда конечно множество

$$\{x \in S : (h_1 h_2)(x) \neq 0\} \subseteq H_1 \cup H_2$$

Следовательно, $h_1 h_2 \in G$. Свойства абелевой группы очевидны.                    ⊙

Пусть $k \in Z$, $x \in S$. Рассмотрим отображение $h = k * x$, определённое равенством

$$k * x(y) = \begin{cases} k & y = x \\ 0 & y \neq x \end{cases}$$

Отображение

$$f : x \in S \to 1 * x \in G$$

---

Теорема 3.6.6. *Универсально отталкивающий объект $G$ категории $Ab(S, +)$ существует.* [3.11]

3.6.6.1: $S \subseteq G$.

3.6.6.2: *Множество $S$ порождает аддивную абелеву группу $G$.*

*Абелева группа $G$ называется* **free Abelian group**, *порождённой множеством $S$.*

Доказательство. Пусть $G$ - множество отображений

$$h : S \to Z$$

таких, что множество

$$H = \{x \in S : h(x) \neq 0\}$$

конечно.

Лемма 3.6.7. *Определим на множестве $G$ операцию*

$$(h_1 + h_2)(x) = h_1(x) + h_2(x)$$

*Множество $G$ является аддитивной абелевой группой.*

Доказательство. Пусть $h_1, h_2 \in G$. Согласно определению, множества

$$H_1 = \{x \in S : h_1(x) \neq 0\}$$
$$H_2 = \{x \in S : h_2(x) \neq 0\}$$

конечны. Тогда конечно множество

$$\{x \in S : (h_1 + h_2)(x) \neq 0\} \subseteq H_1 \cup H_2$$

Следовательно, $h_1 + h_2 \in G$. Свойства абелевой группы очевидны.                    ⊙

Пусть $k \in Z$, $x \in S$. Рассмотрим отображение $h = k * x$, определённое равенством

$$k * x(y) = \begin{cases} k & y = x \\ 0 & y \neq x \end{cases}$$

Отображение

$$f : x \in S \to 1 * x \in G$$

инъективно и позволяет отождествить $S$ и образ $f(S) \subseteq G$ (утверждение 3.6.3.1) и пользоваться записью

$$(1*x)^k = k*x$$

Согласно определению множества $G$, любое отображение $h \in G$ может быть записано в виде

$$(3.6.1) \qquad h = \prod_{x \in S} k_x * x$$

где $k_x \in Z$, $x \in S$, и множество $\{x : k_x \neq 0\}$ конечно. Утверждение 3.6.3.2 является следствием теоремы 3.5.3.

**Лемма 3.6.5.** *Представление* (3.6.1) *отображения $h$ единственно.*

Доказательство. Если отображение $h$ допускает представления

$$(3.6.2) \qquad \prod_{x \in S} k_{1x} * x = \prod_{x \in S} k_{2x} * x$$

то равенство

$$(3.6.3) \qquad \prod_{x \in S}(k_{1x} - k_{2x}) * x = e$$

является следствием равенства (3.6.2). $k_{1x} = k_{2x}$, $x \in S$, следует из равенства (3.6.3). $\odot$

Пусть

$$f' : S \to G'$$

отображение множества $S$ в абелевую группу $G'$. Мы определим гомоморфизм

$$g : G \to G'$$

таким образом, чтобы диаграмма

$$(3.6.4)$$

была коммутативной. Из диаграммы (3.6.4) следует, что

$$(3.6.5) \qquad g(1*x) = f'(x)$$

инъективно и позволяет отождествить $S$ и образ $f(S) \subseteq G$ (утверждение 3.6.6.1) и пользоваться записью

$$k(1*x) = k*x$$

Согласно определению множества $G$, любое отображение $h \in G$ может быть записано в виде

$$(3.6.6) \qquad h = \sum_{x \in S} k_x * x$$

где $k_x \in Z$, $x \in S$, и множество $\{x : k_x \neq 0\}$ конечно. Утверждение 3.6.6.2 является следствием теоремы 3.4.3.

**Лемма 3.6.8.** *Представление* (3.6.6) *отображения $h$ единственно.*

Доказательство. Если отображение $h$ допускает представления

$$(3.6.7) \qquad \sum_{x \in S} k_{1x} * x = \sum_{x \in S} k_{2x} x$$

то равенство

$$(3.6.8) \qquad \sum_{x \in S}(k_{1x} - k_{2x}) * x = 0$$

является следствием равенства (3.6.7). $k_{1x} = k_{2x}$, $x \in S$, следует из равенства (3.6.8). $\odot$

Пусть

$$f' : S \to G'$$

отображение множества $S$ в абелевую группу $G'$. Мы определим гомоморфизм

$$g : G \to G'$$

таким образом, чтобы диаграмма

$$(3.6.9)$$

была коммутативной. Из диаграммы (3.6.9) следует, что

$$(3.6.10) \qquad g(1*x) = f'(x)$$

Согласно теореме 3.5.5, так как отоб-
ражение $g$ - гомоморфизм абелевой груп-
пы, то равенство

$$g(h) = g\left(\prod_{x \in S} k_x * x\right) = \prod_{x \in S} g(k_x * x)$$
$$= \prod_{x \in S} g(1 * x)^{k_x} = \prod_{x \in S} f'(x)^{k_x}$$

для любого отображения $h \in G$ явля-
ется следствием равенств (3.6.5), (3.6.1).
Следовательно, гомоморфизм $g$ опре-
делён однозначно. Согласно определе-
ниям 2.2.2, 3.6.1, абелевая группа $G$
является свободной группой.         □

Согласно теореме 3.4.5, так как отоб-
ражение $g$ - гомоморфизм абелевой груп-
пы, то равенство

$$g(h) = g\left(\sum_{x \in S} k_x * x\right) = \sum_{x \in S} g(k_x * x)$$
$$= \sum_{x \in S} k_x g(1 * x) = \sum_{x \in S} k_x f'(x)$$

для любого отображения $h \in G$ явля-
ется следствием равенств (3.6.10), (3.6.6).
Следовательно, гомоморфизм $g$ опре-
делён однозначно. Согласно определе-
ниям 2.2.2, 3.6.2, абелевая группа $G$
является свободной группой.         □

Теорема 3.6.9. *Мультиплика-
тивная абелевая группа $G$ свободна
тогда и только тогда, когда абелевая
группа $G$ имеет базис $S = \{s_i : i \in
I\}$ и множество $G$-чисел $S$ линейно
независимо.*

Теорема 3.6.10. *Аддитивная
абелевая группа $G$ свободна тогда и
только тогда, когда абелевая группа
$G$ имеет базис $S = \{s_i : i \in I\}$
и множество $G$-чисел $S$ линейно
независимо.*

Доказательство. Согласно опре-
делению 3.6.1 и доказательству теоре-
мы 3.6.3, существует множество $S \subseteq G$
такое, что диаграмма

Доказательство. Согласно опре-
делению 3.6.2 и доказательству теоре-
мы 3.6.6, существует множество $S \subseteq G$
такое, что диаграмма

$$\begin{array}{ccc} G & \xrightarrow{\ g\ } & G' \\ {\scriptstyle f}\nwarrow & & \nearrow{\scriptstyle f'} \\ & S & \end{array}$$

$$\begin{array}{ccc} G & \xrightarrow{\ g\ } & G' \\ {\scriptstyle f}\nwarrow & & \nearrow{\scriptstyle f'} \\ & S & \end{array}$$

коммутативна                        □

коммутативна                        □

Теорема 3.6.11. *Пусть $A$ - абелевый аддитивный моноид. Рассмотрим
категорию $\mathcal{K}(A)$ гомоморфизмов моноида $A$ в абелевые группы. Если $B$ и $C$ -
абелевые группы, то для гомоморфизмов*

$$f : A \to B$$
$$g : A \to C$$

*морфизм категории $\mathcal{K}(A)$ - это гомоморфизм*

$$h : B \to C$$

*такой, что диаграмма*

$$\begin{array}{ccc} A & \xrightarrow{\ f\ } & B \\ {\scriptstyle g}\searrow & & \swarrow{\scriptstyle h} \\ & C & \end{array}$$

*коммутативна. Существует абелевая группа $\mathcal{K}(A)$ и универсально оттал-*

кивающий гомоморфизм

$$\gamma : A \to K(A)$$

*Группа $K(A)$ называется* **группа Гротендика**, *а отображение $\gamma$ называется* **каноническим отображением**.

ЗАМЕЧАНИЕ 3.6.12. *Согласно определению 2.2.2, если $B$ - абелева группа и*

$$f : A \to B$$

*гомоморфизм, то существует единственный гомоморфизм*

$$g : K(A) \to B$$

*такой, что диаграмма*

*коммутативна.* ⊙

ДОКАЗАТЕЛЬСТВО. Пусть $F(A)$ - свободная абелева группа, порождённая множеством $A$. Обозначим $[a]$ образующую группы $F(A)$, порождённую $A$-числом $a$. Пусть $B$ - подгруппа, порождённая всеми $F(A)$-числами вида

$$[a + b] - [a] - [b] \qquad a, b \in A$$

Положим $K(A) = F(A)/B$. Мы определим гомоморфизм $\gamma$ таким образом, что следующая диаграмма

$$f : a \in A \to [a] \in F(A)$$

$$g : a \in F(A) \to aB \in K(A)$$

коммутативна. □

ОПРЕДЕЛЕНИЕ 3.6.13. *В аддитивном моноиде $A$ выполняется* **закон сокращения**, *если, для любых $A$-чисел $a$, $b$, $c$, из равенства $a + c = b + c$ следует $a = b$.* □

ТЕОРЕМА 3.6.14. *Если в аддитивном моноиде $A$ выполняется закон сокращения, то каноническое отображение*

$$\gamma : A \to K(A)$$

*инъективно.*

ДОКАЗАТЕЛЬСТВО. Рассмотрим кортежи $(a, b)$, $a$, $b \in A$. Кортеж $(a, b)$ эквивалентен кортежу $(a', b')$, если

$$a + b' = a' + b$$

Мы определим сложение равенством

$$(a, b) + (c, d) = (a + c, b + d)$$

Классы эквивалентности кортежей образуют группу, нулевым элементом которой является класс кортежа $(0,0)$. Кортеж $(b,a)$ является противоположным для кортежа $(a,b)$. Из закона сокращения следует, что гомоморфизм

$$a \to \text{класс кортежа } (0,a)$$

инъективен.                                                                            □

Пример 3.6.15. *Определим на множестве чисел*

$$A = \{0, 1, ...\}$$

*сложение согласно равенству*

$$(3.6.11) \qquad\qquad a + b = \max(a, b)$$

Лемма 3.6.16. *Сложение* (3.6.11) *коммутативно.*

Доказательство. *Лемма является следствием определения* (3.6.11). ⊙

Лемма 3.6.17. *Сложение* (3.6.11) *ассоциативно.*

Доказательство. *Для доказательства равенства*

$$(3.6.12) \qquad\qquad a + (b + c) = (a + b) + c$$

*достаточно рассмотреть следующие соотношения.*

- *Если* $a \leq \max(b,c)$, *то либо* $a \leq b$, *либо* $a \leq c$.
  - *Пусть* $b \leq c$. *Тогда*

$$(3.6.13) \qquad\qquad a + (b + c) = a + c = c$$

$$(3.6.14) \qquad\qquad (a + b) + c = \max(a, b) + c = c$$

    *Равенство* (3.6.12) *является следствием равенств* (3.6.13), (3.6.14).
  - *Пусть* $c \leq b$. *Тогда*

$$(3.6.15) \qquad\qquad a + (b + c) = a + b = b$$

$$(3.6.16) \qquad\qquad (a + b) + c = b + c = b$$

    *Равенство* (3.6.12) *является следствием равенств* (3.6.15), (3.6.16).
- *Если* $a > \max(b,c)$, *то* $a > b$, $a > c$. *Следовательно*

$$(3.6.17) \qquad\qquad a + (b + c) = a$$

$$(3.6.18) \qquad\qquad (a + b) + c = a + c = a$$

    *Равенство* (3.6.12) *является следствием равенств* (3.6.17), (3.6.18).

*Лемма является следствием равенства* (3.6.12).                          ⊙

Лемма 3.6.18. *A-число* 0 *является нулевым элементом сложения* (3.6.11).

Доказательство. *Так как* $\forall a \in A$, $0 \leq a$, *лемма является следствием определения* (3.6.11).                                                          ⊙

Теорема 3.6.19. *Множество $A$, оснащённое сложением (3.6.11), является абелевым аддитивным моноидом.*

Доказательство. *Теорема является следствием лемм 3.6.16, 3.6.17, 3.6.18 и определения 3.1.1.* ⊙

Теорема 3.6.20. *В аддитивном моноиде $A$ закон сокращения не выполняется.*

Доказательство. *Для доказательства теоремы достаточно убедиться, что*

$$2 + 9 = 9$$
$$4 + 9 = 9$$

*Однако, $2 \neq 4$.* ⊙ □

## 3.7. Прямая сумма абелевых групп

Определение 3.7.1. *Копроизведение в категории мультипликативных абелевых групп $Ab(*)$ называется **прямой суммой**. [3.12] Мы будем пользоваться записью $A \oplus B$ для прямой суммы абелевых групп $A$ и $B$.* □

Определение 3.7.2. *Копроизведение в категории аддитивных абелевых групп $Ab(+)$ называется **прямой суммой**. [3.12] Мы будем пользоваться записью $A \oplus B$ для прямой суммы абелевых групп $A$ и $B$.* □

Теорема 3.7.3. *Пусть $\{A_i, i \in I\}$ - семейство мультипликативных абелевых групп. Пусть*

$$A \subseteq \prod_{i \in I} A_i$$

*такое множество, что*

$$A = \{(x_i, i \in I) : |\{i : x_i \neq e\}| < \infty\}$$

*Тогда* [3.13]

$$A = \bigoplus_{i \in I} A_i$$

Теорема 3.7.4. *Пусть $\{A_i, i \in I\}$ - семейство аддитивных абелевых групп. Пусть*

$$A \subseteq \prod_{i \in I} A_i$$

*такое множество, что*

$$A = \{(x_i, i \in I) : |\{i : x_i \neq 0\}| < \infty\}$$

*Тогда* [3.13]

$$A = \bigoplus_{i \in I} A_i$$

Доказательство. Согласно теореме 2.2.6, существует абелева группа $B = \prod A_i$. Согласно утверждению 2.2.6.3, мы можем представить $B$-число $a$ в виде кортежа $\{a_i, i \in I\}$ $A_i$-

Доказательство. Согласно теореме 2.2.6, существует абелева группа $B = \prod A_i$. Согласно утверждению 2.2.6.3, мы можем представить $B$-число $a$ в виде кортежа $\{a_i, i \in I\}$ $A_i$-

---

[3.12] Определение дано согласно [1], страница 55.
[3.13] Смотри также предложение [1]-10, страница 55.

чисел. Согласно утверждению 2.2.6.4, произведение в группе $B$ определено равенством

$$(3.7.1) \qquad (a_i, i \in I)(b_i, i \in I) \\ = (a_i b_i, i \in I)$$

Согласно построению $A$ является подгруппой абелевой группы $B$. Отображение

$$\lambda_j : A_j \to A$$

определённое равенством

$$(3.7.2) \qquad \lambda_j(x) = (\delta_j^i x, i \in I)$$

является инъективным гомоморфизмом.

Пусть

$$\{f_i : A_i \to , i \in I\}$$

семейство гомоморфизмов в мультипликативную абелеву группу $C$. Мы определим отображение

$$f : A \to C$$

пользуясь равенством

$$(3.7.3) \qquad f\left(\bigoplus_{i \in I} x_i\right) = \prod_{i \in I} f_i(x_i)$$

Произведение в правой части равенства (3.7.3) конечно, так как все сомножители, кроме конечного числа, равны $e$. Из равенства

$$f_i(x_i y_i) = f_i(x_i) f_i(y_i)$$

и равенства (3.7.3) следует, что

$$f(\bigoplus_{i \in I} (x_i y_i)) = \prod_{i \in I} f_i(x_i y_i) \\ = \prod_{i \in I} (f_i(x_i) f_i(y_i)) \\ = \prod_{i \in I} f_i(x_i) \prod_{i \in I} f_i(y_i) \\ = f(\bigoplus_{i \in I} x_i) f(\bigoplus_{i \in I} y_i)$$

Следовательно, отображение $f$ является гомоморфизмом абелевой группы. Равенство

$$f \circ \lambda_j(x) = \prod_{i \in I} f_i(\delta_j^i x) = f_j(x)$$

чисел. Согласно утверждению 2.2.6.4, сумма в группе $B$ определена равенством

$$(3.7.4) \qquad (a_i, i \in I) + (b_i, i \in I) \\ = (a_i + b_i, i \in I)$$

Согласно построению $A$ является подгруппой абелевой группы $B$. Отображение

$$\lambda_j : A_j \to A$$

определённое равенством

$$(3.7.5) \qquad \lambda_j(x) = (\delta_j^i x, i \in I)$$

является инъективным гомоморфизмом.

Пусть

$$\{f_i : A_i \to +, i \in I\}$$

семейство гомоморфизмов в аддитивную абелеву группу $C$. Мы определим отображение

$$f : A \to C$$

пользуясь равенством

$$(3.7.6) \qquad f\left(\bigoplus_{i \in I} x_i\right) = \sum_{i \in I} f_i(x_i)$$

Сумма в правой части равенства (3.7.6) конечна, так как все слагаемые, кроме конечного числа, равны 0. Из равенства

$$f_i(x_i + y_i) = f_i(x_i) + f_i(y_i)$$

и равенства (3.7.6) следует, что

$$f(\bigoplus_{i \in I} (x_i + y_i)) = \sum_{i \in I} f_i(x_i + y_i) \\ = \sum_{i \in I} (f_i(x_i) + f_i(y_i)) \\ = \sum_{i \in I} f_i(x_i) + \sum_{i \in I} f_i(y_i) \\ = f(\bigoplus_{i \in I} x_i) + f(\bigoplus_{i \in I} y_i)$$

Следовательно, отображение $f$ является гомоморфизмом абелевой группы. Равенство

$$f \circ \lambda_j(x) = \sum_{i \in I} f_i(\delta_j^i x) = f_j(x)$$

является следствием равенств (3.7.2), (3.7.3). Так как отображение $\lambda_i$ инъективно, то отображение $f$ определено однозначно. Следовательно, теорема является следствием определений 2.2.12, 3.7.1. □

является следствием равенств (3.7.5), (3.7.6). Так как отображение $\lambda_i$ инъективно, то отображение $f$ определено однозначно. Следовательно, теорема является следствием определений 2.2.12, 3.7.2. □

ТЕОРЕМА 3.7.5. *Пусть $B$, $C$ - подгруппы мультипликавной абелевой группы $A$. Пусть $A = BC$, $B \cap C = e$. Тогда* [3.14] *$A = B \oplus C$.*

ТЕОРЕМА 3.7.6. *Пусть $B$, $C$ - подгруппы аддивной абелевой группы $A$. Пусть $A = B + C$, $B \cap C = 0$. Тогда* [3.14] *$A = B \oplus C$.*

Доказательство. Отображение

(3.7.7)   $f : (x, y) \in B \times C \to xy \in A$

является гомоморфизмом групп. Это утверждение следует из равенств

$$f(x_1 x_2, y_1 y_2) = x_1 x_2 y_1 y_2$$
$$= x_1 y_1 x_2 y_2$$
$$= f(x_1, y_1) f(x_2, y_2)$$

Пусть $f(x, y) = e$. Из равенства (3.7.7) следует, что

(3.7.8)   $xy = e \quad x \in B \quad y \in C$

Из равенства (3.7.8) следует, что $x = y^{-1} \in C$. Следовательно, $x = y = e$. Следовательно, отображение $f$ является изоморфизмом групп.

Согласно теореме 3.7.3

(3.7.9)   $B \times C = B \oplus C$

Теорема является следствием равенства (3.7.9). □

Доказательство. Отображение

(3.7.10)   $f : (x, y) \in B \times C \to x + y \in A$

является гомоморфизмом групп. Это утверждение следует из равенств

$$f(x_1 + x_2, y_1 + y_2) = x_1 + x_2 + y_1 + y_2$$
$$= x_1 + y_1 + x_2 + y_2$$
$$= f(x_1, y_1) + f(x_2, y_2)$$

Пусть $f(x, y) = 0$. Из равенства (3.7.10) следует, что

(3.7.11)   $x + y = 0 \quad x \in B \quad y \in C$

Из равенства (3.7.11) следует, что $x = -y \in C$. Следовательно, $x = y = 0$. Следовательно, отображение $f$ является изоморфизмом групп.

Согласно теореме 3.7.4

(3.7.12)   $B \times C = B \oplus C$

Теорема является следствием равенства (3.7.12). □

ПРИМЕР 3.7.7. *Пусть абелевая группа $H$ является подгруппой мультипликативной абелевой группы $G$. Найдём $H \oplus G$. Мы будем пользоваться теоремой 3.7.3. Положим $A = H \times G$. Мы определим умножение равенством*

$$(a, b)(c, d) = (ac, bd)$$

*Следовательно, если $a \in H$, то кортежи $(a, 0)$ и $(0, a)$ представляют различные $(H \oplus G)$-числа.* □

ПРИМЕР 3.7.8. *Пусть абелевая группа $H$ является подгруппой аддивной абелевой группы $G$. Найдём $H \oplus G$. Мы будем пользоваться теоремой 3.7.4. Положим $A = H \times G$. Мы определим сложение равенством*

$$(a, b) + (c, d) = (a + c, b + d)$$

*Следовательно, если $a \in H$, то кортежи $(a, 0)$ и $(0, a)$ представляют различные $(H \oplus G)$-числа.* □

---

[3.14]  Смотри также замечание в [1], страница 56.

Теорема 3.7.9. *Пусть*

$$f : A \to B$$

*сюръективный гомоморфизм абелевых групп и* $C = \ker f$. *Если* $B$ - *свободная группа, то существует подгруппа* $D$ *группы* $A$ *такая, что ограничение* $f$ *на* $D$ *индуцирует изоморфизм* $D$ *на* $B$ *и* $A = C \oplus D$.

Доказательство.

Пусть $A$, $B$ - мультипликативные абелевые группы. Пусть $(b_i, i \in I)$ - базис группы $B$. Для каждого $i \in I$, пусть $a_i$ - $A$-число такое, что

$$(3.7.13) \qquad f(a_i) = b_i$$

Пусть $D$ - подгруппа в $A$, порождённая $A$-числами $a_i$, $i \in I$. Пусть

$$(3.7.14) \qquad \prod_{i \in I} a_i^{n^i} = e$$

для целых $n^i$ таких, что множество $\{i \in I : n^i \neq 0\}$ конечно. Равенство

$$(3.7.15) \qquad e = \prod_{i \in I} f(a_i)^{n^i} = \prod_{i \in I} b_i^{n^i}$$

следует из равенств (3.7.13), (3.7.14) и теорем 3.3.1, 3.5.5. Так как $(b_i, i \in I)$ - базис, то согласно определениям 3.5.9, 3.5.10, утверждение $n_i = 0$, $i \in I$, следует из равенства (3.7.15). Следовательно, множество $(a_i, i \in I)$ - базис подгруппы $D$.

Пусть $a \in D$ и $f(a) = e$. Так как множество $(a_i, i \in I)$ - базис подгруппы $D$, то существуют целые числа $n^i$ такие, что

$$(3.7.16) \qquad a = a_i^{n^i}$$

и множество $\{i \in I : n^i \neq 0\}$ конечно. Следовательно, равенство

$$(3.7.17) \qquad e = f(a) = f(a_i)^{n^i} = b_i^{n^i}$$

следует из равенства (3.7.16) и теорем 3.3.1, 3.5.5. Так как множество $(b_i, i \in I)$ является базисом группы $B$, $n^i = 0$, $i \in I$, является следствием равенства (3.7.17). Следовательно, $a = e$ и $D \cap C = \{e\}$.

Пусть $a \in A$. Так как $f(x) \in B$ и множество $(b_i, i \in I)$ является базисом группы $B$, то существуют целые

Пусть $A$, $B$ - аддитивные абелевые группы. Пусть $(b_i, i \in I)$ - базис группы $B$. Для каждого $i \in I$, пусть $a_i$ - $A$-число такое, что

$$(3.7.21) \qquad f(a_i) = b_i$$

Пусть $D$ - подгруппа в $A$, порождённая $A$-числами $a_i$, $i \in I$. Пусть

$$(3.7.22) \qquad \sum_{i \in I} n^i a_i = 0$$

для целых $n^i$ таких, что множество $\{i \in I : n^i \neq 0\}$ конечно. Равенство

$$(3.7.23) \qquad 0 = \sum_{i \in I} n^i f(a_i) = \sum_{i \in I} n^i b_i$$

следует из равенств (3.7.21), (3.7.22) и теорем 3.3.2, 3.4.5. Так как $(b_i, i \in I)$ - базис, то согласно определениям 3.4.9, 3.4.10, утверждение $n_i = 0$, $i \in I$, следует из равенства (3.7.23). Следовательно, множество $(a_i, i \in I)$ - базис подгруппы $D$.

Пусть $a \in D$ и $f(a) = 0$. Так как множество $(a_i, i \in I)$ - базис подгруппы $D$, то существуют целые числа $n^i$ такие, что

$$(3.7.24) \qquad a = n^i a_i$$

и множество $\{i \in I : n^i \neq 0\}$ конечно. Следовательно, равенство

$$(3.7.25) \qquad 0 = f(a) = n^i f(a_i) = n^i b_i$$

следует из равенства (3.7.24) и теорем 3.3.2, 3.4.5. Так как множество $(b_i, i \in I)$ является базисом группы $B$, $n^i = 0$, $i \in I$, является следствием равенства (3.7.25). Следовательно, $a = 0$ и $D \cap C = \{0\}$.

Пусть $a \in A$. Так как $f(x) \in B$ и множество $(b_i, i \in I)$ является базисом группы $B$, то существуют целые

числа $n^i$ такие, что

$$(3.7.18) \qquad f(a) = b_i^{n^i}$$

и множество $\{i \in I : n^i \neq 0\}$ конечно. Равенство

$$(3.7.19) \qquad \begin{aligned} f(xa_i^{-n^i}) &= f(x)f(a_i)^{-n^i} \\ &= f(x)b_i^{-n^i} = e \end{aligned}$$

является следствием равенств (3.7.13), (3.7.18). Из равенства (3.7.19) следует, что

$$(3.7.20) \qquad b = aa_i^{-n^i} \in C = \ker f$$

Из утверждения (3.7.20) следует, что $a \in C + D$. Согласно теореме 3.7.5, $A = C \oplus D$.

числа $n^i$ такие, что

$$(3.7.26) \qquad f(a) = n^i b_i$$

и множество $\{i \in I : n^i \neq 0\}$ конечно. Равенство

$$(3.7.27) \qquad \begin{aligned} f(x - n^i a_i) &= f(x) - n^i f(a_i) \\ &= f(x) - n^i b_i = 0 \end{aligned}$$

является следствием равенств (3.7.21), (3.7.26). Из равенства (3.7.27) следует, что

$$(3.7.28) \qquad b = a - n^i a_i \in C = \ker f$$

Из утверждения (3.7.28) следует, что $a \in C + D$. Согласно теореме 3.7.6, $A = C \oplus D$.

$\square$

**Теорема 3.7.10.** *Пусть $A$ - подгруппа абелевой группы $B$. Если группа $B/A$ - свободна, то*

$$(3.7.29) \qquad B = A \oplus B/A$$

Доказательство.

Пусть $B$ - мультипликативная абелева группа. Согласно примерам 3.3.15, 3.3.17, последовательность гомоморфизмов

$$e \longrightarrow A \overset{f}{\longrightarrow} B \overset{g}{\longrightarrow} B/A \longrightarrow e$$

является точной. Согласно примеру 3.3.15, гомоморфизм $f$ является мономорфизмом. Согласно определению 3.3.14, $A = \ker g$. Согласно примеру 3.3.15, гомоморфизм $g$ является эпиморфизмом. Согласно теореме 3.7.9, существует подгруппа $D$ группы $B$, изоморфная группе $B/A$, такая, что

$$(3.7.30) \qquad B = A \oplus D$$

Если отождествить $D$-число $d$ и множество $dA$, то равенство (3.7.29) является следствием равенства (3.7.30).

Пусть $B$ - аддитивная абелева группа. Согласно примерам 3.3.18, 3.3.20, последовательность гомоморфизмов

$$0 \longrightarrow A \overset{f}{\longrightarrow} B \overset{g}{\longrightarrow} B/A \longrightarrow 0$$

является точной. Согласно примеру 3.3.18, гомоморфизм $f$ является мономорфизмом. Согласно определению 3.3.14, $A = \ker g$. Согласно примеру 3.3.18, гомоморфизм $g$ является эпиморфизмом. Согласно теореме 3.7.9, существует подгруппа $D$ группы $B$, изоморфная группе $B/A$, такая, что

$$(3.7.31) \qquad B = A \oplus D$$

Если отождествить $D$-число $d$ и множество $dA$, то равенство (3.7.29) является следствием равенства (3.7.31). $\square$

Теорема 3.7.11. *Пусть*

$$f : A \to B$$

*гомоморфизм абелевых групп. Если* Im $f$ *- свободная группа, то*

(3.7.32)                          $A = \ker f \oplus \operatorname{Im} f$

Доказательство. Поскольку $A/\ker f = \operatorname{Im} f$, равенство (3.7.32) является следствием равенства (3.7.29).                                           □

## 3.8. Унитальное расширение кольца

Пусть $D$ - кольцо.

Лемма 3.8.1. *Пусть* $n \in Z$. *Отображение*

(3.8.1)                          $n : d \in D \to nd \in D$

*является эндоморфизмом абелевой группы* $D$.

Доказательство. Утверждение

$$nd \in D$$

является следствием определения 3.4.1. Согласно теореме 3.4.2 и определению 2.3.1, отображение (3.8.1) является эндоморфизмом абелевой группы $D$.                                           □

Лемма 3.8.2. *Пусть* $a \in D$. *Отображение*

(3.8.2)                          $a : d \in D \to ad \in D$

*является эндоморфизмом абелевой группы* $D$.

Доказательство. Утверждение

$$ad \in D$$

является следствием определения 3.2.1. Согласно утверждению 3.2.1.3, и теоремам 3.2.3, 3.3.2, отображение (3.8.2) является эндоморфизмом абелевой группы $D$.                                           □

Лемма 3.8.3. *Пусть* $n \in Z$, $a \in D$. *Отображение*

(3.8.3)                  $a \oplus n : d \in D \to (a \oplus n)d \in D$

*определённое равенством*

(3.8.4)                          $(a \oplus n)d = ad + nd$

*является эндоморфизмом абелевой группы* $D$.

Доказательство. Так как $D$ является абелевой группой, то утверждение

$$nd + ad \in D \quad n \in Z \quad a \in D$$

является следствием лемм 3.8.1, 3.8.2. Следовательно, для любого $Z$-числа $n$ и любого $D$-числа $a$, мы определили отображение (3.8.3). Согласно леммам 3.8.1,

3.8.2 и теореме 3.3.4, отображение (3.8.3) является эндоморфизмом абелевой группы $V$. □

**Теорема 3.8.4.** *Множество отображений* $a \oplus n$ *порождает*[3.15]*кольцо* $D_{(1)}$ *где сложение определено равенством*

$$(3.8.5) \qquad (a_1 \oplus n_1) + (a_2 \oplus n_2) = (a_1 + a_2) \oplus (n_1 + n_2)$$

*и произведение определено равенством*

$$(3.8.6) \qquad (a_1 \oplus n_1) \circ (a_2 \oplus n_2) = (a_1 a_2 + n_2 a_1 + n_1 a_2) \oplus (n_1 n_2)$$

3.8.4.1: *Если кольцо* $D$ *имеет единицу, то* $Z \subseteq D$, $D_{(1)} = D$.
3.8.4.2: *Если кольцо* $D$ *является идеалом* $Z$, *то* $D \subseteq Z$, $D_{(1)} = Z$.
3.8.4.3: *В противном случае* $D_{(1)} = D \oplus Z$.
3.8.4.4: *кольцо* $D$ *является идеалом кольца* $D_{(1)}$.
*Кольцо* $D_{(1)}$ *называется* **унитальным расширением** *кольца* $D$.

Доказательство. Согласно теореме 3.3.4, равенство

$$((a_1 \oplus n_1) + (a_2 \oplus n_2))d$$
$$= (a_1 \oplus n_1)d + (a_2 \oplus n_2)d$$
$$(3.8.7) \qquad = a_1 d + n_1 d + a_2 d + n_2 d$$
$$= a_1 d + a_2 d + n_1 d + n_2 d$$
$$= (a_1 + a_2)d + (n_1 + n_2)d$$
$$= (a_1 + a_2) \oplus (n_1 + n_2)d$$

является следствием равенства (3.8.4). Равенство (3.8.5) является следствием равенства (3.8.7).

Согласно теореме 2.1.10, равенство

$$((a_1 \oplus n_1) \circ (a_2 \oplus n_2))d$$
$$= (a_1 \oplus n_1)((a_2 \oplus n_2)d)$$
$$= (a_1 \oplus n_1)(a_2 d + n_2 d)$$
$$(3.8.8) \qquad = (a_1 \oplus n_1)(a_2 d)$$
$$+ (a_1 \oplus n_1)(n_2 d)$$
$$= a_1(a_2 d) + n_1(a_2 d)$$
$$+ a_1(n_2 d) + n_1(n_2 d)$$

является следствием равенства (3.8.4). Равенство

$$((a_1 \oplus n_1) \circ (a_2 \oplus n_2))d$$
$$(3.8.9) \qquad = (a_1 a_2)d + (n_1 a_2)d$$
$$+ (n_2 a_1)d + (n_1 n_2)d$$

является следствием равенств (3.4.5), (3.8.8), утверждений 3.1.1.1, 3.2.1.2,

---

[3.15] Смотри определение унитального расширения также на страницах [7]-52, [8]-64.

теоремы 3.4.5,   леммы 3.8.2. Равенство

$$((a_1 \oplus n_1) \circ (a_2 \oplus n_2))d$$
$$=(a_1 a_2 + n_1 a_2 + n_2 a_1 + n_1 n_2)d$$
$$=((a_1 a_2 + n_1 a_2 + n_2 a_1) \oplus n_1 n_2)d$$

а следовательно и равенство (3.8.6), является следствием равенств (3.8.4), (3.8.9) и лемм 3.8.1, 3.8.2.

Пусть кольцо $D$ имеет единицу $e$. Тогда мы можем отождествить $n \in Z$ и $ne \in D$. Равенство

(3.8.10)                          $(a \oplus n)d = ad + nd = (a + ne)d$

является следствием равенства

$$\boxed{(3.8.4) \quad (a \oplus n)d = ad + nd}$$

Утверждение 3.8.4.1 является следствием равенства (3.8.10).

Пусть кольцо $D$ является идеалом $Z$. Тогда $D \subseteq Z$. Равенство

(3.8.11)                          $(a \oplus n)d = ad + nd = (a + n)d$

является следствием равенства (3.8.4). Утверждение 3.8.4.2 является следствием равенства (3.8.11).

Утверждение 3.8.4.4 является следствием леммы 3.8.3.                    $\square$

# $D$-модуль

## 4.1. Модуль над коммутативным кольцом

Согласно теореме 3.3.4, рассматривая представление

$$f : D \longrightarrow V$$

кольца $D$ в абелевой группе $V$, мы полагаем

$$(4.1.1) \qquad (f(a) + f(b))(x) = f(a)(x) + f(b)(x)$$

Согласно определению 2.3.1, отображение $f$ является гомоморфизмом кольца $D$. Следовательно

$$(4.1.2) \qquad f(a + b) = f(a) + f(b)$$

Равенство

$$(4.1.3) \qquad f(a + b)(x) = f(a)(x) + f(b)(x)$$

является следствием равенств (4.1.1), (4.1.2).

ТЕОРЕМА 4.1.1. *Представление*

$$f : D \longrightarrow V$$

*кольца $D$ в абелевой группе $V$ удовлетворяет равенству*

$$(4.1.4) \qquad f(0) = \overline{0}$$

*где*

$$\overline{0} : V \to V$$

*отображение такое, что*

$$(4.1.5) \qquad \overline{0} \circ v = 0$$

ДОКАЗАТЕЛЬСТВО. Равенство

$$(4.1.6) \qquad f(a)(x) = f(a + 0)(x) = f(a)(x) + f(0)(x)$$

является следствием равенства (4.1.3). Равенство

$$(4.1.7) \qquad f(0)(x) = 0$$

является следствием равенства (4.1.6). Равенство (4.1.5) является следствием равенств (4.1.4), (4.1.7). $\qquad \square$

ТЕОРЕМА 4.1.2. *Представление*

$$f : D \longrightarrow V$$

*кольца $D$ в абелевой группе $V$ **эффективно** тогда и только тогда, когда из равенства $f(a) = 0$ следует $a = 0$.*

ДОКАЗАТЕЛЬСТВО. Если $a, b \in D$ порождают одно и то же преобразование, то

$$(4.1.8) \qquad f(a)(m) = f(b)(m)$$

для любого $m \in V$. Из равенств (4.1.3), (4.1.8) следует, что

$$(4.1.9) \qquad\qquad f(a-b)(m) = 0$$

Равенство

$$f(a-b) = \overline{0}$$

является следствием равенств (4.1.5), (4.1.9). Следовательно, представление $f$ эффективно тогда и только тогда, когда $a = b$.  $\square$

---

ОПРЕДЕЛЕНИЕ 4.1.3. *Пусть $D$ - коммутативное кольцо. Эффективное представление*

$$(4.1.10) \qquad\qquad f : D \dashrightarrow V \qquad f(d) : v \to dv$$

*кольца $D$ в абелевой группе $V$ называется* **модулем** *над кольцом $D$. Мы также будем говорить, что абелева группа $V$ является $D$-модулем. $V$-число называется* **вектором**. *Билинейное отображение*

$$(4.1.11) \qquad\qquad (a, v) \in D \times V \to av \in V$$

*порождённое представлением*

$$(4.1.12) \qquad\qquad (a, v) \to av$$

*называется произведением вектора на скаляр.*  $\square$

---

ОПРЕДЕЛЕНИЕ 4.1.4. *Пусть $D$ - поле. Эффективное представление*

$$(4.1.13) \qquad\qquad f : D \dashrightarrow V \qquad f(d) : v \to dv$$

*поля $D$ в абелевой группе $V$ называется* **векторным пространством** *над полем $D$. Мы также будем говорить, что абелева группа $V$ является $D$-векторным пространством. $V$-число называется* **вектором**. *Билинейное отображение*

$$(4.1.14) \qquad\qquad (a, v) \in D \times V \to av \in V$$

*порождённое представлением*

$$(4.1.15) \qquad\qquad (a, v) \to av$$

*называется произведением вектора на скаляр.*  $\square$

---

ТЕОРЕМА 4.1.5. *Следующая диаграмма представлений описывает $D$-модуль $V$*

$$(4.1.16) \qquad\qquad \begin{array}{c} D \xrightarrow{\;g_2\;} V \\ \uparrow{\scriptstyle g_1} \\ Z \end{array}$$

*В диаграмме представлений (4.1.16) верна* **коммутативность представлений** *кольца целых чисел $Z$ и коммутативного кольца $D$ в абелевой группе $V$*

$$(4.1.17) \qquad\qquad a(nv) = n(av)$$

Доказательство. Диаграмма представлений (4.1.16) является следствием определения 4.1.3 и теоремы 3.4.2. Равенство (4.1.17) является следствием утверждения, что преобразование $g_2(a)$ является эндоморфизмом $Z$-модуля $V$. □

Определение 4.1.6. *Подпредставление $D$-модуля $V$ называется* **подмодулем** *$D$-модуля $V$.* □

Теорема 4.1.7. *Пусть $V$ - $D$-модуль. Для любого $D_{(1)}$-числа определён эндоморфизм абелевой группы $V$*

$$a \oplus n : v \in V \to (a \oplus n)v \in V \tag{4.1.18}$$

$$(a \oplus n)v = av + nv \quad a \oplus n \in D_{(1)} \tag{4.1.19}$$

*Множество преобразований (3.8.3) порождает представление кольца $D_{(1)}$ в абелевой группе $V$. Мы будем пользоваться обозначением $D_{(1)}v$ для множества векторов, порождённых вектором $v$.*

Доказательство. Согласно теореме 4.1.5, выражения $av$, $nv$ являются $V$-числами. Поскольку $V$ является абелевой группой, выражение (4.1.19) также является $V$-числом. □

Теорема 4.1.8. *Пусть $V$ является $D$-модулем. $V$-числа удовлетворяют соотношениям*

4.1.8.1: **закон коммутативности**

$$v + w = w + v \tag{4.1.20}$$

4.1.8.2: **закон ассоциативности**

$$(pq)v = p(qv) \tag{4.1.21}$$

4.1.8.3: **закон дистрибутивности**

$$p(v + w) = pv + pw \tag{4.1.22}$$

$$(p + q)v = pv + qv \tag{4.1.23}$$

4.1.8.4: **закон унитарности**

$$(0 \oplus 1)v = v \tag{4.1.24}$$

*Если кольцо имеет единицу, то равенство (4.1.24) имеет вид*

$$1v = v \tag{4.1.25}$$

*для любых $p$, $q \in D_{(1)}$, $v$, $w \in V$.*

Доказательство. Равенство (4.1.20) является следствием определения 4.1.3.

Равенство (4.1.22) является следствием теорем 3.3.2, 4.1.7.

Равенство (4.1.23) является следствием теоремы 3.3.4.

Равенство (4.1.21) является следствием теоремы 2.1.10.

Согласно теоремам 4.1.5, 4.1.7, равенство (4.1.24) является следствием ра-

венства

$$(0 \oplus 1)v = 0v + 1v = v$$

и равенства (3.4.4). Если кольцо имеет единицу, то равенство (4.1.25) является следствием теоремы 3.4.12.                                                                    □

Теорема 4.1.9. *Пусть* $V$ - $D$-*модуль. Для любого множества $V$-чисел*

$$(4.1.26) \qquad v = (v(i) \in V, i \in I)$$

*вектор, порождённый диаграммой представлений* (4.1.16), *имеет следующий вид*[4.1]

$$(4.1.27) \qquad J(v) = \big\{ w : w = c(i)v(i), c(i) \in D_{(1)}, |\{i : c(i) \neq 0\}| < \infty \big\}$$

*где кольцо* $D_{(1)}$ *является унитальным расширением кольца* $D$.

Доказательство. Мы докажем теорему по индукции, опираясь на теорему 2.7.5.

Для произвольного $v(k) \in J(v)$ положим

$$c(i) = \begin{cases} 0 \oplus 1, & i = k \\ 0 \end{cases}$$

Тогда равенство

$$(4.1.28) \qquad v(k) = \sum_{i \in I} c(i)v(i) \quad c(i) \in D_{(1)}$$

является следствием равенства (4.1.19). Из равенств (4.1.27), (4.1.28) следует, что теорема верна на множестве $X_0 = v \subseteq J(v)$.

Пусть $X_k \subseteq J(v)$. Согласно определению 4.1.3 и теоремам 2.7.5, 4.1.7, если $w \in X_k$, то либо $w = w_1 + w_2$, $w_1$, $w_2 \in X_{k-1}$, либо $w = aw_1$, $a \in D_{(1)}$, $w_1 \in X_{k-1}$.

Лемма 4.1.10. *Пусть* $w = w_1 + w_2$, $w_1$, $w_2 \in X_{k-1}$. *Тогда* $w \in J(v)$.

Доказательство. Согласно равенству (4.1.27), существуют $D_{(1)}$-числа $w_1(i)$, $w_2(i)$, $i \in I$, такие, что

$$(4.1.29) \qquad w_1 = \sum_{i \in I} w_1(i)v(i) \quad w_2 = \sum_{i \in I} w_2(i)v(i)$$

где множества

$$(4.1.30) \qquad H_1 = \{i \in I : w_1(i) \neq 0\} \quad H_2 = \{i \in I : w_2(i) \neq 0\}$$

конечны. Так как $V$ - $D$-модуль, то из равенств (4.1.23), (4.1.29) следует, что

$$(4.1.31) \quad \begin{aligned} w_1 + w_2 &= \sum_{i \in I} w_1(i)v(i) + \sum_{i \in I} w_2(i)v(i) = \sum_{i \in I}(w_1(i)v(i) + w_2(i)v(i)) \\ &= \sum_{i \in I}(w_1(i) + w_2(i))v(i) \end{aligned}$$

---

[4.1] Для множества $A$, мы обозначим $|A|$ мощность множества $A$. Запись $|A| < \infty$ означает, что множество $A$ конечно.

Из равенства (4.1.30) следует, что множество

$$\{i \in I : w_1(i) + w_2(i) \neq 0\} \subseteq H_1 \cup H_2$$

конечно. ⊙

**Лемма 4.1.11.** *Пусть* $w = aw_1, a \in D_{(1)}, w_1 \in X_{k-1}$. *Тогда* $w \in J(v)$.

Доказательство. Согласно утверждению 2.7.5.4, для любого $D_{(1)}$-числа $a$,

$$(4.1.32) \qquad\qquad aw_1 \in X_k$$

Согласно равенству (4.1.27), существуют $D_{(1)}$-числа $w(i)$, $i \in I$, такие, что

$$(4.1.33) \qquad\qquad w_1 = \sum_{i \in I} w_1(i) v(i)$$

где

$$(4.1.34) \qquad\qquad |\{i \in I : w_1(i) \neq 0\}| < \infty$$

Из равенства (4.1.33) следует, что

$$(4.1.35) \qquad aw_1 = a \sum_{i \in I} w_1(i) v(i) = \sum_{i \in I} a(w_1(i) v(i)) = \sum_{i \in I} (aw_1(i)) v(i)$$

Из утверждения (4.1.34) следует, что множество $\{i \in I : aw(i) \neq 0\}$ конечно. ⊙

Из лемм 4.1.10, 4.1.11 следует, что $X_k \subseteq J(v)$. □

**Соглашение 4.1.12.** *Мы будем пользоваться соглашением о сумме, в котором повторяющийся индекс в линейной комбинации подразумевает сумму по повторяющемуся индексу. В этом случае предполагается известным множество индекса суммирования и знак суммы опускается*

$$c(i)v(i) = \sum_{i \in I} c(i)v(i)$$

*Я буду явно указывать множество индексов, если это необходимо.* □

**Определение 4.1.13.** *Пусть $V$ - модуль . Пусть*

$$v = (v(i) \in V, i \in I)$$

*- множество векторов. Выражение $c(i)v(i)$, $c(i) \in D_{(1)}$, называется **линейной комбинацией** векторов $v(i)$. Вектор $w = c(i)v(i)$ называется **линейно зависимым** от векторов $v(i)$.* □

**Теорема 4.1.14.** *Пусть $v = (v_i \in V, i \in I)$ - множество векторов $D$-модуля $V$. Если векторы $v_i$, $i \in I$, принадлежат подмодулю $V'$ $D$-модуля $V$, то линейная комбинация векторов $v_i$, $i \in I$, принадлежит подмодулю $V'$.*

Доказательство. Теорема является следствием теоремы 4.1.9 и определений 4.1.13, 4.1.6.                                                                              □

Теорема 4.1.15. *Пусть $D$ - поле. Если уравнение*

$$(4.1.36) \qquad\qquad w(i)v(i) = 0$$

*предполагает существования индекса $i = j$ такого, что $w(j) \neq 0$, то вектор $v(j)$ линейно зависит от остальных векторов $v$.*

Доказательство. Теорема является следствием равенства

$$v(j) = \sum_{i \in I \setminus \{j\}} \frac{w(i)}{w(j)} v(i)$$

и определения 4.1.13.                                                                              □

Очевидно, что для любого множества векторов $v(i)$

$$w(i) = 0 \Rightarrow w(i)v(i) = 0$$

Определение 4.1.16. *Множество векторов[4.2] $v(i)$, $i \in I$, $D$-модуля $V$* **линейно независимо**, *если $c(i) = 0$, $i \in I$, следует из уравнения*

$$(4.1.37) \qquad\qquad c(i)v(i) = 0$$

*В противном случае, множество векторов $v(i)$, $i \in I$,* **линейно зависимо**.
                                                                              □

Следующее определение является следствием теорем 2.7.5, 4.1.9 и определения 2.7.4.

Определение 4.1.17. *$J(v)$ называется* **подмодулем**, *порождённым множеством $v$, а $v$ -* **множеством образующих** *подмодуля $J(v)$. В частности,* **множеством образующих** *$D$-модуля $V$ будет такое подмножество $X \subset V$, что $J(X) = V$.*                                                                              □

Следующее определение является следствием теорем 2.7.5, 4.1.9 и определения 2.7.14.

Определение 4.1.18. *Если множество $X \subset V$ является множеством образующих $D$-модуля $V$, то любое множество $Y$, $X \subset Y \subset V$ также является множеством образующих $D$-модуля $V$. Если существует минимальное множество $X$, порождающее $D$-модуль $V$, то такое множество $X$ называется* **квазибазисом** *$D$-модуля $V$.*                                                                              □

Определение 4.1.19. *Пусть $\overline{\overline{e}}$ - квазибазис $D$-модуля $V$, и вектор $\overline{v} \in V$ имеет разложение*

$$(4.1.38) \qquad\qquad \overline{v} = v(i)e(i)$$

---

[4.2] Я следую определению на страница [1]-100.

*относительно квазибазиса $\overline{\overline{e}}$. $D_{(1)}$-числа $v(i)$ называются* **координатами** *вектора $\overline{v}$ относительно квазибазиса $\overline{\overline{e}}$. Матрица $D_{(1)}$-чисел $v = (v(i), i \in I)$ называется* **координатной матрицей вектора** $\overline{v}$ *в квазибазисе $\overline{\overline{e}}$.* $\square$

Теорема 4.1.20. *Множество векторов $\overline{\overline{e}} = (e(i), i \in I)$ является квазибазисом $D$-модуля $V$, если верны следующие утверждения.*

4.1.20.1: *Произвольный вектор $v \in V$ является линейной комбинацией векторов множества $\overline{\overline{e}}$.*

4.1.20.2: *Вектор $e(i)$ нельзя представить в виде линейной комбинации остальных векторов множества $\overline{\overline{e}}$.*

Доказательство. Согласно утверждению 4.1.20.1, теореме 4.1.9 и определению 4.1.13, множество $\overline{\overline{e}}$ порождает $D$-модуль $V$ (определение 4.1.17). Согласно утверждению 4.1.20.2, множество $\overline{\overline{e}}$ является минимальным множеством, порождающим $D$-модуль $V$. Согласно определению 4.1.17, множество $\overline{\overline{e}}$ является квазибазисом $D$-модуля $V$. $\square$

Теорема 4.1.21. *Пусть $D$ - ассоциативная $D$-алгебра. Пусть $\overline{\overline{e}}$ - квазибазис $D$-модуля $V$. Пусть*

$$(4.1.39) \qquad\qquad c(i)e(i) = 0$$

*линейная зависимость векторов квазибазиса $\overline{\overline{e}}$. Тогда*

4.1.21.1: *$D_{(1)}$-число $c(i)$, $i \in I$, не имеет обратного элемента в кольце $D_{(1)}$.*

4.1.21.2: *Множество $V'$ матриц $c = (c(i), i \in I)$ порождает $D$-модуль $V'$.*

Доказательство. Допустим существует матрица $c = (c(i), i \in I)$ такая, что равенство (4.1.39) верно и существует индекс $i = j$ такой, что $c(j) \neq 0$. Если мы положим, что $D$-число $c(j)$ имеет обратный, то равенство

$$e(j) = \sum_{i \in I \setminus \{j\}} \frac{c(i)}{c(j)} e(i)$$

является следствием равенства (4.1.39). Следовательно вектор $e(j)$ является линейной комбинацией остальных векторов множества $\overline{\overline{e}}$ и множество $\overline{\overline{e}}$ не является квазибазисом. Следовательно, наше предположение неверно, и $D_{(1)}$-число $c(j)$ не имеет обратного.

Пусть матрицы $b = (b(i), i \in I) \in D'$, $c = (c(i), i \in I) \in D'$. Из равенств

$$b(i)e(i) = 0$$
$$c(i)e(i) = 0$$

следует

$$(b(i) + c(i))e(i) = 0$$

Следовательно, множество $D'$ является абелевой группой.

Пусть матрица $c = (c(i), i \in I) \in D'$ и $a \in D$. Из равенства

$$c(i)e(i) = 0$$

следует

$$(ac(i))e(i) = 0$$

Следовательно, абелевая группа $D'$ является $D$-модулем.                □

ТЕОРЕМА 4.1.22.  *Пусть  $D$-модуль $V$ имеет квазибазис $\overline{\overline{e}}$ такой, что в равенстве*

$$(4.1.40) \qquad\qquad\qquad c(i)e(i) = 0$$

*существует индекс  $i = j$  такой, что  $c(j) \neq 0$.  Тогда*

4.1.22.1:  *Матрица  $c = (c(i), i \in I)$   определяет координаты вектора $0 \in V$ относительно квазибазиса $\overline{\overline{e}}$.*

4.1.22.2:  *Координаты вектора $\overline{v}$ относительно квазибазиса $\overline{\overline{e}}$ определены однозначно с точностью до выбора координат вектора $0 \in V$.*

ДОКАЗАТЕЛЬСТВО.  Утверждение 4.1.22.1 является следствием равенства (4.1.40) и определения 4.1.19.

Пусть вектор $\overline{v}$ имеет разложение

$$(4.1.41) \qquad\qquad\qquad \overline{v} = v(i)e(i)$$

относительно квазибазиса $\overline{\overline{e}}$. Равенство

$$(4.1.42) \qquad \overline{v} = \overline{v} + 0 = v(i)e(i) + c(i)e(i) = (v(i) + c(i))e(i)$$

является следствием равенств (4.1.40), (4.1.41). Утверждение 4.1.22.2 является следствием равенств (4.1.41), (4.1.42) и определения 4.1.19.                □

ТЕОРЕМА 4.1.23.  *Пусть множество векторов квазибазиса $\overline{\overline{e}}$ $D$-модуля $V$ линейно независимо. Тогда квазибазис $\overline{\overline{e}}$ является базисом $D$-модуля $V$.*

ДОКАЗАТЕЛЬСТВО.  Теорема является следствием определения 2.7.14 и теоремы 4.1.22.                □

ОПРЕДЕЛЕНИЕ 4.1.24.  *$D$-модуль $V$ - **свободный  $D$-модуль**,*[4.3]*если  $D$-модуль $V$ имеет базис.*                □

ТЕОРЕМА 4.1.25.  *$D$-векторное пространство является свободным $D$-модулем.*

ДОКАЗАТЕЛЬСТВО.  Пусть множество векторов $e(i)$, $i \in I$,  линейно зависимо. Тогда в равенстве

$$w(i)e(i) = 0$$

существует индекс  $i = j$  такой, что  $w(j) \neq 0$. Согласно теореме 4.1.15, вектор $e(j)$ линейно зависит от остальных векторов множества $\overline{\overline{e}}$. Согласно определению 4.1.18, множество векторов $e(i)$, $i \in I$,  не является базисом $D$-векторного пространства $V$.

Следовательно, если множество векторов $e(i)$, $i \in I$,  является базисом, то эти векторы линейно независимы. Так как произвольный вектор $v \in V$ является линейной комбинацией векторов $e(i)$, $i \in I$,  то множество векторов $v, e(i)$, $i \in I$,  не является линейно независимым.                □

---

[4.3]   Я следую определению в [1], страница 103.

Теорема 4.1.26. *Координаты вектора* $v \in V$ *относительно базиса* $\overline{\overline{e}}$ $D$-*модуля* $V$ *определены однозначно. Из равенства*

$$(4.1.43) \qquad ve = we \qquad v(i)e(i) = w(i)e(i)$$

*следует, что*

$$v = w \qquad v(i) = w(i)$$

Доказательство. Теорема является следствием теоремы 4.1.22 и определений 4.1.16, 4.1.24. $\qquad \square$

Пример 4.1.27. *Из теорем 3.4.2, 3.6.10 и определения 4.1.3 следует, что свободная абелевая группа* $G$ *является модулем над кольцом целых чисел* $Z$. $\qquad \square$

## 4.2. Тип $D$-модуля

**4.2.1. Концепция обобщённого индекса.** Единственное преимущество выражения $c(i)v(i)$, $c(i) \in D_{(1)}$, линейной комбинации векторов

$$v = (v(i) \in V, i \in I)$$

- это его универсальность. Однако для того чтобы нагляднее видеть операции над векторами, удобно группировать координаты вектора в таблицу, называемую матрицей. Поэтому вместо аргумента, записанного в скобках, мы будем пользоваться индексами, которые говорят как координаты вектора организованы в матрице.

Организацию координат вектора в матрице мы будем называть типом $D$-модуля.

В этом разделе мы рассмотрим вектор-столбец и вектор-строку. Очевидно, что существуют и другие формы представления вектора. Например, мы можем записать координаты вектора в виде $n \times m$ матрицы или в виде треугольной матриц. Формат представления зависит от рассматриваемой задачи.

Если утверждение зависит от формата представления вектора, мы будем указывать либо тип $D$-модуля, либо формат представления базиса. Оба способа указания типа $D$-модуля эквивалентны.

Если рассматривать выражение $i$ в записи $v(i)$ координат вектора как обобщённый индекс, то мы увидим, что поведение вектора не зависит от типа $D$-модуля. Мы можем расширить определение обобщённого индекса и рассматривать индекс $i$ в выражении $a^i$ как обобщённый индекс. Я пользуюсь символом $\cdot$ перед обобщённым индексом, когда мне необходимо описать его структуру. Я помещаю символ $'-'$ на месте индекса, позиция которого изменилась. Например, если исходное выражение было $a_{ij}$, я пользуюсь записью $a_{i\cdot\_}{}^{j}$ вместо записи $a_i^j$. Это позволяет записать равенство вида

$$(4.2.1) \qquad b'^i = b'^{-}{}_{k} = f^{\cdot}{}_{k}{}^{-}{}_{\_}{}^{l} b^{-}{}_{l} = f^i_j b^j$$

Равенство (4.2.1) предполагает равенство индексов ${}^i = {}^{\cdot}{}^{-}{}_{k}$.

Хотя структура обобщённого индекса произвольна, мы будем предполагать, что существует взаимно однозначное отображение отрезка натуральных чисел $1$, ..., $n$ на множество значений индекса. Пусть $I$ - множество значений индекса $i$. Мы будем обозначать мощность этого множества символом $|I|$ и будем полагать $|I| = n$. Если нам надо перечислить элементы $a_i$, мы будем пользоваться обозначением $a_1$, ..., $a_n$.

В геометрии существует традиция записывать координаты вектора с индексом вверху

$$v(i) = v^i$$

и перечислять векторы базиса, используя индекс внизу

$$e(i) = e_i$$

Мы также будем пользоваться соглашением 4.2.1.

Соглашение 4.2.1. *Мы будем пользоваться соглашением Эйнштейна о сумме, в котором повторяющийся индекс (один вверху и один внизу) подразумевает сумму по повторяющемуся индексу. В этом случае предполагается известным множество индекса суммирования и знак суммы опускается*

$$c^i v_i = \sum_{i \in I} c^i v_i$$

*Я буду явно указывать множество индексов, если это необходимо.*                □

Если мы запишем систему линейных уравнений

$$a_1^1 x^1 + ... + a_n^1 x^n = 0$$

(4.2.2)                                     ......

$$a_1^n x^1 + ... + a_n^n x^n = 0$$

пользуясь матрицами

(4.2.3)           $$\begin{pmatrix} a_1^1 & ... & a_n^1 \\ ... & ... & ... \\ a_1^n & ... & a_n^n \end{pmatrix} \begin{pmatrix} a^1 \\ ... \\ a^n \end{pmatrix} = \begin{pmatrix} 0 \\ ... \\ 0 \end{pmatrix}$$

то мы можем отождествить решение системы линейных уравнений (4.2.2) и столбец матрицы. В тоже время, множество решений системы линейных уравнений (4.2.2) является $D$-модулем.

**4.2.2. Обозначение элементов матрицы.** Представление координат вектора в форме матрицы позволяет сделать запись более компактной. Вопрос о представлении вектора как строка или столбец матрицы является вопросом соглашения. Мы можем распространить концепцию обобщённого индекса на элементы матрицы. Матрица - это двумерная таблица, строки и столбцы которой занумерованы обобщёнными индексами. Так как мы пользуемся обобщёнными индексами, мы не можем сказать, нумерует ли индекс $a$ строки матрицы или столбцы, до тех пор, пока мы не знаем структуры индекса. Однако как мы увидим ниже для нас несущественна форма представления матрицы. Для того, чтобы обозначения, предлагаемые ниже, были согласованы с традиционными, мы будем предполагать, что матрица представлена в виде

$$a = \begin{pmatrix} a_1^1 & ... & a_n^1 \\ ... & ... & ... \\ a_1^m & ... & a_n^m \end{pmatrix}$$

Верхний индекс перечисляет строки, нижний индекс перечисляет столбцы.

Определение 4.2.2. *Я использую следующие имена и обозначения различных* **подматриц** *матрицы* $a$

- **столбец матрицы** *с индексом* $i$

$$a_i = \begin{pmatrix} a_i^1 \\ ... \\ a_i^m \end{pmatrix}$$

*Верхний индекс перечисляет элементы столбца, нижний индекс перечисляет столбцы.*

$a_T$ : *подматрица, полученная из матрицы* $a$ *выбором столбцов с индексом из множества* $T$

$a_{[i]}$ : *подматрица, полученный из матрицы* $a$ *удалением столбца* $a_i$

$a_{[T]}$ : *подматрица, полученный из матрицы* $a$ *удалением столбцов с индексом из множества* $T$

- **строка матрицы** *с индексом* $j$

$$a^j = \begin{pmatrix} a_1^j & ... & a_n^j \end{pmatrix}$$

*Нижний индекс перечисляет элементы строки, верхний индекс перечисляет строки.*

$a^S$ : *подматрица, полученный из матрицы* $a$ *выбором строк с индексом из множества* $S$

$a^{[j]}$ : *подматрица, полученный из матрицы* $a$ *удалением строки* $a^j$

$a^{[S]}$ : *подматрица, полученный из матрицы* $a$ *удалением строк с индексом из множества* $S$

□

Замечание 4.2.3. *Мы будем комбинировать запись индексов. Так* $a_i^j$ *является* $1 \times 1$ *подматрицей. Одновременно, это обозначение элемента матрицы. Это позволяет отождествить* $1 \times 1$ *матрицу и её элемент. Индекс* $a$ *является номером столбца матрицы и индекс* $b$ *является номером строки матрицы.* □

Пусть $I$, $|I| = n$, - множество индексов. **Символ Кронекера** определён равенством

$$\delta_j^i = \begin{cases} 1 & i = j \\ 0 & i \neq j \end{cases} \quad i, j \in I$$

Пусть $k$, $l$, $k \le l$, - целые числа. Мы будем рассматривать множество

$$[k, l] = \{i : k \le i \le l\}$$

### 4.2.3. $D$-модуль столбцов.

ОПРЕДЕЛЕНИЕ 4.2.4. *Мы представили множество векторов* $v(1) = v_1$, *...,* $v(m) = v_m$ *в виде строки матрицы*

$$(4.2.4) \qquad v = \begin{pmatrix} v_1 & ... & v_m \end{pmatrix}$$

*и множество* $D_{(1)}$*-чисел* $c(1) = c^1$, *...,* $c(m) = c^m$ *в виде столбца матрицы*

$$(4.2.5) \qquad c = \begin{pmatrix} c^1 \\ ... \\ c^m \end{pmatrix}$$

*Соответствующее представление* $D$*-модуля* $V$ *мы будем называть* $D$**-модулем столбцов**, *а* $V$*-число мы будем называть* **вектор-столбец** .  $\square$

ТЕОРЕМА 4.2.5. *Мы можем записать линейную комбинацию*

$$\sum_{i=1}^{m} c(i)v(i) = c^i v_i$$

*векторов* $v_1$, *...,* $v_m$ *в виде произведения матриц*

$$(4.2.6) \qquad cv = \begin{pmatrix} c^1 \\ ... \\ c^m \end{pmatrix} \begin{pmatrix} v_1 & ... & v_m \end{pmatrix} = c^i v_i$$

Доказательство. Теорема является следствием определений 4.1.13, 4.2.4.  $\square$

ТЕОРЕМА 4.2.6. *Если мы запишем векторы базиса* $\overline{\overline{e}}$ *в виде строки матрицы*

$$(4.2.7) \qquad e = \begin{pmatrix} e_1 & ... & e_n \end{pmatrix}$$

*и координаты вектора* $\overline{w} = w^i e_i$ *относительно базиса* $\overline{\overline{e}}$ *в виде столбца матрицы*

$$(4.2.8) \qquad w = \begin{pmatrix} w^1 \\ ... \\ w^n \end{pmatrix}$$

*то мы можем представить вектор* $\overline{w}$ *в виде произведения матриц*

$$(4.2.9) \qquad \overline{w} = we = \begin{pmatrix} w^1 \\ ... \\ w^n \end{pmatrix} \begin{pmatrix} e_1 & ... & e_n \end{pmatrix} = w^i e_i$$

Доказательство. Теорема является следствием теоремы 4.2.5. □

Теорема 4.2.7. *Координаты вектора $v \in V$ относительно базиса $\overline{\overline{e}}$ $D$-модуля $V$ определены однозначно. Из равенства*

$$(4.2.10) \qquad ve = we \qquad v^i e_i = w^i e_i$$

*следует, что*

$$v = w \qquad v^i = w^i$$

Доказательство. Теорема является следствием теоремы 4.1.26. □

Теорема 4.2.8. *Пусть множество векторов $\overline{\overline{e}}_A = (e_{A\,k}, k \in K)$ является базисом $D$-алгебры столбцов $A$. Пусть $\overline{\overline{e}}_V = (e_{V\,i}, i \in I)$ базис $A$-модуля столбцов $V$. Тогда множество $V$-чисел*

$$(4.2.11) \qquad \overline{\overline{e}}_{VA} = (e_{VA\,ki} = e_{A\,k} e_{V\,i})$$

*является базисом $D$-модуля $V$. Базис $\overline{\overline{e}}_{VA}$ называется расширением базиса $(\overline{\overline{e}}_A, \overline{\overline{e}}_V)$.*

Лемма 4.2.9. *Множество $V$-чисел $e_{2ik}$ порождает $D$-модуль $V$.*

Доказательство. Пусть

$$(4.2.12) \qquad a = a^i e_{1i}$$

произвольное $V$-число. Для любого $A$-числа $a^i$, согласно определению 5.1.1 и соглашению 5.2.1, существуют $D$-числа $a^{ik}$ такие, что

$$(4.2.13) \qquad a^i = a^{ik} e_k$$

Равенство

$$(4.2.14) \qquad a = a^{ik} e_k e_{1i} = a^{ik} e_{2ik}$$

является следствием равенств (4.2.11), (4.2.12), (4.2.13). Согласно теореме 4.1.9, лемма является следствием равенства (4.2.14). ⊙

Лемма 4.2.10. *Множество $V$-чисел $e_{2ik}$ линейно независимо над кольцом $D$.*

Доказательство. Пусть

$$(4.2.15) \qquad a^{ik} e_{2ik} = 0$$

Равенство

$$(4.2.16) \qquad a^{ik} e_k e_{1i} = 0$$

является следствием равенств (4.2.11), (4.2.15). Согласно теореме 4.2.7, равенство

$$(4.2.17) \qquad a^{ik} e_k = 0 \qquad k \in K$$

является следствием равенства (4.2.16). Согласно теореме 4.2.7, равенство

$$(4.2.18) \qquad a^{ik} = 0 \qquad k \in K \qquad i \in I$$

является следствием равенства (4.2.17). Следовательно, множество $V$-чисел $e_{2ik}$ линейно независимо над кольцом $D$.                                    ⊙

Доказательство теоремы 4.2.8.    Теорема является следствием теоремы 4.2.7, и лемм 4.2.9, 4.2.10.                                    □

### 4.2.4. $D$-модуль строк.

Определение 4.2.11. *Мы представили множество векторов $v(1) = v^1$, ..., $v(m) = v^m$ в виде столбца матрицы*

$$(4.2.19) \qquad v = \begin{pmatrix} v^1 \\ ... \\ v^m \end{pmatrix}$$

*и множество $D_{(1)}$-чисел $c(1) = c_1$, ..., $c(m) = c_m$ в виде строки матрицы*

$$(4.2.20) \qquad c = \begin{pmatrix} c_1 & ... & c_m \end{pmatrix}$$

*Соответствующее представление $D$-модуля $V$ мы будем называть $D$-**модулем строк**, а $V$-число мы будем называть **вектор-строка** .*                                    □

Теорема 4.2.12. *Мы можем записать линейную комбинацию*

$$\sum_{i=1}^{m} c(i)v(i) = c_i v^i$$

*векторов $v^1$, ..., $v^m$ в виде произведения матриц*

$$(4.2.21) \qquad vc = \begin{pmatrix} v^1 \\ ... \\ v^m \end{pmatrix} \begin{pmatrix} c_1 & ... & c_m \end{pmatrix} = c_i v^i$$

Доказательство. Теорема является следствием определений 4.1.13, 4.2.11.                                    □

Теорема 4.2.13. *Если мы запишем векторы базиса $\overline{\overline{e}}$ в виде столбца матрицы*

$$(4.2.22) \qquad e = \begin{pmatrix} e^1 \\ ... \\ e^n \end{pmatrix}$$

*и координаты вектора $\overline{w} = w_i e^i$ относительно базиса $\overline{\overline{e}}$ в виде строки матрицы*

$$(4.2.23) \qquad w = \begin{pmatrix} w_1 & ... & w_n \end{pmatrix}$$

*то мы можем представить вектор $\overline{w}$ в виде произведения матриц*

$$(4.2.24) \qquad \overline{w} = ew = \begin{pmatrix} e^1 \\ \dots \\ e^n \end{pmatrix} \begin{pmatrix} w_1 & \dots & w_n \end{pmatrix} = w_i e^i$$

Доказательство. Теорема является следствием теоремы 4.2.12. $\quad\square$

Теорема 4.2.14. *Координаты вектора $v \in V$ относительно базиса $\overline{\overline{e}}$ $D$-модуля $V$ определены однозначно. Из равенства*

$$(4.2.25) \qquad ve = we \qquad v_i e^i = w_i e^i$$

*следует, что*

$$v = w \qquad v_i = w_i$$

Доказательство. Теорема является следствием теоремы 4.1.26. $\quad\square$

Теорема 4.2.15. *Пусть множество векторов $\overline{\overline{e}}_A = (e_A^k, k \in K)$ является базисом $D$-алгебры строк $A$. Пусть $\overline{\overline{e}}_V = (e_V^i, i \in I)$ базис $A$-модуля строк $V$. Тогда множество $V$-чисел*

$$(4.2.26) \qquad \overline{\overline{e}}_{VA} = (e_{VA}^{ki} = e_A^k e_V^i)$$

*является базисом $D$-модуля $V$. Базис $\overline{\overline{e}}_{VA}$ называется расширением базиса $(\overline{\overline{e}}_A, \overline{\overline{e}}_V)$.*

Лемма 4.2.16. *Множество $V$-чисел $e_2^{ik}$ порождает $D$-модуль $V$.*

Доказательство. Пусть

$$(4.2.27) \qquad a = a_i e_1^i$$

произвольное $V$-число. Для любого $A$-числа $a_i$, согласно определению 5.1.1 и соглашению 5.2.1, существуют $D$-числа $a_{ik}$ такие, что

$$(4.2.28) \qquad a_i = a_{ik} e^k$$

Равенство

$$(4.2.29) \qquad a = a_{ik} e^k e_1^i = a_{ik} e_2^{ik}$$

является следствием равенств (4.2.26), (4.2.27), (4.2.28). Согласно теореме 4.1.9, лемма является следствием равенства (4.2.29). $\quad\odot$

Лемма 4.2.17. *Множество $V$-чисел $e_2^{ik}$ линейно независимо над кольцом $D$.*

Доказательство. Пусть

$$(4.2.30) \qquad a_{ik} e_2^{ik} = 0$$

Равенство

$$(4.2.31) \qquad a_{ik} e^k e_1^i = 0$$

является следствием равенств (4.2.26), (4.2.30). Согласно теореме 4.2.14, равенство

$$(4.2.32) \qquad a_{ik}e^k = 0 \quad k \in K$$

является следствием равенства (4.2.31). Согласно теореме 4.2.14, равенство

$$(4.2.33) \qquad a_{ik} = 0 \quad k \in K \quad i \in I$$

является следствием равенства (4.2.32). Следовательно, множество $V$-чисел $e_2^{ik}$ линейно независимо над кольцом $D$.                                                        ⊙

Доказательство теоремы 4.2.15.  Теорема является следствием теоремы 4.2.14, и лемм 4.2.16, 4.2.17.                                                        □

## 4.3. Линейное отображение $D$-модуля

### 4.3.1. Общее определение.

Определение 4.3.1. *Морфизм представлений*

$$(4.3.1) \qquad (h : D_1 \to D_2 \quad f : V_1 \to V_2)$$

*$D_1$-модуля $V_1$ в $D_2$-модуль $V_2$ называется гомоморфизмом или* **линейным отображением** *$D_1$-модуля $V_1$ в $D_2$-модуль $V_2$. Обозначим* $\mathcal{L}(D_1 \to D_2; V_1 \to V_2)$ *множество линейных отображений $D_1$-модуля $V_1$ в $D_2$-модуль $V_2$.*                                                        □

Если отображение (4.3.1) является линейным отображением $D_1$-модуля $V_1$ в $D_2$-модуль $V_2$, то, согласно определению 2.5.1, я пользуюсь обозначением

$$f \circ a = f(a)$$

для образа отображения $f$.

Теорема 4.3.2. *Линейное отображение*

$$(4.3.1) \qquad (h : D_1 \to D_2 \quad f : V_1 \to V_2)$$

*$D_1$-модуля $V_1$ в $D_2$-модуль $V_2$ удовлетворяет равенствам* [4.4]

$$(4.3.2) \qquad h(p + q) = h(p) + h(q)$$

$$(4.3.3) \qquad h(pq) = h(p)h(q)$$

$$(4.3.4) \qquad f \circ (a + b) = f \circ a + f \circ b$$

$$(4.3.5) \qquad f \circ (pa) = h(p)(f \circ a)$$

$$p, q \in D_1 \quad v, w \in V_1$$

Доказательство.  Согласно определению 2.3.9, отображение $h$ является гомоморфизмом кольца $D_1$ в кольцо $D_2$. Равенства (4.3.2), (4.3.3). являются следствием этого утверждения.  Из определения 2.3.9, следует, что отображение $f$ является гомоморфизмом абелевой группы $V_1$ в абелеву группу $V_2$

---

[4.4]   В некоторых книгах (например, на странице [1]-94) теорема 4.3.2 рассматривается как определение.

(равенство (4.3.4)). Равенство (4.3.5) является следствием равенства (2.3.6). $\square$

ТЕОРЕМА 4.3.3. *Линейное отображение*[4.5]

(4.3.6)   $(h : D_1 \to D_2 \quad \overline{f} : V_1 \to V_2)$

$D_1$-*модуля столбцов* $A_1$ *в* $D_2$-*модуль столбцов* $A_2$ *имеет представление*

(4.3.7)   $\overline{f} \circ (e_1 a) = e_2 f h(a)$

(4.3.8)   $\overline{f} \circ (a^i e_{1i}) = h(a^i) f_i^k e_{2k}$

(4.3.9)   $b = f h(a)$

*относительно заданных базисов. Здесь*

- $a$ - *координатная матрица* $V_1$-*числа* $\overline{a}$ *относительно базиса* $\overline{\overline{e}}_1$ *($D$)*

(4.3.10)   $\overline{a} = e_1 a$

- $h(a) = (h(a_i), i \in I)$ - *матрица* $D_2$-*чисел.*
- $b$ - *координатная матрица* $V_2$-*числа*

$$\overline{b} = \overline{f} \circ \overline{a}$$

*относительно базиса* $\overline{\overline{e}}_2$

(4.3.11)   $\overline{b} = e_2 b$

- $f$ - *координатная матрица множества* $V_2$-*чисел* $(\overline{f} \circ e_{1i}, i \in I)$ *относительно базиса* $\overline{\overline{e}}_2$.

ТЕОРЕМА 4.3.4. *Линейное отображение*[4.6]

(4.3.12)   $(h : D_1 \to D_2 \quad \overline{f} : V_1 \to V_2)$

$D_1$-*модуля строк* $A_1$ *в* $D_2$-*модуль строк* $A_2$ *имеет представление*

(4.3.13)   $\overline{f} \circ (ae_1) = h(a) f e_2$

(4.3.14)   $\overline{f} \circ (a_i e_1^i) = h(a_i) f_k^i e_2^k$

(4.3.15)   $b = h(a) f$

*относительно заданных базисов. Здесь*

- $a$ - *координатная матрица* $V_1$-*числа* $\overline{a}$ *относительно базиса* $\overline{\overline{e}}_1$ *($D$)*

(4.3.16)   $\overline{a} = a e_1$

- $h(a) = (h(a^i), i \in I)$ - *матрица* $D_2$-*чисел.*
- $b$ - *координатная матрица* $V_2$-*числа*

$$\overline{b} = \overline{f} \circ \overline{a}$$

*относительно базиса* $\overline{\overline{e}}_2$

(4.3.17)   $\overline{b} = b e_2$

- $f$ - *координатная матрица множества* $V_2$-*чисел* $(\overline{f} \circ e_1^i, i \in I)$ *относительно базиса* $\overline{\overline{e}}_2$.

ДОКАЗАТЕЛЬСТВО ТЕОРЕМЫ 4.3.3.   Так как

$$\boxed{(4.3.6) \quad (h : D_1 \to D_2 \quad \overline{f} : V_1 \to V_2)}$$

линейное отображение, то равенство

(4.3.18)   $\overline{b} = \overline{f} \circ \overline{a} = \overline{f} \circ (e_1 a) = (\overline{f} \circ e_1) h(a)$

является следствием равенств

$$\boxed{(4.3.5) \quad f \circ (pa) = h(p)(f \circ a)} \qquad \boxed{(4.3.10) \quad \overline{a} = e_1 a}$$

$V_2$-число $\overline{f} \circ e_{1i}$   имеет разложение

(4.3.19)   $\overline{f} \circ e_{1i} = e_2 f_i = e_{2j} f_i^j$

---

[4.5]   В теоремах 4.3.3, 4.3.5, мы опираемся на следующее соглашение. Пусть $\overline{\overline{e}}_1 = (e_{1i}, i \in I)$ базис $D_1$-модуля столбцов $V_1$. Пусть $\overline{\overline{e}}_2 = (e_{2j}, j \in J)$ базис $D_2$-модуля столбцов $V_2$.

[4.6]   В теоремах 4.3.4, 4.3.6, мы опираемся на следующее соглашение. Пусть $\overline{\overline{e}}_1 = (e_1^i, i \in I)$ базис $D_1$-модуля строк $V_1$. Пусть $\overline{\overline{e}}_2 = (e_2^j, j \in J)$ базис $D_2$-модуля строк $V_2$.

относительно базиса $\overline{\overline{e}}_2$. Комбинируя (4.3.18) и (4.3.19), мы получаем

(4.3.20) $$\overline{b} = e_2 fh(a)$$

(4.3.9) следует из сравнения

$$\boxed{(4.3.11) \quad \overline{b} = e_2 b}$$

и (4.3.20) и теоремы 4.2.7.                                                                    □

Доказательство теоремы 4.3.4.   Так как

$$\boxed{(4.3.12) \quad (h : D_1 \to D_2 \quad \overline{f} : V_1 \to V_2)}$$

линейное отображение, то равенство

(4.3.21) $$\overline{b} = \overline{f} \circ \overline{a} = \overline{f} \circ (ae_1) = h(a)(\overline{f} \circ e_1)$$

является следствием равенств

$$\boxed{(4.3.5) \quad f \circ (pa) = h(p)(f \circ a)} \qquad \boxed{(4.3.16) \quad \overline{a} = ae_1}$$

$V_2$-число $\overline{f} \circ e_1^i$  имеет разложение

(4.3.22) $$\overline{f} \circ e_1 = f_i e_2 = f_j^i e_2^j$$

относительно базиса $\overline{\overline{e}}_2$. Комбинируя (4.3.21) и (4.3.22), мы получаем

(4.3.23) $$\overline{b} = h(a)fe_2$$

(4.3.15) следует из сравнения

$$\boxed{(4.3.17) \quad \overline{b} = be_2}$$

и (4.3.23) и теоремы 4.2.14.                                                                    □

---

Теорема 4.3.5. *Пусть отображение*

$$h : D_1 \to D_2$$

*является гомоморфизмом кольца $D_1$ в кольцо $D_2$.   Пусть*

$$f = (f_j^i, i \in I, j \in J)$$

*матрица $D_2$-чисел.   Отображение*[4.5]

$$\boxed{(4.3.6) \quad (h : D_1 \to D_2 \quad \overline{f} : V_1 \to V_2)}$$

*определённое равенством*

$$\boxed{(4.3.7) \quad \overline{f} \circ (e_1 a) = e_2 fh(a)}$$

*является гомоморфизмом $D_1$-модуля столбцов $V_1$ в $D_2$-модуль столбцов $V_2$.   Гомоморфизм (4.3.6), который имеет данную матрицу $f$, определён однозначно.*

Теорема 4.3.6. *Пусть отображение*

$$h : D_1 \to D_2$$

*является гомоморфизмом кольца $D_1$ в кольцо $D_2$.   Пусть*

$$f = (f_i^j, i \in I, j \in J)$$

*матрица $D_2$-чисел.   Отображение*[4.6]

$$\boxed{(4.3.12) \quad (h : D_1 \to D_2 \quad \overline{f} : V_1 \to V_2)}$$

*определённое равенством*

$$\boxed{(4.3.13) \quad \overline{f} \circ (ae_1) = h(a)fe_2}$$

*является гомоморфизмом $D_1$-модуля строк $V_1$ в $D_2$-модуль строк $V_2$.   Гомоморфизм (4.3.12), который имеет данную матрицу $f$, определён однозначно.*

Доказательство. Равенство

$$(4.3.24) \quad \begin{aligned} &\overline{f} \circ (e_1(v+w)) \\ &= e_2 f h(v+w) \\ &= e_2 f(h(v) + h(w)) \\ &= e_2 f h(v) + e_2 f h(w) \\ &= \overline{f} \circ (e_1 v) + \overline{f} \circ (e_1 w) \end{aligned}$$

является следствием равенств

$$(4.3.4) \quad f \circ (a+b) = f \circ a + f \circ b$$

$$(4.3.5) \quad f \circ (pa) = h(p)(f \circ a)$$

Из равенства (4.3.24) следует, что отображение $\overline{f}$ является гомоморфизмом абелевой группы. Равенство

$$(4.3.25) \quad \begin{aligned} &\overline{f} \circ (e_1(av)) \\ &= e_2 f h(av) \\ &= h(a)(e_2 f h(v)) \\ &= h(a)(\overline{f} \circ (e_1 v)) \end{aligned}$$

является следствием равенства (4.3.5). Из равенства (4.3.25) и определений 2.8.5, 4.3.1 следует, что отображение

$$(4.3.6) \quad (h : D_1 \to D_2 \quad \overline{f} : V_1 \to V_2)$$

является гомоморфизмом $D_1$-модуля столбцов $V_1$ в $D_2$-модуль столбцов $V_2$.

Пусть $f$ - матрица гомоморфизмов $\overline{f}$, $\overline{g}$ относительно базисов $\overline{\overline{e}}_1$, $\overline{\overline{e}}_2$. Равенство

$$\overline{f} \circ (e_1 v) = e_2 f h(v) = \overline{g} \circ (e_1 v)$$

является следствием теоремы 4.3.2. Следовательно, $\overline{f} = \overline{g}$. $\qquad \square$

Доказательство. Равенство

$$(4.3.26) \quad \begin{aligned} &\overline{f} \circ ((v+w)e_1) \\ &= h(v+w)f e_2 \\ &= (h(v) + h(w))f e_2 \\ &= h(v)f e_2 + h(w)f e_2 \\ &= \overline{f} \circ (v e_1) + \overline{f} \circ (w e_1) \end{aligned}$$

является следствием равенств

$$(4.3.4) \quad f \circ (a+b) = f \circ a + f \circ b$$

$$(4.3.5) \quad f \circ (pa) = h(p)(f \circ a)$$

Из равенства (4.3.26) следует, что отображение $\overline{f}$ является гомоморфизмом абелевой группы. Равенство

$$(4.3.27) \quad \begin{aligned} &\overline{f} \circ ((av)e_1) \\ &= h(av)f e_2 \\ &= h(a)(h(v)f e_2) \\ &= h(a)(\overline{f} \circ (v e_1)) \end{aligned}$$

является следствием равенства (4.3.5). Из равенства (4.3.27) и определений 2.8.5, 4.3.1 следует, что отображение

$$(4.3.12) \quad (h : D_1 \to D_2 \quad \overline{f} : V_1 \to V_2)$$

является гомоморфизмом $D_1$-модуля строк $V_1$ в $D_2$-модуль строк $V_2$.

Пусть $f$ - матрица гомоморфизмов $\overline{f}$, $\overline{g}$ относительно базисов $\overline{\overline{e}}_1$, $\overline{\overline{e}}_2$. Равенство

$$\overline{f} \circ (v e_1) = h(v)f e_2 = \overline{g} \circ (v e_1)$$

является следствием теоремы 4.3.2. Следовательно, $\overline{f} = \overline{g}$. $\qquad \square$

На основании теорем 4.3.3, 4.3.5, а также на основании теорем 4.3.4, 4.3.6, мы идентифицируем линейное отображение

$$(h : D_1 \to D_2 \quad \overline{f} : V_1 \to V_2)$$

и матрицу $f$ его представления.

Теорема 4.3.7. *Отображение*

$$(I_{(1)} : d \in D \to d + 0 \in D_{(1)} \quad id : v \in V \to v \in V)$$

*является гомоморфизмом $D$-модуля $V$ в $D_{(1)}$-модуль $V$.*

Доказательство. Равенство

$$(4.3.28) \quad (d_1 + 0) + (d_2 + 0) = (d_1 + d_2) + (0 + 0) = (d_1 + d_2) + 0$$

является следствием равенства (3.8.5). Равенство

$$(4.3.29) \qquad (d_1 + 0)(d_2 + 0) = (d_1 d_2 + 0 d_1 + 0 d_2) + (0\,0) = (d_1 d_2) + 0$$

является следствием равенства (3.8.6). Из равенств (4.3.28), (4.3.29) следует, что отображение $I_{(1)}$ является гомоморфизмом кольца $D$ в кольцо $D_{(1)}$.

Равенство

$$(4.3.30) \qquad (d + 0)v = dv + 0v = dv + 0 = dv$$

является следствием равенства

$$\boxed{(4.1.23) \quad (p + q)v = pv + qv}$$

Теорема является следствием равенства (4.3.30) и определения 2.3.9.      □

Соглашение 4.3.8. *Согласно теоремам 4.1.9, 4.3.7, множество слов, порождающих $D$-модуль $V$, совпадает с множеством слов, порождающих $D_{(1)}$-модуль $V$. Поэтому, не нарушая общности, мы будем предполагать, что кольцо $D$ имеет единицу.*      □

### 4.3.2. Линейное отображение, когда кольца $D_1 = D_2 = D$.

Определение 4.3.9. *Морфизм представлений*

$$(4.3.31) \qquad f : V_1 \to V_2$$

*$D$-модуля $V_1$ в $D$-модуль $V_2$ называется гомоморфизмом или **линейным отображением** $D$-модуля $V_1$ в $D$-модуль $V_2$. Обозначим $\mathcal{L}(D; V_1 \to V_2)$ множество линейных отображений $D$-модуля $V_1$ в $D$-модуль $V_2$.*      □

Если отображение (4.3.31) является линейным отображением $D$-модуля $V_1$ в $D$-модуль $V_2$, то, согласно определению 2.5.1, я пользуюсь обозначением

$$f \circ a = f(a)$$

для образа отображения $f$.

Теорема 4.3.10. *Линейное отображение*

$$\boxed{(4.3.31) \quad f : V_1 \to V_2}$$

*$D$-модуля $V_1$ в $D$-модуль $V_2$ удовлетворяет равенствам*[4.7]

$$(4.3.32) \qquad f \circ (a + b) = f \circ a + f \circ b$$

$$(4.3.33) \qquad f \circ (pa) = p(f \circ a)$$

$$p \in D \quad v, w \in V_1$$

Доказательство.      Из определения 2.3.9, следует, что отображение $f$ является гомоморфизмом абелевой группы $V_1$ в абелеву группу $V_2$ (равенство

---

[4.7]    В некоторых книгах (например, на странице [1]-94) теорема 4.3.10 рассматривается как определение.

(4.3.32)). Равенство (4.3.33) является следствием равенства (2.3.6). □

ТЕОРЕМА 4.3.11. *Линейное отоб-ражение*[4.8]

$$(4.3.34) \qquad \overline{f} : V_1 \to V_2$$

*$D$-модуля столбцов $A_1$ в $D$-модуль столбцов $A_2$ имеет представление*

$$(4.3.35) \qquad \overline{f} \circ (e_1 a) = e_2 f a$$

$$(4.3.36) \qquad \overline{f} \circ (a^i e_{1i}) = a^i f_i^k e_{2k}$$

$$(4.3.37) \qquad b = f a$$

*относительно заданных базисов. Здесь*

- *$a$ - координатная матрица $V_1$-числа $\overline{a}$ относительно базиса $\overline{\overline{e}}_1$ (D)*

$$(4.3.38) \qquad \overline{a} = e_1 a$$

- *$b$ - координатная матрица $V_2$-числа*

$$\overline{b} = \overline{f} \circ \overline{a}$$

*относительно базиса $\overline{\overline{e}}_2$*

$$(4.3.39) \qquad \overline{b} = e_2 b$$

- *$f$ - координатная матрица множества $V_2$-чисел ($\overline{f} \circ e_{1i}, i \in I$) относительно базиса $\overline{\overline{e}}_2$.*

ТЕОРЕМА 4.3.12. *Линейное отоб-ражение*[4.9]

$$(4.3.40) \qquad \overline{f} : V_1 \to V_2$$

*$D$-модуля строк $A_1$ в $D$-модуль строк $A_2$ имеет представление*

$$(4.3.41) \qquad \overline{f} \circ (a e_1) = a f e_2$$

$$(4.3.42) \qquad \overline{f} \circ (a_i e_1^i) = a_i f_k^i e_2^k$$

$$(4.3.43) \qquad b = a f$$

*относительно заданных базисов. Здесь*

- *$a$ - координатная матрица $V_1$-числа $\overline{a}$ относительно ба-зиса $\overline{\overline{e}}_1$ (D)*

$$(4.3.44) \qquad \overline{a} = a e_1$$

- *$b$ - координатная матрица $V_2$-числа*

$$\overline{b} = \overline{f} \circ \overline{a}$$

*относительно базиса $\overline{\overline{e}}_2$*

$$(4.3.45) \qquad \overline{b} = b e_2$$

- *$f$ - координатная матрица множества $V_2$-чисел ($\overline{f} \circ e_1^i, i \in I$) относительно ба-зиса $\overline{\overline{e}}_2$.*

ДОКАЗАТЕЛЬСТВО ТЕОРЕМЫ 4.3.11. Так как

$$\boxed{(4.3.34) \quad \overline{f} : V_1 \to V_2}$$

линейное отображение, то равенство

$$(4.3.46) \qquad \overline{b} = \overline{f} \circ \overline{a} = \overline{f} \circ (e_1 a) = (\overline{f} \circ e_1) a$$

является следствием равенств

$$\boxed{(4.3.33) \quad f \circ (pa) = p(f \circ a)} \qquad \boxed{(4.3.38) \quad \overline{a} = e_1 a}$$

$V_2$-число $\overline{f} \circ e_{1i}$ имеет разложение

$$(4.3.47) \qquad \overline{f} \circ e_{1i} = e_2 f_i = e_{2j} f_i^j$$

относительно базиса $\overline{\overline{e}}_2$. Комбинируя (4.3.46) и (4.3.47), мы получаем

$$(4.3.48) \qquad \overline{b} = e_2 f a$$

---

[4.8] В теоремах 4.3.11, 4.3.13, мы опираемся на следующее соглашение. Пусть $\overline{\overline{e}}_1 = (e_{1i}, i \in I)$ базис $D$-модуля столбцов $V_1$. Пусть $\overline{\overline{e}}_2 = (e_{2j}, j \in J)$ базис $D$-модуля столбцов $V_2$.

[4.9] В теоремах 4.3.12, 4.3.14, мы опираемся на следующее соглашение. Пусть $\overline{\overline{e}}_1 = (e_1^i, i \in I)$ базис $D$-модуля строк $V_1$. Пусть $\overline{\overline{e}}_2 = (e_2^j, j \in J)$ базис $D$-модуля строк $V_2$.

(4.3.37) следует из сравнения

$$\boxed{(4.3.39) \quad \overline{b} = e_2 b}$$

и (4.3.48) и теоремы 4.2.7.                                                    □

Доказательство теоремы 4.3.12.    Так как

$$\boxed{(4.3.40) \quad \overline{f} : V_1 \to V_2}$$

линейное отображение, то равенство

$$(4.3.49) \qquad \overline{b} = \overline{f} \circ \overline{a} = \overline{f} \circ (a e_1) = a(\overline{f} \circ e_1)$$

является следствием равенств

$$\boxed{(4.3.33) \quad f \circ (pa) = p(f \circ a)} \qquad \boxed{(4.3.44) \quad \overline{a} = a e_1}$$

$V_2$-число  $\overline{f} \circ e_1^i$  имеет разложение

$$(4.3.50) \qquad \overline{f} \circ e_1^i = f_i e_2 = f_j^i e_2^j$$

относительно базиса $\overline{\overline{e}}_2$. Комбинируя (4.3.49) и (4.3.50), мы получаем

$$(4.3.51) \qquad \overline{b} = a f e_2$$

(4.3.43) следует из сравнения

$$\boxed{(4.3.45) \quad \overline{b} = b e_2}$$

и (4.3.51) и теоремы 4.2.14.                                                    □

---

**Теорема 4.3.13.**  *Пусть*

$$f = (f_j^i, i \in I, j \in J)$$

*матрица $D$-чисел.    Отображение*[4.8]

$$\boxed{(4.3.34) \quad \overline{f} : V_1 \to V_2} \quad \textit{определённое}$$

*равенством*

$$\boxed{(4.3.35) \quad \overline{f} \circ (e_1 a) = e_2 f a}$$

*является гомоморфизмом $D$-модуля столбцов $V_1$ в $D$-модуль столбцов $V_2$. Гомоморфизм (4.3.34), который имеет данную матрицу $f$, определён однозначно.*

Доказательство. Равенство

$$(4.3.52) \quad
\begin{aligned}
& \overline{f} \circ (e_1(v + w)) \\
&= e_2 f(v + w) \\
&= e_2 f v + e_2 f w \\
&= \overline{f} \circ (e_1 v) + \overline{f} \circ (e_1 w)
\end{aligned}$$

является следствием равенств

$$\boxed{(4.3.32) \quad f \circ (a + b) = f \circ a + f \circ b}$$

$$\boxed{(4.3.33) \quad f \circ (pa) = p(f \circ a)}$$

Из равенства (4.3.52) следует, что отображение $\overline{f}$ является гомоморфизмом абелевой группы. Равенство

---

**Теорема 4.3.14.**  *Пусть*

$$f = (f_i^j, i \in I, j \in J)$$

*матрица $D$-чисел.    Отображение*[4.9]

$$\boxed{(4.3.40) \quad \overline{f} : V_1 \to V_2} \quad \textit{определённое}$$

*равенством*

$$\boxed{(4.3.41) \quad \overline{f} \circ (a e_1) = a f e_2}$$

*является гомоморфизмом $D$-модуля строк $V_1$ в $D$-модуль строк $V_2$. Гомоморфизм (4.3.40), который имеет данную матрицу $f$, определён однозначно.*

Доказательство. Равенство

$$(4.3.54) \quad
\begin{aligned}
& \overline{f} \circ ((v + w)e_1) \\
&= (v + w) f e_2 \\
&= v f e_2 + w f e_2 \\
&= \overline{f} \circ (v e_1) + \overline{f} \circ (w e_1)
\end{aligned}$$

является следствием равенств

$$\boxed{(4.3.32) \quad f \circ (a + b) = f \circ a + f \circ b}$$

$$\boxed{(4.3.33) \quad f \circ (pa) = p(f \circ a)}$$

Из равенства (4.3.54) следует, что отображение $\overline{f}$ является гомоморфизмом абелевой группы. Равенство

$$(4.3.53) \quad \begin{aligned} &\overline{f} \circ (e_1(av)) \\ &= e_2 f(av) \\ &= a(e_2 fv) \\ &= a(\overline{f} \circ (e_1 v)) \end{aligned} \qquad (4.3.55) \quad \begin{aligned} &\overline{f} \circ ((av)e_1) \\ &= (av)fe_2 \\ &= a(vfe_2) \\ &= a(\overline{f} \circ (ve_1)) \end{aligned}$$

является следствием равенства (4.3.33). Из равенства (4.3.53) и определений 2.8.5, 4.3.9 следует, что отображение

$$(4.3.34) \quad \boxed{\overline{f} : V_1 \to V_2}$$

является гомоморфизмом $D$-модуля столбцов $V_1$ в $D$-модуль столбцов $V_2$.

Пусть $f$ - матрица гомоморфизмов $\overline{f}$, $\overline{g}$ относительно базисов $\overline{\overline{e}}_1$, $\overline{\overline{e}}_2$. Равенство

$$\overline{f} \circ (e_1 v) = e_2 fv = \overline{g} \circ (e_1 v)$$

является следствием теоремы 4.3.10. Следовательно, $\overline{f} = \overline{g}$. □

является следствием равенства (4.3.33). Из равенства (4.3.55) и определений 2.8.5, 4.3.9 следует, что отображение

$$(4.3.40) \quad \boxed{\overline{f} : V_1 \to V_2}$$

является гомоморфизмом $D$-модуля строк $V_1$ в $D$-модуль строк $V_2$.

Пусть $f$ - матрица гомоморфизмов $\overline{f}$, $\overline{g}$ относительно базисов $\overline{\overline{e}}_1$, $\overline{\overline{e}}_2$. Равенство

$$\overline{f} \circ (ve_1) = vfe_2 = \overline{g} \circ (ve_1)$$

является следствием теоремы 4.3.10. Следовательно, $\overline{f} = \overline{g}$. □

На основании теорем 4.3.11, 4.3.13, а также на основании теорем 4.3.12, 4.3.14, мы идентифицируем линейное отображение

$$\overline{f} : V_1 \to V_2$$

и матрицу $f$ его представления.

ТЕОРЕМА 4.3.15. *Пусть отображение*

$$f : A_1 \to A_2$$

*является линейным отображением $D$-модуля $A_1$ в $D$-модуль $A_2$. Тогда*

$$f \circ 0 = 0$$

ДОКАЗАТЕЛЬСТВО. Следствие равенства

$$f \circ (a + 0) = f \circ a + f \circ 0$$

□

ОПРЕДЕЛЕНИЕ 4.3.16. *Линейное отображение*

$$f : V \to V$$

*$D$-модуля $V$ называется* **эндоморфизмом** *$D$-модуля $V$. Обозначим* $\mathrm{End}(D, V)$ *множество эндоморфизмов $D$-модуля $V$.* □

## 4.4. Полилинейное отображение $D$-модуля

ОПРЕДЕЛЕНИЕ 4.4.1. *Пусть $D$ - коммутативное кольцо. Приведенный полиморфизм $D$-модулей $A_1$, ..., $A_n$ в $D$-модуль $S$*

$$f : A_1 \times ... \times A_n \to S$$

называется **полилинейным отображением** $D$-модулей $A_1$, ..., $A_n$ в $D$-модуль $S$.

*Обозначим* $\mathcal{L}(D; A_1 \times ... \times A_n \to S)$ *множество полилинейных отображений* $D$-модулей $A_1$, ..., $A_n$ в $D$-модуль $S$. Обозначим $\mathcal{L}(D; A^n \to S)$ *множество n-линейных отображений* $D$-модуля $A$ $(A_1 = ... = A_n = A)$ в $D$-модуль $S$. □

**Теорема 4.4.2.** *Пусть* $D$ - *коммутативное кольцо. Полилинейное отображение* $D$-модулей $A_1$, ..., $A_n$ в $D$-модуль $S$

$$f : A_1 \times ... \times A_n \to S$$

*удовлетворяет равенствам*

$$f \circ (a_1, ..., a_i + b_i, ..., a_n) = f \circ (a_1, ..., a_i, ..., a_n) + f \circ (a_1, ..., b_i, ..., a_n)$$

$$f \circ (a_1, ..., pa_i, ..., a_n) = pf \circ (a_1, ..., a_i, ..., a_n)$$

$$1 \leq i \leq n \quad a_i, b_i \in A_i \quad p \in D$$

Доказательство. Теорема является следствием определений 2.6.1, 4.3.9, 4.4.1 и теоремы 4.3.10. □

**Теорема 4.4.3.** *Пусть* $D$ - *коммутативное кольцо. Пусть* $A_1$, ..., $A_n$, $S$ - *D-модули. Отображение*

(4.4.1) $\qquad f + g : A_1 \times ... \times A_n \to S \quad f, g \in \mathcal{L}(D; A_1 \times ... \times A_n \to S)$

*определённое равенством*

(4.4.2) $\qquad (f + g) \circ (a_1, ..., a_n) = f \circ (a_1, ..., a_n) + g \circ (a_1, ..., a_n)$

*называется* **суммой полилинейных отображений** $f$ *и* $g$ *и является полилинейным отображением. Множество* $\mathcal{L}(D; A_1 \times ... \times A_n \to S)$ *является абелевой группой относительно суммы отображений.*

Доказательство. Согласно теореме 4.4.2

(4.4.3) $\quad f \circ (a_1, ..., a_i + b_i, ..., a_n) = f \circ (a_1, ..., a_i, ..., a_n) + f \circ (a_1, ..., b_i, ..., a_n)$

(4.4.4) $\qquad\qquad f \circ (a_1, ..., pa_i, ..., a_n) = pf \circ (a_1, ..., a_i, ..., a_n)$

(4.4.5) $\quad g \circ (a_1, ..., a_i + b_i, ..., a_n) = g \circ (a_1, ..., a_i, ..., a_n) + g \circ (a_1, ..., b_i, ..., a_n)$

(4.4.6) $\qquad\qquad g \circ (a_1, ..., pa_i, ..., a_n) = pg \circ (a_1, ..., a_i, ..., a_n)$

Равенство

(4.4.7) $\quad\begin{aligned} & (f + g) \circ (x_1, ..., x_i + y_i, ..., x_n) \\ & = f \circ (x_1, ..., x_i + y_i, ..., x_n) + g \circ (x_1, ..., x_i + y_i, ..., x_n) \\ & = f \circ (x_1, ..., x_i, ..., x_n) + f \circ (x_1, ..., y_i, ..., x_n) \\ & + g \circ (x_1, ..., x_i, ..., x_n) + g \circ (x_1, ..., y_i, ..., x_n) \\ & = (f + g) \circ (x_1, ..., x_i, ..., x_n) + (f + g) \circ (x_1, ..., y_i, ..., x_n) \end{aligned}$

является следствием равенств (4.4.2), (4.4.3), (4.4.5). Равенство

$$
\begin{aligned}
(f + g) \circ (x_1, ..., px_i, ..., x_n) \\
= f \circ (x_1, ..., px_i, ..., x_n) + g \circ (x_1, ..., px_i, ..., x_n) \\
(4.4.8) \quad = pf \circ (x_1, ..., x_i, ..., x_n) + pg \circ (x_1, ..., x_i, ..., x_n) \\
= p(f \circ (x_1, ..., x_i, ..., x_n) + g \circ (x_1, ..., x_i, ..., x_n)) \\
= p(f + g) \circ (x_1, ..., x_i, ..., x_n)
\end{aligned}
$$

является следствием равенств (4.4.2), (4.4.4), (4.4.6). Из равенств (4.4.7), (4.4.8) и теоремы 4.4.2 следует, что отображение (4.4.1) является полилинейным отображением $D$-модулей.

Пусть $f$, $g$, $h \in \mathcal{L}(D; A_1 \times ... \times A_2 \to S)$. Для любого $a = (a_1, ..., a_n)$, $a_1 \in A_1$, ..., $a_n \in A_n$,

$$
\begin{aligned}
(f + g) \circ a &= f \circ a + g \circ a = g \circ a + f \circ a \\
&= (g + f) \circ a \\
((f + g) + h) \circ a &= (f + g) \circ a + h \circ a = (f \circ a + g \circ a) + h \circ a \\
&= f \circ a + (g \circ a + h \circ a) = f \circ a + (g + h) \circ a \\
&= (f + (g + h)) \circ a
\end{aligned}
$$

Следовательно, сумма полилинейных отображений коммутативна и ассоциативна.

Из равенства (4.4.2) следует, что отображение

$$
0 : v \in A_1 \times ... \times A_n \to 0 \in S
$$

является нулём операции сложения

$$
(0 + f) \circ (a_1, ..., a_n) = 0 \circ (a_1, ..., a_n) + f \circ (a_1, ..., a_n) = f \circ (a_1, ..., a_n)
$$

Из равенства (4.4.2) следует, что отображение

$$
-f : (a_1, ..., a_n) \in A_1 \times ... \times A_n \to -(f \circ (a_1, ..., a_n)) \in S
$$

является отображением, обратным отображению $f$

$$
f + (-f) = 0
$$

так как

$$
\begin{aligned}
(f + (-f)) \circ (a_1, ..., a_n) &= f \circ (a_1, ..., a_n) + (-f) \circ (a_1, ..., a_n) \\
&= f \circ (a_1, ..., a_n) - f \circ (a_1, ..., a_n) \\
&= 0 = 0 \circ (a_1, ..., a_n)
\end{aligned}
$$

Из равенства

$$
\begin{aligned}
(f + g) \circ (a_1, ..., a_n) &= f \circ (a_1, ..., a_n) + g \circ (a_1, ..., a_n) \\
&= g \circ (a_1, ..., a_n) + f \circ (a_1, ..., a_n) \\
&= (g + f) \circ (a_1, ..., a_n)
\end{aligned}
$$

следует, что сумма отображений коммутативно. Следовательно, множество $\mathcal{L}(D; A_1 \times ... \times A_n \to S)$ является абелевой группой. $\qquad \square$

Следствие 4.4.4. *Пусть $A_1$, $A_2$ - D-модули. Отображение*

$$(4.4.9) \qquad f + g : A_1 \to A_2 \quad f, g \in \mathcal{L}(D; A_1 \to A_2)$$

*определённое равенством*

$$(4.4.10) \qquad (f + g) \circ x = f \circ x + g \circ x$$

*называется* **суммой отображений** *$f$ и $g$ и является линейным отображением. Множество $\mathcal{L}(D; A_1 \to A_2)$ является абелевой группой относительно суммы отображений.* □

Теорема 4.4.5. *Пусть $D$ - коммутативное кольцо. Пусть $A_1$, ..., $A_n$, $S$ - D-модули. Отображение*

$$(4.4.11) \qquad d \, f : A_1 \times ... \times A_n \to S \quad d \in D \quad f \in \mathcal{L}(D; A_1 \times ... \times A_n \to S)$$

*определённое равенством*

$$(4.4.12) \qquad (d \, f) \circ (a_1, ..., a_n) = d(f \circ (a_1, ..., a_n))$$

*называется* **произведением отображения** *$f$ на скаляр $d$ и является полилинейным отображением. Представление*

$$(4.4.13) \qquad a : f \in \mathcal{L}(D; A_1 \times ... \times A_n \to S) \to af \in \mathcal{L}(D; A_1 \times ... \times A_n \to S)$$

*кольца $D$ в абелевой группе $\mathcal{L}(D; A_1 \times ... \times A_n \to S)$ порождает структуру D-модуля.*

Доказательство. Согласно теореме 4.4.2

$$(4.4.14) \quad f \circ (a_1, ..., a_i + b_i, ..., a_n) = f \circ (a_1, ..., a_i, ..., a_n) + f \circ (a_1, ..., b_i, ..., a_n)$$

$$(4.4.15) \qquad f \circ (a_1, ..., pa_i, ..., a_n) = pf \circ (a_1, ..., a_i, ..., a_n)$$

Равенство

$$
\begin{aligned}
(4.4.16) \quad & (pf) \circ (x_1, ..., x_i + y_i, ..., x_n) \\
& = p \, f \circ (x_1, ..., x_i + y_i, ..., x_n) \\
& = p \, (f \circ (x_1, ..., x_i, ..., x_n) + f \circ (x_1, ..., y_i, ..., x_n)) \\
& = p(f \circ (x_1, ..., x_i, ..., x_n)) + p(f \circ (x_1, ..., y_i, ..., x_n)) \\
& = (pf) \circ (x_1, ..., x_i, ..., x_n) + (pf) \circ (x_1, ..., y_i, ..., x_n)
\end{aligned}
$$

является следствием равенств (4.4.12), (4.4.14). Равенство

$$
\begin{aligned}
(4.4.17) \quad & (pf) \circ (x_1, ..., qx_i, ..., x_n) \\
& = p(f \circ (x_1, ..., qx_i, ..., x_n)) = pq(f \circ (x_1, ..., x_i, ..., x_n)) \\
& = qp(f \circ (x_1, ..., x_n)) = q(pf) \circ (x_1, ..., x_n)
\end{aligned}
$$

является следствием равенств (4.4.12), (4.4.15). Из равенств (4.4.16), (4.4.17) и теоремы 4.4.2 следует, что отображение (4.4.11) является полилинейным отображением D-модулей.

Равенство

$$(4.4.18) \qquad (p + q)f = pf + qf$$

является следствием равенства

$$
\begin{aligned}
((p+q)f) \circ (x_1, ..., x_n) &= (p+q)(f \circ (x_1, ..., x_n)) \\
&= p(f \circ (x_1, ..., x_n)) + q(f \circ (x_1, ..., x_n)) \\
&= (pf) \circ (x_1, ..., x_n) + (qf) \circ (x_1, ..., x_n)
\end{aligned}
$$

Равенство

(4.4.19) $$p(qf) = (pq)f$$

является следствием равенства

$$
\begin{aligned}
(p(qf)) \circ (x_1, ..., x_n) &= p \ (qf) \circ (x_1, ..., x_n) = p \ (q \ f \circ (x_1, ..., x_n)) \\
&= (pq) \ f \circ (x_1, ..., x_n) = ((pq)f) \circ (x_1, ..., x_n)
\end{aligned}
$$

Из равенств (4.4.18), (4.4.19), следует, что отображение (4.4.13) является представлением кольца $D$ в абелевой группе $\mathcal{L}(D; A_1 \times ... \times A_n \to S)$. Так как указанное представление эффективно, то, согласно определению 4.1.3 и теореме 4.4.3, абелева группа $\mathcal{L}(D; A_1 \to A_2)$ является $D$-модулем. $\qquad \square$

**СЛЕДСТВИЕ 4.4.6.** *Пусть $A_1$, $A_2$ - $D$-модули. Отображение*

(4.4.20) $$d \, f : A_1 \to A_2 \quad d \in D \quad f \in \mathcal{L}(D; A_1 \to A_2)$$

*определённое равенством*

(4.4.21) $$(d \, f) \circ x = d(f \circ x)$$

*называется **произведением отображения** $f$ **на скаляр** $d$ и является линейным отображением. Представление*

(4.4.22) $$a : f \in \mathcal{L}(D; A_1 \to A_2) \to af \in \mathcal{L}(D; A_1 \to A_2)$$

*кольца $D$ в абелевой группе $\mathcal{L}(D; A_1 \to A_2)$ порождает структуру $D$-модуля.* $\qquad \square$

**ТЕОРЕМА 4.4.7.** *Множество $\mathrm{End}(D, V)$ эндоморфизмов $D$-модуля $V$ является $D$-модулем.*

ДОКАЗАТЕЛЬСТВО. Теорема является следсвием определения 4.3.16 и следсвия 4.4.6. $\qquad \square$

### 4.5. $D$-модуль $\mathcal{L}(D; A \to B)$

**ТЕОРЕМА 4.5.1.**

(4.5.1) $$\mathcal{L}(D; A^p \to \mathcal{L}(D; A^q \to B)) = \mathcal{L}(D; A^{p+q} \to B)$$

ДОКАЗАТЕЛЬСТВО. $\qquad \square$

**ТЕОРЕМА 4.5.2.** *Пусть*

$$\overline{\overline{e}} = \{e_i : i \in I\}$$

*базис $D$-модуля $A$. Множество*

(4.5.2) $$\overline{\overline{h}} = \{h^i \in \mathcal{L}(D; A \to D) : i \in I, h^i \circ e_j = \delta^i_j\}$$

*является базисом $D$-модуля $\mathcal{L}(D; A \to D)$.*

Доказательство.

Лемма 4.5.3. *Отображения $h^i$ линейно независимы.*

Доказательство. Пусть существуют $D$-числа $c_i$ такие, что

$$c_i h^i = 0$$

Тогда для любого $A$-числа $e_j$

$$0 = c_i h^i \circ e_j = c_i \delta_j^i = c_j$$

Лемма является следствием определения 4.1.16.                                   ⊙

Лемма 4.5.4. *Отображение $f \in \mathcal{L}(D; A \to D)$, является линейной композицией отображений $h^i$.*

Доказательство. Для любого $A$-числа $a$,

$$(4.5.3) \qquad\qquad a = a^i e_i$$

равенство

$$(4.5.4) \qquad h^i \circ a = h^i \circ (a^j e_j) = a^j (h^i \circ e_j) = a^j \delta_j^i = a^i$$

является следствием равенств (4.5.2), (4.5.3) и теоремы 4.3.10. Равенство

$$(4.5.5) \qquad f \circ a = f \circ (a^i e_i) = a^i (f \circ e_i) = (f \circ e_i)(h^i \circ a)$$

является следствием равенства (4.5.4). Равенство

$$f = (f \circ e_i) h^i$$

является следствием равенств (4.4.10), (4.4.21), (4.5.5).                        ⊙

Теорема является следствием лемм 4.5.3, 4.5.4 и определения 4.1.24.              □

Теорема 4.5.5. *Пусть $D$ - коммутативное кольцо. Пусть*

$$\overline{\overline{e}}_i = \{ e_{i \cdot i} : i \in I_i \}$$

*базис $D$-модуля $A_i$, $i = 1, ..., n$. Пусть*

$$\overline{\overline{e}}_B = \{ e_{B \cdot i} : i \in I \}$$

*базис $D$-модуля $B$. Множество*

$$(4.5.6) \qquad \begin{aligned} \overline{\overline{h}} = \{ \, & h_i^{i_1 \ldots i_n} \in \mathcal{L}(D; A_1 \times \ldots \times A_n \to B) : i \in I, i_i \in I_i, i = 1, \ldots, n, \\ & h_i^{i_1 \ldots i_n} \circ (e_{1 \cdot j_1}, \ldots, e_{n \cdot j_n}) = \delta_{j_1}^{i_1} \ldots \delta_{j_n}^{i_n} e_{B \cdot i} \, \} \end{aligned}$$

*является базисом $D$-модуля $\mathcal{L}(D; A_1 \times \ldots \times A_n \to B)$.*

Доказательство.

Лемма 4.5.6. *Отображения $h_i^{i_1 \ldots i_n}$ линейно независимы.*

Доказательство. Пусть существуют $D$-числа $c_{i_1 \ldots i_n}^i$ такие, что

$$c_{i_1 \ldots i_n}^i h_i^{i_1 \ldots i_n} = 0$$

Тогда для любого набора индексов $j_1, ..., j_n$

$$0 = c_{i_1 \ldots i_n}^i h_i^{i_1 \ldots i_n} \circ (e_{1 \cdot j_1}, ..., e_{n \cdot j_n}) = c_{i_1 \ldots i_n}^i \delta_{j_1}^{i_1} \ldots \delta_{j_n}^{i_n} e_{B \cdot i} = c_{j_1 \ldots j_n}^i e_{B \cdot i}$$

Следовательно, $c^i_{j_1 \ldots j_n} = 0$. Лемма является следствием определения 4.1.16. ⊙

Лемма 4.5.7. *Отображение* $f \in \mathcal{L}(D; A_1 \times \ldots \times A_n \to B)$ *является линейной композицией отображений* $h^{i_1 \ldots i_n}_i$.

Доказательство. Для любого $A_1$-числа $a_1$

$$(4.5.7) \qquad a_1 = a^{i_1}_1 e_{1 \cdot i_1}$$

..., для любого $A_n$-числа $a_n$

$$(4.5.8) \qquad a_n = a^{i_n}_n e_{n \cdot i_n}$$

равенство

$$
(4.5.9)
\begin{aligned}
h^{i_1 \ldots i_n}_i \circ (a_1, \ldots, a_n) &= h^{i_1 \ldots i_n}_i \circ (a^{j_1}_1 e_{1 \cdot j_1}, \ldots, a^{j_n}_n e_{n \cdot j_n}) \\
&= a^{j_1}_1 \ldots a^{j_n}_n (h^{i_1 \ldots i_n}_i \circ (e_{1 \cdot j_1}, \ldots, e_{n \cdot j_n})) \\
&= a^{j_1}_1 \ldots a^{j_n}_n \delta^{i_1}_{j_1} \ldots \delta^{i_n}_{j_n} e_{B \cdot i} \\
&= a^{i_1}_1 \ldots a^{i_n}_n e_{B \cdot i}
\end{aligned}
$$

является следствием равенств (4.5.6), (4.5.7), (4.5.8) и теоремы 4.3.10. Равенство

$$
(4.5.10)
\begin{aligned}
f \circ (a_1, \ldots, a_n) &= f \circ (a^{j_1}_1 e_{1 \cdot j_1}, \ldots, a^{j_n}_n e_{n \cdot j_n}) \\
&= a^{j_1}_1 \ldots a^{j_n}_n f \circ (e_{1 \cdot j_1} \ldots e_{n \cdot j_n})
\end{aligned}
$$

является следствием равенств (4.5.7), (4.5.8). Так как

$$f \circ (e_{1 \cdot j_1} \ldots e_{n \cdot j_n}) \in B$$

то

$$(4.5.11) \qquad f \circ (e_{1 \cdot j_1} \ldots e_{n \cdot j_n}) = f^i_{j_1 \ldots j_n} e_i$$

Равенство

$$
(4.5.12)
\begin{aligned}
f \circ (a_1, \ldots, a_n) &= a^{i_1}_1 \ldots a^{i_n}_n f^i_{i_1 \ldots i_n} e_i \\
&= f^i_{i_1 \ldots i_n} h^{i_1 \ldots i_n}_i \circ (a_1, \ldots, a_n)
\end{aligned}
$$

является следствием равенств (4.5.9), (4.5.10), (4.5.11). Равенство

$$f = f^i_{i_1 \ldots i_n} h^{i_1 \ldots i_n}_i$$

является следствием равенств (4.4.2), (4.4.12), (4.5.12). ⊙

Теорема является следствием лемм 4.5.6, 4.5.7 и определения 4.1.24. □

Теорема 4.5.8. *Пусть* $A_1$, ..., $A_n$, $B$ - *свободные модули над коммутативным кольцом* $D$. $D$-*модуль* $\mathcal{L}(D; A_1 \times \ldots \times A_n \to B)$ *является свободным* $D$-*модулем.*

Доказательство. Теорема является следствием теоремы 4.5.5. □

## 4.6. Тензорное произведение $D$-модулей

**Теорема 4.6.1.** *Коммутативное кольцо $D$ является абелевой мульти-пликативной $\Omega$-группой.*

Доказательство. Мы полагаем, что произведение $\circ$ в кольце $D$ определён согласно правилу
$$a \circ b = ab$$
Так как произведение в кольце дистрибутивно относительно сложения, теорема является следствием определений 2.4.5, 2.4.6.                                   $\square$

**Теорема 4.6.2.** *Пусть $A_1, ..., A_n$ - модули над коммутативным кольцом $D$.* **Тензорное произведение**
$$f : A_1 \times ... \times A_n \to A_1 \otimes ... \otimes A_n$$
*$D$-модулей $A_1, ..., A_n$ существует и единственно. Мы пользуемся обозначением*
$$f \circ (a_1, ..., a_n) = a_1 \otimes ... \otimes a_n$$
*для образа отображения $f$.*

Доказательство. Теорема является следствием определения 4.1.3 и теорем 2.6.3, 4.6.1.                                                     $\square$

**Теорема 4.6.3.** *Пусть $D$ - коммутативное кольцо. Пусть $A_1, ..., A_n$ - $D$-модули. Тензорное произведение дистрибутивно относительно сложения*
$$\begin{aligned}&a_1 \otimes ... \otimes (a_i + b_i) \otimes ... \otimes a_n \\ (4.6.1) \quad &= a_1 \otimes ... \otimes a_i \otimes ... \otimes a_n + a_1 \otimes ... \otimes b_i \otimes ... \otimes a_n \\ &a_i, b_i \in A_i\end{aligned}$$
*Представление кольца $D$ в тензорном произведении определено равенством*
$$\begin{aligned}(4.6.2) \quad &a_1 \otimes ... \otimes (ca_i) \otimes ... \otimes a_n = c(a_1 \otimes ... \otimes a_i \otimes ... \otimes a_n) \\ &a_i \in A_i \quad c \in D\end{aligned}$$

Доказательство. Равенство (4.6.1) является следствием равенства (2.6.1). Равенство (4.6.2) является следствием равенства (2.6.2).               $\square$

**Теорема 4.6.4.** *Пусть $A_1, ..., A_n$ - модули над коммутативным кольцом $D$. Пусть*
$$f : A_1 \times ... \times A_n \to A_1 \otimes ... \otimes A_n$$
*полилинейное отображение, определённое равенством*
$$(4.6.3) \quad f \circ (d_1, ..., d_n) = d_1 \otimes ... \otimes d_n$$
*Пусть*
$$g : A_1 \times ... \times A_n \to V$$
*полилинейное отображение в $D$-модуль $V$. Существует линейное отображение*
$$h : A_1 \otimes ... \otimes A_n \to V$$

*такое, что диаграмма*

(4.6.4)

*коммутативна. Отображение $h$ определено равенством*

(4.6.5)
$$h \circ (a_1 \otimes ... \otimes a_n) = g \circ (a_1, ..., a_n)$$

Доказательство. Теорема является следствием теоремы 2.6.5 и определений 4.3.9, 4.4.1. □

Теорема 4.6.5. *Отображение*
$$(v_1, ..., v_n) \in V_1 \times ... \times V_n \to v_1 \otimes ... \otimes v_n \in V_1 \otimes ... \otimes V_n$$
*является полилинейным отображением.*

Доказательство. Теорема является следствием теоремы 2.6.6 и определения 4.4.1. □

Теорема 4.6.6. *Тензорное произведение $A_1 \otimes ... \otimes A_n$ свободных конечномерных модулей $A_1$, ..., $A_n$ над коммутативным кольцом $D$ является свободным конечномерным модулем.*

*Пусть $\overline{\overline{e}}_i$ - базис модуля $A_i$ над кольцом $D$. Произвольный тензор $a \in A_1 \otimes ... \otimes A_n$ можно представить в виде*

(4.6.6)
$$a = a^{i_1...i_n} e_{1 \cdot i_1} \otimes ... \otimes e_{n \cdot i_n}$$

*Выражение $a^{i_1...i_n}$ определено однозначно и называется **стандартной компонентой тензора**.*

Доказательство. Вектор $a_i \in A_i$ имеет разложение
$$a_i = a_i^k \overline{\overline{e}}_{i \cdot k}$$
относительно базиса $\overline{\overline{e}}_i$. Из равенств (4.6.1), (4.6.2) следует
$$a_1 \otimes ... \otimes a_n = a_1^{i_1} ... a_n^{i_n} e_{1 \cdot i_1} \otimes ... \otimes e_{n \cdot i_n}$$
Так как множество тензоров $a_1 \otimes ... \otimes a_n$ является множеством образующих модуля $A_1 \otimes ... \otimes A_n$, то тензор $a \in A_1 \otimes ... \otimes A_n$ можно записать в виде

(4.6.7)
$$a = a^s a_{s \cdot 1}^{i_1} ... a_{s \cdot n}^{i_n} e_{1 \cdot i_1} \otimes ... \otimes e_{n \cdot i_n}$$
где $a^s$, $a_{s \cdot 1}^{i_1}$, ..., $a_{s \cdot n}^{i_n} \in F$. Положим
$$a^s a_{s \cdot 1}^{i_1} ... a_{s \cdot n}^{i_n} = a^{i_1...i_n}$$
Тогда равенство (4.6.7) примет вид (4.6.6).

Следовательно, множество тензоров $e_{1 \cdot i_1} \otimes ... \otimes e_{n \cdot i_n}$ является множеством образующих модуля $A_1 \otimes ... \otimes A_n$. Так как размерность модуля $A_i$, $i = 1$, ..., $n$, конечна, то конечно множество тензоров $e_{1 \cdot i_1} \otimes ... \otimes e_{n \cdot i_n}$. Следовательно, множество тензоров $e_{1 \cdot i_1} \otimes ... \otimes e_{n \cdot i_n}$ содержит базис модуля $A_1 \otimes ... \otimes A_n$, и модуль $A_1 \otimes ... \otimes A_n$ является свободным модулем над кольцом $D$.                □

### 4.7. Прямая сумма $D$-модулей

Пусть $D$-Module - категория $D$-модулей , морфизмы которой - это гомоморфизмы $D$-модуля.

ОПРЕДЕЛЕНИЕ 4.7.1. *Копроизведение в категории $D$-Module называется* **прямой суммой**.[4.10]*Мы будем пользоваться записью* $V \oplus W$ *для прямой суммы $D$-модулей $V$ и $W$.*                □

ТЕОРЕМА 4.7.2. *В категории $D$-Module существует прямая сумма $D$-модулей.*[4.11]*Пусть $\{V_i, i \in I\}$ - семейство $D$-модулей. Тогда представление*

$$D \xrightarrow{\ *\ } \bigoplus_{i \in I} V_i \qquad d\left(\bigoplus_{i \in I} a_i\right) = \bigoplus_{i \in I} da_i$$

*кольца $D$ в прямой сумме абелевых групп*

$$V = \bigoplus_{i \in I} V_i$$

*является прямой суммой $D$-модулей*

$$V = \bigoplus_{i \in I} V_i$$

ДОКАЗАТЕЛЬСТВО. Пусть $V_i, i \in I$, - множество $D$-модулей. Пусть

$$V = \bigoplus_{i \in I} V_i$$

прямая сумма абелевых групп. Пусть $v = (v_i, i \in I) \in V$. Мы определим действие $D$-числа $a$ на $V$-число $v$ согласно равенству

(4.7.1)                     $a(v_i, i \in I) = (av_i, i \in I)$

Согласно теореме 4.1.8,

(4.7.2)                     $(ab)v_i = a(bv_i)$

(4.7.3)                     $a(v_i + w_i) = av_i + aw_i$

(4.7.4)                     $(a + b)v_i = av_i + bv_i$

(4.7.5)                     $1v_i = v_i$

---

4.10   Смотри так же определение прямой суммы модулей в [1], страница 98. Согласно теореме 1 на той же странице, прямая сумма модулей существует.
4.11   Смотри также предложение [1]-1 на странице 98.

для любых $a, b \in D$, $v_i, w_i \in V_i$, $i \in I$. Равенства

$$
\begin{aligned}
(4.7.6) \quad (ab)v &= (ab)(v_i, i \in I) = ((ab)v_i, i \in I) \\
&= (a(bv_i), i \in I) = a(bv_i, i \in I) = a(b(v_i, i \in I)) \\
&= a(bv)
\end{aligned}
$$

$$
\begin{aligned}
(4.7.7) \quad a(v+w) &= a((v_i, i \in I) + (w_i, i \in I)) = a(v_i + w_i, i \in I) \\
&= (a(v_i + w_i), i \in I) = (av_i + aw_i, i \in I) \\
&= (av_i, i \in I) + (aw_i, i \in I) = a(v_i, i \in I) + a(w_i, i \in I) \\
&= av + aw
\end{aligned}
$$

$$
\begin{aligned}
(4.7.8) \quad (a+b)v &= (a+b)(v_i, i \in I) = ((a+b)v_i, i \in I) \\
&= (av_i + bv_i, i \in I) = (av_i, i \in I) + (bv_i, i \in I) \\
&= a(v_i, i \in I) + b(v_i, i \in I) = av + bv
\end{aligned}
$$

$$
1v = 1(v_i, i \in I) = (1v_i, i \in I) = (v_i, i \in I) = v
$$

являются следствием равенств (3.7.4), (4.7.1), (4.7.2), (4.7.3), (4.7.4), (4.7.5) для любых $a, b \in D$, $v = (v_i, i \in I)$, $w = (w_i, i \in I) \in V$.

Из равенства (4.7.7) следует, что отображение

$$
a : v \to av
$$

определённое равенством (4.7.1), является эндоморфизмом абелевой группы $V$. Из равенств (4.7.6), (4.7.8) следует, что отображение (4.7.1) является гомоморфизмом кольца $D$ в множество эндоморфизмов абелевой группы $V$. Следовательно, отображение (4.7.1) является представлением кольца $D$ в абелевой группе $V$.

Пусть $a_1$, $a_2$ - такие $D$-числа, что

$$
(4.7.9) \quad a_1(v_i, i \in I) = a_2(v_i, i \in I)
$$

для произвольного $V$-числа $v = (v_i, i \in I)$. Равенство

$$
(4.7.10) \quad (a_1 v_i, i \in I) = (a_2 v_i, i \in I)
$$

является следствием равенства (4.7.9). Из равенства (4.7.10) следует, что

$$
(4.7.11) \quad \forall i \in I : a_1 v_i = a_2 v_i
$$

Поскольку для любого $i \in I$, $V_i$ - $D$-модуль, то соответствующее представление эффективно согласно определению 4.1.3. Следовательно, равенство $a_1 = a_2$ является следствием утверждения (4.7.11) и отображение (4.7.1) является эффективным представлением кольца $D$ в абелевой группе $V$. Следовательно, абелевая группе $V$ является $D$-модулем.

Пусть

$$
\{f_i : V_i \to W, i \in I\}
$$

семейство линейных отображений в $D$-модуль $W$. Мы определим отображение

$$
f : V \to W
$$

пользуясь равенством

$$
(4.7.12) \quad f\left(\bigoplus_{i \in I} x_i\right) = \sum_{i \in I} f_i(x_i)
$$

Сумма в правой части равенства (4.7.12) конечна, так как все слагаемые, кроме конечного числа, равны 0. Из равенства

$$f_i(x_i + y_i) = f_i(x_i) + f_i(y_i)$$

и равенства (4.7.12) следует, что

$$f(\bigoplus_{i \in I}(x_i + y_i)) = \sum_{i \in I} f_i(x_i + y_i)$$
$$= \sum_{i \in I}(f_i(x_i) + f_i(y_i))$$
$$= \sum_{i \in I} f_i(x_i) + \sum_{i \in I} f_i(y_i)$$
$$= f(\bigoplus_{i \in I} x_i) + f(\bigoplus_{i \in I} y_i)$$

Из равенства

$$f_i(dx_i) = df_i(x_i)$$

и равенства (4.7.12) следует, что

$$f((dx_i, i \in I)) = \sum_{i \in I} f_i(dx_i) = \sum_{i \in I} df_i(x_i) = d\sum_{i \in I} f_i(x_i)$$
$$= df((x_i, i \in I))$$

Следовательно, отображение $f$ является линейным отображением. Равенство

$$f \circ \lambda_j(x) = \sum_{i \in I} f_i(\delta_j^i x) = f_j(x)$$

является следствием равенств (3.7.5), (4.7.12). Так как отображение $\lambda_i$ инъективно, то отображение $f$ определено однозначно. Следовательно, теорема является следствием определений 2.2.12, 3.7.2. $\qquad \square$

**Теорема 4.7.3.** *Прямая сумма $D$-модулей $V_1$, ..., $V_n$ совпадает с их прямым произведением*

$$V_1 \oplus ... \oplus V_n = V_1 \times ... \times V_n$$

Доказательство. Теорема является следствием теоремы 4.7.2. $\qquad \square$

**Теорема 4.7.4.** *Пусть $\overline{\overline{e}}_1 = (e_{1i}, i \in I)$ базис $D$-модуля столбцов V. Тогда мы можем представить $D$-модуль V как прямую сумму*

$$(4.7.13) \qquad V = \bigoplus_{i \in I} De_i$$

Доказательство. Отображение

$$f: v^i e_i \in V \to (v^i, i \in I) \in \bigoplus_{i \in I} De_i$$

является изоморфизмом. $\qquad \square$

**Теорема 4.7.5.** *Пусть $\overline{\overline{e}}_1 = (e_1^i, i \in I)$ базис $D$-модуля строк V. Тогда мы можем представить $D$-модуль V как прямую сумму*

$$(4.7.14) \qquad V = \bigoplus_{i \in I} De^i$$

Доказательство. Отображение

$$f: v_i e^i \in V \to (v_i, i \in I) \in \bigoplus_{i \in I} De^i$$

является изоморфизмом. $\qquad \square$

Теорема 4.7.6. *Пусть* $V^1$, ..., $V^n$ - *$D$-модули и*

$$V = V^1 \oplus ... \oplus V^n$$

*Представим a-число*

$$a = a^1 \oplus ... \oplus a^n$$

*как вектор столбец*

$$a = \begin{pmatrix} a^1 \\ ... \\ a^n \end{pmatrix}$$

*Представим линейное отображение*

$$f : V \to W$$

*как вектор строку*

$$f = \begin{pmatrix} f_1 & ... & f_n \end{pmatrix}$$

$$f_i : V^i \to W$$

*Тогда значение отображения $f$ в $V$-числе $a$ можно представить как произведение матриц*

$$(4.7.15) \qquad f \circ a = \begin{pmatrix} f_1 & ... & f_n \end{pmatrix} \circ^\circ \begin{pmatrix} a^1 \\ ... \\ a^n \end{pmatrix} = f_i \circ a^i$$

Доказательство. Теорема является следствием определения (4.7.12). $\square$

Теорема 4.7.7. *Пусть* $W^1$, ..., $W^m$ - *$D$-модули и*

$$W = W^1 \oplus ... \oplus W^m$$

*Представим $W$-число*

$$b = b^1 \oplus ... \oplus b^m$$

*как вектор столбец*

$$b = \begin{pmatrix} b^1 \\ ... \\ b^m \end{pmatrix}$$

*Тогда линейное отображение*

$$f : V \to W$$

*имеет представление как вектор столбец отображений*

$$f = \begin{pmatrix} f^1 \\ ... \\ f^m \end{pmatrix}$$

*таким образом, что если $b = f \circ a$, то*

$$\begin{pmatrix} b^1 \\ ... \\ b^m \end{pmatrix} = \begin{pmatrix} f^1 \\ ... \\ f^m \end{pmatrix} \circ a = \begin{pmatrix} f^1 \circ a \\ ... \\ f^m \circ a \end{pmatrix}$$

Доказательство. Теорема является следствием теоремы 2.2.8. $\qquad\square$

Теорема 4.7.8. *Пусть* $V^1, ..., V^n, \quad W^1, ..., W^m$ *- $D$-модули и*

$$V = V^1 \oplus ... \oplus V^n$$
$$W = W^1 \oplus ... \oplus W^m$$

*Представим $V$-число*

$$a = a^1 \oplus ... \oplus a^n$$

*как вектор столбец*

$$(4.7.16) \qquad a = \begin{pmatrix} a^1 \\ ... \\ a^n \end{pmatrix}$$

*Представим $W$-число*

$$b = b^1 \oplus ... \oplus b^m$$

*как вектор столбец*

$$(4.7.17) \qquad b = \begin{pmatrix} b^1 \\ ... \\ b^m \end{pmatrix}$$

*Тогда линейное отображение*

$$f : V \to W$$

*имеет представление как матрица отображений*

$$(4.7.18) \qquad f = \begin{pmatrix} f^1_1 & ... & f^1_n \\ ... & ... & ... \\ f^m_1 & ... & f^m_n \end{pmatrix}$$

*таким образом, что если $b = f \circ a$, то*

$$(4.7.19) \qquad \begin{pmatrix} b^1 \\ ... \\ b^m \end{pmatrix} = \begin{pmatrix} f^1_1 & ... & f^1_n \\ ... & ... & ... \\ f^m_1 & ... & f^m_n \end{pmatrix} \circ \begin{pmatrix} a^1 \\ ... \\ a^n \end{pmatrix} = \begin{pmatrix} f^1_i \circ a^i \\ ... \\ f^m_i \circ a^i \end{pmatrix}$$

*Отображение*

$$f^i_j : V^j \to W^i$$

*является линейным отображением и называется* **частным линейным отображением**

Доказательство. Согласно теореме 4.7.7, существует множество линейных отображений

$$f^i : V \to W^i$$

таких, что

$$(4.7.20) \qquad \begin{pmatrix} b^1 \\ ... \\ b^m \end{pmatrix} = \begin{pmatrix} f^1 \\ ... \\ f^m \end{pmatrix} \circ a = \begin{pmatrix} f^1 \circ a \\ ... \\ f^m \circ a \end{pmatrix}$$

Согласно теореме 4.7.6, для каждого $i$, существует множество линейных отображений

$$f^i_j : V^j \to W^i$$

таких, что

$$(4.7.21) \qquad f^i \circ a = \begin{pmatrix} f^i_1 & ... & f^i_n \end{pmatrix} \circ \begin{pmatrix} a^1 \\ ... \\ a^n \end{pmatrix} = f^i_j \circ a^i$$

Если мы отождествим матрицы

$$\begin{pmatrix} \begin{pmatrix} f^1_1 & ... & f^1_n \end{pmatrix} \\ ... \\ \begin{pmatrix} f^m_1 & ... & f^m_n \end{pmatrix} \end{pmatrix} = \begin{pmatrix} f^1_1 & ... & f^1_n \\ ... & ... & ... \\ f^m_1 & ... & f^m_n \end{pmatrix}$$

то равенство (4.7.19) является следствием равенств (4.7.20), (4.7.21). $\qquad \square$

Пусть $B^1$, ..., $B^m$ - $D$-алгебры. Тогда мы можем записать линейное отображение $f^i_j$ с помощью $B^i \otimes B^i$-чисел.

Теорема 4.7.9. *Пусть $D$-модуль $V$ является свободной аддитивной группой. Пусть $U$ - подмодуль $D$-модуля $V$. Существует подмодуль $W$ такой, что*

$$V = U \oplus W$$

Доказательство. Согласно определениям 4.1.3, 4.1.6, $U$ является подгруппой аддитивной группы $V$. Согласно теоремам 3.7.9, 3.7.10, существует подгруппа $W$ аддитивной группы $V$ такая, что группа $W$ изоморфна группе $V/U$ и

$$V = U \oplus W$$

Согласно теореме 4.7.2, действие $D$-числа $a$ на произвольное $V$-число

$$v = u \oplus w$$

имеет вид

$$av = (au) \oplus (aw)$$

Следовательно, $W$ является подмодулем $D$-модуля $V$. $\qquad \square$

## 4.8. Ядро линейного отображения

Теорема 4.8.1. *Пусть отображение*

(4.3.12)  $\boxed{(h : D_1 \to D_2 \quad f : V_1 \to V_2)}$

*является линейным отображением $D_1$-модуля $V_1$ в $D_2$-модуль $V_2$. Множество*

(4.8.1)                      $\ker(h, f) = \{v \in V : f \circ v = 0\}$

*называется* **ядром линейного отображения** $(h, f)$. *Ядро линейного отображения $(h, f)$ является подмодулем модуля $V_1$.*

Доказательство. Пусть $v, w \in V_1$. Тогда

(4.8.2)                      $f \circ v = 0 \quad f \circ w = 0$

Равенство

$$f \circ (v + w) = f \circ v + f \circ w = 0$$

является следствием равенств (4.3.4), (4.8.2). Следовательно, множество $\ker(h, f)$ является подгруппой аддитивной группы $V_1$. Равенство

$$f \circ (pv) = h(p)(f \circ v) = 0$$

является следствием равенств (4.3.5), (4.8.2). Следовательно, определено представление кольца $D_1$ на множестве $\ker(h, f)$. Согласно определению 4.1.6, множество $\ker(h, f)$ является подмодулем $D_1$-модуля $V_1$.                    $\square$

Определение 4.8.2. *Пусть $V_i$, $i = 1, ..., n$, - $D_i$-модули. Последовательность гомоморфизмов*

$$V_1 \xrightarrow{f_1} V_2 \xrightarrow{f_2} ... \xrightarrow{f_{n-1}} V_n$$
$$\uparrow \qquad \uparrow \qquad \qquad \uparrow$$
$$D_1 \xrightarrow{h_1} D_2 \xrightarrow{h_2} ... \xrightarrow{h_{n-1}} D_n$$

*называется* **точной**, *если $\forall i$, $i = 1, ..., n - 2$, $\operatorname{Im} f_i = \ker(h_{i+1}, f_{i+1})$.*                    $\square$

Теорема 4.8.3. *Пусть отображение*

(4.3.40)  $\boxed{f : V_1 \to V_2}$

*является линейным отображением $D$-модуля $V_1$ в $D$-модуль $V_2$. Множество*

(4.8.3)                      $\ker f = \{v \in V : f \circ v = 0\}$

*называется* **ядром линейного отображения** $f$. *Ядро линейного отображения $f$ является подмодулем модуля $V_1$.*

Доказательство. Пусть $v, w \in V_1$. Тогда

(4.8.4)                      $f \circ v = 0 \quad f \circ w = 0$

Равенство

$$f \circ (v + w) = f \circ v + f \circ w = 0$$

является следствием равенств (4.3.32), (4.8.4). Следовательно, множество $\ker f$ является подгруппой аддитивной группы $V_1$. Равенство

$$f \circ (pv) = p(f \circ v) = 0$$

является следствием равенств (4.3.5), (4.8.4). Следовательно, определено представление кольца $D$ на множестве $\ker f$. Согласно определению 4.1.6, множество $\ker f$ является подмодулем $D$-модуля $V_1$. □

ОПРЕДЕЛЕНИЕ 4.8.4. *Пусть* $V_i$, $i = 1, \ldots, n$, *- $D$-модули. Последовательность гомоморфизмов*

$$V_1 \xrightarrow{f_1} V_2 \xrightarrow{f_2} \ldots \xrightarrow{f_{n-1}} V_n$$

*называется* **точной**, *если* $\forall i$, $i = 1, \ldots, n-2$, $\operatorname{Im} f_i = \ker f_{i+1}$. □

ТЕОРЕМА 4.8.5. *Пусть $D$-модуль $V_1$ является свободной аддитивной группой. Пусть $U$ - подмодуль $D$-модуля $V_1$. Существует гомоморфизм*

$$(4.3.40) \quad f : V_1 \to V_2$$

*такой, что*

$$\ker f = U$$

ДОКАЗАТЕЛЬСТВО. Согласно теореме 4.7.9, существует подмодуль $W$ такой, что

$$V_1 = U \oplus W$$

Отображение

$$f : u \oplus w \in U \oplus W \to w \in W$$

является гомоморфизмом и

$$\ker f = U$$

□

# $D$-алгебра

## 5.1. Алгебра над коммутативным кольцом

ОПРЕДЕЛЕНИЕ 5.1.1. *Пусть $D$ - коммутативное кольцо. $D$-модуль $A$ называется **алгеброй над кольцом** $D$ или $D$-**алгеброй**, если определена операция произведения [5.1] в $A$*

(5.1.1) $$v\,w = C \circ (v, w)$$

*где $C$ - билинейное отображение*

$$C : A \times A \to A$$

*Если $A$ является свободным $D$-модулем, то $A$ называется **свободной алгеброй над кольцом** $D$.* □

ТЕОРЕМА 5.1.2. *Произведение в алгебре $A$ дистрибутивно по отношению к сложению*

(5.1.2) $$(a + b)c = ac + bc$$

(5.1.3) $$a(b + c) = ab + ac$$

Доказательство. Утверждение теоремы следует из цепочки равенств

$$(a + b)c = f \circ (a + b, c) = f \circ (a, c) + f \circ (b, c) = ac + bc$$
$$a(b + c) = f \circ (a, b + c) = f \circ (a, b) + f \circ (a, c) = ab + ac$$

□

ОПРЕДЕЛЕНИЕ 5.1.3. *Если произведение в $D$-алгебре $A$ имеет единичный элемент, то $D$-алгебра $A$ называется **унитальной алгеброй**.* [5.2] □

Произведение в алгебре может быть ни коммутативным, ни ассоциативным. Следующие определения основаны на определениях, данным в [23], с. 13.

ОПРЕДЕЛЕНИЕ 5.1.4. **Коммутатор**

$$[a, b] = ab - ba$$

*служит мерой коммутативности в $D$-алгебре $A$. $D$-алгебра $A$ называется **коммутативной**, если*

$$[a, b] = 0$$

□

---

[5.1] Я следую определению, приведенному в [23], страница 1, [17], страница 4. Утверждение, верное для произвольного $D$-модуля, верно также для $D$-алгебры.

[5.2] Смотри определение унитальной алгебры также на страницах [7]-137.

Определение 5.1.5. **Ассоциатор**

(5.1.4) $$(a, b, c) = (ab)c - a(bc)$$

*служит мерой ассоциативности в D-алгебре A. D-алгебра A называется* **ассоциативной**, *если*

$$(a, b, c) = 0$$

□

Теорема 5.1.6. *Пусть A - алгебра над коммутативным кольцом D.*[5.3]

(5.1.5) $$a(b, c, d) + (a, b, c)d = (ab, c, d) - (a, bc, d) + (a, b, cd)$$

*для любых a, b, c, d ∈ A.*

Доказательство. Равенство (5.1.5) следует из цепочки равенств

$$a(b, c, d) + (a, b, c)d = a((bc)d - b(cd)) + ((ab)c - a(bc))d$$
$$= a((bc)d) - a(b(cd)) + ((ab)c)d - (a(bc))d$$
$$= ((ab)c)d - (ab)(cd) + (ab)(cd)$$
$$+ a((bc)d) - a(b(cd)) - (a(bc))d$$
$$= (ab, c, d) - (a(bc))d + a((bc)d) + (ab)(cd) - a(b(cd))$$
$$= (ab, c, d) - (a, (bc), d) + (a, b, cd)$$

□

Теорема 5.1.7. *Пусть $\overline{\overline{e}}$ - базис D-алгебры A. Тогда*

(5.1.6) $$(a, b, c) = A_{ijl}^k a^i b^j c^l e_k$$

*где* **координаты ассоциатора** *определены равенством*

(5.1.7) $$A_{ijl}^k = C_{pl}^k C_{ij}^p - C_{ip}^k C_{jl}^p$$

Доказательство. Равенство

(5.1.8) $$(a, b, c) = (C_{pl}^k (ab)^p c^l - C_{ip}^k a^i (bc)^p) e_k$$
$$= (C_{pl}^k C_{ij}^p - C_{ip}^k C_{jl}^p) a^i b^j c^l e_k$$

*является следствием равенства* (5.1.4). *Равенство* (5.1.7) *является следствием равенства* (5.1.8).

□

Определение 5.1.8. **Ядро** *D-алгебры A - это множество*[5.4]

$$N(A) = \{a \in A : \forall b, c \in A, (a, b, c) = (b, a, c) = (b, c, a) = 0\}$$

□

---

Определение 5.1.9. **Центр** $D$**-алгебры** $A$ - *это множество*[5.5]
$$Z(A) = \{a \in A : a \in N(A), \forall b \in A, ab = ba\}$$
□

| Определение 5.1.10. **Левый идеал** [5.6] $A_1$ $D$-*алгебры* $A$ - *это подгруппа аддитивной группы* $A$ *такая, что* $$AA_1 \subseteq A_1$$ □ | Определение 5.1.12. **Правый идеал** [5.6] $A_1$ $D$-*алгебры* $A$ - *это подгруппа аддитивной группы* $A$ *такая, что* $$A_1A \subseteq A_1$$ □ |
|---|---|
| Теорема 5.1.11. *Если* $A_1$ *является левым идеалом унитальной* $D$-*алгебры* $A$, *то* $$AA_1 = A_1$$ | Теорема 5.1.13. *Если* $A_1$ *является правым идеалом унитальной* $D$-*алгебры* $A$, *то* $$A_1A = A_1$$ |
| Доказательство. Теорема является следствием определения 5.1.10 и равенства $$eA_1 = A_1$$ □ | Доказательство. Теорема является следствием определения 5.1.12 и равенства $$A_1e = A_1$$ □ |

Определение 5.1.14. **Идеал** [5.6] - *это подгруппа аддитивной группы* $A$, *которая одновременно является левым и правым идеалами.* □

Если $A$ - $D$-алгебра, то множество $\{0\}$ и само $A$ являются идеалами, которые мы называем тривиальными идеалами.

Теорема 5.1.15. *Пусть* $D$ - *коммутативное кольцо. Пусть* $A$ - *ассоциативная* $D$-*алгебра с единицей* $1_A$.

5.1.15.1: *Отображение*
$$f : d \in D \to d1_A \in A$$
*является гомоморфизмом кольца* $D$ *в центр алгебры* $A$ *и позволяет отождествить* $D$-*число* $d$ *и* $A$-*число* $d\,1_A$.

5.1.15.2: *Кольцо* $D$ *является подалгеброй алгебры* $A$.

5.1.15.3: $D \subset Z(A)$.

Доказательство. Пусть $d_1, d_2 \in D$. Согласно утверждению 4.1.8.3,

(5.1.9)                                $d_1 1_A + d_2 1_A = (d_1 + d_2) 1_A$

---

Согласно определению 5.1.1,

$$(5.1.10) \qquad (d_1 1_A)(d_2 1_A) = d_1(1_A(d_2 1_A)) = d_1(d_2(1_A 1_A))$$
$$= d_1(d_2 1_A) = (d_1 d_2)1_A$$

Утверждение 5.1.15.1 является следствием равенств (5.1.9), (5.1.10). Утверждение 5.1.15.2 является следствием утверждения 5.1.15.1.

Согласно определению 5.1.1,

$$(5.1.11) \qquad a(d1_A) = d(a1_A) = d(1_A a) = (d1_A)a$$

Утверждение 5.1.15.3 следует из равенства (5.1.11). $\qquad\square$

---

**Лемма 5.1.16.** *Пусть* $d \in D$. *Отображение*

$$(5.1.12) \qquad d : b \in A \to db \in A$$

*является эндоморфизмом абелевой группы* $A$.

**Доказательство.** Утверждение

$$db \in A$$

является следствием определений 4.1.3, 5.1.1. Согласно равенству (4.1.22) отображение (5.1.12) является эндоморфизмом абелевой группы $A$. $\qquad\square$

**Лемма 5.1.17.** *Пусть* $a \in A$. *Отображение*

$$(5.1.13) \qquad a : b \in A \to ab \in A$$

*является эндоморфизмом абелевой группы* $A$.

**Доказательство.** Утверждение

$$ab \in A$$

является следствием определения 5.1.1. Согласно теореме 5.1.2 отображение (5.1.13) является эндоморфизмом абелевой группы $A$. $\qquad\square$

**Лемма 5.1.20.** *Пусть* $d \in D$, $a \in A$. *Отображение*

$$(5.1.16) \quad a \oplus d : b \in A \to (a \oplus d)b \in A$$

*определённое равенством*

$$(5.1.17) \qquad (a \oplus d)b = ab + db$$

*является эндоморфизмом абелевой группы* $A$.

---

**Лемма 5.1.18.** *Пусть* $d \in D$. *Отображение*

$$(5.1.14) \qquad d : b \in A \to bd \in A$$

*является эндоморфизмом абелевой группы* $A$.

**Доказательство.** Утверждение

$$bd \in A$$

является следствием определений 4.1.3, 5.1.1. Согласно равенству (4.1.22) отображение (5.1.14) является эндоморфизмом абелевой группы $A$. $\qquad\square$

**Лемма 5.1.19.** *Пусть* $a \in A$. *Отображение*

$$(5.1.15) \qquad a : b \in A \to ba \in A$$

*является эндоморфизмом абелевой группы* $A$.

**Доказательство.** Утверждение

$$ba \in A$$

является следствием определения 5.1.1. Согласно теореме 5.1.2 отображение (5.1.15) является эндоморфизмом абелевой группы $A$. $\qquad\square$

**Лемма 5.1.21.** *Пусть* $d \in D$, $a \in A$. *Отображение*

$$(5.1.18) \quad a \oplus d : b \in A \to b(a \oplus d) \in A$$

*определённое равенством*

$$(5.1.19) \qquad b(a \oplus d) = ba + bd$$

*является эндоморфизмом абелевой группы* $A$.

Доказательство. Так как $A$ является абелевой группой, то утверждение

$$db + ab \in A \quad d \in D \quad a \in A$$

является следствием лемм 5.1.16, 5.1.17. Следовательно, для любого $D$-числа $d$ и любого $D$-числа $a$, мы определили отображение (5.1.16). Согласно леммам 5.1.16, 5.1.17 и теореме 3.3.4, отображение (5.1.16) является эндоморфизмом абелевой группы $V$. □

Доказательство. Так как $A$ является абелевой группой, то утверждение

$$bd + ba \in A \quad d \in D \quad a \in A$$

является следствием лемм 5.1.18, 5.1.19. Следовательно, для любого $D$-числа $d$ и любого $D$-числа $a$, мы определили отображение (5.1.18). Согласно леммам 5.1.18, 5.1.19 и теореме 3.3.4, отображение (5.1.18) является эндоморфизмом абелевой группы $V$. □

Теорема 5.1.22. *Пусть $A$ - ассоциативная $D$-алгебра. Множество отображений $a \oplus d$ порождает[5.7]$D$-алгебру $A_{(1)}$ где сложение определено равенством*

$$(5.1.20) \qquad (a_1 \oplus d_1) + (a_2 \oplus d_2) = (a_1 + a_2) \oplus (d_1 + d_2)$$

*и произведение определено равенством*

$$(5.1.21) \qquad (a_1 \oplus d_1) \circ (a_2 \oplus d_2) = (a_1 a_2 + d_2 a_1 + d_1 a_2) \oplus (d_1 d_2)$$

5.1.22.1: *Если $D$-алгебра $A$ имеет единицу, то $D \subseteq A$, $A_{(1)} = A$.*
5.1.22.2: *Если $D$-алгебра $A$ является идеалом $D$, то $A \subseteq D$, $A_{(1)} = D$.*
5.1.22.3: *В противном случае $A_{(1)} = A \oplus D$.*
5.1.22.4: *$D$-алгебра $A$ является идеалом $D$-алгебры $A_{(1)}$.*

*$D$-алгебра $A_{(1)}$ называется* **унитальным расширением** *$D$-алгебры $A$.*

Доказательство. Выражение для операций на множестве $A_{(1)}$ не зависит от того, действует ли эндоморфизм $a \oplus d$ на $A$-числа слева или справа. Однако текст доказательства слегка отличается. Поэтому мы рассмотрим оба варианта доказательства.

Согласно теореме 3.3.4, равенство

$$(5.1.22) \quad \begin{aligned} &((a_1 \oplus d_1) + (a_2 \oplus d_2))b \\ =&(a_1 \oplus d_1)b + (a_2 \oplus d_2)b \\ =&a_1 b + d_1 b + a_2 b + d_2 b \\ =&a_1 b + a_2 b + d_1 b + d_2 b \\ =&(a_1 + a_2)b + (d_1 + d_2)b \\ =&(a_1 + a_2)b \oplus (d_1 + d_2)b \end{aligned}$$

является следствием равенства (5.1.17). Равенство (5.1.20) является следствием равенства (5.1.22).

Согласно теореме 3.3.4, равенство

$$(5.1.27) \quad \begin{aligned} &b((a_1 \oplus d_1) + (a_2 \oplus d_2)) \\ =&b(a_1 \oplus d_1) + b(a_2 \oplus d_2) \\ =&ba_1 + bd_1 + ba_2 + bd_2 \\ =&ba_1 + ba_2 + bd_1 + bd_2 \\ =&b(a_1 + a_2) + b(d_1 + d_2) \\ =&b(a_1 + a_2) \oplus (d_1 + d_2) \end{aligned}$$

является следствием равенства (5.1.19). Равенство (5.1.20) является следствием равенства (5.1.27).

---

[5.7] Смотри определение унитального расширения также на страницах [7]-52, [8]-64.

Согласно теореме 2.1.10, равенство

$$((a_1 \oplus d_1) \circ (a_2 \oplus d_2))b$$

$(5.1.23)$
$$\begin{aligned}
&= (a_1 \oplus d_1)((a_2 \oplus d_2)b)\\
&= (a_1 \oplus d_1)(a_2 b + d_2 b)\\
&= (a_1 \oplus d_1)(a_2 b)\\
&+ (a_1 \oplus d_1)(d_2 b)\\
&= a_1(a_2 b) + d_1(a_2 b)\\
&+ a_1(d_2 b) + d_1(d_2 b)
\end{aligned}$$

является следствием равенства (5.1.17). Равенство

$$((a_1 \oplus d_1) \circ (a_2 \oplus d_2))b$$

$(5.1.24)$
$$\begin{aligned}
&= (a_1 a_2)b + (d_1 a_2)b\\
&+ (d_2 a_1)b + (d_1 d_2)b
\end{aligned}$$

является следствием определений 5.1.1, 5.1.5, равенства (5.1.23), леммы 5.1.17. Равенство

$$((a_1 \oplus d_1) \circ (a_2 \oplus d_2))b$$
$$= (a_1 a_2 + d_1 a_2 + d_2 a_1 + d_1 d_2)b$$
$$= ((a_1 a_2 + d_1 a_2 + d_2 a_1) \oplus d_1 d_2)b$$

а следовательно и равенство (5.1.21), является следствием равенств (5.1.17), (5.1.24) и лемм 5.1.16, 5.1.17.

Пусть $D$-алгебра $A$ имеет единицу $e$. Тогда мы можем отождествить $d \in D$ и $de \in A$. Равенство

$(5.1.25)$ $\quad (a \oplus d)b = ab + db = (a + de)b$

является следствием равенства

$\boxed{(5.1.17) \quad (a \oplus d)b = ab + db}$

Утверждение 5.1.22.1 является следствием равенства (5.1.25).

Пусть $D$-алгебра $A$ является идеалом $D$. Тогда $A \subseteq D$. Равенство

$(5.1.26)$ $\quad (a \oplus d)b = ab + db = (a + d)b$

является следствием равенства (5.1.17). Утверждение 5.1.22.2 является следствием равенства (5.1.26).

Утверждение 5.1.22.4 является следствием определения 5.1.14 и лемм 5.1.20, 5.1.21.

---

Согласно теореме 2.1.10, равенство

$$b((a_2 \oplus d_2) \circ (a_1 \oplus d_1))$$

$(5.1.28)$
$$\begin{aligned}
&= (b(a_2 \oplus d_2))(a_1 \oplus d_1)\\
&= (ba_2 + bd_2)(a_1 \oplus d_1)\\
&= (ba_2)(a_1 \oplus d_1)\\
&+ (bd_2)(a_1 \oplus d_1)\\
&= (ba_2)a_1 + (ba_2)d_1\\
&+ (bd_2)a_1 + (bd_2)d_1
\end{aligned}$$

является следствием равенства (5.1.19). Равенство

$$b((a_2 \oplus d_2) \circ (a_1 \oplus d_1))$$

$(5.1.29)$
$$\begin{aligned}
&= b(a_2 a_1) + b(a_2 d_1)\\
&+ b(a_1 d_2) + b(d_2 d_1)
\end{aligned}$$

является следствием определений 5.1.1, 5.1.5, равенства (5.1.28), леммы 5.1.19. Равенство

$$b((a_2 \oplus d_2) \circ (a_1 \oplus d_1))$$
$$= b(a_2 a_1 + a_2 d_1 + a_1 d_2 + d_2 d_1)$$
$$= b((a_2 a_1 + a_2 d_1 + a_1 d_2) \oplus d_2 d_1)$$

а следовательно и равенство (5.1.21), является следствием равенств (5.1.19), (5.1.29) и лемм 5.1.18, 5.1.19.

Пусть $D$-алгебра $A$ имеет единицу $e$. Тогда мы можем отождествить $d \in D$ и $ed \in A$. Равенство

$(5.1.30)$ $\quad b(a \oplus d) = ba + bd = b(a + ed)$

является следствием равенства

$\boxed{(5.1.19) \quad b(a \oplus d) = ba + bd}$

Утверждение 5.1.22.1 является следствием равенства (5.1.30).

Пусть $D$-алгебра $A$ является идеалом $D$. Тогда $A \subseteq D$. Равенство

$(5.1.31)$ $\quad b(a \oplus d) = ba + bd = b(a + d)$

является следствием равенства (5.1.19). Утверждение 5.1.22.2 является следствием равенства (5.1.31).

□

---

Теорема 5.1.23. *Пусть $A$ - некоммутативная $D$-алгебра. Для любого $b \in A$ существует подалгебра $Z(A, b)$ $D$-алгебры $A$ такая, что*

$(5.1.32)$ $\qquad c \in Z(A, b) \Leftrightarrow cb = bc$

$D$-*алгебра* $Z(A, b)$ *называется* **центром** $A$-**числа** $b$.

Доказательство. Подмножество $Z(A, b)$ не пусто так как $0 \in Z(A, b)$, $b \in Z(A, b)$. Пусть $d_1, d_2 \in D$, $c_1, c_2 \in Z(A, b)$. Тогда

$$(d_1 c_1 + d_2 c_2)b = d_1 c_1 b + d_2 c_2 b = b d_1 c_1 + b d_2 c_2$$
$$= b(d_1 c_1 + d_2 c_2)$$

Откуда следует $d_1 c_1 + d_2 c_2 \in Z(A, b)$. Следовательно, множество $Z(A, b)$ является $D$-модулем.

Пусть $c_1, c_2 \in Z(A, b)$. Тогда

(5.1.33)          $b(c_1 c_2) = (b c_1) c_2 = (c_1 b) c_2 = c_1 (b c_2) = c_1 (c_2 b) = (c_1 c_2)b$

Из равенства (5.1.33) следует, что $c_1 c_2 \in Z(A, b)$. Следовательно, $D$-модуль $Z(A, b)$ является $D$-алгеброй.                                                         □

Лемма 5.1.24. *Пусть $A$ - некоммутативная $D$-алгебра. Для любого $a \in A$, если $p, c \in Z(A, a)$, то*

(5.1.34)                          $p c^n \in Z(A, a)$

Доказательство. Согласно теореме 5.1.23,

$$p = p c^0 \in Z(A, a)$$

Пусть утверждение (5.1.34) верно для $n = k - 1$, $k > 0$,

(5.1.35)                          $p c^{k-1} \in Z(A, a)$

Тогда утверждение

$$p c^k = p c^{k-1} c \in Z(A, a)$$

является следсвием утверждения (5.1.35) и теоремы 5.1.23. Согласно принципу математической индукции, утверждение (5.1.34) верно для любого $n \geq 0$.     □

Теорема 5.1.25. *Пусть $A$ - некоммутативная $D$-алгебра. Для любого $a \in A$, если $c \in Z(A, a)$, то*

(5.1.36)                          $p(c) \in Z(A, a)$

*для любого многочлена*

(5.1.37)                          $p(x) = p_0 + p_1 x + ... + p_n x^n$
$$p_0, \ ..., \ p_n \in D$$

Доказательство. Теорема является следствием леммы 5.1.24 и теоремы 5.1.23.                                                         □

Теорема 5.1.26. *Если $a \in Z(A, b)$, то $b \in Z(A, a)$.*

Доказательство. Пусть $a \in Z(A, b)$. Равенство

(5.1.38)                          $ab = ba$

является следствием утверждения (5.1.32). Теорема является следстием равенства (5.1.38).                                                         □

Теорема 5.1.27. *Если* $a \in Z(A, b)$ *и* $c \in Z(A)$, *то* $a \in Z(A, bc)$.

Доказательство. Равенство

$$(5.1.39) \qquad a(bc) = (ab)c = (ba)c = b(ac) = b(ca) = (bc)a$$

является следствием определений 5.1.8, 5.1.9 и теоремы 5.1.23. Теоремя является следствием равенства (5.1.39) и теоремы 5.1.23. □

Определение 5.1.28. *$D$-алгебра $A$ называется* **алгеброй с делением**, *если для любого $A$-числа $a \neq 0$ существует $A$-число $a^{-1}$.* □

Теорема 5.1.29. *Пусть $D$-алгебра $A$ является алгеброй с делением. $D$-алгебра $A$ имеет единицу.*

Доказательство. Теорема является следствием утверждения, что уравнение

$$ax = a$$

имеет решение для любого $a \in A$. □

Теорема 5.1.30. *Пусть $D$-алгебра $A$ является алгеброй с делением. Кольцо $D$ является полем и подмножеством центра $D$-алгебры $A$.*

Доказательство. Пусть $e$ - единица $D$-алгебры $A$. Так как отображение

$$d \in D \rightarrow d\, e \in A$$

является вложением кольца $D$ в $D$-алгебру $A$, то кольцо $D$ является подмножеством центра $D$-алгебры $A$. □

Теорема 5.1.31. *Пусть $A$ - ассоциативная $D$-алгебра с делением. Из утверждения*

$$(5.1.40) \qquad ab = ac \qquad a \neq 0$$

*следует, что $b = c$.*

Доказательство. Равенство

$$(5.1.41) \qquad b = (a^{-1}a)b = a^{-1}(ab) = a^{-1}(ac) = (a^{-1}a)c = c$$

является следствием утверждения (5.1.41). □

Определение 5.1.32. *$A$-числа $a$ и $b$ называются подобными, если существует $A$-число $c$ такое, что*

$$(5.1.42) \qquad a = c^{-1}bc$$

□

Теорема 5.1.33. *Пусть $A$ - алгебра над коммутативным кольцом $D$. Для любого множества $B$ множество отображений $A^B$ является абелевой группой относительно операции*

$$(f + g)(x) = f(x) + g(x)$$

Доказательство. Теорема является следствием определений 4.1.3, 5.1.1.                                                                      □

## 5.2. Тип $D$-алгебры

Соглашение 5.2.1. *Согласно определениям 4.1.24, 5.1.1 свободная $D$-алгебра $A$ - это $D$-модуль $A$, который имеет базис $\overline{\overline{e}}$. Однако, вообще говоря, $\overline{\overline{e}}$ не является базисом $D$-алгебры $A$, так как в $D$-алгебре $A$ определено произведение. Например, в алгебре кватернионов произвольный кватернион*

$$a = a^0 + a^1 i + a^2 j + a^3 k$$

*может быть представлен в виде*

$$a = a^0 + a^1 i + a^2 j + a^3 ij$$

*Однако при решении большинства задач использование базиса $D$-модуля $A$ проще, чем использование базиса $D$-алгебры $A$. Например, легче определить координаты $H$-числа относительно базиса $(1, i, j, k)$ $R$-векторного пространства $H$, чем определить координаты $H$-числа относительно базиса $(1, i, j)$ $R$-алгебры $H$. Поэтому фраза "мы рассматриваем базис $D$-алгебры $A$" означает, что мы рассматриваем $D$-алгебру $A$ и базис $D$-модуля $A$.*                     □

Соглашение 5.2.2. *Пусть $A$ - свободная алгебра с конечным или счётным базисом. При разложении элемента алгебры $A$ относительно базиса $\overline{\overline{e}}$ мы пользуемся одной и той же корневой буквой для обозначения этого элемента и его координат. В выражении $a^2$ не ясно - это компонента разложения элемента $a$ относительно базиса или это операция возведения в степень. Для облегчения чтения текста мы будем индекс элемента алгебры выделять цветом. Например,*

$$a = a^i e_i$$
                                                                      □

$D$-алгебра $A$ наследует тип $D$-модуля $A$.

Определение 5.2.3. *$D$-алгебра $A$ называется $D$-алгеброй столбцов, если $D$-модуль $A$ является $D$-модулем столбцов.*                     □

Определение 5.2.5. *$D$-алгебра $A$ называется $D$-алгеброй строк, если $D$-модуль $A$ является $D$-модулем строк.*                     □

Соглашение 5.2.4. *Пусть $\overline{\overline{e}}$ - базис $D$-алгебры столбцов $A$. Если $D$-алгебра $A$ имеет единицу, положим $e_0$ - единица $D$-алгебры $A$.*                     □

Соглашение 5.2.6. *Пусть $\overline{\overline{e}}$ - базис $D$-алгебры строк $A$. Если $D$-алгебра $A$ имеет единицу, положим $e^0$ - единица $D$-алгебры $A$.*                     □

Теорема 5.2.7. *Пусть множество векторов* $\overline{\overline{e}} = (e_i, i \in I)$ *является базисом $D$-алгебры столбцов $A$. Пусть*

$$a = a^i e_i \quad b = b^i e_i \quad a, b \in A$$

*Произведение $a$, $b$ можно получить согласно правилу*

$$(5.2.1) \qquad (ab)^k = C^k_{ij} a^i b^j$$

*где* $C^k_{ij}$ *- **структурные константы** алгебры $A$ над кольцом $D$. Произведение базисных векторов в алгебре $A$ определено согласно правилу*

$$(5.2.2) \qquad e_i e_j = C^k_{ij} e_k$$

Доказательство. Равенство (5.2.2) является следствием утверждения, что $\overline{\overline{e}}$ является базисом $D$-алгебры $A$. Так как произведение в алгебре является билинейным отображением, то произведение $a$ и $b$ можно записать в виде

$$(5.2.3) \qquad ab = a^i b^j e_i e_j$$

Из равенств (5.2.2), (5.2.3), следует

$$(5.2.4) \qquad ab = a^i b^j C^k_{ij} e_k$$

Так как $\overline{\overline{e}}$ является базисом алгебры $A$, то равенство (5.2.1) следует из равенства (5.2.4). $\qquad\square$

Теорема 5.2.9. *Если алгебра $A$ коммутативна, то*

$$(5.2.9) \qquad C^p_{ij} = C^p_{ji}$$

*Если алгебра $A$ ассоциативна, то*

$$(5.2.10) \qquad C^p_{ij} C^q_{pk} = C^q_{ip} C^p_{jk}$$

Доказательство. Для коммутативной алгебры, равенство (5.2.9) следует из равенства

$$e_i e_j = e_j e_i$$

Для ассоциативной алгебры, равенство (5.2.10) следует из равенства

$$(e_i e_j) e_k = e_i (e_j e_k)$$

Теорема 5.2.8. *Пусть множество векторов* $\overline{\overline{e}} = (e^i, i \in I)$ *является базисом $D$-алгебры строк $A$. Пусть*

$$a = a_i e^i \quad b = b_i e^i \quad a, b \in A$$

*Произведение $a$, $b$ можно получить согласно правилу*

$$(5.2.5) \qquad (ab)_k = C^{ij}_k a_i b_j$$

*где* $C^{ij}_k$ *- **структурные константы** алгебры $A$ над кольцом $D$. Произведение базисных векторов в алгебре $A$ определено согласно правилу*

$$(5.2.6) \qquad e^i e^j = C^{ij}_k e^k$$

Доказательство. Равенство (5.2.6) является следствием утверждения, что $\overline{\overline{e}}$ является базисом $D$-алгебры $A$. Так как произведение в алгебре является билинейным отображением, то произведение $a$ и $b$ можно записать в виде

$$(5.2.7) \qquad ab = a_i b_j e^i e^j$$

Из равенств (5.2.6), (5.2.7), следует

$$(5.2.8) \qquad ab = a_i b_j C^{ij}_k e^k$$

Так как $\overline{\overline{e}}$ является базисом алгебры $A$, то равенство (5.2.5) следует из равенства (5.2.8). $\qquad\square$

Теорема 5.2.10. *Если алгебра $A$ коммутативна, то*

$$(5.2.11) \qquad C^{ij}_p = C^{ji}_p$$

*Если алгебра $A$ ассоциативна, то*

$$(5.2.12) \qquad C^{ij}_p C^{pk}_q = C^{ip}_q C^{jk}_p$$

Доказательство. Для коммутативной алгебры, равенство (5.2.11) следует из равенства

$$e^i e^j = e^j e^i$$

Для ассоциативной алгебры, равенство (5.2.12) следует из равенства

$$(e^i e^j) e^k = e^i (e^j e^k)$$

□  |                                                      □

## 5.3. Диаграмма представлений

Теорема 5.3.1. *Представление*

(5.3.1) $$g_{23} : A \dashrightarrow A$$

$D$-*модуля* $A$ *в* $D$-*модуле* $A$ *эквивалентно структуре* $D$-*алгебры* $A$.

Теорема 5.3.2. *Пусть* $D$ - *коммутативное кольцо и* $A$ - *абелева группа.* *Диаграмма представлений*

$$A \xrightarrow{g_{23}} A \qquad g_{12}(d) : v \to d\,v$$
$$\nwarrow{}_{g_{12}} \quad \nearrow{}_{g_{12}} \qquad g_{23}(v) : w \to C \circ (v, w)$$
$$D \qquad\qquad C \in \mathcal{L}(D; A^2 \to A)$$

*порождает структуру* $D$-*алгебры* $A$.

Доказательство. Структура $D$-модуля $A$ порождена эффективным представлением

$$g_{12} : D \dashrightarrow A$$

кольца $D$ в абелевой группе $A$.

Лемма 5.3.3. *Пусть в* $D$-*модуле* $A$ *определена структура* $D$-*алгебры* $A$, *порождённая произведением*

$$v\,w = C \circ (v, w)$$

**Левый сдвиг** $D$-**модуля** $A$, *определённый равенством*

(5.3.2) $$l \circ v : w \in A \to v\,w \in A$$

*порождает представление*

$$A \xrightarrow{g_{23}} A \qquad g_{23} : v \to l \circ v$$
$$\qquad\qquad g_{23} \circ v : w \to (l \circ v) \circ w$$

$D$-*модуля* $A$ *в* $D$-*модуле* $A$

Доказательство. Согласно определениям 5.1.1 и 4.4.1, левый сдвиг $D$-модуля $A$ является линейным отображением. Согласно определению 4.3.9, отображение $l \circ v$ является эндоморфизмом $D$-модуля $A$. Равенство

(5.3.3)   $$(l \circ (v_1 + v_2)) \circ w = (v_1 + v_2)w = v_1 w + v_2 w = (l \circ v_1) \circ w + (l \circ v_2) \circ w$$

является следствием определения 4.4.1 и равенства (5.3.2).      Согласно следствию 4.4.4,  равенство

(5.3.4) $$l \circ (v_1 + v_2) = l \circ v_1 + l \circ v_2$$

является следствием равенства (5.3.3). Равенство

(5.3.5) $$(l \circ (dv)) \circ w = (dv)w = d(vw) = d((l \circ v) \circ w)$$

является следствием определения 4.4.1 и равенства (5.3.2). Согласно следствию 4.4.4, равенство

$$(5.3.6) \qquad l \circ (dv) = d(l \circ v)$$

является следствием равенства (5.3.5). Лемма является следствием равенств (5.3.4), (5.3.6). ⊙

Лемма 5.3.4. *Представление*

$$A \xrightarrow{\ g_{23}\ } A \qquad \begin{aligned} g_{23} : v \to l \circ v \\ g_{23} \circ v : w \to (l \circ v) \circ w \end{aligned}$$

$D$-*модуля* $A$ *в* $D$-*модуле* $A$ *определяет произведение в* $D$-*модуле* $A$ *согласно правилу*

$$ab = (g_{23} \circ a) \circ b$$

Доказательство. Поскольку отображение $g_{23} \circ v$ является эндоморфизмом $D$-модуля $A$, то

$$(5.3.7) \qquad \begin{aligned} (g_{23} \circ v)(w_1 + w_2) = (g_{23} \circ v) \circ w_1 + (g_{23} \circ v) \circ w_2 \\ (g_{23} \circ v) \circ (dw) = d((g_{23} \circ v) \circ w) \end{aligned}$$

Поскольку отображение $g_{23}$ является линейным отображением

$$g_{23} : A \to \mathcal{L}(D; A \to A)$$

то, согласно следствиям 4.4.4, 4.4.6,

$$(5.3.8) \quad (g_{23} \circ (v_1 + v_2)) \circ w = (g_{23} \circ v_1 + g_{23} \circ v_2)(w) = (g_{23} \circ v_1) \circ w + (g_{23} \circ v_2) \circ w$$

$$(5.3.9) \qquad (g_{23} \circ (d\,v)) \circ w = (d\,(g_{23} \circ v)) \circ w = d\,((g_{23} \circ v) \circ w)$$

Из равенств (5.3.7), (5.3.8), (5.3.9) и определения 4.4.1, следует, что отображение $g_{23}$ является билинейным отображением. Следовательно, отображение $g_{23}$ определяет произведение в $D$-модуле $A$ согласно правилу

$$ab = (g_{23} \circ a) \circ b$$

⊙

Теорема является следствием лемм 5.3.3, 5.3.4. □

## 5.4. Гомоморфизм $D$-алгебры

### 5.4.1. Общее определение.

Определение 5.4.1. *Пусть* $A_i$, $i = 1, 2$, *- алгебра над коммутативным кольцом* $D_i$. *Пусть*

$$(5.4.1) \qquad \begin{aligned} A_1 \xrightarrow{\ h_{1.23}\ } A_1 \qquad & h_{1.12}(d) : v \to d\,v \\ \underset{h_{1.12}}{\swarrow \ast} \quad \underset{\ast \ h_{1.12}}{\searrow} \qquad & h_{1.23}(v) : w \to C_1 \circ (v, w) \\ D_1 \qquad\qquad & C_1 \in \mathcal{L}(D_1; A_1^2 \to A_1) \end{aligned}$$

*диаграмма представлений, описывающая $D_1$-алгебру $A_1$. Пусть*

$$(5.4.2)$$

$$h_{2.12}(d) : v \to d\, v$$
$$h_{2.23}(v) : w \to C_2 \circ (v, w)$$
$$C_2 \in \mathcal{L}(D_2; A_2^2 \to A_2)$$

*диаграмма представлений, описывающая $D_2$-алгебру $A_2$. Морфизм*

$$(5.4.3) \qquad (h : D_1 \to D_2 \quad f : A_1 \to A_2)$$

*диаграммы представлений* (5.4.1) *в диаграмму представлений* (5.4.2) *называется* **гомоморфизмом** *$D_1$-алгебры $A_1$ в $D_2$-алгебру $A_2$. Обозначим* $\mathrm{Hom}(D_1 \to D_2; A_1 \to A_2)$ *множество гомоморфизмов $D_1$-алгебры $A_1$ в $D_2$-алгебру $A_2$.* □

Теорема 5.4.2. *Гомоморфизм*

$$\boxed{(5.4.3) \quad (h : D_1 \to D_2 \quad f : A_1 \to A_2)}$$

*$D_1$-алгебры $A_1$ в $D_2$-алгебру $A_2$ - это линейное отображение* (5.4.3) *$D_1$-модуля $A_1$ в $D_2$-модуль $A_2$ такое, что*

$$(5.4.4) \qquad f \circ (ab) = (f \circ a)(f \circ b)$$

*и удовлетворяет следующим равенствам*

$$(5.4.5) \qquad h(p + q) = h(p) + h(q)$$

$$(5.4.6) \qquad h(pq) = h(p)h(q)$$

$$(5.4.7) \qquad f \circ (a + b) = f \circ a + f \circ b$$

$$\boxed{(5.4.4) \quad f \circ (ab) = (f \circ a)(f \circ b)}$$

$$(5.4.8) \qquad f \circ (pa) = h(p)(f \circ a)$$

$$p, q \in D_1 \quad a, b \in A_1$$

Доказательство. Согласно определениям 2.3.9, 2.8.5, отображение (5.4.3) являетля морфизмом представления $h_{1.12}$, описывающего $D_1$-модуль $A_1$, в представления $h_{2.12}$, описывающего $D_2$-модуль $A_2$. Следовательно, отображение $f$ является линейным отображением $D_1$-модуля $A_1$ в $D_2$-модуль $A_2$ и равенства (5.4.5), (5.4.6), (5.4.7), (5.4.8) являются следствием теоремы 4.3.2.

Согласно определениям 2.3.9, 2.8.5, отображение $f$ являетля морфизмом представления $h_{1.23}$ в представления $h_{2.23}$. Равенство

$$(5.4.9) \qquad f \circ (h_{1.23}(a)(b)) = h_{2.23}(f \circ a)(f \circ b)$$

следует из равенства (2.3.2). Из равенств (5.4.1), (5.4.2), (5.4.9), следует

$$(5.4.10) \qquad f \circ (C_1 \circ (a, b)) = C_2 \circ (f \circ a, f \circ b)$$

Равенство (5.4.4) следует из равенств (5.1.1), (5.4.10). □

---

**Теорема 5.4.3.** *Гомоморфизм* [5.8]

$$(h : D_1 \to D_2 \quad \overline{f} : A_1 \to A_2)$$

*$D_1$-алгебры столбцов $A_1$ в $D_2$-алгебру столбцов $A_2$ имеет представление*

(5.4.11) $$b = fh(a)$$

(5.4.12) $$\overline{f} \circ (a^i e_{1i}) = h(a^i) f_i^k e_{2k}$$

(5.4.13) $$\overline{f} \circ (e_1 a) = e_2 f h(a)$$

*относительно заданных базисов. Здесь*

- *$a$ - координатная матрица $A_1$-числа $\overline{a}$ относительно базиса $\overline{\overline{e}}_1$ (D)*

(5.4.14) $$\overline{a} = e_1 a$$

- *$h(a) = (h(a_i), i \in I)$ - матрица $D_2$-чисел.*
- *$b$ - координатная матрица $A_2$-числа*

$$\overline{b} = \overline{f} \circ \overline{a}$$

*относительно базиса $\overline{\overline{e}}_2$*

(5.4.15) $$\overline{b} = e_2 b$$

- *$f$ - координатная матрица множества $A_2$-чисел ($\overline{f} \circ e_{1i}, i \in I$) относительно базиса $\overline{\overline{e}}_2$.*
- *Матрица гомоморфизма и структурные константы связаны соотношением*

(5.4.16) $$h(C_{1ij}^{\;k}) f_k^l = f_i^p f_j^q C_{2pq}^{\;l}$$

**Доказательство.** Равенства

| (5.4.11) $b = fh(a)$ |
|---|

| (5.4.12) $\overline{f} \circ (a^i e_{1i}) = h(a^i) f_i^k e_{2k}$ |
|---|

---

**Теорема 5.4.4.** *Гомоморфизм* [5.9]

$$(h : D_1 \to D_2 \quad \overline{f} : A_1 \to A_2)$$

*$D_1$-алгебры строк $A_1$ в $D_2$-алгебру строк $A_2$ имеет представление*

(5.4.22) $$b = h(a)f$$

(5.4.23) $$\overline{f} \circ (a_i e_1^i) = h(a_i) f_k^i e_2^k$$

(5.4.24) $$\overline{f} \circ (ae_1) = h(a) f e_2$$

*относительно заданных базисов. Здесь*

- *$a$ - координатная матрица $A_1$-числа $\overline{a}$ относительно базиса $\overline{\overline{e}}_1$ (D)*

(5.4.25) $$\overline{a} = ae_1$$

- *$h(a) = (h(a^i), i \in I)$ - матрица $D_2$-чисел.*
- *$b$ - координатная матрица $A_2$-числа*

$$\overline{b} = \overline{f} \circ \overline{a}$$

*относительно базиса $\overline{\overline{e}}_2$*

(5.4.26) $$\overline{b} = be_2$$

- *$f$ - координатная матрица множества $A_2$-чисел ($\overline{f} \circ e_1^i, i \in I$) относительно базиса $\overline{\overline{e}}_2$.*
- *Матрица гомоморфизма и структурные константы связаны соотношением*

(5.4.27) $$h(C_{1k}^{\;ij}) f_l^k = f_p^i f_q^j C_{2l}^{\;pq}$$

**Доказательство.** Равенства

| (5.4.22) $b = h(a)f$ |
|---|

| (5.4.23) $\overline{f} \circ (a_i e_1^i) = h(a_i) f_k^i e_2^k$ |
|---|

---

[5.8] В теоремах 5.4.3, 5.4.5, мы опираемся на следующее соглашение. Пусть множество векторов $\overline{\overline{e}}_1 = (e_{1i}, i \in I)$ является базисом и $C_{1ij}^{\;k}$, $k$, $i, j \in I$, - структурные константы $D_1$-алгебры столбцов $A_1$. Пусть множество векторов $\overline{\overline{e}}_2 = (e_{2j}, j \in J)$ является базисом и $C_{2ij}^{\;k}$, $k$, $i, j \in J$, - структурные константы $D_2$-алгебры столбцов $A_2$.

[5.9] В теоремах 5.4.4, 5.4.6, мы опираемся на следующее соглашение. Пусть множество векторов $\overline{\overline{e}}_1 = (e_1^i, i \in I)$ является базисом и $C_{1k}^{\;ij}$, $k$, $i, j \in I$, - структурные константы $D_1$-алгебры строк $A_1$. Пусть множество векторов $\overline{\overline{e}}_2 = (e_2^j, j \in J)$ является базисом и $C_{2k}^{\;ij}$, $k$, $i, j \in J$, - структурные константы $D_2$-алгебры строк $A_2$.

$(5.4.13)$　$\boxed{\ \overline{f}\circ(e_1 a)=e_2 fh(a)\ }$

являются следствием теорем 4.3.3, 5.4.2.

Равенство

$(5.4.17)$　$\overline{f}\circ(e_{1i}e_{1j})=h(C_{1\,ij}^{\ \ k})\overline{f}\circ e_{1k}$

является следствием равенств

$\boxed{(5.2.2)\quad e_i e_j=C_{ij}^k e_k}$

$\boxed{(5.4.8)\quad f\circ(pa)=h(p)(f\circ a)}$

Равенство

$(5.4.18)$　$\overline{f}\circ(e_{1i}e_{1j})=h(C_{1\,ij}^{\ \ k})f_k^l e_{2l}$

является следствием равенства

$\boxed{(4.3.19)\quad \overline{f}\circ e_{1i}=e_2 f_i=e_{2j}f_i^j}$

и равенства $(5.4.17)$. Равенство

$(5.4.19)$
$$\begin{aligned}&\overline{f}\circ(e_{1i}e_{1j})\\&=(\overline{f}\circ e_{1i})(\overline{f}\circ e_{1j})\\&=f_i^p e_{2p}f_j^q e_{2q}\end{aligned}$$

является следствием равенства

$\boxed{(5.4.4)\quad f\circ(ab)=(f\circ a)(f\circ b)}$

Равенство

$(5.4.20)$　$\overline{f}\circ(e_{1i}e_{2j})=f_i^p f_j^q C_{2\,pq}^{\ \ l}e_{2l}$

является следствием равенства

$\boxed{(5.2.2)\quad e_i e_j=C_{ij}^k e_k}$

и равенства $(5.4.19)$. Равенство

$(5.4.21)$　$h(C_{1\,ij}^{\ \ k})f_k^l e_{2l}=f_i^p f_j^q C_{2\,pq}^{\ \ l}e_{2l}$

является следствием равенств $(5.4.18)$, $(5.4.20)$. Равенство

$\boxed{(5.4.16)\quad h(C_{1\,ij}^{\ \ k})f_k^l=f_i^p f_j^q C_{2\,pq}^{\ \ l}}$

следует из равенства $(5.4.21)$ и теоремы 4.2.7.　□

---

$(5.4.24)$　$\boxed{\ \overline{f}\circ(ae_1)=h(a)fe_2\ }$

являются следствием теорем 4.3.4, 5.4.2.

Равенство

$(5.4.28)$　$\overline{f}\circ(e_1^i e_1^j)=h(C_{1\ k}^{ij})\overline{f}\circ e_1^k$

является следствием равенств

$\boxed{(5.2.6)\quad e^i e^j=C_k^{ij}e^k}$

$\boxed{(5.4.8)\quad f\circ(pa)=h(p)(f\circ a)}$

Равенство

$(5.4.29)$　$\overline{f}\circ(e_1^i e_1^j)=h(C_{1\ k}^{ij})f_k^l e_2^l$

является следствием равенства

$\boxed{(4.3.22)\quad \overline{f}\circ e_1^i=f_i e_2=f_j^i e_2^j}$

и равенства $(5.4.28)$. Равенство

$(5.4.30)$
$$\begin{aligned}&\overline{f}\circ(e_1^i e_1^j)\\&=(\overline{f}\circ e_1^i)(\overline{f}\circ e_1^j)\\&=f_p^i e_2^p f_q^j e_2^q\end{aligned}$$

является следствием равенства

$\boxed{(5.4.4)\quad f\circ(ab)=(f\circ a)(f\circ b)}$

Равенство

$(5.4.31)$　$\overline{f}\circ(e_1^i e_2^j)=f_p^i f_q^j C_{2\ l}^{pq}e_2^l$

является следствием равенства

$\boxed{(5.2.6)\quad e^i e^j=C_k^{ij}e^k}$

и равенства $(5.4.30)$. Равенство

$(5.4.32)$　$h(C_{1\ k}^{ij})f_k^l e_2^l=f_p^i f_q^j C_{2\ l}^{pq}e_2^l$

является следствием равенств $(5.4.29)$, $(5.4.31)$. Равенство

$\boxed{(5.4.27)\quad h(C_{1\ k}^{ij})f_k^l=f_p^i f_q^j C_{2\ l}^{pq}}$

следует из равенства $(5.4.32)$ и теоремы 4.2.14.　□

---

**Теорема 5.4.5.** *Пусть отображение*

$$h:D_1\to D_2$$

*является гомоморфизмом кольца $D_1$ в кольцо $D_2$. Пусть*

$$f=(f_j^i,\ i\in I, j\in J)$$

*матрица $D_2$-чисел, которая удовлетворяют равенству*

$\boxed{(5.4.16)\quad h(C_{1\,ij}^{\ \ k})f_k^l=f_i^p f_j^q C_{2\,pq}^{\ \ l}}$

*Отображение*[5.8]

---

**Теорема 5.4.6.** *Пусть отображение*

$$h:D_1\to D_2$$

*является гомоморфизмом кольца $D_1$ в кольцо $D_2$. Пусть*

$$f=(f_i^j,\ i\in I, j\in J)$$

*матрица $D_2$-чисел, которая удовлетворяют равенству*

$\boxed{(5.4.27)\quad h(C_{1\ k}^{ij})f_k^l=f_p^i f_q^j C_{2\ l}^{pq}}$

*Отображение*[5.9]

$(4.3.6)$   $(h : D_1 \to D_2$   $\overline{f} : V_1 \to V_2)$

*определённое равенством*

$(4.3.7)$   $\overline{f} \circ (e_1 a) = e_2 f h(a)$

*является гомоморфизмом $D_1$-алгебры столбцов $A_1$ в $D_2$-алгебру столбцов $A_2$. Гомоморфизм $(4.3.6)$, который имеет данную матрицу $f$, определён однозначно.*

Доказательство. Согласно теореме 4.3.5, отображение $(4.3.6)$ является гомоморфизмом $D_1$-модуля $A_1$ в $D_2$-модуль $A_2$ и гомоморфизм $(4.3.6)$ определён однозначно.

Пусть $a = a^i e_{1i}$, $b = b^i e_{1i}$. Равенство

$$(5.4.33) \quad \begin{aligned} &h(C_{1\,ij}^{\;\;k}) f_k^l h(a^i) h(b^j) e_{2l} \\ &= f_i^p f_j^q C_{2\,pq}^{\;\;l} h(a^i) h(b^j) e_{2l} \end{aligned}$$

является следствием равенства $(5.4.16)$. Так как отображение $h$ является гомоморфизмом кольца $D_1$ в кольцо $D_2$, то

$$(5.4.34) \quad \begin{aligned} &h(C_{1\,ij}^{\;\;k}) f_k^l h(a^i) h(b^j) e_{2l} \\ &= h(C_{1\,ij}^{\;\;k} a^i b^j) f_k^l e_{2l} \end{aligned}$$

Равенство

$$(5.4.35) \quad \begin{aligned} &h(C_{1\,ij}^{\;\;k}) f_k^l h(a^i) h(b^j) e_{2l} \\ &= h((ab)^k) f_k^l e_{2l} \end{aligned}$$

является следствием равенства $(5.4.34)$ и равенства

$(5.2.1)$   $(ab)^k = C_{ij}^k a^i b^j$

Равенство

$$(5.4.36) \quad \begin{aligned} &h(C_{1\,ij}^{\;\;k}) f_k^l h(a^i) h(b^j) e_{2l} \\ &= f \circ (ab) \end{aligned}$$

является следствием равенства $(5.4.35)$ и равенства

$(5.4.12)$   $\overline{f} \circ (a^i e_{1i}) = h(a^i) f_i^k e_{2k}$

Равенство

$$(5.4.37) \quad \begin{aligned} &f_i^p f_j^q C_{2\,pq}^{\;\;l} h(a^i) h(b^j) e_{2l} \\ &= (h(a^i) f_i^p e_{2p})(h(b^j) f_j^q e_{2q}) \end{aligned}$$

является следствием равенства

$(5.2.2)$   $e_i e_j = C_{ij}^k e_k$

$(4.3.12)$   $(h : D_1 \to D_2$   $\overline{f} : V_1 \to V_2)$

*определённое равенством*

$(4.3.13)$   $\overline{f} \circ (ae_1) = h(a) f e_2$

*является гомоморфизмом $D_1$-алгебры строк $A_1$ в $D_2$-алгебру строк $A_2$. Гомоморфизм $(4.3.12)$, который имеет данную матрицу $f$, определён однозначно.*

Доказательство. Согласно теореме 4.3.6, отображение $(4.3.12)$ является гомоморфизмом $D_1$-модуля $A_1$ в $D_2$-модуль $A_2$ и гомоморфизм $(4.3.12)$ определён однозначно.

Пусть $a = a_i e_1^i$, $b = b_i e_1^i$. Равенство

$$(5.4.39) \quad \begin{aligned} &h(C_{1\,k}^{\;\;ij}) f_l^k h(a_i) h(b_j) e_2^l \\ &= f_p^i f_q^j C_{2l}^{\;\;pq} h(a_i) h(b_j) e_2^l \end{aligned}$$

является следствием равенства $(5.4.27)$. Так как отображение $h$ является гомоморфизмом кольца $D_1$ в кольцо $D_2$, то

$$(5.4.40) \quad \begin{aligned} &h(C_{1\,k}^{\;\;ij}) f_l^k h(a_i) h(b_j) e_2^l \\ &= h(C_{1\,k}^{\;\;ij} a_i b_j) f_l^k e_2^l \end{aligned}$$

Равенство

$$(5.4.41) \quad \begin{aligned} &h(C_{1\,k}^{\;\;ij}) f_l^k h(a_i) h(b_j) e_2^l \\ &= h((ab)_k) f_l^k e_2^l \end{aligned}$$

является следствием равенства $(5.4.40)$ и равенства

$(5.2.5)$   $(ab)_k = C_k^{ij} a_i b_j$

Равенство

$$(5.4.42) \quad \begin{aligned} &h(C_{1\,k}^{\;\;ij}) f_l^k h(a_i) h(b_j) e_2^l \\ &= f \circ (ab) \end{aligned}$$

является следствием равенства $(5.4.41)$ и равенства

$(5.4.23)$   $\overline{f} \circ (a_i e_1^i) = h(a_i) f_k^i e_2^k$

Равенство

$$(5.4.43) \quad \begin{aligned} &f_p^i f_q^j C_{2l}^{\;\;pq} h(a_i) h(b_j) e_2^l \\ &= (h(a_i) f_p^i e_2^p)(h(b_j) f_q^j e_2^q) \end{aligned}$$

является следствием равенства

$(5.2.6)$   $e^i e^j = C_k^{ij} e^k$

Равенство

$$(5.4.38) \quad f_i^p f_j^q C_{2\,pq}{}^l h(a^i) h(b^j) e_{2l}$$
$$= (f \circ a)(f \circ b)$$

является следствием равенства (5.4.37) и равенства

$$\boxed{(5.4.12) \quad \overline{f} \circ (a^i e_{1i}) = h(a^i) f_i^k e_{2k}}$$

Равенство

$$\boxed{(5.4.4) \quad f \circ (ab) = (f \circ a)(f \circ b)}$$

является следствием равенств (5.4.33), (5.4.36), (5.4.38). Согласно теореме 5.4.2, отображение (4.3.6) является гомоморфизмом $D_1$-алгебры $A_1$ в $D_2$-алгебру $A_2$. $\quad\square$

Равенство

$$(5.4.44) \quad f_p^i f_q^j C_{2l}{}^{pq} h(a_i) h(b_j) e_2^l$$
$$= (f \circ a)(f \circ b)$$

является следствием равенства (5.4.43) и равенства

$$\boxed{(5.4.23) \quad \overline{f} \circ (a_i e_1^i) = h(a_i) f_k^i e_2^k}$$

Равенство

$$\boxed{(5.4.4) \quad f \circ (ab) = (f \circ a)(f \circ b)}$$

является следствием равенств (5.4.39), (5.4.42), (5.4.44). Согласно теореме 5.4.2, отображение (4.3.12) является гомоморфизмом $D_1$-алгебры $A_1$ в $D_2$-алгебру $A_2$. $\quad\square$

## 5.4.2. Гомоморфизм, когда кольца $D_1 = D_2 = D$.

Определение 5.4.7. *Пусть $A_i$, $i = 1, 2$, - алгебра над коммутативным кольцом $D$. Пусть*

$$(5.4.45)$$

$$h_{1.12}(d) : v \to d\,v$$
$$h_{1.23}(v) : w \to C_1 \circ (v, w)$$
$$C_1 \in \mathcal{L}(D; A_1^2 \to A_1)$$

*диаграмма представлений, описывающая $D_1$-алгебру $A_1$. Пусть*

$$(5.4.46)$$

$$h_{2.12}(d) : v \to d\,v$$
$$h_{2.23}(v) : w \to C_2 \circ (v, w)$$
$$C_2 \in \mathcal{L}(D; A_2^2 \to A_2)$$

*диаграмма представлений, описывающая $D_2$-алгебру $A_2$. Морфизм*

$$(5.4.47) \quad f : A_1 \to A_2$$

*диаграммы представлений (5.4.45) в диаграмму представлений (5.4.46) называется* **гомоморфизмом** *$D$-алгебры $A_1$ в $D$-алгебру $A_2$. Обозначим* $\mathrm{Hom}(D; A_1 \to A_2)$ *множество гомоморфизмов $D$-алгебры $A_1$ в $D$-алгебру $A_2$.* $\quad\square$

Теорема 5.4.8. *Гомоморфизм*

$$\boxed{(5.4.47) \quad f : A_1 \to A_2}$$

*$D$-алгебры $A_1$ в $D$-алгебру $A_2$ - это линейное отображение (5.4.47) $D$-модуля $A_1$ в $D$-модуль $A_2$ такое, что*

$$(5.4.48) \quad f \circ (ab) = (f \circ a)(f \circ b)$$

*и удовлетворяют следующим равенствам*

$$(5.4.49) \quad f \circ (a + b) = f \circ a + f \circ b$$

$$\boxed{(5.4.48) \quad f \circ (ab) = (f \circ a)(f \circ b)}$$

$$(5.4.50) \qquad f \circ (pa) = p(f \circ a)$$

$$p \in D \quad a, b \in A_1$$

Доказательство. Согласно определениям 2.3.9, 2.8.5, отображение (5.4.47) является морфизмом представления $h_{1.12}$, описывающего $D$-модуль $A_1$, в представления $h_{2.12}$, описывающего $D$-модуль $A_2$. Следовательно, отображение $f$ является линейным отображением $D$-модуля $A_1$ в $D$-модуль $A_2$ и равенства (5.4.49), (5.4.50) являются следствием теоремы 4.3.10.

Согласно определениям 2.3.9, 2.8.5, отображение $f$ является морфизмом представления $h_{1.23}$ в представления $h_{2.23}$. Равенство

$$(5.4.51) \qquad f \circ (h_{1.23}(a)(b)) = h_{2.23}(f \circ a)(f \circ b)$$

следует из равенств (2.3.2). Из равенств (5.4.45), (5.4.46), (5.4.51), следует

$$(5.4.52) \qquad f \circ (C_1 \circ (a, b)) = C_2 \circ (f \circ a, f \circ b)$$

Равенство (5.4.48) следует из равенств (5.1.1), (5.4.52). $\qquad\square$

| | |
|---|---|
| ТЕОРЕМА 5.4.9. *Гомоморфизм* [5.10] | ТЕОРЕМА 5.4.10. *Гомоморфизм* [5.11] |
| $$\overline{f} : A_1 \to A_2$$ | $$\overline{f} : A_1 \to A_2$$ |
| *$D$-алгебры столбцов $A_1$ в $D$-алгебру столбцов $A_2$ имеет представление* | *$D$-алгебры строк $A_1$ в $D$-алгебру строк $A_2$ имеет представление* |
| $(5.4.53) \qquad b = fa$ | $(5.4.64) \qquad b = af$ |
| $(5.4.54) \qquad \overline{f} \circ (a^i e_{1i}) = a^i f_i^k e_{2k}$ | $(5.4.65) \qquad \overline{f} \circ (a_i e_1^i) = a_i f_k^i e_2^k$ |
| $(5.4.55) \qquad \overline{f} \circ (e_1 a) = e_2 fa$ | $(5.4.66) \qquad \overline{f} \circ (ae_1) = afe_2$ |
| *относительно заданных базисов. Здесь* | *относительно заданных базисов. Здесь* |
| • *$a$ - координатная матрица $A_1$-числа $\overline{a}$ относительно базиса $\overline{\overline{e}}_1$ (D)* | • *$a$ - координатная матрица $A_1$-числа $\overline{a}$ относительно базиса $\overline{\overline{e}}_1$ (D)* |
| $(5.4.56) \qquad \overline{a} = e_1 a$ | $(5.4.67) \qquad \overline{a} = ae_1$ |
| • *$b$ - координатная матрица $A_2$-числа* | • *$b$ - координатная матрица $A_2$-числа* |
| $$\overline{b} = \overline{f} \circ \overline{a}$$ | $$\overline{b} = \overline{f} \circ \overline{a}$$ |

[5.10] В теоремах 5.4.9, 5.4.11, мы опираемся на следующее соглашение. Пусть множество векторов $\overline{\overline{e}}_1 = (e_{1i}, i \in I)$ является базисом и $C_{1\,j}^{\;\;k}$, $k, i, j \in I$, - структурные константы $D$-алгебры столбцов $A_1$. Пусть множество векторов $\overline{\overline{e}}_2 = (e_{2j}, j \in J)$ является базисом и $C_{2\,j}^{\;\;k}$, $k, i, j \in J$, - структурные константы $D$-алгебры столбцов $A_2$.

[5.11] В теоремах 5.4.10, 5.4.12, мы опираемся на следующее соглашение. Пусть множество векторов $\overline{\overline{e}}_1 = (e_1^i, i \in I)$ является базисом и $C_1{}_{\;k}^{ij}$, $k, i, j \in I$, - структурные константы $D$-алгебры строк $A_1$. Пусть множество векторов $\overline{\overline{e}}_2 = (e_2^j, j \in J)$ является базисом и $C_2{}_{\;k}^{ij}$, $k, i, j \in J$, - структурные константы $D$-алгебры строк $A_2$.

относительно базиса $\overline{\overline{e}}_2$

(5.4.57)            $\overline{b} = e_2 b$

- $f$ - координатная матрица множества $A_2$-чисел ($\overline{f} \circ e_{1i}, i \in I$) относительно базиса $\overline{\overline{e}}_2$.
- Матрица гомоморфизма и структурные константы связаны соотношением

(5.4.58)            $C_{1\,ij}^{\ \ k} f_k^l = f_i^p f_j^q C_{2\,pq}^{\ \ l}$

ДОКАЗАТЕЛЬСТВО. Равенства

$\boxed{(5.4.53) \quad b = fa}$

$\boxed{(5.4.54) \quad \overline{f} \circ (a^i e_{1i}) = a^i f_i^k e_{2k}}$

$\boxed{(5.4.55) \quad \overline{f} \circ (e_1 a) = e_2 f a}$

являются следствием теорем 4.3.11, 5.4.8.

Равенство

(5.4.59)            $\overline{f} \circ (e_{1i} e_{1j}) = C_{1\,ij}^{\ \ k} \overline{f} \circ e_{1k}$

является следствием равенств

$\boxed{(5.2.2) \quad e_i e_j = C_{ij}^k e_k}$

$\boxed{(5.4.50) \quad f \circ (pa) = p(f \circ a)}$

Равенство

(5.4.60)            $\overline{f} \circ (e_{1i} e_{1j}) = C_{1\,ij}^{\ \ k} f_k^l e_{2l}$

является следствием равенства

$\boxed{(4.3.47) \quad \overline{f} \circ e_{1i} = e_2 f_i = e_{2j} f_i^j}$

и равенства (5.4.59). Равенство

(5.4.61)            $\begin{aligned} &\overline{f} \circ (e_{1i} e_{1j}) \\ &= (\overline{f} \circ e_{1i})(\overline{f} \circ e_{1j}) \\ &= f_i^p e_{2p} f_j^q e_{2q} \end{aligned}$

является следствием равенства

$\boxed{(5.4.48) \quad f \circ (ab) = (f \circ a)(f \circ b)}$

Равенство

(5.4.62)            $\overline{f} \circ (e_{1i} e_{2j}) = f_i^p f_j^q C_{2\,pq}^{\ \ l} e_{2l}$

является следствием равенства

$\boxed{(5.2.2) \quad e_i e_j = C_{ij}^k e_k}$

и равенства (5.4.61). Равенство

(5.4.63)            $C_{1\,ij}^{\ \ k} f_k^l e_{2l} = f_i^p f_j^q C_{2\,pq}^{\ \ l} e_{2l}$

относительно базиса $\overline{\overline{e}}_2$

(5.4.68)            $\overline{b} = b e_2$

- $f$ - координатная матрица множества $A_2$-чисел ($\overline{f} \circ e_1^i, i \in I$) относительно базиса $\overline{\overline{e}}_2$.
- Матрица гомоморфизма и структурные константы связаны соотношением

(5.4.69)            $C_{1\,k}^{\ \ ij} f_l^k = f_p^i f_q^j C_{2\,l}^{\ \ pq}$

ДОКАЗАТЕЛЬСТВО. Равенства

$\boxed{(5.4.64) \quad b = af}$

$\boxed{(5.4.65) \quad \overline{f} \circ (a_i e_1^i) = a_i f_k^i e_2^k}$

$\boxed{(5.4.66) \quad \overline{f} \circ (ae_1) = af e_2}$

являются следствием теорем 4.3.12, 5.4.8.

Равенство

(5.4.70)            $\overline{f} \circ (e_1^i e_1^j) = C_{1\,k}^{\ \ ij} \overline{f} \circ e_1^k$

является следствием равенств

$\boxed{(5.2.6) \quad e^i e^j = C_k^{ij} e^k}$

$\boxed{(5.4.50) \quad f \circ (pa) = p(f \circ a)}$

Равенство

(5.4.71)            $\overline{f} \circ (e_1^i e_1^j) = C_{1\,k}^{\ \ ij} f_l^k e_2^l$

является следствием равенства

$\boxed{(4.3.50) \quad \overline{f} \circ e_1^i = f_i e_2 = f_j^i e_2^j}$

и равенства (5.4.70). Равенство

(5.4.72)            $\begin{aligned} &\overline{f} \circ (e_1^i e_1^j) \\ &= (\overline{f} \circ e_1^i)(\overline{f} \circ e_1^j) \\ &= f_p^i e_2^p f_q^j e_2^q \end{aligned}$

является следствием равенства

$\boxed{(5.4.48) \quad f \circ (ab) = (f \circ a)(f \circ b)}$

Равенство

(5.4.73)            $\overline{f} \circ (e_1^i e_1^j) = f_p^i f_q^j C_{2\,l}^{\ \ pq} e_2^l$

является следствием равенства

$\boxed{(5.2.6) \quad e^i e^j = C_k^{ij} e^k}$

и равенства (5.4.72). Равенство

(5.4.74)            $C_{1\,k}^{\ \ ij} f_l^k e_2^l = f_p^i f_q^j C_{2\,l}^{\ \ pq} e_2^l$

является следствием равенств (5.4.60), (5.4.62). Равенство

$$(5.4.58) \quad C_{1\,ij}^{\ \ k} f_k^l = f_i^p f_j^q C_{2\,pq}^{\ \ l}$$

следует из равенства (5.4.63) и теоремы 4.2.7. □

---

**ТЕОРЕМА 5.4.11.** *Пусть*

$$f = (f_j^i, i \in I, j \in J)$$

*матрица $D$-чисел, которая удовлетворяют равенству*

$$(5.4.58) \quad C_{1\,ij}^{\ \ k} f_k^l = f_i^p f_j^q C_{2\,pq}^{\ \ l}$$

*Отображение* [5.10]

$$(4.3.34) \quad \overline{f} : V_1 \to V_2 \quad \text{определённое}$$

*равенством*

$$(4.3.35) \quad \overline{f} \circ (e_1 a) = e_2 f a$$

*является гомоморфизмом $D$-алгебры столбцов $A_1$ в $D$-алгебру столбцов $A_2$. Гомоморфизм (4.3.34), который имеет данную матрицу $f$, определён однозначно.*

ДОКАЗАТЕЛЬСТВО. Согласно теореме 4.3.13, отображение (4.3.34) является гомоморфизмом $D$-модуля $A_1$ в $D$-модуль $A_2$ и гомоморфизм (4.3.34) определён однозначно.

Пусть $a = a^i e_{1i}$, $b = b^i e_{1i}$. Равенство

$$(5.4.75) \quad \begin{aligned} & C_{1\,ij}^{\ \ k} f_k^l a^i b^j e_{2l} \\ & = f_i^p f_j^q C_{2\,pq}^{\ \ l} a^i b^j e_{2l} \end{aligned}$$

является следствием равенства (5.4.58). Равенство

$$(5.4.76) \quad \begin{aligned} & C_{1\,ij}^{\ \ k} f_k^l a^i b^j e_{2l} \\ & = (ab)^k f_k^l e_{2l} \end{aligned}$$

является следствием равенства

$$(5.2.1) \quad (ab)^k = C_{ij}^k a^i b^j$$

Равенство

$$(5.4.77) \quad \begin{aligned} & C_{1\,ij}^{\ \ k} f_k^l a^i b^j e_{2l} \\ & = f \circ (ab) \end{aligned}$$

является следствием равенства (5.4.76) и равенства

$$(5.4.54) \quad \overline{f} \circ (a^i e_{1i}) = a^i f_i^k e_{2k}$$

---

является следствием равенств (5.4.71), (5.4.73). Равенство

$$(5.4.69) \quad C_{1\,k}^{\ ij} f_l^k = f_p^i f_q^j C_{2\,l}^{\ pq}$$

следует из равенства (5.4.74) и теоремы 4.2.14. □

---

**ТЕОРЕМА 5.4.12.** *Пусть*

$$f = (f_i^j, i \in I, j \in J)$$

*матрица $D$-чисел, которая удовлетворяют равенству*

$$(5.4.69) \quad C_{1\,k}^{\ ij} f_l^k = f_p^i f_q^j C_{2\,l}^{\ pq}$$

*Отображение* [5.11]

$$(4.3.40) \quad \overline{f} : V_1 \to V_2 \quad \text{определённое}$$

*равенством*

$$(4.3.41) \quad \overline{f} \circ (ae_1) = afe_2$$

*является гомоморфизмом $D$-алгебры строк $A_1$ в $D$-алгебру строк $A_2$. Гомоморфизм (4.3.40), который имеет данную матрицу $f$, определён однозначно.*

ДОКАЗАТЕЛЬСТВО. Согласно теореме 4.3.14, отображение (4.3.40) является гомоморфизмом $D$-модуля $A_1$ в $D$-модуль $A_2$ и гомоморфизм (4.3.40) определён однозначно.

Пусть $a = a_i e_1^i$, $b = b_i e_1^i$. Равенство

$$(5.4.80) \quad \begin{aligned} & C_{1\,k}^{\ ij} f_l^k a_i b_j e_2^l \\ & = f_p^i f_q^j C_{2\,l}^{\ pq} a_i b_j e_2^l \end{aligned}$$

является следствием равенства (5.4.69). Равенство

$$(5.4.81) \quad \begin{aligned} & C_{1\,k}^{\ ij} f_l^k a_i b_j e_2^l \\ & = (ab)_k f_l^k e_2^l \end{aligned}$$

является следствием равенства

$$(5.2.5) \quad (ab)_k = C_k^{ij} a_i b_j$$

Равенство

$$(5.4.82) \quad \begin{aligned} & C_{1\,k}^{\ ij} f_l^k a_i b_j e_2^l \\ & = f \circ (ab) \end{aligned}$$

является следствием равенства (5.4.81) и равенства

$$(5.4.65) \quad \overline{f} \circ (a_i e_1^i) = a_i f_k^i e_2^k$$

Равенство

$$(5.4.78) \qquad f_i^p f_j^q C_{2pq}{}^l a^i b^j e_{2l}$$
$$= (a^i f_i^p e_{2p})(b^j f_j^q e_{2q})$$

является следствием равенства

$$\boxed{(5.2.2) \quad e_i e_j = C_{ij}^k e_k}$$

Равенство

$$(5.4.79) \qquad f_i^p f_j^q C_{2pq}{}^l a^i b^j e_{2l}$$
$$= (f \circ a)(f \circ b)$$

является следствием равенства (5.4.78) и равенства

$$\boxed{(5.4.54) \quad \overline{f} \circ (a^i e_{1i}) = a^i f_i^k e_{2k}}$$

Равенство

$$\boxed{(5.4.48) \quad f \circ (ab) = (f \circ a)(f \circ b)}$$

является следствием равенств (5.4.75), (5.4.77), (5.4.79). Согласно теореме 5.4.8, отображение (4.3.34) является гомоморфизмом $D$-алгебры $A_1$ в $D$-алгебру $A_2$. $\qquad\square$

Равенство

$$(5.4.83) \qquad f_p^i f_q^j C_{2l}{}^{pq} a_i b_j e_2^l$$
$$= (a_i f_p^i e_2^p)(b_j f_q^j e_2^q)$$

является следствием равенства

$$\boxed{(5.2.6) \quad e^i e^j = C_k^{ij} e^k}$$

Равенство

$$(5.4.84) \qquad f_p^i f_q^j C_{2l}{}^{pq} a_i b_j e_2^l$$
$$= (f \circ a)(f \circ b)$$

является следствием равенства (5.4.83) и равенства

$$\boxed{(5.4.65) \quad \overline{f} \circ (a_i e_1^i) = a_i f_k^i e_2^k}$$

Равенство

$$\boxed{(5.4.48) \quad f \circ (ab) = (f \circ a)(f \circ b)}$$

является следствием равенств (5.4.80), (5.4.82), (5.4.84). Согласно теореме 5.4.8, отображение (4.3.40) является гомоморфизмом $D$-алгебры $A_1$ в $D$-алгебру $A_2$. $\qquad\square$

# Линейное отображение алгебры

## 6.1. Линейное отображение $D$-алгебры

ОПРЕДЕЛЕНИЕ 6.1.1. *Пусть $A_1$ и $A_2$ - алгебры над коммутативным кольцом $D$. Линейное отображение $D$-модуля $A_1$ в $D$-модуль $A_2$ называется **линейным отображением** $D$-алгебры $A_1$ в $D$-алгебру $A_2$.*
*Обозначим $\mathcal{L}(D; A_1 \to A_2)$ множество линейных отображений $D$-алгебры $A_1$ в $D$-алгебру $A_2$.* ☐

ОПРЕДЕЛЕНИЕ 6.1.2. *Пусть $A_1$, ..., $A_n$, $S$ - $D$-алгебры. Полилинейное отображение*

$$f : A_1 \times ... \times A_n \to S$$

*$D$-модулей $A_1$, ..., $A_n$ в $D$-модуль $S$, называется **полилинейным отображением** $D$-алгебр $A_1$, ..., $A_n$ в $D$-алгебру $S$. Обозначим $\mathcal{L}(D; A_1 \times ... \times A_n \to S)$ множество полилинейных отображений $D$-модулей $A_1$, ..., $A_n$ в $D$-модуль $S$. Обозначим $\mathcal{L}(D; A^n \to S)$ множество n-линейных отображений $D$-модуля $A$ ($A_1 = ... = A_n = A$) в $D$-модуль $S$.* ☐

ТЕОРЕМА 6.1.3. *Пусть $A_1$, ..., $A_n$ - $D$-алгебры. Тензорное произведение $A_1 \otimes ... \otimes A_n$ $D$-модулей $A_1$, ..., $A_n$ является $D$-алгеброй, если мы определим произведение согласно правилу*

$$(a_1, ..., a_n) * (b_1, ..., b_n) = (a_1 b_1) \otimes ... \otimes (a_n b_n)$$

ДОКАЗАТЕЛЬСТВО. Согласно определению 5.1.1 и теореме 4.6.2, тензорное произведение $A_1 \otimes ... \otimes A_n$ $D$-алгебр $A_1$, ..., $A_n$ является $D$-модулем.

Рассмотрим отображение

(6.1.1) $$* : (A_1 \times ... \times A_n) \times (A_1 \times ... \times A_n) \to A_1 \otimes ... \otimes A_n$$

определённое равенством

$$(a_1, ..., a_n) * (b_1, ..., b_n) = (a_1 b_1) \otimes ... \otimes (a_n b_n)$$

Для заданных значений переменных $b_1$, ..., $b_n$, отображение (6.1.1) полилинейно по переменным $a_1$, ..., $a_n$. Согласно теореме 4.6.4, существует линейное отображение

(6.1.2) $$*(b_1, ..., b_n) : A_1 \otimes ... \otimes A_n \to A_1 \otimes ... \otimes A_n$$

определённое равенством

(6.1.3) $$(a_1 \otimes ... \otimes a_n) * (b_1, ..., b_n) = (a_1 b_1) \otimes ... \otimes (a_n b_n)$$

Так как произвольный тензор $a \in A_1 \otimes ... \otimes A_n$ можно представить в виде суммы тензоров вида $a_1 \otimes ... \otimes a_n$, то для заданного тензора $a \in A_1 \otimes ... \otimes A_n$ отображение (6.1.2) является полилинейным отображением переменных $b_1$, ..., $b_n$. Согласно теореме 4.6.4, существует линейное отображение

(6.1.4) $$*(a) : A_1 \otimes ... \otimes A_n \to A_1 \otimes ... \otimes A_n$$

определённое равенством

(6.1.5)        $(a_1 \otimes ... \otimes a_n) * (b_1 \otimes ... \otimes b_n) = (a_1 b_1) \otimes ... \otimes (a_n b_n)$

Следовательно, равенство (6.1.5) определяет билинейное отображение

(6.1.6)        $* : (A_1 \otimes ... \otimes A_n) \times (A_1 \otimes ... \otimes A_n) \to A_1 \otimes ... \otimes A_n$

Билинейное отображение (6.1.6) определяет произведение в $D$-модуле $A_1 \otimes ... \otimes A_n$. $\qquad\square$

В случае тензорного произведения $D$-алгебр $A_1$, $A_2$ мы будем рассматривать произведение, определённое равенством

(6.1.7)        $(a_1 \otimes a_2) \circ (b_1 \otimes b_2) = (a_1 b_1) \otimes (b_2 a_2)$

Теорема 6.1.4.  *Пусть $\overline{\overline{e}}_i$ - базис алгебры $A_i$ над кольцом $D$. Пусть $B_{i \cdot kl}^{\cdot \cdot j}$ - структурные константы алгебры $A_i$ относительно базиса $\overline{\overline{e}}_i$. Структурные константы тензорного произведение $A_1 \otimes ... \otimes A_n$ относительно базиса $e_{1 \cdot i_l} \otimes ... \otimes e_{n \cdot i_n}$ имеют вид*

(6.1.8)        $C_{\cdot k_l ... k_n \cdot l_l ... l_n}^{\cdot j_l ... j_n} = C_{1 \cdot k_l l_l}^{\cdot j_l} ... C_{n \cdot k_n l_n}^{\cdot j_n}$

Доказательство.  Непосредственное перемножение тензоров $e_{1 \cdot i_l} \otimes ... \otimes e_{n \cdot i_n}$ имеет вид

$$
\begin{aligned}
& (e_{1 \cdot k_l} \otimes ... \otimes e_{n \cdot k_n})(e_{1 \cdot l_l} \otimes ... \otimes e_{n \cdot l_n}) \\
=& (\overline{e}_{1 \cdot k_l} \overline{e}_{1 \cdot l_l}) \otimes ... \otimes (\overline{e}_{n \cdot k_n} \overline{e}_{n \cdot l_n}) \\
=& (\overline{e}_{1 \cdot k_l} \overline{e}_{1 \cdot l_l}) \otimes ... \otimes (\overline{e}_{n \cdot k_n} \overline{e}_{n \cdot l_n}) \\
=& (C_{1 \cdot k_l l_l}^{\cdot j_l} \overline{e}_{1 \cdot j_l}) \otimes ... \otimes (C_{n \cdot k_n l_n}^{\cdot j_n} \overline{e}_{n \cdot j_n}) \\
=& C_{1 \cdot k_l l_l}^{\cdot j_l} ... C_{n \cdot k_n l_n}^{\cdot j_n} \overline{e}_{1 \cdot j_l} \otimes ... \otimes \overline{e}_{n \cdot j_n}
\end{aligned}
$$

(6.1.9)

Согласно определению структурных констант

(6.1.10)    $(e_{1 \cdot k_l} \otimes ... \otimes e_{n \cdot k_n})(e_{1 \cdot l_l} \otimes ... \otimes e_{n \cdot l_n}) = C_{\cdot k_l ... k_n \cdot l_l ... l_n}^{\cdot j_l ... j_n} (e_{1 \cdot j_l} \otimes ... \otimes e_{n \cdot j_n})$

Равенство (6.1.8) следует из сравнения (6.1.9), (6.1.10).

Из цепочки равенств

$$
\begin{aligned}
& (a_1 \otimes ... \otimes a_n)(b_1 \otimes ... \otimes b_n) \\
=& (a_1^{k_l} \overline{e}_{1 \cdot k_l} \otimes ... \otimes a_n^{k_n} \overline{e}_{n \cdot k_n})(b_1^{l_l} \overline{e}_{1 \cdot l_l} \otimes ... \otimes b_n^{l_n} \overline{e}_{n \cdot l_n}) \\
=& a_1^{k_l} ... a_n^{k_n} b_1^{l_l} ... b_n^{l_n} (\overline{e}_{1 \cdot k_l} \otimes ... \otimes \overline{e}_{n \cdot k_n})(\overline{e}_{1 \cdot l_l} \otimes ... \otimes \overline{e}_{n \cdot l_n}) \\
=& a_1^{k_l} ... a_n^{k_n} b_1^{l_l} ... b_n^{l_n} C_{\cdot k_l ... k_n \cdot l_l ... l_n}^{\cdot j_l ... j_n} (e_{1 \cdot j_l} \otimes ... \otimes e_{n \cdot j_n}) \\
=& a_1^{k_l} ... a_n^{k_n} b_1^{l_l} ... b_n^{l_n} C_{1 \cdot k_l l_l}^{\cdot j_l} ... C_{n \cdot k_n l_n}^{\cdot j_n} (e_{1 \cdot j_l} \otimes ... \otimes e_{n \cdot j_n}) \\
=& (a_1^{k_l} b_1^{l_l} C_{1 \cdot k_l l_l}^{\cdot j_l} \overline{e}_{1 \cdot j_l}) \otimes ... \otimes (a_n^{k_n} b_n^{l_n} C_{n \cdot k_n l_n}^{\cdot j_n} \overline{e}_{n \cdot j_n}) \\
=& (a_1 b_1) \otimes ... \otimes (a_n b_n)
\end{aligned}
$$

следует, что определение произведения (6.1.10) со структурными константами (6.1.8) согласовано с определением произведения (6.1.5). $\qquad\square$

Теорема 6.1.5.  *Для тензоров $a$, $b \in A_1 \otimes ... \otimes A_n$ стандартные компоненты произведения удовлетворяют равенству*

(6.1.11)        $(ab)^{j_l ... j_n} = C_{\cdot k_l ... k_n \cdot l_l ... l_n}^{\cdot j_l ... j_n} a^{k_l ... k_n} b^{l_l ... l_n}$

Доказательство. Согласно определению

$$(6.1.12) \qquad ab = (ab)^{j_1 \ldots j_n} e_{1 \cdot j_1} \otimes \ldots \otimes e_{n \cdot j_n}$$

В тоже время

$$(6.1.13) \qquad ab = a^{k_1 \ldots k_n} e_{1 \cdot k_1} \otimes \ldots \otimes e_{n \cdot k_n} b^{l_1 \ldots l_n} e_{1 \cdot l_1} \otimes \ldots \otimes e_{n \cdot l_n}$$

$$= a^{k_1 \ldots k_n} b^{l_1 \ldots l_n} C^{\cdot j_1 \ldots j_n}_{\cdot k_1 \ldots k_n \cdot l_1 \ldots l_n} e_{1 \cdot j_1} \otimes \ldots \otimes e_{n \cdot j_n}$$

Равенство (6.1.11) следует из равенств (6.1.12), (6.1.13). $\square$

Теорема 6.1.6. *Если алгебра $A_i$, $i = 1$, ..., $n$, ассоциативна, то тензорное произведение $A_1 \otimes \ldots \otimes A_n$ - ассоциативная алгебра.*

Доказательство. Поскольку

$$((e_{1 \cdot i_1} \otimes \ldots \otimes e_{n \cdot i_n})(e_{1 \cdot j_1} \otimes \ldots \otimes e_{n \cdot j_n}))(e_{1 \cdot k_1} \otimes \ldots \otimes e_{n \cdot k_n})$$

$$= ((\overline{e}_{1 \cdot i_1} \overline{e}_{1 \cdot j_1}) \otimes \ldots \otimes (\overline{e}_{n \cdot i_n} \overline{e}_{1 \cdot j_n}))(e_{1 \cdot k_1} \otimes \ldots \otimes e_{n \cdot k_n})$$

$$= ((\overline{e}_{1 \cdot i_1} \overline{e}_{1 \cdot j_1}) \overline{e}_{1 \cdot k_1}) \otimes \ldots \otimes ((\overline{e}_{n \cdot i_n} \overline{e}_{1 \cdot j_n}) \overline{e}_{1 \cdot k_n})$$

$$= (\overline{e}_{1 \cdot i_1} (\overline{e}_{1 \cdot j_1} \overline{e}_{1 \cdot k_1})) \otimes \ldots \otimes (\overline{e}_{n \cdot i_n} (\overline{e}_{1 \cdot j_n} \overline{e}_{1 \cdot k_n}))$$

$$= (e_{1 \cdot i_1} \otimes \ldots \otimes e_{n \cdot i_n})((\overline{e}_{1 \cdot j_1} \overline{e}_{1 \cdot k_1}) \otimes \ldots \otimes (\overline{e}_{1 \cdot j_n} \overline{e}_{1 \cdot k_n}))$$

$$= (e_{1 \cdot i_1} \otimes \ldots \otimes e_{n \cdot i_n})((e_{1 \cdot j_1} \otimes \ldots \otimes e_{n \cdot j_n})(e_{1 \cdot k_1} \otimes \ldots \otimes e_{n \cdot k_n}))$$

то

$$(ab)c = a^{i_1 \ldots i_n} b^{j_1 \ldots j_n} c^{k_1 \ldots k_n}$$

$$((e_{1 \cdot i_1} \otimes \ldots \otimes e_{n \cdot i_n})(e_{1 \cdot j_1} \otimes \ldots \otimes e_{n \cdot j_n}))(e_{1 \cdot k_1} \otimes \ldots \otimes e_{n \cdot k_n})$$

$$= a^{i_1 \ldots i_n} b^{j_1 \ldots j_n} c^{k_1 \ldots k_n}$$

$$(e_{1 \cdot i_1} \otimes \ldots \otimes e_{n \cdot i_n})((e_{1 \cdot j_1} \otimes \ldots \otimes e_{n \cdot j_n})(e_{1 \cdot k_1} \otimes \ldots \otimes e_{n \cdot k_n}))$$

$$= a(bc)$$

$\square$

Теорема 6.1.7. *Пусть $A$ - алгебра над коммутативным кольцом $D$. Существует линейное отображение*

$$h : a \otimes b \in A \otimes A \to ab \in A$$

Доказательство. Теорема является следствием определения 5.1.1 и теоремы 4.6.4. $\square$

Теорема 6.1.8. *Пусть отображение*

$$f : A_1 \to A_2$$

*является линейным отображением $D$-алгебры $A_1$ в $D$-алгебру $A_2$. Тогда отображения $a \circ f$, $b \star f$, $a$, $b \in A_2$, определённые равенствами*

$$(a \circ f) \circ x = a(f \circ x)$$

$$(b \star f) \circ x = (f \circ x)b$$

*также являются линейными.*

Доказательство. Утверждение теоремы следует из цепочек равенств

$$(af) \circ (x+y) = a(f \circ (x+y)) = a(f \circ x + f \circ y) = a(f \circ x) + a(f \circ y)$$
$$= (af) \circ x + (af) \circ y$$
$$(af) \circ (px) = a(f \circ (px)) = ap(f \circ x) = pa(f \circ x)$$
$$= p((af) \circ x)$$
$$(fb) \circ (x+y) = (f \circ (x+y))b = (f \circ x + f \circ y)\, b = (f \circ x)b + (f \circ y)b$$
$$= (fb) \circ x + (fb) \circ y$$
$$(fb) \circ (px) = (f \circ (px))b = p(f \circ x)b$$
$$= p((fb) \circ x)$$

$\square$

## 6.2. Алгебра $\mathcal{L}(D; A \to A)$

Теорема 6.2.1. *Пусть $A$, $B$, $C$ - алгебры над коммутативным кольцом $D$. Пусть $f$ - линейное отображение из $D$-алгебры $A$ в $D$-алгебру $B$. Пусть $g$ - линейное отображение из $D$-алгебры $B$ в $D$-алгебру $C$. Отображение $g \circ f$, определённое диаграммой*

(6.2.1)

*является линейным отображением из $D$-алгебры $A$ в $D$-алгебру $C$.*

Доказательство. Доказательство теоремы следует из цепочек равенств

$$(g \circ f) \circ (a+b) = g \circ (f \circ (a+b)) = g \circ (f \circ a + f \circ b)$$
$$= g \circ (f \circ a) + g \circ (f \circ b) = (g \circ f) \circ a + (g \circ f) \circ b$$
$$(g \circ f) \circ (pa) = g \circ (f \circ (pa)) = g \circ (p\, f \circ a) = p\, g \circ (f \circ a)$$
$$= p\,(g \circ f) \circ a$$

$\square$

Теорема 6.2.2. *Пусть $A$, $B$, $C$ - алгебры над коммутативным кольцом $D$. Пусть $f$ - линейное отображение из $D$-алгебры $A$ в $D$-алгебру $B$. Отображение $f$ порождает линейное отображение*

(6.2.2)           $f^* : g \in \mathcal{L}(D; B \to C) \to g \circ f \in \mathcal{L}(D; A \to C)$

(6.2.3)

Доказательство. Доказательство теоремы следует из цепочек равенств [6.1]

$$((g_1 + g_2) \circ f) \circ a = (g_1 + g_2) \circ (f \circ a) = g_1 \circ (f \circ a) + g_2 \circ (f \circ a)$$
$$= (g_1 \circ f) \circ a + (g_2 \circ f) \circ a$$
$$= (g_1 \circ f + g_2 \circ f) \circ a$$
$$((pg) \circ f) \circ a = (pg) \circ (f \circ a) = p \ g \circ (f \circ a) = p \ (g \circ f) \circ a$$
$$= (p(g \circ f)) \circ a$$

$\square$

Теорема 6.2.3. *Пусть $A$, $B$, $C$ - алгебры над коммутативным кольцом $D$. Пусть $g$ - линейное отображение из $D$-алгебры $B$ в $D$-алгебру $C$. Отображение $g$ порождает линейное отображение*

(6.2.6) $\qquad g_* : f \in \mathcal{L}(D; A \to B) \to g \circ f \in \mathcal{L}(D; A \to C)$

(6.2.7)

Доказательство. Доказательство теоремы следует из цепочек равенств [6.2]

$$(g \circ (f_1 + f_2)) \circ a = g \circ ((f_1 + f_2) \circ a) = g \circ (f_1 \circ a + f_2 \circ a)$$
$$= g \circ (f_1 \circ a) + g \circ (f_2 \circ a) = (g \circ f_1) \circ a + (g \circ f_2) \circ a$$
$$= (g \circ f_1 + g \circ f_2) \circ a$$
$$(g \circ (pf)) \circ a = g \circ ((pf) \circ a) = g \circ (p \ (f \circ a)) = p \ g \circ (f \circ a)$$
$$= p \ (g \circ f) \circ a = (p(g \circ f)) \circ a$$

$\square$

Теорема 6.2.4. *Пусть $A$, $B$, $C$ - алгебры над коммутативным кольцом $D$. Отображение*

(6.2.10) $\qquad \circ : (g, f) \in \mathcal{L}(D; B \to C) \times \mathcal{L}(D; A \to B) \to g \circ f \in \mathcal{L}(D; A \to C)$

*является билинейным отображением.*

Доказательство. Теорема является следствием теорем 6.2.2, 6.2.3. $\square$

---

[6.1]Мы пользуемся следующими определениями операций над отображениями

(6.2.4) $\qquad\qquad\qquad\qquad (f + g) \circ a = f \circ a + g \circ a$

(6.2.5) $\qquad\qquad\qquad\qquad (pf) \circ a = p \ f \circ a$

[6.2]Мы пользуемся следующими определениями операций над отображениями

(6.2.8) $\qquad\qquad\qquad\qquad (f + g) \circ a = f \circ a + g \circ a$

(6.2.9) $\qquad\qquad\qquad\qquad (pf) \circ a = p \ f \circ a$

Теорема 6.2.5. *Пусть $A$ - алгебра над коммутативным кольцом $D$. $D$-модуль $\mathcal{L}(D; A \to A)$, оснащённый произведением*

$$(6.2.11) \qquad \circ : (g, f) \in \mathcal{L}(D; A \to A) \times \mathcal{L}(D; A \to A) \to g \circ f \in \mathcal{L}(D; A \to A)$$

*является алгеброй над $D$.*

Доказательство. Теорема является следствием определения 5.1.1 и теоремы 6.2.4. $\qquad \square$

## 6.3. Линейное отображение в ассоциативную алгебру

Теорема 6.3.1. *Рассмотрим $D$-алгебры $A_1$ и $A_2$. Для заданного отображения $f \in \mathcal{L}(D; A_1 \to A_2)$ отображение*

$$g : A_2 \times A_2 \to \mathcal{L}(D; A_1 \to A_2)$$

$$g(a, b) \circ f = afb$$

*является билинейным отображением.*

Доказательство. Утверждение теоремы следует из цепочек равенств

$$((a_1 + a_2)fb) \circ x = (a_1 + a_2) \ f \circ x \ b = a_1 \ f \circ x \ b + a_2 \ f \circ x \ b$$
$$= (a_1 fb) \circ x + (a_2 fb) \circ x = (a_1 fb + a_2 fb) \circ x$$
$$((pa)fb) \circ x = (pa) \ f \circ x \ b = p(a \ f \circ x \ b) = p((afb) \circ x) = (p(afb)) \circ x$$
$$(af(b_1 + b_2)) \circ x = a \ f \circ x \ (b_1 + b_2) = a \ f \circ x \ b_1 + a \ f \circ x \ b_2$$
$$= (afb_1) \circ x + (afb_2) \circ x = (afb_1 + afb_2) \circ x$$
$$(af(pb)) \circ x = a \ f \circ x \ (pb) = p(a \ f \circ x \ b) = p((afb) \circ x) = (p(afb)) \circ x$$

$\qquad \square$

Теорема 6.3.2. *Для заданного отображения $f \in \mathcal{L}(D; A_1 \to A_2)$ $D$-алгебры $A_1$ в $D$-алгебру $A_2$ существует линейное отображение*

$$h : A_2 \otimes A_2 \to \mathcal{L}(D; A_1 \to A_2)$$

*определённое равенством*

$$(6.3.1) \qquad (a \otimes b) \circ f = b \star (a \circ f)$$

$$((a \otimes b) \circ f) \circ x = a(f \circ x)b$$

Доказательство. Утверждение теоремы является следствием теорем 4.6.4, 6.3.1. $\qquad \square$

Теорема 6.3.3. *Пусть $A_1$ и $A_2$ являются $D$-алгебрами. Пусть произведение в алгебре $A \otimes A$ определено согласно правилу*

$$(6.3.2) \qquad (p_0 \otimes p_1) \circ (q_0 \otimes q_1) = (p_0 q_0) \otimes (q_1 p_1)$$

*Линейное отображение*

$$(6.3.3) \qquad h : A_2 \otimes A_2 \overset{*}{\longrightarrow} \mathcal{L}(D; A_1 \to A_2)$$

*определённое равенством*

(6.3.4)        $(a \otimes b) \circ f = afb \quad a, b \in A_2 \quad f \in \mathcal{L}(D; A_1 \to A_2)$

*является представлением[6.3]алгебры $A_2 \otimes A_2$ в модуле $\mathcal{L}(D; A_1 \to A_2)$.*

Доказательство. Согласно теореме 6.1.8, отображение (6.3.4) является преобразованием модуля $\mathcal{L}(D; A_1 \to A_2)$. Для данного тензора $c \in A_2 \otimes A_2$ преобразование $h(c)$ является линейным преобразованием модуля $\mathcal{L}(D; A_1 \to A_2)$, так как

$$\begin{aligned}
((a \otimes b) \circ (f_1 + f_2)) \circ x &= (a(f_1 + f_2)b) \circ x = a((f_1 + f_2) \circ x)b \\
&= a(f_1 \circ x + f_2 \circ x)b = a(f_1 \circ x)b + a(f_2 \circ x)b \\
&= (af_1b) \circ x + (af_2b) \circ x \\
&= (a \otimes b) \circ f_1 \circ x + (a \otimes b) \circ f_2 \circ x \\
&= ((a \otimes b) \circ f_1 + (a \otimes b) \circ f_2) \circ x \\
((a \otimes b) \circ (pf)) \circ x &= (a(pf)b) \circ x = a((pf) \circ x)b \\
&= a(p \ f \circ x)b = pa(f \circ x)b \\
&= p \ (afb) \circ x = p \ ((a \otimes b) \circ f) \circ x \\
&= (p((a \otimes b) \circ f)) \circ x
\end{aligned}$$

Согласно теореме 6.3.2, отображение (6.3.4) является линейным отображением.

Пусть $f \in \mathcal{L}(D; A_1 \to A_2)$, $a \otimes b$, $c \otimes d \in A_2 \otimes A_2$. Согласно теореме 6.3.2

$$(a \otimes b) \circ f = afb \in \mathcal{L}(D; A_1 \to A_2)$$

Следовательно, согласно теореме 6.3.2

$$(c \otimes d) \circ ((a \otimes b) \circ f) = c(afb)d$$

Поскольку произведение в алгебре $A_2$ ассоциативно, то

$$(c \otimes d) \circ ((a \otimes b) \circ f) = c(afb)d = (ca)f(bd) = (ca \otimes bd) \circ f$$

Следовательно, если мы определим произведение в алгебре $A_2 \otimes A_2$ согласно равенству (6.3.2), то отображение (6.3.3) является морфизмом алгебр. Согласно определению 2.5.1, отображение (6.3.4) является представлением алгебры $A_2 \otimes A_2$ в модуле $\mathcal{L}(D; A_1 \to A_2)$.        □

Теорема 6.3.4. *Пусть $A$ является $D$-алгеброй. Пусть произведение в $D$-модуле $A \otimes A$ определено согласно правилу*

$$(p_0 \otimes p_1) \circ (q_0 \otimes q_1) = (p_0 q_0) \otimes (q_1 p_1)$$

*Представление*

$$h : A \otimes A \relbar\joinrel\relbar\joinrel\twoheadrightarrow \mathcal{L}(D; A \to A) \qquad h(p) : g \to p \circ g$$

*$D$-алгебры $A \otimes A$ в модуле $\mathcal{L}(D; A \to A)$, определённое равенством*

$$(a \otimes b) \circ g = agb \quad a, b \in A \quad g \in \mathcal{L}(D; A \to A)$$

---

[6.3] Определение представления $\Omega$-алгебры дано в определении 2.5.1.

*позволяет отождествить тензор $d \in A \times A$ с линейным отображением $d \circ \delta \in \mathcal{L}(D; A \to A)$, где $\delta \in \mathcal{L}(D; A \to A)$ - тождественное отображение. Линейное отображение $(a \otimes b) \circ \delta$ имеет вид*

$$(6.3.5) \qquad (a \otimes b) \circ c = acb$$

Доказательство. Согласно теореме 6.3.2, отображение $f \in \mathcal{L}(D; A \to A)$ и тензор $d \in A \otimes A$ порождают отображение

$$(6.3.6) \qquad x \to (d \circ f) \circ x$$

Если мы положим $f = \delta$, $d = a \otimes b$, то равенство (6.3.6) приобретает вид

$$(6.3.7) \qquad ((a \otimes b) \circ \delta) \circ x = (a\delta b) \circ x = a \, (\delta \circ x) \, b = axb$$

Если мы положим

$$(6.3.8) \qquad ((a \otimes b) \circ \delta) \circ x = (a \otimes b) \circ (\delta \circ x) = (a \otimes b) \circ x$$

то сравнение равенств (6.3.7) и (6.3.8) даёт основание отождествить действие тензора $a \otimes b$ с преобразованием $(a \otimes b) \circ \delta$. □

Из теоремы 6.3.4 следует, что отображение (6.3.4) можно рассматривать как произведение отображений $a \otimes b$ и $f$. Тензор $a \in A_2 \otimes A_2$ **невырожден**, если существует тензор $b \in A_2 \otimes A_2$ такой, что $a \circ b = 1 \otimes 1$.

ОПРЕДЕЛЕНИЕ 6.3.5. *Рассмотрим[6.4] представление алгебры $A_2 \otimes A_2$ в модуле $\mathcal{L}(D; A_1 \to A_2)$.* **Орбитой линейного отображения** $f \in \mathcal{L}(D; A_1 \to A_2)$ *называется множество*

$$(A_2 \otimes A_2) \circ f = \{g = d \circ f : d \in A_2 \otimes A_2\}$$

□

ТЕОРЕМА 6.3.6. *Рассмотрим $D$-алгебру $A_1$ и ассоциативную $D$-алгебру $A_2$. Рассмотрим представление алгебры $A_2 \otimes A_2$ в модуле $\mathcal{L}(D; A_1 \to A_2)$. Отображение*

$$h : A_1 \to A_2$$

*порождённое отображением*

$$f : A_1 \to A_2$$

*имеет вид*

$$(6.3.9) \qquad h = (a_{s \cdot 0} \otimes a_{s \cdot 1}) \circ f = a_{s \cdot 0} f a_{s \cdot 1}$$

Доказательство. Произвольный тензор $a \in A_2 \otimes A_2$ можно представить в виде

$$a = a_{s \cdot 0} \otimes a_{s \cdot 1}$$

Согласно теореме 6.3.3, отображение (6.3.4) линейно. Это доказывает утверждение теоремы. □

ОПРЕДЕЛЕНИЕ 6.3.7. *Выражение $a_{s \cdot p}$, $p = 0, 1$, в равенстве (6.3.9) называется* **компонентой линейного отображения** $h$. □

---

[6.4]Определение дано по аналогии с определением 2.5.11.

Теорема 6.3.8. *Пусть $A_2$ - алгебра с единицей $e$. Пусть $a \in A_2 \otimes A_2$ - невырожденный тензор. Орбиты линейных отображений $f \in \mathcal{L}(D; A_1 \to A_2)$ и $g = a \circ f$ совпадают*

$$(6.3.10) \qquad (A_2 \otimes A_2) \circ f = (A_2 \otimes A_2) \circ g$$

Доказательство. Если $h \in (A_2 \otimes A_2) \circ g$, то существует $b \in A_2 \otimes A_2$ такое, что $h = b \circ g$. Тогда

$$(6.3.11) \qquad h = b \circ (a \circ f) = (b \circ a) \circ f$$

Следовательно, $h \in (A_2 \otimes A_2) \circ f$,

$$(6.3.12) \qquad (A_2 \otimes A_2) \circ g \subset (A_2 \otimes A_2) \circ f$$

Так как $a$ - невырожденный тензор, то

$$(6.3.13) \qquad f = a^{-1} \circ g$$

Если $h \in (A_2 \otimes A_2) \circ f$, то существует $b \in A_2 \otimes A_2$ такое, что

$$(6.3.14) \qquad h = b \circ f$$

Из равенств (6.3.13), (6.3.14), следует, что

$$h = b \circ (a^{-1} \circ g) = (b \circ a^{-1}) \circ g$$

Следовательно, $h \in (A_2 \otimes A_2) \circ g$,

$$(6.3.15) \qquad (A_2 \otimes A_2) \circ f \subset (A_2 \otimes A_2) \circ g$$

(6.3.10) следует из равенств (6.3.12), (6.3.15). $\qquad\square$

Из теоремы 6.3.8 также следует, что, если $g = a \circ f$ и $a \in A_2 \otimes A_2$ - вырожденный тензор, то отношение (6.3.12) верно. Однако основной результат теоремы 6.3.8 состоит в том, что представления алгебры $A_2 \otimes A_2$ в модуле $\mathcal{L}(D; A_1 \to A_2)$ порождает отношение эквивалентности в модуле $\mathcal{L}(D; A_1 \to A_2)$. Если удачно выбрать представители каждого класса эквивалентности, то полученное множество будет множеством образующих рассматриваемого представления. [6.5]

## 6.4. Полилинейное отображение в ассоциативную алгебру

Теорема 6.4.1. *Пусть $A_1$, ..., $A_n$, $A$ - ассоциативные $D$-алгебры. Пусть*

$$f_i \in \mathcal{L}(D; A_i \to A) \quad i = 1, ..., n$$

$$a_j \in A \quad j = 0, ..., n$$

*Для заданной перестановки $\sigma$ $n$ переменных отображение*

$$((a_0, ..., a_n, \sigma) \circ (f_1, ..., f_n)) \circ (x_1, ..., x_n)$$
$$(6.4.1) \qquad = (a_0 \sigma(f_1) a_1 ... a_{n-1} \sigma(f_n) a_n) \circ (x_1, ..., x_n)$$
$$= a_0 \sigma(f_1 \circ x_1) a_1 ... a_{n-1} \sigma(f_n \circ x_n) a_n$$

*является $n$-линейным отображением в алгебру $A$.*

---

[6.5] Множество образующих представления определено в определении 2.7.4.

Доказательство. Утверждение теоремы следует из цепочек равенств

$$((a_0, ..., a_n, \sigma) \circ (f_1, ..., f_n)) \circ (x_1, ..., x_i + y_i, ..., x_n)$$
$$= a_0\sigma(f_1 \circ x_1)a_1...\sigma(f_i \circ (x_i + y_i))...a_{n-1}\sigma(f_n \circ x_n)a_n$$
$$= a_0\sigma(f_1 \circ x_1)a_1...\sigma(f_i \circ x_i + f_i \circ y_i)...a_{n-1}\sigma(f_n \circ x_n)a_n$$
$$= a_0\sigma(f_1 \circ x_1)a_1...\sigma(f_i \circ x_i)...a_{n-1}\sigma(f_n \circ x_n)a_n$$
$$+ a_0\sigma(f_1 \circ x_1)a_1...\sigma(f_i \circ y_i)...a_{n-1}\sigma(f_n \circ x_n)a_n$$
$$= ((a_0, ..., a_n, \sigma) \circ (f_1, ..., f_n)) \circ (x_1, ..., x_i, ..., x_n)$$
$$+ ((a_0, ..., a_n, \sigma) \circ (f_1, ..., f_n)) \circ (x_1, ..., y_i, ..., x_n)$$

$$((a_0, ..., a_n, \sigma) \circ (f_1, ..., f_n)) \circ (x_1, ..., px_i, ..., x_n)$$
$$= a_0\sigma(f_1 \circ x_1)a_1...\sigma(f_i \circ (px_i))...a_{n-1}\sigma(f_n \circ x_n)a_n$$
$$= a_0\sigma(f_1 \circ x_1)a_1...\sigma(p(f_i \circ x_i))...a_{n-1}\sigma(f_n \circ x_n)a_n$$
$$= p(a_0\sigma(f_1 \circ x_1)a_1...\sigma(f_i \circ x_i)...a_{n-1}\sigma(f_n \circ x_n)a_n)$$
$$= p(((a_0, ..., a_n, \sigma) \circ (f_1, ..., f_n)) \circ (x_1, ..., x_i, ..., x_n))$$

$\square$

В равенстве (6.4.1), также как и в других выражениях полилинейного отображения, принято соглашение, что отображение $f_i$ имеет своим аргументом переменную $x_i$.

Теорема 6.4.2. *Пусть $A_1$, ..., $A_n$, $A$ - ассоциативные $D$-алгебры. Для заданного семейства отображений*

$$f_i \in \mathcal{L}(D; A_i \to A) \quad i = 1, ..., n$$

*отображение*

$$h : A^{n+1} \to \mathcal{L}(D; A_1 \times ... \times A_n \to A)$$

*определённое равенством*

$$(a_0, ..., a_n, \sigma) \circ (f_1, ..., f_n) = a_0\sigma(f_1)a_1...a_{n-1}\sigma(f_n)a_n$$

*является $n+1$-линейным отображением в $D$-модуль $\mathcal{L}(D; A_1 \times ... \times A_n \to A)$.*

Доказательство. Утверждение теоремы следует из цепочек равенств

$$((a_0, ..., a_i + b_i, ...a_n, \sigma) \circ (f_1, ..., f_n)) \circ (x_1, ..., x_n)$$
$$= a_0\sigma(f_1 \circ x_1)a_1...(a_i + b_i)...a_{n-1}\sigma(f_n \circ x_n)a_n$$
$$= a_0\sigma(f_1 \circ x_1)a_1...a_i...a_{n-1}\sigma(f_n \circ x_n)a_n + a_0\sigma(f_1 \circ x_1)a_1...b_i...a_{n-1}\sigma(f_n \circ x_n)a_n$$
$$= ((a_0, ..., a_i, ..., a_n, \sigma) \circ (f_1, ..., f_n)) \circ (x_1, ..., x_n)$$
$$+ ((a_0, ..., b_i, ..., a_n, \sigma) \circ (f_1, ..., f_n)) \circ (x_1, ..., x_n)$$
$$= ((a_0, ..., a_i, ..., a_n, \sigma) \circ (f_1, ..., f_n) + (a_0, ..., b_i, ..., a_n, \sigma) \circ (f_1, ..., f_n)) \circ (x_1, ..., x_n)$$

$$((a_0, ..., pa_i, ...a_n, \sigma) \circ (f_1, ..., f_n)) \circ (x_1, ..., x_n)$$
$$= a_0\sigma(f_1 \circ x_1)a_1...pa_i...a_{n-1}\sigma(f_n \circ x_n)a_n$$
$$= p(a_0\sigma(f_1 \circ x_1)a_1...a_i...a_{n-1}\sigma(f_n \circ x_n)a_n)$$
$$= p(((a_0, ..., a_i, ..., a_n, \sigma) \circ (f_1, ..., f_n)) \circ (x_1, ..., x_n))$$
$$= (p((a_0, ..., a_i, ..., a_n, \sigma) \circ (f_1, ..., f_n))) \circ (x_1, ..., x_n)$$

$\square$

Теорема 6.4.3. *Пусть* $A_1$, ..., $A_n$, $A$ *- ассоциативные* $D$*-алгебры. Для заданного семейства отображений*

$$f_i \in \mathcal{L}(D; A_i \to A) \quad i = 1, ..., n$$

*существует линейное отображение*

$$h : A^{n+1\otimes} \times S_n \to \mathcal{L}(D; A_1 \times ... \times A_n \to A)$$

*определённое равенством*

$$(6.4.2) \qquad (a_0 \otimes ... \otimes a_n, \sigma) \circ (f_1, ..., f_n) = (a_0, ..., a_n, \sigma) \circ (f_1, ..., f_n)$$
$$= a_0 \sigma(f_1) a_1 ... a_{n-1} \sigma(f_n) a_n$$

Доказательство. Утверждение теоремы является следствием теорем 4.6.4, 6.4.2. $\square$

Теорема 6.4.4. *Пусть* $A_1$, ..., $A_n$, $A$ *- ассоциативные* $D$*-алгебры. Для заданного тензора* $a \in A^{n+1\otimes}$ *и заданной перестановки* $\sigma \in S_n$ *отображение*

$$h : \prod_{i=1}^{n} \mathcal{L}(D; A_i \to A) \to \mathcal{L}(D; A_1 \times ... \times A_n \to A)$$

*определённое равенством*

$$(a_0 \otimes ... \otimes a_n, \sigma) \circ (f_1, ..., f_n) = a_0 \sigma(f_1) a_1 ... a_{n-1} \sigma(f_n) a_n$$

*является n-линейным отображением в* $D$*-модуль* $\mathcal{L}(D; A_1 \times ... \times A_n \to A)$.

Доказательство. Утверждение теоремы следует из цепочек равенств

$$((a_0 \otimes ... \otimes a_n, \sigma) \circ (f_1, ..., f_i + g_i, ..., f_n)) \circ (x_1, ..., x_n)$$

$$= (a_0 \sigma(f_1) a_1 ... \sigma(f_i + g_i) ... a_{n-1} \sigma(f_n) a_n) \circ (x_1, ..., x_n)$$

$$= a_0 \sigma(f_1 \circ x_1) a_1 ... \sigma((f_i + g_i) \circ x_i) ... a_{n-1} \sigma(f_n \circ x_n) a_n$$

$$= a_0 \sigma(f_1 \circ x_1) a_1 ... \sigma(f_i \circ x_i + g_i \circ x_i) ... a_{n-1} \sigma(f_n \circ x_n) a_n$$

$$= a_0 \sigma(f_1 \circ x_1) a_1 ... \sigma(f_i \circ x_i) ... a_{n-1} \sigma(f_n \circ x_n) a_n$$

$$+ a_0 \sigma(f_1 \circ x_1) a_1 ... \sigma(g_i \circ x_i) ... a_{n-1} \sigma(f_n \circ x_n) a_n$$

$$= (a_0 \sigma(f_1) a_1 ... \sigma(f_i) ... a_{n-1} \sigma(f_n) a_n) \circ (x_1, ..., x_n)$$

$$+ (a_0 \sigma(f_1) a_1 ... \sigma(g_i) ... a_{n-1} \sigma(f_n) a_n) \circ (x_1, ..., x_n)$$

$$= ((a_0 \otimes ... \otimes a_n, \sigma) \circ (f_1, ..., f_i, ..., f_n)) \circ (x_1, ..., x_n)$$

$$+ ((a_0 \otimes ... \otimes a_n, \sigma) \circ (f_1, ..., g_i, ..., f_n)) \circ (x_1, ..., x_n)$$

$$= ((a_0 \otimes ... \otimes a_n, \sigma) \circ (f_1, ..., f_i, ..., f_n)$$

$$+ (a_0 \otimes ... \otimes a_n, \sigma) \circ (f_1, ..., g_i, ..., f_n)) \circ (x_1, ..., x_n)$$

$$((a_0 \otimes ... \otimes a_n, \sigma) \circ (f_1, ..., pf_i, ..., f_n)) \circ (x_1, ..., x_n)$$

$$= (a_0 \sigma(f_1) a_1, ... \sigma(pf_i) ... a_{n-1} \sigma(f_n) a_n) \circ (x_1, ..., x_n)$$

$$= a_0 \sigma(f_1 \circ x_1) a_1, ... \sigma((pf_i) \circ x_i) ... a_{n-1} \sigma(f_n \circ x_n) a_n$$

$$= a_0 \sigma(f_1 \circ x_1) a_1, ... \sigma(p(f_i \circ x_i)) ... a_{n-1} \sigma(f_n \circ x_n) a_n$$

$$= p(a_0 \sigma(f_1 \circ x_1) a_1, ... \sigma(f_i \circ x_i) ... a_{n-1} \sigma(f_n \circ x_n) a_n)$$

$$= p(((a_0 \otimes ... \otimes a_n, \sigma) \circ (f_1, ..., f_i, ..., f_n)) \circ (x_1, ..., x_n))$$

$$= (p((a_0 \otimes ... \otimes a_n, \sigma) \circ (f_1, ..., f_i, ..., f_n))) \circ (x_1, ..., x_n)$$

$\square$

Теорема 6.4.5. *Пусть* $A_1$, ..., $A_n$, $A$ - *ассоциативные* $D$-*алгебры. Для заданного тензора* $a \in A^{n+1\otimes}$ *и заданной перестановки* $\sigma \in S_n$ *существует линейное отображение*

$$h : \mathcal{L}(D; A_1 \to A) \otimes ... \otimes \mathcal{L}(D; A_n \to A) \to \mathcal{L}(D; A_1 \times ... \times A_n \to A)$$

*определённое равенством*

$$(6.4.3) \qquad (a_0 \otimes ... \otimes a_n, \sigma) \circ (f_1 \otimes ... \otimes f_n) = (a_0 \otimes ... \otimes a_n, \sigma) \circ (f_1, ..., f_n)$$

Доказательство. Утверждение теоремы является следствием теорем <span style="color:magenta">4.6.4</span>, <span style="color:magenta">6.4.4</span>.                    $\square$

# Векторное пространство над алгеброй с делением

## 7.1. Эффективное представление $D$-алгебры

Согласно теореме 3.3.4, рассматривая представление

$$f : A \longrightarrow V$$

$D$-алгебры $A$ в абелевой группе $V$, мы полагаем

(7.1.1) $$(f(a) + f(b))(x) = f(a)(x) + f(b)(x)$$

Согласно определению 2.3.1, отображение $f$ является гомоморфизмом $D$-алгебры $A$. Следовательно

(7.1.2) $$f(a + b) = f(a) + f(b)$$

Равенство

(7.1.3) $$f(a + b)(x) = f(a)(x) + f(b)(x)$$

является следствием равенств (7.1.1), (7.1.2).

**Теорема 7.1.1.** *Представление*

$$f : A \longrightarrow V$$

*$D$-алгебры $A$ в абелевой группе $V$ удовлетворяет равенству*

(7.1.4) $$f(0) = \overline{0}$$

*где*

$$\overline{0} : V \to V$$

*отображение такое, что*

(7.1.5) $$\overline{0} \circ v = 0$$

Доказательство. Равенство

(7.1.6) $$f(a)(x) = f(a + 0)(x) = f(a)(x) + f(0)(x)$$

является следствием равенства (7.1.3). Равенство

(7.1.7) $$f(0)(x) = 0$$

является следствием равенства (7.1.6). Равенство (7.1.5) является следствием равенств (7.1.4), (7.1.7). $\square$

**Теорема 7.1.2.** *Представление*

$$f : A \longrightarrow V$$

*D-алгебры A в абелевой группе V* **эффективно** *тогда и только тогда, когда из равенства f(a) = 0 следует a = 0.*

Доказательство. Если $a, b \in A$ порождают одно и то же преобразование, то

$$(7.1.8) \qquad\qquad f(a)(m) = f(b)(m)$$

для любого $m \in V$. Из равенств (7.1.3), (7.1.8) следует, что

$$(7.1.9) \qquad\qquad f(a - b)(m) = 0$$

Равенство

$$f(a - b) = \overline{0}$$

является следствием равенств (7.1.5), (7.1.9). Следовательно, представление $f$ эффективно тогда и только тогда, когда $a = b$. □

## 7.2. Левое векторное пространство над алгеброй с делением

Определение 7.2.1. *Пусть A - ассоциативная D-алгебра с делением. Пусть V - D-векторное пространство. Пусть в D-векторном пространстве* $\mathrm{End}(D, V)$ *определено произведение эндоморфизмов как композиция отображений. Пусть определён гомоморфизм*

$$(7.2.1) \qquad\qquad g_{34} : A \to \mathrm{End}(D, V)$$

*D-алгебры A в D-алгебру* $\mathrm{End}(D, V)$.

*Эффективное левостороннее представление*

$$(7.2.2) \qquad g_{34} : A \dashrightarrow V \quad g_{34}(a) : v \in V \to av \in V \quad a \in A$$

*D-алгебры A в D-векторном пространстве V называется* **левым векторным пространством** *над D-алгеброй A. Мы также будем говорить, что D-векторное пространство V является* **левым A-векторным пространством**. *V-число называется* **вектором**. *Билинейное отображение*

$$(7.2.3) \qquad\qquad (a, v) \in A \times V \to av \in V$$

*порождённое левосторонним представлением*

$$(7.2.4) \qquad\qquad (a, v) \to av$$

*называется левосторонним произведением вектора на скаляр.* □

Теорема 7.2.2. *Пусть D - коммутативное кольцо с единицей. Пусть A - ассоциативная D-алгебра с делением. Представление, порождённое левым сдвигом,*

$$f : A \dashrightarrow A \quad f(a) : b \in A \to ab \in A$$

*является левым векторным пространством* $\mathcal{L}A$ *над D-алгеброй A.*

Доказательство. Согласно определениям 7.2.1, 5.1.1, $D$-алгебра $A$ является абелевой группой. Из равенства

$$a(b + c) = ab + ac$$

и определений 2.1.7, 2.1.9 следует, что отображение

$$(7.2.5) \qquad f(a) : b \in A \to ab \in A$$

является эндоморфизмом абелевой группы $A$. Из равенств

$$(a_1 + a_2)b = a_1b + a_2b$$

$$(a_2a_1)b = a_2(a_1b)$$

и определений 2.1.7, 2.3.1 следует, что отображение $f$ является представлением $D$-алгебры $A$ в абелевой группе $A$. Теорема является следствием определения 7.2.1. □

ТЕОРЕМА 7.2.3. *Пусть $D$ - коммутативное кольцо с единицей. Пусть $A$ - ассоциативная $D$-алгебра с делением. Для любого множества $B$ представление*

$$f : A \longrightarrow A^B \qquad f(a) : g \in A^B \to ag \in A^B \qquad (ag)(b) = ag(b)$$

*является левым векторным пространством $\mathcal{L}A^B$ над $D$-алгеброй $A$.*

ДОКАЗАТЕЛЬСТВО. Согласно теореме 5.1.33, множество отображений $A^B$ является абелевой группой. Из равенства

$$a(h(b) + g(b)) = ah(b) + ag(b)$$

и определений 2.1.7, 2.1.9 следует, что отображение $f(a)$ является эндоморфизмом абелевой группы $A^B$. Из равенств

$$(a_1 + a_2)h(b) = a_1h(b) + a_2h(b)$$

$$(a_2a_1)h(b) = a_2(a_1h(b))$$

и определений 2.1.7, 2.3.1 следует, что отображение $f$ является представлением $D$-алгебры $A$ в абелевой группы $A^B$. Теорема является следствием определения 7.2.1. □

ТЕОРЕМА 7.2.4. *Следующая диаграмма представлений описывает левое $A$-векторное пространство $V$*

$$(7.2.6)$$

$$g_{12}(d) : a \to d\,a$$
$$g_{23}(v) : w \to C(w, v)$$
$$C \in \mathcal{L}(A^2 \to A)$$
$$g_{34}(a) : v \to a\,v$$
$$g_{14}(d) : v \to d\,v$$

*В диаграмме представлений* (7.2.6) *верна* **коммутативность представлений** *коммутативного кольца $D$ и $D$-алгебры $A$ в абелевой группе $V$*

$$(7.2.7) \qquad a(dv) = d(av)$$

ДОКАЗАТЕЛЬСТВО. Диаграмма представлений (7.2.6) является следствием определения 7.2.1 и теоремы 5.3.2. Равенство (7.2.7) является следствием утверждения, что левое преобразование $g_{3,4}(a)$ является эндоморфизмом $D$-векторного пространства $V$. □

Определение 7.2.5. *Подпредставление левого $A$-векторного простран­ства $V$ называется* **подпространством** *левого $A$-векторного пространства $V$.* □

Теорема 7.2.6. *Пусть $V$ является левым $A$-векторным пространством. $V$-числа удовлетворяют соотношениям*

**7.2.6.1: закон коммутативности**

$$(7.2.8) \qquad v + w = w + v$$

**7.2.6.2: закон ассоциативности**

$$(7.2.9) \qquad (pq)v = p(qv)$$

**7.2.6.3: закон дистрибутивности**

$$(7.2.10) \qquad p(v + w) = pv + pw$$

$$(7.2.11) \qquad (p + q)v = pv + qv$$

**7.2.6.4: закон унитарности**

$$(7.2.12) \qquad 1v = v$$

*для любых $p, q \in A, v, w \in V$.*

Доказательство. Равенство (7.2.8) является следствием определения 7.2.1.

Равенство (7.2.10) является следствием теоремы 3.3.2.

Равенство (7.2.11) является следствием теоремы 3.3.4.

Равенство (7.2.9) является следствием теоремы 2.1.10.

The equality (7.2.12) является следствием теорем 3.4.12, 7.2.4 так как пред­ставление $g_{3,4}$ является левосторонним представлением мультипликативной груп­пы $D$-алгебры $A$. □

### 7.3. Базис левого $A$-векторного пространства

Теорема 7.3.1. *Пусть $V$ - левое $A$-векторное пространство. Для любого множества $V$-чисел*

$$(7.3.1) \qquad v = (v(i) \in V, i \in I)$$

*вектор, порождённый диаграммой представлений (7.2.6), имеет следующий вид[7.1]*

$$(7.3.2) \qquad J(v) = \{w : w = c(i)v(i), c(i) \in A, |\{i : c(i) \neq 0\}| < \infty\}$$

Доказательство. Мы докажем теорему по индукции, опираясь на теоре­му 2.7.5.

---

[7.1] Для множества $A$, мы обозначим $|A|$ мощность множества $A$. Запись $|A| < \infty$ означает, что множество $A$ конечно.

Для произвольного $v(k) \in J(v)$ положим

$$c(i) = \begin{cases} 1, & i = k \\ 0 \end{cases}$$

Тогда

(7.3.3) $$v(k) = \sum_{i \in I} c(i)v(i) \quad c(i) \in A$$

Из равенств (7.3.2), (7.3.3). следует, что теорема верна на множестве $X_0 = v \subseteq J(v)$.

Пусть $X_{k-1} \subseteq J(v)$. Согласно определению 7.2.1 и теореме 2.7.5, если $w \in X_k$, то либо $w = w_1 + w_2$, $w_1$, $w_2 \in X_{k-1}$, либо $w = aw_1$, $a \in A$, $w_1 \in X_{k-1}$.

> **Лемма 7.3.2.** *Пусть* $w = w_1 + w_2$, $w_1$, $w_2 \in X_{k-1}$. *Тогда* $w \in J(v)$.

Доказательство. Согласно равенству (7.3.2), существуют $A$-числа $w_1(i)$, $w_2(i)$, $i \in I$, такие, что

(7.3.4) $$w_1 = \sum_{i \in I} w_1(i)v(i) \quad w_2 = \sum_{i \in I} w_2(i)v(i)$$

где множества

(7.3.5) $$H_1 = \{i \in I : w_1(i) \neq 0\} \quad H_2 = \{i \in I : w_2(i) \neq 0\}$$

конечны. Так как $V$ - левое $A$-векторное пространство, то из равенств (7.2.11), (7.3.4) следует, что

(7.3.6)
$$w_1 + w_2 = \sum_{i \in I} w_1(i)v(i) + \sum_{i \in I} w_2(i)v(i) = \sum_{i \in I}(w_1(i)v(i) + w_2(i)v(i))$$
$$= \sum_{i \in I}(w_1(i) + w_2(i))v(i)$$

Из равенства (7.3.5) следует, что множество

$$\{i \in I : w_1(i) + w_2(i) \neq 0\} \subseteq H_1 \cup H_2$$

конечно. $\odot$

> **Лемма 7.3.3.** *Пусть* $w = aw_1$, $a \in A$, $w_1 \in X_{k-1}$. *Тогда* $w \in J(v)$.

Доказательство. Согласно утверждению 2.7.5.4, для любого $A$-числа $a$,

(7.3.7) $$aw_1 \in X_k$$

Согласно равенству (7.3.2), существуют $A$-числа $w(i)$, $i \in I$, такие, что

(7.3.8) $$w_1 = \sum_{i \in I} w_1(i)v(i)$$

где

(7.3.9) $$|\{i \in I : w_1(i) \neq 0\}| < \infty$$

Из равенства (7.3.8) следует, что

$$(7.3.10) \qquad aw_1 = a\sum_{i \in I} w_1(i)v(i) = \sum_{i \in I} a(w_1(i)v(i)) = \sum_{i \in I}(aw_1(i))v(i)$$

Из утверждения (7.3.9) следует, что множество $\{i \in I : aw(i) \neq 0\}$ конечно.                                                                                   ⊙

Из лемм 7.3.2, 7.3.3 следует, что $X_k \subseteq J(v)$.                                □

СОГЛАШЕНИЕ 7.3.4. *Мы будем пользоваться соглашением о сумме, в котором повторяющийся индекс в линейной комбинации подразумевает сумму по повторяющемуся индексу. В этом случае предполагается известным множество индекса суммирования и знак суммы опускается*

$$c(i)v(i) = \sum_{i \in I} c(i)v(i)$$

*Я буду явно указывать множество индексов, если это необходимо.*        □

ОПРЕДЕЛЕНИЕ 7.3.5. *Пусть V - левостороннее векторное пространство . Пусть*

$$v = (v(i) \in V, i \in I)$$

*- множество векторов. Выражение $c(i)v(i)$, $c(i) \in A$,  называется **линейной комбинацией** векторов  $v(i)$.  Вектор  $w = c(i)v(i)$  называется **линейно зависимым** от векторов  $v(i)$.*                                         □

СОГЛАШЕНИЕ 7.3.6. *Если необходимо явно указать, что мы умножаем вектор $v(i)$ на $A$-число $c(i)$ слева, то мы будем называть выражение $c(i)v(i)$ **левой линейной комбинацией**. Мы будем применять это соглашение к аналогичным терминам. Например, мы будем говорить, что вектор  $w = c(i)v(i)$  линейно зависит от векторов $v(i)$, $i \in I$,  слева.*

ТЕОРЕМА 7.3.7. *Пусть  $v = (v_i \in V, i \in I)$  - множество векторов левого $A$-векторного пространства $V$. Если векторы $v_i$, $i \in I$,  принадлежат подпространству $V'$ левого $A$-векторного пространства $V$, то линейная комбинация векторов $v_i$, $i \in I$,  принадлежит подпространству $V'$.*

ДОКАЗАТЕЛЬСТВО. Теорема является следствием теоремы 7.3.1 и определений 7.3.5, 7.2.5.                                                                            □

ОПРЕДЕЛЕНИЕ 7.3.8. *$J(v)$ называется векторным подпространством, порождённым множеством $v$, а $v$ - **множеством образующих** векторного подпространства $J(v)$. В частности, **множеством образующих** левого $A$-векторного пространства $V$ будет такое подмножество $X \subset V$, что $J(X) = V$.*                                                                          □

Определение 7.3.9. *Если множество $X \subset V$ является множеством образующих левого $A$-векторного пространства $V$, то любое множество $Y$, $X \subset Y \subset V$ также является множеством образующих левого $A$-векторного пространства $V$. Если существует минимальное множество $X$, порождающее левое $A$-векторное пространство $V$, то такое множество $X$ называется **квазибазисом** левого $A$-векторного пространства $V$.* ☐

Определение 7.3.10. *Пусть $\overline{\overline{e}}$ - квазибазис левого $A$-векторного пространства $V$, и вектор $\overline{v} \in V$ имеет разложение*

$$(7.3.11) \qquad \overline{v} = v(i)e(i)$$

*относительно квазибазиса $\overline{\overline{e}}$. $A$-числа $v(i)$ называются **координатами** вектора $\overline{v}$ относительно квазибазиса $\overline{\overline{e}}$. Матрица $A$-чисел $v = (v(i), i \in I)$ называется **координатной матрицей вектора** $\overline{v}$ в квазибазисе $\overline{\overline{e}}$.* ☐

Теорема 7.3.11. *Пусть $A$ - ассоциативная $D$-алгебра с делением. Если уравнение*

$$(7.3.12) \qquad w(i)v(i) = 0$$

*предполагает существования индекса $i = j$ такого, что $w(j) \neq 0$, то вектор $v(j)$ линейно зависит от остальных векторов $v$.*

Доказательство. Теорема является следствием равенства

$$v(j) = \sum_{i \in I \setminus \{j\}} (w(j))^{-1} w(i) v(i)$$

и определения 7.3.5. ☐

Очевидно, что для любого множества векторов $v(i)$

$$w(i) = 0 \Rightarrow w(i)v(i) = 0$$

Определение 7.3.12. *Множество векторов [7.2] $v(i)$, $i \in I$, левого $A$-векторного пространства $V$ **линейно независимо**, если $c(i) = 0$, $i \in I$, следует из уравнения*

$$(7.3.13) \qquad c(i)v(i) = 0$$

*В противном случае, множество векторов $v(i)$, $i \in I$, **линейно зависимо**.* ☐

Теорема 7.3.13. *Пусть $A$ - ассоциативная $D$-алгебра с делением. Множество векторов $\overline{\overline{e}} = (e(i), i \in I)$ является **базисом левого $A$-векторного пространства** $V$, если векторы $e(i)$ линейно независимы и любой вектор $v \in V$ линейно зависит от векторов $e(i)$.*

---

[7.2] Я следую определению на странице [1]-100.

Доказательство. Пусть множество векторов $e(i)$, $i \in I$, линейно зависимо. Тогда в равенстве

$$w(i)e(i) = 0$$

существует индекс $i = j$ такой, что $w(j) \neq 0$. Согласно теореме 7.3.11, вектор $e(j)$ линейно зависит от остальных векторов множества $\overline{\overline{e}}$. Согласно определению 7.3.9, множество векторов $e(i)$, $i \in I$, не является базисом левого $A$-векторного пространства $V$.

Следовательно, если множество векторов $e(i)$, $i \in I$, является базисом, то эти векторы линейно независимы. Так как произвольный вектор $v \in V$ является линейной комбинацией векторов $e(i)$, $i \in I$, то множество векторов $v, e(i)$, $i \in I$, не является линейно независимым. $\square$

ТЕОРЕМА 7.3.14. *Координаты вектора $v \in V$ относительно базиса $\overline{\overline{e}}$ левого $A$-векторного пространства $V$ определены однозначно. Из равенства*

$$(7.3.14) \qquad ve = we \qquad v(i)e(i) = w(i)e(i)$$

*следует, что*

$$v = w \qquad v(i) = w(i)$$

Доказательство. Согласно теореме 7.3.13, система векторов $\overline{v}, e(i)$, $i \in I$, линейно зависима и в уравнении

$$(7.3.15) \qquad b\overline{v} + c(i)e(i) = 0$$

по крайней мере $b$ отлично от 0. Тогда равенство

$$(7.3.16) \qquad \overline{v} = (-b^{-1}c(i))e(i)$$

следует из (7.3.15). Равенство

$$(7.3.17) \qquad \overline{v} = v(i)e(i) \qquad v(i) = -b^{-1}c(i)$$

является следствием из равенства (7.3.16).

Допустим мы имеем другое разложение

$$(7.3.18) \qquad \overline{v} = v'(i)e(i)$$

Вычтя (7.3.17) из (7.3.18), мы получим

$$(7.3.19) \qquad 0 = (v'(i) - v(i))e(i)$$

Так как векторы $e(i)$ линейно независимы, то равенство

$$(7.3.20) \qquad v'(i) - v(i) = 0 \qquad i \in I$$

является следствием равенства (7.3.19). Утверждение теоремы является следствием равенства (7.3.20). $\square$

ТЕОРЕМА 7.3.15. *Мы можем рассматривать $D$-алгебру $A$ с делением как левое $A$-векторное пространство с базисом $\overline{\overline{e}} = \{1\}$.*

Доказательство. Теорема является следствием равенства (7.2.12). $\square$

## 7.4. Правое векторное пространство над алгеброй с делением

Определение 7.4.1. *Пусть $A$ - ассоциативная $D$-алгебра с делением. Пусть $V$ - $D$-векторное пространство. Пусть в $D$-векторном пространстве $\operatorname{End}(D, V)$ определено произведение эндоморфизмов как композиция отображений. Пусть определён гомоморфизм*

$$(7.4.1) \qquad g_{34} : A \to \operatorname{End}(D, V)$$

*$D$-алгебры $A$ в $D$-алгебру $\operatorname{End}(D, V)$.*

*Эффективное правостороннее представление*

$$(7.4.2) \qquad g_{34} : A \multimap V \quad g_{34}(a) : v \in V \to va \in V \quad a \in A$$

*$D$-алгебры $A$ в $D$-векторном пространстве $V$ называется* **правым векторным пространством** *над $D$-алгеброй $A$. Мы также будем говорить, что $D$-векторное пространство $V$ является* **правым $A$-векторным пространством**. *$V$-число называется* **вектором**. *Билинейное отображение*

$$(7.4.3) \qquad (v, a) \in V \times A \to va \in V$$

*порождённое правосторонним представлением*

$$(7.4.4) \qquad (v, a) \to va$$

*называется правосторонним произведением вектора на скаляр.* □

Теорема 7.4.2. *Пусть $D$ - коммутативное кольцо с единицей. Пусть $A$ - ассоциативная $D$-алгебра с делением. Представление, порождённое правым сдвигом,*

$$f : A \multimap A \quad f(a) : b \in A \to ba \in A$$

*является правым векторным пространством* $\mathcal{R}A$ *над $D$-алгеброй $A$.*

Доказательство. Согласно определениям 7.4.1, 5.1.1, $D$-алгебра $A$ является абелевой группой. Из равенства

$$(b + c)a = ba + ca$$

и определений 2.1.7, 2.1.9 следует, что отображение

$$(7.4.5) \qquad f(a) : b \in A \to ba \in A$$

является эндоморфизмом абелевой группы $A$. Из равенств

$$b(a_1 + a_2) = ba_1 + ba_2$$

$$b(a_1 a_2) = (ba_1)a_2$$

и определений 2.1.7, 2.3.1 следует, что отображение $f$ является представлением $D$-алгебры $A$ в абелевой группе $A$. Теорема является следствием определения 7.4.1. □

Теорема 7.4.3. *Пусть $D$ - коммутативное кольцо с единицей. Пусть $A$ - ассоциативная $D$-алгебра с делением. Для любого множества $B$ представление*

$$f : A \multimap A^B \quad f(a) : g \in A^B \to ga \in A^B \quad (ga)(b) = g(b)a$$

*является правым векторным пространством* $\mathcal{R}A^B$ *над D-алгеброй A.*

Доказательство. Согласно теореме 5.1.33, множество отображений $A^B$ является абелевой группой. Из равенства

$$(h(b) + g(b))a = h(b)a + g(b)a$$

и определений 2.1.7, 2.1.9 следует, что отображение $f(a)$ является эндоморфизмом абелевой группы $A^B$. Из равенств

$$h(b)(a_1 + a_2) = h(b)a_1 + h(b)a_2$$

$$h(b)(a_1 a_2) = (h(b)a_1)a_2$$

и определений 2.1.7, 2.3.1 следует, что отображение $f$ является представлением $D$-алгебры $A$ в абелевой группы $A^B$. Теорема является следствием определения 7.4.1.                                                 $\square$

Теорема 7.4.4. *Следующая диаграмма представлений описывает правое A-векторное пространство V*

(7.4.6)

$$g_{12}(d) : a \to d\,a$$
$$g_{23}(v) : w \to C(w, v)$$
$$C \in \mathcal{L}(A^2 \to A)$$
$$g_{34}(a) : v \to v\,a$$
$$g_{14}(d) : v \to v\,d$$

*В диаграмме представлений* (7.4.6) *верна* **коммутативность представлений** *коммутативного кольца D и D-алгебры A в абелевой группе V*

(7.4.7)                                          $$(vd)a = (va)d$$

Доказательство. Диаграмма представлений (7.4.6) является следствием определения 7.4.1 и теоремы 5.3.2. Равенство (7.4.7) является следствием утверждения, что правое преобразование $g_{3,4}(a)$ является эндоморфизмом $D$-векторного пространства $V$.                                                 $\square$

Определение 7.4.5. *Подпредставление правого A-векторного пространства V называется* **подпространством** *правого A-векторного пространства V.*                                                 $\square$

Теорема 7.4.6. *Пусть V является правым A-векторным пространством. V-числа удовлетворяют соотношениям*

7.4.6.1: **закон коммутативности**

(7.4.8)                                          $$v + w = w + v$$

7.4.6.2: **закон ассоциативности**

(7.4.9)                                          $$v(pq) = (vp)q$$

### 7.4.6.3: закон дистрибутивности

(7.4.10) $$(v + w)p = vp + wp$$

(7.4.11) $$v(p + q) = vp + vq$$

### 7.4.6.4: закон унитарности

(7.4.12) $$v1 = v$$

*для любых $p, q \in A$, $v, w \in V$.*

Доказательство. Равенство (7.4.8) является следствием определения 7.4.1.

Равенство (7.4.10) является следствием  теоремы 3.3.2.

Равенство (7.4.11) является следствием теоремы 3.3.4.

Равенство (7.4.9) является следствием теоремы 2.1.10.

The equality (7.4.12) является следствием теорем 3.4.12, 7.4.4 так как представление $g_{3,4}$ является правосторонним представлением мультипликативной группы $D$-алгебры $A$.  □

## 7.5. Базис правого $A$-векторного пространства

Теорема 7.5.1. *Пусть $V$ - правое $A$-векторное пространство. Для любого множества $V$-чисел*

(7.5.1) $$v = (v(i) \in V, i \in I)$$

*вектор, порождённый диаграммой представлений (7.4.6), имеет следующий вид[7.3]*

(7.5.2) $$J(v) = \{w : w = v(i)c(i), c(i) \in A, |\{i : c(i) \neq 0\}| < \infty\}$$

Доказательство. Мы докажем теорему по индукции, опираясь на теорему 2.7.5.

Для произвольного $v(k) \in J(v)$  положим

$$c(i) = \begin{cases} 1, & i = k \\ 0 & \end{cases}$$

Тогда

(7.5.3) $$v(k) = \sum_{i \in I} v(i)c(i) \quad c(i) \in A$$

Из равенств (7.5.2), (7.5.3). следует, что теорема верна на множестве $X_0 = v \subseteq J(v)$.

Пусть $X_{k-1} \subseteq J(v)$. Согласно определению  7.4.1  и теореме 2.7.5, если $w \in X_k$, то либо  $w = w_1 + w_2$, $w_1, w_2 \in X_{k-1}$,  либо  $w = w_1 a$, $a \in A$, $w_1 \in X_{k-1}$.

Лемма 7.5.2. *Пусть  $w = w_1 + w_2$, $w_1, w_2 \in X_{k-1}$.  Тогда $w \in J(v)$.*

---

[7.3] Для множества $A$, мы обозначим $|A|$ мощность множества $A$. Запись $|A| < \infty$ означает, что множество $A$ конечно.

Доказательство. Согласно равенству (7.5.2), существуют $A$-числа $w_1(i)$, $w_2(i)$, $i \in I$, такие, что

(7.5.4) $$w_1 = \sum_{i \in I} v(i) w_1(i) \quad w_2 = \sum_{i \in I} v(i) w_2(i)$$

где множества

(7.5.5) $$H_1 = \{ i \in I : w_1(i) \neq 0 \} \quad H_2 = \{ i \in I : w_2(i) \neq 0 \}$$

конечны. Так как $V$ - правое $A$-векторное пространство, то из равенств (7.4.11), (7.5.4) следует, что

(7.5.6) $$w_1 + w_2 = \sum_{i \in I} v(i) w_1(i) + \sum_{i \in I} v(i) w_2(i) = \sum_{i \in I} (v(i) w_1(i) + v(i) w_2(i))$$
$$= \sum_{i \in I} v(i) (w_1(i) + w_2(i))$$

Из равенства (7.5.5) следует, что множество

$$\{ i \in I : w_1(i) + w_2(i) \neq 0 \} \subseteq H_1 \cup H_2$$

конечно.                                                                      ⊙

**Лемма 7.5.3.** *Пусть* $w = w_1 a$, $a \in A$, $w_1 \in X_{k-1}$. *Тогда* $w \in J(v)$.

Доказательство. Согласно утверждению 2.7.5.4, для любого $A$-числа $a$,

(7.5.7) $$w_1 a \in X_k$$

Согласно равенству (7.5.2), существуют $A$-числа $w(i)$, $i \in I$, такие, что

(7.5.8) $$w_1 = \sum_{i \in I} v(i) w_1(i)$$

где

(7.5.9) $$|\{ i \in I : w_1(i) \neq 0 \}| < \infty$$

Из равенства (7.5.8) следует, что

(7.5.10) $$w a_1 = \left( \sum_{i \in I} v(i) w_1(i) \right) a = \sum_{i \in I} (v(i) w_1(i)) a = \sum_{i \in I} v(i) (w_1(i) a)$$

Из утверждения (7.5.9) следует, что множество $\{ i \in I : w(i) a \neq 0 \}$ конечно.                                                                      ⊙

Из лемм 7.5.2, 7.5.3 следует, что $X_k \subseteq J(v)$.                    □

**Соглашение 7.5.4.** *Мы будем пользоваться соглашением о сумме, в котором повторяющийся индекс в линейной комбинации подразумевает сумму по повторяющемуся индексу. В этом случае предполагается известным множество индекса суммирования и знак суммы опускается*

$$v(i) c(i) = \sum_{i \in I} v(i) c(i)$$

*Я буду явно указывать множество индексов, если это необходимо.*      □

Определение 7.5.5. *Пусть $V$ - правостороннее векторное пространство. Пусть*

$$v = (v(i) \in V, i \in I)$$

*- множество векторов. Выражение $v(i)c(i)$, $c(i) \in A$, называется* **линейной комбинацией** *векторов $v(i)$. Вектор $w = v(i)c(i)$ называется* **линейно зависимым** *от векторов $v(i)$.* □

Соглашение 7.5.6. *Если необходимо явно указать, что мы умножаем вектор $v(i)$ на $A$-число $c(i)$ справа, то мы будем называть выражение $v(i)c(i)$* **правой линейной комбинацией**. *Мы будем применять это соглашение к аналогичным терминам. Например, мы будем говорить, что вектор $w = v(i)c(i)$ линейно зависит от векторов $v(i)$, $i \in I$, справа.* □

Теорема 7.5.7. *Пусть $v = (v_i \in V, i \in I)$ - множество векторов правого $A$-векторного пространства $V$. Если векторы $v_i$, $i \in I$, принадлежат подпространству $V'$ правого $A$-векторного пространства $V$, то линейная комбинация векторов $v_i$, $i \in I$, принадлежит подпространству $V'$.*

Доказательство. Теорема является следствием теоремы 7.5.1 и определений 7.5.5, 7.4.5. □

Определение 7.5.8. *$J(v)$ называется векторным подпространством, порождённым множеством $v$, а $v$ -* **множеством образующих** *векторного подпространства $J(v)$. В частности,* **множеством образующих** *правого $A$-векторного пространства $V$ будет такое подмножество $X \subset V$, что $J(X) = V$.* □

Определение 7.5.9. *Если множество $X \subset V$ является множеством образующих правого $A$-векторного пространства $V$, то любое множество $Y$, $X \subset Y \subset V$ также является множеством образующих правого $A$-векторного пространства $V$. Если существует минимальное множество $X$, порождающее правое $A$-векторное пространство $V$, то такое множество $X$ называется* **квазибазисом** *правого $A$-векторного пространства $V$.* □

Определение 7.5.10. *Пусть $\overline{\overline{e}}$ - квазибазис правого $A$-векторного пространства $V$, и вектор $\overline{v} \in V$ имеет разложение*

$$(7.5.11) \qquad \overline{v} = e(i)v(i)$$

*относительно квазибазиса $\overline{\overline{e}}$. $A$-числа $v(i)$ называются* **координатами** *вектора $\overline{v}$ относительно квазибазиса $\overline{\overline{e}}$. Матрица $A$-чисел $v = (v(i), i \in I)$ называется* **координатной матрицей вектора** $\overline{v}$ *в квазибазисе $\overline{\overline{e}}$.* □

Теорема 7.5.11. *Пусть $A$ - ассоциативная $D$-алгебра с делением. Если уравнение*

$$(7.5.12) \qquad\qquad v(i)w(i) = 0$$

*предполагает существования индекса $i = j$ такого, что $w(j) \neq 0$, то вектор $v(j)$ линейно зависит от остальных векторов $v$.*

Доказательство. Теорема является следствием равенства

$$v(j) = \sum_{i \in I \setminus \{j\}} v(i)w(i)(w(j))^{-1}$$

и определения 7.5.5. $\qquad\qquad\qquad\qquad\qquad\qquad\qquad\qquad\qquad\qquad\square$

Очевидно, что для любого множества векторов $v(i)$

$$w(i) = 0 \Rightarrow v(i)w(i) = 0$$

Определение 7.5.12. *Множество векторов [7.4] $v(i)$, $i \in I$, правого $A$-векторного пространства $V$ **линейно независимо**, если $c(i) = 0$, $i \in I$, следует из уравнения*

$$(7.5.13) \qquad\qquad v(i)c(i) = 0$$

*В противном случае, множество векторов $v(i)$, $i \in I$, **линейно зависимо**.*
$\qquad\qquad\qquad\qquad\qquad\qquad\qquad\qquad\qquad\qquad\qquad\qquad\qquad\qquad\square$

Теорема 7.5.13. *Пусть $A$ - ассоциативная $D$-алгебра с делением. Множество векторов $\overline{\overline{e}} = (e(i), i \in I)$ является **базисом правого $A$-векторного пространства** $V$, если векторы $e(i)$ линейно независимы и любой вектор $v \in V$ линейно зависит от векторов $e(i)$.*

Доказательство. Пусть множество векторов $e(i)$, $i \in I$, линейно зависимо. Тогда в равенстве

$$e(i)w(i) = 0$$

существует индекс $i = j$ такой, что $w(j) \neq 0$. Согласно теореме 7.5.11, вектор $e(j)$ линейно зависит от остальных векторов множества $\overline{\overline{e}}$. Согласно определению 7.5.9, множество векторов $e(i)$, $i \in I$, не является базисом правого $A$-векторного пространства $V$.

Следовательно, если множество векторов $e(i)$, $i \in I$, является базисом, то эти векторы линейно независимы. Так как произвольный вектор $v \in V$ является линейной комбинацией векторов $e(i)$, $i \in I$, то множество векторов $v, e(i), i \in I$, не является линейно независимым. $\qquad\square$

Теорема 7.5.14. *Координаты вектора $v \in V$ относительно базиса $\overline{\overline{e}}$ правого $A$-векторного пространства $V$ определены однозначно. Из равенства*

$$(7.5.14) \qquad\qquad ev = ew \quad e(i)v(i) = e(i)w(i)$$

*следует, что*

$$v = w \quad v(i) = w(i)$$

---

[7.4] Я следую определению на страница [1]-100.

Доказательство. Согласно теореме 7.5.13, система векторов $\overline{v}$, $e(i)$, $i \in I$, линейно зависима и в уравнении

$$(7.5.15) \qquad \overline{v}b + e(i)c(i) = 0$$

по крайней мере $b$ отлично от 0. Тогда равенство

$$(7.5.16) \qquad \overline{v} = e(i)(-c(i)b^{-1})$$

следует из (7.5.15). Равенство

$$(7.5.17) \qquad \overline{v} = e(i)v(i) \quad v(i) = -c(i)b^{-1}$$

является следствием из равенства (7.5.16).

Допустим мы имеем другое разложение

$$(7.5.18) \qquad \overline{v} = e(i)v'(i)$$

Вычтя (7.5.17) из (7.5.18), мы получим

$$(7.5.19) \qquad 0 = e(i)(v'(i) - v(i))$$

Так как векторы $e(i)$ линейно независимы, то равенство

$$(7.5.20) \qquad v'(i) - v(i) = 0 \quad i \in I$$

является следствием равенства (7.5.19). Утверждение теоремы является следствием равенства (7.5.20). $\qquad\qquad\square$

# Тип векторного пространства

## 8.1. Бикольцо

Мы будем рассматривать матрицы, элементы которых принадлежат ассоциативной $D$-алгебре с делением $A$.

Согласно традиции произведение матриц $a$ и $b$ определено как произведение строк матрицы $a$ и столбцов матрицы $b$. В некоммутативной алгебре этой операции умножения недостаточно для решения некоторых задач.

Пример 8.1.1. Мы представим множество векторов $e(1) = e_1$, ..., $e(n) = e_n$ базиса $\overline{\overline{e}}$ левого векторного пространства $V$ над $D$-алгеброй $A$ (смотри определение 7.2.1 и теорему 7.3.13) в виде строки матрицы

$$(8.1.1) \qquad e = \begin{pmatrix} e_1 & ... & e_n \end{pmatrix}$$

Мы представим координаты $v(1) = v^1$, ..., $v(n) = v^n$ вектора $\overline{v} = v^i e_i$ в виде столбца матрицы

$$(8.1.2) \qquad v = \begin{pmatrix} v^1 \\ ... \\ v^n \end{pmatrix}$$

Поэтому мы можем представить вектор $\overline{v}$ как традиционное произведение матриц

$$(8.1.3)$$
$$v = \begin{pmatrix} e_1 & ... & e_n \end{pmatrix} \begin{pmatrix} v^1 \\ ... \\ v^n \end{pmatrix} = e_i v^i$$

□

Пример 8.1.2. Мы представим множество векторов $e(1) = e_1$, ..., $e(n) = e_n$ базиса $\overline{\overline{e}}$ правого векторного пространства $V$ над $D$-алгеброй $A$ (смотри определение 7.4.1 и теорему 7.3.13) в виде строки матрицы

$$(8.1.4) \qquad e = \begin{pmatrix} e_1 & ... & e_n \end{pmatrix}$$

Мы представим координаты $v(1) = v^1$, ..., $v(n) = v^n$ вектора $\overline{v} = v^i e_i$ в виде столбца матрицы

$$(8.1.5) \qquad v = \begin{pmatrix} v^1 \\ ... \\ v^n \end{pmatrix}$$

Однако мы не можем представить вектор $\overline{v}$ как традиционное произведение матриц
$$(8.1.6)$$
$$v = \begin{pmatrix} v^1 \\ ... \\ v^n \end{pmatrix} \qquad e = \begin{pmatrix} e_1 & ... & e_n \end{pmatrix}$$

так как это произведение не определено.

□

Из примеров 8.1.1, 8.1.2 следует, что мы не можем ограничиться традиционным произведением матриц и нам нужно определить два вида произведения матриц. Чтобы различать эти произведения, мы вводим новые обозначения.

ОПРЕДЕЛЕНИЕ 8.1.3. *Пусть число столбцов матрицы a равно числу строк матрицы b.* $_*$**-произведение** *матриц a и b имеет вид*

$$(8.1.7) \qquad a_*{}^*b = \left( a_k^i b_j^k \right)$$

$$(8.1.8) \qquad (a_*{}^*b)_j^i = a_k^i b_j^k$$

$$(8.1.9) \quad \begin{pmatrix} a_1^1 & ... & a_p^1 \\ ... & ... & ... \\ a_1^n & ... & a_p^n \end{pmatrix} {}_*{}^* \begin{pmatrix} b_1^1 & ... & b_m^1 \\ ... & ... & ... \\ b_1^p & ... & b_m^p \end{pmatrix} = \begin{pmatrix} a_k^1 b_1^k & ... & a_k^1 b_m^k \\ ... & ... & ... \\ a_k^n b_1^k & ... & a_k^n b_m^k \end{pmatrix}$$

$$= \begin{pmatrix} (a_*{}^*b)_1^1 & ... & (a_*{}^*b)_m^1 \\ ... & ... & ... \\ (a_*{}^*b)_1^n & ... & (a_*{}^*b)_m^n \end{pmatrix}$$

$_*$*-произведение может быть выражено как произведение строк матрицы a и столбцов матрицы b.* $\square$

ОПРЕДЕЛЕНИЕ 8.1.4. *Пусть число строк матрицы a равно числу столбцов матрицы b.* $^*{}_*$**-произведение** *матриц a и b имеет вид*

$$(8.1.10) \qquad a^*{}_*b = \left( a_i^k b_k^j \right)$$

$$(8.1.11) \qquad (a^*{}_*b)_j^i = a_i^k b_k^j$$

$$(8.1.12) \quad \begin{pmatrix} a_1^1 & ... & a_m^1 \\ ... & ... & ... \\ a_1^p & ... & a_m^p \end{pmatrix} {}^*{}_* \begin{pmatrix} b_1^1 & ... & b_p^1 \\ ... & ... & ... \\ b_1^n & ... & b_p^n \end{pmatrix} = \begin{pmatrix} a_1^k b_k^1 & ... & a_m^k b_k^1 \\ ... & ... & ... \\ a_1^k b_k^n & ... & a_m^k b_k^n \end{pmatrix}$$

$$= \begin{pmatrix} (a^*{}_*b)_1^1 & ... & (a^*{}_*b)_m^1 \\ ... & ... & ... \\ (a^*{}_*b)_1^n & ... & (a^*{}_*b)_m^n \end{pmatrix}$$

$^*{}_*$*-произведение может быть выражено как произведение столбцов матрицы a и строк матрицы b.* $\square$

Мы так же определим следующие операции на множестве матриц.

ОПРЕДЕЛЕНИЕ 8.1.5. *Транспонирование $a^T$ матрицы a меняет местами строки и столбцы*

$$(8.1.13) \qquad \left( a^T \right)_j^i = a_i^j$$

$\square$

ТЕОРЕМА 8.1.6. *Двойная транспозиция есть исходная матрица*

$$(8.1.14) \qquad (a^T)^T = a$$

ДОКАЗАТЕЛЬСТВО. Теорема является следствием определения 8.1.5.  □

ОПРЕДЕЛЕНИЕ 8.1.7. *Сумма матриц a и b определена равенством*

$$(8.1.15) \qquad (a+b)^i_j = a^i_j + b^i_j$$

□

ЗАМЕЧАНИЕ 8.1.8. *Мы будем пользоваться символом $_*{}^*$- или $^*{}_*$- в имени свойств каждого произведения и в обозначениях. Мы можем читать символ $_*{}^*$ как rc-произведение (произведение строки на столбец) и символ $^*{}_*$ как cr-произведение (произведение столбца на строку). Для совместимости обозначений с существующими мы будем иметь в виду $_*{}^*$-произведение, когда нет явных обозначений. Это правило мы распространим на последующую терминологию. Символ произведения сформирован из двух символов операции произведения, которые записываются на месте индекса суммирования. Например, если произведение A-чисел имеет вид $a \circ b$, то $_*{}^*$-произведение матриц a и b имеет вид $a_* \circ^* b$ и $^*{}_*$-произведение матриц a и b имеет вид $a^* \circ_* b$.* □

ТЕОРЕМА 8.1.9. *$_*{}^*$-Произведение матриц* **дистрибутивно** *относительно сложения*

$$(8.1.16) \quad a_*{}^*(b_1 + b_2) = a_*{}^*b_1 + a_*{}^*b_2$$

$$(8.1.17) \quad (b_1 + b_2)_*{}^*a = b_1{}_*{}^*a + b_2{}_*{}^*a$$

ДОКАЗАТЕЛЬСТВО. Равенство

$$(a_*{}^*(b_1 + b_2))^i_j$$
$$= a^i_k(b_1 + b_2)^k_j$$
$$(8.1.20) \qquad = a^i_k(b_1{}^k_j + b_2{}^k_j)$$
$$= a^i_k b_1{}^k_j + a^i_k b_2{}^k_j$$
$$= (a_*{}^*b_1)^i_j + (a_*{}^*b_2)^i_j$$

является следствием равенств (5.1.3), (8.1.7), (8.1.15). Равенство (8.1.16) является следствием равенства (8.1.20).

ТЕОРЕМА 8.1.10. *$^*{}_*$-Произведение матриц* **дистрибутивно** *относительно сложения*

$$(8.1.18) \quad a^*{}_*(b_1 + b_2) = a^*{}_*b_1 + a^*{}_*b_2$$

$$(8.1.19) \quad (b_1 + b_2)^*{}_*a = b_1{}^*{}_*a + b_2{}^*{}_*a$$

ДОКАЗАТЕЛЬСТВО. Равенство

$$(a^*{}_*(b_1 + b_2))^i_j$$
$$= a^k_j(b_1 + b_2)^i_k$$
$$(8.1.22) \qquad = a^k_j(b_1{}^i_k + b_2{}^i_k)$$
$$= a^k_j b_1{}^i_k + a^k_j b_2{}^i_k$$
$$= (a^*{}_*b_1)^i_j + (a^*{}_*b_2)^i_j$$

является следствием равенств (5.1.3), (8.1.10), (8.1.15). Равенство (8.1.18) является следствием равенства

Равенство

$$((b_1 + b_2)_*{}^*a)_j^i$$

$$= (b_1 + b_2)_k^i a_j^k$$

(8.1.21) $\qquad = (b_{1k}^i + b_{2k}^i)a_j^k$

$$= b_{1k}^i a_j^k + b_{2k}^i a_j^k$$

$$= (b_{1*}{}^*a)_j^i + (b_{2*}{}^*a)_j^i$$

является следствием равенств (5.1.2), (8.1.7), (8.1.15). Равенство (8.1.17) является следствием равенства (8.1.21).  □

(8.1.22). Равенство

$$((b_1 + b_2)^*{}_*a)_j^i$$

$$= (b_1 + b_2)_j^k a_k^i$$

(8.1.23) $\qquad = (b_{1j}^k + b_{2j}^k)a_k^i$

$$= b_{1j}^k a_k^i + b_{2j}^k a_k^i$$

$$= (b_1{}^*{}_*a)_j^i + (b_2{}^*{}_*a)_j^i$$

является следствием равенств (5.1.2), (8.1.10), (8.1.15). Равенство (8.1.19) является следствием равенства (8.1.23).  □

Теорема 8.1.11.

(8.1.24) $\qquad (a_*{}^*b)^T = a^{T}{}^*{}_*b^T$

Теорема 8.1.12.

(8.1.25) $\qquad (a^*{}_*b)^T = (a^T)_*{}^*(b^T)$

Доказательство теоремы 8.1.11. Цепочка равенств

(8.1.26) $\qquad ((a_*{}^*b)^T)_j^i = (a_*{}^*b)_i^j = a_i^k b_j^k = (a^T)_i^k (b^T)_k^j = ((a^T)^*{}_*(b^T))_i^j$

следует из (8.1.13), (8.1.7) и (8.1.10). Равенство (8.1.24) следует из (8.1.26).  □

Доказательство теоремы 8.1.12. Мы можем доказать (8.1.25) в случае матриц тем же образом, что мы доказали (8.1.24). Тем не менее для нас более важно показать, что (8.1.25) следует непосредственно из (8.1.24).

Применяя (8.1.14) к каждому слагаемому в левой части (8.1.25), мы получим

(8.1.27) $\qquad (a^*{}_*b)^T = ((a^T)^{T}{}^*{}_*(b^T)^T)^T$

Из (8.1.27) и (8.1.24) следует, что

(8.1.28) $\qquad (a^*{}_*b)^T = ((a^{T}{}^*{}_*b^T)^T)^T$

(8.1.25) следует из (8.1.28) и (8.1.14).  □

Определение 8.1.13. *Пусть*

$$a = \begin{pmatrix} a_1^1 & ... & a_n^1 \\ ... & ... & ... \\ a_1^n & ... & a_n^n \end{pmatrix}$$

*матрица A-чисел. Мы называем матрицу*[8.1]

(8.1.29) $\qquad \mathcal{H}a = \begin{pmatrix} (a_1^1)^{-1} & ... & (a_1^n)^{-1} \\ ... & ... & ... \\ (a_n^1)^{-1} & ... & (a_n^n)^{-1} \end{pmatrix}$

(8.1.30) $\qquad (\mathcal{H}a)_j^i = ((a^T)_j^i)^{-1} = (a_i^j)^{-1}$

**обращением Адамара матрицы** *a ([6]-page 4)*.  □

Множество $n \times n$ матриц замкнуто относительно $_*{}^*$-произведения и $^*{}_*$-произведения, а также относительно суммы.

**Определение 8.1.14.** **Бикольцо** $\mathcal{A}$ - *это множество матриц, на котором мы определили унарную операцию, называемую транспозицией, и три бинарных операции, называемые $_*{}^*$-произведение, $^*{}_*$-произведение и сумма, такие что*

- *$_*{}^*$-произведение и сумма определяют структуру кольца на $\mathcal{A}$*
- *$^*{}_*$-произведение и сумма определяют структуру кольца на $\mathcal{A}$*
- *оба произведения имеют общую единицу*

$$E_n = \left(\delta_j^i\right) \quad 1 \le i, j \le n$$

$\square$

**Теорема 8.1.15.** *Транспозиция единицы есть единица*

$$(8.1.31) \qquad E_n^T = E_n$$

**Доказательство.** Теорема следует из определений 8.1.5, 8.1.14. $\square$

**Определение 8.1.16.** *Мы определим $_*{}^*$-степень $n \times n$ матрицы a, пользуясь рекурсивным правилом*

$$(8.1.32) \qquad a^{0{}_*{}^*} = E_n$$

$$(8.1.33) \qquad a^{k{}_*{}^*} = a^{k-1{}_*{}^*}{}_*{}^* a$$

$\square$

**Определение 8.1.17.** *Мы определим $^*{}_*$-степень $n \times n$ матрицы a, пользуясь рекурсивным правилом*

$$(8.1.34) \qquad a^{0{}^*{}_*} = E_n$$

$$(8.1.35) \qquad a^{k{}^*{}_*} = a^{k-1{}^*{}_*}{}^*{}_* a$$

$\square$

**Теорема 8.1.18.**

$$(8.1.36) \qquad a^{1{}_*{}^*} = a$$

**Теорема 8.1.19.**

$$(8.1.38) \qquad a^{1{}^*{}_*} = a$$

**Доказательство.** Равенство

$$(8.1.37) \quad a^{1{}_*{}^*} = a^{0{}_*{}^*}{}_*{}^* a = E_n{}_*{}^* a = a$$

является следствием равенств (8.1.32), (8.1.33). Равенство (8.1.36) является следствием равенства (8.1.37). $\square$

**Доказательство.** Равенство

$$(8.1.39) \quad a^{1{}^*{}_*} = a^{0{}^*{}_*}{}^*{}_* a = E_n{}^*{}_* a = a$$

является следствием равенств (8.1.34), (8.1.35). Равенство (8.1.38) является следствием равенства (8.1.39). $\square$

**Теорема 8.1.20.**

$$(8.1.40) \qquad (a^T)^{n{}_*{}^*} = (a^{n{}^*{}_*})^T$$

**Теорема 8.1.21.**

$$(8.1.45) \qquad (a^T)^{n{}^*{}_*} = (a^{n{}_*{}^*})^T$$

---

8.1   Запись в равенстве (8.1.30) означает, что при обращении Адамара столбцы и строки меняются местами. Мы можем формально записать правую часть равенства (8.1.30) следующим образом

$$(a_i^j)^{-1} = \frac{1}{a_j^i}$$

Доказательство. Мы проведём доказательство индукцией по $n$.

При $n = 0$ утверждение непосредственно следует из равенств

| | |
|---|---|
| (8.1.32) | $a^{0\cdot^*} = E_n$ |
| (8.1.34) | $a^{0^*\cdot} = E_n$ |
| (8.1.31) | $E_n^T = E_n$ |

Допустим утверждение справедливо при $n = k - 1$

$$(8.1.41) \qquad (a^T)^{k-1\cdot^*} = (a^{k-1^*\cdot})^T$$

Из

| | |
|---|---|
| (8.1.33) | $a^{k\cdot^*} = a^{k-1\cdot^*}\cdot^* a$ |

следует

$$(8.1.42) \qquad (a^T)^{k\cdot^*} = (a^T)^{k-1\cdot^*}\cdot_* a^T$$

Из (8.1.42) и (8.1.41) следует

$$(8.1.43) \qquad (a^T)^{k\cdot^*} = (a^{k-1^*\cdot})^T\cdot_* a^T$$

Из (8.1.43) и (8.1.25) следует

$$(8.1.44) \qquad (a^T)^{k\cdot^*} = (a^{k-1^*\cdot}\cdot_* a)^T$$

Из (8.1.42) и (8.1.35) следует (8.1.40).  □

---

Определение 8.1.22. *Пусть $a$ - $n \times n$ матрица. Если существует $n \times n$ матрица $b$ такая, что*

$$(8.1.50) \qquad a_*\cdot b = E_n$$

*то матрица $b = a^{-1\cdot^*}$ называется $_*\cdot$-обратной матрицей матрицы $a$*

$$(8.1.51) \qquad a_*\cdot a^{-1\cdot^*} = E_n$$

□

---

Определение 8.1.23. *Если $n \times n$ матрица $a$ имеет $_*\cdot$-обратную матрицу, мы будем называть такую матрицу $_*\cdot$-невырожденной матрицей. В противном случае, мы будем называть такую матрицу $_*\cdot$-вырожденной матрицей.*  □

---

Теорема 8.1.26. *Предположим, что $n \times n$ матрица $a$ имеет $_*\cdot$-обратную матрицу. Тогда транспонирован-*

---

Доказательство. Мы проведём доказательство индукцией по $n$.

При $n = 0$ утверждение непосредственно следует из равенств

| | |
|---|---|
| (8.1.31) | $E_n^T = E_n$ |

Допустим утверждение справедливо при $n = k - 1$

$$(8.1.46) \qquad (a^T)^{k-1^*\cdot} = (a^{k-1\cdot^*})^T$$

Из

| | |
|---|---|
| (8.1.35) | $a^{k^*\cdot} = a^{k-1^*\cdot}\cdot_* a$ |

следует

$$(8.1.47) \qquad (a^T)^{k^*\cdot} = (a^T)^{k-1^*\cdot}\cdot_* a^T$$

Из (8.1.47) и (8.1.46) следует

$$(8.1.48) \qquad (a^T)^{k^*\cdot} = (a^{k-1\cdot^*})^T\cdot_* a^T$$

Из (8.1.48) и (8.1.25) следует

$$(8.1.49) \qquad (a^T)^{k^*\cdot} = (a^{k-1\cdot^*}\cdot_* a)^T$$

Из (8.1.47) и (8.1.33) следует (8.1.45).  □

---

Определение 8.1.24. *Пусть $a$ - $n \times n$ матрица. Если существует $n \times n$ матрица $b$ такая, что*

$$(8.1.52) \qquad a^*\cdot_* b = E_n$$

*то матрица $b = a^{-1^*\cdot}$ называется $^*\cdot_*$-обратной матрицей матрицы $a$*

$$(8.1.53) \qquad a^*\cdot_* a^{-1^*\cdot} = E_n$$

□

---

Определение 8.1.25. *Если $n \times n$ матрица $a$ имеет $^*\cdot_*$-обратную матрицу, мы будем называть такую матрицу $^*\cdot_*$-невырожденной матрицей. В противном случае, мы будем называть такую матрицу $^*\cdot_*$-вырожденной матрицей.*  □

---

Теорема 8.1.27. *Предположим, что $n \times n$ матрица $a$ имеет $^*\cdot_*$-обратную матрицу. Тогда транспонирован-*

ная матрица $a^T$ имеет $^*_*$-обратную матрицу и эти матрицы удовлетворяют равенству

(8.1.54)        $(a^T)^{-1^*_*} = (a^{-1^*_*})^T$

Доказательство. Если мы возьмём транспонирование обеих частей (8.1.51) и применим

(8.1.31)   $E_n^T = E_n$

мы получим

$(a_*^*a^{-1^*_*})^T = E_n^T = E_n$

Применяя

(8.1.24)   $(a_*^*b)^T = a^T{}^*_*b^T$

мы получим

(8.1.55)        $E_n = a^T{}^*_*(a^{-1^*_*})^T$

(8.1.54) следует из сравнения (8.1.53) и (8.1.55).   □

ная матрица $a^T$ имеет $_*^*$-обратную матрицу и эти матрицы удовлетворяют равенству

(8.1.56)        $(a^T)^{-1^*_*} = (a^{-1^*_*})^T$

Доказательство. Если мы возьмём транспонирование обеих частей (8.1.53) и применим

(8.1.31)   $E_n^T = E_n$

мы получим

$(a^*_*a^{-1^*_*})^T = E_n^T = E_n$

Применяя

(8.1.25)   $(a^*_*b)^T = (a^T)_*^*(b^T)$

мы получим

(8.1.57)        $E_n = a^T{}_*^*(a^{-1^*_*})^T$

(8.1.56) следует из сравнения (8.1.53) и (8.1.57).   □

Теоремы 8.1.11, 8.1.12, 8.1.20, 8.1.21, 8.1.26, 8.1.27 показывают, что существует двойственность между $_*^*$-произведением и $^*_*$-произведением. Мы можем объединить эти утверждения.

ТЕОРЕМА 8.1.28. **Принцип двойственности для бикольца.** *Пусть $\mathfrak{A}$ - истинное утверждение о бикольце $\mathcal{A}$. Если мы заменим одновременно*

- *$a \in \mathcal{A}$ и $a^T$*
- *$_*^*$-произведение и $^*_*$-произведение*

*то мы снова получим истинное утверждение.*

ТЕОРЕМА 8.1.29. **Принцип двойственности для бикольца матриц.** *Пусть $\mathcal{A}$ является бикольцом матриц. Пусть $\mathfrak{A}$ - истинное утверждение о матрицах. Если мы заменим одновременно*

- *столбцы и строки всех матриц*
- *$_*^*$-произведение и $^*_*$-произведение*

*то мы снова получим истинное утверждение.*

Доказательство. Непосредственное следствие теоремы 8.1.28.   □

ЗАМЕЧАНИЕ 8.1.30. *В выражении*

(8.1.58)        $a_*^*b_*^*c$

*мы выполняем операцию умножения слева направо. Однако мы можем выполнять операцию умножения справа налево. Тогда выражение (8.1.58) примет вид*

$$c^*_*b^*_*a$$

*Мы сохраним правило, что показатель степени и индексы записываются спра-*

*ва от выражения. Например, если исходное выражение имеет вид*

$$a^{-1^{*^*}}{}_**b^i$$

*то выражение, читаемое справа налево, примет вид*

$$b^i{}_**a^{-1^*{}^*}$$

Если задать порядок, в котором мы записываем индексы, то мы можем *утверждать, что мы читаем выражение сверху вниз, читая сперва верхние индексы, потом нижние. Мы можем прочесть выражение снизу вверх. При этом мы дополним правило, что символы операции также читаются в том же направлении, что и индексы. Например, выражение*

$$a^i{}_**b^{-1^*{}^*} = c^i$$

*прочтённое снизу вверх, имеет вид*

$$a_i{}^**b^{-1^*{}^*} = c_i$$

*Согласно принципу двойственности, если верно одно утверждение, то верно и другое.* □

---

ТЕОРЕМА 8.1.31. *Если матрица $a$ имеет $_*^*$-обратную матрицу, то для любых матриц $b$ и $c$ из равенства*

(8.1.59)      $b_*{}^*a = c_*{}^*a$

*следует равенство*

(8.1.60)          $b = c$

ДОКАЗАТЕЛЬСТВО. Равенство (8.1.60) следует из (8.1.59), если обе части равенства (8.1.59) умножить на $a^{-1^*{}^*}$. □

ТЕОРЕМА 8.1.33. *Если $n \times n$ матрицы $a$, $b$ имеют $_*^*$-обратную матрицу, то матрица $a_*{}^*b$ также имеет $_*^*$-обратную матрицу и*

(8.1.63)   $(a_*{}^*b)^{-1^*{}^*} = b^{-1^*{}^*}{}_**a^{-1^*{}^*}$

ДОКАЗАТЕЛЬСТВО. Равенство

$$(a_*{}^*b)_*{}^*(b^{-1^*{}^*}{}_**a^{-1^*{}^*})$$

(8.1.64)
$$= a_*{}^*(b_*{}^*b^{-1^*{}^*})_**a^{-1^*{}^*}$$
$$= a_*{}^*E_n{}_**a^{-1^*{}^*} = a_*{}^*a^{-1^*{}^*}$$
$$= E_n$$

является следствием равенства

(8.1.51)   $a_*{}^*a^{-1^*{}^*} = E_n$

и определения 8.1.14. Согласно определению 8.1.22 и равенству (8.1.64), мат-

---

ТЕОРЕМА 8.1.32. *Если матрица $a$ имеет $^*_*$-обратную матрицу, то для любых матриц $b$ и $c$ из равенства*

(8.1.61)      $b^*{}_*a = c^*{}_*a$

*следует равенство*

(8.1.62)          $b = c$

ДОКАЗАТЕЛЬСТВО. Равенство (8.1.62) следует из (8.1.61), если обе части равенства (8.1.61) умножить на $a^{-1^*{}^*}$. □

ТЕОРЕМА 8.1.34. *Если $n \times n$ матрицы $a$, $b$ имеют $^*_*$-обратную матрицу, то матрица $a^*{}_*b$ также имеет $^*_*$-обратную матрицу и*

(8.1.66)   $(a^*{}_*b)^{-1^*{}^*} = b^{-1^*{}^*{}^*}{}_*a^{-1^*{}^*}$

ДОКАЗАТЕЛЬСТВО. Равенство

$$(a^*{}_*b)^*{}_*(b^{-1^*{}^*{}^*}{}_*a^{-1^*{}^*})$$

(8.1.67)
$$= a^*{}_*(b^*{}_*b^{-1^*{}^*})^*{}_*a^{-1^*{}^*}$$
$$= a^*{}_*E_n{}^*{}_*a^{-1^*{}^*} = a^*{}_*a^{-1^*{}^*}$$
$$= E_n$$

является следствием равенства

(8.1.53)   $a^*{}_*a^{-1^*{}^*} = E_n$

и определения 8.1.14. Согласно определению 8.1.24 и равенству (8.1.67), мат-рица $a_*{}^*b$ имеет $_*{}^*$-обратную матрицу

$$(8.1.65) \qquad (a_*{}^*b)_*{}^* \left(a_*{}^*b\right)^{-1_*{}^*} = E_n$$

Равенство (8.1.63) является следствием равенств (8.1.64), (8.1.65). □

рица $a^*{}_*b$ имеет $^*{}_*$-обратную матрицу

$$(8.1.68) \qquad (a^*{}_*b)^*{}_* \left(a^*{}_*b\right)^{-1^*{}_*} = E_n$$

Равенство (8.1.66) является следствием равенств (8.1.67), (8.1.68). □

## 8.2. $_*{}^*$-Квазидетерминант

Теорема 8.2.1. *Предположим, что* $n \times n$ *матрица a имеет* $_*{}^*$*-обратную матрицу.* [8.2] *Тогда* $k \times k$ *подматрица* $_*{}^*$*-обратной матрицы удовлетворяет равенству*

$$(8.2.1) \qquad \left(\left(a^{-1_*{}^*}\right)^I_J\right)^{-1_*{}^*} = a^J_I - a^J_{[I]}{}_*{}^* \left(a^{[J]}_{[I]}\right)^{-1_*{}^*}{}_*{}^* a^{[J]}_I$$

$|I| = |J| = k < n$.

Доказательство. Мы можем записать определение

$$\boxed{(8.1.51) \quad a_*{}^* a^{-1_*{}^*} = E_n}$$

$_*{}^*$-обратной матрицы следующим образом

$$(8.2.2) \qquad \begin{pmatrix} a^J_I & a^J_{[I]} \\ a^{[J]}_I & a^{[J]}_{[I]} \end{pmatrix} {}_*{}^* \begin{pmatrix} \left(a^{-1_*{}^*}\right)^I_J & \left(a^{-1_*{}^*}\right)^I_{[J]} \\ \left(a^{-1_*{}^*}\right)^{[I]}_J & \left(a^{-1_*{}^*}\right)^{[I]}_{[J]} \end{pmatrix} = \begin{pmatrix} E_k & 0 \\ 0 & E_{n-k} \end{pmatrix}$$

Равенство $\left(a^{[J]}{}_*{}^*(a^{-1_*{}^*})_J = 0\right)$

$$(8.2.3) \qquad a^{[J]}_I{}_*{}^*\left(a^{-1_*{}^*}\right)^I_J + a^{[J]}_{[I]}{}_*{}^*\left(a^{-1_*{}^*}\right)^{[I]}_J = 0$$

и равенство $\left(a^J{}_*{}^*(a^{-1_*{}^*})_J = E_k\right)$

$$(8.2.4) \qquad a^J_I{}_*{}^*\left(a^{-1_*{}^*}\right)^I_J + a^J_{[I]}{}_*{}^*\left(a^{-1_*{}^*}\right)^{[I]}_J = E_k$$

являются следствием равенства (8.2.2) Мы умножим равенство (8.2.3) на $\left(a^{[J]}_{[I]}\right)^{-1_*{}^*}$

$$(8.2.5) \qquad \left(a^{[J]}_{[I]}\right)^{-1_*{}^*}{}_*{}^* a^{[J]}_I{}_*{}^*\left(a^{-1_*{}^*}\right)^I_J + \left(a^{-1_*{}^*}\right)^{[I]}_J = 0$$

Равенство

$$(8.2.6) \qquad a^J_I{}_*{}^*\left(a^{-1_*{}^*}\right)^I_J - a^J_{[I]}{}_*{}^* \left(a^{[J]}_{[I]}\right)^{-1_*{}^*}{}_*{}^* a^{[J]}_I{}_*{}^*\left(a^{-1_*{}^*}\right)^I_J = E_k$$

является следствием равенств (8.2.5), (8.2.4). Равенство (8.2.1) следует из (8.2.6), если мы обе части равенства (8.2.6) умножим на $\left(\left(a^{-1_*{}^*}\right)^I_J\right)^{-1_*{}^*}$. □

Теорема 8.2.2. *Предположим, что* $n \times n$ *матрица a имеет* $_*{}^*$*-обратную матрицу. Тогда элементы* $_*{}^*$*-обратной матрицы удовлетворяют равенству*

$$(8.2.7) \qquad \left(a^{-1_*{}^*}\right)^i_j = \left(a^j_i - a^j_{[i]}{}_*{}^* \left(a^{[j]}_{[i]}\right)^{-1_*{}^*}{}_*{}^* a^{[j]}_i\right)^{-1}$$

---

[8.2] Это утверждение и его доказательство основаны на утверждении 1.2.1 на странице [5]-page 8 для матриц над свободным кольцом с делением.

(8.2.8)
$$\left(\mathcal{H}a^{-1_{*}{}^{*}}\right)_{i}^{j} = a_{i}^{j} - a_{[i]*}^{j}{}^{*}\left(a_{[i]}^{[j]}\right)^{-1_{*}{}^{*}}{}_{*}{}^{*}a_{i}^{[j]}$$

Доказательство. Равенство (8.2.7) является следствием равенства (8.2.1), если мы положим $I = i$, $J = j$. Равенство (8.2.8) является следствием равенств (8.1.30), (8.2.7). $\qquad\square$

Определение 8.2.3. $\binom{j}{i}$-${}_{*}{}^{*}$-**квазидетерминант** $n \times n$ матрицы $a$ - это *формальное выражение*[8.3]

(8.2.9)
$$\det({}_{*}{}^{*})_{i}^{j}\,a = a_{i}^{j} - a_{[i]*}^{j}{}^{*}\left(a_{[i]}^{[j]}\right)^{-1_{*}{}^{*}}{}_{*}{}^{*}a_{i}^{[j]} = ((a^{-1_{*}{}^{*}})_{j}^{i})^{-1}$$

*Мы можем рассматривать $\binom{j}{i}$-${}_{*}{}^{*}$-квазидетерминант как элемент матрицы*

(8.2.10)
$$\det({}_{*}{}^{*})\,a = \begin{pmatrix} \det({}_{*}{}^{*})_{1}^{1}a & \dots & \det({}_{*}{}^{*})_{n}^{1}a \\ \dots & \dots & \dots \\ \det({}_{*}{}^{*})_{1}^{n}a & \dots & \det({}_{*}{}^{*})_{n}^{n}a \end{pmatrix}$$
$$= \begin{pmatrix} ((a^{-1_{*}{}^{*}})_{1}^{1})^{-1} & \dots & ((a^{-1_{*}{}^{*}})_{1}^{n})^{-1} \\ \dots & \dots & \dots \\ ((a^{-1_{*}{}^{*}})_{n}^{1})^{-1} & \dots & ((a^{-1_{*}{}^{*}})_{n}^{n})^{-1} \end{pmatrix}$$

*которую мы будем называть ${}_{*}{}^{*}$-**квазидетерминантом**.* $\qquad\square$

Теорема 8.2.4. *Выражение для ${}_{*}{}^{*}$-обратной матрицы имеет вид*

(8.2.11)
$$a^{-1_{*}{}^{*}} = \mathcal{H}\det({}_{*}{}^{*})a$$
$$\det({}_{*}{}^{*})\,a = \mathcal{H}a^{-1_{*}{}^{*}}$$

Доказательство. Равенство (8.2.11) следует из равенства

(8.2.12)
$$\det({}_{*}{}^{*})_{i}^{j}\,a = (\mathcal{H}a^{-1_{*}{}^{*}})_{i}^{j} = ((a^{-1_{*}{}^{*}})_{j}^{i})^{-1}$$
$\qquad\square$

Теорема 8.2.5. *Выражение для $\binom{j}{i}$-${}_{*}{}^{*}$-квазидетерминанта имеет вид*

(8.2.13)
$$\det({}_{*}{}^{*})_{i}^{j}\,a = a_{i}^{j} - a_{[i]*}^{j}{}^{*}\mathcal{H}\det({}_{*}{}^{*})_{[i]}^{[j]}\,a_{*}{}^{*}a_{i}^{[j]}$$

Доказательство. Равенство (8.2.13) является следствием равенств (8.2.8), (8.2.9). $\qquad\square$

Теорема 8.2.6.

(8.2.14)
$$\begin{pmatrix} \det({}_{*}{}^{*})_{1}^{1}a & \dots & \det({}_{*}{}^{*})_{n}^{1}a \\ \dots & \dots & \dots \\ \det({}_{*}{}^{*})_{1}^{n}a & \dots & \det({}_{*}{}^{*})_{n}^{n}a \end{pmatrix} = \begin{pmatrix} ((a^{-1_{*}{}^{*}})_{1}^{1})^{-1} & \dots & ((a^{-1_{*}{}^{*}})_{1}^{n})^{-1} \\ \dots & \dots & \dots \\ ((a^{-1_{*}{}^{*}})_{n}^{1})^{-1} & \dots & ((a^{-1_{*}{}^{*}})_{n}^{n})^{-1} \end{pmatrix}$$

---

[8.3] В определении 8.2.3, я следую определению [5]-1.2.2 на странице 9.

$$(8.2.15) \quad \begin{pmatrix} (\det({}_*{}^*)^1_1\, a)^{-1} & ... & (\det({}_*{}^*)^n_1\, a)^{-1} \\ ... & ... & ... \\ (\det({}_*{}^*)^1_n\, a)^{-1} & ... & (\det({}_*{}^*)^n_n\, a)^{-1} \end{pmatrix} = \begin{pmatrix} (a^{-1{}_*{}^*})^1_1 & ... & (a^{-1{}_*{}^*})^1_n \\ ... & ... & ... \\ (a^{-1{}_*{}^*})^n_1 & ... & (a^{-1{}_*{}^*})^n_n \end{pmatrix}$$

Доказательство. Равенство (8.2.14) является следствием равенства (8.2.9). Равенство (8.2.15) является следствием равенства (8.2.14). □

Теорема 8.2.7. *Рассмотрим матрицу*

$$a = \begin{pmatrix} a^1_1 & a^1_2 \\ a^2_1 & a^2_2 \end{pmatrix}$$

*Тогда*

$$(8.2.16) \qquad \det({}_*{}^*)a = \begin{pmatrix} a^1_1 - a^1_2(a^2_2)^{-1}a^2_1 & a^1_2 - a^1_1(a^2_1)^{-1}a^2_2 \\ a^2_1 - a^2_2(a^1_2)^{-1}a^1_1 & a^2_2 - a^2_1(a^1_1)^{-1}a^1_2 \end{pmatrix}$$

$$(8.2.17) \qquad a^{-1{}_*{}^*} = \begin{pmatrix} (a^1_1 - a^1_2(a^2_2)^{-1}a^2_1)^{-1} & (a^1_2 - a^1_1(a^2_1)^{-1}a^2_2)^{-1} \\ (a^2_1 - a^2_2(a^1_2)^{-1}a^1_1)^{-1} & (a^2_2 - a^2_1(a^1_1)^{-1}a^1_2)^{-1} \end{pmatrix}$$

Доказательство. Согласно равенству (8.2.9)

$$(8.2.18) \qquad\qquad \det({}_*{}^*)^1_1\, a = a^1_1 - a^1_2(a^2_2)^{-1}a^2_1$$

$$(8.2.19) \qquad\qquad \det({}_*{}^*)^2_1\, a = a^1_2 - a^1_1(a^2_1)^{-1}a^2_2$$

$$(8.2.20) \qquad\qquad \det({}_*{}^*)^1_2\, a = a^2_1 - a^2_2(a^1_2)^{-1}a^1_1$$

$$(8.2.21) \qquad\qquad \det({}_*{}^*)^2_2\, a = a^2_2 - a^2_1(a^1_1)^{-1}a^1_2$$

Равенство (8.2.16) является следствием равенств (8.2.18), (8.2.19), (8.2.20), (8.2.21). □

Теорема 8.2.8.

$$(8.2.22) \qquad\qquad (ma)^{-1{}_*{}^*} = a^{-1{}_*{}^*}m^{-1}$$

$$(8.2.23) \qquad\qquad (am)^{-1{}_*{}^*} = m^{-1}a^{-1{}_*{}^*}$$

Доказательство. Мы докажем равенство (8.2.22) индукцией по размеру матрицы.

Для $1 \times 1$ матрицы утверждение очевидно, так как

$$(ma)^{-1{}_*{}^*} = ((ma)^{-1}) = (a^{-1}m^{-1}) = (a^{-1})\, m^{-1} = a^{-1{}_*{}^*}m^{-1}$$

Допустим утверждение справедливо для $(n-1) \times (n-1)$ матрицы. Тогда из равенства (8.2.1) следует

$$(((ma)^{-1_*}{}^*)^I_J)^{-1_*}{}^* = (ma)^J_I - (ma)^J_{[I]}{}_*^* \left((ma)^{[J]}_{[I]}\right)^{-1_*}{}^* {}_*^* (ma)^{[J]}_I$$

$$= ma^J_I - m\, a^J_{[I]}{}_*^* \left(a^{[J]}_{[I]}\right)^{-1_*}{}^* m^{-1}{}_*^* m\, a^{[J]}_I$$

$$= m\, a^J_I - m\, a^J_{[I]}{}_*^* \left(a^{[J]}_{[I]}\right)^{-1_*}{}^* {}_*^* a^{[J]}_I$$

(8.2.24) $$(((ma)^{-1_*}{}^*)^I_J)^{-1_*}{}^* = m\, (a^{-1_*}{}^*)^I_J$$

Из равенства (8.2.24) следует равенство (8.2.22). Аналогично доказывается равенство (8.2.23). $\square$

ТЕОРЕМА 8.2.9. *Пусть*

(8.2.25) $$a = \begin{pmatrix} 1 & 0 \\ 0 & 1 \end{pmatrix}$$

*Тогда*

(8.2.26) $$a^{-1_*}{}^* = \begin{pmatrix} 1 & 0 \\ 0 & 1 \end{pmatrix}$$

Доказательство. Из (8.2.18), и (8.2.21) очевидно, что $(a^{-1_*}{}^*)^1_1 = 1$ и $(a^{-1_*}{}^*)^2_2 = 1$. Тем не менее выражение для $(a^{-1_*}{}^*)^2_1$ и $(a^{-1_*}{}^*)^1_2$ не может быть определено из (8.2.20) и (8.2.19) так как $a^2_1 = a^1_2 = 0$. Мы можем преобразовать эти выражения. Например

$$(a^{-1_*}{}^*)^2_1 = (a^1_2 - a^1_1(a^2_2)^{-1}a^2_2)^{-1}$$
$$= (a^1_1((a^1_1)^{-1}a^1_2 - (a^2_1)^{-1}a^2_2))^{-1}$$
$$= ((a^2_1)^{-1}a^1_1(a^1_2(a^1_1)^{-1}a^1_2 - a^2_2))^{-1}$$
$$= (a^1_1(a^1_2(a^1_1)^{-1}a^1_2 - a^2_2))^{-1}a^2_1$$

Мы непосредственно видим, что $(a^{-1_*}{}^*)^2_1 = 0$. Таким же образом мы можем найти, что $(a^{-1_*}{}^*)^1_2 = 0$. Это завершает доказательство (8.2.26). $\square$

## 8.3. *_*-Квазидетерминант

ТЕОРЕМА 8.3.1. *Предположим, что* $n \times n$ *матрица a имеет* *_*-*обратную матрицу.*[8.4]*Тогда* $k \times k$ *подматрица* *_*-*обратной матрицы удовлетворяет равенству*

(8.3.1) $$\left((a^{-1_*}{}^*)^J_I\right)^{-1_*}{}^* = a^I_J - a^{[I]}_J{}_*^* \left(a^{[I]}_{[J]}\right)^{-1_*}{}^* {}_*^* a^I_{[J]}$$

$|I| = |J| = k < n.$

---

Доказательство. Мы можем записать определение

$$(8.1.53) \quad a^* {}_* a^{-1^*} {}_* = E_n$$

$*_*$-обратной матрицы следующим образом

$$(8.3.2) \qquad \begin{pmatrix} a_J^J & a_{[J]}^J \\ a_J^{[I]} & a_{[J]}^{[I]} \end{pmatrix} *_* \begin{pmatrix} (a^{-1^*})_I^J & (a^{-1^*})_{[I]}^J \\ (a^{-1^*})_I^{[J]} & (a^{-1^*})_{[I]}^{[J]} \end{pmatrix} = \begin{pmatrix} E_k & 0 \\ 0 & E_{n-k} \end{pmatrix}$$

Равенство $(a_{[J]} *_* (a^{-1^*})^J = 0)$

$$(8.3.3) \qquad a_{[J]}^I *_* (a^{-1^*})_I^J + a_{[J]}^{[I]} *_* (a^{-1^*})_{[I]}^J = 0$$

и равенство $(a_J *_* (a^{-1^*})^J = E_k)$

$$(8.3.4) \qquad a_J^I *_* (a^{-1^*})_I^J + a_J^{[I]} *_* (a^{-1^*})_{[I]}^J = E_k$$

являются следствием равенства (8.3.2) Мы умножим равенство (8.3.3) на $\left(a_{[J]}^{[I]}\right)^{-1^*}{}_*$

$$(8.3.5) \qquad \left(a_{[J]}^{[I]}\right)^{-1^*}{}_* {}_* a_{[J]}^I *_* (a^{-1^*})_I^J + (a^{-1^*})_{[I]}^J = 0$$

Равенство

$$(8.3.6) \qquad a_J^I *_* (a^{-1^*})_I^J - a_J^{[I]} *_* \left(a_{[J]}^{[I]}\right)^{-1^*}{}_* {}_* a_{[J]}^I *_* (a^{-1^*})_I^J = E_k$$

является следствием равенств (8.3.5), (8.3.4). Равенство (8.3.1) следует из (8.3.6), если мы обе части равенства (8.3.6) умножим на $\left((a^{-1^*})_I^J\right)^{-1^*}{}_*$. $\qquad\square$

Теорема 8.3.2. *Предположим, что $n \times n$ матрица $a$ имеет $*_*$-обратную матрицу. Тогда элементы $*_*$-обратной матрицы удовлетворяют равенству*

$$(8.3.7) \qquad (a^{-1^*})_i^j = \left(a_j^i - a_j^{[i]} *_* \left(a_{[j]}^{[i]}\right)^{-1^*}{}_* {}_* a_{[j]}^i\right)^{-1}$$

$$(8.3.8) \qquad (\mathcal{H} a^{-1^*})_j^i = a_j^i - a_j^{[i]} *_* \left(a_{[j]}^{[i]}\right)^{-1^*}{}_* {}_* a_{[j]}^i$$

Доказательство. Равенство (8.3.7) является следствием равенства (8.3.1), если мы положим $I = i$, $J = j$. Равенство (8.3.8) является следствием равенств (8.1.30), (8.3.7). $\qquad\square$

Определение 8.3.3. $\binom{i}{j}$-$*_*$-**квазидетерминант** $n \times n$ матрицы $a$ - это формальное выражение[8.5]

$$(8.3.9) \qquad \det(*_*)_j^i a = a_j^i - a_j^{[i]} *_* \left(a_{[j]}^{[i]}\right)^{-1^*}{}_* {}_* a_{[j]}^i = ((a^{-1^*})_j^i)^{-1}$$

*Мы можем рассматривать $\binom{i}{j}$-$*_*$-квазидетерминант как элемент матрицы*

$$(8.3.10) \qquad \det(*_*) a = \begin{pmatrix} \det(*_*)_1^1 a & \dots & \det(*_*)_n^1 a \\ \dots & \dots & \dots \\ \det(*_*)_1^n a & \dots & \det(*_*)_n^n a \end{pmatrix}$$

$$= \begin{pmatrix} ((a^{-1^*})_1^1)^{-1} & \dots & ((a^{-1^*})_1^n)^{-1} \\ \dots & \dots & \dots \\ ((a^{-1^*})_n^1)^{-1} & \dots & ((a^{-1^*})_n^n)^{-1} \end{pmatrix}$$

которую мы будем называть $^*_*$-**квазидетерминантом**. □

ТЕОРЕМА 8.3.4. *Выражение для $^*_*$-обратной матрицы имеет вид*

$$(8.3.11) \qquad a^{-1^*_*} = \mathcal{H} \det(^*_*) a$$

$$\det(^*_*)\, a = \mathcal{H} a^{-1^*_*}$$

ДОКАЗАТЕЛЬСТВО. Равенство (8.3.11) следует из равенства

$$(8.3.12) \qquad \det(^*_*)^i_j\, a = (\mathcal{H} a^{-1^*_*})^i_j = ((a^{-1^*_*})^j_i)^{-1}$$

□

ТЕОРЕМА 8.3.5. *Выражение для $\binom{i}{j}$-$^*_*$-квазидетерминанта имеет вид*

$$(8.3.13) \qquad \det(^*_*)^i_j\, a = a^i_j - a^{[i]}_j{}^*_* \mathcal{H}\det(^*_*)^{[i]}_{[j]}\, a^*_*\, a^i_{[j]}$$

ДОКАЗАТЕЛЬСТВО. Равенство (8.3.13) является следствием равенств (8.3.8), (8.3.9). □

ТЕОРЕМА 8.3.6.

$$(8.3.14) \qquad \begin{pmatrix} \det(^*_*)^1_1 a & ... & \det(^*_*)^1_n a \\ ... & ... & ... \\ \det(^*_*)^n_1 a & ... & \det(^*_*)^n_n a \end{pmatrix} = \begin{pmatrix} ((a^{-1^*_*})^1_1)^{-1} & ... & ((a^{-1^*_*})^n_1)^{-1} \\ ... & ... & ... \\ ((a^{-1^*_*})^1_n)^{-1} & ... & ((a^{-1^*_*})^n_n)^{-1} \end{pmatrix}$$

$$(8.3.15) \qquad \begin{pmatrix} (\det(^*_*)^1_1 a)^{-1} & ... & (\det(^*_*)^n_1 a)^{-1} \\ ... & ... & ... \\ (\det(^*_*)^1_n a)^{-1} & ... & (\det(^*_*)^n_n a)^{-1} \end{pmatrix} = \begin{pmatrix} (a^{-1^*_*})^1_1 & ... & (a^{-1^*_*})^1_n \\ ... & ... & ... \\ (a^{-1^*_*})^n_1 & ... & (a^{-1^*_*})^n_n \end{pmatrix}$$

ДОКАЗАТЕЛЬСТВО. Равенство (8.3.14) является следствием равенства (8.3.9). Равенство (8.3.15) является следствием равенства (8.3.14). □

ТЕОРЕМА 8.3.7. *Рассмотрим матрицу*

$$a = \begin{pmatrix} a^1_1 & a^1_2 \\ a^2_1 & a^2_2 \end{pmatrix}$$

*Тогда*

$$(8.3.16) \qquad \det(^*_*) a = \begin{pmatrix} a^1_1 - a^2_1 (a^2_2)^{-1} a^1_2 & a^1_2 - a^2_2 (a^2_1)^{-1} a^1_1 \\ a^2_1 - a^1_1 (a^1_2)^{-1} a^2_2 & a^2_2 - a^1_2 (a^1_1)^{-1} a^2_1 \end{pmatrix}$$

$$(8.3.17) \qquad a^{-1^*_*} = \begin{pmatrix} (a^1_1 - a^2_1 (a^2_2)^{-1} a^1_2)^{-1} & (a^1_2 - a^1_1 (a^2_1)^{-1} a^2_2)^{-1} \\ (a^2_1 - a^2_2 (a^1_2)^{-1} a^1_1)^{-1} & (a^2_2 - a^1_2 (a^1_1)^{-1} a^2_1)^{-1} \end{pmatrix}$$

---

[8.5] В определении 8.3.3, я следую определению [5]-1.2.2 на странице 9.

Доказательство. Согласно равенству (8.3.9)

$$(8.3.18) \qquad \det(^*_*)^1_1 \, a = a^1_1 - a^2_1(a^2_2)^{-1}a^1_2$$

$$(8.3.19) \qquad \det(^*_*)^2_1 \, a = a^1_2 - a^2_2(a^2_1)^{-1}a^1_1$$

$$(8.3.20) \qquad \det(^*_*)^1_2 \, a = a^2_1 - a^1_1(a^1_2)^{-1}a^2_2$$

$$(8.3.21) \qquad \det(^*_*)^2_2 \, a = a^2_2 - a^1_2(a^1_1)^{-1}a^2_1$$

Равенство (8.3.16) является следствием равенств (8.3.18), (8.3.19), (8.3.20), (8.3.21).                                                                          □

Теорема 8.3.8.

$$(8.3.22) \qquad (ma)^{-1^*{}^*} = a^{-1^*{}^*}m^{-1}$$

$$(8.3.23) \qquad (am)^{-1^*{}^*} = m^{-1}a^{-1^*{}^*}$$

Доказательство. Мы докажем равенство (8.3.22) индукцией по размеру матрицы.

Для $1 \times 1$ матрицы утверждение очевидно, так как

$$(ma)^{-1^*{}^*} = \big((ma)^{-1}\big) = \big(a^{-1}m^{-1}\big) = \big(a^{-1}\big)\,m^{-1} = a^{-1^*{}^*}m^{-1}$$

Допустим утверждение справедливо для $(n-1) \times (n-1)$ матриц. Тогда из равенства (8.3.1) следует

$$(((ma)^{-1^*{}^*})^J_I)^{-1^*{}^*} = (ma)^I_J - (ma)^{[I]}_J{}_*\left((ma)^{[I]}_{[J]}\right)^{-1^*{}^*}{}_*(ma)^I_{[J]}$$

$$= ma^I_J - m\,a^{[I]}_J{}_*\left(a^{[I]}_{[J]}\right)^{-1^*{}^*}m^{-1}{}_*m\,a^I_{[J]}$$

$$= m\,a^I_J - m\,a^{[I]}_J{}_*\left(a^{[I]}_{[J]}\right)^{-1^*{}^*}{}_*a^I_{[J]}$$

$$(8.3.24) \qquad (((ma)^{-1^*{}^*})^J_I)^{-1^*{}^*} = m\,(a^{-1^*{}^*})^J_I$$

Из равенства (8.3.24) следует равенство (8.3.22). Аналогично доказывается равенство (8.3.23).                                                                          □

Теорема 8.3.9. Пусть

$$(8.3.25) \qquad a = \begin{pmatrix} 1 & 0 \\ 0 & 1 \end{pmatrix}$$

Тогда

$$(8.3.26) \qquad a^{-1^*{}^*} = \begin{pmatrix} 1 & 0 \\ 0 & 1 \end{pmatrix}$$

Доказательство. Из (8.3.18), и (8.3.21) очевидно, что $\left(a^{-1^*{}^*}\right)^1_1 = 1$ и $\left(a^{-1^*{}^*}\right)^2_2 = 1$. Тем не менее выражение для $\left(a^{-1^*{}^*}\right)^2_1$ и $\left(a^{-1^*{}^*}\right)^1_2$ не может быть определено из ($8.3.20$) и ($8.3.19$) так как $a_1^2 = a_2^1 = 0$. Мы можем преобразовать эти выражения. Например

$$(a^{-1*^*})_1^2 = (a_2^1 - a_2^2(a_1^2)^{-1}a_1^1)^{-1}$$
$$= (a_1^1((a_1^1)^{-1}a_2^1 - (a_1^2)^{-1}a_2^2))^{-1}$$
$$= ((a_1^2)^{-1}a_1^1(a_1^2(a_1^1)^{-1}a_2^1 - a_2^2))^{-1}$$
$$= (a_1^1(a_1^2(a_1^1)^{-1}a_2^1 - a_2^2))^{-1}a_1^2$$

Мы непосредственно видим, что $(a^{-1*^*})_1^2 = 0$. Таким же образом мы можем найти, что $(a^{-1*^*})_2^1 = 0$. Это завершает доказательство ($8.3.26$). □

Теорема 8.3.10.
$$(8.3.27) \qquad \det({}_*^*)_j^i \, a^T = \det({}^*_*)_i^j \, a$$

Доказательство. Согласно ($8.2.9$) и ($8.1.29$)
$$\det({}_*^*)_j^i \, a^T = (((a^T)^{-1*^*}){}_i^{\,-\,}{}_{\cdot\,-}^{\,j})^{-1}$$

Пользуясь теоремой $8.1.27$, мы получаем
$$\det({}_*^*)_j^i \, a^T = (((a^{-1*^*})^T){}_i^{\,-\,}{}_{\cdot\,-}^{\,j})^{-1}$$

Пользуясь ($8.1.13$), мы имеем
$$(8.3.28) \qquad \det({}_*^*)_j^i \, a^T = ((a^{-1*^*}){}_{\cdot}^{\,i\,}{}_{-\,}^{\,-}{}_j)^{-1}$$

Пользуясь ($8.3.28$), ($8.3.9$), мы получим ($8.3.27$). □

Теорема $8.3.10$ расширяет принцип двойственности, теорема $8.1.29$, на утверждения о квазидетерминантах и утверждает, что одно и тоже выражение является ${}_*^*$-квазидетерминантом матрицы $a$ и ${}^*_*$-квазидетерминантом матрицы $a^T$.

Теорема 8.3.11. **принцип двойственности.** *Пусть $\mathfrak{A}$ - истинное утверждение о бикольце матриц. Если мы одновременно заменим*
- *строку и столбец*
- *${}_*^*$-квазидетерминант и ${}^*_*$-квазидетерминант*

*то мы снова получим истинное утверждение.*

## 8.4. Тип векторного пространства

### 8.4.1. Левое $A$-векторное пространство столбцов.

Определение 8.4.1. *Мы представили множество векторов $v(1) = v_1$, ..., $v(m) = v_m$ в виде строки матрицы*
$$(8.4.1) \qquad v = \begin{pmatrix} v_1 & ... & v_m \end{pmatrix}$$
*и множество $A$-чисел $c(1) = c^1$, ..., $c(m) = c^m$ в виде столбца матрицы*
$$(8.4.2) \qquad c = \begin{pmatrix} c^1 \\ ... \\ c^m \end{pmatrix}$$

*Соответствующее представление левого $A$-векторного пространства $V$ мы будем называть* **левым $A$-векторным пространством столбцов**, *а $V$-число мы будем называть* **вектор-столбец** . □

**Теорема 8.4.2.** *Мы можем записать линейную комбинацию*

$$\sum_{i=1}^{m} c(i)v(i) = c^i v_i$$

*векторов $v_1$, ..., $v_m$ в виде ${}^*_*$-произведения матриц*

$$(8.4.3) \qquad c^*{}_*v = \begin{pmatrix} c^1 \\ ... \\ c^m \end{pmatrix} {}^*_* \begin{pmatrix} v_1 & ... & v_m \end{pmatrix} = c^i v_i$$

Доказательство. Теорема является следствием определений 7.3.5, 8.4.1.
□

**Теорема 8.4.3.** *Если мы запишем векторы базиса $\overline{\overline{e}}$ в виде строки матрицы*

$$(8.4.4) \qquad e = \begin{pmatrix} e_1 & ... & e_n \end{pmatrix}$$

*и координаты вектора $\overline{w} = w^i e_i$ относительно базиса $\overline{\overline{e}}$ в виде столбца матрицы*

$$(8.4.5) \qquad w = \begin{pmatrix} w^1 \\ ... \\ w^n \end{pmatrix}$$

*то мы можем представить вектор $\overline{w}$ в виде ${}^*_*$-произведения матриц*

$$(8.4.6) \qquad \overline{w} = w^*{}_*e = \begin{pmatrix} w^1 \\ ... \\ w^n \end{pmatrix} {}^*_* \begin{pmatrix} e_1 & ... & e_n \end{pmatrix} = w^i e_i$$

Доказательство. Теорема является следствием теоремы 8.4.2.     □

**Теорема 8.4.4.** *Координаты вектора $v \in V$ относительно базиса $\overline{\overline{e}}$ левого $A$-векторного пространства $V$ определены однозначно. Из равенства*

$$(8.4.7) \qquad v^*{}_*e = w^*{}_*e \qquad v^i e_i = w^i e_i$$

*следует, что*

$$v = w \qquad v^i = w^i$$

Доказательство. Теорема является следствием теоремы 7.3.14.     □

Теорема 8.4.5. *Пусть множество векторов* $\overline{\overline{e}}_A = (e_{Ak}, k \in K)$ *является базисом $D$-алгебры столбцов $A$. Пусть* $\overline{\overline{e}}_V = (e_{Vi}, i \in I)$ *базис левого $A$-векторного пространства столбцов $V$. Тогда множество $V$-чисел*

$$(8.4.8) \qquad \overline{\overline{e}}_{VA} = (e_{VAki} = e_{Ak}e_{Vi})$$

*является базисом $D$-векторного пространства $V$. Базис $\overline{\overline{e}}_{VA}$ называется расширением базиса* $(\overline{e}_A, \overline{e}_V)$.

Лемма 8.4.6. *Множество $V$-чисел $e_{2ik}$ порождает $D$-векторное пространство $V$.*

Доказательство. Пусть

$$(8.4.9) \qquad a = a^i e_{1i}$$

произвольное $V$-число. Для любого $A$-числа $a^i$, согласно определению 5.1.1 и соглашению 5.2.1, существуют $D$-числа $a^{ik}$ такие, что

$$(8.4.10) \qquad a^i = a^{ik}e_k$$

Равенство

$$(8.4.11) \qquad a = a^{ik}e_k e_{1i} = a^{ik}e_{2ik}$$

является следствием равенств (8.4.8), (8.4.9), (8.4.10). Согласно теореме 4.1.9, лемма является следствием равенства (8.4.11). $\qquad \odot$

Лемма 8.4.7. *Множество $V$-чисел $e_{2ik}$ линейно независимо над кольцом $D$.*

Доказательство. Пусть

$$(8.4.12) \qquad a^{ik}e_{2ik} = 0$$

Равенство

$$(8.4.13) \qquad a^{ik}e_k e_{1i} = 0$$

является следствием равенств (8.4.8), (8.4.12). Согласно теореме 8.4.4, равенство

$$(8.4.14) \qquad a^{ik}e_k = 0 \quad k \in K$$

является следствием равенства (8.4.13). Согласно теореме 4.2.7, равенство

$$(8.4.15) \qquad a^{ik} = 0 \quad k \in K \quad i \in I$$

является следствием равенства (8.4.14). Следовательно, множество $V$-чисел $e_{2ik}$ линейно независимо над кольцом $D$. $\qquad \odot$

Доказательство теоремы 8.4.5. Теорема является следствием теоремы 4.2.7, и лемм 8.4.6, 8.4.7. $\qquad \square$

### 8.4.2. Левое $A$-векторное пространство строк.

Определение 8.4.8. *Мы представили множество векторов* $v(1) = v^1$, ..., $v(m) = v^m$ *в виде столбца матрицы*

$$(8.4.16) \qquad v = \begin{pmatrix} v^1 \\ ... \\ v^m \end{pmatrix}$$

*и множество $A$-чисел* $c(1) = c_1$, ..., $c(m) = c_m$ *в виде строки матрицы*

$$(8.4.17) \qquad c = \begin{pmatrix} c_1 & ... & c_m \end{pmatrix}$$

*Соответствующее представление левого $A$-векторного пространства $V$ мы будем называть* **левым $A$-векторным пространством строк**, *а $V$-число мы будем называть* **вектор-строка** . $\square$

Теорема 8.4.9. *Мы можем записать линейную комбинацию*

$$\sum_{i=1}^{m} c(i)v(i) = c_i v^i$$

*векторов* $v^1$, ..., $v^m$ *в виде $_*{}^*$-произведения матриц*

$$(8.4.18) \qquad c_*{}^* v = \begin{pmatrix} c_1 & ... & c_m \end{pmatrix} {}_*{}^* \begin{pmatrix} v^1 \\ ... \\ v^m \end{pmatrix} = c_i v^i$$

Доказательство. Теорема является следствием определений 7.3.5, 8.4.8. $\square$

Теорема 8.4.10. *Если мы запишем векторы базиса $\overline{\overline{e}}$ в виде столбца матрицы*

$$(8.4.19) \qquad e = \begin{pmatrix} e^1 \\ ... \\ e^n \end{pmatrix}$$

*и координаты вектора* $\overline{w} = w_i e^i$ *относительно базиса $\overline{\overline{e}}$ в виде строки матрицы*

$$(8.4.20) \qquad w = \begin{pmatrix} w_1 & ... & w_n \end{pmatrix}$$

*то мы можем представить вектор $\overline{w}$ в виде $_*{}^*$-произведения матриц*

$$(8.4.21) \qquad \overline{w} = w_*{}^* e = \begin{pmatrix} w_1 & ... & w_n \end{pmatrix} {}_*{}^* \begin{pmatrix} e^1 \\ ... \\ e^n \end{pmatrix} = w_i e^i$$

Доказательство. Теорема является следствием теоремы 8.4.9. $\square$

Теорема 8.4.11. *Координаты вектора $v \in V$ относительно базиса $\overline{\overline{e}}$ левого $A$-векторного пространства $V$ определены однозначно. Из равенства*

$$(8.4.22) \qquad v_*{}^*e = w_*{}^*e \qquad v_i e^i = w_i e^i$$

*следует, что*

$$v = w \qquad v_i = w_i$$

Доказательство. Теорема является следствием теоремы 7.3.14. □

Теорема 8.4.12. *Пусть множество векторов $\overline{\overline{e}}_A = (e_A^k, k \in K)$ является базисом $D$-алгебры строк $A$. Пусть $\overline{\overline{e}}_V = (e_V^i, i \in I)$ базис левого $A$-векторного пространства строк $V$. Тогда множество $V$-чисел*

$$(8.4.23) \qquad \overline{\overline{e}}_{VA} = (e_{VA}^{ki} = e_A^k e_V^i)$$

*является базисом $D$-векторного пространства $V$. Базис $\overline{\overline{e}}_{VA}$ называется расширением базиса $(\overline{\overline{e}}_A, \overline{\overline{e}}_V)$.*

Лемма 8.4.13. *Множество $V$-чисел $e_2^{ik}$ порождает $D$-векторное пространство $V$.*

Доказательство. Пусть

$$(8.4.24) \qquad a = a_i e_1^i$$

произвольное $V$-число. Для любого $A$-числа $a_i$, согласно определению 5.1.1 и соглашению 5.2.1, существуют $D$-числа $a_{ik}$ такие, что

$$(8.4.25) \qquad a_i = a_{ik} e^k$$

Равенство

$$(8.4.26) \qquad a = a_{ik} e^k e_1^i = a_{ik} e_2^{ik}$$

является следствием равенств (8.4.23), (8.4.24), (8.4.25). Согласно теореме 4.1.9, лемма является следствием равенства (8.4.26). ⊙

Лемма 8.4.14. *Множество $V$-чисел $e_2^{ik}$ линейно независимо над кольцом $D$.*

Доказательство. Пусть

$$(8.4.27) \qquad a_{ik} e_2^{ik} = 0$$

Равенство

$$(8.4.28) \qquad a_{ik} e^k e_1^i = 0$$

является следствием равенств (8.4.23), (8.4.27). Согласно теореме 8.4.11, равенство

$$(8.4.29) \qquad a_{ik} e^k = 0 \quad k \in K$$

является следствием равенства (8.4.28). Согласно теореме 4.2.14, равенство

$$(8.4.30) \qquad a_{ik} = 0 \quad k \in K \quad i \in I$$

является следствием равенства (8.4.29). Следовательно, множество $V$-чисел $e_2^{ik}$ линейно независимо над кольцом $D$.      $\odot$

Доказательство теоремы 8.4.12. Теорема является следствием теоремы 4.2.14, и лемм 8.4.13, 8.4.14.      $\square$

### 8.4.3. Правое $A$-векторное пространство столбцов.

Определение 8.4.15. *Мы представили множество векторов* $v(1) = v_1$, ..., $v(m) = v_m$ *в виде строки матрицы*

$$(8.4.31) \qquad v = \begin{pmatrix} v_1 & ... & v_m \end{pmatrix}$$

*и множество $A$-чисел* $c(1) = c^1$, ..., $c(m) = c^m$ *в виде столбца матрицы*

$$(8.4.32) \qquad c = \begin{pmatrix} c^1 \\ ... \\ c^m \end{pmatrix}$$

*Соответствующее представление правого $A$-векторного пространства $V$ мы будем называть* **правым $A$-векторным пространством столбцов**, *а $V$-число мы будем называть* **вектор-столбец** *.*      $\square$

Теорема 8.4.16. *Мы можем записать линейную комбинацию*

$$\sum_{i=1}^{m} v(i)c(i) = v_i c^i$$

*векторов* $v_1$, ..., $v_m$ *в виде* $_*{}^*$-*произведения матриц*

$$(8.4.33) \qquad v_* c = \begin{pmatrix} v_1 & ... & v_m \end{pmatrix} {}_*{}^* \begin{pmatrix} c^1 \\ ... \\ c^m \end{pmatrix} = v_i c^i$$

Доказательство. Теорема является следствием определений 7.5.5, 8.4.15.      $\square$

Теорема 8.4.17. *Если мы запишем векторы базиса $\overline{\overline{e}}$ в виде строки матрицы*

$$(8.4.34) \qquad e = \begin{pmatrix} e_1 & ... & e_n \end{pmatrix}$$

*и координаты вектора* $\overline{w} = e_i w^i$ *относительно базиса $\overline{\overline{e}}$ в виде столбца матрицы*

$$(8.4.35) \qquad w = \begin{pmatrix} w^1 \\ ... \\ w^n \end{pmatrix}$$

*то мы можем представить вектор $\overline{w}$ в виде $_*{}^*$-произведения матриц*

$$(8.4.36) \qquad \overline{w} = e_*{}^* w = \begin{pmatrix} e_1 & ... & e_n \end{pmatrix}{}_*{}^* \begin{pmatrix} w^1 \\ ... \\ w^n \end{pmatrix} = e_i w^i$$

Доказательство. Теорема является следствием теоремы 8.4.16.    □

Теорема 8.4.18. *Координаты вектора $v \in V$ относительно базиса $\overline{\overline{e}}$ правого $A$-векторного пространства $V$ определены однозначно. Из равенства*

$$(8.4.37) \qquad e_*{}^* v = e_*{}^* w \qquad e_i v^i = e_i w^i$$

*следует, что*

$$v = w \qquad v^i = w^i$$

Доказательство. Теорема является следствием теоремы 7.5.14.    □

Теорема 8.4.19. *Пусть множество векторов $\overline{\overline{e}}_A = (e_{Ak}, k \in K)$ является базисом $D$-алгебры столбцов $A$. Пусть $\overline{\overline{e}}_V = (e_{Vi}, i \in I)$ базис правого $A$-векторного пространства столбцов $V$. Тогда множество $V$-чисел*

$$(8.4.38) \qquad \overline{\overline{e}}_{VA} = (e_{VAki} = e_{Vi} e_{Ak})$$

*является базисом $D$-векторного пространства $V$. Базис $\overline{\overline{e}}_{VA}$ называется расширением базиса $(\overline{\overline{e}}_A, \overline{\overline{e}}_V)$.*

Лемма 8.4.20. *Множество $V$-чисел $e_{2ik}$ порождает $D$-векторное пространство $V$.*

Доказательство. Пусть

$$(8.4.39) \qquad a = e_{1i} a^i$$

произвольное $V$-число. Для любого $A$-числа $a^i$, согласно определению 5.1.1 и соглашению 5.2.1, существуют $D$-числа $a^{ik}$ такие, что

$$(8.4.40) \qquad a^i = e_k a^{ik}$$

Равенство

$$(8.4.41) \qquad a = e_{1i} e_k a^{ik} = e_{2ik} a^{ik}$$

является следствием равенств (8.4.38), (8.4.39), (8.4.40). Согласно теореме 4.1.9, лемма является следствием равенства (8.4.41).    ⊙

Лемма 8.4.21. *Множество $V$-чисел $e_{2ik}$ линейно независимо над кольцом $D$.*

Доказательство. Пусть

$$(8.4.42) \qquad e_{2ik} a^{ik} = 0$$

Равенство

$$(8.4.43) \qquad\qquad\qquad e_{1i}e_k a^{ik} = 0$$

является следствием равенств (8.4.38), (8.4.42). Согласно теореме 8.4.18, равенство

$$(8.4.44) \qquad\qquad\qquad e_k a^{ik} = 0 \quad k \in K$$

является следствием равенства (8.4.43). Согласно теореме 4.2.7, равенство

$$(8.4.45) \qquad\qquad\qquad a^{ik} = 0 \quad k \in K \quad i \in I$$

является следствием равенства (8.4.44). Следовательно, множество $V$-чисел $e_{2ik}$ линейно независимо над кольцом $D$.                                            ⊙

Доказательство теоремы 8.4.19.   Теорема является следствием теоремы 4.2.7, и лемм 8.4.20, 8.4.21.                                                    □

### 8.4.4. Правое $A$-векторное пространство строк.

Определение 8.4.22.  *Мы представили множество векторов* $v(1) = v^1$, ..., $v(m) = v^m$  *в виде столбца матрицы*

$$(8.4.46) \qquad\qquad\qquad v = \begin{pmatrix} v^1 \\ ... \\ v^m \end{pmatrix}$$

*и множество $A$-чисел* $c(1) = c_1$, ..., $c(m) = c_m$  *в виде строки матрицы*

$$(8.4.47) \qquad\qquad\qquad c = \begin{pmatrix} c_1 & ... & c_m \end{pmatrix}$$

*Соответствующее представление правого $A$-векторного пространства $V$ мы будем называть* **правым $A$-векторным пространством строк**, *а $V$-число мы будем называть* **вектор-строка** .                                            □

Теорема 8.4.23.  *Мы можем записать линейную комбинацию*

$$\sum_{i=1}^{m} v(i)c(i) = v^i c_i$$

*векторов* $v^1$, ..., $v^m$  *в виде* $^*_*$-*произведения матриц*

$$(8.4.48) \qquad\qquad v {^*_*} c = \begin{pmatrix} v^1 \\ ... \\ v^m \end{pmatrix} {^*_*} \begin{pmatrix} c_1 & ... & c_m \end{pmatrix} = v^i c_i$$

Доказательство. Теорема является следствием определений 7.5.5, 8.4.22.
                                                                              □

Теорема 8.4.24. *Если мы запишем векторы базиса $\overline{\overline{e}}$ в виде столбца матрицы*

$$(8.4.49) \qquad e = \begin{pmatrix} e^1 \\ ... \\ e^n \end{pmatrix}$$

*и координаты вектора $\overline{w} = e^i w_i$ относительно базиса $\overline{\overline{e}}$ в виде строки матрицы*

$$(8.4.50) \qquad w = \begin{pmatrix} w_1 & ... & w_n \end{pmatrix}$$

*то мы можем представить вектор $\overline{w}$ в виде $^*_*$-произведения матриц*

$$(8.4.51) \qquad \overline{w} = e^*_* w = \begin{pmatrix} e^1 \\ ... \\ e^n \end{pmatrix} {}^*_* \begin{pmatrix} w_1 & ... & w_n \end{pmatrix} = e^i w_i$$

Доказательство. Теорема является следствием теоремы 8.4.23. □

Теорема 8.4.25. *Координаты вектора $v \in V$ относительно базиса $\overline{\overline{e}}$ правого $A$-векторного пространства $V$ определены однозначно. Из равенства*

$$(8.4.52) \qquad e^*_* v = e^*_* w \qquad e^i v_i = e^i w_i$$

*следует, что*

$$v = w \qquad v_i = w_i$$

Доказательство. Теорема является следствием теоремы 7.5.14. □

Теорема 8.4.26. *Пусть множество векторов $\overline{\overline{e}}_A = (e_A^k, k \in K)$ является базисом $D$-алгебры строк $A$. Пусть $\overline{\overline{e}}_V = (e_V^i, i \in I)$ базис правого $A$-векторного пространства строк $V$. Тогда множество $V$-чисел*

$$(8.4.53) \qquad \overline{\overline{e}}_{VA} = (e_{VA}^{ki} = e_V^i e_A^k)$$

*является базисом $D$-векторного пространства $V$. Базис $\overline{\overline{e}}_{VA}$ называется расширением базиса $(\overline{e}_A, \overline{e}_V)$.*

Лемма 8.4.27. *Множество $V$-чисел $e_2^{ik}$ порождает $D$-векторное пространство $V$.*

Доказательство. Пусть

$$(8.4.54) \qquad a = e_1^i a_i$$

произвольное $V$-число. Для любого $A$-числа $a_i$, согласно определению 5.1.1 и соглашению 5.2.1, существуют $D$-числа $a_{ik}$ такие, что

$$(8.4.55) \qquad a_i = e^k a_{ik}$$

Равенство

(8.4.56) $$a = e_1^i e^k a_{ik} = e_2^{ik} a_{ik}$$

является следствием равенств (8.4.53), (8.4.54), (8.4.55). Согласно теореме 4.1.9, лемма является следствием равенства (8.4.56). ⊙

**Лемма 8.4.28.** *Множество $V$-чисел $e_2^{ik}$ линейно независимо над кольцом $D$.*

Доказательство. Пусть

(8.4.57) $$e_2^{ik} a_{ik} = 0$$

Равенство

(8.4.58) $$e_1^i e^k a_{ik} = 0$$

является следствием равенств (8.4.53), (8.4.57). Согласно теореме 8.4.25, равенство

(8.4.59) $$e^k a_{ik} = 0 \quad k \in K$$

является следствием равенства (8.4.58). Согласно теореме 4.2.14, равенство

(8.4.60) $$a_{ik} = 0 \quad k \in K \quad i \in I$$

является следствием равенства (8.4.59). Следовательно, множество $V$-чисел $e_2^{ik}$ линейно независимо над кольцом $D$. ⊙

Доказательство теоремы 8.4.26. Теорема является следствием теоремы 4.2.14, и лемм 8.4.27, 8.4.28. □

# Система линейных уравнений в $A$-векторном пространстве

## 9.1. Линейная оболочка в левом $A$-векторном пространстве

### 9.1.1. Левое $A$-векторное пространство столбцов.

ОПРЕДЕЛЕНИЕ 9.1.1. *Пусть $V$ - левое $A$-векторное пространство столбцов и $\{\overline{a}_i \in V, i \in I\}$ - множество векторов.* **Линейная оболочка** *в левом $A$-векторном пространстве столбцов - это множество $span(\overline{a}_i, i \in I)$ векторов, линейно зависимых от векторов $\overline{a}_i$.* $\square$

ТЕОРЕМА 9.1.2. *Пусть $V$ - левое $A$-векторное пространство столбцов и строка*

$$\overline{a} = \begin{pmatrix} \overline{a}_1 & ... & \overline{a}_m \end{pmatrix} = (\overline{a}_i, i \in I)$$

*задаёт множество векторов. Вектор $\overline{b}$ линейно зависит от векторов $\overline{a}_i$*

$$\overline{b} \in span(\overline{a}_i, i \in I)$$

*если существует столбец $A$-чисел*

$$x = \begin{pmatrix} x^1 \\ ... \\ x^m \end{pmatrix}$$

*такой, что верно равенство*

$$(9.1.1) \qquad \overline{b} = x^* {}_* \overline{a} = x^i \overline{a}_i$$

ДОКАЗАТЕЛЬСТВО. Равенство (9.1.1) является следствием определения 9.1.1 и теоремы 8.4.2. $\square$

ТЕОРЕМА 9.1.3. *$span(\overline{a}_i, i \in I)$ - подпространство левого $A$-векторного пространства $V$ столбцов.*

ДОКАЗАТЕЛЬСТВО. Предположим, что

$$\overline{b} \in span(\overline{a}_i, i \in I)$$
$$\overline{c} \in span(\overline{a}_i, i \in I)$$

Согласно теореме 9.1.2

$$(9.1.2) \qquad \overline{b} = b^* {}_* \overline{a} = b^i \overline{a}_i$$

$$(9.1.3) \qquad \overline{c} = c^* {}_* \overline{a} = c^i \overline{a}_i$$

Тогда

$$\overline{b} + \overline{c} = a^* {}_* b + a^* {}_* c = a^* {}_* (b + c) \in span(\overline{a}_i, i \in I)$$
$$\overline{b}k = (a^* {}_* b)k = a^* {}_* (bk) \in span(\overline{a}_i, i \in I)$$

Это доказывает утверждение.                                                      □

Теорема 9.1.4.  *Пусть* $\overline{\overline{e}} = (e_j, j \in J)jJ$ *- базис. Пусть векторы* $\overline{b}, \overline{a}_i$
*имеют разложение*

$$(9.1.4) \qquad \overline{b} = b^* {}_* e = b^j e_j \qquad b = \begin{pmatrix} b^1 \\ ... \\ b^n \end{pmatrix}$$

$$(9.1.5) \qquad \overline{a_i} = a_i {}^* {}_* e = a_i^j e_j \qquad a_i = \begin{pmatrix} a_i^1 \\ ... \\ a_i^n \end{pmatrix}$$

*Вектор* $\overline{b}$ *линейно зависит от векторов* $\overline{a}_i$

$$\overline{b} \in span(\overline{a}_i, i \in I)$$

*если существует столбец $A$-чисел*

$$(9.1.6) \qquad x = \begin{pmatrix} x^1 \\ ... \\ x^m \end{pmatrix}$$

*которые удовлетворяют* **системе линейных уравнений**

$$x^1 a_1^1 + ... + x^m a_m^1 = b^1$$
$$(9.1.7) \qquad\qquad ......$$
$$x^1 a_1^n + ... + x^m a_m^n = b^n$$

Доказательство. Равенство

$$(9.1.8) \qquad b^j e_j = x^i(a_i^j e_j) = (x^i a_i^j)e_j$$

является следствием равенств (9.1.1), (9.1.4), (9.1.5). Равенство

$$(9.1.9) \qquad b^j = x^i a_i^j$$

является следствием равенства (9.1.8) и теоремы 8.4.4. Равенство (9.1.7) явля-
ется следствием равенства (9.1.9).                                               □

Теорема 9.1.5.  *Рассмотрим матрицу*

$$(9.1.10) \qquad a = \begin{pmatrix} \begin{pmatrix} a_1^1 \\ ... \\ a_1^n \end{pmatrix} & ... & \begin{pmatrix} a_m^1 \\ ... \\ a_m^n \end{pmatrix} \end{pmatrix} = \begin{pmatrix} a_1^1 & ... & a_m^1 \\ ... & ... & ... \\ a_1^n & ... & a_m^n \end{pmatrix}$$

*Тогда система линейных уравнений* (9.1.7) *может быть записана в следую-*

*щей форме*

$$(9.1.11) \qquad \begin{pmatrix} x^1 \\ ... \\ x^m \end{pmatrix} {}^*_* \begin{pmatrix} a^1_1 & ... & a^1_m \\ ... & ... & ... \\ a^n_1 & ... & a^n_m \end{pmatrix} = \begin{pmatrix} b^1 \\ ... \\ b^n \end{pmatrix}$$

$$(9.1.12) \qquad x^* {}_* a = b$$

ДОКАЗАТЕЛЬСТВО. Равенство (9.1.11) является следствием равенств (9.1.4), (9.1.6), (9.1.7) (9.1.10) и определения 8.1.3. □

ОПРЕДЕЛЕНИЕ 9.1.6. *Предположим, что $a$ - $^*_*$-невырожденная матрица. Мы будем называть соответствующую систему линейных уравнений*

$$\boxed{(9.1.12) \quad x^* {}_* a = b}$$

**невырожденной системой линейных уравнений**. □

ТЕОРЕМА 9.1.7. *Решение невырожденной системы линейных уравнений*

$$\boxed{(9.1.12) \quad x^* {}_* a = b}$$

*определено однозначно и может быть записано в любой из следующих форм*[9.1]

$$(9.1.13) \qquad x = b^* {}_* a^{-1^*_*}$$

$$(9.1.14) \qquad x = b^* {}_* \mathcal{H} \det({}^*_*) a$$

ДОКАЗАТЕЛЬСТВО. Умножая обе части равенства (9.1.12) слева на $a^{-1^*_*}$, мы получим (9.1.13). Пользуясь определением

$$\boxed{(8.3.11) \quad a^{-1^*_*} = \mathcal{H} \det({}^*_*) a}$$

мы получим (9.1.14). Решение системы единственно в силу теоремы 8.1.31. □

## 9.1.2. Левое $A$-векторное пространство строк.

ОПРЕДЕЛЕНИЕ 9.1.8. *Пусть $V$ - левое $A$-векторное пространство строк и $\{\overline{a}^i \in V, i \in I\}$ - множество векторов.* **Линейная оболочка** *в левом $A$-векторном пространстве строк - это множество* $span(\overline{a}^i, i \in I)$ *векторов, линейно зависимых от векторов $\overline{a}^i$.* □

ТЕОРЕМА 9.1.9. *Пусть $V$ - левое $A$-векторное пространство строк и столбец*

$$\overline{a} = \begin{pmatrix} \overline{a}^1 \\ ... \\ \overline{a}^m \end{pmatrix} = (\overline{a}^i, i \in I)$$

---

[9.1] Мы можем найти решение системы (9.1.12) в теореме [5]-1.6.1. Я повторяю это утверждение, так как я слегка изменил обозначения.

*задаёт множество векторов. Вектор $\overline{b}$ линейно зависит от векторов $\overline{a}^i$*

$$\overline{b} \in span(\overline{a}^i, i \in I)$$

*если существует строка $A$-чисел*

$$x = \begin{pmatrix} x_1 & ... & x_m \end{pmatrix}$$

*такая, что верно равенство*

(9.1.15) $$\overline{b} = x_* {}^* \overline{a} = x_i \overline{a}^i$$

Доказательство. Равенство (9.1.15) является следствием определения 9.1.8 и теоремы 8.4.9. $\qquad \square$

Теорема 9.1.10. $span(\overline{a}^i, i \in I)$ - подпространство левого $A$-векторного пространства $V$ строк.

Доказательство. Предположим, что

$$\overline{b} \in \mathrm{span}(\overline{a}^i, i \in I)$$
$$\overline{c} \in \mathrm{span}(\overline{a}^i, i \in I)$$

Согласно теореме 9.1.9

(9.1.16) $$\overline{b} = b_* {}^* \overline{a} = b_i \overline{a}^i$$

(9.1.17) $$\overline{c} = c_* {}^* \overline{a} = c_i \overline{a}^i$$

Тогда

$$\overline{b} + \overline{c} = a_* {}^* b + a_* {}^* c = a_* {}^* (b + c) \in \mathrm{span}(\overline{a}^i, i \in I)$$
$$\overline{b} k = (a_* {}^* b)k = a_* {}^* (bk) \in \mathrm{span}(\overline{a}^i, i \in I)$$

Это доказывает утверждение. $\qquad \square$

Теорема 9.1.11. *Пусть $\overline{\overline{e}} = (e^j, j \in J)jJ$ - базис. Пусть векторы $\overline{b}$, $\overline{a}^i$ имеют разложение*

(9.1.18) $$\overline{b} = b_* {}^* e = b_j e^j \qquad b = \begin{pmatrix} b_1 & ... & b_n \end{pmatrix}$$

(9.1.19) $$\overline{a^i} = a^i {}_* {}^* e = a^i_j e^j \qquad a^i = \begin{pmatrix} a^i_1 & ... & a^i_n \end{pmatrix}$$

*Вектор $\overline{b}$ линейно зависит от векторов $\overline{a}^i$*

$$\overline{b} \in span(\overline{a}^i, i \in I)$$

*если существует строка $A$-чисел*

(9.1.20) $$x = \begin{pmatrix} x_1 & ... & x_m \end{pmatrix}$$

*которые удовлетворяют* **системе линейных уравнений**

$$x_1 a^1_1 + ... + x_m a^m_1 = b_1$$

(9.1.21) $$......$$

$$x_1 a^1_n + ... + x_m a^m_n = b_n$$

Доказательство. Равенство

$$(9.1.22) \qquad b_j e^j = x_i(a_j^i e^j) = (x_i a_j^i)e^j$$

является следствием равенств (9.1.15), (9.1.18), (9.1.19). Равенство

$$(9.1.23) \qquad b_j = x_i a_j^i$$

является следствием равенства (9.1.22) и теоремы 8.4.11. Равенство (9.1.21) является следствием равенства (9.1.23). □

ТЕОРЕМА 9.1.12. *Рассмотрим матрицу*

$$(9.1.24) \qquad a = \begin{pmatrix} \begin{pmatrix} a_1^1 & ... & a_n^1 \end{pmatrix} \\ ... \\ \begin{pmatrix} a_1^m & ... & a_n^m \end{pmatrix} \end{pmatrix} = \begin{pmatrix} a_1^1 & ... & a_n^1 \\ ... & ... & ... \\ a_1^m & ... & a_n^m \end{pmatrix}$$

*Тогда система линейных уравнений (9.1.21) может быть записана в следующей форме*

$$(9.1.25) \qquad \begin{pmatrix} x_1 & ... & x_m \end{pmatrix} {}_*^* \begin{pmatrix} a_1^1 & ... & a_n^1 \\ ... & ... & ... \\ a_1^m & ... & a_n^m \end{pmatrix} = \begin{pmatrix} b^1 \\ ... \\ b^n \end{pmatrix}$$

$$(9.1.26) \qquad x_*{}^* a = b$$

Доказательство. Равенство (9.1.25) является следствием равенств (9.1.18), (9.1.20), (9.1.21) (9.1.24) и определения 8.1.4. □

ОПРЕДЕЛЕНИЕ 9.1.13. *Предположим, что $a$ - ${}_*^*$-невырожденная матрица. Мы будем называть соответствующую систему линейных уравнений*

$$\boxed{(9.1.26) \quad x_*{}^* a = b}$$

**невырожденной системой линейных уравнений**. □

ТЕОРЕМА 9.1.14. *Решение невырожденной системы линейных уравнений*

$$\boxed{(9.1.26) \quad x_*{}^* a = b}$$

*определено однозначно и может быть записано в любой из следующих форм* [9.2]

$$(9.1.27) \qquad x = b_*{}^* a^{-1}{}_*^*$$

$$(9.1.28) \qquad x = b_*{}^* \mathcal{H}\det({}_*^*)a$$

Доказательство. Умножая обе части равенства (9.1.26) слева на $a^{-1}{}_*^*$, мы получим (9.1.27). Пользуясь определением

$$\boxed{(8.2.11) \quad a^{-1}{}_*^* = \mathcal{H}\det({}_*^*)a}$$

мы получим (9.1.28). Решение системы единственно в силу теоремы 8.1.31. □

---

[9.2] Мы можем найти решение системы (9.1.26) в теореме [5]-1.6.1. Я повторяю это утверждение, так как я слегка изменил обозначения.

## 9.2. Линейная оболочка в правом $A$-векторном пространстве

### 9.2.1. Правое $A$-векторное пространство столбцов.

Определение 9.2.1. *Пусть $V$ - правое $A$-векторное пространство столбцов и $\{\overline{a}_i \in V, i \in I\}$ - множество векторов.* **Линейная оболочка** *в правом $A$-векторном пространстве столбцов - это множество $span(\overline{a}_i, i \in I)$ векторов, линейно зависимых от векторов $\overline{a}_i$.* □

Теорема 9.2.2. *Пусть $V$ - правое $A$-векторное пространство столбцов и строка*

$$\overline{a} = \begin{pmatrix} \overline{a}_1 & ... & \overline{a}_m \end{pmatrix} = (\overline{a}_i, i \in I)$$

*задаёт множество векторов. Вектор $\overline{b}$ линейно зависит от векторов $\overline{a}_i$*

$$\overline{b} \in span(\overline{a}_i, i \in I)$$

*если существует столбец $A$-чисел*

$$x = \begin{pmatrix} x^1 \\ ... \\ x^m \end{pmatrix}$$

*такой, что верно равенство*

$$(9.2.1) \qquad \overline{b} = \overline{a}_* {}^* x = \overline{a}_i x^i$$

Доказательство. Равенство (9.2.1) является следствием определения 9.2.1 и теоремы 8.4.16. □

Теорема 9.2.3. *$span(\overline{a}_i, i \in I)$ - подпространство правого $A$-векторного пространства $V$ столбцов.*

Доказательство. Предположим, что

$$\overline{b} \in span(\overline{a}_i, i \in I)$$
$$\overline{c} \in span(\overline{a}_i, i \in I)$$

Согласно теореме 9.2.2

$$(9.2.2) \qquad \overline{b} = \overline{a}_* {}^* b = \overline{a}_i b^i$$

$$(9.2.3) \qquad \overline{c} = \overline{a}_* {}^* c = \overline{a}_i c^i$$

Тогда

$$\overline{b} + \overline{c} = \overline{a}_* {}^* b + \overline{a}_* {}^* c = \overline{a}_* {}^* (b + c) \in span(\overline{a}_i, i \in I)$$

$$\overline{b} k = (\overline{a}_* {}^* b) k = \overline{a}_* {}^* (bk) \in span(\overline{a}_i, i \in I)$$

Это доказывает утверждение. □

Теорема 9.2.4. *Пусть* $\overline{\overline{e}} = (e_j, j \in J) jJ$ *- базис. Пусть векторы* $\overline{b}$, $\overline{a}_i$ *имеют разложение*

$$(9.2.4) \qquad \overline{b} = e_* {}^* b = e_j b^j \qquad b = \begin{pmatrix} b^1 \\ ... \\ b^n \end{pmatrix}$$

$$(9.2.5) \qquad \overline{a_i} = e_* {}^* a_i = e_j a_i^j \qquad a_i = \begin{pmatrix} a_i^1 \\ ... \\ a_i^n \end{pmatrix}$$

*Вектор* $\overline{b}$ *линейно зависит от векторов* $\overline{a}_i$

$$\overline{b} \in span(\overline{a}_i, i \in I)$$

*если существует столбец $A$-чисел*

$$(9.2.6) \qquad x = \begin{pmatrix} x^1 \\ ... \\ x^m \end{pmatrix}$$

*которые удовлетворяют* **системе линейных уравнений**

$$(9.2.7) \qquad \begin{array}{c} a_1^1 x^1 + ... + a_m^1 x^m = b^1 \\ ...... \\ a_1^n x^1 + ... + a_m^n x^m = b^n \end{array}$$

Доказательство. Равенство

$$(9.2.8) \qquad e_j b^j = (e_j a_i^j) x^i = e_j (a_i^j x^i)$$

является следствием равенств (9.2.1), (9.2.4), (9.2.5). Равенство

$$(9.2.9) \qquad b^j = a_i^j x^i$$

является следствием равенства (9.2.8) и теоремы 8.4.18. Равенство (9.2.7) является следствием равенства (9.2.9). $\square$

Теорема 9.2.5. *Рассмотрим матрицу*

$$(9.2.10) \qquad a = \begin{pmatrix} \begin{pmatrix} a_1^1 \\ ... \\ a_1^n \end{pmatrix} & ... & \begin{pmatrix} a_m^1 \\ ... \\ a_m^n \end{pmatrix} \end{pmatrix} = \begin{pmatrix} a_1^1 & ... & a_m^1 \\ ... & ... & ... \\ a_1^n & ... & a_m^n \end{pmatrix}$$

*Тогда система линейных уравнений* (9.2.7) *может быть записана в следующей форме*

$$(9.2.11) \qquad \begin{pmatrix} a_1^1 & ... & a_m^1 \\ ... & ... & ... \\ a_1^n & ... & a_m^n \end{pmatrix} * \begin{pmatrix} x^1 \\ ... \\ x^m \end{pmatrix} = \begin{pmatrix} b^1 \\ ... \\ b^n \end{pmatrix}$$

(9.2.12)                                $a_{*}{}^{*}x = b$

Доказательство. Равенство (9.2.11) является следствием равенств (9.2.4), (9.2.6), (9.2.7) (9.2.10) и определения 8.1.3.                                $\square$

Определение 9.2.6. *Предположим, что $a - {}_{*}{}^{*}$-невырожденная матрица. Мы будем называть соответствующую систему линейных уравнений*

$$\boxed{(9.2.12) \quad a_{*}{}^{*}x = b}$$

**невырожденной системой линейных уравнений**.                                $\square$

Теорема 9.2.7. *Решение невырожденной системы линейных уравнений*

$$\boxed{(9.2.12) \quad a_{*}{}^{*}x = b}$$

*определено однозначно и может быть записано в любой из следующих форм*[9.3]

(9.2.13)                                $x = a^{-1}{}_{*}{}^{*}{}_{*}{}^{*}b$

(9.2.14)                                $x = \mathcal{H}\det({}_{*}{}^{*})a_{*}{}^{*}b$

Доказательство. Умножая обе части равенства (9.2.12) слева на $a^{-1}{}_{*}{}^{*}$, мы получим (9.2.13). Пользуясь определением

$$\boxed{(8.2.11) \quad a^{-1}{}_{*}{}^{*} = \mathcal{H}\det({}_{*}{}^{*})a}$$

мы получим (9.2.14). Решение системы единственно в силу теоремы 8.1.31.
                                                                           $\square$

### 9.2.2. Правое $A$-векторное пространство строк.

Определение 9.2.8. *Пусть $V$ - правое $A$-векторное пространство строк и $\{\overline{a}^i \in V, i \in I\}$ - множество векторов.* **Линейная оболочка** *в правом $A$-векторном пространстве строк - это множество* $span(\overline{a}^i, i \in I)$ *векторов, линейно зависимых от векторов* $\overline{a}^i$.                                $\square$

Теорема 9.2.9. *Пусть $V$ - правое $A$-векторное пространство строк и столбец*

$$\overline{a} = \begin{pmatrix} \overline{a}^1 \\ ... \\ \overline{a}^m \end{pmatrix} = (\overline{a}^i, i \in I)$$

*задаёт множество векторов. Вектор $\overline{b}$ линейно зависит от векторов $\overline{a}^i$*

$$\overline{b} \in span(\overline{a}^i, i \in I)$$

*если существует строка $A$-чисел*

$$x = \begin{pmatrix} x_1 & ... & x_m \end{pmatrix}$$

---

[9.3]   Мы можем найти решение системы (9.2.12) в теореме [5]-1.6.1. Я повторяю это утверждение, так как я слегка изменил обозначения.

*такая, что верно равенство*

$$(9.2.15) \qquad \overline{b} = \overline{a}^*{}_*x = \overline{a}^i x_i$$

Доказательство. Равенство (9.2.15) является следствием определения 9.2.8 и теоремы 8.4.23. $\qquad \square$

Теорема 9.2.10. *$\mathrm{span}(\overline{a}^i, i \in I)$ - подпространство правого $A$-векторного пространства $V$ строк.*

Доказательство. Предположим, что

$$\overline{b} \in \mathrm{span}(\overline{a}^i, i \in I)$$
$$\overline{c} \in \mathrm{span}(\overline{a}^i, i \in I)$$

Согласно теореме 9.2.9

$$(9.2.16) \qquad \overline{b} = \overline{a}^*{}_*b = \overline{a}^i b_i$$

$$(9.2.17) \qquad \overline{c} = \overline{a}^*{}_*c = \overline{a}^i c_i$$

Тогда

$$\overline{b} + \overline{c} = \overline{a}^*{}_*b + \overline{a}^*{}_*c = \overline{a}^*{}_*(b + c) \in \mathrm{span}(\overline{a}^i, i \in I)$$
$$\overline{b}k = (\overline{a}^*{}_*b)k = \overline{a}^*{}_*(bk) \in \mathrm{span}(\overline{a}^i, i \in I)$$

Это доказывает утверждение. $\qquad \square$

Теорема 9.2.11. *Пусть $\overline{\overline{e}} = (e^j, j \in J)jJ$ - базис. Пусть векторы $\overline{b}$, $\overline{a}^i$ имеют разложение*

$$(9.2.18) \qquad \overline{b} = e^*{}_*b = e^j b_j \qquad b = \begin{pmatrix} b_1 & ... & b_n \end{pmatrix}$$

$$(9.2.19) \qquad \overline{a^i} = e^*{}_*a^i = e^j a^i_j \qquad a^i = \begin{pmatrix} a^i_1 & ... & a^i_n \end{pmatrix}$$

*Вектор $\overline{b}$ линейно зависит от векторов $\overline{a}^i$*

$$\overline{b} \in span(\overline{a}^i, i \in I)$$

*если существует строка $A$-чисел*

$$(9.2.20) \qquad x = \begin{pmatrix} x_1 & ... & x_m \end{pmatrix}$$

*которые удовлетворяют* **системе линейных уравнений**

$$(9.2.21) \qquad \begin{aligned} a^1_1 x_1 + ... + a^m_1 x_m &= b_1 \\ &...... \\ a^1_n x_1 + ... + a^m_n x_m &= b_n \end{aligned}$$

Доказательство. Равенство

$$(9.2.22) \qquad e^j b_j = (e^j a^i_j)x_i = e^j(a^i_j x_i)$$

является следствием равенств (9.2.15), (9.2.18), (9.2.19). Равенство

$$(9.2.23) \qquad b_j = a^i_j x_i$$

является следствием равенства (9.2.22) и теоремы 8.4.25. Равенство (9.2.21) является следствием равенства (9.2.23).                                         □

ТЕОРЕМА 9.2.12. *Рассмотрим матрицу*

$$(9.2.24) \qquad a = \begin{pmatrix} \begin{pmatrix} a_1^1 & ... & a_n^1 \end{pmatrix} \\ ... \\ \begin{pmatrix} a_1^m & ... & a_n^m \end{pmatrix} \end{pmatrix} = \begin{pmatrix} a_1^1 & ... & a_n^1 \\ ... & ... & ... \\ a_1^m & ... & a_n^m \end{pmatrix}$$

*Тогда система линейных уравнений* (9.2.21) *может быть записана в следующей форме*

$$(9.2.25) \qquad \begin{pmatrix} a_1^1 & ... & a_n^1 \\ ... & ... & ... \\ a_1^m & ... & a_n^m \end{pmatrix} {}^*{}_* \begin{pmatrix} x_1 & ... & x_m \end{pmatrix} = \begin{pmatrix} b_1 & ... & b_n \end{pmatrix}$$

$$(9.2.26) \qquad\qquad a^*{}_* x = b$$

ДОКАЗАТЕЛЬСТВО. Равенство (9.2.25) является следствием равенств (9.2.18), (9.2.20), (9.2.21) (9.2.24) и определения 8.1.4.                                  □

ОПРЕДЕЛЕНИЕ 9.2.13. *Предположим, что $a$ - ${}^*{}_*$-невырожденная матрица. Мы будем называть соответствующую систему линейных уравнений*

$$\boxed{(9.2.26) \quad a^*{}_* x = b}$$

**невырожденной системой линейных уравнений**.                              □

ТЕОРЕМА 9.2.14. *Решение невырожденной системы линейных уравнений*

$$\boxed{(9.2.26) \quad a^*{}_* x = b}$$

*определено однозначно и может быть записано в любой из следующих форм*[9.4]

$$(9.2.27) \qquad\qquad x = a^{-1\,*}{}^*{}_* b$$

$$(9.2.28) \qquad\qquad x = \mathcal{H} \det({}^*{}_*) a^*{}_* b$$

ДОКАЗАТЕЛЬСТВО. Умножая обе части равенства (9.2.26) слева на $a^{-1\,*}{}^*$, мы получим (9.2.27). Пользуясь определением

$$\boxed{(8.3.11) \quad a^{-1\,*}{}^* = \mathcal{H} \det({}^*{}_*) a}$$

мы получим (9.2.28). Решение системы единственно в силу теоремы 8.1.31.
□

## 9.3. Ранг матрицы

Мы делаем следующие предположения в этом разделе

- $i \in M$, $|M| = m$, $j \in N$, $|N| = n$.
- $a = (a_j^i)$ - произвольная матрица.

---

Мы можем найти решение системы (9.2.26) в теореме [5]-1.6.1. Я повторяю это утверждение, так как я слегка изменил обозначения.

- $k$, $s \in S \supseteq M$, $l$, $t \in T \supseteq N$, $k = |S| = |T|$.
- $p \in M \setminus S$, $r \in N \setminus T$.

Матрица $a_T^S$ является $k \times k$ матрицей.

ОПРЕДЕЛЕНИЕ 9.3.1. *Если под-матрица $a_T^S$ - $_*{}^*$-невырожденная матрица, то мы будем говорить, что $_*{}^*$-ранг матрицы а не мень-ше, чем $k$. $_*{}^*$-**ранг матрицы** $a$, rank$(_*{}^*)a$, - это максимальное значение $k$. Мы будем называть соответствующую подматрицу $_*{}^*$-**главной подматрицей**.* □

ОПРЕДЕЛЕНИЕ 9.3.2. *Если под-матрица $a_T^S$ - $^*{}_*$-невырожденная матрица, то мы будем говорить, что $^*{}_*$-ранг матрицы а не мень-ше, чем $k$. $^*{}_*$-**ранг матрицы** $a$, rank$(^*{}_*)a$, - это максимальное значение $k$. Мы будем называть соответствующую подматрицу $^*{}_*$-**главной подматрицей**.* □

ТЕОРЕМА 9.3.3. *Пусть матрица а - $_*{}^*$-вырожденная матрица и под-матрица $a_T^S$ - $_*{}^*$-главная подматри-ца, тогда*

$$\det(_*{}^*)_r^p \, a_{T \cup \{r\}}^{S \cup \{p\}} = 0$$

ТЕОРЕМА 9.3.4. *Пусть матрица а - $^*{}_*$-вырожденная матрица и под-матрица $a_T^S$ - $^*{}_*$-главная подматри-ца, тогда*

$$\det(^*{}_*)_r^p \, a_{T \cup \{r\}}^{S \cup \{p\}} = 0$$

ДОКАЗАТЕЛЬСТВО ТЕОРЕМЫ 9.3.3. Чтобы понять, почему подматрица

$$(9.3.1) \qquad b = a_{T \cup \{r\}}^{S \cup \{p\}}$$

не имеет $_*{}^*$-обратной матрицы, [9.5] мы предположим, что существует $_*{}^*$-обрат-ная матрица $b^{-1}{}_*{}^*$. Запишем систему линейных уравнений

$$(8.2.3) \quad a_I^{[J]}{}_*{}^*(a^{-1}{}_*{}^*)_J^I + a_{[I]}^{[J]}{}_*{}^*(a^{-1}{}_*{}^*)_J^{[I]} = 0$$

$$(8.2.4) \quad a_I^J{}_*{}^*(a^{-1}{}_*{}^*)_J^I + a_{[I]}^J{}_*{}^*(a^{-1}{}_*{}^*)_J^{[I]} = E_k$$

полагая $I = \{r\}$, $J = \{p\}$ (в этом случае $[I] = T$, $[J] = S$ )

$$(9.3.2) \qquad b_r^S \, b^{-1}{}_*{}^{*\,r}{}_p + b_{T*}^S{}^*b^{-1}{}_*{}^{*\,T}{}_p = 0$$

$$(9.3.3) \qquad b_r^p \, b^{-1}{}_*{}^{*\,r}{}_p + b_{T*}^p{}^*b^{-1}{}_*{}^{*\,T}{}_p = 1$$

и попробуем решить эту систему. Мы умножим (9.3.2) на $(b_T^S)^{-1}{}_*{}^*$

$$(9.3.4) \qquad (b_T^S)^{-1}{}_*{}^*{}_*b_r^S \, b^{-1}{}_*{}^{*\,r}{}_p + b^{-1}{}_*{}^{*\,T}{}_p = 0$$

Теперь мы можем подставить (9.3.4) в (9.3.3)

$$(9.3.5) \qquad b_r^p \, b^{-1}{}_*{}^{*\,r}{}_p - b_{T*}^p{}^*(b_T^S)^{-1}{}_*{}^*{}_*b_r^S \, b^{-1}{}_*{}^{*\,r}{}_p = 1$$

---

[9.5]Естественно ожидать связь между $_*{}^*$-вырожденностью матрицы и её $_*{}^*$-квазидетерми-нантом, подобную связи, известной в коммутативном случае. Однако $_*{}^*$-квазидетерминант определён не всегда. Например, если $_*{}^*$-обратная матрица имеет слишком много элементов, равных 0. Как следует из этой теоремы, не определён $_*{}^*$-квазидетерминант также в случае, когда $_*{}^*$-ранг матрицы меньше $n - 1$.

Из (9.3.5) следует

$$(9.3.6) \qquad (b_r^p - b_{T*}^p {}^*(b_T^S)^{-1*}{}^*{}^*b_r^S) \; b^{-1*}{}_p^{*\,r} = 1$$

Равенство

$$(9.3.7) \qquad (\det({}_*{}^*)_r^p b) b^{-1*}{}_p^{*\,r} = 1$$

является следствием равенства

$$(8.2.9) \qquad \boxed{\det({}_*{}^*)_i^j \, a = a_i^j - a_{[i]}^j {}^* \left( a_{[i]}^{[j]} \right)^{-1*}{}^*{}^*a_i^{[j]} = ((a^{-1*})_j^i)^{-1}}$$

и равенства (9.3.6). Тем самым мы доказали, что квазидетерминант $\det({}_*{}^*)_r^p \, b$ определён и условие его обращения в 0 необходимое и достаточное условие вырожденности матрицы $b$.                     □

ДОКАЗАТЕЛЬСТВО ТЕОРЕМЫ 9.3.4.     Чтобы понять, почему подматрица

$$(9.3.8) \qquad b = a_{T \cup \{r\}}^{S \cup \{p\}}$$

не имеет ${}_*$-обратной матрицы, [9.6] мы предположим, что существует ${}_*$-обратная матрица $b^{-1*}{}^*$. Запишем систему линейных уравнений

$$(8.3.3) \qquad \boxed{a_{[J]}^I {}^*{}^*(a^{-1*}{}^*)_I^J + a_{[J]}^{[I]} {}^*{}^*(a^{-1*}{}^*)_{[I]}^J = 0}$$

$$(8.3.4) \qquad \boxed{a_J^I {}^*{}^*(a^{-1*}{}^*)_I^J + a_J^{[I]} {}^*{}^*(a^{-1*}{}^*)_{[I]}^J = E_k}$$

полагая $I = \{p\}$, $J = \{r\}$  (в этом случае $[I] = S$, $[J] = T$ )

$$(9.3.9) \qquad b_T^p \; b^{-1*}{}_p^{*\,r} + b_T^{S*}{}_*b^{-1*}{}_S^{*\,r} = 0$$

$$(9.3.10) \qquad b_r^p \; b^{-1*}{}_p^{*\,r} + b_r^{S*}{}_*b^{-1*}{}_S^{*\,r} = 1$$

и попробуем решить эту систему. Мы умножим (9.3.9) на $(b_T^S)^{-1*}{}^*$

$$(9.3.11) \qquad (b_T^S)^{-1*}{}^*{}^*b_T^p \; b^{-1*}{}_p^{*\,r} + b^{-1*}{}_S^{*\,r} = 0$$

Теперь мы можем подставить (9.3.11) в (9.3.10)

$$(9.3.12) \qquad b_r^p \; b^{-1*}{}_p^{*\,r} - b_r^{S*}{}_*(b_T^S)^{-1*}{}^*{}^*b_T^p \; b^{-1*}{}_p^{*\,r} = 1$$

Из (9.3.12) следует

$$(9.3.13) \qquad (b_r^p - b_r^{S*}{}_*(b_T^S)^{-1*}{}^*{}^*b_T^p) \; b^{-1*}{}_p^{*\,r} = 1$$

Равенство

$$(9.3.14) \qquad (\det({}^*{}_*)_r^p b) b^{-1*}{}_p^{*\,r} = 1$$

является следствием равенства

$$(8.3.9) \qquad \boxed{\det({}^*{}_*)_j^i \, a = a_j^i - a_j^{[i]} {}_* \left( a_{[j]}^{[i]} \right)^{-1*}{}^*{}_*a_{[j]}^i = ((a^{-1*})_j^i)^{-1}}$$

и равенства (9.3.13). Тем самым мы доказали, что квазидетерминант $\det({}^*{}_*)_r^p \, b$ определён и условие его обращения в 0 необходимое и достаточное условие вырожденности матрицы $b$.                     □

---

[9.6]Естественно ожидать связь между ${}_*$-вырожденностью матрицы и её ${}_*$-квазидетерминантом, подобную связи, известной в коммутативном случае. Однако ${}_*$-квазидетерминант определён не всегда. Например, если ${}_*$-обратная матрица имеет слишком много элементов, равных 0. Как следует из этой теоремы, не определён ${}_*$-квазидетерминант также в случае, когда ${}_*$-ранг матрицы меньше $n-1$.

Теорема 9.3.5. *Предположим, что a - матрица,*

$$\mathrm{rank}(^*_*)a = k < n$$

*и $a_T^S$ - $^*_*$-главная подматрица. Тогда столбец $a_r$ является левой линейной комбинацией столбцов $a_T$*

$$(9.3.15) \qquad a_{N\setminus T} = R^*{}_*a_T$$

$$(9.3.16) \qquad a_r = R_r{}^*_*a_T$$

$$(9.3.17) \qquad a_r^t = R_r^t a_t^t$$

Доказательство. Если матрица $a$ имеет $k$ строк, то, полагая, что столбец $a_r$ - левая линейная комбинация (9.3.16) столбцов $a_T$ с коэффициентами $R_r^t$, мы получим систему левых линейных уравнений (9.3.17). Мы положим, что $x^t = R_r^t$ - это неизвестные переменные в системе левых линейных уравнений (9.3.17). Согласно теореме 9.1.7, система левых линейных уравнений (9.3.17) имеет единственное решение и это решение нетривиально потому, что все $^*_*$-квазидетерминанты отличны от 0.

Остаётся доказать утверждение в случае, когда число строк матрицы $a$ больше чем $k$. Пусть нам даны столбец $a_r$ и строка $a^p$. Согласно предположению, подматрица $a_{T\cup\{r\}}^{S\cup\{p\}}$ - $^*_*$-вырожденная матрица и её $^*_*$-квазидетерминант

$$(9.3.18) \qquad \det(^*_*)_r^p a_{T\cup\{r\}}^{S\cup\{p\}} = 0$$

Согласно равенству (8.3.9) равенство (9.3.18) имеет вид

$$a_r^p - a_r^{S*}((a_T^S)^{-1*})^*{}_*a_T^p = 0$$

Матрица

$$(9.3.19) \qquad R_r = a_r^{S*}((a_T^S)^{-1*})$$

не зависит от $p$. Следовательно, для любых $p \in M \setminus S$

$$(9.3.20) \qquad a_r^p = R_r{}^*_*a_T^p$$

Теорема 9.3.6. *Предположим, что a - матрица,*

$$\mathrm{rank}(_*^*)a = k < n$$

*и $a_T^S$ - $_*^*$-главная подматрица. Тогда столбец $a_r$ является правой линейной комбинацией столбцов $a_T$*

$$(9.3.23) \qquad a_{N\setminus T} = a_T{}_*^*R$$

$$(9.3.24) \qquad a_r = a_T{}_*^*R_r$$

$$(9.3.25) \qquad a_r^a = a_t^t\ R_r^t$$

Доказательство. Если матрица $a$ имеет $k$ строк, то, полагая, что столбец $a_r$ - правая линейная комбинация (9.3.24) столбцов $a_T$ с коэффициентами $R_r^t$, мы получим систему правых линейных уравнений (9.3.25). Мы положим, что $x^t = R_r^t$ - это неизвестные переменные в системе правых линейных уравнений (9.3.25). Согласно теореме 9.2.7, система правых линейных уравнений (9.3.25) имеет единственное решение и это решение нетривиально потому, что все $_*^*$-квазидетерминанты отличны от 0.

Остаётся доказать утверждение в случае, когда число строк матрицы $a$ больше чем $k$. Пусть нам даны столбец $a_r$ и строка $a^p$. Согласно предположению, подматрица $a_{T\cup\{r\}}^{S\cup\{p\}}$ - $_*^*$-вырожденная матрица и её $_*^*$-квазидетерминант

$$(9.3.26) \qquad \det(_*^*)_r^p a_{T\cup\{r\}}^{S\cup\{p\}} = 0$$

Согласно равенству (8.2.9) равенство (9.3.26) имеет вид

$$a_r^p - a_T^p{}_*^*((a_T^S)^{-1*})_*{}^*a_r^S = 0$$

Матрица

$$(9.3.27) \qquad R_r = ((a_T^S)^{-1*})_*{}^*a_r^S$$

не зависит от $p$. Следовательно, для любых $p \in M \setminus S$

$$(9.3.28) \qquad a_r^p = a_T^p{}_*^*R_r$$

Чтобы доказать равенство (9.3.17), необходимо доказать равенство

$$(9.3.21) \qquad a_r^k = R_r{}^*{}_* a_T^k$$

Из равенства

$$((a_T^S)^{-1^*})^*{}_* a_T^k = \delta_S^k$$

следует, что

$$(9.3.22) \qquad \begin{aligned} a_r^k &= a_r^{S*}{}_* \delta_S^k \\ &= a_r^{S*}{}_* ((a_T^S)^{-1^*})^*{}_* a_T^k \end{aligned}$$

Равенство (9.3.21) является следствием равенств (9.3.19), (9.3.22). Равенство (9.3.17) является следствием равенств (9.3.20), (9.3.21). □

ТЕОРЕМА 9.3.7. *Предположим, что $a$ - матрица,*

$$\operatorname{rank}({}^*{}_*)a = k < m$$

*и $a_T^S$ - ${}^*{}_*$-главная подматрица. Тогда строка $a^p$ является правой линейной комбинацией строк $a^S$*

$$(9.3.31) \qquad a^{M \setminus S} = a^{S*}{}_* R$$

$$(9.3.32) \qquad a^p = a^{S*}{}_* R^p$$

$$(9.3.33) \qquad a_p^b = a_s^b R_p^s$$

ДОКАЗАТЕЛЬСТВО. Если матрица $a$ имеет $k$ столбцов, то, полагая, что строка $a^p$ - правая линейная комбинация (9.3.32) строк $a^S$ с коэффициентами $R_s^p$, мы получим систему правых линейных уравнений (9.3.33). Мы положим, что $x_s = R_p^s$ - это неизвестные переменные в системе правых линейных уравнений (9.3.33). Согласно теореме 9.2.14, система правых линейных уравнений (9.3.33) имеет единственное решение и это решение нетривиально потому, что все ${}_*$-квазидетерминанты отличны от 0.

Остаётся доказать утверждение в случае, когда число столбцов матрицы $a$ больше чем $k$. Пусть нам даны строка $a^p$ и столбец $a_r$. Согласно предположению, подматрица $a_{T \cup \{r\}}^{S \cup \{p\}}$ - ${}^*{}_*$-вырожденная матрица и её ${}^*{}_*$-квази-

Чтобы доказать равенство (9.3.25), необходимо доказать равенство

$$(9.3.29) \qquad a_r^k = a_T^k{}_*{}^* R_r$$

Из равенства

$$a_T^k{}_*{}^* ((a_T^S)^{-1^*}) = \delta_S^k$$

следует, что

$$(9.3.30) \qquad \begin{aligned} a_r^k &= \delta_{S*}^k{}^* a_r^S \\ &= a_T^k{}_*{}^* ((a_T^S)^{-1^*})_*{}^* a_r^S \end{aligned}$$

Равенство (9.3.29) является следствием равенств (9.3.27), (9.3.30). Равенство (9.3.25) является следствием равенств (9.3.28), (9.3.29). □

ТЕОРЕМА 9.3.8. *Предположим, что $a$ - матрица,*

$$\operatorname{rank}({}_*{}^*)a = k < m$$

*и $a_T^S$ - ${}_*{}^*$-главная подматрица. Тогда строка $a^p$ является левой линейной комбинацией строк $a^S$*

$$(9.3.39) \qquad a^{M \setminus S} = R_*{}^* a^S$$

$$(9.3.40) \qquad a^p = R^p{}_*{}^* a^S$$

$$(9.3.41) \qquad a_b^p = R_s^p a_b^s$$

ДОКАЗАТЕЛЬСТВО. Если матрица $a$ имеет $k$ столбцов, то, полагая, что строка $a^p$ - левая линейная комбинация (9.3.40) строк $a^S$ с коэффициентами $R_s^p$, мы получим систему левых линейных уравнений (9.3.41). Мы положим, что $x_s = R_s^p$ - это неизвестные переменные в системе левых линейных уравнений (9.3.41). Согласно теореме 9.1.14, система левых линейных уравнений (9.3.41) имеет единственное решение и это решение нетривиально потому, что все ${}_*{}^*$-квазидетерминанты отличны от 0.

Остаётся доказать утверждение в случае, когда число столбцов матрицы $a$ больше чем $k$. Пусть нам даны строка $a^p$ и столбец $a_r$. Согласно предположению, подматрица $a_{T \cup \{r\}}^{S \cup \{p\}}$ - ${}_*{}^*$-вырожденная матрица и её ${}_*{}^*$-квази-

детерминант

$$(9.3.34) \qquad \det({}^{*}{}_{*})^{p}_{r}\, a^{S \cup \{p\}}_{T \cup \{r\}} = 0$$

Согласно равенству (8.3.9) равенство (9.3.34) имеет вид

$$a^{p}_{r} - a^{S*}_{r}{}_{*}((a^{S}_{T})^{-1^{**}})^{*}{}_{*}a^{p}_{T} = 0$$

Матрица

$$(9.3.35) \qquad R^{p} = ((a^{S}_{T})^{-1^{**}})^{*}{}_{*}a^{p}_{T}$$

не зависит от $r$. Следовательно, для любых $r \in N \setminus T$

$$(9.3.36) \qquad a^{p}_{r} = a^{S*}_{r}{}_{*}R^{p}$$

Чтобы доказать равенство (9.3.33), необходимо доказать равенство

$$(9.3.37) \qquad a^{p}_{l} = a^{S*}_{l}{}_{*}R^{p}$$

Из равенства

$$a^{S*}_{l}{}_{*}((a^{S}_{T})^{-1^{**}}) = \delta^{T}_{l}$$

следует, что

$$(9.3.38) \qquad a^{p}_{l} = \delta^{T}_{l}{}^{*}{}_{*}a^{p}_{T}$$
$$= a^{S*}_{l}{}_{*}((a^{S}_{T})^{-1^{**}})^{*}{}_{*}a^{p}_{T}$$

Равенство (9.3.37) является следствием равенств (9.3.35), (9.3.38). Равенство (9.3.33) является следствием равенств (9.3.36), (9.3.37).                  □

детерминант

$$(9.3.42) \qquad \det({}_{*}{}^{*})^{p}_{r}\, a^{S \cup \{p\}}_{T \cup \{r\}} = 0$$

Согласно равенству (8.2.9) равенство (9.3.42) имеет вид

$$a^{p}_{r} - a^{p}_{T}{}_{*}{}^{*}((a^{S}_{T})^{-1^{*}})_{*}{}^{*}a^{S}_{r} = 0$$

Матрица

$$(9.3.43) \qquad R^{p} = a^{p}_{T}{}_{*}{}^{*}((a^{S}_{T})^{-1^{*}})$$

не зависит от $r$. Следовательно, для любых $r \in N \setminus T$

$$(9.3.44) \qquad a^{p}_{r} = R^{p}{}_{*}{}^{*}a^{S}_{r}$$

Чтобы доказать равенство (9.3.41), необходимо доказать равенство

$$(9.3.45) \qquad a^{p}_{l} = R^{p}{}_{*}{}^{*}a^{S}_{l}$$

Из равенства

$$((a^{S}_{T})^{-1^{*}})_{*}{}^{*}a^{S}_{l} = \delta^{T}_{l}$$

следует, что

$$(9.3.46) \qquad a^{p}_{l} = a^{p}_{T}{}_{*}{}^{*}\delta^{T}_{l}$$
$$= a^{p}_{T}{}_{*}{}^{*}((a^{S}_{T})^{-1^{*}})_{*}{}^{*}a^{S}_{l}$$

Равенство (9.3.45) является следствием равенств (9.3.43), (9.3.46). Равенство (9.3.41) является следствием равенств (9.3.44), (9.3.45).                  □

---

Теорема 9.3.9. *Предположим, что матрица a имеет $n$ столбцов. Если*

$$\mathrm{rank}({}_{*}{}^{*})a = k < n$$

*то столбцы матрицы линейно зависимы справа*

$$a_{*}{}^{*}\lambda = 0$$

Теорема 9.3.10. *Предположим, что матрица a имеет $n$ столбцов. Если*

$$\mathrm{rank}({}^{*}{}_{*})a = k < n$$

*то столбцы матрицы линейно зависимы слева*

$$\lambda^{*}{}_{*}a = 0$$

---

Доказательство. Предположим, что столбец $a_r$ является правой линейной комбинацией столбцов (9.3.24). Мы положим $\lambda^{r} = -1$, $\lambda^{t} = R^{t}_{r}$ и остальные $\lambda^{c} = 0$.                  □

Доказательство. Предположим, что столбец $a_r$ является левой линейной комбинацией столбцов (9.3.16). Мы положим $\lambda^{r} = -1$, $\lambda^{t} = R^{t}_{r}$ и остальные $\lambda^{c} = 0$.                  □

---

Теорема 9.3.11. *Предположим, что матрица a имеет $m$ строк. Если*

$$\mathrm{rank}({}_{*}{}^{*})a = k < m$$

*то строки матрицы линейно зависи-*

Теорема 9.3.12. *Предположим, что матрица a имеет $m$ строк. Если*

$$\mathrm{rank}({}^{*}{}_{*})a = k < m$$

*то строки матрицы линейно зависи-*

| *мы слева* | *мы справа* |
|---|---|
| $\lambda_* {}^* a = 0$ | $a^* {}_* \lambda = 0$ |

**Доказательство.** Предположим, что строка $a^p$ является левой линейной комбинацией строк (9.3.40). Мы положим $\lambda_p = -1$, $\lambda_s = R^p_s$ и остальные $\lambda_c = 0$. ☐

**Доказательство.** Предположим, что строка $a^p$ является правой линейной комбинацией строк (9.3.32). Мы положим $\lambda_p = -1$, $\lambda_s = R^p_s$ и остальные $\lambda_c = 0$. ☐

**Теорема 9.3.13.** *Пусть $V$ - правое $A$-векторное пространство столбцов. Пусть $(\overline{a}^i, i \in M, |M| = m)$ семейство векторов, линейно независимых слева. Тогда $_*$-ранг их координатной матрицы равен* $m$.

**Теорема 9.3.14.** *Пусть $V$ - левое $A$-векторное пространство столбцов. Пусть $(\overline{a}^i, i \in M, |M| = m)$ семейство векторов, линейно независимых слева. Тогда $^*_*$-ранг их координатной матрицы равен* $m$.

**Доказательство.** Пусть $\overline{\overline{e}}$ - базис правого $A$-векторного пространства $V$. Согласно теореме 8.4.17, координатная матрица семейства векторов $(\overline{a}^i)$ относительно базиса $\overline{\overline{e}}$ состоит из столбцов, являющихся координатными матрицами векторов $\overline{a}^i$ относительно базиса $\overline{\overline{e}}$. Поэтому $_*$-ранг этой матрицы не может превышать $m$.

Допустим $_*$-ранг координатной матрицы меньше $m$. Согласно теореме 9.3.9, столбцы матрицы линейно зависимы справа

(9.3.47)     $a_* {}^* \lambda = 0$

Положим $c = a_* {}^* \lambda$. Из равенства (9.3.47) следует, что правая линейная комбинация строк

$$\overline{e}_* {}^* c = 0$$

векторов базиса равна 0. Это противоречит утверждению, что векторы $\overline{e}$ образуют базис. Утверждение теоремы доказано. ☐

**Доказательство.** Пусть $\overline{\overline{e}}$ - базис левого $A$-векторного пространства $V$. Согласно теореме 8.4.3, координатная матрица семейства векторов $(\overline{a}^i)$ относительно базиса $\overline{\overline{e}}$ состоит из столбцов, являющихся координатными матрицами векторов $\overline{a}^i$ относительно базиса $\overline{\overline{e}}$. Поэтому $^*_*$-ранг этой матрицы не может превышать $m$.

Допустим $^*_*$-ранг координатной матрицы меньше $m$. Согласно теореме 9.3.10, столбцы матрицы линейно зависимы слева

(9.3.48)     $\lambda^* {}_* a = 0$

Положим $c = \lambda^* {}_* a$. Из равенства (9.3.48) следует, что левая линейная комбинация строк

$$c^* {}_* \overline{e} = 0$$

векторов базиса равна 0. Это противоречит утверждению, что векторы $\overline{e}$ образуют базис. Утверждение теоремы доказано. ☐

**Теорема 9.3.15.** *Пусть $V$ - правое $A$-векторное пространство строк. Пусть $(\overline{a}^i, i \in M, |M| = m)$ семейство векторов, линейно независимых слева. Тогда $_*$-ранг их координатной матрицы равен* $m$.

**Теорема 9.3.16.** *Пусть $V$ - левое $A$-векторное пространство строк. Пусть $(\overline{a}^i, i \in M, |M| = m)$ семейство векторов, линейно независимых слева. Тогда $^*_*$-ранг их координатной матрицы равен* $m$.

Доказательство. Пусть $\overline{\overline{e}}$ - базис правого $A$-векторного пространства $V$. Согласно теореме 8.4.24, координатная матрица семейства векторов $(\overline{a}{}^i)$ относительно базиса $\overline{\overline{e}}$ состоит из строк, являющихся координатными матрицами векторов $\overline{a}{}^i$ относительно базиса $\overline{\overline{e}}$. Поэтому ${}^*_*$-ранг этой матрицы не может превышать $m$.

Допустим ${}^*_*$-ранг координатной матрицы меньше $m$. Согласно теореме 9.3.12, строки матрицы линейно зависимы справа

$$(9.3.49) \qquad a^*{}_* \lambda = 0$$

Положим $c = a^*{}_* \lambda$. Из равенства (9.3.49) следует, что правая линейная комбинация строк

$$\overline{e}^*{}_* c = 0$$

векторов базиса равна 0. Это противоречит утверждению, что векторы $\overline{e}$ образуют базис. Утверждение теоремы доказано. $\qquad\square$

Доказательство. Пусть $\overline{\overline{e}}$ - базис левого $A$-векторного пространства $V$. Согласно теореме 8.4.10, координатная матрица семейства векторов $(\overline{a}{}^i)$ относительно базиса $\overline{\overline{e}}$ состоит из строк, являющихся координатными матрицами векторов $\overline{a}{}^i$ относительно базиса $\overline{\overline{e}}$. Поэтому ${}_*{}^*$-ранг этой матрицы не может превышать $m$.

Допустим ${}_*{}^*$-ранг координатной матрицы меньше $m$. Согласно теореме 9.3.11, строки матрицы линейно зависимы слева

$$(9.3.50) \qquad \lambda_*{}^* a = 0$$

Положим $c = \lambda_*{}^* a$. Из равенства (9.3.50) следует, что левая линейная комбинация строк

$$c_*{}^* \overline{e} = 0$$

векторов базиса равна 0. Это противоречит утверждению, что векторы $\overline{e}$ образуют базис. Утверждение теоремы доказано. $\qquad\square$

## 9.4. Система линейных уравнений $a_* {}^* x = b$

Определение 9.4.1. *Предположим, что $a$ - матрица системы левых линейных уравнений*

$$x_1 a_1^1 + ... + x_m a_1^m = b_1$$
$$(9.4.1) \qquad\qquad ......$$
$$x_1 a_n^1 + ... + x_m a_n^m = b_n$$

*Мы будем называть матрицу*

$$(9.4.2) \qquad \begin{pmatrix} a_i^j \\ b_i \end{pmatrix} = \begin{pmatrix} a_1^1 & ... & a_n^1 \\ ... & ... & ... \\ a_1^m & ... & a_n^m \\ b_1 & ... & b_n \end{pmatrix}$$

**расширенной матрицей** *системы левых линейных уравнений* (9.4.1). $\qquad\square$

Определение 9.4.2. *Предположим, что $a$ - матрица системы правых линейных уравнений*

$$a_1^1 x^1 \quad + ... \quad + a_n^1 x^n \quad = b^1$$
$$(9.4.3) \qquad ... \qquad ... \qquad ... \qquad ...$$
$$a_1^m x^1 \quad + ... \quad + a_n^m x^n \quad = b^m$$

*Мы будем называть матрицу*

$$(9.4.4) \qquad \begin{pmatrix} a_i^j & b^j \end{pmatrix} = \begin{pmatrix} a_1^1 & ... & a_n^1 & b^1 \\ ... & ... & ... & ... \\ a_1^m & ... & a_n^m & b^m \end{pmatrix}$$

**расширенной матрицей** *системы правых линейных уравнений* (9.4.3). $\square$

Теорема 9.4.3. *Система правых линейных уравнений* (9.4.3) *имеет решение тогда и только тогда, когда*

$$(9.4.5) \qquad \mathrm{rank}(_*{}^*)(a_i^j) = \mathrm{rank}(_*{}^*)\begin{pmatrix} a_i^j & b^j \end{pmatrix}$$

Доказательство. Пусть $a_T^S$ - $_*{}^*$-главная подматрица матрицы $a$.

Пусть система правых линейных уравнений (9.4.3) имеет решение $x^i = d^i$. Тогда

$$(9.4.6) \qquad a_*{}^* d = b$$

Уравнение (9.4.6) может быть записано в форме

$$(9.4.7) \qquad a_{T*}{}^* d^T + a_{N \setminus T*}{}^* d^{N \setminus T} = b$$

Подставляя (9.3.23) в (9.4.7), мы получим

$$(9.4.8) \qquad a_{T*}{}^* d^T + a_{T*}{}^* r_*{}^* d^{N \setminus T} = b$$

Из (9.4.8) следует, что $A$-столбец $b$ является $_*{}^*D$-линейной комбинацией $A$-столбцов $a_T$

$$a_{T*}{}^*(d^T + r_*{}^* d^{N \setminus T}) = b$$

Это эквивалентно уравнению (9.4.5).

Нам осталось доказать, что существование решения системы правых линейных уравнений (9.4.3) следует из (9.4.5). Истинность (9.4.5) означает, что $a_T^S$ - так же $_*{}^*$-главная подматрица расширенной матрицы. Из теоремы 9.3.6 следует, что столбец $b$ является правой линейной композицией столбцов $a_T$

$$b = a_{T*}{}^* R^T$$

Полагая $R^r = 0$, мы получим

$$b = a_*{}^* R$$

Следовательно, мы нашли, по крайней мере, одно решение системы правых линейных уравнений (9.4.3). $\square$

Теорема 9.4.4. *Предположим, что* (9.4.3) - *система правых линейных уравнений, удовлетворяющих* (9.4.5). *Если* $\mathrm{rank}(_*{}^*)a = k \le m$, *то решение системы зависит от произвольных значений* $m - k$ *переменных, не включённых в* $_*{}^*$-*главную подматрицу.*

Доказательство. Пусть $a_T^S$ - $_*{}^*$-главная подматрица матрицы $a$. Предположим, что

$$(9.4.9) \qquad a^p{}_*{}^* x = b^p$$

уравнение с номером $p$. Применяя теорему 9.3.8 к расширенной матрице (9.4.2), мы получим

$$(9.4.10) \qquad a^p = R^p{}_*{}^*a^S$$

$$(9.4.11) \qquad b^p = R^p{}_*{}^*b^S$$

Подставляя (9.4.10) и (9.4.11) в (9.4.9), мы получим

$$(9.4.12) \qquad R^p{}_*{}^*a^S{}_*{}^*x = R^p{}_*{}^*b^S$$

(9.4.12) означает, что мы можем исключить уравнение (9.4.9) из системы (9.4.3) и новая система эквивалентна старой. Следовательно, число уравнений может быть уменьшено до $k$.

В этом случае у нас есть два варианта. Если число переменных также равно $k$, то согласно теореме 9.2.14 система имеет единственное решение

$$\boxed{(9.2.28) \quad x = \mathcal{H}\det({}^*_*)a^*{}_*b}$$

Если число переменных $m > k$, то мы можем передвинуть $m - k$ переменных, которые не включены в ${}_*{}^*$-главную подматрицу в правой части. Присваивая произвольные значения этим переменным, мы определяем значение правой части и для этого значения мы получим единственное решение согласно теореме 9.2.14. $\qquad \Box$

СЛЕДСТВИЕ 9.4.5. *Система правых линейных уравнений* (9.4.3) *имеет единственное решение тогда и только тогда, когда её матрица невырожденная.* $\qquad \Box$

ТЕОРЕМА 9.4.6. *Решения однородной системы правых линейных уравнений*

$$(9.4.13) \qquad a_*{}^*x = 0$$

*порождают правое $A$-векторное пространство столбцов.*

ДОКАЗАТЕЛЬСТВО. Пусть $X$ - множество решений системы правых $A$-линейных уравнений (9.4.13). Предположим, что $x = (x^a) \in X$ и $y = (y^a) \in X$ Тогда

$$a^i_j x^j = 0$$
$$a^i_j y^j = 0$$

Следовательно,

$$[a^i_j(x^j + y^j) = [a^i_j x^j + [a^i_j y^j = 0$$
$$x + y = (x^j + y^j) \in X$$

Таким же образом мы видим

$$a^i_j(x^j b) = (a^i_j x^j)b = 0$$
$$xb = (x^j b) \in X$$

Согласно определению 7.4.1, $X$ является правым $A$-векторным пространством столбцов. $\qquad \Box$

## 9.5. Невырожденная матрица

Предположим, что нам дана $n \times n$ матрица $a$. Теоремы 9.3.11 и 9.3.9 говорят нам, что если $\operatorname{rank}(_*{}^*)a < n$, то строки линейно зависимы слева и столбцы линейно зависимы справа. [7]

**ТЕОРЕМА 9.5.1.** *Пусть $a$ - $n \times n$ матрица и столбец $a_r$ - левая линейная комбинация других столбцов. Тогда $\operatorname{rank}(_*{}^*)a < n$.*

ДОКАЗАТЕЛЬСТВО. Утверждение, что столбец $a_r$ является правой линейной комбинацией других столбцов, означает, что система линейных уравнений

$$a_r = a_{[r]*}{}^*\lambda$$

имеет, по крайней мере, одно решение. Согласно теореме 9.4.3

$$\operatorname{rank}(_*{}^*)a = \operatorname{rank}(_*{}^*)a_{[r]}$$

Так как число столбцов меньше, чем $n$, то $\operatorname{rank}(_*{}^*)a_{[r]} < n$. □

**ТЕОРЕМА 9.5.2.** *Пусть $a$ - $n \times n$ матрица и строка $a^p$ - правая линейная комбинация других строк. Тогда $\operatorname{rank}(_*{}^*)a < n$.*

ДОКАЗАТЕЛЬСТВО. Доказательство утверждения похоже на доказательство теоремы 9.5.1. □

**ТЕОРЕМА 9.5.3.** *Предположим, что $a$ и $b$ - $n \times n$ матрицы и*

$$(9.5.1) \qquad c = a_*{}^*b$$

*$c$ - $_*{}^*$-вырожденная матрица тогда и только тогда, когда либо матрица $a$, либо матрица $b$ - $_*{}^*$-вырожденная матрица.*

ДОКАЗАТЕЛЬСТВО. Предположим, что матрица $b$ - $_*{}^*$-вырожденная. Согласно теореме 9.3.6 столбцы матрицы $b$ линейно зависимы справа. Следовательно,

$$(9.5.2) \qquad 0 = b_*{}^*\lambda$$

где $\lambda \neq 0$. Из (9.5.1) и (9.5.2) следует, что

$$c_*{}^*\lambda = a_*{}^*b_*{}^*\lambda = 0$$

Согласно теореме 9.5.1 матрица $c$ является $_*{}^*$-вырожденной.

Предположим, что матрица $b$ - не $_*{}^*$-вырожденная, но матрица $a$ - $_*{}^*$-вырожденная. Согласно теореме 9.3.8, столбцы матрицы $a$ линейно зависимы справа. Следовательно,

$$(9.5.3) \qquad 0 = a_*{}^*\mu$$

где $\mu \neq 0$. Согласно теореме 9.2.14, система линейных уравнений

$$b_*{}^*\lambda = \mu$$

имеет единственное решение, где $\lambda \neq 0$. Следовательно,

$$c_*{}^*\lambda = a_*{}^*b_*{}^*\lambda = a_*{}^*\mu = 0$$

---

[7] Это утверждение похоже на утверждение [6]-1.2.5.

Согласно теореме 9.5.1 матрица $c$ является $_*^*$-вырожденной.

Предположим, что матрица $c$ - $_*^*$-вырожденная матрица. Согласно теореме 9.3.8 столбцы матрицы $c$ линейно зависимы справа. Следовательно,

$$(9.5.4) \qquad 0 = c_*{}^*\lambda$$

где $\lambda \neq 0$. Из (9.5.1) и (9.5.4) следует, что

$$0 = a_*{}^*b_*{}^*\lambda$$

Если

$$0 = b_*{}^*\lambda$$

выполнено, то матрица $b$ - $_*^*$-вырожденная. Предположим, что матрица $b$ не $_*^*$-вырожденная. Положим

$$\mu = b_*{}^*\lambda$$

где $\mu \neq 0$. Тогда

$$(9.5.5) \qquad 0 = a_*{}^*\mu$$

Из (9.5.5) следует, что матрица $a$ - $_*^*$-вырожденная. $\qquad\square$

Согласно теореме 8.3.11, мы можем записать подобные утверждения для левой линейной комбинации строк или правой линейной комбинации столбцов и $^*_*$-квазидетерминанта.

**Теорема 9.5.4.** *Пусть $a$ - $n \times n$ матрица и столбец $a_r$ - правая линейная комбинация других столбцов. Тогда $\mathrm{rank}(^*_*)a < n$.* $\qquad\square$

**Теорема 9.5.5.** *Пусть $a$ - $n \times n$ матрица и строка $a^p$ - левая линейная комбинация других строк. Тогда $\mathrm{rank}(^*_*)a < n$.* $\qquad\square$

**Теорема 9.5.6.** *Предположим, что $a$ и $b$ - $n \times n$ матрицы и*

$$(9.5.6) \qquad c = a^*{}_*b$$

*$c$ - $^*_*$-вырожденная матрица тогда и только тогда, когда либо матрица $a$, либо матрица $b$ - $^*_*$-вырожденная матрица.* $\qquad\square$

**Определение 9.5.7.** $_*^*$-**матричная группа** $GL(n_*{}^*A)$ - это группа $_*^*$-невырожденных матриц, где мы определяем $_*^*$-произведение матриц (8.1.7) и $_*^*$-обратную матрица $a^{-1}{}_*^*$. $\qquad\square$

**Определение 9.5.8.** $^*_*$-**матричная группа** $GL(n^*{}_*A)$ - это группа $^*_*$-невырожденных матриц где мы определяем $^*_*$-произведение матриц (8.1.10) и $^*_*$-обратную матрицу $a^{-1}{}^*_*$. $\qquad\square$

**Теорема 9.5.9.**

$$GL(n_*{}^*A) \neq GL(n^*{}_*A)$$

**Замечание 9.5.10.** *Из теорем 8.2.9, 8.3.9 следует, что существуют матрицы, которые $^*_*$-невырожденны и $_*^*$-невырожденны. Теорема 9.5.9 означает, что множества $^*_*$-невырожденных матриц и $_*^*$-невырожденных матриц не* *совпадают. Например, существует такая $_*^*$-невырожденная матрица, которая $^*_*$-вырожденная матрица.*  □

Доказательство. Это утверждение достаточно доказать для $n = 2$. Предположим, что каждая $^*_*$-вырожденная матрица

$$(9.5.7) \qquad a = \begin{pmatrix} a_1^1 & a_2^1 \\ a_1^2 & a_2^2 \end{pmatrix}$$

является $_*^*$-вырожденной матрицей. Из теорем 9.3.8, 9.3.6 следует, что $_*^*$-вырожденная матрица удовлетворяет условию

$$(9.5.8) \qquad a_2^1 = b a_1^1$$

$$(9.5.9) \qquad a_2^2 = b a_1^2$$

$$(9.5.10) \qquad a_1^2 = a_1^1 c$$

$$(9.5.11) \qquad a_2^2 = a_2^1 c$$

Если мы подставим (9.5.10) в (9.5.9), мы получим

$$a_2^2 = b \; a_1^1 \; c$$

$b$ и $c$ - произвольные $A$-числа и $_*^*$-вырожденная матрица (9.5.7) имеет вид $(d = a_1^1)$

$$(9.5.12) \qquad a = \begin{pmatrix} d & dc \\ bd & bdc \end{pmatrix}$$

Подобным образом мы можем показать, что $^*_*$-вырожденная матрица имеет вид

$$(9.5.13) \qquad a = \begin{pmatrix} d & c'd \\ db' & c'db' \end{pmatrix}$$

Из предположения следует, что (9.5.13) и (9.5.12) представляют одну и ту же матрицу. Сравнивая (9.5.13) и (9.5.12), мы получим, что для любых $d, c \in D$ существует такое $c' \in D$, которое не зависит от $d$ и удовлетворяет уравнению

$$dc = c'd$$

Это противоречит утверждению, что $A$ - алгебра с делением.  □

Пример 9.5.11. *В случае алгебры кватернионов мы положим $b = 1 + k$, $c = j$, $d = k$. Тогда мы имеем*

$$a = \begin{pmatrix} k & kj \\ (1+k)k & (1+k)kj \end{pmatrix} = \begin{pmatrix} k & -i \\ k-1 & -i-j \end{pmatrix}$$

$$\det(^*_*)_2^2 \, a = a_2^2 - a_2^1 (a_1^1)^{-1} a_1^2$$
$$= -i - j - (k-1)(k)^{-1}(-i) = -i - j - (k-1)(-k)(-i)$$
$$= -i - j - kki + ki = -i - j + i + j = 0$$

$$\det({}_*{}^*)^{\underline{1}}_{\underline{1}}\,a = a^1_1 - a^2_1(a^2_2)^{-1}a^1_2$$

$$= k - (k-1)(-i-j)^{-1}(-i) = k - (k-1)\frac{1}{2}(i+j)(-i)$$

$$= k + \frac{1}{2}((k-1)i + (k-1)j)i = k + \frac{1}{2}(ki - i + kj - j)i$$

$$= k + \frac{1}{2}(j - i - i - j)i = k - ii = k + 1$$

$$\det({}_*{}^*)^{\underline{2}}_{\underline{1}}\,a = a^2_1 - a^1_1(a^1_2)^{-1}a^2_2$$

$$= k - 1 - k(-i)^{-1}(-i-j) = k - 1 + ki(i+j)$$

$$= k - 1 + ji(i+j) = k - 1 + ji + jj$$

$$= k - 1 - k - 1 = -2$$

$$\det({}_*{}^*)^{\underline{1}}_{\underline{2}}\,a = a^1_2 - a^2_2(a^2_1)^{-1}a^1_1$$

$$= (-i) - (-i-j)(k-1)^{-1}k = -i + (i+j)\frac{1}{2}(-k-1)k$$

$$= -i - \frac{1}{2}(i+j)(k+1)k = -i - \frac{1}{2}(ik + i + jk + j)k$$

$$= -i - \frac{1}{2}(-j + i + i + j)k = -i - ik = -i + j$$

$$\det({}^*{}_*)^{\underline{2}}_{\underline{2}}\,a = a^2_2 - a^1_2(a^1_1)^{-1}a^2_1$$

$$= -i - j - (-i)(k)^{-1}(k-1) = -i - j + i(-k)(k-1)$$

$$= -i - j + j(k-1) = -i - j + jk - j = -i - j + i - j = -2j$$

*Система линейных уравнений*

$$(9.5.14) \qquad \begin{pmatrix} k & -i \\ k-1 & -i-j \end{pmatrix}{}_*{}^* \begin{pmatrix} x^1 \\ x^2 \end{pmatrix} = \begin{pmatrix} b^1 \\ b^2 \end{pmatrix}$$

*имеет ${}_*{}^*$-вырожденную матрицу. Мы можем записать систему линейных уравнений* (9.5.14) *в виде*

$$\begin{cases} k\ x^1 - \quad\ i\ x^2 = b^1 \\ (k-1)\ x^1 - (i+j)\ x^2 = b^2 \end{cases}$$

*Система линейных уравнений*

$$(9.5.15) \qquad \begin{pmatrix} k & -i \\ k-1 & -i-j \end{pmatrix}{}^*{}_* \begin{pmatrix} x_1 & x_2 \end{pmatrix} = \begin{pmatrix} b_1 & b_2 \end{pmatrix}$$

*имеет ${}^*{}_*$-невырожденную матрицу. Мы можем записать систему линейных уравнений* (9.5.15) *в виде*

$$\begin{cases} k\ x_1 - \quad\ i\ x_1 \\ + (k-1)\ x_2 - (i+j)\ x_2 \\ = \quad\ b_1 = \quad\ b_2 \end{cases} \qquad \begin{cases} k\ x_1 + (k-1)\ x_2 = b_1 \\ -i\ x_1 - (i+j)\ x_2 = b_2 \end{cases}$$

*Система линейных уравнений*

$$(9.5.16) \qquad \begin{pmatrix} x_1 & x_2 \end{pmatrix} {}_*{}^* \begin{pmatrix} k & -i \\ k-1 & -i-j \end{pmatrix} = \begin{pmatrix} b_1 & b_2 \end{pmatrix}$$

*имеет $_*{}^*$-вырожденную матрицу. Мы можем записать систему линейных уравнений (9.5.16) в виде*

$$\begin{cases} x_1 k & -x_1 i \\ + x_2(k-1) & -x_2(i+j) \\ = b_1 & = b_2 \end{cases} \qquad \begin{cases} x_1 k & +x_2(k-1) = b_1 \\ -x_1 i & -x_2(i+j) = b_2 \end{cases}$$

*Система линейных уравнений*

$$(9.5.17) \qquad \begin{pmatrix} x^1 \\ x^2 \end{pmatrix} {}^*{}_* \begin{pmatrix} k & -i \\ k-1 & -i-j \end{pmatrix} = \begin{pmatrix} b^1 \\ b^2 \end{pmatrix}$$

*имеет $^*{}_*$-невырожденную матрицу. Мы можем записать систему линейных уравнений (9.5.17) в виде*

$$\begin{cases} x^1 \, k & - x^2 \, i & = b^1 \\ x^1 \, (k-1) - x^2 \, (i+j) & = b^2 \end{cases}$$

$\square$

## 9.6. Обратная матрица

**Определение 9.6.1.** *Пусть $a$ - $n \times m$ матрица*

$$a = \begin{pmatrix} a_1^1 & \dots & a_n^1 \\ \dots & \dots & \dots \\ a_1^m & \dots & a_n^m \end{pmatrix}$$

*и*

$$k = \min(n, m)$$

*Матрица $a$ называется $^*{}_*$-невырожденной слева, если существует $m \times n$ матрица $b$*

$$b = \begin{pmatrix} b_1^1 & \dots & b_m^1 \\ \dots & \dots & \dots \\ b_1^n & \dots & b_m^n \end{pmatrix}$$

*такая, что $^*{}_*$-произведение*

$$(9.6.1) \qquad c = b^*{}_* a$$

**Определение 9.6.2.** *Пусть $a$ - $n \times m$ матрица*

$$a = \begin{pmatrix} a_1^1 & \dots & a_n^1 \\ \dots & \dots & \dots \\ a_1^m & \dots & a_n^m \end{pmatrix}$$

*и*

$$k = \min(n, m)$$

*Матрица $a$ называется $_*{}^*$-невырожденной справа, если существует $m \times n$ матрица $b$*

$$b = \begin{pmatrix} b_1^1 & \dots & b_m^1 \\ \dots & \dots & \dots \\ b_1^n & \dots & b_m^n \end{pmatrix}$$

*такая, что $_*{}^*$-произведение*

$$(9.6.2) \qquad c = a_*{}^* b$$

содержит подматрицу

$$c_J^I = E_k$$

и остальные элементы матрицы c равны 0. Матрица b называется левой ${}^*_*$-обратной матрице a.  □

ОПРЕДЕЛЕНИЕ 9.6.3. *Пусть* a - $n \times m$ *матрица*

$$a = \begin{pmatrix} a_1^1 & ... & a_n^1 \\ ... & ... & ... \\ a_1^m & ... & a_n^m \end{pmatrix}$$

*и*

$$k = \min(n, m)$$

*Матрица* a *называется* ${}^*_*$*-невырожденной слева, если существует* $m \times n$ *матрица* b

$$b = \begin{pmatrix} b_1^1 & ... & b_m^1 \\ ... & ... & ... \\ b_1^n & ... & b_m^n \end{pmatrix}$$

*такая, что* ${}^*_*$*-произведение*

(9.6.3)          $c = b {}^*_* a$

*содержит подматрицу*

$$c_J^I = E_k$$

*и остальные элементы матрицы* c *равны* 0. *Матрица* b *называется левой* ${}^*_*$*-обратной матрице* a.  □

ТЕОРЕМА 9.6.5. *Пусть*

$$k = \min(n, m)$$

*Матрица* a *имеет левую* ${}^*_*$*-обратную матрицу* b *тогда и только тогда, когда*

(9.6.5)          $\operatorname{rank}({}^*_*)a = k$

ДОКАЗАТЕЛЬСТВО. Пусть

$$J = [1, k]$$

Пусть I - множество индексов линейно независимых слева строк матрицы a. Если $|I| < k$, то мы дополним

---

содержит подматрицу

$$c_J^I = E_k$$

и остальные элементы матрицы c равны 0. Матрица b называется правой ${}^*_*$-обратной матрице a.  □

ОПРЕДЕЛЕНИЕ 9.6.4. *Пусть* a - $n \times m$ *матрица*

$$a = \begin{pmatrix} a_1^1 & ... & a_n^1 \\ ... & ... & ... \\ a_1^m & ... & a_n^m \end{pmatrix}$$

*и*

$$k = \min(n, m)$$

*Матрица* a *называется* ${}^*_*$*-невырожденной справа, если существует* $m \times n$ *матрица* b

$$b = \begin{pmatrix} b_1^1 & ... & b_m^1 \\ ... & ... & ... \\ b_1^n & ... & b_m^n \end{pmatrix}$$

*такая, что* ${}^*_*$*-произведение*

(9.6.4)          $c = a {}^*_* b$

*содержит подматрицу*

$$c_J^I = E_k$$

*и остальные элементы матрицы* c *равны* 0. *Матрица* b *называется правой* ${}^*_*$*-обратной матрице* a.  □

ТЕОРЕМА 9.6.6. *Пусть*

$$k = \min(n, m)$$

*Матрица* a *имеет правую* ${}^*_*$*-обратную матрицу* b *тогда и только тогда, когда*

(9.6.8)          $\operatorname{rank}({}^*_*)a = k$

ДОКАЗАТЕЛЬСТВО. Пусть

$$J = [1, k]$$

Пусть I - множество индексов линейно независимых справа строк матрицы a. Если $|I| < k$, то мы дополним множество $I$ индексами оставшихся строк матрицы $a$ таким образом, чтобы получить $|I| = k$. Согласно определению 8.1.4, равенство

$$(9.6.6) \qquad c_J^I = b_J {}^* {}_* a^I$$

является следствием равенства

$$\boxed{(9.6.1) \quad c = b^* {}_* a}$$

Согласно определениям 8.1.4, 9.6.1, для каждого индекса $j \in J$, система линейных уравнений

$$b_j^1 a_{1}^{i_1} + ... + b_j^n a_n^{i_1} = \delta_j^1$$
$$(9.6.7) \qquad\qquad ......$$
$$b_j^1 a_{1}^{i_k} + ... + b_j^n a_n^{i_k} = \delta_j^k$$

является следствием равенства (9.6.6). Мы положим, что неизвестные переменные здесь - это $x^i = b_j^i$.

Если мы рассмотрим все система линейных уравнений (9.6.7), то их расширенная матрица будет иметь вид

$$\begin{pmatrix} a_1^{i_1} & ... & a_n^{i_1} & \delta_1^1 & ... & \delta_k^1 \\ ... & ... & ... & ... & ... & ... \\ a_1^{i_k} & ... & a_n^{i_k} & \delta_1^k & ... & \delta_k^k \end{pmatrix}$$

Ранг расширенной матрицы равен $k$. Согласно теореме 9.4.3, системы линейных уравнений (9.6.7) имеют решение тогда и только тогда, когда выполнено условие (9.6.5). В частности, множество строк $a^i$, $i \in I$, линейно независимо слева. □

---

множество $I$ индексами оставшихся строк матрицы $a$ таким образом, чтобы получить $|I| = k$. Согласно определению 8.1.3, равенство

$$(9.6.9) \qquad c_J^I = a^I {}_* {}^* b_J$$

является следствием равенства

$$\boxed{(9.6.2) \quad c = a_* {}^* b}$$

Согласно определениям 8.1.3, 9.6.2, для каждого индекса $j \in J$, система линейных уравнений

$$a_1^{i_1} b_j^1 + ... + a_n^{i_1} b_j^n = \delta_j^1$$
$$(9.6.10) \qquad\qquad ......$$
$$a_1^{i_k} b_j^1 + ... + a_n^{i_k} b_j^n = \delta_j^k$$

является следствием равенства (9.6.9). Мы положим, что неизвестные переменные здесь - это $x^i = b_j^i$.

Если мы рассмотрим все система линейных уравнений (9.6.10), то их расширенная матрица будет иметь вид

$$\begin{pmatrix} a_1^{i_1} & ... & a_n^{i_1} & \delta_1^1 & ... & \delta_k^1 \\ ... & ... & ... & ... & ... & ... \\ a_1^{i_k} & ... & a_n^{i_k} & \delta_1^k & ... & \delta_k^k \end{pmatrix}$$

Ранг расширенной матрицы равен $k$. Согласно теореме 9.4.3, системы линейных уравнений (9.6.10) имеют решение тогда и только тогда, когда выполнено условие (9.6.8). В частности, множество строк $a^i$, $i \in I$, линейно независимо справа. □

---

Теорема 9.6.7. *Пусть*

$$k = \min(n, m)$$

*Матрица a имеет левую $_*{}^*$-обратную матрицу b тогда и только тогда, когда*

$$(9.6.11) \qquad \mathrm{rank}(_*{}^*)a = k$$

Доказательство. Пусть

$$I = [1, k]$$

Пусть $J$ - множество индексов линейно независимых слева столб-

---

Теорема 9.6.8. *Пусть*

$$k = \min(n, m)$$

*Матрица a имеет правую $^*{}_*$-обратную матрицу b тогда и только тогда, когда*

$$(9.6.14) \qquad \mathrm{rank}(^*{}_*)a = k$$

Доказательство. Пусть

$$I = [1, k]$$

Пусть $J$ - множество индексов линейно независимых справа столб-цов матрицы $a$. Если $|J| < k$, то мы дополним множество $J$ индексами оставшихся столбцов матрицы $a$ таким образом, чтобы получить $|J| = k$. Согласно определению 8.1.3, равенство

$$(9.6.12) \qquad c_J^I = b^I {}_* * a_J$$

является следствием равенства

$$\boxed{(9.6.3) \quad c = b_* * a}$$

Согласно определениям 8.1.3, 9.6.3, для каждого индекса $i \in I$, система линейных уравнений

$$b_1^i a_{j_1}^1 + \ldots + b_m^i a_{j_1}^m = \delta_1^i$$
$$(9.6.13) \qquad \ldots\ldots$$
$$b_1^i a_{j_k}^1 + \ldots + b_m^i a_{j_k}^m = \delta_k^i$$

является следствием равенства (9.6.12). Мы положим, что неизвестные переменные здесь - это $x_j = b_j^i$.

Если мы рассмотрим все система линейных уравнений (9.6.13), то их расширенная матрица будет иметь вид

$$\begin{pmatrix} a_{j_1}^1 & \ldots & a_{j_k}^1 \\ \ldots & \ldots & \ldots \\ a_{j_1}^m & \ldots & a_{j_k}^m \\ \delta_1^1 & \ldots & \delta_k^1 \\ \ldots & \ldots & \ldots \\ \delta_1^k & \ldots & \delta_k^k \end{pmatrix}$$

Ранг расширенной матрицы равен $k$. Согласно теореме 9.4.3, системы линейных уравнений (9.6.13) имеют решение тогда и только тогда, когда выполнено условие (9.6.11). В частности, множество столбцов $a_j$, $j \in J$, линейно независимо слева. $\qquad \square$

цов матрицы $a$. Если $|J| < k$, то мы дополним множество $J$ индексами оставшихся столбцов матрицы $a$ таким образом, чтобы получить $|J| = k$. Согласно определению 8.1.4, равенство

$$(9.6.15) \qquad c_J^I = a_J * * b^I$$

является следствием равенства

$$\boxed{(9.6.4) \quad c = a^* * b}$$

Согласно определениям 8.1.4, 9.6.4, для каждого индекса $i \in I$, система линейных уравнений

$$a_{j_1}^1 b_1^i + \ldots + a_{j_1}^m b_m^i = \delta_1^i$$
$$(9.6.16) \qquad \ldots\ldots$$
$$a_{j_k}^1 b_1^i + \ldots + a_{j_k}^m b_m^i = \delta_k^i$$

является следствием равенства (9.6.15). Мы положим, что неизвестные переменные здесь - это $x_j = b_j^i$.

Если мы рассмотрим все система линейных уравнений (9.6.16), то их расширенная матрица будет иметь вид

$$\begin{pmatrix} a_{j_1}^1 & \ldots & a_{j_k}^1 \\ \ldots & \ldots & \ldots \\ a_{j_1}^m & \ldots & a_{j_k}^m \\ \delta_1^1 & \ldots & \delta_k^1 \\ \ldots & \ldots & \ldots \\ \delta_1^k & \ldots & \delta_k^k \end{pmatrix}$$

Ранг расширенной матрицы равен $k$. Согласно теореме 9.4.3, системы линейных уравнений (9.6.16) имеют решение тогда и только тогда, когда выполнено условие (9.6.14). В частности, множество столбцов $a_j$, $j \in J$, линейно независимо справа. $\qquad \square$

Если $m = n$, то левая $^*_*$-обратная матрица определена однозначно. Вообще говоря, это неверно, если $m \neq n$. Рассмотрим систему линейных уравнений (9.6.7).

- Если $m < n$, то система линейных уравнений (9.6.7) содержит неизвестных больше, чем число уравнений. Решение системы линейных уравнений (9.6.7) определено однозначно с точностью до выбора значений избыточных переменных.

- Если $m > n$, то система линейных уравнений (9.6.7) содержит $k$ неизвестных и $k$ уравнений. Следовательно, система линейных уравнений (9.6.7) имеет единственное решение. Однако выбор множества $I$ неоднозначен.

ОПРЕДЕЛЕНИЕ 9.6.9. *Множество* $left(a^{-1^*}{}_*)$ *- это множество левых* $^*_*$-*обратных матриц матрице* $a$. □

ОПРЕДЕЛЕНИЕ 9.6.10. *Множество* $right(a^{-1^*}{}_*)$ *- это множество правых* $^*_*$-*обратных матриц матрице* $a$. □

ОПРЕДЕЛЕНИЕ 9.6.11. *Множество* $left(a^{-1}{}_*^*)$ *- это множество левых* $_*^*$-*обратных матриц матрице* $a$. □

ОПРЕДЕЛЕНИЕ 9.6.12. *Множество* $right(a^{-1}{}_*^*)$ *- это множество правых* $_*^*$-*обратных матриц матрице* $a$. □

ТЕОРЕМА 9.6.13. $a \in left(b^{-1^*}{}_*)$ *тогда и только тогда, когда* $b \in right(a^{-1^*}{}_*)$.

ТЕОРЕМА 9.6.14. $a \in right(b^{-1^*}{}_*)$ *тогда и только тогда, когда* $b \in left(a^{-1^*}{}_*)$.

ДОКАЗАТЕЛЬСТВО. Теорема является следствием определений 9.6.1, 9.6.4, 9.6.9, 9.6.12. □

ДОКАЗАТЕЛЬСТВО. Теорема является следствием определений 9.6.2, 9.6.3, 9.6.10, 9.6.11. □

# Гомоморфизм левого $A$-векторного пространства

## 10.1. Общее определение

Пусть $A_i$, $i = 1, 2$, - алгебра с делением над коммутативным кольцом $D_i$. Пусть $V_i$, $i = 1, 2$, - левое $A_i$-модуль.

ОПРЕДЕЛЕНИЕ 10.1.1. *Пусть*

(10.1.1)

$$g_{1.12}(d) : a \to d\,a$$
$$g_{1.23}(v) : w \to C_1(w, v)$$
$$C_1 \in \mathcal{L}(A_1^2 \to A_1)$$
$$g_{1.34}(a) : v \to a\,v$$
$$g_{1.14}(d) : v \to d\,v$$

*диаграмма представлений, описывающая левое $A_1$-векторное пространство $V_1$. Пусть*

(10.1.2)

$$g_{2.12}(d) : a \to d\,a$$
$$g_{2.23}(v) : w \to C_2(w, v)$$
$$C_2 \in \mathcal{L}(A_2^2 \to A_2)$$
$$g_{2.34}(a) : v \to a\,v$$
$$g_{2.14}(d) : v \to d\,v$$

*диаграмма представлений, описывающая левое $A_2$-векторное пространство $V_2$. Морфизм*

(10.1.3) $$(h : D_1 \to D_2 \quad g : A_1 \to A_2 \quad f : V_1 \to V_2)$$

*диаграммы представлений* (10.1.1) *в диаграмму представлений* (10.1.2) *называется* **гомоморфизмом** *левого $A_1$-векторного пространства $V_1$ в левое $A_2$-векторное пространство $V_2$. Обозначим* $\operatorname{Hom}(D_1 \to D_2; A_1 \to A_2; *V_1 \to *V_2)$ *множество гомоморфизмов левого $A_1$-векторного пространства $V_1$ в левое $A_2$-векторное пространство $V_2$.* $\square$

Мы будем пользоваться записью

$$f \circ a = f(a)$$

для образа гомоморфизма $f$.

ТЕОРЕМА 10.1.2. *Гомоморфизм*

(10.1.3) $$(h : D_1 \to D_2 \quad g : A_1 \to A_2 \quad f : V_1 \to V_2)$$

*левого $A_1$-векторного пространства $V_1$ в левое $A_2$-векторное пространство $V_2$ удовлетворяет следующим равенствам*

(10.1.4) $$h(p + q) = h(p) + h(q)$$

(10.1.5) $$h(pq) = h(p)h(q)$$

$$(10.1.6) \qquad g \circ (a + b) = g \circ a + g \circ b$$

$$(10.1.7) \qquad g \circ (ab) = (g \circ a)(g \circ b)$$

$$(10.1.8) \qquad g \circ (pa) = h(p)(g \circ a)$$

$$(10.1.9) \qquad f \circ (u + v) = f \circ u + f \circ v$$

$$(10.1.10) \qquad f \circ (av) = (g \circ a)(f \circ v)$$

$$p, q \in D_1 \quad a, b \in A_1 \quad u, v \in V_1$$

Доказательство. Согласно определениям 2.3.9, 2.8.5, 7.2.1, 11 отображение

$$(h : D_1 \to D_2 \quad g : A_1 \to A_2)$$

является гомоморфизмом $D_1$-алгебры $A_1$ в $D_2$-алгебру $A_2$. Следовательно, равенства (10.1.4), (10.1.5), (10.1.6), (10.1.7), (10.1.8) являются следствием теоремы 4.3.2. Равенство (10.1.9) является следствием определения 10.1.1, так как, согласно определению 2.3.9, отображение $f$ является гомоморфизмом абелевой группы. Равенство (10.1.10) является следствием равенства (2.3.6), так как отображение

$$(g : A_1 \to A_2 \quad f : V_1 \to V_2)$$

является морфизмом представления $g_{1.34}$ в представление $g_{2.34}$.　　□

Определение 10.1.3. *Гомоморфизм* [10.1]

$$(h : D_1 \to D_2 \quad g : A_1 \to A_2 \quad f : V_1 \to V_2)$$

*называется* **изоморфизмом** *между левым $A_1$-векторным пространством $V_1$ и левым $A_2$-векторным пространством $V_2$, если существует отображение*

$$(h^{-1} : D_2 \to D_1 \quad g^{-1} : A_2 \to A_1 \quad f^{-1} : V_2 \to V_1)$$

*которое является гомоморфизмом.*　　□

## 10.1.1. Гомоморфизм левого $A$-векторного пространства столбцов.

Теорема 10.1.4. *Гомоморфизм* [10.2]

$$\boxed{(10.1.3) \quad (h : D_1 \to D_2 \quad \overline{g} : A_1 \to A_2 \quad \overline{f} : V_1 \to V_2)}$$

*левого $A_1$-векторного пространства столбцов $V_1$ в левое $A_2$-векторное пространство столбцов $V_2$ имеет представление*

$$(10.1.11) \qquad b = gh(a)$$

$$(10.1.12) \qquad \overline{g} \circ (a^i e_{A_1 i}) = h(a^i) g_i^k e_{A_2 k}$$

---

[10.1] Я следую определению на странице [18]-63.

[10.2] В теоремах 10.1.4, 10.1.5, мы опираемся на следующее соглашение. Пусть $i = 1$, 2. Пусть множество векторов $\overline{\overline{e}}_{A_i} = (e_{A_i k}, k \in K_i)$ является базисом и $C_{i \cdot ij}^{\cdot \cdot k}$, $k$, $i$, $j \in K_i$, - структурные константы $D_i$-алгебры столбцов $A_i$. Пусть множество векторов $\overline{\overline{e}}_{V_i} = (e_{V_i i}, i \in I_i)$ является базисом левого $A_i$-векторного пространства $V_i$.

$$(10.1.13) \qquad \overline{g} \circ (e_{A_1} a) = e_{A_2} gh(a)$$

$$(10.1.14) \qquad w = (\overline{g} \circ v)^*{}_* f$$

$$(10.1.15) \qquad \overline{f} \circ (v^i e_{V_1 i}) = (\overline{g} \circ v^i) f_i^k e_{V_2 k}$$

$$(10.1.16) \qquad \overline{f} \circ (v^*{}_* e_{V_1}) = (\overline{g} \circ v)^*{}_* (f^*{}_* e_{V_2})$$

*относительно выбранных базисов. Здесь*

- *$a$ - координатная матрица $A_1$-числа $\overline{a}$ относительно базиса $\overline{\overline{e}}_{A_1}$ (A)*

$$\overline{a} = e_{A_1} a$$

- *$h(a) = (h(a_k), k \in K)$ - матрица $D_2$-чисел.*
- *$b$ - координатная матрица $A_2$-числа*

$$\overline{b} = \overline{g} \circ \overline{a}$$

  *относительно базиса $\overline{\overline{e}}_{A_2}$*

$$\overline{b} = e_{A_2} b$$

- *$g$ - координатная матрица множества $A_2$-чисел $(\overline{g} \circ e_{A_1 k}, k \in K)$ относительно базиса $\overline{\overline{e}}_{A_2}$.*
- *Матрица гомоморфизма и структурные константы связаны соотношением*

$$(10.1.17) \qquad h(C_{1\,ij}^{\ k}) g_k^l = g_i^p g_j^q C_{2\,pq}^{\ \ l}$$

- *$v$ - координатная матрица $V_1$-числа $\overline{v}$ относительно базиса $\overline{\overline{e}}_{V_1}$ (A)*

$$\overline{v} = v^*{}_* e_{V_1}$$

- *$g(v) = (g(v_i), i \in I)$ - матрица $A_2$-чисел.*
- *$w$ - координатная матрица $V_2$-числа*

$$\overline{w} = \overline{f} \circ \overline{v}$$

  *относительно базиса $\overline{\overline{e}}_{V_2}$*

$$\overline{w} = w^*{}_* e_{V_2}$$

- *$f$ - координатная матрица множества $V_2$-чисел $(\overline{f} \circ e_{V_1 i}, i \in I)$ относительно базиса $\overline{\overline{e}}_{V_2}$.*

*Для заданного гомоморфизма*

$$\boxed{(10.1.3) \quad (h : D_1 \to D_2 \quad \overline{g} : A_1 \to A_2 \quad \overline{f} : V_1 \to V_2)}$$

*множество матриц $(g, f)$ определено однозначно и называется* **координатами гомоморфизма** (10.1.3) *относительно базисов* $(\overline{\overline{e}}_{A_1}, \overline{\overline{e}}_{V_1})$, $(\overline{\overline{e}}_{A_2}, \overline{\overline{e}}_{V_2})$.

**Доказательство.** Согласно определениям 2.3.9, 2.8.5, 5.4.1, 10.1.1, отображение

$$(h : D_1 \to D_2 \quad \overline{g} : A_1 \to A_2)$$

является гомоморфизмом $D_1$-алгебры $A_1$ в $D_2$-алгебру $A_2$. Следовательно, равенства (10.1.11), (10.1.12), (10.1.13), (10.1.17) следуют из равенств

$$\boxed{(5.4.11) \quad b = fh(a)}$$

$$\boxed{(5.4.12) \quad \overline{f} \circ (a^i e_{1i}) = h(a^i) f_i^k e_{2k}}$$

$$(5.4.13) \quad \overline{f} \circ (e_1 a) = e_2 f h(a)$$

$$(5.4.16) \quad h(C_1{}_{ij}^{\,k}) f_k^l = f_i^p f_j^q C_2{}_{pq}^{\,l}$$

Из теоремы 4.3.3, следует что матрица $g$ определена однозначно.

Вектор $\overline{v} \in V_1$ имеет разложение

$$(10.1.18) \qquad \overline{v} = v^* {}_* e_{V_1}$$

относительно базиса $\overline{\overline{e}}_{V_1}$. Вектор $\overline{w} \in V_2$ имеет разложение

$$(10.1.19) \qquad \overline{w} = w^* {}_* e_{V_2}$$

относительно базиса $\overline{\overline{e}}_{V_2}$. Так как $\overline{f}$ - гомоморфизм, то равенство

$$(10.1.20) \qquad \overline{w} = \overline{f} \circ \overline{v} = \overline{f} \circ (v^* {}_* e_{V_1})$$
$$= (\overline{g} \circ v)^* {}_* (\overline{f} \circ e_{V_1})$$

является следствием равенств (10.1.9), (10.1.10). $V_2$-число $\overline{f} \circ e_{V_1 i}$ имеет разложение

$$(10.1.21) \qquad \overline{f} \circ e_{V_1 i} = f_i^* {}_* e_{V_2} = f_i^j e_{V_2 j}$$

относительно базиса $\overline{\overline{e}}_{V_2}$. Равенство

$$(10.1.22) \qquad \overline{w} = (\overline{g} \circ v)^* {}_* f^* {}_* e_{V_2}$$

является следствием равенств (10.1.20), (10.1.21). Равенство (10.1.14) следует из сравнения (10.1.19) и (10.1.22) и теоремы 8.4.4. Из равенства (10.1.21) и теоремы 8.4.4 следует что матрица $f$ определена однозначно. $\qquad \square$

ТЕОРЕМА 10.1.5. *Пусть отображение*

$$h : D_1 \to D_2$$

*является гомоморфизмом кольца $D_1$ в кольцо $D_2$.*     *Пусть*

$$g = (g_l^k, k \in K, l \in L)$$

*матрица $D_2$-чисел, которая удовлетворяют равенству*

$$(10.1.17) \quad h(C_1{}_{ij}^{\,k}) g_k^l = g_i^p g_j^q C_2{}_{pq}^{\,l}$$

*Пусть*

$$f = (f_j^i, i \in I, j \in J)$$

*матрица $A_2$-чисел.*     *Отображение* [10.2]

$$(10.1.3) \quad (h : D_1 \to D_2 \quad \overline{g} : A_1 \to A_2 \quad \overline{f} : V_1 \to V_2)$$

*определённое равенствами*

$$(10.1.13) \quad \overline{g} \circ (e_{A_1} a) = e_{A_2} g h(a)$$

$$(10.1.16) \quad \overline{f} \circ (v^* {}_* e_{V_1}) = (\overline{g} \circ v)^* {}_* (f^* {}_* e_{V_2})$$

*является гомоморфизмом левого $A_1$-векторного пространства столбцов $V_1$ в левое $A_2$-векторное пространство столбцов $V_2$.*     *Гомоморфизм*

$$(10.1.3) \quad (h : D_1 \to D_2 \quad \overline{g} : A_1 \to A_2 \quad \overline{f} : V_1 \to V_2)$$

*который имеет данное множество матриц $(g, f)$, определён однозначно.*

ДОКАЗАТЕЛЬСТВО. Согласно теореме 5.4.5, отображение

$$(10.1.23) \qquad (h : D_1 \to D_2 \quad \overline{g} : A_1 \to A_2)$$

является гомоморфизмом $D_1$-алгебры $A_1$ в $D_2$-алгебру $A_2$ и гомоморфизм (10.1.23) определён однозначно.

Равенство

$$(10.1.24) \quad \begin{aligned} &\overline{f} \circ ((v+w)^*{}_* e_{V_1}) \\ &= (g \circ (v+w))^*{}_* f^*{}_* e_{V_2} \\ &= (g \circ v)^*{}_* f^*{}_* e_{V_2} \\ &+ (g \circ w)^*{}_* f^*{}_* e_{V_2} \\ &= \overline{f} \circ (v^*{}_* e_{V_1}) + \overline{f} \circ (w^*{}_* e_{V_1}) \end{aligned}$$

является следствием равенств

| | |
|---|---|
| (8.1.19) $\quad (b_1+b_2)^*{}_* a = b_1{}^*{}_* a + b_2{}^*{}_* a$ | (10.1.6) $\quad g \circ (a+b) = g \circ a + g \circ b$ |
| (10.1.10) $\quad f \circ (av) = (g \circ a)(f \circ v)$ | (10.1.16) $\quad \overline{f} \circ (v^*{}_* e_{V_1}) = (\overline{g} \circ v)^*{}_* (f^*{}_* e_{V_2})$ |

Из равенства (10.1.24) следует, что отображение $\overline{f}$ является гомоморфизмом абелевой группы. Равенство

$$(10.1.25) \quad \begin{aligned} &\overline{f} \circ ((av)^*{}_* e_{V_1}) \\ &= (\overline{g} \circ (av))^*{}_* f^*{}_* e_{V_2} \\ &= (\overline{g} \circ a)((\overline{g} \circ v)^*{}_* f^*{}_* e_{V_2}) \\ &= (\overline{g} \circ a)(\overline{f} \circ (v^*{}_* e_{V_1})) \end{aligned}$$

является следствием равенств

| |
|---|
| (10.1.10) $\quad f \circ (av) = (g \circ a)(f \circ v)$ |
| (5.4.4) $\quad g \circ (ab) = (g \circ a)(g \circ b)$ |
| (10.1.16) $\quad \overline{f} \circ (v^*{}_* e_{V_1}) = (\overline{g} \circ v)^*{}_* (f^*{}_* e_{V_2})$ |

Из равенства (10.1.25), и определений 2.8.5, 10.1.1 следует, что отображение

| |
|---|
| (10.1.3) $\quad (h : D_1 \to D_2 \quad \overline{g} : A_1 \to A_2 \quad \overline{f} : V_1 \to V_2)$ |

является гомоморфизмом $A_1$-векторного пространства столбцов $V_1$ в $A_2$-векторное пространство столбцов $V_2$.

Пусть $f$ - матрица гомоморфизмов $\overline{f}$, $\overline{g}$ относительно базисов $\overline{\overline{e}}_V$, $\overline{\overline{e}}_W$. Равенство

$$\overline{f} \circ (v^*{}_* e_V) = v^*{}_* f^*{}_* e_W = \overline{g} \circ (v^*{}_* e_V)$$

является следствием теоремы 10.3.6. Следовательно, $\overline{f} = \overline{g}$. $\qquad\square$

На основании теорем 10.1.4, 10.1.5 мы идентифицируем гомоморфизм

| |
|---|
| (10.1.3) $\quad (h : D_1 \to D_2 \quad \overline{g} : A_1 \to A_2 \quad \overline{f} : V_1 \to V_2)$ |

левого $A$-векторного пространства $V$ столбцов и координаты его представления

| |
|---|
| (10.1.13) $\quad \overline{g} \circ (e_{A_1} a) = e_{A_2} gh(a)$ |
| (10.1.16) $\quad \overline{f} \circ (v^*{}_* e_{V_1}) = (\overline{g} \circ v)^*{}_* (f^*{}_* e_{V_2})$ |

## 10.1.2. Гомоморфизм левого $A$-векторного пространства строк.

Теорема 10.1.6. *Гомоморфизм* [10.3]

$$(10.1.3) \quad (h : D_1 \to D_2 \quad \overline{g} : A_1 \to A_2 \quad \overline{f} : V_1 \to V_2)$$

*левого $A_1$-векторного пространства строк $V_1$ в левое $A_2$-векторное пространство строк $V_2$ имеет представление*

$$(10.1.26) \qquad b = h(a)g$$

$$(10.1.27) \qquad \overline{g} \circ (a_i e^i_{A_1}) = h(a_i) g^i_k e^k_{A_2}$$

$$(10.1.28) \qquad \overline{g} \circ (a e_{A_1}) = h(a) g e_{A_2}$$

$$(10.1.29) \qquad w = (\overline{g} \circ v)_*{}^* f$$

$$(10.1.30) \qquad \overline{f} \circ (v_i e^i_{V_1}) = (\overline{g} \circ v_i) f^i_k e^k_{V_2}$$

$$(10.1.31) \qquad \overline{f} \circ (v_*{}^* e_{V_1}) = (\overline{g} \circ v)_*{}^* f_*{}^* e_{V_2}$$

*относительно выбранных базисов. Здесь*

- *$a$ - координатная матрица $A_1$-числа $\overline{a}$ относительно базиса $\overline{\overline{e}}_{A_1}$ (A)*

$$\overline{a} = a e_{A_1}$$

- *$h(a) = (h(a^k), k \in K)$ - матрица $D_2$-чисел.*
- *$b$ - координатная матрица $A_2$-числа*

$$\overline{b} = \overline{g} \circ \overline{a}$$

  *относительно базиса $\overline{\overline{e}}_{A_2}$*

$$\overline{b} = b e_{A_2}$$

- *$g$ - координатная матрица множества $A_2$-чисел $(\overline{g} \circ e^k_{A_1}, k \in K)$ относительно базиса $\overline{\overline{e}}_{A_2}$.*
- *Матрица гомоморфизма и структурные константы связаны соотношением*

$$(10.1.32) \qquad h(C_{1\,k}^{\;ij}) g^k_l = g^i_p g^j_q C_{2\,l}^{\;pq}$$

- *$v$ - координатная матрица $V_1$-числа $\overline{v}$ относительно базиса $\overline{\overline{e}}_{V_1}$ (A)*

$$\overline{v} = v_*{}^* e_{V_1}$$

- *$g(v) = (g(v^i), i \in I)$ - матрица $A_2$-чисел.*
- *$w$ - координатная матрица $V_2$-числа*

$$\overline{w} = \overline{f} \circ \overline{v}$$

  *относительно базиса $\overline{\overline{e}}_{V_2}$*

$$\overline{w} = w_*{}^* e_{V_2}$$

---

[10.3] В теоремах 10.1.6, 10.1.7, мы опираемся на следующее соглашение. Пусть $i = 1, 2$. Пусть множество векторов $\overline{\overline{e}}_{A_i} = (e^k_{A_i}, k \in K_i)$ является базисом и $C_{i\,k}^{\;ij}$, $k$, $i$, $j \in K_i$, - структурные константы $D_i$-алгебры строк $A_i$. Пусть множество векторов $\overline{\overline{e}}_{V_i} = (e^i_{V_i}, i \in I_i)$ является базисом левого $A_i$-векторного пространства $V_i$.

- $f$ - координатная матрица множества $V_2$-чисел  $\left(\overline{f} \circ e_{V_1}^i, i \in I\right)$  относительно базиса  $\overline{\overline{e}}_{V_2}$.

*Для заданного гомоморфизма*

(10.1.3)   $(h : D_1 \to D_2 \quad \overline{g} : A_1 \to A_2 \quad \overline{f} : V_1 \to V_2)$

*множество матриц* $(g, f)$ *определено однозначно и называется* **координатами гомоморфизма** (10.1.3)   *относительно базисов* $(\overline{\overline{e}}_{A_1}, \overline{\overline{e}}_{V_1})$,   $(\overline{\overline{e}}_{A_2}, \overline{\overline{e}}_{V_2})$.

Доказательство.   Согласно определениям   2.3.9, 2.8.5, 5.4.1, 10.1.1, отображение

$$(h : D_1 \to D_2 \quad \overline{g} : A_1 \to A_2)$$

является гомоморфизмом $D_1$-алгебры $A_1$ в $D_2$-алгебру $A_2$. Следовательно, равенства   (10.1.26), (10.1.27), (10.1.28), (10.1.32)   следуют из равенств

(5.4.22)   $b = h(a)f$

(5.4.23)   $\overline{f} \circ (a_i e_1^i) = h(a_i)f_k^i e_2^k$

(5.4.24)   $\overline{f} \circ (ae_1) = h(a)fe_2$

(5.4.27)   $h(C_{1k}^{ij})f_l^k = f_p^i f_q^j C_{2l}^{pq}$

Из теоремы 4.3.4, следует что матрица $g$ определена однозначно.

Вектор  $\overline{v} \in V_1$  имеет разложение

(10.1.33)                         $\overline{v} = v_*{}^* e_{V_1}$

относительно базиса  $\overline{\overline{e}}_{V_1}$. Вектор  $\overline{w} \in V_2$  имеет разложение

(10.1.34)                         $\overline{w} = w_*{}^* e_{V_2}$

относительно базиса  $\overline{\overline{e}}_{V_2}$. Так как  $\overline{f}$ - гомоморфизм, то равенство

(10.1.35)       $\overline{w} = \overline{f} \circ \overline{v} = \overline{f} \circ (v_*{}^* e_{V_1})$

$$= (\overline{g} \circ v)_*{}^* (\overline{f} \circ e_{V_1})$$

является следствием равенств (10.1.9), (10.1.10). $V_2$-число  $\overline{f} \circ e_{V_1}^i$  имеет разложение

(10.1.36)                         $\overline{f} \circ e_{V_1}^i = f^i{}_*{}^* e_{V_2} = f_j^i e_{V_2}^j$

относительно базиса  $\overline{\overline{e}}_{V_2}$. Равенство

(10.1.37)                         $\overline{w} = (\overline{g} \circ v)_*{}^* f_*{}^* e_{V_2}$

является следствием равенств (10.1.35), (10.1.36). Равенство (10.1.29) следует из сравнения (10.1.34) и (10.1.37) и теоремы 8.4.11. Из равенства (10.1.36) и теоремы 8.4.11 следует что матрица $f$ определена однозначно.   $\square$

Теорема 10.1.7. *Пусть отображение*

$$h : D_1 \to D_2$$

*является гомоморфизмом кольца $D_1$ в кольцо $D_2$.*       *Пусть*

$$g = (g_k^l, k \in K, l \in L)$$

*матрица $D_2$-чисел, которая удовлетворяют равенству*

$$(10.1.32) \quad h(C_{1\,k}^{ij})g_l^k = g_p^i g_q^j C_{2l}^{pq}$$

*Пусть*

$$f = (f_i^j, i \in I, j \in J)$$

*матрица $A_2$-чисел.          Отображение* [10.3]

$$(10.1.3) \quad (h : D_1 \to D_2 \quad \overline{g} : A_1 \to A_2 \quad \overline{f} : V_1 \to V_2)$$

*определённое равенствами*

$$(10.1.28) \quad \overline{g} \circ (ae_{A_1}) = h(a)ge_{A_2}$$

$$(10.1.31) \quad \overline{f} \circ (v_*{}^*e_{V_1}) = (\overline{g} \circ v)_*{}^* f_*{}^* e_{V_2}$$

*является гомоморфизмом левого $A_1$-векторного пространства строк $V_1$ в ле-
вое $A_2$-векторное пространство строк $V_2$.    Гомоморфизм*

$$(10.1.3) \quad (h : D_1 \to D_2 \quad \overline{g} : A_1 \to A_2 \quad \overline{f} : V_1 \to V_2)$$

*который имеет  данное множество матриц $(g, f)$,  определён однозначно.*

Доказательство.   Согласно теореме [5.4.6], отображение

$$(10.1.38) \quad (h : D_1 \to D_2 \quad \overline{g} : A_1 \to A_2)$$

является гомоморфизмом $D_1$-алгебры $A_1$ в $D_2$-алгебру $A_2$ и гомоморфизм
(10.1.38) определён однозначно.

Равенство

$$
\begin{aligned}
(10.1.39) \quad & \overline{f} \circ ((v+w)_*{}^*e_{V_1}) \\
& = (g \circ (v+w))_*{}^* f_*{}^* e_{V_2} \\
& = (g \circ v)_*{}^* f_*{}^* e_{V_2} \\
& + (g \circ w)_*{}^* f_*{}^* e_{V_2} \\
& = \overline{f} \circ (v_*{}^*e_{V_1}) + \overline{f} \circ (w_*{}^*e_{V_1})
\end{aligned}
$$

является следствием равенств

$$(8.1.17) \quad (b_1 + b_2)_*{}^*a = b_{1*}{}^*a + b_{2*}{}^*a \qquad (10.1.6) \quad g \circ (a+b) = g \circ a + g \circ b$$

$$(10.1.10) \quad f \circ (av) = (g \circ a)(f \circ v) \qquad (10.1.31) \quad \overline{f} \circ (v_*{}^*e_{V_1}) = (\overline{g} \circ v)_*{}^* f_*{}^* e_{V_2}$$

Из равенства (10.1.39) следует, что отображение $\overline{f}$ является гомоморфизмом
абелевой группы. Равенство

$$
\begin{aligned}
(10.1.40) \quad & \overline{f} \circ ((av)_*{}^*e_{V_1}) \\
& = (\overline{g} \circ (av))_*{}^* f^*{}_* e_{V_2} \\
& = (\overline{g} \circ a)((\overline{g} \circ v)_*{}^* f_*{}^* e_{V_2}) \\
& = (\overline{g} \circ a)(\overline{f} \circ (e_{V_1}{}^*{}_* v))
\end{aligned}
$$

является следствием равенств

$$(10.1.10) \quad f \circ (av) = (g \circ a)(f \circ v)$$

$$(5.4.4) \quad g \circ (ab) = (g \circ a)(g \circ b)$$

$$(10.1.31) \quad \overline{f} \circ (v_*{}^*e_{V_1}) = (\overline{g} \circ v)_*{}^* f_*{}^* e_{V_2}$$

Из равенства (10.1.40), и определений 2.8.5, 10.1.1 следует, что отображение

$$(10.1.3) \quad (h : D_1 \to D_2 \quad \overline{g} : A_1 \to A_2 \quad \overline{f} : V_1 \to V_2)$$

является гомоморфизмом $A_1$-векторного пространства строк $V_1$ в $A_2$-векторное пространство строк $V_2$.

Пусть $f$ - матрица гомоморфизмов $\overline{f}$, $\overline{g}$ относительно базисов $\overline{\overline{e}}_V$, $\overline{\overline{e}}_W$. Равенство

$$\overline{f} \circ (v_* {}^* e_V) = v_* {}^* f_* {}^* e_W = \overline{g} \circ (v_* {}^* e_V)$$

является следствием теоремы 10.3.7. Следовательно, $\overline{f} = \overline{g}$. $\qquad\square$

На основании теорем 10.1.6, 10.1.7 мы идентифицируем гомоморфизм

$$(10.1.3) \quad (h : D_1 \to D_2 \quad \overline{g} : A_1 \to A_2 \quad \overline{f} : V_1 \to V_2)$$

левого $A$-векторного пространства $V$ строк и координаты его представления

$$(10.1.28) \quad \overline{g} \circ (a e_{A_1}) = h(a) g e_{A_2}$$

$$(10.1.31) \quad \overline{f} \circ (v_* {}^* e_{V_1}) = (\overline{g} \circ v)_* {}^* f_* {}^* e_{V_2}$$

## 10.2. Гомоморфизм, когда кольца $D_1 = D_2 = D$

Пусть $A_i$, $i = 1, 2$, - алгебра с делением над коммутативным кольцом $D$. Пусть $V_i$, $i = 1, 2$, - левое $A_i$-модуль.

ОПРЕДЕЛЕНИЕ 10.2.1. *Пусть*

$$
\begin{aligned}
&g_{1.12}(d) : a \to d\,a \\
&g_{1.23}(v) : w \to C_1(w, v) \\
&\qquad C_1 \in \mathcal{L}(A_1^2 \to A_1) \\
&g_{1.34}(a) : v \to a\,v \\
&g_{1.14}(d) : v \to d\,v
\end{aligned}
$$

(10.2.1)

*диаграмма представлений, описывающая левое $A_1$-векторное пространство $V_1$. Пусть*

$$
\begin{aligned}
&g_{2.12}(d) : a \to d\,a \\
&g_{2.23}(v) : w \to C_2(w, v) \\
&\qquad C_2 \in \mathcal{L}(A_2^2 \to A_2) \\
&g_{2.34}(a) : v \to a\,v \\
&g_{2.14}(d) : v \to d\,v
\end{aligned}
$$

(10.2.2)

*диаграмма представлений, описывающая левое $A_2$-векторное пространство $V_2$. Морфизм*

$$(10.2.3) \qquad\qquad (g : A_1 \to A_2 \quad f : V_1 \to V_2)$$

*диаграммы представлений* (10.2.1) *в диаграмму представлений* (10.2.2) *называется* **гомоморфизмом** *левого $A_1$-векторного пространства $V_1$ в левое $A_2$-* *векторное пространство $V_2$. Обозначим* $\mathrm{Hom}(D; A_1 \to A_2; *V_1 \to *V_2)$ *множество гомоморфизмов левого $A_1$-векторного пространства $V_1$ в левое $A_2$-векторное пространство $V_2$.* □

Мы будем пользоваться записью

$$f \circ a = f(a)$$

для образа гомоморфизма $f$.

Теорема 10.2.2. *Гомоморфизм*

$$\boxed{(10.2.3) \quad (g : A_1 \to A_2 \quad f : V_1 \to V_2)}$$

*левого $A_1$-векторного пространства $V_1$ в левое $A_2$-векторное пространство $V_2$ удовлетворяет следующим равенствам*

$$(10.2.4) \qquad\qquad g \circ (a + b) = g \circ a + g \circ b$$

$$(10.2.5) \qquad\qquad g \circ (ab) = (g \circ a)(g \circ b)$$

$$(10.2.6) \qquad\qquad g \circ (pa) = p(g \circ a)$$

$$(10.2.7) \qquad\qquad f \circ (u + v) = f \circ u + f \circ v$$

$$(10.2.8) \qquad\qquad f \circ (av) = (g \circ a)(f \circ v)$$

$$p \in D \quad a, b \in A_1 \quad u, v \in V_1$$

Доказательство. Согласно определениям 2.3.9, 2.8.5, 7.2.1, 1 отображение

$$g : A_1 \to A_2$$

является гомоморфизмом $D$-алгебры $A_1$ в $D$-алгебру $A_2$. Следовательно, равенства (10.2.4), (10.2.5), (10.2.6) являются следствием теоремы 4.3.10. Равенство (10.2.7) является следствием определения 10.2.1, так как, согласно определению 2.3.9, отображение $f$ является гомоморфизмом абелевой группы. Равенство (10.2.8) является следствием равенства (2.3.6), так как отображение

$$(g : A_1 \to A_2 \quad f : V_1 \to V_2)$$

является морфизмом представления $g_{1.34}$ в представление $g_{2.34}$. □

Определение 10.2.3. *Гомоморфизм* [10.4]

$$(g : A_1 \to A_2 \quad f : V_1 \to V_2)$$

*называется* **изоморфизмом** *между левым $A_1$-векторным пространством $V_1$ и левым $A_2$-векторным пространством $V_2$, если существует отображение*

$$(g^{-1} : A_2 \to A_1 \quad f^{-1} : V_2 \to V_1)$$

*которое является гомоморфизмом.* □

---

[10.4] Я следую определению на странице [18]-63.

**Теорема 10.2.4.** *Гомоморфизм* [10.5]

$$(10.2.3) \quad (\overline{g}: A_1 \to A_2 \quad \overline{f}: V_1 \to V_2)$$

*левого $A_1$-векторного пространства столбцов $V_1$ в левое $A_2$-векторное пространство столбцов $V_2$ имеет представление*

$$(10.2.9) \quad b = ga$$

$$(10.2.10) \quad \overline{g} \circ (a^i e_{A_1 i}) = a^i g_i^k e_{A_2 k}$$

$$(10.2.11) \quad \overline{g} \circ (e_{A_1} a) = e_{A_2} ga$$

$$(10.2.12) \quad w = (\overline{g} \circ v)_{*}{}^{*} f$$

$$(10.2.13) \quad \overline{f} \circ (v^i e_{V_1 i}) = (\overline{g} \circ v^i) f_i^k e_{V_2 k}$$

$$(10.2.14)$$
$$\overline{f} \circ (v^{*} e_{V_1}) = (\overline{g} \circ v)^{*} (f^{*}{}_{*} e_{V_2})$$

*относительно выбранных базисов. Здесь*

- *$a$ - координатная матрица $A_1$-числа $\overline{a}$ относительно базиса $\overline{\overline{e}}_{A_1}$ (A)*

  $$\overline{a} = e_{A_1} a$$

- *$b$ - координатная матрица $A_2$-числа*

  $$\overline{b} = \overline{g} \circ \overline{a}$$

  *относительно базиса $\overline{\overline{e}}_{A_2}$*

  $$\overline{b} = e_{A_2} b$$

- *$g$ - координатная матрица множества $A_2$-чисел ($\overline{g} \circ e_{A_1 k}, k \in K$) относительно базиса $\overline{\overline{e}}_{A_2}$.*

**Теорема 10.2.5.** *Гомоморфизм* [10.6]

$$(10.2.3) \quad (\overline{g}: A_1 \to A_2 \quad \overline{f}: V_1 \to V_2)$$

*левого $A_1$-векторного пространства строк $V_1$ в левое $A_2$-векторное пространство строк $V_2$ имеет представление*

$$(10.2.16) \quad b = ag$$

$$(10.2.17) \quad \overline{g} \circ (a_i e_{A_1}^i) = a_i g_k^i e_{A_2}^k$$

$$(10.2.18) \quad \overline{g} \circ (ae_{A_1}) = age_{A_2}$$

$$(10.2.19) \quad w = (\overline{g} \circ v)_{*}{}^{*} f$$

$$(10.2.20) \quad \overline{f} \circ (v_i e_{V_1}^i) = (\overline{g} \circ v_i) f_k^i e_{V_2}^k$$

$$(10.2.21) \quad \overline{f} \circ (v_{*} e_{V_1}) = (\overline{g} \circ v)_{*} f_{*}{}^{*} e_{V_2}$$

*относительно выбранных базисов. Здесь*

- *$a$ - координатная матрица $A_1$-числа $\overline{a}$ относительно базиса $\overline{\overline{e}}_{A_1}$ (A)*

  $$\overline{a} = ae_{A_1}$$

- *$b$ - координатная матрица $A_2$-числа*

  $$\overline{b} = \overline{g} \circ \overline{a}$$

  *относительно базиса $\overline{\overline{e}}_{A_2}$*

  $$\overline{b} = be_{A_2}$$

- *$g$ - координатная матрица множества $A_2$-чисел ($\overline{g} \circ e_{A_1}^k, k \in K$) относительно базиса $\overline{\overline{e}}_{A_2}$.*

---

[10.5] В теоремах 10.2.4, 10.2.6, мы опираемся на следующее соглашение. Пусть $i = 1$, 2. Пусть множество векторов $\overline{\overline{e}}_{A_i} = (e_{A_i k}, k \in K_i)$ является базисом и $C_{ij}^k$, $k$, $i$, $j \in K_i$, - структурные константы $D$-алгебры столбцов $A_i$. Пусть множество векторов $\overline{\overline{e}}_{V_i} = (e_{V_i i}, i \in I_i)$ является базисом левого $A_i$-векторного пространства $V_i$.

[10.6] В теоремах 10.2.5, 10.2.7, мы опираемся на следующее соглашение. Пусть $i = 1, 2$. Пусть множество векторов $\overline{\overline{e}}_{A_i} = (e_{A_i}^k, k \in K_i)$ является базисом и $C_{ik}^{ij}$, $k$, $i$, $j \in K_i$, - структурные константы $D$-алгебры строк $A_i$. Пусть множество векторов $\overline{\overline{e}}_{V_i} = (e_{V_i}^i, i \in I_i)$ является базисом левого $A_i$-векторного пространства $V_i$.

- *Матрица гомоморфизма и структурные константы связаны соотношением*

$$(10.2.15) \qquad C_{1\,ij}^{\phantom{1}k} g_k^l = g_i^p g_j^q C_{2\,pq}^{\phantom{2}l}$$

- *$v$ - координатная матрица $V_1$-числа $\overline{v}$ относительно базиса $\overline{\overline{e}}_{V_1}$ ($A$)*

$$\overline{v} = v^*{}_*e_{V_1}$$

- *$g(v) = (g(v_i),\ i \in I)$ - матрица $A_2$-чисел.*

- *$w$ - координатная матрица $V_2$-числа*

$$\overline{w} = \overline{f} \circ \overline{v}$$

*относительно базиса $\overline{\overline{e}}_{V_2}$*

$$\overline{w} = w^*{}_*e_{V_2}$$

- *$f$ - координатная матрица множества $V_2$-чисел ($\overline{f} \circ e_{V_1 i},\ i \in I$) относительно базиса $\overline{\overline{e}}_{V_2}$.*

*Для заданного гомоморфизма*

$$\boxed{(10.2.3) \quad (\overline{g} : A_1 \to A_2 \quad \overline{f} : V_1 \to V_2)}$$

*множество матриц $(g, f)$ определено однозначно и называется* **координатами гомоморфизма** (10.2.3) *относительно базисов $(\overline{\overline{e}}_{A_1}, \overline{\overline{e}}_{V_1})$, $(\overline{\overline{e}}_{A_2}, \overline{\overline{e}}_{V_2})$.*

Доказательство. Согласно определениям 2.3.9, 2.8.5, 5.4.7, 10.2.1, отображение

$$\overline{g} : A_1 \to A_2$$

является гомоморфизмом $D$-алгебры $A_1$ в $D$-алгебру $A_2$. Следовательно, равенства (10.2.9), (10.2.10), (10.2.11), (10.2.15) следуют из равенств

$$\boxed{(5.4.53) \quad b = fa}$$

$$\boxed{(5.4.54) \quad \overline{f} \circ (a^i e_{1i}) = a^i f_i^k e_k}$$

$$\boxed{(5.4.55) \quad \overline{f} \circ (e_{1a}) = e_2 fa}$$

$$\boxed{(5.4.58) \quad C_{1\,ij}^{\phantom{1}k} f_k^l = f_i^p f_j^q C_{2\,pq}^{\phantom{2}l}}$$

Из теоремы 4.3.11, следует что матрица $g$ определена однозначно.

- *Матрица гомоморфизма и структурные константы связаны соотношением*

$$(10.2.22) \qquad C_{1\,k}^{\phantom{1}ij} g_l^k = g_p^i g_q^j C_{2\,l}^{\phantom{2}pq}$$

- *$v$ - координатная матрица $V_1$-числа $\overline{v}$ относительно базиса $\overline{\overline{e}}_{V_1}$ ($A$)*

$$\overline{v} = v_*{}^*e_{V_1}$$

- *$g(v) = (g(v^i),\ i \in I)$ - матрица $A_2$-чисел.*

- *$w$ - координатная матрица $V_2$-числа*

$$\overline{w} = \overline{f} \circ \overline{v}$$

*относительно базиса $\overline{\overline{e}}_{V_2}$*

$$\overline{w} = w_*{}^*e_{V_2}$$

- *$f$ - координатная матрица множества $V_2$-чисел ($\overline{f} \circ e_{V_1}^i,\ i \in I$) относительно базиса $\overline{\overline{e}}_{V_2}$.*

*Для заданного гомоморфизма*

$$\boxed{(10.2.3) \quad (\overline{g} : A_1 \to A_2 \quad \overline{f} : V_1 \to V_2)}$$

*множество матриц $(g, f)$ определено однозначно и называется* **координатами гомоморфизма** (10.2.3) *относительно базисов $(\overline{\overline{e}}_{A_1}, \overline{\overline{e}}_{V_1})$, $(\overline{\overline{e}}_{A_2}, \overline{\overline{e}}_{V_2})$.*

Доказательство. Согласно определениям 2.3.9, 2.8.5, 5.4.7, 10.2.1, отображение

$$\overline{g} : A_1 \to A_2$$

является гомоморфизмом $D$-алгебры $A_1$ в $D$-алгебру $A_2$. Следовательно, равенства (10.2.16), (10.2.17), (10.2.18), (10.2.22) следуют из равенств

$$\boxed{(5.4.64) \quad b = af}$$

$$\boxed{(5.4.65) \quad \overline{f} \circ (a_i e_1^i) = a_i f_k^i e^k}$$

$$\boxed{(5.4.66) \quad \overline{f} \circ (ae_1) = afe_2}$$

$$\boxed{(5.4.69) \quad C_{1\,k}^{\phantom{1}ij} f_l^k = f_p^i f_q^j C_{2\,l}^{\phantom{2}pq}}$$

Из теоремы 4.3.12, следует что матри-

Вектор $\overline{v} \in V_1$ имеет разложение

$$(10.2.23) \qquad \overline{v} = v^* {}_* e_{V_1}$$

относительно базиса $\overline{\overline{e}}_{V_1}$. Вектор $\overline{w} \in V_2$ имеет разложение

$$(10.2.24) \qquad \overline{w} = w^* {}_* e_{V_2}$$

относительно базиса $\overline{\overline{e}}_{V_2}$. Так как $\overline{f}$ - гомоморфизм, то равенство

$$(10.2.25) \quad \begin{aligned} \overline{w} = \overline{f} \circ \overline{v} &= \overline{f} \circ (v^* {}_* e_{V_1}) \\ &= (\overline{g} \circ v)^* {}_* (\overline{f} \circ e_{V_1}) \end{aligned}$$

является следствием равенств (10.2.7), (10.2.8). $V_2$-число $\overline{f} \circ e_{V_1 i}$ имеет разложение

$$(10.2.26) \quad \overline{f} \circ e_{V_1 i} = f_i{}^* {}_* e_{V_2} = f_i^j e_{V_2 j}$$

относительно базиса $\overline{\overline{e}}_{V_2}$. Равенство

$$(10.2.27) \qquad \overline{w} = (\overline{g} \circ v)^* {}_* f^* {}_* e_{V_2}$$

является следствием равенств (10.2.25), (10.2.26). Равенство (10.2.12) следует из сравнения (10.2.24) и (10.2.27) и теоремы 8.4.4. Из равенства (10.2.26) и теоремы 8.4.4 следует что матрица $f$ определена однозначно. $\qquad \square$

ца $g$ определена однозначно.

Вектор $\overline{v} \in V_1$ имеет разложение

$$(10.2.28) \qquad \overline{v} = v_* {}^* e_{V_1}$$

относительно базиса $\overline{\overline{e}}_{V_1}$. Вектор $\overline{w} \in V_2$ имеет разложение

$$(10.2.29) \qquad \overline{w} = w_* {}^* e_{V_2}$$

относительно базиса $\overline{\overline{e}}_{V_2}$. Так как $\overline{f}$ - гомоморфизм, то равенство

$$(10.2.30) \quad \begin{aligned} \overline{w} = \overline{f} \circ \overline{v} &= \overline{f} \circ (v_* {}^* e_{V_1}) \\ &= (\overline{g} \circ v)_* {}^* (\overline{f} \circ e_{V_1}) \end{aligned}$$

является следствием равенств (10.2.7), (10.2.8). $V_2$-число $\overline{f} \circ e_{V_1}^i$ имеет разложение

$$(10.2.31) \quad \overline{f} \circ e_{V_1}^i = f^i {}_* {}^* e_{V_2} = f_j^i e_{V_2}^j$$

относительно базиса $\overline{\overline{e}}_{V_2}$. Равенство

$$(10.2.32) \qquad \overline{w} = (\overline{g} \circ v)_* {}^* f_* {}^* e_{V_2}$$

является следствием равенств (10.2.30), (10.2.31). Равенство (10.2.19) следует из сравнения (10.2.29) и (10.2.32) и теоремы 8.4.11. Из равенства (10.2.31) и теоремы 8.4.11 следует что матрица $f$ определена однозначно. $\qquad \square$

---

**ТЕОРЕМА 10.2.6.** *Пусть*

$$g = (g_l^k, k \in K, l \in L)$$

*матрица $D$-чисел, которая удовлетворяют равенству*

$$(10.2.15) \quad C_1{}^k_{ij} g_k^l = g_i^p g_j^q C_2{}^l_{pq}$$

*Пусть*

$$f = (f_j^i, i \in I, j \in J)$$

*матрица $A_2$-чисел.*    *Отображение* 10.5

$$(10.2.3) \quad (\overline{g} : A_1 \to A_2 \quad \overline{f} : V_1 \to V_2)$$

*определённое равенствами*

$$(10.2.11) \quad \overline{f} \circ (e_{A_1} a) = e_{A_2} f a$$

$$(10.2.14) \quad \overline{f} \circ (v^* {}_* e_{V_1}) = (\overline{g} \circ v)^* {}_* (f^* {}_* e_{V_1})$$

*является гомоморфизмом левого $A_1$-векторного пространства столбцов $V_1$ в левое $A_2$-векторное пространство столбцов $V_2$. Гомоморфизм*

**ТЕОРЕМА 10.2.7.** *Пусть*

$$g = (g_k^l, k \in K, l \in L)$$

*матрица $D$-чисел, которая удовлетворяют равенству*

$$(10.2.22) \quad C_1{}^{ij}_k g_l^k = g_p^i g_l^j C_2{}^{pq}_l$$

*Пусть*

$$f = (f_i^j, i \in I, j \in J)$$

*матрица $A_2$-чисел.*    *Отображение* 10.6

$$(10.2.3) \quad (\overline{g} : A_1 \to A_2 \quad \overline{f} : V_1 \to V_2)$$

*определённое равенствами*

$$(10.2.18) \quad \overline{f} \circ (a e_{A_1}) = a f e_{A_2}$$

$$(10.2.21) \quad \overline{f} \circ (v_* {}^* e_{V_1}) = (\overline{g} \circ v)_* {}^* f_* {}^* e_{V_2}$$

*является гомоморфизмом левого $A_1$-векторного пространства строк $V_1$ в левое $A_2$-векторное пространство строк $V_2$. Гомоморфизм*

| |
|---|
| (10.2.3)   $(\overline{g}: A_1 \to A_2 \quad \overline{f}: V_1 \to V_2)$ |

который имеет данное множество матриц $(g, f)$, определён однозначно.

Доказательство. Согласно теореме 5.4.11, отображение

$$(10.2.33) \qquad \overline{g}: A_1 \to A_2$$

является гомоморфизмом $D$-алгебры $A_1$ в $D$-алгебру $A_2$ и гомоморфизм (10.2.33) определён однозначно.

Равенство

$$\overline{f} \circ ((v + w)^* {}_* e_{V_1})$$
$$= (g \circ (v + w))^* {}_* f^* {}_* e_{V_2}$$
$$(10.2.34) = (g \circ v)^* {}_* f^* {}_* e_{V_2}$$
$$+ (g \circ w)^* {}_* f^* {}_* e_{V_2}$$
$$= \overline{f} \circ (v^* {}_* e_{V_1}) + \overline{f} \circ (w^* {}_* e_{V_1})$$

является следствием равенств

| |
|---|
| (8.1.19)   $(b_1 + b_2)^* {}_* a = b_1{}^* {}_* a + b_2{}^* {}_* a$ |
| (10.1.6)   $g \circ (a + b) = g \circ a + g \circ b$ |
| (10.1.10)   $f \circ (av) = (g \circ a)(f \circ v)$ |
| (10.1.16)   $\overline{f} \circ (v^* {}_* e_{V_1}) = (\overline{g} \circ v)^* {}_* (f^* {}_* e_{V_2})$ |

Из равенства (10.2.34) следует, что отображение $\overline{f}$ является гомоморфизмом абелевой группы. Равенство

$$\overline{f} \circ ((av)^* {}_* e_{V_1})$$
$$= (\overline{g} \circ (av))^* {}_* f^* {}_* e_{V_2}$$
$$(10.2.35) = (\overline{g} \circ a)((\overline{g} \circ v)^* {}_* f^* {}_* e_{V_2})$$
$$= (\overline{g} \circ a)(\overline{f} \circ (v^* {}_* e_{V_1}))$$

является следствием равенств

| |
|---|
| (10.1.10)   $f \circ (av) = (g \circ a)(f \circ v)$ |
| (5.4.4)   $g \circ (ab) = (g \circ a)(g \circ b)$ |
| (10.1.16)   $\overline{f} \circ (v^* {}_* e_{V_1}) = (\overline{g} \circ v)^* {}_* (f^* {}_* e_{V_2})$ |

Из равенства (10.2.35), и определений 2.8.5, 10.2.1 следует, что отображение

| |
|---|
| (10.2.3)   $(\overline{g}: A_1 \to A_2 \quad \overline{f}: V_1 \to V_2)$ |

---

Доказательство. Согласно теореме 5.4.12, отображение

$$(10.2.36) \qquad \overline{g}: A_1 \to A_2$$

является гомоморфизмом $D$-алгебры $A_1$ в $D$-алгебру $A_2$ и гомоморфизм (10.2.36) определён однозначно.

Равенство

$$\overline{f} \circ ((v + w)_* {}^* e_{V_1})$$
$$= (g \circ (v + w))_* {}^* f_* {}^* e_{V_2}$$
$$(10.2.37) = (g \circ v)_* {}^* f_* {}^* e_{V_2}$$
$$+ (g \circ w)_* {}^* f_* {}^* e_{V_2}$$
$$= \overline{f} \circ (v_* {}^* e_{V_1}) + \overline{f} \circ (w_* {}^* e_{V_1})$$

является следствием равенств

| |
|---|
| (8.1.17)   $(b_1 + b_2)_* {}^* a = b_1{}_* {}^* a + b_2{}_* {}^* a$ |
| (10.1.6)   $g \circ (a + b) = g \circ a + g \circ b$ |
| (10.1.10)   $f \circ (av) = (g \circ a)(f \circ v)$ |
| (10.1.31)   $\overline{f} \circ (v_* {}^* e_{V_1}) = (\overline{g} \circ v)_* {}^* f_* {}^* e_{V_2}$ |

Из равенства (10.2.37) следует, что отображение $\overline{f}$ является гомоморфизмом абелевой группы. Равенство

$$\overline{f} \circ ((av)_* {}^* e_{V_1})$$
$$= (\overline{g} \circ (av))_* {}^* f_* {}^* e_{V_2}$$
$$(10.2.38) = (\overline{g} \circ a)((\overline{g} \circ v)_* {}^* f_* {}^* e_{V_2})$$
$$= (\overline{g} \circ a)(\overline{f} \circ (e_{V_1} {}^* v))$$

является следствием равенств

| |
|---|
| (10.1.10)   $f \circ (av) = (g \circ a)(f \circ v)$ |
| (5.4.4)   $g \circ (ab) = (g \circ a)(g \circ b)$ |
| (10.1.31)   $\overline{f} \circ (v_* {}^* e_{V_1}) = (\overline{g} \circ v)_* {}^* f_* {}^* e_{V_2}$ |

Из равенства (10.2.38), и определений 2.8.5, 10.2.1 следует, что отображение

| |
|---|
| (10.2.3)   $(\overline{g}: A_1 \to A_2 \quad \overline{f}: V_1 \to V_2)$ |

является гомоморфизмом $A_1$-векторного пространства столбцов $V_1$ в $A_2$-векторное пространство столбцов $V_2$.

Пусть $f$ - матрица гомоморфизмов $\overline{f}$, $\overline{g}$ относительно базисов $\overline{\overline{e}}_V$, $\overline{\overline{e}}_W$. Равенство

$$\overline{f} \circ (v^* {}_* e_V) = v^* {}_* f^* {}_* e_W = \overline{g} \circ (v^* {}_* e_V)$$

является следствием теоремы 10.3.6. Следовательно, $\overline{f} = \overline{g}$.   □

На основании теорем 10.2.4, 10.2.6 мы идентифицируем гомоморфизм

$$(10.2.3) \quad (\overline{g} : A_1 \to A_2 \quad \overline{f} : V_1 \to V_2)$$

левого $A$-векторного пространства $V$ столбцов и координаты его представления

$$(10.2.11) \quad \overline{f} \circ (e_{A_1} a) = e_{A_2} f a$$

$$(10.2.14) \quad \overline{f} \circ (v^* {}_* e_{V_1}) = (\overline{g} \circ v)^* {}_* (f^* {}_* e_{V_2})$$

является гомоморфизмом $A_1$-векторного пространства строк $V_1$ в $A_2$-векторное пространство строк $V_2$.

Пусть $f$ - матрица гомоморфизмов $\overline{f}$, $\overline{g}$ относительно базисов $\overline{\overline{e}}_V$, $\overline{\overline{e}}_W$. Равенство

$$\overline{f} \circ (v^* {}_* e_V) = v^* {}_* f^* {}_* e_W = \overline{g} \circ (v^* {}_* e_V)$$

является следствием теоремы 10.3.7. Следовательно, $\overline{f} = \overline{g}$.   □

На основании теорем 10.2.5, 10.2.7 мы идентифицируем гомоморфизм

$$(10.2.3) \quad (\overline{g} : A_1 \to A_2 \quad \overline{f} : V_1 \to V_2)$$

левого $A$-векторного пространства $V$ строк и координаты его представления

$$(10.2.18) \quad \overline{f} \circ (a e_{A_1}) = a f e_{A_2}$$

$$(10.2.21) \quad \overline{f} \circ (v_* {}^* e_{V_1}) = (\overline{g} \circ v)_* {}^* f_* {}^* e_{V_2}$$

## 10.3. Гомоморфизм, когда $D$-алгебры $A_1 = A_2 = A$

Пусть $A$ - алгебра с делением над коммутативным кольцом $D$. Пусть $V_i$, $i = 1, 2$, - левое $A$-векторное пространство.

ОПРЕДЕЛЕНИЕ 10.3.1. *Пусть*

$$(10.3.1)$$

$$g_{1.12}(d) : a \to d\,a$$
$$g_{1.23}(v) : w \to C(w, v)$$
$$C \in \mathcal{L}(A^2 \to A)$$
$$g_{1.34}(a) : v \to a\,v$$
$$g_{1.14}(d) : v \to d\,v$$

*диаграмма представлений, описывающая левое $A$-векторное пространство $V_1$. Пусть*

$$(10.3.2)$$

$$g_{2.12}(d) : a \to d\,a$$
$$g_{2.23}(v) : w \to C(w, v)$$
$$C \in \mathcal{L}(A^2 \to A)$$
$$g_{2.34}(a) : v \to a\,v$$
$$g_{2.14}(d) : v \to d\,v$$

*диаграмма представлений, описывающая левое $A$-векторное пространство*

$V_2$. *Морфизм*

(10.3.3)                    $$f : V_1 \to V_2$$

*диаграммы представлений* (10.3.1) *в диаграмму представлений* (10.3.2) *называется* **гомоморфизмом** *левого $A$-векторного пространства $V_1$ в левое $A$-векторное пространство $V_2$. Обозначим* $\mathrm{Hom}(D; A; *V_1 \to *V_2)$ *множество гомоморфизмов левого $A$-векторного пространства $V_1$ в левое $A$-векторное пространство $V_2$.* □

Мы будем пользоваться записью

$$f \circ a = f(a)$$

для образа гомоморфизма $f$.

Теорема 10.3.2. *Гомоморфизм*

$\boxed{(10.3.3) \quad f : V_1 \to V_2}$

*левого $A$-векторного пространства $V_1$ в левое $A$-векторное пространство $V_2$ удовлетворяет следующим равенствам*

(10.3.4)                    $$f \circ (u + v) = f \circ u + f \circ v$$

(10.3.5)                    $$f \circ (av) = a(f \circ v)$$

$$a \in A \quad u, v \in V_1$$

Доказательство. Равенство (10.3.4) является следствием определения 10.3.1, так как, согласно определению 2.3.9, отображение $f$ является гомоморфизмом абелевой группы. Равенство (10.3.5) является следствием равенства (2.3.6), так как отображение

$$f : V_1 \to V_2$$

является морфизмом представления $g_{1.34}$ в представление $g_{2.34}$. □

Определение 10.3.3. *Гомоморфизм* [10.7]

$$f : V_1 \to V_2$$

*называется* **изоморфизмом** *между левым $A$-векторным пространством $V_1$ и левым $A$-векторным пространством $V_2$, если существует отображение*

$$f^{-1} : V_2 \to V_1$$

*которое является гомоморфизмом.* □

Определение 10.3.4. *Гомоморфизм* [10.7]

$$f : V \to V$$

*источником и целью которого является одно и тоже левое $A$-векторное пространство, называется* **эндоморфизмом**. *Эндоморфизм*

$$f : V \to V$$

---

[10.7] Я следую определению на странице [18]-63.

*левого $A$-векторного пространства $V$ называется **автоморфизмом**, если существует отображение $f^{-1}$ которое является эндоморфизмом.*   $\square$

**Теорема 10.3.5.** *Множество $GL(V)$ автоморфизмов левого $A$-векторного пространства $V$ является группой.*

---

**Теорема 10.3.6.** *Гомоморфизм* [10.8]

$$\boxed{(10.3.3) \quad \overline{f} : V_1 \to V_2}$$

*левого $A$-векторного пространства столбцов $V_1$ в левое $A$-векторное пространство столбцов $V_2$ имеет представление*

$$(10.3.6) \qquad w = v^*{}_* f$$

$$(10.3.7) \qquad \overline{f} \circ (v^i e_{V_1 i}) = v^i f_i^k e_{V_2 k}$$

$$(10.3.8) \qquad \overline{f} \circ (v^*{}_* e_{V_1}) = v^*{}_* f^*{}_* e_{V_2}$$

*относительно выбранных базисов. Здесь*

- $v$ - *координатная матрица $V_1$-числа $\overline{v}$ относительно базиса $\overline{\overline{e}}_{V_1}$ $(A)$*

$$\overline{v} = v^*{}_* e_{V_1}$$

- $w$ - *координатная матрица $V_2$-числа*

$$\overline{w} = \overline{f} \circ \overline{v}$$

  *относительно базиса $\overline{\overline{e}}_{V_2}$*

$$\overline{w} = w^*{}_* e_{V_2}$$

- $f$ - *координатная матрица множества $V_2$-чисел ($\overline{f} \circ e_{V_1 i}, i \in I$) относительно базиса $\overline{\overline{e}}_{V_2}$.*

*Для заданного гомоморфизма $\overline{f}$ матрица $f$ определена однозначно и называется **матрицей гомоморфизма** $\overline{f}$ относительно базисов $\overline{\overline{e}}_1$, $\overline{\overline{e}}_2$.*

---

**Теорема 10.3.7.** *Гомоморфизм* [10.9]

$$\boxed{(10.3.3) \quad \overline{f} : V_1 \to V_2}$$

*левого $A$-векторного пространства строк $V_1$ в левое $A$-векторное пространство строк $V_2$ имеет представление*

$$(10.3.9) \qquad w = v_*{}^* f$$

$$(10.3.10) \qquad \overline{f} \circ (v_i e_{V_1}^i) = v_i f_k^i e_{V_2}^k$$

$$(10.3.11) \qquad \overline{f} \circ (v_*{}^* e_{V_1}) = v_*{}^* f_*{}^* e_{V_2}$$

*относительно выбранных базисов. Здесь*

- $v$ - *координатная матрица $V_1$-числа $\overline{v}$ относительно базиса $\overline{\overline{e}}_{V_1}$ $(A)$*

$$\overline{v} = v_*{}^* e_{V_1}$$

- $w$ - *координатная матрица $V_2$-числа*

$$\overline{w} = \overline{f} \circ \overline{v}$$

  *относительно базиса $\overline{\overline{e}}_{V_2}$*

$$\overline{w} = w_*{}^* e_{V_2}$$

- $f$ - *координатная матрица множества $V_2$-чисел ($\overline{f} \circ e_{V_1}^i, i \in I$) относительно базиса $\overline{\overline{e}}_{V_2}$.*

*Для заданного гомоморфизма $\overline{f}$ матрица $f$ определена однозначно и называется **матрицей гомоморфизма** $\overline{f}$ относительно базисов $\overline{\overline{e}}_1$, $\overline{\overline{e}}_2$.*

---

[10.8] В теоремах 10.3.6, 10.3.8, мы опираемся на следующее соглашение. Пусть $i = 1$, 2. Пусть множество векторов $\overline{\overline{e}}_{V_i} = (e_{V_i i}, i \in I_i)$ является базисом левого $A$-векторного пространства $V_i$.

[10.9] В теоремах 10.3.7, 10.3.9, мы опираемся на следующее соглашение. Пусть $i = 1$, 2. Пусть множество векторов $\overline{\overline{e}}_{V_i} = (e_{V_i}^i, i \in I_i)$ является базисом левого $A$-векторного пространства $V_i$.

**Доказательство.** Вектор $\overline{v} \in V_1$ имеет разложение

$$(10.3.12) \qquad \overline{v} = v^*{}_* e_{V_1}$$

относительно базиса $\overline{\overline{e}}_{V_1}$. Вектор $\overline{w} \in V_2$ имеет разложение

$$(10.3.13) \qquad \overline{w} = w^*{}_* e_{V_2}$$

относительно базиса $\overline{\overline{e}}_{V_2}$. Так как $\overline{f}$ - гомоморфизм, то равенство

$$(10.3.14) \quad \begin{aligned} \overline{w} = \overline{f} \circ \overline{v} &= \overline{f} \circ (v^*{}_* e_{V_1}) \\ &= v^*{}_* (\overline{f} \circ e_{V_1}) \end{aligned}$$

является следствием равенств (10.3.4), (10.3.5). $V_2$-число $\overline{f} \circ e_{V_1 i}$ имеет разложение

$$(10.3.15) \quad \overline{f} \circ e_{V_1 i} = f_i{}^*{}_* e_{V_2} = f_i^j e_{V_2 j}$$

относительно базиса $\overline{\overline{e}}_{V_2}$. Равенство

$$(10.3.16) \qquad \overline{w} = v^*{}_* f^*{}_* e_{V_2}$$

является следствием равенств (10.3.14), (10.3.15). Равенство (10.3.6) следует из сравнения (10.3.13) и (10.3.16) и теоремы 8.4.4. Из равенства (10.3.15) и теоремы 8.4.4 следует что матрица $f$ определена однозначно. $\qquad \square$

**Доказательство.** Вектор $\overline{v} \in V_1$ имеет разложение

$$(10.3.17) \qquad \overline{v} = v_*{}^* e_{V_1}$$

относительно базиса $\overline{\overline{e}}_{V_1}$. Вектор $\overline{w} \in V_2$ имеет разложение

$$(10.3.18) \qquad \overline{w} = w_*{}^* e_{V_2}$$

относительно базиса $\overline{\overline{e}}_{V_2}$. Так как $\overline{f}$ - гомоморфизм, то равенство

$$(10.3.19) \quad \begin{aligned} \overline{w} = \overline{f} \circ \overline{v} &= \overline{f} \circ (v_*{}^* e_{V_1}) \\ &= v_*{}^* (\overline{f} \circ e_{V_1}) \end{aligned}$$

является следствием равенств (10.3.4), (10.3.5). $V_2$-число $\overline{f} \circ e_{V_1}^i$ имеет разложение

$$(10.3.20) \quad \overline{f} \circ e_{V_1}^i = f^i{}_*{}^* e_{V_2} = f_i^j e_{V_2}^j$$

относительно базиса $\overline{\overline{e}}_{V_2}$. Равенство

$$(10.3.21) \qquad \overline{w} = v_*{}^* f_*{}^* e_{V_2}$$

является следствием равенств (10.3.19), (10.3.20). Равенство (10.3.9) следует из сравнения (10.3.18) и (10.3.21) и теоремы 8.4.11. Из равенства (10.3.20) и теоремы 8.4.11 следует что матрица $f$ определена однозначно. $\qquad \square$

**Теорема 10.3.8.** *Пусть*

$$f = (f_j^i, i \in I, j \in J)$$

*матрица $A$-чисел. Отображение* [10.8]

$$(10.3.3) \quad \boxed{\overline{f} : V_1 \to V_2}$$

*определённое равенством*

$$(10.3.8) \quad \boxed{\overline{f} \circ (v^*{}_* e_{V_1}) = v^*{}_* f^*{}_* e_{V_2}}$$

*является гомоморфизмом левого $A$-векторного пространства столбцов $V_1$ в левое $A$-векторное пространство столбцов $V_2$. Гомоморфизм*

$$(10.3.3) \quad \boxed{\overline{f} : V_1 \to V_2}$$

*который имеет данную матрицу $f$, определён однозначно.*

**Теорема 10.3.9.** *Пусть*

$$f = (f_i^j, i \in I, j \in J)$$

*матрица $A$-чисел. Отображение* [10.9]

$$(10.3.3) \quad \boxed{\overline{f} : V_1 \to V_2}$$

*определённое равенством*

$$(10.3.11) \quad \boxed{\overline{f} \circ (v_*{}^* e_{V_1}) = v_*{}^* f_*{}^* e_{V_2}}$$

*является гомоморфизмом левого $A$-векторного пространства строк $V_1$ в левое $A$-векторное пространство строк $V_2$. Гомоморфизм*

$$(10.3.3) \quad \boxed{\overline{f} : V_1 \to V_2}$$

*который имеет данную матрицу $f$, определён однозначно.*

Доказательство. Равенство

$$(10.3.22) \quad \begin{aligned} &\overline{f} \circ ((v+w)^*{}_*e_{V_1}) \\ &= (v+w)^*{}_*f^*{}_*e_{V_2} \\ &= v^*{}_*f^*{}_*e_{V_2} + w^*{}_*f^*{}_*e_{V_2} \\ &= \overline{f} \circ (v^*{}_*e_{V_1}) + \overline{f} \circ (w^*{}_*e_{V_1}) \end{aligned}$$

является следствием равенств

| | |
|---|---|
| (8.1.19) | $(b_1 + b_2)^*{}_*a = b_1{}^*{}_*a + b_2{}^*{}_*a$ |
| (10.2.8) | $f \circ (av) = a(f \circ v)$ |
| (10.3.8) | $\overline{f} \circ (v^*{}_*e_{V_1}) = v^*{}_*f^*{}_*e_{V_2}$ |

Из равенства (10.3.22) следует, что отображение $\overline{f}$ является гомоморфизмом абелевой группы. Равенство

$$(10.3.23) \quad \begin{aligned} &\overline{f} \circ ((av)^*{}_*e_{V_1}) \\ &= (av)^*{}_*f^*{}_*e_{V_2} \\ &= a(v^*{}_*f^*{}_*e_{V_2}) \\ &= a(\overline{f} \circ (v^*{}_*e_{V_1})) \end{aligned}$$

является следствием равенств

| | |
|---|---|
| (10.2.8) | $f \circ (av) = a(f \circ v)$ |
| (10.3.8) | $\overline{f} \circ (v^*{}_*e_{V_1}) = v^*{}_*f^*{}_*e_{V_2}$ |

Из равенства (10.3.23), и определений 2.8.5, 10.3.1 следует, что отображение

| | |
|---|---|
| (10.3.3) | $\overline{f} : V_1 \to V_2$ |

является гомоморфизмом $A$-векторного пространства столбцов $V_1$ в $A$-векторное пространство столбцов $V_2$.

Пусть $f$ - матрица гомоморфизмов $\overline{f}$, $\overline{g}$ относительно базисов $\overline{\overline{e}}_V$, $\overline{\overline{e}}_W$. Равенство

$$\overline{f} \circ (v^*{}_*e_V) = v^*{}_*f^*{}_*e_W = \overline{g} \circ (v^*{}_*e_V)$$

является следствием теоремы 10.3.6. Следовательно, $\overline{f} = \overline{g}$. □

На основании теорем 10.3.6, 10.3.8 мы идентифицируем гомоморфизм

| | |
|---|---|
| (10.3.3) | $\overline{f} : V_1 \to V_2$ |

левого $A$-векторного пространства $V$ столбцов и координаты его представления

| | |
|---|---|
| (10.3.8) | $\overline{f} \circ (v^*{}_*e_{V_1}) = v^*{}_*f^*{}_*e_{V_2}$ |

Доказательство. Равенство

$$(10.3.24) \quad \begin{aligned} &\overline{f} \circ ((v+w)_*{}^*e_{V_1}) \\ &= (v+w)_*{}^*f_*{}^*e_{V_2} \\ &= v_*{}^*f_*{}^*e_{V_2} + w_*{}^*f_*{}^*e_{V_2} \\ &= \overline{f} \circ (v_*{}^*e_{V_1}) + \overline{f} \circ (w_*{}^*e_{V_1}) \end{aligned}$$

является следствием равенств

| | |
|---|---|
| (8.1.17) | $(b_1 + b_2)_*{}^*a = b_1{}_*{}^*a + b_2{}_*{}^*a$ |
| (10.2.8) | $f \circ (av) = a(f \circ v)$ |
| (10.3.11) | $\overline{f} \circ (v_*{}^*e_{V_1}) = v_*{}^*f_*{}^*e_{V_2}$ |

Из равенства (10.3.24) следует, что отображение $\overline{f}$ является гомоморфизмом абелевой группы. Равенство

$$(10.3.25) \quad \begin{aligned} &\overline{f} \circ ((av)_*{}^*e_{V_1}) \\ &= (av)_*{}^*f^*{}_*e_{V_2} \\ &= a(v_*{}^*f_*{}^*e_{V_2}) \\ &= a(\overline{f} \circ (v_*{}^*e_{V_1})) \end{aligned}$$

является следствием равенств

| | |
|---|---|
| (10.2.8) | $f \circ (av) = a(f \circ v)$ |
| (10.3.11) | $\overline{f} \circ (v_*{}^*e_{V_1}) = v_*{}^*f_*{}^*e_{V_2}$ |

Из равенства (10.3.25), и определений 2.8.5, 10.3.1 следует, что отображение

| | |
|---|---|
| (10.3.3) | $\overline{f} : V_1 \to V_2$ |

является гомоморфизмом $A$-векторного пространства строк $V_1$ в $A$-векторное пространство строк $V_2$.

Пусть $f$ - матрица гомоморфизмов $\overline{f}$, $\overline{g}$ относительно базисов $\overline{\overline{e}}_V$, $\overline{\overline{e}}_W$. Равенство

$$\overline{f} \circ (v_*{}^*e_V) = v_*{}^*f_*{}^*e_W = \overline{g} \circ (v_*{}^*e_V)$$

является следствием теоремы 10.3.7. Следовательно, $\overline{f} = \overline{g}$. □

На основании теорем 10.3.7, 10.3.9 мы идентифицируем гомоморфизм

| | |
|---|---|
| (10.3.3) | $\overline{f} : V_1 \to V_2$ |

левого $A$-векторного пространства $V$ строк и координаты его представления

| | |
|---|---|
| (10.3.11) | $\overline{f} \circ (v_*{}^*e_{V_1}) = v_*{}^*f_*{}^*e_{V_2}$ |

## 10.4. Сумма гомоморфизмов $A$-векторного пространства

Теорема 10.4.1. *Пусть $V$, $W$ - левые $A$-векторные пространства и отображения $g : V \to W$, $h : V \to W$ - гомоморфизмы левого $A$-векторного пространства. Пусть отображение $f : V \to W$ определено равенством*

$$(10.4.1) \qquad \forall v \in V : f \circ v = g \circ v + h \circ v$$

*Отображение $f$ является гомоморфизмом левого $A$-векторного пространства и называется* **суммой**

$$f = g + h$$

**гомоморфизмов** $g$ *and* $h$.

Доказательство. Согласно теореме 7.2.6 (равенство (7.2.8)), равенство

$$(10.4.2) \qquad \begin{aligned} & f \circ (u + v) \\ & = g \circ (u + v) + h \circ (u + v) \\ & = g \circ u + g \circ v \\ & + h \circ u + h \circ u \\ & = f \circ u + f \circ v \end{aligned}$$

является следствием равенства

$$\boxed{(10.3.4) \quad f \circ (u + v) = f \circ u + f \circ v}$$

и равенства (10.4.1). Согласно теореме 7.2.6 (равенство (7.2.10)), равенство

$$(10.4.3) \qquad \begin{aligned} & f \circ (av) \\ & = g \circ (av) + h \circ (av) \\ & = a(g \circ v) + a(h \circ v) \\ & = a(g \circ v + h \circ v) \\ & = a(f \circ v) \end{aligned}$$

является следствием равенства

$$\boxed{(10.3.5) \quad f \circ (av) = a(f \circ v)}$$

и равенства (10.4.1).

Теорема является следствием теоремы 10.3.2 и равенств (10.4.2), (10.4.3). $\square$

Теорема 10.4.2. *Пусть $V$, $W$ - левые $A$-векторные пространства столбцов. Пусть $\overline{\overline{e}}_V$ базис левого $A$-векторного пространства $V$. Пусть $\overline{\overline{e}}_W$ базис левого $A$-векторного пространства $W$. Гомоморфизм $\overline{f} : V \to W$ является суммой гомоморфизмов $\overline{g} : V \to W$, $\overline{h} : V \to W$ тогда и только тогда, когда матрица $f$ гомоморфизма $\overline{f}$ относительно базисов $\overline{\overline{e}}_V$, $\overline{\overline{e}}_W$ равна сумме матрицы $g$ гомоморфизма $\overline{g}$ относительно базисов*

Теорема 10.4.3. *Пусть $V$, $W$ - левые $A$-векторные пространства строк. Пусть $\overline{\overline{e}}_V$ базис левого $A$-векторного пространства $V$. Пусть $\overline{\overline{e}}_W$ базис левого $A$-векторного пространства $W$. Гомоморфизм $\overline{f} : V \to W$ является суммой гомоморфизмов $\overline{g} : V \to W$, $\overline{h} : V \to W$ тогда и только тогда, когда матрица $f$ гомоморфизма $\overline{f}$ относительно базисов $\overline{\overline{e}}_V$, $\overline{\overline{e}}_W$ равна сумме матрицы $g$ гомоморфизма $\overline{g}$ относительно базисов $\overline{\overline{e}}_V$, $\overline{\overline{e}}_W$* $\overline{\overline{e}}_V$, $\overline{\overline{e}}_W$ и матрицы $h$ гомоморфизма $\overline{h}$ относительно базисов $\overline{\overline{e}}_V$, $\overline{\overline{e}}_W$.

и матрицы $h$ гомоморфизма $\overline{h}$ относительно базисов $\overline{\overline{e}}_V$, $\overline{\overline{e}}_W$.

**Доказательство.** Согласно теореме [10.3.6](#)

$$(10.4.4) \qquad \overline{f} \circ (a^*{}_*e_U) = a^*{}_*f^*{}_*e_V$$

$$(10.4.5) \qquad \overline{g} \circ (a^*{}_*e_U) = a^*{}_*g^*{}_*e_V$$

$$(10.4.6) \qquad \overline{h} \circ (a^*{}_*e_U) = a^*{}_*h^*{}_*e_V$$

Согласно теореме [10.4.1](#), равенство

$$(10.4.7) \qquad \begin{aligned} &\overline{f} \circ (v^*{}_*e_U) \\ &= \overline{g} \circ (v^*{}_*e_U) + \overline{h} \circ (v^*{}_*e_U) \\ &= v^*{}_*g^*{}_*e_V + v^*{}_*h^*{}_*e_V \\ &= v^*{}_*(g+h)^*{}_*e_V \end{aligned}$$

является следствием равенств ([8.1.18](#)), ([8.1.19](#)), ([10.4.5](#)), ([10.4.6](#)). Равенство

$$f = g + h$$

является следствием равенств ([10.4.4](#)), ([10.4.7](#)) и теоремы [10.3.6](#).

Пусть
$$f = g + h$$

Согласно теореме [10.3.8](#), существуют гомоморфизмы $\overline{f} : U \to V$, $\overline{g} : U \to V$, $\overline{h} : U \to V$ такие, что

$$(10.4.8) \qquad \begin{aligned} &\overline{f} \circ (v^*{}_*e_U) \\ &= v^*{}_*f^*{}_*e_V \\ &= v^*{}_*(g+h)^*{}_*e_V \\ &= v^*{}_*g^*{}_*e_V + v^*{}_*h^*{}_*e_V \\ &= \overline{g} \circ (v^*{}_*e_U) + \overline{h} \circ (v^*{}_*e_U) \end{aligned}$$

Из равенства ([10.4.8](#)) и теоремы [10.4.1](#) следует, что гомоморфизм $\overline{f}$ является суммой гомоморфизмов $\overline{g}$, $\overline{h}$. $\square$

**Доказательство.** Согласно теореме [10.3.7](#)

$$(10.4.9) \qquad \overline{f} \circ (a_*{}^*e_U) = a_*{}^*f_*{}^*e_V$$

$$(10.4.10) \qquad \overline{g} \circ (a_*{}^*e_U) = a_*{}^*g_*{}^*e_V$$

$$(10.4.11) \qquad \overline{h} \circ (a_*{}^*e_U) = a_*{}^*h_*{}^*e_V$$

Согласно теореме [10.4.1](#), равенство

$$(10.4.12) \qquad \begin{aligned} &\overline{f} \circ (v_*{}^*e_U) \\ &= \overline{g} \circ (v_*{}^*e_U) + \overline{h} \circ (v_*{}^*e_U) \\ &= v_*{}^*g_*{}^*e_V + v_*{}^*h_*{}^*e_V \\ &= v_*{}^*(g+h)_*{}^*e_V \end{aligned}$$

является следствием равенств ([8.1.16](#)), ([8.1.17](#)), ([10.4.10](#)), ([10.4.11](#)). Равенство

$$f = g + h$$

является следствием равенств ([10.4.9](#)), ([10.4.12](#)) и теоремы [10.3.7](#).

Пусть
$$f = g + h$$

Согласно теореме [10.3.9](#), существуют гомоморфизмы $\overline{f} : U \to V$, $\overline{g} : U \to V$, $\overline{h} : U \to V$ такие, что

$$(10.4.13) \qquad \begin{aligned} &\overline{f} \circ (v_*{}^*e_U) \\ &= v_*{}^*f_*{}^*e_V \\ &= v_*{}^*(g+h)_*{}^*e_V \\ &= v_*{}^*g_*{}^*e_V + v_*{}^*h_*{}^*e_V \\ &= \overline{g} \circ (v_*{}^*e_U) + \overline{h} \circ (v_*{}^*e_U) \end{aligned}$$

Из равенства ([10.4.13](#)) и теоремы [10.4.1](#) следует, что гомоморфизм $\overline{f}$ является суммой гомоморфизмов $\overline{g}$, $\overline{h}$. $\square$

## 10.5. Произведение гомоморфизмов $A$-векторного пространства

**Теорема 10.5.1.** *Пусть $U$, $V$, $W$ - левые $A$-векторные пространства. Предположим, что мы имеем коммутативную диаграмму отображений*

$$(10.5.1)$$

$$(10.5.2) \qquad \forall u \in U : f \circ u = (h \circ g) \circ u$$
$$= h \circ (g \circ u)$$

где отображения $g$, $h$ являются гомоморфизмами левого $A$-векторного пространства. Отображение $f = h \circ g$ является гомоморфизмом левого $A$-векторного пространства и называется **произведением гомоморфизмов** $h$, $g$.

Доказательство. Равенство ($10.5.2$) следует из коммутативности диаграммы ($10.5.1$). Равенство

$$(10.5.3) \qquad \begin{aligned} &f \circ (u + v) \\ =\ & h \circ (g \circ (u + v)) \\ =\ & h \circ (g \circ u + g \circ v)) \\ =\ & h \circ (g \circ u) + h \circ (g \circ v) \\ =\ & f \circ u + f \circ v \end{aligned}$$

является следствием равенства

$$\boxed{(10.3.4) \quad f \circ (u + v) = f \circ u + f \circ v}$$

и равенства ($10.5.2$). Равенство

$$(10.5.4) \qquad \begin{aligned} f \circ (av) &= h \circ (g \circ (av)) \\ &= h \circ (a(g \circ v)) \\ &= a(h \circ (g \circ v)) \\ &= a(f \circ v) \end{aligned}$$

является следствием равенства

$$\boxed{(10.3.5) \quad f \circ (av) = a(f \circ v)}$$

и равенства ($10.4.1$). Теорема является следствием теоремы $10.3.2$ и равенств ($10.5.3$), ($10.5.4$). $\qquad\square$

Теорема 10.5.2. *Пусть $U$, $V$, $W$ - левые $A$-векторные пространства столбцов. Пусть $\overline{\overline{e}}_U$ базис левого $A$-векторного пространства $U$. Пусть $\overline{\overline{e}}_V$ базис левого $A$-векторного пространства $V$. Пусть $\overline{\overline{e}}_W$ базис левого $A$-векторного пространства $W$. Гомоморфизм $\overline{f} : U \to W$ является произведением гомоморфизмов $\overline{h} : V \to W$, $\overline{g} : U \to V$ тогда и только тогда, когда матрица $f$ гомоморфизма $\overline{f}$ относительно базисов $\overline{\overline{e}}_U$, $\overline{\overline{e}}_W$ равна $*_*$-произведению матрицы $g$ гомоморфизма $\overline{g}$ относительно базисов $\overline{\overline{e}}_U$, $\overline{\overline{e}}_V$ на матрицу $h$ гомоморфизма*

Теорема 10.5.3. *Пусть $U$, $V$, $W$ - левые $A$-векторные пространства строк. Пусть $\overline{\overline{e}}_U$ базис левого $A$-векторного пространства $U$. Пусть $\overline{\overline{e}}_V$ базис левого $A$-векторного пространства $V$. Пусть $\overline{\overline{e}}_W$ базис левого $A$-векторного пространства $W$. Гомоморфизм $\overline{f} : U \to W$ является произведением гомоморфизмов $\overline{h} : V \to W$, $\overline{g} : U \to V$ тогда и только тогда, когда матрица $f$ гомоморфизма $\overline{f}$ относительно базисов $\overline{\overline{e}}_U$, $\overline{\overline{e}}_W$ равна $*^*$-произведению матрицы $h$ гомоморфизма $\overline{h}$ относительно базисов $\overline{\overline{e}}_V$, $\overline{\overline{e}}_W$ на матрицу $g$ гомоморфизма $\overline{g}$*

$\overline{h}$ *относительно базисов* $\overline{\overline{e}}_V$, $\overline{\overline{e}}_W$

(10.5.5) $\qquad f = g^*{}_* h$

**Доказательство.** Согласно теореме 10.3.6

(10.5.7) $\qquad \overline{f} \circ (u^*{}_* e_U) = u^*{}_* f^*{}_* e_W$

(10.5.8) $\qquad \overline{g} \circ (u^*{}_* e_U) = u^*{}_* g^*{}_* e_V$

(10.5.9) $\qquad \overline{h} \circ (v^*{}_* e_V) = v^*{}_* h^*{}_* e_W$

Согласно теореме 10.5.1, равенство

$$
\begin{aligned}
& \overline{f} \circ (u^*{}_* e_U) \\
(10.5.10) \quad & = \overline{h} \circ (\overline{g} \circ (u^*{}_* e_U)) \\
& = \overline{h} \circ (u^*{}_* g^*{}_* e_V) \\
& = u^*{}_* g^*{}_* h^*{}_* e_W
\end{aligned}
$$

является следствием равенств (10.5.8), (10.5.9). Равенство (10.5.5) является следствием равенств (10.5.7), (10.5.10) и теоремы 10.3.6.

Пусть $f = g^*{}_* h$. Согласно теореме 10.3.8, существуют гомоморфизмы $\overline{f}: U \to V$, $\overline{g}: U \to V$, $\overline{h}: U \to V$ такие, что

$$
\begin{aligned}
& \overline{f} \circ (u^*{}_* e_U) = u^*{}_* f^*{}_* e_W \\
& = u^*{}_* g^*{}_* h^*{}_* e_W \\
(10.5.11) \quad & = \overline{h} \circ (u^*{}_* g^*{}_* e_V) \\
& = \overline{h} \circ (\overline{g} \circ (u^*{}_* e_U)) \\
& = (\overline{h} \circ \overline{g}) \circ (u^*{}_* e_U)
\end{aligned}
$$

Из равенства (10.5.11) и теоремы 10.4.1 следует, что гомоморфизм $\overline{f}$ является произведением гомоморфизмов $\overline{h}$, $\overline{g}$. $\qquad\square$

*относительно базисов* $\overline{\overline{e}}_U$, $\overline{\overline{e}}_V$

(10.5.6) $\qquad f = h_*{}^* g$

**Доказательство.** Согласно теореме 10.3.7

(10.5.12) $\qquad \overline{f} \circ (u_*{}^* e_U) = u_*{}^* f_*{}^* e_W$

(10.5.13) $\qquad \overline{g} \circ (u_*{}^* e_U) = u_*{}^* g_*{}^* e_V$

(10.5.14) $\qquad \overline{h} \circ (v_*{}^* e_V) = v_*{}^* h_*{}^* e_W$

Согласно теореме 10.5.1, равенство

$$
\begin{aligned}
& \overline{f} \circ (u_*{}^* e_U) \\
(10.5.15) \quad & = \overline{h} \circ (\overline{g} \circ (u_*{}^* e_U)) \\
& = \overline{h} \circ (u_*{}^* g_*{}^* e_V) \\
& = u_*{}^* g_*{}^* h_*{}^* e_W
\end{aligned}
$$

является следствием равенств (10.5.13), (10.5.14). Равенство (10.5.6) является следствием равенств (10.5.12), (10.5.15) и теоремы 10.3.7.

Пусть $f = h_*{}^* g$. Согласно теореме 10.3.9, существуют гомоморфизмы $\overline{f}: U \to V$, $\overline{g}: U \to V$, $\overline{h}: U \to V$ такие, что

$$
\begin{aligned}
& \overline{f} \circ (u_*{}^* e_U) = u_*{}^* f_*{}^* e_W \\
& = u_*{}^* g_*{}^* h_*{}^* e_W \\
(10.5.16) \quad & = \overline{h} \circ (u_*{}^* g_*{}^* e_V) \\
& = \overline{h} \circ (\overline{g} \circ (u_*{}^* e_U)) \\
& = (\overline{h} \circ \overline{g}) \circ (u_*{}^* e_U)
\end{aligned}
$$

Из равенства (10.5.16) и теоремы 10.4.1 следует, что гомоморфизм $\overline{f}$ является произведением гомоморфизмов $\overline{h}$, $\overline{g}$. $\qquad\square$

## 10.6. Прямая сумма левых $A$-векторных пространств

Пусть Left $A$-Vector - категория левых $A$-векторных пространств, морфизмы которой - это гомоморфизмы $A$-модуля.

**Определение 10.6.1.** *Копроизведение в категории* Left $A$-Vector *называется* **прямой суммой**. [10.10] *Мы будем пользоваться записью* $V \oplus W$ *для прямой суммы левых $A$-векторных пространств* $V$ *и* $W$. $\qquad\square$

---

[10.10] Смотри так же определение прямой суммы модулей в [1], страница 98. Согласно теореме 1 на той же странице, прямая сумма модулей существует.

Теорема 10.6.2. *В категории* Left $A$-Vector *существует прямая сумма левых $A$-векторных пространств.*[10.11]*Пусть* $\{V_i, i \in I\}$ *- семейство левых $A$-векторных пространств. Тогда представление*

$$A \xrightarrow{\ *\ } \bigoplus_{i \in I} V_i \qquad a\left(\bigoplus_{i \in I} v_i\right) = \bigoplus_{i \in I} av_i$$

*$D$-алгебры $A$ в прямой сумме абелевых групп*

$$V = \bigoplus_{i \in I} V_i$$

*является прямой суммой левых $A$-векторных пространств*

$$V = \bigoplus_{i \in I} V_i$$

Доказательство. Пусть $V_i, i \in I$, - множество левых $A$-векторных пространств. Пусть

$$V = \bigoplus_{i \in I} V_i$$

прямая сумма абелевых групп. Пусть $v = (v_i, i \in I) \in V$. Мы определим действие $A$-числа $a$ на $V$-число $v$ согласно равенству

$$(10.6.1) \qquad a(v_i, i \in I) = (av_i, i \in I)$$

Согласно теореме 7.2.6,

$$(10.6.2) \qquad (ab)v_i = a(bv_i)$$

$$(10.6.3) \qquad a(v_i + w_i) = av_i + aw_i$$

$$(10.6.4) \qquad (a + b)v_i = av_i + bv_i$$

$$(10.6.5) \qquad 1v_i = v_i$$

для любых $a, b \in A$, $v_i, w_i \in V_i$, $i \in I$. Равенства

$$(10.6.6) \quad \begin{aligned} (ab)v &= (ab)(v_i, i \in I) = ((ab)v_i, i \in I) \\ &= (a(bv_i), i \in I) = a(bv_i, i \in I) = a(b(v_i, i \in I)) \\ &= a(bv) \end{aligned}$$

$$(10.6.7) \quad \begin{aligned} a(v + w) &= a((v_i, i \in I) + (w_i, i \in I)) = a(v_i + w_i, i \in I) \\ &= (a(v_i + w_i), i \in I) = (av_i + aw_i, i \in I) \\ &= (av_i, i \in I) + (aw_i, i \in I) = a(v_i, i \in I) + a(w_i, i \in I) \\ &= av + aw \end{aligned}$$

$$(10.6.8) \quad \begin{aligned} (a + b)v &= (a + b)(v_i, i \in I) = ((a + b)v_i, i \in I) \\ &= (av_i + bv_i, i \in I) = (av_i, i \in I) + (bv_i, i \in I) \\ &= a(v_i, i \in I) + b(v_i, i \in I) = av + bv \end{aligned}$$

$$1v = 1(v_i, i \in I) = (1v_i, i \in I) = (v_i, i \in I) = v$$

---

[10.11]  Смотри также предложение [1]-1 на странице 98.

являются следствием равенств $(3.7.4)$, $(10.6.1)$, $(10.6.2)$, $(10.6.3)$, $(10.6.4)$, $(10.6.5)$
для любых $a$, $b \in A$, $v = (v_i, i \in I)$, $w = (w_i, i \in I) \in V$.

Из равенства $(10.6.7)$ следует, что отображение

$$a : v \to av$$

определённое равенством $(10.6.1)$, является эндоморфизмом абелевой группы
$V$. Из равенств $(10.6.6)$, $(10.6.8)$ следует, что отображение $(10.6.1)$ является
гомоморфизмом $D$-алгебры $A$ в множество эндоморфизмов абелевой группы
$V$. Следовательно, отображение $(10.6.1)$ является представлением $D$-алгебры
$A$ в абелевой группе $V$.

Пусть $a_1$, $a_2$ - такие $A$-числа, что

$$(10.6.9) \qquad a_1(v_i, i \in I) = a_2(v_i, i \in I)$$

для произвольного $V$-числа $v = (v_i, i \in I)$. Равенство

$$(10.6.10) \qquad (a_1 v_i, i \in I) = (a_2 v_i, i \in I)$$

является следствием равенства $(10.6.9)$. Из равенства $(10.6.10)$ следует, что

$$(10.6.11) \qquad \forall i \in I : a_1 v_i = a_2 v_i$$

Поскольку для любого $i \in I$, $V_i$ - $A$-векторное пространство, то соответствую-
щее представление эффективно согласно определению $7.2.1$  Следовательно,
равенство $a_1 = a_2$ является следствием утверждения $(10.6.11)$ и отображе-
ние $(10.6.1)$ является эффективным представлением $D$-алгебры $A$ в абелевой
группе $V$. Следовательно, абелевая группе $V$ является левым $A$-векторным
пространством.

Пусть

$$\{f_i : V_i \to W, i \in I\}$$

семейство линейных отображений в  левое $A$-векторное пространство $W$. Мы
определим отображение

$$f : V \to W$$

пользуясь равенством

$$(10.6.12) \qquad f\left(\bigoplus_{i \in I} x_i\right) = \sum_{i \in I} f_i(x_i)$$

Сумма в правой части равенства $(10.6.12)$ конечна, так как все слагаемые,
кроме конечного числа, равны 0. Из равенства

$$f_i(x_i + y_i) = f_i(x_i) + f_i(y_i)$$

и равенства $(10.6.12)$ следует, что

$$\begin{aligned}
f(\bigoplus_{i \in I}(x_i + y_i)) &= \sum_{i \in I} f_i(x_i + y_i) \\
&= \sum_{i \in I} (f_i(x_i) + f_i(y_i)) \\
&= \sum_{i \in I} f_i(x_i) + \sum_{i \in I} f_i(y_i) \\
&= f(\bigoplus_{i \in I} x_i) + f(\bigoplus_{i \in I} y_i)
\end{aligned}$$

Из равенства

$$f_i(dx_i) = df_i(x_i)$$

и равенства (10.6.12) следует, что

$$f((dx_i, i \in I)) = \sum_{i \in I} f_i(dx_i) = \sum_{i \in I} df_i(x_i) = d \sum_{i \in I} f_i(x_i)$$
$$= df((x_i, i \in I))$$

Следовательно, отображение $f$ является линейным отображением. Равенство

$$f \circ \lambda_j(x) = \sum_{i \in I} f_i(\delta_j^i x) = f_j(x)$$

является следствием равенств (3.7.5), (10.6.12). Так как отображение $\lambda_i$ инъективно, то отображение $f$ определено однозначно. Следовательно, теорема является следствием определений 2.2.12, 3.7.2.    □

**Определение 10.6.3.** *Пусть $D$ - коммутативное кольцо с единицей. Пусть $A$ - ассоциативная $D$-алгебра с делением. Для любого целого $n > 0$ мы определим левое $A$-векторное пространство $\mathcal{L}^n A$ используя определение по индукции*

$$\mathcal{L}^1 A = \mathcal{L}A$$
$$\mathcal{L}^k A = \mathcal{L}^{k-1} A \oplus \mathcal{L}A$$

□

**Определение 10.6.4.** *Пусть $D$ - коммутативное кольцо с единицей. Пусть $A$ - ассоциативная $D$-алгебра с делением. Для любого множества $B$ и любого целого $n > 0$ мы определим левое $A$-векторное пространство $\mathcal{L}^n A^B$ используя определение по индукции*

$$\mathcal{L}^1 A^B = \mathcal{L}A^B$$
$$\mathcal{L}^k A^B = \mathcal{L}^{k-1} A^B \oplus \mathcal{L}A^B$$

□

**Теорема 10.6.5.** *Прямая сумма левых $A$-векторных пространств $V_1$, ..., $V_n$ совпадает с их прямым произведением*

$$V_1 \oplus ... \oplus V_n = V_1 \times ... \times V_n$$

**Доказательство.** Теорема является следствием теоремы 10.6.2.    □

**Теорема 10.6.6.** *Пусть $\overline{\overline{e}}_1 = (e_1{}_i, i \in I)$ базис левого $A$-векторного пространства столбцов $V$. Тогда мы можем представить $A$-модуль $V$ как прямую сумму*

$$(10.6.13) \qquad V = \bigoplus_{i \in I} A e_i$$

**Теорема 10.6.7.** *Пусть $\overline{\overline{e}}_1 = (e_1^i, i \in I)$ базис левого $A$-векторного пространства строк $V$. Тогда мы можем представить $A$-модуль $V$ как прямую сумму*

$$(10.6.14) \qquad V = \bigoplus_{i \in I} A e^i$$

Доказательство. Отображение

$$f : v^i e_i \in V \to (v^i, i \in I) \in \bigoplus_{i \in I} A e_i$$

является изоморфизмом.          $\square$

Доказательство. Отображение

$$f : v_i e^i \in V \to (v_i, i \in I) \in \bigoplus_{i \in I} A e^i$$

является изоморфизмом.          $\square$

# Гомоморфизм правого $A$-векторного пространства

## 11.1. Общее определение

Пусть $A_i$, $i = 1, 2$, - алгебра с делением над коммутативным кольцом $D_i$. Пусть $V_i$, $i = 1, 2$, - правое $A_i$-модуль.

ОПРЕДЕЛЕНИЕ 11.1.1. *Пусть*

$$(11.1.1) \quad \begin{array}{c} A_1 \xrightarrow{g_{1.23}} A_1 \xrightarrow{g_{1.34}} V_1 \\ \scriptstyle{g_{1.12}} \quad \scriptstyle{g_{1.12}} \quad \scriptstyle{g_{1.14}} \\ D_1 \end{array} \qquad \begin{array}{l} g_{1.12}(d) : a \to d\,a \\ g_{1.23}(v) : w \to C_1(w, v) \\ \quad C_1 \in \mathcal{L}(A_1^2 \to A_1) \\ g_{1.34}(a) : v \to v\,a \\ g_{1.14}(d) : v \to v\,d \end{array}$$

*диаграмма представлений, описывающая правое $A_1$-векторное пространство $V_1$. Пусть*

$$(11.1.2) \quad \begin{array}{c} A_2 \xrightarrow{g_{2.23}} A_2 \xrightarrow{g_{2.34}} V_2 \\ \scriptstyle{g_{2.12}} \quad \scriptstyle{g_{2.12}} \quad \scriptstyle{g_{2.14}} \\ D_2 \end{array} \qquad \begin{array}{l} g_{2.12}(d) : a \to d\,a \\ g_{2.23}(v) : w \to C_2(w, v) \\ \quad C_2 \in \mathcal{L}(A_2^2 \to A_2) \\ g_{2.34}(a) : v \to v\,a \\ g_{2.14}(d) : v \to v\,d \end{array}$$

*диаграмма представлений, описывающая правое $A_2$-векторное пространство $V_2$. Морфизм*

$$(11.1.3) \qquad (h : D_1 \to D_2 \quad g : A_1 \to A_2 \quad f : V_1 \to V_2)$$

*диаграммы представлений* (11.1.1) *в диаграмму представлений* (11.1.2) *называется* **гомоморфизмом** *правого $A_1$-векторного пространства $V_1$ в правое $A_2$-векторное пространство $V_2$. Обозначим* $\mathrm{Hom}(D_1 \to D_2; A_1 \to A_2; V_1* \to V_2*)$ *множество гомоморфизмов правого $A_1$-векторного пространства $V_1$ в правое $A_2$-векторное пространство $V_2$.* □

Мы будем пользоваться записью

$$f \circ a = f(a)$$

для образа гомоморфизма $f$.

ТЕОРЕМА 11.1.2. *Гомоморфизм*

$$(11.1.3) \quad (h : D_1 \to D_2 \quad g : A_1 \to A_2 \quad f : V_1 \to V_2)$$

*правого $A_1$-векторного пространства $V_1$ в правое $A_2$-векторное пространство*

*$V_2$ удовлетворяет следующим равенствам*

$$(11.1.4) \qquad h(p + q) = h(p) + h(q)$$

$$(11.1.5) \qquad h(pq) = h(p)h(q)$$

$$(11.1.6) \qquad g \circ (a + b) = g \circ a + g \circ b$$

$$(11.1.7) \qquad g \circ (ab) = (g \circ a)(g \circ b)$$

$$(11.1.8) \qquad g \circ (ap) = (g \circ a)h(p)$$

$$(11.1.9) \qquad f \circ (u + v) = f \circ u + f \circ v$$

$$(11.1.10) \qquad f \circ (va) = (f \circ v)(g \circ a)$$

$$p, q \in D_1 \qquad a, b \in A_1 \qquad u, v \in V_1$$

Доказательство. Согласно определениям 2.3.9, 2.8.5, 7.4.1, 11 отображение

$$(h : D_1 \to D_2 \quad g : A_1 \to A_2)$$

является гомоморфизмом $D_1$-алгебры $A_1$ в $D_2$-алгебру $A_2$. Следовательно, равенства (11.1.4), (11.1.5), (11.1.6), (11.1.7), (11.1.8) являются следствием теоремы 4.3.2. Равенство (11.1.9) является следствием определения 11.1.1, так как, согласно определению 2.3.9, отображение $f$ является гомоморфизмом абелевой группы. Равенство (11.1.10) является следствием равенства (2.3.6), так как отображение

$$(g : A_1 \to A_2 \quad f : V_1 \to V_2)$$

является морфизмом представления $g_{1.34}$ в представление $g_{2.34}$. □

Определение 11.1.3. *Гомоморфизм* [11.1]

$$(h : D_1 \to D_2 \quad g : A_1 \to A_2 \quad f : V_1 \to V_2)$$

*называется* **изоморфизмом** *между правым $A_1$-векторным пространством $V_1$ и правым $A_2$-векторным пространством $V_2$, если существует отображение*

$$(h^{-1} : D_2 \to D_1 \quad g^{-1} : A_2 \to A_1 \quad f^{-1} : V_2 \to V_1)$$

*которое является гомоморфизмом.* □

### 11.1.1. Гомоморфизм правого $A$-векторного пространства столбцов.

Теорема 11.1.4. *Гомоморфизм* [11.2]

$$(11.1.3) \qquad (h : D_1 \to D_2 \quad \overline{g} : A_1 \to A_2 \quad \overline{f} : V_1 \to V_2)$$

---

[11.1] Я следую определению на странице [18]-63.

[11.2] В теоремах 11.1.4, 11.1.5, мы опираемся на следующее соглашение. Пусть $i = 1$, 2. Пусть множество векторов $\overline{\overline{e}}_{A_i} = (e_{A_i k}, k \in K_i)$ является базисом и $C_{ijk}^{~k}$, $k$, $i$, $j \in K_i$, - структурные константы $D_i$-алгебры столбцов $A_i$. Пусть множество векторов $\overline{\overline{e}}_{V_i} = (e_{V_i i}, i \in I_i)$ является базисом правого $A_i$-векторного пространства $V_i$.

*правого $A_1$-векторного пространства столбцов $V_1$ в правое $A_2$-векторное пространство столбцов $V_2$ имеет представление*

$$(11.1.11) \qquad b = gh(a)$$

$$(11.1.12) \qquad \overline{g} \circ (a^i e_{A_1 i}) = h(a^i)g_i^k e_{A_2 k}$$

$$(11.1.13) \qquad \overline{g} \circ (e_{A_1} a) = e_{A_2} gh(a)$$

$$(11.1.14) \qquad w = f_*{}^* (\overline{g} \circ v)$$

$$(11.1.15) \qquad \overline{f} \circ (e_{V_1 i} v^i) = e_{V_2 k} f_i^k (\overline{g} \circ v^i)$$

$$(11.1.16) \qquad \overline{f} \circ (e_{V_1 *}{}^* v) = e_{V_2 *}{}^* f_*{}^* (\overline{g} \circ v)$$

*относительно выбранных базисов. Здесь*

- *$a$ - координатная матрица $A_1$-числа $\overline{a}$ относительно базиса $\overline{\overline{e}}_{A_1}$ (A)*

$$\overline{a} = e_{A_1} a$$

- *$h(a) = (h(a_k), k \in K)$  - матрица $D_2$-чисел.*
- *$b$ - координатная матрица $A_2$-числа*

$$\overline{b} = \overline{g} \circ \overline{a}$$

  *относительно базиса $\overline{\overline{e}}_{A_2}$*

$$\overline{b} = e_{A_2} b$$

- *$g$ - координатная матрица множества $A_2$-чисел $(\overline{g} \circ e_{A_1 k}, k \in K)$ относительно базиса $\overline{\overline{e}}_{A_2}$.*
- *Матрица гомоморфизма и структурные константы связаны соотношением*

$$(11.1.17) \qquad h(C_{1\,ij}^{\phantom{1ij}k})g_k^l = g_i^p g_j^q C_{2\,pq}^{\phantom{2pq}l}$$

- *$v$ - координатная матрица $V_1$-числа $\overline{v}$ относительно базиса $\overline{\overline{e}}_{V_1}$ (A)*

$$\overline{v} = e_{V_1 *}{}^* v$$

- *$g(v) = (g(v_i), i \in I)$  - матрица $A_2$-чисел.*
- *$w$ - координатная матрица $V_2$-числа*

$$\overline{w} = \overline{f} \circ \overline{v}$$

  *относительно базиса $\overline{\overline{e}}_{V_2}$*

$$\overline{w} = e_{V_2 *}{}^* w$$

- *$f$ - координатная матрица множества $V_2$-чисел $(\overline{f} \circ e_{V_1 i}, i \in I)$ относительно базиса $\overline{\overline{e}}_{V_2}$.*

*Для заданного гомоморфизма*

$$\boxed{(11.1.3) \quad (h : D_1 \to D_2 \quad \overline{g} : A_1 \to A_2 \quad \overline{f} : V_1 \to V_2)}$$

*множество матриц $(g, f)$ определено однозначно и называется **координатами гомоморфизма** (11.1.3)  относительно базисов $(\overline{\overline{e}}_{A_1}, \overline{\overline{e}}_{V_1})$, $(\overline{\overline{e}}_{A_2}, \overline{\overline{e}}_{V_2})$.*

Доказательство.　Согласно определениям　2.3.9, 2.8.5, 5.4.1, 11.1.1, отображение

$$(h : D_1 \to D_2 \quad \overline{g} : A_1 \to A_2)$$

является гомоморфизмом $D_1$-алгебры $A_1$ в $D_2$-алгебру $A_2$. Следовательно, равенства (11.1.11), (11.1.12), (11.1.13), (11.1.17) следуют из равенств

$(5.4.11)$   $b = fh(a)$

$(5.4.12)$   $\overline{f} \circ (a^i e_{1i}) = h(a^i) f_i^k e_{2k}$

$(5.4.13)$   $\overline{f} \circ (e_1 a) = e_2 fh(a)$

$(5.4.16)$   $h(C_{1\,ij}^{\ \ k}) f_k^l = f_i^p f_j^q C_{2\,pq}^{\ \ l}$

Из теоремы 4.3.3, следует что матрица $g$ определена однозначно.

Вектор $\overline{v} \in V_1$ имеет разложение

$(11.1.18)$ $$\overline{v} = e_{V_1 *}{}^* v$$

относительно базиса $\overline{\overline{e}}_{V_1}$. Вектор $\overline{w} \in V_2$ имеет разложение

$(11.1.19)$ $$\overline{w} = e_{V_2 *}{}^* w$$

относительно базиса $\overline{\overline{e}}_{V_2}$. Так как $\overline{f}$ - гомоморфизм, то равенство

$(11.1.20)$ $$\overline{w} = \overline{f} \circ \overline{v} = \overline{f} \circ (e_{V_1 *}{}^* v)$$
$$= (\overline{f} \circ e_{V_1})_*{}^* (\overline{g} \circ v)$$

является следствием равенств (11.1.9), (11.1.10). $V_2$-число $\overline{f} \circ e_{V_1 i}$ имеет разложение

$(11.1.21)$ $$\overline{f} \circ e_{V_1 i} = e_{V_2 *}{}^* f_i = e_{V_2}^j f_i^j$$

относительно базиса $\overline{\overline{e}}_{V_2}$. Равенство

$(11.1.22)$ $$\overline{w} = e_{V_2 *}{}^* f_*{}^* (\overline{g} \circ v)$$

является следствием равенств (11.1.20), (11.1.21). Равенство (11.1.14) следует из сравнения (11.1.19) и (11.1.22) и теоремы 8.4.18. Из равенства (11.1.21) и теоремы 8.4.18 следует что матрица $f$ определена однозначно. □

ТЕОРЕМА 11.1.5. *Пусть отображение*
$$h : D_1 \to D_2$$
*является гомоморфизмом кольца $D_1$ в кольцо $D_2$.*     *Пусть*
$$g = (g_l^k, k \in K, l \in L)$$
*матрица $D_2$-чисел, которая удовлетворяют равенству*

$(11.1.17)$   $h(C_{1\,ij}^{\ \ k}) g_k^l = g_i^p g_j^q C_{2\,pq}^{\ \ l}$

*Пусть*
$$f = (f_j^i, i \in I, j \in J)$$
*матрица $A_2$-чисел.*     *Отображение* [11.2]

$(11.1.3)$   $(h : D_1 \to D_2 \quad \overline{g} : A_1 \to A_2 \quad \overline{f} : V_1 \to V_2)$

*определённое равенствами*

$(11.1.13)$   $\overline{g} \circ (e_{A_1} a) = e_{A_2} gh(a)$

$(11.1.16)$   $\overline{f} \circ (e_{V_1 *}{}^* v) = e_{V_2 *}{}^* f_*{}^* (\overline{g} \circ v)$

*является гомоморфизмом правого $A_1$-векторного пространства столбцов $V_1$ в правое $A_2$-векторное пространство столбцов $V_2$.*  *Гомоморфизм*

$$\boxed{(11.1.3) \quad (h:D_1 \to D_2 \quad \overline{g}:A_1 \to A_2 \quad \overline{f}:V_1 \to V_2}$$
*который имеет данное множество матриц $(g,f)$, определён однозначно.*

Доказательство. Согласно теореме 5.4.5, отображение

$$(11.1.23) \qquad (h:D_1 \to D_2 \quad \overline{g}:A_1 \to A_2)$$

является гомоморфизмом $D_1$-алгебры $A_1$ в $D_2$-алгебру $A_2$ и гомоморфизм (11.1.23) определён однозначно.

Равенство

$$(11.1.24) \quad
\begin{aligned}
&\overline{f} \circ (e_{V_1*}{}^*(v+w)) \\
&= e_{V_2*}{}^* f_*{}^*(g \circ (v+w)) \\
&= e_{V_2*}{}^* f_*{}^*(g \circ v) \\
&+ e_{V_2*}{}^* f_*{}^*(g \circ w) \\
&= \overline{f} \circ (e_{V_1*}{}^*v) + \overline{f} \circ (e_{V_1*}{}^*w)
\end{aligned}$$

является следствием равенств

$$\boxed{(8.1.16) \quad a_*{}^*(b_1+b_2) = a_*{}^*b_1 + a_*{}^*b_2} \qquad \boxed{(11.1.6) \quad g \circ (a+b) = g \circ a + g \circ b}$$

$$\boxed{(11.1.10) \quad f \circ (va) = (f \circ v)(g \circ a)} \qquad \boxed{(11.1.16) \quad \overline{f} \circ (e_{V_1*}{}^*v) = e_{V_2*}{}^* f_*{}^*(\overline{g} \circ v)}$$

Из равенства (11.1.24) следует, что отображение $\overline{f}$ является гомоморфизмом абелевой группы. Равенство

$$(11.1.25) \quad
\begin{aligned}
&\overline{f} \circ (e_{V_1*}{}^*(va)) \\
&= e_{V_2*}{}^* f_*{}^*(\overline{g} \circ (va)) \\
&= (e_{V_2*}{}^* f_*{}^*(\overline{g} \circ v))(\overline{g} \circ a) \\
&= (\overline{f} \circ (e_{V_1*}{}^*v))(\overline{g} \circ a)
\end{aligned}$$

является следствием равенств

$$\boxed{(11.1.10) \quad f \circ (va) = (f \circ v)(g \circ a)}$$

$$\boxed{(5.4.4) \quad g \circ (ab) = (g \circ a)(g \circ b)}$$

$$\boxed{(11.1.16) \quad \overline{f} \circ (e_{V_1*}{}^*v) = e_{V_2*}{}^* f_*{}^*(\overline{g} \circ v)}$$

Из равенства (11.1.25), и определений 2.8.5, 11.1.1 следует, что отображение

$$\boxed{(11.1.3) \quad (h:D_1 \to D_2 \quad \overline{g}:A_1 \to A_2 \quad \overline{f}:V_1 \to V_2)}$$

является гомоморфизмом $A_1$-векторного пространства столбцов $V_1$ в $A_2$-векторное пространство столбцов $V_2$.

Пусть $f$ - матрица гомоморфизмов $\overline{f}$, $\overline{g}$ относительно базисов $\overline{\overline{e}}_V$, $\overline{\overline{e}}_W$. Равенство

$$\overline{f} \circ (e_{V*}{}^*v) = e_{W*}{}^* f_*{}^*v = \overline{g} \circ (e_{V*}{}^*v)$$

является следствием теоремы 11.3.6. Следовательно, $\overline{f} = \overline{g}$. $\qquad\square$

На основании теорем 11.1.4, 11.1.5 мы идентифицируем гомоморфизм

$$\boxed{(11.1.3) \quad (h:D_1 \to D_2 \quad \overline{g}:A_1 \to A_2 \quad \overline{f}:V_1 \to V_2)}$$

правого $A$-векторного пространства $V$ столбцов и координаты его представления

$$\boxed{(11.1.13) \quad \overline{g} \circ (e_{A_1}a) = e_{A_2}gh(a)}$$

$$\boxed{(11.1.16) \quad \overline{f} \circ (e_{V_1*}v) = e_{V_2*}f_*(\overline{g} \circ v)}$$

## 11.1.2. Гомоморфизм правого $A$-векторного пространства строк.

Теорема 11.1.6. *Гомоморфизм* [11.3]

$$\boxed{(11.1.3) \quad (h : D_1 \to D_2 \quad \overline{g} : A_1 \to A_2 \quad \overline{f} : V_1 \to V_2)}$$

*правого $A_1$-векторного пространства строк $V_1$ в правое $A_2$-векторное пространство строк $V_2$ имеет представление*

$$(11.1.26) \qquad b = h(a)g$$

$$(11.1.27) \qquad \overline{g} \circ (a_i e_{A_1}^i) = h(a_i)g_k^i e_{A_2}^k$$

$$(11.1.28) \qquad \overline{g} \circ (a e_{A_1}) = h(a)g e_{A_2}$$

$$(11.1.29) \qquad w = f_*(\overline{g} \circ v)$$

$$(11.1.30) \qquad \overline{f} \circ (e_{V_1}^i v_i) = e_{V_2}^k f_k^i(\overline{g} \circ v_i)$$

$$(11.1.31) \qquad \overline{f} \circ (e_{V_1*}v) = e_{V_2*}f_*(\overline{g} \circ v)$$

*относительно выбранных базисов. Здесь*

- $a$ - *координатная матрица $A_1$-числа $\overline{a}$ относительно базиса $\overline{\overline{e}}_{A_1}$ (A)*

$$\overline{a} = a e_{A_1}$$

- $h(a) = (h(a^k), k \in K)$ - *матрица $D_2$-чисел.*
- $b$ - *координатная матрица $A_2$-числа*

$$\overline{b} = \overline{g} \circ \overline{a}$$

*относительно базиса $\overline{\overline{e}}_{A_2}$*

$$\overline{b} = b e_{A_2}$$

- $g$ - *координатная матрица множества $A_2$-чисел $(\overline{g} \circ e_{A_1}^k, k \in K)$ относительно базиса $\overline{\overline{e}}_{A_2}$.*
- *Матрица гомоморфизма и структурные константы связаны соотношением*

$$(11.1.32) \qquad h(C_{1k}^{ij})g_l^k = g_p^i g_q^j C_{2l}^{pq}$$

- $v$ - *координатная матрица $V_1$-числа $\overline{v}$ относительно базиса $\overline{\overline{e}}_{V_1}$ (A)*

$$\overline{v} = e_{V_1*}v$$

- $g(v) = (g(v^i), i \in I)$ - *матрица $A_2$-чисел.*

---

[11.3] В теоремах [11.1.6], [11.1.7], мы опираемся на следующее соглашение. Пусть $i = 1, 2$. Пусть множество векторов $\overline{\overline{e}}_{A_i} = (e_{A_i}^k, k \in K_i)$ является базисом и $C_{ik}^{ij}$, $k, i, j \in K_i$, - структурные константы $D_i$-алгебры строк $A_i$. Пусть множество векторов $\overline{\overline{e}}_{V_i} = (e_{V_i}^i, i \in I_i)$ является базисом правого $A_i$-векторного пространства $V_i$.

- $w$ - координатная матрица $V_2$-числа

$$\overline{w} = \overline{f} \circ \overline{v}$$

относительно базиса $\overline{\overline{e}}_{V_2}$

$$\overline{w} = e_{V_2}{}^*{}_*w$$

- $f$ - координатная матрица множества $V_2$-чисел $\left(\overline{f} \circ e^i_{V_1}, i \in I\right)$ относительно базиса $\overline{\overline{e}}_{V_2}$.

*Для заданного гомоморфизма*

$$\boxed{(11.1.3) \quad (h : D_1 \to D_2 \quad \overline{g} : A_1 \to A_2 \quad \overline{f} : V_1 \to V_2)}$$

*множество матриц* $(g, f)$ *определено однозначно и называется* **координатами гомоморфизма** (11.1.3) *относительно базисов* $(\overline{\overline{e}}_{A_1}, \overline{\overline{e}}_{V_1})$, $(\overline{\overline{e}}_{A_2}, \overline{\overline{e}}_{V_2})$.

Доказательство. Согласно определениям 2.3.9, 2.8.5, 5.4.1, 11.1.1, отображение

$$(h : D_1 \to D_2 \quad \overline{g} : A_1 \to A_2)$$

является гомоморфизмом $D_1$-алгебры $A_1$ в $D_2$-алгебру $A_2$. Следовательно, равенства (11.1.26), (11.1.27), (11.1.28), (11.1.32) следуют из равенств

$$\boxed{(5.4.22) \quad b = h(a)f}$$

$$\boxed{(5.4.23) \quad \overline{f} \circ (a_i e^i_1) = h(a_i)f^i_k e^k_2}$$

$$\boxed{(5.4.24) \quad \overline{f} \circ (a e_1) = h(a)f e_2}$$

$$\boxed{(5.4.27) \quad h(C^{ij}_{1k})f^k_l = f^i_p f^j_q C^{pq}_{2l}}$$

Из теоремы 4.3.4, следует что матрица $g$ определена однозначно.

Вектор $\overline{v} \in V_1$ имеет разложение

$$(11.1.33) \qquad \overline{v} = e_{V_1}{}^*{}_*v$$

относительно базиса $\overline{\overline{e}}_{V_1}$. Вектор $\overline{w} \in V_2$ имеет разложение

$$(11.1.34) \qquad \overline{w} = e_{V_2}{}^*{}_*w$$

относительно базиса $\overline{\overline{e}}_{V_2}$. Так как $\overline{f}$ - гомоморфизм, то равенство

$$(11.1.35) \qquad \overline{w} = \overline{f} \circ \overline{v} = \overline{f} \circ (e_{V_1}{}^*{}_*v)$$
$$= (\overline{f} \circ e_{V_1})^*{}_*(\overline{g} \circ v)$$

является следствием равенств (11.1.9), (11.1.10). $V_2$-число $\overline{f} \circ e^i_{V_1}$ имеет разложение

$$(11.1.36) \qquad \overline{f} \circ e^i_{V_1} = e_{V_2}{}^*{}_*f^i = e_{V_2 j}f^i_j$$

относительно базиса $\overline{\overline{e}}_{V_2}$. Равенство

$$(11.1.37) \qquad \overline{w} = e_{V_2}{}^*{}_*f^*{}_*(\overline{g} \circ v)$$

является следствием равенств (11.1.35), (11.1.36). Равенство (11.1.29) следует из сравнения (11.1.34) и (11.1.37) и теоремы 8.4.25. Из равенства (11.1.36) и теоремы 8.4.25 следует что матрица $f$ определена однозначно. $\square$

Теорема 11.1.7. *Пусть отображение*

$$h : D_1 \to D_2$$

*является гомоморфизмом кольца $D_1$ в кольцо $D_2$. Пусть*

$$g = (g_k^l, k \in K, l \in L)$$

*матрица $D_2$-чисел, которая удовлетворяют равенству*

$$(11.1.32) \quad h(C_{1k}^{ij})g_l^k = g_p^i g_q^j C_{2l}^{pq}$$

*Пусть*

$$f = (f_i^j, i \in I, j \in J)$$

*матрица $A_2$-чисел. Отображение*[11.3]

$$(11.1.3) \quad (h : D_1 \to D_2 \quad \overline{g} : A_1 \to A_2 \quad \overline{f} : V_1 \to V_2)$$

*определённое равенствами*

$$(11.1.28) \quad \overline{g} \circ (ae_{A_1}) = h(a)ge_{A_2}$$

$$(11.1.31) \quad \overline{f} \circ (e_{V_1}{}^*{}_*v) = e_{V_2}{}^*{}_*f^*{}_*(\overline{g} \circ v)$$

*является гомоморфизмом правого $A_1$-векторного пространства строк $V_1$ в правое $A_2$-векторное пространство строк $V_2$. Гомоморфизм*

$$(11.1.3) \quad (h : D_1 \to D_2 \quad \overline{g} : A_1 \to A_2 \quad \overline{f} : V_1 \to V_2)$$

*который имеет данное множество матриц $(g, f)$, определён однозначно.*

Доказательство. Согласно теореме 5.4.6, отображение

$$(11.1.38) \quad (h : D_1 \to D_2 \quad \overline{g} : A_1 \to A_2)$$

является гомоморфизмом $D_1$-алгебры $A_1$ в $D_2$-алгебру $A_2$ и гомоморфизм (11.1.38) определён однозначно.

Равенство

$$(11.1.39) \quad \begin{aligned} &\overline{f} \circ (e_{V_1}{}^*{}_*(v + w)) \\ &= e_{V_2}{}^*{}_*f^*{}_*(g \circ (v + w)) \\ &= e_{V_2}{}^*{}_*f^*{}_*(g \circ v) \\ &+ e_{V_2}{}^*{}_*f^*{}_*(g \circ w) \\ &= \overline{f} \circ (e_{V_1}{}^*{}_*v) + \overline{f} \circ (e_{V_1}{}^*{}_*w) \end{aligned}$$

является следствием равенств

$$(8.1.18) \quad a^*{}_*(b_1 + b_2) = a^*{}_*b_1 + a^*{}_*b_2 \qquad (11.1.6) \quad g \circ (a + b) = g \circ a + g \circ b$$

$$(11.1.10) \quad f \circ (va) = (f \circ v)(g \circ a) \qquad (11.1.31) \quad \overline{f} \circ (e_{V_1}{}^*{}_*v) = e_{V_2}{}^*{}_*f^*{}_*(\overline{g} \circ v)$$

Из равенства (11.1.39) следует, что отображение $\overline{f}$ является гомоморфизмом абелевой группы. Равенство

$$(11.1.40) \quad \begin{aligned} &\overline{f} \circ (e_{V_1}{}^*{}_*(va)) \\ &= e_{V_2}{}^*{}_*f^*{}_*(\overline{g} \circ (va)) \\ &= (e_{V_2}{}^*{}_*f^*{}_*(\overline{g} \circ v))(\overline{g} \circ a) \\ &= (\overline{f} \circ (e_{V_1}{}^*{}_*v))(\overline{g} \circ a) \end{aligned}$$

является следствием равенств

$$(11.1.10) \quad f \circ (va) = (f \circ v)(g \circ a)$$

$$(5.4.4) \quad g \circ (ab) = (g \circ a)(g \circ b)$$

$$(11.1.31) \quad \overline{f} \circ (e_{V_1}{}^*{}_*v) = e_{V_2}{}^*{}_* f^*{}_* (\overline{g} \circ v)$$

Из равенства (11.1.40), и определений 2.8.5, 11.1.1 следует, что отображение

$$(11.1.3) \quad (h : D_1 \to D_2 \quad \overline{g} : A_1 \to A_2 \quad \overline{f} : V_1 \to V_2)$$

является гомоморфизмом $A_1$-векторного пространства строк $V_1$ в $A_2$-векторное пространство строк $V_2$.

Пусть $f$ - матрица гомоморфизмов $\overline{f}$, $\overline{g}$ относительно базисов $\overline{\overline{e}}_V$, $\overline{\overline{e}}_W$. Равенство

$$\overline{f} \circ (e_V{}^*{}_*v) = e_W{}^*{}_* f^*{}_* v = \overline{g} \circ (e_V{}^*{}_*v)$$

является следствием теоремы 11.3.7. Следовательно, $\overline{f} = \overline{g}$. $\qquad \square$

На основании теорем 11.1.6, 11.1.7 мы идентифицируем гомоморфизм

$$(11.1.3) \quad (h : D_1 \to D_2 \quad \overline{g} : A_1 \to A_2 \quad \overline{f} : V_1 \to V_2)$$

правого $A$-векторного пространства $V$ строк и координаты его представления

$$(11.1.28) \quad \overline{g} \circ (ae_{A_1}) = h(a)ge_{A_2}$$

$$(11.1.31) \quad \overline{f} \circ (e_{V_1}{}^*{}_*v) = e_{V_2}{}^*{}_* f^*{}_* (\overline{g} \circ v)$$

## 11.2. Гомоморфизм, когда кольца $D_1 = D_2 = D$

Пусть $A_i$, $i = 1, 2$, - алгебра с делением над коммутативным кольцом $D$. Пусть $V_i$, $i = 1, 2$, - правое $A_i$-модуль.

Определение 11.2.1. *Пусть*

$$(11.2.1)$$

$$A_1 \xrightarrow{g_{1.23}} A_1 \xrightarrow{g_{1.34}} V_1$$

$g_{1.12}(d) : a \to d\,a$

$g_{1.23}(v) : w \to C_1(w, v)$

$C_1 \in \mathcal{L}(A_1^2 \to A_1)$

$g_{1.34}(a) : v \to v\,a$

$g_{1.14}(d) : v \to v\,d$

*диаграмма представлений, описывающая  правое $A_1$-векторное пространство $V_1$. Пусть*

$$(11.2.2)$$

$$A_2 \xrightarrow{g_{2.23}} A_2 \xrightarrow{g_{2.34}} V_2$$

$g_{2.12}(d) : a \to d\,a$

$g_{2.23}(v) : w \to C_2(w, v)$

$C_2 \in \mathcal{L}(A_2^2 \to A_2)$

$g_{2.34}(a) : v \to v\,a$

$g_{2.14}(d) : v \to v\,d$

*диаграмма представлений, описывающая  правое $A_2$-векторное пространство*

$V_2$. *Морфизм*

$$(11.2.3) \qquad (g : A_1 \to A_2 \quad f : V_1 \to V_2)$$

*диаграммы представлений* (11.2.1) *в диаграмму представлений* (11.2.2) *называется* **гомоморфизмом** *правого $A_1$-векторного пространства $V_1$ в правое $A_2$-векторное пространство $V_2$. Обозначим* $\mathrm{Hom}(D; A_1 \to A_2; V_1* \to V_2*)$ *множество гомоморфизмов правого $A_1$-векторного пространства $V_1$ в правое $A_2$-векторное пространство $V_2$.* □

Мы будем пользоваться записью

$$f \circ a = f(a)$$

для образа гомоморфизма $f$.

Теорема 11.2.2. *Гомоморфизм*

$$\boxed{(11.2.3) \quad (g : A_1 \to A_2 \quad f : V_1 \to V_2)}$$

*правого $A_1$-векторного пространства $V_1$ в правое $A_2$-векторное пространство $V_2$ удовлетворяет следующим равенствам*

$$(11.2.4) \qquad g \circ (a + b) = g \circ a + g \circ b$$

$$(11.2.5) \qquad g \circ (ab) = (g \circ a)(g \circ b)$$

$$(11.2.6) \qquad g \circ (ap) = (g \circ a)p$$

$$(11.2.7) \qquad f \circ (u + v) = f \circ u + f \circ v$$

$$(11.2.8) \qquad f \circ (va) = (f \circ v)(g \circ a)$$

$$p \in D \quad a, b \in A_1 \quad u, v \in V_1$$

Доказательство. Согласно определениям 2.3.9, 2.8.5, 7.4.1, 1 отображение

$$g : A_1 \to A_2$$

является гомоморфизмом $D$-алгебры $A_1$ в $D$-алгебру $A_2$. Следовательно, равенства (11.2.4), (11.2.5), (11.2.6) являются следствием теоремы 4.3.10. Равенство (11.2.7) является следствием определения 11.2.1, так как, согласно определению 2.3.9, отображение $f$ является гомоморфизмом абелевой группы. Равенство (11.2.8) является следствием равенства (2.3.6), так как отображение

$$(g : A_1 \to A_2 \quad f : V_1 \to V_2)$$

является морфизмом представления $g_{1.34}$ в представление $g_{2.34}$. □

Определение 11.2.3. *Гомоморфизм* [11.4]

$$(g : A_1 \to A_2 \quad f : V_1 \to V_2)$$

*называется* **изоморфизмом** *между правым $A_1$-векторным пространством $V_1$ и правым $A_2$-векторным пространством $V_2$, если существует отображение*

$$(g^{-1} : A_2 \to A_1 \quad f^{-1} : V_2 \to V_1)$$

---

11.4 Я следую определению на странице [18]-63.

*которое является гомоморфизмом.*          □

| ТЕОРЕМА 11.2.4. *Гомоморфизм*[11.5] | ТЕОРЕМА 11.2.5. *Гомоморфизм*[11.6] |
|---|---|
| (11.2.3)   $(\overline{g}: A_1 \to A_2 \quad \overline{f}: V_1 \to V_2)$ | (11.2.3)   $(\overline{g}: A_1 \to A_2 \quad \overline{f}: V_1 \to V_2)$ |

*правого $A_1$-векторного пространства столбцов $V_1$ в правое $A_2$-векторное пространство столбцов $V_2$ имеет представление*

$$(11.2.9) \qquad b = ga$$

$$(11.2.10) \qquad \overline{g} \circ (a^i e_{A_1 i}) = a^i g_i^k e_{A_2 k}$$

$$(11.2.11) \qquad \overline{g} \circ (e_{A_1} a) = e_{A_2} ga$$

$$(11.2.12) \qquad w = f_* {}^* (\overline{g} \circ v)$$

$$(11.2.13) \qquad \overline{f} \circ (e_{V_1 i} v^i) = e_{V_2 k} f_i^k (\overline{g} \circ v^i)$$

$$(11.2.14) \qquad \overline{f} \circ (e_{V_1 *} {}^* v) = e_{V_2 *} {}^* f_* {}^* (\overline{g} \circ v)$$

*относительно выбранных базисов. Здесь*

- $a$ - *координатная матрица $A_1$-числа $\overline{a}$ относительно базиса* $\overline{\overline{e}}_{A_1}$ *(A)*

  $$\overline{a} = e_{A_1} a$$

- $b$ - *координатная матрица $A_2$-числа*

  $$\overline{b} = \overline{g} \circ \overline{a}$$

  *относительно базиса* $\overline{\overline{e}}_{A_2}$

  $$\overline{b} = e_{A_2} b$$

- $g$ - *координатная матрица множества $A_2$-чисел* $(\overline{g} \circ e_{A_1 k}, k \in K)$ *относительно базиса* $\overline{\overline{e}}_{A_2}$.

*правого $A_1$-векторного пространства строк $V_1$ в правое $A_2$-векторное пространство строк $V_2$ имеет представление*

$$(11.2.16) \qquad b = ag$$

$$(11.2.17) \qquad \overline{g} \circ (a_i e_{A_1}^i) = a_i g_k^i e_{A_2}^k$$

$$(11.2.18) \qquad \overline{g} \circ (ae_{A_1}) = age_{A_2}$$

$$(11.2.19) \qquad w = f^*{}_* (\overline{g} \circ v)$$

$$(11.2.20) \qquad \overline{f} \circ (e_{V_1}^i v_i) = e_{V_2}^k f_k^i (\overline{g} \circ v_i)$$

$$(11.2.21) \qquad \overline{f} \circ (e_{V_1}^* {}^* v) = e_{V_2}^* {}^* f^*{}_* (\overline{g} \circ v)$$

*относительно выбранных базисов. Здесь*

- $a$ - *координатная матрица $A_1$-числа $\overline{a}$ относительно базиса* $\overline{\overline{e}}_{A_1}$ *(A)*

  $$\overline{a} = ae_{A_1}$$

- $b$ - *координатная матрица $A_2$-числа*

  $$\overline{b} = \overline{g} \circ \overline{a}$$

  *относительно базиса* $\overline{\overline{e}}_{A_2}$

  $$\overline{b} = be_{A_2}$$

- $g$ - *координатная матрица множества $A_2$-чисел* $(\overline{g} \circ e_{A_1}^k, k \in K)$ *относительно базиса* $\overline{\overline{e}}_{A_2}$.

---

[11.5] В теоремах 11.2.4, 11.2.6, мы опираемся на следующее соглашение. Пусть $i = 1$, 2. Пусть множество векторов $\overline{\overline{e}}_{A_i} = (e_{A_i k}, k \in K_i)$ является базисом и $C_{i\,j}^{\quad k}$, $k$, $i$, $j \in K_i$, - структурные константы $D$-алгебры столбцов $A_i$. Пусть множество векторов $\overline{\overline{e}}_{V_i} = (e_{V_i i}, i \in I_i)$ является базисом правого $A_i$-векторного пространства $V_i$.

[11.6] В теоремах 11.2.5, 11.2.7, мы опираемся на следующее соглашение. Пусть $i = 1, 2$. Пусть множество векторов $\overline{\overline{e}}_{A_i} = (e_{A_i}^k, k \in K_i)$ является базисом и $C_{i\,k}^{\quad ij}$, $k$, $i$, $j \in K_i$, - структурные константы $D$-алгебры строк $A_i$. Пусть множество векторов $\overline{\overline{e}}_{V_i} = (e_{V_i}^i, i \in I_i)$ является базисом правого $A_i$-векторного пространства $V_i$.

- *Матрица гомоморфизма и структурные константы связаны соотношением*

$$(11.2.15) \qquad C_{1\ ij}^{\ k} g_k^l = g_i^p g_j^q C_{2\ pq}^{\ l}$$

- *$v$ - координатная матрица $V_1$-числа $\overline{v}$ относительно базиса $\overline{\overline{e}}_{V_1}$ (A)*

$$\overline{v} = e_{V_1 *}{}^* v$$

- *$g(v) = (g(v_i), i \in I)$ - матрица $A_2$-чисел.*

- *$w$ - координатная матрица $V_2$-числа*

$$\overline{w} = \overline{f} \circ \overline{v}$$

*относительно базиса $\overline{\overline{e}}_{V_2}$*

$$\overline{w} = e_{V_2 *}{}^* w$$

- *$f$ - координатная матрица множества $V_2$-чисел ($\overline{f} \circ e_{V_1 i}, i \in I$) относительно базиса $\overline{\overline{e}}_{V_2}$.*

*Для заданного гомоморфизма*

$$\boxed{(11.2.3) \quad (\overline{g} : A_1 \to A_2 \quad \overline{f} : V_1 \to V_2)}$$

*множество матриц $(g, f)$ определено однозначно и называется* **координатами гомоморфизма** (11.2.3) *относительно базисов* $(\overline{\overline{e}}_{A_1}, \overline{\overline{e}}_{V_1}), \quad (\overline{\overline{e}}_{A_2}, \overline{\overline{e}}_{V_2})$.

Доказательство. Согласно определениям 2.3.9, 2.8.5, 5.4.7, 11.2.1, отображение

$$\overline{g} : A_1 \to A_2$$

является гомоморфизмом $D$-алгебры $A_1$ в $D$-алгебру $A_2$. Следовательно, равенства (11.2.9), (11.2.10), (11.2.11), (11.2.15) следуют из равенств

$$\boxed{(5.4.53) \quad b = fa}$$

$$\boxed{(5.4.54) \quad \overline{f} \circ (a^i e_{1i}) = a^i f_i^k e_k}$$

$$\boxed{(5.4.55) \quad \overline{f} \circ (e_{1a}) = e_2 fa}$$

$$\boxed{(5.4.58) \quad C_{1\ ij}^{\ k} f_k^l = f_i^p f_j^q C_{2\ pq}^{\ l}}$$

Из теоремы 4.3.11, следует что матрица $g$ определена однозначно.

- *Матрица гомоморфизма и структурные константы связаны соотношением*

$$(11.2.22) \qquad C_{1\ k}^{\ ij} g_l^k = g_p^i g_q^j C_{2\ l}^{\ pq}$$

- *$v$ - координатная матрица $V_1$-числа $\overline{v}$ относительно базиса $\overline{\overline{e}}_{V_1}$ (A)*

$$\overline{v} = e_{V_1}{}^*{}_* v$$

- *$g(v) = (g(v^i), i \in I)$ - матрица $A_2$-чисел.*

- *$w$ - координатная матрица $V_2$-числа*

$$\overline{w} = \overline{f} \circ \overline{v}$$

*относительно базиса $\overline{\overline{e}}_{V_2}$*

$$\overline{w} = e_{V_2}{}^*{}_* w$$

- *$f$ - координатная матрица множества $V_2$-чисел ($\overline{f} \circ e_{V_1}^i, i \in I$) относительно базиса $\overline{\overline{e}}_{V_2}$.*

*Для заданного гомоморфизма*

$$\boxed{(11.2.3) \quad (\overline{g} : A_1 \to A_2 \quad \overline{f} : V_1 \to V_2)}$$

*множество матриц $(g, f)$ определено однозначно и называется* **координатами гомоморфизма** (11.2.3) *относительно базисов* $(\overline{\overline{e}}_{A_1}, \overline{\overline{e}}_{V_1}), \quad (\overline{\overline{e}}_{A_2}, \overline{\overline{e}}_{V_2})$.

Доказательство. Согласно определениям 2.3.9, 2.8.5, 5.4.7, 11.2.1, отображение

$$\overline{g} : A_1 \to A_2$$

является гомоморфизмом $D$-алгебры $A_1$ в $D$-алгебру $A_2$. Следовательно, равенства (11.2.16), (11.2.17), (11.2.18), (11.2.22) следуют из равенств

$$\boxed{(5.4.64) \quad b = af}$$

$$\boxed{(5.4.65) \quad \overline{f} \circ (a_i e_1^i) = a_i f_k^i e^k}$$

$$\boxed{(5.4.66) \quad \overline{f} \circ (ae_1) = afe_2}$$

$$\boxed{(5.4.69) \quad C_{1\ k}^{\ ij} f_l^k = f_p^i f_q^j C_{2\ l}^{\ pq}}$$

Из теоремы 4.3.12, следует что матри-

Вектор $\overline{v} \in V_1$ имеет разложение

$$(11.2.23) \qquad \overline{v} = e_{V_1 *}{}^* v$$

относительно базиса $\overline{\overline{e}}_{V_1}$. Вектор $\overline{w} \in V_2$ имеет разложение

$$(11.2.24) \qquad \overline{w} = e_{V_2 *}{}^* w$$

относительно базиса $\overline{\overline{e}}_{V_2}$. Так как $\overline{f}$ - гомоморфизм, то равенство

$$(11.2.25) \quad \begin{aligned} \overline{w} &= \overline{f} \circ \overline{v} = \overline{f} \circ (e_{V_1 *}{}^* v) \\ &= (\overline{f} \circ e_{V_1})_* {}^* (\overline{g} \circ v) \end{aligned}$$

является следствием равенств (11.2.7), (11.2.8). $V_2$-число $\overline{f} \circ e_{V_1 i}$ имеет разложение

$$(11.2.26) \quad \overline{f} \circ e_{V_1 i} = e_{V_2 *}{}^* f_i = e_{V_2}^j f_i^j$$

относительно базиса $\overline{\overline{e}}_{V_2}$. Равенство

$$(11.2.27) \qquad \overline{w} = e_{V_2 *}{}^* f_*{}^* (\overline{g} \circ v)$$

является следствием равенств (11.2.25), (11.2.26). Равенство (11.2.12) следует из сравнения (11.2.24) и (11.2.27) и теоремы 8.4.18. Из равенства (11.2.26) и теоремы 8.4.18 следует что матрица $f$ определена однозначно. $\square$

ца $g$ определена однозначно.

Вектор $\overline{v} \in V_1$ имеет разложение

$$(11.2.28) \qquad \overline{v} = e_{V_1}{}^*{}_* v$$

относительно базиса $\overline{\overline{e}}_{V_1}$. Вектор $\overline{w} \in V_2$ имеет разложение

$$(11.2.29) \qquad \overline{w} = e_{V_2}{}^*{}_* w$$

относительно базиса $\overline{e}_{V_2}$. Так как $\overline{f}$ - гомоморфизм, то равенство

$$(11.2.30) \quad \begin{aligned} \overline{w} &= \overline{f} \circ \overline{v} = \overline{f} \circ (e_{V_1}{}^*{}_* v) \\ &= (\overline{f} \circ e_{V_1})^* {}_* (\overline{g} \circ v) \end{aligned}$$

является следствием равенств (11.2.7), (11.2.8). $V_2$-число $\overline{f} \circ e_{V_1}^i$ имеет разложение

$$(11.2.31) \quad \overline{f} \circ e_{V_1}^i = e_{V_2}{}^*{}_* f^i = e_{V_2 j} f_j^i$$

относительно базиса $\overline{e}_{V_2}$. Равенство

$$(11.2.32) \qquad \overline{w} = e_{V_2}{}^*{}_* f^*{}_* (\overline{g} \circ v)$$

является следствием равенств (11.2.30), (11.2.31). Равенство (11.2.19) следует из сравнения (11.2.29) и (11.2.32) и теоремы 8.4.25. Из равенства (11.2.31) и теоремы 8.4.25 следует что матрица $f$ определена однозначно. $\square$

---

**Теорема 11.2.6.** *Пусть*

$$g = (g_l^k, k \in K, l \in L)$$

*матрица $D$-чисел, которая удовлетворяют равенству*

$$(11.2.15) \quad \boxed{C_{1\,ij}^{\;k} g_k^l = g_i^p g_j^q C_{2\,pq}^l}$$

*Пусть*

$$f = (f_j^i, i \in I, j \in J)$$

*матрица $A_2$-чисел.     Отображение* 11.5

$$(11.2.3) \quad \boxed{(\overline{g} : A_1 \to A_2 \quad \overline{f} : V_1 \to V_2)}$$

*определённое равенствами*

$$(11.2.11) \quad \boxed{\overline{f} \circ (e_{A_1} a) = e_{A_2} fa}$$

$$(11.2.14) \quad \boxed{\overline{f} \circ (e_{V_1 *}{}^* v) = e_{V_2 *}{}^* f_*{}^* (\overline{g} \circ v)}$$

*является   гомоморфизмом правого $A_1$-векторного пространства столбцов $V_1$ в правое $A_2$-векторное*

**Теорема 11.2.7.** *Пусть*

$$g = (g_k^l, k \in K, l \in L)$$

*матрица $D$-чисел, которая удовлетворяют равенству*

$$(11.2.22) \quad \boxed{C_{1\,k}^{\;ij} g_l^k = g_p^i g_q^j C_{2\,l}^{pq}}$$

*Пусть*

$$f = (f_i^j, i \in I, j \in J)$$

*матрица $A_2$-чисел.     Отображение* 11.6

$$(11.2.3) \quad \boxed{(\overline{g} : A_1 \to A_2 \quad \overline{f} : V_1 \to V_2)}$$

*определённое равенствами*

$$(11.2.18) \quad \boxed{\overline{f} \circ (a e_{A_1}) = a f e_{A_2}}$$

$$(11.2.21) \quad \boxed{\overline{f} \circ (e_{V_1}{}^*{}_* v) = e_{V_2}{}^*{}_* f^*{}_* (\overline{g} \circ v)}$$

*является   гомоморфизмом право-го $A_1$-векторного пространства строк $V_1$ в правое $A_2$-векторное пространство строк $V_2$.   Гомоморфизм* пространство столбцов $V_2$. Гомоморфизм

$$(11.2.3) \quad (\overline{g} : A_1 \to A_2 \quad \overline{f} : V_1 \to V_2)$$

который имеет данное множество матриц $(g, f)$, определён однозначно.

Доказательство. Согласно теореме 5.4.11, отображение

$$(11.2.33) \quad \overline{g} : A_1 \to A_2$$

является гомоморфизмом $D$-алгебры $A_1$ в $D$-алгебру $A_2$ и гомоморфизм (11.2.33) определён однозначно.

Равенство

$$(11.2.34) \quad \begin{aligned} & \overline{f} \circ (e_{V_1 *}{}^*(v + w)) \\ &= e_{V_2 *}{}^* f_*{}^*(g \circ (v + w)) \\ &= e_{V_2 *}{}^* f_*{}^*(g \circ v) \\ &\quad + e_{V_2 *}{}^* f_*{}^*(g \circ w) \\ &= \overline{f} \circ (e_{V_1 *}{}^* v) + \overline{f} \circ (e_{V_1 *}{}^* w) \end{aligned}$$

является следствием равенств

$$(8.1.16) \quad a_*{}^*(b_1 + b_2) = a_*{}^* b_1 + a_*{}^* b_2$$

$$(11.1.6) \quad g \circ (a + b) = g \circ a + g \circ b$$

$$(11.1.10) \quad f \circ (va) = (f \circ v)(g \circ a)$$

$$(11.1.16) \quad \overline{f} \circ (e_{V_1 *}{}^* v) = e_{V_2 *}{}^* f_*{}^*(\overline{g} \circ v)$$

Из равенства (11.2.34) следует, что отображение $\overline{f}$ является гомоморфизмом абелевой группы. Равенство

$$(11.2.35) \quad \begin{aligned} & \overline{f} \circ (e_{V_1 *}{}^*(va)) \\ &= e_{V_2 *}{}^* f_*{}^*(\overline{g} \circ (va)) \\ &= (e_{V_2 *}{}^* f_*{}^*(\overline{g} \circ v))(\overline{g} \circ a) \\ &= (\overline{f} \circ (e_{V_1 *}{}^* v))(\overline{g} \circ a) \end{aligned}$$

является следствием равенств

$$(11.1.10) \quad f \circ (va) = (f \circ v)(g \circ a)$$

$$(5.4.4) \quad g \circ (ab) = (g \circ a)(g \circ b)$$

$$(11.1.16) \quad \overline{f} \circ (e_{V_1 *}{}^* v) = e_{V_2 *}{}^* f_*{}^*(\overline{g} \circ v)$$

Из равенства (11.2.35), и определений 2.8.5, 11.2.1 следует, что отображение

$$(11.2.3) \quad (\overline{g} : A_1 \to A_2 \quad \overline{f} : V_1 \to V_2)$$

который имеет данное множество матриц $(g, f)$, определён однозначно.

Доказательство. Согласно теореме 5.4.12, отображение

$$(11.2.36) \quad \overline{g} : A_1 \to A_2$$

является гомоморфизмом $D$-алгебры $A_1$ в $D$-алгебру $A_2$ и гомоморфизм (11.2.36) определён однозначно.

Равенство

$$(11.2.37) \quad \begin{aligned} & \overline{f} \circ (e_{V_1}{}^*{}_*(v + w)) \\ &= e_{V_2}{}^*{}_* f^*{}_*(g \circ (v + w)) \\ &= e_{V_2}{}^*{}_* f^*{}_*(g \circ v) \\ &\quad + e_{V_2}{}^*{}_* f^*{}_*(g \circ w) \\ &= \overline{f} \circ (e_{V_1}{}^*{}_* v) + \overline{f} \circ (e_{V_1}{}^*{}_* w) \end{aligned}$$

является следствием равенств

$$(8.1.18) \quad a^*{}_*(b_1 + b_2) = a^*{}_* b_1 + a^*{}_* b_2$$

$$(11.1.6) \quad g \circ (a + b) = g \circ a + g \circ b$$

$$(11.1.10) \quad f \circ (va) = (f \circ v)(g \circ a)$$

$$(11.1.31) \quad \overline{f} \circ (e_{V_1}{}^*{}_* v) = e_{V_2}{}^*{}_* f^*{}_*(\overline{g} \circ v)$$

Из равенства (11.2.37) следует, что отображение $\overline{f}$ является гомоморфизмом абелевой группы. Равенство

$$(11.2.38) \quad \begin{aligned} & \overline{f} \circ (e_{V_1}{}^*{}_*(va)) \\ &= e_{V_2}{}^*{}_* f^*{}_*(\overline{g} \circ (va)) \\ &= (e_{V_2}{}^*{}_* f^*{}_*(\overline{g} \circ v))(\overline{g} \circ a) \\ &= (\overline{f} \circ (e_{V_1}{}^*{}_* v))(\overline{g} \circ a) \end{aligned}$$

является следствием равенств

$$(11.1.10) \quad f \circ (va) = (f \circ v)(g \circ a)$$

$$(5.4.4) \quad g \circ (ab) = (g \circ a)(g \circ b)$$

$$(11.1.31) \quad \overline{f} \circ (e_{V_1}{}^*{}_* v) = e_{V_2}{}^*{}_* f^*{}_*(\overline{g} \circ v)$$

Из равенства (11.2.38), и определений 2.8.5, 11.2.1 следует, что отображение

$$(11.2.3) \quad (\overline{g} : A_1 \to A_2 \quad \overline{f} : V_1 \to V_2)$$

$$\boxed{(11.2.3) \quad (\overline{g} : A_1 \to A_2 \quad \overline{f} : V_1 \to V_2)}$$

является гомоморфизмом $A_1$-векторного пространства столбцов $V_1$ в $A_2$-векторное пространство столбцов $V_2$.

Пусть $f$ - матрица гомоморфизмов $\overline{f}$, $\overline{g}$ относительно базисов $\overline{\overline{e}}_V$, $\overline{\overline{e}}_W$. Равенство

$$\overline{f} \circ (e_{V*}{}^* v) = e_{W*}{}^* f_* v = \overline{g} \circ (e_{V*}{}^* v)$$

является следствием теоремы 11.3.6. Следовательно, $\overline{f} = \overline{g}$. □

На основании теорем 11.2.4, 11.2.6 мы идентифицируем гомоморфизм

$$\boxed{(11.2.3) \quad (\overline{g} : A_1 \to A_2 \quad \overline{f} : V_1 \to V_2)}$$

правого $A$-векторного пространства $V$ столбцов и координаты его представления

$$\boxed{(11.2.11) \quad \overline{f} \circ (e_{A_1} a) = e_{A_2} fa}$$

$$\boxed{(11.2.14) \quad \overline{f} \circ (e_{V_1*}{}^* v) = e_{V_2*}{}^* f_*{}^* (\overline{g} \circ v)}$$

является гомоморфизмом $A_1$-векторного пространства строк $V_1$ в $A_2$-векторное пространство строк $V_2$.

Пусть $f$ - матрица гомоморфизмов $\overline{f}$, $\overline{g}$ относительно базисов $\overline{\overline{e}}_V$, $\overline{\overline{e}}_W$. Равенство

$$\overline{f} \circ (e_{V*}{}^* v) = e_{W*}{}^* f^*{}_* v = \overline{g} \circ (e_{V*}{}^* v)$$

является следствием теоремы 11.3.7. Следовательно, $\overline{f} = \overline{g}$. □

На основании теорем 11.2.5, 11.2.7 мы идентифицируем гомоморфизм

$$\boxed{(11.2.3) \quad (\overline{g} : A_1 \to A_2 \quad \overline{f} : V_1 \to V_2)}$$

правого $A$-векторного пространства $V$ строк и координаты его представления

$$\boxed{(11.2.18) \quad \overline{f} \circ (a e_{A_1}) = a f e_{A_2}}$$

$$\boxed{(11.2.21) \quad \overline{f} \circ (e_{V_1*}{}^* v) = e_{V_2*}{}^* f^*{}_* (\overline{g} \circ v)}$$

## 11.3. Гомоморфизм, когда $D$-алгебры $A_1 = A_2 = A$

Пусть $A$ - алгебра с делением над коммутативным кольцом $D$. Пусть $V_i$, $i = 1, 2$, - правое $A$-векторное пространство.

ОПРЕДЕЛЕНИЕ 11.3.1. *Пусть*

$$(11.3.1)$$

$$g_{1.12}(d) : a \to d\,a$$
$$g_{1.23}(v) : w \to C(w, v)$$
$$C \in \mathcal{L}(A^2 \to A)$$
$$g_{1.34}(a) : v \to v\,a$$
$$g_{1.14}(d) : v \to v\,d$$

*диаграмма представлений, описывающая правое $A$-векторное пространство $V_1$. Пусть*

$$(11.3.2)$$

$$g_{2.12}(d) : a \to d\,a$$
$$g_{2.23}(v) : w \to C(w, v)$$
$$C \in \mathcal{L}(A^2 \to A)$$
$$g_{2.34}(a) : v \to v\,a$$
$$g_{2.14}(d) : v \to v\,d$$

*диаграмма представлений, описывающая  правое $A$-векторное пространство $V_2$. Морфизм*

$$(11.3.3) \qquad f : V_1 \to V_2$$

*диаграммы представлений* (11.3.1) *в диаграмму представлений* (11.3.2) *называется* **гомоморфизмом** *правого $A$-векторного пространства $V_1$ в правое $A$-векторное пространство $V_2$. Обозначим* $\mathrm{Hom}(D; A; V_1* \to V_2*)$  *множество гомоморфизмов правого $A$-векторного пространства $V_1$ в  правое $A$-векторное пространство $V_2$.* □

Мы будем пользоваться записью

$$f \circ a = f(a)$$

для образа гомоморфизма $f$.

ТЕОРЕМА 11.3.2. *Гомоморфизм*

$$\boxed{(11.3.3) \quad f : V_1 \to V_2}$$

*правого $A$-векторного пространства $V_1$ в  правое $A$-векторное пространство $V_2$ удовлетворяет следующим равенствам*

$$(11.3.4) \qquad f \circ (u + v) = f \circ u + f \circ v$$

$$(11.3.5) \qquad f \circ (va) = (f \circ v)a$$

$$a \in A \quad u, v \in V_1$$

ДОКАЗАТЕЛЬСТВО.  Равенство (11.3.4) является следствием определения 11.3.1, так как, согласно определению 2.3.9, отображение $f$ является гомоморфизмом абелевой группы. Равенство (11.3.5) является следствием равенства (2.3.6), так как отображение

$$f : V_1 \to V_2$$

является морфизмом представления $g_{1.34}$ в представление $g_{2.34}$. □

ОПРЕДЕЛЕНИЕ 11.3.3. *Гомоморфизм* [11.7]

$$f : V_1 \to V_2$$

*называется* **изоморфизмом** *между правым $A$-векторным пространством $V_1$ и правым $A$-векторным пространством $V_2$, если существует отображение*

$$f^{-1} : V_2 \to V_1$$

*которое является гомоморфизмом.* □

ОПРЕДЕЛЕНИЕ 11.3.4. *Гомоморфизм* [11.7]

$$f : V \to V$$

*источником и целью которого является одно и тоже  правое $A$-векторное пространство, называется* **эндоморфизмом**. *Эндоморфизм*

$$f : V \to V$$

---

[11.7] Я следую определению на странице [18]-63.

правого $A$-векторного пространства $V$ называется **автоморфизмом**, если существует отображение $f^{-1}$ которое является эндоморфизмом. $\square$

---

Теорема 11.3.5. *Множество $GL(V)$ автоморфизмов правого $A$-векторного пространства $V$ является группой.*

---

Теорема 11.3.6. *Гомоморфизм* [11.8]

$$\boxed{\text{(11.3.3)}\quad \overline{f} : V_1 \to V_2}$$

*правого $A$-векторного пространства столбцов $V_1$ в правое $A$-векторное пространство столбцов $V_2$ имеет представление*

(11.3.6)     $w = f_*{}^* v$

(11.3.7)     $\overline{f} \circ (v i^{e_{V_1 i}}) = e_{V_2 k} f_i^k v^i$

(11.3.8)     $\overline{f} \circ (e_{V_1 *}{}^* v) = e_{V_2 *}{}^* f_* {}^* v$

*относительно выбранных базисов. Здесь*

- $v$ - *координатная матрица $V_1$-числа $\overline{v}$ относительно базиса $\overline{\overline{e}}_{V_1}$ (A)*

$$\overline{v} = e_{V_1 *}{}^* v$$

- $w$ - *координатная матрица $V_2$-числа*

$$\overline{w} = \overline{f} \circ \overline{v}$$

*относительно базиса $\overline{\overline{e}}_{V_2}$*

$$\overline{w} = e_{V_2 *}{}^* w$$

- $f$ - *координатная матрица множества $V_2$-чисел ($\overline{f} \circ e_{V_1 i}, i \in I$) относительно базиса $\overline{\overline{e}}_{V_2}$.*

*Для заданного гомоморфизма $\overline{f}$ матрица $f$ определена однозначно и называется* **матрицей гомоморфизма $\overline{f}$** *относительно базисов $\overline{\overline{e}}_1$, $\overline{\overline{e}}_2$.*

---

Теорема 11.3.7. *Гомоморфизм* [11.9]

$$\boxed{\text{(11.3.3)}\quad \overline{f} : V_1 \to V_2}$$

*правого $A$-векторного пространства строк $V_1$ в правое $A$-векторное пространство строк $V_2$ имеет представление*

(11.3.9)     $w = f^* {}_* v$

(11.3.10)     $\overline{f} \circ (v i_{e_{V_1}^i}) = e_{V_2}^k f_k^i v_i$

(11.3.11)     $\overline{f} \circ (e_{V_1}^* {}_* v) = e_{V_2}^* {}_* f^* {}_* v$

*относительно выбранных базисов. Здесь*

- $v$ - *координатная матрица $V_1$-числа $\overline{v}$ относительно базиса $\overline{\overline{e}}_{V_1}$ (A)*

$$\overline{v} = e_{V_1}^* {}_* v$$

- $w$ - *координатная матрица $V_2$-числа*

$$\overline{w} = \overline{f} \circ \overline{v}$$

*относительно базиса $\overline{\overline{e}}_{V_2}$*

$$\overline{w} = e_{V_2}^* {}_* w$$

- $f$ - *координатная матрица множества $V_2$-чисел ($\overline{f} \circ e_{V_1}^i, i \in I$) относительно базиса $\overline{\overline{e}}_{V_2}$.*

*Для заданного гомоморфизма $\overline{f}$ матрица $f$ определена однозначно и называется* **матрицей гомоморфизма $\overline{f}$** *относительно базисов $\overline{\overline{e}}_1$, $\overline{\overline{e}}_2$.*

---

[11.8] В теоремах 11.3.6, 11.3.8, мы опираемся на следующее соглашение. Пусть $i = 1$, 2. Пусть множество векторов $\overline{\overline{e}}_{V_i} = (e_{V_i i}, i \in I_i)$ является базисом правого $A$-векторного пространства $V_i$.

[11.9] В теоремах 11.3.7, 11.3.9, мы опираемся на следующее соглашение. Пусть $i = 1$, 2. Пусть множество векторов $\overline{\overline{e}}_{V_i} = (e_{V_i}^i, i \in I_i)$ является базисом правого $A$-векторного пространства $V_i$.

**Доказательство.** Вектор $\overline{v} \in V_1$ имеет разложение

$$(11.3.12) \qquad \overline{v} = e_{V_1 *}{}^* v$$

относительно базиса $\overline{\overline{e}}_{V_1}$. Вектор $\overline{w} \in V_2$ имеет разложение

$$(11.3.13) \qquad \overline{w} = e_{V_2 *}{}^* w$$

относительно базиса $\overline{\overline{e}}_{V_2}$. Так как $\overline{f}$ - гомоморфизм, то равенство

$$(11.3.14) \qquad \begin{aligned} \overline{w} &= \overline{f} \circ \overline{v} = \overline{f} \circ (e_{V_1 *}{}^* v) \\ &= (\overline{f} \circ e_{V_1})_*{}^* v \end{aligned}$$

является следствием равенств (11.3.4), (11.3.5). $V_2$-число $\overline{f} \circ e_{V_1 i}$ имеет разложение

$$(11.3.15) \quad \overline{f} \circ e_{V_1 i} = e_{V_2 *}{}^* f_i = e_{V_2}^j f_i^j$$

относительно базиса $\overline{\overline{e}}_{V_2}$. Равенство

$$(11.3.16) \qquad \overline{w} = e_{V_2 *}{}^* f_*{}^* v$$

является следствием равенств (11.3.14), (11.3.15). Равенство (11.3.6) следует из сравнения (11.3.13) и (11.3.16). Из равенства (11.3.15) и теоремы 8.4.18 следует что матрица $f$ определена однозначно. $\square$

---

**Доказательство.** Вектор $\overline{v} \in V_1$ имеет разложение

$$(11.3.17) \qquad \overline{v} = e_{V_1}{}_*{}^* v$$

относительно базиса $\overline{\overline{e}}_{V_1}$. Вектор $\overline{w} \in V_2$ имеет разложение

$$(11.3.18) \qquad \overline{w} = e_{V_2 *}{}^* w$$

относительно базиса $\overline{\overline{e}}_{V_2}$. Так как $\overline{f}$ - гомоморфизм, то равенство

$$(11.3.19) \qquad \begin{aligned} \overline{w} &= \overline{f} \circ \overline{v} = \overline{f} \circ (e_{V_1}{}_*{}^* v) \\ &= (\overline{f} \circ e_{V_1})^*{}_* v \end{aligned}$$

является следствием равенств (11.3.4), (11.3.5). $V_2$-число $\overline{f} \circ e_{V_1}^i$ имеет разложение

$$(11.3.20) \quad \overline{f} \circ e_{V_1}^i = e_{V_2 *}{}^* f^i = e_{V_2 j} f_j^i$$

относительно базиса $\overline{\overline{e}}_{V_2}$. Равенство

$$(11.3.21) \qquad \overline{w} = e_{V_2}{}_*{}^* f^*{}_* v$$

является следствием равенств (11.3.19), (11.3.20). Равенство (11.3.9) следует из сравнения (11.3.18) и (11.3.21) и теоремы 8.4.25. Из равенства (11.3.20) и теоремы 8.4.25 следует что матрица $f$ определена однозначно. $\square$

---

**Теорема 11.3.8.** *Пусть*

$$f = (f_j^i, i \in I, j \in J)$$

*матрица $A$-чисел.*     *Отображение* [11.8]

$$\boxed{(11.3.3) \quad \overline{f} : V_1 \to V_2}$$

*определённое равенством*

$$\boxed{(11.3.8) \quad \overline{f} \circ (e_{V_1 *}{}^* v) = e_{V_2 *}{}^* f_*{}^* v}$$

*является гомоморфизмом правого $A$-векторного пространства столбцов $V_1$ в правое $A$-векторное пространство столбцов $V_2$. Гомоморфизм*

$$\boxed{(11.3.3) \quad \overline{f} : V_1 \to V_2}$$

*который имеет данную матрицу $f$, определён однозначно.*

---

**Теорема 11.3.9.** *Пусть*

$$f = (f_i^j, i \in I, j \in J)$$

*матрица $A$-чисел.*     *Отображение* [11.9]

$$\boxed{(11.3.3) \quad \overline{f} : V_1 \to V_2}$$

*определённое равенством*

$$\boxed{(11.3.11) \quad \overline{f} \circ (e_{V_1}{}_*{}^* v) = e_{V_2 *}{}^* f^*{}_* v}$$

*является гомоморфизмом правого $A$-векторного пространства строк $V_1$ в правое $A$-векторное пространство строк $V_2$. Гомоморфизм*

$$\boxed{(11.3.3) \quad \overline{f} : V_1 \to V_2}$$

*который имеет данную матрицу $f$, определён однозначно.*

Доказательство. Равенство

$$(11.3.22) \quad \begin{aligned} &\overline{f} \circ (e_{V_1 *}{}^*(v + w)) \\ &= e_{V_2 *}{}^* f_*{}^*(v + w) \\ &= e_{V_2 *}{}^* f_*{}^* v + e_{V_2 *}{}^* f_*{}^* w \\ &= \overline{f} \circ (e_{V_1 *}{}^* v) + \overline{f} \circ (e_{V_1 *}{}^* w) \end{aligned}$$

является следствием равенств

$(8.1.16)$   $a_*(b_1 + b_2) = a_* b_1 + a_* b_2$

$(11.2.8)$   $f \circ (va) = (f \circ v)a$

$(11.3.8)$   $\overline{f} \circ (e_{V_1 *}{}^* v) = e_{V_2 *}{}^* f_*{}^* v$

Из равенства $(11.3.22)$ следует, что отображение $\overline{f}$ является гомоморфизмом абелевой группы. Равенство

$$(11.3.23) \quad \begin{aligned} &\overline{f} \circ (e_{V_1 *}{}^*(va)) \\ &= e_{V_2 *}{}^* f_*{}^*(va) \\ &= (e_{V_2 *}{}^* f_*{}^* v)a \\ &= (\overline{f} \circ (e_{V_1 *}{}^* v))a \end{aligned}$$

является следствием равенств

$(11.2.8)$   $f \circ (va) = (f \circ v)a$

$(11.3.8)$   $\overline{f} \circ (e_{V_1 *}{}^* v) = e_{V_2 *}{}^* f_*{}^* v$

Из равенства $(11.3.23)$, и определений $2.8.5$, $11.3.1$ следует, что отображение

$(11.3.3)$   $\overline{f} : V_1 \to V_2$

является гомоморфизмом $A$-векторного пространства столбцов $V_1$ в $A$-векторное пространство столбцов $V_2$.

Пусть $f$ - матрица гомоморфизмов $\overline{f}$, $\overline{g}$ относительно базисов $\overline{\overline{e}}_V$, $\overline{\overline{e}}_W$. Равенство

$$\overline{f} \circ (e_{V *}{}^* v) = e_{W *}{}^* f_*{}^* v = \overline{g} \circ (e_{V *}{}^* v)$$

является следствием теоремы $11.3.6$. Следовательно, $\overline{f} = \overline{g}$.   □

На основании теорем $11.3.6$, $11.3.8$ мы идентифицируем гомоморфизм

$(11.3.3)$   $\overline{f} : V_1 \to V_2$

правого $A$-векторного пространства $V$ столбцов и координаты его представления

$(11.3.8)$   $\overline{f} \circ (e_{V_1 *}{}^* v) = e_{V_2 *}{}^* f_*{}^* v$

---

Доказательство. Равенство

$$(11.3.24) \quad \begin{aligned} &\overline{f} \circ (e_{V_1}{}^*{}_*(v + w)) \\ &= e_{V_2}{}^*{}_* f^*{}_*(v + w) \\ &= e_{V_2}{}^*{}_* f^*{}_* v + e_{V_2}{}^*{}_* f^*{}_* w \\ &= \overline{f} \circ (e_{V_1}{}^*{}_* v) + \overline{f} \circ (e_{V_1}{}^*{}_* w) \end{aligned}$$

является следствием равенств

$(8.1.18)$   $a^*(b_1 + b_2) = a^* b_1 + a^* b_2$

$(11.2.8)$   $f \circ (va) = (f \circ v)a$

$(11.3.11)$   $\overline{f} \circ (e_{V_1}{}^*{}_* v) = e_{V_2}{}^*{}_* f^*{}_* v$

Из равенства $(11.3.24)$ следует, что отображение $\overline{f}$ является гомоморфизмом абелевой группы. Равенство

$$(11.3.25) \quad \begin{aligned} &\overline{f} \circ (e_{V_1}{}^*{}_*(va)) \\ &= e_{V_2}{}^*{}_* f^*{}_*(va) \\ &= (e_{V_2}{}^*{}_* f^*{}_* v)a \\ &= (\overline{f} \circ (e_{V_1}{}^*{}_* v))a \end{aligned}$$

является следствием равенств

$(11.2.8)$   $f \circ (va) = (f \circ v)a$

$(11.3.11)$   $\overline{f} \circ (e_{V_1}{}^*{}_* v) = e_{V_2}{}^*{}_* f^*{}_* v$

Из равенства $(11.3.25)$, и определений $2.8.5$, $11.3.1$ следует, что отображение

$(11.3.3)$   $\overline{f} : V_1 \to V_2$

является гомоморфизмом $A$-векторного пространства строк $V_1$ в $A$-векторное пространство строк $V_2$.

Пусть $f$ - матрица гомоморфизмов $\overline{f}$, $\overline{g}$ относительно базисов $\overline{\overline{e}}_V$, $\overline{\overline{e}}_W$. Равенство

$$\overline{f} \circ (e_{V}{}^*{}_* v) = e_{W}{}^*{}_* f^*{}_* v = \overline{g} \circ (e_{V}{}^*{}_* v)$$

является следствием теоремы $11.3.7$. Следовательно, $\overline{f} = \overline{g}$.   □

На основании теорем $11.3.7$, $11.3.9$ мы идентифицируем гомоморфизм

$(11.3.3)$   $\overline{f} : V_1 \to V_2$

правого $A$-векторного пространства $V$ строк и координаты его представления

$(11.3.11)$   $\overline{f} \circ (e_{V_1}{}^*{}_* v) = e_{V_2}{}^*{}_* f^*{}_* v$

## 11.4. Сумма гомоморфизмов $A$-векторного пространства

Теорема 11.4.1. *Пусть $V$, $W$ - правые $A$-векторные пространства и отображения $g : V \to W$, $h : V \to W$ - гомоморфизмы правого $A$-векторного пространства. Пусть отображение $f : V \to W$ определено равенством*

$$\text{(11.4.1)} \qquad \forall v \in V : f \circ v = g \circ v + h \circ v$$

*Отображение $f$ является гомоморфизмом правого $A$-векторного пространства и называется* **суммой**

$$f = g + h$$

**гомоморфизмов** $g$ *and* $h$.

Доказательство. Согласно теореме 7.4.6 (равенство (7.4.8)), равенство

$$
\begin{aligned}
\text{(11.4.2)} \qquad & f \circ (u + v) \\
& = g \circ (u + v) + h \circ (u + v) \\
& = g \circ u + g \circ v \\
& + h \circ u + h \circ u \\
& = f \circ u + f \circ v
\end{aligned}
$$

является следствием равенства

$$\boxed{\text{(11.3.4)} \quad f \circ (u + v) = f \circ u + f \circ v}$$

и равенства (11.4.1). Согласно теореме 7.4.6 (равенство (7.4.10)), равенство

$$
\begin{aligned}
\text{(11.4.3)} \qquad & f \circ (va) \\
& = g \circ (va) + h \circ (va) \\
& = (g \circ v)a + (h \circ v)a \\
& = (g \circ v + h \circ v)a \\
& = (f \circ v)a
\end{aligned}
$$

является следствием равенства

$$\boxed{\text{(11.3.5)} \quad f \circ (va) = (f \circ v)a}$$

и равенства (11.4.1).

Теорема является следствием теоремы 11.3.2 и равенств (11.4.2), (11.4.3). $\square$

---

Теорема 11.4.2. *Пусть $V$, $W$ - правые $A$-векторные пространства столбцов. Пусть $\overline{\overline{e}}_V$ базис правого $A$-векторного пространства $V$. Пусть $\overline{\overline{e}}_W$ базис правого $A$-векторного пространства $W$. Гомоморфизм $\overline{f} : V \to W$ является суммой гомоморфизмов $\overline{g} : V \to W$, $\overline{h} : V \to W$ тогда и только тогда, когда матрица $f$ гомоморфизма $\overline{f}$ относительно базисов $\overline{\overline{e}}_V$, $\overline{\overline{e}}_W$ равна сумме матрицы $g$ гомоморфизма $\overline{g}$ относительно базисов*

Теорема 11.4.3. *Пусть $V$, $W$ - правые $A$-векторные пространства строк. Пусть $\overline{\overline{e}}_V$ базис правого $A$-векторного пространства $V$. Пусть $\overline{\overline{e}}_W$ базис правого $A$-векторного пространства $W$. Гомоморфизм $\overline{\overline{f}} : V \to W$ является суммой гомоморфизмов $\overline{g} : V \to W$, $\overline{h} : V \to W$ тогда и только тогда, когда матрица $f$ гомоморфизма $\overline{\overline{f}}$ относительно базисов $\overline{\overline{e}}_V$, $\overline{\overline{e}}_W$ равна сумме матрицы $g$ гомоморфизма $\overline{g}$ относительно базисов $\overline{\overline{e}}_V$, $\overline{\overline{e}}_W$* $\overline{\overline{e}}_V,\quad \overline{\overline{e}}_W$   и матрицы $h$ гомоморфизма $\overline{h}$ относительно базисов   $\overline{\overline{e}}_V,\quad \overline{\overline{e}}_W$.

и матрицы $h$ гомоморфизма $\overline{h}$ относительно базисов   $\overline{\overline{e}}_V,\quad \overline{\overline{e}}_W$.

**Доказательство.** Согласно теореме 11.3.6

$$(11.4.4)\qquad \overline{f}\circ(e_{U*}{}^*a)=e_{V*}{}^*f_*{}^*a$$

$$(11.4.5)\qquad \overline{g}\circ(e_{U*}{}^*a)=e_{V*}{}^*g_*{}^*a$$

$$(11.4.6)\qquad \overline{h}\circ(e_{U*}{}^*a)=e_{V*}{}^*h_*{}^*a$$

Согласно теореме 11.4.1, равенство

$$(11.4.7)\quad \begin{aligned}&\overline{f}\circ(e_{U*}{}^*v)\\ =&\overline{g}\circ(e_{U*}{}^*v)+\overline{h}\circ(e_{U*}{}^*v)\\ =&e_{V*}{}^*g_*{}^*v+e_{V*}{}^*h_*{}^*v\\ =&e_{V*}{}^*(g+h)_*{}^*v\end{aligned}$$

является следствием равенств (8.1.16), (8.1.17), (11.4.5), (11.4.6). Равенство

$$f=g+h$$

является следствием равенств (11.4.4), (11.4.7) и теоремы 11.3.6.

   Пусть

$$f=g+h$$

Согласно теореме 11.3.8, существуют гомоморфизмы   $\overline{f}:U\to V,\quad \overline{g}:U\to V,\quad \overline{h}:U\to V$ такие, что

$$(11.4.8)\quad \begin{aligned}&\overline{f}\circ(e_{U*}{}^*v)\\ =&e_{V*}{}^*f_*{}^*v\\ =&e_{V*}{}^*(g+h)_*{}^*v\\ =&e_{V*}{}^*g_*{}^*v+e_{V*}{}^*h_*{}^*v\\ =&\overline{g}\circ(e_{U*}{}^*v)+\overline{h}\circ(e_{U*}{}^*v)\end{aligned}$$

Из равенства (11.4.8) и теоремы 11.4.1 следует, что гомоморфизм $\overline{f}$ является суммой гомоморфизмов $\overline{g},\ \overline{h}$.   $\square$

**Доказательство.** Согласно теореме 11.3.7

$$(11.4.9)\qquad \overline{f}\circ(e_U{}^*{}_*a)=e_V{}^*{}_*f^*{}_*a$$

$$(11.4.10)\qquad \overline{g}\circ(e_U{}^*{}_*a)=e_V{}^*{}_*g^*{}_*a$$

$$(11.4.11)\qquad \overline{h}\circ(e_U{}^*{}_*a)=e_V{}^*{}_*h^*{}_*a$$

Согласно теореме 11.4.1, равенство

$$(11.4.12)\quad \begin{aligned}&\overline{f}\circ(e_U{}^*{}_*v)\\ =&\overline{g}\circ(e_U{}^*{}_*v)+\overline{h}\circ(e_U{}^*{}_*v)\\ =&e_V{}^*{}_*g^*{}_*v+e_V{}^*{}_*h^*{}_*v\\ =&e_V{}^*{}_*(g+h)^*{}_*v\end{aligned}$$

является следствием равенств (8.1.18), (8.1.19), (11.4.10), (11.4.11). Равенство

$$f=g+h$$

является следствием равенств (11.4.9), (11.4.12) и теоремы 11.3.7.

   Пусть

$$f=g+h$$

Согласно теореме 11.3.9, существуют гомоморфизмы   $\overline{f}:U\to V,\quad \overline{g}:U\to V,\quad \overline{h}:U\to V$ такие, что

$$(11.4.13)\quad \begin{aligned}&\overline{f}\circ(e_U{}^*{}_*v)\\ =&e_V{}^*{}_*f^*{}_*v\\ =&e_V{}^*{}_*(g+h)^*{}_*v\\ =&e_V{}^*{}_*g^*{}_*v+e_V{}^*{}_*h^*{}_*v\\ =&\overline{g}\circ(e_U{}^*{}_*v)+\overline{h}\circ(e_U{}^*{}_*v)\end{aligned}$$

Из равенства (11.4.13) и теоремы 11.4.1 следует, что гомоморфизм $\overline{f}$ является суммой гомоморфизмов $\overline{g},\ \overline{h}$.   $\square$

## 11.5. Произведение гомоморфизмов $A$-векторного пространства

**Теорема 11.5.1.** *Пусть $U,\ V,\ W$ - правые $A$-векторные пространства. Предположим, что мы имеем коммутативную диаграмму отображений*

$$(11.5.1)\qquad \begin{array}{ccc}U & \xrightarrow{\ f\ } & W\\ & \searrow{\scriptstyle g}\quad\nearrow{\scriptstyle h} & \\ & V & \end{array}$$

$$(11.5.2) \qquad \forall u \in U : f \circ u = (h \circ g) \circ u$$
$$= h \circ (g \circ u)$$

*где отображения $g$, $h$ являются гомоморфизмами правого $A$-векторного пространства. Отображение $f = h \circ g$ является гомоморфизмом правого $A$-векторного пространства и называется* **произведением гомоморфизмов** *$h$, $g$.*

Доказательство. Равенство (11.5.2) следует из коммутативности диаграммы (11.5.1). Равенство

$$(11.5.3) \qquad \begin{aligned} &f \circ (u + v) \\ &= h \circ (g \circ (u + v)) \\ &= h \circ (g \circ u + g \circ v)) \\ &= h \circ (g \circ u) + h \circ (g \circ v) \\ &= f \circ u + f \circ v \end{aligned}$$

является следствием равенства

$$\boxed{(11.3.4) \quad f \circ (u + v) = f \circ u + f \circ v}$$

и равенства (11.5.2). Равенство

$$(11.5.4) \qquad \begin{aligned} f \circ (va) &= h \circ (g \circ (va)) \\ &= h \circ ((g \circ v)a) \\ &= (h \circ (g \circ v))a \\ &= (f \circ v)a \end{aligned}$$

является следствием равенства

$$\boxed{(11.3.5) \quad f \circ (va) = (f \circ v)a}$$

и равенства (11.4.1). Теорема является следствием теоремы 11.3.2 и равенств (11.5.3), (11.5.4). $\qquad\square$

Теорема 11.5.2. *Пусть $U$, $V$, $W$ - правые $A$-векторные пространства столбцов. Пусть $\overline{\overline{e}}_U$ базис правого $A$-векторного пространства $U$. Пусть $\overline{\overline{e}}_V$ базис правого $A$-векторного пространства $V$. Пусть $\overline{\overline{e}}_W$ базис правого $A$-векторного пространства $W$. Гомоморфизм $\overline{f} : U \to W$ является произведением гомоморфизмов $\overline{h} : V \to W$, $\overline{g} : U \to V$ тогда и только тогда, когда матрица $f$ гомоморфизма $\overline{f}$ относительно базисов $\overline{\overline{e}}_U$, $\overline{\overline{e}}_W$ равна $_*$\*-произведению матрицы $g$ гомоморфизма $\overline{g}$ относительно базисов $\overline{\overline{e}}_U$, $\overline{\overline{e}}_V$ на матрицу $h$ гомоморфизма*

Теорема 11.5.3. *Пусть $U$, $V$, $W$ - правые $A$-векторные пространства строк. Пусть $\overline{\overline{e}}_U$ базис правого $A$-векторного пространства $U$. Пусть $\overline{\overline{e}}_V$ базис правого $A$-векторного пространства $V$. Пусть $\overline{\overline{e}}_W$ базис правого $A$-векторного пространства $W$. Гомоморфизм $\overline{f} : U \to W$ является произведением гомоморфизмов $\overline{h} : V \to W$, $\overline{g} : U \to V$ тогда и только тогда, когда матрица $f$ гомоморфизма $\overline{f}$ относительно базисов $\overline{\overline{e}}_U$, $\overline{\overline{e}}_W$ равна $^*_*$-произведению матрицы $h$ гомоморфизма $\overline{h}$ относительно базисов $\overline{\overline{e}}_V$, $\overline{\overline{e}}_W$ на матрицу $g$ гомоморфизма $\overline{g}$*

$\overline{h}$ *относительно базисов*  $\overline{\overline{e}}_V$, $\overline{\overline{e}}_W$

(11.5.5) $\qquad f = g_*{}^*h$

*относительно базисов*  $\overline{\overline{e}}_U$, $\overline{\overline{e}}_V$

(11.5.6) $\qquad f = h^*{}_*g$

**Доказательство.** Согласно теореме 11.3.6

(11.5.7) $\qquad \overline{f} \circ (e_{U*}{}^*u) = e_{W*}{}^*f_*{}^*u$

(11.5.8) $\qquad \overline{g} \circ (e_{U*}{}^*u) = e_{V*}{}^*g_*{}^*u$

(11.5.9) $\qquad \overline{h} \circ (e_{V*}{}^*v) = e_{W*}{}^*h_*{}^*v$

Согласно теореме 11.5.1, равенство

$$
\begin{aligned}
(11.5.10) \quad & \overline{f} \circ (e_{U*}{}^*u) \\
& = \overline{h} \circ (\overline{g} \circ (e_{U*}{}^*u)) \\
& = \overline{h} \circ (u_* g_* e_V) \\
& = e_{W*}{}^*h_*{}^*g_*{}^*u
\end{aligned}
$$

является следствием равенств (11.5.8), (11.5.9). Равенство (11.5.5) является следствием равенств (11.5.7), (11.5.10) и теоремы 11.3.6.

Пусть $f = g_*{}^*h$. Согласно теореме 11.3.8, существуют гомоморфизмы $\overline{f}: U \to V$, $\overline{g}: U \to V$, $\overline{h}: U \to V$ такие, что

$$
\begin{aligned}
(11.5.11) \quad & \overline{f} \circ (e_{U*}{}^*u) = e_{W*}{}^*f_*{}^*u \\
& = e_{W*}{}^*h_*{}^*g_*{}^*u \\
& = \overline{h} \circ (e_{V*}{}^*g_*{}^*u) \\
& = \overline{h} \circ (\overline{g} \circ (e_{U*}{}^*u)) \\
& = (\overline{h} \circ \overline{g}) \circ (e_{U*}{}^*u)
\end{aligned}
$$

Из равенства (11.5.11) и теоремы 11.4.1 следует, что гомоморфизм $\overline{f}$ является произведением гомоморфизмов $\overline{h}$, $\overline{g}$. □

**Доказательство.** Согласно теореме 11.3.7

(11.5.12) $\qquad \overline{f} \circ (e_U{}^*{}_*u) = e_W{}^*{}_*f^*{}_*u$

(11.5.13) $\qquad \overline{g} \circ (e_U{}^*{}_*u) = e_V{}^*{}_*g^*{}_*u$

(11.5.14) $\qquad \overline{h} \circ (e_V{}^*{}_*v) = e_W{}^*{}_*h^*{}_*v$

Согласно теореме 11.5.1, равенство

$$
\begin{aligned}
(11.5.15) \quad & \overline{f} \circ (e_U{}^*{}_*u) \\
& = \overline{h} \circ (\overline{g} \circ (e_U{}^*{}_*u)) \\
& = \overline{h} \circ (u^*{}_*g^*{}_*e_V) \\
& = u^*{}_*g^*{}_*h^*{}_*e_W
\end{aligned}
$$

является следствием равенств (11.5.13), (11.5.14). Равенство (11.5.6) является следствием равенств (11.5.12), (11.5.15) и теоремы 11.3.7.

Пусть $f = h^*{}_*g$. Согласно теореме 11.3.9, существуют гомоморфизмы $\overline{f}: U \to V$, $\overline{g}: U \to V$, $\overline{h}: U \to V$ такие, что

$$
\begin{aligned}
(11.5.16) \quad & \overline{f} \circ (e_U{}^*{}_*u) = e_W{}^*{}_*f^*{}_*u \\
& = e_W{}^*{}_*h^*{}_*g^*{}_*u \\
& = \overline{h} \circ (e_V{}^*{}_*g^*{}_*u) \\
& = \overline{h} \circ (\overline{g} \circ (e_U{}^*{}_*u)) \\
& = (\overline{h} \circ \overline{g}) \circ (e_U{}^*{}_*u)
\end{aligned}
$$

Из равенства (11.5.16) и теоремы 11.4.1 следует, что гомоморфизм $\overline{f}$ является произведением гомоморфизмов $\overline{h}$, $\overline{g}$. □

## 11.6. Прямая сумма левых $A$-векторных пространств

Пусть Right $A$-Vector - категория правых $A$-векторных пространств, морфизмы которой - это гомоморфизмы $A$-модуля.

Определение 11.6.1. *Копроизведение в категории* Right $A$-Vector *называется* **прямой суммой**.[11.10] *Мы будем пользоваться записью* $V \oplus W$ *для прямой суммы правых $A$-векторных пространств $V$ и $W$.* □

---

11.10 Смотри так же определение прямой суммы модулей в [1], страница 98. Согласно теореме 1 на той же странице, прямая сумма модулей существует.

**Теорема 11.6.2.** *В категории* Right $A$-Vector *существует прямая сумма правых $A$-векторных пространств.*[11.11] *Пусть* $\{V_i, i \in I\}$ *- семейство правых $A$-векторных пространств. Тогда представление*

$$A \xrightarrow{\ast} \bigoplus_{i \in I} V_i \qquad \left(\bigoplus_{i \in I} v_i\right) a = \bigoplus_{i \in I} v_i a$$

*$D$-алгебры $A$ в прямой сумме абелевых групп*

$$V = \bigoplus_{i \in I} V_i$$

*является прямой суммой правых $A$-векторных пространств*

$$V = \bigoplus_{i \in I} V_i$$

**Доказательство.** Пусть $V_i$, $i \in I$, - множество правых $A$-векторных пространств. Пусть

$$V = \bigoplus_{i \in I} V_i$$

прямая сумма абелевых групп. Пусть $v = (v_i, i \in I) \in V$. Мы определим действие $A$-числа $a$ на $V$-число $v$ согласно равенству

$$(11.6.1) \qquad a(v_i, i \in I) = (av_i, i \in I)$$

Согласно теореме 7.4.6,

$$(11.6.2) \qquad v_i(ba) = (v_i b)a$$

$$(11.6.3) \qquad (v_i + w_i)a = v_i a + w_i a$$

$$(11.6.4) \qquad v_i(a + b) = v_i a + v_i b$$

$$(11.6.5) \qquad v_i 1 = v_i$$

для любых $a, b \in A$, $v_i, w_i \in V_i$, $i \in I$. Равенства

$$v(ba) = (v_i, i \in I)(ba) = (v_i(ba), i \in I)$$
$$(11.6.6) \qquad = ((v_i b)a, i \in I) = (v_i b, i \in I)a = ((v_i, i \in I)b)a$$
$$= (vb)a$$

$$(v + w)a = ((v_i, i \in I) + (w_i, i \in I))a = (v_i + w_i, i \in I)a$$
$$(11.6.7) \qquad = ((v_i + w_i)a, i \in I) = (v_i a + w_i a, i \in I)$$
$$= (v_i a, i \in I) + (w_i a, i \in I) = (v_i, i \in I)a + (w_i, i \in I)a$$
$$= va + wa$$

$$v(a + b) = (v_i, i \in I)(a + b) = (v_i(a + b), i \in I)$$
$$(11.6.8) \qquad = (v_i a + v_i b, i \in I) = (v_i a, i \in I) + (v_i b, i \in I)$$
$$= (v_i, i \in I)a + (v_i, i \in I)b = va + vb$$

$$v1 = (v_i, i \in I)1 = (v_i 1, i \in I) = (v_i, i \in I) = v$$

---

[11.11]   Смотри также предложение [1]-1 на странице 98.

являются следствием равенств (3.7.4), (11.6.1), (11.6.2), (11.6.3), (11.6.4), (11.6.5) для любых $a$, $b \in A$, $v = (v_i, i \in I)$, $w = (w_i, i \in I) \in V$.

Из равенства (11.6.7) следует, что отображение

$$a : v \to va$$

определённое равенством (11.6.1), является эндоморфизмом абелевой группы $V$. Из равенств (11.6.6), (11.6.8) следует, что отображение (11.6.1) является гомоморфизмом $D$-алгебры $A$ в множество эндоморфизмов абелевой группы $V$. Следовательно, отображение (11.6.1) является представлением $D$-алгебры $A$ в абелевой группе $V$.

Пусть $a_1$, $a_2$ - такие $A$-числа, что

(11.6.9)          $$a_1(v_i, i \in I) = a_2(v_i, i \in I)$$

для произвольного $V$-числа $v = (v_i, i \in I)$. Равенство

(11.6.10)          $$(a_1 v_i, i \in I) = (a_2 v_i, i \in I)$$

является следствием равенства (11.6.9). Из равенства (11.6.10) следует, что

(11.6.11)          $$\forall i \in I : a_1 v_i = a_2 v_i$$

Поскольку для любого $i \in I$, $V_i$ - $A$-векторное пространство, то соответствующее представление эффективно согласно определению 7.4.1 Следовательно, равенство $a_1 = a_2$ является следствием утверждения (11.6.11) и отображение (11.6.1) является эффективным представлением $D$-алгебры $A$ в абелевой группе $V$. Следовательно, абелевая группе $V$ является правым $A$-векторным пространством.

Пусть

$$\{f_i : V_i \to W, i \in I\}$$

семейство линейных отображений в правое $A$-векторное пространство $W$. Мы определим отображение

$$f : V \to W$$

пользуясь равенством

(11.6.12)          $$f\left(\bigoplus_{i \in I} x_i\right) = \sum_{i \in I} f_i(x_i)$$

Сумма в правой части равенства (11.6.12) конечна, так как все слагаемые, кроме конечного числа, равны 0. Из равенства

$$f_i(x_i + y_i) = f_i(x_i) + f_i(y_i)$$

и равенства (11.6.12) следует, что

$$\begin{aligned}
f(\bigoplus_{i \in I}(x_i + y_i)) &= \sum_{i \in I} f_i(x_i + y_i) \\
&= \sum_{i \in I}(f_i(x_i) + f_i(y_i)) \\
&= \sum_{i \in I} f_i(x_i) + \sum_{i \in I} f_i(y_i) \\
&= f(\bigoplus_{i \in I} x_i) + f(\bigoplus_{i \in I} y_i)
\end{aligned}$$

Из равенства
$$f_i(dx_i) = df_i(x_i)$$
и равенства (11.6.12) следует, что
$$f((dx_i, i \in I)) = \sum_{i \in I} f_i(dx_i) = \sum_{i \in I} df_i(x_i) = d \sum_{i \in I} f_i(x_i)$$
$$= df((x_i, i \in I))$$
Следовательно, отображение $f$ является линейным отображением. Равенство
$$f \circ \lambda_j(x) = \sum_{i \in I} f_i(\delta_j^i x) = f_j(x)$$
является следствием равенств (3.7.5), (11.6.12). Так как отображение $\lambda_i$ инъективно, то отображение $f$ определено однозначно. Следовательно, теорема является следствием определений 2.2.12, 3.7.2.    $\square$

ОПРЕДЕЛЕНИЕ 11.6.3. *Пусть $D$ - коммутативное кольцо с единицей. Пусть $A$ - ассоциативная $D$-алгебра с делением. Для любого целого $n > 0$ мы определим правое $A$-векторное пространство $\mathcal{R}^n A$ используя определение по индукции*
$$\mathcal{R}^1 A = \mathcal{R} A$$
$$\mathcal{R}^k A = \mathcal{R}^{k-1} A \oplus \mathcal{R} A$$
$\square$

ОПРЕДЕЛЕНИЕ 11.6.4. *Пусть $D$ - коммутативное кольцо с единицей. Пусть $A$ - ассоциативная $D$-алгебра с делением. Для любого множества $B$ и любого целого $n > 0$ мы определим правое $A$-векторное пространство $\mathcal{R}^n A^B$ используя определение по индукции*
$$\mathcal{R}^1 A^B = \mathcal{R} A^B$$
$$\mathcal{R}^k A^B = \mathcal{R}^{k-1} A^B \oplus \mathcal{R} A^B$$
$\square$

ТЕОРЕМА 11.6.5. *Прямая сумма правых $A$-векторных пространств $V_1$, ..., $V_n$ совпадает с их прямым произведением*
$$V_1 \oplus ... \oplus V_n = V_1 \times ... \times V_n$$

Доказательство. Теорема является следствием теоремы 11.6.2.    $\square$

ТЕОРЕМА 11.6.6. *Пусть $\overline{\overline{e}}_1 = (e_{1i}, i \in I)$ базис правого $A$-векторного пространства столбцов $V$. Тогда мы можем представить $A$-модуль $V$ как прямую сумму*

(11.6.13) $$V = \bigoplus_{i \in I} e_i A$$

ТЕОРЕМА 11.6.7. *Пусть $\overline{\overline{e}}_1 = (e_1^i, i \in I)$ базис правого $A$-векторного пространства строк $V$. Тогда мы можем представить $A$-модуль $V$ как прямую сумму*

(11.6.14) $$V = \bigoplus_{i \in I} e^i A$$

Доказательство. Отображение

$$f : e_i v^i \in V \to (v^i, i \in I) \in \bigoplus_{i \in I} e_i A$$

является изоморфизмом.                     $\square$

Доказательство. Отображение

$$f : e^i v_i \in V \to (v_i, i \in I) \in \bigoplus_{i \in I} e^i A$$

является изоморфизмом.                     $\square$

# Многообразие базисов $A$-векторного пространства

## 12.1. Размерность левого $A$-векторного пространства

**Теорема 12.1.1.** *Пусть $V$ - левое $A$-векторное пространство столбцов. Пусть множество векторов $\overline{\overline{e}} = (e_i, i \in I)$ является базисом левого $A$-векторного пространства $V$. Пусть множество векторов $\overline{\overline{e}} = (e_j, j \in J)$ является базисом левого $A$-векторного пространства $V$. Если $|I|$ и $|J|$ - конечные числа, то $|I| = |J|$.*

**Теорема 12.1.2.** *Пусть $V$ - левое $A$-векторное пространство строк. Пусть множество векторов $\overline{\overline{e}} = (e^i, i \in I)$ является базисом левого $A$-векторного пространства $V$. Пусть множество векторов $\overline{\overline{e}} = (e^j, j \in J)$ является базисом левого $A$-векторного пространства $V$. Если $|I|$ и $|J|$ - конечные числа, то $|I| = |J|$.*

**Доказательство.** Предположим, что $|I| = m$ и $|J| = n$. Предположим, что

(12.1.1) $\qquad m < n$

Так как $\overline{\overline{e}}_1$ - базис, любой вектор $e_{2j}$, $j \in J$, имеет разложение

$$e_{2j} = a_j{}^*{}_* e_1$$

Так как $\overline{\overline{e}}_2$ - базис,

(12.1.2) $\qquad \lambda = 0$

следует из

(12.1.3) $\qquad \lambda^*{}_* e_2 = \lambda^*{}_* a^*{}_* e_1 = 0$

Так как $\overline{\overline{e}}_1$ - базис, мы получим

(12.1.4) $\qquad \lambda^*{}_* a = 0$

Согласно (12.1.1), $\mathrm{rank}(^*{}_*)a \leq m$ и система линейных уравнений (12.1.4) имеет больше переменных, чем уравнений. Согласно теореме 9.4.4, $\lambda \neq 0$, что противоречит утверждению (12.1.2). Следовательно, утверждение $m < n$ неверно. Таким же образом мы можем доказать, что утверждение $n < m$ неверно. Это завершает доказательство теоремы. □

**Доказательство.** Предположим, что $|I| = m$ и $|J| = n$. Предположим, что

(12.1.5) $\qquad m < n$

Так как $\overline{\overline{e}}_1$ - базис, любой вектор $e^j_2$, $j \in J$, имеет разложение

$$e^j_2 = a^j{}_* {}^* e_1$$

Так как $\overline{\overline{e}}_2$ - базис,

(12.1.6) $\qquad \lambda = 0$

следует из

(12.1.7) $\qquad \lambda_*{}^* e_2 = \lambda_*{}^* a_*{}^* e_1 = 0$

Так как $\overline{\overline{e}}_1$ - базис, мы получим

(12.1.8) $\qquad \lambda_*{}^* a = 0$

Согласно (12.1.5), $\mathrm{rank}(_*{}^*)a \leq m$ и система линейных уравнений (12.1.8) имеет больше переменных, чем уравнений. Согласно теореме 9.4.4, $\lambda \neq 0$, что противоречит утверждению (12.1.6). Следовательно, утверждение $m < n$ неверно. Таким же образом мы можем доказать, что утверждение $n < m$ неверно. Это завершает доказательство теоремы. □

**Определение 12.1.3.** *Мы будем называть* **размерностью левого $A$-векторного пространства** *число векторов в базисе.* □

Теорема 12.1.4. *Если левое $A$-векторное пространство $V$ имеет конечную размерность, то из утверждения*

$$(12.1.9) \quad \forall v \in V, av = bv \quad a, b \in A$$

*следует, что $a = b$.*

Доказательство. Пусть множество векторов $\overline{\overline{e}} = (e_j, j \in J)$ является базисом левого $A$-векторного пространства $V$. Согласно теоремам 7.2.6, 8.4.4, равенство

$$(12.1.11) \qquad av^i = bv^i$$

является следствием равенства (12.1.9). Теорема следует из равенства (12.1.11), если мы положим $v = e_I$. $\square$

---

Теорема 12.1.5. *Если левое $A$-векторное пространство $V$ имеет конечную размерность, то из утверждения*

$$(12.1.10) \quad \forall v \in V, av = bv \quad a, b \in A$$

*следует, что $a = b$.*

Доказательство. Пусть множество векторов $\overline{\overline{e}} = (e_j, j \in J)$ является базисом левого $A$-векторного пространства $V$. Согласно теоремам 7.2.6, 8.4.11, равенство

$$(12.1.12) \qquad av_i = bv_i$$

является следствием равенства (12.1.10). Теорема следует из равенства (12.1.12), если мы положим $v = e_I$. $\square$

---

Теорема 12.1.6. *Если левое $A$-векторное пространство $V$ имеет конечную размерность, то соответствующее представление свободно.*

Доказательство. Теорема является следствием определения 2.3.3 и теорем 12.1.4, 12.1.5. $\square$

---

Теорема 12.1.7. *Пусть $V$ - левое $A$-векторное пространство столбцов. Координатная матрица базиса $\overline{\overline{g}}$ относительно базиса $\overline{\overline{e}}$ левого $A$-векторного пространства $V$ является $*_*$-невырожденной матрицей.*

Доказательство. Согласно теореме 9.3.5, $*_*$-ранг координатной матрицы $\overline{\overline{g}}$ относительно базиса $\overline{\overline{e}}$ равен размерности левого $A$-пространства столбцов, откуда следует утверждение теоремы. $\square$

---

Теорема 12.1.8. *Пусть $V$ - левое $A$-векторное пространство строк. Координатная матрица базиса $\overline{\overline{g}}$ относительно базиса $\overline{\overline{e}}$ левого $A$-векторного пространства $V$ является $_*^*$-невырожденной матрицей.*

Доказательство. Согласно теореме 9.3.8, $_*^*$-ранг координатной матрицы базиса $\overline{\overline{g}}$ относительно базиса $\overline{\overline{e}}$ равен размерности левого $A$-пространства столбцов, откуда следует утверждение теоремы. $\square$

---

Теорема 12.1.9. *Пусть $V$ - левое $A$-векторное пространство столбцов. Пусть $\overline{\overline{e}}$ - базис левого $A$-векторного пространства $V$. Тогда любой автоморфизм $\overline{f}$ левого $A$-векторного пространства $V$ имеет вид*

$$(12.1.13) \qquad v' = v *_* f$$

*где $f$ - $*_*$-невырожденная матрица.*

---

Теорема 12.1.10. *Пусть $V$ - левое $A$-векторное пространство строк. Пусть $\overline{\overline{e}}$ - базис левого $A$-векторного пространства $V$. Тогда любой автоморфизм $\overline{f}$ левого $A$-векторного пространства $V$ имеет вид*

$$(12.1.14) \qquad v' = v_* ^* f$$

где $f$ - $_*{}^*$-невырожденная матрица.

**Доказательство.** Равенство (12.1.13) следует из теоремы 10.3.6. Так как $\overline{f}$ - изоморфизм, то для каждого вектора $\overline{v}'$ существует единственный вектор $\overline{v}$ такой, что $\overline{v}' = \overline{f} \circ \overline{v}$. Следовательно, система линейных уравнений (12.1.13) имеет единственное решение. Согласно следствию 9.4.5 матрица $f$ $_*{}^*$-невырожденна. □

**Доказательство.** Равенство (12.1.14) следует из теоремы 10.3.7. Так как $\overline{f}$ - изоморфизм, то для каждого вектора $\overline{v}'$ существует единственный вектор $\overline{v}$ такой, что $\overline{v}' = \overline{f} \circ \overline{v}$. Следовательно, система линейных уравнений (12.1.14) имеет единственное решение. Согласно следствию 9.4.5 матрица $f$ $_*{}^*$-невырожденна. □

**Теорема 12.1.11.** *Матрицы автоморфизмов левого $A$-пространства столбцов порождают группу $GL(n_*{}^*A)$.*

**Теорема 12.1.12.** *Матрицы автоморфизмов левого $A$-пространства строк порождают группу $GL(n_*{}^*A)$.*

**Доказательство.** Если даны два автоморфизма $\overline{f}$ и $\overline{g}$, то мы можем записать

$$v' = v_*{}^* f$$

$$v'' = v'_*{}^* g = v_*{}^* f_*{}^* g$$

Следовательно, результирующий автоморфизм имеет матрицу $f_*{}^* g$. □

**Доказательство.** Если даны два автоморфизма $\overline{f}$ и $\overline{g}$, то мы можем записать

$$v' = v_*{}^* f$$

$$v'' = v'_*{}^* g = v_*{}^* f_*{}^* g$$

Следовательно, результирующий автоморфизм имеет матрицу $f_*{}^* g$. □

**Теорема 12.1.13.** *Пусть $\overline{\overline{e}}_{V1}$, $\overline{\overline{e}}_{V2}$ - базисы левого $A$-векторного пространства $V$ столбцов. Пусть $\overline{\overline{e}}_{W1}$, $\overline{\overline{e}}_{W2}$ - базисы левого $A$-векторного пространства $W$ столбцов. Пусть $a_1$ - матрица гомоморфизма*

$$(12.1.15) \qquad a : V \to W$$

*относительно базисов $\overline{\overline{e}}_{V1}$, $\overline{\overline{e}}_{W1}$ и $a_2$ - матрица гомоморфизма (12.1.15) относительно базисов $\overline{\overline{e}}_{V2}$, $\overline{\overline{e}}_{W2}$. Если базис $\overline{\overline{e}}_{V1}$ имеет координатную матрицу $b$ относительно базиса $\overline{\overline{e}}_{V2}$*

$$\overline{\overline{e}}_{V1} = b_*{}^* \overline{\overline{e}}_{V2}$$

*и $\overline{\overline{e}}_{W1}$ имеет координатную матрицу $c$ относительно базиса $\overline{\overline{e}}_{W2}$*

$$(12.1.16) \qquad \overline{\overline{e}}_{W1} = c_*{}^* \overline{\overline{e}}_{W2}$$

*то матрицы $a_1$ и $a_2$ связаны соотношением*

$$a_1 = b_*{}^* a_2{}_*{}^* c^{-1}{}_*{}^*$$

**Теорема 12.1.14.** *Пусть $\overline{\overline{e}}_{V1}$, $\overline{\overline{e}}_{V2}$ - базисы левого $A$-векторного пространства $V$ столбцов. Пусть $\overline{\overline{e}}_{W1}$, $\overline{\overline{e}}_{W2}$ - базисы левого $A$-векторного пространства $W$ столбцов. Пусть $a_1$ - матрица гомоморфизма*

$$(12.1.17) \qquad a : V \to W$$

*относительно базисов $\overline{\overline{e}}_{V1}$, $\overline{\overline{e}}_{W1}$ и $a_2$ - матрица гомоморфизма (12.1.17) относительно базисов $\overline{\overline{e}}_{V2}$, $\overline{\overline{e}}_{W2}$. Если базис $\overline{\overline{e}}_{V1}$ имеет координатную матрицу $b$ относительно базиса $\overline{\overline{e}}_{V2}$*

$$\overline{\overline{e}}_{V1} = b_*{}^* \overline{\overline{e}}_{V2}$$

*и $\overline{\overline{e}}_{W1}$ имеет координатную матрицу $c$ относительно базиса $\overline{\overline{e}}_{W2}$*

$$(12.1.18) \qquad \overline{\overline{e}}_{W1} = c_*{}^* \overline{\overline{e}}_{W2}$$

*то матрицы $a_1$ и $a_2$ связаны соотношением*

$$a_1 = b_*{}^* a_2{}_*{}^* c^{-1}{}_*{}^*$$

Доказательство. Вектор $\overline{v} \in V$ имеет разложение

$$\overline{v} = v^*{}_* e_{V1} = v^*{}_* b^*{}_* e_{V2}$$

относительно базисов $\overline{\overline{e}}_{V1}$, $\overline{\overline{e}}_{V2}$. Так как $a$ - гомоморфизм, то мы можем записать его в форме

(12.1.19)    $\overline{w} = v^*{}_* a_1{}^*{}_* e_{W1}$

относительно базисов $\overline{\overline{e}}_{V1}$, $\overline{\overline{e}}_{W1}$ и

(12.1.20)    $\overline{w} = v_*{}^* b^*{}_* a_2{}^*{}_* e_{W2}$

относительно базисов $\overline{\overline{e}}_{V2}$, $\overline{\overline{e}}_{W2}$. В силу теоремы 12.1.7 матрица $c$ имеет $*{}_*$-обратную и из равенства (12.1.16) следует

(12.1.21)    $\overline{\overline{e}}_{W2} = c^{-1^*}{}_*{}_* \overline{\overline{e}}_{W1}$

Равенство

(12.1.22)    $\overline{w} = v^*{}_* b^*{}_* a_2{}^*{}_* c^{-1^*}{}_*{}_* e_{W1}$

является следствием равенств (12.1.20), (12.1.21). В силу теоремы 8.4.4 сравнение равенств (12.1.19) и (12.1.22) даёт

(12.1.23)    $v^*{}_* a_1 = v^*{}_* b^*{}_* a_2{}^*{}_* c^{-1^*}$

Так как вектор $a$ - произвольный вектор, утверждение теоремы является следствием теоремы 12.2.1 и равенства (12.1.23).    □

---

Доказательство. Вектор $\overline{v} \in V$ имеет разложение

$$\overline{v} = v_*{}^* e_{V1} = v_*{}^* b_*{}^* e_{V2}$$

относительно базисов $\overline{\overline{e}}_{V1}$, $\overline{\overline{e}}_{V2}$. Так как $a$ - гомоморфизм, то мы можем записать его в форме

(12.1.24)    $\overline{w} = v_*{}^* a_{1*}{}^* e_{W1}$

относительно базисов $\overline{\overline{e}}_{V1}$, $\overline{\overline{e}}_{W1}$ и

(12.1.25)    $\overline{w} = v^*{}_* b_*{}^* a_2{}_*{}^* e_{W2}$

относительно базисов $\overline{\overline{e}}_{V2}$, $\overline{\overline{e}}_{W2}$. В силу теоремы 12.1.8 матрица $c$ имеет $_*{}^*$-обратную и из равенства (12.1.18) следует

(12.1.26)    $\overline{\overline{e}}_{W2} = c^{-1^*}{}_*{}^* \overline{\overline{e}}_{W1}$

Равенство

(12.1.27)    $\overline{w} = v_*{}^* b_*{}^* a_2{}_*{}^* c^{-1^*}{}_*{}^* e_{W1}$

является следствием равенств (12.1.25), (12.1.26). В силу теоремы 8.4.11 сравнение равенств (12.1.24) и (12.1.27) даёт

(12.1.28)    $v^*{}_* a_1 = v_*{}^* b_*{}^* a_2{}_*{}^* c^{-1^*}$

Так как вектор $a$ - произвольный вектор, утверждение теоремы является следствием теоремы 12.2.2 и равенства (12.1.28).    □

---

Теорема 12.1.15. *Пусть $\overline{a}$ - автоморфизм левого A-векторного пространства $V$ столбцов. Пусть $a_1$ - матрица этого автоморфизма, заданная относительно базиса $\overline{\overline{e}}_1$, и $a_2$ - матрица того же автоморфизма, заданная относительно базиса $\overline{\overline{e}}_2$. Если базис $\overline{\overline{e}}_1$ имеет координатную матрицу $b$ относительно базиса $\overline{\overline{e}}_2$*

$$\overline{\overline{e}}_1 = b^*{}_* \overline{\overline{e}}_2$$

*то матрицы $a_1$ и $a_2$ связаны соотношением*

$$a_1 = b^*{}_* a_2{}^*{}_* b^{-1^*}$$

Доказательство. Утверждение является следствием теоремы 12.1.13, так

---

Теорема 12.1.16. *Пусть $\overline{a}$ - автоморфизм левого A-векторного пространства $V$ столбцов. Пусть $a_1$ - матрица этого автоморфизма, заданная относительно базиса $\overline{\overline{e}}_1$, и $a_2$ - матрица того же автоморфизма, заданная относительно базиса $\overline{\overline{e}}_2$. Если базис $\overline{\overline{e}}_1$ имеет координатную матрицу $b$ относительно базиса $\overline{\overline{e}}_2$*

$$\overline{\overline{e}}_1 = b_*{}^* \overline{\overline{e}}_2$$

*то матрицы $a_1$ и $a_2$ связаны соотношением*

$$a_1 = b_*{}^* a_{2*}{}^* b^{-1^*}$$

Доказательство. Утверждение является следствием теоремы 12.1.14, так как в данном случае $c = b$. □ | как в данном случае $c = b$. □

## 12.2. Многообразие базисов левого $A$-векторного пространства

Пусть $V$ - левое $A$-векторное пространство столбцов. Мы доказали в теореме 12.1.9, что если задан базис левого $A$-векторного пространства $V$, то любой автоморфизм $\overline{f}$ левого $A$-векторного пространства $V$ можно отождествить с $_*$-невырожденной матрицей $f$. Соответствующее преобразование координат вектора

$$v' = v^*{}_* f$$

называется **линейным преобразованием**.

Пусть $V$ - левое $A$-векторное пространство строк. Мы доказали в теореме 12.1.10, что если задан базис левого $A$-векторного пространства $V$, то любой автоморфизм $\overline{f}$ левого $A$-векторного пространства $V$ можно отождествить с $_*$-невырожденной матрицей $f$. Соответствующее преобразование координат вектора

$$v' = v_*{}^* f$$

называется **линейным преобразованием**.

Теорема 12.2.1. *Пусть $V$ - левое $A$-векторное пространство столбцов и $\overline{e}$ - базис левого $A$-векторного пространства $V$. Автоморфизмы левого $A$-векторного пространства столбцов порождают правостороннее линейное эффективное $GL(n^*{}_*A)$-представление.*

Теорема 12.2.2. *Пусть $V$ - левое $A$-векторное пространство строк и $\overline{e}$ - базис левого $A$-векторного пространства $V$. Автоморфизмы левого $A$-векторного пространства строк порождают правостороннее линейное эффективное $GL(n_*{}^*A)$-представление.*

Доказательство. Пусть $a$, $b$ - матрицы автоморфизмов $\overline{a}$ и $\overline{b}$ относительно базиса $\overline{e}$. Согласно теореме 10.3.6 преобразование координат имеет вид

$$(12.2.1) \qquad v' = v^*{}_* a$$

$$(12.2.2) \qquad v'' = v'^*{}_* b$$

Равенство

$$v'' = v^*{}_* a^*{}_* b$$

является следствием равенств (12.2.1), (12.2.2). Согласно теореме 10.5.2, произведение автоморфизмов $\overline{a}$ и $\overline{b}$ имеет матрицу $a^*{}_* b$. Следовательно, автоморфизмы левого $A$-векторного пространства столбцов порождают правостороннее линейное $GL(n^*{}_*A)$-представление.

Нам осталось показать, что ядро неэффективности состоит только из единичного элемента. Тождественное пре-

Доказательство. Пусть $a$, $b$ - матрицы автоморфизмов $\overline{a}$ и $\overline{b}$ относительно базиса $\overline{e}$. Согласно теореме 10.3.7 преобразование координат имеет вид

$$(12.2.4) \qquad v' = v_*{}^* a$$

$$(12.2.5) \qquad v'' = v'_*{}^* b$$

Равенство

$$v'' = v_*{}^* a_*{}^* b$$

является следствием равенств (12.2.4), (12.2.5). Согласно теореме 10.5.3, произведение автоморфизмов $\overline{a}$ и $\overline{b}$ имеет матрицу $a_*{}^* b$. Следовательно, автоморфизмы левого $A$-векторного пространства строк порождают правостороннее линейное $GL(n_*{}^*A)$-представление.

Нам осталось показать, что ядро неэффективности состоит только из единичного элемента. Тождественное пре-образование удовлетворяет равенству

$$v^i = v^j a_j^i$$

Выбирая значения координат в форме $v^i = \delta_k^i$ , где $k$ задано, мы получим

(12.2.3)     $$\delta_k^i = \delta_k^j a_j^i$$

Из (12.2.3) следует

$$\delta_k^i = a_k^i$$

Поскольку $k$ произвольно, мы приходим к заключению $a = \delta$.     □

образование удовлетворяет равенству

$$v_i = v_j a_i^j$$

Выбирая значения координат в форме $v_i = \delta_i^k$ , где $k$ задано, мы получим

(12.2.6)     $$\delta_i^k = \delta_j^k a_i^j$$

Из (12.2.6) следует

$$\delta_i^k = a_i^k$$

Поскольку $k$ произвольно, мы приходим к заключению $a = \delta$.     □

<div style="background-color:#eef5ee">

Теорема 12.2.3. *Автоморфизм $a$, действуя на каждый вектор базиса в левом $A$-векторном пространстве столбцов , отображает базис в другой базис.*

Теорема 12.2.4. *Автоморфизм $a$, действуя на каждый вектор базиса в левом $A$-векторном пространстве строк , отображает базис в другой базис.*

</div>

Доказательство. Пусть $\overline{\overline{e}}$ - базис в левом $A$-векторном пространстве $V$ столбцов . Согласно теореме 12.1.9, вектор $e_i$ отображается в вектор $e_i'$

(12.2.7)     $$e_i' = e_i{}^*{}_* a$$

Если векторы $e_i'$ линейно зависимы, то в линейной комбинации

(12.2.8)     $$\lambda^*{}_* e' = 0$$

$\lambda \neq 0$. Из равенств (12.2.7) и (12.2.8) следует, что

$$\lambda^*{}_* e'^*{}_* a^{-1^*} = \lambda^*{}_* e = 0$$

и $\lambda \neq 0$. Это противоречит утверждению, что векторы $e_i$ линейно независимы. Следовательно, векторы $e_i'$ линейно независимы и порождают базис.     □

Доказательство. Пусть $\overline{\overline{e}}$ - базис в левом $A$-векторном пространстве $V$ строк . Согласно теореме 12.1.10, вектор $e^i$ отображается в вектор $e'^i$

(12.2.9)     $$e'^i = e^i{}_*{}^* a$$

Если векторы $e'^i$ линейно зависимы, то в линейной комбинации

(12.2.10)     $$\lambda_*{}^* e' = 0$$

$\lambda \neq 0$. Из равенств (12.2.9) и (12.2.10) следует, что

$$\lambda_*{}^* e'_*{}^* a^{-1^*} = \lambda_*{}^* e = 0$$

и $\lambda \neq 0$. Это противоречит утверждению, что векторы $e^i$ линейно независимы. Следовательно, векторы $e'^i$ линейно независимы и порождают базис.     □

Таким образом, мы можем распространить правостороннее линейное $GL(V_*)$-представление в левом $A$-векторном пространстве $V_*$ на множество базисов левого $A$-векторного пространства $V$.     Мы будем называть преобразование этого правостороннего представления на множестве базисов левого $A$-векторного пространства $V$ **активным преобразованием** потому, что гоморфизм левого $A$-векторного пространства породил это преобразование (Смотри также определение в разделе [4]-14.1-3, а также определение на странице [2]-214).

Соответственно определению мы будем записывать действие активного преобразования $a \in GL(V_*)$ на базис $\overline{\overline{e}}$ в форме $\overline{\overline{e}}{}_* a$. Рассмотрим равенство

(12.2.11)     $$v^*{}_* \overline{\overline{e}}{}^*{}_* a = v^*{}_* \overline{\overline{e}}{}^*{}_* a$$

Соответственно определению мы будем записывать действие активного преобразования $a \in GL(V_*)$ на базис $\overline{\overline{e}}$ в форме $\overline{\overline{e}}{}_* a$. Рассмотрим равенство

(12.2.12)     $$v_*{}^* \overline{\overline{e}}{}_*{}^* a = v_*{}^* \overline{\overline{e}}{}_*{}^* a$$

Выражение $\overline{\overline{e}}{}^*{}_*a$ в левой части равенства (12.2.11) является образом базиса $\overline{\overline{e}}$ при активном преобразовании $a$. Выражение $v^*{}_*\overline{e}$ в правой части равенства (12.2.11) является разложением вектора $\overline{v}$ относительно базиса $\overline{e}$. Следовательно, выражение в правой части равенства (12.2.11) является образом вектора $\overline{v}$ относительно эндоморфизма $a$ и выражение в левой части равенства

$$(12.2.11) \quad v^*{}_*\overline{\overline{e}}{}^*{}_*a = v^*{}_*\overline{\overline{e}}{}^*{}_*a$$

является разложением образа вектора $\overline{v}$ относительно образа базиса $\overline{\overline{e}}$. Следовательно, из равенства (12.2.11) следует, что эндоморфизм $a$ левого $A$-векторного пространства и соответствующее активное преобразование $a$ действуют синхронно и координаты $a \circ \overline{v}$ образа вектора $\overline{v}$ относительно образа $\overline{e}{}^*{}_*a$ базиса $\overline{e}$ совпадают с координатами вектора $\overline{v}$ относительно базиса $\overline{\overline{e}}$.

Выражение $\overline{\overline{e}}{}_*{}^*a$ в левой части равенства (12.2.12) является образом базиса $\overline{\overline{e}}$ при активном преобразовании $a$. Выражение $v_*{}^*\overline{e}$ в правой части равенства (12.2.12) является разложением вектора $\overline{v}$ относительно базиса $\overline{e}$. Следовательно, выражение в правой части равенства (12.2.12) является образом вектора $\overline{v}$ относительно эндоморфизма $a$ и выражение в левой части равенства

$$(12.2.12) \quad v_*{}^*\overline{\overline{e}}{}_*{}^*a = v_*{}^*\overline{\overline{e}}{}_*{}^*a$$

является разложением образа вектора $\overline{v}$ относительно образа базиса $\overline{\overline{e}}$. Следовательно, из равенства (12.2.12) следует, что эндоморфизм $a$ левого $A$-векторного пространства и соответствующее активное преобразование $a$ действуют синхронно и координаты $a \circ \overline{v}$ образа вектора $\overline{v}$ относительно образа $\overline{e}{}_*{}^*a$ базиса $\overline{e}$ совпадают с координатами вектора $\overline{v}$ относительно базиса $\overline{\overline{e}}$.

ТЕОРЕМА 12.2.5. *Активное правостороннее* $GL(n^*{}_*A)$-*представление на множестве базисов однотранзитивно. Множество базисов, отождествлённое с траекторией* $\overline{\overline{e}}{}^*{}_*GL(n^*{}_*A)(V)$ *активного правостороннего* $GL(n^*{}_*A)$-*представления называется* **многообразием базисов** *левого $A$-векторного пространства $V$.*

ТЕОРЕМА 12.2.6. *Активное правостороннее* $GL(n_*{}^*A)$-*представление на множестве базисов однотранзитивно. Множество базисов, отождествлённое с траекторией* $\overline{\overline{e}}{}_*{}^*GL(n_*{}^*A)(V)$ *активного правостороннего* $GL(n_*{}^*A)$-*представления называется* **многообразием базисов** *левого $A$-векторного пространства $V$.*

Чтобы доказать теоремы 12.2.5, 12.2.6, достаточно показать, что для любых двух базисов определено по крайней мере одно преобразование правостороннего представления и это преобразование единственно.

ДОКАЗАТЕЛЬСТВО. Гомоморфизм $a$, действуя на базис $\overline{e}$ имеет вид

$$e'_i = e_i{}^*{}_*a$$

где $e'_i$ - координатная матрица вектора $\overline{e}'_i$ относительно базиса $\overline{\overline{h}}$ и $e_i$ - координатная матрица вектора $\overline{e}_i$ относительно базиса $\overline{\overline{h}}$. Следовательно, координатная матрица образа базиса равна ${}^*{}_*$-произведению координатной матрицы исходного базиса и матрицы

ДОКАЗАТЕЛЬСТВО. Гомоморфизм $a$, действуя на базис $\overline{e}$ имеет вид

$$e'^i = e^i{}_*{}^*a$$

где $e'^i$ - координатная матрица вектора $\overline{e}'^i$ относительно базиса $\overline{\overline{h}}$ и $e^i$ - координатная матрица вектора $\overline{e}^i$ относительно базиса $\overline{\overline{h}}$. Следовательно, координатная матрица образа базиса равна ${}_*{}^*$-произведению координатной матрицы исходного базиса и матрицы

автоморфизма

$$e' = e^*_* g$$

В силу теоремы 12.1.7, матрицы $g$ и $e$ невырождены. Следовательно, матрица

$$a = e^{-1*}_{*}e'$$

является матрицей автоморфизма, отображающего базис $\overline{\overline{e}}$ в базис $\overline{\overline{e}}'$.

Допустим элементы $g_1$, $g_2$ группы $G$ и базис $\overline{\overline{e}}$ таковы, что

$$e^*_* g_1 = e^*_* g_2$$

В силу теорем 12.1.7 и 8.1.32 справедливо равенство $g_1 = g_2$. Отсюда следует утверждение теоремы.    □

автоморфизма

$$e' = e_{*}^{*} g$$

В силу теоремы 12.1.8, матрицы $g$ и $e$ невырождены. Следовательно, матрица

$$a = e^{-1\,*^*}_{\quad *}e'$$

является матрицей автоморфизма, отображающего базис $\overline{\overline{e}}$ в базис $\overline{\overline{e}}'$.

Допустим элементы $g_1$, $g_2$ группы $G$ и базис $\overline{\overline{e}}$ таковы, что

$$e_{*}^{*} g_1 = e_{*}^{*} g_2$$

В силу теорем 12.1.8 и 8.1.32 справедливо равенство $g_1 = g_2$. Отсюда следует утверждение теоремы.    □

Если в левом $A$-векторном пространстве $V$ определена дополнительная структура, не всякий автоморфизм сохраняет свойства заданной структуры. Например, если мы определим норму в левом $A$-векторном пространстве $V$, то для нас интересны автоморфизмы, которые сохраняют норму вектора.

ОПРЕДЕЛЕНИЕ 12.2.7. *Нормальная подгруппа $G(V)$ группы $GL(V)$, которая порождает автоморфизмы, сохраняющие свойства заданной структуры называется* **группой симметрии**.

*Не нарушая общности, мы будем отождествлять элемент $g$ группы $G(V)$ с соответствующим автоморфизмом и записывать его действие на вектор $v \in V$ в виде $v^*_* g$.*    □

ОПРЕДЕЛЕНИЕ 12.2.8. *Нормальная подгруппа $G(V)$ группы $GL(V)$, которая порождает автоморфизмы, сохраняющие свойства заданной структуры называется* **группой симметрии**.

*Не нарушая общности, мы будем отождествлять элемент $g$ группы $G(V)$ с соответствующим автоморфизмом и записывать его действие на вектор $v \in V$ в виде $v_{*}^{*} g$.*    □

Мы будем называть правостороннее представление группы $G$ в левом $A$-векторном пространстве **линейным $G$-представлением**.

ОПРЕДЕЛЕНИЕ 12.2.9. *Мы будем называть правостороннее представление группы $G(V)$ на множестве базисов левого $A$-векторного пространства $V$* **активным левосторонним $G$-представлением**.    □

Если базис $\overline{\overline{e}}$ задан, то мы можем отождествить автоморфизм левого $A$-векторного пространства $V$ и его координаты относительно базиса $\overline{\overline{e}}$. Множество $G(V_*)$ координат автоморфизмов относительно базиса $\overline{\overline{e}}$ является группой, изоморфной группе $G(V)$.

Не всякие два базиса могут быть связаны преобразованием группы симметрии потому, что не всякое невырожденное линейное преобразование принадлежит представлению группы $G(V)$. Таким образом, множество базисов можно представить как объединение орбит группы $G(V)$. В частности, если базис $\overline{\overline{e}} \in G(V)$, то орбита базиса $\overline{\overline{e}}$ совпадает с группой $G(V)$.

| | |
|---|---|
| **Определение 12.2.10.** *Мы будем называть орбиту* $\overline{\overline{e}}*_*G(V)$ *выбранного базиса* $\overline{\overline{e}}$ **многообразием базисов** *левого $A$-векторного пространства $V$ столбцов .* □ | **Определение 12.2.11.** *Мы будем называть орбиту* $\overline{\overline{e}}_**G(V)$ *выбранного базиса* $\overline{\overline{e}}$ **многообразием базисов** *левого $A$-векторного пространства $V$ строк .* □ |

**Теорема 12.2.12.** *Активное правостороннее $G(V)$-представление на многообразии базисов однотранзитивно.*

Доказательство. Теорема является следствием теоремы 12.2.5, и определений 12.2.10, а также теорема является следствием теоремы 12.2.6 и определений 12.2.11. □

| | |
|---|---|
| **Теорема 12.2.13.** *На многообразии базисов* $\overline{\overline{e}}*_*G(V)$ *левого $A$-векторного пространства столбцов существует однотранзитивное левостороннее $G(V)$-представление, перестановочное с активным.* | **Теорема 12.2.14.** *На многообразии базисов* $\overline{\overline{e}}_**G(V)$ *левого $A$-векторного пространства строк существует однотранзитивное левостороннее $G(V)$-представление, перестановочное с активным.* |
| Доказательство. Из теоремы 12.2.12 следует, что многообразие базисов $\overline{\overline{e}}*_*G(V)$ является однородным пространством группы $G$. Теорема является следствием теоремы 2.5.17. □ | Доказательство. Из теоремы 12.2.12 следует, что многообразие базисов $\overline{\overline{e}}_**G(V)$ является однородным пространством группы $G$. Теорема является следствием теоремы 2.5.17. □ |

Как мы видим из замечания 2.5.18 преобразование левостороннего $G(V)$-представления отличается от активного преобразования и не может быть сведено к преобразованию пространства $V$.

| | |
|---|---|
| **Определение 12.2.15.** *Преобразование левостороннего $G(V)$-представления называется* **пассивным преобразованием** *многообразия базисов* $\overline{\overline{e}}*_*G(V)$ *левого $A$-векторного пространства столбцов , а левостороннее $G(V)$-представление называется* **пассивным левосторонним $G(V)$-представлением**. *Согласно определению мы будем записывать пассивное преобразование базиса $\overline{\overline{e}}$, порождённое элементом $a \in G(V)$, в форме* $a*_*\overline{\overline{e}}$. □ | **Определение 12.2.16.** *Преобразование левостороннего $G(V)$-представления называется* **пассивным преобразованием** *многообразия базисов* $\overline{\overline{e}}_**G(V)$ *левого $A$-векторного пространства строк , а левостороннее $G(V)$-представление называется* **пассивным левосторонним $G(V)$-представлением**. *Согласно определению мы будем записывать пассивное преобразование базиса $\overline{\overline{e}}$, порождённое элементом $a \in G(V)$, в форме* $a_**\overline{\overline{e}}$. □ |
| **Теорема 12.2.17.** *Координатная матрица базиса $\overline{\overline{e}}'$ относительно ба-* | **Теорема 12.2.18.** *Координатная матрица базиса $\overline{\overline{e}}'$ относительно ба-* |

зиса $\overline{\overline{e}}$ левого $A$-векторного пространства $V$ столбцов совпадает с матрицей пассивного преобразования, отображающего базис $\overline{\overline{e}}$ в базис $\overline{\overline{e}}'$.

Доказательство. Согласно теореме 8.4.3, координатная матрица базиса $\overline{\overline{e}}'$ относительно базиса $\overline{\overline{e}}$ состоит из столбцов, являющихся координатными матрицами векторов $\overline{e}'_i$ относительно базиса $\overline{\overline{e}}$. Следовательно,

$$\overline{e}'_i = e'^*_i {}_* \overline{e}$$

В тоже время пассивное преобразование $a$, связывающее два базиса, имеет вид

$$\overline{e}'_i = a_i {}^*_* \overline{e}$$

В силу теоремы 8.4.4,

$$e'_i = a_i$$

для любого $i$. Это доказывает теорему. □

Замечание 12.2.19. Согласно теоремам 2.5.15, 12.2.17, мы можем отождествить базис $\overline{\overline{e}}'$ многообразия базисов $\overline{\overline{e}}^*_* G(V)$ левого $A$-векторного пространства столбцов и матрицу $g$ координат базиса $\overline{\overline{e}}'$ относительно базиса $\overline{\overline{e}}$

$$(12.2.13) \qquad \overline{\overline{e}}' = g^*_* \overline{\overline{e}}$$

Если $\overline{\overline{e}}' = \overline{\overline{e}}$, то равенство (12.2.13) принимает вид

$$(12.2.14) \qquad \overline{\overline{e}} = E_n {}^*_* \overline{\overline{e}}$$

Опираясь на равенство (12.2.14) мы будем пользоваться записью $\overline{\overline{E_n}} {}^*_* G(V)$ для многообразия базисов $\overline{\overline{e}}^*_* G(V)$. □

зиса $\overline{\overline{e}}$ левого $A$-векторного пространства $V$ строк совпадает с матрицей пассивного преобразования, отображающего базис $\overline{\overline{e}}$ в базис $\overline{\overline{e}}'$.

Доказательство. Согласно теореме 8.4.10, координатная матрица базиса $\overline{\overline{e}}'$ относительно базиса $\overline{\overline{e}}$ состоит из строк, являющихся координатными матрицами векторов $\overline{e}'_i$ относительно базиса $\overline{\overline{e}}$. Следовательно,

$$\overline{e}'^i = e'^i {}_* \overline{e}$$

В тоже время пассивное преобразование $a$, связывающее два базиса, имеет вид

$$\overline{e}'^i = a^i {}^*_* \overline{e}$$

В силу теоремы 8.4.11,

$$e'^i = a^i$$

для любого $i$. Это доказывает теорему. □

Замечание 12.2.20. Согласно теоремам 2.5.15, 12.2.18, мы можем отождествить базис $\overline{\overline{e}}'$ многообразия базисов $\overline{\overline{e}}_* {}^* G(V)$ левого $A$-векторного пространства столбцов и матрицу $g$ координат базиса $\overline{\overline{e}}'$ относительно базиса $\overline{\overline{e}}$

$$(12.2.15) \qquad \overline{\overline{e}}' = g_* {}^* \overline{\overline{e}}$$

Если $\overline{\overline{e}}' = \overline{\overline{e}}$, то равенство (12.2.15) принимает вид

$$(12.2.16) \qquad \overline{\overline{e}} = E_n {}_* {}^* \overline{\overline{e}}$$

Опираясь на равенство (12.2.16) мы будем пользоваться записью $\overline{\overline{E_n}}_* {}^* G(V)$ для многообразия базисов $\overline{\overline{e}}_* {}^* G(V)$. □

## 12.3. Пассивное преобразование в левом $A$-векторном пространстве

Активное преобразование изменяет базисы и векторы согласовано и координаты вектора относительно базиса не меняются. Пассивное преобразование меняет только базис, и это ведёт к преобразованию координат вектора относительно базиса.

Теорема 12.3.1. Пусть $V$ - левое $A$-векторное пространство столбцов.

Теорема 12.3.2. Пусть $V$ - левое $A$-векторное пространство строк.

Пусть пассивное преобразование $g \in G(V_*)$ отображает базис $\overline{\overline{e}}_1$ в базис $\overline{\overline{e}}_2$

$$(12.3.1) \qquad \overline{\overline{e}}_2 = g^*{}_*\overline{\overline{e}}_1$$

Пусть

$$(12.3.2) \qquad v_i = \begin{pmatrix} v_i^1 \\ ... \\ v_i^n \end{pmatrix}$$

матрица координат вектора $\overline{v}$ относительно базиса $\overline{\overline{e}}_i$, $i = 1, 2$. Преобразования координат

$$(12.3.3) \qquad v_1 = v_2{}^*{}_*g$$

$$(12.3.4) \qquad v_2 = v_1{}^*{}_*g^{-1*}{}^*$$

не зависят от вектора $\overline{v}$ или базиса $\overline{\overline{e}}$, а определенно исключительно координатами вектора $\overline{v}$ относительно базиса $\overline{\overline{e}}$.

**Доказательство.** Согласно теореме 8.4.3,

$$(12.3.8) \qquad \overline{v} = v_1{}^*{}_*\overline{\overline{e}}_1 = v_2{}^*{}_*\overline{\overline{e}}_2$$

Равенство

$$(12.3.9) \qquad v_1{}^*{}_*\overline{\overline{e}}_1 = v_2{}^*{}_*g^*{}_*\overline{\overline{e}}_1$$

является следствием равенств (12.3.1), (12.3.8). Согласно теореме 8.4.4, преобразование координат (12.3.3) является следствием равенства (12.3.9). Так как $g$ - $^*{}_*$-невырожденная матрица, то преобразование координат (12.3.4) является следствием равенства (12.3.3). $\square$

---

**Теорема 12.3.3.** *Преобразования координат*

$$\boxed{(12.3.4) \quad v_2 = v_1{}^*{}_*g^{-1*}{}_*}$$

*порождают эффективное линейное правостороннее контравариантное $G$-представление, называемое* **координатным представлением в левом $A$-векторном пространстве** *столбцов.*

---

Пусть пассивное преобразование $g \in G(V_*)$ отображает базис $\overline{\overline{e}}_1$ в базис $\overline{\overline{e}}_2$

$$(12.3.5) \qquad \overline{\overline{e}}_2 = g_*{}^*\overline{\overline{e}}_1$$

Пусть

$$v_i = \begin{pmatrix} v_{i1} & ... & v_{in} \end{pmatrix}$$

матрица координат вектора $\overline{v}$ относительно базиса $\overline{\overline{e}}_i$, $i = 1, 2$. Преобразования координат

$$(12.3.6) \qquad v_1 = v_{2*}{}^*g$$

$$(12.3.7) \qquad v_2 = v_{1*}{}^*g^{-1*}{}^*$$

не зависят от вектора $\overline{v}$ или базиса $\overline{\overline{e}}$, а определенно исключительно координатами вектора $\overline{v}$ относительно базиса $\overline{\overline{e}}$.

**Доказательство.** Согласно теореме 8.4.10,

$$(12.3.10) \qquad \overline{v} = v_{1*}{}^*\overline{\overline{e}}_1 = v_{2*}{}^*\overline{\overline{e}}_2$$

Равенство

$$(12.3.11) \qquad v_{1*}{}^*\overline{\overline{e}}_1 = v_{2*}{}^*g_*{}^*\overline{\overline{e}}_1$$

является следствием равенств (12.3.5), (12.3.10). Согласно теореме 8.4.11, преобразование координат (12.3.6) является следствием равенства (12.3.11). Так как $g$ - $_*{}^*$-невырожденная матрица, то преобразование координат (12.3.7) является следствием равенства (12.3.6). $\square$

---

**Теорема 12.3.4.** *Преобразования координат*

$$\boxed{(12.3.7) \quad v_2 = v_{1*}{}^*g^{-1*}{}^*}$$

*порождают эффективное линейное правостороннее контравариантное $G$-представление, называемое* **координатным представлением в левом $A$-векторном пространстве** *строк.*

Доказательство.   Пусть

$$(12.3.12) \qquad v_i = \begin{pmatrix} v_i^1 \\ ... \\ v_i^n \end{pmatrix}$$

матрица координат вектора $\overline{v}$ относительно базиса $\overline{\overline{e}}_i$, $i = 1, 2, 3$.   Тогда

$$(12.3.13) \qquad \overline{v} = v_1{}^*{}_*\overline{\overline{e}}_1 = v_2{}^*{}_*\overline{\overline{e}}_2$$
$$= v_3{}^*{}_*\overline{\overline{e}}_3$$

Пусть пассивное преобразование $g_1$ отображает базис $\overline{\overline{e}}_1 \in \overline{e}^*{}_*G(V)$ в базис $\overline{\overline{e}}_2 \in \overline{e}^*{}_*G(V)$

$$(12.3.14) \qquad \overline{\overline{e}}_2 = g_1{}^*{}_*\overline{\overline{e}}_1$$

Пусть пассивное преобразование $g_2$ отображает базис $\overline{\overline{e}}_2 \in \overline{e}^*{}_*G(V)$ в базис $\overline{\overline{e}}_3 \in \overline{e}^*{}_*G(V)$

$$(12.3.15) \qquad \overline{\overline{e}}_3 = g_2{}^*{}_*\overline{\overline{e}}_2$$

Преобразование

$$(12.3.16) \qquad \overline{\overline{e}}_3 = (g_2{}^*{}_*g_1)^*{}_*\overline{\overline{e}}_1$$

является произведением преобразований (12.3.14), (12.3.15). Преобразование координат

$$(12.3.17) \qquad v_2 = v_1{}^*{}_*g_1^{-1}{}^*{}_*$$

соответствует пассивному преобразованию $g_1$ и является следствием равенств (12.3.13), (12.3.14) и теоремы 12.3.1.   Аналогично, согласно теореме 12.3.1, преобразование координат

$$(12.3.18) \qquad v_3 = v_2{}^*{}_*g_2^{-1}{}^*{}_*$$

соответствует пассивному преобразованию $g_2$ и является следствием равенств (12.3.13), (12.3.15). Произведение преобразований координат (12.3.17), (12.3.18) имеет вид

$$(12.3.19) \qquad v_3 = v_1{}^*{}_*g_1^{-1}{}^*{}_*{}^*{}_*g_2^{-1}{}^*{}_*$$

и является координатным преобразованием, соответствующим пассивному преобразованию (12.3.16). Равенство

$$(12.3.20) \qquad v_3 = v_1{}^*{}_*(g_2{}^*{}_*g_1)^{-1}{}^*{}_*$$

Доказательство.   Пусть

$$v_i = \begin{pmatrix} v_{i\,1} & ... & v_{i\,n} \end{pmatrix}$$

матрица координат вектора $\overline{v}$ относительно базиса $\overline{e}_i$, $i = 1, 2, 3$.   Тогда

$$(12.3.21) \qquad \overline{v} = v_{1*}{}^*\overline{e}_1 = v_{2*}{}^*\overline{e}_2$$
$$= v_{3*}{}^*\overline{e}_3$$

Пусть пассивное преобразование $g_1$ отображает базис $\overline{\overline{e}}_1 \in \overline{e}_*{}^*G(V)$ в базис $\overline{\overline{e}}_2 \in \overline{e}_*{}^*G(V)$

$$(12.3.22) \qquad \overline{\overline{e}}_2 = g_{1*}{}^*\overline{\overline{e}}_1$$

Пусть пассивное преобразование $g_2$ отображает базис $\overline{\overline{e}}_2 \in \overline{e}_*{}^*G(V)$ в базис $\overline{\overline{e}}_3 \in \overline{e}_*{}^*G(V)$

$$(12.3.23) \qquad \overline{\overline{e}}_3 = g_{2*}{}^*\overline{\overline{e}}_2$$

Преобразование

$$(12.3.24) \qquad \overline{\overline{e}}_3 = (g_{2*}{}^*g_1)_*{}^*\overline{\overline{e}}_1$$

является произведением преобразований (12.3.22), (12.3.23). Преобразование координат

$$(12.3.25) \qquad v_2 = v_{1*}{}^*g_1^{-1}{}_*{}^*$$

соответствует пассивному преобразованию $g_1$ и является следствием равенств (12.3.21), (12.3.22) и теоремы 12.3.2.   Аналогично, согласно теореме 12.3.2, преобразование координат

$$(12.3.26) \qquad v_3 = v_{2*}{}^*g_2^{-1}{}_*{}^*$$

соответствует пассивному преобразованию $g_2$ и является следствием равенств (12.3.21), (12.3.23). Произведение преобразований координат (12.3.17), (12.3.18) имеет вид

$$(12.3.27) \qquad v_3 = v_{1*}{}^*g_1^{-1}{}_*{}^*{}^*g_2^{-1}{}_*{}^*$$

и является координатным преобразованием, соответствующим пассивному преобразованию (12.3.24). Равенство

$$(12.3.28) \qquad v_3 = v_{1*}{}^*(g_{2*}{}^*g_1)^{-1}{}_*{}^*$$

является следствием равенств (8.1.66), (12.3.19). Согласно определению 2.5.9, преобразования координат порождают линейное правостороннее контравариантное $G$-представление. Если координатное преобразование не изменяет векторы $\delta_k$, то ему соответствует единица группы $G$, так как пассивное представление однотранзитивно. Следовательно, координатное представление эффективно. $\square$

является следствием равенств (8.1.66), (12.3.27). Согласно определению 2.5.9, преобразования координат порождают линейное правостороннее контравариантное $G$-представление. Если координатное преобразование не изменяет векторы $\delta_k$, то ему соответствует единица группы $G$, так как пассивное представление однотранзитивно. Следовательно, координатное представление эффективно. $\square$

ТЕОРЕМА 12.3.5. *Пусть $V$ - левое $A$-векторное пространство столбцов и $\overline{\overline{e}}_1$, $\overline{\overline{e}}_2$ - базисы в левом $A$-векторном пространстве $V$. Пусть базисы $\overline{\overline{e}}_1$, $\overline{\overline{e}}_2$ связаны пассивным преобразованием $g$*

$$\overline{\overline{e}}_2 = g^*{}_*\overline{\overline{e}}_1$$

*Пусть $f$ - эндоморфизм левого $A$-векторного пространства $V$. Пусть $f_i$, $i = 1$, $2$, - матрица эндоморфизма $f$ относительно базиса $\overline{\overline{e}}_i$. Тогда*

(12.3.29) $$f_2 = g^*{}_* f_1{}^*{}_* g^{-1^*{}^*}$$

ТЕОРЕМА 12.3.6. *Пусть $V$ - левое $A$-векторное пространство строк и $\overline{\overline{e}}_1$, $\overline{\overline{e}}_2$ - базисы в левом $A$-векторном пространстве $V$. Пусть базисы $\overline{\overline{e}}_1$, $\overline{\overline{e}}_2$ связаны пассивным преобразованием $g$*

$$\overline{\overline{e}}_2 = g_{**}\overline{\overline{e}}_1$$

*Пусть $f$ - эндоморфизм левого $A$-векторного пространства $V$. Пусть $f_i$, $i = 1$, $2$, - матрица эндоморфизма $f$ относительно базиса $\overline{\overline{e}}_i$. Тогда*

(12.3.30) $$f_2 = g_{**} f_1{}_{**} g^{-1_{*}{}^{*}}$$

ДОКАЗАТЕЛЬСТВО ТЕОРЕМЫ 12.3.5. Пусть эндоморфизм $f$ отображает вектор $v$ в вектор $w$

(12.3.31) $$w = f \circ v$$

Пусть $v_i$, $i = 1$, $2$, - матрица вектора $v$ относительно базиса $\overline{\overline{e}}_i$. Тогда

(12.3.32) $$v_1 = v_2{}^*{}_* g$$

следует из теоремы 12.3.1.

Пусть $w_i$, $i = 1$, $2$, - матрица вектора $w$ относительно базиса $\overline{\overline{e}}_i$. Тогда

(12.3.33) $$w_1 = w_2{}^*{}_* g$$

следует из теоремы 12.3.1.

Согласно теореме 10.3.6, равенства

(12.3.34) $$w_1 = v_1{}^*{}_* f_1$$

(12.3.35) $$w_2 = v_2{}^*{}_* f_2$$

являются следствием равенства (12.3.31). Равенство

(12.3.36) $$w_2{}^*{}_* g = v_1{}^*{}_* f_1 = v_2{}^*{}_* g^*{}_* f_1$$

является следствием равенств (12.3.32), (12.3.34), (12.3.33). Так как $g$ - невырожденная матрица, то равенство

(12.3.37) $$w_2 = v_1{}^*{}_* f_1 = v_2{}^*{}_* g^*{}_* f_1{}^*{}_* g^{-1^*{}^*}$$

является следствием равенства (12.3.36). Равенство

$$(12.3.38) \qquad v_2{}^*{}_*f_2 = v_2{}^*{}_*g^*{}_*f_1{}^*{}_*g^{-1^*}{}^*$$

является следствием равенств (12.3.35), (12.3.37). Равенство (12.3.29) является следствием равенства (12.3.38). □

Доказательство теоремы 12.3.6. Пусть эндоморфизм $f$ отображает вектор $v$ в вектор $w$

$$(12.3.39) \qquad w = f \circ v$$

<table>
<tr><td>

Пусть $v_i$, $i = 1, 2,$ - матрица вектора $v$ относительно базиса $\overline{\overline{e}}_i$. Тогда

$$(12.3.40) \qquad v_1 = v_2{}_*{}^*g$$

следует из теоремы 12.3.2.

</td><td>

Пусть $w_i$, $i = 1, 2,$ - матрица вектора $w$ относительно базиса $\overline{\overline{e}}_i$. Тогда

$$(12.3.41) \qquad w_1 = w_2{}_*{}^*g$$

следует из теоремы 12.3.2.

</td></tr>
</table>

Согласно теореме 10.3.6, равенства

$$(12.3.42) \qquad w_1 = v_1{}_*{}^*f_1$$

$$(12.3.43) \qquad w_2 = v_2{}_*{}^*f_2$$

являются следствием равенства (12.3.39). Равенство

$$(12.3.44) \qquad w_2{}_*{}^*g = v_1{}_*{}^*f_1 = v_2{}_*{}^*g_*{}^*f_1$$

является следствием равенств (12.3.40), (12.3.42), (12.3.41). Так как $g$ - невырожденная матрица, то равенство

$$(12.3.45) \qquad w_2 = v_1{}_*{}^*f_1 = v_2{}_*{}^*g_*{}^*f_1{}_*{}^*g^{-1^*}$$

является следствием равенства (12.3.44). Равенство

$$(12.3.46) \qquad v_2{}_*{}^*f_2 = v_2{}_*{}^*g_*{}^*f_1{}_*{}^*g^{-1^*}$$

является следствием равенств (12.3.43), (12.3.45). Равенство (12.3.30) является следствием равенства (12.3.46). □

<table>
<tr><td>

Теорема 12.3.7. *Пусть $V$ - левое $A$-векторное пространство столбцов. Пусть эндоморфизм $\overline{f}$ левого $A$-векторного пространства $V$ является суммой эндоморфизмов $\overline{f}_1$, $\overline{f}_2$. Матрица $f$ эндоморфизма $\overline{f}$ равна сумме матриц $f_1$, $f_2$ эндоморфизмов $\overline{f}_1$, $\overline{f}_2$ и это равенство не зависит от выбора базиса.*

Доказательство. Пусть $\overline{\overline{e}}$, $\overline{\overline{e}}'$ - базисы левого $A$-векторного пространства $V$ и связаны пассивным преобразованием $g$

$$e' = e^*{}_*g$$

</td><td>

Теорема 12.3.8. *Пусть $V$ - левое $A$-векторное пространство строк. Пусть эндоморфизм $\overline{f}$ левого $A$-векторного пространства $V$ является суммой эндоморфизмов $\overline{f}_1$, $\overline{f}_2$. Матрица $f$ эндоморфизма $\overline{f}$ равна сумме матриц $f_1$, $f_2$ эндоморфизмов $\overline{f}_1$, $\overline{f}_2$ и это равенство не зависит от выбора базиса.*

Доказательство. Пусть $\overline{\overline{e}}$, $\overline{\overline{e}}'$ - базисы левого $A$-векторного пространства $V$ и связаны пассивным преобразованием $g$

$$e' = e_*{}^*g$$

</td></tr>
</table>

Пусть $f$ - матрица эндоморфизма $\overline{f}$ относительно базиса $\overline{e}$. Пусть $f'$ - матрица эндоморфизма $\overline{f}$ относительно базиса $\overline{e}'$. Пусть $f_i$, $i = 1, 2$, - матрица эндоморфизма $\overline{f}_i$ относительно базиса $\overline{e}$. Пусть $f'_i$, $i = 1, 2$, - матрица эндоморфизма $\overline{f}_i$ относительно базиса $\overline{e}'$. Согласно теореме 10.4.2,

$$(12.3.47) \qquad f = f_1 + f_2$$

$$(12.3.48) \qquad f' = f'_1 + f'_2$$

Согласно теореме 12.3.5,

$$(12.3.49) \qquad f' = g^*{}_* f^*{}_* g^{-1^*}{}^*$$

$$(12.3.50) \qquad f'_1 = g^*{}_* f_1^*{}_* g^{-1^*}{}^*$$

$$(12.3.51) \qquad f'_2 = g^*{}_* f_2^*{}_* g^{-1^*}{}^*$$

Равенство

$$
\begin{aligned}
f' &= g^*{}_* f^*{}_* g^{-1^*}{}^* \\
&= g^*{}_* (f_1 + f_2)^*{}_* g^{-1^*}{}^* \\
(12.3.52) \qquad &= g^*{}_* f_1^*{}_* g^{-1^*}{}^* \\
&\quad + g^*{}_* f_2^*{}_* g^{-1^*}{}^* \\
&= f'_1 + f'_2
\end{aligned}
$$

является следствием равенств (12.3.47), (12.3.49), (12.3.50), (12.3.51). Из равенств (12.3.48), (12.3.52) следует, что равенство (12.3.47) ковариантно относительно выбора базиса левого $A$-векторного пространства $V$. $\square$

---

Пусть $f$ - матрица эндоморфизма $\overline{f}$ относительно базиса $\overline{e}$. Пусть $f'$ - матрица эндоморфизма $\overline{f}$ относительно базиса $\overline{e}'$. Пусть $f_i$, $i = 1, 2$, - матрица эндоморфизма $\overline{f}_i$ относительно базиса $\overline{e}$. Пусть $f'_i$, $i = 1, 2$, - матрица эндоморфизма $\overline{f}_i$ относительно базиса $\overline{e}'$. Согласно теореме 10.4.3,

$$(12.3.53) \qquad f = f_1 + f_2$$

$$(12.3.54) \qquad f' = f'_1 + f'_2$$

Согласно теореме 12.3.6,

$$(12.3.55) \qquad f' = g_*{}^* f_*{}^* g^{-1}{}_*{}^*$$

$$(12.3.56) \qquad f'_1 = g_*{}^* f_{1*}{}^* g^{-1}{}_*{}^*$$

$$(12.3.57) \qquad f'_2 = g_*{}^* f_{2*}{}^* g^{-1}{}_*{}^*$$

Равенство

$$
\begin{aligned}
f' &= g_*{}^* f_*{}^* g^{-1}{}_*{}^* \\
&= g_*{}^* (f_1 + f_2)_*{}^* g^{-1}{}_*{}^* \\
(12.3.58) \qquad &= g_*{}^* f_{1*}{}^* g^{-1}{}_*{}^* \\
&\quad + g_*{}^* f_{2*}{}^* g^{-1}{}_*{}^* \\
&= f'_1 + f'_2
\end{aligned}
$$

является следствием равенств (12.3.53), (12.3.55), (12.3.56), (12.3.57). Из равенств (12.3.54), (12.3.58) следует, что равенство (12.3.53) ковариантно относительно выбора базиса левого $A$-векторного пространства $V$. $\square$

---

Теорема 12.3.9. *Пусть $V$ - левое $A$-векторное пространство столбцов. Пусть эндоморфизм $\overline{f}$ левого $A$-векторного пространства $V$ является произведением эндоморфизмов $\overline{f}_1$, $\overline{f}_2$. Матрица $f$ эндоморфизма $\overline{f}$ равна $^*{}_*$-произведению матриц $f_1$, $f_2$ эндоморфизмов $\overline{f}_1$, $\overline{f}_2$ и это равенство не зависит от выбора базиса.*

Доказательство. Пусть $\overline{e}$, $\overline{e}'$ - базисы левого $A$-векторного пространства $V$ и связаны пассивным преобра-

---

Теорема 12.3.10. *Пусть $V$ - левое $A$-векторное пространство строк. Пусть эндоморфизм $\overline{f}$ левого $A$-векторного пространства $V$ является произведением эндоморфизмов $\overline{f}_1$, $\overline{f}_2$. Матрица $f$ эндоморфизма $\overline{f}$ равна $_*{}^*$-произведению матриц $f_1$, $f_2$ эндоморфизмов $\overline{f}_1$, $\overline{f}_2$ и это равенство не зависит от выбора базиса.*

Доказательство. Пусть $\overline{e}$, $\overline{e}'$ - базисы левого $A$-векторного пространства $V$ и связаны пассивным преобра-

зованием $g$

$$e' = e^*{}_*g$$

Пусть $f$ - матрица эндоморфизма $\overline{f}$ относительно базиса $\overline{\overline{e}}$. Пусть $f'$ - матрица эндоморфизма $\overline{f}$ относительно базиса $\overline{\overline{e}}'$. Пусть $f_i$, $i = 1, 2$, - матрица эндоморфизма $\overline{f}_i$ относительно базиса $\overline{\overline{e}}$. Пусть $f'_i$, $i = 1, 2$, - матрица эндоморфизма $\overline{f}_i$ относительно базиса $\overline{\overline{e}}'$. Согласно теореме 10.5.2,

(12.3.59)     $f = f_{1*}{}^*f_2$

(12.3.60)     $f' = f'_{1*}{}^*f'_2$

Согласно теореме 12.3.5,

(12.3.61)     $f' = g^*{}_*f^*{}_*g^{-1^*}{}^*$

(12.3.62)     $f'_1 = g^*{}_*f_1^*{}_*g^{-1^*}{}^*$

(12.3.63)     $f'_2 = g^*{}_*f_2^*{}_*g^{-1^*}{}^*$

Равенство

$$f' = g^*{}_*f^*{}_*g^{-1^*}{}^*$$
$$= g^*{}_*f_1^*{}_*f_2^*{}_*g^{-1^*}{}^*$$

(12.3.64)     $= g^*{}_*f_1^*{}_*g^{-1^*}{}^*$

$$^*{}_*g^*{}_*f_2^*{}_*g^{-1^*}{}^*$$
$$= f'^*_1{}_*f'_2$$

является следствием равенств (12.3.59), (12.3.61), (12.3.62), (12.3.63). Из равенств (12.3.60), (12.3.64) следует, что равенство (12.3.59) ковариантно относительно выбора базиса левого $A$-векторного пространства $V$.     □

зованием $g$

$$e' = e_*{}^*g$$

Пусть $f$ - матрица эндоморфизма $\overline{f}$ относительно базиса $\overline{\overline{e}}$. Пусть $f'$ - матрица эндоморфизма $\overline{f}$ относительно базиса $\overline{\overline{e}}'$. Пусть $f_i$, $i = 1, 2$, - матрица эндоморфизма $\overline{f}_i$ относительно базиса $\overline{\overline{e}}$. Пусть $f'_i$, $i = 1, 2$, - матрица эндоморфизма $\overline{f}_i$ относительно базиса $\overline{\overline{e}}'$. Согласно теореме 10.5.3,

(12.3.65)     $f = f_1{}^*{}_*f_2$

(12.3.66)     $f' = f'_1{}^*{}_*f'_2$

Согласно теореме 12.3.6,

(12.3.67)     $f' = g_*{}^*f_*{}^*g^{-1_*}{}^*$

(12.3.68)     $f'_1 = g_*{}^*f_{1*}{}^*g^{-1_*}{}^*$

(12.3.69)     $f'_2 = g_*{}^*f_{2*}{}^*g^{-1_*}{}^*$

Равенство

$$f' = g_*{}^*f_*{}^*g^{-1_*}{}^*$$
$$= g_*{}^*f_{1*}{}^*f_{2*}{}^*g^{-1_*}{}^*$$

(12.3.70)     $= g_*{}^*f_{1*}{}^*g^{-1_*}{}^*$

$$_*{}^*g_*{}^*f_{2*}{}^*g^{-1_*}{}^*$$
$$= f'_{1*}{}^*f'_2$$

является следствием равенств (12.3.65), (12.3.67), (12.3.68), (12.3.69). Из равенств (12.3.66), (12.3.70) следует, что равенство (12.3.65) ковариантно относительно выбора базиса левого $A$-векторного пространства $V$.     □

## 12.4. Размерность правого $A$-векторного пространства

**Теорема 12.4.1.** *Пусть $V$ - правое $A$-векторное пространство столбцов. Пусть множество векторов $\overline{\overline{e}} = (e_i, i \in I)$ является базисом правого $A$-векторного пространства $V$. Пусть множество векторов $\overline{\overline{e}} = (e_j, j \in J)$ является базисом правого $A$-векторного пространства $V$. Если $|I|$ и $|J|$ - конечные числа, то $|I| = |J|$.*

**Теорема 12.4.2.** *Пусть $V$ - правое $A$-векторное пространство строк. Пусть множество векторов $\overline{\overline{e}} = (e^i, i \in I)$ является базисом правого $A$-векторного пространства $V$. Пусть множество векторов $\overline{\overline{e}} = (e^j, j \in J)$ является базисом правого $A$-векторного пространства $V$. Если $|I|$ и $|J|$ - конечные числа, то $|I| = |J|$.*

Доказательство. Предположим, что $|I| = m$ и $|J| = n$. Предположим, что

$$(12.4.1) \qquad m < n$$

Так как $\overline{\overline{e}}_1$ - базис, любой вектор $e_{2j}$, $j \in J$, имеет разложение

$$e_{2j} = e_{1*}{}^*a_j$$

Так как $\overline{\overline{e}}_2$ - базис,

$$(12.4.2) \qquad \lambda = 0$$

следует из

$$(12.4.3) \qquad e_{2*}{}^*\lambda = e_{1*}{}^*a_*{}^*\lambda = 0$$

Так как $\overline{\overline{e}}_1$ - базис, мы получим

$$(12.4.4) \qquad a_*{}^*\lambda = 0$$

Согласно (12.4.1), $\mathrm{rank}(_*{}^*)a \le m$ и система линейных уравнений (12.4.4) имеет больше переменных, чем уравнений. Согласно теореме 9.4.4, $\lambda \ne 0$, что противоречит утверждению (12.4.2). Следовательно, утверждение $m < n$ неверно. Таким же образом мы можем доказать, что утверждение $n < m$ неверно. Это завершает доказательство теоремы. $\qquad\square$

Доказательство. Предположим, что $|I| = m$ и $|J| = n$. Предположим, что

$$(12.4.5) \qquad m < n$$

Так как $\overline{\overline{e}}_1$ - базис, любой вектор $e_2^j$, $j \in J$, имеет разложение

$$e_2^j = e_1{}^*{}_*a^j$$

Так как $\overline{\overline{e}}_2$ - базис,

$$(12.4.6) \qquad \lambda = 0$$

следует из

$$(12.4.7) \qquad e_2{}^*{}_*\lambda = e_1{}^*{}_*a^*{}_*\lambda = 0$$

Так как $\overline{\overline{e}}_1$ - базис, мы получим

$$(12.4.8) \qquad a^*{}_*\lambda = 0$$

Согласно (12.4.5), $\mathrm{rank}(^*{}_*)a \le m$ и система линейных уравнений (12.4.8) имеет больше переменных, чем уравнений. Согласно теореме 9.4.4, $\lambda \ne 0$, что противоречит утверждению (12.4.6). Следовательно, утверждение $m < n$ неверно. Таким же образом мы можем доказать, что утверждение $n < m$ неверно. Это завершает доказательство теоремы. $\qquad\square$

---

Определение 12.4.3. *Мы будем называть* **размерностью правого $A$-векторного пространства** *число векторов в базисе.* $\qquad\square$

---

Теорема 12.4.4. *Если правое $A$-векторное пространство $V$ имеет конечную размерность, то из утверждения*

$$(12.4.9) \qquad \forall v \in V, va = vb \quad a, b \in A$$

*следует, что $a = b$.*

Теорема 12.4.5. *Если правое $A$-векторное пространство $V$ имеет конечную размерность, то из утверждения*

$$(12.4.10) \qquad \forall v \in V, va = vb \quad a, b \in A$$

*следует, что $a = b$.*

Доказательство. Пусть множество векторов $\overline{e} = (e_j, j \in J)$ является базисом левого $A$-векторного пространство $V$. Согласно теоремам 7.4.6, 8.4.18, равенство

$$(12.4.11) \qquad v^i a = v^i b$$

является следствием равенства (12.4.9). Теорема следует из равенства (12.4.11), если мы положим $v = e_1$.

Доказательство. Пусть множество векторов $\overline{e} = (e_j, j \in J)$ является базисом левого $A$-векторного пространство $V$. Согласно теоремам 7.4.6, 8.4.25, равенство

$$(12.4.12) \qquad v_i a = v_i b$$

является следствием равенства (12.4.10). Теорема следует из равенства (12.4.12), если мы положим

☐ | $v = e_1$.                                        ☐

**Теорема 12.4.6.** *Если правое $A$-векторное пространство $V$ имеет конечную размерность, то соответствующее представление свободно.*

Доказательство. Теорема является следствием определения 2.3.3 и теорем 12.4.4, 12.4.5.                                        ☐

**Теорема 12.4.7.** *Пусть $V$ - правое $A$-векторное пространство столбцов. Координатная матрица $\overline{\overline{g}}$ базиса $\overline{\overline{g}}$ относительно базиса $\overline{\overline{e}}$ правого $A$-векторного пространства $V$ является $_*{}^*$-невырожденной матрицей.*

Доказательство. Согласно теореме 9.3.6, $_*{}^*$-ранг координатной матрицы базиса $\overline{\overline{g}}$ относительно базиса $\overline{\overline{e}}$ равен размерности левого $A$-пространства столбцов, откуда следует утверждение теоремы.                ☐

**Теорема 12.4.9.** *Пусть $V$ - правое $A$-векторное пространство столбцов. Пусть $\overline{\overline{e}}$ - базис правого $A$-векторного пространства $V$. Тогда любой автоморфизм $\overline{f}$ правого $A$-векторного пространства $V$ имеет вид*

(12.4.13)          $v' = f_*{}^* v$

*где $f$ - $_*{}^*$-невырожденная матрица.*

Доказательство. Равенство (12.4.13) следует из теоремы 11.3.6. Так как $\overline{f}$ - изоморфизм, то для каждого вектора $\overline{v}'$ существует единственный вектор $\overline{v}$ такой, что $\overline{v}' = \overline{f} \circ \overline{v}$. Следовательно, система линейных уравнений (12.4.13) имеет единственное решение. Согласно следствию 9.4.5 матрица $f$ $_*{}^*$-невырожденна.                ☐

**Теорема 12.4.11.** *Матрицы автоморфизмов правого $A$-пространства столбцов порождают группу $GL(n_* {}^* A)$.*

**Теорема 12.4.8.** *Пусть $V$ - правое $A$-векторное пространство строк. Координатная матрица базиса $\overline{\overline{g}}$ относительно базиса $\overline{\overline{e}}$ правого $A$-векторного пространства $V$ является $_*{}^*$-невырожденной матрицей.*

Доказательство. Согласно теореме 9.3.7, $_*{}^*$-ранг координатной матрицы базиса $\overline{\overline{g}}$ относительно базиса $\overline{\overline{e}}$ равен размерности левого $A$-пространства столбцов, откуда следует утверждение теоремы.                ☐

**Теорема 12.4.10.** *Пусть $V$ - правое $A$-векторное пространство строк. Пусть $\overline{\overline{e}}$ - базис правого $A$-векторного пространства $V$. Тогда любой автоморфизм $\overline{f}$ правого $A$-векторного пространства $V$ имеет вид*

(12.4.14)          $v' = f^* {}_* v$

*где $f$ - $^*{}_*$-невырожденная матрица.*

Доказательство. Равенство (12.4.14) следует из теоремы 11.3.7. Так как $\overline{f}$ - изоморфизм, то для каждого вектора $\overline{v}'$ существует единственный вектор $\overline{v}$ такой, что $\overline{v}' = \overline{f} \circ \overline{v}$. Следовательно, система линейных уравнений (12.4.14) имеет единственное решение. Согласно следствию 9.4.5 матрица $f$ $^*{}_*$-невырожденна.                ☐

**Теорема 12.4.12.** *Матрицы автоморфизмов правого $A$-пространства строк порождают группу $GL(n^* {}_* A)$.*

Доказательство. Если даны два автоморфизма $\overline{f}$ и $\overline{g}$, то мы можем записать

$$v' = f_*{}^* v$$

$$v'' = g_*{}^* v' = g_*{}^* f_*{}^* v$$

Следовательно, результирующий автоморфизм имеет матрицу $f_*{}^* g$. $\quad\square$

Теорема 12.4.13. *Пусть* $\overline{\overline{e}}_{V1}$ *,* $\overline{\overline{e}}_{V2}$ *- базисы правого $A$-векторного пространства $V$ столбцов. Пусть* $\overline{\overline{e}}_{W1}$ *,* $\overline{\overline{e}}_{W2}$ *- базисы правого $A$-векторного пространства $W$ столбцов. Пусть $a_1$ - матрица гомоморфизма*

$$(12.4.15) \qquad a : V \to W$$

*относительно базисов* $\overline{\overline{e}}_{V1}$*,* $\overline{\overline{e}}_{W1}$ *и $a_2$ - матрица гомоморфизма* (12.4.15) *относительно базисов* $\overline{\overline{e}}_{V2}$*,* $\overline{\overline{e}}_{W2}$*. Если базис* $\overline{\overline{e}}_{V1}$ *имеет координатную матрицу $b$ относительно базиса* $\overline{\overline{e}}_{V2}$

$$\overline{\overline{e}}_{V1} = \overline{\overline{e}}_{V2*}{}^* b$$

*и* $\overline{\overline{e}}_{W1}$ *имеет координатную матрицу $c$ относительно базиса* $\overline{\overline{e}}_{W2}$

$$(12.4.16) \qquad \overline{\overline{e}}_{W1} = \overline{\overline{e}}_{W2*}{}^* c$$

*то матрицы $a_1$ и $a_2$ связаны соотношением*

$$a_1 = c^{-1}{}_*{}^* {}_* a_{2*}{}^* b$$

Доказательство. Вектор $\overline{v} \in V$ имеет разложение

$$\overline{v} = e_{V1}{}^*{}_* v = e_{V2}{}^*{}_* b^*{}_* v$$

относительно базисов $\overline{\overline{e}}_{V1}$, $\overline{\overline{e}}_{V2}$. Так как $a$ - гомоморфизм, то мы можем записать его в форме

$$(12.4.19) \qquad \overline{w} = e_{W1}{}^*{}_* a_1{}^*{}_* v$$

относительно базисов $\overline{\overline{e}}_{V1}$, $\overline{\overline{e}}_{W1}$ и

$$(12.4.20) \qquad \overline{w} = e_{W2*}{}^* a_{2*}{}^* b^*{}_* v$$

относительно базисов $\overline{\overline{e}}_{V2}$, $\overline{\overline{e}}_{W2}$. В силу теоремы 12.4.7 матрица $c$ имеет ${}_*{}^*$-обратную и из равенства (12.4.16) следует

$$(12.4.21) \qquad \overline{\overline{e}}_{W2} = e_{W1}{}^*{}_* \overline{\overline{c}}^{-1}{}_*{}^*$$

Доказательство. Если даны два автоморфизма $\overline{f}$ и $\overline{g}$, то мы можем записать

$$v' = f^*{}_* v$$

$$v'' = g^*{}_* v' = g^*{}_* f^*{}_* v$$

Следовательно, результирующий автоморфизм имеет матрицу $f^*{}_* g$. $\quad\square$

Теорема 12.4.14. *Пусть* $\overline{\overline{e}}_{V1}$ *,* $\overline{\overline{e}}_{V2}$ *- базисы правого $A$-векторного пространства $V$ столбцов. Пусть* $\overline{\overline{e}}_{W1}$ *,* $\overline{\overline{e}}_{W2}$ *- базисы правого $A$-векторного пространства $W$ столбцов. Пусть $a_1$ - матрица гомоморфизма*

$$(12.4.17) \qquad a : V \to W$$

*относительно базисов* $\overline{\overline{e}}_{V1}$*,* $\overline{\overline{e}}_{W1}$ *и $a_2$ - матрица гомоморфизма* (12.4.17) *относительно базисов* $\overline{\overline{e}}_{V2}$*,* $\overline{\overline{e}}_{W2}$*. Если базис* $\overline{\overline{e}}_{V1}$ *имеет координатную матрицу $b$ относительно базиса* $\overline{\overline{e}}_{V2}$

$$\overline{\overline{e}}_{V1} = \overline{\overline{e}}_{V2}{}^*{}_* b$$

*и* $\overline{\overline{e}}_{W1}$ *имеет координатную матрицу $c$ относительно базиса* $\overline{\overline{e}}_{W2}$

$$(12.4.18) \qquad \overline{\overline{e}}_{W1} = \overline{\overline{e}}_{W2}{}^*{}_* c$$

*то матрицы $a_1$ и $a_2$ связаны соотношением*

$$a_1 = c^{-1}{}^*{}_*{}^* a_2{}^*{}_* b$$

Доказательство. Вектор $\overline{v} \in V$ имеет разложение

$$\overline{v} = e_{V1*}{}^* v = e_{V2*}{}^* b_*{}^* v$$

относительно базисов $\overline{\overline{e}}_{V1}$, $\overline{\overline{e}}_{V2}$. Так как $a$ - гомоморфизм, то мы можем записать его в форме

$$(12.4.24) \qquad \overline{w} = e_{W1*}{}^* a_1{}_*{}^* v$$

относительно базисов $\overline{\overline{e}}_{V1}$, $\overline{\overline{e}}_{W1}$ и

$$(12.4.25) \qquad \overline{w} = e_{W2*}{}^* a_{2*}{}^* b_*{}^* v$$

относительно базисов $\overline{\overline{e}}_{V2}$, $\overline{\overline{e}}_{W2}$. В силу теоремы 12.4.8 матрица $c$ имеет ${}_*{}^*$-обратную и из равенства (12.4.18) следует

$$(12.4.26) \qquad \overline{\overline{e}}_{W2} = e_{W1*}{}^* \overline{\overline{c}}^{-1}{}^*{}_*$$

Равенство

$$(12.4.22) \quad \overline{w} = e_{W1*}{}^* c^{-1*}{}^*{}_*{}^* a_{2*}{}^* b_*{}^* v$$

является следствием равенств (12.4.20), (12.4.21). В силу теоремы 8.4.18 сравнение равенств (12.4.19) и (12.4.22) даёт

$$(12.4.23) \quad v^*{}_* a_1 = c^{-1*}{}^*{}_*{}^* a_{2*}{}^* b_*{}^* v$$

Так как вектор $a$ - произвольный вектор, утверждение теоремы является следствием теоремы 12.5.1 и равенства (12.4.23).

<div style="text-align:right">□</div>

---

**Теорема 12.4.15.** *Пусть $\overline{a}$ - автоморфизм левого $A$-векторного пространства $V$ столбцов. Пусть $a_1$ - матрица этого автоморфизма, заданная относительно базиса $\overline{\overline{e}}_1$, и $a_2$ - матрица того же автоморфизма, заданная относительно базиса $\overline{\overline{e}}_2$. Если базис $\overline{\overline{e}}_1$ имеет координатную матрицу $b$ относительно базиса $\overline{\overline{e}}_2$*

$$\overline{\overline{e}}_1 = \overline{\overline{e}}_{2*}{}^* b$$

*то матрицы $a_1$ и $a_2$ связаны соотношением*

$$a_1 = b^{-1*}{}^*{}_*{}^* a_{2*}{}^* b$$

Доказательство. Утверждение является следствием теоремы 12.4.13, так как в данном случае $c = b$. □

---

Равенство

$$(12.4.27) \quad \overline{w} = e_{W1*}{}^* c^{-1*}{}^*{}_*{}^* a_{2}{}^*{}_* b^*{}_* v$$

является следствием равенств (12.4.25), (12.4.26). В силу теоремы 8.4.25 сравнение равенств (12.4.24) и (12.4.27) даёт

$$(12.4.28) \quad v^*{}_* a_1 = c^{-1*}{}^*{}_*{}^* a_2{}^*{}_* b^*{}_* v$$

Так как вектор $a$ - произвольный вектор, утверждение теоремы является следствием теоремы 12.5.2 и равенства (12.4.28).

<div style="text-align:right">□</div>

---

**Теорема 12.4.16.** *Пусть $\overline{a}$ - автоморфизм левого $A$-векторного пространства $V$ столбцов. Пусть $a_1$ - матрица этого автоморфизма, заданная относительно базиса $\overline{\overline{e}}_1$, и $a_2$ - матрица того же автоморфизма, заданная относительно базиса $\overline{\overline{e}}_2$. Если базис $\overline{\overline{e}}_1$ имеет координатную матрицу $b$ относительно базиса $\overline{\overline{e}}_2$*

$$\overline{\overline{e}}_1 = \overline{\overline{e}}_2{}^*{}_* b$$

*то матрицы $a_1$ и $a_2$ связаны соотношением*

$$a_1 = b^{-1*}{}^*{}_*{}^* a_2{}^*{}_* b$$

Доказательство. Утверждение является следствием теоремы 12.4.14, так как в данном случае $c = b$. □

---

## 12.5. Многообразие базисов правого $A$-векторного пространства

Пусть $V$ - правое $A$-векторное пространство столбцов. Мы доказали в теореме 12.4.9, что если задан базис правого $A$-векторного пространства $V$, то любой автоморфизм $\overline{f}$ правого $A$-векторного пространства $V$ можно отождествить с $_*{}^*$-невырожденной матрицей $f$. Соответствующее преобразование координат вектора

$$v' = f_*{}^* v$$

называется **линейным преобразованием**.

Пусть $V$ - правое $A$-векторное пространство строк. Мы доказали в теореме 12.4.10, что если задан базис правого $A$-векторного пространства $V$, то любой автоморфизм $\overline{f}$ правого $A$-векторного пространства $V$ можно отождествить с $^*{}_*$-невырожденной матрицей $f$. Соответствующее преобразование координат вектора

$$v' = f^*{}_* v$$

называется **линейным преобразованием**.

Теорема 12.5.1. *Пусть $V$ - правое $A$-векторное пространство столбцов и $\overline{\overline{e}}$ - базис правого $A$-векторного пространства $V$. Автоморфизмы правого $A$-векторного пространства столбцов порождают левостороннее линейное эффективное $GL(n_{*}A)$-представление.*

Доказательство. Пусть $a, b$ - матрицы автоморфизмов $\overline{a}$ и $\overline{b}$ относительно базиса $\overline{\overline{e}}$. Согласно теореме 11.3.6 преобразование координат имеет вид

$$(12.5.1) \qquad v' = a_{*}{}^{*}v$$

$$(12.5.2) \qquad v'' = b_{*}{}^{*}v'$$

Равенство

$$v'' = b_{*}{}^{*}a_{*}{}^{*}v$$

является следствием равенств (12.5.1), (12.5.2). Согласно теореме 11.5.2, произведение автоморфизмов $\overline{a}$ и $\overline{b}$ имеет матрицу $b_{*}{}^{*}a$. Следовательно, автоморфизмы правого $A$-векторного пространства столбцов порождают левостороннее линейное $GL(n_{*}A)$-представление.

Нам осталось показать, что ядро неэффективности состоит только из единичного элемента. Тождественное преобразование удовлетворяет равенству

$$v^i = a^i_j v^j$$

Выбирая значения координат в форме $v^i = \delta^i_k$, где $k$ задано, мы получим

$$(12.5.3) \qquad \delta^i_k = a^i_j \delta^j_k$$

Из (12.5.3) следует

$$\delta^i_k = a^i_k$$

Поскольку $k$ произвольно, мы приходим к заключению $a = \delta$. $\qquad \square$

Теорема 12.5.3. *Автоморфизм $a$, действуя на каждый вектор базиса в правом $A$-векторном пространстве столбцов , отображает базис в другой базис.*

Теорема 12.5.2. *Пусть $V$ - правое $A$-векторное пространство строк и $\overline{\overline{e}}$ - базис правого $A$-векторного пространства $V$. Автоморфизмы правого $A$-векторного пространства строк порождают левостороннее линейное эффективное $GL(n^{*}{}_{*}A)$-представление.*

Доказательство. Пусть $a, b$ - матрицы автоморфизмов $\overline{a}$ и $\overline{b}$ относительно базиса $\overline{\overline{e}}$. Согласно теореме 11.3.7 преобразование координат имеет вид

$$(12.5.4) \qquad v' = a^{*}{}_{*}v$$

$$(12.5.5) \qquad v'' = b^{*}{}_{*}v'$$

Равенство

$$v'' = b^{*}{}_{*}a^{*}{}_{*}v$$

является следствием равенств (12.5.4), (12.5.5). Согласно теореме 11.5.3, произведение автоморфизмов $\overline{a}$ и $\overline{b}$ имеет матрицу $b^{*}{}_{*}a$. Следовательно, автоморфизмы правого $A$-векторного пространства строк порождают левостороннее линейное $GL(n^{*}{}_{*}A)$-представление.

Нам осталось показать, что ядро неэффективности состоит только из единичного элемента. Тождественное преобразование удовлетворяет равенству

$$v_i = a^j_i v_j$$

Выбирая значения координат в форме $v_i = \delta^k_i$, где $k$ задано, мы получим

$$(12.5.6) \qquad \delta^i_k = a^i_j \delta^k_j$$

Из (12.5.6) следует

$$\delta^k_i = a^k_i$$

Поскольку $k$ произвольно, мы приходим к заключению $a = \delta$. $\qquad \square$

Теорема 12.5.4. *Автоморфизм $a$, действуя на каждый вектор базиса в правом $A$-векторном пространстве строк , отображает базис в другой базис.*

Доказательство. Пусть $\overline{\overline{e}}$ - базис в правом $A$-векторном пространстве $V$ столбцов . Согласно теореме 12.4.9, вектор $e_i$ отображается в вектор $e'_i$

(12.5.7)          $e'_i = a_* {}^* e_i$

Если векторы $e'_i$ линейно зависимы, то в линейной комбинации

(12.5.8)          $e'_* {}^* \lambda = 0$

$\lambda \neq 0$. Из равенств (12.5.7) и (12.5.8) следует, что

$$a^{-1 *} {}_* {}^* e'_* {}^* \lambda = e_* {}^* \lambda = 0$$

и $\lambda \neq 0$. Это противоречит утверждению, что векторы $e_i$ линейно независимы. Следовательно, векторы $e'_i$ линейно независимы и порождают базис. $\square$

Доказательство. Пусть $\overline{\overline{e}}$ - базис в правом $A$-векторном пространстве $V$ строк . Согласно теореме 12.4.10, вектор $e^i$ отображается в вектор $e'^i$

(12.5.9)          $e'^i = a^* {}_* e^i$

Если векторы $e'^i$ линейно зависимы, то в линейной комбинации

(12.5.10)          $e'^* {}_* \lambda = 0$

$\lambda \neq 0$. Из равенств (12.5.9) и (12.5.10) следует, что

$$a^{-1 *} {}^* {}_* e'^* {}_* \lambda = e^* {}_* \lambda = 0$$

и $\lambda \neq 0$. Это противоречит утверждению, что векторы $e^i$ линейно независимы. Следовательно, векторы $e'^i$ линейно независимы и порождают базис. $\square$

Таким образом, мы можем распространить левостороннее линейное $GL(V_*)$-представление в правом $A$-векторном пространстве $V_*$ на множество базисов правого $A$-векторного пространства $V$. Мы будем называть преобразование этого левостороннего представления на множестве базисов правого $A$-векторного пространства $V$ **активным преобразованием** потому, что гоморфизм правого $A$-векторного пространства породил это преобразование (Смотри также определение в разделе [4]-14.1-3, а также определение на странице [2]-214).

Соответственно определению мы будем записывать действие активного преобразования $a \in GL(V_*)$ на базис $\overline{\overline{e}}$ в форме $a_* {}^* \overline{\overline{e}}$. Рассмотрим равенство

(12.5.11)          $a_* {}^* \overline{\overline{e}}_* {}^* v = a_* {}^* \overline{\overline{e}}_* {}^* v$

Выражение $a_* {}^* \overline{\overline{e}}$ в левой части равенства (12.5.11) является образом базиса $\overline{\overline{e}}$ при активном преобразовании $a$. Выражение $\overline{\overline{e}}_* {}^* v$ в правой части равенства (12.5.11) является разложением вектора $\overline{v}$ относительно базиса $\overline{\overline{e}}$. Следовательно, выражение в правой части равенства (12.5.11) является образом вектора $\overline{v}$ относительно эндоморфизма $a$ и выражение в левой части равенства

Соответственно определению мы будем записывать действие активного преобразования $a \in GL(V_*)$ на базис $\overline{\overline{e}}$ в форме $a^* {}_* \overline{\overline{e}}$. Рассмотрим равенство

(12.5.12)          $a^* {}_* \overline{\overline{e}}^* {}_* v = a^* {}_* \overline{\overline{e}}^* {}_* v$

Выражение $a^* {}_* \overline{\overline{e}}$ в левой части равенства (12.5.12) является образом базиса $\overline{\overline{e}}$ при активном преобразовании $a$. Выражение $\overline{\overline{e}}^* {}_* v$ в правой части равенства (12.5.12) является разложением вектора $\overline{v}$ относительно базиса $\overline{\overline{e}}$. Следовательно, выражение в правой части равенства (12.5.12) является образом вектора $\overline{v}$ относительно эндоморфизма $a$ и выражение в левой части равенства

| (12.5.11)          $a_* {}^* \overline{\overline{e}}_* {}^* v = a_* {}^* \overline{\overline{e}}_* {}^* v$ |
|---|

является разложением образа вектора $\overline{v}$ относительно образа базиса $\overline{\overline{e}}$. Следовательно, из равенства (12.5.11) следует, что эндоморфизм $a$ правого

| (12.5.12)          $a^* {}_* \overline{\overline{e}}^* {}_* v = a^* {}_* \overline{\overline{e}}^* {}_* v$ |
|---|

является разложением образа вектора $\overline{v}$ относительно образа базиса $\overline{\overline{e}}$. Следовательно, из равенства (12.5.12) следует, что эндоморфизм $a$ правого $A$-векторного пространства и соответствующее активное преобразование $a$ действуют синхронно и координаты $a \circ \overline{v}$ образа вектора $\overline{v}$ относительно образа $a_* \overline{\overline{e}}$ базиса $\overline{\overline{e}}$ совпадают с координатами вектора $\overline{v}$ относительно базиса $\overline{\overline{e}}$.

$A$-векторного пространства и соответствующее активное преобразование $a$ действуют синхронно и координаты $a \circ \overline{v}$ образа вектора $\overline{v}$ относительно образа $a^*_* \overline{\overline{e}}$ базиса $\overline{\overline{e}}$ совпадают с координатами вектора $\overline{v}$ относительно базиса $\overline{\overline{e}}$.

**Теорема 12.5.5.** *Активное левостороннее $GL(n_* {}_* A)$-представление на множестве базисов однотранзитивно. Множество базисов, отождествлённое с траекторией $GL(n_* {}_* A)(V)_* {}^* \overline{\overline{e}}$ активного левостороннего $GL(n_* {}_* A)$-представления называется* **многообразием базисов** *правого $A$-векторного пространства $V$.*

**Теорема 12.5.6.** *Активное левостороннее $GL(n^* {}_* A)$-представление на множестве базисов однотранзитивно. Множество базисов, отождествлённое с траекторией $GL(n^* {}_* A)(V)^* {}_* \overline{\overline{e}}$ активного левостороннего $GL(n^* {}_* A)$-представления называется* **многообразием базисов** *правого $A$-векторного пространства $V$.*

Чтобы доказать теоремы 12.5.5, 12.5.6, достаточно показать, что для любых двух базисов определено по крайней мере одно преобразование левостороннего представления и это преобразование единственно.

**Доказательство.** Гомоморфизм $a$, действуя на базис $\overline{e}$ имеет вид

$$e'_i = a_* {}_* e_i$$

где $e'_i$ - координатная матрица вектора $\overline{e}'_i$ относительно базиса $\overline{\overline{h}}$ и $e_i$ - координатная матрица вектора $\overline{e}_i$ относительно базиса $\overline{\overline{h}}$. Следовательно, координатная матрица образа базиса равна $_*{}^*$-произведению матрицы автоморфизма и координатной матрицы исходного базиса

$$e' = g_* {}^* e$$

В силу теоремы 12.4.7, матрицы $g$ и $e$ невырождены. Следовательно, матрица

$$a = e'_* {}^* e^{-1_*{}^*}$$

является матрицей автоморфизма, отображающего базис $\overline{\overline{e}}$ в базис $\overline{\overline{e}}'$.

Допустим элементы $g_1$, $g_2$ группы $G$ и базис $\overline{\overline{e}}$ таковы, что

$$g_{1*} {}^* e = g_{2*} {}^* e$$

В силу теорем 12.4.7 и 8.1.32 справедливо равенство $g_1 = g_2$. Отсюда следует утверждение теоремы. □

**Доказательство.** Гомоморфизм $a$, действуя на базис $\overline{e}$ имеет вид

$$e'^i = a^* {}_* e^i$$

где $e'^i$ - координатная матрица вектора $\overline{e}'^i$ относительно базиса $\overline{\overline{h}}$ и $e^i$ - координатная матрица вектора $\overline{e}^i$ относительно базиса $\overline{\overline{h}}$. Следовательно, координатная матрица образа базиса равна $_*{}^*$-произведению матрицы автоморфизма и координатной матрицы исходного базиса

$$e' = g^* {}_* e$$

В силу теоремы 12.4.8, матрицы $g$ и $e$ невырождены. Следовательно, матрица

$$a = e'^* {}_* e^{-1^* {}_*}$$

является матрицей автоморфизма, отображающего базис $\overline{\overline{e}}$ в базис $\overline{\overline{e}}'$.

Допустим элементы $g_1$, $g_2$ группы $G$ и базис $\overline{\overline{e}}$ таковы, что

$$g_1 {}^* e = g_2 {}^* e$$

В силу теорем 12.4.8 и 8.1.32 справедливо равенство $g_1 = g_2$. Отсюда следует утверждение теоремы. □

Если в правом $A$-векторном пространстве $V$ определена дополнительная структура, не всякий автоморфизм сохраняет свойства заданной структуры. Например, если мы определим норму в правом $A$-векторном пространстве $V$, то для нас интересны автоморфизмы, которые сохраняют норму вектора.

ОПРЕДЕЛЕНИЕ 12.5.7. *Нормальная подгруппа* $G(V)$ *группы* $GL(V)$, *которая порождает автоморфизмы, сохраняющие свойства заданной структуры называется* **группой симметрии**.

*Не нарушая общности, мы будем отождествлять элемент $g$ группы $G(V)$ с соответствующим автоморфизмом и записывать его действие на вектор $v \in V$ в виде $g_* {}^* v$.* □

ОПРЕДЕЛЕНИЕ 12.5.8. *Нормальная подгруппа* $G(V)$ *группы* $GL(V)$, *которая порождает автоморфизмы, сохраняющие свойства заданной структуры называется* **группой симметрии**.

*Не нарушая общности, мы будем отождествлять элемент $g$ группы $G(V)$ с соответствующим автоморфизмом и записывать его действие на вектор $v \in V$ в виде $g^* {}_* v$.* □

Мы будем называть левостороннее представление группы $G$ в правом $A$-векторном пространстве **линейным $G$-представлением**.

ОПРЕДЕЛЕНИЕ 12.5.9. *Мы будем называть левостороннее представление группы $G(V)$ на множестве базисов правого $A$-векторного пространства $V$* **активным правосторонним $G$-представлением**. □

Если базис $\overline{\overline{e}}$ задан, то мы можем отождествить автоморфизм правого $A$-векторного пространства $V$ и его координаты относительно базиса $\overline{\overline{e}}$. Множество $G(V_*)$ координат автоморфизмов относительно базиса $\overline{\overline{e}}$ является группой, изоморфной группе $G(V)$.

Не всякие два базиса могут быть связаны преобразованием группы симметрии потому, что не всякое невырожденное линейное преобразование принадлежит представлению группы $G(V)$. Таким образом, множество базисов можно представить как объединение орбит группы $G(V)$. В частности, если базис $\overline{\overline{e}} \in G(V)$, то орбита базиса $\overline{\overline{e}}$ совпадает с группой $G(V)$.

ОПРЕДЕЛЕНИЕ 12.5.10. *Мы будем называть орбиту $G(V)_* {}^* \overline{\overline{e}}$ выбранного базиса $\overline{\overline{e}}$* **многообразием базисов** *правого $A$-векторного пространства $V$ столбцов* . □

ОПРЕДЕЛЕНИЕ 12.5.11. *Мы будем называть орбиту $G(V)^* {}_* \overline{\overline{e}}$ выбранного базиса $\overline{\overline{e}}$* **многообразием базисов** *правого $A$-векторного пространства $V$ строк* . □

ТЕОРЕМА 12.5.12. *Активное левостороннее $G(V)$-представление на многообразии базисов однотранзитивно.*

ДОКАЗАТЕЛЬСТВО. Теорема является следствием теоремы 12.5.5, и определений 12.5.10, а также теорема является следствием теоремы 12.5.6 и определений 12.5.11. □

Теорема 12.5.13. *На многообразии базисов $G(V)_*{}^{\overline{\overline{e}}}$ правого $A$-векторного пространства столбцов существует однотранзитивное правостороннее $G(V)$-представление, перестановочное с активным.*

Теорема 12.5.14. *На многообразии базисов $G(V)^*{}_{\overline{\overline{e}}}$ правого $A$-векторного пространства строк существует однотранзитивное правостороннее $G(V)$-представление, перестановочное с активным.*

Доказательство. Из теоремы 12.5.12 следует, что многообразие базисов $G(V)_*{}^{\overline{\overline{e}}}$ является однородным пространством группы $G$. Теорема является следствием теоремы 2.5.17. $\square$

Доказательство. Из теоремы 12.5.12 следует, что многообразие базисов $G(V)^*{}_{\overline{\overline{e}}}$ является однородным пространством группы $G$. Теорема является следствием теоремы 2.5.17. $\square$

Как мы видим из замечания 2.5.18 преобразование правостороннего $G(V)$-представления отличается от активного преобразования и не может быть сведено к преобразованию пространства $V$.

Определение 12.5.15. *Преобразование правостороннего $G(V)$-представления называется* **пассивным преобразованием** *многообразия базисов $G(V)_*{}^{\overline{\overline{e}}}$ правого $A$-векторного пространства столбцов , а правостороннее $G(V)$-представление называется* **пассивным правосторонним** $G(V)$**-представлением**. *Согласно определению мы будем записывать пассивное преобразование базиса $\overline{e}$, порождённое элементом $a \in G(V)$, в форме $\overline{\overline{e}}_*{}^*a$.* $\square$

Определение 12.5.16. *Преобразование правостороннего $G(V)$-представления называется* **пассивным преобразованием** *многообразия базисов $G(V)^*{}_{\overline{\overline{e}}}$ правого $A$-векторного пространства строк , а правостороннее $G(V)$-представление называется* **пассивным правосторонним** $G(V)$**-представлением**. *Согласно определению мы будем записывать пассивное преобразование базиса $\overline{e}$, порождённое элементом $a \in G(V)$, в форме $\overline{\overline{e}}_*{}^*a$.* $\square$

Теорема 12.5.17. *Координатная матрица базиса $\overline{\overline{e}}'$ относительно базиса $\overline{\overline{e}}$ правого $A$-векторного пространства $V$ столбцов совпадает с матрицей пассивного преобразования, отображающего базис $\overline{\overline{e}}$ в базис $\overline{\overline{e}}'$.*

Теорема 12.5.18. *Координатная матрица базиса $\overline{\overline{e}}'$ относительно базиса $\overline{\overline{e}}$ правого $A$-векторного пространства $V$ строк совпадает с матрицей пассивного преобразования, отображающего базис $\overline{\overline{e}}$ в базис $\overline{\overline{e}}'$.*

Доказательство. Согласно теореме 8.4.17, координатная матрица базиса $\overline{\overline{e}}'$ относительно базиса $\overline{\overline{e}}$ состоит из столбцов, являющихся координатными матрицами векторов $\overline{e}'_i$ относительно базиса $\overline{\overline{e}}$. Следовательно,

$$\overline{e}'_i = \overline{e}_*{}^*e'_i$$

В тоже время пассивное преобразование $a$, связывающее два базиса, имеет

Доказательство. Согласно теореме 8.4.24, координатная матрица базиса $\overline{\overline{e}}'$ относительно базиса $\overline{\overline{e}}$ состоит из строк, являющихся координатными матрицами векторов $\overline{e}'_i$ относительно базиса $\overline{\overline{e}}$. Следовательно,

$$\overline{e}'^i = \overline{e}_*{}^*e'^i$$

В тоже время пассивное преобразование $a$, связывающее два базиса, имеет

вид

$$\overline{e}'_i = \overline{e}_{*}{}^{*} a_i$$

В силу теоремы 8.4.18,

$$e'_i = a_i$$

для любого $i$. Это доказывает теорему. □

Замечание 12.5.19. *Согласно теоремам 2.5.15, 12.5.17, мы можем отождествить базис $\overline{\overline{e}}'$ многообразия базисов $G(V)_{*}{}_{*}\overline{\overline{e}}$ правого $A$-векторного пространства столбцов и матрицу $g$ координат базиса $\overline{\overline{e}}'$ относительно базиса $\overline{\overline{e}}$*

(12.5.13)        $\overline{\overline{e}}' = \overline{\overline{e}}_{*}{}^{*} g$

*Если $\overline{\overline{e}}' = \overline{\overline{e}}$, то равенство (12.5.13) принимает вид*

(12.5.14)        $\overline{\overline{e}} = \overline{\overline{e}}_{*}{}^{*} E_n$

*Опираясь на равенство (12.5.14) мы будем пользоваться записью $G(V)_{*}{}^{*}\overline{\overline{E_n}}$ для многообразия базисов $G(V)_{*}{}^{*}\overline{\overline{e}}$.* □

вид

$$\overline{e}'{}^{i} = \overline{e}^{*}{}_{*} a^{i}$$

В силу теоремы 8.4.25,

$$e'{}^{i} = a^{i}$$

для любого $i$. Это доказывает теорему. □

Замечание 12.5.20. *Согласно теоремам 2.5.15, 12.5.18, мы можем отождествить базис $\overline{\overline{e}}'$ многообразия базисов $G(V)^{*}{}_{*}\overline{\overline{e}}$ правого $A$-векторного пространства столбцов и матрицу $g$ координат базиса $\overline{\overline{e}}'$ относительно базиса $\overline{\overline{e}}$*

(12.5.15)        $\overline{\overline{e}}' = \overline{\overline{e}}^{*}{}_{*} g$

*Если $\overline{\overline{e}}' = \overline{\overline{e}}$, то равенство (12.5.15) принимает вид*

(12.5.16)        $\overline{\overline{e}} = \overline{\overline{e}}^{*}{}_{*} E_n$

*Опираясь на равенство (12.5.16) мы будем пользоваться записью $G(V)^{*}{}_{*}\overline{\overline{E_n}}$ для многообразия базисов $G(V)^{*}{}_{*}\overline{\overline{e}}$.* □

## 12.6. Пассивное преобразование в правом $A$-векторном пространстве

Активное преобразование изменяет базисы и векторы согласовано и координаты вектора относительно базиса не меняются. Пассивное преобразование меняет только базис, и это ведёт к преобразованию координат вектора относительно базиса.

Теорема 12.6.1. *Пусть $V$ - правое $A$-векторное пространство столбцов. Пусть пассивное преобразование $g \in G(V_{*})$ отображает базис $\overline{\overline{e}}_1$ в базис $\overline{\overline{e}}_2$*

(12.6.1)        $\overline{\overline{e}}_2 = \overline{\overline{e}}_{1*}{}^{*} g$

*Пусть*

(12.6.2)        $v_i = \begin{pmatrix} v_i^1 \\ ... \\ v_i^n \end{pmatrix}$

*матрица координат вектора $\overline{v}$ относительно базиса $\overline{\overline{e}}_i$, $i = 1, 2$. Преоб-*

Теорема 12.6.2. *Пусть $V$ - правое $A$-векторное пространство строк. Пусть пассивное преобразование $g \in G(V_{*})$ отображает базис $\overline{\overline{e}}_1$ в базис $\overline{\overline{e}}_2$*

(12.6.5)        $\overline{\overline{e}}_2 = \overline{\overline{e}}_1{}^{*}{}_{*} g$

*Пусть*

$$v_i = \begin{pmatrix} v_{i1} & ... & v_{in} \end{pmatrix}$$

*матрица координат вектора $\overline{v}$ относительно базиса $\overline{\overline{e}}_i$, $i = 1, 2$. Преобразования координат*

(12.6.6)        $v_1 = g^{*}{}_{*} v_2$

*разование координат*

$$(12.6.3) \qquad v_1 = g_*{}^* v_2$$

$$(12.6.4) \qquad v_2 = g^{-1}{}_*{}^*{}^* v_1$$

*не зависят от вектора $\overline{v}$ или базиса $\overline{\overline{e}}$, а определенно исключительно координатами вектора $\overline{v}$ относительно базиса $\overline{\overline{e}}$.*

ДОКАЗАТЕЛЬСТВО. Согласно теореме 8.4.17,

$$(12.6.8) \qquad \overline{v} = \overline{\overline{e}}_{1*}{}^* v_1 = \overline{\overline{e}}_{2*}{}^* v_2$$

Равенство

$$(12.6.9) \qquad \overline{\overline{e}}_{1*}{}^* v_1 = \overline{\overline{e}}_{1*}{}^* g_*{}^* v_2$$

является следствием равенств (12.6.1), (12.6.8). Согласно теореме 8.4.18, преобразование координат (12.6.3) является следствием равенства (12.6.9). Так как $g$ - $_*$*-невырожденная матрица, то преобразование координат (12.6.4) является следствием равенства (12.6.3). □

ТЕОРЕМА 12.6.3. *Преобразования координат*

$$\boxed{(12.6.4) \quad v_2 = g^{-1}{}_*{}^*{}^* v_1}$$

*порождают эффективное линейное левостороннее контравариантное $G$-представление, называемое* **координатным представлением в правом $A$-векторном пространстве** *столбцов.*

ДОКАЗАТЕЛЬСТВО. Пусть

$$(12.6.12) \qquad v_i = \begin{pmatrix} v_i^1 \\ ... \\ v_i^n \end{pmatrix}$$

матрица координат вектора $\overline{v}$ относительно базиса $\overline{\overline{e}}_i$, $i = 1, 2, 3$. Тогда

$$(12.6.13) \qquad \begin{aligned} \overline{v} &= \overline{\overline{e}}_{1*}{}^* v_1 = \overline{\overline{e}}_{2*}{}^* v_2 \\ &= \overline{\overline{e}}_{3*}{}^* v_3 \end{aligned}$$

Пусть пассивное преобразование $g_1$ отображает базис $\overline{\overline{e}}_1 \in G(V)_*{}^* \overline{\overline{e}}$ в ба-

---

$$(12.6.7) \qquad v_2 = g^{-1}{}^*{}^*{}_* v_1$$

*не зависят от вектора $\overline{v}$ или базиса $\overline{\overline{e}}$, а определенно исключительно координатами вектора $\overline{v}$ относительно базиса $\overline{\overline{e}}$.*

ДОКАЗАТЕЛЬСТВО. Согласно теореме 8.4.24,

$$(12.6.10) \qquad \overline{v} = \overline{\overline{e}}_1{}^*{}_* v_1 = \overline{\overline{e}}_2{}^*{}_* v_2$$

Равенство

$$(12.6.11) \qquad \overline{\overline{e}}_1{}^*{}_* v_1 = \overline{\overline{e}}_1{}^*{}_* g^*{}_* v_2$$

является следствием равенств (12.6.5), (12.6.10). Согласно теореме 8.4.25, преобразование координат (12.6.6) является следствием равенства (12.6.11). Так как $g$ - $^*{}_*$-невырожденная матрица, то преобразование координат (12.6.7) является следствием равенства (12.6.6). □

ТЕОРЕМА 12.6.4. *Преобразования координат*

$$\boxed{(12.6.7) \quad v_2 = g^{-1}{}^*{}^*{}_* v_1}$$

*порождают эффективное линейное левостороннее контравариантное $G$-представление, называемое* **координатным представлением в правом $A$-векторном пространстве** *строк.*

ДОКАЗАТЕЛЬСТВО. Пусть

$$v_i = \begin{pmatrix} v_{i1} & ... & v_{in} \end{pmatrix}$$

матрица координат вектора $\overline{v}$ относительно базиса $\overline{\overline{e}}_i$, $i = 1, 2, 3$. Тогда

$$(12.6.21) \qquad \begin{aligned} \overline{v} &= \overline{\overline{e}}_1{}^*{}_* v_1 = \overline{\overline{e}}_2{}^*{}_* v_2 \\ &= \overline{\overline{e}}_3{}^*{}_* v_3 \end{aligned}$$

Пусть пассивное преобразование $g_1$ отображает базис $\overline{\overline{e}}_1 \in G(V)^*{}_* \overline{\overline{e}}$ в ба-зис $\overline{\overline{e}}_2 \in G(V)_*{}^*\overline{\overline{e}}$

$$(12.6.14) \qquad \overline{\overline{e}}_2 = \overline{\overline{e}}_{1*}{}^* g_1$$

Пусть пассивное преобразование $g_2$ отображает базис $\overline{\overline{e}}_2 \in G(V)_*{}^*\overline{\overline{e}}$ в базис $\overline{\overline{e}}_3 \in G(V)_*{}^*\overline{\overline{e}}$

$$(12.6.15) \qquad \overline{\overline{e}}_3 = \overline{\overline{e}}_{2*}{}^* g_2$$

Преобразование

$$(12.6.16) \qquad \overline{\overline{e}}_3 = \overline{\overline{e}}_{1*}{}^* (g_{1*}{}^* g_2)$$

является произведением преобразований (12.6.14), (12.6.15). Преобразование координат

$$(12.6.17) \qquad v_2 = g_1^{-1*}{}^*{}_* v_1$$

соответствует пассивному преобразованию $g_1$ и является следствием равенств (12.6.13), (12.6.14) и теоремы 12.6.1. Аналогично, согласно теореме 12.6.1, преобразование координат

$$(12.6.18) \qquad v_3 = g_2^{-1*}{}^*{}_* v_2$$

соответствует пассивному преобразованию $g_2$ и является следствием равенств (12.6.13), (12.6.15). Произведение преобразований координат (12.3.17), (12.3.18) имеет вид

$$(12.6.19) \qquad v_3 = g_2^{-1*}{}^*{}_* g_1^{-1*}{}^*{}_* v_1$$

и является координатным преобразованием, соответствующим пассивному преобразованию (12.6.16). Равенство

$$(12.6.20) \qquad v_3 = (g_{2*}{}^* g_1)^{-1*}{}^*{}_* v_1$$

является следствием равенств (8.1.66), (12.6.19). Согласно определению 2.5.9, преобразования координат порождают линейное правостороннее контравариантное $G$-представление. Если координатное преобразование не изменяет векторы $\delta_k$, то ему соответствует единица группы $G$, так как пассивное представление однотранзитивно. Следовательно, координатное представление эффективно. □

зис $\overline{\overline{e}}_2 \in G(V)^*{}_*\overline{\overline{e}}$

$$(12.6.22) \qquad \overline{\overline{e}}_2 = \overline{\overline{e}}_1{}^*{}_* g_1$$

Пусть пассивное преобразование $g_2$ отображает базис $\overline{\overline{e}}_2 \in G(V)^*{}_*\overline{\overline{e}}$ в базис $\overline{\overline{e}}_3 \in G(V)^*{}_*\overline{\overline{e}}$

$$(12.6.23) \qquad \overline{\overline{e}}_3 = \overline{\overline{e}}_2{}^*{}_* g_2$$

Преобразование

$$(12.6.24) \qquad \overline{\overline{e}}_3 = \overline{\overline{e}}_1{}^*{}_* (g_1{}^*{}_* g_2)$$

является произведением преобразований (12.6.22), (12.6.23). Преобразование координат

$$(12.6.25) \qquad v_2 = g_1^{-1*}{}^*{}_* v_1$$

соответствует пассивному преобразованию $g_1$ и является следствием равенств (12.6.21), (12.6.22) и теоремы 12.6.2. Аналогично, согласно теореме 12.6.2, преобразование координат

$$(12.6.26) \qquad v_3 = g_2^{-1*}{}^*{}_* v_2$$

соответствует пассивному преобразованию $g_2$ и является следствием равенств (12.6.21), (12.6.23). Произведение преобразований координат (12.3.17), (12.3.18) имеет вид

$$(12.6.27) \qquad v_3 = g_2^{-1*}{}^*{}_* g_1^{-1*}{}^*{}_* v_1$$

и является координатным преобразованием, соответствующим пассивному преобразованию (12.6.24). Равенство

$$(12.6.28) \qquad v_3 = (g_2{}^*{}_* g_1)^{-1*}{}^*{}_* v_1$$

является следствием равенств (8.1.66), (12.6.27). Согласно определению 2.5.9, преобразования координат порождают линейное правостороннее контравариантное $G$-представление. Если координатное преобразование не изменяет векторы $\delta_k$, то ему соответствует единица группы $G$, так как пассивное представление однотранзитивно. Следовательно, координатное представление эффективно. □

Теорема 12.6.5. *Пусть* $V$ *- правое* $A$*-векторное пространство*

Теорема 12.6.6. *Пусть* $V$ *- правое* $A$*-векторное пространство строк*

столбцов и $\overline{\overline{e}}_1$, $\overline{\overline{e}}_2$ - базисы в правом $A$-векторном пространстве $V$. Пусть базисы $\overline{\overline{e}}_1$, $\overline{\overline{e}}_2$ связаны пассивным преобразованием $g$

$$\overline{\overline{e}}_2 = \overline{\overline{e}}_{1*}{}^* g$$

Пусть $f$ - эндоморфизм правого $A$-векторного пространства $V$. Пусть $f_i$, $i = 1, 2$, - матрица эндоморфизма $f$ относительно базиса $\overline{\overline{e}}_i$. Тогда

$$(12.6.29) \qquad f_2 = g^{-1*}{}^*{}_* {}^* f_{1*}{}^* g$$

и $\overline{\overline{e}}_1$, $\overline{\overline{e}}_2$ - базисы в правом $A$-векторном пространстве $V$. Пусть базисы $\overline{\overline{e}}_1$, $\overline{\overline{e}}_2$ связаны пассивным преобразованием $g$

$$\overline{\overline{e}}_2 = \overline{\overline{e}}_1{}^*{}_* g$$

Пусть $f$ - эндоморфизм правого $A$-векторного пространства $V$. Пусть $f_i$, $i = 1, 2$, - матрица эндоморфизма $f$ относительно базиса $\overline{\overline{e}}_i$. Тогда

$$(12.6.30) \qquad f_2 = g^{-1*}{}^*{}_*{}^* f_1{}^*{}_* g$$

Доказательство теоремы 12.6.5. Пусть эндоморфизм $f$ отображает вектор $v$ в вектор $w$

$$(12.6.31) \qquad w = f \circ v$$

Пусть $v_i$, $i = 1, 2$, - матрица вектора $v$ относительно базиса $\overline{\overline{e}}_i$. Тогда

$$(12.6.32) \qquad v_1 = g_*{}^* v_2$$

следует из теоремы 12.6.1.

Пусть $w_i$, $i = 1, 2$, - матрица вектора $w$ относительно базиса $\overline{\overline{e}}_i$. Тогда

$$(12.6.33) \qquad w_1 = g_*{}^* w_2$$

следует из теоремы 12.6.1.

Согласно теореме 11.3.6, равенства

$$(12.6.34) \qquad w_1 = f_{1*}{}^* v_1$$

$$(12.6.35) \qquad w_2 = f_{2*}{}^* v_2$$

являются следствием равенства (12.6.31). Равенство

$$(12.6.36) \qquad g_*{}^* w_2 = f_{1*}{}^* v_1 = f_{1*}{}^* g_*{}^* v_2$$

является следствием равенств (12.6.32), (12.6.34), (12.6.33). Так как $g$ - невырожденная матрица, то равенство

$$(12.6.37) \qquad w_2 = f_{1*}{}^* v_1 = g^{-1*}{}^*{}_*{}^* f_{1*}{}^* g_*{}^* v_2$$

является следствием равенства (12.6.36). Равенство

$$(12.6.38) \qquad f_{2*}{}^* v_2 = g^{-1*}{}^*{}_*{}^* f_{1*}{}^* g_*{}^* v_2$$

является следствием равенств (12.6.35), (12.6.37). Равенство (12.6.29) является следствием равенства (12.6.38). $\qquad\square$

Доказательство теоремы 12.6.6. Пусть эндоморфизм $f$ отображает вектор $v$ в вектор $w$

$$(12.6.39) \qquad w = f \circ v$$

Пусть $v_i$, $i = 1, 2$, - матрица вектора $v$ относительно базиса $\overline{\overline{e}}_i$. Тогда

$$(12.6.40) \qquad v_1 = g^*{}_* v_2$$

следует из теоремы 12.6.2.

Пусть $w_i$, $i = 1, 2$, - матрица вектора $w$ относительно базиса $\overline{\overline{e}}_i$. Тогда

$$(12.6.41) \qquad w_1 = g^*{}_* w_2$$

следует из теоремы 12.6.2.

Согласно теореме 11.3.6, равенства

$$(12.6.42) \qquad w_1 = f_1{}^*{}_* v_1$$

(12.6.43) $$w_2 = {f_2}^*{}_* v_2$$

являются следствием равенства (12.6.39). Равенство

(12.6.44) $$g^*{}_* w_2 = {f_1}^*{}_* v_1 = {f_1}^*{}_* g^*{}_* v_2$$

является следствием равенств (12.6.40), (12.6.42), (12.6.41). Так как $g$ - невырожденная матрица, то равенство

(12.6.45) $$w_2 = {f_1}^*{}_* v_1 = g^{-1^*}{}^*{}_* {f_1}^*{}_* g^*{}_* v_2$$

является следствием равенства (12.6.44). Равенство

(12.6.46) $${f_2}^*{}_* v_2 = g^{-1^*}{}^*{}_* {f_1}^*{}_* g^*{}_* v_2$$

является следствием равенств (12.6.43), (12.6.45). Равенство (12.6.30) является следствием равенства (12.6.46).            $\square$

| | |
|---|---|
| **Теорема 12.6.7.** *Пусть $V$ - правое $A$-векторное пространство столбцов. Пусть эндоморфизм $\overline{f}$ правого $A$-векторного пространства $V$ является суммой эндоморфизмов $\overline{f}_1$, $\overline{f}_2$. Матрица $f$ эндоморфизма $\overline{f}$ равна сумме матриц $f_1$, $f_2$ эндоморфизмов $\overline{f}_1$, $\overline{f}_2$ и это равенство не зависит от выбора базиса.* | **Теорема 12.6.8.** *Пусть $V$ - правое $A$-векторное пространство строк. Пусть эндоморфизм $\overline{f}$ правого $A$-векторного пространства $V$ является суммой эндоморфизмов $\overline{f}_1$, $\overline{f}_2$. Матрица $f$ эндоморфизма $\overline{f}$ равна сумме матриц $f_1$, $f_2$ эндоморфизмов $\overline{f}_1$, $\overline{f}_2$ и это равенство не зависит от выбора базиса.* |

**Доказательство.** Пусть $\overline{e}$, $\overline{e}'$ - базисы правого $A$-векторного пространства $V$ и связаны пассивным преобразованием $g$

$$e' = g_*{}^* e$$

Пусть $f$ - матрица эндоморфизма $\overline{f}$ относительно базиса $\overline{e}$. Пусть $f'$ - матрица эндоморфизма $\overline{f}$ относительно базиса $\overline{e}'$. Пусть $f_i$, $i = 1, 2$, - матрица эндоморфизма $\overline{f}_i$ относительно базиса $\overline{e}$. Пусть $f_i'$, $i = 1, 2$, - матрица эндоморфизма $\overline{f}_i$ относительно базиса $\overline{e}'$. Согласно теореме 11.4.2,

(12.6.47) $$f = f_1 + f_2$$

(12.6.48) $$f' = f_1' + f_2'$$

Согласно теореме 12.6.5,

(12.6.49) $$f' = g^{-1^*}{}^*{}_* f_*{}^* g$$

(12.6.50) $$f_1' = g^{-1^*}{}^*{}_* f_{1*}{}^* g$$

(12.6.51) $$f_2' = g^{-1^*}{}^*{}_* f_{2*}{}^* g$$

**Доказательство.** Пусть $\overline{e}$, $\overline{e}'$ - базисы правого $A$-векторного пространства $V$ и связаны пассивным преобразованием $g$

$$e' = g^*{}_* e$$

Пусть $f$ - матрица эндоморфизма $\overline{f}$ относительно базиса $\overline{e}$. Пусть $f'$ - матрица эндоморфизма $\overline{f}$ относительно базиса $\overline{e}'$. Пусть $f_i$, $i = 1, 2$, - матрица эндоморфизма $\overline{f}_i$ относительно базиса $\overline{e}$. Пусть $f_i'$, $i = 1, 2$, - матрица эндоморфизма $\overline{f}_i$ относительно базиса $\overline{e}'$. Согласно теореме 11.4.3,

(12.6.53) $$f = f_1 + f_2$$

(12.6.54) $$f' = f_1' + f_2'$$

Согласно теореме 12.6.6,

(12.6.55) $$f' = g^{-1^*}{}^*{}_* f^*{}_* g$$

(12.6.56) $$f_1' = g^{-1^*}{}^*{}_* f_1{}^*{}_* g$$

(12.6.57) $$f_2' = g^{-1^*}{}^*{}_* f_2{}^*{}_* g$$

Равенство

$$(12.6.52) \quad \begin{aligned} f' &= g^{-1*}{}_*{}^*f_*{}^*g \\ &= g^{-1*}{}_*{}^*(f_1+f_2)_*{}^*g \\ &= g^{-1*}{}_*{}^*f_{1*}{}^*g \\ &+ g^{-1*}{}_*{}^*f_{2*}{}^*g \\ &= f_1' + f_2' \end{aligned}$$

является следствием равенств (12.6.47), (12.6.49), (12.6.50), (12.6.51). Из равенств (12.6.48), (12.6.52) следует, что равенство (12.6.47) ковариантно относительно выбора базиса правого $A$-векторного пространства $V$.  □

---

Равенство

$$(12.6.58) \quad \begin{aligned} f' &= g^{-1*}{}^*{}_*f^*{}_*g \\ &= g^{-1*}{}^*{}_*(f_1+f_2)^*{}_*g \\ &= g^{-1*}{}^*{}_*f_1^*{}_*g \\ &+ g^{-1*}{}^*{}_*f_2^*{}_*g \\ &= f_1' + f_2' \end{aligned}$$

является следствием равенств (12.6.53), (12.6.55), (12.6.56), (12.6.57). Из равенств (12.6.54), (12.6.58) следует, что равенство (12.6.53) ковариантно относительно выбора базиса правого $A$-векторного пространства $V$.  □

---

**Теорема 12.6.9.** *Пусть $V$ - правое $A$-векторное пространство столбцов. Пусть эндоморфизм $\overline{f}$ правого $A$-векторного пространства $V$ является произведением эндоморфизмов $\overline{f}_1$, $\overline{f}_2$. Матрица $f$ эндоморфизма $\overline{f}$ равна $_*{}^*$-произведению матриц $f_2$, $f_1$ эндоморфизмов $\overline{f}_2$, $\overline{f}_1$ и это равенство не зависит от выбора базиса.*

Доказательство. Пусть $\overline{\overline{e}}$, $\overline{\overline{e}}'$ - базисы правого $A$-векторного пространства $V$ и связаны пассивным преобразованием $g$

$$e' = g_*{}^*e$$

Пусть $f$ - матрица эндоморфизма $\overline{f}$ относительно базиса $\overline{\overline{e}}$. Пусть $f'$ - матрица эндоморфизма $\overline{f}$ относительно базиса $\overline{\overline{e}}'$. Пусть $f_i$, $i = 1, 2$, - матрица эндоморфизма $\overline{f}_i$ относительно базиса $\overline{\overline{e}}$. Пусть $f_i'$, $i = 1, 2$, - матрица эндоморфизма $\overline{f}_i$ относительно базиса $\overline{\overline{e}}'$. Согласно теореме 11.5.2,

$$(12.6.59) \quad f = f_{1*}{}^*f_2$$

$$(12.6.60) \quad f' = f_{1*}'{}^*f_2'$$

Согласно теореме 12.6.5,

$$(12.6.61) \quad f' = g^{-1*}{}_*{}^*f_*{}^*g$$

$$(12.6.62) \quad f_1' = g^{-1*}{}_*{}^*f_{1*}{}^*g$$

---

**Теорема 12.6.10.** *Пусть $V$ - правое $A$-векторное пространство строк. Пусть эндоморфизм $\overline{f}$ правого $A$-векторного пространства $V$ является произведением эндоморфизмов $\overline{f}_1$, $\overline{f}_2$. Матрица $f$ эндоморфизма $\overline{f}$ равна $^*{}_*$-произведению матриц $f_2$, $f_1$ эндоморфизмов $\overline{f}_2$, $\overline{f}_1$ и это равенство не зависит от выбора базиса.*

Доказательство. Пусть $\overline{\overline{e}}$, $\overline{\overline{e}}'$ - базисы правого $A$-векторного пространства $V$ и связаны пассивным преобразованием $g$

$$e' = g^*{}_*e$$

Пусть $f$ - матрица эндоморфизма $\overline{f}$ относительно базиса $\overline{\overline{e}}$. Пусть $f'$ - матрица эндоморфизма $\overline{f}$ относительно базиса $\overline{\overline{e}}'$. Пусть $f_i$, $i = 1, 2$, - матрица эндоморфизма $\overline{f}_i$ относительно базиса $\overline{\overline{e}}$. Пусть $f_i'$, $i = 1, 2$, - матрица эндоморфизма $\overline{f}_i$ относительно базиса $\overline{\overline{e}}'$. Согласно теореме 11.5.3,

$$(12.6.65) \quad f = f_{1*}{}^*f_2$$

$$(12.6.66) \quad f' = f_{1*}'{}^*f_2'$$

Согласно теореме 12.6.6,

$$(12.6.67) \quad f' = g^{-1*}{}^*{}_*f^*{}_*g$$

$$(12.6.68) \quad f_1' = g^{-1*}{}^*{}_*f_{1*}^*{}_*g$$

(12.6.63) $\qquad f'_2 = g^{-1}{}_*{}^*{}_* f_{2*}{}^* g$

Равенство

$$f' = g^{-1}{}_*{}^*{}_* f_*{}^* g$$
$$= g^{-1}{}_*{}^*{}_* f_{1*} f_{2*}{}^* g$$
(12.6.64)
$$= g^{-1}{}_*{}^*{}_* f_{1*}{}^* g$$
$$\quad {}_*{}^* g^{-1}{}_*{}^*{}_* f_{2*}{}^* g$$
$$= f'_{1*}{}^* f'_2$$

является следствием равенств (12.6.59), (12.6.61), (12.6.62), (12.6.63). Из равенств (12.6.60), (12.6.64) следует, что равенство (12.6.59) ковариантно относительно выбора базиса правого $A$-векторного пространства $V$. $\qquad\square$

(12.6.69) $\qquad f'_2 = g^{-1}{}^*{}_*{}_* f_2{}^*{}_* g$

Равенство

$$f' = g^{-1}{}^*{}_* f^*{}_* g$$
$$= g^{-1}{}^*{}_* f_1{}^*{}_* f_2{}^*{}_* g$$
(12.6.70)
$$= g^{-1}{}^*{}_* f_1{}^*{}_* g$$
$$\quad {}^*{}_* g^{-1}{}^*{}_* f_2{}^*{}_* g$$
$$= f'_1{}^*{}_* f'_2$$

является следствием равенств (12.6.65), (12.6.67), (12.6.68), (12.6.69). Из равенств (12.6.66), (12.6.70) следует, что равенство (12.6.65) ковариантно относительно выбора базиса правого $A$-векторного пространства $V$. $\qquad\square$

# Ковариантность в $A$-векторном пространстве

### 13.1. Геометрический объект в левом $A$-векторном пространстве

Предположим, что $V$, $W$ - левые $A$-векторные пространства столбцов . Пусть $G(V_*)$ - группа симметрий левого $A$-векторного пространства $V$. Гомоморфизм

$$(13.1.1) \qquad F: G(V_*) \to GL(W_*)$$

отображает пассивное преобразование $g \in G(V_*)$

$$(13.1.2) \qquad \overline{\overline{e}}_{V2} = g^*{}_* \overline{\overline{e}}_{V1}$$

левого $A$-векторного пространства $V$ в пассивное преобразование $F(g) \in GL(W_*)$

$$(13.1.3) \qquad \overline{\overline{e}}_{W2} = F(g)^*{}_* \overline{\overline{e}}_{W1}$$

левого $A$-векторного пространства $W$.

Тогда координатное преобразование в левом $A$-векторном пространстве $W$ принимает вид

$$(13.1.4) \qquad w_2 = w_1{}^*{}_* F(g)^{-1^*{}^*}$$

Следовательно, отображение

$$(13.1.5) \qquad F_1(g) = F(g)^{-1^*{}^*}$$

является правосторонним представлением

$$F_1: G(V_*) \longrightarrow W_*$$

группы $G(V_*)$ на множестве $W_*$.

---

ОПРЕДЕЛЕНИЕ 13.1.1. *Мы будем называть орбиту*

$$(13.1.11) \quad \begin{aligned} &\mathcal{O}(V, W, \overline{\overline{e}}_V, w) \\ &= (w^*{}_* F(G)^{-1^*{}^*}, G^*{}_* \overline{\overline{e}}_V) \end{aligned}$$

*представления $F_1$ многообразием* **координат геометрического объекта** *в левом $A$-векторном пространстве $V$ столбцов. Для любого базиса*

$(13.1.2)$ $\quad \overline{\overline{e}}_{V2} = g^*{}_* \overline{\overline{e}}_{V1}$

*соответствующая точка*

$(13.1.4)$ $\quad w_2 = w_1{}^*{}_* F(g)^{-1^*{}^*}$

*орбиты определяет* **координаты**

---

Предположим, что $V$, $W$ - левые $A$-векторные пространства строк . Пусть $G(V_*)$ - группа симметрий левого $A$-векторного пространства $V$. Гомоморфизм

$$(13.1.6) \qquad F: G(V_*) \to GL(W_*)$$

отображает пассивное преобразование $g \in G(V_*)$

$$(13.1.7) \qquad \overline{\overline{e}}_{V2} = g_*{}^* \overline{\overline{e}}_{V1}$$

левого $A$-векторного пространства $V$ в пассивное преобразование $F(g) \in GL(W_*)$

$$(13.1.8) \qquad \overline{\overline{e}}_{W2} = F(g)_*{}^* \overline{\overline{e}}_{W1}$$

левого $A$-векторного пространства $W$.

Тогда координатное преобразование в левом $A$-векторном пространстве $W$ принимает вид

$$(13.1.9) \qquad w_2 = w_1{}_* F(g)^{-1^*{}^*}$$

Следовательно, отображение

$$(13.1.10) \qquad F_1(g) = F(g)^{-1^*{}^*}$$

является правосторонним представлением

$$F_1: G(V_*) \longrightarrow W_*$$

группы $G(V_*)$ на множестве $W_*$.

---

ОПРЕДЕЛЕНИЕ 13.1.2. *Мы будем называть орбиту*

$$(13.1.12) \quad \begin{aligned} &\mathcal{O}(V, W, \overline{\overline{e}}_V, w) \\ &= (w_*{}^* F(G)^{-1^*{}^*}, G_*{}^* \overline{\overline{e}}_V) \end{aligned}$$

*представления $F_1$ многообразием* **координат геометрического объекта** *в левом $A$-векторном пространстве $V$ строк. Для любого базиса*

$(13.1.7)$ $\quad \overline{\overline{e}}_{V2} = g_*{}^* \overline{\overline{e}}_{V1}$

*соответствующая точка*

$(13.1.9)$ $\quad w_2 = w_1{}_* F(g)^{-1^*{}^*}$

*орбиты определяет* **координаты**

геометрического объекта *в ко-ординатном левом $A$-векторном пространстве относительно базиса* $\overline{\overline{e}}_{V2}$ .          □

ОПРЕДЕЛЕНИЕ 13.1.3. *Допустим даны координаты* $w_1$ *вектора* $\overline{w}$ *относительно базиса* $\overline{\overline{e}}_{W1}$. *Мы будем называть множество векторов*

$$(13.1.13) \qquad \overline{w}_2 = w_2{}^*{}_*e_{W2}$$

**геометрическим объектом**, *определённым в левом $A$-векторном пространстве $V$ столбцов. Для любого базиса* $\overline{\overline{e}}_{W2}$, *соответствующая точка*

$$\boxed{(13.1.4) \quad w_2 = w_1{}^*{}_*F(g)^{-1}{}^*{}_*}$$

*многообразия координат определяет вектор*

$$(13.1.14) \qquad \overline{w}_2 = w_2{}^*{}_*e_{W2}$$

*который называется* **представителем геометрического объекта** *в левом $A$-векторном пространстве $V$ в базисе* $\overline{\overline{e}}_{V2}$.          □

ОПРЕДЕЛЕНИЕ 13.1.4. *Допустим даны координаты* $w_1$ *вектора* $\overline{w}$ *относительно базиса* $\overline{\overline{e}}_{W1}$. *Мы будем называть множество векторов*

$$(13.1.15) \qquad \overline{w}_2 = w_{2*}{}^*e_{W2}$$

**геометрическим объектом**, *определённым в левом $A$-векторном пространстве $V$ строк. Для любого базиса* $\overline{\overline{e}}_{W2}$, *соответствующая точка*

$$\boxed{(13.1.9) \quad w_2 = w_{1*}{}^*F(g)^{-1}{}_*{}^*}$$

*многообразия координат определяет вектор*

$$(13.1.16) \qquad \overline{w}_2 = w_{2*}{}^*e_{W2}$$

*который называется* **представителем геометрического объекта** *в левом $A$-векторном пространстве $V$ в базисе* $\overline{\overline{e}}_{V2}$.          □

Так как геометрический объект - это орбита представления, то согласно теореме 2.5.13 определение геометрического объекта корректно.

Мы будем также говорить, что $\overline{w}$ - это **геометрический объект типа** $F$.

Определения 13.1.1, 13.1.2 строят геометрический объект в координатном пространстве. Определения 13.1.3, 13.1.4 предполагают, что мы выбрали базис в векторном пространстве $W$. Это позволяет использовать представитель геометрического объекта вместо его координат.

ТЕОРЕМА 13.1.5. **(Принцип ковариантности)**. *Представитель геометрического объекта не зависит от выбора базиса* $\overline{\overline{e}}_{V2}$.

ДОКАЗАТЕЛЬСТВО. Чтобы определить представителя геометрического объекта, мы должны выбрать базис $\overline{\overline{e}}_V$, базис $\overline{\overline{e}}_W$ и координаты геометрического объекта $w^\alpha$. Соответствующий представитель геометрического объекта имеет вид

$$\overline{w} = w^*{}_*e_W$$

ДОКАЗАТЕЛЬСТВО. Чтобы определить представителя геометрического объекта, мы должны выбрать базис $\overline{\overline{e}}_V$, базис $\overline{\overline{e}}_W$ и координаты геометрического объекта $w^\alpha$. Соответствующий представитель геометрического объекта имеет вид

$$\overline{w} = w_*{}^*e_W$$

Предположим базис $\overline{\overline{e}}'_V$ связан с базисом $\overline{\overline{e}}_V$ пассивным преобразованием

$$\overline{\overline{e}}'_V = g^*{}_*\overline{\overline{e}}_V$$

Согласно построению это порождает пассивное преобразование

$$(13.1.3) \quad \overline{\overline{e}}_{W2} = F(g)^*{}_*\overline{\overline{e}}_{W1}$$

и координатное преобразование

$$(13.1.4) \quad w_2 = w_1{}^*{}_*F(g)^{-1*}{}_*$$

Соответствующий представитель геометрического объекта имеет вид

$$\overline{w}' = w'^*{}_*e'_W$$
$$= w^*{}_*F(a)^{-1*}{}^*{}_*F(a)^*{}_*e_W$$
$$= w^*{}_*e_W = \overline{w}$$

Следовательно, представитель геометрического объекта инвариантен относительно выбора базиса. □

---

ОПРЕДЕЛЕНИЕ 13.1.6. *Пусть $V$ - левое $A$-векторное пространство столбцов и*

$$\overline{w}_1 = w_1{}^*{}_*e_W$$
$$\overline{w}_2 = w_2{}^*{}_*e_W$$

*геометрические объекты одного и того же типа, определённые в левом $A$-векторном пространстве $V$. Геометрический объект*

$$\overline{w} = (w_1 + w_2)^*{}_*e_W$$

*называется* **суммой**

$$\overline{w} = \overline{w}_1 + \overline{w}_2$$

**геометрических объектов** $\overline{w}_1$ *и* $\overline{w}_2$. □

---

ОПРЕДЕЛЕНИЕ 13.1.8. *Пусть $V$ - левое $A$-векторное пространство столбцов и*

$$\overline{w}_1 = w_1{}^*{}_*e_W$$

*геометрический объект, определённый в левом $A$-векторном пространстве $V$. Геометрический объект*

$$\overline{w}_2 = (kw_1)^*{}_*e_W$$

---

Предположим базис $\overline{\overline{e}}'_V$ связан с базисом $\overline{\overline{e}}_V$ пассивным преобразованием

$$\overline{\overline{e}}'_V = g_*{}^*\overline{\overline{e}}_V$$

Согласно построению это порождает пассивное преобразование

$$(13.1.8) \quad \overline{\overline{e}}_{W2} = F(g)_*{}^*\overline{\overline{e}}_{W1}$$

и координатное преобразование

$$(13.1.9) \quad w_2 = w_1{}_*{}^*F(g)^{-1}{}_*{}^*$$

Соответствующий представитель геометрического объекта имеет вид

$$\overline{w}' = w'_*{}^*e'_W$$
$$= w_*{}^*F(a)^{-1}{}_*{}^*{}_*{}^*F(a)_*{}^*e_W$$
$$= w_*{}^*e_W = \overline{w}$$

Следовательно, представитель геометрического объекта инвариантен относительно выбора базиса. □

---

ОПРЕДЕЛЕНИЕ 13.1.7. *Пусть $V$ - левое $A$-векторное пространство строк и*

$$\overline{w}_1 = w_1{}_*{}^*e_W$$
$$\overline{w}_2 = w_2{}_*{}^*e_W$$

*геометрические объекты одного и того же типа, определённые в левом $A$-векторном пространстве $V$. Геометрический объект*

$$\overline{w} = (w_1 + w_2)_*{}^*e_W$$

*называется* **суммой**

$$\overline{w} = \overline{w}_1 + \overline{w}_2$$

**геометрических объектов** $\overline{w}_1$ *и* $\overline{w}_2$. □

---

ОПРЕДЕЛЕНИЕ 13.1.9. *Пусть $V$ - левое $A$-векторное пространство строк и*

$$\overline{w}_1 = w_1{}_*{}^*e_W$$

*геометрический объект, определённый в левом $A$-векторном пространстве $V$. Геометрический объект*

$$\overline{w}_2 = (kw_1)_*{}^*e_W$$

называется **произведением**

$$\overline{w}_2 = k\overline{w}_1$$

**геометрического объекта** $\overline{w}_1$ **и константы** $k \in A$.          □

---

называется **произведением**

$$\overline{w}_2 = k\overline{w}_1$$

**геометрического объекта** $\overline{w}_1$ **и константы** $k \in A$.          □

---

Теорема 13.1.10. *Геометрические объекты типа $F$, определённые в левом $A$-векторном пространстве $V$ столбцов, образуют левое $A$-векторное пространство столбцов.*

---

Теорема 13.1.11. *Геометрические объекты типа $F$, определённые в левом $A$-векторном пространстве $V$ строк, образуют левое $A$-векторное пространство строк.*

---

Доказательство. Утверждение теорем 13.1.10, 13.1.11 следует из непосредственной проверки свойств векторного пространства.          □

## 13.2. Геометрический объект в правом $A$-векторном пространстве

Предположим, что $V$, $W$ - правые $A$-векторные пространства столбцов . Пусть $G(V_*)$ - группа симметрий правого $A$-векторного пространства $V$. Гомоморфизм

$$(13.2.1) \qquad F : G(V_*) \to GL(W_*)$$

отображает пассивное преобразование $g \in G(V_*)$

$$(13.2.2) \qquad \overline{\overline{e}}_{V2} = \overline{\overline{e}}_{V1*}{}^{*}g$$

правого $A$-векторного пространства $V$ в пассивное преобразование $F(g) \in GL(W_*)$

$$(13.2.3) \qquad \overline{\overline{e}}_{W2} = \overline{\overline{e}}_{W1*}{}^{*}F(g)$$

правого $A$-векторного пространства $W$.

Тогда координатное преобразование в правом $A$-векторном пространстве $W$ принимает вид

$$(13.2.4) \qquad w_2 = F(g)^{-1*}{}_*^*w_1$$

Следовательно, отображение

$$(13.2.5) \qquad F_1(g) = F(g)^{-1*}$$

является левосторонним представлением

$$F_1 : G(V_*) \longrightarrow_* W_*$$

группы $G(V_*)$ на множестве $W_*$.

---

Предположим, что $V$, $W$ - правые $A$-векторные пространства строк . Пусть $G(V_*)$ - группа симметрий правого $A$-векторного пространства $V$. Гомоморфизм

$$(13.2.6) \qquad F : G(V_*) \to GL(W_*)$$

отображает пассивное преобразование $g \in G(V_*)$

$$(13.2.7) \qquad \overline{\overline{e}}_{V2} = \overline{\overline{e}}_{V1}{}^{*}{}_*g$$

правого $A$-векторного пространства $V$ в пассивное преобразование $F(g) \in GL(W_*)$

$$(13.2.8) \qquad \overline{\overline{e}}_{W2} = \overline{\overline{e}}_{W1}{}^{*}{}_*F(g)$$

правого $A$-векторного пространства $W$.

Тогда координатное преобразование в правом $A$-векторном пространстве $W$ принимает вид

$$(13.2.9) \qquad w_2 = F(g)^{-1*}{}^*{}_*w_1$$

Следовательно, отображение

$$(13.2.10) \qquad F_1(g) = F(g)^{-1*}{}^{*}$$

является левосторонним представлением

$$F_1 : G(V_*) \longrightarrow_* W_*$$

группы $G(V_*)$ на множестве $W_*$.

---

Определение 13.2.1. *Мы будем*

---

Определение 13.2.2. *Мы будем*

*называть орбиту*

(13.2.11)
$$\mathcal{O}(V, W, \overline{\overline{e}}_V, w)$$
$$= (F(G)^{-1}{}_*{}^*{}_*w, \overline{\overline{e}}_{V*}{}^*G)$$

*представления $F_1$ многообразием* **координат геометрического объекта** *в правом $A$-векторном пространстве $V$ столбцов. Для любого базиса*

(13.2.2)    $\overline{\overline{e}}_{V2} = \overline{\overline{e}}_{V1*}{}^*g$

*соответствующая точка*

(13.2.4)    $w_2 = F(g)^{-1}{}_*{}^*{}_*w_1$

*орбиты определяет* **координаты геометрического объекта** *в координатном правом $A$-векторном пространстве относительно базиса $\overline{\overline{e}}_{V2}$ .*                                                   □

*называть орбиту*

(13.2.12)
$$\mathcal{O}(V, W, \overline{\overline{e}}_V, w)$$
$$= (F(G)^{-1}{}^*{}_*{}_*w, \overline{\overline{e}}_V{}^*{}_*G)$$

*представления $F_1$ многообразием* **координат геометрического объекта** *в правом $A$-векторном пространстве $V$ строк. Для любого базиса*

(13.2.7)    $\overline{\overline{e}}_{V2} = \overline{\overline{e}}_{V1}{}^*{}_*g$

*соответствующая точка*

(13.2.9)    $w_2 = F(g)^{-1}{}^*{}_*{}_*w_1$

*орбиты определяет* **координаты геометрического объекта** *в координатном правом $A$-векторном пространстве относительно базиса $\overline{\overline{e}}_{V2}$ .*                                                   □

**Определение 13.2.3.** *Допустим даны координаты $w_1$ вектора $\overline{w}$ относительно базиса $\overline{\overline{e}}_{W1}$. Мы будем называть множество векторов*

(13.2.13)          $\overline{w}_2 = e_{W2*}{}^*w_2$

**геометрическим объектом**, *определённым в правом $A$-векторном пространстве $V$ столбцов. Для любого базиса $\overline{\overline{e}}_{W2}$, соответствующая точка*

(13.2.4)    $w_2 = F(g)^{-1}{}_*{}^*{}_*w_1$

*многообразия координат определяет вектор*

(13.2.14)          $\overline{w}_2 = e_{W2*}{}^*w_2$

*который называется* **представителем геометрического объекта** *в правом $A$-векторном пространстве $V$ в базисе $\overline{\overline{e}}_{V2}$.*                                    □

**Определение 13.2.4.** *Допустим даны координаты $w_1$ вектора $\overline{w}$ относительно базиса $\overline{\overline{e}}_{W1}$. Мы будем называть множество векторов*

(13.2.15)          $\overline{w}_2 = e_{W2}{}^*{}_*w_2$

**геометрическим объектом**, *определённым в правом $A$-векторном пространстве $V$ строк. Для любого базиса $\overline{\overline{e}}_{W2}$, соответствующая точка*

(13.2.9)    $w_2 = F(g)^{-1}{}^*{}_*{}_*w_1$

*многообразия координат определяет вектор*

(13.2.16)          $\overline{w}_2 = e_{W2}{}^*{}_*w_2$

*который называется* **представителем геометрического объекта** *в правом $A$-векторном пространстве $V$ в базисе $\overline{\overline{e}}_{V2}$.*                                    □

Так как геометрический объект - это орбита представления, то согласно теореме 2.5.13 определение геометрического объекта корректно.

Мы будем также говорить, что $\overline{w}$ - это **геометрический объект типа $F$**.

Определения 13.2.1, 13.2.2 строят геометрический объект в координатном пространстве. Определения 13.2.3, 13.2.4 предполагают, что мы выбрали базис в векторном пространстве $W$. Это позволяет использовать представитель геометрического объекта вместо его координат.

Теорема 13.2.5. **(Принцип ковариантности).** *Представитель геометрического объекта не зависит от выбора базиса* $\overline{\overline{e}}_{V2}$.

Доказательство. Чтобы определить представителя геометрического объекта, мы должны выбрать базис $\overline{\overline{e}}_V$, базис $\overline{\overline{e}}_W$ и координаты геометрического объекта $w^\alpha$. Соответствующий представитель геометрического объекта имеет вид

$$\overline{w} = e_{W*}{}^* w$$

Предположим базис $\overline{\overline{e}}'_V$ связан с базисом $\overline{\overline{e}}_V$ пассивным преобразованием

$$\overline{\overline{e}}'_V = \overline{\overline{e}}_{V*}{}^* g$$

Согласно построению это порождает пассивное преобразование

$$\boxed{(13.2.3)\qquad \overline{\overline{e}}_{W2} = \overline{\overline{e}}_{W1*}{}^* F(g)}$$

и координатное преобразование

$$\boxed{(13.2.4)\quad w_2 = F(g)^{-1*}{}^*{}_* w_1}$$

Соответствующий представитель геометрического объекта имеет вид

$$\overline{w}' = e'_{W*}{}^* w'$$
$$= e_{W*}{}^* F(a)_*{}^* F(a)^{-1*}{}^*{}_* w$$
$$= e_{W*}{}^* w = \overline{w}$$

Следовательно, представитель геометрического объекта инвариантен относительно выбора базиса. □

Определение 13.2.6. *Пусть $V$ - правое $A$-векторное пространство столбцов и*

$$\overline{w}_1 = e_{W*}{}^* w_1$$
$$\overline{w}_2 = e_{W*}{}^* w_2$$

*геометрические объекты одного и того же типа, определённые в правом $A$-векторном пространстве $V$. Геометрический объект*

$$\overline{w} = e_{W*}{}^* (w_1 + w_2)$$

*называется* **суммой**

$$\overline{w} = \overline{w}_1 + \overline{w}_2$$

**геометрических объектов** $\overline{w}_1$ *и* $\overline{w}_2$. □

Доказательство. Чтобы определить представителя геометрического объекта, мы должны выбрать базис $\overline{\overline{e}}_V$, базис $\overline{\overline{e}}_W$ и координаты геометрического объекта $w^\alpha$. Соответствующий представитель геометрического объекта имеет вид

$$\overline{w} = e_W{}^*{}_* w$$

Предположим базис $\overline{\overline{e}}'_V$ связан с базисом $\overline{\overline{e}}_V$ пассивным преобразованием

$$\overline{\overline{e}}'_V = \overline{\overline{e}}_V{}^*{}_* g$$

Согласно построению это порождает пассивное преобразование

$$\boxed{(13.2.8)\qquad \overline{\overline{e}}_{W2} = \overline{\overline{e}}_{W1}{}^*{}_* F(g)}$$

и координатное преобразование

$$\boxed{(13.2.9)\quad w_2 = F(g)^{-1*}{}^*{}_* w_1}$$

Соответствующий представитель геометрического объекта имеет вид

$$\overline{w}' = e'_W{}^*{}_* w'$$
$$= e_W{}^*{}_* F(a)^*{}_* F(a)^{-1*}{}^*{}_* w$$
$$= e_W{}^*{}_* w = \overline{w}$$

Следовательно, представитель геометрического объекта инвариантен относительно выбора базиса. □

Определение 13.2.7. *Пусть $V$ - правое $A$-векторное пространство строк и*

$$\overline{w}_1 = e_W{}^*{}_* w_1$$
$$\overline{w}_2 = e_W{}^*{}_* w_2$$

*геометрические объекты одного и того же типа, определённые в правом $A$-векторном пространстве $V$. Геометрический объект*

$$\overline{w} = e_W{}^*{}_* (w_1 + w_2)$$

*называется* **суммой**

$$\overline{w} = \overline{w}_1 + \overline{w}_2$$

**геометрических объектов** $\overline{w}_1$ *и* $\overline{w}_2$. □

Определение 13.2.8. *Пусть* $V$ - *правое* $A$-*векторное пространство столбцов и*

$$\overline{w}_1 = e_{W*}{}^* w_1$$

*геометрический объект, определённый в правом* $A$-*векторном пространстве* $V$. *Геометрический объект*

$$\overline{w}_2 = e_{W*}{}^* (w_1 k)$$

*называется* **произведением**

$$\overline{w}_2 = \overline{w}_1 k$$

**геометрического объекта** $\overline{w}_1$ **и константы** $k \in A$. □

Определение 13.2.9. *Пусть* $V$ - *правое* $A$-*векторное пространство строк и*

$$\overline{w}_1 = e_{W}{}^*{}_* w_1$$

*геометрический объект, определённый в правом* $A$-*векторном пространстве* $V$. *Геометрический объект*

$$\overline{w}_2 = e_{W}{}^*{}_* (w_1 k)$$

*называется* **произведением**

$$\overline{w}_2 = \overline{w}_1 k$$

**геометрического объекта** $\overline{w}_1$ **и константы** $k \in A$. □

Теорема 13.2.10. *Геометрические объекты типа* $F$, *определённые в правом* $A$-*векторном пространстве* $V$ *столбцов, образуют правое* $A$-*векторное пространство столбцов.*

Теорема 13.2.11. *Геометрические объекты типа* $F$, *определённые в правом* $A$-*векторном пространстве* $V$ *строк, образуют правое* $A$-*векторное пространство строк.*

Доказательство. Утверждение теорем 13.2.10, 13.2.11 следует из непосредственной проверки свойств векторного пространства. □

## 13.3. Итоги

Вопрос как велико разнообразие геометрических объектов хорошо изучен в случае векторных пространств. Однако он не столь очевиден в случае левого или правого $A$-векторных пространств. Как видно из таблицы 13.3.1 левое $A$-векторное пространство столбцов и правое $A$-векторное пространство строк имеют общую группу симметрии $GL(n^*{}_* A)$. Это позволяет, рассматривая пассивное представление в левом $A$-векторном пространстве столбцов, изучать геометрический объект в правом $A$-векторном пространстве строк. Можем ли мы одновременно изучать геометрический объект в левом $A$-векторном пространстве строк? На первый взгляд, в силу теоремы 9.5.9 ответ отрицательный. Однако, равенство (8.1.24) устанавливает искомое линейное отображение между $GL(n_*{}^* A)$ и $GL(n^*{}_* A)$.

Теорема 13.3.1. *Многообразие базисов левого* $A$-*векторного пространства* $V$ *и многообразие базисов правого* $A$-*векторного пространства* $V$ *отличны*

$$GL(n_*{}^* A)_*{}^* E_n \neq E_{n}{}^*{}_* GL(n^*{}_* A)$$

Доказательство.

Лемма 13.3.2. *Мы можем отождествить многообразие базисов* $GL(n_*{}^* A)(V)_*{}^* \overline{E_n}$ *и множество* $_*{}^*$-*невырожденных матриц* $GL(n_*{}^* A)$.

Лемма 13.3.3. *Мы можем отождествить многообразие базисов* $\overline{E_n}{}^*{}_* GL(n^*{}_* A)(V)$ *и множество* $^*{}_*$-*невырожденных матриц* $GL(n^*{}_* A)$.

Таблица 13.3.1. Активное и пассивное представления
в $A$-векторном пространстве

| векторное пространство | группа представления | активное представление | пассивное представление |
|---|---|---|---|
| левое $A$-векторное пространство строк | $GL(n_*{}^*A)$ | правостороннее | левостороннее |
| левое $A$-векторное пространство столбцов | $GL(n^*{}_*A)$ | правостороннее | левостороннее |
| правое $A$-векторное пространство строк | $GL(n^*{}_*A)$ | левостороннее | правостороннее |
| правое $A$-векторное пространство столбцов | $GL(n_*{}^*A)$ | левостороннее | правостороннее |

Доказательство. Согласно замечанию 12.5.20, мы можем отождествить базис $\overline{\overline{e}}$ правого $A$-векторного пространства $V$ и матрицу $e$ координат базиса $\overline{\overline{e}}$ относительно базиса $\overline{\overline{E_n}}$

$$\overline{\overline{e}} = e^*{}_* \overline{\overline{E_n}}$$

⊙

Доказательство. Согласно замечанию 12.2.20, мы можем отождествить базис $\overline{\overline{e}}$ левого $A$-векторного пространства $V$ и матрицу $e$ координат базиса $\overline{\overline{e}}$ относительно базиса $\overline{\overline{E_n}}$

$$\overline{\overline{e}} = \overline{\overline{E_n}}{}_*{}^* e$$

⊙

Из лемм 13.3.2, 13.3.3 следует, что если множество векторов $(f_i, i \in I, |I| = n)$ порождают базис левого $A$-векторного пространства столбцов и базис правого $A$-векторного пространства столбцов, то их координатная матрица

$$(13.3.1) \qquad f = \begin{pmatrix} f_1^1 & ... & f_1^n \\ ... & ... & ... \\ f_n^1 & ... & f_n^n \end{pmatrix}$$

является $_*{}^*$-невырожденной и $^*{}_*$-невырожденной матрицей. Утверждение теоремы следует из теоремы 9.5.9. □

# Линейное отображение $A$-векторного пространства

## 14.1. Линейное отображение

ОПРЕДЕЛЕНИЕ 14.1.1. *Пусть $A_i$, $i = 1, 2$, - алгебра с делением над коммутативным кольцом $D_i$. Пусть $V_i$, $i = 1, 2$, - $A_i$-векторное пространство. Морфизм диаграммы представлений*

$$A_1 \xleftarrow{\quad*\quad} D_1 \xrightarrow{\quad*\quad} V_1$$

*в диаграмму представлений*

$$A_2 \xleftarrow{\quad*\quad} D_2 \xrightarrow{\quad*\quad} V_2$$

*называется* **линейным отображением** *$A_1$-векторного пространства $V_1$ в $A_2$-векторное пространство $V_2$. Обозначим* $\mathcal{L}(D_1 \to D_2; A_1 \to A_2; V_1 \to V_2)$ *множество линейных отображений $A_1$-векторного пространства $V_1$ в $A_2$-векторное пространство $V_2$.* $\square$

Из определения 14.1.1 следует, что линейное отображение $A_1$-векторного пространства $V_1$ в $A_2$-векторное пространство $V_2$ - это линейное отображение $D_1$-векторного пространства $V_1$ в $D_2$-векторное пространство $V_2$.

ТЕОРЕМА 14.1.2. *Пусть $V^1, ..., V^n$ - $A$-векторные пространства и*

$$V = V^1 \oplus ... \oplus V^n$$

*Представим $v$-число*

$$v = v^1 \oplus ... \oplus v^n$$

*как вектор столбец*

$$v = \begin{pmatrix} v^1 \\ ... \\ v^n \end{pmatrix}$$

*Представим линейное отображение*

$$f : V \to W$$

*как вектор строку*

$$f = \begin{pmatrix} f_1 & ... & f_n \end{pmatrix}$$

$$f_i : V^i \to W$$

*Тогда значение отображения $f$ в $V$-числе $v$ можно представить как произведение матриц*

$$(14.1.1) \qquad f \circ v = \begin{pmatrix} f_1 & ... & f_n \end{pmatrix} \circ^\circ \begin{pmatrix} v^1 \\ ... \\ v^n \end{pmatrix} = f_i \circ v^i$$

Доказательство. Теорема является следствием определений (10.6.12), (11.6.12). □

ТЕОРЕМА 14.1.3. *Пусть $W^1$, ..., $W^m$ - $A$-векторные пространства и*
$$W = W^1 \oplus ... \oplus W^m$$

*Представим $W$-число*
$$w = w^1 \oplus ... \oplus w^m$$

*как вектор столбец*
$$w = \begin{pmatrix} w^1 \\ ... \\ w^m \end{pmatrix}$$

*Тогда линейное отображение*
$$f : V \to W$$

*имеет представление как вектор столбец отображений*
$$f = \begin{pmatrix} f^1 \\ ... \\ f^m \end{pmatrix}$$

*таким образом, что если $w = f \circ v$, то*
$$\begin{pmatrix} w^1 \\ ... \\ w^m \end{pmatrix} = \begin{pmatrix} f^1 \\ ... \\ f^m \end{pmatrix} \circ v = \begin{pmatrix} f^1 \circ v \\ ... \\ f^m \circ v \end{pmatrix}$$

Доказательство. Теорема является следствием теоремы 2.2.8. □

ТЕОРЕМА 14.1.4. *Пусть $V^1$, ..., $V^n$, $W^1$, ..., $W^m$ - $A$-векторные пространства и*
$$V = V^1 \oplus ... \oplus V^n$$
$$W = W^1 \oplus ... \oplus W^m$$

*Представим $V$-число*
$$v = v^1 \oplus ... \oplus v^n$$

*как вектор столбец*

$$(14.1.2) \qquad v = \begin{pmatrix} v^1 \\ ... \\ v^n \end{pmatrix}$$

*Представим $W$-число*

$$w = w^1 \oplus ... \oplus w^m$$

*как вектор столбец*

$$(14.1.3) \qquad w = \begin{pmatrix} w^1 \\ ... \\ w^m \end{pmatrix}$$

*Тогда линейное отображение*

$$f : V \to W$$

*имеет представление как матрица отображений*

$$(14.1.4) \qquad f = \begin{pmatrix} f^1_1 & ... & f^1_n \\ ... & ... & ... \\ f^m_1 & ... & f^m_n \end{pmatrix}$$

*таким образом, что если $w = f \circ v$, то*

$$(14.1.5) \qquad \begin{pmatrix} w^1 \\ ... \\ w^m \end{pmatrix} = \begin{pmatrix} f^1_1 & ... & f^1_n \\ ... & ... & ... \\ f^m_1 & ... & f^m_n \end{pmatrix} \circ \begin{pmatrix} v^1 \\ ... \\ v^n \end{pmatrix} = \begin{pmatrix} f^1_i \circ v^i \\ ... \\ f^m_i \circ v^i \end{pmatrix}$$

*Отображение*

$$f^i_j : V^j \to W^i$$

*является линейным отображением и называется* **частным линейным отображением**

Доказательство. Согласно теореме 14.1.3, существует множество линейных отображений

$$f^i : V \to W^i$$

таких, что

$$(14.1.6) \qquad \begin{pmatrix} w^1 \\ ... \\ w^m \end{pmatrix} = \begin{pmatrix} f^1 \\ ... \\ f^m \end{pmatrix} \circ v = \begin{pmatrix} f^1 \circ v \\ ... \\ f^m \circ v \end{pmatrix}$$

Согласно теореме 14.1.2, для каждого $i$, существует множество линейных отображений

$$f^i_j : V^j \to W^i$$

таких, что

$$(14.1.7) \qquad f^i \circ v = \begin{pmatrix} f^i_1 & ... & f^i_n \end{pmatrix} \circ^\circ \begin{pmatrix} v^1 \\ ... \\ v^n \end{pmatrix} = f^i_j \circ v^j$$

Если мы отождествим матрицы

$$\begin{pmatrix} \begin{pmatrix} f^1_1 & ... & f^1_n \end{pmatrix} \\ ... \\ \begin{pmatrix} f^m_1 & ... & f^m_n \end{pmatrix} \end{pmatrix} = \begin{pmatrix} f^1_1 & ... & f^1_n \\ ... & ... & ... \\ f^m_1 & ... & f^m_n \end{pmatrix}$$

то равенство (14.1.5) является следствием равенств (14.1.6), (14.1.7). $\qquad \square$

## 14.2. Матрица линейного отображения

Определение 14.2.1. *Пусть* $A_i$, $i = 1, ..., n$, *- множество* $D$-*алгебр. Матрица*

$$f = \begin{pmatrix} f^1_1 & ... & f^1_n \\ ... & ... & ... \\ f^m_1 & ... & f^m_n \end{pmatrix}$$

*элементы которой являются отображениями*

$$f^i_j : A_i \to A_j$$

*называется* **матрицей отображений**. $\qquad \square$

На множестве матриц отображений определены следующие операции

* сложение

$$\begin{pmatrix} f^1_1 & ... & f^1_n \\ ... & ... & ... \\ f^m_1 & ... & f^m_n \end{pmatrix} + \begin{pmatrix} g^1_1 & ... & g^1_n \\ ... & ... & ... \\ g^m_1 & ... & g^m_n \end{pmatrix} = \begin{pmatrix} f^1_1 + g^1_1 & ... & f^1_n + g^1_n \\ ... & ... & ... \\ f^m_1 + g^m_1 & ... & f^m_n + g^m_n \end{pmatrix}$$

* $\circ^\circ$-произведение

$$\begin{pmatrix} f^1_1 & ... & f^1_k \\ ... & ... & ... \\ f^m_1 & ... & f^m_k \end{pmatrix} \circ^\circ \begin{pmatrix} g^1_1 & ... & g^1_n \\ ... & ... & ... \\ g^k_1 & ... & g^k_n \end{pmatrix} = \begin{pmatrix} f^1_i \circ g^i_1 & ... & f^1_i \circ g^i_n \\ ... & ... & ... \\ f^m_i \circ g^i_1 & ... & f^m_i \circ g^i_n \end{pmatrix}$$

* $_\circ\circ$-произведение

$$\begin{pmatrix} f^1_1 & ... & f^1_n \\ ... & ... & ... \\ f^k_1 & ... & f^k_n \end{pmatrix} {}_\circ\circ \begin{pmatrix} g^1_1 & ... & g^1_k \\ ... & ... & ... \\ g^m_1 & ... & g^m_k \end{pmatrix} = \begin{pmatrix} f^i_1 \circ g^1_i & ... & f^i_n \circ g^1_i \\ ... & ... & ... \\ f^i_1 \circ g^m_i & ... & f^i_n \circ g^m_i \end{pmatrix}$$

# Линейное отображение левого $A$-векторного пространства

## 15.1. Общее определение

Пусть $A_i$, $i = 1, 2$, - алгебра с делением над коммутативным кольцом $D_i$. Пусть $V_i$, $i = 1, 2$, - левое $A_i$-модуль.

ТЕОРЕМА 15.1.1. *Линейное отображение*

$$(15.1.1) \qquad (h : D_1 \to D_2 \quad g : A_1 \to A_2 \quad f : V_1 \to V_2)$$

*левого $A_2$-векторного пространства $V_1$ в левое $A_2$-векторное пространство $V_2$ удовлетворяет равенствам*

$$(15.1.2) \qquad h(p + q) = h(p) + h(q)$$

$$(15.1.3) \qquad h(pq) = h(p)h(q)$$

$$(15.1.4) \qquad g \circ (a + b) = g \circ a + g \circ b$$

$$(15.1.5) \qquad g \circ (ab) = (g \circ a)(g \circ b)$$

$$(15.1.6) \qquad g \circ (pa) = h(p)(g \circ a)$$

$$(15.1.7) \qquad f \circ (u + v) = f \circ u + f \circ v$$

$$(15.1.8) \qquad f \circ (pv) = h(p)(f \circ v)$$

$$p, q \in D_1 \quad a, b \in A_1 \quad u, v \in V_1$$

ДОКАЗАТЕЛЬСТВО. Согласно определениям 2.3.9, 2.8.5, 14.1.1, 11 отображение

$$(h : D_1 \to D_2 \quad g : A_1 \to A_2)$$

является гомоморфизмом $D_1$-алгебры $A_1$ в $D_2$-алгебру $A_2$. Следовательно, равенства (15.1.2), (15.1.3), (15.1.4), (15.1.5), (15.1.6) являются следствием теоремы 4.3.2. Равенство (15.1.7) является следствием определения 10.1.1, так как, согласно определению 2.3.9, отображение $f$ является гомоморфизмом абелевой группы. Равенство (15.1.8) является следствием равенства (2.3.6), так как отображение

$$(h : D_1 \to D_2 \quad f : V_1 \to V_2)$$

является морфизмом представления $g_{1.34}$ в представление $g_{2.34}$. $\square$

---

[15.1] В теореме 15.1.2, мы опираемся на следующее соглашение. Пусть $i = 1, 2$. Пусть множество векторов $\overline{\overline{e}}_{A_i} = (e_{A_i k}, k \in K_i)$ является базисом и $C_{ij}^k$, $k$, $i$, $j \in K_i$, - структурные константы $D_i$-алгебры столбцов $A_i$. Пусть множество векторов $\overline{\overline{e}}_{V_i} = (e_{V_i i}, i \in$

Теорема 15.1.2. *Линейное отображение* [15.1]

$$(15.1.1) \qquad (h : D_1 \to D_2 \quad \overline{g} : A_1 \to A_2 \quad \overline{f} : V_1 \to V_2)$$

*имеет представление*

$$(15.1.9) \qquad\qquad b = gh(a)$$

$$(15.1.10) \qquad\qquad \overline{g} \circ (a^i e_{A_1 i}) = h(a^i) g_i^k e_{A_2 k}$$

$$(15.1.11) \qquad\qquad \overline{g} \circ (e_{A_1} a) = e_{A_2} gh(a)$$

$$(15.1.12) \qquad\qquad f \circ (v^*{}_* e_{V_1}) = (f_i^j \circ g \circ v^i) e_{V_2 \cdot j}$$

$$w = f_\circ{}^\circ (g \circ v)$$

*относительно заданных базисов. Здесь*

- $a$ - *координатная матрица $A_1$-числа $\overline{a}$ относительно базиса* $\overline{\overline{e}}_{A_1}$ *(A)*

$$\overline{a} = e_{A_1} a$$

- $h(a) = (h(a_k), k \in \mathbf{K})$ - *матрица $D_2$-чисел.*
- $b$ - *координатная матрица $A_2$-числа*

$$\overline{b} = \overline{g} \circ \overline{a}$$

  *относительно базиса* $\overline{\overline{e}}_{A_2}$

$$\overline{b} = e_{A_2} b$$

- $g$ - *координатная матрица множества $A_2$-чисел* $(\overline{g} \circ e_{A_1 k}, k \in \mathbf{K})$ *относительно базиса* $\overline{\overline{e}}_{A_2}$.
- *Матрица гомоморфизма и структурные константы связаны соотношением*

$$(15.1.13) \qquad\qquad h(C_{1\,ij}^{\ \ k}) g_k^l = g_i^p g_j^q C_{2\,pq}^{\ \ l}$$

- $v$ - *координатная матрица $V_1$-числа $\overline{v}$ относительно базиса* $\overline{\overline{e}}_{V_1}$ *(A)*

$$\overline{v} = v^*{}_* e_{V_1}$$

- $g(v) = (g(v_i), i \in \mathbf{I})$ - *матрица $A_2$-чисел.*
- $w$ - *координатная матрица $V_2$-числа*

$$\overline{w} = \overline{f} \circ \overline{v}$$

  *относительно базиса* $\overline{\overline{e}}_{V_2}$

$$\overline{w} = w^*{}_* e_{V_2}$$

- $f$ - *координатная матрица множества $V_2$-чисел* $(\overline{f} \circ e_{V_1 i}, i \in \mathbf{I})$ *относительно базиса* $\overline{\overline{e}}_{V_2}$.

*Отображение*

$$f_j^i : A_2(g \circ e_{V_1 \cdot j}) \to A_2 e_i$$

*является линейным отображением и называется* **частным линейным отображением**.

---

$\mathbf{I}_i$) является базисом левого $A_i$-векторного пространства $V_i$.

Доказательство. Согласно определениям   2.3.9, 2.8.5, 14.1.1,  отображение

$$(15.1.14) \qquad (h : D_1 \to D_2 \quad \overline{g} : A_1 \to A_2)$$

является гомоморфизмом $D_1$-алгебры $A_1$ в $D_2$-алгебру $A_2$. Следовательно, равенства   (15.1.9), (15.1.10), (15.1.11), (15.1.13)  следуют из равенств

$$(5.4.11) \quad b = fh(a)$$

$$(5.4.12) \quad \overline{f} \circ (a^i e_{1i}) = h(a^i) f_i^k e_{2k}$$

$$(5.4.13) \quad \overline{f} \circ (e_1 a) = e_2 fh(a)$$

$$(5.4.16) \quad h(C_{1\,ij}^{\ \ k}) f_k^l = f_i^p f_j^q C_{2\,pq}^{\ \ l}$$

Из теоремы 4.3.3, следует что матрица $g$ определена однозначно.

$A_1$-модуль $V_1$ является прямой суммой

$$V_1 = \bigoplus_{k \in K} A_1 e_{V_1 \cdot k}$$

$A_2$-модуль $V_2$ является прямой суммой

$$V_2 = \bigoplus_{l \in L} A_2 e_{V_2 \cdot l}$$

Равенство (15.1.12) является следствием теоремы 14.1.4.                    $\square$

## 15.2. Линейное отображение, когда кольца  $D_1 = D_2 = D$

Пусть  $A_i$, $i = 1,\,2$,  - алгебра с делением над коммутативным кольцом $D$. Пусть  $V_i$, $i = 1,\,2$,  -  левое $A_i$-модуль.

Теорема 15.2.1. *Линейное отображение*

$$(15.2.1) \qquad (g : A_1 \to A_2 \quad f : V_1 \to V_2)$$

*левого $A_2$-векторного пространства $V_1$ в  левое $A_2$-векторное пространство $V_2$ удовлетворяет равенствам*

$$(15.2.2) \qquad g \circ (a + b) = g \circ a + g \circ b$$

$$(15.2.3) \qquad g \circ (ab) = (g \circ a)(g \circ b)$$

$$(15.2.4) \qquad g \circ (pa) = p(g \circ a)$$

$$(15.2.5) \qquad f \circ (u + v) = f \circ u + f \circ v$$

$$(15.2.6) \qquad f \circ (pv) = p(f \circ v)$$

$$p \in D \quad a, b \in A_1 \quad u, v \in V_1$$

Доказательство. Согласно определениям   2.3.9, 2.8.5, 14.1.1,  1 отображение

$$g : A_1 \to A_2$$

является гомоморфизмом $D$-алгебры $A_1$ в $D$-алгебру $A_2$. Следовательно, равенства  (15.2.2), (15.2.3), (15.2.4)  являются следствием теоремы 4.3.10.  Равенство (15.2.5) является следствием определения 10.2.1, так как, согласно опре-делению [2.3.9](#), отображение $f$ является гомоморфизмом абелевой группы. Равенство [(15.2.6)](#) является следствием равенства [(2.3.6)](#), так как отображение

$$f : V_1 \to V_2$$

является морфизмом представления $g_{1.34}$ в представление $g_{2.34}$.    □

**Теорема 15.2.2.** *Линейное отображение* [15.2](#)

$$\boxed{(15.2.1) \quad (\overline{g} : A_1 \to A_2 \quad \overline{f} : V_1 \to V_2)}$$

*имеет представление*

$$(15.2.7) \qquad b = ga$$

$$(15.2.8) \qquad \overline{g} \circ (a^i e_{A_1 i}) = a^i g_i{}^k e_{A_2 k}$$

$$(15.2.9) \qquad \overline{g} \circ (e_{A_1} a) = e_{A_2} ga$$

$$(15.2.10) \qquad f \circ (v^*{}_* e_{V_1}) = (f_i{}^j \circ g \circ v^i) e_{V_2 \cdot j}$$

$$w = f_\circ{}^\circ (g \circ v)$$

*относительно заданных базисов. Здесь*

- $a$ - *координатная матрица $A_1$-числа $\overline{a}$ относительно базиса* $\overline{\overline{e}}_{A_1}$ *(A)*

$$\overline{a} = e_{A_1} a$$

- $b$ - *координатная матрица $A_2$-числа*

$$\overline{b} = \overline{g} \circ \overline{a}$$

*относительно базиса* $\overline{\overline{e}}_{A_2}$

$$\overline{b} = e_{A_2} b$$

- $g$ - *координатная матрица множества $A_2$-чисел* $(\overline{g} \circ e_{A_1 k}, k \in K)$ *относительно базиса* $\overline{\overline{e}}_{A_2}$.
- *Матрица гомоморфизма и структурные константы связаны соотношением*

$$(15.2.11) \qquad C_1{}^k_{ij} g_k{}^l = g_i{}^p g_j{}^q C_2{}^l_{pq}$$

- $v$ - *координатная матрица $V_1$-числа $\overline{v}$ относительно базиса* $\overline{\overline{e}}_{V_1}$ *(A)*

$$\overline{v} = v^*{}_* e_{V_1}$$

- $g(v) = (g(v_i), i \in I)$ - *матрица $A_2$-чисел.*
- $w$ - *координатная матрица $V_2$-числа*

$$\overline{w} = \overline{f} \circ \overline{v}$$

*относительно базиса* $\overline{\overline{e}}_{V_2}$

$$\overline{w} = w^*{}_* e_{V_2}$$

- $f$ - *координатная матрица множества $V_2$-чисел* $(\overline{f} \circ e_{V_1 i}, i \in I)$ *относительно базиса* $\overline{\overline{e}}_{V_2}$.

---

15.2  В теореме [15.2.2](#), мы опираемся на следующее соглашение. Пусть $i = 1, 2$. Пусть множество векторов $\overline{\overline{e}}_{A_i} = (e_{A_i k}, k \in K_i)$ является базисом и $C_i{}^k_{ij}$, $k, i, j \in K_i$, - структурные константы $D$-алгебры столбцов $A_i$. Пусть множество векторов $\overline{\overline{e}}_{V_i} = (e_{V_i i}, i \in I_i)$ является базисом левого $A_i$-векторного пространства $V_i$.

*Отображение*

$$f_j^i : A_2(g \circ e_{V_1 \cdot j}) \to A_2 e_i$$

*является линейным отображением и называется* **частным линейным отображением**.

Доказательство. Согласно определениям 2.3.9, 2.8.5, 14.1.1, отображение

$$(15.2.12) \qquad \overline{g} : A_1 \to A_2$$

является гомоморфизмом $D$-алгебры $A_1$ в $D$-алгебру $A_2$. Следовательно, равенства (15.2.7), (15.2.8), (15.2.9), (15.2.11) следуют из равенств

$$\boxed{(5.4.53) \quad b = fa}$$

$$\boxed{(5.4.54) \quad \overline{f} \circ (a^i e_{1i}) = a^i f_i^k e_k}$$

$$\boxed{(5.4.55) \quad \overline{f} \circ (e_1 a) = e_2 fa}$$

$$\boxed{(5.4.58) \quad C_{1\,ij}^{\ \ k} f_k^l = f_i^p f_j^q C_{2\,pq}^{\ \ l}}$$

Из теоремы 4.3.11, следует что матрица $g$ определена однозначно.

$A_1$-модуль $V_1$ является прямой суммой

$$V_1 = \bigoplus_{k \in K} A_1 e_{V_1 \cdot k}$$

$A_2$-модуль $V_2$ является прямой суммой

$$V_2 = \bigoplus_{l \in L} A_2 e_{V_2 \cdot l}$$

Равенство (15.2.10) является следствием теоремы 14.1.4. $\qquad\qquad \square$

## 15.3. Линейное отображение, когда $D$-алгебры $A_1 = A_2 = A$

Пусть $A$ - алгебра с делением над коммутативным кольцом $D$. Пусть $V_i$, $i = 1, 2$, - левое $A$-векторное пространство.

Теорема 15.3.1. *Линейное отображение*

$$(15.3.1) \qquad f : V_1 \to V_2$$

*левого $A$-векторного пространства $V_1$ в левое $A$-векторное пространство $V_2$ удовлетворяет равенствам*

$$(15.3.2) \qquad f \circ (u + v) = f \circ u + f \circ v$$

$$(15.3.3) \qquad f \circ (av) = a(f \circ v)$$

$$a \in A \quad u, v \in V_1$$

Доказательство. Равенство (15.3.2) является следствием определения 10.3.1, так как, согласно определению 2.3.9, отображение $f$ является гомоморфизмом абелевой группы. Равенство (15.3.3) является следствием равенства (2.3.6), так как отображение

$$f : V_1 \to V_2$$

является морфизмом представления $g_{1.34}$ в представление $g_{2.34}$.          □

Теорема 15.3.2. *Линейное отображение* [15.3]

$$\boxed{(15.3.1) \quad \overline{f} : V_1 \to V_2}$$

*имеет представление*

$$(15.3.4) \qquad f \circ (v^*{}_* e_{V_1}) = (f_i^j \circ v^i) e_{V_2 \cdot j}$$

$$w = f_\circ{}^\circ v$$

*относительно заданных базисов. Здесь*

- $v$ - *координатная матрица $V_1$-числа $\overline{v}$ относительно базиса* $\overline{\overline{e}}_{V_1}$ *(A)*

$$\overline{v} = v^*{}_* e_{V_1}$$

- $w$ - *координатная матрица $V_2$-числа*

$$\overline{w} = \overline{f} \circ \overline{v}$$

*относительно базиса* $\overline{\overline{e}}_{V_2}$

$$\overline{w} = w^*{}_* e_{V_2}$$

- $f$ - *координатная матрица множества $V_2$-чисел* $(\overline{f} \circ e_{V_1 i}, i \in I)$
  *относительно базиса* $\overline{\overline{e}}_{V_2}$.

*Отображение*

$$f_j^i : A_2(g \circ e_{V_1 \cdot j}) \to A_2 e_i$$

*является линейным отображением и называется* **частным линейным отображением**.

Доказательство.
$A_1$-модуль $V_1$ является прямой суммой

$$V_1 = \bigoplus_{k \in K} A_1 e_{V_1 \cdot k}$$

$A_2$-модуль $V_2$ является прямой суммой

$$V_2 = \bigoplus_{l \in L} A_2 e_{V_2 \cdot l}$$

Равенство (15.3.4) является следствием теоремы 14.1.4.          □

---

# Линейное отображение правого $A$-векторного пространства

## 16.1. Общее определение

Пусть $A_i$, $i = 1, 2$, - алгебра с делением над коммутативным кольцом $D_i$. Пусть $V_i$, $i = 1, 2$, - правое $A_i$-модуль.

ТЕОРЕМА 16.1.1. *Линейное отображение*

$$(16.1.1) \qquad (h : D_1 \to D_2 \quad g : A_1 \to A_2 \quad f : V_1 \to V_2)$$

*правого $A_2$-векторного пространства $V_1$ в правое $A_2$-векторное пространство $V_2$ удовлетворяет равенствам*

$$(16.1.2) \qquad h(p + q) = h(p) + h(q)$$

$$(16.1.3) \qquad h(pq) = h(p)h(q)$$

$$(16.1.4) \qquad g \circ (a + b) = g \circ a + g \circ b$$

$$(16.1.5) \qquad g \circ (ab) = (g \circ a)(g \circ b)$$

$$(16.1.6) \qquad g \circ (ap) = (g \circ a)h(p)$$

$$(16.1.7) \qquad f \circ (u + v) = f \circ u + f \circ v$$

$$(16.1.8) \qquad f \circ (vp) = (f \circ v)h(p)$$

$$p, q \in D_1 \quad a, b \in A_1 \quad u, v \in V_1$$

ДОКАЗАТЕЛЬСТВО. Согласно определениям 2.3.9, 2.8.5, 14.1.1, 11 отображение

$$(h : D_1 \to D_2 \quad g : A_1 \to A_2)$$

является гомоморфизмом $D_1$-алгебры $A_1$ в $D_2$-алгебру $A_2$. Следовательно, равенства (16.1.2), (16.1.3), (16.1.4), (16.1.5), (16.1.6) являются следствием теоремы 4.3.2. Равенство (16.1.7) является следствием определения 11.1.1, так как, согласно определению 2.3.9, отображение $f$ является гомоморфизмом абелевой группы. Равенство (16.1.8) является следствием равенства (2.3.6), так как отображение

$$(h : D_1 \to D_2 \quad f : V_1 \to V_2)$$

является морфизмом представления $g_{1.34}$ в представление $g_{2.34}$. $\qquad \square$

---

16.1  В теореме 16.1.2, мы опираемся на следующее соглашение. Пусть $i = 1, 2$. Пусть множество векторов $\overline{\overline{e}}_{A_i} = (e_{A_i k}, k \in K_i)$ является базисом и $C_{ij}^{k}$, $k$, $i$, $j \in K_i$, - структурные константы $D_i$-алгебры столбцов $A_i$. Пусть множество векторов $\overline{\overline{e}}_{V_i} = (e_{V_i i}, i \in$

Теорема 16.1.2. *Линейное отображение* <span style="color:magenta">16.1</span>

$$\boxed{(16.1.1) \qquad (h : D_1 \to D_2 \quad \overline{g} : A_1 \to A_2 \quad \overline{f} : V_1 \to V_2)}$$

*имеет представление*

$$(16.1.9) \qquad b = gh(a)$$

$$(16.1.10) \qquad \overline{g} \circ (a^i e_{A_1 \, i}) = h(a^i) g_i^{\,k} e_{A_2 \, k}$$

$$(16.1.11) \qquad \overline{g} \circ (e_{A_1} a) = e_{A_2} gh(a)$$

$$(16.1.12) \qquad f \circ (e_{V_1 *}{}^* v) = e_{V_2 \cdot j}(f_i^{\,j} \circ g \circ v^i)$$

$$w = f_\circ{}^\circ (g \circ v)$$

*относительно заданных базисов. Здесь*

- $a$ - *координатная матрица $A_1$-числа $\overline{a}$ относительно базиса* $\overline{\overline{e}}_{A_1}$ *(A)*

$$\overline{a} = e_{A_1} a$$

- $h(a) = (h(a_k), k \in K)$ - *матрица $D_2$-чисел.*
- $b$ - *координатная матрица $A_2$-числа*

$$\overline{b} = \overline{g} \circ \overline{a}$$

*относительно базиса* $\overline{\overline{e}}_{A_2}$

$$\overline{b} = e_{A_2} b$$

- $g$ - *координатная матрица множества $A_2$-чисел* $(\overline{g} \circ e_{A_1 \, k}, k \in K)$ *относительно базиса* $\overline{\overline{e}}_{A_2}$.
- *Матрица гомоморфизма и структурные константы связаны соотношением*

$$(16.1.13) \qquad h(C_{1 \, ij}^{\ \ k}) g_k^{\,l} = g_i^{\,p} g_j^{\,q} C_{2 \, pq}^{\ \ l}$$

- $v$ - *координатная матрица $V_1$-числа $\overline{v}$ относительно базиса* $\overline{\overline{e}}_{V_1}$ *(A)*

$$\overline{v} = e_{V_1 *}{}^* v$$

- $g(v) = (g(v_i), i \in I)$ - *матрица $A_2$-чисел.*
- $w$ - *координатная матрица $V_2$-числа*

$$\overline{w} = \overline{f} \circ \overline{v}$$

*относительно базиса* $\overline{\overline{e}}_{V_2}$

$$\overline{w} = e_{V_2 *}{}^* w$$

- $f$ - *координатная матрица множества $V_2$-чисел* $(\overline{f} \circ e_{V_1 \, i}, i \in I)$ *относительно базиса* $\overline{\overline{e}}_{V_2}$.

*Отображение*

$$f_j^{\,i} : (g \circ e_{V_1 \cdot j}) A_2 \to e_i A_2$$

*является линейным отображением и называется* **частным линейным отображением**.

---

<span style="color:magenta">$I_i$</span>) является базисом правого $A_i$-векторного пространства $V_i$.

Доказательство. Согласно определениям 2.3.9, 2.8.5, 14.1.1, отображение

$$(16.1.14) \qquad (h : D_1 \to D_2 \quad \overline{g} : A_1 \to A_2)$$

является гомоморфизмом $D_1$-алгебры $A_1$ в $D_2$-алгебру $A_2$. Следовательно, равенства (16.1.9), (16.1.10), (16.1.11), (16.1.13) следуют из равенств

$$\boxed{(5.4.11) \quad b = fh(a)}$$

$$\boxed{(5.4.12) \quad \overline{f} \circ (a^i e_{1i}) = h(a^i) f_i^k e_{2k}}$$

$$\boxed{(5.4.13) \quad \overline{f} \circ (e_1 a) = e_2 fh(a)}$$

$$\boxed{(5.4.16) \quad h(C_{1\,ij}^{\;\;k}) f_k^l = f_i^p f_j^q C_{2\,pq}^{\;\;l}}$$

Из теоремы 4.3.3, следует что матрица $g$ определена однозначно.

$A_1$-модуль $V_1$ является прямой суммой

$$V_1 = \bigoplus_{k \in K} A_1 e_{V_1 \cdot k}$$

$A_2$-модуль $V_2$ является прямой суммой

$$V_2 = \bigoplus_{l \in L} A_2 e_{V_2 \cdot l}$$

Равенство (16.1.12) является следствием теоремы 14.1.4.     $\square$

## 16.2. Линейное отображение, когда кольца $D_1 = D_2 = D$

Пусть $A_i$, $i = 1, 2$, - алгебра с делением над коммутативным кольцом $D$. Пусть $V_i$, $i = 1, 2$, - правое $A_i$-модуль.

Теорема 16.2.1. *Линейное отображение*

$$(16.2.1) \qquad (g : A_1 \to A_2 \quad f : V_1 \to V_2)$$

*правого $A_2$-векторного пространства $V_1$ в правое $A_2$-векторное пространство $V_2$ удовлетворяет равенствам*

$$(16.2.2) \qquad g \circ (a + b) = g \circ a + g \circ b$$

$$(16.2.3) \qquad g \circ (ab) = (g \circ a)(g \circ b)$$

$$(16.2.4) \qquad g \circ (ap) = (g \circ a)p$$

$$(16.2.5) \qquad f \circ (u + v) = f \circ u + f \circ v$$

$$(16.2.6) \qquad f \circ (vp) = (f \circ v)p$$

$$p \in D \quad a, b \in A_1 \quad u, v \in V_1$$

Доказательство. Согласно определениям 2.3.9, 2.8.5, 14.1.1, 1 отображение

$$g : A_1 \to A_2$$

является гомоморфизмом $D$-алгебры $A_1$ в $D$-алгебру $A_2$. Следовательно, равенства (16.2.2), (16.2.3), (16.2.4) являются следствием теоремы 4.3.10. Равенство (16.2.5) является следствием определения 11.2.1, так как, согласно опре-делению 2.3.9, отображение $f$ является гомоморфизмом абелевой группы. Равенство (16.2.6) является следствием равенства (2.3.6), так как отображение

$$f : V_1 \to V_2$$

является морфизмом представления $g_{1.34}$ в представление $g_{2.34}$.    □

ТЕОРЕМА 16.2.2. *Линейное отображение* [16.2]

$$\boxed{(16.2.1)\quad (\overline{g} : A_1 \to A_2 \quad \overline{f} : V_1 \to V_2)}$$

*имеет представление*

$$(16.2.7)\qquad\qquad b = ga$$

$$(16.2.8)\qquad\qquad \overline{g} \circ (a^i e_{A_1 i}) = a^i g_i^k e_{A_2 k}$$

$$(16.2.9)\qquad\qquad \overline{g} \circ (e_{A_1} a) = e_{A_2} ga$$

$$(16.2.10)\qquad\qquad f \circ (e_{V_1 *}{}^* v) = e_{V_2 \cdot j}(f_i^j \circ g \circ v^i)$$

$$w = f_\circ{}^\circ (g \circ v)$$

*относительно заданных базисов. Здесь*

- $a$ - *координатная матрица $A_1$-числа $\overline{a}$ относительно базиса* $\overline{\overline{e}}_{A_1}$ *(A)*

$$\overline{a} = e_{A_1} a$$

- $b$ - *координатная матрица $A_2$-числа*

$$\overline{b} = \overline{g} \circ \overline{a}$$

*относительно базиса* $\overline{\overline{e}}_{A_2}$

$$\overline{b} = e_{A_2} b$$

- $g$ - *координатная матрица множества $A_2$-чисел*  $(\overline{g} \circ e_{A_1 k}, k \in K)$ *относительно базиса* $\overline{\overline{e}}_{A_2}$.
- *Матрица гомоморфизма и структурные константы связаны соотношением*

$$(16.2.11)\qquad\qquad C_1{}_{ij}^{\ k} g_k^l = g_i^p g_j^q C_2{}_{pq}^{\ l}$$

- $v$ - *координатная матрица $V_1$-числа $\overline{v}$ относительно базиса* $\overline{\overline{e}}_{V_1}$ *(A)*

$$\overline{v} = e_{V_1 *}{}^* v$$

- $g(v) = (g(v_i), i \in I)$  - *матрица $A_2$-чисел.*
- $w$ - *координатная матрица $V_2$-числа*

$$\overline{w} = \overline{f} \circ \overline{v}$$

*относительно базиса* $\overline{\overline{e}}_{V_2}$

$$\overline{w} = e_{V_2 *}{}^* w$$

- $f$ - *координатная матрица множества $V_2$-чисел*  $(\overline{f} \circ e_{V_1 i}, i \in I)$ *относительно базиса* $\overline{\overline{e}}_{V_2}$.

---

[16.2]  В теореме 16.2.2, мы опираемся на следующее соглашение.   Пусть  $i = 1, 2$.   Пусть множество векторов $\overline{\overline{e}}_{A_i} = (e_{A_i k}, k \in K_i)$  является базисом и  $C_i{}_{ij}^{\ k}$,    $k, i, j \in K_i$, - структурные константы $D$-алгебры столбцов $A_i$.   Пусть множество векторов $\overline{\overline{e}}_{V_i} = (e_{V_i i}, i \in I_i)$  является базисом правого $A_i$-векторного пространства $V_i$.

*Отображение*

$$f^i_j : (g \circ e_{V_1 \cdot j}) A_2 \to e_i A_2$$

*является линейным отображением и называется* **частным линейным отображением**.

ДОКАЗАТЕЛЬСТВО.   Согласно определениям   2.3.9, 2.8.5, 14.1.1,   отображение

$$(16.2.12) \qquad\qquad \overline{g} : A_1 \to A_2$$

является гомоморфизмом $D$-алгебры $A_1$ в $D$-алгебру $A_2$. Следовательно, равенства   (16.2.7), (16.2.8), (16.2.9), (16.2.11)   следуют из равенств

$$(5.4.53) \quad b = fa$$

$$(5.4.54) \quad \overline{f} \circ (a^i e_{1i}) = a^i f^k_i e_k$$

$$(5.4.55) \quad \overline{f} \circ (e_1 a) = e_2 fa$$

$$(5.4.58) \quad C_{1\,ij}^{\ \ k} f^l_k = f^p_i f^q_j C_{2\,pq}^{\ \ l}$$

Из теоремы 4.3.11, следует что матрица $g$ определена однозначно.

$A_1$-модуль $V_1$ является прямой суммой

$$V_1 = \bigoplus_{k \in K} A_1 e_{V_1 \cdot k}$$

$A_2$-модуль $V_2$ является прямой суммой

$$V_2 = \bigoplus_{l \in L} A_2 e_{V_2 \cdot l}$$

Равенство (16.2.10) является следствием теоремы 14.1.4.    $\square$

## 16.3. Линейное отображение, когда $D$-алгебры   $A_1 = A_2 = A$

Пусть $A$ - алгебра с делением над коммутативным кольцом $D$.   Пусть   $V_i$, $i = 1, 2$,  -  правое $A$-векторное пространство.

ТЕОРЕМА 16.3.1.   *Линейное отображение*

$$(16.3.1) \qquad\qquad f : V_1 \to V_2$$

*правого $A$-векторного пространства $V_1$ в  правое $A$-векторное пространство $V_2$ удовлетворяет равенствам*

$$(16.3.2) \qquad\qquad f \circ (u + v) = f \circ u + f \circ v$$

$$(16.3.3) \qquad\qquad f \circ (va) = (f \circ v)a$$

$$a \in A \quad u, v \in V_1$$

ДОКАЗАТЕЛЬСТВО.   Равенство (16.3.2) является следствием определения 11.3.1, так как, согласно определению 2.3.9, отображение $f$ является гомоморфизмом абелевой группы. Равенство (16.3.3) является следствием равенства (2.3.6), так как отображение

$$f : V_1 \to V_2$$

является морфизмом представления $g_{1.34}$ в представление $g_{2.34}$.        □

Теорема 16.3.2. *Линейное отображение* 16.3

$$\boxed{(16.3.1) \quad \overline{f} : V_1 \to V_2}$$

*имеет представление*

$$(16.3.4) \qquad f \circ (e_{V_1*}{}^*v) = e_{V_2 \cdot j}(f_i^j \circ v^i)$$

$$w = f \circ {}^\circ v$$

*относительно заданных базисов. Здесь*

- $v$ - *координатная матрица* $V_1$*-числа* $\overline{v}$ *относительно базиса* $\overline{\overline{e}}_{V_1}$ *(A)*

$$\overline{v} = e_{V_1*}{}^*v$$

- $w$ - *координатная матрица* $V_2$*-числа*

$$\overline{w} = \overline{f} \circ \overline{v}$$

  *относительно базиса* $\overline{\overline{e}}_{V_2}$

$$\overline{w} = e_{V_2*}{}^*w$$

- $f$ - *координатная матрица множества* $V_2$*-чисел* $(\overline{f} \circ e_{V_1 i}, i \in I)$ *относительно базиса* $\overline{\overline{e}}_{V_2}$.

*Отображение*

$$f_j^i : (g \circ e_{V_1 \cdot j})A_2 \to e_i A_2$$

*является линейным отображением и называется* **частным линейным отображением**.

Доказательство.
$A_1$-модуль $V_1$ является прямой суммой

$$V_1 = \bigoplus_{k \in K} A_1 e_{V_1 \cdot k}$$

$A_2$-модуль $V_2$ является прямой суммой

$$V_2 = \bigoplus_{l \in L} A_2 e_{V_2 \cdot l}$$

Равенство (16.3.4) является следствием теоремы 14.1.4.        □

---

16.3  В теореме 16.3.2, мы опираемся на следующее соглашение.  Пусть $i = 1, 2$.  Пусть множество векторов $\overline{\overline{e}}_{V_i} = (e_{V_i i}, i \in I_i)$ является базисом правого $A$-векторного пространства $V_i$.

# Преобразование подобия

Пусть $A$ - ассоциативная $D$-алгебра с делением.

## 17.1. Преобразование подобия в левом $A$-векторном пространстве

**Теорема 17.1.1.** *Пусть $V$ - левое $A$-векторное пространство столбцов и $\overline{\overline{e}}$ - базис левого $A$-векторного пространства $V$. Пусть $\dim V = $ $n$ и $E_n$ - $n \times n$ единичная матрица. Для любого $A$-числа $a$, существует эндоморфизм $\overline{\overline{e}}a$ левого $A$-векторного пространства $V$, который имеет матрицу $E_n a$ относительно базиса $\overline{\overline{e}}$.*

**Теорема 17.1.2.** *Пусть $V$ - левое $A$-векторное пространство строк и $\overline{\overline{e}}$ - базис левого $A$-векторного пространства $V$. Пусть $\dim V = $ $n$ и $E_n$ - $n \times n$ единичная матрица. Для любого $A$-числа $a$, существует эндоморфизм $\overline{\overline{e}}a$ левого $A$-векторного пространства $V$, который имеет матрицу $E_n a$ относительно базиса $\overline{\overline{e}}$.*

**Доказательство.** Пусть $\overline{\overline{e}}$ - базис левого $A$-векторного пространства столбцов. Согласно теореме 10.3.8 и определению 10.3.4, отображение

$$\overline{\overline{e}}a : v^*{}_*e \in V \to v^*{}_*(E_n a)^*{}_*e \in V$$

является эндоморфизмом левого $A$-векторного пространства $V$. $\square$

**Доказательство.** Пусть $\overline{\overline{e}}$ - базис левого $A$-векторного пространства столбцов. Согласно теореме 10.3.9 и определению 10.3.4, отображение

$$\overline{\overline{e}}a : v_*{}^*e \in V \to v_*{}^*(E_n a)_*{}^*e \in V$$

является эндоморфизмом левого $A$-векторного пространства $V$. $\square$

**Определение 17.1.3.** *Эндоморфизм $\overline{\overline{e}}a$ левого $A$-векторного пространства $V$ называется* **преобразованием подобия** *относительно базиса $\overline{\overline{e}}$.* $\square$

**Теорема 17.1.4.** *Пусть $V$ - левое $A$-векторное пространство столбцов и $\overline{\overline{e}}_1$, $\overline{\overline{e}}_2$ - базисы левого $A$-векторного пространства $V$. Пусть базисы $\overline{\overline{e}}_1$, $\overline{\overline{e}}_2$ связаны пассивным преобразованием $g$*

$$\overline{\overline{e}}_2 = g^*{}_* \overline{\overline{e}}_1$$

*Преобразование подобия $\overline{\overline{e}}_1 a$ имеет матрицу*

$$(17.1.1) \qquad g a^*{}_* g^{-1}{}^*{}_*$$

**Теорема 17.1.5.** *Пусть $V$ - левое $A$-векторное пространство строк и $\overline{\overline{e}}_1$, $\overline{\overline{e}}_2$ - базисы левого $A$-векторного пространства $V$. Пусть базисы $\overline{\overline{e}}_1$, $\overline{\overline{e}}_2$ связаны пассивным преобразованием $g$*

$$\overline{\overline{e}}_2 = g_*{}^* \overline{\overline{e}}_1$$

*Преобразование подобия $\overline{\overline{e}}_1 a$ имеет матрицу*

$$(17.1.2) \qquad g a_*{}^* g^{-1}{}_*{}^*$$

*относительно базиса $\overline{\overline{e}}_2$.*

| | |
|---|---|
| Доказательство. Согласно теореме 12.3.5, преобразование подобия $\overline{\overline{e}}_1 a$ имеет матрицу | Доказательство. Согласно теореме 12.3.6, преобразование подобия $\overline{\overline{e}}_1 a$ имеет матрицу |
| (17.1.3) $\qquad g^*{}_* E_n a^*{}_* g^{-1^*{}^*}$ | (17.1.4) $\qquad g_* E_n a_*{}^* g^{-1_*{}^*}$ |
| относительно базиса $\overline{\overline{e}}_2$. Теорема является следствием выражения (17.1.3). $\qquad\square$ | относительно базиса $\overline{\overline{e}}_2$. Теорема является следствием выражения (17.1.4). $\qquad\square$ |

Теорема 17.1.6. *Пусть $V$ - левое $A$-векторное пространство столбцов и $\overline{\overline{e}}$ - базис левого $A$-векторного пространства $V$. Множество преобразований подобия $\overline{\overline{e}}a$ левого $A$-векторного пространства $V$ является правосторонним эффективным представлением $D$-алгебры $A$ в левом $A$-векторном пространстве $V$.*

Доказательство. Рассмотрим множество матриц $E_n a$ преобразований подобия $\overline{\overline{e}}a$ относительно базиса $\overline{\overline{e}}$. Пусть $V^*$ - множество координат $V$-чисел относительно базиса $\overline{\overline{e}}$. Согласно теоремам 13.1.10, 13.1.11, множество $V^*$ является левым $A$-векторным пространством. Согласно теореме 17.1.1 отображение $\overline{\overline{\overline{e}}a}$ является эндоморфизмом левого $A$-векторного пространства столбцов.

| | |
|---|---|
| Лемма 17.1.7. *Пусть $V$ - левое $A$-векторное пространство столбцов. Отображение* | Лемма 17.1.8. *Пусть $V$ - левое $A$-векторное пространство строк. Отображение* |
| (17.1.5) $\quad a \in A \to E_n a \in \operatorname{End}(A, V^*)$ | (17.1.6) $\quad a \in A \to E_n a \in \operatorname{End}(A, V^*)$ |
| *является гомоморфизмом $D$-алгебры $A$.* | *является гомоморфизмом $D$-алгебры $A$.* |
| Доказательство. Согласно теореме 5.4.8, лемма является следствием равенств | Доказательство. Согласно теореме 5.4.8, лемма является следствием равенств |
| $$E_n(a+b) = E_n a + E_n b$$ $$E_n(ba) = b(E_n a)$$ $$E_n(ab) = (E_n a)^*{}_*(E_n b)$$ | $$E_n(a+b) = E_n a + E_n b$$ $$E_n(ba) = b(E_n a)$$ $$E_n(ab) = (E_n a)_*{}^*(E_n b)$$ |
| $\qquad\qquad\qquad\qquad\qquad\odot$ | $\qquad\qquad\qquad\qquad\qquad\odot$ |
| Согласно определению 2.5.1 и лемме 17.1.7, отображение (17.1.5) является правосторонним представлением $D$-алгебры $A$ в левом $A$-векторном пространстве $V^*$. Согласно теореме 5.1.31, из равенства | Согласно определению 2.5.1 и лемме 17.1.8, отображение (17.1.6) является правосторонним представлением $D$-алгебры $A$ в левом $A$-векторном пространстве $V^*$. Согласно теореме 5.1.31, из равенства |
| $$E_n a = E_n b$$ | $$E_n a = E_n b$$ |
| следует, что $a = b$. Согласно определению 2.3.2, представление (17.1.5) эф- | следует, что $a = b$. Согласно определению 2.3.2, представление (17.1.6) эф- |

фективно.

Согласно теоремам 12.3.7, 12.3.9, представление (17.1.5) ковариантно относительно выбора базиса левого $A$-векторного пространства $V$. Согласно теореме 17.1.1 множество отображений $\overline{\overline{e}}a$ является правосторонним представлением $D$-алгебры $A$ в левом $A$-векторном пространстве $V$.

фективно.

Согласно теоремам 12.3.8, 12.3.10, представление (17.1.6) ковариантно относительно выбора базиса левого $A$-векторного пространства $V$. Согласно теореме 17.1.2 множество отображений $\overline{\overline{e}}a$ является правосторонним представлением $D$-алгебры $A$ в левом $A$-векторном пространстве $V$. $\qquad\square$

## 17.2. Преобразование подобия в правом $A$-векторном пространстве

Теорема 17.2.1. *Пусть $V$ - правое $A$-векторное пространство столбцов и $\overline{\overline{e}}$ - базис правого $A$-векторного пространства $V$. Пусть $\dim V = n$ и $E_n$ - $n \times n$ единичная матрица. Для любого $A$-числа $a$, существует эндоморфизм $a\overline{\overline{e}}$ правого $A$-векторного пространства $V$, который имеет матрицу $aE_n$ относительно базиса $\overline{\overline{e}}$.*

Теорема 17.2.2. *Пусть $V$ - правое $A$-векторное пространство строк и $\overline{\overline{e}}$ - базис правого $A$-векторного пространства $V$. Пусть $\dim V = n$ и $E_n$ - $n \times n$ единичная матрица. Для любого $A$-числа $a$, существует эндоморфизм $a\overline{\overline{e}}$ правого $A$-векторного пространства $V$, который имеет матрицу $aE_n$ относительно базиса $\overline{\overline{e}}$.*

Доказательство. Пусть $\overline{\overline{e}}$ - базис правого $A$-векторного пространства столбцов. Согласно теореме 11.3.8 и определению 11.3.4, отображение

$$a\overline{\overline{e}} : e_* {}^*v \in V \to e_* {}^*(aE_n)_* {}^*v \in V$$

является эндоморфизмом правого $A$-векторного пространства $V$. $\qquad\square$

Доказательство. Пусть $\overline{\overline{e}}$ - базис правого $A$-векторного пространства столбцов. Согласно теореме 11.3.9 и определению 11.3.4, отображение

$$a\overline{\overline{e}} : e^* {}_*v \in V \to e^* {}_*(aE_n)^* {}_*v \in V$$

является эндоморфизмом правого $A$-векторного пространства $V$. $\qquad\square$

Определение 17.2.3. *Эндоморфизм $a\overline{\overline{e}}$ правого $A$-векторного пространства $V$ называется* **преобразованием подобия** *относительно базиса $\overline{\overline{e}}$.* $\qquad\square$

Теорема 17.2.4. *Пусть $V$ - правое $A$-векторное пространство столбцов и $\overline{\overline{e}}_1$, $\overline{\overline{e}}_2$ - базисы правого $A$-векторного пространства $V$. Пусть базисы $\overline{\overline{e}}_1$, $\overline{\overline{e}}_2$ связаны пассивным преобразованием $g$*

$$\overline{\overline{e}}_2 = \overline{\overline{e}}_{1*} {}^*g$$

*Преобразование подобия $a\overline{\overline{e}}_1$ имеет матрицу*

$$(17.2.1) \qquad g^{-1*} {}^*{}^*ag$$

Теорема 17.2.5. *Пусть $V$ - правое $A$-векторное пространство строк и $\overline{\overline{e}}_1$, $\overline{\overline{e}}_2$ - базисы правого $A$-векторного пространства $V$. Пусть базисы $\overline{\overline{e}}_1$, $\overline{\overline{e}}_2$ связаны пассивным преобразованием $g$*

$$\overline{\overline{e}}_2 = \overline{\overline{e}}_1 {}^* {}_*g$$

*Преобразование подобия $a\overline{\overline{e}}_1$ имеет матрицу*

$$(17.2.2) \qquad g^{-1*} {}^* {}_*ag$$

*относительно базиса $\overline{\overline{e}}_2$.*

*относительно базиса $\overline{\overline{e}}_2$.*

Доказательство. Согласно теореме 12.6.5, преобразование подобия $\overline{a}\overline{\overline{e}}_1$ имеет матрицу

$$(17.2.3) \qquad g^{-1}{}_{*}{}^{*}{}_{*}{}^{*}aE_{n*}{}^{*}g$$

относительно базиса $\overline{\overline{e}}_2$. Теорема является следствием выражения (17.2.3).  $\square$

Доказательство. Согласно теореме 12.6.6, преобразование подобия $\overline{a}\overline{\overline{e}}_1$ имеет матрицу

$$(17.2.4) \qquad g^{-1}{}^{*}{}_{*}{}^{*}{}_{*}aE_{n}{}^{*}{}_{*}g$$

относительно базиса $\overline{\overline{e}}_2$. Теорема является следствием выражения (17.2.4).  $\square$

ТЕОРЕМА 17.2.6. *Пусть $V$ - правое $A$-векторное пространство столбцов и $\overline{\overline{e}}$ - базис правого $A$-векторного пространства $V$. Множество преобразований подобия $\overline{a}\overline{\overline{e}}$ правого $A$-векторного пространства $V$ является левосторонним эффективным представлением $D$-алгебры $A$ в правом $A$-векторном пространстве $V$.*

Доказательство. Рассмотрим множество матриц $aE_n$ преобразований подобия $\overline{a}\overline{\overline{e}}$ относительно базиса $\overline{\overline{e}}$. Пусть $V^*$ - множество координат $V$-чисел относительно базиса $\overline{\overline{e}}$. Согласно теоремам 13.2.10, 13.2.11, множество $V^*$ является правым $A$-векторным пространством. Согласно теореме 17.2.1 отображение $\overline{a}\overline{\overline{e}}$ является эндоморфизмом правого $A$-векторного пространства столбцов.

ЛЕММА 17.2.7. *Пусть $V$ - правое $A$-векторное пространство столбцов. Отображение*

$$(17.2.5) \quad a \in A \to aE_n \in \operatorname{End}(A, V^*)$$

*является гомоморфизмом $D$-алгебры $A$.*

ЛЕММА 17.2.8. *Пусть $V$ - правое $A$-векторное пространство строк. Отображение*

$$(17.2.6) \quad a \in A \to aE_n \in \operatorname{End}(A, V^*)$$

*является гомоморфизмом $D$-алгебры $A$.*

Доказательство. Согласно теореме 5.4.8, лемма является следствием равенств

$$(a + b)E_n = aE_n + bE_n$$
$$(ba)E_n = b(aE_n)$$
$$((ab)E_n) = (aE_n)_{*}{}^{*}(bE_n)$$

$\odot$

Доказательство. Согласно теореме 5.4.8, лемма является следствием равенств

$$(a + b)E_n = aE_n + bE_n$$
$$(ba)E_n = b(aE_n)$$
$$((ab)E_n) = (aE_n)^{*}{}_{*}(bE_n)$$

$\odot$

Согласно определению 2.5.4 и лемме 17.2.7, отображение (17.2.5) является левосторонним представлением $D$-алгебры $A$ в правом $A$-векторном пространстве $V^*$. Согласно теореме 5.1.31, из равенства

$$aE_n = bE_n$$

следует, что $a = b$. Согласно определе-

Согласно определению 2.5.4 и лемме 17.2.8, отображение (17.2.6) является левосторонним представлением $D$-алгебры $A$ в правом $A$-векторном пространстве $V^*$. Согласно теореме 5.1.31, из равенства

$$aE_n = bE_n$$

следует, что $a = b$. Согласно определе-нию [2.3.2](#), представление [(17.2.5)](#) эффективно.

Согласно теоремам [12.6.7](#), [12.6.9](#), представление [(17.2.5)](#) ковариантно относительно выбора базиса правого $A$-векторного пространства $V$. Согласно теореме [17.2.1](#) множество отображений $\overline{a\overline{e}}$ является левосторонним представлением $D$-алгебры $A$ в правом $A$-векторном пространстве $V$.

нию [2.3.2](#), представление [(17.2.6)](#) эффективно.

Согласно теоремам [12.6.8](#), [12.6.10](#), представление [(17.2.6)](#) ковариантно относительно выбора базиса правого $A$-векторного пространства $V$. Согласно теореме [17.2.2](#) множество отображений $\overline{a\overline{e}}$ является левосторонним представлением $D$-алгебры $A$ в правом $A$-векторном пространстве $V$. $\qquad\square$

## 17.3. Парные представления $D$-алгебры с делением

ОПРЕДЕЛЕНИЕ 17.3.1. *Пусть $V$ - правое $A$-векторное пространство и $\overline{e}$ - базис правого $A$-векторного пространства $V$. Билинейное отображение*

$$(17.3.1) \qquad (a,v) \in A \times V \to av \in V$$

*порождённое левосторонним представлением*

$$(17.3.2) \qquad (a,v) \to \overline{a\overline{e}} \circ v$$

*называется* **левосторонним произведением** *вектора на скаляр.* $\qquad\square$

ОПРЕДЕЛЕНИЕ 17.3.2. *Пусть $V$ - левое $A$-векторное пространство и $\overline{e}$ - базис левого $A$-векторного пространства $V$. Билинейное отображение*

$$(17.3.3) \qquad (v,a) \in V \times A \to va \in V$$

*порождённое правосторонним представлением*

$$(17.3.4) \qquad (v,a) \to \overline{\overline{e}a} \circ v$$

*называется* **правосторонним произведением** *вектора на скаляр.* $\qquad\square$

ТЕОРЕМА 17.3.3. $\forall a, b \in A,\ v \in V$

$$(17.3.5) \qquad (av)b = a(vb)$$

Доказательство. Согласно теореме [17.2.1](#) и определению [17.3.1](#), равенство

$$(17.3.6) \qquad \begin{aligned} a(vb) &= \overline{a\overline{e}} \circ (vb) \\ &= (\overline{a\overline{e}} \circ v)b = (av)b \end{aligned}$$

является следствием равенства

$$\boxed{(11.3.5) \quad f \circ (va) = (f \circ v)a}$$

Равенство [(17.3.5)](#) является следствием равенства [(17.3.6)](#). $\qquad\square$

Векторы $v$ и $av$ называются **колинеарными** слева. Из утверждения слева, что векторы $v$ и $w$ колинеарны слева не следует, что векторы $v$ и $w$ линейно зависимы справа.

ТЕОРЕМА 17.3.4. $\forall a, b \in A,\ v \in V$

$$(17.3.7) \qquad (av)b = a(vb)$$

Доказательство. Согласно теореме [17.1.1](#) и определению [17.3.2](#), равенство

$$(17.3.8) \qquad \begin{aligned} (av)b &= \overline{\overline{e}b} \circ (av) \\ &= a(\overline{\overline{e}b} \circ v) = a(vb) \end{aligned}$$

является следствием равенства

$$\boxed{(10.3.5) \quad f \circ (av) = a(f \circ v)}$$

Равенство [(17.3.7)](#) является следствием равенства [(17.3.8)](#). $\qquad\square$

Векторы $v$ и $va$ называются **колинеарными** справа. Из утверждения, что векторы $v$ и $w$ колинеарны справа не следует, что векторы $v$ и $w$ линейно зависимы слева.

Теорема 17.3.5. *Пусть $V$ - правое $A$-векторное пространство*

$$f : A \longrightarrow\!\!\!* \longrightarrow V$$

*Эффективное левостороннее представление*

$$g(\overline{e}) : A \longrightarrow\!\!\!* \longrightarrow V$$

*$D$-алгебры $A$ в правом $A$-векторном пространстве $V$ зависит от базиса $\overline{e}$ правого $A$-векторного пространства $V$. Следующая диаграмма*

(17.3.9)

$$
\begin{array}{ccc}
V & \xrightarrow{f(a)} & V \\
\downarrow{\scriptstyle g(\overline{e})(b)} & & \downarrow{\scriptstyle g(\overline{e})(b)} \\
V & \xrightarrow{f(a)} & V
\end{array}
$$

*коммутативна для любых $a, b \in A$.*

Теорема 17.3.6. *Пусть $V$ - левое $A$-векторное пространство*

$$f : A \longrightarrow\!\!\!* \longrightarrow V$$

*Эффективное правостороннее представление*

$$g(\overline{\overline{e}}) : A \longrightarrow\!\!\!* \longrightarrow V$$

*$D$-алгебры $A$ в левом $A$-векторном пространстве $V$ зависит от базиса $\overline{\overline{e}}$ левого $A$-векторного пространства $V$. Следующая диаграмма*

(17.3.10)

$$
\begin{array}{ccc}
V & \xrightarrow{f(a)} & V \\
\downarrow{\scriptstyle g(\overline{\overline{e}})(b)} & & \downarrow{\scriptstyle g(\overline{\overline{e}})(b)} \\
V & \xrightarrow{f(a)} & V
\end{array}
$$

*коммутативна для любых $a, b \in A$.*

Доказательство. Коммутативность диаграммы (17.3.9) является следствием равенства (17.3.5). □

Доказательство. Коммутативность диаграммы (17.3.10) является следствием равенства (17.3.7). □

Определение 17.3.7. *Мы будем называть представления $f$ и $g(\overline{e})$, рассмотренные в теореме 17.3.5,* **парными представлениями** *$D$-алгебры $A$. Равенство* (17.3.5) *представляет* **закон ассоциативности** *для парных представлений. Это позволяет нам писать подобное выражение не пользуясь скобками.* □

Определение 17.3.8. *Мы будем называть представления $f$ и $g(\overline{\overline{e}})$, рассмотренные в теореме 17.3.6,* **парными представлениями** *$D$-алгебры $A$. Равенство* (17.3.7) *представляет* **закон ассоциативности** *для парных представлений. Это позволяет нам писать подобное выражение не пользуясь скобками.* □

## 17.4. $A \otimes A$-модуль

Теорема 17.4.1. *Пусть $f$, $g(\overline{e})$ парные представления $D$-алгебры $A$ в $D$-модуле $V$ и $\overline{e}$ - базис правого $A$-векторного пространства $V$. Пусть произведение в $D$-модуле $A \otimes A$ определено согласно правилу*

$$(p_0 \otimes p_1) \circ (q_0 \otimes q_1) = (p_0 q_0) \otimes (q_1 p_1)$$

Теорема 17.4.2. *Пусть $f$, $g(\overline{\overline{e}})$ парные представления $D$-алгебры $A$ в $D$-модуле $V$ и $\overline{\overline{e}}$ - базис левого $A$-векторного пространства $V$. Пусть произведение в $D$-модуле $A \otimes A$ определено согласно правилу*

$$(p_0 \otimes p_1) \circ (q_0 \otimes q_1) = (p_0 q_0) \otimes (q_1 p_1)$$

| | |
|---|---|
| *Тогда представление* | *Тогда представление* |
| $$h : A \otimes A \longrightarrow V$$ | $$h : A \otimes A \longrightarrow V$$ |
| *D-алгебры $A \otimes A$ в D-модуле $V$ опре-делёное равенством* | *D-алгебры $A \otimes A$ в D-модуле $V$ опре-делёное равенством* |
| $$h(a \otimes b) \circ v = avb$$ | $$h(a \otimes b) \circ v = avb$$ |
| *является левым $A \otimes A$-модулем.* | *является левым $A \otimes A$-модулем.* |

Доказательство. Согласно теоремам 17.2.6, 17.1.6, для заданного $v \in V$, отображение

$$p : (a, b) \in A \times A \to avb \in V$$

билинейно. Согласно теореме 4.6.4, для заданного $v \in V$, существует отображение

$$h(a \otimes b) = p(a, b)$$
$$h : a \otimes b \in A \otimes A \to avb \in V$$

Так как $v \in V$ - произвольно, то для любых $a$, $b$  $h(a \otimes b)$  является отображением

$$h(a \otimes b) : v \in V \to avb \in V$$

Согласно теоремам 17.2.1, 17.1.1, (а также 17.2.2, 17.1.2)

$$h(a \otimes b) \in \mathcal{L}(D; V \to V)$$

Если в $D$-модуле $A \otimes A$ определить произведение

$$(a_0 \otimes a_1) \circ (b_0 \otimes b_1) = (a_0 b_0) \otimes (b_1 a_1)$$

то $D$-модуле $A \otimes A$ является $D$-алгеброй. Из равенств

$$(a_0 \otimes a_1 + b_0 \otimes b_1) \circ v = a_0 v a_1 + b_0 v b_1 = (a_0 \otimes a_1) \circ v + (b_0 \otimes b_1) \circ v$$
$$(pa_0 \otimes a_1) \circ v = pa_0 v a_1 = (p(a_0 \otimes a_1)) \circ v$$
$$(a_0 \otimes a_1) \circ (b_0 \otimes b_1) \circ v = (a_0 \otimes a_1) \circ (b_0 v b_1) = a_0 b_0 v b_1 a_1 = (a_0 b_0 \otimes (b))$$

следует, что отображение

$$h : A \otimes A \to \mathcal{L}(D, V \to V)$$

является гомоморфизмом. Следовательно, отображение $h$ является представлением. □

Согласно теоремам 17.4.2, 17.4.1, существует три различных представления $D$-алгебры $A$ в $D$-модуле $V$. Выбор того или иного представления зависит от рассматриваемой задачи.

Например, множество решений системы линейных уравнений

$$5ix^1 + 6jx^2 + 3kx^3 + 7x^4 = i + k$$
$$6kx^1 + 7jx^2 + 9kx^3 + jx^4 = j + 2k$$

в алгебре кватернионов является правым $H$-векторным пространством. Левая линейная комбинация решений вообще говоря решением не является. Однако есть задачи, в которых мы вынуждены отказаться от простой модели и одновременно рассматривать обе структуры векторного пространства.

Пусть $V$ - левый $A \otimes A$-модуль. Если $v$, $w \in V$, то для любых $a$, $b$, $c$, $d \in A$

(17.4.1)                          $$avb + cwd \in V$$

Выражение (17.4.1) называется **линейной комбинацией векторов** левого $A \otimes A$-модуля $V$.

| | |
|---|---|
| Определение 17.4.3. *Пусть* $V$ *левый* $A \otimes A$-*модуль столбцов. Множество векторов* $a_i$, $i \in I$, *левого* $A \otimes A$-*модуля* $V$ **линейно независимо**, *если из равенства* $$b^i a_i c^i = 0$$ *следует* $b^i c^i = 0$, $i = 1$, ..., $n$. *В противном случае, множество векторов* $a_i$ **линейно зависимо**.    $\square$ | Определение 17.4.4. *Пусть* $V$ *левый* $A \otimes A$-*модуль строк. Множество векторов* $a^i$, $i \in I$, *левого* $A \otimes A$-*модуля* $V$ **линейно независимо**, *если из равенства* $$b_i a^i c_i = 0$$ *следует* $b_i c_i = 0$, $i = 1$, ..., $n$. *В противном случае, множество векторов* $a_i$ **линейно зависимо**.    $\square$ |
| Определение 17.4.5. *Мы называем множество векторов* $$\overline{\overline{e}} = (e_i, i \in I)$$ **базисом** $A \otimes A$-*модуля* $V$, *если множество векторов* $e_i$ *линейно независимо и, добавив к этой системе любой другой вектор, мы получим новое множество, которое является линейно зависимым.*    $\square$ | Определение 17.4.6. *Мы называем множество векторов* $$\overline{\overline{e}} = (e^i, i \in I)$$ **базисом** $A \otimes A$-*модуля* $V$, *если множество векторов* $e^i$ *линейно независимо и, добавив к этой системе любой другой вектор, мы получим новое множество, которое является линейно зависимым.*    $\square$ |

## 17.5. Собственное значение эндоморфизма

| | |
|---|---|
| Теорема 17.5.1. *Пусть* $V$ - *левое* $A$-*векторное пространство и* $\overline{\overline{e}}$ - *базис левого* $A$-*векторного пространства* $V$. *Если* $b \notin Z(A)$, *то не существует эндоморфизма* $\overline{f}$ *такого, что* $$(17.5.1) \qquad \overline{f} \circ v = bv$$ | Теорема 17.5.2. *Пусть* $V$ - *правое* $A$-*векторное пространство и* $\overline{\overline{e}}$ - *базис правого* $A$-*векторного пространства* $V$. *Если* $b \notin Z(A)$, *то не существует эндоморфизма* $\overline{f}$ *такого, что* $$(17.5.4) \qquad \overline{f} \circ v = vb$$ |
| Доказательство. Допустим, что существует эндоморфизм $\overline{f}$ такой, что равенство (17.5.1) верно. Равенство $$(ba)v = b(av) = \overline{f} \circ (av)$$ $$(17.5.2) \qquad = a(\overline{f} \circ v) = a(bv)$$ $$= (ab)v$$ является следствием равенств | Доказательство. Допустим, что существует эндоморфизм $\overline{f}$ такой, что равенство (17.5.4) верно. Равенство $$v(ab) = (va)b = \overline{f} \circ (va)$$ $$(17.5.5) \qquad = (\overline{f} \circ v)a = (vb)a$$ $$= v(ba)$$ является следствием равенств |
| $\boxed{(7.2.9) \quad (pq)v = p(qv)}$ $\boxed{(10.3.5) \quad f \circ (av) = a(f \circ v)}$ и равенства (17.5.1). Согласно теореме | $\boxed{(7.4.9) \quad v(pq) = (vp)q}$ $\boxed{(11.3.5) \quad f \circ (va) = (f \circ v)a}$ и равенства (17.5.4). Согласно теореме |

<table>
<tr><td>

[12.1.4](), равенство

$$(17.5.3) \qquad ab = ba$$

для любого $a \in A$ является следствием равенства [(17.5.2)](). Равенство [(17.5.3)]() для любого $a \in A$ противоречит утверждению $b \notin Z(A)$. Следовательно, рассматриваемый эндоморфизм $\overline{f}$ не существует. $\qquad \square$

Согласно теореме [17.5.1](), в отличие от правостороннего произведения вектора на скаляр левостороннее произведение вектора на скаляр не является эндоморфизмом левого $A$-векторного пространства. Из сравнения определения [17.1.3]() и теоремы [17.5.1]() следует, что преобразование подобия играет в некоммутативной алгебре такую же роль, что преобразование подобия играет в коммутативной алгебре.

</td><td>

[12.4.4](), равенство

$$(17.5.6) \qquad ab = ba$$

для любого $a \in A$ является следствием равенства [(17.5.5)](). Равенство [(17.5.6)]() для любого $a \in A$ противоречит утверждению $b \notin Z(A)$. Следовательно, рассматриваемый эндоморфизм $\overline{f}$ не существует. $\qquad \square$

Согласно теореме [17.5.2](), в отличие от левостороннего произведения вектора на скаляр правостороннее произведение вектора на скаляр не является эндоморфизмом правого $A$-векторного пространства. Из сравнения определения [17.2.3]() и теоремы [17.5.2]() следует, что преобразование подобия играет в некоммутативной алгебре такую же роль, что преобразование подобия играет в коммутативной алгебре.

</td></tr>
<tr><td>

ОПРЕДЕЛЕНИЕ 17.5.3. *Пусть $V$ - левое $A$-векторное пространство и $\overline{\overline{e}}$ - базис левого $A$-векторного пространства $V$. Вектор $v \in V$ называется* **собственным вектором** *эндоморфизма*

$$\overline{f} : V \to V$$

*относительно базиса $\overline{\overline{e}}$, если существует $b \in A$ такое, что*

$$(17.5.7) \qquad \overline{f} \circ v = \overline{\overline{e}}b \circ v$$

*$A$-число $b$ называется* **собственным значением** *эндоморфизма $f$ относительно базиса $\overline{\overline{e}}$.* $\qquad \square$

</td><td>

ОПРЕДЕЛЕНИЕ 17.5.4. *Пусть $V$ - правое $A$-векторное пространство и $\overline{\overline{e}}$ - базис правого $A$-векторного пространства $V$. Вектор $v \in V$ называется* **собственным вектором** *эндоморфизма*

$$\overline{f} : V \to V$$

*относительно базиса $\overline{\overline{e}}$, если существует $b \in A$ такое, что*

$$(17.5.8) \qquad \overline{f} \circ v = b\overline{\overline{e}} \circ v$$

*$A$-число $b$ называется* **собственным значением** *эндоморфизма $f$ относительно базиса $\overline{\overline{e}}$.* $\qquad \square$

</td></tr>
<tr><td>

ТЕОРЕМА 17.5.5. *Пусть $V$ - левое $A$-векторное пространство столбцов и $\overline{\overline{e}}_1$, $\overline{\overline{e}}_2$ - базисы левого $A$-векторного пространства $V$. Пусть базисы $\overline{\overline{e}}_1$, $\overline{\overline{e}}_2$ связаны пассивным преобразованием $g$*

$$\overline{\overline{e}}_2 = g^*{}_* \overline{\overline{e}}_1$$

*Пусть $f$ - матрица эндоморфизма $\overline{f}$ и $v$ - матрица вектора $\overline{v}$ относительно базиса $\overline{\overline{e}}_2$. Вектор $\overline{v}$ является собственным вектором эндоморфизма $\overline{f}$ относительно базиса $\overline{\overline{e}}_1$ тогда и*

</td><td>

ТЕОРЕМА 17.5.7. *Пусть $V$ - правое $A$-векторное пространство столбцов и $\overline{\overline{e}}_1$, $\overline{\overline{e}}_2$ - базисы правого $A$-векторного пространства $V$. Пусть базисы $\overline{\overline{e}}_1$, $\overline{\overline{e}}_2$ связаны пассивным преобразованием $g$*

$$\overline{\overline{e}}_2 = \overline{\overline{e}}_{1*}{}^* g$$

*Пусть $f$ - матрица эндоморфизма $\overline{f}$ и $v$ - матрица вектора $\overline{v}$ относительно базиса $\overline{\overline{e}}_2$. Вектор $\overline{v}$ является собственным вектором эндоморфиз-*

</td></tr>
</table>

только тогда, когда существует A-число $b$ такое, что система линейных уравнений

$$(17.5.9) \qquad v^*{}_*(f - gb^*{}_*g^{-1*}{}_*) = 0$$

имеет нетривиальное решение.

Доказательство. Согласно теоремам 10.3.6, 17.1.4, равенство

$$(17.5.10) \quad \begin{aligned} v^*{}_*f^*{}_*e &= \overline{f} \circ \overline{v} = \overline{\overline{b}\overline{e}} \circ \overline{v} \\ &= v^*{}_*gb^*{}_*g^{-1*}{}_*{}_*e \end{aligned}$$

является следствием равенства (17.5.7). Равенство (17.5.9) является следствием равенства (17.5.10) и теоремы 8.4.4. Так как равенство $\overline{v} = 0$ для нас не интересно, система линейных уравнений (17.5.9) должна иметь нетривиальное решение. □

Теорема 17.5.6. *Пусть $V$ - левое $A$-векторное пространство столбцов и $\overline{\overline{e}}_1, \overline{\overline{e}}_2$ - базисы левого $A$-векторного пространства $V$. Пусть базисы $\overline{\overline{e}}_1, \overline{\overline{e}}_2$ связаны пассивным преобразованием $g$*

$$\overline{\overline{e}}_2 = g^*{}_*\overline{\overline{e}}_1$$

*Пусть $f$ - матрица эндоморфизма $\overline{f}$ относительно базиса $\overline{\overline{e}}_2$. $A$-число $b$ является собственным значением эндоморфизма $\overline{f}$ тогда и только тогда, когда матрица*

$$(17.5.11) \qquad f - gb^*{}_*g^{-1*}{}_*$$

$*{}_*$-*вырождена.*

Доказательство. Теорема является следствием теорем 9.4.3, 17.5.5. □

Теорема 17.5.9. *Пусть $V$ - левое $A$-векторное пространство строк и $\overline{\overline{e}}_1, \overline{\overline{e}}_2$ - базисы левого $A$-векторного пространства $V$. Пусть базисы $\overline{\overline{e}}_1, \overline{\overline{e}}_2$ связаны пассивным преобразованием $g$*

$$\overline{\overline{e}}_2 = g_*{}^*\overline{\overline{e}}_1$$

ма $\overline{f}$ относительно базиса $\overline{e}_1$ тогда и только тогда, когда существует A-число $b$ такое, что система линейных уравнений

$$(17.5.12) \qquad (f - g^{-1*}{}^*{}_*bg)_*{}^*v = 0$$

имеет нетривиальное решение.

Доказательство. Согласно теоремам 11.3.6, 17.2.4, равенство

$$(17.5.13) \quad \begin{aligned} e_*{}^*f_*{}^*v &= \overline{f} \circ \overline{v} = \overline{b\overline{e}} \circ \overline{v} \\ &= e_*{}^*g^{-1*}{}^*{}_*bg_*{}^*v \end{aligned}$$

является следствием равенства (17.5.8). Равенство (17.5.12) является следствием равенства (17.5.13) и теоремы 8.4.18. Так как равенство $\overline{v} = 0$ для нас не интересно, система линейных уравнений (17.5.12) должна иметь нетривиальное решение. □

Теорема 17.5.8. *Пусть $V$ - правое $A$-векторное пространство столбцов и $\overline{e}_1, \overline{e}_2$ - базисы правого $A$-векторного пространства $V$. Пусть базисы $\overline{e}_1, \overline{e}_2$ связаны пассивным преобразованием $g$*

$$\overline{e}_2 = \overline{e}_{1*}{}^*g$$

*Пусть $f$ - матрица эндоморфизма $\overline{f}$ относительно базиса $\overline{e}_2$. $A$-число $b$ является собственным значением эндоморфизма $\overline{f}$ тогда и только тогда, когда матрица*

$$(17.5.14) \qquad f - g^{-1*}{}^*{}_*bg$$

$*{}^*$-*вырождена.*

Доказательство. Теорема является следствием теорем 9.4.3, 17.5.7. □

Теорема 17.5.11. *Пусть $V$ - правое $A$-векторное пространство строк и $\overline{e}_1, \overline{e}_2$ - базисы правого $A$-векторного пространства $V$. Пусть базисы $\overline{e}_1, \overline{e}_2$ связаны пассивным преобразованием $g$*

$$\overline{e}_2 = \overline{e}_1{}^*{}_*g$$

*Пусть $f$ - матрица эндоморфизма $\overline{f}$ и $v$ - матрица вектора $\overline{v}$ относительно базиса $\overline{\overline{e}}_2$. Вектор $\overline{v}$ является собственным вектором эндоморфизма $\overline{f}$ относительно базиса $\overline{\overline{e}}_1$ тогда и только тогда, когда существует $A$-число $b$ такое, что система линейных уравнений*

$$(17.5.15) \qquad v_*{}^*(f - gb_*{}^*g^{-1*^*}) = 0$$

*имеет нетривиальное решение.*

Доказательство. Согласно теоремам 10.3.7, 17.1.5, равенство

$$(17.5.16) \qquad \begin{aligned} v_*{}^*f_*{}^*e &= \overline{f} \circ \overline{v} = \overline{b\overline{\overline{e}}} \circ \overline{v} \\ &= v_*{}^*gb_*{}^*g^{-1*^*}{}_*{}^*e \end{aligned}$$

является следствием равенства (17.5.7). Равенство (17.5.15) является следствием равенства (17.5.16) и теоремы 8.4.11. Так как равенство $\overline{v} = 0$ для нас не интересно, система линейных уравнений (17.5.15) должна иметь нетривиальное решение. □

Теорема 17.5.10. *Пусть $V$ - левое $A$-векторное пространство строк и $\overline{\overline{e}}_1, \overline{\overline{e}}_2$ - базисы левого $A$-векторного пространства $V$. Пусть базисы $\overline{\overline{e}}_1, \overline{\overline{e}}_2$ связаны пассивным преобразованием $g$*

$$\overline{\overline{e}}_2 = g_*{}^*\overline{\overline{e}}_1$$

*Пусть $f$ - матрица эндоморфизма $\overline{f}$ относительно базиса $\overline{\overline{e}}_2$. $A$-число $b$ является собственным значением эндоморфизма $\overline{f}$ тогда и только тогда, когда матрица*

$$(17.5.17) \qquad f - gb_*{}^*g^{-1*^*}$$

$_*{}^*$-*вырождена.*

Доказательство. Теорема является следствием теорем 9.4.3, 17.5.9. □

Теорема 17.5.13. *Пусть $V$ - левое $A$-векторное пространство столбцов и $\overline{\overline{e}}_1, \overline{\overline{e}}_2$ - базисы левого $A$-векторного пространства $V$.*

*Пусть $f$ - матрица эндоморфизма $\overline{f}$ и $v$ - матрица вектора $\overline{v}$ относительно базиса $\overline{\overline{e}}_2$. Вектор $\overline{v}$ является собственным вектором эндоморфизма $\overline{f}$ относительно базиса $\overline{\overline{e}}_1$ тогда и только тогда, когда существует $A$-число $b$ такое, что система линейных уравнений*

$$(17.5.18) \qquad (f - g^{-1*^*}{}_*{}^*bg)_*{}^*v = 0$$

*имеет нетривиальное решение.*

Доказательство. Согласно теоремам 11.3.7, 17.2.5, равенство

$$(17.5.19) \qquad \begin{aligned} e_*{}^*f_*{}^*v &= \overline{f} \circ \overline{v} = \overline{\overline{\overline{e}}b} \circ \overline{v} \\ &= e_*{}^*g^{-1*^*}{}_*{}^*bg_*{}^*v \end{aligned}$$

является следствием равенства (17.5.8). Равенство (17.5.18) является следствием равенства (17.5.19) и теоремы 8.4.25. Так как равенство $\overline{v} = 0$ для нас не интересно, система линейных уравнений (17.5.18) должна иметь нетривиальное решение. □

Теорема 17.5.12. *Пусть $V$ - правое $A$-векторное пространство строк и $\overline{\overline{e}}_1, \overline{\overline{e}}_2$ - базисы правого $A$-векторного пространства $V$. Пусть базисы $\overline{\overline{e}}_1, \overline{\overline{e}}_2$ связаны пассивным преобразованием $g$*

$$\overline{\overline{e}}_2 = \overline{\overline{e}}_1{}_*{}^*g$$

*Пусть $f$ - матрица эндоморфизма $\overline{f}$ относительно базиса $\overline{\overline{e}}_2$. $A$-число $b$ является собственным значением эндоморфизма $\overline{f}$ тогда и только тогда, когда матрица*

$$(17.5.20) \qquad f - g^{-1*^*}{}_*{}^*bg$$

$_*{}^*$-*вырождена.*

Доказательство. Теорема является следствием теорем 9.4.3, 17.5.11. □

Теорема 17.5.20. *Пусть $V$ - правое $A$-векторное пространство столбцов и $\overline{\overline{e}}_1, \overline{\overline{e}}_2$ - базисы правого $A$-векторного пространства $V$.*

УТВЕРЖДЕНИЕ 17.5.14. *Собственное значение b эндоморфизма $\overline{f}$ относительно базиса $\overline{\overline{e}}_1$ не зависит от выбора базиса $\overline{\overline{e}}_2$.* ⊙

УТВЕРЖДЕНИЕ 17.5.15. *Собственный вектор $\overline{v}$ эндоморфизма $\overline{f}$, соответствующий собственному значению b, не зависит от выбора базиса $\overline{\overline{e}}_2$.* ⊙

Доказательство. Пусть $f_i$, $i = 1, 2$, - матрица эндоморфизма f относительно базиса $\overline{\overline{e}}_i$. Пусть b - собственное значение эндоморфизма $\overline{f}$ относительно базиса $\overline{\overline{e}}_1$ и $\overline{v}$ собственный вектор эндоморфизма $\overline{f}$, соответствующее собственному значению b. Пусть $v_i$, $i = 1, 2$, - матрица вектора $\overline{v}$ относительно базиса $\overline{\overline{e}}_i$.

Предположим сперва, что $\overline{\overline{e}}_1 = \overline{\overline{e}}_2$. Тогда базисы $\overline{\overline{e}}_1$, $\overline{\overline{e}}_2$ связаны пассивным преобразованием $E_n$

$$\overline{\overline{e}}_2 = E_n{}^*{}_*\overline{\overline{e}}_1$$

Согласно теореме 17.1.4, преобразование подобия $\overline{\overline{e}}_1 b$ имеет матрицу

$$E_n b^*{}_* E_n^{-1*} = E_n b$$

относительно базиса $\overline{\overline{e}}_2$. Следовательно, согласно теореме 17.5.6, мы доказали следующие утверждения.

ЛЕММА 17.5.16. *Матрица $f_1 - E_n b$ $_*{}^*$-вырождена* ⊙

ЛЕММА 17.5.17. *Вектор столбец $v_1$ удовлетворяет ситеме линейных уравнений*

(17.5.21) $$v_1{}^*{}_*(f_1 - E_n b) = 0$$ ⊙

Предположим, что $\overline{\overline{e}}_1 \neq \overline{\overline{e}}_2$. Согласно теореме 12.2.13, существует единственное пассивное преобразование

$$\overline{\overline{e}}_2 = g^*{}_*\overline{\overline{e}}_1$$

УТВЕРЖДЕНИЕ 17.5.21. *Собственное значение b эндоморфизма $\overline{f}$ относительно базиса $\overline{\overline{e}}_1$ не зависит от выбора базиса $\overline{\overline{e}}_2$.* ⊙

УТВЕРЖДЕНИЕ 17.5.22. *Собственный вектор $\overline{v}$ эндоморфизма $\overline{f}$, соответствующий собственному значению b, не зависит от выбора базиса $\overline{\overline{e}}_2$.* ⊙

Доказательство. Пусть $f_i$, $i = 1, 2$, - матрица эндоморфизма f относительно базиса $\overline{\overline{e}}_i$. Пусть b - собственное значение эндоморфизма $\overline{f}$ относительно базиса $\overline{\overline{e}}_1$ и $\overline{v}$ собственный вектор эндоморфизма $\overline{f}$, соответствующее собственному значению b. Пусть $v_i$, $i = 1, 2$, - матрица вектора $\overline{v}$ относительно базиса $\overline{\overline{e}}_i$.

Предположим сперва, что $\overline{\overline{e}}_1 = \overline{\overline{e}}_2$. Тогда базисы $\overline{\overline{e}}_1$, $\overline{\overline{e}}_2$ связаны пассивным преобразованием $E_n$

$$\overline{\overline{e}}_2 = \overline{\overline{e}}_1{}_*{}^* E_n$$

Согласно теореме 17.2.4, преобразование подобия $b\overline{\overline{e}}_1$ имеет матрицу

$$E_n^{-1*}{}^*{}_* b E_n = b E_n$$

относительно базиса $\overline{\overline{e}}_2$. Следовательно, согласно теореме 17.5.8, мы доказали следующие утверждения.

ЛЕММА 17.5.23. *Матрица $f_1 - b E_n$ $_*{}^*$-вырождена* ⊙

ЛЕММА 17.5.24. *Вектор столбец $v_1$ удовлетворяет ситеме линейных уравнений*

(17.5.25) $$(f_1 - b E_n)_*{}^* v_1 = 0$$ ⊙

Предположим, что $\overline{\overline{e}}_1 \neq \overline{\overline{e}}_2$. Согласно теореме 12.5.13, существует единственное пассивное преобразование

$$\overline{\overline{e}}_2 = \overline{\overline{e}}_1{}_*{}^* g$$

и $g$ - $^*{}_*$-невырожденая матрица. Согласно теореме 12.3.5,

(17.5.22) $\qquad f_2 = g^*{}_* f_1{}^*{}_* g^{-1^*}{}_*$

Согласно теореме 17.5.6, если $b$ - собственное значение эндоморфизма $\overline{f}$, то мы должны доказать лемму 17.5.18.

**Лемма 17.5.18.** *Матрица*
$$f_2 - g b^*{}_* g^{-1^*}{}_*$$
$^*{}_*$*-вырождена* ⊙

Так как равенство
$$f_2 - g b^*{}_* g^{-1^*}{}_*$$
$$= g^*{}_* f_1{}^*{}_* g^{-1^*}{}_* - g b^*{}_* g^{-1^*}{}_*$$
$$= g^*{}_* (f_1 - E_n b)^*{}_* g^{-1^*}{}_*$$
следует из равенства (17.5.22), то лемма 17.5.18 следует из леммы 17.5.16. Следовательно, мы доказали утверждение 17.5.14.

Согласно теореме 17.5.5, если $\overline{v}$ - собственный вектор эндоморфизма $\overline{f}$, то мы должны доказать лемму 17.5.19.

**Лемма 17.5.19.** *Вектор столбец $v_2$ удовлетворяет ситеме линейных уравнений*

(17.5.23) $\quad v_2{}^*{}_* (f_2 - g b^*{}_* g^{-1^*}{}_*) = 0$

Доказательство. Согласно теореме 12.3.1,

(17.5.24) $\qquad v_2 = v_1{}^*{}_* g^{-1^*}{}_*$

Так как равенство
$$v_2{}^*{}_* (f_2 - g b^*{}_* g^{-1^*}{}_*)$$
$$= v_2{}^*{}_* (g^*{}_* f_1{}^*{}_* g^{-1^*}{}_* - g b^*{}_* g^{-1^*}{}_*)$$
$$= v_1{}^*{}_* g^{-1^*}{}_*{}_* g^*{}_* (f_1 - E_n b)^*{}_* g^{-1^*}{}_*$$
следует из равенств (17.5.22), (17.5.24), то лемма 17.5.19 следует из леммы 17.5.16.

Следовательно, мы доказали утверждение 17.5.15. ⊡

ТЕОРЕМА 17.5.27. *Пусть $V$ - левое $A$-векторное пространство строк*

---

и $g$ - $^*{}_*$-невырожденая матрица. Согласно теореме 12.6.5,

(17.5.26) $\qquad f_2 = g^{-1^*}{}_*{}^* f_1{}_* g$

Согласно теореме 17.5.8, если $b$ - собственное значение эндоморфизма $\overline{f}$, то мы должны доказать лемму 17.5.25.

**Лемма 17.5.25.** *Матрица*
$$f_2 - g^{-1^*}{}_*{}^* b g$$
$^*{}_*$*-вырождена* ⊙

Так как равенство
$$f_2 - g^{-1^*}{}_*{}^* b g$$
$$= g^{-1^*}{}_*{}^* f_1{}_* g - g^{-1^*}{}_*{}^* b g$$
$$= g^{-1^*}{}_*{}^* (f_1 - b E_n){}_* g$$
следует из равенства (17.5.26), то лемма 17.5.25 следует из леммы 17.5.23. Следовательно, мы доказали утверждение 17.5.21.

Согласно теореме 17.5.7, если $\overline{v}$ - собственный вектор эндоморфизма $\overline{f}$, то мы должны доказать лемму 17.5.26.

**Лемма 17.5.26.** *Вектор столбец $v_2$ удовлетворяет ситеме линейных уравнений*

(17.5.27) $\quad (f_2 - g^{-1^*}{}_*{}^* b g){}_* v_2 = 0$

Доказательство. Согласно теореме 12.6.1,

(17.5.28) $\qquad v_2 = g^{-1^*}{}_*{}^* v_1$

Так как равенство
$$(f_2 - g^{-1^*}{}_*{}^* b g){}_* v_2$$
$$= (g^{-1^*}{}_*{}^* f_1{}_* g - g^{-1^*}{}_*{}^* b g){}_* v_2$$
$$= g^{-1^*}{}_*{}^* (f_1 - b E_n){}_* g_* g^{-1^*}{}_*{}^* v_1$$
следует из равенств (17.5.26), (17.5.28), то лемма 17.5.26 следует из леммы 17.5.23.

Следовательно, мы доказали утверждение 17.5.22. ⊡

ТЕОРЕМА 17.5.30. *Пусть $V$ - правое $A$-векторное пространство строк*

и $\overline{\overline{e}}_1, \overline{\overline{e}}_2$ - базисы левого $A$-векторно-го пространства $V$.

**Утверждение 17.5.28.** *Собственное значение $b$ эндоморфизма $\overline{f}$ относительно базиса $\overline{\overline{e}}_1$ не зависит от выбора базиса $\overline{\overline{e}}_2$.* ⊙

**Утверждение 17.5.29.** *Собственный вектор $\overline{v}$ эндоморфизма $\overline{f}$, соответствующий собственному значению $b$, не зависит от выбора базиса $\overline{\overline{e}}_2$.* ⊙

Доказательство. Пусть $f_i$, $i = 1$, $2$, - матрица эндоморфизма $f$ относительно базиса $\overline{\overline{e}}_i$. Пусть $b$ - собственное значение эндоморфизма $\overline{f}$ относительно базиса $\overline{\overline{e}}_1$ и $\overline{v}$ собственный вектор эндоморфизма $\overline{f}$, соответствующее собственному значению $b$. Пусть $v_i$, $i = 1$, $2$, - матрица вектора $\overline{v}$ относительно базиса $\overline{\overline{e}}_i$.

Предположим сперва, что $\overline{\overline{e}}_1 = \overline{\overline{e}}_2$. Тогда базисы $\overline{\overline{e}}_1, \overline{\overline{e}}_2$ связаны пассивным преобразованием $E_n$

$$\overline{\overline{e}}_2 = E_{n*}{}^*\overline{\overline{e}}_1$$

Согласно теореме 17.1.5, преобразование подобия $\overline{\overline{\overline{e}}_1}b$ имеет матрицу

$$E_n b_*{}^* E_n^{-1*} = E_n b$$

относительно базиса $\overline{\overline{e}}_2$. Следовательно, согласно теореме 17.5.10, мы доказали следующие утверждения.

**Лемма 17.5.33.** *Матрица $f_1 - E_n b_{*}{}^*$-вырождена* ⊙

**Лемма 17.5.34.** *Вектор строка $v_1$ удовлетворяет ситеме линейных уравнений*

(17.5.29)     $v_{1*}{}^*(f_1 - E_n b) = 0$

⊙

Предположим, что $\overline{\overline{e}}_1 \neq \overline{\overline{e}}_2$. Согласно теореме 12.2.14, существует един-и $\overline{\overline{e}}_1, \overline{\overline{e}}_2$ - базисы правого $A$-векторного пространства $V$.

**Утверждение 17.5.31.** *Собственное значение $b$ эндоморфизма $\overline{f}$ относительно базиса $\overline{\overline{e}}_1$ не зависит от выбора базиса $\overline{\overline{e}}_2$.* ⊙

**Утверждение 17.5.32.** *Собственный вектор $\overline{v}$ эндоморфизма $\overline{f}$, соответствующий собственному значению $b$, не зависит от выбора базиса $\overline{\overline{e}}_2$.* ⊙

Доказательство. Пусть $f_i$, $i = 1$, $2$, - матрица эндоморфизма $f$ относительно базиса $\overline{\overline{e}}_i$. Пусть $b$ - собственное значение эндоморфизма $\overline{f}$ относительно базиса $\overline{\overline{e}}_1$ и $\overline{v}$ собственный вектор эндоморфизма $\overline{f}$, соответствующее собственному значению $b$. Пусть $v_i$, $i = 1$, $2$, - матрица вектора $\overline{v}$ относительно базиса $\overline{\overline{e}}_i$.

Предположим сперва, что $\overline{\overline{e}}_1 = \overline{\overline{e}}_2$. Тогда базисы $\overline{\overline{e}}_1, \overline{\overline{e}}_2$ связаны пассивным преобразованием $E_n$

$$\overline{\overline{e}}_2 = \overline{\overline{e}}_1{}^*{}_* E_n$$

Согласно теореме 17.2.5, преобразование подобия $b\overline{\overline{e}}_1$ имеет матрицу

$$E_n^{-1*}{}^*{}_* b E_n = b E_n$$

относительно базиса $\overline{\overline{e}}_2$. Следовательно, согласно теореме 17.5.12, мы доказали следующие утверждения.

**Лемма 17.5.37.** *Матрица $f_1 - b E_n {}_*{}^*$-вырождена* ⊙

**Лемма 17.5.38.** *Вектор строка $v_1$ удовлетворяет ситеме линейных уравнений*

(17.5.33)     $(f_1 - b E_n)^*{}_* v_1 = 0$

⊙

Предположим, что $\overline{\overline{e}}_1 \neq \overline{\overline{e}}_2$. Согласно теореме 12.5.14, существует един-ственное пассивное преобразование

$$\overline{\overline{e}}_2 = g_*{}^*\overline{\overline{e}}_1$$

и $g$ - $_*{}^*$-невырожденая матрица. Согласно теореме 12.3.6,

$$(17.5.30) \qquad f_2 = g_*{}^*f_{1*}g^{-1*}$$

Согласно теореме 17.5.10, если $b$ - собственное значение эндоморфизма $\overline{f}$, то мы должны доказать лемму 17.5.35.

> **Лемма 17.5.35.** *Матрица*
> $$f_2 - gb_*{}^*g^{-1*}$$
> $_*{}^*$-*вырождена* ⊙

Так как равенство

$$f_2 - gb_*{}^*g^{-1*}$$
$$= g_*{}^*f_{1*}g^{-1*} - gb_*{}^*g^{-1*}$$
$$= g_*{}^*(f_1 - E_n b)_*{}^*g^{-1*}$$

следует из равенства (17.5.30), то лемма 17.5.35 следует из леммы 17.5.33. Следовательно, мы доказали утверждение 17.5.28.

Согласно теореме 17.5.9, если $\overline{v}$ - собственный вектор эндоморфизма $\overline{f}$, то мы должны доказать лемму 17.5.36.

> **Лемма 17.5.36.** *Вектор строка $v_2$ удовлетворяет ситеме линейных уравнений*
> $$(17.5.31) \quad v_{2*}{}^*(f_2 - gb_*{}^*g^{-1*}) = 0$$

**Доказательство.** Согласно теореме 12.3.2,

$$(17.5.32) \qquad v_2 = v_{1*}{}^*g^{-1*}$$

Так как равенство

$$v_{2*}{}^*(f_2 - gb_*{}^*g^{-1*})$$
$$= v_{2*}{}^*(g_*{}^*f_{1*}g^{-1*} - gb_*{}^*g^{-1*})$$
$$= v_{1*}{}^*g^{-1*}{}_*{}^*g_*{}^*(f_1 - E_n b)_*{}^*g^{-1*}$$

следует из равенств (17.5.30), (17.5.32), то лемма 17.5.36 следует из леммы 17.5.33. ⊙

Следовательно, мы доказали утверждение 17.5.29. □

---

ственное пассивное преобразование

$$\overline{\overline{e}}_2 = \overline{\overline{e}}_1{}_*{}^*g$$

и $g$ - $_*{}^*$-невырожденая матрица. Согласно теореме 12.6.6,

$$(17.5.34) \qquad f_2 = g^{-1*}{}^*{}_*f_1{}^*{}_*g$$

Согласно теореме 17.5.12, если $b$ - собственное значение эндоморфизма $\overline{f}$, то мы должны доказать лемму 17.5.39.

> **Лемма 17.5.39.** *Матрица*
> $$f_2 - g^{-1*}{}^*{}_*bg$$
> $_*{}^*$-*вырождена* ⊙

Так как равенство

$$f_2 - g^{-1*}{}^*{}_*bg$$
$$= g^{-1*}{}^*{}_*f_1{}^*{}_*g - g^{-1*}{}^*{}_*bg$$
$$= g^{-1*}{}^*{}_*(f_1 - bE_n){}^*{}_*g$$

следует из равенства (17.5.34), то лемма 17.5.39 следует из леммы 17.5.37. Следовательно, мы доказали утверждение 17.5.31.

Согласно теореме 17.5.11, если $\overline{v}$ - собственный вектор эндоморфизма $\overline{f}$, то мы должны доказать лемму 17.5.40.

> **Лемма 17.5.40.** *Вектор строка $v_2$ удовлетворяет ситеме линейных уравнений*
> $$(17.5.35) \quad (f_2 - g^{-1*}{}^*{}_*bg){}^*{}_*v_2 = 0$$

**Доказательство.** Согласно теореме 12.6.2,

$$(17.5.36) \qquad v_2 = g^{-1*}{}^*{}_*v_1$$

Так как равенство

$$(f_2 - g^{-1*}{}^*{}_*bg){}^*{}_*v_2$$
$$= (g^{-1*}{}^*{}_*f_1{}^*{}_*g - g^{-1*}{}^*{}_*bg){}^*{}_*v_2$$
$$= g^{-1*}{}^*{}_*(f_1 - bE_n){}^*{}_*g{}^*{}_*g^{-1*}{}^*{}_*v_1$$

следует из равенств (17.5.34), (17.5.36), то лемма 17.5.40 следует из леммы 17.5.37. ⊙

Следовательно, мы доказали утверждение 17.5.32. □

## 17.6. Собственное значение матрицы

Определение 17.5.3 не зависит от структуры координат левого $A$-векторного пространства $V$. Когда мы выбираем базис $\overline{\overline{e}}$ левого $A$-векторного пространства $V$, мы видим дополнительную структуру, которая отражена в теоремах 17.5.5, 17.5.9.

Определение 17.5.4 не зависит от структуры координат правого $A$-векторного пространства $V$. Когда мы выбираем базис $\overline{\overline{e}}$ правого $A$-векторного пространства $V$, мы видим дополнительную структуру, которая отражена в теоремах 17.5.7, 17.5.11.

Одна и таже матрица может соответствовать различным эндоморфизмам в зависимости рассматриваем ли мы левое или правое $A$-векторного пространство. Поэтому мы должны рассмотреть все определения собственных векторов, которые могут соответствовать данной матрице. Из теорем 17.5.5, 17.5.7, 17.5.9, 17.5.11 следует, что для данной матрицы существует четыре вида собственных векторов и четыре вида собственных чисел. Пара собственный вектор и собственное значение определённого вида используется в решении задачи, где возникает соответствующий вид $A$-векторного пространства. Это утверждение отражено в названии собственного вектора и собственного значения.

Собственный вектор эндоморфизма зависит от выбора двух базисов. Однако когда мы рассматриваем собственный вектор матрицы, для нас важен не выбор базисов, а матрица пассивного преобразования между этими базисами.

---

ОПРЕДЕЛЕНИЕ 17.6.1. $A$-число $b$ называется $*_*$-**собственным значением** пары матриц $(f,g)$, где $g$ - $*_*$-невырожденная матрица, если матрица

$$f - gb^* {_*}g^{-1}{^*}{_*}$$

$*_*$-вырожденная матрица.                □

ОПРЕДЕЛЕНИЕ 17.6.2. $A$-число $b$ называется $_*{^*}$-**собственным значением** пары матриц $(f,g)$, где $g$ - $_*{^*}$-невырожденная матрица, если матрица

$$f - g^{-1}{_*}{^*}{_*}{^*}bg$$

$_*{^*}$-вырожденная матрица.                □

---

ТЕОРЕМА 17.6.3. *Пусть $f$ - $n \times n$ матрица $A$-чисел. Пусть $g$ - $*_*$-невырожденная $n \times n$ матрица $A$-чисел. Тогда ненулевые $*_*$-собственные значения пары матриц $(f,g)$ являются корнями любого уравнения*

$$\det({^*}{_*})^i_j \,(f - gb^* {_*}g^{-1}{^*}{_*}) = 0$$

ТЕОРЕМА 17.6.4. *Пусть $f$ - $n \times n$ матрица $A$-чисел. Пусть $g$ - $_*{^*}$-невырожденная $n \times n$ матрица $A$-чисел. Тогда ненулевые $_*{^*}$-собственные значения пары матриц $(f,g)$ являются корнями любого уравнения*

$$\det(_*{^*})^i_j \,(f - g^{-1}{_*}{^*}{_*}{^*}bg) = 0$$

---

ДОКАЗАТЕЛЬСТВО. Теорема является следствием теоремы 9.3.4 и определения 17.6.1.                □

ДОКАЗАТЕЛЬСТВО. Теорема является следствием теоремы 9.3.3 и определения 17.6.2.                □

---

ОПРЕДЕЛЕНИЕ 17.6.5. *Пусть $A$-число $b$ является $*_*$-собственным значением пары матриц $(f,g)$, где $g$ - $*_*$-невырожденная матрица. Столбец $v$ называется* **собственным столбцом** *пары матриц $(f,g)$,*

ОПРЕДЕЛЕНИЕ 17.6.7. *Пусть $A$-число $b$ является $_*{^*}$-собственным значением пары матриц $(f,g)$, где $g$ - $_*{^*}$-невырожденная матрица. Столбец $v$ называется* **собственным столбцом** *пары матриц $(f,g)$,* соответствующим $*_*$-*собственному значению* $b$, *если верно равенство*

$$v^*{}_*f = v^*{}_*gb^*{}_*g^{-1^*{}_*}$$

□

Определение 17.6.5 является следствием определения 17.6.1 и теоремы 17.5.5.

ТЕОРЕМА 17.6.6. *Пусть A-число* $b$ *является* $*_*$-*собственным значением пары матриц* $(f, g)$, *где* $g$ - $*_*$-*невырожденная матрица. Множество собственных столбцов пары матриц* $(f, g)$, *соответствующих* $*_*$-*собственному значению* $b$, *является левым A-векторным пространством столбцов.*

ДОКАЗАТЕЛЬСТВО. Согласно определениям 17.6.1, 17.6.5, координаты собственного столбца являются решением системы линейных уравнений

$$(17.6.1) \qquad v^*{}_*(f - gb^*{}_*g^{-1^*{}_*}) = 0$$

с $*_*$-вырожденной матрицей

$$f - gb^*{}_*g^{-1^*{}_*}$$

Согласно теореме 9.4.6, множество решений системы линейных уравнений (17.6.1) является левым пространством столбцов. □

ОПРЕДЕЛЕНИЕ 17.6.9. *Пусть A-число* $b$ *является* $*_*$-*собственным значением пары матриц* $(f, g^{-1^*{}_*})$, *где* $g$ - $*_*$-*невырожденная матрица. Строка* $v$ *называется* **собственной строкой** *пары матриц* $(f, g)$, *соответствующей* $*_*$-*собственному значению* $b$, *если верно равенство*

$$f^*{}_*v = g^{-1^*{}_*{}^*{}_*}bg^*{}_*v$$

□

Поскольку

$$g^{-1^*{}_*{}^*{}_*}bg = (g^{-1^*{}_*}b)^*{}_*g$$
$$= (g^{-1^*{}_*}b)^*{}_*(g^{-1^*{}_*})^{-1^*{}_*}$$

соответствующим $_*{}^*$-*собственному значению* $b$, *если верно равенство*

$$f_*{}^*v = g^{-1_*{}^*}bg_*{}^*v$$

□

Определение 17.6.7 является следствием определения 17.6.2 и теоремы 17.5.7.

ТЕОРЕМА 17.6.8. *Пусть A-число* $b$ *является* $_*{}^*$-*собственным значением пары матриц* $(f, g)$, *где* $g$ - $_*{}^*$-*невырожденная матрица. Множество собственных столбцов пары матриц* $(f, g)$, *соответствующих* $_*{}^*$-*собственному значению* $b$, *является правым A-векторным пространством столбцов.*

ДОКАЗАТЕЛЬСТВО. Согласно определениям 17.6.2, 17.6.7, координаты собственного столбца являются решением системы линейных уравнений

$$(17.6.2) \qquad (f - g^{-1_*{}^*}{}_*{}^*bg)_*{}^*v = 0$$

с $_*{}^*$-вырожденной матрицей

$$f - g^{-1_*{}^*}{}_*{}^*bg$$

Согласно теореме 9.4.6, множество решений системы линейных уравнений (17.6.2) является правым пространством столбцов. □

ОПРЕДЕЛЕНИЕ 17.6.11. *Пусть A-число* $b$ *является* $_*{}^*$-*собственным значением пары матриц* $(f, g^{-1_*{}^*})$, *где* $g$ - $_*{}^*$-*невырожденная матрица. Строка* $v$ *называется* **собственной строкой** *пары матриц* $(f, g)$, *соответствующей* $_*{}^*$-*собственному значению* $b$, *если верно равенство*

$$v_*{}^*f = v_*{}^*gb_*{}^*g^{-1_*{}^*}$$

□

Поскольку

$$gb_*{}^*g^{-1_*{}^*} = g_*{}^*(bg^{-1_*{}^*})$$
$$= (g^{-1_*{}^*})^{-1_*{}^*}{}_*{}^*(bg^{-1_*{}^*})$$

то определение 17.6.9 является следствием определения 17.6.1 и теоремы 17.5.11.

то определение 17.6.11 является следствием определения 17.6.2 и теоремы 17.5.9.

**ТЕОРЕМА 17.6.10.** *Пусть A-число $b$ является ${}_*^*$-собственным значением пары матриц $(f, g^{-1_*})$, где $g$ - ${}_*^*$-невырожденная матрица. Множество собственных строк пары матриц $(f, g)$, соответствующих ${}_*^*$-собственному значению $b$, является правым A-векторным пространством строк.*

**ТЕОРЕМА 17.6.12.** *Пусть A-число $b$ является ${}_*$-собственным значением пары матриц $(f, g^{-1^*})$, где $g$ - ${}_*$-невырожденная матрица. Множество собственных строк пары матриц $(f, g)$, соответствующих ${}_*$-собственному значению $b$, является левым A-векторным пространством строк.*

Доказательство. Согласно определениям 17.6.1, 17.6.9, координаты собственной строки являются решением системы линейных уравнений

$$(17.6.3) \qquad (f - g^{-1^*} {}_* {}^* bg)^*_* v = 0$$

с ${}_*^*$-вырожденной матрицей

$$f - g^{-1^*} {}^* {}_* bg$$

Согласно теореме 9.4.6, множество решений системы линейных уравнений (17.6.3) является правым пространством строк. □

Доказательство. Согласно определениям 17.6.2, 17.6.11, координаты собственной строки являются решением системы линейных уравнений

$$(17.6.4) \qquad v_* {}^* (f - g b_* {}^* g^{-1^*}) = 0$$

с ${}_*$-вырожденной матрицей

$$f - g b_* {}^* g^{-1^*}$$

Согласно теореме 9.4.6, множество решений системы линейных уравнений (17.6.4) является левым пространством строк. □

ЗАМЕЧАНИЕ 17.6.13. *Пусть*

$$(17.6.5) \qquad g^i_j \in Z(A) \quad i = 1, ..., n \quad j = 1, ..., n$$

*Тогда*

$$(17.6.6) \qquad g^{-1^*} {}_* {}^* bg = (g^{-1^*} {}_* {}^* g)b = E_n b$$

$$(17.6.7) \qquad g b_* {}^* g^{-1^*} = b(g_* {}^* g^{-1^*}) = b E_n$$

*Из равенств (17.6.6), (17.6.7) следует, что, если выполнено условие (17.6.5), то мы можем упростить определения собственного значения и собственного вектора матрицы и получить определения, которые будут соответствовать определениям в коммутативной алгебре.* □

**ОПРЕДЕЛЕНИЕ 17.6.14.** *A-число $b$ называется ${}_*$-собственным значением матрицы $f$, если матрица $f - bE_n$ - ${}_*$-вырожденная матрица.* □

**ОПРЕДЕЛЕНИЕ 17.6.17.** *A-число $b$ называется ${}_*^*$-собственным значением матрицы $f$, если матрица $f - bE_n$ - ${}_*^*$-вырожденная матрица.* □

**ТЕОРЕМА 17.6.15.** *Пусть $f$ - $n \times n$ матрица A-чисел и*

$${}_*^* \operatorname{rank} a = k < n$$

**ТЕОРЕМА 17.6.18.** *Пусть $f$ - $n \times n$ матрица A-чисел и*

$${}^*_* \operatorname{rank} a = k < n$$

*Тогда 0 является $_*{}^*$-собственным значением кратности $n - k$.*

**Доказательство.** Теорема является следствием равенства

$$a = a - 0E$$

определения 17.6.14 и теоремы 9.4.4. □

---

*Тогда 0 является $^*{}_*$-собственным значением кратности $n - k$.*

**Доказательство.** Теорема является следствием равенства

$$a = a - 0E$$

определения 17.6.17 и теоремы 9.4.4. □

---

**ТЕОРЕМА 17.6.16.** *Пусть $f$ - $n \times n$ матрица $A$-чисел. Тогда ненулевые $_*{}^*$-собственные значения матрицы $f$ являются корнями любого уравнения*

$$\det({}_*{}^*)^i_j \, (a - bE) = 0$$

**Доказательство.** Теорема является следствием теоремы 9.3.3 и определения 17.5.3. □

---

**ТЕОРЕМА 17.6.19.** *Пусть $f$ - $n \times n$ матрица $A$-чисел. Тогда ненулевые $^*{}_*$-собственные значения матрицы $f$ являются корнями любого уравнения*

$$\det({}^*{}_*)^i_j \, (a - bE) = 0$$

**Доказательство.** Теорема является следствием теоремы 9.3.4 и определения 17.5.3. □

---

**ОПРЕДЕЛЕНИЕ 17.6.20.** *Пусть $A$-число $b$ является $_*{}^*$-собственным значением матрицы $f$. Столбец $v$ называется **собственным столбцом** матрицы $f$, соответствующим $_*{}^*$-собственному значению $b$, если верно равенство*

$$(17.6.8) \qquad f_*{}^*v = bv$$

□

---

**ОПРЕДЕЛЕНИЕ 17.6.21.** *Пусть $A$-число $b$ является $^*{}_*$-собственным значением матрицы $f$. Столбец $v$ называется **собственным столбцом** матрицы $f$, соответствующим $^*{}_*$-собственному значению $b$, если верно равенство*

$$(17.6.9) \qquad v^*{}_*f = vb$$

□

---

**ОПРЕДЕЛЕНИЕ 17.6.22.** *Пусть $A$-число $b$ является $_*{}^*$-собственным значением матрицы $f$. Строка $v$ называется **собственной строкой** матрицы $f$, соответствующей $_*{}^*$-собственному значению $b$, если верно равенство*

$$(17.6.10) \qquad v_*{}^*f = vb$$

□

---

**ОПРЕДЕЛЕНИЕ 17.6.23.** *Пусть $A$-число $b$ является $^*{}_*$-собственным значением матрицы $f$. Строка $v$ называется **собственной строкой** матрицы $f$, соответствующей $^*{}_*$-собственному значению $b$, если верно равенство*

$$(17.6.11) \qquad f^*{}_*v = bv$$

□

---

**ТЕОРЕМА 17.6.24.** *Пусть $A$-число $b$ является $_*{}^*$-собственным значением матрицы $a$. Множество собственных столбцов матрицы $a$, соответствующих $_*{}^*$-собственному*

---

**ТЕОРЕМА 17.6.25.** *Пусть $A$-число $b$ является $^*{}_*$-собственным значением матрицы $a$. Множество собственных столбцов матрицы $a$, соответствующих $^*{}_*$-собственному* *значению b, является правым A-век-
торным пространством столбцов.*

Доказательство. Согласно опре-
делениям 17.6.14, 17.6.20, координаты
собственного столбца являются реше-
нием системы линейных уравнений

$$(17.6.12) \qquad (f - bE_n)_*{}^*v = 0$$

с $_*{}^*$-вырожденной матрицей $f - bE_n$.
Согласно теореме 9.4.6, множество ре-
шений системы линейных уравнений
(17.6.12) является правым A-векторным
пространством столбцов.          □

---

*значению b, является левым A-век-
торным пространством столбцов.*

Доказательство. Согласно опре-
делениям 17.6.17, 17.6.21, координаты
собственного столбца являются реше-
нием системы линейных уравнений

$$(17.6.13) \qquad v^*{}_*(f - E_n b) = 0$$

с $^*{}_*$-вырожденной матрицей $f - E_n b$.
Согласно теореме 9.4.6, множество ре-
шений системы линейных уравнений
(17.6.13) является левым A-векторным
пространством столбцов.          □

---

Теорема 17.6.26. *Пусть A-чис-
ло b является $_*{}^*$-собственным зна-
чением матрицы a. Множество соб-
ственных строк матрицы a, соот-
ветствующих $_*{}^*$-собственному зна-
чению b, является левым A-вектор-
ным пространством строк.*

Доказательство.          Согласно
определениям 17.6.14, 17.6.22, коор-
динаты собственной строки являются
решением системы линейных уравне-
ний

$$(17.6.14) \qquad v^*{}^*(f - E_n b) = 0$$

с $_*{}^*$-вырожденной матрицей $f - E_n b$.
Согласно теореме 9.4.6, множество ре-
шений системы линейных уравнений
(17.6.14) является левым A-вектор-
ным пространством строк.          □

---

Теорема 17.6.27. *Пусть A-чис-
ло b является $^*{}_*$-собственным зна-
чением матрицы a. Множество соб-
ственных строк матрицы a, соот-
ветствующих $^*{}_*$-собственному зна-
чению b, является правым A-вектор-
ным пространством строк.*

Доказательство.          Согласно
определениям 17.6.17, 17.6.23, коор-
динаты собственной строки являются
решением системы линейных уравне-
ний

$$(17.6.15) \qquad (f - bE_n)^*{}_*v = 0$$

с $^*{}_*$-вырожденной матрицей $f - bE_n$.
Согласно теореме 9.4.6, множество ре-
шений системы линейных уравнений
(17.6.15) является правым A-векто-
ным пространством строк.          □

# Левые и правые собственные значения

Пусть $A$ - ассоциативная $D$-алгебра с делением.

## 18.1. Левые и правые ${}_*$*-собственные значения

ОПРЕДЕЛЕНИЕ 18.1.1. *Пусть $a_2$ - $n \times n$ матрица, которая ${}_*$*-подобна диагональной матрице $a_1$*

$$a_1 = \mathrm{diag}(b(1), ..., b(n))$$

*Следовательно, существует ${}_*$*-невырожденая матрица $u_2$ такая, что*

(18.1.1) $\qquad u_{2*}{}^* a_{2*}{}^* u_2^{-1*}{}^* = a_1$

(18.1.2) $\qquad u_{2*}{}^* a_2 = a_{1*}{}^* u_2$

*Строка $u_2^i$ матрицы $u_2$ удовлетворяет равенству*

(18.1.3) $\qquad u_2^i{}_* {}^* a_2 = b(i) u_2^i$

*$A$-число $b(i)$ называется левым ${}_*$***-собственным значением**. *Вектор-строка $u_2^i$ называется **собственной строкой**, соответствующей левому ${}_*$*-собственному значению $b(i)$.* $\qquad\square$

ОПРЕДЕЛЕНИЕ 18.1.2. *Пусть $a_2$ - $n \times n$ матрица, которая ${}_*$*-подобна диагональной матрице $a_1$*

$$a_1 = \mathrm{diag}(b(1), ..., b(n))$$

*Следовательно, существует ${}_*$*-невырожденая матрица $u_2$ такая, что*

(18.1.4) $\qquad u_2^{-1*}{}^* {}_* {}^* a_{2*}{}^* u_2 = a_1$

(18.1.5) $\qquad a_{2*}{}^* u_2 = u_{2*}{}^* a_1$

*Столбец $u_{2i}$ матрицы $u_2$ удовлетворяет равенству*

(18.1.6) $\qquad a_{2*}{}^* u_{2i} = u_{2i} b(i)$

*$A$-число $b(i)$ называется правым ${}_*$***-собственным значением**. *Вектор-столбец $u_{2i}$ называется **собственным столбцом**, соответствующим правому ${}_*$*-собственному значению $b(i)$.* $\qquad\square$

ОПРЕДЕЛЕНИЕ 18.1.3. *Множество ${}_*$*-$\mathrm{spec}(a)$ всех левых и правых ${}_*$*-собственных значений называется ${}_*$*-спектром матрицы $a$.* $\qquad\square$

ТЕОРЕМА 18.1.4. *Пусть $A$-число $b$ - левое ${}_*$*-собственное значение матрицы $a_2$.[18.1] Тогда любое $A$-число, $A$-сопряжённое с $A$-числом $b$, является левым ${}_*$*-собственным значением матрицы $a_2$.*

ТЕОРЕМА 18.1.5. *Пусть $A$-число $b$ - правое ${}_*$*-собственное значение матрицы $a_2$.[18.1] Тогда любое $A$-число, $A$-сопряжённое с $A$-числом $b$, является правым ${}_*$*-собственным значением матрицы $a_2$.*

ДОКАЗАТЕЛЬСТВО. Пусть $v$ - собственная строка, соответствующей левому ${}_*$*-собственному значению $b$. Пусть $c \neq 0$ - $A$-число. Утверждение теоремы является следствием равенства

$$cv_*{}^* a_2 = cbv = cbc^{-1} cv$$

ДОКАЗАТЕЛЬСТВО. Пусть $v$ - собственный столбец, соответствующим правому ${}_*$*-собственному значению $b$. Пусть $c \neq 0$ - $A$-число. Утверждение теоремы является следствием равенства

$$a_{2*}{}^* vc = vbc = vc\,c^{-1} bc$$

---

[18.1] Смотри аналогичное утверждение и доказательство на странице [20]-375.

□　　　　　　　　　　　　　　　　　　　　　□

**Теорема 18.1.6.** *Пусть $v$ - собственная строка, соответствующая левому $_*^*$-собственному значению $b$. Вектор $vc$ является собственной строкой , соответствующей левому $_*^*$-собственному значению $b$ тогда и только тогда, когда $A$-число $c$ коммутирует со всеми элементами матрицы $a_2$.*

Доказательство. Вектор $vc$ является собственной строкой, соответствующей левому $_*^*$-собственному значению $b$ тогда и только тогда, когда верно следующее равенство

$$(18.1.7) \qquad bv_ic = v_j a_2{}_i^j c = v_j ca_2{}_i^j$$

Равенство (18.1.7) верно тогда и только тогда, когда $A$-число $c$ коммутирует со всеми элементами матрицы $a_2$. □

Из теорем 18.1.4, 18.1.6 следует, что множество собственных строк, соответствующих данному левому $_*^*$-собственному значению, не является $A$-векторным пространством.

**Теорема 18.1.8.** *Собственный вектор $u_2^i$ матрицы $a_2$, соответствующий левому $_*^*$-собственному значению $b(i)$, не зависит от выбора базиса $\overline{\overline{e}}_2$*

$$(18.1.9) \qquad e_1^i = u_{2*}^i {}^*\overline{\overline{e}}_2$$

Доказательство. Согласно теореме 10.3.9, мы можем рассматривать матрицу $a_2$ как матрицу автоморфизма $a$ левого $A$-векторного пространства относительно базиса $\overline{\overline{e}}_2$. Согласно теореме 12.3.6, матрица $u_2$ является матрицей пассивного преобразования

$$(18.1.10) \qquad \overline{\overline{e}}_1 = u_{2*} {}^*\overline{\overline{e}}_2$$

и матрица $a_1$ является матрицей автоморфизма $a$ левого $A$-векторного пространства относительно базиса $\overline{\overline{e}}_1$.

Матрица координат базиса $\overline{\overline{e}}_1$ относительно базиса $\overline{\overline{e}}_1$ - это единичная

**Теорема 18.1.7.** *Пусть $v$ - собственный столбец, соответствующий правому $_*^*$-собственному значению $b$. Вектор $cv$ является собственным столбцом , соответствующим правому $_*^*$-собственному значению $b$ тогда и только тогда, когда $A$-число $c$ коммутирует со всеми элементами матрицы $a_2$.*

Доказательство. Вектор $cv$ является собственным столбцом, соответствующим правому $_*^*$-собственному значению $b$ тогда и только тогда, когда верно следующее равенство

$$(18.1.8) \qquad cv^i b = ca_2{}_j^i v^i = a_j^i cv^j$$

Равенство (18.1.8) верно тогда и только тогда, когда $A$-число $c$ коммутирует со всеми элементами матрицы $a_2$. □

Из теорем 18.1.5, 18.1.7 следует, что множество собственных столбцов, соответствующих данному правому $_*^*$-собственному значению, не является $A$-векторным пространством.

**Теорема 18.1.9.** *Собственный вектор $u_{2i}$ матрицы $a_2$, соответствующий правому $_*^*$-собственному значению $b(i)$, не зависит от выбора базиса $\overline{\overline{e}}_2$*

$$(18.1.15) \qquad e_{1i} = \overline{\overline{e}}_{2*} {}^*u_{2i}$$

Доказательство. Согласно теореме 11.3.8, мы можем рассматривать матрицу $a_2$ как матрицу автоморфизма $a$ правого $A$-векторного пространства относительно базиса $\overline{\overline{e}}_2$. Согласно теореме 12.6.5, матрица $u_2$ является матрицей пассивного преобразования

$$(18.1.16) \qquad \overline{\overline{e}}_1 = \overline{\overline{e}}_{2*} {}^*u_2$$

и матрица $a_1$ является матрицей автоморфизма $a$ правого $A$-векторного пространства относительно базиса $\overline{\overline{e}}_1$.

Матрица координат базиса $\overline{\overline{e}}_1$ относительно базиса $\overline{\overline{e}}_1$ - это единич-матрица $E_n$. Следовательно, вектор $e_1^i$ является собственной строкой, соответствующей левому ∗∗-собственному значению $b(i)$

$$(18.1.11) \qquad e_{1*}{}^{*}a_1 = b(i)e_1^i$$

Равенство

$$(18.1.12) \quad e_{1*}{}^{*}u_{2*}{}^{*}a_{2*}{}^{*}u_2^{-1*}{}^{*} = b(i)e_1^i$$

является следствием равенств (18.1.1), (18.1.11). Равенство

$$(18.1.13) \qquad e_{1*}{}^{*}u_{2*}{}^{*}a_2 = b(i)e_{1*}^i{}^{*}u_2$$

является следствием равенства (18.1.12). Равенство

$$(18.1.14) \qquad e_{1*}{}^{*}u_2 = u_2^i$$

является следствием равенств (18.1.3), (18.1.13). Теорема является следствием равенств (18.1.10), (18.1.14). □

Равенство (18.1.9) является следствием равенства (18.1.10). Однако смысл этих равенств различный. Равенство (18.1.10) является представлением пассивного преобразования. Мы видим в равенстве (18.1.9) разложение собственного вектора для данного левого ∗∗-собственного значения относительно выбранного базиса. Задача доказательства теоремы 18.1.8 - раскрыть смысл равенства (18.1.9).

**Теорема 18.1.10.** *Любое левое ∗∗-собственное значение $b(i)$ является ∗∗-собственным значением. Собственная строка $u^i$, соответствующая левому ∗∗-собственному значению $b(i)$, является собственной строкой, соответствующей ∗∗-собственному значению $b(i)$.*

Доказательство. Все элементы строки матрицы

$$(18.1.21) \qquad a_1 - b(i)E_n$$

равны нулю. Следовательно, матрица (18.1.21) ∗∗-вырождена и левое ∗∗-собственное значение $b(i)$ является ∗∗-

ная матрица $E_n$. Следовательно, вектор $e_{1i}$ является собственным столбцом, соответствующим правому ∗∗-собственному значению $b(i)$

$$(18.1.17) \qquad a_{1*}{}^{*}e_{1i} = e_{1i}b(i)$$

Равенство

$$(18.1.18) \qquad u_2^{-1*}{}^{*}a_{2*}{}^{*}u_{2*}{}^{*}e_{1i}$$
$$= e_{1i}b(i)$$

является следствием равенств (18.1.4), (18.1.17). Равенство

$$(18.1.19) \quad a_{2*}{}^{*}u_{2*}{}^{*}e_{1i} = u_{2*}{}^{*}e_{1i}b(i)$$

является следствием равенства (18.1.18). Равенство

$$(18.1.20) \qquad u_{2*}{}^{*}e_{1i} = u_{2i}$$

является следствием равенств (18.1.6), (18.1.19). Теорема является следствием равенств (18.1.16), (18.1.20). □

Равенство (18.1.15) является следствием равенства (18.1.16). Однако смысл этих равенств различный. Равенство (18.1.16) является представлением пассивного преобразования. Мы видим в равенстве (18.1.15) разложение собственного вектора для данного правого ∗∗-собственного значения относительно выбранного базиса. Задача доказательства теоремы 18.1.9 - раскрыть смысл равенства (18.1.15).

**Теорема 18.1.11.** *Любое правое ∗∗-собственное значение $b(i)$ является ∗∗-собственным значением. Собственный столбец $u_i$, соответствующий правому ∗∗-собственному значению $b(i)$, является собственным столбцом, соответствующим ∗∗-собственному значению $b(i)$.*

Доказательство. Все элементы строки матрицы

$$(18.1.27) \qquad a_1 - b(i)E_n$$

равны нулю. Следовательно, матрица (18.1.27) ∗∗-вырождена и правое ∗∗-собственное значение $b(i)$ является ∗∗-

собственным значением. Поскольку

$$(18.1.22) \qquad e_1{}^j_i = \delta^j_i$$

то собственный вектор $e_1^i$, соответствующий левому $_*$*-собственному значению $b(i)$, является собственным вектором, соответствующим $_*$*-собственному значению $b(i)$. Равенство

$$(18.1.23) \qquad e_1^i{}_*{}^*a_1 = e_1^i b(i)$$

является следствием определения 17.6.22.

Из равенства (18.1.1) следует, что мы можем записать матрицу (18.1.21) в виде

$$u_{2*}{}^*a_{2*}{}^*u_2^{-1}{}^* - b(i)E_n$$
$$= u_{2*}{}^*(a_2 - u_2^{-1}{}^*{}^*b(i)u_2)_*{}^*u_2^{-1}{}^*$$

Так как матрица (18.1.21) $_*$*-вырождена, то матрица

$$(18.1.24) \qquad a_2 - u_2^{-1}{}^*{}_*{}^*b(i)u_2$$

также $_*$*-вырождена. Следовательно, левое $_*$*-собственное значение $b(i)$ матрицы $a_2$ является $_*$*-собственным значением пары матриц $(a_2, u_2)$.

Равенства

$$(18.1.25) \quad e_1^i{}_*{}^*u_{2*}{}^*a_{2*}{}^*u_2^{-1}{}^* = e_1^i b(i)$$

$$e_1^i{}_*{}^*u_{2*}{}^*a_2$$

$$(18.1.26) \quad = e_1^i b(i)_*{}^*u_2$$
$$= e_1^i{}_*{}^*u_{2*}{}^*u_2^{-1}{}^*{}_*{}^*b(i)u_2$$

являются следствием равенств (18.1.1), (18.1.23). Из равенств (18.1.9), (18.1.26) и из теоремы 18.1.8 следует, что собственная строка $u_2^i$ матрицы $a_2$, соответствующий левому $_*$*-собственному значению $b(i)$, является собственной строкой $u_2^i$ пары матриц $(a_2, u_2)$, соответствующей $_*$*-собственному значению $b(i)$. □

ЗАМЕЧАНИЕ 18.1.12. *Равенство*

$$(18.1.33) \quad u_1^i{}_*{}^*u_2^{-1}{}^*{}_*{}^*b(i)u_2 = b(i)u_1^i$$

*является следствием теоремы*

собственным значением. Поскольку

$$(18.1.28) \qquad e_1{}^j_i = \delta^j_i$$

то собственный вектор $e_{1i}$, соответствующий правому $_*$*-собственному значению $b(i)$, является собственным вектором, соответствующим $_*$*-собственному значению $b(i)$. Равенство

$$(18.1.29) \qquad a_{1*}{}^*e_{1i} = b(i)e_{1i}$$

является следствием определения 17.6.20.

Из равенства (18.1.4) следует, что мы можем записать матрицу (18.1.27) в виде

$$u_2^{-1}{}^*{}_*{}^*a_{2*}{}^*u_2 - b(i)E_n$$
$$= u_2^{-1}{}^*{}_*{}^*(a_2 - u_2 b(i)_*{}^*u_2^{-1}{}^*)_*{}^*u_2$$

Так как матрица (18.1.27) $_*$*-вырождена, то матрица

$$(18.1.30) \qquad a_2 - u_2 b(i)_*{}^*u_2^{-1}{}^*$$

также $_*$*-вырождена. Следовательно, правое $_*$*-собственное значение $b(i)$ матрицы $a_2$ является $_*$*-собственным значением пары матриц $(a_2, u_2)$.

Равенства
$$(18.1.31)$$
$$u_2^{-1}{}^*{}_*{}^*a_{2*}{}^*u_2{}_*{}^*e_{1i} = b(i)e_{1i}$$

$$a_{2*}{}^*u_{2*}{}^*e_{1i}$$

$$(18.1.32) \quad = u_{2*}{}^*b(i)e_{1i}$$
$$= u_2 b(i)_*{}^*u_2^{-1}{}^*{}_*{}^*u_{2*}{}^*e_{1i}$$

являются следствием равенств (18.1.4), (18.1.29). Из равенств (18.1.15), (18.1.32) и из теоремы 18.1.9 следует, что собственный столбец $u_{2i}$ матрицы $a_2$, соответствующий правому $_*$*-собственному значению $b(i)$, является собственным столбцом $u_{2i}$ пары матриц $(a_2, u_2)$, соответствующим $_*$*-собственному значению $b(i)$. □

ЗАМЕЧАНИЕ 18.1.13. *Равенство*
$$(18.1.35)$$
$$u_2 b(i)_*{}^*u_2^{-1}{}^*{}_*{}^*u_{1i} = u_{1i} b(i)$$

*является следствием теоремы*

*18.1.10.* Равенство (18.1.33) кажется необычным. Однако, если мы его слегка изменим

$$(18.1.34) \qquad u_{1*}^i{}^*u_2^{-1}{}^{*}b(i)$$
$$= b(i)u_{1*}^i{}^*u_2^{-1}{}^{*}$$

то равенство (18.1.34) следует из равенства (18.1.9).    □

*18.1.11.* Равенство (18.1.35) кажется необычным. Однако, если мы его слегка изменим

$$(18.1.36) \qquad b(i)u_2^{-1}{}^{*}{}_*u_{1i}$$
$$= u_2^{-1}{}^{*}{}_*u_{1i}b(i)$$

то равенство (18.1.36) следует из равенства (18.1.15).    □

## 18.2. Левые и правые $^*_*$-собственные значения

Определение 18.2.1. *Пусть $a_2$ - $n \times n$ матрица, которая $^*_*$-подобна диагональной матрице $a_1$*

$$a_1 = \mathrm{diag}(b(1), ..., b(n))$$

*Следовательно, существует $^*_*$-невырожденая матрица $u_2$ такая, что*

$$(18.2.1) \qquad u_2^*{}_*a_2{}^*{}_*u_2^{-1}{}^{*} = a_1$$

$$(18.2.2) \qquad u_2^*{}_*a_2 = a_1{}^*{}_*u_2$$

*Столбец $u_{2i}$ матрицы $u_2$ удовлетворяет равенству*

$$(18.2.3) \qquad u_{2i}^*{}_*a_2 = b(i)u_{2i}$$

*A-число $b(i)$ называется левым $^*_*$-***собственным значением***. Вектор-столбец $u_{2i}$ называется ***собственным столбцом***, соответствующим левому $^*_*$-собственному значению $b(i)$.*    □

Определение 18.2.2. *Пусть $a_2$ - $n \times n$ матрица, которая $^*_*$-подобна диагональной матрице $a_1$*

$$a_1 = \mathrm{diag}(b(1), ..., b(n))$$

*Следовательно, существует $^*_*$-невырожденая матрица $u_2$ такая, что*

$$(18.2.4) \qquad u_2^{-1}{}^{*}{}^*{}_*a_2{}^*{}_*u_2 = a_1$$

$$(18.2.5) \qquad a_2{}^*{}_*u_2 = u_2{}^*{}_*a_1$$

*Строка $u_2^i$ матрицы $u_2$ удовлетворяет равенству*

$$(18.2.6) \qquad a_2{}^*{}_*u_2^i = u_2^i b(i)$$

*A-число $b(i)$ называется правым $^*_*$-***собственным значением***. Вектор-строка $u_2^i$ называется ***собственной строкой***, соответствующей правому $^*_*$-собственному значению $b(i)$.*    □

Определение 18.2.3. *Множество $^*_*$-$\mathrm{spec}(a)$ всех левых и правых $^*_*$-собственных значений называется $^*_*$-спектром матрицы $a$.*    □

Теорема 18.2.4. *Пусть A-число $b$ - левое $^*_*$-собственное значение матрицы $a_2$.[18.2] Тогда любое A-число, A-сопряжённое с A-числом $b$, является левым $^*_*$-собственным значением матрицы $a_2$.*

Доказательство. Пусть $v$ - собственный столбец, соответствующим ле-

Теорема 18.2.5. *Пусть A-число $b$ - правое $^*_*$-собственное значение матрицы $a_2$.[18.2] Тогда любое A-чис-ло, A-сопряжённое с A-числом $b$, является правым $^*_*$-собственным значением матрицы $a_2$.*

Доказательство. Пусть $v$ - собственная строка, соответствующей пра-

---

[18.2] Смотри аналогичное утверждение и доказательство на странице [20]-375.

вому $*$-собственному значению $b$. Пусть $c \neq 0$ - $A$-число. Утверждение теоремы является следствием равенства

$$cv^*{}_*a_2 = cbv = cbc^{-1}\,cv$$

$\square$

Теорема 18.2.6. *Пусть $v$ - собственный столбец, соответствующий левому $*_*$-собственному значению $b$. Вектор $vc$ является собственным столбцом , соответствующим левому $*_*$-собственному значению $b$ тогда и только тогда, когда $A$-число $c$ коммутирует со всеми элементами матрицы $a_2$.*

Доказательство. Вектор $vc$ является собственным столбцом, соответствующим левому $*_*$-собственному значению $b$ тогда и только тогда, когда верно следующее равенство

(18.2.7)    $$bv^i c = v^j a_j^i c = v^j c a_j^i$$

Равенство (18.2.7) верно тогда и только тогда, когда $A$-число $c$ коммутирует со всеми элементами матрицы $a_2$.   $\square$

Из теорем 18.2.4, 18.2.6 следует, что множество собственных столбцов, соответствующих данному левому $*_*$-собственному значению, не является $A$-векторным пространством.

Теорема 18.2.8. *Собственный вектор $u_{2i}$ матрицы $a_2$, соответствующий левому $*_*$-собственному значению $b(i)$, не зависит от выбора базиса $\overline{\overline{e}}_2$*

(18.2.9)    $$e_{1i} = u_{2i}{}^*{}_*\overline{\overline{e}}_2$$

Доказательство. Согласно теореме 10.3.8, мы можем рассматривать матрицу $a_2$ как матрицу автоморфизма $a$ левого $A$-векторного пространства относительно базиса $\overline{\overline{e}}_2$. Согласно теореме 12.3.5, матрица $u_2$ является матрицей пассивного преобразования

(18.2.10)    $$\overline{\overline{e}}_1 = u_2{}^*{}_*\overline{\overline{e}}_2$$

вому $*$-собственному значению $b$. Пусть $c \neq 0$ - $A$-число. Утверждение теоремы является следствием равенства

$$a_2{}^*{}_*vc = vbc = vc\,c^{-1}bc$$

$\square$

Теорема 18.2.7. *Пусть $v$ - собственная строка, соответствующая правому $*_*$-собственному значению $b$. Вектор $cv$ является собственной строкой , соответствующей правому $*_*$-собственному значению $b$ тогда и только тогда, когда $A$-число $c$ коммутирует со всеми элементами матрицы $a_2$.*

Доказательство. Вектор $cv$ является собственной строкой, соответствующей правому $*_*$-собственному значению $b$ тогда и только тогда, когда верно следующее равенство

(18.2.8)    $$cv_i b = ca_i^j v_j = a_i^j cv_j$$

Равенство (18.2.8) верно тогда и только тогда, когда $A$-число $c$ коммутирует со всеми элементами матрицы $a_2$.   $\square$

Из теорем 18.2.5, 18.2.7 следует, что множество собственных строк, соответствующих данному правому $*_*$-собственному значению, не является $A$-векторным пространством.

Теорема 18.2.9. *Собственный вектор $u_2^i$ матрицы $a_2$, соответствующий правому $*_*$-собственному значению $b(i)$, не зависит от выбора базиса $\overline{\overline{e}}_2$*

(18.2.15)    $$e_1^i = \overline{\overline{e}}_2{}^*{}_*u_2^i$$

Доказательство. Согласно теореме 11.3.9, мы можем рассматривать матрицу $a_2$ как матрицу автоморфизма $a$ правого $A$-векторного пространства относительно базиса $\overline{\overline{e}}_2$. Согласно теореме 12.6.6, матрица $u_2$ является матрицей пассивного преобразования

(18.2.16)    $$\overline{\overline{e}}_1 = \overline{\overline{e}}_2{}^*{}_*u_2$$

и матрица $a_1$ является матрицей автоморфизма $a$ левого $A$-векторного пространства относительно базиса $\overline{\overline{e}}_1$.

Матрица координат базиса $\overline{e}_1$ относительно базиса $\overline{\overline{e}}_1$ - это единичная матрица $E_n$. Следовательно, вектор $e_{1i}$ является собственным столбцом, соответствующим левому $^*_*$-собственному значению $b(i)$

(18.2.11) $\qquad e_{1i}{}^*{}_*a_1 = b(i)e_{1i}$

Равенство

(18.2.12)
$$e_{1i}{}^*{}_*u_2{}^*{}_*a_2{}^*{}_*u_2^{-1^*} = b(i)e_{1i}$$

является следствием равенств (18.2.1), (18.2.11). Равенство

(18.2.13) $\quad e_{1i}{}^*{}_*u_2{}^*{}_*a_2 = b(i)e_{1i}{}^*{}_*u_2$

является следствием равенства (18.2.12). Равенство

(18.2.14) $\qquad e_{1i}{}^*{}_*u_2 = u_{2i}$

является следствием равенств (18.2.3), (18.2.13). Теорема является следствием равенств (18.2.10), (18.2.14). $\qquad\square$

Равенство (18.2.9) является следствием равенства (18.2.10). Однако смысл этих равенств различный. Равенство (18.2.10) является представлением пассивного преобразования. Мы видим в равенстве (18.2.9) разложение собственного вектора для данного левого $^*_*$-собственного значения относительно выбранного базиса. Задача доказательства теоремы 18.2.8 - раскрыть смысл равенства (18.2.9).

---

ТЕОРЕМА 18.2.10. *Любое левое $^*_*$-собственное значение $b(i)$ является $^*_*$-собственным значением. Собственный столбец $u_i$, соответствующий левому $^*_*$-собственному значению $b(i)$, является собственным столбцом, соответствующим $^*_*$-собственному значению $b(i)$.*

ДОКАЗАТЕЛЬСТВО. Все элементы строки матрицы

(18.2.21) $\qquad a_1 - b(i)E_n$

---

и матрица $a_1$ является матрицей автоморфизма $a$ правого $A$-векторного пространства относительно базиса $\overline{\overline{e}}_1$.

Матрица координат базиса $\overline{\overline{e}}_1$ относительно базиса $\overline{\overline{e}}_1$ - это единичная матрица $E_n$. Следовательно, вектор $e_1^i$ является собственной строкой, соответствующей правому $^*_*$-собственному значению $b(i)$

(18.2.17) $\qquad a_1{}^*{}_*e_1^i = e_1^ib(i)$

Равенство

(18.2.18) $\quad u_2^{-1^*}{}^*{}_*a_2{}^*{}_*u_2{}^*{}_*e_1^i = e_1^ib(i)$

является следствием равенств (18.2.4), (18.2.17). Равенство

(18.2.19) $\qquad a_2{}^*{}_*u_2{}^*{}_*e_1^i = u_2{}^*{}_*e_1^ib(i)$

является следствием равенства (18.2.18). Равенство

(18.2.20) $\qquad u_2{}^*{}_*e_1^i = u_2^i$

является следствием равенств (18.2.6), (18.2.19). Теорема является следствием равенств (18.2.16), (18.2.20). $\qquad\square$

Равенство (18.2.15) является следствием равенства (18.2.16). Однако смысл этих равенств различный. Равенство (18.2.16) является представлением пассивного преобразования. Мы видим в равенстве (18.2.15) разложение собственного вектора для данного правого $^*_*$-собственного значения относительно выбранного базиса. Задача доказательства теоремы 18.2.9 - раскрыть смысл равенства (18.2.15).

---

ТЕОРЕМА 18.2.11. *Любое правое $^*_*$-собственное значение $b(i)$ является $^*_*$-собственным значением. Собственная строка $u^i$, соответствующая правому $^*_*$-собственному значению $b(i)$, является собственной строкой, соответствующей $^*_*$-собственному значению $b(i)$.*

ДОКАЗАТЕЛЬСТВО. Все элементы строки матрицы

(18.2.27) $\qquad a_1 - b(i)E_n$

равны нулю. Следовательно, матрица (18.2.21) $^*_*$-вырождена и левое $^*_*$-собственное значение $b(i)$ является $^*_*$-собственным значением. Поскольку

$$(18.2.22) \qquad e_{1\,i}{}^j = \delta^j_i$$

то собственный вектор $e_{1\,i}$, соответствующий левому $^*_*$-собственному значению $b(i)$, является собственным вектором, соответствующим $^*_*$-собственному значению $b(i)$. Равенство

$$(18.2.23) \qquad e_{1\,i}{}^*_* a_1 = e_{1\,i} b(i)$$

является следствием определения 17.6.21.

Из равенства (18.2.1) следует, что мы можем записать матрицу (18.2.21) в виде

$$u_2{}^*_* a_2{}^*_* u_2^{-1^*} - b(i) E_n$$
$$= u_2{}^*_* (a_2 - u_2^{-1^*}{}^*_* b(i) u_2)^*_* u_2^{-1^*}$$

Так как матрица (18.2.21) $^*_*$-вырождена, то матрица

$$(18.2.24) \qquad a_2 - u_2^{-1^*}{}^*_* b(i) u_2$$

также $^*_*$-вырождена. Следовательно, левое $^*_*$-собственное значение $b(i)$ матрицы $a_2$ является $^*_*$-собственным значением пары матриц $(a_2, u_2)$.

Равенства
$$(18.2.25)$$
$$e_{1\,i}{}^*_* u_2{}^*_* a_2{}^*_* u_2^{-1^*} = e_{1\,i} b(i)$$

$$e_{1\,i}{}^*_* u_2{}^*_* a_2$$
$$(18.2.26) \quad = e_{1\,i} b(i)^*_* u_2$$
$$= e_{1\,i}{}^*_* u_2{}^*_* u_2^{-1^*}{}^*_* b(i) u_2$$

являются следствием равенств (18.2.1), (18.2.23). Из равенств (18.2.9), (18.2.26) и из теоремы 18.2.8 следует, что собственный столбец $u_{2\,i}$ матрицы $a_2$, соответствующий левому $^*_*$-собственному значению $b(i)$, является собственным столбцом $u_{2\,i}$ пары матриц $(a_2, u_2)$, соответствующим $^*_*$-собственному значению $b(i)$.    $\square$

равны нулю. Следовательно, матрица (18.2.27) $^*_*$-вырождена и правое $^*_*$-собственное значение $b(i)$ является $^*_*$-собственным значением. Поскольку

$$(18.2.28) \qquad e_1{}^j_i = \delta^i_j$$

то собственный вектор $e_1^i$, соответствующий правому $^*_*$-собственному значению $b(i)$, является собственным вектором, соответствующим $^*_*$-собственному значению $b(i)$. Равенство

$$(18.2.29) \qquad a_1{}^*_* e_1^i = b(i) e_1^i$$

является следствием определения 17.6.23.

Из равенства (18.2.4) следует, что мы можем записать матрицу (18.2.27) в виде

$$u_2^{-1^*}{}^*_* a_2{}^*_* u_2 - b(i) E_n$$
$$= u_2^{-1^*}{}^*_* (a_2 - u_2 b(i)^*_* u_2^{-1^*})^*_* u_2$$

Так как матрица (18.2.27) $^*_*$-вырождена, то матрица

$$(18.2.30) \qquad a_2 - u_2 b(i)^*_* u_2^{-1^*}$$

также $^*_*$-вырождена. Следовательно, правое $^*_*$-собственное значение $b(i)$ матрицы $a_2$ является $^*_*$-собственным значением пары матриц $(a_2, u_2)$.

Равенства
$$(18.2.31) \quad u_2^{-1^*}{}^*_* a_2{}^*_* u_2{}^*_* e_1^i = b(i) e_1^i$$

$$a_2{}^*_* u_2{}^*_* e_1^i$$
$$(18.2.32) \quad = u_2{}^*_* b(i) e_1^i$$
$$= u_2 b(i)^*_* u_2^{-1^*}{}^*_* u_2{}^*_* e_1^i$$

являются следствием равенств (18.2.4), (18.2.29). Из равенств (18.2.15), (18.2.32) и из теоремы 18.2.9 следует, что собственная строка $u_2^i$ матрицы $a_2$, соответствующий правому $^*_*$-собственному значению $b(i)$, является собственной строкой $u_2^i$ пары матриц $(a_2, u_2)$, соответствующей $^*_*$-собственному значению $b(i)$.    $\square$

ЗАМЕЧАНИЕ 18.2.12. *Равенство*
(18.2.33)

$$u_{1i}{}^*{}_*u_2^{-1^*}{}^*{}_*b(i)u_2 = b(i)u_{1i}$$

*является следствием теоремы 18.2.10. Равенство* (18.2.33) *кажется необычным. Однако, если мы его слегка изменим*

(18.2.34)
$$u_{1i}{}^*{}_*u_2^{-1^*}{}^*{}_*b(i)$$
$$= b(i)u_{1i}{}^*{}_*u_2^{-1^*}{}^*$$

*то равенство* (18.2.34) *следует из равенства* (18.2.9). □

ЗАМЕЧАНИЕ 18.2.13. *Равенство*
(18.2.35)  $u_2 b(i)^*{}_*u_2^{-1^*}{}^*{}_*u_1^i = u_1^i b(i)$

*является следствием теоремы 18.2.11. Равенство* (18.2.35) *кажется необычным. Однако, если мы его слегка изменим*

(18.2.36)
$$b(i)u_2^{-1^*}{}^*{}_*u_1^i$$
$$= u_2^{-1^*}{}^*{}_*u_1^i b(i)$$

*то равенство* (18.2.36) *следует из равенства* (18.2.15). □

## 18.3. Проблема обратной теоремы

ВОПРОС 18.3.1. *Из теоремы [18.3] 18.1.11 следует, что любое правое $_*{}^*$-собственное значение $d(i)$ является $_*{}^*$-собственным значением. Верна ли обратная теорема, а именно, является ли $_*{}^*$-собственное значение правым $_*{}^*$-собственным значением?* □

Получить ответ на вопрос 18.3.1, вообще говоря, непросто. Что произойдёт, если кратность $_*{}^*$-собственного значения будет больше 1? Может ли $n \times n$ матрица иметь меньше чем $n$ собственных значений с учётом кратности?

Может показаться, что, если $n \times n$ матрица имеет $n$ собственных значений кратности 1, то это самый простой случай. Если $n \times n$ матрица с элементами из коммутативной алгебры имеет $n$ собственных значений кратности 1, то множество собственных векторов является базисом векторного пространства размерности $n$. Однако это утверждение неверно, если элементы матрицы принадлежать некоммутативной алгебре.

Рассмотрим теорему 18.3.2.

ТЕОРЕМА 18.3.2. *Пусть $A$ - ассоциативная $D$-алгебра и $V$ - $A$-векторное пространство. Если эндоморфизм имеет $k$ различных собственных значений [18.4] $b(1)$, ..., $b(k)$, то собственные вектора $v(1)$, ..., $v(k)$, относящиеся к различным собственным значениям, линейно независимы.*

Однако теорема 18.3.2 не верна. Мы рассмотрим отдельно доказательства для левого и правого $A$-векторного пространств.

ДОКАЗАТЕЛЬСТВО. Мы будем доказывать теорему индукцией по $k$.

Поскольку собственный вектор отличен от нуля, теорема верна для $k = 1$.

ДОКАЗАТЕЛЬСТВО. Мы будем доказывать теорему индукцией по $k$.

Поскольку собственный вектор отличен от нуля, теорема верна для $k = 1$.

---

[18.3]  Подобные рассуждения верны для левого $^*{}_*$-собственного значения,    а также для левого и правого $^*{}_*$-собственного значения.

[18.4]  Смотри также теорему на странице [3]-210.

**Утверждение 18.3.3.** *Пусть теорема верна для $k = m - 1$.* ⊙

Согласно теоремам 17.5.13, 17.5.27, мы можем положить, что $\overline{\overline{e}}_2 = \overline{\overline{e}}_1$. Пусть $b(1)$, ..., $b(m-1)$, $b(m)$ - различные собственные значения и $v(1)$, ..., $v(m-1)$, $v(m)$ - соответствующие собственные векторы

$$(18.3.1) \qquad f \circ v(i) = v(i)b(i)$$

$$i \neq j \Rightarrow b(i) \neq b(j)$$

Пусть утверждение 18.3.4 верно.

**Утверждение 18.3.4.** *Существует линейная зависимость*

$$a(1)v(1) + ...$$
$$(18.3.2) \qquad + a(m-1)v(m-1)$$
$$+ a(m)v(m) = 0$$

*где* $a(1) \neq 0$. ⊙

Равенство

$$a(1)(f \circ v(1)) + ...$$
$$(18.3.3) \qquad + a(m-1)(f \circ v(m-1))$$
$$+ a(m)(f \circ v(m)) = 0$$

является следствием равенства (18.3.2) и теоремы 10.3.2. Равенство

$$a(1)v(1)b(1) + ...$$
$$+ a(m-1)$$
$$(18.3.4) \qquad * v(m-1)b(m-1)$$
$$+ a(m)v(m)b(m) = 0$$

является следствием равенств (18.3.1), (18.3.3). Так как произведение не коммутативно, мы не можем упростить выражение (18.3.4) до линейной зависимости. □

**Утверждение 18.3.5.** *Пусть теорема верна для $k = m - 1$.* ⊙

Согласно теоремам 17.5.20, 17.5.30, мы можем положить, что $\overline{\overline{e}}_2 = \overline{\overline{e}}_1$. Пусть $b(1)$, ..., $b(m-1)$, $b(m)$ - различные собственные значения и $v(1)$, ..., $v(m-1)$, $v(m)$ - соответствующие собственные векторы

$$(18.3.5) \qquad f \circ v(i) = b(i)v(i)$$

$$i \neq j \Rightarrow b(i) \neq b(j)$$

Пусть утверждение 18.3.6 верно.

**Утверждение 18.3.6.** *Существует линейная зависимость*

$$v(1)a(1) + ...$$
$$(18.3.6) \qquad + v(m-1)a(m-1)$$
$$+ v(m)a(m) = 0$$

*где* $a(1) \neq 0$. ⊙

Равенство

$$(f \circ v(1))a(1) + ...$$
$$(18.3.7) \qquad + (f \circ v(m-1))a(m-1)$$
$$+ (f \circ v(m))a(m) = 0$$

является следствием равенства (18.3.6) и теоремы 11.3.2. Равенство

$$b(1)v(1)a(1) + ...$$
$$+ b(m-1)v(m-1)$$
$$(18.3.8) \qquad * a(m-1)$$
$$+ b(m)v(m)a(m) = 0$$

является следствием равенств (18.3.5), (18.3.7). Так как произведение не коммутативно, мы не можем упростить выражение (18.3.8) до линейной зависимости. □

**Замечание 18.3.7.** *Согласно данному выше определению, множества левых и правых $_*^*$-собственных значений совпадают. Однако мы рассмотрели самый простой случай. Кон рассматривает определения*

$$\boxed{(18.1.5) \quad a_{2*}{}^* u_2 = u_{2*}{}^* a_1} \qquad \boxed{(18.1.2) \quad u_{2*}{}^* a_2 = a_{1*}{}^* u_2}$$

*где $a_2$ - $_*^*$-невырожденная матрица, как начальную точку для исследования.* □

Мы предполагаем, что кратность всех левых $_*^*$-собственных значений $n \times n$ матрица $a_2$ равна 1. Допустим $n \times n$

Мы предполагаем, что кратность всех правых $_*^*$-собственных значений $n \times n$ матрица $a_2$ равна 1. Допустим матрица $a_2$ не имеет $n$ левых $_*$*-собственных значений, но имеет $k$, $k < n$, левых $_*$*-собственных значений. В этом случае, матрица $u_2$ в равенстве

$$(18.1.2) \quad u_{2*}{}^*a_2 = a_{1*}{}^*u_2$$

является $n \times k$ матрицей, и матрица $a_1$ является $k \times k$ матрицей. Однако из равенства (18.1.2) не следует, что

$$\mathrm{rank}(_*{}^*)u_2 = k$$

Если

$$\mathrm{rank}(_*{}^*)u_2 = k$$

то существует $k \times n$ матрица $w \in right(u_2^{-1_*{}^*})$ такая, что

$$(18.3.9) \qquad u_{2*}{}^*w = E_k$$

Равенство

$$(18.3.10) \qquad u_{2*}{}^*a_{2*}{}^*w = a_1$$

является следствием равенств (18.1.2), (18.3.9).

Если

$$\mathrm{rank}(_*{}^*)u_2 < k$$

то между левыми собственными векторами существует линейная зависимость. В этом случае множества левых и правых $_*$*-собственных значений не совпадают.

$n \times n$ матрица $a_2$ не имеет $n$ правых $_*$*-собственных значений, но имеет $k$, $k < n$, правых $_*$*-собственных значений. В этом случае, матрица $u_2$ в равенстве

$$(18.1.5) \quad a_{2*}{}^*u_2 = u_{2*}{}^*a_1$$

является $k \times n$ матрицей, и матрица $a_1$ является $k \times k$ матрицей. Однако из равенства (18.1.5) не следует, что

$$\mathrm{rank}(_*{}^*)u_2 = k$$

Если

$$\mathrm{rank}(_*{}^*)u_2 = k$$

то существует $n \times k$ матрица $w \in left(u_2^{-1_*{}^*})$ такая, что

$$(18.3.11) \qquad w_*{}^*u_2 = E_k$$

Равенство

$$(18.3.12) \qquad w_*{}^*a_{2*}{}^*u_2 = a_1$$

является следствием равенств (18.1.5), (18.3.11).

Если

$$\mathrm{rank}(_*{}^*)u_2 < k$$

то между правыми собственными векторами существует линейная зависимость. В этом случае множества левых и правых $_*$*-собственных значений не совпадают.

## 18.4. Собственное значение матрицы линейного отображения

Пусть

$$a = \begin{pmatrix} a^{1}_{1s0} \otimes a^{1}_{1s1} & ... & a^{1}_{ns0} \otimes a^{1}_{ns1} \\ ... & ... & ... \\ a^{n}_{1s0} \otimes a^{n}_{1s1} & ... & a^{n}_{ns0} \otimes a^{n}_{ns1} \end{pmatrix}$$

матрица линейного отображения

$$\overline{a} : V \to V$$

$A$-векторного пространства $V$.

### 18.4.1. $_\circ{}^\circ$-собственное значение матрицы линейного отображения.

Определение 18.4.1. *$A$-число $b$ называется левым $_\circ{}^\circ$-собственным значением матрицы $a$, если суще-*

Определение 18.4.2. *$A$-число $b$ называется правым $_\circ{}^\circ$-собственным значением матрицы $a$, если суще-*

ствует вектор-столбец $v$ такой, что

(18.4.1)    $$a \circ v = bv$$

*Вектор-столбец $v$ называется* **собственным столбцом,** *соответствующим левому $\circ^\circ$-собственному значению $b$.*    □

ствует вектор-столбец $v$ такой, что

(18.4.2)    $$a \circ v = vb$$

*Вектор-столбец $v$ называется* **собственным столбцом,** *соответствующим правому $\circ^\circ$-собственному значению $b$.*    □

---

Теорема 18.4.3. *Пусть элементы матрицы $a$ удовлетворяют равенству*

(18.4.3)    $$a^i_{j\,s0} \otimes a^i_{j\,s1} = 1 \otimes a^i_{1\,j}$$

*Тогда левое $\circ^\circ$-собственное значение $b$ является левым $^*_*$-собственным значением матрицы $a_1$.*

Теорема 18.4.4. *Пусть элементы матрицы $a$ удовлетворяют равенству*

(18.4.4)    $$a^i_{j\,s0} \otimes a^i_{j\,s1} = a_{0\,j}^{\ i} \otimes 1$$

*Тогда правое $\circ^\circ$-собственное значение $b$ является правым $_*^*$-собственным значением матрицы $a_0$.*

---

Доказательство теоремы 18.4.3.    Равенство

(18.4.5)    $$(a^i_{j\,s0} \otimes a^i_{j\,s1}) \circ v^j = (1 \otimes a^i_{1\,j}) \circ v^j = v^j a^i_{1\,j}$$

является следствием равенства $\boxed{(18.4.3)\quad a^i_{j\,s0} \otimes a^i_{j\,s1} = 1 \otimes a^i_{1\,j}}$ .
Равенство

(18.4.6)

$$\begin{pmatrix} v^1 \\ ... \\ v^n \end{pmatrix} {}^*_* \begin{pmatrix} a_1{}^1_1 & ... & a_1{}^1_n \\ ... & ... & ... \\ a_1{}^n_1 & ... & a_1{}^n_n \end{pmatrix} = \begin{pmatrix} a^1_{1\,s0} \otimes a^1_{1\,s1} & ... & a^1_{n\,s0} \otimes a^1_{n\,s1} \\ ... & ... & ... \\ a^n_{1\,s0} \otimes a^n_{1\,s1} & ... & a^n_{n\,s0} \otimes a^n_{n\,s1} \end{pmatrix} \circ \begin{pmatrix} v^1 \\ ... \\ v^n \end{pmatrix}$$

является следствием равенства (18.4.5). Равенство

(18.4.7)

$$\begin{pmatrix} a^1_{1\,s0} \otimes a^1_{1\,s1} & ... & a^1_{n\,s0} \otimes a^1_{n\,s1} \\ ... & ... & ... \\ a^n_{1\,s0} \otimes a^n_{1\,s1} & ... & a^n_{n\,s0} \otimes a^n_{n\,s1} \end{pmatrix} \circ \begin{pmatrix} v^1 \\ ... \\ v^n \end{pmatrix} = b \begin{pmatrix} v^1 \\ ... \\ v^n \end{pmatrix}$$

является следствием определения 18.4.1 левого $\circ^\circ$-собственного значения матрицы $a$. Равенство

(18.4.8)

$$\begin{pmatrix} v^1 \\ ... \\ v^n \end{pmatrix} {}^*_* \begin{pmatrix} a_1{}^1_1 & ... & a_1{}^1_n \\ ... & ... & ... \\ a_1{}^n_1 & ... & a_1{}^n_n \end{pmatrix} = b \begin{pmatrix} v^1 \\ ... \\ v^n \end{pmatrix}$$

является следствием равенств (18.4.6), (18.4.7). Теорема является следствием равенства (18.4.8) и определения 18.2.1 левого $^*_*$-собственного значения матрицы $a_1$.    □

Доказательство теоремы 18.4.4.    Равенство

(18.4.9)    $$(a^i_{j\,s0} \otimes a^i_{j\,s1}) \circ v^j = (a_{0\,j}^{\ i} \otimes 1) \circ v^j = a_{0\,j}^{\ i} v^j$$

является следствием равенства $\boxed{(18.4.4) \quad a^i_{j\,s0} \otimes a^i_{j\,s1} = a^i_{0\,j} \otimes 1}$.

Равенство

(18.4.10)

$$\begin{pmatrix} a^1_{0\,1} & ... & a^1_{0\,n} \\ ... & ... & ... \\ a^n_{0\,1} & ... & a^n_{0\,n} \end{pmatrix} *^* \begin{pmatrix} v^1 \\ ... \\ v^n \end{pmatrix} = \begin{pmatrix} a^1_{1\,s0} \otimes a^1_{1\,s1} & ... & a^1_{n\,s0} \otimes a^1_{n\,s1} \\ ... & ... & ... \\ a^n_{1\,s0} \otimes a^n_{1\,s1} & ... & a^n_{n\,s0} \otimes a^n_{n\,s1} \end{pmatrix} \circ^\circ \begin{pmatrix} v^1 \\ ... \\ v^n \end{pmatrix}$$

является следствием равенства (18.4.9). Равенство

(18.4.11)

$$\begin{pmatrix} a^1_{1\,s0} \otimes a^1_{1\,s1} & ... & a^1_{n\,s0} \otimes a^1_{n\,s1} \\ ... & ... & ... \\ a^n_{1\,s0} \otimes a^n_{1\,s1} & ... & a^n_{n\,s0} \otimes a^n_{n\,s1} \end{pmatrix} \circ^\circ \begin{pmatrix} v^1 \\ ... \\ v^n \end{pmatrix} = \begin{pmatrix} v^1 \\ ... \\ v^n \end{pmatrix} b$$

является следствием определения 18.4.2 правого $\circ^\circ$-собственного значения матрицы $a$. Равенство

(18.4.12)

$$\begin{pmatrix} a^1_{0\,1} & ... & a^1_{0\,n} \\ ... & ... & ... \\ a^n_{0\,1} & ... & a^n_{0\,n} \end{pmatrix} *^* \begin{pmatrix} v^1 \\ ... \\ v^n \end{pmatrix} = \begin{pmatrix} v^1 \\ ... \\ v^n \end{pmatrix} b$$

является следствием равенств (18.4.10), (18.4.11). Теорема является следствием равенства (18.4.12) и определения 18.1.2 правого $*^*$-собственного значения матрицы $a_0$. □

**ТЕОРЕМА 18.4.5.** *Пусть элементы матрицы $a$ удовлетворяют равенству*

(18.4.13)    $a^i_{j\,s0} \otimes a^i_{j\,s1} = a^i_{0\,j} \otimes 1$

*Тогда левое $\circ^\circ$-собственное значение $b$ является $*^*$-собственным значением матрицы $a_0$.*

**ТЕОРЕМА 18.4.6.** *Пусть элементы матрицы $a$ удовлетворяют равенству*

(18.4.14)    $a^i_{j\,s0} \otimes a^i_{j\,s1} = 1 \otimes a^i_{1\,j}$

*Тогда правое $\circ^\circ$-собственное значение $b$ является $*_*$-собственным значением матрицы $a_1$.*

ДОКАЗАТЕЛЬСТВО ТЕОРЕМЫ 18.4.5. Равенство

(18.4.15)    $(a^i_{j\,s0} \otimes a^i_{j\,s1}) \circ v^j = (a^i_{0\,j} \otimes 1) \circ v^j = a^i_{0\,j} v^j$

является следствием равенства $\boxed{(18.4.13) \quad a^i_{j\,s0} \otimes a^i_{j\,s1} = a^i_{0\,j} \otimes 1}$.

Равенство

(18.4.16)

$$\begin{pmatrix} a^1_{0\,1} & ... & a^1_{0\,n} \\ ... & ... & ... \\ a^n_{0\,1} & ... & a^n_{0\,n} \end{pmatrix} *^* \begin{pmatrix} v^1 \\ ... \\ v^n \end{pmatrix} = \begin{pmatrix} a^1_{1\,s0} \otimes a^1_{1\,s1} & ... & a^1_{n\,s0} \otimes a^1_{n\,s1} \\ ... & ... & ... \\ a^n_{1\,s0} \otimes a^n_{1\,s1} & ... & a^n_{n\,s0} \otimes a^n_{n\,s1} \end{pmatrix} \circ^\circ \begin{pmatrix} v^1 \\ ... \\ v^n \end{pmatrix}$$

является следствием равенства (18.4.15). Равенство

$$(18.4.17) \quad \begin{pmatrix} a_{1\,s0}^{1} \otimes a_{1\,s1}^{1} & ... & a_{n\,s0}^{1} \otimes a_{n\,s1}^{1} \\ ... & ... & ... \\ a_{1\,s0}^{n} \otimes a_{1\,s1}^{n} & ... & a_{n\,s0}^{n} \otimes a_{n\,s1}^{n} \end{pmatrix} \circ^{\circ} \begin{pmatrix} v^{1} \\ ... \\ v^{n} \end{pmatrix} = b \begin{pmatrix} v^{1} \\ ... \\ v^{n} \end{pmatrix}$$

является следствием определения 18.4.1 левого $\circ^{\circ}$-собственного значения матрицы $a$. Равенство

$$(18.4.18) \quad \begin{pmatrix} a_{0\,1}^{1} & ... & a_{0\,n}^{1} \\ ... & ... & ... \\ a_{0\,1}^{n} & ... & a_{0\,n}^{n} \end{pmatrix} *^{*} \begin{pmatrix} v^{1} \\ ... \\ v^{n} \end{pmatrix} = b \begin{pmatrix} v^{1} \\ ... \\ v^{n} \end{pmatrix}$$

является следствием равенств (18.4.16), (18.4.17). Теорема является следствием равенства (18.4.18), определения 17.6.14 $*_*$-собственного значения матрицы $a_0$ и определения 17.6.20 соответствующего собственного вектора $v$. □

Доказательство теоремы 18.4.6. Равенство

$$(18.4.19) \quad (a_{j\,s0}^{i} \otimes a_{j\,s1}^{i}) \circ v^{j} = (1 \otimes a_{1\,j}^{i}) \circ v^{j} = v^{j} a_{1\,j}^{i}$$

является следствием равенства $\boxed{(18.4.14) \quad a_{j\,s0}^{i} \otimes a_{j\,s1}^{i} = 1 \otimes a_{1\,j}^{i}}$ .
Равенство
$$(18.4.20)$$
$$\begin{pmatrix} v^{1} \\ ... \\ v^{n} \end{pmatrix} *_{*} \begin{pmatrix} a_{1\,1}^{1} & ... & a_{1\,n}^{1} \\ ... & ... & ... \\ a_{1\,1}^{n} & ... & a_{1\,n}^{n} \end{pmatrix} = \begin{pmatrix} a_{1\,s0}^{1} \otimes a_{1\,s1}^{1} & ... & a_{n\,s0}^{1} \otimes a_{n\,s1}^{1} \\ ... & ... & ... \\ a_{1\,s0}^{n} \otimes a_{1\,s1}^{n} & ... & a_{n\,s0}^{n} \otimes a_{n\,s1}^{n} \end{pmatrix} \circ^{\circ} \begin{pmatrix} v^{1} \\ ... \\ v^{n} \end{pmatrix}$$

является следствием равенства (18.4.19). Равенство

$$(18.4.21) \quad \begin{pmatrix} a_{1\,s0}^{1} \otimes a_{1\,s1}^{1} & ... & a_{n\,s0}^{1} \otimes a_{n\,s1}^{1} \\ ... & ... & ... \\ a_{1\,s0}^{n} \otimes a_{1\,s1}^{n} & ... & a_{n\,s0}^{n} \otimes a_{n\,s1}^{n} \end{pmatrix} \circ^{\circ} \begin{pmatrix} v^{1} \\ ... \\ v^{n} \end{pmatrix} = \begin{pmatrix} v^{1} \\ ... \\ v^{n} \end{pmatrix} b$$

является следствием определения 18.4.2 правого $\circ^{\circ}$-собственного значения матрицы $a$. Равенство

$$(18.4.22) \quad \begin{pmatrix} v^{1} \\ ... \\ v^{n} \end{pmatrix} *_{*} \begin{pmatrix} a_{1\,1}^{1} & ... & a_{1\,n}^{1} \\ ... & ... & ... \\ a_{1\,1}^{n} & ... & a_{1\,n}^{n} \end{pmatrix} = \begin{pmatrix} v^{1} \\ ... \\ v^{n} \end{pmatrix} b$$

является следствием равенств (18.4.20), (18.4.21). Теорема является следствием равенства (18.4.22), определения 17.6.17 $*_*$-собственного значения матрицы $a_1$ и определения 17.6.21 соответствующего собственного вектора $v$. □

Теорема 18.4.7. *Пусть вектор-столбец $v$ является собственным столбцом, соответствующим лево-му $\circ^{\circ}$-собственному значению $b$ мат-*

Теорема 18.4.8. *Пусть вектор-столбец $v$ является собственным столбцом, соответствующим право-му $\circ^{\circ}$-собственному значению $b$ в мат-* *рицы a. Пусть A-число b удовлетворяет условию*

$$(18.4.23) \qquad b \in \bigcap_{i=1}^{n} \bigcap_{j=1}^{n} Z(A, a_{j\,s0}^{i})$$

*Тогда для любого многочлена*

$$(18.4.24) \qquad p(x) = p_0 + p_1 x + ... + p_n x^n$$
$$p_0, \ ..., \ p_n \in D$$

*вектор-столбец*

$$cv = \begin{pmatrix} cv^1 \\ ... \\ cv^n \end{pmatrix}$$
$$c = p(b)$$

*также является собственным столбцом, соответствующим левому ∘-собственному значению b.*

Доказательство. Согласно теореме 5.1.25, равенство

$$(18.4.25) \qquad \begin{aligned} (a_{j\,s0}^{i} \otimes a_{j\,s1}^{i}) \circ (cv^j) \\ = a_{j\,s0}^{i} cv^j a_{j\,s1}^{i} \\ = c(a_{j\,s0}^{i} v^j a_{j\,s1}^{i}) \\ = cbv^i = b(cv^i) \end{aligned}$$

является следствием условия (18.4.23). Согласно определению 18.4.1, из равенства (18.4.25) следует, что вектор-столбец $cv$ является собственным столбцом, соответствующим левому ∘-собственному значению $b$ матрицы $a$.  □

*рицы a. Пусть A-число b удовлетворяет условию*

$$(18.4.26) \qquad b \in \bigcap_{i=1}^{n} \bigcap_{j=1}^{n} Z(A, a_{j\,s1}^{i})$$

*Тогда для любого многочлена*

$$(18.4.27) \qquad p(x) = p_0 + p_1 x + ... + p_n x^n$$
$$p_0, \ ..., \ p_n \in D$$

*вектор-столбец*

$$vc = \begin{pmatrix} v^1 c \\ ... \\ v^n c \end{pmatrix}$$
$$c = p(b)$$

*также является собственным столбцом, соответствующим правому ∘-собственному значению b.*

Доказательство. Согласно теореме 5.1.25, равенство

$$(18.4.28) \qquad \begin{aligned} (a_{j\,s0}^{i} \otimes a_{j\,s1}^{i}) \circ (v^j c) \\ = a_{j\,s0}^{i} v^j c a_{j\,s1}^{i} \\ = (a_{j\,s0}^{i} v^j a_{j\,s1}^{i})c \\ = v^i bc = (v^i c)b \end{aligned}$$

является следствием условия (18.4.26). Согласно определению 18.4.2, из равенства (18.4.28) следует, что вектор-столбец $vc$ является собственным столбцом, соответствующим правому ∘-собственному значению $b$ матрицы $a$.  □

## 18.4.2. ∘-собственное значение матрицы линейного отображения.

Определение 18.4.9. *A-число b называется левым ∘-собственным значением матрицы a, если существует вектор-строка v такой, что*

$$(18.4.29) \qquad a^{\circ} \circ v = bv$$

*Вектор-строка v называется собственной строкой, соответствующей левому ∘-собственному значению b.*  □

Определение 18.4.10. *A-число b называется правым ∘-собственным значением матрицы a, если существует вектор-строка v такой, что*

$$(18.4.30) \qquad a^{\circ} \circ v = vb$$

*Вектор-строка v называется собственной строкой, соответствующей правому ∘-собственному значению b.*  □

<table>
<tr><td>

Теорема 18.4.11. *Пусть элементы матрицы $a$ удовлетворяют равенству*

$$(18.4.31) \qquad a^j_{i\,s0} \otimes a^j_{i\,s1} = 1 \otimes a^j_{1\,i}$$

*Тогда левое $\circ$-собственное значение $b$ является левым $_*$-собственным значением матрицы $a_1$.*

</td><td>

Теорема 18.4.12. *Пусть элементы матрицы $a$ удовлетворяют равенству*

$$(18.4.32) \qquad a^j_{i\,s0} \otimes a^j_{i\,s1} = a^j_{0\,i} \otimes 1$$

*Тогда правое $\circ$-собственное значение $b$ является правым $^*_*$-собственным значением матрицы $a_0$.*

</td></tr>
</table>

Доказательство теоремы 18.4.11. Равенство

$$(18.4.33) \qquad (a^j_{i\,s0} \otimes a^j_{i\,s1}) \circ v_i = (1 \otimes a^j_{1\,i}) \circ v_i = v_i a^j_{1\,i}$$

является следствием равенства $\boxed{(18.4.31) \quad a^j_{i\,s0} \otimes a^j_{i\,s1} = 1 \otimes a^j_{1\,i}}$. Равенство

$$(18.4.34)$$

$$\begin{pmatrix} v_1 & \dots & v_n \end{pmatrix} {}_*{*} \begin{pmatrix} a^1_{1\,1} & \dots & a^1_{1\,n} \\ \dots & \dots & \dots \\ a^n_{1\,1} & \dots & a^n_{1\,n} \end{pmatrix} = \begin{pmatrix} a^1_{1\,s0} \otimes a^1_{1\,s1} & \dots & a^n_{1\,s0} \otimes a^n_{1\,s1} \\ \dots & \dots & \dots \\ a^1_{n\,s0} \otimes a^1_{n\,s1} & \dots & a^n_{n\,s0} \otimes a^n_{n\,s1} \end{pmatrix} \circ \begin{pmatrix} v_1 & \dots & v_n \end{pmatrix}$$

является следствием равенства (18.4.33). Равенство

$$(18.4.35) \qquad \begin{pmatrix} a^1_{1\,s0} \otimes a^1_{1\,s1} & \dots & a^n_{1\,s0} \otimes a^n_{1\,s1} \\ \dots & \dots & \dots \\ a^1_{n\,s0} \otimes a^1_{n\,s1} & \dots & a^n_{n\,s0} \otimes a^n_{n\,s1} \end{pmatrix} \circ \begin{pmatrix} v_1 & \dots & v_n \end{pmatrix} = b \begin{pmatrix} v_1 & \dots & v_n \end{pmatrix}$$

является следствием определения 18.4.9 левого $\circ$-собственного значения матрицы $a$. Равенство

$$(18.4.36) \qquad \begin{pmatrix} v_1 & \dots & v_n \end{pmatrix} {}_*{*} \begin{pmatrix} a^1_{1\,1} & \dots & a^1_{1\,n} \\ \dots & \dots & \dots \\ a^n_{1\,1} & \dots & a^n_{1\,n} \end{pmatrix} = b \begin{pmatrix} v_1 & \dots & v_n \end{pmatrix}$$

является следствием равенств (18.4.34), (18.4.35). Теорема является следствием равенства (18.4.36) и определения 18.1.1 левого $_*$-собственного значения матрицы $a_1$. $\square$

Доказательство теоремы 18.4.12. Равенство

$$(18.4.37) \qquad (a^j_{i\,s0} \otimes a^j_{i\,s1}) \circ v_i = (a^j_{0\,i} \otimes 1) \circ v_i = a^j_{0\,i} v_i$$

является следствием равенства $\boxed{(18.4.32) \quad a^j_{i\,s0} \otimes a^j_{i\,s1} = a^j_{0\,i} \otimes 1}$. Равенство

$$(18.4.38)$$

$$\begin{pmatrix} a^1_{0\,1} & \dots & a^1_{0\,n} \\ \dots & \dots & \dots \\ a^n_{0\,1} & \dots & a^n_{0\,n} \end{pmatrix} {}_*{*} \begin{pmatrix} v_1 & \dots & v_n \end{pmatrix} = \begin{pmatrix} a^1_{1\,s0} \otimes a^1_{1\,s1} & \dots & a^n_{1\,s0} \otimes a^n_{1\,s1} \\ \dots & \dots & \dots \\ a^1_{n\,s0} \otimes a^1_{n\,s1} & \dots & a^n_{n\,s0} \otimes a^n_{n\,s1} \end{pmatrix} \circ \begin{pmatrix} v_1 & \dots & v_n \end{pmatrix}$$

является следствием равенства (18.4.37). Равенство

$$(18.4.39) \quad \begin{pmatrix} a_{1\,s0}^{1} \otimes a_{1\,s1}^{1} & \dots & a_{1\,s0}^{n} \otimes a_{1\,s1}^{n} \\ \dots & \dots & \dots \\ a_{n\,s0}^{1} \otimes a_{n\,s1}^{1} & \dots & a_{n\,s0}^{n} \otimes a_{n\,s1}^{n} \end{pmatrix} \circ_{\circ} \begin{pmatrix} v_{1} & \dots & v_{n} \end{pmatrix} = \begin{pmatrix} v_{1} & \dots & v_{n} \end{pmatrix} b$$

является следствием определения 18.4.10 правого $^{\circ}{}_{\circ}$-собственного значения матрицы $a$. Равенство

$$(18.4.40) \quad \begin{pmatrix} a_{0\,1}^{1} & \dots & a_{0\,n}^{1} \\ \dots & \dots & \dots \\ a_{0\,1}^{n} & \dots & a_{0\,n}^{n} \end{pmatrix} {}^{*}_{*} \begin{pmatrix} v_{1} & \dots & v_{n} \end{pmatrix} = \begin{pmatrix} v_{1} & \dots & v_{n} \end{pmatrix} b$$

является следствием равенств (18.4.38), (18.4.39). Теорема является следствием равенства (18.4.40) и определения 18.2.2 правого $^{*}{}_{*}$-собственного значения матрицы $a_0$. □

| | |
|---|---|
| Теорема 18.4.13. *Пусть элементы матрицы a удовлетворяют равенству* $$(18.4.41) \quad a_{i\,s0}^{j} \otimes a_{i\,s1}^{j} = a_{i}^{j} \otimes 1$$ *Тогда левое $^{\circ}{}_{\circ}$-собственное значение b является $^{*}{}_{*}$-собственным значением матрицы $a_0$.* | Теорема 18.4.14. *Пусть элементы матрицы a удовлетворяют равенству* $$(18.4.42) \quad a_{i\,s0}^{j} \otimes a_{i\,s1}^{j} = 1 \otimes a_{1\,i}^{j}$$ *Тогда правое $^{\circ}{}_{\circ}$-собственное значение b является $_{*}{}^{*}$-собственным значением матрицы $a_1$.* |

Доказательство теоремы 18.4.13. Равенство

$$(18.4.43) \quad (a_{i\,s0}^{j} \otimes a_{i\,s1}^{j}) \circ v_i = (a_{0\,i}^{j} \otimes 1) \circ v_i = a_{0\,i}^{j} v_i$$

является следствием равенства $\boxed{(18.4.41) \quad a_{i\,s0}^{j} \otimes a_{i\,s1}^{j} = a_{0\,i}^{j} \otimes 1}$.
Равенство
$(18.4.44)$

$$\begin{pmatrix} a_{0\,1}^{1} & \dots & a_{0\,n}^{1} \\ \dots & \dots & \dots \\ a_{0\,1}^{n} & \dots & a_{0\,n}^{n} \end{pmatrix} {}^{*}_{*} \begin{pmatrix} v_{1} & \dots & v_{n} \end{pmatrix} = \begin{pmatrix} a_{1\,s0}^{1} \otimes a_{1\,s1}^{1} & \dots & a_{1\,s0}^{n} \otimes a_{1\,s1}^{n} \\ \dots & \dots & \dots \\ a_{n\,s0}^{1} \otimes a_{n\,s1}^{1} & \dots & a_{n\,s0}^{n} \otimes a_{n\,s1}^{n} \end{pmatrix} \circ_{\circ} \begin{pmatrix} v_{1} & \dots & v_{n} \end{pmatrix}$$

является следствием равенства (18.4.43). Равенство

$$(18.4.45) \quad \begin{pmatrix} a_{1\,s0}^{1} \otimes a_{1\,s1}^{1} & \dots & a_{1\,s0}^{n} \otimes a_{1\,s1}^{n} \\ \dots & \dots & \dots \\ a_{n\,s0}^{1} \otimes a_{n\,s1}^{1} & \dots & a_{n\,s0}^{n} \otimes a_{n\,s1}^{n} \end{pmatrix} \circ_{\circ} \begin{pmatrix} v_{1} & \dots & v_{n} \end{pmatrix} = b \begin{pmatrix} v_{1} & \dots & v_{n} \end{pmatrix}$$

является следствием определения 18.4.9 левого $^\circ_\circ$-собственного значения матрицы $a$. Равенство

$$(18.4.46) \qquad \begin{pmatrix} a_{0\,1}^{\phantom{0}1} & \dots & a_{0\,n}^{\phantom{0}1} \\ \dots & \dots & \dots \\ a_{0\,1}^{\phantom{0}n} & \dots & a_{0\,n}^{\phantom{0}n} \end{pmatrix} {}_*^* \begin{pmatrix} v_1 & \dots & v_n \end{pmatrix} = b \begin{pmatrix} v_1 & \dots & v_n \end{pmatrix}$$

является следствием равенств (18.4.44), (18.4.45). Теорема является следствием равенства (18.4.46), определения 17.6.17 $_*^*$-собственного значения матрицы $a_0$ и определения 17.6.23 соответствующего собственного вектора $v$. $\qquad\square$

Доказательство теоремы 18.4.14. Равенство

$$(18.4.47) \qquad (a_{i\,s0}^j \otimes a_{i\,s1}^j) \circ v_i = (1 \otimes a_{1\,i}^j) \circ v_i = v_i a_{1\,i}^j$$

является следствием равенства $\boxed{(18.4.42) \quad a_{i\,s0}^j \otimes a_{i\,s1}^j = 1 \otimes a_{1\,i}^j}$ .
Равенство
(18.4.48)

$$\begin{pmatrix} v_1 & \dots & v_n \end{pmatrix} {}_*^* \begin{pmatrix} a_{1\,1}^{\phantom{1}1} & \dots & a_{1\,n}^{\phantom{1}1} \\ \dots & \dots & \dots \\ a_{1\,1}^{\phantom{1}n} & \dots & a_{1\,n}^{\phantom{1}n} \end{pmatrix} = \begin{pmatrix} a_{1\,s0}^1 \otimes a_{1\,s1}^1 & \dots & a_{n\,s0}^1 \otimes a_{n\,s1}^1 \\ \dots & \dots & \dots \\ a_{1\,s0}^n \otimes a_{1\,s1}^n & \dots & a_{n\,s0}^n \otimes a_{n\,s1}^n \end{pmatrix} {}_\circ^\circ \begin{pmatrix} v_1 & \dots & v_n \end{pmatrix}$$

является следствием равенства (18.4.47). Равенство

$$(18.4.49) \qquad \begin{pmatrix} a_{1\,s0}^1 \otimes a_{1\,s1}^1 & \dots & a_{n\,s0}^1 \otimes a_{n\,s1}^1 \\ \dots & \dots & \dots \\ a_{1\,s0}^n \otimes a_{1\,s1}^n & \dots & a_{n\,s0}^n \otimes a_{n\,s1}^n \end{pmatrix} {}_\circ^\circ \begin{pmatrix} v_1 & \dots & v_n \end{pmatrix} = \begin{pmatrix} v_1 & \dots & v_n \end{pmatrix} b$$

является следствием определения 18.4.10 правого $^\circ_\circ$-собственного значения матрицы $a$. Равенство

$$(18.4.50) \qquad \begin{pmatrix} v_1 & \dots & v_n \end{pmatrix} {}_*^* \begin{pmatrix} a_{1\,1}^{\phantom{1}1} & \dots & a_{1\,n}^{\phantom{1}1} \\ \dots & \dots & \dots \\ a_{1\,1}^{\phantom{1}n} & \dots & a_{1\,n}^{\phantom{1}n} \end{pmatrix} = \begin{pmatrix} v_1 & \dots & v_n \end{pmatrix} b$$

является следствием равенств (18.4.48), (18.4.49). Теорема является следствием равенства (18.4.50), определения 17.6.14 $_*^*$-собственного значения матрицы $a_1$ и определения 17.6.22 соответствующего собственного вектора $v$. $\qquad\square$

Теорема 18.4.15. *Пусть вектор-строка $v$ является собственной строкой, соответствующей левому $^\circ_\circ$-собственному значению $b$ матрицы $a$. Пусть $A$-число $b$ удовлетворяет условию*

$$(18.4.51) \qquad b \in \bigcap_{i=1}^n \bigcap_{j=1}^n Z(A, a_{i\,s0}^j)$$

Теорема 18.4.16. *Пусть вектор-строка $v$ является собственной строкой, соответствующей правому $^\circ_\circ$-собственному значению $b$ матрицы $a$. Пусть $A$-число $b$ удовлетворяет условию*

$$(18.4.54) \qquad b \in \bigcap_{i=1}^n \bigcap_{j=1}^n Z(A, a_{i\,s1}^j)$$

*Тогда для любого многочлена*

$$(18.4.52) \qquad p(x) = p_0 + p_1 x + ... + p_n x^n$$

$$p_0, \ ..., \ p_n \in D$$

*вектор-строка*

$$cv = \begin{pmatrix} cv_1 & ... & cv_n \end{pmatrix}$$

$$c = p(b)$$

*также является собственной строкой, соответствующей левому °ₒ-собственному значению b.*

---

ДОКАЗАТЕЛЬСТВО. Согласно теореме 5.1.25, равенство

$$(18.4.53) \qquad \begin{aligned} & (a^j_{i\,s0} \otimes a^j_{i\,s1}) \circ (cv_i) \\ &= a^j_{i\,s0} cv_i a^j_{i\,s1} \\ &= c(a^j_{i\,s0} v_i a^j_{i\,s1}) \\ &= cbv_j = b(cv_j) \end{aligned}$$

является следствием условия (18.4.51). Согласно определению 18.4.9, из равенства (18.4.53) следует, что вектор-строка $cv$ является собственной строкой, соответствующей левому °ₒ-собственному значению $b$ матрицы $a$.  □

---

*Тогда для любого многочлена*

$$(18.4.55) \qquad p(x) = p_0 + p_1 x + ... + p_n x^n$$

$$p_0, \ ..., \ p_n \in D$$

*вектор-строка*

$$vc = \begin{pmatrix} v_1 c & ... & v_n c \end{pmatrix}$$

$$c = p(b)$$

*также является собственной строкой, соответствующей правому °ₒ-собственному значению b.*

---

ДОКАЗАТЕЛЬСТВО. Согласно теореме 5.1.25, равенство

$$(18.4.56) \qquad \begin{aligned} & (a^j_{i\,s0} \otimes a^j_{i\,s1}) \circ (v_i c) \\ &= a^j_{i\,s0} v_i c a^j_{i\,s1} \\ &= (a^j_{i\,s0} v_i a^j_{i\,s1}) c \\ &= v_j b c = (v_f c) b \end{aligned}$$

является следствием условия (18.4.54). Согласно определению 18.4.10, из равенства (18.4.56) следует, что вектор-строка $vc$ является собственной строкой, соответствующей правому °ₒ-собственному значению $b$ матрицы $a$.  □

# Список литературы

# Предметный указатель

# Специальные символы и обозначения